%% file: essai_livre_sites_topos.tex
\documentclass{book}
\usepackage{amsmath,amssymb,amscd,stmaryrd,bm,rotating,trfsigns,setspace,enumitem,adjustbox}
\usepackage{amsthm} 
\usepackage{mathtools,extarrows}
\usepackage[french]{babel}
\usepackage[all]{xy}
\usepackage{enumitem} 
\usepackage{xcolor}
\usepackage{hyperref}
\def\un{{\rm 1\mkern-4mu I}}
\makeatletter
\def\revddots{\mathinner{\mkern1mu\raise\p@
		\vbox{\kern7\p@\hbox{.}}\mkern2mu
		\raise4\p@\hbox{.}\mkern2mu\raise7\p@\hbox{.}\mkern1mu}}
\makeatother

\makeatletter
\newcommand*{\relrelbarsep}{.386ex}
\newcommand*{\relrelbar}{%
	\mathrel{%
		\mathpalette\@relrelbar\relrelbarsep
	}%
}
\newcommand*{\@relrelbar}[2]{%
	\raise#2\hbox to 0pt{$\m@th#1\relbar$\hss}%
	\lower#2\hbox{$\m@th#1\relbar$}%
}
\providecommand*{\rightrightarrowsfill@}{%
	\arrowfill@\relrelbar\relrelbar\rightrightarrows
}
\providecommand*{\leftleftarrowsfill@}{%
	\arrowfill@\leftleftarrows\relrelbar\relrelbar
}
\providecommand*{\xrightrightarrows}[2][]{%
	\ext@arrow 0359\rightrightarrowsfill@{#1}{#2}%
}
\providecommand*{\xleftleftarrows}[2][]{%
	\ext@arrow 3095\leftleftarrowsfill@{#1}{#2}%
}
\makeatother

\usepackage{scalerel,stackengine}
\stackMath
\newcommand\reallywidehat[1]{%
	\savestack{\tmpbox}{\stretchto{%
			\scaleto{%
				\scalerel*[\widthof{\ensuremath{#1}}]{\kern-.6pt\bigwedge\kern-.6pt}%
				{\rule[-\textheight/2]{1ex}{\textheight}}
			}{\textheight}%
		}{0.5ex}}%
	\stackon[1pt]{#1}{\tmpbox}%
}

\def\dar[#1]{\ar@<2pt>[#1]\ar@<-2pt>[#1]}

\newtheoremstyle{propo}
{1.5em}
{1.5em}
{\itshape
	\interlinepenalty=10
	\predisplaypenalty=100
	\postdisplaypenalty=100
}
{}
{\bfseries}
{.\,--}
{\newline}
{}
\theoremstyle{propo}

\newcounter{thm}[section]
\renewcommand{\thethm}{\Roman{chapter}.\arabic{section}.\arabic{thm}}
\newenvironment{thm}[1][]{\refstepcounter{thm}\par\medskip
	\noindent \textbf{Th\'{e}or\`{e}me~\thethm. -- #1} 
	{\color{white}Grothendieck}
	
	\begin{em}}{\end{em}\medskip}
\newenvironment{lem}[1][]{\refstepcounter{thm}\par\medskip
	\noindent \textbf{Lemme~\thethm. -- #1} 
	{\color{white}Grothendieck}
	
	\begin{em}}{\end{em}\medskip}
\newenvironment{defn}[1][]{\refstepcounter{thm}\par\medskip
	\noindent \textbf{D\'{e}finition~\thethm. -- #1} 
	{\color{white}Grothendieck}
	
	\begin{em}}{\end{em}\medskip}
\newenvironment{prop}[1][]{\refstepcounter{thm}\par\medskip
	\noindent \textbf{Proposition~\thethm. -- #1} 
	{\color{white}Grothendieck}
	
	\begin{em}}{\end{em}\medskip}

\newenvironment{cor}[1][]{\refstepcounter{thm}\par\medskip
	\noindent \textbf{Corollaire~\thethm. -- #1} 
	{\color{white}Grothendieck}
	
	\begin{em}}{\end{em}\medskip}

\newenvironment{remark}[1][]{\par\medskip
	\noindent \textbf{Remarque~ : #1} 
	{\color{white}Grothendieck}
	
	\begin{rm}}{\end{rm}\medskip}

\newenvironment{remarks}[1][]{\par\medskip
	\noindent \textbf{Remarques~ : #1} 
	{\color{white}Grothendieck}
	
	\begin{rm}}{\end{rm}\medskip}

\newenvironment{remarkqed}[1][]{\par\medskip
	\noindent \textbf{Remarque~ : #1} 
	{\color{white}Grothendieck}
	
	\begin{rm}}{\end{rm}{\color{white}Grothendieck}\null\hfill\qedsymbol\medskip}

\newenvironment{remarksqed}[1][]{\par\medskip
	\noindent \textbf{Remarques~ : #1} 
	{\color{white}Grothendieck}
	
	\begin{rm}}{\end{rm}{\color{white}Grothendieck}\nopagebreak\null\hfill\qedsymbol\medskip}

\newenvironment{demo}[1][]{\par\medskip
	\noindent \textbf{D\'{e}monstration~ : #1} 
	{\color{white}Grothendieck}
	
	\begin{rm}}{\end{rm}{\color{white}Grothendieck}\nopagebreak\null\hfill\qedsymbol\medskip}

\newenvironment{demolem}[1][]{\par\medskip
	\noindent \textbf{D\'{e}monstration du lemme~ : #1} 
	{\color{white}Grothendieck}
	
	\begin{rm}}{\end{rm}{\color{white}Grothendieck}\nopagebreak\null\hfill\qedsymbol\medskip}

\newenvironment{demothm}[1][]{\par\medskip
	\noindent \textbf{D\'{e}monstration du th\'{e}or\`{e}me~ : #1} 
	{\color{white}Grothendieck}
	
	\begin{rm}}{\end{rm}{\color{white}Grothendieck}\nopagebreak\null\hfill\qedsymbol\medskip}

\newenvironment{demopart}[1][]{\par\medskip
	\noindent \textbf{Démonstration partielle~ : #1} 
	{\color{white}Grothendieck}
	
	\begin{rm}}{\end{rm}{\color{white}Grothendieck}\nopagebreak\null\hfill\qedsymbol\medskip}

\newenvironment{democor}[1][]{\par\medskip
	\noindent \textbf{Démonstration du corollaire~ : #1} 
	{\color{white}Grothendieck}
	
	\begin{rm}}{\end{rm}{\color{white}Grothendieck}\nopagebreak\null\hfill\qedsymbol\medskip}

\newenvironment{demoprop}[1][]{\par\medskip
	\noindent \textbf{Démonstration de la proposition~ : #1} 
	{\color{white}Grothendieck}
	
	\begin{rm}}{\end{rm}{\color{white}Grothendieck}\nopagebreak\null\hfill\qedsymbol\medskip}

\newenvironment{demosansqed}[1][]{\par\medskip
	\noindent \textbf{Démonstration~ : #1} 
	{\color{white}Grothendieck}
	
	\begin{rm}}{\end{rm}\medskip}

\newenvironment{demopropsansqed}[1][]{\par\medskip
	\noindent \textbf{Démonstration de la proposition~ : #1} 
	{\color{white}Grothendieck}
	
	\begin{rm}}{\end{rm}\medskip}

\newenvironment{commdemo}[1][]{\par\medskip
	\noindent \textbf{Commencement de la d\'emonstration:~ : #1} 
	{\color{white}Grothendieck}
	
	\begin{rm}}{\end{rm}{\color{white}Grothendieck}\nopagebreak\null\hfill\qedsymbol\medskip}

\newenvironment{listeisansmarge}[1][]{\begin{enumerate}[label=(\roman*),leftmargin=0pt]}{\end{enumerate}}

\newenvironment{listeimarge}[1][]{\begin{enumerate}[label=(\roman*)]}{\end{enumerate}}

\allowdisplaybreaks

\title{Sites et topos de Grothendieck: une introduction}

	\author{\Large Olivia Caramello$^{(*)}$ et Laurent Lafforgue$^{(**)}$}

\date{
	(*) {\normalsize Università di Insubria, Como, Italie} \\
	{\normalsize Istituto Grothendieck, Corso Statuto 24, Mondovi CN, Italie} \\
	{\normalsize \&} \\
	(**) {\normalsize Centre Lagrange, 103 rue de Grenelle, 75007 Paris, France} \\
	{\normalsize Centre de recherche de Huawei, 18 quai du Point du Jour,} \\
	{\normalsize  92100 Boulogne-Billancourt, France}
}

\begin{document}

\maketitle

\section*{Introduction}

La théorie des topos occupe une position singulière dans le paysage des mathématiques contemporaines : née et développée pour les besoins de la géométrie algébrique d'Alexandre Grothendieck (\cite{SGA4tome1}), elle a progressivement révélé sa puissance comme langage unificateur de la géométrie, de la topologie, de l'algèbre, et même de la logique mathématique. Ce livre propose une introduction structurée et progressive à cette théorie, en partant des théories classiques ou antérieures sur la base desquelles elle a été développée (groupes et actions de groupes, espaces topologiques, catégories, faisceaux) pour aboutir aux fondements avancés de la géométrie des topos et à ses liens profonds avec la logique géométrique qui s’organisent autour du théorème de construction des topos classifiants.

L'ambition est double. D'une part, montrer comment la notion de topos surgit naturellement dès qu'on pousse ``jusqu'au bout'' (pour utiliser les mots de Grothendieck) l'idée que les faisceaux codent ``l'expression'' d'objets géométriques tels que les espaces topologiques classiques, permettant alors de voir et étudier comme des espaces n’importe quelles catégories munies de ``topologies de Grothendieck''. D'autre part, mettre au jour la richesse de cette théorie d'un point de vue interne : un topos n'est pas qu'une ``catégorie de faisceaux'' abstraite, c'est aussi un univers mathématique où s'interprètent structures, équations, preuves et modèles, qui possède toutes les propriétés ``constructives'' de la catégorie des ensembles qui s’expriment en termes de limites et colimites, et admet des exponentielles arbitraires et un ``classificateur des sous-objets''.

Ce texte propose une entrée progressive, ambitieuse et unifiée dans cette théorie à la manière de Grothendieck. L'architecture générale suit un fil conducteur très clair : partir des notions mathématiques familières (groupes, espaces, faisceaux), en abstraire le langage (catégories, foncteurs, adjonctions, limites), formaliser ce que signifie ``faire de la géométrie par les faisceaux'' (sites et topologies de Grothendieck), puis dégager la notion abstraite de topos et ses propriétés structurelles fondamentales, avant d'en explorer la géométrie (morphismes, points, sous-topos, localisations) et enfin d'en faire la passerelle naturelle vers la logique (théories, modèles, catégories syntactiques, topos classifiants). Voici, chapitre par chapitre, la trame et les idées essentielles.

\textbf{I. Avant les topos : groupes, espaces topologiques, catégories}

Le livre commence en terrain familier. D'abord les groupes et monoïdes : on n'en donne pas seulement les définitions classiques, mais on insiste sur leur ``expression'' par leurs actions. Cette perspective est fondamentale : un groupe se révèle par ses actions sur des objets (ensembles, modules, espaces), et des notions comme celles d’orbites, de stabilisateurs ou de quotients surgissent naturellement dans ce cadre. Cette grammaire des actions prépare l'idée qu'une structure se révèle par les foncteurs qu'elle engendre.

Viennent ensuite les espaces topologiques et leurs morphismes (applications continues), puis surtout les faisceaux sur un espace : l'accent est mis sur le passage des préfaisceaux aux faisceaux via le recollement local-global, et sur les premiers exemples clefs (fonctions continues, différentiables ou holomorphes). Les schémas apparaissent très tôt comme une synthèse puissante : on explique le spectre d'un anneau, sa topologie de Zariski, puis on construit le faisceau de structure et les faisceaux quasi-cohérents associés à des modules. Cela constitue un modèle concret et particulièrement important de ce que fera la théorie des sites et des topos à un niveau d'abstraction supérieur.

Le cœur conceptuel de ce chapitre est la théorie des catégories : objets, flèches, composition, sous-catégories pleines. On détaille les premiers exemples utiles (la catégorie des ensembles, celle des espaces topologiques, celle des modules sur un anneau, la catégorie des petites catégories, les catégories des actions d'un groupe, les catégories de préfaisceaux et de faisceaux). Le lemme de Yoneda se déploie comme le pivot de la théorie des catégories et prend un sens presque philosophique : un objet mathématique d’un certain type est entièrement déterminé par la manière dont tous les objets du même type pointent vers lui. Du lemme de Yoneda découlent les notions de foncteur représentable et d'élément universel, et donc la méthode qui consiste à introduire des objets mathématiques comme les solutions de ``problèmes universels'' (produits tensoriels, puissances symétriques et alternées, espaces affines, groupes linéaires,orthogonaux ou symplectiques). Les notions de limite et de colimite sont ensuite présentées comme les formalisations mathématiques de l’idée de ``propriété universelle'' au plus haut niveau de généralité.

On introduit enfin la notion fondamentale d'adjonction entre foncteurs, et les premiers exemples importants de foncteurs adjoints (images directes et réciproques, faisceautisation, prolongements par zéro). Le chapitre se conclut par un panorama de cas notables de limites et colimites, tels que les produits fibrés et les sommes amalgamées, et par les constructions des catégories relatives ou fléchées, qui permettent de formaliser le fait de raisonner ``au-dessus'' d'un objet.

En filigrane, ce chapitre donne déjà un aperçu de la dualité entre local et global dont la théorie des topos constituera une expression achevée : on apprend à passer d'un calcul local (sur des ouverts, sur des présentations) à des notions globales par recollement universel.

\textbf{II. Topologies de Grothendieck, sites et faisceaux sur un site}

Ce chapitre aborde le pas conceptuel majeur dans l'introduction des topos par Grothendieck, celui qui consiste à libérer la notion de recouvrement du cadre topologique ordinaire. Au sens de Grothendieck, une topologie n'est plus une collection d'ouverts d’un ensemble mais une notion de crible couvrant sur les objets d'une catégorie, dont on demande seulement qu’elle comprenne les cribles maximaux, soit stable par image réciproque et soit transitive. Un ``site'' est une catégorie munie d'une telle topologie. Cette définition abstraite formalise et capture, dans un langage purement catégorique, l’intuition naturelle de recouvrement dans un cadre géométrique.

Le livre explique longuement ce qu'est-ce qu'un crible et comment transporter des cribles, et il donne une multitude d’exemples concrets (topologie usuelle d'un espace, topologie canonique d'un treillis, topologies de densité, topologies atomiques, topologies induites par des ensembles ordonnés ou par des mesures, topologies définies à partir d'une classe géométrique de morphismes et d'une notion de recouvrement). On voit ainsi se dégager une large famille de \emph{sites} : Zariski, étale, fppf, fpqc et, du côté différentiable, la topologie lisse (submersions).

On refait donc la théorie des faisceaux ``sur un site'' : les mêmes axiomes de séparation et recollement apparaissent, mais cette fois pour des cribles couvrants. On construit les foncteurs de faisceautisation comme les adjoints à gauche des foncteurs de plongement des faisceaux dans les préfaisceaux, on introduit les ``foncteurs canoniques'' composés des foncteurs de faisceautisation avec les foncteurs de Yoneda, et on montre la stabilité des constructions par passage au site localisé au-dessus d'un objet. Puis les sites annelés entrent en scène : un faisceau d'anneaux sur un site joue le rôle de faisceau de structures, tandis que les faisceaux de modules, les faisceaux localement libres et les faisceaux cohérents permettent de faire de l'algèbre dans ce contexte. On montre en détail comment les opérations algébriques habituelles (modules des homomorphismes de modules, produits tensoriels) se transplantent dans ce cadre général. Le chapitre bâtit un dictionnaire : presque tout ce qu'on sait faire sur les ouverts d'un espace ou sur les ouverts affines d'un schéma se fait, \emph{mutatis mutandis}, sur un site.

Enfin, on relie ces topologies ``abstraites'' aux géométries concrètes, notamment celles qui surgissent dans le contexte des schémas. Les exemples illustrent l’idée générale qu'une topologie de Grothendieck encode exactement ce qui suffit pour donner un contenu géométrique à une situation mathématique. 

\textbf{III. Définition et propriétés catégoriques des topos}

Ayant étudié au chapitre précédent la notion de faisceau sur un site, on monte au niveau supérieur dans l’abstraction : un topos est, par définition, une catégorie équivalente à la catégorie des faisceaux sur un site. Cette définition constructive des topos est équivalente d’après le théorème de Giraud à une définition axiomatique. On peut ``reconnaître'' si une catégorie est un topos par ses propriétés, semblables à celles de la catégorie des ensembles. Cela explique pourquoi les topos forment une classe robuste d’univers catégoriques dans lesquelles on peut développer des mathématiques comme dans l’univers familier des ensembles. 

Le chapitre démontre d’abord les propriétés générales partagées par tous les topos qui, selon l’expression de Grothendieck, en font des ``pastiches'' de la catégorie des ensembles : un topos est complet et cocomplet au sens que les limites et colimites y sont toujours bien définies, les foncteurs de produits fibrés respectent toutes les colimites, les foncteurs de colimites filtrantes respectent les limites finies, les sommes y sont disjointes, toutes les relation d’équivalence sont effectives au sens qu’elles définissent des quotients et réciproquement tous les quotients sont définis par des relations d’équivalence, les exponentielles existent ainsi qu’un ``classificateur des sous-objets'' qui représente le foncteur qui associe à tout objet l’ensemble de ses sous-objets, ce qui permet de parler de la ``logique interne'' propre à n’importe quel topos et de l’étudier comme une structure algébrique. On y montre notamment que les exponentielles se calculent explicitement dans un topos de faisceaux comme des faisceaux d’homomorphismes locaux, et que les sous-objets de chaque objet forment un treillis riche, muni en particulier d’opérateurs de pseudo-négation et d’implication qui sont respectées par les foncteurs d’images réciproques et induisent des opérations algébriques internes du classificateur des sous-objets .

La profondeur vient aussi de deux principes techniques essentiels. D'abord, le ``lemme de comparaison'' : pour n’importe quelle sous-catégorie pleine et dense d'un site (par exemple celle définie par une base d'ouverts d’un espace topologique), on ne change pas le topos des faisceaux associé si l’on remplace le site d’origine par cette sous-catégorie munie de la topologie induite. Ensuite, la ``présentation par objets représentables''b: tout faisceau est colimite d'objets représentables, ce qui permet souvent de ramener des preuves globales à des vérifications locales sur des ``cartes''.
Le chapitre développe enfin la théorie des objets linéaires internes, c’est-à-dire des faisceaux de modules sur un faisceau d'anneaux. Il montre en particulier l'existence d'assez d'objets ``injectifs'' et les conditions d'exactitude nécessaires au développement de la cohomologie – dans le cadre unifié et très général constitué par les sites. C’est un retour sur la genèse des topos, puisque c’est en dégageant le cadre catégorique le plus général permettant de faire de la cohomologie et en cherchant à unifier les classes d’exemples connus qui rentraient dans ce cadre – d’une part les catégories de faisceaux sur des espaces topologiques et d’autre part les catégories d’actions de groupes ou plus généralement de ``diagrammes'', c’est-à-dire de préfaisceaux – que Grothendieck a été conduit à introduire la notion générale de site, celle de faisceau sur un site et celle de topos, dont il a vérifié que les catégories d’objets linéaires satisfaisaient toujours les propriétés générales.

Enfin, le théorème de Giraud donne la caractérisation axiomatique des topos : toute catégorie localement petite, avec limites finies et colimites arbitraires, quotients effectifs des relations d'équivalence, colimites respectées par changements de base et une famille séparante d'objets, est équivalente à un topos. Mieux, pour toute famille séparante d’une telle catégorie, le site constitué par la sous-catégorie pleine sur cette famille et la topologie des familles globalement épimorphiques définit un topos des faisceaux qui est équivalent à la catégorie considérée. Comme corollaire, on voit que tout topos admet une infinie diversité de représentations comme topos des faisceaux sur des sites.

\textbf{IV. Géométrie des topos}

Une fois les topos définis et leurs propriétés internes dégagées, on fait de la géométrie. Le chapitre présente la théorie des morphismes de topos, dont la définition catégorique générale est dictée par le traduction de la notion d’application continue entre deux espaces topologiques en termes de foncteurs entre les catégories de faisceaux sur ces espaces. Un morphisme entre deux topos est défini comme une paire de foncteurs adjoints, appelés composante d’image réciproque et composante d’image directe, dont on demande que la première respecte aussi les limites finies.
Ces foncteurs permettent de transporter les objets et les structures internes d’un topos à un autre, pour la composante d’image directe toutes les structures qui s’expriment en termes de limites arbitraires et pour celle d’image d’image réciproque toutes les structures qui s’expriment en termes de colimites arbitraires et de limites finies. En particulier, toutes les structures de type algébrique, comme celles d’anneau interne ou de module interne, sont transportées dans les deux sens. Les morphismes de topos se composent. Un point d’un topos est par définition un morphisme vers ce topos défini sur le ``topos ponctuel'' constitué de la catégorie des ensembles.  Dans le cas des espaces topologiques ``sobres'' (comme sont tous les espaces que l’on rencontre dans la pratique mathématique usuelle, en particulier tous les espaces séparés et les espaces sous-jacents des schémas), les points des espaces correspondent bijectivement aux points des topos de faisceaux sur ces espaces, et les applications continues aux morphismes de topos associés.

On étudie particulièrement l’action des morphismes de topos sur les structures algébriques, la notion de morphisme de topos annelés et les actions d’un tel morphisme sur les catégories de modules internes. Cela prépare le terrain pour la cohomologie relative dans des cadres très généraux.

L’explicitation du lien profond entre les topos et la logique qui sera donnée au chapitre suivant est préparée par ``l'équivalence de Diaconescu'' qui fait, pour n’importe quel site, la théorie des morphismes des topos vers le topos des faisceaux sur ce site. Elle montre que les catégories de tels morphismes sont équivalentes aux catégories des foncteurs de la catégorie sous-jacente à ce site vers ces topos qui satisfont deux propriétés formulées au premier ordre : la propriété de ``platitude'' et celle de ``continuité'' qui signifie que le foncteur transforme les familles couvrantes de morphismes en des familles globalement épimorphiques des topos. L’équivalence de Diaconescu généralise un résultat de Grothendieck pour les sites dont la catégorie sous-jacente est ``cartésienne'', c’est-à-dire a toutes les limites finies : dans ce cas, la propriété de platitude d’un foncteur d’une telle catégorie sous-jacente vers un topos signifie qu’il respecte les limites finies. L’équivalence de Diaconescu donne notamment une description concrète des points d'un topos et de leurs morphismes.

On présente ensuite les foncteurs entre les catégories sous-jacentes de sites qui induisent des morphismes entre les topos associés. Suivant que le morphisme de topos induit va dans le sens contraire ou dans le même sens que le foncteur, on parle de morphisme ou de comorphisme de sites. On étudie la dualité des morphismes et des comorphismes de sites, et on montre comment passer d’une représentation d’un morphisme de topos par un morphisme de sites à une représentation par un comorphisme, ou l’inverse.

Puis est présentée la théorie des sous-topos, qui est un analogue topossique plus riche de la théorie des sous-espaces. On étudie les opérations internes d’intersection et de réunion des sous-topos, ainsi que celle de différence de deux sous-topos donc en particulier de passage au complémentaire d’un sous-topos, et les opérations externes d’image réciproque et d’image directe des sous-topos par un morphisme de topos. On montre que les sous-topos du topos des faisceaux sur un site correspondent bijectivement aux topologies sur la catégorie sous-jacente au site qui raffinent la topologie donnée. En termes intrinsèques, sans référence à aucune représentation par un site, les sous-topos d’un topos correspondent aux opérateurs internes du classificateur des sous-objets qui sont des ``opérateurs de fermeture'' au sens qu’ils satisfont deux propriétés de croissance et d’idempotence. 
On montre encore que toute localisation d’un topos relativement à l’un de ses objets est encore un topos et l’on étudie les morphismes entre topos localisés induits par les morphismes entre les objets qui les définissent. On appelle ouverts d’un topos ses sous-topos qui sont aussi des localisations, et on appelle fermés les sous-topos complémentaires des ouverts.

En somme, ce chapitre généralise aux topos toutes les définitions, notions et opérations fondamentales de la topologie et en donne des expressions concrètes, ce qui rend la théorie opérationnelle et propre au calcul.

\textbf{V. Théories, modèles, catégories syntaxiques et topos classifiants}

Le dernier chapitre expose comment les topos relient la topologie et la logique du premier ordre, principalement à travers la notion de topos classifiant découverte dans sa plus grande généralité par l'école de logique catégorique de Montréal (Makkai et Reyes \cite{CategoricalLogic}) à partir de premiers exemples donnés dans la thèse de Monique Hakim (\cite{TASR}) réalisée sous la direction de Grothendieck.
On introduit la notion de langage ou ``signature'' du premier ordre, constitué de ``sortes'', ``symboles de fonctions'' et ``symboles de relations'', c’est-à-dire de noms d’objets, de noms de morphismes d’un produit fini d’objets qui ont des noms dans un objet qui a un nom et de noms de sous-objets dans un produit fini d’objets qui ont des noms.  Puis on définit la notion sémantique de ``structure'' ou interprétation d’un tel langage dans un topos, qui consiste en des objets, des morphismes et des sous-objets nommés par les éléments de vocabulaire de ce langage. Une formule d’un tel langage est une écriture symbolique constituée récursivement à partir de ses éléments de vocabulaires et de variables formelles associées à ses ``sortes'' ou noms d’objets, par application répétée de substitutions et d’opérations logiques telles que la conjonction, la disjonction ou la quantification existentielle sur une partie des variables. Pour toute structure ou interprétation d’un tel langage dans un topos, toute formule de ce langage en une suite finie de variables associées à des sortes définit un sous-objet du produit des objets nommés par ces sortes. Deux formules en la même suite finie de variables définissent deux sous-objets du même objet produit. On dit que le ``séquent'' que forment ces deux formules est satisfait par la structure considérée si le sous-objet défini par la première formule est contenu dans celui défini par la seconde. Une théorie du premier ordre consiste en un langage et une famille de séquents, c’est-à-dire de paires de formules en les mêmes variables, qui sont appelés les axiomes de la théorie. Une structure ou interprétation d’un tel langage dans un topos est appelée un ``modèle'' de la théorie définie sur ce langage par une certaine famille d’axiomes si elle satisfait tous ces axiomes. Les modèles d’une théorie du premier ordre dans n’importe quel topos forment une catégorie. Si cette théorie est ``géométrique'', c’est-à-dire si les formules constitutives de ses axiomes ne font pas appel à d’autres symboles logiques que les conjonctions finies, les disjonctions et les quantifications existentielles, les composantes d’image réciproque des morphismes de topos respectent les interprétations de ces formules, qui se réalisent en termes de limites finies et de colimites arbitraires, et donc elles induisent des foncteurs entre les catégories de modèles d’une telle théorie dans les différents topos.

Ainsi, toute théorie géométrique du premier ordre définit un ``foncteur des modèles'' qui est un foncteur contravariant de la 2-catégorie des topos vers la 2-catégorie des catégories localement petites. Le théorème d’existence des topos classifiants dit que, pour toute telle théorie, son ``foncteur des modèles'' est représentable par un topos , unique à équivalence près, qui est appelé son ``topos classifiant''. Il est ainsi caractérisé par le fait qu’il possède un ``modèle universel'' de la théorie, au sens que la construction des images réciproques de ce ``modèle universel'' par les morphismes de n’importe quel topos vers le ``topos classifiant'' définit une équivalence de ca catégorie de ces morphismes vers celle des modèles de la théorie dans le topos source des morphismes.

Le ``topos classifiant'' d’une théorie géométrique du premier ordre est construit comme le topos des faisceaux sur une catégorie dite ``syntactique'' dont les objets sont les formules du langage de la théorie et dont les morphismes sont les formules dites ``démontrablement fonctionnelles'' entre les formules de source et de but, car leurs interprétations dans n’importe quels modèles de la théorie dans n’importe quel topos définissent des sous-objets qui sont automatiquement les graphes de morphismes entre les objets définis par les formules de source et de but. La catégorie syntactique est munie d’une topologie dite syntactique telle que les familles couvrantes sont celles qui s’interprètent dans tous les modèles de tous les topos comme des familles globalement épimorphiques.
Si l’on considère des théories dont les formules constitutives des axiomes appartiennent à des classes plus restrictives que la classe ``géométrique'' (telles que les classes ``cartésiennes'', ``régulières'' ou ``cohérentes''), leur topos classifiant se construit aussi comme topos des faisceaux sur des ``catégories syntactiques'' plus petites, constituées des seules formules de la classe considérée.
En particulier, si une théorie du premier ordre est algébrique ou plus généralement ``cartésienne'', comme par exemple la théorie des catégories, son topos classifiant se construit comme le topos des préfaisceaux sur sa ``catégorie syntactique cartésienne'', sans topologie.
La démonstration de ce que le topos des faisceaux sur la ``catégorie syntactique'' d’une théorie munie de sa ``topologie syntactique'' est un ``topos classifiant'' de cette théorie repose sur l’équivalence de Diaconescu, ou plus précisément sur son cas cartésien dû à Grothendieck puisque les catégories syntactiques ont toujours des limites finies arbitraires.
L’équivalence de Diaconescu montre aussi que, réciproquement, toute présentation d’un topos comme topos des faisceaux sur un site définit une théorie géométrique du premier ordre associée à ce site dont le topos classifiant s’identifie au topos considéré.
Ainsi, tout topos peut être présenté comme le topos classifiant d’une infinie diversité de théories géométriques du premier ordre. Deux théories sont dites ``sémantiquement équivalentes''  ou ``équivalentes au sens de Morita'' lorsque leurs topos classifiants sont équivalents ou, ce qui revient au même, lorsque les foncteurs des modèles qu’ils définissent sont équivalents.

Suivant le livre \cite{TST}, on montre encore que pour tout topos présenté comme le topos classifiant d’une certaine théorie, ses sous-topos correspondent bijectivement aux classes de ``théories quotient'' de la théorie considérée, c’est-à-dire aux théories définies dans le même langage en ajoutant des axiomes supplémentaires, considérées à équivalence démontrable près. Ce théorème implique que tout problème de démontrabilité en logique géométrique du premier ordre peut être traduit en un problème d’engendrement de topologies de Grothendieck par des cribles, sur la catégorie sous-jacente de n’importe quel site dont le topos des faisceaux contient comme sous-topos le topos classifiant de la théorie dans laquelle on se place.
Suivant toujours le livre \cite{TST}, on étudie particulièrement les théories ``de type préfaisceaux'', c’est-à-dire celles dont le topos classifiant peut être présenté comme le topos des préfaisceaux sur une catégorie, sans topologie. On montre en particulier que ces théories se caractérisent par une équivalence de leur syntaxe et de leur sémantique.
Dans une future version de ce manuscrit, on présentera la technique des ``topos comme ponts'' de \cite{TST}, qui se fonde sur la notion d'équivalence sémantique entre théories pour établir des connexions - souvent profondes et surprenantes – entre des contextes mathématiques différents, ainsi que des résultats récents sur les topos relatifs.

\vspace{0.5cm}
\textbf{Public visé et fil méthodologique}

L'ouvrage s'adresse à des lectrices et lecteurs ayant une familiarité raisonnable avec l'algèbre (groupes, anneaux, modules), la topologie générale, et les rudiments de théorie des catégories. Les parties I–II peuvent être lues comme un ``long échauffement'' : elles construisent progressivement la compréhension des notions essentielles (recouvrements abstraits, faisceaux, sites) et l’apprentissage de l'outillage (adjoints, limites et colimites, lemme de Yoneda) qui rendront naturelles les propriétés des topos (III) puis leur géométrie (IV) et leur logique (V).
La pédagogie s'appuie sur:
\begin{itemize}
	\item[$\bullet$] une montée en abstraction guidée par des exemples canoniques (topologie usuelle, topologie de Zariski, schémas affines et projectifs, recouvrements lisses ou étales, sites fpqc ou fppf),
	
	\item[$\bullet$] des résultats-clefs (lemme de comparaison, équivalence de Diaconescu, axiomes de Giraud),  
	
	\item[$\bullet$] des constructions universelles (faisceautisation, classificateur des sous-objets, exponentielles, topos classifiants).
\end{itemize}

L'accent est mis sur la réversibilité des points de vue : externe (faisceaux sur un site donné) et interne (logique du topos et structures internes), syntactique (catégories et topologies construites à partir d'une théorie) et sémantique (modèles dans des topos), local (changements de base, recouvrements, sites) et global (topos, invariants). Ces dualités, omniprésentes, constituent sans doute l'apport conceptuel le plus marquant de la théorie des topos.

\vspace{0.5cm}
\textbf{Pourquoi lire ce livre?}

Parce qu'il donne une vision unifiée où ``faire de la géométrie'' signifie organiser des principes de localisation (recouvrements), de cohérence (recollement), et de stabilité (changement de base), puis reconnaître les topos comme le cadre naturel où ces opérations deviennent structurelles. Parce qu'il montre que la logique géométrique n'est pas un appendice, mais une grammaire inhérente aux topos, apte à classifier des théories, à comparer des modèles, et à révéler des équivalences de nature sémantique.
Enfin, parce que la théorie des topos est un langage : après l'avoir apprise, on parle autrement des recouvrements, des techniques de descente, des morphismes de schémas, des modules, des espaces de solutions, des théories logiques. On y gagne une carte conceptuelle où les ``espaces'' sont aussi des catégories de faisceaux et où les ``théories'' vivent dans des topos qui les classifient. Ce livre en propose les clefs, en alternant concepts, exemples et principes universels, de sorte que la lectrice ou le lecteur puisse s'approprier ce langage et l'employer avec profit dans ses propres contextes mathématiques familiers et pour ses propres questionnements.

\tableofcontents

\input{Chapitre1_num.tex}

\input{Chapitre2_num.tex}

\input{Chapitre3_num.tex}

\input{Chapitre4_num.tex}

\input{Chapitre5_num.tex}

\vglue 1cm

\end{document}

%% file: Chapitre1_num.tex






\vglue 15mm

\chapter{Avant les topos: groupes, espaces topologiques, cat\'egories}\label{chap1}

\section{Les groupes}\label{sec11}

\subsection{La notion de groupe (et celle de mono{\"\i}de)}\label{subsec111}

On rappelle la d\'efinition suivante, qui est l'une des plus importantes des math\'ematiques:

\begin{defn}\label{defI11}
\begin{listeimarge}
\item Un groupe est un ensemble $G$ muni de

\medskip

$\left\{ \begin{matrix}
\bullet &\mbox{une loi de multiplication} \hfill \\
{ \ } \\
&\begin{matrix}
G \times G &\xrightarrow{ \ m \ } &G \, , \hfill \\
\hfill (g,g') &\xmapsto{ \ \ \ \ } &m(g,g') = g \cdot g' \, ,
\end{matrix} \\
{ \ } \\
\bullet &\mbox{un \'el\'ement ``unit\'e''} \hfill \\
{ \ } \\
&1 = 1_G \in G \, , \\
{ \ } \\
\bullet &\mbox{une application de passage \`a un ``inverse''} \hfill \\
{ \ } \\
&\begin{matrix}
G &\xrightarrow{ \, i = (\bullet)^{-1} \, } &G \, , \hfill \\
\hfill g &\xmapsto{ \ \ \ \ \ \ \ \ \ } &g^{-1} \, ,
\end{matrix}
\end{matrix} \right.
$

\medskip

\noindent tels que soient v\'erifi\'es les axiomes suivants:
\begin{enumerate}[label=(\Alph*)]
\item {\it Associativit\'e:}
$$
(g_1 \cdot g_2) \cdot g_3 = g_1 \cdot (g_2 \cdot g_3) \, , \quad \forall \, g_1 , g_2 , g_3 \in G \, .
$$
\item {\it Neutralit\'e:}
$$
1 \cdot g = g = g \cdot 1 \, , \quad \forall \, g \in G \, .
$$
\item {\it Inverse:}
$$
g \cdot g^{-1} = 1 = g^{-1} \cdot g \, , \qquad \forall \, g \in G \, .
$$
\end{enumerate}

\medskip

\item Un morphisme (ou homomorphisme) d'un groupe $G_1$ dans un groupe $G_2$ est une application
$$
\rho : G_1 \longrightarrow G_2
$$
telle que
\begin{eqnarray}
\rho (g \cdot g') &= &\rho (g) \cdot \rho (g') \, , \quad \forall \, g,g' \in G \, , \nonumber \\
\hfill \rho (1_{G_1}) &= &1_{G_2} \, . \hfill \nonumber
\end{eqnarray}
\end{listeimarge}
\end{defn}

\begin{remarkqed}

Un morphisme de groupes $\rho : G_1 \to G_2$ respecte \'egalement les passages \`a un inverse au sens que
$$
\rho (g^{-1}) = \rho (g)^{-1} \, , \quad \forall \, g \in G_1 \, .
$$
En effet, on a
$$
\rho (g^{-1}) \cdot \rho (g) = \rho (g^{-1} \cdot g) = \rho (1_{G_1}) = 1_{G_2}
$$
et donc
$$
\rho (g^{-1}) = \rho (g^{-1}) \cdot (\rho (g) \cdot \rho (g)^{-1}) = (\rho (g^{-1}) \cdot \rho (g)) \cdot \rho (g)^{-1} = \rho (g)^{-1} \, .
$$
\end{remarkqed}

\medskip

Les groupes sont des cas particuliers de mono{\"\i}des au sens suivant:

\begin{defn}\label{defI12}
\begin{listeimarge}
\item Un mono{\"\i}de est un ensemble $M$ muni de

\medskip

$\left\{ \begin{matrix}
\bullet &\mbox{une loi de multiplication} \hfill \\
{ \ } \\
&\begin{matrix}
M \times M &\xrightarrow{ \  \ } &M \, , \hfill \\
\hfill (m,m') &\xmapsto{ \ \ } &m \cdot m' \, ,
\end{matrix} \\
{ \ } \\
\bullet &\mbox{un \'el\'ement ``unit\'e'' $1 = 1_M \in M$,} \hfill
\end{matrix} \right.
$

\medskip

\noindent v\'erifiant les axiomes {\rm (A)} d'associativit\'e et {\rm (B)} de neutralit\'e de la d\'efinition~\ref{defI11}.

\item Un morphisme d'un mono{\"\i}de $M_1$ dans un mono{\"\i}de $M_2$ est une application $\rho : M_1 \to M_2$ telle que
\begin{eqnarray}
\rho (m \cdot m') &= &\rho (m) \cdot \rho (m') \, , \quad \forall \, m,m' \in M \, , \nonumber \\
\hfill \rho (1_{M_1}) &= &1_{M_2} \, . \hfill \nonumber
\end{eqnarray}
\end{listeimarge}
\end{defn}

\begin{remarksqed}
\begin{listeisansmarge}
\item Tout groupe est un mono{\"\i}de.

\item Pour tout \'el\'ement $m$ d'un mono{\"\i}de $M$, il existe au plus un \'el\'ement $m'$ de $M$ tel que $m \cdot m' = 1 = m' \cdot m$. 

\smallskip

En effet, si $m''$ est un autre \'el\'ement v\'erifiant $m \cdot m'' = 1$, on a
$$
m'' = 1 \cdot m'' = (m' \cdot m) \cdot m'' = m' \cdot (m \cdot m'') = m' \cdot 1 = m' \, .
$$

S'il existe un tel \'el\'ement $m'$, on dit que $m$ est un \'el\'ement inversible de $M$ d'inverse $m' = m^{-1}$.

\smallskip

Les \'el\'ements inversibles de $M$ forment un groupe not\'e $M^{\times}$.

\medskip

\item Etant donn\'e deux groupes $G_1$ et $G_2$, les morphismes de mono{\"\i}des de $G_1$ dans $G_2$ sont les morphismes de groupes de $G_1$ dans $G_2$.

\medskip

\item Tout morphisme de mono{\"\i}des $\rho : M_1 \to M_2$ induit un morphisme de groupes $\rho : M_1^{\times} \to M_2^{\times}$.

\medskip

\item Le compos\'e de deux morphismes de mono{\"\i}des $M_1 \to M_2 \to M_3$ est un morphisme de mono{\"\i}des $M_1 \to M_3$.

\smallskip

En particulier, le compos\'e de deux morphismes de groupes $G_1 \to G_2 \to G_3$ est un morphisme de groupes $G_1 \to G_3$. 
\end{listeisansmarge}
\end{remarksqed}

\subsection{Exemples de mono{\"\i}des}\label{subsec112}

\noindent $\bullet$ {\bf Les mono{\"\i}des additifs et multiplicatifs:}

\smallskip

Dans un anneau $A$, toute partie contenant $0$ qui est stable par l'addition est un mono{\"\i}de. De m\^eme, toute partie contenant $1$ qui est stable par la multiplication est un mono{\"\i}de.

\smallskip

Par exemple, dans $A = {\mathbb Z}$, ${\mathbb N}$ est un mono{\"\i}de additif. Et, pour tout ensemble $P$ de nombres premiers, l'ensemble des entiers dont la d\'ecomposition en facteurs premiers ne fait appara{\^\i}tre que des \'el\'ements de $P$ est un mono{\"\i}de multiplicatif.

\medskip

\noindent $\bullet$ {\bf Les mono{\"\i}des de matrices:}

\smallskip

Pour tout anneau $A$, l'ensemble $M_n (A)$ des applications $A$-lin\'eaires $A^n \to A^n$ muni de la loi de composition est un mono{\"\i}de.

\medskip

\noindent $\bullet$ {\bf Les mono{\"\i}des d'applications:}

\smallskip

Pour tout ensemble $I$, l'ensemble ${\rm End} (I)$ des applications $I \to I$ muni de la loi de composition est un mono{\"\i}de.

\subsection{Exemples de groupes}\label{subsec113}

La plupart des groupes naturels des math\'ematiques apparaissent comme groupes de sym\'etries, c'est-\`a-dire de transformations de divers objets qui respectent leurs structures.

\medskip

\noindent $\bullet$ {\bf Les groupes de sym\'etries d'un plan affine:}

\smallskip

Consid\'erant l'ensemble ${\mathcal P}$ des points et celui ${\mathcal D}$ des droites d'un plan affine, munis de la relation $\in$ d'appartenance d'un point \`a une droite et de la relation $\sslash$ de parall\'elisme entre droites, on dispose des groupes de sym\'etries en les sept sens suivants:
\begin{enumerate}
\item[--] les ``transformations affines'' d\'efinies comme les paires de bijections de ${\mathcal P}$ et ${\mathcal D}$ qui respectent les relations $\in$ et $\sslash$ vues comme des sous-ensembles de ${\mathcal P} \times {\mathcal D}$ et ${\mathcal D} \times {\mathcal D}$,
\item[--] les ``d\'eplacements'' d\'efinis comme les transformations affines qui transforment toute droite $D$ en une droite parall\`ele \`a $D$,
\item[--] les ``translations de direction $D$'' d\'efinies comme les d\'eplacements qui transforment tout point $A$ en un point $A'$ de l'unique droite passant par $A$ et parall\`ele \`a une droite $D$ fix\'ee,
\item[--] les ``translations'' d\'efinies comme les d\'eplacements qui commutent avec les ``translations de direction $D$'' pour toute droite $D$,
\item[--] les ``homoth\'eties de centre $P$'' d\'efinies comme les d\'eplacements qui fixent un \'el\'ement $P$ de ${\mathcal P}$,
\item[--] les ``dilatations'' d\'efinies comme les homoth\'eties de centre un \'el\'ement $O$ de ${\mathcal P}$ appel\'e le ``point origine'',
\item[--] les ``transformations lin\'eaires'' d\'efinies comme les transformations affines qui commutent avec les dilatations.
\end{enumerate}

\medskip

\noindent $\bullet$ {\bf Les groupes de sym\'etries d'un plan euclidien:}

\smallskip

Consid\'erant un plan euclidien comme un plan affine muni de surcro{\^\i}t d'une relation $\perp$ d'orthogonalit\'e entre droites, on dispose des groupes suppl\'ementaires de sym\'etries en les sens suivants:
\begin{enumerate}
\item[--] les similitudes d\'efinies comme les transformations affines qui respectent la relation $\perp$ vue comme un sous-ensemble de ${\mathcal D} \times {\mathcal D}$,
\item[--] les similitudes lin\'eaires d\'efinies comme les transformations lin\'eaires (relativement \`a un point origine $O$ fix\'e) qui sont des similitudes,
\item[--] les transformations orthogonales d\'efinies comme les similitudes lin\'eaires qui respectent n'importe quel produit scalaire qui d\'efinit la relation d'orthogonalit\'e $\perp$ ou, ce qui revient au m\^eme, dont le d\'eterminant est $1$ ou $-1$,
\item[--] les rotations d\'efinies comme les transformations orthogonales dont le d\'eterminant est $1$.
\end{enumerate}

\medskip

\noindent $\bullet$ {\bf Les groupes additifs et multiplicatifs:}

\smallskip

\`A tout anneau $A$ sont associ\'es les groupes additifs
$$
{\mathbb A}^n (A) = (A^n , +) \, , \quad n \geq 1 \, ,
$$
ainsi que le groupe multiplicatif
$$
A^{\times} = {\mathbb G}_m (A) = {\rm GL}_1 (A) \, .
$$

\medskip

\noindent $\bullet$ {\bf Les groupes lin\'eaires:}

\smallskip

\`A tout anneau $A$ sont associ\'es les groupes de matrices inversibles
$$
{\rm GL}_n (A) = M_n (A)^{\times} \, , \quad n \geq 1 \, .
$$

\medskip

\noindent $\bullet$ {\bf Les groupes orthogonaux et symplectiques:}

\smallskip

Si $A$ est un anneau commutatif et $A^n$ [resp. $A^{2n}$] est muni d'un produit scalaire sym\'etrique [resp. antisym\'etrique au sens qu'il s'annule sur la diagonale de $A^{2n} \times A^{2n}$] et non d\'eg\'en\'er\'e (au sens que sa matrice est inversible), les applications $A$-lin\'eaires $A^n \to A^n$ [resp. $A^{2n} \to A^{2n}$] qui respectent ce produit scalaire forment un groupe not\'e
$$
O_n (A) \qquad \mbox{[resp. $Sp_{2n} (A)$]}
$$
et appel\'e groupe orthogonal [resp. groupe symplectique].

\medskip

\noindent $\bullet$ {\bf Les groupes de bijections:}

\smallskip

Pour tout ensemble $I$, les bijections $I \xrightarrow{ \ \sim \ } I$ forment un groupe que l'on peut noter ${\mathfrak S}_I$.

\smallskip

En particulier, le groupe ${\mathfrak S}_n$ des bijections $\{1,2,\cdots , n\} \to \{1,2,\cdots , n\}$ est appel\'e le groupe sym\'etrique de degr\'e $n$.

\medskip

\noindent $\bullet$ {\bf Les groupes de Galois:}

\smallskip

\`A toute extension alg\'ebrique $L$ d'un corps $K$ (c'est-\`a-dire tout corps commutatif $L$ contenant $K$ et dont tout \'el\'ement est racine d'un polyn\^ome non nul \`a coefficients dans $K$) est associ\'e son ``groupe de Galois''
$$
{\rm Gal} (L/K) = {\rm Aut}_K(L) = \left\{ \sigma : L \to L \ \Biggl\vert \ \begin{matrix}
\sigma (\ell_1 + \ell_2) = \sigma (\ell_1)  + \sigma (\ell_2) \, , &\forall \, \ell_1 , \ell_2 \in L \, , \\
\sigma (\ell_1 \, \ell_2) = \sigma(\ell_1) \, \sigma(\ell_2) \, , \hfill &\forall \, \ell_1 , \ell_2 \in L \, , \\
\sigma (k) = k \, , \hfill &\forall \, k \in K \hfill
\end{matrix} \right\}.
$$

\medskip

\noindent $\bullet$ {\bf Les groupes fondamentaux de Poincar\'e:}

\smallskip

Dans un espace topologique $X$, un lacet d'un point $x$ \`a un point $y$ est une application continue $\sigma : [0,1] \to X$ telle que $\sigma (0) = x$ et $\sigma (1) = y$. Deux lacets $\sigma_0 , \sigma_1$ de $x$ \`a $y$ sont dits ``homotopes'' s'il existe une application continue $h : [0,1] \times [0,1] \to X$ telle que $\sigma_0 (t) = h(0,t)$, $\sigma_1 (t) = h(1,t)$, $\forall \, t$ et $h(s,0) = x$, $h(s,1) = y$, $\forall \, s$. La relation d'homotopie est une relation d'\'equivalence si bien que l'on peut introduire l'ensemble $\pi_X (x,y)$ des classes d'homotopie de lacets de $x$ \`a $y$.

\smallskip

Le compos\'e d'un lacet $\sigma$ de $x$ \`a $y$ et d'un lacet $\tau$ de $y$ \`a $z$ est d\'efini comme le lacet de $x$ \`a $z$
$$
\tau \circ \sigma : t \longmapsto \left\{ \begin{matrix}
\sigma (2t) \hfill &\mbox{si} &0 \leq t \leq 1/2 \ , \\
\tau (2t-1) &\mbox{si} &1/2 \leq t \leq 1 \, .
\end{matrix} \right.
$$

La composition est compatible avec les relations d'homotopie et d\'efinit une application
\begin{eqnarray}
\pi_X (y,z) \times \pi_X (x,y) &\longrightarrow &\pi_X (x,z) \, , \nonumber \\
(\tau , \sigma) &\longmapsto &\tau \circ \sigma \, .  \nonumber
\end{eqnarray}
En particulier, pour tout point $x$ de $X$, l'ensemble
$$
\pi_X (x,x)
$$
est muni d'une loi de composition interne.

\smallskip

On v\'erifie qu'elle est associative, poss\`ede un \'el\'ement neutre repr\'esent\'e par le lacet constant $t \mapsto x$ et que tout \'el\'ement repr\'esent\'e par un lacet $\sigma$ y poss\`ede un inverse repr\'esent\'e par le lacet $t \mapsto \sigma (1-t)$.

\smallskip

Ainsi, $\pi_X (x,x)$ est un groupe appel\'e le groupe fondamental, ou groupe de Poincar\'e, de $X$ au point $x$.

\section{Expression des groupes: leurs actions}\label{sec12}

\subsection{La notion d'action d'un groupe (ou d'un mono{\"\i}de)}\label{subsec121}

\begin{defn}\label{defI21}

Soit $G$ un groupe (ou un mono{\"\i}de).

\begin{listeimarge}

\item Une action de $G$, ou un $G$-ensemble, consiste en un ensemble $X$ muni d'une application
\begin{eqnarray}
G \times X &\longrightarrow &X \, , \nonumber \\
(g,x) &\longmapsto &g \cdot x \nonumber
\end{eqnarray}
v\'erifiant les axiomes

\medskip

$\left\{ \begin{matrix}
\bullet &g_1 \cdot (g_2 \cdot x) = (g_1 \cdot g_2) \cdot x \, , &\forall \, g_1 , g_2 \in G \, , \ \forall \, x \in X \, , \hfill \\
\bullet &1 \cdot x = x \, , \hfill &\forall \, x \in X \, . \hfill
\end{matrix} \right.
$

\medskip

\item Un morphisme, ou application $G$-\'equivariante, entre deux actions $X,Y$ est une application
$$
u : X \longrightarrow Y
$$
telle que
$$
u(g \cdot x) = g \cdot u(x) \, , \qquad \forall \, g \in G \, , \ \forall \, x \in X \, .
$$
\end{listeimarge}
\end{defn}

\begin{remarkqed}

La compos\'ee de deux applications $G$-\'equivariantes
$$
X_1 \longrightarrow X_2 \longrightarrow X_3
$$
est $G$-\'equivariante.
\end{remarkqed}

\begin{lem}\label{lemI22}

Soit $G$ un groupe agissant sur un ensemble $X$. Alors:

\begin{listeimarge}

\item Pour tout \'el\'ement $x$ de $X$,
$$
G_x = \{ g \in G \mid g \cdot x = x \}
$$
est un sous-groupe de $G$, appel\'e le fixateur de l'\'el\'ement $x$.

\medskip

\item La relation entre paires d'\'el\'ements de $X$
$$
\{ (x_1 , x_2) \in X \times X \mid \exists \, g \in G \, , \ g \cdot x_1 = x_2 \}
$$
est une relation d'\'equivalence.

\smallskip

Les classes de $X$ pour cette relation d'\'equivalence sont appel\'ees les orbites de $X$ sous l'action de $G$.

\smallskip

Leur ensemble peut \^etre not\'e $G \backslash X$ et appel\'e le quotient de $X$ sous l'action de $G$.
\end{listeimarge}
\end{lem}

\begin{remarksqed}
\begin{listeisansmarge}
\item Une action est dite ``transitive'' si elle a une seule orbite.


\item Toute application $G$-\'equivariante
$$
u : X \longrightarrow Y
$$
induit une aplication
$$
G \backslash X \longrightarrow G \backslash Y \, .
$$
De plus, pour tout \'el\'ement $x$ de $X$, on a l'inclusion entre fixateurs
$$
G_x \subset G_{u(x)} \, .
$$
\end{listeisansmarge}
\end{remarksqed}

\subsection{Exemples d'actions de groupes (ou de mono{\"\i}des)}\label{subsec122}

L'int\'er\^et des groupes r\'eside dans leurs actions et c'est toujours pour elles qu'ils sont \'etudi\'es. De plus, tout groupe d\'efini comme le groupe des sym\'etries d'un certain objet est muni d'une action canonique sur cet objet.

\medskip

\noindent $\bullet$ {\bf L'action multiplicative d'un anneau sur lui-m\^eme:}

\smallskip

Un anneau $A$ consid\'er\'e comme un mono{\"\i}de multiplicatif agit sur lui-m\^eme par la multiplication.

\smallskip

Cette action respecte la loi d'addition de $A$, son \'el\'ement neutre $0$ et donc aussi le passage \`a l'oppos\'e $a \mapsto -a$. C'est l'axiome de ``distributivit\'e'' de la multiplication par rapport \`a l'addition.

\medskip

\noindent $\bullet$ {\bf L'action d'un groupe sur lui-m\^eme par conjugaison:}

\smallskip

Un groupe $G$ agit sur lui-m\^eme par
\begin{eqnarray}
G \times G &\longrightarrow &G \, , \nonumber \\
(g,h) &\longmapsto &g \cdot h \cdot g^{-1} \, . \nonumber
\end{eqnarray}
Cette action respecte la structure de groupe de $G$.

\smallskip

Les orbites de cette action sont appel\'ees classes de conjugaison et ses fixateurs sont les sous-groupes de commutateurs
$$
G_h = \{ g \in G \mid g \cdot h = h \cdot g \} \, .
$$

\medskip

\noindent $\bullet$ {\bf L'action canonique des applications et des bijections:}

\smallskip

Pour tout ensemble $I$, le mono{\"\i}de ${\rm End} (I)$ des applications $I \to I$ et le groupe ${\mathfrak S}_I = {\rm End} (I)^{\times}$ de ses bijections agissent sur $I$.

\smallskip

Ils agissent aussi sur l'ensemble ${\mathcal P} (I)$ des parties de $I$ par
$$
(f,J) \longmapsto f(J) = \{ i \in I \mid \exists \, j \in J \, , \ f(j) = i \} \, .
$$

\medskip

\noindent $\bullet$ {\bf L'action sur un plan affine de ses sym\'etries:}

\smallskip

Tous les groupes de sym\'etries d'un plan affine agissent sur l'ensemble ${\mathcal P}$ de ses points et l'ensemble ${\mathcal D}$ de ses droites. Ils sont m\^eme tous d\'efinis comme des fixateurs relativement \`a des actions d\'eriv\'ees de ces deux actions:
\begin{enumerate}
\item[--] les ``transformations affines'' sont les \'el\'ements de ${\mathfrak S}_{\mathcal P} \times {\mathfrak S}_{\mathcal D}$ qui fixent les parties $\in$ et $\sslash$ de ${\mathcal P} \times {\mathcal D}$ et ${\mathcal D} \times {\mathcal D}$,
\item[--] les ``d\'eplacements" sont les transformations affines qui fixent les classes d'\'equivalence de ${\mathcal D}$ pour la relation de parall\'elisme,
\item[--] les ``translations de direction $D$'' sont les d\'eplacements qui fixent les droites parall\`eles \`a $D$,
\item[--] les ``translations'' sont les d\'eplacements qui fixent les graphes des ``translations de direction $D$'' pour toute droite $D$,
\item[--] les ``homoth\'eties de centre $P$'' [resp. les dilatations] sont les d\'eplacements qui fixent un \'el\'ement $P$ de ${\mathcal P}$ [resp. le point origine $O$],
\item[--] les ``transformations lin\'eaires'' sont les transformations affines qui fixent les graphes des dilatations.
\end{enumerate}

\medskip

\noindent $\bullet$ {\bf L'action sur un plan euclidien de ses sym\'etries:}

\smallskip

Les transformations affines d'un plan affine agissent sur l'ensemble de ses relations d'orthogonalit\'e, et ses transformations lin\'eaires agissent sur l'ensemble de ses produits scalaires sym\'etriques non d\'eg\'en\'er\'es.

\smallskip

Les similitudes sont les transformations affines qui fixent une relation d'orthogonalit\'e $\perp$.

\smallskip

De m\^eme, les transformations orthogonales sont les transformations lin\'eaires qui fixent un produit scalaire de d\'efinition de $\perp$.

\smallskip

Elles agissent sur l'ensemble \`a deux \'el\'ements des orientations du plan euclidien.

\smallskip

Celles qui fixent les orientations sont les rotations.

\medskip

\noindent $\bullet$ {\bf Les actions canoniques des groupes lin\'eaires, orthogonaux et symplectiques:}

\smallskip

Pour tout anneau $A$, le mono{\"\i}de $M_n (A)$ et le groupe ${\rm GL}_n (A) = M_n (A)^{\times}$ agissent sur $A^n$.

\smallskip

Si $A$ est commutatif [resp. et si $n$ est pair], ${\rm GL}_n (A)$ agit sur l'ensemble des produits scalaires sym\'etriques [resp. antisym\'etriques] non d\'eg\'en\'er\'es de $A^n$. Les sous-groupes fixateurs de ces produits scalaires sont les groupes orthogonaux $O_n (A)$ [resp. les groupes symplectiques $Sp_n (A)$].

\medskip

\noindent $\bullet$ {\bf Les actions d'un groupe de Galois sur les ensembles de racines d'un polyn\^ome:}

\smallskip

Si $L$ est une extension alg\'ebrique d'un corps commutatif $K$, et que $P \in K[X]$ est un polyn\^ome \`a coefficients dans $K$, le groupe de Galois ${\rm Gal} (L/K) = {\rm Aut}_K(L)$ agit sur l'ensemble des \'el\'ements de $L$ qui sont racines de $P$.

\medskip

\noindent $\bullet$ {\bf L'action d'un groupe de Poincar\'e sur les rev\^etements:}

\smallskip

Un rev\^etement d'un espace topologique $X$ est un espace topologique $X'$ muni d'une application continue $p : X' \to X$ telle que $X$ admette un recouvrement par des ouverts $U_i$ au-dessus desquels chaque $p^{-1} (U_i)$ muni de $p$ soit hom\'eomorphe au produit de $U_i$ et d'un ensemble discret (autrement dit, \`a une somme disjointe de copies de $U_i$).

\smallskip

Si $X' \xrightarrow{ \ p \ } X$ est un tel rev\^etement, $[0,1] \xrightarrow{ \ \sigma \ } X$ est un lacet d'un point $x$ de $X$ \`a un point $y$, et $x'$ est un point de la fibre $p^{-1} (x)$ de $X'$ au-dessus de $x$, il existe un unique lacet $[0,1] \xrightarrow{ \ \sigma' \ } X'$ tel que $p\circ \sigma' = \sigma$ et $\sigma' (0) = x'$. De plus, le point $\sigma' (1) = y'$ de $p^{-1} (y)$ ne d\'epend que de la classe d'homotopie de $\sigma$.

\smallskip

Ainsi, tout \'el\'ement de $\pi_X (x,y)$ d\'efinit une application de la fibre $p^{-1} (x)$ dans la fibre $p^{-1} (y)$.

\smallskip

Si $x=y$, cela d\'efinit une action du groupe de Poincar\'e
$$
\pi_X (x,x)
$$
sur la fibre $p^{-1} (x)$.

\section{Les espaces topologiques}\label{sec13}

\subsection{La notion d'espace topologique}\label{subsec131}

\medskip

On rappelle maintenant la d\'efinition suivante, \'egalement l'une des plus importantes des math\'ematiques:

\begin{defn}\label{defI31}

\begin{listeimarge}

\item Un espace topologique est un ensemble $X$ muni d'une famille de sous-ensembles de $X$, appel\'es ses ``ouverts'', telle que:

\medskip

\noindent $\left\{ \begin{matrix}
\bullet &\mbox{le sous-ensemble vide $\emptyset$ et le sous-ensemble total $X$ de $X$ sont ouverts,} \hfill \\
\bullet &\mbox{pour toute famille $(U_i)_{i \in I}$ d'ouverts de $X$, leur r\'eunion $\underset{i \in I}{\bigcup} U_i$ est un ouvert de $X$,} \hfill \\
\bullet &\mbox{pour toute famille finie d'ouverts $U_1 , \cdots , U_n$ de $X$, leur intersection $U_1 \cap \cdots \cap U_n$ est un ouvert de $X$.} \hfill
\end{matrix} \right.
$

\medskip

\item Une application continue (ou morphisme d'espaces topologiques) entre deux espaces topologiques $X$ et $Y$ est une application
$$
f : X \longrightarrow Y
$$
telle que, pour tout ouvert $V$ de $Y$, son image r\'eciproque
$$
f^{-1} V = \{ x \in X \mid f(x) \in V \}
$$
soit un ouvert de $X$.
\end{listeimarge}
\end{defn}

\begin{remarks}
\begin{listeisansmarge}
\item La compos\'ee de deux applications continues
$$
X_1 \longrightarrow X_2 \longrightarrow X_3
$$
est une application continue.

\medskip

\item Un sous-ensemble $Z$ d'un espace topologique $X$ est appel\'e un ``ferm\'e'' si son compl\'ementaire
$$
X \backslash Z = X - Z \quad \mbox{est un ouvert.}
$$
Ainsi, $X$ et $\emptyset$ sont \`a la fois ouverts et ferm\'es, toute intersection de ferm\'es est un ferm\'e, et toute r\'eunion finie de ferm\'es est un ferm\'e.

\medskip

\item Par cons\'equent, toute partie $X'$ d'un espace topologique $X$ est contenue dans un plus petit ferm\'e $\overline{X'}$ de $X$, appel\'e l'adh\'erence de $X'$ dans $X$, qui est l'intersection de tous les ferm\'es de $X$ contenant $X'$.

\medskip

\item Une structure d'espace topologique sur un ensemble donn\'e $X$ est appel\'ee une topologie de $X$.

\smallskip

L'ensemble des topologies sur un ensemble $X$ est ordonn\'e par l'inclusion.

\smallskip

Toute intersection de topologies sur un ensemble $X$ est encore une topologie.

\smallskip

Par cons\'equent, toute famille de topologies $({\mathcal T}_i)_{i \in I}$ sur $X$ engendre une topologie ${\mathcal T}$, qui est l'intersection de toutes les topologies qui contiennent chaque topologie ${\mathcal T}_i$, $i \in I$.

\medskip

\item Si $Y \xrightarrow{ \ f \ } X$ est une application d'un ensemble $Y$ dans un espace topologique $X$, on peut munir $Y$ de la plus petite topologie telle que
$$
f : Y \longrightarrow X \quad \mbox{soit une application continue.}
$$
Dans cette topologie, les ouverts de $Y$ sont les images r\'eciproques $f^{-1} U$ des ouverts $U$ de $X$.

\smallskip

Si $Y \xrightarrow{ \ f \ } X$ est injective, $Y$ est appel\'e un sous-espace de $X$.

\medskip

\item Si $X \xrightarrow{ \ f \ } Y$ est une application d'un espace topologique $X$ dans un ensemble $Y$, on peut munir $Y$ de la plus grande topologie pour laquelle
$$
f : X \longrightarrow Y \quad \mbox{soit une application continue.}
$$
Dans cette topologie, une partie $V$ de $Y$ est un ouvert si et seulement si $f^{-1} V$ est un ouvert de $X$.

\smallskip

Si $X \xrightarrow{ \ f \ } Y$ est surjective, $Y$ est appel\'e un espace quotient de $X$.
\end{listeisansmarge}
\end{remarks}

\subsection{Exemples d'espaces topologiques}\label{subsec132}

\noindent $\bullet$ {\bf Les espaces m\'etriques:}

\smallskip

On rappelle qu'une ``m\'etrique'' ou fonction ``distance'' sur un ensemble $X$ est une application
$$
d : X \times X \longrightarrow {\mathbb R}_+
$$
telle que
$$
\left\{ \begin{matrix}
d(x_1 , x_2) = d(x_2 , x_1) \, , \hfill &\forall \, x_1 , x_2 \in X \hfill \\
d(x_1 , x_3) \leq d (x_1 , x_2) + d(x_2 , x_3) \, , \hfill &\forall \, x_1 , x_2 , x_3 \in X \, , \hfill \\
d(x_1 , x_2) = 0 \Leftrightarrow x_1 = x_2 \, , \hfill &\forall \, x_1 , x_2 \in X \, . \hfill
\end{matrix} \right.
$$

\smallskip

\noindent Toute m\'etrique d\'efinit une topologie:

\begin{defn}\label{defI32}

Une partie $U$ d'un espace m\'etrique $(X,d)$ est d\'eclar\'ee ouverte pour la topologie induite si, pour tout $x \in U$, il existe un nombre r\'eel $\varepsilon > 0$ tel que
$$
\forall \, x' \, , \quad d(x,x') < \varepsilon \Rightarrow x' \in U \, .
$$
\end{defn}

\begin{remarksqed}
\begin{listeisansmarge}
\item Etant donn\'es deux espaces m\'etriques $(X,d)$ et $(Y,d')$, une application
$$
f : X \longrightarrow Y
$$
est continue pour les topologies induites si et seulement si, pour tout $x \in X$ et tout $\varepsilon > 0$, il existe $\eta > 0$ tel que
$$
d(x,x') < \eta \Rightarrow d' (f(x) , f(x')) < \varepsilon \, , \quad \forall \, x,x' \in X \, .
$$

\item La topologie d'un sous-espace $X'$ d'un espace m\'etrique $X$ est induite par la restriction \`a $X'$ de la m\'etrique de $X$.

\smallskip

En particulier, tout sous-espace d'un espace m\'etrique est encore un espace m\'etrique.

\medskip

\item La topologie d'un espace m\'etrique est toujours ``s\'epar\'ee'' au sens que pour tous \'el\'ements $x_1 \ne x_2$, il existe des ouverts $U_1 , U_2$ tels que
$$
x_1 \in U_1 \, , \ x_2 \in U_2 \quad \mbox{et} \quad U_1 \cap U_2 = \emptyset \, .
$$

\item Le quotient d'un espace m\'etrique n'est pas n\'ecessairement s\'epar\'e, par exemple celui d\'efini par l'application surjective
\begin{eqnarray}
{\mathbb R} &\longrightarrow &\{0,1\} \, , \nonumber \\
t &\longmapsto &\left\{ \begin{matrix}
0 &\mbox{si} &t=0 \, , \\
1 &\mbox{si} &t \ne 0 \, .
\end{matrix} \right. \nonumber
\end{eqnarray}

A fortiori, le quotient d'un espace m\'etrique n'est pas n\'ecessairement un espace m\'etrique.
\end{listeisansmarge}
\end{remarksqed}

\medskip

\noindent $\bullet$ {\bf Les espaces d'Alexandrov:}

\smallskip

On rappelle qu'un ordre sur un ensemble $X$ est une relation binaire $\leq$ qui est transitive (au sens que $x_1 \leq x_3$ si $x_1 \leq x_2$ et $x_2 \leq x_3$), r\'eflexive (au sens que $x \leq x$, $\forall \, x$) et anti-sym\'etrique (au sens que $x_1 = x_2$ si $x_1 \leq x_2$ et $x_2 \leq x_1$).

\smallskip

Un pr\'e-ordre est une relation binaire $\leq$ qui est \'egalement transitive et r\'eflexive mais pas n\'ecessairement anti-sym\'etrique.

\begin{defn}\label{defI33}

Un espace topologique $X$ est appel\'e espace d'Alexandrov s'il v\'erifie les deux propri\'et\'es \'equivalentes suivantes:

\begin{enumerate}[label=(\arabic*)]

\item Toute intersection d'ouverts (pas n\'ecessairement finie) est un ouvert.


\item Il existe un (unique) pr\'e-ordre $\leq$ sur $X$ tel qu'une partie $U$ de $X$ soit un ouvert si et seulement si, pour tout $x \in U$, 
$$
x \leq x' \Rightarrow x' \in U \, .
$$
\end{enumerate}
\end{defn}

\begin{remarksqed}
\begin{listeisansmarge}
\item Le pr\'e-ordre de $X$ est d\'efini \`a partir de sa topologie en d\'ecidant que
$$
x \leq x' \quad \mbox{si et seulement si} \quad x \in \overline{x'} \, ,
$$
c'est-\`a-dire si et seulement si tout ouvert qui contient $x$ contient \'egalement $x'$.

\medskip

\item Une application $f : X \to Y$ entre deux espaces d'Alexandrov est continue si et seulement si $x \leq x' \Rightarrow f(x) \leq f(x')$. 
\end{listeisansmarge}
\end{remarksqed}

\medskip

\noindent $\bullet$ {\bf Le spectre d'un anneau commutatif:}

\smallskip

On rappelle qu'une partie $I$ d'un anneau commutatif $A$ est appel\'ee un id\'eal si

\medskip

$\left\{ \begin{matrix}
\bullet &a_1 , a_2 \in I \Rightarrow a_1 + a_2 \in I \, , \hfill \\
\bullet &a \in I \, , \ b \in A \Rightarrow ab \in I \, . \hfill
\end{matrix} \right.
$

\medskip

\noindent Un id\'eal $I$ est dit premier si $1 \notin I$ et que, pour tous $a,b \in A$, 
$$
ab \in I
$$
si et seulement si $a \in I$ ou $b \in I$ (ou, ce qui revient au m\^eme, si l'anneau quotient $A/I$ est int\`egre au sens que les produits d'\'el\'ements non nuls sont toujours non nuls).

\newpage

\begin{defn}\label{defI34}

\begin{listeimarge}
\item Pour tout anneau commutatif $A$, on appelle ``spectre'' de $A$ l'ensemble ${\rm Spec} (A)$ de ses id\'eaux premiers.

\smallskip

Une partie $U$ de ${\rm Spec} (A)$ est appel\'ee un ouvert de la topologie de Zariski s'il existe une famille $(f_i)_{i \in I}$ d'\'el\'ements de $A$ telle que, pour tout id\'eal premier $p$,
$$
p \in U \Leftrightarrow \exists \, i \in U \, , \ f_i \notin p \, .
$$

\item Pour tout morphisme d'anneaux commutatifs
$$
u : A \longrightarrow B \, ,
$$
on note ${\rm Spec} (u) : {\rm Spec} (B) \to {\rm Spec} (A)$ l'application continue
$$
q \longmapsto u^{-1} (q) \, .
$$
\end{listeimarge}
\end{defn}

\begin{remarksqed}
\begin{listeisansmarge}
	
\item Par d\'efinition, les ouverts de ${\rm Spec} (A)$ sont les r\'eunions arbitraires de parties de la forme
$$
{\rm Spec} (A)_f = \{ p \in {\rm Spec} (A) \mid f \notin p \} \, , \quad f \in A \, .
$$

Leurs intersections finies sont encore des ouverts car, pour tous \'el\'ements $f_1 , \cdots , f_n \in A$, on a
$$
{\rm Spec} (A)_{f_1} \cap \cdots \cap {\rm Spec} (A)_{f_n} = {\rm Spec} (A)_{f_1 \cdots f_n} \, .
$$

\item Une partie $Z$ de ${\rm Spec} (A)$ est ferm\'ee si et seulement si il existe une famille $F$ d'\'el\'ements $f \in A$ telle que
$$
p \in Z \Leftrightarrow f \in p \, , \quad \forall \, f \in F \, .
$$
En notant $I$ l'id\'eal de $A$ engendr\'e par la famille $F$, c'est \'equivalent \`a
$$
I \subseteq p
$$
et $Z$ s'identifie \`a l'ensemble ${\rm Spec} (A/I)$ des id\'eaux premiers de l'anneau quotient.

\smallskip

Le ferm\'e $Z$ est \'egalement d\'efini par l'id\'eal au moins aussi grand
$$
I_Z = \underset{p \in Z}{\bigcap} \ p \, .
$$

\item Si $u : A \to B$ est un morphisme d'anneaux, on a pour tout $f \in A$
$$
{\rm Spec} (u)^{-1} ({\rm Spec} (A)_f) = {\rm Spec} (B)_{u(f)} \, .
$$

Par cons\'equent, si $Z$ est un ferm\'e de ${\rm Spec} (A)$ d\'efini par un id\'eal $I$ de $A$, ${\rm Spec} (u)^{-1} (Z)$ est le ferm\'e de ${\rm Spec} (B)$ d\'efini par l'id\'eal $J$ engendr\'e par $u(I)$. 
\end{listeisansmarge}
\end{remarksqed}

\medskip

Par exemple, si $A = {\mathbb Z}$, ${\rm Spec} ({\mathbb Z})$ est constitu\'e de l'id\'eal $\{0\}$ et des nombres premiers. Ses ouverts non vides sont les sous-ensembles constitu\'es de $\{0\}$ et de tous les nombres premiers sauf un nombre fini.

\smallskip

Si $K$ est un corps commutatif, ${\rm Spec} (K)$ est r\'eduit \`a un seul \'el\'ement.

\smallskip

D'autre part, l'ensemble ${\rm Spec} (K[X])$ est constitu\'e de l'id\'eal $\{0\}$ et des polyn\^omes irr\'eductibles unitaires dans $K[X]$. Ses ouverts non vides sont les sous-ensembles constitu\'es de $\{0\}$ et de tous les polyn\^omes irr\'eductibles unitaires sauf un nombre fini.

\bigskip

Citons quelques propri\'et\'es importantes de la topologie des espaces de la forme ${\rm Spec} (A)$:

\begin{defn}\label{defI35}

Soit $A$ un anneau commutatif.

\begin{listeimarge}

\item Pour tout \'el\'ement $f \in A$ et $A_f = A[X] / (f \cdot X-1)$, le morphisme
$$
A \longrightarrow A_f
$$
induit une bijection canonique
$$
{\rm Spec} (A_f) \xrightarrow{ \ \sim \ } {\rm Spec} (A)_f = \{ p \in {\rm Spec} (A) \mid f \notin p \} \, .
$$

\item L'espace ${\rm Spec} (A)$ est la r\'eunion d'une famille d'ouverts ${\rm Spec} (A_{f_i})$ d\'efinis par des \'el\'ements $f_i \in A$, $i \in I$, si et seulement si les $f_i$ engendrent l'id\'eal total $A$, c'est-\`a-dire si $1$ est une combinaison lin\'eaire (finie) des $f_i$.

\medskip

\item Un ferm\'e $Z$ de ${\rm Spec} (A)$ est irr\'eductible (au sens que deux ouverts non vides de $Z$ ont une intersection non vide) si et seulement si l'id\'eal
$$
I_Z = \underset{p \in Z}{\bigcap} \ p
$$
est premier.
\end{listeimarge}
\end{defn}

\begin{remarksqed}
\begin{listeisansmarge}
\item Il r\'esulte de (ii) que tout recouvrement ouvert de ${\rm Spec} (A)$ poss\`ede un sous-recouvrement fini.

\smallskip

Autrement dit, l'espace ${\rm Spec} (A)$ est ``quasi-compact''.

\medskip

\item Il r\'esulte de (iii) que l'espace topologique ${\rm Spec} (A)$ est ``sobre'' au sens que tout ferm\'e irr\'eductible est l'adh\'erence d'un point uniquement d\'etermin\'e. 
\end{listeisansmarge}
\end{remarksqed}

\section{Expression des espaces topologiques: les faisceaux}\label{sec14}

\subsection{La notion de faisceau (et celle de pr\'efaisceau)}\label{subsec141}

\begin{defn}\label{defI41}

Soit $X$ un espace topologique.

\begin{listeimarge}

\item Un ``pr\'efaisceau'' sur $X$ est une double application $P$ qui associe

\medskip

$\left\{\begin{matrix}
\bullet &\mbox{\`a tout ouvert $U$ de $X$} \hfill \\
&\mbox{un ensemble $P(U)$} \hfill \\
&\mbox{appel\'e ``ensemble des sections de $P$ sur $U$'',} \hfill \\
{ \ } \\
\bullet &\mbox{\`a toute inclusion $U_1 \subset U_2$ entre deux ouverts de $X$} \hfill \\
&\mbox{une application $r_{U_1 , U_2} : P(U_2) \to P(U_1)$} \hfill \\
&\mbox{appel\'ee ``application de restriction des sections de $U_2$ \`a $U_1$'',} \hfill
\end{matrix} \right.$

\medskip

et telle que

\medskip

$\left\{\begin{matrix}
\bullet &\mbox{pour toute inclusion compos\'ee $U_1 \subset U_2 \subset U_3$,} \hfill \\
{ \ } \\
&r_{U_1 , U_3} = r_{U_1 , U_2} \circ r_{U_2 , U_3} : P(U_3) \longrightarrow P(U_2) \longrightarrow P(U_1) \, , \\
{ \ } \\
\bullet &\mbox{pour tout ouvert $U$,} \hfill \\
{ \ } \\
&r_{U,U} = {\rm id} : P(U) \longrightarrow P(U) \, .
\end{matrix} \right.$


\item Un morphisme $u : P_1 \to P_2$ entre deux pr\'efaisceaux $P_1 , P_2$ sur $X$ est une famille d'applications
$$
u_U : P_1 (U) \longrightarrow P_2 (U)
$$
index\'ees par les ouverts $U$ de $X$, et telle que, pour toute inclusion $U \subset V$, le carr\'e
$$
\xymatrix{
P_1(V) \ar[d]_{r_{U,V}} \ar[r]^{u_V} &P_2(V) \ar[d]^{r_{U,V}} \\
P_1 (U) \ar[r]^{u_U} &P_2 (U)
}
$$
soit commutatif (au sens que $r_{U,V} \circ u_V = u_U \circ r_{U,V}$).

\end{listeimarge}
\end{defn}
\vspace{-0.4cm}
\begin{remarksqed} 
\begin{listeisansmarge}
\item Le compos\'e de deux morphismes de pr\'efaisceaux

\vspace{-0.3cm}
$$
P_1 \xrightarrow{ \ u \ } P_2 \xrightarrow{ \ v \ } P_3
$$

\vspace{-0.3cm}
est d\'efini comme la famille des applications compos\'ees
\vspace{-0.3cm}
$$
v_U \circ u_U : P_1 (U) \longrightarrow P_2(U) \longrightarrow P_3(U) \, .
$$

\vspace{-0.4cm}
C'est un morphisme de pr\'efaisceaux.

\medskip

\item Si $P$ est un pr\'efaisceau, que chaque ensemble $P(U)$ est muni d'une structure de groupe [resp. d'anneau] et que les applications de restriction $r_{U_1 , U_2} : P(U_2) \to P(U_1)$ sont des morphismes de groupes [resp. d'anneaux], alors $P$ est appel\'e un pr\'efaisceau de groupes [resp. d'anneaux].

\medskip

\item Un morphisme entre deux pr\'efaisceaux de groupes [resp. d'anneaux] $P_1$ et $P_2$ est un morphisme de pr\'efaisceaux
\vspace{-0.3cm}
$$
u : P_1 \longrightarrow P_2
$$
\vspace{-0.3cm}
dont toutes les composantes
$$
u_U : P_1 (U) \longrightarrow P_2(U)
$$
sont des morphismes de groupes [resp. d'anneaux].

\medskip

\item Si $A$ est un pr\'efaisceau d'anneaux sur $X$, un pr\'efaisceau de modules sur $A$ est un pr\'efaisceau de groupes ab\'eliens $M$ muni d'applications
$$
A(U) \times M(U) \longrightarrow M(U)
$$
qui font de chaque $M(U)$ un module sur $A(U)$ et telles que, pour toute inclusion d'ouverts $U \subset V$, le carr\'e

\vspace{-0.4cm}
$$
\xymatrix{
A(V) \times M(V) \ar[d] \ar[r] &M(V) \ar[d] \\
A(U) \times M(U) \ar[r] &M(U)
}
$$

\vspace{-0.3cm}
soit commutatif.


Un morphisme de pr\'efaisceaux de modules sur $A$
$$
M \longrightarrow N
$$
est un morphisme de pr\'efaisceaux dont chaque composante
$$
M(U) \longrightarrow N(U)
$$
\vspace{-0.2cm}
est un morphisme de modules sur $A(U)$. 
\end{listeisansmarge}
\end{remarksqed}
\pagebreak 
 
Pour d\'efinir les faisceaux \`a partir des pr\'efaisceaux, on a besoin de la notion de recouvrement ouvert:

\begin{defn}\label{defI42}

Soit $X$ un espace topologique.

\begin{listeimarge}

\item Un recouvrement d'un ouvert $U$ de $X$ est une famille d'ouverts $(U_i)_{i \in I}$ contenus dans $U$ telle que
$$
U = \underset{i \in I}{\bigcup} \ U_i \, .
$$

\item On dit que deux sections d'un pr\'efaisceau $P$ sur un ouvert $U$ ``co{\"\i}ncident localement'' si elles ont m\^eme image dans
$$
\prod_{i \in I} P(U_i)
$$
pour un recouvrement ouvert $(U_i)_{i \in I}$ de $U$. 
\end{listeimarge}
\end{defn}

Un faisceau est un pr\'efaisceau dont les sections peuvent \^etre d\'efinies localement, par recollement le long de n'importe quel recouvrement ouvert:

\begin{defn}\label{defI43}

Soit $X$ un espace topologique.

\begin{listeimarge}

\item Un pr\'efaisceau $P$ sur $X$ est dit ``s\'epar\'e'' si toute paire de sections de $P$ sur un ouvert qui co{\"\i}ncident localement sont n\'ecessairement \'egales.

\medskip

\item Un faisceau sur $X$ est un pr\'efaisceau s\'epar\'e $F$ tel que, pour tout recouvrement $(U_i)_{i \in I}$ de tout ouvert $U$, l'image de l'application injective
$$
F(U) \xhookrightarrow{ \ \ \ \ } \prod_{i \in I} F(U_i)
$$
consiste en les familles de sections
$$
(s_i \in F(U_i))_{i \in I}
$$
qui satisfont la condition de ``recollement''
$$
r_{U_i \cap U_j , U_i} (s_i) = r_{U_i \cap U_j , U_j} (s_j) \, , \quad \forall \, i,j \in I \, .
$$

\item Un morphisme de faisceaux $F_1 \to F_2$ est un morphisme de $F_1$ dans $F_2$ consid\'er\'es comme des pr\'efaisceaux.
\end{listeimarge}
\end{defn}

\begin{remarksqed}
	\begin{listeisansmarge}
\item Un faisceau de groupes [resp. d'anneaux, resp. de modules sur un faisceau d'anneaux $A$] est un pr\'efaisceau de groupes [resp. d'anneaux, resp. de modules] qui est un faisceau.

\medskip

\item Un morphisme de faisceaux de groupes [resp. d'anneaux, resp. de modules sur un faisceau d'anneaux $A$]
$$
G_1 \longrightarrow G_2 \qquad \mbox{[resp.} \ A_1 \longrightarrow A_2 \, , \ \mbox{resp.} \ M \longrightarrow N \mbox{]}
$$
est un morphisme entre eux consid\'er\'es comme des pr\'efaisceaux de groupes [resp. d'anneaux, resp. de modules].

\medskip

\item Consid\'erant l'immersion $i : U \hookrightarrow X$ d'un ouvert $U$ de $X$, la restriction d'un pr\'efaisceau $P$ sur $X$ aux ouverts contenus dans $U$ d\'efinit un pr\'efaisceau $i^* P$ sur $U$. L'op\'eration $P \mapsto i^* P$ transforme les faisceaux sur $X$ en faisceaux sur $U$ et les morphismes de pr\'efaisceaux (ou de faisceaux) sur $X$ en morphismes de pr\'efaisceaux (ou de faisceaux) sur $U$. Elle respecte les structures de pr\'efaisceaux (ou de faisceaux) de groupes, d'anneaux ou de modules.

\medskip

\item De m\^eme, pour toute application continue $f : X \to Y$ et tout pr\'efaisceau $P$ sur $X$, associer \`a tout ouvert $V$ de $Y$ l'ensemble $P(f^{-1} V)$ et \`a toute inclusion $V_1 \subset V_2$ entre ouverts de $Y$ l'application $P (f^{-1} V_2) \to P(f^{-1} V_1)$ d\'efinit un pr\'efaisceau sur $Y$ not\'e $f_* P$.

\smallskip

C'est un faisceau sur $Y$ si $P$ est un faisceau sur $X$.

\smallskip

L'op\'eration $P \mapsto j_* P$ transforme les morphismes de pr\'efaisceaux (ou de faisceaux) en morphismes de pr\'efaisceaux (ou de faisceaux) et respecte les structures de pr\'efaisceaux (ou de faisceaux) de groupes, d'anneaux ou de modules. 
\end{listeisansmarge}
\end{remarksqed}

\subsection{Exemples de faisceaux}\label{subsec142}

\noindent $\bullet$ {\bf Le faisceau des applications d'un ensemble dans un autre:}

\smallskip

Soient $X$ et $Y$ deux ensembles.

\smallskip

Munissons $X$ de la ``topologie discr\`ete'' pour laquelle toute partie $U$ de $X$ est un ensemble ouvert.

\smallskip

Alors associer \`a toute partie $U$ de $X$ l'ensemble
$$
{\rm Hom} (U,Y)
$$
des applications de $U$ dans $Y$ d\'efinit un faisceau sur $X$.

\medskip

\noindent $\bullet$ {\bf Le faisceau des applications continues d'un espace dans un autre:}

\smallskip

Soient $X$ et $Y$ deux espaces topologiques.

\smallskip

Alors associer \`a tout ouvert $U$ de $X$ l'ensemble
$$
{\rm Hom} (U,Y)
$$
des applications continues de $U$ dans $Y$ d\'efinit un faisceau sur $X$.

\medskip

\noindent $\bullet$ {\bf Les faisceaux d'applications diff\'erentiables:}

\smallskip

Si $U$ est un ouvert de ${\mathbb R}^n$ et $k \geq 1$ un entier ou le symbole $\infty$, associer \`a tout ouvert $V$ de $U$ l'anneau
$$
C^k(V)
$$
des applications de classe $C^k$ de $V$ dans ${\mathbb R}$ d\'efinit un faisceau d'anneaux sur $U$.

\smallskip

De m\^eme, si $U$ est un ouvert de ${\mathbb C}^n$, associer \`a tout ouvert $V$ de $U$ l'anneau
$$
{\rm Hol} (V)
$$
des fonctions holomorphes de $V$ dans ${\mathbb C}$ d\'efinit un faisceau d'anneaux sur $U$.

\medskip

\noindent $\bullet$ {\bf La notion de vari\'et\'e diff\'erentielle:}

\begin{defn}\label{defI44}

Soit $k \geq 1$ un entier ou le symbole $\infty$.

\smallskip

Une vari\'et\'e diff\'erentielle de classe $C^k$ et de dimension $n$ est un espace topologique $V$ muni d'un sous-faisceau ${\mathcal O}_V$ du faisceau d'anneaux
$$
(U \subset V) \longmapsto \{\mbox{applications continues} \ U \to {\mathbb R}\}
$$
et tel qu'existe un recouvrement ouvert
$$
V = \bigcup_{i \in I} V_i
$$
et des hom\'eomorphismes $\sigma_i$ des $V_i$ vers des ouverts $U_i$ de ${\mathbb R}^n$ qui, par composition avec $\sigma_i^{-1}$, identifient la restriction de ${\mathcal O}_V$ \`a $V_i$ au sous-faisceau
$$
\begin{matrix}
&(U \subset U_i) &\longmapsto &C^k (U) = \{\mbox{applications $U \to {\mathbb R}$ de classe} \ C^k\} \hfill \\
\mbox{de} &(U \subset U_i) &\longmapsto &\{\mbox{applications continues} \ U \to {\mathbb R}\} \, . \hfill
\end{matrix}
$$
\end{defn}

\begin{remarkqed}

Pour tout ouvert $U$ d'une telle vari\'et\'e $V$, ${\mathcal O}_V (U)$ peut \^etre appel\'e l'anneau des fonctions $U \to {\mathbb R}$ de classe $C^k$. 
\end{remarkqed}

\medskip

\noindent $\bullet$ {\bf La notion de vari\'et\'e analytique:}

\begin{defn}\label{defI45}

Une vari\'et\'e analytique de dimension $n$ est une space topologique $V$ muni d'un sous-faisceau ${\mathcal O}_V$ du faisceau d'anneaux
$$
(U \subset V) \longmapsto \{\mbox{applications continues} \ U \to {\mathbb C}\}
$$
et tel qu'existe un recouvrement ouvert
$$
V = \bigcup_{i \in I} V_i
$$
et des hom\'eomorphismes $\sigma_i$ des $V_i$ vers des ouverts $U_i$ de ${\mathbb C}^n$ qui, par composition avec $\sigma_i^{-1}$, identifient la restriction de ${\mathcal O}_V$ \`a $V_i$ au sous-faisceau
$$
\begin{matrix}
&(U \subset U_i) &\longmapsto &{\rm Hol} (U) = \{\mbox{fonctions holomorphes} \ U \to {\mathbb C}\} \hfill \\
\mbox{de} &(U \subset U_i) &\longmapsto &\{\mbox{fonctions continues} \ U \to {\mathbb C}\} \, . \hfill
\end{matrix}
$$
\end{defn}

\begin{remarkqed}

Pour tout ouvert $U$ d'une telle vari\'et\'e analytique $V$, ${\mathcal O}_V (U)$ peut \^etre appel\'e l'anneau des fonctions holomorphes $U \to {\mathbb C}$.

\end{remarkqed}

\subsection{La notion de sch\'ema affine et celle de sch\'ema}\label{subsec143}

\smallskip

La construction des sch\'emas affines associ\'es aux anneaux commutatifs fournit des exemples plus sophistiqu\'es de faisceaux d'anneaux (et de modules) sur des espaces topologiques.

\begin{lem}\label{lemI46}

Soit $A$ un anneau commutatif.

\begin{listeimarge}

\item Pour tout \'el\'ement $f \in A$, un morphisme d'anneaux commutatifs
$$
u : A \longrightarrow B
$$
se factorise, de mani\`ere n\'ecessairement unique, \`a travers le morphisme canonique
$$
A \longrightarrow A_f = A[X] / (f \cdot X-1)
$$
si et seulement si $u(f)$ est inversible dans $B$.

\smallskip

En particulier, si $f' = f a$ est un multiple de $f$ dans $A$, on a un morphisme induit
$$
A_f \longrightarrow A_{fa} \, .
$$

\item Pour tout $A$-module $M$ et tout $f \in A$, les \'el\'ements de $A_f \otimes_A M = M_f$ peuvent \^etre \'ecrits $f^{-n} \cdot m$ avec $n \in {\mathbb N}$, $m \in M$. Deux \'el\'ements $f^{-n} \cdot m$ et $f^{-n'} \cdot m'$ sont \'egaux dans $M_f$ si et seulement si il existe $N \in {\mathbb N}$ tel que
$$
f^N \cdot (f^{n'} \cdot m - f^n \cdot m') = 0 \quad \mbox{dans} \quad M \, .
$$

\item Pour tous \'el\'ements $f_i$, $i \in I$, de $A$ tels que $\underset{i}{\sum} \, f_i \cdot A = A$ ou, ce qui revient au m\^eme, ${\rm Spec} (A) = \underset{i \in I}{\bigcup} {\rm Spec} (A_{f_i})$, et pour tout $A$-module $M$, les \'el\'ements de $M$ s'identifient aux familles d'\'el\'ements
$$
m_i \in M_{f_i} \, , \quad i \in I \, ,
$$
telles que, pour tous $i,j$, $m_i$ et $m_j$ aient la m\^eme image dans $M_{f_i f_j}$.
\end{listeimarge}
\end{lem}

\begin{demo}
\begin{listeisansmarge}
\item est \'evident.

\medskip

\item Le $A_f$-module $M_f = A_f \otimes_A M$ est le quotient du $A[X]$-module
$$
A[X] \otimes_A M = \bigoplus_{n \in {\mathbb N}} X^n \otimes M
$$
par le sous-module $(f \cdot X - 1) \cdot A[X] \otimes_A M$.

\smallskip

Tout \'el\'ement de $M_f$ peut \^etre repr\'esent\'e par une expression
$$
P = 1 \otimes m_0 + X \otimes m_1 + \cdots + X^n \otimes m_n
$$
avec $m_0 , m_1 , \cdots , m_n \in M$. Alors $f^n \cdot P$ est aussi repr\'esent\'e par
$$
f^n \cdot m_0 + f^{n-1} \cdot m_1 + \cdots + f \cdot m_{n-1} + m_0 \in M
$$
puisque $f^k \cdot X^k = 1$ dans $A_f$ pour tout $k \in {\mathbb N}$.

\smallskip

Si un \'el\'ement $m \in M$ est $0$ dans $M_f$, il existe une expression
$$
P = 1 \otimes m_0 + X \otimes m_1 + \cdots + X^n \otimes m_n \in A[X] \otimes_A M
$$
telle que
$$
m = (f \cdot X - 1) \cdot P \quad \mbox{dans} \quad A[X] \otimes_A M \, .
$$
Cela implique
$$
m = m_0 \, , \ f \cdot m_0 = m_1 , \cdots , f \cdot m_{n-1} = m_n \, , \ f \cdot m_n = 0
$$
et donc
$$
f^{n+1} \cdot m = 0 \, .
$$

\medskip

\item L'\'egalit\'e $\underset{i \in I}{\sum} \, f_i \cdot A = A$ est \'equivalente \`a $1 \in \underset{i\in I}{\sum} \, f_i \cdot A$. Donc on peut supposer que $I$ est fini et \'egal \`a $\{ 1,\ldots , k \}$.

\smallskip

Elle est aussi \'equivalente \`a ce que, pour tout id\'eal premier $p$ de $A$, il existe $i$ tel que $f_i \notin p$. Donc chaque $f_i$ peut \^etre remplac\'e par une puissance arbitraire de $f_i$.

\smallskip

Consid\'erons un \'el\'ement $m \in M$ dont l'image dans chaque $M_{f_i}$ est $0$. Alors il existe des entiers $n_i \geq 1$ tels que
$$
f_i^{n_i} \cdot m = 0 \quad \mbox{dans $M$ pour tout $i$.}
$$
Comme il existe des \'el\'ements $a_i \in A$ tels que
$$
a_1 \, f_1^{n_1} + \cdots + a_k \, f_k^{n_k} = 1 \, ,
$$
on conclut
$$
m = a_1 \, f_1^{n_1} \cdot m + \cdots + a_k \, f_k^{n_k} \cdot m = 0 \quad \mbox{dans} \quad M \, .
$$
Ainsi, le morphisme canonique
$$
M \longrightarrow \prod_{1 \leq i \leq k} M_{f_i}
$$
est injectif.

\smallskip

Consid\'erons enfin une famille d'\'el\'ements $f_i^{-n_i} \cdot m_i \in M_{f_i}$, $1 \leq i \leq k$, telle que, pour tous $i,j$, $f_i^{-n_i} \cdot m_i = f_j^{-n_j} \cdot m_j$ dans $M_{f_i f_j}$. On peut supposer que tous les entiers $n_i$ sont \'egaux \`a un m\^eme entier $n$.

\smallskip

Alors il existe un entier $N \in {\mathbb N}$ tel que, pour tous $i,j$,
$$
(f_i \, f_j)^N \, f_j^n \cdot m_i = (f_i \, f_j)^N \, f_i^n \cdot m_j \quad \mbox{dans} \quad M \, .
$$
Posant $m'_i = f_i^N \cdot m_i$ et $f'_i = f_i^{N+n}$ pour tout $i$, on a
$$
f_i^{-n_i} \cdot m_i = f'^{-1}_i \cdot m'_i \quad \mbox{dans chaque} \quad M_{f_i} = M_{f'_i} \, ,
$$
et ces \'el\'ements sont reli\'es par les \'egalit\'es
$$
f'_j \cdot m'_i = f'_i \cdot m'_j \quad \mbox{dans $M$ pour tous $i,j$.}
$$
Choisissant des \'el\'ements $a_i \in A$ tels que $a_1 \, f'_1 + \cdots + a_k \, f'_k = 1$, on d\'efinit l'\'el\'ement
$$
m = a_1 \cdot m'_1 + \cdots + a_k \cdot m'_k \quad \mbox{dans} \quad M \, .
$$
Pour tout indice $i$, on a dans $M$ l'\'egalit\'e
$$
f'_i \cdot m= \sum_j a_j \, f'_i \cdot m'_j = \sum_j a_j \, f'_j \cdot m'_i = m'_i \, ,
$$
ce qui signifie que $f'^{-1}_i \cdot m'_i = m$ dans chaque $M_{f'_i} = M_{f_i}$.
\end{listeisansmarge}
\end{demo}

\medskip

Ce lemme permet de poser la d\'efinition suivante:

\begin{defn}\label{defI47}

Le ``sch\'ema affine'' associ\'e \`a un anneau commutatif $A$ est l'espace topologique ${\rm Spec} (A)$ muni de l'unique faisceau d'anneaux commutatifs ${\mathcal O}_A$ tel que
$$
{\mathcal O}_A ({\rm Spec} (A)_f) = A_f
$$
pour tout \'el\'ement $f \in A$.

\smallskip

On l'appelle le ``spectre'' de $A$.
\end{defn}
\pagebreak
\begin{remarksqed}
\begin{listeisansmarge}
\item L'unicit\'e du faisceau ${\mathcal O}_A$ r\'esulte de ce que tout ouvert de ${\rm Spec} (A)$ est r\'eunion d'ouverts de la forme ${\rm Spec} (A)_f$.

\smallskip

Son existence r\'esulte du lemme \ref{lemI46}~(iii) en d\'ecidant que, pour tout ouvert $U$ de ${\rm Spec} (A)$, ${\mathcal O}_A (U)$ est l'anneau commutatif des familles d'\'el\'ements 

\smallskip

$a_f \in A_f$ index\'es par les $f \in A$ tels que ${\rm Spec} (A)_f \subset U$, 

\smallskip

\noindent telles que pour tout multiple $f' = fa$ d'un tel $f$, 

\smallskip

$a_f$ s'envoie sur $a_{f'}$ par $A_f \to A_{f'}$.

\medskip

\item Consid\'erons l'application continue
$$
p = {\rm Spec} (u) : {\rm Spec}(B) \longrightarrow {\rm Spec}(A)
$$
induite par un morphisme d'anneaux commutatifs $u : A \to B$.

\smallskip

Pour tout \'el\'ement $f \in A$, on a un morphisme induit
$$
\xymatrix{
A_f \ar@{=}[d] \ar[r] &B_{u(f)} \ar@{=}[d] \\
{\mathcal O}_A ({\rm Spec} (A)_f) &{\mathcal O}_B ({\rm Spec} (B)_{u(f)})
}
$$
avec ${\rm Spec} (B)_{u(f)} = p^{-1} ({\rm Spec} (A)_f)$.

\smallskip

Ainsi, $u$ induit un morphisme de faisceaux d'anneaux commutatifs sur ${\rm Spec} (A)$
$$
{\mathcal O}_A \longrightarrow p_* {\mathcal O}_B \, .
$$
La composante de ce morphisme pour les sections de ces faisceaux au-dessus de ${\rm Spec} (A)$ tout entier n'est autre que:
$$
\xymatrix{
A \ar@{=}[d] \ar[r]^u &B \ar@{=}[d] \\
{\mathcal O}_A ({\rm Spec} (A)) &{\mathcal O}_B ({\rm Spec} (B)) = p_* {\mathcal O}_B ({\rm Spec} (A))
}
$$

\end{listeisansmarge}
\end{remarksqed}

\medskip

De la m\^eme fa\c con, le lemme permet de poser:

\begin{defn}\label{defI48}

Soient $A$ un anneau commutatif et ${\rm Spec} (A)$ son spectre muni de son faisceau de structure ${\mathcal O}_A$.

\smallskip

Alors tout $A$-module $M$ d\'efinit un unique faisceau de modules $\widetilde M$ sur ${\mathcal O}_A$ tel que
$$
\widetilde M ({\rm Spec} (A)_f) = M_f = A_f \otimes_A M
$$
pour tout \'el\'ement $f$ de $A$.
\end{defn}

\begin{remarksqed} 
\begin{listeisansmarge}
\item Un faisceau de modules sur ${\mathcal O}_A$ associ\'e de cette fa\c con \`a un module sur $A$ (qui est n\'ecessairement le module de ses sections sur ${\rm Spec} (A)$ tout entier) est appel\'e un faisceau ``quasi-coh\'erent''. Il est dit ``coh\'erent'' s'il est associ\'e \`a un $A$-module $M$ qui est de pr\'esentation finie sur $A$ (ou, ce qui revient au m\^eme, de type fini quand $A$ est un anneau n\oe th\'erien).

\medskip

\item Si $M$ est un $A$-module, $\widetilde M$ le faisceau quasi-coh\'erent associ\'e et $U$ un ouvert de ${\rm Spec} (A)$, $\widetilde M (U)$ est le module des familles d'\'el\'ements

\smallskip

$m_f \in M_f$ index\'es par les $f \in A$ tels que ${\rm Spec} (A)_f \subset U$,

\smallskip

\noindent telles que pour tout multiple $f' = fa$ d'un tel $f$,

\smallskip

$m_f$ s'envoie sur $m_{f'}$ par le morphisme
$$
M_f = A_f \otimes_A M \longrightarrow A_{f'} \otimes_A M = M_{f'} \, .
$$

\item Un morphisme de $A$-modules $M \to N$ induit pour tout $f \in A$ un morphisme de $A_f$-modules
$$
\xymatrix{
M_f \ar@{=}[d] \ar[r] &N_f \ar@{=}[d] \\
\widetilde M ({\rm Spec} (A)_f) &\widetilde N ({\rm Spec} (A)_f)
}
$$
et donc un morphisme de faisceaux de modules sur ${\mathcal O}_A$
$$
\widetilde M \longrightarrow \widetilde N \, .
$$
La composante de ce morphisme pour les sections de ces faisceaux au-dessus de ${\rm Spec} (A)$ tout entier n'est autre que le morphisme de $A$-modules de d\'epart
$$
M \longrightarrow N \, .
$$

\end{listeisansmarge}
\end{remarksqed}

\medskip
\vspace{-0.5cm}
On peut enfin poser:

\begin{defn} \label{defI49}
\begin{listeimarge}
\item Un ``espace annel\'e'' est un espace topologique $X$ muni d'un faisceau d'anneaux ${\mathcal O}_X$.

\medskip

\item Un sch\'ema est un espace annel\'e $(X,{\mathcal O}_X)$ qui admet un recouvrement ouvert
$$
X = \bigcup_{i \in I} U_i
$$
par des ouverts $U_i$ qui, munis des restrictions ${\mathcal O}_{U_i}$ aux $U_i$ du faisceau de structure ${\mathcal O}_X$, sont isomorphes aux sch\'emas affines $({\rm Spec} (A_i) , {\mathcal O}_{A_i})$ associ\'es \`a des anneaux commutatifs $A_i$.
\end{listeimarge}
\end{defn}

\begin{remarksqed} 
\begin{listeisansmarge}
\item Un isomorphisme entre deux espaces annel\'es $(X,{\mathcal O}_X)$ et $(Y,{\mathcal O}_Y)$ consiste en un hom\'eomorphisme entre les deux espaces topologiques sous-jacents $\sigma : X \xrightarrow{ \ \sim \ }Y$ compl\'et\'e par un isomorphisme ${\mathcal O}_Y \xrightarrow{ \ \sim \ } \sigma_* {\mathcal O}_X$ entre faisceaux d'anneaux sur $Y$.

\medskip

\item La d\'efinition \ref{defI44} [resp. \ref{defI45}] a d\'efini une vari\'et\'e diff\'erentielle de classe $C^k$ [resp. une vari\'et\'e analytique] de dimension $n$ comme un espace annel\'e qui est localement isomorphe \`a des ouverts de ${\mathbb R}^n$ [resp. ${\mathbb C}^n$] munis du faisceau des fonctions de classe $C^k$ [resp. holomorphes] \`a valeurs dans ${\mathbb R}$ [resp. ${\mathbb C}$]. 
\end{listeisansmarge}
\end{remarksqed}

\section{Cat\'egories et foncteurs}\label{sec15}

\subsection{La notion de cat\'egorie}\label{subsec151}

Chaque fois que nous avons d\'efini la notion de groupe, de mono{\"\i}de, d'action d'un groupe ou d'un mono{\"\i}de, d'espace topologique, de pr\'efaisceau ou de faisceau sur un espace topologique, nous avons \'egalement d\'efini la notion de morphisme entre groupes, entre mono{\"\i}des, entre ensembles munis d'une action d'un groupe ou d'un mono{\"\i}de, entre espaces topologiques, entre pr\'efaisceaux ou entre faisceaux sur un tel espace. Et chaque fois nous avons observ\'e que ces morphismes se composaient.

\smallskip

Cela conduit \`a introduire la d\'efinition g\'en\'erale suivante:

\begin{defn}\label{defI51}

Une cat\'egorie ${\mathcal C}$ consiste en

\medskip

$\left\{ \begin{matrix}
\bullet &\mbox{une collection ${\rm Ob} ({\mathcal C})$ dont les \'el\'ements sont appel\'es les objets de ${\mathcal C}$,} \hfill \\
{ \ } \\
\bullet &\mbox{pour tous objets $X,Y$ de ${\mathcal C}$, une collection ${\rm Hom}_{\mathcal C} (X,Y) = {\rm Hom} (X,Y)$ dont les \'el\'ements} \hfill \\
&\mbox{sont appel\'es les morphismes, ou les fl\`eches, de $X$ dans $Y$ et sont not\'es $u : X \to Y$,} \hfill \\
{ \ } \\
\bullet &\mbox{pour tous objets $X,Y,Z$, une loi de composition} \hfill \\
{ \ } \\
&\begin{matrix}
{\rm Hom} (Y,Z) \times {\rm Hom} (X,Y) &\longrightarrow &{\rm Hom} (X,Z) \, , \hfill \\
(v: Y \to Z , u : X \to Y)) &\longmapsto &(v \circ u : X \to Z),
\end{matrix}
\end{matrix} \right.
$

\medskip

satisfaisant les conditions suivantes:

\medskip

$\left\{\begin{matrix}
\bullet &\mbox{la loi de composition est associative au sens que,} \hfill \\
&\mbox{pour tous objets $X_1 , X_2 , X_3 , X_4$ et tous morphismes} \hfill \\
{ \ } \\
&u : X_1 \to X_2 \, , \ v : X_2 \to X_3 \, , \ w : X_3 \to X_4 \, , \\
&\mbox{on a} \hfill \\
&w \circ (v \circ u) = (w \circ v) \circ u : X_1 \to X_4 \, , \\
{ \ } \\
\bullet &\mbox{pour tout objet $X$, il existe un morphisme (n\'ecessairement unique)} \hfill \\
&{\rm id}_X : X \to X \ \mbox{tel que, pour tout objet $Y$,} \hfill \\
{ \ } \\
&{\rm id}_X \circ u = u \, , \quad \forall \, u \in {\rm Hom} (Y,X) \\
&\mbox{et} \hfill \\
&u \circ {\rm id}_X = u \, , \quad \forall \, u \in {\rm Hom} (X,Y) \, .
\end{matrix} \right.$
\end{defn}

\begin{remarksqed}
\begin{listeisansmarge}
\item On peut aussi consid\'erer la collection ${\rm Hom} ({\mathcal C})$ de toutes les fl\`eches d'une cat\'egorie ${\mathcal C}$. Chaque fl\`eche $u : X \to Y$ a une ``origine'' $X$ et un ``but'' $Y$.

\medskip

\item Si, pour tous objets $X,Y$ de ${\mathcal C}$, la collection ${\rm Hom} (X,Y)$ est un ensemble, la cat\'egorie ${\mathcal C}$ est dite ``localement petite''.

\medskip

\item Si ${\mathcal C}$ est une cat\'egorie localement petite et que, de plus, la collection ${\rm Ob} ({\mathcal C})$ est un ensemble, la cat\'egorie ${\mathcal C}$ est dite petite.

\medskip

\item Pour tout objet $X$ d'une cat\'egorie ${\mathcal C}$, les \'el\'ements de ${\rm Hom} (X,X) = {\rm End} (X)$ sont appel\'es les endomorphismes de ${\mathcal C}$.

\smallskip

Si ${\mathcal C}$ est localement petite, chaque ${\rm End} (X)$ est un mono{\"\i}de dont l'\'el\'ement unit\'e est ${\rm id}_X$, l'endomorphisme ``identit\'e'' de $X$.

\medskip

\item On appelle ``triangle'' [resp. ``carr\'e''] un diagramme de fl\`eches de ${\mathcal C}$ de la forme
$$
\xymatrix{
X \ar[rr]^f \ar[rd]_p &&Y \ar[ld]^q \\
&S
} \qquad \xymatrix{
X \ar[r]^f \ar[d]_{\mbox{[resp.} \qquad p} &Y \ar[d]^{q \quad \mbox{].}} \\
S' \ar[r]^s &S
}
$$
Un tel diagramme est dit ``commutatif'' si
$$
q \circ f = p \qquad \mbox{[resp.} \quad q \circ f = s \circ p \ \mbox{].}
$$

\end{listeisansmarge}
\end{remarksqed}

\begin{defn}\label{defI52}
\begin{listeimarge}

\item Une sous-cat\'egorie ${\mathcal C}'$ d'une cat\'egorie ${\mathcal C}$ consiste en une sous-collection ${\rm Ob} ({\mathcal C}')$ de ${\rm Ob} ({\mathcal C})$ et, pour tous $X,Y \in {\rm Ob} ({\mathcal C}')$, une sous-collection ${\rm Hom}_{{\mathcal C}'} (X,Y)$ de ${\rm Hom}_{\mathcal C} (X,Y)$ v\'erifiant

\medskip

$\left\{\begin{matrix}
\bullet &\mbox{pour tout $X \in {\rm Ob} ({\mathcal C}')$, l'\'el\'ement identit\'e ${\rm id}_X$ de ${\rm Hom}_{\mathcal C} (X,X)$ est dans ${\rm Hom}_{{\mathcal C}'} (X,X)$,} \hfill \\
{ \ } \\
\bullet &\mbox{pour tous $X,Y,Z \in {\rm Ob} ({\mathcal C}')$, la loi de multiplication} \hfill \\
{ \ } \\
&{\rm Hom}_{\mathcal C} (Y,Z) \times {\rm Hom}_{\mathcal C} (X,Y) \longrightarrow {\rm Hom}_{\mathcal C} (X,Z) \\
{ \ } \\
&\mbox{envoie ${\rm Hom}_{{\mathcal C}'} (Y,Z) \times {\rm Hom}_{{\mathcal C}'} (X,Y)$ dans ${\rm Hom}_{{\mathcal C}'} (X,Z)$.} \hfill
\end{matrix} \right.
$

\medskip

\item Une telle sous-cat\'egorie est dite ``pleine'' si, pour tous $X,Y \in {\rm Ob} ({\mathcal C}')$, on a ${\rm Hom}_{{\mathcal C}'} (X,Y) = {\rm Hom}_{\mathcal C} (X,Y)$.
\end{listeimarge}
\end{defn}

\begin{remarkqed}

Se donner une sous-cat\'egorie pleine de ${\mathcal C}$ \'equivaut \`a se donner une sous-collection de ${\rm Ob} ({\mathcal C})$.
 
\end{remarkqed}

\subsection{Premiers exemples de cat\'egories}\label{subsec152}

\medskip

\noindent $\bullet$ {\bf La cat\'egorie des ensembles:}

\smallskip

Les ensembles et les applications entre ensembles munies de la loi de composition des applications constituent une cat\'egorie localement petite not\'ee ${\rm Ens}$.

\medskip

\noindent $\bullet$ {\bf La cat\'egorie des mono{\"\i}des et celle des groupes:}

\smallskip

Les mono{\"\i}des et les morphismes de mono{\"\i}des forment une cat\'egorie localement petite.

\smallskip

Les groupes forment une sous-cat\'egorie pleine de celle des mono{\"\i}des.

\medskip

\noindent $\bullet$ {\bf La cat\'egorie des anneaux et celle des corps:}

\smallskip

Les anneaux et les morphismes d'anneaux forment une cat\'egorie localement petite.

\smallskip

Les corps forment une sous-cat\'egorie pleine de celle des anneaux.

\medskip

\noindent $\bullet$ {\bf Les cat\'egories de modules:}

\smallskip

Pour tout anneau $A$, les modules sur $A$ et les applications $A$-lin\'eaires entre modules forment une cat\'egorie localement petite que l'on peut noter ${\rm Mod}_A$.

\smallskip

En particulier, pour tout corps $K$, les espaces vectoriels sur $K$ forment une cat\'egorie ${\rm Mod}_K = {\rm Vec}_K$.

\newpage

\noindent $\bullet$ {\bf Les cat\'egories d'actions d'un groupe ou d'un mono{\"\i}de:}

\smallskip

Pour tout groupe (ou mono{\"\i}de) $G$, les ensembles munis d'une action de $G$ et les applications $G$-equivariantes entre tels ensembles forment une cat\'egorie localement petite not\'ee $BG$ et appel\'ee le ``classifiant'' de $G$.

\medskip

\noindent $\bullet$ {\bf Les groupes ou mono{\"\i}des comme cas particulier de cat\'egories:}

\smallskip

Tout mono{\"\i}de (en particulier, tout groupe) $M$ peut \^etre vu comme une petite cat\'egorie avec un unique objet $X$ et ${\rm End} (X) = M$.

\medskip

\noindent $\bullet$ {\bf Les ensembles ordonn\'es comme cas particuliers de cat\'egories:}

\smallskip

Tout ensemble $O$ muni d'un ordre partiel $\leq$ peut \^etre vu comme une petite cat\'egorie dont les objets sont les \'el\'ements $x \in O$ et dont les ensembles de fl\`eches ${\rm Hom} (x_1 , x_2)$ ont un unique \'el\'ement si $x_1 \leq x_2$ et sont vides sinon.

\smallskip

En particulier, pour tout espace topologique $X$, l'ensemble $O(X)$ de ses ouverts ordonn\'e par la relation d'inclusion $\subset$ peut \^etre vu comme une petite cat\'egorie.

\medskip

\noindent $\bullet$ {\bf La cat\'egorie des espaces topologiques:}

\smallskip

Les espaces topologiques et les applications continues entre de tels espaces forment une cat\'egorie localement petite que l'on peut noter ${\rm Top}$.

\medskip

\noindent $\bullet$ {\bf Les cat\'egories de pr\'efaisceaux et de faisceaux sur un espace:}

\smallskip

Pour tout espace topologique $X$, les pr\'efaisceaux sur $X$ et les morphismes de pr\'efaisceaux forment une cat\'egorie localement petite que l'on  notera $\widehat{O(X)}$.

\smallskip

Les faisceaux sur $X$ forment une sous-cat\'egorie pleine de $\widehat{O(X)}$ que l'on notera ${\mathcal E}_X$.

\medskip

\noindent $\bullet$ {\bf Les cat\'egories de lacets:}

\smallskip

Tout espace topologique $X$ d\'efinit une petite cat\'egorie $\pi_X$ dont les objets sont les \'el\'ements $x$ de $X$ et dont les morphismes $x \to y$ sont les classes d'homotopie de ``lacets'' c'est-\`a-dire d'applications continues
$$
u : [0,1] \longrightarrow X
$$
dont l'origine est $u(0) = x$ et dont le but est $u(1) = y$.

\subsection{La notion de foncteur}\label{subsec153}

\begin{defn}\label{defI53}
\begin{listeimarge}
\item Un foncteur d'une cat\'egorie ${\mathcal C}$ dans une cat\'egorie ${\mathcal D}$
$$
F : {\mathcal C} \longrightarrow {\mathcal D}
$$
consiste en une application
\begin{eqnarray}
{\rm Ob} ({\mathcal C}) &\longrightarrow &{\rm Ob} ({\mathcal D}) \, , \nonumber \\
X &\longmapsto &F(X) \nonumber
\end{eqnarray}
et une collection d'applications index\'ees par les objets $X,Y$ de ${\mathcal C}$
\begin{eqnarray}
{\rm Hom}_{\mathcal C} (X,Y) &\longrightarrow &{\rm Hom} _{\mathcal D} (F(X) , F(Y)) \nonumber \\
(X \xrightarrow{ \ u \ } Y ) &\longmapsto &(F(X) \xrightarrow{ \, F(u) \ } F(Y)), \nonumber
\end{eqnarray}
telles que:
\begin{enumerate}
\item[$\bullet$] pour tout objet $X$ de ${\mathcal C}$
$$
F({\rm id}_X) = {\rm id}_{F(X)} \, ,
$$
\item[$\bullet$] pour tous morphismes $u : X \to Y$ et $v : Y \to Z$ de ${\mathcal C}$,
$$
F(v \circ u) = F(v) \circ F(u) \, .
$$
\end{enumerate}

\item Etant donn\'es deux cat\'egories ${\mathcal C} , {\mathcal D}$ et deux foncteurs
$$
F,G : {\mathcal C} \rightrightarrows {\mathcal D}
$$
de ${\mathcal C}$ dans ${\mathcal D}$, un morphisme (ou transformation naturelle) du foncteur $F$ dans le foncteur $G$ est une collection de morphismes de ${\mathcal D}$
$$
u_X : F(X) \longrightarrow G(X)
$$
index\'es par les objets $X$ de ${\mathcal C}$, telle que, pour tout morphisme $u : X \to Y$ de ${\mathcal C}$, le carr\'e de ${\mathcal D}$
$$
\xymatrix{
F(X) \ar[d]_{u_X} \ar[r]^{F(u)} &F(Y) \ar[d]^{u_Y} \\
G(X) \ar[r]^{G(u)} &G(Y)
}
$$
soit commutatif.

\end{listeimarge}
\end{defn}

\begin{remarksqed}
\begin{listeisansmarge}
\item Toute cat\'egorie ${\mathcal C}$ poss\`ede un ``foncteur identit\'e'' ${\rm id}_{\mathcal C}$ d\'efini par
$$
\begin{matrix}
{\rm id}_{\mathcal C} (X) = X \, , &\forall \, X \in {\rm Ob} ({\mathcal C}) \, , \hfill \\
{\rm id}_{\mathcal C} (u) = u \, , \hfill &\forall \, u \in {\rm Hom} (X,Y) \, , \ \forall \, X,Y \in {\rm Ob} ({\mathcal C}) \, .
\end{matrix}
$$

\item Deux foncteurs $F : {\mathcal C} \to {\mathcal D}$ et $G : {\mathcal D} \to {\mathcal E}$ ont un ``foncteur compos\'e'' d\'efini comme
$$
\begin{matrix}
G \circ F : &\hfill X &\longmapsto &G(F(X)) \, , \hfill \\
&(X \xrightarrow{ \ u \ } Y) &\longmapsto &(G(F(X)) \xrightarrow{ \ G(F(u)) \ } G(F(Y))) \, .
\end{matrix}
$$
La loi de composition des foncteurs est associative au sens que, pour tous foncteurs ${\mathcal C}_1 \xrightarrow{ \ F \ } {\mathcal C}_2 \xrightarrow{ \ G \ } {\mathcal C}_3 \xrightarrow{ \ H \ } {\mathcal C}_4$, on a
$$
H \circ (G \circ F) = (H \circ G) \circ F \, .
$$

\item Par cons\'equent, les petites cat\'egories et les foncteurs entre petites cat\'egories forment une cat\'egorie localement petite not\'ee ${\rm Cat}$.

\medskip

\item Tout foncteur $F : {\mathcal C} \to {\mathcal D}$ entre deux cat\'egories ${\mathcal C}$ et ${\mathcal D}$ a un ``endomorphisme identit\'e'' ${\rm id}_F$ d\'efini comme
$$
{\rm Ob} ({\mathcal C}) \ni X \longmapsto {\rm id}_{F(X)} \, .
$$

\item Etant donn\'es trois foncteurs $F,G,H$ d'une cat\'egorie ${\mathcal C}$ dans une cat\'egorie ${\mathcal D}$, deux morphismes de foncteurs
$$
u : F \longrightarrow G \qquad \mbox{et} \qquad v : G \longrightarrow H
$$
ont un morphisme compos\'e $v \circ u : F \to H$ d\'efini comme
$$
{\rm Ob} ({\mathcal C}) \ni X \longmapsto (v \circ u)_X = (v_X \circ u_X : F(X) \to G(X) \to H(X)) \, .
$$
La loi de composition des morphismes de foncteurs est associative au sens que, pour tous morphismes
$$
F_1 \xrightarrow{ \ u \ } F_2 \xrightarrow{ \ v \ } F_3 \xrightarrow{ \ w \ } F_4
$$
entre foncteurs $F_1 , F_2 , F_3 , F_4$ de ${\mathcal C}$ dans ${\mathcal D}$, on a
$$
w \circ (v \circ u) = (w \circ v) \circ u \, .
$$

\item Par cons\'equent, pour toutes cat\'egories ${\mathcal C}$ et ${\mathcal D}$, les foncteurs ${\mathcal C} \to {\mathcal D}$ et leurs morphismes forment une cat\'egorie not\'ee $[{\mathcal C} , {\mathcal D}]$.

\smallskip

Cette cat\'egorie est localement petite [resp. petite] si ${\mathcal C}$ est une petite cat\'egorie et ${\mathcal D}$ est localement petite [resp. petite].

\medskip

\item Si $M$ est un mono{\"\i}de (en particulier un groupe) consid\'er\'e comme une cat\'egorie \`a un seul objet, la cat\'egorie $BM$ des actions de $M$ n'est autre que la cat\'egorie $[M,{\rm Ens}]$ des foncteurs $M \to {\rm Ens}$.

\medskip

\item Pour toutes cat\'egories ${\mathcal C} , {\mathcal D}_1 , {\mathcal D}_2$, tout foncteur $G : {\mathcal D}_1 \to {\mathcal D}_2$ d\'efinit un foncteur
$$
[{\mathcal C} , {\mathcal D}_1] \longrightarrow [{\mathcal C} , {\mathcal D}_2]
$$
qui associe \`a tout foncteur $F : {\mathcal C} \to {\mathcal D}_1$ le compos\'e
$$
G \circ F : {\mathcal C} \longrightarrow {\mathcal D}_2
$$
et \`a toute transformation naturelle de foncteurs $F_1 , F_2 : {\mathcal C} \to {\mathcal D}_1$
$$
(u : F_1 \longrightarrow F_2) = (F_1 (X) \xrightarrow{ \ u_X \ } F_2 (X))_{X \in {\rm Ob} ({\mathcal C})}
$$
la transformation naturelle $G(u) : G \circ F_1 \to G \circ F_2$ d\'efinie par
$$
G(u)_X = (G(F_1 (X)) \xrightarrow{ \ G(u_X) \ } G(F_2(X))) \, , \quad \forall \, X \in {\rm Ob} ({\mathcal C}) \, .
$$

\item De m\^eme, pour toutes cat\'egories ${\mathcal C}_1 , {\mathcal C}_2 , {\mathcal D}$, tout foncteur $F : {\mathcal C}_1 \to {\mathcal C}_2$ d\'efinit par composition un foncteur
$$
\begin{matrix}
[{\mathcal C}_2 , {\mathcal D}] &\longrightarrow &[{\mathcal C}_1 , {\mathcal D}] \, , \\
\hfill G &\longmapsto &G \circ F \, . \hfill
\end{matrix}
$$

\end{listeisansmarge}
\end{remarksqed}

La d\'efinition des foncteurs est compl\'et\'ee par celle des ``foncteurs contravariants'' qui passe par la remarque que toute cat\'egorie a une ``cat\'egorie oppos\'ee'':

\begin{defn}\label{defI54}
\begin{listeimarge}
\item La ``cat\'egorie oppos\'ee'' ${\mathcal C}^{\rm op}$ d'une cat\'egorie ${\mathcal C}$ est la cat\'egorie d\'efinie par
\begin{enumerate}
\item[$\bullet$] ${\rm Ob} ({\mathcal C}^{\rm op}) = {\rm Ob} ({\mathcal C})$,
\item[$\bullet$] pour tous objets $X,Y$ de ${\mathcal C}$ ou ${\mathcal C}^{\rm op}$, les fl\`eches $X \to Y$ de ${\mathcal C}^{\rm op}$ sont les fl\`eches $Y \to X$ de ${\mathcal C}$, soit
$$
{\rm Hom}_{{\mathcal C}^{\rm op}} (X,Y) = {\rm Hom}_{\mathcal C} (Y,X) \, ,
$$
\item[$\bullet$] pour tous objets $X,Y,Z$, la loi de composition de ${\mathcal C}^{\rm op}$
\begin{eqnarray}
{\rm Hom}_{{\mathcal C}^{\rm op}} (Y,Z) \times {\rm Hom}_{{\mathcal C}^{\rm op}} (X,Y) &\longrightarrow &{\rm Hom}_{{\mathcal C}^{\rm op}} (X,Z) \, , \nonumber \\
(Y \xrightarrow{ \ v \ } Z , X \xrightarrow{ \ u \ } Y) &\longmapsto &(X \xrightarrow{ \ v \circ u \ } Z) \nonumber
\end{eqnarray}
\end{enumerate}
est la loi de composition de ${\mathcal C}$
\begin{eqnarray}
{\rm Hom}_{\mathcal C} (Z,Y) \times {\rm Hom}_{\mathcal C} (Y,X) &\longrightarrow &{\rm Hom}_{\mathcal C} (Z,X) \, , \nonumber \\
(Z \xrightarrow{ \ v \ } Y , Y \xrightarrow{ \ u \ } X) &\longmapsto &(Z \xrightarrow{ \ u \circ v \ } X) \, . \nonumber
\end{eqnarray}

\item Un ``foncteur contravariant'' d'une cat\'egorie ${\mathcal C}$ dans une cat\'egorie ${\mathcal D}$ est un foncteur
$$
{\mathcal C}^{\rm op} \longrightarrow {\mathcal D} \, .
$$
\end{listeimarge}
\end{defn}

\begin{remarksqed}
\begin{listeisansmarge}
\item Tout objet $X$ d'une cat\'egorie ${\mathcal C}$ a le m\^eme endomorphisme identit\'e ${\rm id}_X$ dans ${\mathcal C}$ et dans ${\mathcal C}^{\rm op}$.

\medskip

\item Pour toute cat\'egorie ${\mathcal C}$, $({\mathcal C}^{\rm op})^{\rm op} = {\mathcal C}$.

\medskip

\item Une cat\'egorie ${\mathcal C}$ est localement petite [resp. petite] si et seulement si ${\mathcal C}$ l'est.

\medskip

\item Si ${\mathcal C}$ est un mono{\"\i}de [resp. un groupe], son oppos\'e ${\mathcal C}^{\rm op}$ est un mono{\"\i}de [resp. un groupe].

\medskip

\item Un foncteur ${\mathcal C} \to {\mathcal D}$ est la m\^eme chose qu'un foncteur ${\mathcal C}^{\rm op} \to {\mathcal D}^{\rm op}$.

\smallskip

En particulier, un foncteur contravariant ${\mathcal C}^{\rm op} \to {\mathcal D}$ peut aussi \^etre vu comme un foncteur ${\mathcal C} \to {\mathcal D}^{\rm op}$.

\medskip

\item Un foncteur contravariant $F : {\mathcal C}^{\rm op} \to {\mathcal D}$ consiste \`a associer \`a tout objet $X$ de ${\mathcal C}$ un objet $F(X)$ de ${\mathcal D}$ et \`a toute fl\`eche $u : X \to Y$ de ${\mathcal C}$ une fl\`eche en sens inverse
$$
F(u) : F(Y) \longrightarrow F(X) \quad \mbox{de} \quad {\mathcal C} \, ,
$$
de telle fa\c con que $F({\rm id}_X) = {\rm id}_{F(X)}$, $\forall \, X \in {\rm Ob} ({\mathcal C})$, et que, pour toutes fl\`eches $X \xrightarrow{ \ u \ } Y \xrightarrow{ \ v \ } Z$ de ${\mathcal C}$,
$$
F(v \circ u) = F(u) \circ F(v) : F(Z) \longrightarrow F(Y) \longrightarrow F(X) \, .
$$
\end{listeisansmarge}
\end{remarksqed}

On pose:

\begin{defn}\label{defI55}

Pour toute cat\'egorie localement petite ${\mathcal C}$, on note
$$
\widehat{\mathcal C} = [{\mathcal C}^{\rm op} , {\rm Ens}]
$$
la cat\'egorie des foncteurs contravariants ${\mathcal C}^{\rm op} \to {\rm Ens}$, et $y$ le ``foncteur de Yoneda''
\begin{eqnarray}
{\mathcal C} &\longrightarrow &\widehat{\mathcal C} \nonumber \\
X &\longmapsto &y(X) = {\rm Hom} (\bullet , X) = \left\{ \begin{matrix}
{\mathcal C}^{\rm op} &\longrightarrow &{\rm Ens} \, , \hfill \\
\hfill Y &\longmapsto &{\rm Hom} (Y,X) \, .
\end{matrix} \right. \nonumber
\end{eqnarray}
\end{defn}

\begin{remarksqed} 
\begin{listeisansmarge}
\item La cat\'egorie $\widehat{\mathcal C}$ est appel\'ee la cat\'egorie des pr\'efaisceaux sur ${\mathcal C}$.

\smallskip

En effet, si $X$ est un espace topologique et $O(X)$ l'ensemble ordonn\'e de ses ouverts consid\'er\'e comme une cat\'egorie, la cat\'egorie des pr\'efaisceaux sur $X$ n'est autre que
$$
\widehat{O(X)} = [O(X)^{\rm op} , {\rm Ens}] \, .
$$

\item Pour tout objet $X$ de ${\mathcal C}$, $y(X) = {\rm Hom} (\bullet , X)$ transforme toute fl\`eche $v : Y_1 \to Y_2$ en l'application de composition
\begin{eqnarray}
{\rm Hom} (Y_2 , X) &\longrightarrow &{\rm Hom} (Y_1 , X) \, , \nonumber \\
(f : Y_2 \to X) &\longmapsto &(f \circ v : Y_1 \to Y_2 \to X) \, . \nonumber
\end{eqnarray}
C'est un foncteur parce que la loi de composition des fl\`eches de ${\mathcal C}$ est associative.

\medskip

\item Le foncteur de Yoneda $y$ transforme toute fl\`eche $u : X_1 \to X_2$ en la transformation naturelle
$$
{\rm Hom} (\bullet , X_1) \longrightarrow {\rm Hom} (\bullet , X_2)
$$
qui consiste en la famille des applications
\begin{eqnarray}
{\rm Hom} (Y,X_1) &\longrightarrow &{\rm Hom} (Y,X_2) \, , \nonumber \\
(f : Y \to X_1) &\longmapsto &(u \circ f : Y \to X_1 \to X_2) \, . \nonumber
\end{eqnarray}
Pour toute fl\`eche $v : Y_1 \to Y_2$, le carr\'e
$$
\xymatrix{
{\rm Hom} (Y_2 , X_1) \ar[d]_{\bullet \circ v} \ar[rr]^{u \circ \bullet} &&{\rm Hom} (Y_2 , X_2) \ar[d]^{\bullet \circ v} \\
{\rm Hom} (Y_1 , X_1) \ar[rr]^{u \circ \bullet} &&{\rm Hom} (Y_1 , X_2)
}
$$
est commutatif car, par associativit\'e,
$$
(u \circ f) \circ v = u \circ (f \circ v)
$$
pour toute fl\`eche $f \in {\rm Hom} (Y_2 , X_1)$.

\end{listeisansmarge}
\end{remarksqed}

\subsection{Exemples simples de foncteurs}\label{subsec154}

\noindent $\bullet$ {\bf Les foncteurs d'oubli:}

\smallskip

Associer \`a un groupe [resp. mono{\"\i}de, resp. anneau, resp. corps, resp. module sur un anneau, resp. espace vectoriel sur un corps, resp. espace topologique] son ensemble sous-jacent d\'efinit un foncteur, appel\'e ``foncteur d'oubli'' (au sens d'oubli de la structure), de la cat\'egorie des groupes [resp. mono{\"\i}des, resp. anneaux, resp. corps, resp. modules sur un anneau fix\'e, resp. espaces vectoriels sur un corps fix\'e, resp. espaces topologiques] dans la cat\'egorie ${\rm Ens}$ des ensembles.

\medskip

\noindent $\bullet$ {\bf Les foncteurs de substitution d'une action par une autre:}

\smallskip

Un morphisme de groupes ou plus g\'en\'eralement de mono{\"\i}des
$$
\rho : M_1 \longrightarrow M_2
$$
peut \^etre vu comme un foncteur.

\smallskip

La composition d'un tel foncteur $\rho$ avec les foncteurs
$$
M_2 \longrightarrow {\rm Ens}
$$
d\'efinit un foncteur induit
$$
\xymatrix{
\rho^* : BM_2 \ar@{=}[d] \ar[r] &BM_1 \ar@{=}[d] \\
[M_2 , {\rm Ens}] &[M_1 , {\rm Ens}]
}
$$
qui consiste \`a associer \`a tout ensemble muni d'une action de $M_2$ le m\^eme ensemble muni de l'action induite de $M_1$.

\smallskip

En particulier, pour tout groupe ou plus g\'en\'eralement tout mono{\"\i}de $M$, l'unique morphisme
$$
\rho : \{1\} \longrightarrow M
$$
d\'efinit le foncteur d'oubli de l'action de $M$
$$
\rho^* : BM \longrightarrow {\rm Ens} \, .
$$
Le foncteur $\rho^*$ est un objet de la cat\'egorie $[BM,{\rm Ens}]$ et, en tant que tel, il a des endomorphismes.

\smallskip

La proposition suivante montre que l'objet $\rho^*$ de $[BM,{\rm Ens}]$ a $M$ pour mono{\"\i}de des endomorphismes et donc que le mono{\"\i}de $M$ peut \^etre reconstruit \`a partir de la cat\'egorie $BM$ munie du foncteur d'oubli
$$
BM \longrightarrow {\rm Ens} \, .
$$

\begin{prop}\label{propI56}

Soient $M$ un groupe ou plus g\'en\'eralement un mono{\"\i}de, $BM$ la cat\'egorie des actions de $M$ et $\rho^* : BM \to {\rm Ens}$ le foncteur d'oubli des actions.

\smallskip

Alors:

\begin{listeimarge}

\item Tout \'el\'ement de $M$ d\'efinit un endomorphisme de $\rho^*$.


\item La loi de multiplication de $M$ correspond \`a la loi de composition des endomorphismes de $\rho^*$.


\item Tout endomorphisme de $\rho^*$ provient d'un unique \'el\'ement de $M$.
\end{listeimarge}
\end{prop}

\begin{demo}
\begin{listeisansmarge}
\item Un \'el\'ement $m \in M$ d\'efinit une application
\begin{eqnarray}
X &\longrightarrow &X \nonumber \\
x &\longmapsto &m \cdot x \nonumber
\end{eqnarray}
pour tout objet $X$ de $BM$. 

\smallskip

De plus, pour toute application $M$-\'equivariante $X_1 \xrightarrow{ \ f \ } X_2$, le carr\'e
$$
\xymatrix{
X_1 \ar[dd]_{\begin{matrix} x_1 \\ \mbox{\rotatebox{-90}{$\mapsto$}} \\ m \cdot x_1 \end{matrix}} \ar[rr]^f &&X_2 \ar[dd]^{\begin{matrix} x_2 \\ \mbox{\rotatebox{-90}{$\mapsto$}} \\ m \cdot x_2 \end{matrix}} \\
{ \ } \\
X_1 \ar[rr]^f &&X_2
}
$$
est commutatif.

\medskip

\item r\'esulte de ce que, par d\'efinition des actions, on a pour tout objet $X$ de $BM$
$$
m_1 \cdot (m_2 \cdot x) = (m_1 \cdot m_2) \cdot x \, , \quad \forall x \in X \, , \ \forall \, m_1 , m_2 \in M \, .
$$

\item Soit $u$ un endomorphisme de l'objet $\rho^*$ de $[BM,{\rm Ens}]$.

\smallskip

Soit $X_0$ l'objet de $BM$ constitu\'e de l'ensemble $M$ muni de la multiplication \`a gauche $M \times M \to M$ par les \'el\'ements de $M$.

\smallskip

L'endomorphisme $u$ induit une application
$$
u_{X_0} : X_0 \longrightarrow X_0
$$
et on peut noter $m = u_{X_0} (1)$.

\smallskip

Pour tout objet $X$ de $BM$ et tout \'el\'ement $x \in X$, l'application
\begin{eqnarray}
f : X_0 = M &\longrightarrow &X \, , \nonumber \\
m &\longmapsto &m \cdot x \nonumber
\end{eqnarray}
est l'unique morphisme $X_0 \to X$ de $BM$ qui envoie l'\'el\'ement $1 \in M = X_0$ sur l'\'el\'ement $x \in X$.

\smallskip

La commutativit\'e du carr\'e
$$
\xymatrix{
X_0 \ar[d]_{u_{X_0}} \ar[r]^f &X \ar[d]^{u_X} \\
X_0 \ar[r]^f &X
}
$$
montre que $u_X (x) = m \cdot x$. 

\end{listeisansmarge}
\end{demo}

\noindent $\bullet$ {\bf Les foncteurs de composition entre cat\'egories de pr\'efaisceaux:}

\smallskip

Tout foncteur entre deux cat\'egories localement petites
$$
\rho : {\mathcal C}_1 \longrightarrow {\mathcal C}_2
$$
d\'efinit par composition un foncteur induit:
$$
\xymatrix{
\widehat{{\mathcal C}_2} \ar@{=}[d] \ar[r]^{\rho^*} &\widehat{{\mathcal C}_1} \ar@{=}[d] \\
[{\mathcal C}_2^{\rm op} , {\rm Ens}] &[{\mathcal C}_1^{\rm op} , {\rm Ens}]
}
$$
On remarque que, pour toute paire de foncteurs
$$
{\mathcal C}_1 \xrightarrow{ \ \rho \ } {\mathcal C}_2 \xrightarrow{ \ \rho' \ } {\mathcal C}_3 \, ,
$$
on a $(\rho' \circ \rho)^* = \rho^* \circ \rho'^*$.

\smallskip

En particulier, si $O_1$ et $O_2$ sont deux ensembles (partiellement) ordonn\'es, toute application
$$
f : O_1 \longrightarrow O_2
$$
qui respecte les relations d'ordre d\'efinit un foncteur
$$
\xymatrix{
f^* : &\widehat{O_2} \ar@{=}[d] \ar[r] &\widehat{O_1} \ar@{=}[d] \\
&[O_2^{\rm op} , {\rm Ens}] &[O_1^{\rm op} , {\rm Ens}]
}
$$
et, pour toute paire de telles applications
$$
O_1 \xrightarrow{ \ f \ } O_2 \xrightarrow{ \ f' \ } O_3 \, ,
$$
on a
$$
(f' \circ f)^* = f^* \circ f'^* \, .
$$

Cette construction g\'en\'eralise celle du foncteur
$$
\xymatrix{
f_* : &\widehat{O(X)} \ar@{=}[d] \ar[r] &\widehat{O(Y)} \ar@{=}[d] \\
&[O(X)^{\rm op} , {\rm Ens}] \ar[r] &[O(Y)^{\rm op} , {\rm Ens}]
}
$$
associ\'e \`a toute application continue entre espaces topologiques
$$
f : X \longrightarrow Y
$$
via l'application induite entre les ensembles d'ouverts
$$
f^{-1} : O(Y) \longrightarrow O(X) \, .
$$
Pour toute paire d'applications continues
$$
X \xrightarrow{ \ f \ } Y \xrightarrow{ \ g \ } Z \, ,
$$
on a donc
$$
(g \circ f)_* = g_* \circ f_* \, .
$$

\bigskip

\noindent $\bullet$ {\bf Les foncteurs d'image directe entre cat\'egories de faisceaux:}

\smallskip

On a d\'ej\`a mentionn\'e que pour tous espaces topologiques $X,Y$ reli\'es par une application continue
$$
f : X \longrightarrow Y \, ,
$$
le foncteur induit
$$
f_* : \widehat{O(X)} \longrightarrow \widehat{O(Y)}
$$
se restreint en un foncteur entre les sous-cat\'egories pleines des faisceaux sur $X$ et $Y$
$$
f_* : {\mathcal E}_X \longrightarrow {\mathcal E}_Y
$$
appel\'e foncteur d'image directe.

\smallskip

Cela r\'esulte de ce que, pour tout recouvrement ouvert
$$
V = \bigcup_{i \in I} \, V_i
$$
d'un ouvert $V$ de $Y$, on a
$$
f^{-1} V = \bigcup_{i \in I} \, f^{-1} V_i
$$
et
$$
f^{-1} V_i \cap f^{-1} V_j = f^{-1} (V_i \cap V_j) \, , \quad \forall \, i,j \in I \, .
$$

Pour toute paire d'applications continues
$$
X \xrightarrow{ \ f \ } Y \xrightarrow{ \ g \ } Z \, ,
$$
on a encore
$$
(g \circ f)_* = g_* \circ f_* : {\mathcal E}_X \longrightarrow {\mathcal E}_Y \longrightarrow {\mathcal E}_Z \, .
$$

\medskip

\noindent $\bullet$ {\bf Les foncteurs d'image directe en les points:}

\smallskip

Tout point $x$ d'un espace topologique $X$ peut \^etre vu comme une application continue
$$
x : \{\bullet\} \longrightarrow X
$$
et induit donc un foncteur d'image directe
$$
x_* : {\rm Ens} = {\mathcal E}_{\{\bullet\}} \longrightarrow {\mathcal E}_X \, .
$$
Pour toute application continue $f : X \to Y$ avec $f(x) = y$, on a
$$
y_* = f_* \circ x_* : {\rm Ens} \longrightarrow {\mathcal E}_X \longrightarrow {\mathcal E}_Y \, .
$$


\noindent $\bullet$ {\bf Les foncteurs de restriction \`a un ouvert:}

\smallskip

Pour tout ouvert $U \xhookrightarrow{ \, i \, } X$ d'un espace topologique $X$, l'inclusion $O(U) \hookrightarrow O(X)$ d\'efinit un foncteur
$$
i^* : \widehat{O(X)} \longrightarrow \widehat{O(U)}
$$
qui induit le foncteur
$$
i^* : {\mathcal E}_X \longrightarrow {\mathcal E}_U
$$
de restriction des faisceaux de $X$ \`a $U$.

\bigskip

\noindent $\bullet$ {\bf Les foncteurs de sections des pr\'efaisceaux ou des faisceaux:}

\smallskip

Pour toute cat\'egorie localement petite ${\mathcal C}$, tout objet $Y$ de ${\mathcal C}$ peut \^etre vu comme un foncteur
$$
1 \longrightarrow {\mathcal C}
$$
(si $1$ d\'esigne ici la cat\'egorie \`a un objet et un morphisme ${\rm id}$) et d\'efinit donc un foncteur
\begin{eqnarray}
\widehat{\mathcal C} = [{\mathcal C}^{\rm op} , {\rm Ens}] &\longrightarrow &{\rm Ens} \, , \nonumber \\
P &\longmapsto &P(Y) \, , \nonumber
\end{eqnarray}
appel\'e le foncteur d'\'evaluation en $Y$ ou le foncteur des sections sur $Y$.

\smallskip

En particulier, tout ouvert $U$ d'un espace topologique $X$ d\'efinit un foncteur
\begin{eqnarray}
\widehat{O(X)} &\longrightarrow &{\rm Ens} \, , \nonumber \\
P &\longmapsto &P(U) \nonumber
\end{eqnarray}
et sa restriction \`a la sous-cat\'egorie pleine des faisceaux
\begin{eqnarray}
{\mathcal E}_X &\longrightarrow &{\rm Ens} \, , \nonumber \\
F &\longmapsto &F(U) \, . \nonumber
\end{eqnarray}
Ce sont les foncteurs des sections sur l'ouvert $U$.

\smallskip

Si $U = X$ est l'espace tout entier, ces foncteurs co{\"\i}ncident avec les foncteurs d'image directe
\begin{eqnarray}
p_* : \widehat{O(X)} &\longrightarrow &{\rm Ens} \, , \nonumber \\
{\mathcal E}_X &\longmapsto &{\rm Ens} \nonumber
\end{eqnarray}
associ\'es \`a l'unique application continue $p :  X \to \{\bullet\}$.

\subsection{Espaces annel\'es et morphismes d'espaces annel\'es}\label{subsec155}

Le foncteur d'image directe des faisceaux associ\'e \`a une application continue $f : X \to Y$
$$
\begin{matrix}
f_* : &\hfill {\mathcal E}_X &\longrightarrow &{\mathcal E}_Y \, , \hfill \\
&\hfill F &\longmapsto &f_* F = [V \mapsto F(f^{-1} V)] \, , \hfill \\
&(F_1 \xrightarrow{u} F_2) &\longmapsto &f_* u = [V \mapsto (F_1 (f^{-1} V) \xrightarrow{u} F_2 (f^{-1} V)] 
\end{matrix}
$$
transforme les faisceaux d'anneaux en faisceaux d'anneaux et les morphismes de faisceaux d'anneaux en morphismes de faisceaux d'anneaux.

\smallskip

Cela permet de poser:

\begin{defn}\label{defI57}
\begin{listeimarge}
\item On appelle ``espace annel\'e'' la donn\'ee d'un espace topologique $X$ et d'un faisceau d'anneaux ${\mathcal O}_X$ sur $X$, dit son ``faisceau de structure''.

\medskip

\item On appelle morphisme d'espaces annel\'es
$$
(X, {\mathcal O}_X) \longrightarrow (Y,{\mathcal O}_Y)
$$
la donn\'ee d'une application continue $f : X \to Y$ et d'un morphisme de faisceaux d'anneaux sur $Y$
$$
{\mathcal O}_Y \longrightarrow f_* {\mathcal O}_X \, .
$$

\item Le compos\'e de deux morphismes d'espaces annel\'es
$$
(X , {\mathcal O}_X) \xrightarrow{ \ f \ } (Y,{\mathcal O}_Y) \xrightarrow{ \ g \ } (Z,{\mathcal O}_z)
$$
est d\'efini comme la paire compos\'ee de l'application continue $g \circ f : X \to Z$ et du morphisme de faisceaux d'anneaux sur $Z$ compos\'e
$$
{\mathcal O}_Z \longrightarrow g_* {\mathcal O}_Y \longrightarrow g_* f_* {\mathcal O}_X \, .
$$
\end{listeimarge}
\end{defn}

\begin{remarksqed} 
\begin{listeisansmarge}	

\item Les espaces annel\'es forment une cat\'egorie.

\medskip

\item Un espace annel\'e $(X,{\mathcal O}_X)$ est dit commutatif si son faisceau de structure ${\mathcal O}_X$ est un faisceau d'anneaux commutatifs.

\smallskip

Les espaces annel\'es commutatifs forment une sous-cat\'egorie pleine de celle des espaces annel\'es.
\end{listeisansmarge}
\end{remarksqed}


On peut appeler ``cat\'egories g\'eom\'etriques'' les sous-cat\'egories de la cat\'egorie des espaces annel\'es dont les objets et les morphismes sont caract\'eris\'es localement:

\begin{defn}\label{defI58}

Une sous-cat\'egorie (pas n\'ecessairement pleine) ${\mathcal G}$ de la cat\'egorie ${\rm Top}_{\rm an}$ des espaces annel\'es sera dite g\'eom\'etrique si:

\begin{listeimarge}

\item Pour tout objet $(X,{\mathcal O}_X)$ de ${\mathcal G}$, tout ouvert $U \xhookrightarrow{ \ i \ } X$ de $X$ muni de la restriction ${\mathcal O}_U = i^* {\mathcal O}_{X}$ du faisceau de structure ${\mathcal O}_X$ est un objet de ${\mathcal G}$ et le morphisme $(U \xhookrightarrow{ \ i \ } X , {\mathcal O}_X \to i_* {\mathcal O}_U)$ est un morphisme de ${\mathcal G}$ (appel\'ee une immersion ouverte).

\medskip

\item R\'eciproquement, un objet $(X,{\mathcal O}_X)$ de ${\rm Top}_{\rm an}$ est dans ${\mathcal G}$ si $X$ admet un recouvrement par des ouverts $U_i$ qui, munis des restrictions ${\mathcal O}_{U_i}$ du faisceau de structure ${\mathcal O}_X$, sont des objets de ${\mathcal G}$.

\medskip

\item Si $(X,{\mathcal O}_X)$ et $(Y,{\mathcal O}_Y)$ sont deux objets de ${\mathcal G}$ et les $(U_i , {\mathcal O}_{U_i})$ forment un recouvrement ouvert de $(X,{\mathcal O}_X)$, alors un morphisme d'espaces annel\'es $(X,{\mathcal O}_X) \to (Y,{\mathcal O}_Y)$ est dans ${\mathcal G}$ si (et seulement si) ses restrictions $(U_i , {\mathcal O}_{U_i}) \to (Y , {\mathcal O}_Y)$ sont dans ${\mathcal G}$.

\medskip

\item Si $(X,{\mathcal O}_X)$ et $(Y , {\mathcal O}_Y)$ sont deux objets de ${\mathcal G}$, et $V$ un ouvert de $Y$ muni de la restriction ${\mathcal O}_V$ de ${\mathcal O}_Y$, un morphisme d'espaces annel\'es $(X,{\mathcal O}_X) \to (V,{\mathcal O}_V)$ est dans ${\mathcal G}$ si (et seulement si) le compos\'e $(X,{\mathcal O}_X) \to (Y , {\mathcal O}_Y)$ est dans ${\mathcal G}$. 
\end{listeimarge}

\end{defn}

\subsection{Exemples de cat\'egories g\'eom\'etriques}\label{subsec156}

\noindent $\bullet$ {\bf Les vari\'et\'es diff\'erentielles de classe $C^k$:}

\smallskip

Pour tout $k \geq 1$ qui est un entier ou le symbole $\infty$, les vari\'et\'es de classe $C^k$ et de dimension $n$ sont, d'apr\`es la d\'efinition \ref{defI44}, des espaces annel\'es $(X,{\mathcal O}_X)$ qui admettent un recouvrement par des ouverts $(U_i , {\mathcal O}_{U_i})$ tels qu'existent des hom\'eomorphismes $\sigma_i : U_i \xrightarrow{ \ \sim \ } U'_i$ vers des ouverts $U'_i$ de ${\mathbb R}^n$ qui identifient, par composition avec $\sigma_i^{-1}$, le faisceau ${\mathcal O}_{U_i}$ au faisceau des fonctions de classe $C^k$ sur $U'_i$ \`a valeurs dans ${\mathbb R}$. 

\smallskip

On obtient une cat\'egorie g\'eom\'etrique en compl\'etant cette d\'efinition par celle des morphismes entre vari\'et\'es de classe $C^k$:

\begin{defn}\label{defI59}

Etant donn\'ees deux vari\'et\'es diff\'erentielles $(X,{\mathcal O}_X)$ et $(Y,{\mathcal O}_Y)$ de classe $C^k$, un morphisme entre elles est un morphisme d'espaces annel\'es
$$
(X \xrightarrow{ \ f \ } Y \, , \ {\mathcal O}_Y \longrightarrow f_* {\mathcal O}_X)
$$
tel que, pour tout ouvert $V$ de $Y$, le morphisme d'anneaux
$$
{\mathcal O}_Y (V) \longrightarrow f_* {\mathcal O}_X (V) = {\mathcal O}_X (f^{-1} V)
$$
est la composition des applications $V \to {\mathbb R}$ avec $f : f^{-1} V \to V$.
\end{defn}

\begin{remarkqed}

Cette d\'efinition \'equivaut \`a demander que, pour tout point $x$ de $X$, le morphisme d'espaces annel\'es $(X,{\mathcal O}_X) \to (Y,{\mathcal O}_Y)$ s'identifie dans des voisinages ouverts de $x$ et $f(x)$ \`a une application de classe $C^k$ d'un ouvert de ${\mathbb R}^n$ dans un ouvert d'un ${\mathbb R}^m$.

\end{remarkqed}

\smallskip

\noindent $\bullet$ {\bf Les vari\'et\'es analytiques:}

\smallskip

Les vari\'et\'es analytiques de dimension $n$ sont, d'apr\`es la d\'efinition \ref{defI45}, des espaces annel\'es $(X,{\mathcal O}_X)$ qui admettent un recouvrement par des ouverts $(U_i , {\mathcal O}_{U_i})$ tels qu'existent des hom\'eomorphismes $\sigma_i : U_i \xrightarrow{ \ \sim \ } U'_i$ vers des ouverts $U'_i$ de ${\mathbb C}^n$ qui identifient, par composition avec $\sigma_i^{-1}$, le faisceau ${\mathcal O}_{U_i}$ au faisceau des fonctions holomorphes sur $U'_i$ \`a valeurs dans ${\mathbb C}$.

\smallskip

On obtient une cat\'egorie g\'eom\'etrique en compl\'etant cette d\'efinition par celle des morphismes entre vari\'et\'es analytiques:

\begin{defn}\label{defI510}

Etant donn\'ees deux vari\'et\'es analytiques $(X,{\mathcal O}_X)$ et $(Y,{\mathcal O}_Y)$, un morphisme entre elles est un morphisme d'espaces annel\'es
$$
(X \xrightarrow{ \ f \ } Y \, , \ {\mathcal O}_Y \longrightarrow f_* {\mathcal O}_X)
$$
tel que, pour tout ouvert $V$ de $Y$, le morphisme d'anneaux
$$
{\mathcal O}_Y (V) \longrightarrow f_* {\mathcal O}_X (V) = {\mathcal O}_X (f^{-1} V)
$$
est la composition des applications $V \to {\mathbb C}$ avec $f : f^{-1} V \to V$.
\end{defn}

\begin{remarkqed}

Cette d\'efinition \'equivaut \`a demander que, pour tout point $x$ de $X$, le morphisme d'espaces annel\'es $(X,{\mathcal O}_X) \to (Y,{\mathcal O}_Y)$ s'identifie dans des voisinages ouverts de $x$ et $f(x)$ \`a une application holomorphe d'un ouvert de ${\mathbb C}^n$ dans un ouvert d'un ${\mathbb C}^m$.

\end{remarkqed}


\noindent $\bullet$ {\bf Les sch\'emas:}

\smallskip

Dans la d\'efinition \ref{defI47}, on a associ\'e \`a tout anneau commutatif $A$ un espace annel\'e $({\rm Spec} (A) , {\mathcal O}_A)$, souvent not\'e simplement ${\rm Spec} (A)$, appel\'e le ``spectre'' de $A$.

\smallskip

On appelle sch\'emas affines les espaces annel\'es de la forme ${\rm Spec} (A)$.

\smallskip

Dans la d\'efinition \ref{defI49} (ii), on a appel\'e ``sch\'emas'' les espaces annel\'es $(X , {\mathcal O}_X)$ qui admettent un recouvrement par des ouverts $(U_i , {\mathcal O}_{U_i})$ isomorphes \`a des sch\'emas affines ${\rm Spec} (A_i)$.

\smallskip

On obtient une cat\'egorie g\'eom\'etrique en compl\'etant cette d\'efinition par celle des morphismes de sch\'emas:

\begin{defn}\label{defI511}
\begin{listeimarge}
\item Un morphisme de sch\'emas affines
$$
{\rm Spec} (B) \longrightarrow {\rm Spec} (A)
$$
est un morphisme d'espaces annel\'es
$$
({\rm Spec} (B) \xrightarrow{ \ x \ } {\rm Spec} (A) \, , \ {\mathcal O}_A \longrightarrow x_* {\mathcal O}_B)
$$
induit par un morphisme d'anneaux
$$
u : A \longrightarrow B
$$
au sens que

\medskip

$\left\{ \begin{matrix}
\bullet &\mbox{pour tout id\'eal premier $q \in {\rm Spec} (B)$ de $B$,} \hfill \\
&\mbox{$x(q)$ est l'id\'eal premier $p = u^{-1} (q) \in {\rm Spec} (A)$ de $A$,} \hfill \\
{ \ } \\
\bullet &\mbox{pour tout \'el\'ement $f \in A$ d\'efinissant l'ouvert} \hfill \\
{ \ } \\
&{\rm Spec} (A)_f \quad \mbox{de} \quad {\rm Spec} (A) \\
{ \ } \\
&\mbox{et son image r\'eciproque} \hfill \\
{ \ } \\
&x^{-1} ({\rm Spec} (A)_f) = {\rm Spec} (B)_{u(f)} \quad \mbox{dans} \quad {\rm Spec} (B) \, , \\
{ \ } \\
&\mbox{le morphisme d'anneaux} \hfill \\
{ \ } \\
&A_f = {\mathcal O}_A ({\rm Spec} (A)_f) \longrightarrow x_* {\mathcal O}_B ({\rm Spec} (A)_f) = B_{u(f)} \\
{ \ } \\
&\mbox{est celui induit par $u : A \to B$.} \hfill
\end{matrix} \right.
$

\medskip

\item Un morphisme de sch\'emas
$$
X \longrightarrow Y
$$
est un morphisme d'espaces annel\'es
$$
(X \xrightarrow{ \ f \ } Y \, , \ {\mathcal O}_Y \longrightarrow f_* {\mathcal O}_X)
$$
tel que, pour tout point $x \in X$ d'image $f(x) \in Y$, il existe des voisinages ouverts $(U , {\mathcal O}_U)$ et $(V , {\mathcal O}_V)$ de $x$ et $f(x)$ dans les sch\'emas $X$ et $Y$, et des isomorphismes de ces voisinages vers des sch\'emas affines ${\rm Spec} (B)$ et ${\rm Spec} (A)$ induisant un carr\'e commutatif de la cat\'egorie des espaces annel\'es
$$
\xymatrix{
{\rm Spec} (B) \ar[d] \ar[r] &X \ar[d] \\
{\rm Spec} (A) \ar[r] &Y
}
$$
tel que le morphisme ${\rm Spec} (B) \to {\rm Spec} (A)$ soit un morphisme de sch\'emas affines.

\end{listeimarge}
\end{defn}


\begin{remarksqed}
\begin{listeisansmarge}
\item La cat\'egorie des sch\'emas est une sous-cat\'egorie g\'eom\'etrique de celle des espaces annel\'es, mais ce n'est pas le cas de la cat\'egorie des sch\'emas affines: en effet, un ouvert d'un sch\'ema affine n'est pas n\'ecessairement un sch\'ema affine.

\medskip

\item Cependant, associer \`a un anneau commutatif $A$ son spectre ${\rm Spec} (A)$ d\'efinit un foncteur contravariant de la cat\'egorie des anneaux commutatifs vers celle des espaces annel\'es.
\end{listeisansmarge}
\end{remarksqed}

\pagebreak

On note:

\begin{defn}\label{defI512}
\begin{listeimarge}
\item Pour tout sch\'ema $X$ et tout sch\'ema affine ${\rm Spec} (A)$, se donner un morphisme de sch\'emas
$$
X \longrightarrow {\rm Spec} (A)
$$
\'equivaut \`a se donner un morphisme d'anneaux commutatifs
$$
A \longrightarrow {\mathcal O}_X (X) \, .
$$

\item En particulier, pour tous sch\'emas affines ${\rm Spec} (A)$ et ${\rm Spec} (B)$, se donner un morphisme de sch\'emas
$$
{\rm Spec} (B) \longrightarrow {\rm Spec} (A)
$$
\'equivaut \`a se donner un morphisme d'anneaux commutatifs
$$
A \longrightarrow B \, .
$$

\end{listeimarge}
\end{defn}

\begin{remark}
\begin{listeisansmarge}
\item[(ii)] signifie que la cat\'egorie des sch\'emas affines est une sous-cat\'egorie pleine de celle des sch\'emas.
\end{listeisansmarge}
\end{remark}

\begin{demo}
\begin{listeisansmarge}
\item[(ii)] est un cas particulier de (i) puisque
$$
{\mathcal O}_B ({\rm Spec} (B)) = B
$$
pour tout anneau commutatif $B$.

\medskip

\item[(i)] Consid\'erons n'importe quel recouvrement du sch\'ema $X$ par des ouverts $U_i$ isomorphes en tant qu'espaces annel\'es \`a des sch\'emas affines ${\rm Spec} (B_i)$.

\smallskip

Par d\'efinition, se donner une famille de morphismes de sch\'emas affines
$$
{\rm Spec} (B_i) \longrightarrow {\rm Spec} (A)
$$
\'equivaut \`a se donner une famille de morphismes d'anneaux commutatifs
$$
A \longrightarrow B_i \, .
$$

De plus, demander que les morphismes correspondants d'espaces annel\'es
$$
(U_i , {\mathcal O}_{U_i}) \longrightarrow {\rm Spec} (A)
$$
co{\"\i}ncident sur les intersections $U_i \cap U_j$ \'equivaut \`a demander que les morphismes compos\'es
$$
A \longrightarrow B_i \xrightarrow{ \ \sim \ } {\mathcal O}_{U_i} (U_i) = {\mathcal O}_X (U_i)
$$
s'inscrivent dans des carr\'es commutatifs
$$
\xymatrix{
A \ar[d] \ar[r] &{\mathcal O}_X (U_i) \ar[d] \\
{\mathcal O}_X (U_j) \ar[r] &{\mathcal O}_X (U_i \cap U_j)
}
$$
index\'es par les paires d'indices $i,j$.

\smallskip

Comme ${\mathcal O}_X$ est un faisceau, cela \'equivaut \`a demander que les morphismes $A \to {\mathcal O}_X (U_i)$ se rel\`event en un unique morphisme $A \to {\mathcal O}_X(X)$. 

\end{listeisansmarge}
\end{demo}

\section{Isomorphismes, monomorphismes, \'epimorphismes}\label{sec16}

\subsection{Traduction cat\'egorique des notions de bijection, d'injection et de surjection}\label{subsec161}

Les notions d'applications bijectives, injectives ou surjectives jouent un tr\`es grand r\^ole dans les math\'emati\-ques formul\'ees dans le langage des ensembles.

\smallskip

Les exprimer dans le langage des cat\'egories, c'est-\`a-dire en termes de fl\`eches et non pas de points, permet de d\'efinir des notions qui ont un sens dans le contexte de n'importe quelle cat\'egorie.

\smallskip

Ainsi, les notions ensemblistes de bijection, d'injection et de surjection sont g\'en\'eralis\'ees en les notions cat\'egoriques d'isomorphisme, de monomorphisme et d'\'epimorphisme:

\begin{defn}\label{defI61}

Soit $u : X \to Y$ un morphisme d'une cat\'egorie ${\mathcal C}$. Alors:

\begin{listeimarge}

\item Le morphisme $u$ est appel\'e un ``isomorphisme'' s'il existe un morphisme en sens inverse (n\'ecessairement unique)
$$
u^{-1} : Y \longrightarrow X \, , \ \mbox{appel\'e l'inverse de $u$},
$$
tel que
$$
u^{-1} \circ u = {\rm id}_X \quad \mbox{et} \quad u \circ u^{-1} = {\rm id}_Y \, .
$$

\item Il est appel\'e un ``monomorphisme'' si, pour tout objet $X'$ et tous morphismes $f,g : X' \rightrightarrows X$, on a l'implication
$$
u \circ f = u \circ g \Rightarrow f=g \, .
$$

\item Il est appel\'e un ``\'epimorphisme'' si, pour tout objet $Y'$ et tous morphismes $f,g : Y \rightrightarrows Y'$, on a l'implication
$$
f \circ u = g \circ u \Rightarrow f=g \, .
$$
\end{listeimarge}
\end{defn}

\begin{remarks}
\begin{listeisansmarge}
\item Un morphisme $u$ d'une cat\'egorie ${\mathcal C}$ est un isomorphisme de ${\mathcal C}$ si et seulement si c'est un isomorphisme de ${\mathcal C}^{\rm op}$.

\smallskip

D'autre part, c'est un monomorphisme [resp. \'epimorphisme] de ${\mathcal C}$ si et seulement si c'est un \'epimorphisme [resp. monomorphisme] de ${\mathcal C}^{\rm op}$.

\smallskip

On dit que la propri\'et\'e d'\^etre un isomorphisme est auto-duale et que les propri\'et\'es d'\^etre un monomorphisme ou un \'epimorphisme sont duales l'une de l'autre.

\medskip

\item Le compos\'e de deux isomorphismes $u : X \to Y$ et $v : Y \to Z$ est un isomorphisme et on a
$$
(v \circ u)^{-1} = u^{-1} \circ v^{-1} \, .
$$

Le compos\'e de deux monomorphismes [resp. \'epimorphismes] est un \'epimorphisme.

\medskip

\item Tout isomorphisme d'une cat\'egorie ${\mathcal C}$ est \`a la fois un monomorphisme et un \'epimorphisme.

\smallskip

La r\'eciproque n'est pas vraie en g\'en\'eral.

\smallskip

Une cat\'egorie ${\mathcal C}$ est appel\'ee ``balanc\'ee'' si la r\'eciproque est vraie, c'est-\`a-dire si tout morphisme de ${\mathcal C}$ qui est \`a la fois un monomorphisme et un \'epimorphisme est un isomorphisme.

\medskip

\item Pour tout objet $X$ d'une cat\'egorie ${\mathcal C}$, les isomorphismes $X \xrightarrow{ \ \sim \ } X$ sont appel\'es les automorphismes de $X$.

\medskip

\item Si ${\mathcal C}$ est une cat\'egorie localement petite, les automorphismes d'un objet $X$ de ${\mathcal C}$ forment un groupe not\'e ${\rm Aut} (X)$.

\smallskip

Ce groupe n'est autre que le groupe ${\rm End} (X)^{\times}$ des \'el\'ements inversibles du mono{\"\i}de ${\rm End} (X)$ des endomorphismes de $X$.

\medskip

\item Un foncteur $F : {\mathcal C} \to {\mathcal D}$ transforme tout isomorphisme $u : X \to Y$ de ${\mathcal C}$ en un isomorphisme $F(u) : F(X) \to F(Y)$ de ${\mathcal D}$ avec
$$
F(u)^{-1} = F(u^{-1}) \, .
$$

On dit qu'un tel foncteur $F$ est ``conservatif'' si tout morphisme $u : X \to Y$ de ${\mathcal C}$ tel que $F(u) : F(X) \to F(Y)$ soit un isomorphisme de ${\mathcal D}$ est un isomorphisme de ${\mathcal C}$.

\medskip

\item En g\'en\'eral, le transform\'e par un foncteur $F : {\mathcal C} \to {\mathcal D}$ d'un monomorphisme [resp. \'epimorphisme] de ${\mathcal C}$ n'est pas un monomorphisme [resp. \'epimorphisme] de ${\mathcal D}$.
\end{listeisansmarge}
\end{remarks}

\subsection{Exemples de traductions des mots ``isomorphisme'', ``monomorphisme'' et ``\'epimorphisme'' dans diff\'erentes cat\'egories}\label{subsec162}

\noindent $\bullet$ {\bf La cat\'egorie des ensembles et les cat\'egories de pr\'efaisceaux:}

\smallskip

Les notions d'isomorphisme, de monomorphisme et d'\'epimorphisme ont \'et\'e d\'efinies de telle fa\c con qu'un morphisme $X \to Y$ de la cat\'egorie ${\rm Ens}$ des ensembles, c'est-\`a-dire une application, soit un isomorphisme [resp. un monomorphisme, resp. un \'epimorphisme] si et seulement si c'est une bijection [resp. une injection, resp. une surjection].

\smallskip

Plus g\'en\'eralement, pour toute cat\'egorie localement petite ${\mathcal C}$, un morphisme de la cat\'egorie $\widehat{\mathcal C} = [{\mathcal C}^{\rm op} , {\rm Ens}]$ des pr\'efaisceaux sur ${\mathcal C}$
$$
P_1 \longrightarrow P_2
$$
est un isomorphisme [resp. un monomorphisme, resp. un \'epimorphisme] si et seulement si, pour tout objet $X$ de ${\mathcal C}$, l'application
$$
P_1 (X) \longrightarrow P_2(X)
$$
est bijective [resp. injective, resp. surjective].

\smallskip

En particulier, les cat\'egories de pr\'efaisceaux $\widehat{\mathcal C}$ sont toujours balanc\'ees.

\smallskip

Le cas o\`u ${\mathcal C}$ a un unique objet est celui des cat\'egories $BM$ d'actions d'un mono{\"\i}de $M$.

\smallskip

Un morphisme de $BM$, c'est-\`a-dire une application $M$-\'equivariante
$$
X_1 \longrightarrow X_2
$$
est un isomorphisme [resp. un monomorphisme, resp. un \'epimorphisme] si et seulement si l'application sous-jacente est bijective [resp. injective, resp. surjective].

\smallskip

En particulier, le foncteur d'oubli de l'action de $M$
$$
BM \longrightarrow {\rm Ens}
$$
est conservatif.

\medskip

\noindent $\bullet$ {\bf Les cat\'egories de faiseaux:}

\smallskip

Un morphisme de la cat\'egorie ${\mathcal E}_X$ des faisceaux sur un espace topologique $X$
$$
u : F_1 \longrightarrow F_2
$$
est un isomorphisme [resp. un monomorphisme] si et seulement si, pour tout ouvert $U$ de $X$, l'application
$$
u_U : F_1 (U) \longrightarrow F_2 (U)
$$
est bijective [resp. injective].

\smallskip

C'est un \'epimorphisme si, pour tout ouvert $U$ de $X$ et tout \'el\'ement $s \in F_2 (U)$, il existe un recouvrement de $U$ par des ouverts $U_i$, $i \in I$, et pour tout \'el\'ement $i \in I$ un \'el\'ement $s_i \in F_1 (U_i)$ tel que
$$
u_{U_i} (s_i) = r_{U_i,U} (s) \quad \mbox{dans} \quad F_2 (U_i) \, .
$$

En particulier, pour tout espace topologique $X$, la cat\'egorie ${\mathcal E}_X$ des faisceaux sur $X$ est balanc\'ee et le foncteur de plongement des faisceaux dans les pr\'efaisceaux
$$
{\mathcal E}_X \longrightarrow \widehat{O(X)}
$$
est conservatif.

\medskip

\noindent $\bullet$ {\bf Les cat\'egories d'ensembles munis d'une structure alg\'ebrique:}

\smallskip

Si ${\mathcal C}$ est la cat\'egorie des ensembles munis d'un certain type de structure alg\'ebrique -- telle que celle de groupe, de mono{\"\i}de, d'anneau, d'espace vectoriel sur un corps ou de module sur un anneau --, alors un morphisme de ${\mathcal C}$ est un isomorphisme [resp. un monomorphisme, resp. un \'epimorphisme] si et seulement si l'application sous-jacente est bijective [resp. injective, resp. surjective].

\smallskip

En particulier, une telle cat\'egorie ${\mathcal C}$ est balanc\'ee et son foncteur d'oubli de la structure alg\'ebrique
$$
{\mathcal C} \longrightarrow {\rm Ens}
$$
est conservatif.

\medskip

\noindent $\bullet$ {\bf La cat\'egorie des espaces topologiques:}

\smallskip

Un morphisme de la cat\'egorie ${\rm Top}$, c'est-\`a-dire une application continue entre deux espaces topologiques
$$
f : X \longrightarrow Y
$$
est un isomorphisme si et seulement si elle est bijective et identifie les topologies de $X$ et $Y$, c'est-\`a-dire est un hom\'eomorphisme.

\smallskip

Elle est un monomorphisme [resp. un \'epimorphisme] si et seulement si elle est injective [resp. surjective].

\smallskip

En particulier, la cat\'egorie ${\rm Top}$ n'est pas balanc\'ee et son foncteur d'oubli de la structure topologique
$$
{\rm Top} \longrightarrow {\rm Ens}
$$
n'est pas conservatif.

\medskip

\noindent $\bullet$ {\bf La sous-cat\'egorie des espaces topologiques s\'epar\'es}

\smallskip

Un espace topologique $Y$ est s\'epar\'e si et seulement si pour toute paire d'applications continues
$$
f,g : X \rightrightarrows Y
$$
d'un espace topologique $X$ dans $Y$, le sous-ensemble
$$
\{x \in X \mid f(x) = g(x)\}
$$
est un ferm\'e de $X$.

\smallskip

On en d\'eduit que, dans la sous-cat\'egorie pleine ${\rm Top}_s$ de ${\rm Top}$ constitu\'ee des espaces topologiques s\'epar\'es, un morphisme c'est-\`a-dire une application continue
$$
f : X \longrightarrow Y
$$
est un \'epimorphisme si (et seulement si) l'image de $X$ dans $Y$ est dense.

\smallskip

D'autre part, $f : X \to Y$ est un isomorphisme [resp. un monomorphisme] de ${\rm Top}_s$ si et seulement si c'en est un dans ${\rm Top}$.

\smallskip

Ainsi, le foncteur de plongement
$$
{\rm Top}_s \longrightarrow {\rm Top}
$$
est conservatif et pr\'eserve les monomorphismes mais il ne pr\'eserve pas les \'epimorphismes.

\subsection{Sous-objets et quotients dans une cat\'egorie}\label{subsec163}

\smallskip

Les notions de sous-ensemble et d'ensemble quotient peuvent s'exprimer dans le langage des cat\'egories et donc se g\'en\'eraliser de la cat\'egorie des ensembles \`a toutes les cat\'egories:

\begin{defn}\label{defI62}

Soit $X$ un objet d'une cat\'egorie ${\mathcal C}$.

\begin{listeimarge}

\item Deux monomorphismes dans $X$
$$
i_1 : X_1 \xhookrightarrow{ \ \ } X \qquad \mbox{et} \qquad i_2 : X_2 \xhookrightarrow{ \ \ } X 
$$
sont dits \'equivalents s'il existe un isomorphisme de ${\mathcal C}$
$$
u : X_1 \xrightarrow{ \ \sim \ } X_2
$$
tel que $i_1 = i_2 \circ u $.

\smallskip

Une classe d'\'equivalence de monomorphismes $X' \hookrightarrow X$ pour cette relation est appel\'ee un ``sous-objet'' de $X$.

\medskip

\item Deux \'epimorphismes de $X$
$$
p_1 : X \twoheadrightarrow X_1 \qquad \mbox{et} \qquad p_2 : X \twoheadrightarrow X_2
$$
sont dits \'equivalents s'il existe un isomorphisme de ${\mathcal C}$
$$
u : X_1 \xrightarrow{ \ \sim \ } X_2
$$
tel que $u \circ p_1 = p_2$.

\smallskip

Une classe d'\'equivalence d'\'epimorphismes $X \twoheadrightarrow X'$ pour cette relation est appel\'ee un ``quotient'' de $X$.
\end{listeimarge}
\end{defn}

\begin{remarksqed}
\begin{listeisansmarge}
\item Etant donn\'e un objet $X$ d'une cat\'egorie ${\mathcal C}$, les sous-objets [resp. quotients] de $X$ dans ${\mathcal C}$ sont les quotients [resp. sous-objets] de $X$ dans ${\mathcal C}^{\rm op}$.

\smallskip

Autrement dit, les notions de quotient et de sous-objet sont duales l'une de l'autre.

\medskip

\item Etant donn\'es deux sous-objets $\overline X_1$ et $\overline X_2$ d'un objet $X$ repr\'esent\'es par deux monomorphismes $i_1 : X_1 \hookrightarrow X$ et $i_2 : X_2 \hookrightarrow X$, il existe au plus un morphisme $i : X_1 \to X_2$ de ${\mathcal C}$ tel que
$$
i_1 = i_2 \circ i
$$
et, s'il existe, c'est un monomorphisme.

\smallskip

On dit alors que $\overline X_1 \leq \overline X_2$.

\smallskip

Cela d\'efinit une relation d'ordre sur les sous-objets de $X$.

\smallskip

Les sous-objets d'un sous-objet $\overline X_1$ de $X$ s'identifient aux sous-objets $\overline X_2$ de $\overline X$ tels que $\overline X_2 \leq \overline X_1$.

\medskip

\item De m\^eme, il y a une relation d'ordre sur les quotients d'un objet $X$.

\smallskip

Elle est telle que les quotients d'un quotient $\overline X_1$ de $X$ s'identifient aux quotients $\overline X_2$ de $X$ tels que $\overline X_2 \leq \overline X_1$. 
\end{listeisansmarge}
\end{remarksqed}

\subsection{Exemples de sous-objets et de quotients dans des cat\'egories}\label{subsec164}

\medskip

\noindent $\bullet$ {\bf La cat\'egorie des ensembles et les cat\'egories de pr\'efaisceaux:}

\smallskip

Les notions cat\'egoriques de sous-objet et de quotient ont \'et\'e d\'efinies de telle fa\c con qu'un sous-objet [resp. un quotient] d'un objet $E$ de la cat\'egorie ${\rm Ens}$ est un sous-ensemble [resp. un ensemble quotient, d\'efini par une relation d'\'equivalence] de l'ensemble $E$.

\smallskip

Plus g\'en\'eralement, si ${\mathcal C}$ est une cat\'egorie localement petite, un sous-objet [resp. un quotient] d'un objet $P$ de la cat\'egorie $\widehat{\mathcal C} = [{\mathcal C}^{\rm op} , {\rm Ens}]$ des pr\'efaisceaux sur ${\mathcal C}$ est une collection de sous-ensembles [resp. d'ensembles quotients] $P'(X)$ des $P(X)$, $X \in {\rm Ob} ({\mathcal C})$, telle que, pour toute fl\`eche $X_1 \xrightarrow{ \ u \ } X_2$ de ${\mathcal C}$, l'application $P(u) : P(X_2) \to P(X_1)$ induise une application $P' (X_2) \to P'(X_1)$.

\smallskip

Dans le cas des cat\'egories $BM$ des actions de mono{\"\i}des $M$, les sous-objets [resp. les quotients] d'un objet $X$ sont les sous-ensembles de $X$ qui sont respect\'es par l'action de $M$ [resp. correspondent aux relations d'\'equivalences de $X$ respect\'ees par $M$] .

\medskip

\noindent $\bullet$ {\bf Les cat\'egories de faisceaux:}

\smallskip

Un sous-objet d'un objet $F$ de la cat\'egorie ${\mathcal E}_X$ des faisceaux sur un espace topologique $X$ est un sous-pr\'efaisceau $F'$ de $F$ consid\'er\'e comme pr\'efaisceau qui est un faisceau. Ainsi, tout sous-objet de $F$ dans la cat\'egorie ${\mathcal E}_X$ est un sous-objet de $F$ dans la cat\'egorie $\widehat{O(X)}$ des pr\'efaisceaux sur $X$, mais la r\'eciproque n'est pas vraie en g\'en\'eral.

\smallskip

Un quotient d'un objet $F$ de ${\mathcal E}_X$ est une classe d'isomorphie de faisceaux $F'$ munis d'un morphisme $F \to F'$ tel que toute section $s$ de $F'$ sur un ouvert $U$ de $X$ se rel\`eve en des sections $s_i \in F(U_i)$ de $F$ sur des ouverts $U_i$ qui recouvrent $U$. Pour tout tel quotient $F'$ et tout ouvert $U$, on peut noter $R(U)$ le sous-ensemble de $F(U) \times F(U)$ constitu\'e des paires de sections qui ont m\^eme image dans $F'(U)$. C'est une relation d'\'equivalence dans $F(U)$. De plus, la famille des $R(U)$, $U \in O(X)$, forme un sous-faisceau de $F \times F$.

\smallskip

On v\'erifiera plus loin (et dans le cadre plus g\'en\'eral des topos) que, r\'eciproquement, tout sous-faisceau $R$ de $F \times F$ tel que chaque $R(U)$ soit une relation d'\'equivalence dans $F(U)$, correspond \`a un unique quotient de $F$ dans la cat\'egorie ${\mathcal E}_X$.

\medskip

\noindent $\bullet$ {\bf Les cat\'egories d'ensembles munis d'une structure alg\'ebrique:}

\smallskip

Si ${\mathcal C}$ est la cat\'egorie des ensembles munis d'un certain type de structure alg\'ebrique -- telle que celle de groupe, de mono{\"\i}de, d'anneau ou de module sur un anneau --, alors un sous-objet $X'$ d'un objet $X$ est un sous-ensemble [resp. un ensemble quotient] de $X$ sur lequel la structure de $X$ induit une structure de m\^eme type.

\smallskip

On retrouve les notions usuelles de sous-groupes [resp. de groupe quotient par un sous-groupe distingu\'e], de sous-mono{\"\i}de [resp. de mono{\"\i}de quotient], de sous-anneau [resp. d'anneau quotient par un id\'eal], de sous-module [resp. de module quotient par un sous-module].


\noindent $\bullet$ {\bf La cat\'egorie des petites cat\'egories:}

\smallskip

Un sous-objet d'un objet ${\mathcal C}$ de la cat\'egorie ${\rm Cat}$ des petites cat\'egories est une sous-cat\'egorie de ${\mathcal C}$ au sens de la d\'efinition \ref{defI52}.

\smallskip

Un quotient ${\mathcal C}'$ d'un objet ${\mathcal C}$ de ${\rm Cat}$ tel que ${\rm Ob} ({\mathcal C}') = {\rm Ob} ({\mathcal C})$ est une cat\'egorie quotient de ${\mathcal C}$ au sens suivant:

\begin{defn}\label{defI63}
Soit ${\mathcal C}$ une cat\'egorie localement petite.

\smallskip

Une cat\'egorie quotient ${\mathcal C}'$ de ${\mathcal C}$ (sous-entendu: ayant m\^emes objets) consiste en une famille d'ensembles quotients
$$
{\rm Hom}_{{\mathcal C}'} (X,Y) \quad \mbox{des} \quad {\rm Hom}_{\mathcal C} (X,Y) \, , \qquad X,Y \in {\rm Ob} ({\mathcal C}) = {\rm Ob} ({\mathcal C}') \, ,
$$
telle que, pour tous objets $X,Y,Z$, la loi de composition
\begin{eqnarray}
{\rm Hom}_{\mathcal C} (Y,Z) \times {\rm Hom}_{\mathcal C} (X,Y) &\longrightarrow &{\rm Hom}_{\mathcal C} (X,Z) \, , \nonumber \\
(v,u) &\longmapsto &v \circ u \nonumber
\end{eqnarray}
induise par passage aux quotients une application
$$
{\rm Hom}_{{\mathcal C}'} (Y,Z) \times {\rm Hom}_{{\mathcal C}'} (X,Y) \longrightarrow {\rm Hom}_{{\mathcal C}'} (X,Z) \, .
$$
\end{defn}

\begin{remarkqed}

Ainsi, ${\mathcal C}'$ est munie de l'unique structure de cat\'egorie telle que les applications de passage au quotient
$$
{\rm Hom}_{\mathcal C} (X,Y) \longrightarrow {\rm Hom}_{{\mathcal C}'} (X,Y) \, , \qquad X,Y \in {\rm Ob} ({\mathcal C}) = {\rm Ob} ({\mathcal C}') \, ,
$$
d\'efinissent un foncteur
$$
F : {\mathcal C} \longrightarrow {\mathcal C}'
$$
avec $F(X) = X$, $\forall \, X \in {\rm Ob} ({\mathcal C}) = {\rm Ob} ({\mathcal C}')$. 
\end{remarkqed}
\medskip

\noindent $\bullet$ {\bf La cat\'egorie des espaces topologiques:}

\smallskip

Un sous-objet [resp. un quotient] d'un objet $X$ de la cat\'egorie ${\rm Top}$ des espaces topologiques est un sous-ensemble [resp. un ensemble quotient] de $X$ muni d'une topologie qui contient [resp. contenue dans] la topologie induite par $X$.

\smallskip

La notion de sous-espace topologique [resp. d'espace topologique quotient] d'un espace topologique $X$ introduite dans la remarque (v) [resp. (vi)] suivant la d\'efinition \ref{defI31} est donc plus stricte que celle de sous-objet [resp. de quotient] de $X$ dans ${\rm Top}$.

\medskip

\noindent $\bullet$ {\bf La cat\'egorie des espaces topologiques s\'epar\'es:}

\smallskip

Les sous-objets d'un objet $X$ de la cat\'egorie ${\rm Top}_s$ des espaces topologiques s\'epar\'es sont les sous-objets de $X$ dans ${\rm Top}$.

\smallskip

En revanche, les quotients de $X$ dans ${\rm Top}_s$ sont diff\'erents: ce sont les classes d'isomorphie d'espaces topologiques s\'epar\'es $X'$ munis d'une application continue $p : X \to X'$ telle que $p(X)$ soit dense dans $X'$.

\section{Foncteurs fid\`eles et pleinement fid\`eles, \'equivalences de cat\'egories}\label{sec17}

\subsection{Les notions de foncteur fid\`ele ou pleinement fid\`ele}\label{subsec171}

\begin{defn}\label{defI71}
Un foncteur $F : {\mathcal C} \to {\mathcal D}$ entre deux cat\'egories est appel\'e
\begin{enumerate}
\item[$\bullet$] ``fid\`ele'' si les applications
\begin{eqnarray}
{\rm Hom}_{\mathcal C} (X_1 , X_2) &\longrightarrow & {\rm Hom}_{\mathcal D} (F(X_1) , F(X_2)) \, , \quad X_1 , X_2 \in {\rm Ob} ({\mathcal C}) \, , \nonumber \\
(u : X_1 \to X_2) &\longmapsto &(F(u) : F (X_1) \to F(X_2)) \nonumber
\end{eqnarray}
sont toutes injectives,
\item[$\bullet$] ``pleinement fid\`ele'' si ces applications sont toutes bijectives.
\end{enumerate}
\end{defn}

\begin{remarksqed}
\begin{listeisansmarge}
\item Tout foncteur pleinement fid\`ele est fid\`ele.

\medskip

\item Le foncteur de plongement dans une cat\'egorie ${\mathcal C}$ d'une sous-cat\'egorie ${\mathcal C}'$ de ${\mathcal C}$ est fid\`ele.

\smallskip

Il est pleinement fid\`ele si et seulement si ${\mathcal C}'$ est une sous-cat\'egorie pleine de ${\mathcal C}$.

\medskip

\item Le compos\'e ${\mathcal C}_1 \to {\mathcal C}_2 \to {\mathcal C}_3$ de deux foncteurs fid\`eles [resp. pleinement fid\`eles] ${\mathcal C}_1 \to {\mathcal C}_2$ et ${\mathcal C}_2 \to {\mathcal C}_3$ est fid\`ele [resp. pleinement fid\`ele]. 
\end{listeisansmarge}
\end{remarksqed}

\subsection{Exemples de foncteurs fid\`eles ou pleinement fid\`eles}\label{subsec172}

\smallskip

\noindent $\bullet$ {\bf Les foncteurs d'oubli:}

\smallskip

Si ${\mathcal C}$ est la cat\'egorie des ensembles munis d'un certain type de structure alg\'ebrique -- telle que celle de groupe, de mono{\"\i}de,  d'anneau, d'espace vectoriel sur un corps, de module sur un anneau --, alors le foncteur d'oubli de la structure
$$
{\mathcal C} \longrightarrow {\rm Ens}
$$
est fid\`ele.

\smallskip

De m\^eme, le foncteur d'oubli des structures topologiques
$$
{\rm Top} \longrightarrow {\rm Ens}
$$
est fid\`ele.

\medskip

\noindent $\bullet$ {\bf Les foncteurs de composition \`a gauche:}

\smallskip

Soit $F : {\mathcal D}_1 \to {\mathcal D}_2$ un foncteur fid\`ele.

\smallskip

Alors, pour toute cat\'egorie ${\mathcal C}$, le foncteur
\begin{eqnarray}
[{\mathcal C} , {\mathcal D}_1] &\longrightarrow &[{\mathcal C} , {\mathcal D}_2] \, , \nonumber \\
G &\longmapsto &F \circ G \nonumber
\end{eqnarray}
est fid\`ele.

\medskip

\noindent $\bullet$ {\bf Les foncteurs de composition \`a droite:}

\smallskip

Soit $F : {\mathcal C}_1 \to {\mathcal C}_2$ un foncteur tel que l'application
\begin{eqnarray}
{\rm Ob} ({\mathcal C}_1) &\longrightarrow &{\rm Ob} ({\mathcal C}_2) \, , \nonumber \\
X &\longmapsto &F(X) \nonumber
\end{eqnarray}
soit surjective ou, plus g\'en\'eralement, tel que tout objet $X_2$ de ${\mathcal C}_2$ soit isomorphe \`a l'image $F(X_1)$ d'un objet $X_1$ de ${\mathcal C}_1$.

\smallskip

Alors, pour toute cat\'egorie ${\mathcal D}$, le foncteur
\begin{eqnarray}
[{\mathcal C}_2 , {\mathcal D}] &\longrightarrow &[{\mathcal C}_1 , {\mathcal D}] \, , \nonumber \\
G &\longmapsto &G \circ F \nonumber
\end{eqnarray}
est fid\`ele.

\smallskip

En particulier, si cette condition est v\'erifi\'ee et que les cat\'egories ${\mathcal C}_1 , {\mathcal C}_2$ sont localement petites, le foncteur
$$
\widehat{\mathcal C}_2 \longrightarrow \widehat{\mathcal C}_1
$$
de composition avec $F : {\mathcal C}_1 \to {\mathcal C}_2$ est fid\`ele.

\smallskip

Dans le cas o\`u ${\mathcal C}_1$ et ${\mathcal C}_2$ sont deux mono{\"\i}des $M_1$ et $M_2$ reli\'es par un morphisme $\rho : M_1 \to M_2$, on obtient que le foncteur de substitution des actions de $M_2$ par des actions de $M_1$ induit par $\rho$
$$
\rho^* : BM_2 \longrightarrow BM_1
$$
est fid\`ele.

\medskip

\noindent $\bullet$ {\bf Les foncteurs d'extension des faisceaux d'un sous-espace \`a un espace:}

\smallskip

Si $f : X \to Y$ est une application continue entre deux espaces topologiques, les foncteurs d'image directe des pr\'efaisceaux ou des faisceaux
\begin{eqnarray}
f_* : \widehat{O(X)} &\longrightarrow &\widehat{O(Y)} \, , \nonumber \\
f_* : {\mathcal E}_X &\longrightarrow &{\mathcal E}_Y \nonumber
\end{eqnarray}
sont pleinement fid\`eles si tout ouvert de $X$ est l'image r\'eciproque $f^{-1} V$ d'un ouvert $V$ de $Y$.

\smallskip

C'est le cas en particulier si $X$ est un sous-espace de $Y$, c'est-\`a-dire un sous-ensemble muni de la topologie induite par celle de $Y$.

\subsection{Le lemme de Yoneda et la notion de foncteur repr\'esentable}\label{subsec173}

Le r\'esultat le plus important de la th\'eorie des cat\'egories est le lemme simple suivant:

\begin{lem}[\bf (Yoneda)]\label{lemI72}

Soit ${\mathcal C}$ une cat\'egorie localement petite.

\smallskip

Alors:

\begin{listeimarge}

\item Le foncteur de Yoneda
\begin{eqnarray}
y : {\mathcal C} &\longrightarrow &\widehat{\mathcal C} = [{\mathcal C}^{\rm op} , {\rm Ens}] \, , \nonumber \\
X &\longmapsto &y(X) = {\rm Hom} (\bullet , X) = \left\{ \begin{matrix}
\hfill Y &\longmapsto &{\rm Hom} (Y,X) \, , \hfill \\
(Y_1 \xrightarrow{ \ u \ } Y_2) &\longmapsto &\left[ \begin{matrix}
{\rm Hom} (Y_2 , X) &\to &{\rm Hom} (Y_1 , X) \, , \\
\hfill f &\mapsto &f \circ u \hfill
\end{matrix} \right]
\end{matrix} \right. \nonumber
\end{eqnarray}
est pleinement fid\`ele.

\medskip

\item Plus g\'en\'eralement, pour tout pr\'efaisceau $P : {\mathcal C}^{\rm op} \to {\rm Ens}$ et tout objet $X$ de ${\mathcal C}$, se donner un \'el\'ement de $P(X)$ \'equivaut \`a se donner un morphisme de $\widehat{\mathcal C}$
$$
{\rm Hom} (\bullet , X) \longrightarrow P \, .
$$
\end{listeimarge}
\end{lem}

\begin{remarks}
\begin{listeisansmarge}
\item En d'autres termes, se donner un morphisme $u : X_1 \to X_2$ entre deux objets $X_1 , X_2$ de ${\mathcal C}$ revient \`a se donner une famille d'applications
$$
u_Y : {\rm Hom} (Y,X_1) \longrightarrow {\rm Hom} (Y,X_2)
$$
index\'ees par les objets $Y$ de ${\mathcal C}$, telle que pour tout morphisme $v : Y \to Y'$ de ${\mathcal C}$, le carr\'e
$$
\xymatrix{
{\rm Hom} (Y' , X_1) \ar[d]_{\bullet \circ v} \ar[rr]^{u_{Y'}} &&{\rm Hom} (Y' , X_2) \ar[d]^{\bullet \circ v} \\
{\rm Hom} (Y,X_1) \ar[rr]^{u_Y} &&{\rm Hom} (Y,X_2)
}
$$
soit commutatif.


\item Rempla\c cant ${\mathcal C}$ par la cat\'egorie oppos\'ee ${\mathcal C}^{\rm op}$, on obtient que le foncteur
\begin{eqnarray}
y : {\mathcal C}^{\rm op} &\longrightarrow &\widehat{{\mathcal C}^{\rm op}} = [{\mathcal C} , {\rm Ens}] \, , \nonumber \\
X &\longmapsto &{\rm Hom} (X, \bullet) = \left\{ \begin{matrix}
\hfill Y &\longmapsto &{\rm Hom} (X,Y) \, , \hfill \\
(Y_1 \xrightarrow{ \ u \ } Y_2) &\longmapsto &\left[ \begin{matrix}
{\rm Hom} (X,Y_1) &\to &{\rm Hom} (X,Y_2) \, , \\
\hfill f &\mapsto &u \circ f \hfill
\end{matrix} \right.
\end{matrix} \right. \nonumber
\end{eqnarray}
est aussi pleinement fid\`ele.
\end{listeisansmarge}
\end{remarks}

\medskip

\noindent {\bf D\'emonstration du lemme de Yoneda:}

\begin{listeisansmarge}

\item est le cas particulier de (ii) o\`u $P$ est un pr\'efaisceau de la forme $y(X_2) = {\rm Hom} (\bullet , X_2)$.

\medskip

\item Consid\'erons un pr\'efaisceau $P$ sur ${\mathcal C}$, un objet $X$ de ${\mathcal C}$ et le pr\'efaisceau $y(X) = {\rm Hom} (\bullet , X)$ qu'il repr\'esente.

\smallskip

On dispose de l'application
$$
P(X) \longrightarrow {\rm Hom} (y(X) , P)
$$
qui associe \`a tout \'el\'ement $x \in P(X)$ la collection des applications
$$
\begin{matrix}
{\rm Hom} (Y,X) &\longrightarrow &P(Y) \, , \hfill \\
\hfill (Y \xrightarrow{u} X) &\longmapsto &P(u)(x)
\end{matrix}
$$
index\'ees par les objets $Y$ de ${\mathcal C}$.

\smallskip

Pour tout morphisme $Y_1 \xrightarrow{v} Y_2$ de ${\mathcal C}$, le carr\'e
$$
\xymatrix{
{\rm Hom} (Y_2, X) \ar[d]_{\bullet \circ v} \ar[r] &P(Y_2) \ar[d]^{P(v)} \\
{\rm Hom} (Y_1,X) \ar[r] &P(Y_1)
}
$$
est commutatif car, pour tout morphisme $u : Y_2 \to X$, on a
$$
P (u \circ v) = P(v) \circ P(u) \, .
$$

R\'eciproquement, consid\'erons un morphisme de $\widehat{\mathcal C}$
$$
\alpha : {\rm Hom} (\bullet , X) \longrightarrow P \, .
$$
Il consiste en une collection d'applications
$$
\alpha_Y : {\rm Hom} (Y,X) \longrightarrow P(Y) \, , \qquad Y \in {\rm Ob} ({\mathcal C}) \, ,
$$
rendant commutatif le carr\'e
$$
\xymatrix{
{\rm Hom} (Y_2, X) \ar[d]_{\bullet \circ v} \ar[rr]^{\alpha_{Y_2}} &&P(Y_2) \ar[d]^{P(v)} \\
{\rm Hom} (Y_1,X) \ar[rr]^{\alpha_{Y_1}} &&P(Y_1)
}
$$
associ\'e \`a tout morphisme $v : Y_1 \to Y_2$ de ${\mathcal C}$.

\smallskip

Prenant $Y = X$, on dispose de l'application
$$
\alpha_X : {\rm Hom} (X,X) \longrightarrow P(X)
$$
et on peut introduire l'\'el\'ement $x \in P(X)$ image par cette application de l'\'el\'ement identit\'e ${\rm id}_X \in {\rm Hom} (X,X)$.

\smallskip

Pour tout objet $Y$ de ${\mathcal C}$ et tout morphisme $u : Y \to X$, la commutativit\'e du carr\'e
$$
\xymatrix{
{\rm Hom} (X , X) \ar[d]_{\bullet \circ u} \ar[rr]^{\alpha_X} &&P(X) \ar[d]^{P(u)} \\
{\rm Hom} (Y,X) \ar[rr]^{\alpha_Y} &&P(Y)
}
$$
implique que 
$$
\alpha_Y (u) = P(u)(x) \, .
$$
Cela signifie que $\alpha$ est l'image de $x \in P(X)$ par l'application
$$
P(X) \longrightarrow {\rm Hom} (y(X),P) \, .
$$\hfill $\Box$

\end{listeisansmarge}

Le lemme de Yoneda conduit \`a poser la d\'efinition suivante:

\begin{defn}\label{defI73}
Soit ${\mathcal C}$ une cat\'egorie localement petite.

\smallskip

Un pr\'efaisceau sur ${\mathcal C}$ c'est-\`a-dire un foncteur contravariant
$$
F : {\mathcal C}^{\rm op} \longrightarrow {\rm Ens}
$$
est dit ``repr\'esentable'' s'il existe un objet $X$ de ${\mathcal C}$ et un isomorphisme
$$
y(X) \xrightarrow{ \ \sim \ } F \qquad \mbox{dans} \qquad \widehat{\mathcal C} = [{\mathcal C}^{\rm op} , {\rm Ens}] \, .
$$
On dit alors qu'un tel  objet $X$ repr\'esente le foncteur $F$ ou que $F$ est repr\'esent\'e par $X$.
\end{defn}


\begin{remarksqed}
\begin{listeisansmarge}
\item Il r\'esulte du lemme de Yoneda que si un tel objet $X$ de ${\mathcal C}$ existe, il est unique \`a unique isomorphisme pr\`es.

\smallskip

En effet, si $X_1$ et $X_2$ sont deux tels objets, l'isomorphisme compos\'e
$$
y(X_1) \xrightarrow{ \ \sim \ } F \xrightarrow{ \ \sim \ } y(X_2)
$$
provient d'apr\`es le lemme de Yoneda d'un unique isomorphisme de ${\mathcal C}$
$$
X_1 \xrightarrow{ \ \sim \ }  X_2 \, .
$$

\item S'il existe un tel isomorphisme
$$
u : y (X) = {\rm Hom} (\bullet , X) \xrightarrow{ \ \sim \ }  F \, ,
$$
on peut consid\'erer l'\'el\'ement
$$
u_X ({\rm id}_X) = s \in F(X) \, .
$$
Alors, pour tout objet $Y$ de ${\mathcal C}$, la bijection
$$
u_Y : {\rm Hom} (Y,X) \xrightarrow{ \ \sim \ }  F(Y)
$$
n'est autre que l'application
$$
(Y \xrightarrow{ \ f \ } X) \longmapsto F(f)(s) \, .
$$
Pour cette raison, cet \'el\'ement $s \in F(X)$ est appel\'e ``l'\'el\'ement universel'' (ou la ``section universelle'') du foncteur $F$.


\item Rempla\c cant la cat\'egorie ${\mathcal C}$ par son oppos\'ee ${\mathcal C}^{\rm op}$, on dit qu'un foncteur
$$
F : {\mathcal C} \longrightarrow {\rm Ens}
$$
est ``repr\'esentable'' s'il existe un objet $X$ de ${\mathcal C}$ et un isomorphisme
$$
{\rm Hom} (X,\bullet) \xrightarrow{ \ \sim \ }  F \qquad \mbox{dans} \qquad [{\mathcal C} , {\rm Ens}] \, .
$$
On dit alors qu'un tel objet $X$ repr\'esente le foncteur $F$ ou que $F$ est repr\'esent\'e par $X$.

\smallskip

Si $F$ est repr\'esent\'e par un objet $X$ de ${\mathcal C}$, celui-ci est unique \`a unique isomorphisme pr\`es.
\end{listeisansmarge}
\end{remarksqed}

\subsection{Premiers exemples de foncteurs repr\'esentables}\label{subsec174}

\smallskip

Le lemme de Yoneda et la d\'efinition \ref{defI73} qui le suit ouvrent la possibilit\'e de construire des objets de cat\'egories ${\mathcal C}$ localement petites en introduisant d'abord des foncteurs
$$
F : {\mathcal C}^{\rm op} \longrightarrow {\rm Ens}
$$
puis en montrant que ces foncteurs sont repr\'esentables.

\smallskip

Beaucoup des objets les plus int\'eressants des math\'ematiques, particuli\`erement en alg\`ebre et en g\'eom\'etrie, sont construits de cette fa\c con. Certains sont tr\`es sophistiqu\'es mais il existe de nombreux exemples simples et importants.

\medskip

\noindent $\bullet$ {\bf Le classificateur des sous-ensembles:}

\smallskip

Le foncteur contravariant
$$
\begin{matrix}
{\mathcal P} : &\hfill {\rm Ens}^{\rm op} &\longrightarrow &{\rm Ens} \, , \hfill \\
&\hfill I &\longmapsto &{\mathcal P} (I) = \mbox{ensemble des parties $J$ de $I$} \, , \\
&(I_1 \xrightarrow{ \ f \ } I_2) &\longmapsto &\left( (J \subseteq I_2) \longmapsto f^{-1} J \subseteq I_1 \right) \hfill
\end{matrix}
$$
est repr\'esentable par l'ensemble $\Omega = \{0,1\}$.

\smallskip

L'isomorphisme
$$
{\mathcal P} \xrightarrow{ \ \sim \ } {\rm Hom} (\bullet , \Omega)
$$
consiste \`a associer \`a toute partie $J$ d'un ensemble $I$ sa fonction caract\'eristique
\begin{eqnarray}
\un_J : I &\longrightarrow &\{0,1 \} \, , \nonumber \\
i &\longmapsto &\left\{ \begin{matrix}
1 &\mbox{si $i \in J$,} \\
0 &\mbox{sinon.} \hfill
\end{matrix} \right. \nonumber
\end{eqnarray}

\noindent $\bullet$ {\bf Le classifiant des ouverts:}

\smallskip

Le foncteur contravariant
$$
\begin{matrix}
O : &\hfill {\rm Top}^{\rm op} &\longrightarrow &{\rm Ens} \, , \hfill \\
&\hfill X &\longmapsto &O(X) = \mbox{ensemble des ouverts de $X$,} \\
&(X_1 \xrightarrow{ \ f \ } X_2) &\longmapsto &\left( O(X_2) \ni V \longmapsto f^{-1} V \in O(X_1)\right) \hfill
\end{matrix}
$$
est repr\'esentable par ``l'espace de Sierpinski ${\mathbb S}$'' constitu\'e de l'ensemble $\{0,1\}$ muni de la topologie dont le seul ouvert non trivial est $\{1\}$. Autrement dit, ${\mathbb S}$ est l'espace d'Alexandrov associ\'e \`a l'ensemble ordonn\'e $\{0,1\}$.

\smallskip

L'isomorphisme
$$
O \xrightarrow{ \ \sim \ } {\rm Hom} (\bullet , {\mathbb S})
$$
consiste \`a associer \`a tout ouvert $U$ d'un espace topologique $X$ l'unique application continue $\un_U : X \to {\mathbb S}$ telle que
$$
U = \un_U^{-1} \left(\{1\}\right) \, .
$$

\medskip

\noindent $\bullet$ {\bf Produits tensoriels, puissances sym\'etriques et puissances altern\'ees:}

\smallskip

On rappelle qu'\'etant donn\'es des modules $M_1 , \cdots , M_n$ et $N$ sur un anneau commutatif $A$, une application
$$
u : M_1 \times \cdots \times M_n \longrightarrow A
$$
est dite $n$-lin\'eaire si elle est $A$-lin\'eaire en chacune des $n$ variables. Autrement dit, $u$ est $n$-lin\'eaire si, pour tout indice $i$ et tous $m_1 \in M_1 , \cdots , m_{i-1} \in M_{i-1}$, $m_{i+1} \in M_{i+1} , \cdots , m_n \in M_n$, l'application induite
\begin{eqnarray}
M_i &\longrightarrow &N \, , \nonumber \\
m_i &\longmapsto &u(m_1 , \cdots , m_i , \cdots , m_n) \nonumber
\end{eqnarray}
est un morphisme de $A$-modules.

\smallskip

De plus, si $M_1 = \cdots = M_n = M$, une telle application $n$-lin\'eaire est dite sym\'etrique [resp. altern\'ee] si elle est invariante par permutation des $n$ coordonn\'ees $(m_1 , \cdots , m_n)$ [resp. si elle s'annule sur tous les \'el\'ements $(m_1 , \cdots , m_n)$ tels que $m_i = m_j$ pour au moins deux indices $i \ne j$].

\smallskip

On rappelle:

\begin{prop}\label{propI74}

Soit $A$ un anneau commutatif.

\begin{listeimarge}

\item Pour tous $A$-modules $M_1 , \ldots , M_n$, le foncteur
\begin{eqnarray}
{\rm Mod}_A &\longrightarrow &{\rm Ens} \nonumber \\
N &\longmapsto &\left\{ \begin{matrix}
u : &M_1 \times \cdots \times M_n \longrightarrow N \\
&\mbox{$n$-lin\'eaire} \hfill
\end{matrix} \right\} \nonumber
\end{eqnarray}
est repr\'esentable par un objet de ${\rm Mod}_A$ not\'e
$$
M_1 \otimes_A \cdots \otimes_A M_n
$$
et appel\'e le ``produit tensoriel'' de $M_1 , \cdots , M_n$ sur $A$

\medskip

\item Pour tout $A$-module $M$ et tout entier $n \geq 1$, le foncteur
\begin{eqnarray}
{\rm Mod}_A &\longrightarrow &{\rm Ens} \, , \nonumber \\
N &\longmapsto &\left\{ \begin{matrix}
u : &\underbrace{M \times \cdots \times M}_{\mbox{$n$ fois}} \longrightarrow N \qquad \mbox{$n$-lin\'eaire sym\'etrique}
\end{matrix} \right\} \nonumber
\end{eqnarray}
est repr\'esentable par un objet de ${\rm Mod}_A$ not\'e
$$
{\rm Sym}^n M
$$
et appel\'e la ``puissance sym\'etrique de degr\'e $n$'' de $M$ sur $A$.

\medskip

\item Pour tout $A$-module $M$ et tout entier $n \geq 1$, le foncteur
\begin{eqnarray}
{\rm Mod}_A &\longrightarrow &{\rm Ens} \nonumber \\
N &\longmapsto &\left\{ \begin{matrix}
u : &\underbrace{M \times \cdots \times M}_{\mbox{$n$ fois}} \longrightarrow N \qquad \mbox{$n$-lin\'eaire altern\'ee}
\end{matrix} \right\} \nonumber
\end{eqnarray}
est repr\'esentable par un objet de ${\rm Mod}_A$ not\'e
$$
\Lambda^n M
$$
et appel\'e la ``puissance altern\'ee de degr\'e $n$'' de $M$ sur $A$.
\end{listeimarge}
\end{prop}

\bigskip

\begin{remarks}
\begin{listeisansmarge}
\item Ainsi, $M_1 \otimes_A \cdots \otimes_A M_n$ est muni d'une application $n$-lin\'eaire canonique
\begin{eqnarray}
M_1 \times \cdots \times M_n &\longrightarrow &M_1 \otimes_A \cdots \otimes_A M_n \, , \nonumber \\
(m_1 , \cdots, m_n) &\longmapsto &m_1 \otimes \cdots \otimes m_n \nonumber
\end{eqnarray}
\`a travers laquelle toute application $n$-lin\'eaire
$$
u : M_1 \times \cdots \times M_n \longrightarrow N
$$
se factorise de mani\`ere unique en une application $A$-lin\'eaire
$$
M_1 \otimes_A \cdots \otimes_A M_n \longrightarrow N \, .
$$

\item De m\^eme, ${\rm Sym}^n M$ [resp. $\Lambda^n M$] est muni d'une application $n$-lin\'eaire sym\'etrique [resp. altern\'ee] canonique
\begin{eqnarray}
M \times \cdots \times M &\longrightarrow &{\rm Sym}^n M \, , \nonumber \\
(m_1,\cdots , m_n) &\longmapsto &m_1 \otimes \cdots \otimes m_n \nonumber \\
\mbox{[resp.} \qquad M \times \cdots \times M &\longrightarrow &\Lambda^n M \, , \nonumber \\
(m_1,\cdots , m_n) &\longmapsto &m_1 \wedge \cdots \wedge m_n \quad \mbox{ ]} \nonumber
\end{eqnarray}
\`a travers laquelle toute application $n$-lin\'eaire sym\'etrique [resp. altern\'ee]
$$
u : M \times \cdots \times M \longrightarrow N
$$
se factorise de mani\`ere unique en une application $A$-lin\'eaire
\begin{eqnarray}
{\rm Sym}^n M &\longrightarrow &N \nonumber \\
\mbox{[resp.} \qquad \Lambda^n M &\longrightarrow &N \quad \mbox{]}. \nonumber
\end{eqnarray}
\end{listeisansmarge}
\end{remarks}

\noindent $\bullet$ {\bf Les espaces affines:}

\smallskip

On rappelle que tout morphisme d'espaces annel\'es
$$
(X,{\mathcal O}_X) \longrightarrow (Y , {\mathcal O}_Y)
$$
consiste en une application continue $f : X \to Y$ compl\'et\'ee par un morphisme de faisceaux d'anneaux ${\mathcal O}_Y \to f_* {\mathcal O}_Y$.

\smallskip

En particulier, il induit un morphisme d'anneaux
$$
{\mathcal O}_Y (Y) \longrightarrow f_* {\mathcal O}_X (Y) = {\mathcal O}_X (X) \, .
$$

Cela conduit \`a poser:

\begin{defn}\label{defI75}

Soit ${\mathcal C}$ une sous-cat\'egorie de la cat\'egorie ${\rm Top}_{\rm an}$ des espaces annel\'es.

\smallskip

On dit que l'espace affine ${\mathbb A}^n$ de dimension $n$ est bien d\'efini dans ${\mathcal C}$ si le foncteur contravariant
\begin{eqnarray}
{\mathcal C}^{\rm op} &\longrightarrow &{\rm Ens} \, , \nonumber \\
(X,{\mathcal O}_X) &\longmapsto &{\mathcal O}_X (X)^n \nonumber
\end{eqnarray}
est repr\'esentable par un objet ${\mathbb A}^n$ de ${\mathcal C}$.
\end{defn}

\noindent {\bf Remarque:}

\smallskip

Dans le cas $n=1$, un objet ${\mathbb A}^1$ de ${\mathcal C}$ qui repr\'esente le foncteur $(X,{\mathcal O}_X) \longmapsto {\mathcal O}_X(X)$ est appel\'e la droite affine de ${\mathcal C}$. 
\hfill $\Box$

\bigskip

Ainsi, la cat\'egorie des vari\'et\'es diff\'erentielles de classe $C^k$ [resp. analytiques] contient un espace affine de chaque dimension $n$ qui est ${\mathbb R}^n$ [resp. ${\mathbb C}^n$] vu comme une vari\'et\'e diff\'erentielle [resp. analytique].

\smallskip

De m\^eme, la cat\'egorie ${\rm Sch}$ des sch\'emas contient un espace affine ${\mathbb A}^n$ de chaque dimension $n$ qui est le sch\'ema affine
$$
{\mathbb A}^n = {\rm Spec} ({\mathbb Z} [X_1 , \cdots , X_n]) \, .
$$

\pagebreak

\noindent $\bullet$ {\bf Le groupe multiplicatif et les groupes lin\'eaires:}

\smallskip

On pose encore:

\begin{defn}\label{defI76}

Soit ${\mathcal C}$ une sous-cat\'egorie de la cat\'egorie ${\rm Top}_{\rm an}$ des espaces annel\'es.

\smallskip

On dit que le groupe multiplicatif ${\mathbb G}_m$ [resp. le groupe lin\'eaire ${\rm GL}_r$ de rang $r \geq 1$] est bien d\'efini dans ${\mathcal C}$ si le foncteur contravariant
\begin{eqnarray}
{\mathcal C}^{\rm op} &\longrightarrow &{\rm Ens} \, , \nonumber \\
(X,{\mathcal O}_X) &\longmapsto &{\mathcal O}_X (X)^{\times} \quad \mbox{[resp.} \quad {\rm GL}_r ({\mathcal O}_X (X)) \ \mbox{]} \nonumber
\end{eqnarray}
est repr\'esentable par un objet ${\mathbb G}_m$ [resp. ${\rm GL}_r$] de ${\mathcal C}$.
\end{defn}

\begin{remarkqed}

On a ${\mathbb G}_m = {\rm GL}_1$ s'il existe. 
\end{remarkqed}

\bigskip

Ainsi, la cat\'egorie des vari\'et\'es diff\'erentielles de classe $C^k$ [resp. analytiques] contient un groupe multiplicatif ${\mathbb G}_m = {\mathbb R}^{\times} = {\mathbb R} - \{0\}$ [resp. ${\mathbb C}^{\times} = {\mathbb C} - \{0\}$] et un groupe lin\'eaire ${\rm GL}_r$ de chaque rang $r \geq 1$.

\smallskip

De m\^eme, la cat\'egorie ${\rm Sch}$ des sch\'emas contient un groupe multiplicatif qui est le sch\'ema affine
$$
{\mathbb G}_m = {\rm Spec} ({\mathbb Z} [X,Y] / (XY-1))
$$
et un groupe lin\'eaire ${\rm GL}_r$ de chaque rang $r \geq 1$ qui est le sch\'ema affine
$$
{\rm GL}_r = {\rm Spec} \left({\mathbb Z} [(X_{i,j})_{1 \leq i,j \leq r} , Y]/(\det (X_{i,j}) \cdot Y-1)\right).
$$

\medskip

\noindent $\bullet$ {\bf Les groupes orthogonaux et symplectiques:}

\smallskip

Soit $A$ un anneau commutatif.

\smallskip

On note $A$-alg la cat\'egorie des $A$-alg\`ebres: ses objets sont les anneaux commutatifs $B$ munis d'un morphisme $A \to B$ et ses morphismes $(A \to B_1) \to (A \to B_2)$ sont les morphismes d'anneaux $B_1 \to B_2$ qui rendent commutatif le triangle:
$$
\xymatrix{
B_1 \ar[rr] &&B_2 \\
&A \ar[ru] \ar[lu]
}
$$

Soient $n$ un entier [resp. un entier pair] et $\langle \bullet , \bullet \rangle$ une forme bilin\'eaire sym\'etrique [resp. altern\'ee]
$$
\langle \bullet , \bullet \rangle : A^n \times A^n \longrightarrow A
$$
qui est non d\'eg\'en\'er\'ee au sens que sa matrice
$$
P = (p_{i,j} = \langle e_i , e_j \rangle)_{1 \leq i,j \leq n}
$$
dans la base canonique de $A^n$ est inversible.

\smallskip

On observe que pour tout objet $B$ de $A$-alg, $\langle \bullet , \bullet \rangle$ s'\'etend en une forme bilin\'eaire sym\'etrique [resp. altern\'ee]
$$
\langle \bullet , \bullet \rangle : B^n \times B^n \longrightarrow B \, .
$$
On peut consid\'erer son groupe de sym\'etrie $O_n (B)$ [resp. ${\rm Sp}_n (B)$] d\'efini comme $\{g \in {\rm GL}_n (B) \mid \langle g(\bullet) , g(\bullet) \rangle = \langle \bullet , \bullet \rangle \}$.

\smallskip

On a:

\begin{prop}\label{propI77}

Dans les conditions ci-dessus, le foncteur
\begin{eqnarray}
\mbox{$A$-{\rm alg}} &\longrightarrow &{\rm Ens} \, , \nonumber \\
B &\longmapsto &O_n (B) \qquad \mbox{[resp.} \quad {\rm Sp}_n (B) \ \mbox{]} \nonumber
\end{eqnarray}
est repr\'esentable par un objet $A_{O_n}$ [resp. $A_{{\rm Sp}_n}$] qui est le quotient de la $A$-alg\`ebre
$$
A [ (X_{i,j})_{1 \leq i,j \leq n} ]
$$
par son id\'eal engendr\'e par les polyn\^omes
$$
\sum_{i',j'} X_{i,i'} \cdot p_{i',j'} \cdot X_{j',j} - p_{i,j} \, , \qquad 1 \leq i,j \leq n \, .
$$
\end{prop}

\begin{remarkqed}

Supposons que ${\mathcal C}$ est une sous-cat\'egorie de ${\rm Top}_{\rm an}$ dont tout objet $(X,{\mathcal O}_X)$ est muni d'un morphisme canonique
$$
A \longrightarrow {\mathcal O}_X(X)
$$
et dont tout morphisme $(X,{\mathcal O}_X) \to (Y , {\mathcal O}_Y)$ rend commutatif le triangle:
$$
\xymatrix{
{\mathcal O}_Y (Y) \ar[rr] &&{\mathcal O}_X(X) \\
&A \ar[ru] \ar[lu]
}
$$
On dispose alors du foncteur contravariant
\begin{eqnarray}
{\mathcal C}^{\rm op} &\longrightarrow &{\rm Ens} \, , \nonumber \\
(X,{\mathcal O}_X) &\longmapsto &O_n ({\mathcal O}_X(X)) \qquad \mbox{[resp.} \quad {\rm Sp}_n ({\mathcal O}_X(X)) \ \mbox{]}. \nonumber
\end{eqnarray}
Si ce foncteur est repr\'esentable dans ${\mathcal C}$, on dit que le groupe orthogonal $O_n$ [resp. le groupe symplectique ${\rm Sp}_n$] est bien d\'efini dans ${\mathcal C}$.

\smallskip

Il en est toujours ainsi dans le cas o\`u $A = {\mathbb R}$ [resp. $A = {\mathbb C}$] et ${\mathcal C}$ est la cat\'egorie des vari\'et\'es diff\'erentielles de classe $C^k$ [resp. analytiques].

\smallskip

De m\^eme, si $A = {\mathbb Z}$ et ${\mathcal C}$ est la cat\'egorie ${\rm Sch}$ des sch\'emas, ce foncteur est repr\'esentable par le sch\'ema affine
$$
O_n = {\rm Spec} (A_{O_n}) \qquad \mbox{[resp.} \quad {\rm Sp}_n = {\rm Spec} (A_{{\rm Sp}_n}) \ \mbox{]}.
$$
\end{remarkqed}

\subsection{La notion d'\'equivalence de cat\'egories}\label{subsec175}

\smallskip

Deux cat\'egories seront consid\'er\'ees comme \'equivalentes lorsqu'elles sont reli\'ees par une ``\'equivalence faible'' au sens suivant:

\begin{defn}\label{defI78}

Un foncteur $F : {\mathcal C} \to {\mathcal D}$ entre deux cat\'egories est appel\'e
\begin{enumerate}
\item[$\bullet$] une ``\'equivalence faible'' s'il est pleinement fid\`ele et si tout objet $Y$ de ${\mathcal D}$ est isomorphe \`a l'image $F(X)$ d'un objet $X$ de ${\mathcal C}$,
\item[$\bullet$] une ``\'equivalence'' s'il existe un foncteur en sens inverse
$$
G : {\mathcal D} \longrightarrow {\mathcal C}
$$
et deux isomorphismes de foncteurs
$$
\begin{matrix}
G \circ F &\xrightarrow{ \ \sim \ } &{\rm id}_{\mathcal C} &\mbox{dans} &[{\mathcal C} , {\mathcal C}] \, , \\
F \circ G &\xrightarrow{ \ \sim \ } &{\rm id}_{\mathcal D} &\mbox{dans} &[{\mathcal D} , {\mathcal D}] \, .
\end{matrix}
$$
\end{enumerate}
\end{defn}

\begin{remarksqed}
\begin{listeisansmarge}
\item Toute \'equivalence est une \'equivalence faible.

\smallskip

R\'eciproquement, une \'equivalence faible $F : {\mathcal C} \to {\mathcal D}$ est une \'equivalence si et seulement si il est possible de choisir uniform\'ement pour tout objet $Y$ de ${\mathcal D}$ un objet $X$ de ${\mathcal C}$ et un isomorphisme $F(X) \xrightarrow{ \ \sim \ } Y$.

\smallskip

C'est toujours possible si ${\mathcal C}$ et ${\mathcal D}$ sont des petites cat\'egories, d\`es lors qu'on a suppos\'e que les ensembles satisfont l'axiome du choix.

\medskip

\item Un foncteur $F : {\mathcal C} \to {\mathcal D}$ est pleinement fid\`ele si et seulement si il se factorise \`a travers une \'equivalence faible
$$
{\mathcal C} \longrightarrow {\mathcal D}'
$$
vers une sous-cat\'egorie pleine ${\mathcal D}'$ de ${\mathcal D}$. 
\end{listeisansmarge}
\end{remarksqed}

\medskip

La notion d'\'equivalence de cat\'egories est tellement importante qu'on est amen\'e \`a poser:

\begin{defn}\label{defI79}

Une propri\'et\'e des cat\'egories [resp. une construction sur les cat\'egories] est appel\'ee ``invariante'' si, pour toute \'equivalence faible
$$
F : {\mathcal C} \longrightarrow {\mathcal D} \, ,
$$
elle est v\'erifi\'ee par ${\mathcal D}$ d\`es lors qu'elle l'est par ${\mathcal C}$ [resp. elle donne le m\^eme r\'esultat sur ${\mathcal C}$ et sur ${\mathcal D}$].
\end{defn}

\bigskip

\begin{remarkqed}

Toutes les propri\'et\'es ou constructions cat\'egoriques consid\'er\'ees en pratique sont des propri\'et\'es invariantes. 

\end{remarkqed}

\medskip

Si une propri\'et\'e des cat\'egories n'est pas invariante, il convient de la remplacer par la propri\'et\'e invariante qu'elle engendre.

\smallskip

C'est ainsi que l'on pose:

\begin{defn}\label{defI710}
Une cat\'egorie ${\mathcal C}$ est dite ``essentiellement petite'' s'il existe une petite cat\'egorie ${\mathcal C}'$ et une \'equivalence faible
$$
{\mathcal C}' \longrightarrow {\mathcal C} \, .
$$
\end{defn}

\begin{remarksqed} 
\begin{listeisansmarge}
\item La propri\'et\'e pour une cat\'egorie d'\^etre essentiellement petite est invariante, tout comme la propri\'et\'e d'\^etre localement petite.

\medskip

\item Toute cat\'egorie essentiellement petite est localement petite. 
\end{listeisansmarge}
\end{remarksqed}

\newpage

\subsection{ Exemples d'\'equivalences de cat\'egories}\label{subsec176}

\medskip

\noindent $\bullet$ {\bf L'\'equivalence fondamentale de l'alg\`ebre lin\'eaire:}

\smallskip

Etant donn\'e un corps $K$, consid\'erons la cat\'egorie ${\rm Vect}_K^f$ des espaces vectoriels de dimension finie sur $K$ et des applications $K$-lin\'eaires entre de tels espaces.

\smallskip

Toute l'alg\`ebre lin\'eaire est fond\'ee sur le r\'esultat suivant:

\begin{prop}\label{propI711}

Soit ${\rm Mat}_K$ la cat\'egorie dont

\bigskip

$\left\{ \begin{matrix}
\bullet &\mbox{les objets sont les entiers $n \in {\mathbb N}$,} \hfill \\
{ \ } \\
\bullet &\mbox{les morphismes $n \to m$ sont les matrices \`a $m$ lignes et $n$ colonnes \`a coefficients dans le corps $K$,} \hfill \\
{ \ } \\
\bullet &\mbox{la loi de composition des morphismes est celle des matrices} \hfill \\
{ \ } \\
&\begin{matrix}
(N \ , \ M) &\xmapsto{ \ \ \ \ \ \ \ } &N \cdot M \, . \\
\Vert \ \ \ \ \Vert &&\Vert \\
(n_{i,j})_{1 \leq i \leq m \atop 1 \leq j \leq n} \ (m_{j,k})_{1 \leq j \leq n \atop 1 \leq k \leq r} && \ \ \ \left( \underset{1 \leq j \leq n}{\sum} \, n_{i,j} \, m_{j,k} \right)_{1 \leq i \leq m \atop 1 \leq k \leq r}
\end{matrix}
\end{matrix} \right.
$

\bigskip

Alors le foncteur

$$\left\{ \begin{matrix}
{\rm Mat}_K &\xrightarrow{ \ \ \ \ \ \ \ } &{\rm Vect}_K^f \, , \hfill \\
n &\xmapsto{ \ \ \ \ \ \ \ } &K^n \, , \hfill \\
M &\xmapsto{ \ \ \ \ \ \ \ } &\mbox{application lin\'eaire $K^n \to K^m$} \hfill \\
\Vert &&\mbox{d\'efinie par la matrice $M$} \hfill \\
\mbox{matrice} \ (m_{ij})_{1 \leq i \leq m \atop 1 \leq j \leq n} &&(k_j)_{1 \leq j \leq n} \mapsto \left( \underset{1 \leq j \leq n}{\sum} \, m_{i,j} \, k_j \right)_{1 \leq i \leq m}
\end{matrix} \right.
$$
est une \'equivalence faible. \hfill $\Box$
\end{prop}

\noindent $\bullet$ {\bf L'\'equivalence fondamentale de la th\'eorie de Galois:}

\smallskip

On rappelle qu'un corps commutatif $L$ contenant un corps $K$ est une extension finie de $K$, autrement dit est de dimension finie en tant qu'espace vectoriel sur $K$, si et seulement si $L$ est engendr\'e sur $K$ par un nombre fini de racines de polyn\^omes \`a coefficients dans $K$.

\smallskip

D'autre part, un polyn\^ome $P = a_n X^n + \cdots + a_1 X + a_0$ \`a coefficients dans $K$ est dit s\'eparable s'il est premier avec son polyn\^ome d\'eriv\'e $P' = n \, a_n X^{n-1} + \cdots + a_2 X + a_1$ ou, ce qui revient au m\^eme, si ses racines dans n'importe quelle extension finie de $K$ o\`u il est scind\'e (c'est-\`a-dire s'\'ecrit comme un produit de  polyn\^omes de degr\'e $1$) sont deux \`a deux distinctes.

\smallskip

Le r\'esultat central de la th\'eorie de Galois peut \^etre formul\'e de la mani\`ere suivante:

\begin{thm}\label{thmI712}

Soient $P$ un polyn\^ome s\'eparable \`a coefficients dans $K$ et $L$ une extension finie de $K$ dans laquelle $P$ est scind\'e et qui est engendr\'ee sur $K$ par les racines de $P$.

\smallskip

Soit $G = {\rm Aut}_K (L) = {\rm Gal} (L/K)$ le groupe de Galois de $L$ sur $K$, c'est-\`a-dire le groupe des automorphismes $L \xrightarrow{ \ \sim \ } L$ qui fixent les \'el\'ements de $K$.

\smallskip

Soit ${\rm Ext}_K^L$ la cat\'egorie dont

\medskip

$\left\{\begin{matrix}
\bullet &\mbox{les objets sont les sous-corps $L_1$ de $L$ qui contiennent $K$,} \hfill \\
\bullet &\mbox{les morphismes sont les morphismes de corps $L_1 \to L_2$ qui fixent les \'el\'ements de $K$.} \hfill
\end{matrix} \right.
$

\medskip

Alors le foncteur
$$\left\{ \begin{matrix}
({\mathcal E}\!xt_K^L)^{\rm op} &\xrightarrow{ \ \ \ \ \ \ \ } &BG \, , \hfill \\
L_1 &\xmapsto{ \ \ \ \ \ \ \ } &{\rm Hom} (L_1 , L) \, , \hfill \\
L_1 \xrightarrow{\iota} L_2 &\xmapsto{ \ \ \ \ \ \ \ } &\left[ \begin{matrix}
{\rm Hom} (L_2 , L) &\longrightarrow &{\rm Hom} (L_1 , L) \\
\hfill \iota_2 &\longmapsto &\iota_2 \circ \iota \hfill
\end{matrix} \right] \hfill \\
\end{matrix} \right.
$$
induit une \'equivalence de $({\mathcal E}\!xt_K^L)^{\rm op}$ vers la sous-cat\'egorie pleine de $BG$ constitu\'ee des ensembles finis $X$ sur lesquels le groupe fini $G$ agit transitivement.
\hfill $\Box$
\end{thm}

\medskip

\noindent $\bullet$ {\bf L'\'equivalence fondamentale de la th\'eorie des rev\^etements:}

\smallskip

Un rev\^etement d'un espace topologique $X$ est une paire $(X' , p)$ constitu\'ee d'un espace topologique $X'$ et d'une application continue
$$
p : X' \longrightarrow X
$$
telle qu'existent un recouvrement de $X$ par des ouverts $U_i$ et des hom\'eomorphismes (compatibles aux projections sur les $U_i$) entre les images r\'eciproques $p^{-1} (U_i)$ et des r\'eunions disjointes de copies de chaque $U_i$.

\smallskip

Un morphisme de rev\^etements
$$
(X'_1 , p_1) \longrightarrow (X'_2 , p_2)
$$
est une application continue $X'_1 \to X'_2$ qui rend commutatif le triangle:
$$
\xymatrix{
X'_1 \ar[rd] \ar[rr] &&X'_2 \ar[ld] \\
&X
}
$$

On note ${\mathcal C}\!ov_X$ la cat\'egorie dont les objets sont les rev\^etements de $X$ et dont les morphismes sont les morphismes de rev\^etements.

\smallskip

Rappelons qu'un espace topologique $X$ est dit localement contractile si tout voisinage ouvert $U$ d'un point $x$ de $X$ contient un voisinage $U'$ de $x$ qui admet une application continue
$$
h : U' \times [0,1] \longrightarrow U'
$$
telle que
$$
\begin{matrix}
&h(y,0) = x \, , &\forall \, y \in U' \, , \\
\mbox{et} &h(y,1) = y \, , &\forall \, y \in U' \, .
\end{matrix}
$$

Ceci \'etant rappel\'e, on peut formuler le r\'esultat central de la th\'eorie des rev\^etements:

\begin{thm}\label{thmI713}

Soit $X$ un espace topologique connexe et localement contractile.

\smallskip

Soient $x$ un \'el\'ement de $X$ et $G = \pi_X (x,x)$ son groupe fondamental de Poincar\'e.

\smallskip

Alors le foncteur
$$
\left\{ \begin{matrix}
\hfill {\mathcal R}ev_X &\xrightarrow{ \ \ \ \ \ \ \ } &BG \, , \hfill \\
\hfill (X' \xrightarrow{ \ p \ } X) &\xmapsto{ \ \ \ \ \ \ \ } &p^{-1} (x) \ \mbox{muni de l'action de} \ G = \pi_X (x,x) \, , \\
\left( (X'_1 \xrightarrow{ \ p_1 \ } X) \xrightarrow{f} (X'_2 \xrightarrow{ \ p_2 \ } X)\right) &\xmapsto{ \ \ \ \ \ \ \ } &\left( p_1^{-1} (x) \xrightarrow{ \ f \ } p_2^{-1} (x)\right), \hfill
\end{matrix} \right.
$$
est une \'equivalence de la cat\'egorie des rev\^etements de $X$ vers la cat\'egorie des actions du groupe fondamental $\pi_X (x,x)$.
\end{thm}

\bigskip

\noindent {\bf Remarque:}

\smallskip

Grothendieck a montr\'e que ce th\'eor\`eme de topologie et le th\'eor\`eme pr\'ec\'edent sur l'\'equivalence de Galois en alg\`ebre peuvent \^etre consid\'er\'es comme deux cas particuliers d'un unique th\'eor\`eme qui caract\'erise toutes les cat\'egories \'equivalentes \`a des cat\'egories $BG$ d'actions de groupe $G$. \hfill $\Box$

\bigskip

\noindent $\bullet$ {\bf Quelques exemples de cat\'egories essentiellement petites:}

\smallskip

La cat\'egorie des ensembles finis est essentiellement petite. En effet, sa sous-cat\'egorie pleine constitu\'ee des ensembles $\{0,1,\cdots,n\}$, $n \in {\mathbb N}$, lui est faiblement \'equivalente.

\smallskip

La cat\'egorie des ensembles munis d'un certain type de structure alg\'ebrique -- telle que celle de groupe, de mono{\"\i}de, de corps, d'anneau, d'espace vectoriel sur un corps fix\'e ou de module sur un anneau fix\'e -- et qui sont ``de pr\'esentation finie'', c'est-\`a-dire sont d\'efinis par une famille finie de g\'en\'erateurs et de relations, est essentiellement petite.

\smallskip

La cat\'egorie des vari\'et\'es diff\'erentielles de classe $C^k$ [resp. analytiques] qui sont ``d\'enombrables \`a l'infini'', c'est-\`a-dire sont r\'eunions de familles d\'enombrables d'ouverts isomorphes \`a des ouverts de ${\mathbb R}^n$ [resp. ${\mathbb C}^n$], est essentiellement petite.

\smallskip

Il en va encore de m\^eme de la cat\'egorie des sch\'emas de pr\'esentation finie sur un sch\'ema fix\'e au sens suivant:

\begin{defn}\label{defI714}
\begin{listeimarge}
\item Un anneau commutatif $B$ sur un anneau commutatif $A$ est dit ``de pr\'esentation finie sur $A$'' [resp. ``de type fini sur $A$''] si, avec son morphisme de structure $A \to B$, il est isomorphe au quotient d'une alg\`ebre de polyn\^omes $A[X_1 , \cdots , X_n]$ par un id\'eal engendr\'e par un nombre fini d'\'el\'ements [resp. par un id\'eal].

\smallskip

On dit alors que le sch\'ema affine ${\rm Spec} (B)$ est de pr\'esentation finie [resp. de type fini] sur le sch\'ema affine ${\rm Spec} (A)$.

\medskip

\item Un sch\'ema $X$ muni d'un morphisme $X \xrightarrow{ \ p \ } S$ vers un sch\'ema $S$ est dit ``de pr\'esentation finie [resp. de type fini] sur $S$'' si, pour tout ouvert affine $U \cong {\rm Spec} (A)$ de $S$, son image r\'eciproque $p^{-1} (U)$ est une r\'eunion finie d'ouverts affines $V_i \cong {\rm Spec} (B_i)$ qui sont de pr\'esentation finie [resp. de type fini] sur ${\rm Spec} (A)$.
\end{listeimarge}
\end{defn}

\begin{remarkqed}

On rappelle qu'un anneau commutatif $A$ est dit ``n\oe th\'erien'' si tout id\'eal de $A$ est engendr\'e par un nombre fini d'\'el\'ements.

\smallskip

On montre que si $A$ est n\oe th\'erien, il en est de m\^eme de chaque $A[X_1 , \cdots , X_n]$.

\smallskip

Cela entra{\^\i}ne que toute $A$-alg\`ebre $B$ de type fini est n\oe th\'erienne et de pr\'esentation finie sur $B$. 
\end{remarkqed}

\section{Foncteurs adjoints}\label{sec18}

\subsection{Les notions de foncteur adjoint \`a gauche ou \`a droite}\label{subsec181}

\begin{defn}\label{defI81}

Soient
$$
F : {\mathcal C} \longrightarrow {\mathcal D} \qquad \mbox{et} \qquad G : {\mathcal D} \longrightarrow {\mathcal C}
$$
deux foncteurs entre deux cat\'egories localement petites ${\mathcal C}$ et ${\mathcal D}$.

\smallskip

La paire $(F,G)$ est appel\'ee une ``paire adjointe'',

\smallskip

\qquad ou $F$ est appel\'e un ``adjoint \`a gauche'' de $G$,

\smallskip

\qquad ou $G$ est appel\'e un ``adjoint \`a droite'' de $F$

\smallskip

\noindent si les conditions \'equivalentes suivantes sont satisfaites:

\begin{enumerate}[label=(\Alph*)]

\item Les deux foncteurs
$$
\begin{matrix}
&{\mathcal C}^{\rm op} \times {\mathcal D} &\longrightarrow &{\rm Ens} \, , \hfill \\
&\hfill (X,Y) &\longmapsto &{\rm Hom} (F(X),Y) \\
\mbox{et} &\hfill (X,Y) &\longmapsto &{\rm Hom} (X,G(Y))
\end{matrix}
$$
sont isomorphes.

\medskip

\item Pour tout objet $Y$ de ${\mathcal D}$, $G(Y)$ repr\'esente le foncteur contravariant
\begin{eqnarray}
{\mathcal C}^{\rm op} &\longrightarrow &{\rm Ens} \, , \nonumber \\
X &\longmapsto &{\rm Hom} (F(X),Y) \nonumber
\end{eqnarray}
et pour tout morphisme $v : Y_1 \to Y_2$ de ${\mathcal D}$, le morphisme $G(v) : G(Y_1) \to G(Y_2)$ correspond au morphisme de foncteurs
$$
{\rm Hom} (F(\bullet) , Y_1) \longrightarrow {\rm Hom} (F(\bullet) , Y_2)
$$
dans $[{\mathcal C}^{\rm op} , {\rm Ens}] = \widehat{\mathcal C}$.

\medskip

\item Pour tout objet $X$ de ${\mathcal C}$, $F(X)$ repr\'esente le foncteur
\begin{eqnarray}
{\mathcal D} &\longrightarrow &{\rm Ens} \, , \nonumber \\
Y &\longmapsto &{\rm Hom} (X , G(Y)) \nonumber
\end{eqnarray}
et pour tout morphisme $u : X_1 \to X_2$ de ${\mathcal C}$, le morphisme $F(u) : F(X_1) \to F(X_2)$ correspond au morphisme de foncteurs
$$
{\rm Hom} (X_2 , G(\bullet)) \longrightarrow {\rm Hom} (X_1 , G(\bullet))
$$
dans $[{\mathcal D} , {\rm Ens}] = \widehat{{\mathcal D}^{\rm op}}$.
\end{enumerate}
\end{defn}

\medskip

\begin{remarksqed}
\begin{listeisansmarge}
\item Un foncteur $F : {\mathcal C} \to {\mathcal D}$ est adjoint \`a gauche de $G : {\mathcal D} \to {\mathcal C}$ si et seulement si $F : {\mathcal C}^{\rm op} \to {\mathcal D}^{\rm op}$ est adjoint \`a droite de $G : {\mathcal D}^{\rm op} \to {\mathcal C}^{\rm op}$.

\medskip

\item Un foncteur $F : {\mathcal C} \to {\mathcal D}$ a un adjoint \`a droite si et seulement si, pour tout objet $Y$ de ${\mathcal D}$, le foncteur
\begin{eqnarray}
{\mathcal C}^{\rm op} &\longrightarrow &{\rm Ens} \, , \nonumber \\
X &\longmapsto &{\rm Hom} (F(X),Y) \nonumber
\end{eqnarray}
est repr\'esentable et qu'il est possible de choisir uniform\'ement en $Y \in {\rm Ob} ({\mathcal D})$ un objet $G(Y)$ qui repr\'esente ce foncteur.

\medskip

\item De m\^eme, un foncteur $G : {\mathcal D} \to {\mathcal C}$ a un adjoint \`a gauche si et seulement si, pour tout objet $X$ de ${\mathcal C}$, le foncteur
\begin{eqnarray}
{\mathcal D} &\longrightarrow &{\rm Ens} \, , \nonumber \\
Y &\longmapsto &{\rm Hom} (X , G(Y)) \nonumber
\end{eqnarray}
est repr\'esentable et qu'il est possible de choisir uniform\'ement en $X \in {\rm Ob} ({\mathcal C})$ un objet $F(X)$ qui repr\'esente ce foncteur.

\medskip

\item Une paire de foncteurs $({\mathcal C} \xrightarrow{ \ F \ } {\mathcal D} , {\mathcal D} \xrightarrow{ \ G \ } {\mathcal C})$ est adjointe si et seulement si il existe un morphisme de foncteurs
$$
{\rm id}_{\mathcal C} \xrightarrow{ \ \alpha \ } G \circ F \qquad \mbox{[resp.} \quad F \circ G \xrightarrow{ \ \beta \ } {\rm id}_{\mathcal D} \ \mbox{]}
$$
tel que, pour tous objets $X$ de ${\mathcal C}$ et $Y$ de ${\mathcal D}$, l'application induite
$$
{\rm Hom} (F(X),Y) \xrightarrow{ \ G(\bullet) \ } {\rm Hom} (G \circ F (X) , G(Y)) \xrightarrow{ \ \bullet \circ \alpha_X \ } {\rm Hom} (X,G(Y))
$$
$$
\mbox{[resp.} \qquad {\rm Hom} (X,G(Y)) \xrightarrow{ F(\bullet) \ } {\rm Hom} (F(X) , F \circ G(Y)) \xrightarrow{ \ \beta_Y \circ \bullet \ } {\rm Hom} (F(X),Y) \ \mbox{]}
$$
est bijective.

\medskip

\item Si $({\mathcal C} \xrightarrow{ \ F \ } {\mathcal D} , {\mathcal D} \xrightarrow{ \ G \ } {\mathcal C})$ et $({\mathcal D} \xrightarrow{ \ F' \ } {\mathcal E} , {\mathcal E} \xrightarrow{ \ G' \ } {\mathcal D})$ sont deux paires adjointes, alors la paire $({\mathcal C} \xrightarrow{ \ F' \circ F \ } {\mathcal E} , {\mathcal E} \xrightarrow{ \ G \circ G' \ } {\mathcal C})$ est adjointe.

\medskip

\item Si deux foncteurs ${\mathcal C} \xrightarrow{ \ F \ } {\mathcal D}$ et ${\mathcal D} \xrightarrow{ \ G \ } {\mathcal C}$ d\'efinissent une \'equivalence de cat\'egories, alors les deux paires $(F,G)$ et $(G,F)$ sont adjointes.

\smallskip

En d'autres termes, la notion de paire de foncteurs adjoints est (beaucoup) plus faible que celle d'\'equivalence de cat\'egories.

\medskip

\item Il r\'esulte de (v) et (vi) que si $({\mathcal C} \xrightarrow{ \ F \ } {\mathcal D} , {\mathcal D} \xrightarrow{ \ G \ } {\mathcal C})$ est une paire adjointe et $F' : {\mathcal C} \to {\mathcal D}$, $G' : {\mathcal D} \to {\mathcal C}$ sont deux foncteurs isomorphes \`a $F$ et $G$, alors la paire $(F',G')$ est encore adjointe.

\smallskip

R\'eciproquement, si deux foncteurs sont adjoints \`a gauche [resp. \`a droite] d'un m\^eme foncteur, alors ils sont canoniquement isomorphes.

\medskip

\item Si deux paires de foncteurs $({\mathcal C} \xrightarrow{ \ F_1 \ } {\mathcal D} , {\mathcal D} \xrightarrow{ \ G_1 \ } {\mathcal C})$ et $({\mathcal C} \xrightarrow{ \ F_2 \ } {\mathcal D} , {\mathcal D} \xrightarrow{ \ G_2 \ } {\mathcal C})$ sont adjointes, alors les morphismes de foncteurs dans $[{\mathcal C}, {\mathcal D}]$
$$
F_1 \longrightarrow F_2
$$
sont en correspondance bijective avec les morphismes de foncteurs dans $[{\mathcal D} , {\mathcal C}]$
$$
G_2 \longrightarrow G_1 \, .
$$

\item Un foncteur
$$
F : {\mathcal C} \longrightarrow {\mathcal D}
$$
peut avoir \`a la fois un adjoint \`a gauche et un adjoint \`a droite.

\smallskip

S'il existe deux tels adjoints de $F$ \`a gauche et \`a droite, ils sont diff\'erents en g\'en\'eral (m\^eme s'il peut arriver dans certains cas qu'ils soient isomorphes). 
\end{listeisansmarge}
\end{remarksqed}

\medskip

Ainsi, lorsqu'un foncteur $F : {\mathcal C} \to {\mathcal D}$ a un adjoint \`a droite [resp. \`a gauche], celui-ci est uniquement d\'etermin\'e \`a unique isomorphisme pr\`es. C'est pourquoi on l'appelle couramment le foncteur adjoint \`a droite [resp. \`a gauche] de $F$.

\smallskip

Il arrive souvent que, partant d'un foncteur $F$ tr\`es simple on puisse lui associer un adjoint \`a droite ou \`a gauche beaucoup plus sophistiqu\'e et int\'eressant.

\smallskip

Cette observation est \`a l'origine d'un grand nombre de constructions math\'ematiques, particuli\`erement en alg\`ebre et en g\'eom\'etrie.

\medskip

Le premier r\'esultat suivant sur les paires de foncteurs adjoints est tr\`es simple mais utile:

\begin{lem}\label{lemI82}

Soit $({\mathcal C} \xrightarrow{ \ F \ } {\mathcal D} , {\mathcal D} \xrightarrow{ \ G \ } {\mathcal C})$ une paire adjointe.

\smallskip

Alors le foncteur $F$ [resp. $G$] est pleinement fid\`ele si et seulement si le morphisme canonique de foncteurs
$$
{\rm id}_{\mathcal C} \longrightarrow G \circ F \qquad \mbox{[resp.} \quad F \circ G \longrightarrow {\rm id}_{\mathcal D} \ \mbox{]}
$$
est un isomorphisme.
\end{lem}

\medskip

\begin{remark}

En particulier, si $F : {\mathcal C} \to {\mathcal D}$ a un adjoint \`a droite $G : {\mathcal D} \to {\mathcal C}$ et un adjoint \`a gauche $H : {\mathcal D} \to {\mathcal C}$, alors $G$ est pleinement fid\`ele si et seulement si $H$ est pleinement fid\`ele.
\end{remark}

\medskip

\begin{demo}

Remplacer les cat\'egories ${\mathcal C} , {\mathcal D}$ par leurs oppos\'ees ${\mathcal C}^{\rm op}$, ${\mathcal D}^{\rm op}$ r\'eduit le second cas du lemme au premier cas.

\smallskip

Pour tous objets $X_1 , X_2$ de ${\mathcal C}$, l'application
$$
{\rm Hom} (X_1 , X_2) \xrightarrow{ \ F(\bullet) \ } {\rm Hom} (F(X_1) , F(X_2))
$$
est une bijection si et seulement si sa compos\'ee
$$
{\rm Hom} (X_1 , X_2) \longrightarrow {\rm Hom} (X_1 , G \circ F(X_2))
$$
avec la bijection d'adjonction
$$
{\rm Hom} (F(X_1) , F(X_2)) \xrightarrow{ \ \sim \ } {\rm Hom} (X_1 , G \circ F(X_2))
$$
est une bijection.

\smallskip

D'apr\`es le lemme de Yoneda, cela \'equivaut \`a demander que, pour tout objet $X_2$ de ${\mathcal C}$, le morphisme canonique
$$
X_2 \longrightarrow G \circ F (X_2)
$$
soit un isomorphisme. 
\end{demo}

\subsection{Premiers exemples de foncteurs adjoints}\label{subsec182}

\medskip

\noindent $\bullet$ {\bf Le foncteur des puissances d'ensembles:}

\smallskip

Pour tout ensemble $A$, le foncteur
\begin{eqnarray}
{\rm Ens} &\longrightarrow &{\rm Ens} \, , \nonumber \\
C &\longmapsto &A \times C \nonumber
\end{eqnarray}
admet pour adjoint \`a droite le foncteur
\begin{eqnarray}
{\rm Ens} &\longrightarrow &{\rm Ens} \, , \nonumber \\
B &\longmapsto &B^A = \{\mbox{applications} \ A \to B \} \, . \nonumber
\end{eqnarray}

\medskip

\noindent $\bullet$ {\bf Le foncteur des parties d'ensembles:}

\smallskip

Le foncteur contravariant
\begin{eqnarray}
{\rm Ens}^{\rm op} &\longrightarrow &{\rm Ens} \, , \nonumber \\
I &\longmapsto &{\mathcal P} (I) = \{\mbox{parties de $I$}\} \nonumber
\end{eqnarray}
est adjoint \`a droite de lui-m\^eme consid\'er\'e comme un foncteur
$$
{\rm Ens} \longrightarrow {\rm Ens}^{\rm op} \, .
$$
En effet, rappelant que le foncteur ${\mathcal P}$ est repr\'esent\'e par l'objet $\Omega = \{0,1\}$, on remarque que pour tous ensembles $I$ et $J$, ${\rm Hom} (I,\Omega^J)$ et ${\rm Hom} (J,\Omega^I)$ s'identifient tous deux \`a l'ensemble ${\rm Hom} (I \times J , \Omega) = {\rm Hom} (J \times I , \Omega)$.

\medskip

\noindent $\bullet$ {\bf Les foncteurs d'engendrement des structures alg\'ebriques libres:}

\smallskip

Soit ${\mathcal C}$ la cat\'egorie des ensembles munis d'un certain type de structure alg\'ebrique, telle que celle de mono{\"\i}de, de groupe, d'anneau, d'anneau commutatif, de corps commutatif ou de module sur un anneau.

\smallskip

Alors le foncteur d'oubli de la structure alg\'ebrique
$$
G : {\mathcal C} \longrightarrow {\rm Ens}
$$
qui associe \`a tout objet de ${\mathcal C}$ son ensemble sous-jacent, admet un adjoint \`a gauche
$$
F : {\rm Ens} \longrightarrow {\mathcal C}
$$
qui associe \`a tout ensemble $I$ la structure alg\'ebrique ``libre'' engendr\'ee par les \'el\'ements de $I$.

\smallskip

Par exemple, si ${\mathcal C}$ est la cat\'egorie des modules sur un anneau $A$, $F(I)$ est le module libre $\underset{i \in I}{\bigoplus} \, A$.

\smallskip

Si ${\mathcal C}$ est la cat\'egorie des anneaux commutatifs, $F(I)$ est l'anneau de polyn\^omes ${\mathbb Z} [(X_i)_{i \in I}]$ en des variables $X_i$ index\'ees par les \'el\'ements de $I$.

\smallskip

Si ${\mathcal C}$ est la cat\'egorie des mono{\"\i}des, $F(I)$ est le mono{\"\i}de constitu\'e d'un \'el\'ement unit\'e $1$ et de suites formelles finies
$$
m_1 \cdots m_k
$$
compos\'ees d'\'el\'ements de $I$. La loi de multiplication de ces suites consiste \`a les juxtaposer.

\smallskip

Si ${\mathcal C}$ est la cat\'egorie des groupes, $F(I)$ est le groupe constitu\'e d'un \'el\'ement unit\'e $1$ et de suites formelles finies
$$
g_1 \cdots g_k
$$
compos\'ees d'\'el\'ements de $I$ ou de tels \'el\'ements affect\'es formellement de l'exposant $-1$. La loi de multiplication de ces suites consiste \`a les juxtaposer, et on consid\`ere deux suites comme \'equivalentes lorsque l'une peut \^etre d\'eduite de l'autre en enlevant une paire d'\'el\'ements juxtapos\'es de la forme $gg^{-1}$ ou $g^{-1} g$.

\medskip

\noindent $\bullet$ {\bf Le foncteur des points fixes et celui des orbites:}

\smallskip

Pour tout groupe $G$, le foncteur induit par le morphisme $G \to \{1\}$
\begin{eqnarray}
{\rm Ens} &\longrightarrow &BG \, , \nonumber \\
X &\longmapsto &\mbox{$X$ muni de l'action triviale de $G$,} \nonumber
\end{eqnarray}
admet pour adjoint \`a gauche le foncteur des points fixes
$$
(X , \mbox{action de $G$}) \longmapsto X_G = \{x \in X \mid g \cdot x = x \, , \quad \forall \, g \in G \}
$$
et pour adjoint \`a droite le foncteur des orbites
$$
(X , \mbox{action de $G$}) \longmapsto G \backslash X = \{\mbox{orbites de $X$ sous l'action de $G$}\} \, .
$$

Ceci fournit une premi\`ere s\'erie d'exemples de foncteurs qui admettent un adjoint \`a gauche et un adjoint \`a droite distincts.

\medskip

\noindent $\bullet$ {\bf Les extensions de Kan de pr\'efaisceaux:}

\smallskip

On verra plus loin que pour tout foncteur
$$
\rho : {\mathcal C} \longrightarrow {\mathcal D}
$$
entre deux cat\'egories essentiellement petites, le foncteur de composition avec $\rho$
$$
\xymatrix{
\rho^* : &\widehat{\mathcal D} \ar@{=}[d] \ar[r] &\widehat{\mathcal C} \ar@{=}[d] \\
&[{\mathcal D}^{\rm op} , {\rm Ens}] &[{\mathcal C}^{\rm op} , {\rm Ens}]
}
$$
admet un adjoint \`a droite 
$$
\rho_* : \widehat{\mathcal C} \longrightarrow \widehat{\mathcal D}
$$
et un adjoint \`a gauche
$$
\rho_! : \widehat{\mathcal C} \longrightarrow \widehat{\mathcal D} \, .
$$

\medskip

\noindent $\bullet$ {\bf Les foncteurs d'actions induites ou co-induites:}

\smallskip

Comme cas particulier d'extensions de Kan, pour tout morphisme de mono{\"\i}des
$$
\rho : M_1 \longrightarrow M_2 \, ,
$$
le foncteur de substitution des actions
$$
\rho^* : BM_2 \longrightarrow BM_1
$$
admet un adjoint \`a droite
$$
\rho_* : BM_1 \longrightarrow BM_2
$$
appel\'e le foncteur des actions induites de $M_1$ \`a $M_2$, et un adjoint \`a gauche
$$
\rho_! : BM_1 \longrightarrow BM_2
$$
appel\'e le foncteur des actions co-induites de $M_1$ \`a $M_2$.

\medskip

\noindent $\bullet$ {\bf Les foncteurs de faisceautisation:}

\smallskip

On verra plus loin que pour tout espace topologique $X$, le foncteur de plongement des faisceaux dans les pr\'efaisceaux
$$
j_* : {\mathcal E}_X \xhookrightarrow{ \ \ } \widehat{O(X)}
$$
admet un adjoint \`a gauche 
$$
j^* : \widehat{O(X)} \longrightarrow {\mathcal E}_X
$$
appel\'e le foncteur de faisceautisation des pr\'efaisceaux.

\smallskip

Comme le foncteur de plongement $j_* : {\mathcal E}_X \hookrightarrow \widehat{O(X)}$ est pleinement fid\`ele, il r\'esulte du lemme \ref{lemI82} que le morphisme canonique dans $[{\mathcal E}_X,{\mathcal E}_X]$
$$
j^* \circ j_* \longrightarrow {\rm id}
$$
est un isomorphisme. Autrement dit, tout faisceau est canoniquement isomorphe \`a sa faisceautisation.

\medskip

\noindent $\bullet$ {\bf Les foncteurs d'images r\'eciproques des faisceaux:}

\smallskip

On verra aussi que, pour toute application continue entre deux espaces topologiques
$$
f : X \longrightarrow Y \, ,
$$
le foncteur d'image directe
$$
f_* : {\mathcal E}_X \longrightarrow {\mathcal E}_Y
$$
admet un adjoint \`a gauche
$$
f^* : {\mathcal E}_Y \longrightarrow {\mathcal E}_X
$$
appel\'e le foncteur d'image r\'eciproque des faisceaux par $f$.

\smallskip

Dans le cas d'un point d'un espace topologique $X$
$$
x : \{\bullet\} \longrightarrow X \, ,
$$
l'adjoint \`a gauche de $x_* : {\rm Ens} = {\mathcal E}_{\{\bullet\}} \to {\mathcal E}_X$
$$
x^* : {\mathcal E}_X \longrightarrow {\rm Ens}
$$
est appel\'e le ``foncteur fibre'' en le point $x$.

\medskip

\noindent $\bullet$ {\bf Les foncteurs de prolongement des faisceaux \`a partir des ouverts:}

\smallskip

Si $U \xhookrightarrow{ \ i \ } X$ est un ouvert d'un espace topologique $X$, le foncteur d'image r\'eciproque
$$
i^* : {\mathcal E}_X \longrightarrow {\mathcal E}_U
$$
n'est autre que le foncteur de restriction des faisceaux de $X$ \`a $U$ d\'efini par composition avec l'inclusion $O(U) \subset O(X)$.

\smallskip

Il admet pour adjoint \`a droite le foncteur d'image directe 
$$
i_* : {\mathcal E}_U \longrightarrow {\mathcal E}_X
$$
et pour adjoint \`a gauche le foncteur
$$
\begin{matrix}
i_! : &{\mathcal E}_U &\longrightarrow &{\mathcal E}_X \, , \hfill \\
&\hfill F &\longmapsto &i_! \, F = \left[ \begin{matrix} 
O(X) &\longrightarrow &{\rm Ens}, \hfill \\
\hfill V &\longmapsto &\left\{ \begin{matrix} F(V) &\mbox{si} \ V \subset U \, , \\ \emptyset \hfill &\mbox{sinon} \hfill \end{matrix}\right.
\end{matrix} \right]
\end{matrix}
$$
que l'on appelle le foncteur de prolongement des faisceaux de $U$ \`a $X$.

\section{Limites et colimites}\label{sec19}

\subsection{La notion de diagramme dans une cat\'egorie}\label{subsec191}

\begin{defn}\label{defI91}
\begin{listeimarge}
\item On appelle ``carquois'' la donn\'ee $D$ de

\bigskip

$\left\{ \begin{matrix}
\bullet &\mbox{un ensemble ${\rm Ob} (D)$ dont les \'el\'ements sont appel\'es les objets de $D$,} \hfill \\
{ \ } \\
\bullet &\mbox{pour toute paire d'objets $d_1 , d_2$ de $D$, un ensemble} \hfill \\
{ \ } \\
&{\rm Hom}_D (d_1 , d_2) = {\rm Hom} (d_1 , d_2) \\
{ \ } \\
&\mbox{de fl\`eches $d_1 \to d_2$ de $d_1$ dans $d_2$ (ou d'origine $d_1$ et de but $d_2$).} \hfill
\end{matrix} \right.
$


\item Pour tout carquois $D$, un $D$-diagramme $X_{\bullet}$ dans une cat\'egorie ${\mathcal C}$ consiste en
\begin{enumerate}
\item[$\bullet$] une famille d'objets $X_d$ de ${\mathcal C}$ index\'es par les objets $d$ de $D$,
\item[$\bullet$] une famille de morphismes $X_u : X_{d_1} \to X_{d_2}$ de ${\mathcal C}$ index\'es par les fl\`eches $u : d_1 \to d_2$ de $D$.
\end{enumerate}

\medskip

\item Un morphisme $\alpha$ entre deux diagrammes $X_{\bullet}$ et $Y_{\bullet}$ dans une cat\'egorie ${\mathcal C}$ est une famille de morphismes.
$$
\alpha_d : X_d \longrightarrow Y_d
$$
index\'es par les objets $d$ de $D$, telle que le carr\'e associ\'e \`a toute fl\`eche $u : d_1 \to d_2$ de $D$
$$
\xymatrix{
X_{d_1} \ar[rr]^{\alpha_{d_1}} \ar[d]_{X_u} &&Y_{d_1} \ar[d]^{Y_u} \\
X_{d_2} \ar[rr]^{\alpha_{d_2}} &&Y_{d_2}
}
$$
soit commutatif.

\end{listeimarge}
\end{defn}

\begin{remarks}
\begin{listeisansmarge}
\item Tout carquois $D$ a un ``carquois oppos\'e'' $D^{\rm op}$ qui a les m\^emes objets et dont les fl\`eches sont donn\'ees par la r\`egle
$$
{\rm Hom}_{D^{\rm op}} (X,Y) = {\rm Hom}_D (Y,X) \, .
$$

\medskip

\item Pour tout carquois $D$, les $D$-diagrammes dans une cat\'egorie ${\mathcal C}$ forment une cat\'egorie
$$
D\mbox{-diag} \, ({\mathcal C}) \, .
$$

Si ${\mathcal C}$ est localement petite [resp. petite], alors la cat\'egorie $D$-diag $({\mathcal C})$ est localement petite [resp. petite].

\medskip

\item Pour tout carquois $D$, tout foncteur entre cat\'egories
$$
F : {\mathcal C}_1 \longrightarrow {\mathcal C}_2
$$
induit un foncteur
$$
\begin{matrix}
D\mbox{-diag} \, ({\mathcal C}_1) &\longrightarrow &D\mbox{-diag} \, ({\mathcal C}_2) \, , \hfill \\
\hfill X_{\bullet} &\longmapsto &F(X_{\bullet}) = \left\{ \begin{matrix}
\hfill d &\longmapsto &F(X_d) \, , \hfill \\
(d_1 \xrightarrow{ \, u \, } d_2) &\longmapsto &(F(X_{d_1}) \xrightarrow{ \, F(u) \, } F(X_{d_2})) \, .
\end{matrix} \right.
\end{matrix}
$$

\item On appelle morphisme entre deux carquois $D$ et $D'$
$$
\alpha : D \longrightarrow D'
$$
une famille d'applications
$$
\alpha : {\rm Ob} (D) \longrightarrow {\rm Ob} (D')
$$
et
$$
\alpha : {\rm Hom}_D (d_1 , d_2) \longrightarrow {\rm Hom}_{D'} (\alpha (d_1) , \alpha (d_2))
$$
index\'ees par les objets $d_1 , d_2$ de $D$.

\smallskip

On d\'efinit le compos\'e de deux morphismes de carquois
$$
D \xrightarrow{ \ \alpha \ } D' \xrightarrow{ \ \beta \ } D''
$$
comme la famille des applications compos\'ees
$$
\beta \circ \alpha : {\rm Ob} (D) \longrightarrow {\rm Ob} (D'')
$$
et
$$
\beta \circ \alpha : {\rm Hom}_D (d_1 , d_2) \longrightarrow {\rm Hom}_{D'} (\alpha (d_1) , \alpha (d_2)) \longrightarrow {\rm Hom}_{D''} (\beta \circ \alpha (d_1) , \beta \circ \alpha (d_2)) \, .
$$

Les carquois et leurs morphismes forment une cat\'egorie localement petite que l'on notera Carq.

\smallskip

Le passage aux carquois oppos\'es d\'efinit un foncteur ${\rm Carq} \to {\rm Carq}$.

\medskip

\item Pour toute cat\'egorie ${\mathcal C}$, tout morphisme de carquois
$$
\alpha : D \longrightarrow D'
$$
d\'efinit un foncteur entre cat\'egories de diagrammes
$$
\begin{matrix}
D'\mbox{-diag} \, ({\mathcal C}) &\longrightarrow &D\mbox{-diag} \, ({\mathcal C}) \, , \hfill \\
\hfill X_{\bullet} &\longmapsto &X_{\bullet} \circ \alpha = \left\{ \begin{matrix}
\hfill {\rm Ob} (D) \ni d &\longmapsto &X_{\alpha(d)} \, , \hfill \\
(d_1 \xrightarrow{ \, u \, } d_2) &\longmapsto &(X_{\alpha(d_1)} \xrightarrow{ \, X_{\alpha(u)} \, } X_{\alpha (d_2)}) \, .
\end{matrix} \right.
\end{matrix}
$$

\item Toute petite cat\'egorie ${\mathcal D}$ peut \^etre consid\'er\'ee comme un carquois en oubliant sa loi de composition des morphismes.

\smallskip

Pour toute cat\'egorie ${\mathcal C}$, la cat\'egorie des foncteurs
$$
[{\mathcal D} , {\mathcal C}]
$$
est une sous-cat\'egorie pleine de la cat\'egorie
$$
{\mathcal D}\mbox{-diag} \, ({\mathcal C})
$$
des ${\mathcal D}$-diagrammes de ${\mathcal C}$.

\medskip

\item L'oubli de la loi de composition des morphismes des petites cat\'egories d\'efinit un foncteur
$$
{\rm Cat} \longrightarrow {\rm Carq}
$$
de la cat\'egorie des petites cat\'egories dans celle des carquois.

\smallskip

Ce foncteur admet un adjoint \`a gauche
$$
{\rm Carq} \longrightarrow {\rm Cat}
$$
qui associe \`a tout diagramme $D$ la petite cat\'egorie ``libre'' ${\mathcal C}_D$ engendr\'ee par $D$ de la mani\`ere suivante:

\bigskip

$
\left\{\begin{matrix}
\bullet &\mbox{les objets de ${\mathcal C}_D$ sont les objets de $D$,} \hfill \\
{ \ } \\
\bullet &\mbox{les morphismes de ${\mathcal C}_D$ sont des fl\`eches identit\'es ${\rm id}_d : d \to d$} \hfill \\
&\mbox{et, pour toute paire d'objets $d,d'$ de $D$, les suites finies} \hfill \\
{ \ } \\
&d = d_0 \xrightarrow{ \ u_1 \ } d_1 \xrightarrow{ \ u_2 \ } \cdots \xrightarrow{ \ u_{k-1} \ } d_{k-1} \xrightarrow{ \ u_k \ } d_k = d' \\
{ \ } \\
&\mbox{de fl\`eches $u_1 , u_2 , \cdots , u_k$ de $D$ telles que le but de chaque $u_i$ soit l'origine de $u_{i+1}$,} \hfill \\
{ \ } \\
\bullet &\mbox{la loi de composition des morphismes est la juxtaposition de ces suites.} \hfill
\end{matrix}
\right.
$
\end{listeisansmarge}
\end{remarks}

\subsection{Exemples de carquois}\label{subsec192}

\medskip

\noindent $\bullet$ {\bf Les carquois sans fl\`eches:}

\smallskip

Tout ensemble $I$ (fini ou infini) d\'efinit un carquois sans fl\`eche dont l'ensemble des objets est $I$.

\smallskip

En particulier, m\^eme l'ensemble vide $\emptyset$ d\'efinit un carquois.

\smallskip

Pour toute cat\'egorie ${\mathcal C}$, $\emptyset$-diag $({\mathcal C})$ est la cat\'egorie \`a un objet et un morphisme.

\medskip

\noindent $\bullet$ {\bf La paire de fl\`eches parall\`eles:}

\smallskip

Le carquois de la paire de fl\`eches parall\`eles est constitu\'e de deux objets $O$ et $B$ et de deux fl\`eches distinctes d'origine $O$ et de but $B$:
$$
\bullet \rightrightarrows \bullet
$$

L'oppos\'e de ce carquois lui est canoniquement isomorphe.

\medskip

\noindent $\bullet$ {\bf Le produit fibr\'e:}

\smallskip

Le carquois de produit fibr\'e (ou de changement de base) est constitu\'e de trois objets $X,S,S'$ et de deux fl\`eches de $X$ vers $S$ et de $S'$ vers $S$:
$$
\xymatrix{&\bullet \ar[d] \\
\bullet \ar[r] &\bullet
}
$$

\noindent $\bullet$ {\bf La somme amalgam\'ee:}

\smallskip

Le carquois de somme amalgam\'ee est l'oppos\'e du pr\'ec\'edent. Il est constitu\'e de trois objets $X,Y,S$ et de deux fl\`eches de $S$ vers $X$ et de $S$ vers $Y$:
$$
\xymatrix{
\bullet \ar[d] \ar[r] &\bullet \\
\bullet
}
$$

\noindent $\bullet$ {\bf Le carr\'e et le triangle:}

\smallskip

Le carquois carr\'e [resp. triangulaire] a d\'ej\`a \'et\'e envisag\'e, tellement il est omnipr\'esent dans l'\'etude des cat\'egories. Il a la forme
$$
\xymatrix{
\bullet \ar[d] \ar[r] &\bullet \ar[d] \\
\bullet \ar[r] &\bullet
} \qquad\qquad \xymatrix{
\bullet \ar[rr] \ar[rd]_{\mbox{[resp. \quad}} &&\bullet \ar[ld]^{\qquad\mbox{].}} \\
&\bullet
}
$$

L'oppos\'e du carquois carr\'e [resp. triangulaire] lui est canoniquement isomorphe.

\medskip

\noindent $\bullet$ {\bf La cha{\^\i}ne:}

\smallskip

Le carquois de cha{\^\i}ne est constitu\'e d'objets index\'es par les entiers $n \in {\mathbb Z}$ et a une fl\`eche $n \to n+1$ pour chaque entier $n$:
$$
\cdots \longrightarrow \underset{-1}{\bullet} \longrightarrow \underset{0}{\bullet} \longrightarrow \underset{1}{\bullet} \longrightarrow \underset{2}{\bullet} \longrightarrow \cdots
$$

Contrairement aux carquois pr\'ec\'edemment donn\'es en exemples, le carquois de cha{\^\i}ne a des automorphismes qui sont les ``translations'' des indices des objets par les entiers $n \in {\mathbb Z}$.

\smallskip

L'oppos\'e du carquois de cha{\^\i}ne lui est isomorphe.

\smallskip

Ce carquois joue lui aussi un grand r\^ole en math\'ematiques, \`a travers les notions de ``complexes de cha{\^\i}nes'' ou de ``complexes de cocha{\^\i}nes'' qui sont les objets des th\'eories ``homologiques'' ou ``cohomologiques'':

\smallskip

Si ${\mathcal A}$ est une cat\'egorie ``lin\'eaire'', telle que la cat\'egorie ${\rm Mod}_A$ des modules sur un anneau $A$, un complexe de cha{\^\i}nes [resp. de cocha{\^\i}nes] de ${\mathcal A}$ est un diagramme
$$
\cdots \longrightarrow X_{-1} \xrightarrow{ \ u_0 \ } X_0 \xrightarrow{ \ u_1 \ } X_1 \xrightarrow{ \ u_2 \ } X_2 \longrightarrow \cdots
$$
$$
\mbox{[resp.} \qquad \cdots \longrightarrow X_1 \xrightarrow{ \ u_1 \ } X_0 \xrightarrow{ \ u_{0} \ } X_{-1} \xrightarrow{ \ u_{-1} \ } X_{-2} \longrightarrow \cdots \ \mbox{]}
$$
tel que
$$
u_{n+1} \circ u_n = 0 \, , \qquad \forall \, n \in {\mathbb Z}
$$
$$
\mbox{[resp.} \qquad u_{n-1} \circ u_n = 0 \, , \qquad \forall \, n \in {\mathbb Z} \ \mbox{]}.
$$
Si $D$ est le carquois de cha{\^\i}ne, les complexes de cha{\^\i}nes [resp. de cocha{\^\i}nes] de ${\mathcal A}$ forment une sous-cat\'egorie pleine de la cat\'egorie $D$-diag $({\mathcal A})$ [resp. $D^{\rm op}$-diag $({\mathcal A})$].

\subsection{Les notions de limite et de colimite}\label{subsec193}

\smallskip

Pour toute petite cat\'egorie ${\mathcal C}$ et tout carquois $D$, on dispose du ``foncteur diagonal''
$$
\begin{matrix}
\Delta_D : &\hfill {\mathcal C} &\longrightarrow &D\mbox{-diag} \, ({\mathcal C}) \, , \hfill \\
&\hfill X &\longmapsto &\Delta_D (X) = \left[ \begin{matrix}
\hfill d &\longmapsto &X \hfill \\
(d_1 \xrightarrow{ \, u \, } d_2) &\longmapsto &{\rm id}_X
\end{matrix} \right] , \\
&(X_1 \xrightarrow{ \ f \ } X_2) &\longmapsto &\Delta_D (f) = [d \longmapsto (X_1 \xrightarrow{ \, f \, } X_2)]. \hfill
\end{matrix}
$$

Cela permet de poser:

\begin{defn}\label{defI92}

Soient ${\mathcal C}$ une cat\'egorie localement petite et $D$ un carquois.

\begin{listeimarge}

\item Un $D$-diagramme $X_{\bullet}$ de ${\mathcal C}$ a une ``limite'' (ou ``limite projective'') dans ${\mathcal C}$ si le foncteur
$$
\begin{matrix}
{\mathcal C}^{\rm op} &\longrightarrow &{\rm Ens} \, , \hfill \\
\hfill X &\longmapsto &{\rm Hom}_{\mbox{\footnotesize $D$-{\rm diag} $({\mathcal C})$}} (\Delta_D (X) , X_{\bullet}) \hfill \\
&&= \left\{ (X \xrightarrow{ \ \alpha_d \ } X_d)_{d \in {\rm Ob} (D)} \Biggl\vert \begin{matrix}
\hfill X_u \circ \alpha_{d_1} &=& \alpha_{d_2} \, , \hfill \\
\forall \, (d_1 \xrightarrow{ \, u \, } d_2) &= &\mbox{fl\`eche de $D$}
\end{matrix}  \right\}
\end{matrix}
$$
est repr\'esentable par un objet de $D$, alors not\'e
$$
\varprojlim_D X_{\bullet} \, .
$$

\smallskip

\item Un $D$-diagramme $X_{\bullet}$ de ${\mathcal C}$ a une ``colimite'' (ou ``limite inductive'') si le foncteur
$$
\begin{matrix}
{\mathcal C} &\longrightarrow &{\rm Ens} \, , \hfill \\
\hfill X &\longmapsto &{\rm Hom}_{\mbox{\footnotesize $D$-{\rm diag} $({\mathcal C})$}} (X_{\bullet} , \Delta_D (X)) \hfill \\
&&= \left\{ (X_d \xrightarrow{ \ \alpha_d \ } X)_{d \in {\rm Ob} (D)} \Biggl\vert \begin{matrix}
\hfill \alpha_{d_2} \circ X_u &=& \alpha_{d_1} \, , \hfill \\
\forall \, (d_1 \xrightarrow{ \, u \, } d_2) &= &\mbox{fl\`eche de $D$}
\end{matrix}  \right\}
\end{matrix}
$$
est repr\'esentable par un objet de ${\mathcal C}$, alors not\'e
$$
\varinjlim_D X_{\bullet} \, .
$$
\end{listeimarge}
\end{defn}

\begin{remarksqed}
\begin{listeisansmarge}
\item Si un $D$-diagramme $X_{\bullet}$ de ${\mathcal C}$ admet une limite [resp. une colimite] dans ${\mathcal C}$, celle-ci est uniquement d\'etermin\'ee \`a unique isomorphisme pr\`es.

\medskip

\item Si $X_{\bullet}$ et $Y_{\bullet}$ sont deux $D$-diagrammes de ${\mathcal C}$ qui admettent des limites [resp. des colimites] dans ${\mathcal C}$, tout morphisme de $D$-diag $({\mathcal C})$
$$
X_{\bullet} \longrightarrow Y_{\bullet}
$$
induit un morphisme de ${\mathcal C}$
$$
\varprojlim_D X_{\bullet} \longrightarrow \varprojlim_D Y_{\bullet}
$$
$$
\mbox{[resp.} \qquad \varinjlim_D X_{\bullet} \longrightarrow \varinjlim_D Y_{\bullet} \ \mbox{]}.
$$

\end{listeisansmarge}
\end{remarksqed}

\subsection{Exemples de limites et de colimites}\label{subsec194}

\noindent $\bullet$ {\bf Les produits et les sommes:}

\smallskip

Dans le cas d'un carquois $D$ sans fl\`eche, qui donc est un ensemble $I$ d'objets, un $D$-diagramme dans une cat\'egorie ${\mathcal C}$ est une famille $(X_i)_{i \in I}$ d'objets.

\smallskip

Si elle existe, sa limite [resp. sa colimite] est not\'ee
$$
\prod_{i \in I} X_i \qquad \mbox{[resp.} \quad \coprod_{i \in I} X_i \ \mbox{]}
$$
et appel\'ee son ``produit'' [resp. sa ``somme''].

\smallskip

Elle est munie d'une famille de morphismes canoniques
$$
\prod_{i \in I} X_i \longrightarrow X_j \qquad \mbox{[resp.} \quad X_j \longrightarrow \prod_{i \in I} X_i \ \mbox{]} \, , \ j \in I \, ,
$$
qui induit pour tout objet $X$ de ${\mathcal C}$ une bijection
$$
{\rm Hom} \left( X , \prod_{i \in I} X_i \right) \xrightarrow{ \ \sim \ } \prod_{i \in I} {\rm Hom} (X,X_i)
$$
$$
\mbox{[resp.} \qquad {\rm Hom} \left( \coprod_{i \in I} X_i , X \right) \xrightarrow{ \ \sim \ } \prod_{i \in I} {\rm Hom} (X_i,X) \ \mbox{]}.
$$

Autrement dit, se donner une famille de morphismes
$$
(X \longrightarrow X_i)_{i \in I} \qquad \mbox{[resp.} \quad (X_i \longrightarrow X)_{i \in I} \ \mbox{]}
$$
\'equivaut \`a se donner un morphisme
$$
X \longrightarrow \prod_{i \in I} X_i \qquad \mbox{[resp.} \quad \coprod_{i \in I} X_i \longrightarrow X \ \mbox{]}.
$$

On note que les produits et les sommes sont \'echang\'es si l'on passe de ${\mathcal C}$ \`a son oppos\'ee ${\mathcal C}^{\rm op}$. En d'autres termes, les notions de produit et de somme sont duales l'une de l'autre.

\smallskip

Enfin, si $I$ est l'ensemble fini $\{1,\cdots , n \}$, on note
$$
\prod_{1 \leq i \leq n} X_i = X_1 \times \ldots \times X_n \qquad \mbox{[resp.} \quad \coprod_{1 \leq i \leq n} X_i = X_1 \amalg \cdots \amalg X_n \ \mbox{]}.
$$

\medskip

\noindent $\bullet$ {\bf Objet terminal et objet initial:}

\smallskip

Si $D=I$ est le carquois vide $\emptyset$, la limite [resp. la colimite] est not\'ee $1$ [resp. $0$] quand elle existe.

\smallskip

Elle est alors appel\'ee l'objet ``terminal'' [resp. l'objet ``initial''] de ${\mathcal C}$ car elle est caract\'eris\'ee par la propri\'et\'e que tout objet $X$ de ${\mathcal C}$ a un unique morphisme
$$
X \longrightarrow 1 \qquad \mbox{[resp.} \quad 0 \longrightarrow X \ \mbox{]}.
$$

Les notions d'objet terminal et d'objet initial sont duales l'une de l'autre.

\medskip

\noindent $\bullet$ {\bf Egalisateurs et co\'egalisateurs:}

\smallskip

Si le carquois $D$ est la paire de fl\`eches parall\`eles, un $D$-diagramme d'une cat\'egorie ${\mathcal C}$ a la forme
$$
X \overset{f}{\underset{g}{\rightrightarrows}} Y \, .
$$

Quand elle existe, sa limite [resp. sa colimite] est appel\'ee l'\'egalisateur [resp. le co\'egalisateur] de $f$ et $g$, et elle est not\'ee
$$
{\rm eg} \, (X \overset{f}{\underset{g}{\rightrightarrows}} Y) = {\rm eg} \, (f,g) \qquad \mbox{[resp.} \quad {\rm coeg} \, (X \overset{f}{\underset{g}{\rightrightarrows}} Y) = {\rm coeg} \, (f,g) \ \mbox{]}.
$$

Elle est caract\'eris\'ee par la propri\'et\'e que se donner un morphisme de ${\mathcal C}$
$$
Z \longrightarrow {\rm eg} \, (f,g) \qquad \mbox{[resp.} \quad {\rm coeg} \, (f,g) \longrightarrow Z \ \mbox{]}
$$
\'equivaut \`a se donner un morphisme
$$
h : Z \longrightarrow X \qquad \mbox{[resp.} \quad h : Y \longrightarrow Z \ \mbox{]}
$$
tel que
$$
f \circ h = g \circ h \qquad \mbox{[resp.} \quad h \circ f = h \circ g \ \mbox{]} .
$$

Les notions d'\'egalisateur et de co\'egalisateur sont duales l'une de l'autre.

\medskip

\noindent $\bullet$ {\bf Produits fibr\'es et sommes amalgam\'ees:}

\smallskip

Si $D$ est le carquois des produits fibr\'es [resp. des sommes amalgam\'ees], un $D$-diagramme de ${\mathcal C}$ a la forme
$$
\xymatrix{
&Y \ar[d]^g \\
X \ar[r]^f &S
} \qquad \xymatrix{
S \ar[r]^f \ar[d]_{\mbox{[resp.} \quad g} &X \ar@{}[d]^{\mbox{\quad ].}} \\
Y &
}
$$

Quand elle existe, sa limite [resp. sa colimite] est appel\'ee le produit fibr\'e de $X$ et $Y$ sur $S$ [resp. la somme amalgam\'ee de $X$ et $Y$ le long de $S$] et not\'ee
$$
X \times_S Y \qquad \mbox{[resp.} \quad X \amalg_S Y \ \mbox{].}
$$

Elle est caract\'eris\'ee par la propri\'et\'e que se donner un morphisme
$$
Z \longrightarrow X \times_S Y \qquad \mbox{[resp.} \quad X \amalg_S Y \longrightarrow Z \ \mbox{]}
$$
\'equivaut \`a se donner une paire de morphismes
$$
Z \xrightarrow{ \ f' \ } X , Z \xrightarrow{ \ g' \ } Y \qquad \mbox{[resp.} \quad X \xrightarrow{ \ f' \ } Z , Y \xrightarrow{ \ g' \ } Z \ \mbox{]}
$$
tel que $f \circ f' = g \circ g'$ [resp. $f' \circ f = g' \circ g$].

\smallskip

Autrement dit, $X \times_S Y$ [resp. $X \amalg_S Y$] est muni de deux morphismes canoniques vers $X$ et $Y$ [resp. depuis $X$ et $Y$] qui s'inscrivent dans un carr\'e commutatif
$$
\xymatrix{
X \times_S Y \ar[d] \ar[r]&Y \ar[d]^g \\
X \ar[r]^f &S
} \qquad \xymatrix{
S \ar[r]^f \ar[d]_{\mbox{[resp.} \quad g} &X \ar[d]^{\mbox{\qquad ]}} \\
Y \ar[r]&X \amalg_S Y
}
$$
tel que toute compl\'etion du diagramme de d\'epart en un carr\'e commutatif
$$
\xymatrix{
Z \ar[d] \ar[r]&Y \ar[d]^g \\
X \ar[r]^f &S
} \qquad \xymatrix{
S \ar[r]^f \ar[d]_{\mbox{[resp.} \quad g} &X \ar[d]^{\mbox{\qquad ]}} \\
Y \ar[r]&Z
}
$$
se factorise en un unique morphisme
$$
Z \longrightarrow X \times_S Y \qquad \mbox{[resp.} \quad X \amalg_S Y \longrightarrow Z \ \mbox{].}
$$

\subsection{Foncteurs limites et colimites}\label{subsec195}

\smallskip

Faisant varier les $D$-diagrammes d'une cat\'egorie ${\mathcal C}$, on peut encore demander:

\begin{defn}\label{defI93}

Soit ${\mathcal C}$ une cat\'egorie localement petite.

\smallskip

Si $D$ est un carquois, on dit que les $D$-limites [resp. les $D$-colimites] sont bien d\'efinies dans ${\mathcal C}$ si le foncteur diagonal
$$
\Delta_D : {\mathcal C} \longrightarrow D\mbox{\rm -diag} \, ({\mathcal C})
$$
a un adjoint \`a droite [resp. \`a gauche] not\'e
$$
\varprojlim_{D} : D\mbox{-{\rm diag}} \, ({\mathcal C}) \longrightarrow {\mathcal C} \qquad \mbox{[resp.} \quad \varinjlim_D : D\mbox{{\rm -diag}} \, ({\mathcal C}) \longrightarrow {\mathcal C} \ \mbox{]}.
$$
\end{defn}


\begin{remarksqed}
\begin{listeisansmarge}
\item  Les $D$-limites [resp. les $D$-colimites] sont bien d\'efinies dans ${\mathcal C}$ si tout $D$-diagramme $X_{\bullet}$ de ${\mathcal C}$ admet une limite [resp. une colimite] et que celle-ci peut \^etre choisie uniform\'ement en $X_{\bullet}$.

\medskip

\item Si les $D$-limites [resp. $D$-colimites] sont bien d\'efinies pour tout carquois $D$ sans fl\`eche et fini, on dit que ${\mathcal C}$ a des produits finis [resp. des sommes finies] arbitraires.

\smallskip

Dans ce cas, ${\mathcal C}$ a en particulier un objet terminal $1$ [resp. un objet initial $0$].

\medskip

\item Si cette propri\'et\'e est v\'erifi\'ee pour tout carquois $D$ sans fl\`eche, on dit que ${\mathcal C}$ a des produits [resp. des sommes] arbitraires.

\medskip

\item Si cette propri\'et\'e est v\'erifi\'ee pour tout carquois fini $D$ (avec ou sans fl\`eche), on dit que ${\mathcal C}$ a des limites finies [resp. des colimites finies] arbitraires.

\medskip

\item Enfin, si cette propri\'et\'e est v\'erifi\'ee pour tout carquois $D$ sans restriction, on dit que ${\mathcal C}$ a des limites [resp. des colimites] arbitraires, ou encore que ${\mathcal C}$ est une cat\'egorie compl\`ete [resp. cocompl\`ete]. 
\end{listeisansmarge}
\end{remarksqed}
\pagebreak

On remarque que les limites [resp. les colimites] sont bien d\'efinies s'il en est ainsi des produits [resp. des sommes] et des \'egalisateurs [resp. des co\'egalisateurs]:

\begin{lem}\label{lemI94}

Soit ${\mathcal C}$ une cat\'egorie localement petite.

\smallskip

Soient $D$ un carquois, ${\rm Hom} (D)$ l'ensemble de ses fl\`eches et $o,b : {\rm Hom} (D) \to {\rm Ob} (D)$ les deux applications qui associent \`a toute fl\`eche de $D$ son origine et son but.

\begin{listeimarge}

\item Etant donn\'e un $D$-diagramme $X_{\bullet}$ de ${\mathcal C}$, supposons que les deux produits [resp. les deux sommes]
$$
\prod_{d \in {\rm Ob} (D)} X_d \qquad \mbox{et} \quad \prod_{(u : d_1 \to d_2) \in {\rm Hom} (D)} X_{d_2}
$$
$$
\mbox{[resp.} \qquad \coprod_{d \in {\rm Ob} (D)} X_d \qquad \mbox{et} \quad \coprod_{(u : d_1 \to d_2) \in {\rm Hom} (D)} X_{d_1} \ \mbox{]}
$$
sont bien d\'efinis dans ${\mathcal C}$.

\smallskip

Alors un objet de ${\mathcal C}$ est une limite [resp. une colimite] du diagramme $X_{\bullet}$ si et seulement si il est un \'egalisateur [resp. un co\'egalisateur] des deux fl\`eches
$$
\prod_{d \in {\rm Ob} (D)} X_d \rightrightarrows \prod_{(u : d_1 \to d_2) \in {\rm Hom} (D)} X_{d_2}
$$
$$
\mbox{[resp.} \qquad \coprod_{(u : d_1 \to d_2) \in {\rm Hom} (D)} X_{d_1} \rightrightarrows \coprod_{d \in {\rm Ob} (D)} X_d \ \mbox{]}
$$
qui correspondent d'apr\`es le lemme de Yoneda aux deux familles d'applications index\'ees par les objets $Y$ de ${\mathcal C}$
$$
\begin{matrix}
&\displaystyle \prod_d {\rm Hom} (Y,X_d) &\rightrightarrows &\displaystyle \prod_{(u : d_1 \to d_2)} {\rm Hom} (Y,X_{d_2}) \, , \\
&\hfill (f_d : Y \to X_d)_d &\longmapsto &(f_{d_2} : Y \to X_{d_2})_u \hfill \\
\mbox{et} &(f_d : Y \to X_d)_d &\longmapsto &(X_u \circ f_{d_1} : Y \to X_{d_2})_u
\end{matrix}
$$

$$
\begin{matrix}
\mbox{[resp.} &\displaystyle \prod_d {\rm Hom} (X_d,Y) &\rightrightarrows &\displaystyle \prod_{(u : d_1 \to d_2)} {\rm Hom} (X_{d_1}, Y) \, , \\
&\hfill (f_d : X_d \to Y)_d &\longmapsto &(f_{d_1} : X_{d_1} \to Y)_u \hfill \\
\mbox{et} &(f_d : X_d \to Y)_d &\longmapsto &(f_{d_2} \circ X_u : X_{d_1} \to Y)_u \ \mbox{]}.
\end{matrix}
$$

\item Pour que le foncteur de limite $\underset{D}{\varprojlim}$ [resp. de colimite $\underset{D}{\varinjlim}$\,] soit bien d\'efini dans ${\mathcal C}$, il suffit que les foncteurs de produits [resp. de sommes]
$$
\begin{matrix}
&\displaystyle \prod_d &\mbox{et} &\displaystyle \prod_{(u : d_1 \to d_2)} \\
\mbox{[resp.} \quad &\displaystyle \coprod_d &\mbox{et} &\displaystyle \coprod_{(u : d_1 \to d_2)} \quad \mbox{]}
\end{matrix}
$$
soient bien d\'efinis dans ${\mathcal C}$, ainsi que le foncteur d'\'egalisation ${\rm eg}$ [resp. de co\'egalisation ${\rm coeg}$] associ\'e au carquois de la paire de fl\`eches.
\end{listeimarge}
\end{lem}

\begin{remarks}
\begin{listeisansmarge}
\item Il r\'esulte de ce lemme qu'une cat\'egorie localement petite ${\mathcal C}$ est compl\`ete [resp. cocompl\`ete] si et seulement si elle a des produits [resp. des sommes] arbitraires et un foncteur d'\'egalisation [resp. de co\'egalisation].

\medskip

\item De m\^eme, une cat\'egorie localement petite a des limites finies [resp. des colimites finies] arbitraires si et seulement si elle a des produits finis [resp. des sommes finies] 	arbitraires et un foncteur d'\'egalisation [resp. de co\'egalisation].

\smallskip

Cela \'equivaut aussi \`a demander que cette cat\'egorie ait un objet terminal $1$ [resp. un objet initial $0$] et un foncteur de produit fibr\'e [resp. de somme amalgam\'ee] associ\'e au carquois
$$
\xymatrix{
&\bullet \ar[d] \\
\bullet \ar[r] &\bullet
} \qquad\qquad \xymatrix{
\bullet \ar[r] \ar[d]_{\mbox{[resp. \quad}} &\bullet \ar@{}[d]^{\qquad\mbox{].}} \\
\bullet &
}
$$

\end{listeisansmarge}
\end{remarks}
\medskip

\begin{demolem}
\begin{listeisansmarge}
\item est une cons\'equence formelle de la d\'efinition des limites [resp. des colimites].

\medskip

\item r\'esulte de (i). 
\end{listeisansmarge}
\end{demolem}

\medskip

Il r\'esulte de ce lemme que la notion de produit fibr\'e [resp. de somme amalgam\'ee] est particuli\`erement importante.

\smallskip

C'est pourquoi on pose encore la d\'efinition suivante:

\begin{defn}\label{defI95}

Soit ${\mathcal C}$ une cat\'egorie localement petite.

\begin{listeimarge}

\item Un carr\'e commutatif de ${\mathcal C}$
$$
\xymatrix{
Z \ar[d] \ar[r] &X \ar[d] \\
Y \ar[r] &S
}
$$
est dit ``cart\'esien'' [resp. ``cocart\'esien''] s'il induit un isomorphisme
$$
Z \longrightarrow X \times_S Y \qquad \mbox{[resp.} \quad X \amalg_Z Y \longrightarrow S \ \mbox{]}.
$$

\item Un morphisme de ${\mathcal C}$
$$
X \longrightarrow S \qquad \mbox{[resp.} \quad Z \longrightarrow X \ \mbox{]}
$$
est dit ``carrable'' [resp. cocarrable''] si pour tout morphisme
$$
Y \xrightarrow{ \ v \ } S \qquad \mbox{[resp.} \quad Z \xrightarrow{ \ v \ } Y \ \mbox{],}
$$
le produit fibr\'e [resp. la somme amalgam\'ee]
$$
X \times_S Y \qquad \mbox{[resp.} \quad X \amalg_Z Y \ \mbox{]}
$$
existe dans ${\mathcal C}$ et que les $X \times_S Y$ [resp. $X \amalg_Z Y$] peuvent \^etre choisis uniform\'ement en $v$.
\end{listeimarge}
\end{defn}

\subsection{Transformation des limites ou des colimites par un foncteur}\label{subsec196}

On a le lemme \'evident:

\begin{lem}\label{lemI96}

Soit $F : {\mathcal C} \to {\mathcal D}$ un foncteur entre deux cat\'egories localement petites.

\smallskip

Soient $D$ un carquois et $X_{\bullet}$ un $D$-diagramme de ${\mathcal C}$ qui admet une limite [resp. une colimite] dans ${\mathcal C}$ et tel que $F(X_{\bullet})$ admette une limite [resp. une colimite] dans ${\mathcal D}$.

\smallskip

Alors ces limites [resp. ces colimites] sont reli\'ees par un morphisme canonique de ${\mathcal D}$
$$
 F (\varprojlim_D X_{\bullet}) \longrightarrow \displaystyle \varprojlim_D F(X_{\bullet}) 
 $$
 $$
\mbox{[resp.} \qquad  \varinjlim_D F(X_{\bullet}) \longrightarrow \displaystyle F (\varinjlim_D X_{\bullet}) \quad \mbox{].}
$$
\end{lem}

\begin{demo}

Quitte \`a remplacer ${\mathcal C}$ par ${\mathcal C}^{\rm op}$, on peut se contenter de traiter le cas des limites.

\smallskip

Si $X = \underset{D}{\varprojlim} X_{\bullet}$, le morphisme canonique de $D$-diag $({\mathcal C})$
$$
\Delta_X \longrightarrow X_{\bullet}
$$
est transform\'e par $F$ en un morphisme de $D$-diag $({\mathcal D})$
$$
\Delta_{F(X)} \longrightarrow F(X_{\bullet})
$$
qui correspond par adjonction \`a un morphisme de ${\mathcal D}$
$$
F(X) \longrightarrow \varprojlim_D F(X_{\bullet}) \, .
$$

\end{demo}


Ce lemme conduit \`a poser la d\'efinition suivante:

\begin{defn}\label{defI97}

Soit $F : {\mathcal C} \to {\mathcal D}$ un foncteur entre deux cat\'egories localement petites.

\smallskip

Soient $D$ un carquois et $X_{\bullet}$ un $D$-diagramme de ${\mathcal C}$ qui a une limite [resp. une colimite] dans ${\mathcal C}$.

\smallskip

On dit que $F$ respecte la limite [resp. la colimite] de $X_{\bullet}$ si, dans la cat\'egorie ${\mathcal D}$,
$$
F(\varprojlim_D X_{\bullet}) \ \mbox{est une limite du $D$-diagramme $F(X_{\bullet})$} 
$$
$$
\mbox{[resp.} \qquad F(\varinjlim_{D} X_{\bullet}) \ \mbox{est une colimite du $D$-diagramme $F(X_{\bullet})$ \ \mbox{]}.}
$$
\hfill $\Box$
\end{defn}

La propri\'et\'e pour un foncteur de respecter les limites finies ou les colimites finies (c'est-\`a-dire associ\'ees \`a des carquois finis) est tellement importante qu'on lui a donn\'e un nom:

\begin{defn}\label{defI98}
\begin{listeimarge}
\item Un foncteur $F : {\mathcal C} \to {\mathcal D}$ entre deux cat\'egories localement petites est appel\'e cart\'esien [resp. cocart\'esien] si les produits finis [resp. les sommes finies] sont toujours bien d\'efinis dans ${\mathcal C}$ et sont respect\'es par $F$.

\medskip

\item Un tel foncteur est dit exact \`a gauche [resp. exact \`a droite] si les limites finies [resp. les colimites finies] sont toujours bien d\'efinies dans ${\mathcal C}$ et sont respect\'ees par $F$.

\medskip

\item Un tel foncteur est dit exact s'il est \`a la fois exact \`a gauche et exact \`a droite.
\end{listeimarge}
\end{defn}

\medskip

\begin{remarkqed}

Il r\'esulte du lemme \ref{lemI94} qu'un foncteur cart\'esien [resp. cocart\'esien] $F : {\mathcal C} \to {\mathcal D}$ est exact \`a gauche [resp. \`a droite] si et seulement si le foncteur eg [resp. coeg] est bien d\'efini dans ${\mathcal C}$ et que $F$ respecte les \'egalisateurs [resp. les co\'egalisateurs]. 

\end{remarkqed}

\medskip

Une large classe de foncteurs qui respectent toutes les limites [resp. toutes les colimites] est fournie par la proposition suivante:

\begin{prop}\label{propI99}

Soit $F : {\mathcal C} \to {\mathcal D}$ un foncteur entre deux cat\'egories localement petites.

\smallskip

Si $F$ admet un adjoint \`a gauche [resp. un adjoint \`a droite], alors il respecte toutes les limites [resp. toutes les colimites].
\end{prop}

\begin{remark}

Cette proposition donne une condition n\'ecessaire pour qu'un foncteur $F : {\mathcal C} \to {\mathcal D}$ admette un adjoint \`a gauche [resp. \`a droite]: respecter toutes les limites [resp. toutes les colimites].

\smallskip

On verra plus loin que cette condition n\'ecessaire est aussi suffisante lorsque la cat\'egorie ${\mathcal C}$ poss\`ede certaines propri\'et\'es g\'en\'erales, en particulier lorsqu'elle est un topos.
\end{remark}

\medskip

\begin{demoprop}

Quitte \`a remplacer ${\mathcal C}$ et ${\mathcal D}$ par les cat\'egories oppos\'ees ${\mathcal C}^{\rm op}$ et ${\mathcal D}^{\rm op}$, on peut se contenter de traiter le cas o\`u $F$ admet un adjoint \`a droite $G : {\mathcal D} \to {\mathcal C}$.

\smallskip

Soient $D$ un carquois et $X_{\bullet}$ un $D$-diagramme dans ${\mathcal C}$ qui admet une colimite $X$ dans ${\mathcal C}$. Il s'agit de prouver que $F(X)$ est une colimite du $D$-diagramme $F(X_{\bullet})$ dans ${\mathcal D}$.

\smallskip

Or, pour tout objet $Y$ de ${\mathcal D}$, on a la suite de bijections canoniques
$$
\begin{matrix}
{\rm Hom} (F(X) , Y) &\cong &{\rm Hom} (X,G(Y)) \hfill &\mbox{par adjonction,} \hfill \\
&\cong &{\rm Hom} (X_{\bullet} , \Delta_{G(Y)}) &\mbox{par d\'efinition de $X$} \hfill \\
&&&\mbox{comme colimite de $X_{\bullet}$,} \\
&\cong &{\rm Hom} (F(X_{\bullet}) , \Delta_Y) &\mbox{par adjonction.} \hfill
\end{matrix}
$$
On remarque en effet que les deux foncteurs induits par $F$ et $G$
$$
\begin{matrix}
F : D\mbox{-diag} ({\mathcal C}) &\longrightarrow &D\mbox{-diag} ({\mathcal D}) \, , \\
G : D\mbox{-diag} ({\mathcal D}) &\longrightarrow &D\mbox{-diag} ({\mathcal C}) \hfill
\end{matrix}
$$
forment une paire de foncteurs adjoints.

\smallskip

Les bijections compos\'ees
$$
{\rm Hom} (F(X),Y) \cong {\rm Hom} (F(X_{\bullet}) , \Delta_Y)
$$
font de $F(X)$ une colimite du $D$-diagramme $F(X_{\bullet})$ dans ${\mathcal D}$. 

\end{demoprop}

\medskip

Cette proposition implique en particulier que les limites respectent les limites et que les colimites respectent les colimites:

\begin{cor}\label{corI910}

Soient ${\mathcal C}$ une cat\'egorie localement petite et $D$ un carquois tel que les limites [resp. les colimites] des $D$-diagrammes soient bien d\'efinies dans ${\mathcal C}$.

\smallskip

Alors le foncteur
$$
\varprojlim_D : D\mbox{\rm -diag} \, ({\mathcal C}) \longrightarrow {\mathcal C} 
$$
$$
\mbox{[resp.} \qquad  \varinjlim_D : D\mbox{\rm -diag} \, ({\mathcal C}) \longrightarrow {\mathcal C}  \quad \mbox{]}
$$
respecte toutes les limites [resp. toutes les colimites].
\end{cor}

\begin{remarkqed}

Explicitons ce que cela signifie, par exemple dans le cas des limites.

\smallskip

Pour tout carquois $D'$, un $D'$-diagramme de $D$-diag $({\mathcal C})$ est un $D' \times D$-diagramme $X_{\bullet , \bullet}$ de ${\mathcal C}$ dont les objets $X_{d',d}$ sont index\'es par les paires d'objets $d'$ de $D'$ et $d$ de $D$, et dont les fl\`eches $X_{u',u}$ sont index\'ees par les paires de fl\`eches $u'$ de $D'$ et $u$ de $D$.

\smallskip

Supposant que pour tout $d$, le diagramme $X_{\bullet , d}$ a une limite $X_d$ dans ${\mathcal C}$, les $X_d$ forment un $D$-diagramme dont la limite est aussi une limite du $D'$-diagramme $d' \mapsto \underset{D}{\varprojlim} \, X_{d',\bullet}$. 

\end{remarkqed}

\subsection{Exemples de constructions de limites ou de colimites}\label{subsec197}

\medskip

\noindent $\bullet$ {\bf Limites et colimites dans les ensembles:}

\smallskip

On a le lemme facile mais tr\`es important:

\begin{lem}\label{lemI911}

La cat\'egorie ${\rm Ens}$ des ensembles est compl\`ete et cocompl\`ete.

\smallskip

Plus pr\'ecis\'ement, pour tout carquois $D$ et tout $D$-diagramme d'ensembles $X_{\bullet}$, on a:

\begin{listeimarge}

\item Le diagramme $X_{\bullet}$ admet pour limite le sous-ensemble de $\underset{d \in {\rm Ob} (D)}{\prod} \, X_d$ constitu\'e des familles d'\'el\'ements
$$
(x_d \in X_d)_{d \in {\rm Ob} (D)}
$$ 
telles que, pour toute fl\`eche $u : d_1 \to d_2$ de $D$, on ait
$$
X_u (x_{d_1}) = x_{d_2} \, .
$$

\item Le diagramme $X_{\bullet}$ admet pour colimite l'ensemble quotient de la r\'eunion disjointe
$$
\coprod_{d \in {\rm Ob} (D)} X_d
$$ 
par la relation d'\'equivalence engendr\'ee par les paires d'\'el\'ements
$$
(x_{d_1} , X_u (x_{d_2}))
$$
associ\'ees \`a tout choix d'une fl\`eche $u : d_1 \to d_2$ de $D$ et d'un \'el\'ement $x_{d_1} \in X_{d_1}$.
\end{listeimarge}
\end{lem}

\begin{remark}

Si $X$ est un ensemble et $R$ un sous-ensemble de $X \times X$, la relation d'\'equivalence $\overline R$ engendr\'ee par $R$ est le plus petit sous-ensemble de $X \times X$ contenant $R$ et la diagonale et qui est sym\'etrique et transitif.

\smallskip

Il est constitu\'e des \'el\'ements
$$
(x,x) \, , \ x \in X \, ,
$$
et des \'el\'ements 
$$
(x,x')
$$
tels qu'existe une suite d'\'el\'ements
$$
x = x_0 , x_1 , \cdots , x_{k-1} , x_k = x'
$$
avec, pour tout $i$, $0 \leq i < k$,
$$
(x_i , x_{i+1}) \in R \qquad \mbox{ou} \qquad (x_{i+1} , x_i) \in R \, .
$$
\end{remark}

\medskip

\begin{demo}
\begin{listeisansmarge}
\item est \'evident sur la d\'efinition des limites.

\medskip

\item r\'esulte de ce que toute application $X \to Y$ d\'efinit une relation d'\'equivalence $X \times_Y X$ sur $X$ et de ce que, r\'eciproquement, toute relation d'\'equivalence $R$ sur $X$ d\'efinit un ensemble quotient $\overline X$ tel que le carr\'e
$$
\xymatrix{
R \ar[d] \ar[r] &X \ar[d] \\
X \ar[r] &\overline X
}
$$

soit \`a la fois cocart\'esien et cart\'esien. 
\end{listeisansmarge}
\end{demo}

\bigskip

\noindent $\bullet$ {\bf Limites et colimites dans les pr\'efaisceaux:}

\smallskip

Des ensembles, on passe automatiquement aux pr\'efaisceaux:

\begin{prop}\label{propI912}

Soit ${\mathcal C}$ une cat\'egorie essentiellement petite.

\smallskip

Alors la cat\'egorie $\widehat{\mathcal C}$ des pr\'efaisceaux sur ${\mathcal C}$ est compl\`ete et cocompl\`ete.

\smallskip

Plus pr\'ecis\'ement, pour tout carquois $D$ et tout $D$-diagramme $P_{\bullet}$ de $\widehat{\mathcal C}$, on a:

\begin{listeimarge}

\item Le diagramme $P_{\bullet}$ admet pour limite le pr\'efaisceau qui associe \`a tout objet $X$ de ${\mathcal C}$ l'ensemble limite
$$
\varprojlim_D P_{\bullet} (X)
$$
et \`a tout morphisme $u : X \to Y$ de ${\mathcal C}$ l'application induite 
$$
\varprojlim_{D} P_{\bullet} (Y) \longrightarrow \varprojlim_D P_{\bullet} (X) \, .
$$

\item Le diagramme $P_{\bullet}$ admet pour colimite le pr\'efaisceau qui associe \`a tout objet $X$ de ${\mathcal C}$ l'ensemble colimite
$$
\varinjlim_D P_{\bullet} (X)
$$
et \`a tout morphisme $u : X \to Y$ de ${\mathcal C}$ l'application induite
$$
\varinjlim_{D} P_{\bullet} (Y) \longrightarrow \varinjlim_D P_{\bullet} (X) \, . 
$$

\end{listeimarge}
\end{prop}

\begin{remark}

Cela implique que pour tout objet $X$ d'une cat\'egorie essentiellement petite ${\mathcal C}$, le foncteur d'\'evaluation en $X$
$$
\begin{matrix}
\widehat{\mathcal C} &\longrightarrow &{\rm Ens} \, , \\
P &\longmapsto &P(X)
\end{matrix}
$$
respecte les limites et les colimites.

\smallskip

En particulier, il est exact.
\end{remark}

\medskip

\begin{demo}

Cela r\'esulte de ce que, pour tout carquois $D$, $\underset{D}{\varprojlim}$ et $\underset{D}{\varinjlim}$ sont bien d\'efinis en tant que foncteurs
$$
D\mbox{-diag} \, ({\rm Ens}) \longrightarrow {\rm Ens} \, .
$$

\end{demo}

\medskip

Cette proposition a une cons\'equence tr\`es importante pour l'existence et le calcul des limites dans une cat\'egorie essentiellement petite:

\begin{cor}\label{corI913}

Soit ${\mathcal C}$ une cat\'egorie essentiellement petite.

\smallskip

Soient $D$ un carquois et $X_{\bullet}$ un $D$-diagramme de ${\mathcal C}$.

\smallskip

Alors un objet $X$ de ${\mathcal C}$ est une limite de $X_{\bullet}$ si et seulement si il repr\'esente le pr\'efaisceau limite
$$
\varprojlim_D y(X_{\bullet})
$$
du $D$-diagramme $y(X_{\bullet})$ d\'eduit de $X_{\bullet}$ par le foncteur de Yoneda
$$
y : {\mathcal C} \longrightarrow \widehat{\mathcal C} \, .
$$
\end{cor}

\begin{remark}

En revanche, il est tr\`es rare que l'image par $y$ d'une colimite d'un $D$-diagramme $X_{\bullet}$ dans ${\mathcal C}$ soit une colimite de $y(X_{\bullet})$ dans $\widehat{\mathcal C}$.

\smallskip

Comme on le verra, la th\'eorie des topos permet cependant de calculer sous certaines conditions le pr\'efaisceau
$$
y ( \varinjlim_D X_{\bullet})
$$
\`a partir du $D$-diagramme $y(X_{\bullet})$ dans $\widehat{\mathcal C}$.

\end{remark}

\medskip

\begin{demo}

Soit $P$ le pr\'efaisceau $\underset{D}{\varprojlim} \, y(X_{\bullet})$.

\smallskip

Pour tout objet $Y$ de ${\mathcal C}$, l'ensemble
$$
\begin{matrix}
P(Y) &= &{\rm Hom}_{\widehat{\mathcal C}} \, (y(Y) , P) \hfill \\
&= &{\rm Hom}_{\mbox{\scriptsize$D$-diag} (\widehat{\mathcal C})} (\Delta_{y(Y)} , y(X_{\bullet}))
\end{matrix}
$$
s'identifie d'apr\`es le lemme de Yoneda \`a l'ensemble
$$
{\rm Hom}_{\mbox{\scriptsize$D$-diag} ({\mathcal C})} (\Delta_Y , X_{\bullet}) \, .
$$
Par cons\'equent, un objet $X$ de ${\mathcal C}$ repr\'esente $P$ si et seulement si c'est une limite du $D$-diagramme $X_{\bullet}$ dans ${\mathcal C}$. 

\end{demo}
\pagebreak

\bigskip

\noindent $\bullet$ {\bf Limites et colimites dans les faisceaux:}

\smallskip

Pour tout espace topologique $X$, le foncteur $j_* : {\mathcal E}_X \hookrightarrow \widehat{O(X)}$ de plongement de la cat\'egorie ${\mathcal E}_X$ des faisceaux sur $X$ comme sous-cat\'egorie pleine de celle $\widehat{O(X)}$ des pr\'efaisceaux admet un adjoint \`a gauche, le foncteur de faisceautisation
$$
j^* : \widehat{O(X)} \longrightarrow {\mathcal E}_X \, .
$$
Consid\'erons un carquois $D$, un $D$-diagramme $F_{\bullet}$ dans ${\mathcal E}_X$ et son image $j_* F_{\bullet}$ dans $\widehat{O(X)}$.

\smallskip

Le $D$-diagramme $j_* F_{\bullet}$ a une limite $F$ et une colimite $P$ dans $\widehat{O(X)}$.

\smallskip

Comme $j^*$ a un adjoint \`a droite, il respecte les colimites et le faisceau $j^* P$ est une colimite du $D$-diagramme $j^* \circ j_* F_{\bullet}$. Or, comme $j_*$ est pleinement fid\`ele, le morphisme canonique $j^* \circ j_* \to {\rm id}$ est un isomorphisme. Donc $j^* P$ est une colimite du $D$-diagramme $F_{\bullet}$ dans ${\mathcal E}_X$.

\smallskip

Le fait que chaque $F_d$, $d \in {\rm Ob} (D)$, est un faisceau signifie que pour tout recouvrement d'un ouvert $U$ de $X$ par des ouverts $U_i$, $i \in I$, l'application
$$
F_d (U) \longrightarrow {\rm eg} \left( \prod_{i \in I} F_d (U_i) \rightrightarrows \prod_{i,j \in I} F_d (U_i \cap U_j) \right)
$$
est bijective.

\smallskip

Comme cette condition ne fait appara{\^\i}tre que des limites de diagrammes d'ensembles et que les limites respectent les limites, elle est \'egalement v\'erifi\'ee par la limite $F$ dans $\widehat{O(X)}$ du $D$-diagramme des $F_d$.

\smallskip

Donc $F$ est un faisceau et c'est une limite du $D$-diagramme $F_{\bullet}$ dans ${\mathcal E}_X$.

\smallskip

Ainsi, la cat\'egorie ${\mathcal E}_X$ des faisceaux sur $X$ est compl\`ete et cocompl\`ete, le foncteur $j_* : {\mathcal E}_X \hookrightarrow \widehat{O(X)}$ respecte les limites et son adjoint \`a gauche $j^* : \widehat{O(X)} \to {\mathcal E}_X$ respecte les colimites. 

\smallskip

On verra plus loin dans le contexte g\'en\'eral des topos que le foncteur de faisceautisation
$$
j^* : \widehat{O(X)} \longrightarrow {\mathcal E}_X
$$
respecte aussi les limites finies. En particulier, il est exact.

\medskip

\noindent $\bullet$ {\bf Limites et colimites d'ensembles munis d'une structure alg\'ebrique:}

\smallskip

Soit ${\mathcal A}$ la cat\'egorie des ensembles munis d'un certain type de structure alg\'ebrique, telle que celle de mono{\"\i}de, de groupe, d'anneau ou de module sur un anneau fix\'e.

\smallskip

Soient $D$ un carquois et $X_{\bullet}$ un $D$-diagramme de ${\mathcal A}$ constitu\'e d'objets $X_d$ reli\'es par des applications $X_u : X_{d_1} \to X_{d_2}$ index\'ees par les fl\`eches $u : d_1 \to d_2$ de $D$, qui respectent les structures des $X_{d_1}$ et $X_{d_2}$.

\smallskip

Alors l'ensemble produit
$$
\prod_{d \in {\rm Ob} (D)} X_d
$$
a une structure alg\'ebrique d\'eduite de celle des $X_d$ et son sous-ensemble d\'efini par la famille d'\'equations
$$
X_u (x_{d_1}) = x_{d_2}
$$
est stable par les op\'erations qui d\'efinissent le type de structure alg\'ebrique consid\'er\'e.

\smallskip

Cela r\'esulte de ce qu'un type de structure alg\'ebrique s'exprime en termes de produits finis, de fl\`eches et de conditions de commutativit\'e de diagrammes (par exemple une op\'eration interne sur un objet $X$ est une fl\`eche $X \times X \to X$ et elle est associative si et seulement si les deux fl\`eches compos\'ees $X \times X \times X \rightrightarrows X$ qui s'en d\'eduisent sont \'egales). En effet, le foncteur de limite $\underset{D}{\varprojlim}$ dans ${\rm Ens}$ respecte toute structure exprim\'ee en de tels termes.

\smallskip

Donc la cat\'egorie ${\mathcal A}$ est compl\`ete et le foncteur d'oubli de la structure alg\'ebrique
$$
{\mathcal A} \longrightarrow {\rm Ens}
$$
respecte toutes les limites.

\smallskip

La cat\'egorie ${\mathcal A}$ est \'egalement cocompl\`ete.

\smallskip

Pour construire la colimite d'un $D$-diagramme $X_{\bullet}$ de ${\mathcal A}$, on commence par former la r\'eunion disjointe
$$
X = \coprod_{d \in {\rm Ob} (D)} X_d
$$
des ensembles sous-jacents des objets $X_d$, $d \in {\rm Ob} (D)$.

\smallskip

Puis on forme l'objet ``libre'' $\widetilde X$ de ${\mathcal A}$ engendr\'e par l'ensemble $X$.

\smallskip

Enfin, la colimite de $X_{\bullet}$ dans ${\mathcal A}$ est l'objet quotient de $\widetilde X$ par les relations alg\'ebriques qui existent dans chaque $X_d$ et par les relations
$$
x_{d_2} = u(x_{d_1})
$$
associ\'ees aux fl\`eches $u : d_1 \to d_2$ de $D$.

\medskip

\noindent $\bullet$ {\bf Limites et colimites d'espaces topologiques:}

\smallskip

Le cas de la cat\'egorie des ensembles entra{\^\i}ne encore:

\begin{prop}\label{propI914}

La cat\'egorie ${\rm Top}$ des espaces topologiques est compl\`ete et cocompl\`ete et le foncteur d'oubli des topologies
$$
{\rm Top} \longrightarrow {\rm Ens}
$$
respecte toutes les limites et toutes les colimites.

\smallskip

Plus pr\'ecis\'ement, on a pour tout carquois $D$ et tout $D$-diagramme $X_{\bullet}$ de ${\rm Top}$:

\begin{listeimarge}

\item Le diagramme $X_{\bullet}$ admet pour limite l'ensemble
$$
\varprojlim_D X_{\bullet} \qquad \mbox{(calcul\'e dans ${\rm Ens}$)}
$$
muni de la topologie la moins fine pour laquelle chacune des projections canoniques
$$
\varprojlim_D X_{\bullet} \longrightarrow X_d \, , \qquad d \in {\rm Ob} \, ,
$$
est continue.

\medskip

\item Le diagramme $X_{\bullet}$ admet pour colimite l'ensemble
$$
\varinjlim_D X_{\bullet} \qquad \mbox{(calcul\'e dans ${\rm Ens}$)}
$$
muni de la topologie dont les ouverts sont les parties $U$ dont l'image r\'eciproque par chaque projection canonique
$$
X_d \longrightarrow \varinjlim_D X_{\bullet} \, , \qquad d \in {\rm Ob} (D) \, ,
$$
est un ouvert de $X_d$.
\end{listeimarge}
\end{prop}

\begin{remark}

En particulier, si $D$ est un carquois sans fl\`eche r\'eduit \`a un ensemble $I$ d'objets,
$$
\prod_{d \in I} X_d \qquad \mbox{et} \qquad \coprod_{d \in I} X_d
$$
sont l'espace topologique produit des $X_d$ et l'espace topologique somme des $X_d$ au sens usuel.

\smallskip

Dans le cas g\'en\'eral,
$$
\varprojlim_{D} X_d \xhookrightarrow{ \ \ \ } \prod_{d \in {\rm Ob} (D)} X_d
$$
est muni de la topologie induite par celle du produit, et $\underset{D}{\varinjlim} \, X_d$ est muni de la topologie induite par l'application surjective
$$
\coprod_{d \in {\rm Ob} (D)} X_d \longrightarrow \varinjlim_D X_{\bullet} \, .
$$
\end{remark}

\medskip

\begin{demo}

Ce sont des cons\'equences formelles de la d\'efinition de la notion d'application continue. 

\end{demo}

\medskip

\noindent $\bullet$ {\bf Limites et colimites de vari\'et\'es diff\'erentielles ou analytiques:}

\smallskip

On a d'abord le lemme \'evident:

\begin{lem}\label{lemI915}

La cat\'egorie ${\mathcal V}$ des vari\'et\'es diff\'erentielles de classe $C^k$ [resp. des vari\'et\'es analytiques] admet des produits finis arbitraires et des sommes arbitraires.

\smallskip

Plus pr\'ecis\'ement:

\begin{listeimarge}

\item Pour toutes vari\'et\'es $X_1 , \cdots , X_k$, leur produit dans ${\mathcal V}$ est l'espace topologique produit $X_1 \times \cdots \times X_k$ muni du faisceau de structure ${\mathcal O}$ dont la restriction \`a tout produit $U_1 \times \cdots \times U_k$ d'ouverts $U_1 , \cdots , U_k$ de $X_1 , \cdots , X_k$ isomorphes \`a des ouverts $U'_1 \subset {\mathbb R}^{n_1} , \cdots , U'_k \subset {\mathbb R}^{n_k}$ [resp. $U'_1 \subset {\mathbb C}^{n_1} , \cdots , U'_k \subset {\mathbb C}^{n_k}$] s'identifie, par composition avec l'hom\'eomorphisme
$$
U'_1 \times \cdots \times U'_k \xrightarrow{ \ \sim \ } U_1 \times \cdots \times U_k \, ,
$$
avec le faisceau des fonctions r\'eelles de classe $C^k$ [resp. des fonctions complexes holomorphes] sur l'ouvert $U'_1 \times \cdots \times U'_k$ de ${\mathbb R}^{n_1} \times \cdots \times {\mathbb R}^{n_k}$ [resp. ${\mathbb C}^{n_1} \times \cdots \times {\mathbb C}^{n_k}$].

\medskip

\item Pour toutes vari\'et\'es $X_i$, $i \in I$, leur somme dans ${\mathcal V}$ est l'espace topologique $\underset{i \in I}{\coprod} \, X_i$ muni du faisceau de structure ${\mathcal O}$ dont la restriction \`a chaque ouvert $X_i$ de $\underset{i \in I}{\coprod} \, X_i$ est le faisceau de structure de $X_i$.
\end{listeimarge}
\end{lem}

\begin{remarksqed}
\begin{listeisansmarge}
\item Ainsi le foncteur d'oubli de la structure diff\'erentielle [resp. analytique]
$$
{\mathcal V} \longrightarrow {\rm Top}
$$
respecte les produits finis et les sommes.

\medskip

\item La sous-cat\'egorie pleine de ${\mathcal V}$ constitu\'ee des vari\'et\'es ``d\'enombrables \`a l'infini'' est stable par produits finis et par sommes finies ou d\'enombrables. 

\end{listeisansmarge}
\end{remarksqed}

\pagebreak

Pour aller plus loin, on a besoin de la d\'efinition suivante:

\begin{defn}\label{defI916}

Soit ${\mathcal V}$ la cat\'egorie des vari\'et\'es diff\'erentielles de classe $C^k$ [resp. des vari\'et\'es analytiques].

\begin{listeimarge}

\item Un morphisme de ${\mathcal V}$
$$
f : X \to Y
$$
est dit ``submersif'' si, pour tout diagramme commutatif
$$
\xymatrix{
U' \ar[d]_{f'} \ar[r]^{\sim} &U \ \ar@{^{(}->}[r] &X \ar[d]^f \\
V' \ar[r]^{\sim} &V \ \ar@{^{(}->}[r] &Y
}
$$
dont les fl\`eches horizontales sont les compos\'ees d'immersions ouvertes $U \hookrightarrow X$, $V \hookrightarrow Y$ et d'isomorphismes $U' \xrightarrow{ \ \sim \ } U$, $V' \xrightarrow{ \ \sim \ } V$ d\'efinis sur des ouverts $U',V'$ de ${\mathbb R}^n$, ${\mathbb R}^m$ [resp. ${\mathbb C}^n$, ${\mathbb C}^m$], la matrice des d\'eriv\'ees partielles de $f' = (f'_1 , \ldots , f'_m)$
$$
\frac{\partial f'}{\partial x} = \left( \frac{\partial f'_i}{\partial x_j} \right)_{1 \leq i \leq m \atop 1 \leq j \leq n}
$$
est de rang $m$ (avec donc $m \leq n$) en tout point de $U'$.

\medskip

\item Un diagramme de ${\mathcal V}$ de la forme
$$
X \underset{g}{\overset{f}{\rightrightarrows}} Y
$$
est dit ``transversal'' si pour tout point $x \in X$ tel que $f(x) = g(x)$ et tous voisinages ouverts $U \hookrightarrow X$ de $x$, $V \hookrightarrow Y$ de $f(x)$ isomorphes \`a des ouverts $U',V'$ de ${\mathbb R}^n , {\mathbb R}^m$ [resp. ${\mathbb C}^n , {\mathbb C}^m$] et tels que $f(U) \subset V$, $g (U) \subset V$, avec donc un diagramme induit
$$
U' \underset{g'}{\overset{f'}{\rightrightarrows}} V' \, ,
$$
la diff\'erence des deux matrices de d\'eriv\'ees partielles en le point $x' \in U'$ correspondant \`a $x$
$$
\frac{\partial f'}{\partial x} (x') - \frac{\partial g'}{\partial x} (x')
$$
est de rang $m$ (avec donc $m \leq n$).

\medskip

\item Pour tout carquois fini $D$, un $D$-diagramme $X_{\bullet}$ de ${\mathcal V}$ est dit ``transversal'' si le diagramme induit
$$
\prod_{d \in {\rm Ob} (D)} X_d \rightrightarrows \prod_{(u : d_1 \to d_2) \in {\rm Hom} (D)} X_{d_2}
$$
est transversal au sens de (ii).
\end{listeimarge}
\end{defn}

\begin{remarksqed}
\begin{listeisansmarge}
\item Toute immersion ouverte est un morphisme submersif. 

\smallskip

Tout compos\'e de morphismes submersifs est submersif.

\medskip

\item Un diagramme de produit fibr\'e
$$
\xymatrix{
{ \ } \ar@{}[d] \ar@{}[r] &X \ar[d]^f \\
Y \ar[r]^g &S
}
$$
est transversal si pour tous points $x \in X$, $y \in Y$ tels que $f(x) = g(y)$ et tous voisinages ouverts $U \hookrightarrow X$ de $x$, $V \hookrightarrow Y$ de $y$, $W \hookrightarrow S$ de $f(x) = g(y)$ isomorphes \`a des ouverts $U' , V' , W'$ de ${\mathbb R}^n , {\mathbb R}^{\ell} , {\mathbb R}^m$ [resp. ${\mathbb C}^n , {\mathbb C}^{\ell} , {\mathbb C}^m$] et tels que $f(U) \subset W$, $g(V) \subset W$, avec donc un diagramme induit
$$
\xymatrix{
{ \ } \ar@{}[d] \ar@{}[r] &U' \ar[d]^{f'} \\
V' \ar[r]^{g'} &W'
}
$$
les images des deux matrices de d\'eriv\'ees partielles en les points $x' \in U'$, $y' \in V'$ correspondant \`a $x$ et $y$
$$
\frac{\partial f'}{\partial x} (x') \, , \quad \frac{\partial g'}{\partial y} (y') \, ,
$$
engendrent l'espace ${\mathbb R}^m$ [resp. ${\mathbb C}^m$] tout entier.

\medskip

\item Cette condition est r\'ealis\'ee en particulier si le morphisme
$$
f : X \to S \qquad \mbox{ou} \qquad g : Y \to S
$$
est submersif. 

\end{listeisansmarge}
\end{remarksqed}

\medskip

On a encore besoin de la d\'efinition suivante:

\begin{defn}\label{defI917}

Soit $X$ une vari\'et\'e diff\'erentielle de classe $C^k$ [resp. une vari\'et\'e analytique].

\smallskip

Un sous-espace ferm\'e $Y$ de l'espace topologique sous-jacent \`a $X$ est dit ``r\'egulier'' si $X$ peut \^etre recouvert par des ouverts $U$ isomorphes \`a des ouverts $U'$ d'un ${\mathbb R}^n$ [resp. ${\mathbb C}^n$] tels que l'image de $Y$ dans chaque $U'$ soit l'intersection de $U'$ et d'un sous-espace affine de ${\mathbb R}^n$ [resp. ${\mathbb C}^n$].

\end{defn}

Tout ferm\'e r\'egulier d'une vari\'et\'e diff\'erentielle [resp. analytique] h\'erite d'une structure de vari\'et\'e:

\begin{lem}\label{lemfI918}

Dans les conditions de la d\'efinition pr\'ec\'edente, le ferm\'e r\'egulier $Y$ de $X$ a une unique structure de vari\'et\'e diff\'erentielle de classe $C^k$ [resp. de vari\'et\'e analytique] telle que:
\begin{enumerate}
\item[$\bullet$] le plongement $Y \hookrightarrow X$ devient un morphisme de vari\'et\'es,
\item[$\bullet$] tout morphisme de vari\'et\'es $Y' \to X$ dont l'image de l'espace sous-jacent est contenue dans $Y$ se factorise de mani\`ere unique en un morphisme de vari\'et\'es
$$
Y' \to Y \, .
$$
\end{enumerate}
\end{lem}

\begin{demo}

Il suffit de d\'emontrer ce lemme localement sur $X$.

\smallskip

Cela ram\`ene au cas \'evident o\`u $X$ est un ouvert $U'$ d'un ${\mathbb R}^n$ [resp. ${\mathbb C}^n$] et $Y$ est l'intersection de cet ouvert et d'un sous-espace affine qu'on peut d'ailleurs supposer de la forme ${\mathbb R}^m \times \{0\}$ [resp. ${\mathbb C}^m \times \{0\}$]. 

\end{demo}

\medskip

On peut maintenant \'enoncer des conditions g\'en\'erales d'existence de limites ou de colimites de diagrammes de vari\'et\'es diff\'erentielles [resp. analytiques]:

\begin{thm}\label{thmI919}

Soit ${\mathcal V}$ la cat\'egorie des vari\'et\'es diff\'erentielles de classe $C^k$ [resp. des vari\'et\'es analytiques].

\begin{listeimarge}

\item Soient $D$ un carquois fini, $X_{\bullet}$ un $D$-diagramme de ${\mathcal V}$ et
$$
\prod_{d \in {\rm Ob} (D)} X_d  \underset{p_2}{\overset{p_1}{\rightrightarrows}} \prod_{(u : d_1 \to d_2) \in {\rm Hom} (D)} X_{d_2}
$$
le diagramme \`a deux fl\`eches parall\`eles qui lui est associ\'e.

\smallskip

Si ce diagramme est transversal, le sous-espace topologique ferm\'e
$$
Y \xhookrightarrow{ \ \ } \prod_{d \in {\rm Ob} (D)} X_d
$$
d\'efini par l'\'equation $p_1 = p_2$ est r\'egulier et, muni de sa structure canonique de vari\'et\'e, il est une limite du $D$-diagramme $X_{\bullet}$ dans la cat\'egorie ${\mathcal V}$.

\medskip

\item Soient $D$ un carquois, $X_{\bullet}$ un $D$-diagramme de ${\mathcal V}$ et
$$
R = \coprod_{(u:d_1 \to d_2) \in {\rm Hom} (D)} X_{d_1}  \underset{q_2}{\overset{q_1}{\rightrightarrows}} \coprod_{d \in {\rm Ob} (D)} X_d = X
$$
le diagramme \`a deux fl\`eches parall\`eles qui lui est associ\'e.

\smallskip

Soit $\overline R$ le plus petit ferm\'e de $X \times X$ qui est une relation d'\'equivalence et qui contient l'image de
$$
R \xrightarrow{ \ (q_1 , q_2) \ } X \times X \, .
$$

Si $\overline R$ est un ferm\'e r\'egulier, l'espace topologique $\overline X$ quotient de $X$ par la relation d'\'equivalence $\overline R$ a une structure canonique de vari\'et\'e diff\'erentielle [resp. analytique] qui en fait une colimite du $D$-diagramme $X_{\bullet}$ dans la cat\'egorie ${\mathcal V}$.

\smallskip

De plus, le morphisme
$$
X \longrightarrow \overline X
$$
est submersif et le carr\'e commutatif
$$
\xymatrix{
\overline R \ar[d] \ar[r] &X \ar[d] \\
X \ar[r] &\overline X
}
$$
est \`a la fois cocart\'esien et cart\'esien dans la cat\'egorie ${\mathcal V}$.

\end{listeimarge}
\end{thm}

\begin{remarks}
\begin{listeisansmarge}
\item En particulier, pour tout morphisme submersif
$$
f : X \longrightarrow S \, ,
$$
tout morphisme $g : Y \to S$ d\'efinit un produit fibr\'e
$$
X \times_S Y
$$
dans ${\mathcal V}$ et le morphisme $X \times_S Y \to Y$ est encore submersif.

\smallskip

Autrement dit, tous les morphismes submersifs sont carrables et leur classe est ``stable par changement de base''.

\medskip

\item La partie (ii) du th\'eor\`eme s'applique en particulier \`a la question de l'existence d'un quotient d'une vari\'et\'e diff\'erentielle [resp. analytique] $X$ par l'action d'un groupe discret $G$ qui agit sur elle.

\smallskip

Le diagramme de fl\`eches parall\`eles associ\'e \`a une telle situation s'\'ecrit
$$
R = \coprod_{\sigma \in G} X \underset{q_2}{\overset{q_1}{\rightrightarrows}} X
$$
o\`u les composantes de $q_1$ sont les automorphismes de $X$ associ\'es aux \'el\'ements de $G$ et celles de $q_2$ sont toutes \'egales \`a ${\rm id}_X$.

\smallskip

Dans ce cas, la relation d'\'equivalence ferm\'ee $\overline R$ est l'adh\'erence dans $X \times X$ de l'image de $R$ par $(q_1 , q_2)$.

\smallskip

Pour que $\overline R$ soit un ferm\'e r\'egulier, il suffit que le groupe $G$ agisse sur $X$ librement et discr\`etement, c'est-\`a-dire que l'application continue
$$
\coprod_{\sigma \in G} X \xrightarrow{ \ (q_1 , q_2) \ } X \times X
$$
soit injective et ait pour image un ferm\'e.
\end{listeisansmarge}
\end{remarks}

\medskip

\begin{demothm}

Elle est laiss\'ee au lecteur.

\smallskip

C'est une cons\'equence du ``th\'eor\`eme des fonctions implicites''. 

\end{demothm}

\bigskip

\noindent $\bullet$ {\bf Limites et colimites de sch\'emas affines:}

\smallskip

La cat\'egorie ${\rm Aff}$ des sch\'emas affines est l'oppos\'ee de la cat\'egorie des anneaux commutatifs. Comme cette derni\`ere est la cat\'egorie des ensembles munis d'une structure d'anneau commutatif et que ce type de structure est alg\'ebrique, la cat\'egorie ${\rm Aff}$ est compl\`ete et cocompl\`ete:

\smallskip

On peut pr\'eciser la forme des limites et des colimites dans ${\rm Aff}$:

\begin{prop}\label{propI920}
\begin{listeimarge}

\item Soient $D$ un carquois et $X_{\bullet}$ un $D$-diagramme de ${\rm Aff}$ constitu\'e de sch\'emas affines $X_d = {\rm Spec} (A_d)$ et de morphismes $X_{d_1} \xrightarrow{ \ X_u \ } X_{d_2}$, $(u : d_1 \to d_2) \in {\rm Hom} (D)$, associ\'es \`a des morphismes d'anneaux commutatifs $A_u : A_{d_2} \to A_{d_1}$.

\smallskip

Alors le $D$-diagramme $X_{\bullet}$ admet pour colimite dans ${\rm Aff}$ le sch\'ema affine
$$
{\rm Spec} (A)
$$
associ\'e au sous-anneau
$$
A \xhookrightarrow{ \ \ } \prod_{d \in {\rm Ob} (D)} A_d
$$
d\'efini comme
$$
A = \{ (a_d)_{d \in {\rm Ob} (D)} \mid a_{d_1} = A_u (a_{d_2}) \, , \ \forall \, (d_1 \xrightarrow{ \, u \, } d_2) \in {\rm Hom} (D) \} \, .
$$

\item Tout diagramme de changement de base dans ${\rm Aff}$
$$
\xymatrix{
{ \ } \ar@{}[d] \ar@{}[r] &{\rm Spec} (B) \ar[d] \\
{\rm Spec} (C) \ar[r] &{\rm Spec}(A)
}
$$
correspondant \`a un diagramme d'anneaux commutatifs
$$
\xymatrix{
A \ar[d] \ar[r] &B \\
C
}
$$
admet pour limite dans ${\rm Aff}$ le sch\'ema affine
$$
{\rm Spec} (B \otimes_A C)
$$
o\`u $B \otimes_A C$ est le $A$-module produit tensoriel de $B$ et $C$ muni de la structure d'anneau commutatif d\'eduite de celles de $B$ et $C$.

\end{listeimarge}
\end{prop}

\begin{remarks}
\begin{listeisansmarge}
\item En particulier, la somme $\underset{i \in I}{\coprod} \, X_i$ dans ${\rm Aff}$ de sch\'emas affines $X_i = {\rm Spec} (A_i)$ est le sch\'ema affine
$$
{\rm Spec} \left( \prod_{i \in I} A_i \right) \, .
$$

\item Si $B = A [X_1 , \cdots , X_n]$ est un anneau de polyn\^omes, le produit fibr\'e
$$
{\rm Spec} (A[X_1 , \cdots  , X_n]) \times_{{\rm Spec} (A)} {\rm Spec} (C)
$$
s'identifie \`a 
$$
{\rm Spec} (C[X_1 , \cdots , X_n]) \, .
$$
En particulier
$$
{\rm Spec} (A[X_1 , \cdots , X_n]) \times_{{\rm Spec} (A)} {\rm Spec} (A[X_{n+1} , \cdots , X_{n+m}])
$$
s'identifie \`a 
$$
{\rm Spec} (A[X_1 , \cdots , X_{n+m}]) \, .
$$
Autrement dit, le produit fibr\'e sur un sch\'ema affine de base des deux espaces affines de dimensions $n$ et $m$ sur cette base est l'espace affine de dimension $n+m$.
\end{listeisansmarge}
\end{remarks}

\medskip

\begin{demo}
\begin{listeisansmarge}
\item est un cas particulier de l'expression g\'en\'erale des limites dans la cat\'egorie des ensembles munis d'un type de structure alg\'ebrique.

\medskip

\item Se donner un carr\'e commutatif d'anneaux commutatifs
$$
\xymatrix{
A \ar[d] \ar[r] &B \ar[d] \\
C \ar[r] &D
}
$$
\'equivaut en effet \`a se donner un morphisme d'anneaux commutatifs
$$
B \otimes_A C \longrightarrow D \, .
$$

\end{listeisansmarge}
\end{demo}

\smallskip

\noindent $\bullet$ {\bf Limites finies et sommes de sch\'emas:}

\smallskip

On remarque d'abord le r\'esultat g\'en\'eral:

\begin{lem}\label{lemI921}

Soit ${\mathcal G}$ une sous-cat\'egorie g\'eom\'etrique de la cat\'egorie ${\rm Top}_{\rm an}$ des espaces annel\'es.

\smallskip

Alors toute famille $(X_i)_{i \in I}$ d'objets de ${\mathcal G}$ a une somme
$$
\coprod_{i \in I} X_i
$$
qui est constitu\'ee de la somme des espaces topologiques sous-jacents aux $X_i$, munie du faisceau de structure ${\mathcal O}$ dont la restriction \`a chaque $X_i$ est le faisceau de structure ${\mathcal O}_{X_i}$ de celui-ci.
\end{lem}

\begin{remark}

On a d\'ej\`a appliqu\'e cette construction dans le cas de la cat\'egorie g\'eom\'etrique des vari\'et\'es diff\'erentielles de classe $C^k$ [resp. des vari\'et\'es analytiques].
\end{remark}

\medskip

\begin{demo}

L'objet $\underset{i \in I}{\coprod} \, X_i$ de ${\rm Top}_{\rm an}$ est dans ${\mathcal G}$ car il est recouvert par les ouverts $X_i$ qui sont des objets de ${\mathcal G}$.

\smallskip

De plus, un morphisme de ${\rm Top}_{\rm an}$
$$
\coprod_{i \in I} X_i \longrightarrow Y
$$
est dans ${\mathcal G}$ si et seulement si ses restrictions $X_i \to Y$ sont des morphismes de ${\mathcal G}$.

\smallskip

D'o\`u la conclusion. 

\end{demo}

\medskip

Ce lemme s'applique en particulier \`a la cat\'egorie g\'eom\'etrique des sch\'emas et on obtient:

\begin{cor}\label{corI922}

La cat\'egorie ${\rm Sch}$ des sch\'emas a des sommes arbitraires. De plus, le foncteur d'inclusion
$$
{\rm Aff} \longrightarrow {\rm Sch}
$$
respecte les sommes finies.
\end{cor}

\begin{remarkqed}

En revanche, le foncteur ${\rm Aff} \to {\rm Sch}$ ne respecte pas les sommes infinies.

\end{remarkqed}

\bigskip

D'autre part, la cat\'egorie des sch\'emas ${\rm Sch}$ poss\`ede la propri\'et\'e fondamentale que tout morphisme y est carrable:

\begin{thm}\label{thmI923}

Toutes les limites finies sont bien d\'efinies dans la cat\'egorie ${\rm Sch}$ des sch\'emas.

\smallskip

De plus, le foncteur d'inclusion
$$
{\rm Aff} \longrightarrow {\rm Sch}
$$
respecte les limites finies.
\end{thm}

\begin{remark}

Pour toute immersion ouverte $V \hookrightarrow Y$ et tout morphisme de sch\'emas $X \xrightarrow{ \ f \ } Y$, le produit fibr\'e $V \times_Y X$ est le sous-sch\'ema ouvert $f^{-1} (V)$ de $X$.

\smallskip

On dit que les immersions ouvertes sont stables par changement de base.
\end{remark}

\medskip

\begin{demosansqed}

Elle r\'esulte facilement du lemme suivant qui repr\'esente une seconde mani\`ere de voir les sch\'emas et dont la preuve est laiss\'ee au lecteur:
\end{demosansqed}

\begin{lem}\label{lemI924}
\begin{listeimarge}
\item Tout sch\'ema $X$ d\'efinit un foncteur contravariant
\begin{eqnarray}
{\rm Aff}^{\rm op} &\longrightarrow &{\rm Ens} \, , \nonumber \\
{\rm Spec} (A) &\longmapsto &{\rm Hom} ({\rm Spec} (A),X) = X(A) \, . \nonumber
\end{eqnarray}

\item Le foncteur ainsi d\'efini
$$
{\rm Sch} \longrightarrow [{\rm Aff}^{\rm op} , {\rm Ens}] = \widehat{\rm Aff}
$$
est pleinement fid\`ele.

\medskip

\item Un foncteur contravariant
$$
X : {\rm Aff}^{\rm op} \longrightarrow {\rm Ens}
$$
est un sch\'ema si et seulement si il existe des morphismes
$$
x_i : {\rm Hom} (\bullet , {\rm Spec} (A_i)) \longrightarrow X
$$
tels que:

\bigskip

$\left\{\begin{matrix}
\bullet &\mbox{les ${\rm Hom} (\bullet , {\rm Spec} (A_i))$ sont des objets repr\'esentables de $\widehat{\rm Aff}$,} \hfill \\
&\mbox{c'est-\`a-dire associ\'es \`a des objets ${\rm Spec} (A_i)$ de ${\rm Aff}$,} \hfill \\
{ \ } \\
\bullet &\mbox{chaque $x_i$ est ouvert au sens que pour tout objet ${\rm Spec} (A)$ de ${\rm Aff}$ et} \hfill \\
&\mbox{tout morphisme ${\rm Hom} (\bullet , {\rm Spec} (A)) \to X$ le produit fibr\'e dans $\widehat{\rm Aff}$} \hfill \\
{ \ } \\
&{\rm Hom} (\bullet , {\rm Spec} (A_i)) \times_X {\rm Hom} (\bullet , {\rm Spec} (A)) \\
{ \ } \\
&\mbox{est repr\'esentable par un sous-sch\'ema ouvert ${\rm Spec} (A_i) \times_X {\rm Spec} (A)$ du sch\'ema affine ${\rm Spec} (A)$,} \hfill \\
{ \ } \\
\bullet &\mbox{la famille des $x_i$ forme un recouvrement au sens que pour tout morphisme ${\rm Hom} (\bullet , {\rm Spec} (A)) \to X$,} \hfill \\
&\mbox{les ouverts} \hfill \\
&{\rm Spec} (A_i) \times_X {\rm Spec} (A) \\
{ \ } \\
&\mbox{de ${\rm Spec} (A)$ forment un recouvrement de celui-ci.} \hfill
\end{matrix} \right.
$
\end{listeimarge}
\end{lem}

\bigskip

\begin{remarkqed}

Pour tout sch\'ema $X$ et pour tout anneau commutatif $A$, l'ensemble
$$
X(A) = {\rm Hom} ({\rm Spec} (A) , X)
$$
est appel\'e l'ensemble des points de $X$ \`a valeurs (ou \`a coefficients) dans $A$.

\smallskip

Ainsi, un morphisme de sch\'emas $X \to Y$ peut \^etre vu comme une famille d'applications
$$
X(A) \longrightarrow Y(A)
$$
telle que, pour tout morphisme d'anneaux commutatifs $A \to B$, le carr\'e
$$
\xymatrix{
X(A) \ar[d] \ar[r] &Y(A) \ar[d] \\
X(B) \ar[r] &Y(B)
}
$$

soit commutatif. 

\end{remarkqed}

\section{Cat\'egories relatives}\label{sec110}

\subsection{La notion de cat\'egorie relative (et celle de cat\'egorie fl\'ech\'ee)}\label{subsec1101}

\begin{defn}\label{defI101}

Si ${\mathcal C}$ est une cat\'egorie et $S$ un objet de ${\mathcal C}$, la cat\'egorie relative \`a droite ${\mathcal C} / S$ [resp. \`a gauche $S \backslash {\mathcal C}$] est d\'efinie de la mani\`ere suivante:
\begin{enumerate}
\item[$\bullet$] Ses objets sont les morphismes de ${\mathcal C}$
$$
X \xrightarrow{ \ p \ } S \qquad \mbox{[resp.} \quad S \xrightarrow{ \ i \ } X \ \mbox{]}
$$
dont le but [resp. l'origine] est l'objet $S$.
\item[$\bullet$] Ses morphismes
$$
(X_1 \xrightarrow{ \ p_1 \ } S) \longrightarrow (X_2 \xrightarrow{ \ p_2 \ } S) \qquad \mbox{[resp.} \quad (S \xrightarrow{ \ i_1 \ } X_1) \longrightarrow (S \xrightarrow{ \ i_2 \ } X_2) \ \mbox{]}
$$
sont les morphismes de ${\mathcal C}$
$$
X_1 \xrightarrow{ \ u \ } X_2
$$
tels que le triangle

$$\xymatrix{
	X_1 \ar[rd]_{p_1} \ar[rr]^u &&X_2 \ar[ld]^{p_2} \\
	&S
} \qquad\qquad \xymatrix{
	&S \ar[ld]_{\mbox{[resp. \quad} i_1} \ar[rd]^{i_2\qquad\mbox{].}} \\
	X_1 \ar[rr]^u &&X_2
}
$$

soit commutatif.

\end{enumerate}
\end{defn}

\pagebreak

\begin{remarks}
\begin{listeisansmarge}
\item Si ${\mathcal C}$ est une cat\'egorie localement petite [resp. petite], il en est de m\^eme de ses cat\'egories relatives ${\mathcal C}/S$ ou $S \backslash {\mathcal C}$.

\medskip

\item Les cat\'egories relatives sont munies d'un foncteur d'oubli
$$
\begin{matrix}
\hfill {\mathcal C} / S &\longrightarrow &{\mathcal C} \, , \\
(X \to S) &\longmapsto &X \hfill
\end{matrix} \qquad \mbox{[resp.} \quad \begin{matrix}
\hfill S \backslash {\mathcal C} &\longrightarrow &{\mathcal C} \, , \\
(S \to X) &\longmapsto &X \hfill
\end{matrix} \ \mbox{]}.
$$

\item Si $S$ est un objet terminal $1$ [resp. un objet initial $0$] de ${\mathcal C}$, le foncteur d'oubli
$$
{\mathcal C} / S \longrightarrow {\mathcal C}  \qquad \mbox{[resp.} \quad S \backslash {\mathcal C} \longrightarrow {\mathcal C} \ \mbox{]}
$$
est un isomorphisme de cat\'egories.

\medskip

\item Toute cat\'egorie relative \`a droite ${\mathcal C} / S$ [resp. \`a gauche $S \backslash {\mathcal C}$] admet pour objet terminal [resp. initial] l'objet $S \xrightarrow{ \ {\rm id} \ } S$.

\medskip

\item Tout morphisme dans une cat\'egorie ${\mathcal C}$
$$
S \xrightarrow{ \ s \ } S'
$$
d\'efinit un foncteur entre cat\'egories relatives \`a droite [resp. \`a gauche]
$$
\begin{matrix}
\hfill {\mathcal C} / S &\longrightarrow &{\mathcal C}/S' \, , \hfill \\
(X \to S) &\longmapsto &(X \to S \xrightarrow{u} S') \hfill
\end{matrix} \qquad \mbox{[resp.} \quad \begin{matrix}
\hfill S' \backslash {\mathcal C} &\longrightarrow &S \backslash{\mathcal C} \, , \hfill \\
(S' \to X) &\longmapsto &(S \to S' \to X) \hfill
\end{matrix} \ \mbox{]}
$$
par composition \`a droite [resp. \`a gauche] avec $u$.

\medskip

\item Tout morphisme d'une cat\'egorie ${\mathcal C}$
$$
S \xrightarrow{ \ s \ } S'
$$
qui est carrable [resp. cocarrable] d\'efinit un foncteur
$$
\begin{matrix}
\hfill {\mathcal C} / S' &\longrightarrow  &{\mathcal C} / S \, , \hfill \\
(X \to S') &\longmapsto &(X \times_{S'} S \to S)
\end{matrix}
$$
$$
\begin{matrix}
\mbox{[resp.} &\hfill S \backslash {\mathcal C} &\longrightarrow &S \backslash {\mathcal C}' \hfill \\
&(S \to X) &\longmapsto &(S' \to S' \amalg_S X) \ \mbox{]}
\end{matrix}
$$
qui est adjoint \`a droite [resp. \`a gauche] du foncteur de la remarque pr\'ec\'edente.
\end{listeisansmarge}
\end{remarks}
\medskip

Plus g\'en\'eralement, on peut encore d\'efinir:

\begin{defn}\label{defI102}

Si $F : {\mathcal C} \to {\mathcal D}$ est un foncteur entre deux cat\'egories et $Y$ est un objet de ${\mathcal D}$, la cat\'egorie fl\'ech\'ee \`a droite ${\mathcal C} / Y$ ou $(F \downarrow Y)$ ou encore ${\mathcal C} {}_{F}/ Y$ [resp. \`a gauche $Y\backslash{\mathcal C}$ ou $(Y \downarrow F)$ ou encore $Y\backslash_{F}\,{\mathcal C}$] est d\'efinie de la mani\`ere suivante:
\begin{enumerate}
\item[$\bullet$] Ses objets sont les paires $(X,v)$ constitu\'ees d'un objet $X$ de ${\mathcal C}$ et d'un morphisme de ${\mathcal D}$
$$
F(X) \xrightarrow{ \ v \ } Y \qquad \mbox{[resp.} \quad Y \xrightarrow{ \ v \ } F(X) \ \mbox{]}.
$$
\item[$\bullet$] Ses morphismes
$$
(X_1 , v_1) \longrightarrow (X_2 , v_2)
$$
sont les morphismes de ${\mathcal C}$
$$
X_1 \xrightarrow{ \ u \ } X_2
$$
tels que le triangle
$$
\xymatrix{
F(X_1) \ar[rd]_{v_1} \ar[rr]^{F(u)} &&F(X_2) \ar[ld]^{v_2} \\
&Y
} \qquad\qquad \xymatrix{
&Y \ar[ld]_{\mbox{[resp. \quad} v_1} \ar[rd]^{v_2\qquad\mbox{]}} \\
F(X_1) \ar[rr]^{F(u)} &&F(X_2)
}
$$
soit commutatif dans ${\mathcal D}$.

\end{enumerate}
\end{defn}

\begin{remarksqed}
\begin{listeisansmarge}
\item Si ${\mathcal C}$ est localement petite, les cat\'egories fl\'ech\'ees ${\mathcal C} / Y$ ou $Y \backslash {\mathcal C}$, $Y \in {\rm Ob} ({\mathcal D})$, sont toujours localement petites.

\smallskip

Si ${\mathcal C}$ est une petite cat\'egorie et ${\mathcal D}$ est localement petite, les cat\'egories fl\'ech\'ees ${\mathcal C}/Y$ ou $Y\backslash {\mathcal C}$, $Y \in {\rm Ob} ({\mathcal D})$, sont toujours petites.

\medskip

\item Les cat\'egories fl\'ech\'ees ${\mathcal C}/Y$ ou $Y \backslash {\mathcal C}$ sont munies d'un foncteur d'oubli
$$
\begin{matrix}
\hfill {\mathcal C} / Y &\longrightarrow &{\mathcal C} \, , \hfill \\
(X , F(X) \to Y) &\longmapsto &X \hfill
\end{matrix} \qquad \mbox{[resp.} \quad \begin{matrix}
\hfill Y \backslash {\mathcal C} &\longrightarrow &{\mathcal C} \, , \hfill \\
(X,Y \to F(X)) &\longmapsto &X \hfill
\end{matrix} \ \mbox{]}.
$$

\item Tout morphisme dans la cat\'egorie ${\mathcal D}$
$$
Y \xrightarrow{ \ y \ } Y'
$$
d\'efinit un foncteur entre cat\'egories fl\'ech\'ees \`a droite [resp. \`a gauche]
$$
\begin{matrix}
\hfill {\mathcal C} / Y &\longrightarrow &{\mathcal C}/Y' \hfill \\
(X , F(X) \to Y) &\longmapsto &(X , F(X) \to Y \xrightarrow{y} Y') \hfill
\end{matrix} \qquad \mbox{[resp.} \quad \begin{matrix}
\hfill Y' \backslash {\mathcal C} &\longrightarrow &Y \backslash {\mathcal C} \, , \hfill \\
(X,Y' \to F(X)) &\longmapsto &(X , Y \xrightarrow{y} Y' \to F(X)) \hfill
\end{matrix} \ \mbox{]}
$$
par composition \`a droite [resp. \`a gauche] avec $y$.

\medskip

\item Si le foncteur $F : {\mathcal C} \to {\mathcal D}$ admet un adjoint \`a droite [resp. \`a gauche]
$$
G : {\mathcal D} \longrightarrow {\mathcal C} \, ,
$$
chaque cat\'egorie fl\'ech\'ee \`a droite [resp. \`a gauche] associ\'ee \`a un objet $Y$ de ${\mathcal D}$
$$
{\mathcal C} / Y \qquad \mbox{[resp.} \quad Y \backslash {\mathcal C} \ \mbox{]}
$$
est isomorphe \`a la cat\'egorie relative
$$
{\mathcal C} / G(Y) \qquad \mbox{[resp.} \quad G(Y) \backslash {\mathcal C} \ \mbox{]}.
$$
En particulier, elle poss\`ede pour objet terminal [resp. initial] la paire
$$
(G(Y) , F \circ G (Y) \longrightarrow Y) \qquad \mbox{[resp.} \quad (G(Y) , Y \longrightarrow F \circ G (Y)) \ \mbox{]}
$$
o\`u
$$
F \circ G(Y) \longrightarrow Y \qquad \mbox{[resp.} \quad Y \longrightarrow F \circ G (Y) \ \mbox{]}
$$
est le morphisme qui correspond par adjonction \`a
$$
G(Y) \xrightarrow{ \ {\rm id} \ } G(Y) \, .
$$

\end{listeisansmarge}
\end{remarksqed}

\subsection{Exemples de cat\'egories relatives}\label{subsec1102}

\noindent $\bullet$ {\bf Les cat\'egories d'alg\`ebres:}

\medskip

Si ${\rm An}$ d\'esigne la cat\'egorie des anneaux [resp. des anneaux commutatifs] et $A$ est un objet de cette cat\'egorie, la cat\'egorie relative \`a gauche
$$
A \backslash {\rm An} = {\rm Alg}_A
$$
est appel\'ee la cat\'egorie des $A$-alg\`ebres [resp. des $A$-alg\`ebres commutatives].

\bigskip

\noindent $\bullet$ {\bf Les cat\'egories d'extensions de corps:}

\medskip

Si ${\rm Co}$ d\'esigne la cat\'egorie des corps commutatifs et $K$ est un objet de cette cat\'egorie, la cat\'egorie relative \`a gauche
$$
K \backslash {\rm Co} = {\rm Ext}_K
$$
est appel\'ee la cat\'egorie des extensions du corps $K$.

\medskip

On consid\`ere souvent ses sous-cat\'egories pleines que sont les cat\'egories des extensions alg\'ebriques [resp. alg\'ebriques s\'eparables, resp. finies, resp. finies s\'eparables] du corps $K$.

\bigskip

\noindent $\bullet$ {\bf Les cat\'egories de sch\'emas sur un sch\'ema de base:}

\smallskip

Si $S$ est un objet de la cat\'egorie ${\rm Sch}$ des sch\'emas, la cat\'egorie relative \`a droite
$$
{\rm Sch} / S = {\rm Sch}_S
$$
est appel\'ee la cat\'egorie des sch\'emas sur $S$.

\medskip

Si $S = {\rm Spec} (A)$ est le spectre d'un anneau commutatif $A$, on la note aussi ${\rm Sch}_A$ et l'appelle la cat\'egorie des sch\'emas sur $A$.

\medskip

Comme tout morphisme
$$
S \xrightarrow{ \ s \ } S'
$$
de ${\rm Sch}$ est carrable, il d\'efinit un foncteur de changement de base
$$
\begin{matrix}
\hfill {\rm Sch}_{S'} &\longrightarrow &{\rm Sch}_S \, , \hfill \\
(X \to S') &\longmapsto &(X \times_{S'} S \to S)
\end{matrix}
$$
qui est adjoint \`a droite du foncteur
$$
\begin{matrix}
\hfill {\rm Sch}_{S} &\longrightarrow &{\rm Sch}_{S'} \, , \hfill \\
(X \to S) &\longmapsto &(X \to S \to S') \, .
\end{matrix}
$$

\bigskip

\noindent $\bullet$ {\bf Les cat\'egories de vari\'et\'es alg\'ebriques:}

\medskip

Pour tout sch\'ema $S$, on distingue dans la cat\'egorie ${\rm Sch}_S = {\rm Sch}/S$ des sch\'emas sur $S$ la sous-cat\'egorie pleine ${\rm Sch}_S^{\rm pf}$ constitut\'ee des morphismes $X \to S$ qui sont de pr\'esentation finie.

\medskip

Si $S = {\rm Spec} (K)$ est le spectre d'un corps $K$, la cat\'egorie ${\rm Sch}_S^{\rm pf} = {\rm Sch}_K^{\rm pf} = {\rm Var}_K$ est appel\'ee la cat\'egorie des vari\'et\'es alg\'ebriques sur $K$.

\newpage

\subsection{La cat\'egorie des \'el\'ements d'un pr\'efaisceau}\label{subsec1103}

\smallskip

Une classe tr\`es importante de cat\'egories fl\'ech\'ees est celle que d\'efinissent les foncteurs de Yoneda:

\begin{defn}\label{defI103}

Si ${\mathcal C}$ est une cat\'egorie localement petite et
$$
y : {\mathcal C} \longrightarrow \widehat{\mathcal C} = [{\mathcal C}^{\rm op} , {\rm Ens}]
$$
son foncteur de Yoneda, alors pour tout objet $P$ de $\widehat{\mathcal C}$, la cat\'egorie fl\'ech\'ee \`a droite
$$
{\mathcal C} / P
$$
est aussi not\'ee $\int\!P$ et appel\'ee la ``cat\'egorie des \'el\'ements de $P$''. Ainsi:
\begin{enumerate}
\item[$\bullet$] Ses objets sont les paires $(X,x)$ constitu\'ees d'un objet $X$ de ${\mathcal C}$ et d'un \'el\'ement $x \in P(X)$.
\item[$\bullet$] Ses morphismes
$$
(X_1 , x_1) \longrightarrow (X_2 , x_2)
$$
sont les morphismes de ${\mathcal C}$
$$
X_1 \xrightarrow{ \ u \ } X_2
$$
tels que
$$
P(u)(x_2) = x_1 \, .
$$
\end{enumerate}
\end{defn}

\begin{remarksqed}
\begin{listeisansmarge}
\item Comme ${\mathcal C}$ est localement petite, les cat\'egories ${\mathcal C}/P = \int\!P$ sont toujours localement petites.

\smallskip

Elles sont petites [resp. essentiellement petites] si la cat\'egorie ${\mathcal C}$ est petite [resp. essentiellement petite].

\medskip

\item Pour tout pr\'efaisceau $P$, on a un foncteur d'oubli
$$
\begin{matrix}
\int\!P &\longrightarrow &{\mathcal C} \, , \\
(X,x) &\longmapsto &X \, .
\end{matrix}
$$

\item Tout morphisme de pr\'efaisceaux
$$
P \xrightarrow{ \ p \ } P'
$$
d\'efinit un foncteur
$$
\begin{matrix}
\int\!P &\longrightarrow &\int\!P' \, , \hfill \\
(X,x) &\longmapsto &(X , p_X (x)) \, .
\end{matrix}
$$

\noindent (iv) Si un pr\'efaisceau $P$ est repr\'esent\'e par un objet $S$ de ${\mathcal C}$, la cat\'egorie des \'el\'ements de $P$
$$
{\textstyle \int}P = {\mathcal C} / P
$$
s'identifie \`a la cat\'egorie relative \`a droite
$$
{\mathcal C} / S \, .
$$
\end{listeisansmarge}
\end{remarksqed}

\pagebreak

Tout pr\'efaisceau sur une petite cat\'egorie s'\'ecrit comme une colimite de pr\'efaisceaux repr\'esentables index\'es par la cat\'egorie de ses \'el\'ement consid\'er\'ee comme un carquois:

\begin{prop}\label{propI104}

Soient ${\mathcal C}$ une petite cat\'egorie, $P$ un objet de $\widehat{\mathcal C}$, $\int\!P$ la cat\'egorie de ses \'el\'ements consid\'er\'ee comme un carquois qui indexe le diagramme de $\widehat{\mathcal C}$
$$
\begin{matrix}
\int\!P &\longrightarrow &\widehat{\mathcal C} \, , \hfill \\
(X,x) &\longmapsto &y(X)  = {\rm Hom} (\bullet , X) \, .
\end{matrix}
$$

Alors $P$ est une colimite de ce diagramme de $\widehat{\mathcal C}$, ce qui s'\'ecrit
$$
P = \varinjlim_{(X,x) \in \int\!P} {\rm Hom} (\bullet , X) \, .
$$
\end{prop}

\medskip

\begin{demo}

Soit $Q$ un pr\'efaisceau arbitraire.

\smallskip

Un morphisme de $\widehat{\mathcal C}$
$$
P \xrightarrow{ \ p \ } Q
$$
consiste \`a associer \`a tout objet $X$ de ${\mathcal C}$ et tout \'el\'ement $x$ de $P(X)$ un \'el\'ement
$$
p_X (x) \in Q(X)
$$
ou, ce qui revient au m\^eme, un morphisme de pr\'efaisceaux
$$
{\rm Hom} (\bullet , X) \longrightarrow Q \, ,
$$
de telle fa\c con que, pour tout morphisme de ${\mathcal C}$
$$
X_1 \xrightarrow{ \ u \ } X_2
$$
et tous \'el\'ements $x_1 \in P(X_1)$, $x_2 \in P(X_2)$ v\'erifiant
$$
x_1 = P(u)(x_2) \, ,
$$
on ait
$$
p_{X_1} (x_1) = Q(u) (p_{X_2} (x_2)) \, .
$$

\smallskip

Cette condition \'equivaut \`a demander que le triangle de $\widehat{\mathcal C}$
$$
\xymatrix{
{\rm Hom} (\bullet , X_1) \ar[rd]_{p_{X_1} (x_1)} \ar[rr]^{u \circ \bullet} &&{\rm Hom} (\bullet , X_2) \ar[ld]^{p_{X_2} (x_2)} \\
&Q
}
$$
soit commutatif.

\smallskip

Cela signifie que $P$ est colimite du diagramme des ${\rm Hom} (\bullet , X)$ index\'e par le carquois des $(X,x) \in \int\!P$. 

\end{demo}

\subsection{Extensions de Kan \`a droite ou \`a gauche}\label{subsec1104}

\smallskip

Les cat\'egories fl\'ech\'ees permettent de construire des extensions de Kan \`a droite et \`a gauche, une notion tr\`es g\'en\'erale et tr\`es importante de la th\'eorie des cat\'egories.

\bigskip

Pr\'esentons d'abord de quoi il s'agit:

\begin{defn}\label{defI105}

Soient $F : {\mathcal C} \to {\mathcal D}$ un foncteur entre deux cat\'egories essentiellement petites et ${\mathcal E}$ une cat\'egorie localement petite.

\smallskip

On appelle extension de Kan \`a droite [resp. \`a gauche] le long de $F$ un foncteur
$$
\begin{matrix}
&F_* = {\rm Ran}_F : [{\mathcal C} , {\mathcal E}] &\longrightarrow &[{\mathcal D} , {\mathcal E}] \hfill \\
\mbox{[resp.} \qquad &\hfill F_! = {\rm Lan}_F : [{\mathcal C} , {\mathcal E}] &\longrightarrow &[{\mathcal D} , {\mathcal E}] \ \mbox{]}
\end{matrix}
$$
qui est adjoint \`a droite [resp. \`a gauche] du foncteur de composition avec $F$
$$
\begin{matrix}
F^* : &[{\mathcal D} , {\mathcal E}] &\longrightarrow &[{\mathcal C} , {\mathcal E}] \, , \hfill \\
&\hfill G &\longmapsto &G \circ F \, .
\end{matrix}
$$
\end{defn}

\begin{remarksqed}
\begin{listeisansmarge}
\item Comme ${\mathcal C}$ et ${\mathcal D}$ sont suppos\'ees essentiellement petites, et ${\mathcal E}$ est suppos\'ee localement petite, les cat\'egories $[{\mathcal C}, {\mathcal E}]$ et $[{\mathcal D} , {\mathcal E}]$ sont localement petites.

\medskip

\item Plus g\'en\'eralement, on dit qu'un foncteur
$$
H : {\mathcal C} \longrightarrow {\mathcal E}
$$
admet pour extension de Kan \`a droite [resp. \`a gauche] un foncteur
$$
G : {\mathcal D} \longrightarrow {\mathcal E}
$$
si $G$ repr\'esente le foncteur contravariant [resp. covariant]
$$
\begin{matrix}
&[{\mathcal D} , {\mathcal E}]^{\rm op} &\longrightarrow &{\rm Ens} \, , \hfill \\
&\hfill G' &\longmapsto &{\rm Hom} (G' \circ F , H) \\
{ \ } \\
\mbox{[resp.} \qquad &[{\mathcal D} , {\mathcal E}] &\longrightarrow &{\rm Ens} \, , \hfill \\
&\hfill G' &\longmapsto &{\rm Hom} (H , G' \circ F) \ \mbox{]}.
\end{matrix}
$$

On peut alors noter
$$
G = {\rm Ran}_F H \qquad \mbox{[resp.} \quad G = {\rm Lan}_F H \ \mbox{]}.
$$

\item Ainsi, $F$ a une extension de Kan \`a droite $F_*$ [resp. une extension de Kan \`a gauche $F_!$] si et seulement si tout foncteur $H : {\mathcal C} \to {\mathcal E}$ a une extension de Kan \`a droite ${\rm Ran}_F H$ [resp. \`a gauche ${\rm Lan}_F H$] et qu'il est possible de choisir ces extensions uniform\'ement en $H$.
\end{listeisansmarge}
\end{remarksqed}

\pagebreak

Ayant d\'efini la notion g\'en\'erale d'extensions de Kan, montrons comment les construire gr\^ace aux cat\'egories fl\'ech\'ees:

\begin{prop}\label{propI106}

Soient $F : {\mathcal C} \to {\mathcal D}$ un foncteur entre deux petites cat\'egories et ${\mathcal E}$ une cat\'egorie localement petite.

\smallskip

Supposons que pour tout objet $Y$ de ${\mathcal D}$, le foncteur
$$
\varprojlim_{Y \backslash {\mathcal C}} \qquad \mbox{[resp.} \quad \varinjlim_{{\mathcal C} / Y} \ \mbox{]}
$$
de limite [resp. colimite] des diagrammes de ${\mathcal E}$ index\'es par la cat\'egorie fl\'ech\'ee $Y \backslash {\mathcal C} = (Y \downarrow F)$ [resp. ${\mathcal C} / Y = (F \downarrow Y)$] soit bien d\'efini.

\smallskip

Alors il existe un foncteur d'extension de Kan \`a droite [resp. \`a gauche] le long de $F$
$$
\begin{matrix}
&F_* : [{\mathcal C} , {\mathcal E}] &\longrightarrow &[{\mathcal D} , {\mathcal E}] \hfill \\
\mbox{[resp.} \qquad &\hfill F_! : [{\mathcal C} , {\mathcal E}] &\longrightarrow &[{\mathcal D} , {\mathcal E}] \ \mbox{]}.
\end{matrix}
$$
Il associe \`a tout foncteur 
$$
H : {\mathcal C} \longrightarrow {\mathcal E}
$$
le foncteur
$$
F_* H = G : {\mathcal D} \longrightarrow {\mathcal E} \qquad \mbox{[resp.} \quad F_! H = G : {\mathcal D} \longrightarrow {\mathcal E} \ \mbox{]}
$$
ainsi d\'efini:

\bigskip

\noindent $\left\{ \begin{matrix}
\bullet &\mbox{Pour tout objet $Y$ de ${\mathcal D}$, on a} \hfill \\
{ \ } \\
&\displaystyle G(Y) = \varprojlim_{Y \backslash {\mathcal C}} H \qquad \mbox{[resp.} \quad G(Y) = \varinjlim_{{\mathcal C}/Y} H \ \mbox{]} \\
{ \ } \\
&\mbox{o\`u l'on note encore $H$ le compos\'e du foncteur} \hfill \\
{ \ } \\
&H : {\mathcal C} \longrightarrow {\mathcal E} \\
&\mbox{et du foncteur d'oubli} \hfill \\
{ \ } \\
&\begin{matrix}
\hfill Y \backslash {\mathcal C} &\longrightarrow &{\mathcal C} \, , \hfill \\
(X , Y \to F(X)) &\longmapsto &X \hfill
\end{matrix} \qquad \mbox{[resp.} \quad \begin{matrix}
\hfill {\mathcal C}/Y &\longrightarrow &{\mathcal C} \, , \hfill \\
(X,F(X) \to Y) &\longmapsto &X \hfill
\end{matrix} \ \mbox{]}. \\
{ \ } \\
\bullet &\mbox{Pour toute fl\`eche $v : Y_1 \to Y_2$ de ${\mathcal D}$, la fl\`eche de ${\mathcal E}$} \hfill \\
{ \ } \\
&G(v) : G(Y_1) \longrightarrow G(Y_2) \\
{ \ } \\
&\mbox{est induite par le foncteur canonique de composition avec $v$} \hfill \\
{ \ } \\
&\begin{matrix}
\hfill Y_2 \backslash {\mathcal C} &\longrightarrow &Y_1 \backslash {\mathcal C} \, , \hfill \\
(X , Y_2 \to F(X)) &\longmapsto &(X, Y_1 \xrightarrow{v} Y_2 \to F(X)) \hfill
\end{matrix} \quad \mbox{[resp.} \quad \begin{matrix}
\hfill {\mathcal C}/Y_1 &\longrightarrow &{\mathcal C} /Y_2 \hfill \\
(X,F(X) \to Y_1) &\longmapsto &(X , F(X) \to Y_1 \xrightarrow{v} Y_2) \hfill
\end{matrix} \ \mbox{]}.
\end{matrix}\right.
$
\end{prop}

\bigskip

\begin{demo}

Traitons le cas des extensions de Kan \`a gauche, celui des extensions de Kan \`a droite s'en d\'eduisant en rempla\c cant ${\mathcal C} , {\mathcal C}'$ et ${\mathcal E}$ par leurs oppos\'ees ${\mathcal C}^{\rm op} , {\mathcal C}'^{\rm op}$ et ${\mathcal E}^{\rm op}$.

\smallskip

Pour tout foncteur $G' : {\mathcal D} \to {\mathcal E}$, une transformation naturelle
$$
\alpha : H \longrightarrow G' \circ F
$$
consiste en une famille de fl\`eches de ${\mathcal E}$
$$
\alpha_X : H(X) \longrightarrow G' \circ F(X) \, , \qquad X \in {\rm Ob} ({\mathcal C}) \, ,
$$
telle que, pour toute fl\`eche de ${\mathcal C}$
$$
u : X_1 \longrightarrow X_2 \, ,
$$
le carr\'e
$$
\xymatrix{
H(X_1) \ar[d]_{H(u)} \ar[rr]^{\alpha_{X_1}} &&G' \circ F(X_1) \ar[d]^{G' \circ F(u)} \\
H(X_2) \ar[rr]^{\alpha_{X_2}} &&G' \circ F(X_2)
}
$$
soit commutatif.

\smallskip

Autrement dit, elle consiste en une famille de fl\`eches de ${\mathcal E}$
$$
\varinjlim_{{\mathcal C} / Y} H \longrightarrow G'(Y) \, , \qquad Y \in {\rm Ob} ({\mathcal D}) \, ,
$$
telle que, pour toute fl\`eche de ${\mathcal D}$
$$
v : Y_1 \longrightarrow Y_2 \, ,
$$
le carr\'e induit
$$
\xymatrix{
\displaystyle\varinjlim_{{\mathcal C} / Y_1} H \ar[d] \ar[r] &G'(Y_1) \ar[d]^{G'(v)} \\
\displaystyle\varinjlim_{{\mathcal C}/Y_2} H \ar[r] &G'(Y_2)
}
$$
soit commutatif.

\smallskip

D'o\`u la conclusion. 

\end{demo}

\subsection{Foncteurs adjoints et extensions de Kan}\label{subsec1105}

\smallskip

La proposition qui pr\'ec\`ede entra{\^\i}ne le corollaire suivant:

\begin{cor}\label{corI107}

Soit $F : {\mathcal C} \to {\mathcal D}$ un foncteur entre deux cat\'egories essentiellement petites.

\begin{listeimarge}

\item Supposons que $F$ admet un adjoint {\bf \`a gauche} [resp. {\bf \`a droite}]
$$
G : {\mathcal D} \longrightarrow {\mathcal C} \, .
$$

Alors, pour toute cat\'egorie localement petite ${\mathcal E}$, le foncteur de composition avec $G$
$$
\begin{matrix}
[{\mathcal C} , {\mathcal E}] &\longrightarrow &[{\mathcal D} , {\mathcal E}] \, , \\
\hfill H &\longmapsto &H \circ G \hfill
\end{matrix}
$$
est un foncteur d'extension de Kan {\bf \`a droite} [resp. {\bf \`a gauche}] le long de $F$.

\smallskip

En particulier, $G$ est une extension de Kan {\bf \`a droite} [resp. {\bf \`a gauche}] de ${\rm id}_{\mathcal C} : {\mathcal C} \to {\mathcal C}$ le long de $F : {\mathcal C} \to {\mathcal D}$.

\medskip

\item R\'eciproquement, supposons que ${\rm id}_{\mathcal C} : {\mathcal C} \to {\mathcal C}$ admette une extension de Kan {\bf \`a droite} [resp. {\bf \`a gauche}] le long de $F$
$$
G : {\mathcal D} \longrightarrow {\mathcal C} \, .
$$
Alors elle est munie d'une transformation naturelle
$$
\alpha : G \circ F \longrightarrow {\rm id}_{\mathcal C} \qquad \mbox{[resp.} \quad \alpha : {\rm id}_{\mathcal C} \longrightarrow G \circ F \mbox{]}.
$$
De plus, elle est un adjoint {\bf \`a gauche} [resp. {\bf \`a droite}] de $F$ si et seulement si, pour tous objets $X$ de ${\mathcal C}$ et $Y$ de ${\mathcal D}$, l'application
$$
\begin{matrix}
&{\rm Hom} (Y , F(X)) &\longrightarrow &{\rm Hom} (G(Y) , X) \, , \hfill \\
&\hfill f &\longmapsto &\alpha_X \circ G(f) \hfill \\
{ \ } \\
\mbox{[resp.} \qquad &{\rm Hom} (F(X),Y) &\longrightarrow &{\rm Hom} (X,G(Y)) \, , \hfill \\
&\hfill g &\longmapsto &G(g) \circ \alpha_X \qquad \mbox{]}
\end{matrix}
$$
est bijective.
\end{listeimarge}
\end{cor}

\begin{remark}

Ainsi, pour tout foncteur $F : {\mathcal C} \to {\mathcal D}$ entre deux cat\'egories essentiellement petites, la notion d'extension de Kan {\bf \`a droite} [resp. {\bf \`a gauche}] de ${\rm id}_{\mathcal C} : {\mathcal C} \to {\mathcal C}$ le long de $F$ est plus faible que celle d'adjoint \`a gauche [resp. \`a droite] de $F$.

\end{remark}

\medskip

\begin{demo}
\begin{listeisansmarge}
\item Cela r\'esulte de la proposition \ref{propI106} combin\'ee avec la remarque (iv) apr\`es la d\'efinition \ref{defI102}.

\smallskip

En effet, si $G$ est un adjoint \`a gauche [resp. \`a droite] de $F : {\mathcal C} \to {\mathcal D}$, alors pour tout objet $Y$ de ${\mathcal D}$ la cat\'egorie fl\'ech\'ee $Y \backslash {\mathcal C}$ [resp. ${\mathcal C}/Y$] admet pour objet initial [resp. terminal] l'objet $G(Y)$ muni du morphisme
$$
Y \longrightarrow F \circ G(Y) \qquad \mbox{[resp.} \quad F \circ G(Y) \longrightarrow Y \ \mbox{]}.
$$
Il en r\'esulte que chaque limite $\underset{Y \backslash {\mathcal C}}{\varprojlim} \, H$ [resp. chaque colimite $\underset{{\mathcal C}/Y}{\varinjlim} \, H$] est bien d\'efinie et vaut $H \circ G(Y)$.

\medskip

\item Si $G$ est une extension de Kan \`a droite [resp. \`a gauche] de ${\rm id}_{\mathcal C}$ le long de $F$, il repr\'esente par d\'efinition le foncteur
$$
\begin{matrix}
&[{\mathcal D} , {\mathcal C}]^{\rm op} &\longrightarrow &{\rm Ens} \, , \hfill \\
&\hfill G' &\longmapsto &{\rm Hom} (G' \circ F , {\rm id}_{\mathcal C}) \\
{ \ } \\
\mbox{[resp.} \qquad &[{\mathcal D} , {\mathcal C}] &\longrightarrow &{\rm Ens} \, , \hfill \\
&\hfill G' &\longmapsto &{\rm Hom} ({\rm id}_{\mathcal C} , G' \circ F) \ \mbox{]}.
\end{matrix}
$$
En particulier, la transformation identique de $G$ correspond \`a une transformation naturelle
$$
\alpha : G \circ F \longrightarrow {\rm id}_{\mathcal C} \qquad \mbox{[resp.} \quad \alpha : {\rm id}_{\mathcal C} \longrightarrow G \circ F \ \mbox{]}.
$$
La seconde assertion est la d\'efinition de la propri\'et\'e d'adjonction entre deux foncteurs.
\end{listeisansmarge}
\end{demo}

\subsection{Un crit\`ere g\'en\'eral d'existence des extensions de Kan}\label{subsec1106}

\smallskip

On d\'eduit aussi de la proposition \ref{propI106} la cons\'equence imm\'ediate:

\begin{cor}\label{corI108}

Soient $F : {\mathcal C} \to {\mathcal D}$ un foncteur entre deux cat\'egories essentiellement petites et ${\mathcal E}$ une cat\'egorie localement petite suppos\'ee compl\`ete [resp. cocompl\`ete].

\smallskip

Alors:

\begin{listeimarge}

\item Il existe un foncteur d'extension de Kan \`a droite [resp. \`a gauche] le long de $F$
$$
\begin{matrix}
&F_* : [{\mathcal C} , {\mathcal E}] &\longrightarrow &[{\mathcal D} , {\mathcal E}] \hfill \\
\mbox{[resp.} \qquad &\hfill F_! : [{\mathcal C} , {\mathcal E}] &\longrightarrow &[{\mathcal D} , {\mathcal E}] \ \mbox{]}.
\end{matrix}
$$

\item Si ${\mathcal C}$ est petite, on peut \'ecrire pour tout foncteur $H : {\mathcal C} \to {\mathcal E}$ et tout objet $Y$ de ${\mathcal D}$
$$
\begin{matrix}
&F_* H(Y) &= &\displaystyle \varprojlim_{Y\backslash {\mathcal C}} H \hfill \\
\mbox{[resp.} \qquad &F_! H(Y) &= &\displaystyle \varinjlim_{{\mathcal C}/Y} H \ \mbox{]}.
\end{matrix}
$$
\end{listeimarge}
\end{cor}

\begin{remarksqed}
\begin{listeisansmarge}
\item Ce corollaire s'applique en particulier si ${\mathcal E}$ est la cat\'egorie des ensembles.

\smallskip

Ainsi, pour tout foncteur $\rho : {\mathcal C} \to {\mathcal D}$ entre cat\'egories essentiellement petites, le foncteur de composition
$$
\begin{matrix}
\rho^* &: &\widehat{\mathcal D} &\longrightarrow &\widehat{\mathcal C} \, , \hfill \\
&&\hfill P &\longmapsto &P \circ \rho
\end{matrix}
$$
admet un adjoint \`a droite
$$
\rho_* : \widehat{\mathcal C} \longrightarrow \widehat{\mathcal D}
$$
et un adjoint \`a gauche
$$
\rho_! : \widehat{\mathcal C} \longrightarrow \widehat{\mathcal D} \, .
$$

\item Le corollaire s'applique aussi si ${\mathcal E}$ est une cat\'egorie de pr\'efaisceaux, la cat\'egorie ${\rm Top}$ des espaces topologiques ou n'importe quelle cat\'egorie constitu\'ee des ensembles munis d'un certain type de structure alg\'ebrique.

\medskip

\item En particulier, le corollaire s'applique si ${\mathcal E}$ est la cat\'egorie ${\rm Mod}_A$ des modules sur un anneau $A$. 
\end{listeisansmarge}
\end{remarksqed}

\subsection{Repr\'esentations induites ou co-induites}\label{subsec1107}

\smallskip

Particuli\`erement int\'eressant est le cas du corollaire pr\'ec\'edent o\`u ${\mathcal C}$ et ${\mathcal C}'$ sont deux mono{\"\i}des $M$ et $M'$ (comme par exemple deux groupes) reli\'es par un foncteur c'est-\`a-dire un morphisme
$$
\rho : M \longrightarrow M' \, .
$$

Pour toute cat\'egorie localement petite ${\mathcal E}$, les cat\'egories $[M,{\mathcal E}]$ et $[M' , {\mathcal E}]$ ont pour objets les objets $V$ de ${\mathcal E}$ munis d'actions $\sigma$ ou $\sigma'$ de $M$ ou $M'$ c'est-\`a-dire de morphismes de mono{\"\i}des
$$
M \xrightarrow{ \ \sigma \ } {\rm End} (V) \qquad \mbox{ou} \qquad M' \xrightarrow{ \ \sigma' \ } {\rm End} (V) \, .
$$
Le foncteur
$$
\rho^* : [M' , {\mathcal E}] \longrightarrow [M,{\mathcal E}]
$$
consiste \`a associer \`a tout objet $V$ de ${\mathcal E}$ muni d'une action $\sigma' : M' \to {\rm End} (M)$ l'objet $V$ muni de l'action compos\'ee $\sigma = \sigma' \circ \rho : M \xrightarrow{ \ \rho \ } M' \xrightarrow{ \ \sigma' \ } {\rm End} (V)$.

\smallskip

Le corollaire pr\'ec\'edent s'\'ecrit dans ce cadre:

\begin{cor}\label{corI109}

Soient $\rho : M \to M'$ un morphisme de mono{\"\i}des et ${\mathcal E}$ une cat\'egorie localement petite suppos\'ee compl\`ete [resp. cocompl\`ete].

\smallskip

Alors le foncteur de composition des actions
$$
\begin{matrix}
\rho^* = {\rm Res}_{\rho} &: &[M' , {\mathcal E}] &\longrightarrow &[M,{\mathcal E}] \, , \hfill \\
&&\hfill (V,\sigma') &\longmapsto &(V , \sigma' \circ \rho)
\end{matrix}
$$
admet pour adjoint \`a droite [resp. \`a gauche] le foncteur
$$
\begin{matrix}
&\rho^* = {\rm Ind}_{\rho} &: &[M , {\mathcal E}] &\longrightarrow &[M',{\mathcal E}] \, , \hfill \\
&&&\hfill (V,\sigma) &\longmapsto &\displaystyle \varprojlim_{(\bullet \downarrow \rho)} V \hfill \\
{} \\
\mbox{[resp.} &\rho_! = {\rm coInd}_{\rho} &: &[M , {\mathcal E}] &\longrightarrow &[M',{\mathcal E}] \, , \hfill \\
&&&(V,\sigma) &\longmapsto &\displaystyle \varinjlim_{(\rho \downarrow \bullet)} V \ \mbox{]}
\end{matrix}
$$
o\`u:
\begin{enumerate}
\item[$\bullet$] la cat\'egorie $(\bullet \downarrow \rho)$ [resp. $(\rho \downarrow \bullet)$] admet pour objets les \'el\'ements $m' \in M'$ et pour fl\`eches
$$
m'_1 \longrightarrow m'_2
$$
les \'el\'ements $m \in M$ tels que
$$
\rho (m) \, m'_1 = m'_2 \qquad \mbox{[resp.} \quad m'_1 = m'_2 \, \rho (m) \ \mbox{]} ,
$$
\item[$\bullet$] la lettre $V$ d\'esigne encore le foncteur qui associe \`a tout objet $m'$ l'objet $V$ de ${\mathcal E}$ et \`a toute fl\`eche $m$ le morphisme
$$
\rho (m) : V \longrightarrow V \, .
$$
\end{enumerate}
\end{cor}

\begin{remarksqed}
\begin{listeisansmarge}
\item Les foncteurs ${\rm Res}_{\rho}$, ${\rm Ind}_{\rho}$ et ${\rm coInd}_{\rho}$ sont appel\'es les foncteurs de restriction, d'induction et de co-induction.

\medskip

\item Si le mono{\"\i}de $M'$ est trivial, le foncteur
$$
{\rm Ind}_{\rho} : [M,{\mathcal E}] \longrightarrow {\mathcal E} \qquad \mbox{[resp.} \quad {\rm coInd}_{\rho} : [M,{\mathcal E}] \longrightarrow {\mathcal E} \ \mbox{]}
$$
associe \`a tout objet $V$ de ${\mathcal E}$ muni d'une action $\sigma : M \to {\rm End} (V)$ de $M$ le sous-objet [resp. l'objet quotient]
$$
V_M \qquad \mbox{[resp.} \quad M \backslash V \ \mbox{]}
$$
des invariants [resp. des co-invariants] de cette action. 
\end{listeisansmarge}
\end{remarksqed}

\bigskip

Si ${\mathcal E} = {\rm Mod}_A$ est la cat\'egorie des modules sur un anneau $A$, les cat\'egories $[M,{\mathcal E}]$ et $[M' , {\mathcal E}]$ sont celles des repr\'esentations $A$-lin\'eaires des mono{\"\i}des $M$ et $M'$.

\smallskip

Le foncteur de composition avec un morphisme $\rho : M \to M'$
$$
{\rm Res}_{\rho} : [M' , {\mathcal E}] \longrightarrow [M,{\mathcal E}]
$$
est le foncteur de restriction des repr\'esentations, et ses adjoints \`a droite et \`a gauche
$$
\begin{matrix}
\hfill{\rm Ind}_{\rho} &: &[M,{\mathcal E}] &\longrightarrow &[M',{\mathcal E}] \, , \\
{\rm coInd}_{\rho} &: &[M,{\mathcal E}] &\longrightarrow &[M' , {\mathcal E}] \hfill
\end{matrix}
$$
associent \`a toute repr\'esentation $(V,\sigma)$ de $M$ sa repr\'esentation induite ${\rm Ind}_{\rho} (V,\sigma)$ et sa repr\'esentation co-induite ${\rm coInd} (V,\sigma)$.

\bigskip

Dans le cadre de la th\'eorie des groupes, il est fr\'equent que les repr\'esentations induites et co-induites co{\"\i}ncident:

\begin{prop}\label{propI1010}

Soient $\rho : G \to G'$ un morphisme de groupes et $A$ un anneau.

\smallskip

Supposons que le sous-groupe image de $\rho$ dans $G'$ est d'indice fini et que son noyau $H$ est fini d'un cardinal $\vert H \vert$ qui est inversible dans l'anneau $A$.

\smallskip

Alors les deux foncteurs
$$
\begin{matrix}
\hfill{\rm Ind}_{\rho} &: &[G,{\rm Mod}_A] &\longrightarrow &[G',{\rm Mod}_A] \, , \\
{\rm coInd}_{\rho} &: &[G,{\rm Mod}_A] &\longrightarrow &[G' , {\rm Mod}_A] \hfill
\end{matrix}
$$
sont canoniquement isomorphes.
\end{prop}

\begin{remarks}
\begin{listeisansmarge}
\item Autrement dit, les foncteurs ${\rm Res}_{\rho}$ et ${\rm Ind}_{\rho} \cong {\rm coInd}_{\rho}$ sont adjoints \`a la fois \`a gauche et \`a droite l'un de l'autre.

\smallskip

En particulier, ils respectent les limites et les colimites arbitraires.

\medskip

\item Cette proposition s'applique en particulier lorsque $A$ est un corps dont la caract\'eristique ne divise pas $\vert H \vert$.
\end{listeisansmarge}
\end{remarks}

\medskip

\begin{demopropsansqed}

Il s'agit de d\'emontrer que ${\rm Ind}_{\rho}$ est aussi un adjoint \`a gauche de ${\rm Res}_{\rho}$.

\smallskip

Comme le morphisme $\rho$ se factorise en
$$
G \twoheadrightarrow {\rm Im} (\rho) \hookrightarrow G' \, ,
$$
il suffit de traiter s\'epar\'ement le cas o\`u $\rho$ est injectif et celui o\`u il est surjectif.

\smallskip

Voyons d'abord le cas o\`u $G$ est un sous-groupe de $G'$ d'indice fini.

\smallskip

Pour toute repr\'esentation $R$-lin\'eaire $(V,\sigma)$ de $G$, ${\rm Ind} (V,\sigma)$ est constitu\'e des familles $(v_{g'})$ d'\'el\'ements $v_{g'} \in V$ index\'ees par les $g' \in G'$ et telles que
$$
v_{gg'} = \sigma (g) (v_{g'}) \, , \qquad \forall \, g' \in G' \, , \ \forall \, g \in G \, .
$$

L'action de $G'$ sur ${\rm Ind} (V,\sigma)$ est donn\'ee par
$$
h \cdot (v_{g'})_{g' \in G'} = (v_{g'h})_{g' \in G'} \, , \qquad \forall \, h \in G' \, .
$$

D'autre part, ${\rm coInd} (V,\sigma)$ est le plus grand quotient de la somme directe
$$
\bigoplus_{g' \in G'} V
$$
sur lequel les actions de $G$ par
$$
g \cdot \left( \bigoplus_{g' \in G'} v_{g'} \right) = \bigoplus_{g' \in G'} v_{gg'}
$$
et
$$
g \cdot \left( \bigoplus_{g' \in G'} v_{g'} \right) = \bigoplus_{g' \in G'} \sigma (g) (v_{g'})
$$
co{\"\i}ncident.

\smallskip

L'action de $G'$ sur ${\rm coInd} (V,\sigma)$ est d\'eduite de celle sur les
$$
\bigoplus_{g' \in G'} v_{g'}
$$
d\'efinie par
$$
h \cdot \left( \bigoplus_{g' \in G'} v_{g'} \right) = \bigoplus_{g' \in G'} v_{g'h} \, , \qquad \forall \, h \in G' \, .
$$

Choisissons une famille d'\'el\'ements $g'_1 , \cdots , g'_n$ de $G'$ qui repr\'esentent l'ensemble $G \backslash G'$ des classes \`a gauche modulo $G$.

\smallskip

L'application $A$-lin\'eaire compos\'ee
$$
{\rm Ind} (V,\sigma) \xhookrightarrow{ \ \ } \prod_{g' \in G'} V \longrightarrow {\rm coInd} (V,\sigma) \, ,
$$
$$
(v_{g'})_{g' \in G'} \longmapsto v_{g'_1} \oplus \cdots \oplus v_{g'_n}
$$
ne d\'epend pas du choix des repr\'esentants $g'_1 , \ldots , g'_n$.

\smallskip

Pour cette raison, elle respecte les actions de $G'$.

\smallskip

De plus, c'est un isomorphisme de $A$-modules.

\smallskip

Enfin, la famille des isomorphismes
$$
{\rm Ind} (V,\sigma) \xrightarrow{ \ \sim \ } {\rm coInd} (V,\sigma)
$$
est compatible par construction avec les morphismes de repr\'esentations $A$-lin\'eaires de $G$
$$
(V_1 , \sigma_1) \longrightarrow (V_2 , \sigma_2)
$$
au sens que les carr\'es induits
$$
\xymatrix{
{\rm Ind} (V_1 , \sigma_1) \ar[d] \ar[r]^-{\sim} &{\rm coInd} (V_1 , \sigma_1) \ar[d] \\
{\rm Ind} (V_2 , \sigma_2) \ar[r]^-{\sim} &{\rm coInd} (V_2 , \sigma_2)
}
$$
sont commutatifs.

\smallskip

Cela conclut le cas o\`u $\rho$ est injectif.

\smallskip

Le cas o\`u $\rho$ est surjectif r\'esulte du lemme suivant dont la d\'emonstration est laiss\'ee au lecteur:
\end{demopropsansqed}

\begin{lem}\label{lemI1011}

Soient $\rho : G \twoheadrightarrow G'$ un \'epimorphisme de groupes et $A$ un anneau.

\smallskip

On suppose que le noyau $H$ de $\rho$ est fini d'un cardinal $\vert H \vert$ qui est inversible dans $A$.

\smallskip

Alors:

\begin{listeimarge}

\item Pour toute repr\'esentation $A$-lin\'eaire $(V,\sigma)$ de $G$, le morphisme
$$
\begin{matrix}
p_V &: &V &\longrightarrow &V \, , \hfill \\
&&v &\longmapsto &\displaystyle \frac1{\vert H \vert} \cdot \displaystyle \sum_{h \in H} \sigma (h)(v)
\end{matrix}
$$
respecte l'action de $G$ sur $V$ et il est idempotent au sens que
$$
p_V \circ p_V = p_V
$$
et donc 
$$
V = {\rm Ker} (p_V) \oplus {\rm Ker} ({\rm Id} - p_V) \, .
$$

\item De plus, on a pour tout $v \in V$
$$
p_V (v) = v \Leftrightarrow \sigma (h)(v) = v \, , \qquad \forall \, h \in H \, ,
$$
et donc
$$
{\rm Ind}_{\rho} (V,\sigma) = V_H = {\rm Ker} ({\rm Id} - p_V) \, .
$$
\end{listeimarge}
\end{lem}

\medskip

\noindent {\bf Fin de la d\'emonstration de la proposition \ref{propI1010}.:}

\smallskip

Soit $(V,\sigma)$ une repr\'esentation $A$-lin\'eaire de $G$.

\smallskip

Comme ${\rm Ind}_{\rho} (V,\sigma) = V_H$, on sait d\'ej\`a que, pour toute repr\'esentation $A$-lin\'eaire $(V',\sigma')$ de $G'$, se donner un morphisme de $[G,{\rm Mod}_A]$
$$
(V' , \sigma' \circ \rho) \longrightarrow (V,\sigma)
$$
\'equivaut \`a se donner un morphisme de $[G',{\rm Mod}_A]$
$$
(V',\sigma') \longrightarrow V_H \, .
$$

Pour conclure, il suffit de montrer que, de m\^eme, se donner un morphisme de $[G , {\rm Mod}_A]$
$$
(V,\sigma) \longrightarrow (V' , \sigma' \circ \rho)
$$
\'equivaut \`a se donner un morphisme de $[G' , {\rm Mod}_A]$
$$
V_H \longrightarrow (V',\sigma') \, .
$$

Consid\'erons la d\'ecomposition de $V$ en somme directe
$$
V = {\rm Ker} (p_V) \oplus V_H
$$
respect\'ee par l'action de $G$.

\smallskip

Comme l'action de $G$ sur $V_H$ se factorise \`a travers son quotient $G'$, un morphisme de $[G , {\rm Mod}_A]$
$$
V_H \longrightarrow (V' , \sigma' \circ \rho)
$$
est la m\^eme chose qu'un morphisme de $[G' , {\rm Mod}_A]$
$$
V_H \longrightarrow (V',\sigma') \, .
$$

D'autre part, comme l'action de $H$ sur $V'$ via $\sigma' \circ \rho$ est triviale, on a dans $V'$ muni de l'action de $G$ par $\sigma' \circ \rho$
$$
{\rm Ker} (p_{V'}) = 0 \, .
$$
Il en r\'esulte que le seul morphisme de $[G,{\rm Mod}_A]$
$$
{\rm Ker} (p_V) \longrightarrow (V' , \sigma' \circ \rho)
$$
est $0$.

\smallskip

D'o\`u la conclusion. \hfill $\Box$


%% file: Chapitre2_num.tex







\vglue 15mm

\chapter{Topologies de Grothendieck, sites et faisceaux sur un site}\label{chap2}

\section{La notion de crible}\label{sec21}

\subsection{Cribles sur un objet d'une cat\'egorie}\label{subsec211}

\begin{defn}\label{defII11}

Soient ${\mathcal C}$ une cat\'egorie et $X$ un objet de ${\mathcal C}$.

\smallskip

Un crible sur $X$ est une collection $S$ de fl\`eches de ${\mathcal C}$
$$
U \longrightarrow X
$$
telle que, pour tout triangle commutatif de ${\mathcal C}$
$$
\xymatrix{
U' \ar[r]^f \ar[rd]_{p'} &U \ar[d]^p \\
&X
}
$$
$p' = p \circ f$ est dans $S$ si $p$ est dans $S$.

\end{defn}

\begin{remarksqed}
\begin{listeisansmarge}
\item La collection vide est un crible de n'importe quel objet $X$, et la collection de toutes les fl\`eches $U \to X$ de but $X$ est un crible de $X$, appel\'e son crible maximal.

\medskip

\item Toute intersection ou toute r\'eunion de cribles sur un objet $X$ de ${\mathcal C}$ est un crible sur $X$.

\medskip

\item Par cons\'equent, toute fl\`eche $U \xrightarrow{ \ p \ } X$ [resp. toute famille de fl\`eches $U_i \xrightarrow{ \ p_i \ } X$, $i \in I$] vers un objet $X$ d'une cat\'egorie ${\mathcal C}$ engendre un crible sur $X$.

\smallskip

Il est constitu\'e des fl\`eches $U' \xrightarrow{ \ p' \ } X$ de ${\mathcal C}$ telles qu'il existe une fl\`eche $f : U' \to U$ [resp. un indice $i \in I$ et une fl\`eche $f : U' \to U_i$] avec
$$
p' = p \circ f \qquad \mbox{[resp.} \quad p' = p_i \circ f \ \mbox{]}.
$$

\item Si ${\mathcal C}$ est une cat\'egorie localement petite, un crible sur $X$ n'est pas autre chose qu'un sous-objet
$$
S \xhookrightarrow{ \ \ \ } y(X) = {\rm Hom} (\bullet , X)
$$
du pr\'efaisceau $y(X) = {\rm Hom} (\bullet , X)$ dans la cat\'egorie $\widehat{\mathcal C} = [{\mathcal C}^{\rm op} , {\rm Ens}]$ des pr\'efaisceaux sur ${\mathcal C}$.

\smallskip

En effet, la condition de stabilit\'e de $S$ par composition \`a droite avec n'importe quelle fl\`eche de ${\mathcal C}$ correspond \`a la propri\'et\'e de $S$ plong\'e dans ${\rm Hom} (\bullet , X)$ d'\^etre un pr\'efaisceau.

\medskip

\item Si $F : {\mathcal C} \to {\mathcal C}'$ est une \'equivalence faible de cat\'egories et $X'$ est un objet de ${\mathcal C}'$ isomorphe \`a l'image $F(X)$ d'un objet $X$ de ${\mathcal C}$, alors les cribles sur $X'$ dans ${\mathcal C}'$ sont en bijection avec les cribles sur $X$ dans ${\mathcal C}$.

\medskip

\item Il r\'esulte de la remarque pr\'ec\'edente que si $X$ est un objet d'une cat\'egorie ${\mathcal C}$ essentiellement petite, alors les cribles sur $X$ dans ${\mathcal C}$ forment un ensemble. 
\end{listeisansmarge}
\end{remarksqed}

\subsection{Exemples de cribles}\label{subsec212}

\noindent $\bullet$ {\bf Les cribles d'un mono{\"\i}de:} 

\smallskip

Si $M$ est un mono{\"\i}de consid\'er\'e comme une cat\'egorie, les cribles de son unique objet sont les parties $S$ de $M$ stables par multiplication \`a droite au sens que
$$
m \in S \, , \ m' \in M \Rightarrow m \cdot m' \in S \, .
$$

On remarque qu'un tel mono{\"\i}de $M$ est un groupe si et seulement si son seul crible non vide est le crible maximal $M$.

\medskip

\noindent $\bullet$ {\bf Les cribles d'un mono{\"\i}de commutatif:}

\smallskip

Si $M$ est un mono{\"\i}de commutatif, ses cribles sont les parties de $M$ stables par multiplication par n'importe quel \'el\'ement de $M$.

\smallskip

Par exemple, le mono{\"\i}de ${\mathbb N}^{\times}$ des entiers $\geq 1$ muni de la multiplication admet pour crible l'ensemble des entiers qui ne sont pas premiers. C'est le crible d'Eratosth\`ene.

\medskip

\noindent $\bullet$ {\bf Les cribles d'un ensemble ordonn\'e:}

\smallskip

Si $(O , \leq)$ est un ensemble ordonn\'e consid\'er\'e comme une cat\'egorie, un crible sur un \'el\'ement $u$ de $O$ peut \^etre vu comme une partie $S$ de $O$ telle que
$$
\left\{\begin{matrix}
\bullet &v \leq u \, , \ \forall \, v \in S \, , \hfill \\
\bullet &v \in S \, , \ v' \leq v \Rightarrow v' \in S \, .
\end{matrix} \right.
$$

Par exemple, dans $({\mathbb R} , \leq)$, les cribles non vides sur un \'el\'ement $x$ sont les intervalles de la forme $] - \infty , a \, [$ et $]-\infty , a \, ]$ avec $a \leq x$.

\medskip

\noindent $\bullet$ {\bf Les cribles d'un espace topologique:}

\smallskip

Si $X$ est un espace topologique et $O(X)$ l'ensemble ordonn\'e de ses ouverts vu comme une cat\'egorie, un crible sur un \'el\'ement $U$ de $O(X)$ est une partie $S$ constitu\'ee d'ouverts $V \subset U$, telle que
$$
V \in S \, , \ V' \subset V \Rightarrow V' \in S \, .
$$

Le crible sur un ouvert $U$ engendr\'e par une famille d'ouverts $U_i \subset U$, $i \in I$, est l'ensemble $S$ des ouverts $V \subset U$ qui sont contenus dans l'un au moins des $U_i$.

\smallskip

Si $(X,d)$ est un espace m\'etrique, alors pour tout nombre r\'eel $\varepsilon > 0$, les ouverts $V$ contenus dans un ouvert $U$ et de diam\`etre $\leq \varepsilon$ au sens que
$$
d(x_1 , x_2) \leq \varepsilon \, , \quad \forall \, x_1 , x_2 \in V \, ,
$$
forment un crible de $U$.

\subsection{Images r\'eciproques des cribles}\label{subsec213}

\begin{defn}\label{defII12}

Soient ${\mathcal C}$ une cat\'egorie et $h : X \to Y$ un morphisme de ${\mathcal C}$.

\smallskip

Pour tout crible $S$ sur $Y$, son image r\'eciproque par $h$ est le crible sur $X$
$$
h^* S = \{ U \xrightarrow{ \ p \ } X \mid h \circ p \in S \} \, .
$$
\end{defn}

\begin{remarksqed}
\begin{listeisansmarge}
\item Si ${\mathcal C}$ est une cat\'egorie localement petite et un crible $S$ sur $Y$ est vu comme un sous-objet
$$
S \xhookrightarrow{ \ \ } y(Y) = {\rm Hom} (\bullet , Y) \quad \mbox{dans} \quad \widehat{\mathcal C} \, ,
$$
alors son image r\'eciproque $h^* S$ par un morphisme $h : X \to Y$ est le sous-objet
$$
S \times_{y(Y)} y(X) \xhookrightarrow{ \ \ } y(X) = {\rm Hom} (\bullet , X)
$$
qui s'en d\'eduit par le changement de base par le morphisme
$$
y(h) : y(X) \longrightarrow y(Y)
$$
de $\widehat{\mathcal C}$.

\medskip

\item Pour tous morphismes $X \xrightarrow{ \ h \ } Y \xrightarrow{ \ g \ } Z$ d'une cat\'egorie ${\mathcal C}$ et tout crible $S$ sur $Z$, on a
$$
h^* (g^* S) = (g \circ h)^* S \, .
$$

\item Si ${\mathcal C}$ est une cat\'egorie essentiellement petite, associer \`a tout objet $X$ de ${\mathcal C}$ l'ensemble $\Omega (X)$ de ses cribles d\'efinit un foncteur contravariant
$$
\begin{matrix}
\Omega &: &\hfill {\mathcal C}^{\rm op} &\longrightarrow &{\rm Ens} \, , \hfill \\
&&\hfill X &\longmapsto &\Omega (X) \, , \hfill \\
&&(X \xrightarrow{ \ h \ } Y) &\longmapsto &(h^* : \Omega (Y) \longrightarrow \Omega (X)) \, .
\end{matrix}
$$

Autrement dit, les cribles sur les objets $X$ de ${\mathcal C}$ sont les \'el\'ements d'un pr\'efaisceau $\Omega$.
\end{listeisansmarge}
\end{remarksqed}

\section{Topologies et sites de Grothendieck}\label{sec22}

\subsection{La notion de topologie de Grothendieck et celle de site}\label{subsec221}

\begin{defn}\label{defII21}

Soit ${\mathcal C}$ une cat\'egorie essentiellement petite.

\smallskip

Une topologie (de Grothendieck) sur ${\mathcal C}$ est une application
$$
J : X \longmapsto J(X)
$$
qui associe \`a tout objet $X$ de ${\mathcal C}$ un sous-ensemble
$$
J(X) \subset \Omega (X)
$$
de l'ensemble $\Omega(X)$ des cribles sur $X$ dans ${\mathcal C}$, et qui satisfait les trois conditions suivantes:

\bigskip

\noindent $\left\{ \begin{matrix}
\mbox{\rm (M)} &\mbox{\bf Maximalit\'e:} \hfill \\
&\mbox{Pour tout objet $X$ de ${\mathcal C}$, $J(X)$ contient le ``crible maximal'' sur $X$,} \hfill \\
&\mbox{constitu\'e de toutes les fl\`eches $U \to X$ de but $X$.} \hfill \\
{ \ } \\
\mbox{\rm (S)} &\mbox{\bf Stabilit\'e:} \hfill \\
&\mbox{Pour tout morphisme $X \xrightarrow{ \ h \ } Y$ de ${\mathcal C}$, l'application} \hfill \\
{ \ } \\
&h^* : \Omega (Y) \longrightarrow \Omega (X) \\
{ \ } \\
&\mbox{envoie le sous-ensemble $J(Y)$ de $\Omega(Y)$ dans le sous-ensemble $J(X)$ de $\Omega (X)$.} \hfill \\
{ \ } \\
\mbox{\rm (T)} &\mbox{\bf Transitivit\'e:} \hfill \\
&\mbox{Pour tout objet $X$ de ${\mathcal C}$ et tout $S \in J(X)$, tout crible $S'$ sur $X$ tel que} \hfill \\
{ \ } \\
&h^* S' \in J(U) \, , \quad \forall \ (U \xrightarrow{ \ h \ } X) \in S \, , \\
{ \ } \\
&\mbox{est n\'ecessairement \'el\'ement de $J(X)$.} \hfill
\end{matrix} \right.
$
\end{defn}

\begin{remarksqed}
\begin{listeisansmarge}
\item Si $J$ est une topologie sur ${\mathcal C}$, les cribles \'el\'ements d'un ensemble $J(X)$ sont appel\'es les cribles $J$-couvrants de $X$ ou, s'il n'y a pas ambigu{\"\i}t\'e, les cribles couvrants.

\medskip

\item Une famille de fl\`eches ${\mathcal C}$
$$
U_i \longrightarrow X \,, \quad i \in I \, ,
$$
est appel\'ee ``$J$-couvrante'' ou ``couvrante'' si le crible sur $X$ qu'elle engendre est couvrant.

\smallskip

Par exemple, l'axiome de maximalit\'e (M) signifie que, pour tout objet $X$ de ${\mathcal C}$, le morphisme
$$
X \xrightarrow{ \ {\rm id}_X \ } X
$$
est couvrant pour n'importe quelle topologie.

\medskip

\item Dans la d\'efinition d'une topologie $J$ sur une cat\'egorie ${\mathcal C}$, on n'exclut pas que certains objets de ${\mathcal C}$ puissent admettre parmi leurs cribles $J$-couvrants le crible vide.

\medskip

\item Si ${\mathcal C}$ est une cat\'egorie essentiellement petite, les topologies $J$ sur ${\mathcal C}$ forment un ensemble ordonn\'e par l'inclusion.

\smallskip

Quand deux topologies $J_1$ et $J_2$ satisfont la relation $J_1 \subset J_2$, on dit que $J_2$ est plus fine que $J_1$ ou que $J_1$ est moins fine que $J_2$.

\medskip

\item Toute famille $(J_i)_{i \in I}$ de topologies sur ${\mathcal C}$ a une intersection
$$
\bigwedge_{i \in I} J_i : X \longmapsto \bigcap_{i \in I} J_i (X)
$$
qui est encore une topologie sur ${\mathcal C}$.

\medskip

\item Par cons\'equent, toute famille de cribles d'objets de ${\mathcal C}$ engendre une topologie, qui est l'intersection de toutes les topologies contenant cette famille.

\smallskip

En particulier, toute famille $(J_i)_{i \in I}$ de topologies de ${\mathcal C}$ engendre une topologie
$$
\bigvee_{i \in I} J_i
$$
caract\'eris\'ee par la propri\'et\'e que, pour toute topologie $J'$ de ${\mathcal C}$
$$
J' \supset \bigvee_{i \in I} J_i \Longleftrightarrow J' \supset J_i \, , \quad \forall \, i \in I \, .
$$

\end{listeisansmarge}
\end{remarksqed}

\medskip

La donn\'ee d'une cat\'egorie essentiellement petite et d'une topologie de Grothendieck sur cette cat\'egorie porte un nom:

\begin{defn}\label{defII22}

On appelle site une paire constitu\'ee d'une cat\'egorie essentiellement petite ${\mathcal C}$ et d'une topologie $J$ sur ${\mathcal C}$.
\end{defn}

\begin{remarkqed}

Si $({\mathcal C} , J)$ est un site, alors pour tout objet $S$ de ${\mathcal C}$, $J$ induit une topologie sur la cat\'egorie relative ${\mathcal C}/S$: on d\'ecide qu'une famille de morphismes $(U_i \to S) \to (X \to S)$, $i \in I$, vers un objet $X \to S$ de ${\mathcal C}/S$ est couvrante si la famille de morphismes $U_i \to X$, $i \in I$, de ${\mathcal C}$ est $J$-couvrante.

\smallskip

Le site $({\mathcal C}/S,J)$ est appel\'e le localis\'e en $S$ de $({\mathcal C},J)$.

\end{remarkqed}

\subsection{Premiers exemples de topologies de Grothendieck}\label{subsec222}

\medskip

\noindent $\bullet$ {\bf La topologie usuelle des espaces topologiques:}

\smallskip

Soient $X$ un espace topologique et $O(X)$ l'ensemble de ses ouverts ordonn\'e par l'inclusion, consid\'er\'e comme une cat\'egorie.

\smallskip

On d\'efinit une topologie de Grothendieck sur $O(X)$ en d\'ecidant qu'une famille de fl\`eches
$$
U_i \subset U \, , \quad i \in I \, ,
$$
est couvrante si et seulement si
$$
\bigcup_{i \in I} U_i = U \, .
$$

C'est la notion usuelle de recouvrement d'un ouvert par des ouverts plus petits.

\medskip

\noindent $\bullet$ {\bf La topologie canonique des treillis:}

\smallskip

L'exemple des espaces topologiques se g\'en\'eralise au cadre plus vaste des treillis:

\begin{defn}\label{defII23}

Un treillis est un ensemble $O$ muni d'une relation d'ordre partiel $\leq$ telle que:
\begin{enumerate}
\item[$\bullet$] toute famille $(u_i)_{i \in I}$ d'\'el\'ements de $O$ admet un supremum $\underset{i \in I}{\bigvee} \, u_i$ caract\'eris\'e par la propri\'et\'e que, pour tout $u' \in O$,
$$
\bigvee_{i \in I} u_i \leq u' \Longleftrightarrow u_i \leq u' \, , \quad \forall \, i \in I \, ,
$$
\item[$\bullet$] toute famille finie $u_1 , \cdots , u_n$ d'\'el\'ements de $O$ admet un infimum $u_1 \wedge \cdots \wedge u_n$ caract\'eris\'e par la propri\'et\'e que, pour tout $u' \in O$,
$$
u' \leq u_1 \wedge \cdots \wedge u_n \Longleftrightarrow u' \leq u_i \, , \quad \forall \, i \in I \, ,
$$
\item[$\bullet$] pour tous \'el\'ements $u$ et $u_i$, $i \in I$, de $O$, on a
$$
u \wedge \left(\bigvee_{i \in I} u_i \right) = \bigvee_{i \in I} (u \wedge u_i) \, .
$$
\end{enumerate}
\end{defn}

\begin{remarkqed}

Pour tout espace topologique $X$, l'ensemble $O(X)$ de ses ouverts ordonn\'e par l'inclusion est un treillis au sens de cette d\'efinition.
 
\end{remarkqed}

\medskip

La topologie canonique d'un treillis $(O,\leq)$ consid\'er\'e comme une cat\'egorie est d\'efinie en d\'ecidant qu'une famille de fl\`eches
$$
u_i \leq u \, , \quad i \in I \, ,
$$
est couvrante si et seulement si
$$
\bigvee_{i \in I} u_i = u \, .
$$

Comme annonc\'e, on retrouve la topologie usuelle d'un espace topologique $X$ dans le cas o\`u $O$ est le treillis $O(X)$ de ses ouverts.

\medskip

\noindent $\bullet$ {\bf La topologie induite par un sous-espace d'un espace topologique:}

\smallskip

Si $X$ est un espace topologique, $O(X)$ le treillis de ses ouverts et $Y \hookrightarrow X$ un sous-ensemble de $X$, on d\'efinit une topologie sur la cat\'egorie $O(X)$ en d\'ecidant qu'une famille de fl\`eches
$$
U_i \subset U \, , \quad i \in I \, ,
$$
est couvrante si et seulement si
$$
\bigcup_{i \in I} \ (U_i \cap Y) = U \cap Y \, .
$$

Si $Y=X$, on retrouve bien s\^ur la topologie usuelle de $O(X)$.

\medskip

\noindent $\bullet$ {\bf La topologie induite par un \'el\'ement d'un treillis:}

\smallskip

Si $(O,\leq)$ est un treillis et $y$ un \'el\'ement de $O$, on d\'efinit une topologie sur la cat\'egorie $O$ en d\'ecidant qu'une famille de fl\`eches
$$
u_i \leq u \, , \quad i \in I \, ,
$$
est couvrante si et seulement si
$$
\bigvee_{i \in I} (u_i \wedge y) =  u \wedge y \, .
$$

Cet exemple recoupe le pr\'ec\'edent lorsque $O$ est le treillis $O(X)$ des ouverts d'un espace topologique $X$ et $y$ est un ouvert $Y$ de $X$.

\medskip

\noindent $\bullet$ {\bf La topologie naturelle d'un espace muni d'une mesure:}

\smallskip

Soit $X$ un ensemble muni d'une mesure $\mu$.

\smallskip

Soit $O$ l'ensemble des parties mesurables de $X$, ordonn\'e par la relation d'inclusion.

\smallskip

On d\'efinit une topologie sur l'ensemble ordonn\'e $O$ en d\'ecidant qu'une famille de fl\`eches
$$
Y_i \subset Y \, , \quad i \in I \, ,
$$
est couvrante si elle contient une sous-famille d\'enombrable
$$
Y_{i_n} \subset Y \, , \quad n \in {\mathbb N} \, ,
$$
telle que
$$
\mu \left( Y - \bigcup_{n \in {\mathbb N}} Y_{i_n} \right) = 0 \, .
$$

\newpage

\noindent $\bullet$ {\bf La topologie de la densit\'e dans un espace topologique:}

\smallskip

Soient $X$ un espace topologique et $O(X)$ le treillis de ses ouverts.

\smallskip

La topologie de la densit\'e sur $O(X)$ consid\'er\'e comme une cat\'egorie est d\'efinie en d\'ecidant qu'une famille de fl\`eches
$$
U_i \subset U \, , \quad i \in I \, ,
$$
est couvrante si et seulement si
$$
\bigcup_{i \in I} U_i \quad \mbox{est dense dans} \quad U \, .
$$
Autrement dit, une famille
$$
U_i \subset U \, , \quad i \in I \, ,
$$
est couvrante si et seulement si, pour tout ouvert
$$
V \subset U
$$
qui est non vide, il existe $i \in I$ tel que l'intersection
$$
V \cap U_i
$$
ne soit pas vide.

\medskip

\noindent $\bullet$ {\bf La topologie de la densit\'e dans une cat\'egorie:}

\smallskip

Soit ${\mathcal C}$ une cat\'egorie essentiellement petite.

\smallskip

La ``topologie de la densit\'e'' de ${\mathcal C}$ est d\'efinie en d\'ecidant qu'un crible $S$ sur un objet $U$ de ${\mathcal C}$ est couvrant si, pour toute fl\`eche de ${\mathcal C}$ de but $U$
$$
h : V \longrightarrow U \, ,
$$
le crible $h^*S$ n'est pas vide, autrement dit il existe une fl\`eche $U' \xrightarrow{ \ g \ } V$ de ${\mathcal C}$ telle que
$$
h \circ g : U' \longrightarrow U
$$
est \'el\'ement de $S$.

\smallskip

Cet exemple g\'en\'eralise celui de la topologie de la densit\'e d'un espace topologique $X$: il suffit en effet de prendre pour ${\mathcal C}$ la sous-cat\'egorie pleine de $O(X)$ constitu\'ee de tous les ouverts non vides.

\medskip

\noindent $\bullet$ {\bf La topologie atomique d'une cat\'egorie codirig\'ee:}

\smallskip

La d\'efinition des topologies de la densit\'e se simplifie dans le cadre des cat\'egories codirig\'ees:

\begin{defn}\label{defII24}
\begin{listeimarge}
\item Une cat\'egorie ${\mathcal C}$ est dite dirig\'ee [resp. codirig\'ee] si tout diagramme de ${\mathcal C}$ de la forme
$$
\xymatrix{
Z \ar[d] \ar[r] &X \\
Y
} \qquad \qquad
\xymatrix{
{ \ } \ar@{}[d]_{\mbox{[resp.} \quad }&X \ar[d]^{\quad \mbox{]}} \\
Y \ar[r] &S
}
$$
peut \^etre compl\'et\'e en un carr\'e commutatif:
$$
\xymatrix{
Z \ar[d] \ar[r] &X \ar[d] \\
Y \ar[r] &S
}
$$

\item La topologie atomique d'une cat\'egorie codirig\'ee essentiellement petite est celle dont les cribles couvrants sont les cribles non vides.
\end{listeimarge}
\end{defn}

\begin{remarkqed}

Dans une cat\'egorie codirig\'ee essentiellement petite, la topologie de la densit\'e est la topologie atomique. 

\end{remarkqed}

\bigskip

\noindent $\bullet$ {\bf Les topologies qui ont des cribles couvrants minimaux:}

\begin{lem}\label{lemII25}

Soient ${\mathcal C}$ une cat\'egorie essentiellement petite et ${\mathcal M}$ une collection de morphismes de ${\mathcal C}$ telle que:

\bigskip

\noindent $\left\{\begin{matrix}
\bullet &\mbox{pour toute fl\`eche $U \xrightarrow{ \ f \ } X$ qui est dans ${\mathcal M}$ et toute fl\`eche $X \xrightarrow{ \ h \ } Y$ de ${\mathcal C}$,} \hfill \\
&\mbox{la compos\'ee $h \circ f : U \to Y$ est dans ${\mathcal M}$,} \hfill \\
{ \ } \\
\bullet &\mbox{toute fl\`eche $U \xrightarrow{ \ f \ } X$ qui est dans ${\mathcal M}$ s'\'ecrit comme la compos\'ee $U \xrightarrow{ \ g \ } V \xrightarrow{ \ h \ } X$} \hfill \\
&\mbox{de deux fl\`eches $g,h$ de ${\mathcal M}$.} \hfill
\end{matrix} \right.$

\bigskip

Alors on d\'efinit une topologie $J_{\mathcal M}$ sur ${\mathcal C}$ en d\'ecidant qu'un crible $S$ sur un objet $X$ de ${\mathcal C}$ est couvrant si et seulement si il contient toutes les fl\`eches $U \to X$ de ${\mathcal M}$ dont le but est $X$.
\end{lem}

\begin{remarks}
\begin{listeisansmarge}
\item Les topologies de la forme $J_{\mathcal M}$ sont les topologies de ${\mathcal C}$ telles que tout objet $X$ de ${\mathcal C}$ admette un crible couvrant minimal.

\medskip

\item Si $O$ est une collection d'objets de ${\mathcal C}$, on peut en particulier prendre pour ${\mathcal M}$ la collection des fl\`eches $U \to X$ dont l'origine $U$ est \'el\'ement de $O$.
\end{listeisansmarge}
\end{remarks}

\bigskip

\begin{demo}

Il faut prouver que $J_{\mathcal M}$ satisfait les trois conditions (M), (S) et (T) qui d\'efinissent la notion de topologie de Grothendieck.

\smallskip

La condition de maximalit\'e (M) est \'evidente.

\smallskip

Si $X \xrightarrow{ \ f \ } Y$ est une fl\`eche de ${\mathcal C}$, $S$ est un crible couvrant de $Y$ et $U \xrightarrow{ \ h \ } X$ est un \'el\'ement de ${\mathcal M}$, alors $f \circ h : U \to Y$ est aussi dans ${\mathcal M}$ ce qui implique que $h \in f^* S$. Donc $f^* S$ est un crible couvrant de $X$ et la condition de stabilit\'e (S) est satisfaite.

\smallskip

Enfin, soit $S$ un crible sur un objet $X$ de ${\mathcal C}$ tel que, pour toute fl\`eche $V \xrightarrow{ \ h \ } X$ qui est dans ${\mathcal M}$, $h^* S$ soit un crible couvrant de $V$. Toute fl\`eche $U \xrightarrow{ \ f \ } X$ de but $X$ qui est dans ${\mathcal M}$ s'\'ecrit comme la compos\'ee
$$
f = h \circ g : U \xrightarrow{ \ g \ } V \xrightarrow{ \ h \ } X
$$
de deux fl\`eches $g,h$ qui sont elles-m\^emes dans ${\mathcal M}$. Comme $h^* S$ est un crible couvrant de $V$, on a n\'ecessairement $g \in h^* S$ soit $h \circ g = f \in S$. Ainsi, $S$ est un crible couvrant de $X$ et la condition de transitivit\'e (T) est satisfaite. 

\end{demo}

\medskip

\noindent $\bullet$ {\bf Les topologies des mono{\"\i}des:}

\smallskip

On a vu qu'un crible d'un mono{\"\i}de $M$ est une partie $S$ de $M$ stable par multiplication \`a droite.

\smallskip

Une topologie de $M$ est une partie $J$ de l'ensemble $\Omega$ des cribles de $M$ telle que, pour tout $S \in J$ et tout \'el\'ement $m \in M$, le crible
$$
m^* S = \{m' \in M \mid m \cdot m' \in S \}
$$
soit \'egalement \'el\'ement de $J$.

\smallskip

Par exemple, si $M = {\mathbb N}^{\times}$ est le mono{\"\i}de des entiers $n \geq 1$ muni de la multiplication, toute partie ${\mathcal P}$ de l'ensemble des nombres premiers d\'efinit une topologie $J_{\mathcal P}$ de ${\mathbb N}^{\times}$ dont les cribles couvrants sont les ensembles de la forme
$$
p_1^{n_1} \cdots p_k^{n_k} \cdot  {\mathbb N}^{\times} \quad \mbox{avec} \quad p_1 , \cdots , p_k \in {\mathcal P} \quad \mbox{et} \quad n_1 , \cdots , n_k \geq 0 \, .
$$

Si $M = G$ est un groupe, les seuls cribles de $G$ sont $\emptyset$ et le crible maximal $G$ si bien que les seules topologies de $G$ sont $\{ G \}$ et $\{ \emptyset , G \}$.

\bigskip

\noindent $\bullet$ {\bf La topologie discr\`ete d'une cat\'egorie:}

\smallskip

La ``topologie discr\`ete'' d'une cat\'egorie essentiellement petite ${\mathcal C}$ est d\'efinie en d\'ecidant que l'unique crible couvrant de tout objet $X$ de ${\mathcal C}$ est le ``crible maximal'' sur $X$ engendr\'e par la fl\`eche ${\rm id}_X : X \to X$.

\smallskip

La topologie discr\`ete de ${\mathcal C}$ est la moins fine de toutes les topologies de ${\mathcal C}$.

\section{Sites g\'eom\'etriques}\label{sec23}

\subsection{Gros sites et petits sites}\label{subsec231}

Introduisons des propri\'et\'es tr\`es g\'en\'erales des classes de morphismes de cat\'egories qui permettent de d\'efinir des topologies de Grothendieck:

\begin{defn}\label{defII31}

Soit ${\mathcal C}$ une cat\'egorie localement petite.

\begin{listeimarge}

\item Une classe ${\mathcal M}$ de morphismes de ${\mathcal C}$ sera dite ``g\'eom\'etrique'' si

\bigskip

\noindent $\left\{\begin{matrix}
\bullet &\mbox{cette classe comprend tous les isomorphismes, en particulier les morphismes ${\rm id}_X : X \to X$,} \hfill \\
{ \ } \\
\bullet &\mbox{elle est stable par composition au sens que pour tous morphismes de ${\mathcal M}$} \hfill \\
{ \ } \\
&f : X \longrightarrow Y \qquad \mbox{et} \qquad g : Y \longrightarrow Z \, , \\
{ \ } \\
&\mbox{leur compos\'e $g \circ f : X \to Z$ est encore dans ${\mathcal M}$,} \hfill \\
{ \ } \\
\bullet &\mbox{elle est stable par changement de base au sens que tout morphisme $f : X \to Y$ de ${\mathcal M}$ est carrable dans ${\mathcal C}$} \hfill \\
&\mbox{et tel que pour tout morphisme $Y' \to Y$ de ${\mathcal C}$, le morphisme induit} \hfill \\
{ \ } \\
&X \times_Y Y' \longrightarrow Y' \\
&\mbox{est encore dans ${\mathcal M}$}. \hfill
\end{matrix} \right.$

\bigskip

\item Dans les conditions de {\rm (i)}, une propri\'et\'e {\rm (R)} des familles de morphismes de ${\mathcal M}$ de m\^eme but
$$
(U_i \longrightarrow X)_{i \in I}
$$
sera appel\'ee une ``notion g\'eom\'etrique de recouvrement'' si
\begin{enumerate}
\item[$\bullet$] toute famille $(U_i \to X)_{i \in I}$ de morphismes de ${\mathcal M}$ qui comprend au moins un isomorphisme $U_{i_0} \xrightarrow{  \sim  } X$ poss\`ede la propri\'et\'e {\rm (R)},
\item[$\bullet$] la propri\'et\'e {\rm (R)} est stable par changement de base au sens que pour toute famille $(U_i \to X)_{i \in I}$ de morphismes de ${\mathcal M}$ qui la poss\`ede, toutes les familles
$$
(U_i \times_X X' \longrightarrow X')_{i \in I}
$$
qui s'en d\'eduisent par des changements de base $X' \to X$ poss\`edent encore la propri\'et\'e ${\rm (R)}$,
\item[$\bullet$] la propri\'et\'e {\rm (R)} est ``locale'' au sens que, r\'eciproquement, si $(U_i \to X)_{i \in I}$ est une famille de morphismes de ${\mathcal M}$ qui poss\`ede la propri\'et\'e {\rm (R)} et $(V_j \to X)_{j \in J}$ est une famille de morphismes de ${\mathcal M}$ telle que, pour tout $i \in I$, la famille induite
$$
(V_j \times_X U_i \longrightarrow U_i)_{i \in I}
$$
poss\`ede la propri\'et\'e {\rm (R)}, alors la famille $(V_j \to X)_{j \in J}$ poss\`ede la propri\'et\'e ${\rm (R)}$.
\end{enumerate}
\end{listeimarge}
\end{defn}

\begin{remarkqed}

Comme cons\'equence des conditions de (ii), une notion g\'eom\'etrique de recouvrement (R) poss\`ede encore les propri\'et\'es suivantes:
\begin{enumerate}
\item[$\bullet$] toute famille $(U_i \to X)_{i \in I}$ de morphismes de ${\mathcal M}$ qui admet une sous-famille poss\'edant la propri\'et\'e (R) poss\`ede a fortiori la propri\'et\'e (R),
\item[$\bullet$] si $(U_i \to X)_{i \in I}$ est une famille de morphismes de ${\mathcal M}$ qui poss\`ede la propri\'et\'e (R), alors pour tous isomorphismes $U'_i \xrightarrow{\sim} U_i$, $i \in I$, la famille induite $(U'_i \to X)_{i \in I}$ poss\`ede encore la propri\'et\'e (R).
\end{enumerate}
\end{remarkqed}

\bigskip

La donn\'ee d'une classe g\'eom\'etrique de morphismes d'une cat\'egorie et d'une notion g\'eom\'etrique de recouvrement permettent de d\'efinir des topologies de Grothendieck et donc des sites:

\begin{defn}\label{defII32}

Soient ${\mathcal C}$ une cat\'egorie essentiellement petite, ${\mathcal M}$ une classe g\'eom\'etrique de morphismes de ${\mathcal C}$ et {\rm (R)} une notion g\'eom\'etrique de recouvrement.

\smallskip

Alors:

\begin{listeimarge}

\item On d\'efinit une topologie $J_R$ sur ${\mathcal C}$ en d\'ecidant qu'un crible sur un objet $X$ de ${\mathcal C}$ est couvrant s'il contient une sous-famille constitu\'ee d'\'el\'ements de ${\mathcal M}$ et qui poss\`ede la propri\'et\'e {\rm (R)}.

\medskip

\item Pour tout objet $X$ de ${\mathcal C}$, on appelle ``gros site'' de $X$ (dans ${\mathcal C}$) d\'efini par $({\mathcal M} , R)$ la paire constitu\'ee de la cat\'egorie relative ${\mathcal C} / X$ et de sa topologie induite par $J_R$.

\medskip

\item Pour tout objet $X$ de ${\mathcal C}$, on appelle ``petit site'' de $X$ (dans ${\mathcal C}$) d\'efini par $({\mathcal M}, R)$ la paire constitu\'ee de la sous-cat\'egorie $({\mathcal C} /X)_{\mathcal M}$ de ${\mathcal C} / X$ telle que

\bigskip

\noindent $\left\{\begin{matrix}
\bullet &\mbox{un objet $X' \xrightarrow{p} X$ de ${\mathcal C}/X$ est dans $({\mathcal C}/X)_{\mathcal M}$ si et seulement si $p$ est \'el\'ement de ${\mathcal M}$,} \hfill \\
{ \ } \\
\bullet &\mbox{un morphisme entre deux objets de $({\mathcal C}/X)_{\mathcal M}$} \hfill \\
{ \ } \\
&\xymatrix{
X'_1 \ar[rd]_{p_1} \ar[rr]^f &&X'_2 \ar[ld]^{p_2} \\
&X
} \\
{ \ } \\
&\mbox{est dans $({\mathcal C}/X)_{\mathcal M}$ si et seulement si $f$ est dans ${\mathcal M}$,} \hfill
\end{matrix}\right.$

\bigskip

\noindent et de la topologie de $({\mathcal C}/X)_{\mathcal M}$ induite par $J_R$ pour laquelle un crible est couvrant quand son image par le foncteur d'oubli
$$
\begin{matrix}
\hfill ({\mathcal C} / X)_{\mathcal M} &\longrightarrow &{\mathcal C} \, , \hfill \\
(X' \xrightarrow{ \ p \ } X) &\longmapsto &X'
\end{matrix}
$$
engendre un crible $J_R$-couvrant c'est-\`a-dire poss\`ede la propri\'et\'e (R).
\end{listeimarge}
\end{defn}

\begin{remarkqed}

La topologie $J_R$ satisfait l'axiome de maximalit\'e (M) car, pour tout objet $X$ de ${\mathcal C}$, le morphisme ${\rm id}_X : X \to X$ est dans ${\mathcal M}$ et poss\`ede la propri\'et\'e (R).

\smallskip

Elle satisfait l'axiome de stabilit\'e (S) car la propri\'et\'e (R) est stable par changement de base.

\smallskip

Enfin, elle satisfait l'axiome de transitivit\'e (T) car la propri\'et\'e (R) est ``locale''. 

\end{remarkqed}

\subsection{La topologie ordinaire des espaces topologiques}\label{subsec232}

\smallskip

Les propri\'et\'es d\'efinissant les ``classes g\'eom\'etriques de morphismes'' et les ``notions g\'eom\'etriques de recouvrement'' ont \'et\'e abstraites \`a partir de celles de la classe des immersions ouvertes et de la notion de recouvrement ouvert dans le cadre des espaces topologiques:

\begin{lem}\label{lemII33}

Soit ${\mathcal C}$ une sous-cat\'egorie pleine de la cat\'egorie ${\rm Top}$ des espaces topologiques qui est essentiellement petite et contient toutes les immersions ouvertes $U \hookrightarrow X$ dans les objets $X$ de ${\mathcal C}$.

\smallskip

Alors:

\begin{listeimarge}

\item La classe ${\mathcal M}$ des compos\'es d'un hom\'eomorphisme et d'une immersion ouverte est un classe g\'eom\'etrique de morphismes de ${\mathcal C}$.

\medskip

\item La propri\'et\'e {\rm (R)} des familles de morphismes de ${\mathcal M}$ de m\^eme but
$$
U_i \xrightarrow{ \ f_i \ } X \, , \qquad i \in I \, ,
$$
qui consiste \`a demander que
$$
X = \bigcup_{i \in I} {\rm Im} (f_i)
$$
est une notion g\'eom\'etrique de recouvrement.
\end{listeimarge}
\end{lem}

\bigskip

\begin{demo}

Ce sont les propri\'et\'es famili\`eres des immersions ouvertes et des recouvrements ouverts.

\smallskip

Elles r\'esultent de ce que, pour toute application continue entre deux espaces topologiques
$$
f : X \longrightarrow Y
$$
et pour toute immersion ouverte $V \hookrightarrow Y$, l'ouvert
$$
f^{-1} V \xhookrightarrow{ \ \ \ } X
$$
est un produit fibr\'e $V \times_Y X$. 

\end{demo}

\pagebreak

Si ${\mathcal G}$ est une sous-cat\'egorie g\'eom\'etrique (au sens de la d\'efinition \ref{defI58}) de la cat\'egorie ${\rm Top}_{\rm an}$ des espaces annel\'es, toutes les immersions ouvertes vers des objets de ${\rm Top}_{\rm an}$ qui sont dans ${\mathcal G}$ sont \'egalement dans ${\mathcal G}$.

\smallskip

On d\'eduit du lemme ci-dessus:

\begin{cor}\label{corII34}

Soit ${\mathcal G}$ une sous-cat\'egorie g\'eom\'etrique de ${\rm Top}_{\rm an}$ qui est essentiellement petite.

\smallskip

Alors:

\begin{listeimarge}

\item La classe ${\mathcal M}$ des compos\'es d'un isomorphisme de ${\mathcal G}$ et d'une immersion ouverte de ${\mathcal G}$ est une classe g\'eom\'etrique de morphismes de ${\mathcal G}$.

\medskip

\item La propri\'et\'e {\rm (R)} des familles de morphismes de ${\mathcal M}$ de m\^eme but
$$
(U_i , {\mathcal O}_{U_i}) \xrightarrow{ \ f_i \ } (X,{\mathcal O}_X) \, , \qquad i \in I \, ,
$$
qui consiste \`a demander que
$$
X = \bigcup_{i \in I} {\rm Im} (f_i)
$$
est une notion g\'eom\'etrique de recouvrement.
\end{listeimarge}
\end{cor}

\begin{remarksqed}
\begin{listeisansmarge}
\item Ce corollaire s'applique en particulier \`a la cat\'egorie g\'eom\'etrique ${\mathcal G}$ des vari\'et\'es diff\'erentielles de classe $C^k$, $k \geq 1$, [resp. des vari\'et\'es analytiques] qui sont d\'enombrables \`a l'infini.

\smallskip

La cat\'egorie ${\mathcal G}$ munie de la ``topologie des recouvrements ouverts'' $J_R$ d\'efinie par la classe ${\mathcal M}$ des immersions ouvertes et la propri\'et\'e (R) de recouvrement ouvert est un site.

\smallskip

Et tout objet $X$ de ${\mathcal G}$ d\'efinit un gros site $({\mathcal G} / X , J_R)$ et un petit site $(({\mathcal G}/X)_{\mathcal M} , J_R)$ pour cette topologie.

\medskip

\item Le corollaire s'applique \'egalement \`a la cat\'egorie g\'eom\'etrique ${\rm Sch}$ des sch\'emas ou plut\^ot \`a n'importe quelle sous-cat\'egorie pleine ${\mathcal G}$ de ${\rm Sch}$ qui est g\'eom\'etrique (donc contient toutes les immersions ouvertes dans des objets de ${\mathcal G}$) et essentiellement petite.

\smallskip

Par exemple, on peut prendre pour ${\mathcal G}$ la sous-cat\'egorie pleine de ${\rm Sch}$ constitu\'ee des sch\'emas qui admettent un morphisme de pr\'esentation finie vers un sch\'ema $S$ donn\'e.

\smallskip

La topologie $J_R$ sur ${\mathcal G}$ d\'efinie par la classe ${\mathcal M}$ des immersions ouvertes et la propri\'et\'e (R) de recouvrement ouvert est appel\'ee la topologie de Zariski de ${\mathcal G}$.

\smallskip

Tout sch\'ema $X$ de ${\mathcal G}$ d\'efinit ainsi un gros site de Zariski $({\mathcal G}/X , J_R)$ et un petit site de Zariski $(({\mathcal G}/X)_{\mathcal M} , J_R)$. 

\end{listeisansmarge}
\end{remarksqed}

\medskip

On observe que dans le cadre des vari\'et\'es diff\'erentielles ou analytiques la classe des morphismes submersifs (au sens de la d\'efinition I.9.16 (i)) est g\'eom\'etrique et qu'il lui est naturellement associ\'e une notion g\'eom\'etrique de recouvrement:

\begin{cor}\label{corII35}

Soit ${\mathcal G}$ la cat\'egorie g\'eom\'etrique (et essentiellement petite) des vari\'et\'es diff\'erentielles de classe $C^k$ [resp. des vari\'et\'es analytiques] qui sont d\'enombrables \`a l'infini.

\smallskip

Alors:

\begin{listeimarge}

\item La classe ${\mathcal M}$ des morphismes submersifs est une classe g\'eom\'etrique de morphismes de ${\mathcal G}$.

\medskip

\item La propri\'et\'e {\rm (R)} des familles de morphismes submersifs de m\^eme but
$$
(U_i \xrightarrow{ \ f_i \ } X)_{i \in I}
$$
qui consiste \`a demander que
$$
X = \bigcup_{i \in I} {\rm Im} (f_i)
$$
est une notion g\'eom\'etrique de recouvrement.
\end{listeimarge}
\end{cor}

\begin{remarks}
\begin{listeisansmarge}
\item La topologie associ\'ee $J_R$ de ${\mathcal G}$ peut \^etre appel\'ee la topologie des submersions ou la topologie lisse.

\medskip

\item Tout objet $X$ de ${\mathcal G}$ d\'efinit ainsi un gros site lisse $({\mathcal G} / X , J_R)$ et un petit site lisse $(({\mathcal G}/X)_{\mathcal M} , J_R)$.

\medskip

\item Si $U \xrightarrow{ \ f \ } X$ est un morphisme submersif et $x$ un \'el\'ement de ${\rm Im} (f) \subset X$, il existe d'apr\`es le th\'eor\`eme des fonctions implicites un voisinage ouvert $V$ de $x$ dans $X$ tel que l'immersion ouverte
$$
V \xhookrightarrow{ \ \ \ } X
$$
se factorise en $V \longrightarrow U \xrightarrow{ \ f \ } X$.

\smallskip

Par cons\'equent, la topologie ordinaire et la topologie lisse co{\"\i}ncident sur la cat\'egorie ${\mathcal G}$.

\smallskip

En revanche, si $X$ est un objet de ${\mathcal G}$, le petit site lisse de $X$ et le petit site de $X$ pour la topologie ordinaire ne co{\"\i}ncident pas.
\end{listeisansmarge}
\end{remarks}

\medskip

\begin{democor}
\begin{listeisansmarge}
\item Il est \'evident que tout isomorphisme de ${\mathcal G}$ est submersif.

\smallskip

Le compos\'e de deux morphismes submersifs est un morphisme submersif car la compos\'ee de deux matrices surjectives est une matrice surjective.

\smallskip

Enfin, on a vu dans la remarque (i) suivant le th\'eor\`eme \ref{thmI919} que tout morphisme submersif est carrable dans ${\mathcal G}$ et que les morphismes qui s'en d\'eduisent par changement de base sont encore submersifs.

\medskip

\item r\'esulte de (i) et de ce que le morphisme d'oubli de la structure diff\'erentielle [resp. analytique]
$$
{\mathcal G} \longrightarrow {\rm Top} \longrightarrow {\rm Ens}
$$
respecte les produits fibr\'es
$$
U \times_X X'
$$
d'un morphisme submersif $U \to X$ avec un morphisme arbitraire $X' \to X$ de ${\mathcal G}$.
\end{listeisansmarge}
\end{democor}

\subsection{Des classes g\'eom\'etriques de morphismes de sch\'emas}\label{subsec233}

\smallskip

Dans la cat\'egorie ${\rm Sch}$ des sch\'emas, tout morphisme est carrable. Cette propri\'et\'e permet de d\'efinir une grande vari\'et\'e de classes de morphismes de sch\'emas qui sont ``g\'eom\'etriques'' au sens de la d\'efinition \ref{defII31}~(i).

\smallskip

Il en est d'abord ainsi de la classe des morphismes de pr\'esentation finie introduite dans la d\'efinition~\ref{defI714}:

\begin{lem}\label{lemII36}

Dans la cat\'egorie ${\rm Sch}$ des sch\'emas, les classes de morphismes suivantes sont g\'eom\'etriques:
\begin{enumerate}[label=(\arabic*)]
\item La classe des morphismes de pr\'esentation finie [resp. de type fini].
\item La classe des morphismes $X \to S$ qui sont ``localement de pr\'esentation finie'' [resp. ``localement de type fini''] au sens que $X$ est r\'eunion d'ouverts $U$ tels que $U \to X$ est de pr\'esentation finie [resp. de type fini].
\item La classe des morphismes $X \xrightarrow{ \ p \ } S$ qui sont ``quasi-compacts'' au sens que pour tout ouvert affine $V \cong {\rm Spec} (A)$ de $S$, $p^{-1} (V)$ est une r\'eunion finie d'ouverts affines de $X$.
\end{enumerate}
\end{lem}

\begin{remarks}
\begin{listeisansmarge}
\item Un morphisme est de pr\'esentation finie [resp. de type fini] si et seulement si il est \`a la fois localement de pr\'esentation finie [resp. localement de type fini] et quasi-compact.

\medskip

\item Dans un sch\'ema affine ${\rm Spec} (A)$, une famille d'ouverts de la forme ${\rm Spec} (A_{f_i})$, $f_i \in A$, $i \in I$, recouvre ${\rm Spec} (A)$ si et seulement si la famille des $f_i$ n'est contenue dans aucun id\'eal premier ou, ce qui revient au m\^eme,
$$
\sum_{i \in I} f_i \cdot A = A \, .
$$
Il en r\'esulte que tout sch\'ema affine ${\rm Spec} (A)$ est quasi-compact au sens que tout recouvrement ouvert de ${\rm Spec} (A)$ admet un sous-recouvrement fini.

\smallskip

Par cons\'equent, un morphisme de sch\'emas $X \xrightarrow{ \ p \ } S$ est quasi-compact si et seulement si l'image r\'eciproque de tout ouvert quasi-compact de $S$ est un ouvert quasi-compact de $X$.
\end{listeisansmarge}
\end{remarks}

\bigskip

\begin{demolem}

Il est \'evident que les classes de (1), (2) et (3) comprennent les isomorphismes.

\smallskip

La stabilit\'e par composition est \'evidente pour (3). Elle r\'esulte pour (1) et (2) de ce que, pour tous morphismes de pr\'esentation finie [resp. de type fini] entre anneaux commutatifs
$$
A \longrightarrow B \qquad \mbox{et} \qquad B \longrightarrow C \, ,
$$
le compos\'e $A \to C$ est encore de pr\'esentation finie [resp. de type fini], avec pour g\'en\'erateurs ou relations les r\'eunions de ceux de $C$ sur $B$ et de ceux de $B$ sur $A$.

\smallskip

La stabilit\'e par changement de base pour (3) r\'esulte de ce que tout produit fibr\'e de sch\'emas affines
$$
{\rm Spec} (A') \times_{{\rm Spec} (A)} {\rm Spec} (B)
$$
est un sch\'ema affine ${\rm Spec} (A' \otimes_A B)$.

\smallskip

La stabilit\'e par changement de base pour (1) et (2) r\'esulte de ce que, de plus, $A' \otimes_A B$ est de pr\'esentation finie [resp. de type fini] sur $A'$ si $B$ est de pr\'esentation finie [resp. de type fini] sur $A$, avec les m\^emes g\'en\'erateurs et les m\^emes relations.

\end{demolem}

\pagebreak

Voici encore d'autres classes g\'eom\'etriques de morphismes de sch\'emas:

\begin{defn}\label{defII37}

Un morphisme $X \to S$ de sch\'emas est dit
\begin{enumerate}
\item[(4)] affine si pour tout morphisme ${\rm Spec} (A) \to S$ depuis un sch\'ema affine, le produit fibr\'e $X \times_S {\rm Spec} (A)$ est affine,
\item[(5)] fini s'il est affine et si pour tout morphisme ${\rm Spec} (A) \to S$ avec $X \times_S {\rm Spec} (A) \cong {\rm Spec} (B)$, $B$ est de type fini comme module (et non alg\`ebre) sur $A$,
\item[(6)] une immersion ferm\'ee s'il est fini et que pour tout ${\rm Spec} (A) \to S$ avec $X \times_S {\rm Spec} (A) \cong {\rm Spec} (B)$, le morphisme $A \to B$ est surjectif,
\item[(7)] une immersion localement ferm\'ee s'il est le compos\'e d'une immersion ouverte et d'une immersion ferm\'ee.
\end{enumerate}
\end{defn}

\begin{remarksqed}
\begin{listeisansmarge}
\item Il est \'evident sur les d\'efinitions que les classes de morphismes (4), (5), (6) et (7) sont g\'eom\'etriques.


\item Pour qu'un morphisme $X \to S$ soit affine [resp. fini, resp. une immersion ferm\'ee], il suffit que $S$ admette un recouvrement par des ouverts affines $V_i \cong {\rm Spec} (A_i)$ tels que les produits fibr\'es $V_i \times_S X$ soient des sch\'emas affines ${\rm Spec} (B_i)$ [resp. avec chaque $B_i$ de type fini comme module sur $A_i$, resp. avec pour chaque $B_i$ un quotient de $A_i$].

\end{listeisansmarge}

La d\'emonstration de ces faits (qui ne sont pas \'evidents) est laiss\'ee au lecteur.\end{remarksqed}


Les morphismes diagonaux $X \to X \times_S X$ associ\'es aux morphismes de sch\'emas $X \to S$ fournissent des familles tr\`es importantes d'immersions ferm\'ees ou localement ferm\'ees:

\begin{lem}\label{lemII38}

Le morphisme diagonal $X \to X \times_S X$ associ\'e \`a un morphisme de sch\'emas $X \to S$ est

\begin{listeimarge}

\item une immersion ferm\'ee si $X$ et $S$ sont affines,

\medskip

\item une immersion localement ferm\'ee dans le cas g\'en\'eral.\end{listeimarge}
\end{lem}

\begin{demo}
\begin{listeisansmarge}
\item En effet, si $S = {\rm Spec} (A)$ et $X = {\rm Spec} (B)$, on a
$$
X \times_S X \cong {\rm Spec} (B \otimes_A B)
$$
et le morphisme diagonal $X \to X \times_S X$ est induit par le morphisme d'anneaux
$$
\begin{matrix}
B \otimes_A B &\longrightarrow &B \, , \hfill \\
\hfill b \otimes b' &\longmapsto &bb'
\end{matrix}
$$
lequel est surjectif.


\item On peut recouvrir $X$ par des ouverts affines $U_i \cong {\rm Spec} (B_i)$ qui s'envoient dans des ouverts affines $V_i \cong{\rm Spec} (A_i)$ de $S$.

\smallskip

Alors le morphisme diagonal $X \to X \times_S X$ se factorise \`a travers l'ouvert de $X \times_S X$ d\'efini comme la r\'eunion des
$$
U_i \times_{V_i} U_i \cong {\rm Spec} (B_i \otimes_{A_i} B_i) \, .
$$
De plus, chaque produit fibr\'e
$$
X \times_{X \times_S X} U_i \times_{V_i} U_i
$$
s'identifie \`a $U_i$ plong\'e diagonalement dans $U_i \times_{V_i} U_i$.

\smallskip

Donc le morphisme de $X$ dans la r\'eunion des $U_i \times_{V_i} U_i$ est une immersion ferm\'ee.\end{listeisansmarge}\end{demo}

\pagebreak

Les morphismes diagonaux sont d'autant plus importants qu'ils peuvent servir de base au d\'eveloppement d'un calcul diff\'erentiel alg\'ebrique sur les sch\'emas:

\begin{prop}\label{propII39}

Soit un morphisme de sch\'emas affines
$$
{\rm Spec} (B) \longrightarrow {\rm Spec} (A) \, .
$$

Soit $I$ le noyau du morphisme surjectif
$$
\begin{matrix}
B \otimes_A B &\longrightarrow &B \, , \hfill \\
\hfill b \otimes b' &\longmapsto &bb' \, . \hfill
\end{matrix}
$$

Alors le $B$-module
$$
I/I^2 = \Omega_{B/A}
$$
repr\'esente le foncteur covariant
$$
\begin{matrix}
{\rm Mod}_B  &\longrightarrow &{\rm Ens} \, , \hfill \\
\hfill M &\longmapsto &{\rm Der}_A (B,M)
\end{matrix}
$$
qui associe \`a tout $B$-module $M$ l'ensemble ${\rm Der}_A (B,M)$ des ``d\'erivations'' de $B$ sur $A$ \`a valeurs dans $M$ c'est-\`a-dire des applications $A$-lin\'eaires
$$
{\rm d} : B \longrightarrow M
$$
telles que
$$
{\rm d}(bb') = b \cdot {\rm d}b' + b' \cdot {\rm d}b \, , \qquad \forall \, b , b' \in B \, .
$$
\end{prop}

\begin{remarks}
\begin{listeisansmarge}
\item Le $B$-module $\Omega_{B/A}$ est appel\'e le module des diff\'erentielles de $B$ sur $A$.

\medskip

\item Une d\'erivation $A$-lin\'eaire ${\rm d} : B \to M$ s'annule n\'ecessairement en l'\'el\'ement $1$ et plus g\'en\'eralement en les images des \'el\'ements de $A$.
\end{listeisansmarge}
\end{remarks}

\medskip

\begin{demo}

Le quotient $I/I^2$ de l'id\'eal $I$ de $B \otimes_A B$ par son carr\'e $I^2$ (qui est l'id\'eal engendr\'e par les produits de paires d'\'el\'ements de $I$) est un module sur l'anneau quotient $B \otimes_A B/I = B$.

\smallskip

L'application
$$
\begin{matrix}
{\rm d} : &B &\longrightarrow &I/I^2 \, , \hfill \\
&b &\longmapsto &b \otimes 1 - 1 \otimes b
\end{matrix}
$$
respecte l'addition, s'annule sur l'image de $A$ dans $B$ et v\'erifie la formule
$$
{\rm d}(bb') = b \cdot {\rm d}b' + b' \cdot {\rm d}b
$$
car
$$
\begin{matrix}
&(bb' \otimes 1-1 \otimes bb') - (bb' \otimes 1-b \otimes b') - (b'b \otimes 1-b' \otimes b) \\
= &b \otimes b' + b' \otimes b - bb' \otimes 1 - 1 \otimes bb' \hfill \\
= &(b \otimes 1-1 \otimes b) \cdot (1 \otimes b'-b' \otimes 1) \hfill
\end{matrix}
$$
est un \'el\'ement de $I^2$.

\smallskip

Il reste \`a prouver que toute d\'erivation $A$-lin\'eaire \`a valeurs dans un $B$-module $M$
$$
{\rm d} : B \longrightarrow M
$$
se factorise \`a travers ${\rm d} : B \to I/I^2$ en un morphisme $B$-lin\'eaire uniquement d\'etermin\'e
$$
I/I^2 \longrightarrow M \, .
$$
L'application
$$
\begin{matrix}
B \times B &\longrightarrow &M \, , \hfill \\
\hfill (b , b') &\longmapsto &-b \cdot {\rm d}b' \hfill
\end{matrix}
$$
est $B$-lin\'eaire en la premi\`ere variable et $A$-lin\'eaire en la seconde.

\smallskip

Elle se factorise en un morphisme $A$-lin\'eaire
$$
B \otimes_A B \longrightarrow M
$$
qui se restreint en
$$
\begin{matrix}
u : &\hfill I &\longrightarrow &M \, , \hfill \\
&\sum (b_i \otimes b'_i) &\longmapsto &-\sum b_i \cdot {\rm d}b'_i = \sum b'_i \cdot {\rm d} b_i
\end{matrix}
$$
puisque tout \'el\'ement $\sum (b_i \otimes b'_i)$ de $I$ v\'erifie
$$
\sum_i b_i \, b'_i = 0 \, .
$$
Donc le morphisme $u$ est $B \otimes_A B$-lin\'eaire. Comme $M$ est un module sur $B = B \otimes_A B/I$, $u$ s'annule sur $I^2$ et d\'efinit un morphisme $B$-lin\'eaire
$$
I/I^2 \longrightarrow M
$$
qui est la factorisation cherch\'ee de ${\rm d} : B \to M$.

\end{demo}

\medskip

Le module des diff\'erentielles d'un anneau commutatif sur un autre se calcule facilement gr\^ace \`a sa caract\'erisation par le foncteur des d\'erivations:

\begin{cor}\label{corII310}

Soit un morphisme de sch\'emas affines
$$
{\rm Spec} (B) \longrightarrow {\rm Spec} (A) \, .
$$
Alors le module des diff\'erentielles relatives
$$
\Omega_{B/A}
$$
est

\begin{listeimarge}

\item le $B$-module libre
$$
\bigoplus_{i \in I} B \cdot {\rm d} X_i
$$
si $B$ est la $A$-alg\`ebre
$$
B = A [(X_i)_{i \in I}]
$$
des polyn\^omes en des variables $X_i$, $i \in I$,

\medskip

\item le quotient
$$
\bigoplus_{i \in I} B \cdot {\rm d}X_i \Biggl/ \left( \sum_i \frac{\partial P_j}{\partial X_i} \cdot {\rm d}X_i \right)_{j\in J}
$$
du $B$-module libre $\underset{i \in I}{\bigoplus} \, B \cdot {\rm d} X_i$ par le sous-module engendr\'e par les relations
$$
\sum_i \frac{\partial P_j}{\partial X_i} \cdot {\rm d} X_i
$$
si $B$ est le quotient
$$
B = A [(X_i)_{i \in I}] \bigl/ (P_j)_{j \in J}
$$
de l'alg\`ebre de polyn\^omes $A[(X_i)_{i \in I}]$ par son id\'eal engendr\'e par les polyn\^omes $P_j$, $j \in J$,

\medskip

\item \'egal \`a $0$ si [resp. si et seulement si] le morphisme diagonal
$$
{\rm Spec} (B) \longrightarrow {\rm Spec} (B) \times_{{\rm Spec} (A)} {\rm Spec} (B)
$$
est aussi une immersion ouverte [resp. d\`es lors que la $A$-alg\`ebre $B$ est de type fini].
\end{listeimarge}
\end{cor}

\begin{remark}

Il r\'esulte de (ii) que le $B$-module $\Omega_{B/A}$ est de type fini [resp. de pr\'esentation finie] si la $A$-alg\`ebre $B$ est de type fini [resp. de pr\'esentation finie].
\end{remark}

\medskip

\begin{demo}
\begin{listeisansmarge}
\item r\'esulte de ce que, si $B = A [(X_i)_{i \in I}]$, il existe pour tout $B$-module $M$ et tous \'el\'ements $m_i \in M$, $i \in I$, une unique d\'erivation $A$-lin\'eaire
$$
{\rm d} : B \longrightarrow M
$$
telle que
$$
{\rm d} X_i = m_i \, , \qquad \forall \, i \in I \, .
$$

\item r\'esulte de ce que, de plus, une d\'erivation $A$-lin\'eaire
$$
{\rm d} : A [(X_i)_{i \in I}] \longrightarrow M
$$
s'annule sur l'id\'eal engendr\'e par les polyn\^omes $P_j$, $j \in J$, si et seulement si les images $m_i = {\rm d} X_i$, $i \in I$, v\'erifient
$$
P_j \cdot m_i = 0 \, , \qquad \forall \, j \, , \ \forall \, i \, ,
$$
et
$$
\sum_i \frac{\partial P_j}{\partial X_i} \cdot m_i = 0 \, , \qquad \forall \, j \, .
$$

\item Si l'immersion ferm\'ee
$$
{\rm Spec} (B) \longrightarrow {\rm Spec} (B) \times_{{\rm Spec} (A)} {\rm Spec} (B)
$$
est aussi une immersion ouverte, il existe un \'el\'ement $f \in B \otimes_A B$ qui co{\"\i}ncide avec $1$ sur l'ouvert ${\rm Spec} (B)$ et avec $0$ sur l'ouvert compl\'ementaire.

\smallskip

Il v\'erifie
$$
f \cdot I = 0
$$
et a fortiori
$$
f \cdot I/I^2 = 0 \, .
$$
Comme $f$ co{\"\i}ncide avec $1$ dans $B$, cela prouve que
$$
I/I^2 = 0 \, .
$$

R\'eciproquement, supposons que $B$ est engendr\'e sur $A$ par un nombre fini d'\'el\'ements $b_1 , \cdots , b_n$ et que
$$
\Omega_{B/A} = I/I^2 = 0 \, .
$$
L'id\'eal $I$ de $B \otimes_A B$ est engendr\'e par la famille des \'el\'ements
$$
b_i \otimes 1 - 1 \otimes b_i = f_i \, , \qquad 1 \leq i \leq n \, .
$$
Comme $I=I^2$, on peut \'ecrire pour tout $i$, $1 \leq i \leq n$,
$$
f_i = f'_{i,1} \cdot f_1 + \cdots + f'_{i,n} \cdot f_n \, ,
$$
avec $f'_{i,j} \in I$, $\forall \, i$, $\forall \, j$.

\smallskip

Le d\'eterminant de la matrice
$$
{\rm Id} - (f'_{i,j})_{1 \leq i,j \leq n}
$$
a la forme $1+f$, avec $f \in I$, et il v\'erifie
$$
(1+f) \cdot f_i = 0 \, , \qquad 1 \leq i \leq n \, ,
$$
soit $(1+f) \cdot I = 0$.

\smallskip

Cela prouve que ${\rm Spec} (B)$ s'identifie \`a l'ouvert de ${\rm Spec} (B) \times_{{\rm Spec} (A)} {\rm Spec} (B)$ d\'efini par l'\'el\'ement $1+f \in B \otimes_A B$. 
\end{listeisansmarge}
\end{demo}

\medskip

La consid\'eration des morphismes diagonaux permet enfin de d\'efinir une autre classe g\'eom\'etrique importante de morphismes de sch\'emas, celle des morphismes s\'epar\'es, et par suite celle des morphismes propres:

\begin{defn}\label{defII311}

Un morphisme de sch\'emas $X \to S$ est dit
\begin{enumerate}
\item[(8)] ``s\'epar\'e'' si le morphisme diagonal associ\'e
$$
X \longrightarrow X \times_S X
$$
est une immersion ferm\'ee,
\item[(9)] ``propre'' s'il est s\'epar\'e, de type fini et ``universellement ferm\'e'' au sens que, pour tout morphisme de changement de base $S' \to S$, le morphisme induit
$$
X \times_S S' \longrightarrow S'
$$
transforme tout ferm\'e de $X \times_S S'$ en un ferm\'e de $S'$.
\end{enumerate}
\end{defn}

\pagebreak

\begin{remarksqed}
\begin{listeisansmarge}
\item La classe des morphismes s\'epar\'es comprend les isomorphismes.

\smallskip

Elle est stable par composition car, pour tous morphismes
$$
X \longrightarrow Y \qquad \mbox{et} \qquad Y \longrightarrow Z \, ,
$$
le morphisme diagonal $X \to X \times_Z X$ est le compos\'e du morphisme diagonal $X \to X \times_Y X$ et du morphisme $X \times_Y X \to X \times_Z X$ qui s'inscrit dans le carr\'e cart\'esien:
$$
\xymatrix{
X \times_Y X \ar[d] \ar[r] &X \times_Z X \ar[d] \\
Y \ar[r] &Y \times_Z Y
}
$$

Elle est stable par changement de base car tout carr\'e cart\'esien
$$
\xymatrix{
X' \ar[d] \ar[r] &X \ar[d] \\
S' \ar[r] &S
}
$$
induit un carr\'e cart\'esien:
$$
\xymatrix{
X' \ar[d] \ar[r] &X' \times_{S'} X' \ar[d] \\
X \ar[r] &X \times_S X
}
$$
C'est pourquoi la classe des morphismes s\'epar\'es est g\'eom\'etrique.

\medskip

\item La classe des morphismes de type fini est g\'eom\'etrique.

\smallskip

La d\'efinition des morphismes universellement ferm\'es rend \'evident que leur classe est aussi g\'eom\'etrique.

\smallskip

Donc la classe des morphismes propres est g\'eom\'etrique parce que celle des morphismes s\'epar\'es est g\'eom\'etrique.

\medskip

\item Il r\'esulte du lemme \ref{lemII38} (i) que tous les morphismes affines de ${\rm Sch}$ sont des morphismes s\'epar\'es.

\medskip

\item Un morphisme $X \to S$ est un monomorphisme si et seulement si le morphisme diagonal
$$
X \longrightarrow X \times_S X
$$
est un isomorphisme. Donc tout monomorphisme de ${\rm Sch}$ est s\'epar\'e.

\smallskip

En particulier, les immersions ouvertes (ou localement ferm\'ees) de ${\rm Sch}$ sont des morphismes s\'epar\'es.

\medskip

\item Toutes les immersions ferm\'ees de ${\rm Sch}$ sont des morphismes propres.

\smallskip

On peut montrer que, plus g\'en\'eralement, tous les morphismes finis de ${\rm Sch}$ sont des morphismes propres.

\medskip

\item La notion de morphisme s\'epar\'e de ${\rm Sch}$ est un analogue relatif \`a une base $S$ de la notion d'espace topologique s\'epar\'e (ou de Hausdorff).

\smallskip

En effet, un espace topologique $X$ est s\'epar\'e (au sens que pour tous $x_1 \ne x_2$ dans $X$, il existe des ouverts $U_1 \ni x_1$, $U_2 \ni x_2$ de $X$ tels que $U_1 \cap U_2 = \emptyset$) si et seulement si la diagonale
$$
X \xhookrightarrow{ \ \ \ } X \times X
$$
est un ferm\'e.

\medskip

\item La notion de morphisme propre de ${\rm Sch}$ est un analogue relatif \`a une base $S$ de la notion d'espace topologique compact.

\smallskip

En effet, un espace topologique s\'epar\'e $X$ est compact (au sens que tout recouvrement ouvert a un sous-recouvrement fini) si et seulement si, pour tout espace topologique $Y$, la projection
$$
X \times Y \longrightarrow Y
$$
transforme tout ferm\'e de $X \times Y$ en un ferm\'e de $Y$. 
\end{listeisansmarge}
\end{remarksqed}

\medskip

La classe la plus importante de morphismes propres est celle des morphismes projectifs.

\smallskip

Leur d\'efinition d\'epend de la construction des espaces projectifs:

\begin{prop}\label{propII312}

Pour tout entier $n \geq 0$, notons
$$
{\mathbb A}^{n+1} - \{0\}
$$
le sous-sch\'ema ouvert de ${\mathbb A}^{n+1} = {\rm Spec} ({\mathbb Z} [X_0 , \cdots , X_n])$ d\'efini comme la r\'eunion des ouverts affines
$$
{\mathbb A}_i^{n+1} = {\rm Spec} ({\mathbb Z} [X_0 , \cdots , X_n , X_i^{-1}]) \, , \qquad 0 \leq i \leq n \, ,
$$
et muni des deux morphismes
$$
{\mathbb G}_m \times ({\mathbb A}^{n+1} - \{0\}) \rightrightarrows {\mathbb A}^{n+1} - \{0\}
$$
dont l'un est la projection sur le facteur ${\mathbb A}^{n+1} - \{0\}$ et l'autre correspond \`a la famille des applications
$$
\begin{matrix}
{\mathbb G}_m (A) \times ({\mathbb A}^{n+1} - \{0\}) (A) &\longrightarrow &({\mathbb A}^{n+1} - \{0\})(A) \, , \hfill  \\
\hfill (a, (a_0 , a_1 , \cdots , a_n)) &\longmapsto &(aa_0 , aa_1 , \cdots , aa_n)
\end{matrix}
$$
index\'ees par les anneaux commutatifs $A$.

\smallskip

Alors:

\begin{listeimarge}

\item Il existe un sch\'ema ${\mathbb P}^n$ qui s'inscrit dans un carr\'e cart\'esien et cocart\'esien:
$$
\xymatrix{
{\mathbb G}_m \times ({\mathbb A}^{n+1} - \{0\}) \ar[d] \ar[r] &{\mathbb A}^{n+1} - \{0\} \ar[d] \\
{\mathbb A}^{n+1} - \{0\} \ar[r] &{\mathbb P}^n
}
$$

\item Ce sch\'ema ${\mathbb P}^n$ est r\'eunion d'ouverts ${\mathbb P}_i^n$, $0 \leq i \leq n$, qui s'inscrivent dans des carr\'es cart\'esiens et cocart\'esiens
$$
\xymatrix{
{\mathbb G}_m \times {\mathbb A}^{n+1}_i \ar[d] \ar[r] &{\mathbb A}^{n+1}_i \ar[d] \\
{\mathbb A}^{n+1}_i \ar[r] &{\mathbb P}_i^n
}
$$
et sont munis d'isomorphismes canoniques
$$
{\mathbb P}_i^n \xrightarrow{ \ \sim \ } {\rm Spec} \left( {\mathbb Z} \left[ \frac{X_0}{X_i} , \cdots , \frac{X_{i-1}}{X_i} , \frac{X_{i+1}}{X_i} , \cdots , \frac{X_n}{X_i} \right]\right).
$$

\end{listeimarge}
\end{prop}

\newpage

\begin{demo}

Pour tout $i$, $0 \leq i \leq n$, le morphisme
$$
\begin{matrix}
{\rm Spec} ({\mathbb Z} [X_0 , \cdots , X_n , X_i^{-1}]) &\longrightarrow &\displaystyle {\rm Spec} \left( {\mathbb Z} \left[ \frac{X_0}{X_i} , \cdots , \frac{X_{i-1}}{X_i} , \frac{X_{i+1}}{X_i} , \cdots , \frac{X_n}{X_i} \right]\right), \\
\hfill (a_0 , a_1 , \cdots , a_n) &\longmapsto &\displaystyle \left( \frac{a_0}{a_i} , \cdots , \frac{a_{i-1}}{a_i} , \frac{a_{i+1}}{a_i} , \cdots , \frac{a_n}{a_i} \right) \hfill
\end{matrix}
$$
d\'efinit un carr\'e \`a la fois cart\'esien et cocart\'esien:
$$
\xymatrix{
{\mathbb G}_m \times {\mathbb A}_i^{n+1} \ar[d] \ar[r] &{\mathbb A}_i^{n+1} \ar[d] \\
{\mathbb A}_i^{n+1} \ar[r] &\displaystyle {\mathbb P}_i^n = {\rm Spec} \left( {\mathbb Z} \left[ \frac{X_0}{X_i} , \cdots , \frac{X_{i-1}}{X_i} , \frac{X_{i+1}}{X_i} , \cdots , \frac{X_n}{X_i} \right]\right)
}
$$

Pour tous indices $i \ne j$, l'intersection ${\mathbb A}_{i,j}^{n+1} = {\mathbb A}_i^{n+1} \cap {\mathbb A}_j^{n+1}$ des deux ouverts ${\mathbb A}_i^{n+1}$ et ${\mathbb A}_j^{n+1}$ de ${\mathbb A}^{n+1} - \{0\}$ est l'image r\'eciproque de deux ouverts
$$
\begin{matrix}
&{\mathbb P}_{i,j}^n &\mbox{de} &{\mathbb P}_i^n \\
\mbox{et} &{\mathbb P}_{j,i}^n &\mbox{de} &{\mathbb P}_j^n
\end{matrix}
$$
reli\'es par un isomorphisme canonique
$$
{\mathbb P}_{i,j}^n \xrightarrow{ \ \sim \ } {\mathbb P}_{j,i}^n \, ,
$$
$$
\left( \frac{a_0}{a_1} , \cdots , \frac{a_{i-1}}{a_i} , \frac{a_{i+1}}{a_i} , \cdots , \frac{a_n}{a_i} \right) \longmapsto \left(\frac{a_0}{a_j} , \cdots , \frac{a_{j-1}}{a_j} , \frac{a_{j+1}}{a_j} , \cdots , \frac{a_n}{a_j} \right)
$$
avec, pour tout indice $k$,
$$
\frac{a_k}{a_j} = \frac{a_k}{a_i} \Bigl/ \frac{a_j}{a_i} \, .
$$

On observe  que pour tous indices $i,j,k$, l'isomorphisme canonique
$$
{\mathbb P}_{i,k}^n \xrightarrow{ \ \sim \ } {\mathbb P}_{k,i}^n
$$
est le compos\'e des isomorphismes induits par
$$
{\mathbb P}_{i,j}^n \xrightarrow{ \ \sim \ } {\mathbb P}_{j,i}^n \qquad \mbox{et} \qquad {\mathbb P}_{j,k}^n \longrightarrow {\mathbb P}_{k,j}^n
$$
sur les ouverts de ${\mathbb P}_i^n$, ${\mathbb P}_j^n$ et ${\mathbb P}_k^n$ dont l'image r\'eciproque dans ${\mathbb A}^{n+1} - \{0\}$ est ${\mathbb A}^{n+1}_i \cap {\mathbb A}^{n+1}_j \cap {\mathbb A}^{n+1}_k$.

\smallskip

Il existe donc un sch\'ema ${\mathbb P}^n$, unique \`a unique isomorphe pr\`es, qui est r\'eunion d'ouverts
$$
{\mathbb P}_i^n \, , \qquad 0 \leq i \leq n \, ,
$$
dont les intersections ${\mathbb P}_i^n \cap {\mathbb P}_j^n$ sont les ouverts
$$
{\mathbb P}_{i,j}^n \subset {\mathbb P}_i^n \qquad \mbox{ou} \qquad {\mathbb P}_{j,i}^n \subset {\mathbb P}_j^n
$$
identifi\'es via les isomorphismes canoniques
$$
{\mathbb P}_{i,j}^n \xrightarrow{ \ \sim \ } {\mathbb P}_{j,i}^n \, .
$$

Le sch\'ema ${\mathbb P}^n$ ainsi construit r\'epond \`a la question pos\'ee. 

\end{demo}

\pagebreak

On a:

\begin{thm}\label{thmII313}

Pour tout entier $n \geq 0$, le morphisme
$$
{\mathbb P}^n \longrightarrow {\rm Spec} ({\mathbb Z})
$$
est propre et de pr\'esentation finie.
\end{thm}

\begin{remark}

Pour $n=0$, le morphisme ${\mathbb P}^0 \to {\rm Spec} ({\mathbb Z})$ est un isomorphisme.

\end{remark}

\medskip

\begin{demopart}

Le sch\'ema ${\mathbb P}^n$ est r\'eunion des ouverts
$$
{\mathbb P}_i^n \cong {\rm Spec} \left( {\mathbb Z} \left[ \frac{X_0}{X_i} , \cdots , \frac{X_{i-1}}{X_i} , \frac{X_{i+1}}{X_i} , \cdots , \frac{X_n}{X_i} \right]\right) , \quad 0 \leq i \leq n \, ,
$$
donc il est de pr\'esentation finie.

\smallskip

Le morphisme de projection
$$
p : {\mathbb A}^{n+1} - \{0\} \longrightarrow {\mathbb P}^n
$$
admet au-dessus de chaque ouvert ${\mathbb P}_i^n$ une section c'est-\`a-dire un morphisme
$$
s_i : {\mathbb P}_i^n \longrightarrow {\mathbb A}_i^{n+1} \xhookrightarrow{ \ \ \ } {\mathbb A}^{n+1} - \{0\}
$$
tel que
$$
p \circ s_i = {\rm id}_{{\mathbb P}_i^n} \, .
$$

Pour montrer que la diagonale
$$
{\mathbb P}^n \longrightarrow {\mathbb P}^n \times {\mathbb P}^n
$$
est ferm\'ee, il suffit donc de d\'emontrer que son image r\'eciproque dans
$$
({\mathbb A}^{n+1} - \{0\}) \times ({\mathbb A}^{n+1} - \{0\}) \xhookrightarrow{ \ \ \ } {\rm Spec} ({\mathbb Z} [X_0 , \cdots , X_n]) \times {\rm Spec} ({\mathbb Z} [Y_0 , \cdots , Y_n])
$$
est ferm\'ee.

\smallskip

Or celle-ci est d\'efinie par les \'equations
$$
X_i \, Y_j = X_j \, Y_i \, , \qquad 0 \leq i,j \leq n \, .
$$

Cela montre que ${\mathbb P}^n \to {\rm Spec} ({\mathbb Z})$ est s\'epar\'e.

\smallskip

La d\'emonstration de ce que
$$
{\mathbb P}^n \longrightarrow {\rm Spec} ({\mathbb Z})
$$
est universellement ferm\'e est laiss\'ee en probl\`eme au lecteur. 

\end{demopart}

\bigskip

Ce th\'eor\`eme implique la partie (ii) du corollaire suivant:

\begin{cor}\label{corII314}

Un morphisme de sch\'emas
$$
X \longrightarrow S
$$
est dit projectif [resp. quasi-projectif] s'il est possible de l'\'ecrire comme le compos\'e d'une immersion ferm\'ee [resp. localement ferm\'ee]
$$
X \xhookrightarrow{ \ \ \ } {\mathbb P}^n \times S
$$
et de la projection canonique associ\'ee \`a un espace projectif ${\mathbb P}^n$
$$
{\mathbb P}^n \times S \longrightarrow S \, .
$$

\pagebreak
Alors:

\begin{listeimarge}

\item La classes des morphismes projectifs [resp. quasi-projectifs] est g\'eom\'etrique.

\medskip

\item Tout morphisme projectif [resp. quasi-projectif] est propre [resp. s\'epar\'e].
\end{listeimarge}
\end{cor}

\begin{demo}
\begin{listeisansmarge}
\item[(ii)] r\'esulte du th\'eor\`eme \ref{thmII313} et des remarques (iv) et (v) apr\`es la d\'efinition \ref{defII311}.

\medskip

\item[(i)] Il est \'evident sur la d\'efinition que la classe des morphismes projectifs [resp. quasi-projectifs] comprend les isomorphismes et qu'elle est stable par changement de base.

\smallskip

Pour v\'erifier qu'elle est stable par composition, il suffit de prouver que le produit ${\mathbb P}^n \times {\mathbb P}^m$ de deux espaces projectifs ${\mathbb P}^n$ et ${\mathbb P}^m$ se plonge comme sous-sch\'ema ferm\'e dans un espace projectif.

\smallskip

Or le morphisme
$$
\begin{matrix}
({\mathbb A}^{n+1} - \{0\}) \times ({\mathbb A}^{m+1} - \{0\}) &\longrightarrow &{\mathbb A}^{(n+1)(m+1)} - \{0\} \, , \\
\hfill ((a_0 , \cdots , a_n) , (b_0 , \cdots , b_m)) &\longmapsto &(a_i \, b_j)_{0 \leq i \leq n \atop 0 \leq j \leq m} \hfill
\end{matrix}
$$
d\'efinit un morphisme
$$
{\mathbb P}^n \times {\mathbb P}^m \longrightarrow {\mathbb P}^{n+m+nm}
$$
qui est une immersion ferm\'ee.

\smallskip

En effet, son produit fibr\'e avec le morphisme
$$
{\mathbb A}^{(n+1)(m+1)} - \{0\} \longrightarrow {\mathbb P}^{n+m+nm}
$$
est le sous-sch\'ema ferm\'e de
$$
{\mathbb A}^{(n+1)(m+1)} - \{0\} \xhookrightarrow{ \ \ \ } {\mathbb A}^{(n+1)(m+1)} = {\rm Spec} \left({\mathbb Z} \left[(X_{i,j})_{0 \leq i \leq n \atop 0 \leq j \leq m}\right]\right)
$$
d\'efini par les \'equations
$$
X_{i,j} \, X_{i',j'} = X_{i,j'} \, X_{i',j} \, , \qquad 0 \leq i,i' \leq n \, , \ 0 \leq j,j' \leq m \, .
$$
\end{listeisansmarge}
\end{demo}

\bigskip

Nous allons enfin introduire trois derni\`eres classes g\'eom\'etriques de morphismes de sch\'emas: celles des morphismes plats, lisses ou \'etales.

\smallskip

Elles sont d\'efinies \`a partir de la notion de platitude d'un module sur un anneau.

\smallskip

Pour tout module $M$ sur un anneau commutatif $A$, le foncteur
$$
\begin{matrix}
M \otimes_A \bullet : &{\rm Mod}_A &\longrightarrow &{\rm Mod}_A \, , \hfill \\
&\hfill N &\longmapsto &M \otimes_A N
\end{matrix}
$$
admet pour adjoint \`a droite le foncteur
$$
\begin{matrix}
{\rm Mod}_A &\longrightarrow &{\rm Mod}_A \, , \hfill \\
\hfill L &\longmapsto &{\rm Hom}_A (M,L)
\end{matrix}
$$
donc il respecte les colimites arbitraires.

\smallskip

En revanche, il ne respecte pas n\'ecessairement les limites, ce qui conduit \`a poser:

\pagebreak

\begin{defn}\label{defII315}

Soit $A$ un anneau commutatif.

\begin{listeimarge}

\item Un $A$-module $M$ est dit ``plat'' si le foncteur
$$
\begin{matrix}
M \otimes_A \bullet : &{\rm Mod}_A &\longrightarrow &{\rm Mod}_A \, , \hfill \\
&\hfill N &\longmapsto &M \otimes_A N
\end{matrix}
$$
est exact, c'est-\`a-dire respecte les limites finies (ou, ce qui revient au m\^eme, les noyaux) en plus des colimites arbitraires.

\medskip

\item Une $A$-alg\`ebre commutative $B$ est dite ``plate'' si elle est plate comme $A$-module, c'est-\`a-dire si le foncteur
$$
\begin{matrix}
{\rm Mod}_A &\longrightarrow &{\rm Mod}_B \, , \hfill \\
\hfill M &\longmapsto &B \otimes_A M
\end{matrix}
$$
est exact.
\end{listeimarge}
\end{defn}

\medskip

\begin{remarksqed}
\begin{listeisansmarge}
\item Tout $A$-module libre, c'est-\`a-dire isomorphe \`a une somme directe $\underset{i \in I}{\bigoplus} A$ est plat sur $A$.

\medskip

\item En particulier, toute $A$-alg\`ebre de polyn\^omes
$$
A[(X_i)_{i\in I}]
$$
est plate sur $A$.

\medskip

\item De m\^eme, pour tout polyn\^ome $P \in A [X]$ dont le coefficient dominant est $1$
$$
P = X^n + a_{n-1} X^{n-1} + \cdots + a_1 X + a_0 \, ,
$$
l'alg\`ebre quotient $A[X]/(P)$ est plate sur $A$.

\medskip

\item Tout espace vectoriel sur un corps $K$ est plat sur $K$.

\smallskip

En particulier, toute $K$-alg\`ebre est plate sur $K$.

\medskip

\item Un ${\mathbb Z}$-module $M$ est plat si et seulement si il est sans torsion au sens que
$$
n \cdot m = 0 \Rightarrow m=0 \, , \qquad \forall \, n \in {\mathbb Z} - \{0\} \, , \ \forall \, m \in M \, .
$$

En particulier, un anneau commutatif $A$ est plat sur ${\mathbb Z}$ si et seulement si
$$
n \, a = 0 \Rightarrow a = 0 \, , \qquad \forall \, n \in {\mathbb Z} - \{0\} \, , \ \forall \, a \in A \, .
$$

\item Si $A \to B$ est un morphisme d'anneaux commutatifs et $M$ un module sur $A$, les deux foncteurs ${\rm Mod}_B \rightrightarrows {\rm Mod}_B$
$$
\begin{matrix}
&N &\longmapsto &(B \otimes_A M) \otimes_B N \\
\mbox{et} &N &\longmapsto &M \otimes_A N \hfill
\end{matrix}
$$
sont isomorphes.

\smallskip

Donc $B \otimes_A M$ est plat sur $B$ si $M$ est plat sur $A$.

\smallskip

En particulier, si $C$ est une $A$-alg\`ebre plate, $B \otimes_A C$ est plate sur $B$.

\medskip

\item Si $A \to B \to C$ sont deux morphismes d'anneaux commutatifs, les deux foncteurs ${\rm Mod}_A \rightrightarrows {\rm Mod}_A$
$$
\begin{matrix}
&M &\longmapsto &C \otimes_A M \hfill  \\
\mbox{et} &M &\longmapsto &C \otimes_B (B \otimes_A M) 
\end{matrix}
$$
sont isomorphes.

\smallskip

Donc $C$ est plat sur $A$ si $C$ est plat sur $B$ et $B$ est plat sur $A$.

\medskip

\item Pour tout \'el\'ement $f$ d'un anneau commutatif $A$, l'anneau $A_f = A[X] / (f \cdot X-1)$ est plat sur $A$.

\smallskip

Cela r\'esulte de la description du foncteur
$$
\begin{matrix}
{\rm Mod}_A &\longrightarrow &{\rm Mod}_{A_f} \, , \hfill \\
\hfill M &\longmapsto &M_f = A_f \otimes_A M
\end{matrix}
$$
dans la partie (ii) du lemme \ref{lemI46}.

\medskip

\item Pour toute famille $(f_i)_{i \in I}$ d'\'el\'ements d'un anneau commutatif $A$ telle que $\underset{i \in I}{\sum} f_i \, A = A$ ou, ce qui revient au m\^eme,
$$
{\rm Spec} (A) = \bigcup_{i \in I} {\rm Spec} (A_f) \, ,
$$
un module $M$ est plat sur $A$ si et seulement si, pour tout $i \in I$, $M_{f_i}$ est plat sur $A_{f_i}$.

\smallskip

L'implication en sens direct r\'esulte de la remarque (vi).

\smallskip

L'implication dans l'autre sens r\'esulte de ce que, d'apr\`es le lemme \ref{lemI46} (iii), les deux foncteurs ${\rm Mod}_A \rightrightarrows {\rm Mod}_A$
$$
\begin{matrix}
&N &\longmapsto &M \otimes_A N \hfill  \\
\mbox{et} &N &\longmapsto &\displaystyle {\rm eg} \left( \prod_{i \in I} M_{f_i} \otimes_{A_{f_i}} N_{f_i} \rightrightarrows \prod_{i,j \in I} M_{f_i f_j} \otimes_{A_{f_i f_j}} N_{f_i f_j} \right)
\end{matrix}
$$
sont isomorphes et du fait que les limites respectent les limites.

\smallskip

En particulier, tout $A$-module $M$ ``localement libre'' (au sens que ${\rm Spec} (A)$ est recouvert par des ouverts ${\rm Spec}(A_{f_i})$ tels que chaque $M_{f_i}$ soit un $A_{f_i}$-module libre) est plat.

\medskip

\noindent (x) Si $A \to B$ est un morphisme d'anneaux commutatifs et $(g_i)_{i \in I}$ est une famille d'\'el\'ements de $B$ telle que
$$
{\rm Spec} (B) = \bigcup_{i\in I} {\rm Spec} (B_{g_i}) \, ,
$$
alors $B$ est plat sur $A$ si et seulement si chaque $B_{g_i}$ est plat sur $A$.

\smallskip

L'implication en sens direct r\'esulte des remarques (viii) et (vii).

\smallskip

L'implication dans l'autre sens r\'esulte \`a nouveau du lemme \ref{lemI46} (iii), qui implique que les deux foncteurs ${\rm Mod}_A \rightrightarrows {\rm Mod}_B$
$$
\begin{matrix}
&M &\longmapsto &B \otimes_A M \hfill  \\
\mbox{et} &M &\longmapsto &\displaystyle {\rm eg} \left( \prod_{i \in I} B_{g_i} \otimes_A M \rightrightarrows \prod_{i,j \in I} B_{g_i g_j} \otimes_A M \right)
\end{matrix}
$$
sont isomorphes, et du fait que les limites respectent les limites.
\end{listeisansmarge}
\end{remarksqed}

\pagebreak

La consid\'eration des diff\'erentielles relatives d'un anneau sur un autre fournit un important crit\`ere de platitude:

\begin{prop}\label{propII316}

Soit $A \to B$ un morphisme plat et de type fini d'anneaux commutatifs.

\smallskip

Si $f_1 , \cdots , f_n$, sont des \'el\'ements de $B$ dont les diff\'erentielles ${\rm d} f_1 , \cdots , {\rm d} f_n \in \Omega_{B/A}$ forment une base du $B$-module $\Omega_{B/A}$, avec donc
$$
\Omega_{B/A} = \bigoplus_{1 \leq i \leq n} B \cdot {\rm d} f_i \, ,
$$
alors l'unique morphisme de $A$-alg\`ebres
$$
A[X_1 , \cdots , X_n] \longrightarrow B
$$
qui envoie les variables $X_1 , \cdots , X_n$ sur $f_1 , \cdots , f_n$ est plat.
\end{prop}

\medskip

\begin{demo}

Elle est laiss\'e au lecteur \`a titre de probl\`eme. 

\end{demo}

\medskip

Nous pouvons maintenant d\'efinir les classes g\'eom\'etriques des morphismes plats, des morphismes lisses et des morphismes \'etales:

\begin{defn}\label{defII317}

Un morphisme $X \to S$ de sch\'emas est dit
\begin{enumerate}
\item[(10)] plat si pour tout ouvert affine $U \cong {\rm Spec} (B)$ de $X$ s'envoyant dans un ouvert affine $V \cong {\rm Spec} (A)$ de $S$, le morphisme induit $A \to B$ est plat,
\item[(11)] lisse s'il est localement de pr\'esentation finie et plat, et si $X$ admet un recouvrement par des ouverts affines
$$
U_i \cong {\rm Spec} (B_i)
$$
s'envoyant dans des ouverts affines de $S$
$$
V_i \cong {\rm Spec} (A_i)
$$
tels que chaque $B_i$-module
$$
\Omega_{B_i / A_i}
$$
soit libre et admette pour base les diff\'erentielles
$$
{\rm d} f_1 , \cdots , {\rm d} f_{n_i}
$$
d'\'el\'ements $f_1 , \cdots , f_{n_i} \in B_i$,
\item[(12)] \'etale s'il est localement de pr\'esentation finie et plat, et si le morphisme diagonal
$$
X \longrightarrow X \times_S X
$$
est une immersion ouverte.
\end{enumerate}
\end{defn}

\pagebreak

\begin{remarksqed}
\begin{listeisansmarge}	
\item Il r\'esulte des remarques (vii), (viii) et (x) suivant la d\'efinition \ref{defII315} qu'un morphisme de sch\'emas
$$
X \longrightarrow S
$$
est plat si $X$ admet un recouvrement par des ouverts affines
$$
U_i \cong {\rm Spec} (B_i)
$$
s'envoyant dans des ouverts affines $V_i \cong {\rm Spec} (A_i)$ de $S$ et tels que chaque morphisme induit $A_i \to B_i$ soit plat.

\smallskip

Il r\'esulte de cette remarque et des remarques (vi) et (vii) suivant la d\'efinition \ref{defII315} que la classe des morphismes plats est g\'eom\'etrique.

\medskip

\item Les m\^emes arguments que ceux de la remarque (i) suivant la d\'efinition \ref{defII311} montrent que les morphismes $X \to S$ dont le morphisme diagonal associ\'e
$$
X \longrightarrow X \times_S X
$$
est une immersion ouverte forment une classe g\'eom\'etrique.

\smallskip

Donc les morphismes \'etales forment une classe g\'eom\'etrique.

\medskip

\item Il r\'esulte du corollaire \ref{corII310} (iii) qu'un morphisme $X \to S$ plat et localement de pr\'esentation finie est \'etale si et seulement si, pour tout ouvert affine $U \cong {\rm Spec} (B)$ de $X$ s'envoyant dans un ouvert affine $V \cong {\rm Spec}(A)$ de $S$, on a
$$
\Omega_{B/A} = 0 \, .
$$

\item Pour tout morphisme d'anneaux commutatifs $A \to B$, tout \'el\'ement $f \in B$ et tout module $M$ sur $B_f$, toute diff\'erentielle $A$-lin\'eaire
$$
{\rm d} : B \longrightarrow M
$$
se prolonge de mani\`ere unique en une diff\'erentielle $A$-lin\'eaire
$$
{\rm d} : B_f \longrightarrow M
$$
par la formule
$$
{\rm d} (f^{-n} b) = -n f^{n-1} b \cdot {\rm d} f + f^{-n} \cdot {\rm d} b \, , \qquad \forall \, n \in {\mathbb N} \, , \ \forall \, b \in B \, .
$$
D'o\`u un isomorphisme canonique
$$
\Omega_{B_f/A} \cong B_f \otimes_B \Omega_{B/A} \, .
$$
Dans la d\'efinition (11), on peut donc remplacer le recouvrement ouvert $(U_i)$ par n'importe quel raffinement.

\medskip

\item Pour tous morphismes d'anneaux commutatifs $A \to B$ et $A \to C$, et pour tout module $M$ sur $C \otimes_A B$, se donner une diff\'erentielle $C$-lin\'eaire
$$
{\rm d} : C \otimes_A B \longrightarrow M
$$
\'equivaut \`a se donner une diff\'erentielle $A$-lin\'eaire
$$
{\rm d} : B \longrightarrow M \, .
$$
D'o\`u un isomorphisme canonique
$$
\Omega_{C \otimes_A B/C} \cong C \otimes_A \Omega_{B/A} \, .
$$
C'est pourquoi la classe des morphismes lisses est stable par changement de base.

\medskip

\item Le compos\'e d'une projection affine
$$
{\mathbb A}_S^n = ({\rm Spec} \, A[X_1 , \cdots , X_n ]) \longrightarrow {\rm Spec} (A) = S
$$
et d'un morphisme \'etale de la forme
$$
X = {\rm Spec} (B) \longrightarrow {\mathbb A}_S^n
$$
est lisse.

\smallskip

En effet, il est d'une part plat et de pr\'esentation finie.

\smallskip

D'autre part, le carr\'e commutatif
$$
\xymatrix{
X \ar[d] \ar[r] &X \times_S X \ar[d] \\
{\mathbb A}_S^n \ar[r] &{\mathbb A}_S^n \times_S {\mathbb A}_S^n
}
$$
identifie $X$ \`a un ouvert du produit fibr\'e
$$
{\mathbb A}_S^n \times_{{\mathbb A}_S^n \times_S {\mathbb A}_S^n} X \times_S X = X \times_{{\mathbb A}_S^n} X \, .
$$
Comme le foncteur
$$
(B \otimes_A B) \otimes_{A[X_1,\cdots , \, X_n] \otimes_A A [X_1 , \cdots , \, X_n]} \bullet
$$
est exact, on en d\'eduit que les noyaux
$$
\begin{matrix}
&I &= &{\rm Ker} (B \otimes_A B \longrightarrow B) \hfill \\
\mbox{et} &J &= &{\rm Ker} (A[X_1 , \cdots , X_n] \otimes_A A [X_1 , \cdots , X_n] \longrightarrow A[X_1 , \cdots , X_n])
\end{matrix}
$$
sont reli\'es par la formule
$$
I/I^2 = (B \otimes_A B) \otimes_{A[X_1 , \cdots , \, X_n] \otimes_A A[X_1,\cdots , \, X_n]} J/J^2 = B \otimes_{A[X_1,\cdots , \, X_n]} J/J^2
$$
et donc que les images $f_1 , \cdots , f_n$ de $X_1 , \cdots , X_n$ dans $B$ satisfont la formule
$$
\Omega_{B/A} = \bigoplus_{1 \leq i \leq n} B \cdot {\rm d} f_i \, .
$$

\item Il r\'esulte de la remarque pr\'ec\'edente qu'un morphisme
$$
X \longrightarrow S
$$
est lisse s'il peut \^etre recouvert par des ouverts $U_i$ tels que chaque morphisme $U_i \to S$ s'\'ecrive comme le compos\'e d'une projection affine
$$
{\mathbb A}^{n_i} \times S \longrightarrow S
$$
et d'un morphisme \'etale
$$
U_i \longrightarrow {\mathbb A}^{n_i} \times S \, .
$$

Si tous les entiers $n_i$ sont \'egaux \`a un m\^eme entier $n$, on dit que $X$ est lisse sur $S$ de dimension relative $n$.

\smallskip

En particulier, les morphismes \'etales sont les morphismes lisses de dimension relative $0$.

\medskip

\item La remarque pr\'ec\'edente montre que la classe des morphismes lisses est stable par composition et donc g\'eom\'etrique.

\smallskip

En effet, si deux morphismes $U \to {\mathbb A}^n \times V$ et $V \to {\mathbb A}^m \times W$ sont \'etales, alors le morphisme compos\'e
$$
U \longrightarrow {\mathbb A}^n \times V \longrightarrow {\mathbb A}^n \times {\mathbb A}^m \times W
$$
est \'etale. 
\end{listeisansmarge}
\end{remarksqed}

\subsection{Des notions g\'eom\'etriques de recouvrement de sch\'emas}\label{subsec234}

\smallskip

En plus de la notion usuelle de recouvrement pour la topologie de Zariski, nous allons introduire plusieurs notions g\'eom\'etriques de recouvrement dans la cat\'egorie ${\rm Sch}$ des sch\'emas: les recouvrements \'etales, lisses, fppf (fid\`element plats de pr\'esentation finie) et fpqc (fid\`element plats quasi-compacts).

\smallskip

Toutes sont fond\'ees sur la notion de morphisme fid\`element plat ou de classe de morphismes fid\`element plate quasi-compacte:

\begin{defn}\label{defII318}
\begin{listeimarge}
\item Un morphisme d'anneaux commutatifs
$$
A \longrightarrow B
$$
est dit ``fid\`element plat'' s'il est plat et que tout id\'eal premier de $A$ est l'image r\'eciproque d'au moins un id\'eal premier de $B$.

\medskip

\item Un morphisme de sch\'emas
$$
X \longrightarrow S
$$
est dit ``fid\`element plat quasi-compact'' (fpqc) s'il est plat, surjectif et quasi-compact.

\medskip

\item Une famille de morphismes de sch\'emas
$$
X_i \longrightarrow S \, , \qquad i \in I \, ,
$$
est dite ``fid\`element plate quasi-compacte'' (fpqc) si le morphisme induit
$$
\coprod_{i \in I} X_i \longrightarrow S
$$
est fid\`element plat quasi-compact, autrement dit si

\bigskip

\noindent $\left\{\begin{matrix}
\bullet &\mbox{chaque $X_i \to S$ est plat et quasi-compact,} \hfill \\
{ \ } \\
\bullet &\mbox{$S$ est la r\'eunion des images des $X_i$,} \hfill \\
{ \ } \\
\bullet &\mbox{la famille des $X_i \to S$ est localement finie au sens que pour tout ouvert affine $V = {\rm Spec} (A)$ de $S$,} \hfill \\
&\mbox{tous les $X_i \times_S V$ sont vides sauf un nombre fini.} \hfill
\end{matrix} \right. $
\end{listeimarge}
\end{defn}

\begin{remarksqed}
\begin{listeisansmarge}
\item La classe des morphismes surjectifs est g\'eom\'etrique. En effet, elle contient les isomorphismes et elle est stable par composition et par changement de base.

\smallskip

Par cons\'equent, la classe des morphismes fid\`element plats est aussi g\'eom\'etrique.

\medskip

\item Un morphisme de sch\'emas $X \to S$ est fid\`element plat quasi-compact si, pour tout ouvert affine $V = {\rm Spec} (A)$ de $S$, son image r\'eciproque $X \times_S V$ est une r\'eunion finie d'ouverts affines $U_i = {\rm Spec} (B_i)$, $1 \leq i \leq n$, tels que le morphisme
$$
A \longrightarrow B_1 \times \cdots \times B_n
$$
soit fid\`element plat. 
\end{listeisansmarge}
\end{remarksqed}

\pagebreak

L'expression ``fid\`element plat'' a \'et\'e choisie pour faire r\'ef\'erence \`a la propri\'et\'e \'enonc\'ee dans la partie (iii) du lemme suivant:

\begin{lem}\label{lemII319}

Soit $A \to B$ un morphisme d'anneaux commutatifs qui est fid\`element plat.

\smallskip

Alors:

\begin{listeimarge}

\item Pour tout module $M$ sur $A$, le morphisme d\'eduit de $A \to B$
$$
M \longrightarrow B \otimes_A M
$$
est injectif.

\smallskip

En particulier, le morphisme $A \to B$ est injectif.

\medskip

\item Un module $M$ sur $A$ est $0$ si (et seulement si) $B \otimes_A M$ est $0$.

\medskip

\item Un complexe de modules sur $A$
$$
M_1 \xrightarrow{ \ u \ } M_2 \xrightarrow{ \ v \ } M_3
$$
v\'erifiant $v \circ u = 0$, est exact au sens que
$$
{\rm Im} (u) = {\rm Ker} (v)
$$
si et seulement si le complexe induit
$$
B \otimes_A M_1 \xrightarrow{ \ {\rm id}_B \otimes u \ } B \otimes_A M_2 \xrightarrow{ \ {\rm id}_B \otimes v \ } B \otimes_A M_3
$$
est exact.

\smallskip

En particulier, un morphisme de modules sur $A$
$$
M \longrightarrow N
$$
est injectif [resp. surjectif] si et seulement si le morphisme induit
$$
B \otimes_A M \longrightarrow B \otimes_A N
$$
est injectif [resp. surjectif].
\end{listeimarge}
\end{lem}

\begin{demo}
\begin{listeisansmarge}
\item Un \'el\'ement non nul $m$ de $M$ peut \^etre vu comme un morphisme de modules sur $A$
$$
A \longrightarrow M
$$
dont le noyau $I$ est un id\'eal autre que $A$. Donc $I$ est contenu dans un id\'eal premier $p$ de $A$ et, par hypoth\`ese, il existe un id\'eal premier $q$ de $B$ dont l'image r\'eciproque par le morphisme $A \to B$ est $p$.

\smallskip

Alors $B/q$ est un quotient de $B \otimes_A (A/I)$.

\smallskip

Comme $B$ est plat sur $A$, le morphisme
$$
B \otimes_A (A/I) \longrightarrow B \otimes_A M
$$
est injectif.

\smallskip

Cela montre que l'image de $m$ dans $B \otimes_A M$ n'est pas $0$.

\medskip

\item Si $M=0$, alors bien s\^ur $B \otimes_A M = 0$.

\smallskip

L'implication en sens inverse r\'esulte de ce que le morphisme $M \longrightarrow B \otimes_A M$ est toujours injectif.

\medskip

\item Comme $v \circ u = 0$, on peut introduire le module quotient
$$
H = {\rm Ker} (v) / {\rm Im} (u) \, .
$$
La platitude de $B$ sur $A$ signifie que le foncteur $B \otimes_A \bullet$ est exact. On en d\'eduit que
$$
B \otimes_A H = {\rm Ker} ({\rm id}_B \otimes v) / {\rm Im} ({\rm id}_B \otimes u) \, .
$$
Or, d'apr\`es (ii), $H$ est $0$ si et seulement si $B \otimes_A H$ est $0$.

\end{listeisansmarge}
\end{demo}

\medskip

On montre encore \`a partir de ce lemme:

\begin{lem}\label{lemII320}

Soit $A \xrightarrow{ \ u \ } B$ un morphisme d'anneaux commutatifs qui est fid\`element plat.

\smallskip

Notons $p_1 , p_2 : B \rightrightarrows B \otimes_A B$ les deux morphismes
$$
u \otimes {\rm id}_B \, , \qquad {\rm id}_B \otimes u \, ,
$$
et $q_1 , q_2 , q_3 : B \overset{\mbox{$\rightarrow$}}{\underset{\mbox{$\rightarrow$}}\rightarrow} B \otimes_A B \otimes_A B$ les trois morphismes
$$
u \otimes {\rm id}_B \otimes {\rm id}_B \, , \quad {\rm id}_B \otimes u \otimes {\rm id}_B \, , \quad  {\rm id}_B \otimes {\rm id}_B \otimes u \, .
$$

Alors:

\begin{listeimarge}

\item Tout $A$-module $M$ est canoniquement isomorphe \`a
$$
{\rm eg} \left( \raisebox{.7ex}{\xymatrix{ B \otimes_A M  \dar[rr]^-{^{^{\mbox{\scriptsize$p_1 \otimes {\rm id}_M$}}}}_-{p_2 \, \otimes \, {\rm id}_M} &&B \otimes_A B \otimes_A M}} \right).
$$

\medskip

\item Soit $M'$ un $B$-module muni d'un isomorphisme $B \otimes_A B$-lin\'eaire
$$
\sigma : (B \otimes_A B) \otimes_{p_1 , B} M' \xrightarrow{ \ \sim \ } (B \otimes_A B) \otimes_{p_2 , B} M'
$$
tel que le triangle induit
$$
\xymatrix{
(B \otimes_A B \otimes_A B) \otimes_{q_1 , B} M' \ar[rd]^{\sim} \ar[rr]^{\sim} &&(B \otimes_A B \otimes_A B) \otimes_{q_3 , B} M' \\
&(B \otimes_A B \otimes_A B) \otimes_{q_2 , B} M' \ar[ru]^{\sim}
}
$$
soit commutatif.

\smallskip

Soit
$$
M = {\rm eg} \left( M' \rightrightarrows B \otimes_A M' = (B \otimes_A B) \otimes_{p_2,B} M \right)
$$
le $A$-module \'egalisateur du morphisme
$$
u \otimes {\rm id}_{M'}
$$
et de son compos\'e avec
$$
B \otimes_A M' = (B \otimes_A B) \otimes_{p_1,B} M' \xrightarrow{ \ {\sigma \atop \sim} \ } (B \otimes_A B) \otimes_{p_2,B} M' = B \otimes_A M' \, .
$$

Alors le morphisme canonique de modules sur $B$
$$
B \otimes_A M \longrightarrow M'
$$
est un isomorphisme.

\medskip

\item Le foncteur qui associe \`a tout module $M$ sur $A$ le module sur $B$
$$
M' = B \otimes_A M
$$
muni de l'isomorphisme canonique
$$
\xymatrix{
(B \otimes_A B) \otimes_{p_1,B} M' \ar@{=}[d] \ar[r]^{\sim} &(B \otimes_A B) \otimes_{p_2,B} M' \ar@{=}[d] \\
(B \otimes_A B) \otimes_A M &(B \otimes_A B) \otimes_A M
}
$$

\noindent est une \'equivalence de la cat\'egorie ${\rm Mod}_A$ des modules sur $A$ vers la cat\'egorie des modules $M'$ sur $B$ munis d'un isomorphisme
$$
\sigma : (B \otimes_A B) \otimes_{p_1 , B} M' \xrightarrow{ \ \sim \ } (B \otimes_A B) \otimes_{p_2 , B} M'
$$
qui satisfait la condition de (ii).
\end{listeimarge}
\end{lem}

\begin{demo}
\begin{listeisansmarge}
\item Comme $B$ est fid\`element plat sur $A$, pour montrer que la suite
$$
0 \longrightarrow M \xrightarrow{ \ u \otimes {\rm id}_M \ } B \otimes_A M \xrightarrow{ \ p_1 \otimes {\rm id}_M - p_2 \otimes {\rm id}_M \ } B \otimes_A B \otimes_A M
$$
est exacte, il suffit de montrer qu'il en est ainsi de la suite
$$
0 \longrightarrow B \otimes_A M \xrightarrow{ \ {\rm id}_B \otimes u \otimes {\rm id}_M \ } B \otimes_A B \otimes_A M \xrightarrow{ \ ({\rm id}_B \otimes p_1 - {\rm id}_B \otimes p_2) \, \otimes \,  {\rm id}_M \ } B \otimes_A B \otimes_A B \otimes_A M \, .
$$
Or le morphisme
$$
\begin{matrix}
\Delta : &B \otimes_A B &\longrightarrow &B \, , \\
&\hfill b \otimes b' &\longmapsto &bb'
\end{matrix}
$$
satisfait la formule
$$
\Delta \circ ({\rm id}_B \otimes u) = {\rm id}_B
$$
et a fortiori
$$
(\Delta \otimes {\rm id}_M) \circ ({\rm id}_B \otimes u \otimes {\rm id}_M) = {\rm id}_{B \otimes_A M} \, .
$$
De plus, si un \'el\'ement de $B \otimes_A B \otimes_A M$
$$
\sum_{1 \leq i \leq n} b_i \otimes b'_i \otimes m_i
$$
est dans le noyau de $({\rm id}_B \otimes p_1 - {\rm id}_B \otimes p_2) \otimes {\rm id}_M$ c'est-\`a-dire v\'erifie la formule
$$
\sum_{1 \leq i \leq n} b_i \otimes 1 \otimes b'_i \otimes m_i = \sum_{1 \leq i \leq n} b_i \otimes b'_i \otimes 1 \otimes m_i
$$
alors, appliquant $\Delta \otimes {\rm id}_B \otimes {\rm id}_M$, il s'\'ecrit
$$
\sum_{1 \leq i \leq n} b_i \otimes b'_i \otimes m_i = \sum_{1 \leq i \leq n} b_i b'_i \otimes 1 \otimes m_i
$$
ce qui signifie qu'il est dans l'image de
$$
{\rm id}_B \otimes u \otimes {\rm id}_M \, .
$$

\item Pour montrer que
$$
B \otimes_A M \longrightarrow M'
$$
est un isomorphisme, il suffit de montrer qu'il en est ainsi du morphisme image par le foncteur $B \otimes_A \bullet$
$$
B \otimes_A B \otimes_A M \longrightarrow B \otimes_A M' \, .
$$
Or celui-ci s'\'ecrit \'egalement
$$
(B \otimes_A B) \otimes_{p_1,B} (B \otimes_A M) \longrightarrow (B \otimes_A B) \otimes_{p_1 , B} M'
$$
avec
$$
B \otimes_A M = {\rm eg} \left( (B \otimes_A B) \otimes_{p_1,B} M' \rightrightarrows (B \otimes_A B) \otimes_{p_1 , B , p_1} (B \otimes_A B) \otimes_{p_1,B} M' \right).
$$
Donc (ii) se d\'eduit de (i) appliqu\'e au $B$-module $M'$ et au morphisme fid\`element plat
$$
p_1 : B \longrightarrow B \otimes_A B \, .
$$
\item r\'esulte de (i) et (ii). 

\end{listeisansmarge}
\end{demo}

\medskip

On d\'eduit de la partie (iii) de ce lemme:

\begin{cor}\label{corII321}

Consid\'erons une famille fid\`element plate quasi-compacte
$$
(S_i \longrightarrow S)_{i \in I}
$$
de morphismes vers un sch\'ema $S$.

\smallskip

Alors:

\begin{listeimarge}

\item Le foncteur
$$
(X \longrightarrow S) \longmapsto (X \times_S S_i \longrightarrow S_i)_{i \in I}
$$
d\'efinit une \'equivalence de la cat\'egorie des sch\'emas sur $S$
$$
X \longrightarrow S
$$
qui sont affines [resp. finis, resp. des sous-sch\'emas ferm\'es] vers la cat\'egorie des familles de sch\'emas sur les $S_i$
$$
(X_i \longrightarrow S_i)_{i \in I}
$$
qui sont affines [resp. finis, resp. des sous-sch\'emas ferm\'es] et sont munis d'isomorphismes au-dessus des $S_i \times_S S_j$
$$
\xymatrix{
X_i \times_{S_i} (S_i \times_S S_j) \ar[rd] \ar[rr]^{\sim} &&(S_i \times_S S_j) \times_{S_j} X_j \ar[ld] \\
&S_i \times_S S_j
}
$$
tels que les triangles
$$
\xymatrix{
X_i \times_{S_i} (S_i \times_S S_j \times_S S_k) \ar[rd]^{\sim} \ar[rr]^{\sim} &&X_k \times_{S_k} (S_i \times_S S_j \times_S S_k) \\
&X_j \times_{S_j} (S_i \times_S S_j \times_S S_k) \ar[ru]^{\sim}
}
$$
soient commutatifs.

\medskip

\item Pour tout sch\'ema $X$ sur $S$ et
$$
\begin{matrix}
\hfill X_i &= &X \times_S S_i \, , \hfill &\forall \, i \, , \hfill \\
X_{i,j} &= &X \times_S (S_i \times_S S_j) \, , &\forall \, i,j \, ,
\end{matrix}
$$
et pour tout ouvert $U$ de $X$ et ses images r\'eciproques
$$
\begin{matrix}
\hfill U_i &\subset &X_i \, , \hfill &i \in I \, , \hfill \\
U_{i,j} &\subset &X_{i,j} \, , &i,j \in I \, ,
\end{matrix}
$$
on a
$$
{\mathcal O}_X (U) = {\rm eg} \left( \prod_{i \in I} {\mathcal O}_{X_i} (U_i) \rightrightarrows \prod_{i,j \in I} {\mathcal O}_{X_{i,j}} (U_{i,j})\right) .
$$
\end{listeimarge}
\end{cor}


\begin{demo}

Il suffit de traiter le cas o\`u $S = {\rm Spec} (A)$ et
$$
\coprod_{i \in I} S_i = {\rm Spec} (B)
$$
sont deux sch\'emas affines reli\'es par un morphisme fid\`element plat
$$
A \longrightarrow B \, .
$$
\begin{listeisansmarge}
\item La cat\'egorie des sch\'emas affines sur $S = {\rm Spec} (A)$ est l'oppos\'ee de celle des $A$-modules $M$ munis d'un morphisme de multiplication
$$
M \times M \longrightarrow M
$$
et d'un \'el\'ement unit\'e $1 \in M$ c'est-\`a-dire d'un morphisme
$$
A \longrightarrow M
$$
v\'erifiant les axiomes de d\'efinition des $A$-alg\`ebres commutatives.

\smallskip

Ces axiomes s'expriment comme des propri\'et\'es de commutativit\'e de diagrammes impliquant les objets $M \times M \times M$, $M \times M , M$ et $A \times M$ donc ils se lisent aussi bien sur leurs transform\'es par le foncteur $B \otimes_A \bullet$.

\smallskip

D'o\`u la conclusion dans le cas de la cat\'egorie des sch\'emas affines sur $S = {\rm Spec} (A)$.

\smallskip

Le cas des sch\'emas finis sur $S = {\rm Spec} (A)$ s'en d\'eduit puisqu'un $A$-module $M$ est de type fini si et seulement si le $B$-module $B \otimes_A M$ est de type fini.

\smallskip

De m\^eme, un morphisme $A \to M$ est surjectif si et seulement si le morphisme induit $B \to B \otimes_A M$ est surjectif, ce qui conclut le cas des immersions ferm\'ees.

\medskip

\item Comme la famille
$$
(U_i \longrightarrow U)_{i \in I}
$$
est fid\`element plate et
$$
U_{i,j} = U_i \times_U U_j \, , \qquad \forall \, i,j \in I \, ,
$$
il suffit de traiter le cas o\`u $U = X = S = {\rm Spec} (A)$ et $\underset{i \in I}{\coprod} \, U_i = {\rm Spec} (B)$.

\smallskip

On est ramen\'e \`a l'assertion (i) du lemme \ref{lemII320} dans le cas particulier $M=A$.\end{listeisansmarge} 
\end{demo}

\pagebreak

On introduit la d\'efinition g\'en\'erale suivante:

\begin{defn}\label{defII322}

Soient ${\mathcal C}$ une cat\'egorie localement petite et ${\mathcal M}$ une classe g\'eom\'etrique de morphismes de ${\mathcal C}$.

\smallskip

Soit {\rm (R)} une propri\'et\'e des familles de morphismes de ${\mathcal C}$ ayant le m\^eme but
$$
(S_i \longrightarrow S)_{i \in I} \, .
$$

On dit que la classe ${\mathcal M}$ est caract\'eris\'ee ${\rm R}$-localement (ou que les propri\'et\'es qui d\'efinissent ${\mathcal M}$ sont ${\rm R}$-locales) si un morphisme carrable de ${\mathcal C}$
$$
X \longrightarrow S
$$
est \'el\'ement de ${\mathcal M}$ d\`es qu'il existe une famille de morphismes de but $S$ v\'erifiant {\rm (R)}
$$
(S_i \longrightarrow S)_{i \in I}
$$
telle que les morphismes d\'eduits de $X \to S$ par les changements de base $S_i \to S$, $i \in I$,
$$
X \times_S S_i \longrightarrow S_i
$$
soient \'el\'ements de ${\mathcal M}$.
\end{defn}

\medskip

\begin{remarkqed}

Dans la pratique, (R) est presque toujours une notion g\'eom\'etrique de recouvrement au sens de la d\'efinition~\ref{defII32}. Elle est alors associ\'ee \`a une classe g\'eom\'etrique ${\mathcal M}'$ de morphismes de ${\mathcal C}$ mais les classes ${\mathcal M}$ et ${\mathcal M}'$ peuvent \^etre compl\`etement diff\'erentes. 
\end{remarkqed}

\medskip

Cette d\'efinition g\'en\'erale \'etant pos\'ee, on peut \'enoncer:

\begin{thm}\label{thmII323}

Dans la cat\'egorie ${\rm Sch}$ des sch\'emas, notons (fpqc) la propri\'et\'e des familles de morphismes de m\^eme but
$$
(S_i \longrightarrow S)_{i \in I}
$$
d'\^etre fid\`element plates quasi-compactes.

\smallskip

Alors les classes g\'eom\'etriques suivantes de morphismes de ${\mathcal M}$ sont (fpqc)-locales:
\begin{enumerate}
\item[(0)] La classe des isomorphismes.
\item[(1)] La classe des morphismes de pr\'esentation finie [resp. de type fini].
\item[(2)] La classe des morphismes localement de pr\'esentation finie [resp. localement de type fini].
\item[(3)] La classe des morphismes quasi-compacts.
\item[(4)] La classe des morphismes affines.
\item[(5)] La classe des morphismes finis.
\item[(6)] La classe des immersions ferm\'ees.
\item[(7)] La classe des immersions localement ferm\'ees [resp. ouvertes].
\item[(8)] La classe des morphismes s\'epar\'es.
\item[(9)]  La classe des morphismes propres.
\item[(10)] La classe des morphismes plats.
\item[(11)] La classe des morphismes lisses.
\item[(12)] La classe des morphismes \'etales.
\end{enumerate}
\end{thm}

\begin{demo}

Dans tous les cas, il suffit de d\'emontrer que si
$$
S' = {\rm Spec} (B) \longrightarrow {\rm Spec} (A) = S
$$
est un morphisme fid\`element plat de sch\'emas affines, et si
$$
X \longrightarrow {\rm Spec} (A) = S
$$
est un morphisme d'un sch\'ema $X$ vers le sch\'ema affine $S$ tel que
$$
X \times_S S' \longrightarrow S'
$$
soit \'el\'ement de l'une de ces classes ${\mathcal M}$, alors le morphisme $X \to S$ est \'el\'ement de la m\^eme classe ${\mathcal M}$.

\begin{listeisansmarge}
\item[(0)] Commen\c cons par le cas o\`u $X' \to S'$ est un isomorphisme.

\smallskip

Montrons d'abord que l'application sous-jacente \`a $X \to S$ est une bijection.

\smallskip

Soient $p \in {\rm Spec} (A) = S$ un point de $S$, $p' \in {\rm Spec} (B) = S'$ un point de $S'$ au-dessus de $p$ et $\kappa_p$, $\kappa_{p'}$ les corps des fractions des anneaux int\`egres $A/p$ et $B/p'$. Ainsi,
$$
X_p = X \times_S {\rm Spec} (\kappa_p) \qquad \mbox{et} \qquad X'_{p'} = X' \times_{S'} {\rm Spec} (\kappa_{p'})
$$
sont les fibres de $X$ et $X'$ au-dessus de $p$ et $p'$.

\smallskip

Le carr\'e
$$
\xymatrix{
X'_{p'} \ar[d] \ar[r] &X_p \ar[d] \\
{\rm Spec} (\kappa_{p'}) \ar[r] &{\rm Spec} (\kappa_p)
}
$$
est encore cart\'esien et ses fl\`eches horizontales sont fid\`element plates, en particulier surjectives.

\smallskip

De plus, ${\rm Spec} (\kappa_p)$ et ${\rm Spec} (\kappa_{p'})$ ont un unique point, et il en est de m\^eme de $X'_{p'}$ puisque $X'_{p'} \to {\rm Spec} (\kappa_{p'})$ est un isomorphisme.

\smallskip

Donc $X_p$ aussi a un unique \'el\'ement.

\smallskip

Ainsi, l'application $X \to S$ est bien une bijection.

\smallskip

Il r\'esulte du corollaire \ref{corII321}~(i) et de l'isomorphisme
$$
X' \xrightarrow{ \ \sim \ } S'
$$
que les sous-sch\'emas ferm\'es de $X$ et de $S$ se correspondent. En particulier, la bijection $X \to S$ est un hom\'eomorphisme.

\smallskip

Enfin, il r\'esulte du corollaire \ref{corII321}~(ii) et de l'isomorphisme $X' \xrightarrow{ \ \sim \ } S'$ que l'hom\'eomorphisme $X \to S$ identifie les faisceaux de structure ${\mathcal O}_X$ et ${\mathcal O}_S$ de $X$ et $S$.

\smallskip

Cela termine de traiter le cas (0) des isomorphismes.

\medskip

\item[(4)(5)(6)] Les cas des morphismes ferm\'es, des morphismes finis et des immersions ferm\'ees r\'esultent du corollaire \ref{corII321}~(i) combin\'e avec le cas (0).

\medskip

\item[(8)] Le cas des morphismes s\'epar\'es se d\'eduit du cas (6) des immersions ferm\'ees puisque le morphisme diagonal
$$
X' \longrightarrow X' \times_{S'} X'
$$
se d\'eduit du morphisme diagonal
$$
X \longrightarrow X \times_S X
$$
par le changement de base fid\`element plat
$$
S' \longrightarrow S
$$
ou
$$
X' \times_{S'} X' \longrightarrow X \times_S X \, .
$$

\item[(7)] Traitons le cas o\`u $X' \to S' = {\rm Spec} (B)$ est une immersion ouverte.

\smallskip

Soient
$$
I = {\rm Ker} \left( A \longrightarrow \prod_{p' \in {\rm Spec} (B) \atop p' \notin X'} B / p' \right)
$$
et
$$
Y = {\rm Spec} (A/I) \xhookrightarrow{ \ \ \ } {\rm Spec} (A) = S
$$
le sous-sch\'ema ferm\'e de $S$ d\'efini par l'id\'eal $I$.

\smallskip

Alors on a
$$
B \otimes_A I = \bigcap_{p' \in {\rm Spec} (B') \atop p' \notin X'} p'
$$
si bien que $X'$ s'identifie \`a l'ouvert compl\'ementaire du sous-sch\'ema ferm\'e de $S'$
$$
Y \times_S S' \xhookrightarrow{ \ \ \ } S' \, .
$$
On conclut d'apr\`es le cas (0) que $X$ s'identifie \`a l'ouvert compl\'ementaire du sous-sch\'ema ferm\'e $Y$ de $S$.

\smallskip

Puis traitons le cas o\`u $X' \to S'$ est une immersion localement ferm\'ee.

\smallskip

Soient
$$
I = {\rm Ker} (A \longrightarrow {\mathcal O}_X (X))
$$
et
$$
Y = {\rm Spec} (A/I) \xhookrightarrow{ \ \ \ } {\rm Spec} (A) = S 
$$
le sous-sch\'ema ferm\'e de $S$ d\'efini par l'id\'eal $I$.

\smallskip

Alors on a
$$
B \otimes_A I = {\rm Ker} (B \longrightarrow {\mathcal O}_{X'} (X'))
$$
si bien que $Y \times_S S'$ est le plus petit sous-sch\'ema ferm\'e de $S'$ \`a travers lequel se factorise $X' \to S'$.

\smallskip

Donc $X' \to Y \times_S S'$ est une immersion ouverte et, par cons\'equent, $X \to Y$ est une immersion ouverte.

\smallskip

Cela ach\`eve de traiter le cas (7) des immersions ouvertes et des immersions localement ferm\'ees.

\medskip

\item[(3)] Traitons le cas des morphismes quasi-compacts.

\smallskip

Si les $U_i$, $i \in I$, sont des ouverts de $X$ et les $U'_i$ sont leurs images r\'eciproques dans $X'$, il r\'esulte de la surjectivit\'e des applications sous-jacentes aux $U'_i \to U_i$ que
$$
X = \bigcup_{i \in I} U_i
$$
si et seulement si
$$
X' = \bigcup_{i \in I} U'_i \, .
$$
Donc $X \to S$ est quasi-compact si $X' \to S'$ est quasi-compact.

\medskip

\item[(2)] Le cas des morphismes localement de pr\'esentation finie [resp. localement de type fini] r\'esulte de ce que une $A$-alg\`ebre commutative $C$ est de pr\'esentation finie [resp. de type fini] si et seulement si la $B$-alg\`ebre $B \otimes_A C$ est de pr\'esentation finie [resp. de type fini].

\medskip

\item[(1)] Le cas des morphismes de pr\'esentation finie [resp. de type fini] r\'esulte de la combinaison des cas (2) et (3).

\medskip

\item[(9)] Traitons le cas des morphismes propres.

\smallskip

Le cas (1) des morphismes de type fini et le cas (8) des morphismes s\'epar\'es \'etant connus, il suffit de montrer que si $X' \to S'$ est universellement ferm\'e, alors $X \to S$ est universellement ferm\'e.

\smallskip

Les morphismes fid\`element plats \'etant stables par changement de base, il suffit de prouver que si $X' \to S'$ transforme tout ferm\'e en un ferm\'e, il en est de m\^eme de $X \to S$.

\smallskip

Soit donc un ferm\'e $Z$ de $X$.

\smallskip

Soit $I$ l'id\'eal de $A$ constitu\'e des \'el\'ements dont l'image dans ${\mathcal O}_X (X)$ s'annule en les points de $Z$.

\smallskip

Ainsi
$$
Y = {\rm Spec} (A/I) \xhookrightarrow{ \ \ \ } {\rm Spec} (A) = S
$$
est le plus petit sous-sch\'ema ferm\'e de $S$ qui contient l'image de $F$.

\smallskip

Alors
$$
Y \times_S S' \xhookrightarrow{ \ \ \ } S'
$$
est le plus petit sous-sch\'ema ferm\'e de $S$ qui contient l'image du ferm\'e $Z'$ de $X'$ image r\'eciproque du ferm\'e $Z$ de $X$.

\smallskip

Il r\'esulte de l'hypoth\`ese que l'application
$$
Z' \longrightarrow Y \times_S S'
$$
est surjective.

\smallskip

Donc l'application
$$
Z \longrightarrow Y
$$
est \'egalement surjective, ce qui signifie comme voulu que l'image de $Z$ dans $S$ est un ferm\'e.

\medskip

\item[(10)] Pour traiter le cas des morphismes plats, on peut supposer que $X = {\rm Spec} (C)$ est affine avec donc $X' = {\rm Spec} (B \otimes_A C)$.

\smallskip

Si $B \otimes_A C$ est plat sur $B$ et donc sur $A$, le foncteur
$$
\begin{matrix}
{\rm Mod}_A &\longrightarrow &{\rm Mod}_{B \otimes_A C} \, , \hfill \\
\hfill M &\longmapsto &((B \otimes_A C) \otimes_A M) = (B \otimes_A C) \otimes_C (C \otimes_A M)
\end{matrix}
$$
est exact donc aussi le foncteur
$$
\begin{matrix}
{\rm Mod}_A &\longrightarrow &{\rm Mod}_C \, , \hfill \\
\hfill M &\longmapsto &C \otimes_A M \hfill 
\end{matrix}
$$
puisque le morphisme $C \to B \otimes_A C$ est fid\`element plat.

\smallskip

Ainsi, $B$ est plat sur $A$.

\medskip

\item[(12)] Traitons le cas o\`u $X' \to S'$ est \'etale.

\smallskip

Le cas (1) des morphismes localement de pr\'esentation finie et le cas (10) des morphismes plats \'etant connus, il reste \`a prouver que si le morphisme diagonal
$$
X' \longrightarrow X' \times_{S'} X'
$$
est une immersion ouverte, il en est de m\^eme de
$$
X \longrightarrow X \times_S X \, .
$$
Or ceci r\'esulte du cas (7).

\medskip

\item[(11)] Supposons enfin que $X' \to S'$ est lisse.

\smallskip

On sait d\'ej\`a que $X \to S$ est plat et localement de pr\'esentation finie.

\smallskip

On peut supposer que $X = {\rm Spec} (C)$ est affine avec donc $X' = {\rm Spec} (B \otimes_A C)$. On a la formule
$$
\Omega_{B \otimes_A C / B} = (B \otimes_A C) \otimes_C \Omega_{C/A} = B \otimes_A \Omega_{C/A} \, .
$$
Comme la $A$-alg\`ebre $C$ est de pr\'esentation finie, le $C$-module $\Omega_{C/A}$ est de pr\'esentation finie.

\smallskip

D'autre part, $X' = {\rm Spec} (B \otimes_A C)$ peut \^etre recouvert par des ouverts affines o\`u le module $B \otimes_A \Omega_{C/A}$ devient libre.

\smallskip

Pour tout point $q$ de $X = {\rm Spec} (C)$, notons $\kappa_q$ le corps des fractions de l'anneau int\`egre $C/q$. On peut choisir des \'el\'ements $f_1 , \cdots , f_n \in C$ dont les diff\'erentielles
$$
{\rm d} f_1 , \cdots , {\rm d} f_n
$$
forment une base du $\kappa_q$-espace vectoriel $\kappa_q \otimes_C \Omega_{C/A}$.

\smallskip

Alors le morphisme
$$
\begin{matrix}
\hfill C^n &\longmapsto &\Omega_{C/A} \, , \hfill \\
(c_1 , \cdots , c_n) &\longmapsto &c_1 \cdot {\rm d} f_1 + \cdots + c_n \cdot {\rm d} f_n
\end{matrix}
$$
est un isomorphisme sur un voisinage ouvert $U$ de $q$ dans ${\rm Spec} (C)$.

\smallskip

En effet, l'id\'eal engendr\'e par son d\'eterminant pour n'importe quel choix d'une base de $B \otimes_A \Omega_{C/A} = \Omega_{B \otimes_A C/B}$ localement sur $X' = {\rm Spec} (B \otimes_A C)$ d\'efinit un sous-sch\'ema ferm\'e
$$
Y' \xhookrightarrow{ \ \ \ } X'
$$
qui ne d\'epend d'aucun choix.

\smallskip

D'apr\`es le corollaire \ref{corII321}~(i), le sous-sch\'ema ferm\'e $Y' \hookrightarrow X'$ provient d'un sous-sch\'ema ferm\'e
$$
Y \xhookrightarrow{ \ \ \ } X
$$
dont le compl\'ementaire $U$ est un voisinage ouvert de $q$ qui r\'epond \`a la question pos\'ee. 
\end{listeisansmarge}
\end{demo}

\medskip

On d\'eduit de ce th\'eor\`eme:

\begin{cor}\label{corII324}

Dans la cat\'egorie {\rm Sch} des sch\'emas, on a:

\begin{listeimarge}

\item Si ${\mathcal M}_{\rm p}$ est la classe g\'eom\'etrique des morphismes plats, la propri\'et\'e (fpqc) [resp. (fppf)] des familles de morphismes de ${\mathcal M}_{\rm p}$ de m\^eme but
$$
(S_i \longrightarrow S)_{i\in I}
$$
de contenir une sous-famille $(S_i \to S)_{i \in I'}$, $I' \subseteq I$, qui soit ``fid\`element plate quasi-compacte'' [resp. et telle que chaque $S_i \to S$, $i \in I'$, soit (localement) de pr\'esentation finie] est une notion g\'eom\'etrique de recouvrement.

\medskip

\item Si ${\mathcal M}_{\rm l}$ [resp. ${\mathcal M}_{\mbox{\scriptsize\rm \'et}}$] est la classe g\'eom\'etrique des morphismes lisses [resp. \'etales], les propri\'et\'es (fpqc) et (fppf) des familles de morphismes de ${\mathcal M}_{\rm l}$ [resp. ${\mathcal M}_{\mbox{\scriptsize\rm \'et}}$] de m\^eme but
$$
(S_i \longrightarrow S)_{i \in I}
$$
co{\"\i}ncident et d\'efinissent une notion g\'eom\'etrique de recouvrement.
\end{listeimarge}
\end{cor}

\begin{demo}
\begin{listeisansmarge}
\item La propri\'et\'e (fpqc) [resp. (fppf)] est v\'erifi\'ee par les familles qui comprennent au moins un isomorphisme.

\smallskip

D'autre part, elle est stable par changement de base puisqu'il en est ainsi des propri\'et\'es pour un morphisme $S' \to S$ d'\^etre plat, surjectif, quasi-compact ou de pr\'esentation finie et que la formation des sommes disjointes $\underset{i \in I}{\coprod}$ est respect\'ee par les foncteurs de produits fibr\'es $X \times_S \bullet$.

\smallskip

Enfin, cette propri\'et\'e est (fpqc)-locale d'apr\`es le th\'eor\`eme \ref{thmII323}. A fortiori, elle est (fppf)-locale.

\medskip

\item Les morphismes lisses (en particulier \'etales) sont localement de pr\'esentation finie donc les propri\'et\'es (fppf) et (fpqc) des familles de morphismes de ${\mathcal M}_{\rm l}$ [resp. ${\mathcal M}_{\mbox{\scriptsize\rm \'et}}$] de m\^eme but co{\"\i}ncident.

\smallskip

D'apr\`es (i), ces propri\'et\'es d\'efinissent une notion g\'eom\'etrique de recouvrement.\end{listeisansmarge}\end{demo}


Les paires constitu\'ees de ces classes g\'eom\'etriques de morphismes et des notions g\'eom\'etriques de recouvrement (fpqc) ou (fppf) d\'efinissent des topologies de Grothendieck et donc des gros et petits sites sur les sch\'emas:

\begin{defn}\label{defII325}

Soit ${\mathcal G}$ une sous-cat\'egorie pleine de {\rm Sch} qui est essentiellement petite, stable par formation des produits fibr\'es, et contient tous les morphismes $X \to S$ vers un objet $S$ de ${\mathcal G}$ qui sont des immersions ouvertes ou sont de pr\'esentation finie.

\smallskip

Alors:

\begin{listeimarge}

\item La paire $({\mathcal M}_{\rm p}, \mbox{fpqc})$ [resp. $({\mathcal M}_{\rm p} , \mbox{fppf} \, )$, resp. $({\mathcal M}_{\rm l} , \mbox{fppf} \, )$, resp. $({\mathcal M}_{\mbox{\scriptsize\rm \'et}} , \mbox{fppf} \, )$] d\'efinit sur ${\mathcal G}$ une topologie de Grothendieck $J_{\rm fpqc}$ [resp. $J_{\rm fppf}$, resp. $J_{\rm l}$, resp. $J_{\mbox{\scriptsize\rm \'et}}$] appel\'ee la topologie fid\`element plate quasi-compacte (fpqc) [resp. fid\`element plate de pr\'esentation finie (fppf), resp. lisse, resp. \'etale].

\medskip

\item A chaque objet $X$ de ${\mathcal G}$ est associ\'e un gros site fpqc $({\mathcal G} / X , J_{\rm fpqc})$ [resp. fppf $({\mathcal G} / X , J_{\rm fppf})$, resp. lisse $({\mathcal G} / X , J_{\rm l})$, resp. \'etale $({\mathcal G}/X ,J_{\mbox{\scriptsize\rm \'et}})$] et un petit site fpqc $(({\mathcal G}/X)_{{\mathcal M}_{\rm p}} , J_{\rm fpqc})$ [resp. fppf $(({\mathcal G}/X)_{{\mathcal M}_{\rm p}} , J_{\rm fppf})$, resp. lisse $(({\mathcal G}/X)_{{\mathcal M}_{\rm l}} , J_{\rm l})$, resp. \'etale $(({\mathcal G}/X)_{{\mathcal M}_{\mbox{\tiny\rm \'et}}} , J_{\mbox{\scriptsize\rm \'et}})$].
\end{listeimarge}
\end{defn}

\begin{remarkqed}

On peut montrer que si
$$
X' \longrightarrow X
$$
est un morphisme lisse et $x$ un point de $X$ contenu dans l'image de $X'$, il existe un voisinage ouvert $U$ de $x$ dans $X$ et un morphisme \'etale surjectif de pr\'esentation finie
$$
U' \longrightarrow U
$$
tel que le compos\'e
$$
U' \longrightarrow U \xhookrightarrow{ \ \ \ } X
$$
se factorise \`a travers $X' \to X$.

\smallskip

Par cons\'equent, la topologie lisse $J_{\rm l}$ et la topologie \'etale $J_{\mbox{\scriptsize\rm \'et}}$ co{\"\i}ncident sur ${\mathcal G}$.

\smallskip

En revanche, si $X$ est un objet de ${\mathcal G}$, le petit site lisse de $X$ et son petit site \'etale ne co{\"\i}ncident pas. 
\end{remarkqed}

\section{Faisceaux sur un site de Grothendieck}\label{sec24}

\subsection{La notion de faisceau sur un site}\label{subsec241}

On rappelle que, par d\'efinition, la cat\'egorie $\widehat{\mathcal C}$ des pr\'efaisceaux sur une cat\'egorie essentiellement petite (ou plus g\'en\'eralement localement petite) ${\mathcal C}$ est la cat\'egorie $[{\mathcal C}^{\rm op} , {\rm Ens}]$ des foncteurs contravariants
$$
F : {\mathcal C}^{\rm op} \longrightarrow {\rm Ens} \, .
$$

La notion de topologie de Grothendieck $J$ sur une cat\'egorie essentiellement petite ${\mathcal C}$ est la plus vaste qui permet de g\'en\'eraliser au cadre des cat\'egories la notion ant\'erieure de faisceau sur un espace topologique:

\begin{defn}\label{defII41}

Soit ${\mathcal C}$ une cat\'egorie essentiellement petite munie d'une topologie de Grothendieck $J$, constituant un site $({\mathcal C} , J)$.

\smallskip

Alors:

\begin{listeimarge}

\item Un pr\'efaisceau sur ${\mathcal C}$
$$
F : {\mathcal C}^{\rm op} \longrightarrow {\rm Ens} 
$$
est appel\'e un faisceau pour la topologie $J$, ou un $J$-faisceau, si, pour tout objet $X$ de ${\mathcal C}$ et tout crible $J$-couvrant de $X$
$$
S \in J(X) \, ,
$$
consid\'er\'e comme la sous-cat\'egorie pleine de ${\mathcal C}/X$ dont les objets sont les \'el\'ements $U \to X$ de $S$, l'application canonique
$$
F(X) \longrightarrow \varprojlim_{(U \to X) \in S} F(U)
$$
est une bijection.

\medskip

\item Un morphisme $F_1 \to F_2$ entre deux $J$-faisceaux $F_1 , F_2$ sur ${\mathcal C}$ est un morphisme de pr\'efaisceaux $F_1 \to F_2$.

\smallskip

Autrement dit, les $J$-faisceaux sur ${\mathcal C}$ d\'efinissent une sous-cat\'egorie pleine de la cat\'egorie $\widehat{\mathcal C}$ des pr\'efaisceaux sur ${\mathcal C}$ que l'on notera
$$
\widehat{\mathcal C}_J \quad \mbox{ou} \quad {\rm Sh} ({\mathcal C},J) \, .
$$
\end{listeimarge}
\end{defn}

\medskip

\begin{remarksqed}
\begin{listeisansmarge}
\item Si $S \in J(X)$, la sous-cat\'egorie pleine $S$ de ${\mathcal C}/X$ n'est pas n\'ecessairement petite (\`a moins que ${\mathcal C}$ ne soit une petite cat\'egorie) mais, comme ${\mathcal C}$, elle est essentiellement petite, et cela suffit pour que le foncteur diagonal 
$$
\Delta : {\rm Ens} \longrightarrow [S^{\rm op} , {\rm Ens}]
$$
admette un adjoint \`a  droite
$$
\varprojlim_{(U \to X) \in S} : [S^{\rm op} , {\rm Ens}] \longrightarrow {\rm Ens} \, .
$$

\item Dans le cas o\`u ${\mathcal C}$ est la cat\'egorie $O(X)$ des ouverts d'un espace topologique $X$ et $J$ est la topologie ordinaire de $O(X)$, la notion de faisceau sur le site $({\mathcal C},J)$ co{\"\i}ncide avec celle de faisceau sur l'espace topologique $X$ introduite au paragraphe \ref{sec14}.

\smallskip

En effet, les cribles couvrants $S$ d'un objet $U$ de ${\mathcal C} = O(X)$ sont les cribles engendr\'es par une famille d'ouverts
$$
U_i \subset U \, , \ i \in I \, ,
$$
telle que
$$
\bigcup_{i \in I} U_i = U \, ,
$$
et on a alors pour tout pr\'efaisceau $F : O(X)^{\rm op} \to {\rm Ens}$, une identification
$$
\varprojlim_{(U' \hookrightarrow U)\in S} F(U') = {\rm eg} \left[ \prod_{i \in I} F(U_i) \rightrightarrows \prod_{i,j \in I} F(U_i \cap U_j) \right].
$$

\item On a 
$$
\widehat{\mathcal C}_J = \widehat{\mathcal C} \, ,
$$
autrement dit tout pr\'efaisceau sur ${\mathcal C}$ est un $J$-faisceau, si et seulement si $J$ est la topologie discr\`ete de ${\mathcal C}$ pour laquelle les seuls cribles couvrants sont les cribles maximaux engendr\'es par les fl\`eches ${\rm id}_X : X \to X$.

\medskip

\item Si $J$ et $J'$ sont deux topologies sur une m\^eme cat\'egorie essentiellement petite ${\mathcal C}$ telles que
$$
J \subseteq J' \, ,
$$
alors tout $J'$-faisceau est a fortiori un $J$-faisceau.

\smallskip

Autrement dit, $\widehat{\mathcal C}_{J'}$ est une sous-cat\'egorie pleine de $\widehat{\mathcal C}_J$.

\smallskip

Ainsi, plus une topologie $J$ sur une cat\'egorie ${\mathcal C}$ fix\'ee est fine, plus la notion de $J$-faisceau est restrictive. 
\end{listeisansmarge}
\end{remarksqed}

\subsection{Un exemple universel: le faisceau des cribles $J$-ferm\'es}\label{subsec242}

\medskip

Nous allons donner une d\'efinition g\'en\'erale qui s'applique \`a  tout site $({\mathcal C},J)$ et fournit un exemple int\'eressant de faisceau. On verra plus loin que ce faisceau d\'efini sur tout site est caract\'eris\'e par une propri\'et\'e cat\'egorique remarquable qui le rend important.

\smallskip

Sa d\'efinition est fond\'ee sur celle de crible ferm\'e relativement \`a  une topologie de Grothendieck:

\begin{defn}\label{defII42}

Soit $({\mathcal C},J)$ un site.

\smallskip

Un crible $S$ d'un objet $X$ de la cat\'egorie ${\mathcal C}$ est dit ferm\'e pour la topologie $J$, ou $J$-ferm\'e, si pour toute fl\`eche de ${\mathcal C}$ de but $X$
$$
p : U \longrightarrow X
$$
telle qu'existe une famille $J$-couvrante de $U$
$$
q_i : U_i \longrightarrow U
$$
dont les compos\'es $p \circ q_i : U_i \to X$ sont tous \'el\'ements de $S$, alors $p$ est n\'ecessairement \'el\'ement de $S$.
\end{defn}

\begin{remarksqed}
\begin{listeisansmarge}

\item Pour toute fl\`eche de ${\mathcal C}$
$$
f : X \to Y \, ,
$$
l'application associ\'ee d'image r\'eciproque des cribles
$$
\begin{matrix}
f^* : \Omega (Y) &\longmapsto &\Omega (X) \, , \\
\hfill S &\longmapsto &f^*S \hfill
\end{matrix}
$$
transforme les cribles $J$-ferm\'es en cribles $J$-ferm\'es.

\medskip

\item Si $J$ est la topologie discr\`ete de ${\mathcal C}$, dont les seuls cribles $J$-couvrants sont les cribles maximaux, alors tout crible sur un objet de ${\mathcal C}$ est $J$-ferm\'e.

\medskip

\item Si $X$ est un espace topologique et ${\mathcal C} = O(X)$ est la cat\'egorie des ouverts de $X$ munie de la topologie usuelle $J$, les cribles $J$-ferm\'es d'un objet $U \subset X$ de ${\mathcal C}$ sont les cribles engendr\'es par un unique ouvert $V \subset U$ de $U$.

\medskip

\item Plus g\'en\'eralement, si $(O,\leq)$ est un treillis muni de sa topologie canonique $J$, les cribles $J$-ferm\'es d'un objet $u$ sont les cribles engendr\'es par un unique objet $v \leq u$. 

\end{listeisansmarge}
\end{remarksqed}


Cette d\'efinition \'etant pos\'ee, on peut \'enoncer:

\begin{lem}\label{lemII43}

Soit $({\mathcal C},J)$ un site.

\smallskip

Associons \`a  tout objet $X$ de ${\mathcal C}$ l'ensemble
$$
\Omega_J (X) \subset \Omega (X)
$$
des cribles $J$-ferm\'es de $X$ et \`a  tout morphisme $f : X \to Y$ de ${\mathcal C}$ l'application
$$
f^* : \Omega_J (Y) \longrightarrow \Omega_J (X)
$$
induite par l'application $\Omega (Y) \to \Omega (X)$ d'image r\'eciproque des cribles.

\smallskip

Alors le pr\'efaisceau
$$
\Omega_J : {\mathcal C}^{\rm op} \longrightarrow {\rm Ens}
$$
est un faisceau pour la topologie $J$.
\end{lem}

\bigskip

\begin{demo}

Etant donn\'e un crible $J$-couvrant $S$ d'un objet $X$ de ${\mathcal C}$, il s'agit de montrer que l'application
$$
\Omega_J (X) \longrightarrow \varprojlim_{(U \to X) \in S} \Omega_J (U)
$$
est bijective.

\smallskip

Consid\'erons donc un \'el\'ement de la limite $\underset{(U \to X) \in S}{\varprojlim} \, \Omega_J (U)$ constitu\'e d'une famille compatible de cribles $J$-ferm\'es
$$
S_{U,u} \in \Omega_J (U)
$$
index\'es par les objets $U \xrightarrow{ \ u \ } X$ du crible $S$.

\smallskip

Pour toute fl\`eche $v : V \to X$ de ${\mathcal C}$ de but $X$, le crible $v^* S$ de $V$ est $J$-couvrant.

\smallskip

Si $T$ est un crible $J$-ferm\'e de $X$, la fl\`eche $v : V \to X$ est \'el\'ement de $T$ si et seulement si les fl\`eches
$$
u = v \circ p \, , \ (U \xrightarrow{ \ p \ } V) \in v^* S
$$
sont \'el\'ements de $T$.

\smallskip

Si $T \in \Omega_J(X)$ s'envoie sur la famille $(S_{U,u})_{(U \xrightarrow{  u  } X) \in S}$, cela \'equivaut \`a  demander que, pour tout $(U \xrightarrow{ \ p \ } V) \in v^* S$, on ait
$$
{\rm id}_U \in S_{U, \, v \, \circ \, p} \, .
$$

Ainsi, la famille $(S_{U,u})_{(U \xrightarrow{  u  } X) \in S}$ a au plus un ant\'ec\'edent dans l'ensemble $\Omega_J (X)$.

\smallskip

R\'eciproquement, posons
$$
T = \left\{ V \xrightarrow{ \ v \ } X \ \bigl\vert \ {\rm id}_U \in S_{U, \, v \, \circ \, p} \, , \ \forall \, (U \xrightarrow{ \ p \ } V) \in v^* S \right\} \, .
$$

Montrons d'abord que $T$ est un crible, autrement dit que pour tout \'el\'ement $V \xrightarrow{ \ v \ } X$ de $T$ et tout morphisme $W \xrightarrow{ \ w \ } V$ de ${\mathcal C}$, le compos\'e $W \xrightarrow{ \ v \, \circ \, w \ } X$ est encore \'el\'ement de $T$. En effet, pour tout \'el\'ement $U \xrightarrow{ \ q \ } W$ de $(v \circ w)^* S$, le compos\'e $U \xrightarrow{ \ w \, \circ \, q \ } V$ est \'el\'ement de $v^* S$ et donc ${\rm id}_U$ est \'el\'ement de $S_{U , \, v \, \circ \, (w \, \circ \, q)} = S_{U,(v \, \circ \, w) \, \circ \, q}$.

\smallskip

Puis montrons que le crible $T$ est $J$-ferm\'e. Consid\'erons donc un morphisme $v : V \to X$ de but $X$ tel qu'existe un crible $J$-couvrant $Q$ de $V$ constitu\'e de morphismes
$$
q : V' \longrightarrow V
$$
dont les compos\'es $v \circ q : V' \to X$ sont \'el\'ements de $T$. Il s'agit de prouver que $v$ est alors \'el\'ement de $T$, autrement dit que pour tout \'el\'ement $U \xrightarrow{ \ p \ } V$ de $v^* S$, ${\rm id}_U$ est \'el\'ement de $S_{U , \, v \, \circ \, p}$. Comme le crible $p^* Q$ de $U$ est $J$-couvrant et que le crible $S_{U, \, v \, \circ \, p}$ de $U$ est $J$-ferm\'e, c'est \'equivalent \`a  demander que tout \'el\'ement $q' : U' \to U$ de $p^* Q$ soit aussi \'el\'ement de $S_{U, \, v \, \circ \, p}$. Or, si $p \circ q' : U' \xrightarrow{ \ q' \ } U \xrightarrow{ \ p \ } V$ est \'el\'ement de $Q$, $v \circ p \circ q' : U' \xrightarrow{ \ q' \ } U \xrightarrow{ \ p \ } V \xrightarrow{ \ v \ } X$ est \'el\'ement de $T$ et ${\rm id}_{U'} \in S_{U' , \, v \, \circ \, p \, \circ \, q'}$ puisque $p \circ q' \in v^* S$. La conclusion r\'esulte de ce que $S_{U' ,  \, v \, \circ \, p \, \circ \,q'} = q'^* S_{U, \, v \, \circ \, p}$.

\smallskip

Enfin, montrons que pour tout \'el\'ement $u : U \to X$ de $S$, on a
$$
u^* T = S_{U,u} \, .
$$

Or, un morphisme $p : V \to U$ est \'el\'ement de $u^* T$ si et seulement si, pour tout \'el\'ement $q : V' \to V$ de $(u \circ p)^* S$, on a
$$
{\rm id}_{V'} \in S_{V' , \, u \, \circ \, p \, \circ \, q} = q^* S_{V, \, u \, \circ \, p}
$$
c'est-\`a -dire
$$
q \in S_{V, \, u \, \circ \, p} \, .
$$

Comme le crible $(u \circ p)^* S$ de $V$ est $J$-couvrant et que le crible $S_{V, \, u \, \circ \, p}$ de $V$ est $J$-ferm\'e, cela \'equivaut \`a 
$$
{\rm id}_V \in S_{V, \, u \, \circ \, p} = p^* S_{U,u}
$$
soit 
$$
p \in S_{U,u}
$$
comme annonc\'e.\end{demo}

\subsection{Faisceaux de fonctions}\label{subsec243}

\medskip

Les fonctions continues vers un espace topologique fix\'e d\'efinissent des faisceaux pour la topologie usuelle des espaces topologiques:

\begin{lem}\label{lemII44}

Soit ${\mathcal C}$ une sous-cat\'egorie pleine et essentiellement petite de la cat\'egorie ${\rm Top}$ des espaces topologiques telle que, pour tout objet $X$ de ${\mathcal C}$, les ouverts de $X$ sont encore des objets de ${\mathcal C}$.

\smallskip

Soit $J$ la topologie usuelle de ${\mathcal C}$ pour laquelle un crible sur un objet $X$ est couvrant s'il contient une famille d'immersions ouvertes $U_i \hookrightarrow X$ dont les images recouvrent $X$ au sens que $X = \underset{i}{\bigcup} \, U_i$.

\smallskip

Soit $C$ un espace topologique.

\smallskip

Alors le pr\'efaisceau
$$
\begin{matrix}
{\mathcal C}^{\rm op} &\longrightarrow &{\rm Ens} \, , \hfill \\
\hfill X &\longmapsto &{\rm Hom}_{\rm Top} (X,C) = \{\mbox{applications continues} \ X \to C\}
\end{matrix}
$$
est un faisceau pour la topologie $J$.
\end{lem}

\begin{remarks}
\begin{listeisansmarge}
\item Si $C$ est un ensemble muni de la topologie discr\`ete, alors pour tout espace topologique $X$
$$
{\rm Hom}_{\rm Top} (X,C)
$$
est l'ensemble des applications localement constantes $X \to C$.

\medskip

\item Si $C = {\mathbb S}$ est l'espace de Sierpinski constitu\'e de l'ensemble $\{0,1\}$ muni de la topologie dont l'unique ouvert non trivial est $\{1\}$, alors pour tout espace topologique $X$
$$
{\rm Hom}_{\rm Top} (X,{\mathbb S})
$$
s'identifie \`a  l'ensemble $O(X)$ des ouverts de $X$.
\end{listeisansmarge}
\end{remarks}

\begin{demo}

En effet, pour tout recouvrement d'un espace topologique $X$ par une famille d'ouverts $U_i$, se donner une application continue
$$
u : X \longrightarrow C
$$
\'equivaut \`a  se donner une famille d'applications continues
$$
u_i : U_i \longrightarrow C
$$
telle que, pour tous indices $i$ et $j$, les applications $u_i$ et $u_j$ co{\"\i}ncident sur l'intersection $U_i \cap U_j$. 
\end{demo}

\medskip

Ce r\'esultat s'\'etend aux sous-cat\'egories de la cat\'egorie ${\rm Top}_{\rm an}$ des espaces annel\'es qui sont ``g\'eom\'etriques'' au sens de la d\'efinition \ref{defI58}:

\begin{prop}\label{propII45}

Soit ${\mathcal G}$ une sous-cat\'egorie ``g\'eom\'etrique'' de la cat\'egorie ${\rm Top}_{\rm an}$ des espaces annel\'es.

\smallskip

Soit ${\mathcal C}$ une sous-cat\'egorie pleine et essentiellement petite de ${\mathcal G}$ telle que, pour tout objet $(X,{\mathcal O}_X)$ de ${\mathcal C}$, ses ouverts $(U \xhookrightarrow{ \ i \ } X , {\mathcal O}_U = i^* {\mathcal O}_X)$ sont encore des objets de ${\mathcal C}$.

\smallskip

Soit $J$ la topologie usuelle de ${\mathcal C}$ pour laquelle un crible sur un objet $(X,{\mathcal O}_X)$ est couvrant s'il contient une famille d'immersions ouvertes $(U_i , {\mathcal O}_{U_i}) \hookrightarrow (X,{\mathcal O}_X)$ dont les images recouvrent $X$ au sens que $X = \underset{i}{\bigcup} \, U_i$.

\smallskip

Alors:

\begin{listeimarge}

\item Pour tout objet $C$ de ${\mathcal G}$, le pr\'efaisceau
$$
\begin{matrix}
{\mathcal C}^{\rm op} &\longrightarrow &{\rm Ens} \, , \hfill \\
\hfill X &\longmapsto &{\rm Hom}_{\mathcal G} (X,C)
\end{matrix}
$$
est un faisceau pour la topologie ${\mathcal G}$.

\medskip

\item Le pr\'efaisceau des ``fonctions de structure''
$$
\begin{matrix}
\hfill {\mathcal C}^{\rm op} &\longrightarrow &{\rm Ens} \, , \hfill \\
(X,{\mathcal O}_X) &\longmapsto &{\mathcal O}_X (X)
\end{matrix}
$$
est un faisceau pour la topologie $J$.
\end{listeimarge}
\end{prop}

\begin{remarks}
\begin{listeisansmarge}
\item La partie (ii) est un cas particulier de (i) si le pr\'efaisceau
$$
\begin{matrix}
\hfill {\mathcal G} &\longrightarrow &{\rm Ens} \hfill \\
(X,{\mathcal O}_X) &\longmapsto &{\mathcal O}_X (X)
\end{matrix}
$$
est repr\'esentable par un objet ${\mathbb A}^1$ de ${\mathcal G}$.

\medskip

\item Cette proposition s'applique en particulier lorsque ${\mathcal G}$ est la cat\'egorie g\'eom\'etrique des vari\'et\'es diff\'erentielles de classe $C^k$, $k \geq 1$, ou bien celle des vari\'et\'es analytiques, ou encore la cat\'egorie ${\rm Sch}$ [resp. ${\rm Sch}/S$] des sch\'emas [resp. des sch\'emas sur un sch\'ema de base $S$].
\end{listeisansmarge}
\end{remarks}

\medskip

\begin{demo}
\begin{listeisansmarge}
\item[(ii)] r\'esulte de ce que, par d\'efinition des espaces annel\'es, le pr\'efaisceau de structure ${\mathcal O}_X$ de tout espace annel\'e $(X,{\mathcal O}_X)$ est un faisceau sur l'espace topologique $X$.

\medskip

\item[(i)] Pour tout recouvrement d'un espace annel\'e $(X,{\mathcal O}_X)$ par des ouverts $(U_i , {\mathcal O}_{U_i})$, se donner un morphisme d'espaces annel\'es
$$
\left( X \xrightarrow{ \ u \ } C \, , \ {\mathcal O}_C \longrightarrow u_* {\mathcal O}_X \right)
$$
\'equivaut \`a  se donner une famille de morphismes d'espaces annel\'es
$$
\left(U_i \xrightarrow{ \ u_i \ } C \, , \ {\mathcal O}_C \longrightarrow (u_i)_* {\mathcal O}_{U_i} \right)
$$
telle que, pour tous indices $i$ et $j$, les morphismes d'espaces annel\'es induits
$$
(U_i \cap U_j , {\mathcal O}_{U_i \, \cap \, U_j}) \rightrightarrows (C,{\mathcal O}_C)
$$
co{\"\i}ncident.

\smallskip

De plus, il r\'esulte des propri\'et\'es (i) et (iii) de la d\'efinition \ref{defI58} de la notion de ``sous-cat\'egorie g\'eom\'etrique'' de ${\rm Top}_{\rm an}$ que, si $(X,{\mathcal O}_X)$ et $(C,{\mathcal O}_C)$ sont des objets d'une sous-cat\'egorie g\'eom\'etrique ${\mathcal G}$, le morphisme
$$
\left( X \xrightarrow{ \ u \ } C \, , \ {\mathcal O}_C \longrightarrow u_* {\mathcal O}_X \right)
$$
est un morphisme de ${\mathcal G}$ si et seulement si il en est ainsi des morphismes induits
$$
\left(U_i \xrightarrow{ \ u_i \ } C \, , \ {\mathcal O}_C \longrightarrow (u_i)_* {\mathcal O}_{U_i} \right) .
$$
\end{listeisansmarge}
\end{demo}

\subsection{Sites annel\'es, faisceaux de modules et faisceaux localement libres}\label{subsec244}

\medskip

La notion d'espace annel\'e se g\'en\'eralise en celle de site annel\'e:

\begin{defn}\label{defII46}

Un site annel\'e est un site $({\mathcal C} , J)$ muni d'un faisceau d'anneaux ${\mathcal O}$, c'est-\`a -dire d'un foncteur contravariant vers la cat\'egorie ${\rm An}$ des anneaux
$$
{\mathcal O} : {\mathcal C}^{\rm op} \longrightarrow {\rm An}
$$
dont le compos\'e avec le foncteur ${\rm An} \to {\rm Ens}$ d'oubli des structures d'anneaux
$$
{\mathcal C}^{\rm op} \xrightarrow{ \ {\mathcal O} \ } {\rm An} \longrightarrow {\rm Ens}
$$
est un faisceau pour la topologie $J$.
\end{defn}

\bigskip

\begin{remarksqed}
\begin{listeisansmarge}
\item Autrement dit, un faisceau d'anneaux sur un site $({\mathcal C} , J)$ est un faisceau d'ensembles ${\mathcal O}$ muni d'une structure d'anneau de ${\mathcal O} (X)$ pour tout objet $X$ de ${\mathcal C}$ et tel que, pour tout morphisme $X \to Y$ de ${\mathcal C}$, l'application
$$
{\mathcal O} (Y) \longrightarrow {\mathcal O} (X)
$$
soit un morphisme d'anneaux. 

\medskip

\item Un site annel\'e $({\mathcal C} , J, {\mathcal O})$ est dit commutatif si, pour tout objet $X$ de ${\mathcal C}$, l'anneau ${\mathcal O}(X)$ est commutatif.

\medskip 

\item La proposition \ref{propII45} (ii) signifie que pour toute sous-cat\'egorie g\'eom\'etrique ${\mathcal G}$ de la cat\'egorie ${\rm Top}_{\rm an}$ des espaces annel\'es, et pour toute sous-cat\'egorie pleine et essentiellement petite ${\mathcal C}$ de ${\mathcal G}$ qui contient les immersions ouvertes vers les objets de ${\mathcal C}$, alors  la cat\'egorie ${\mathcal C}$ munie de la topologie usuelle $J$ et du foncteur contravariant
$$
(X,{\mathcal O}_X) \longmapsto {\mathcal O}_X (X)
$$
est un site annel\'e commutatif.

\smallskip

Il en est en particulier ainsi lorsque ${\mathcal G}$ est la cat\'egorie des vari\'et\'es diff\'erentielles de classe $C^k$, $k \geq 1$, ou bien celle des vari\'et\'es analytiques, ou encore celle des sch\'emas (ou des sch\'emas sur un sch\'ema de base fix\'e). 
\end{listeisansmarge}
\end{remarksqed}

\medskip

On dispose de la notion de faisceau de modules sur un site annel\'e:

\begin{defn}\label{defII47}
\begin{listeimarge}
\item Un faisceau de groupes ab\'eliens sur un site $({\mathcal C} , J)$ est un foncteur contravariant vers la cat\'egorie ${\rm Ab}$ des groupes ab\'eliens
$$
{\mathcal M} : {\mathcal C}^{\rm op} \longrightarrow {\rm Ab}
$$
dont le compos\'e avec le foncteur ${\rm Ab} \to {\rm Ens}$ d'oubli des structures de groupes ab\'eliens
$$
{\mathcal C}^{\rm op} \xrightarrow{ \ {\mathcal M} \ } {\rm Ab} \longrightarrow {\rm Ens}
$$
est un faisceau pour la topologie $J$.

\smallskip

Un morphisme entre deux faisceaux de groupes ab\'eliens sur $({\mathcal C},J)$ est un morphisme dans la cat\'egorie $[{\mathcal C}^{\rm op} , {\rm Ab}]$ des foncteurs contravariants ${\mathcal C}^{\rm op} \to {\rm Ab}$.

\medskip

\item Un faisceau de modules ${\mathcal M}$ sur un site annel\'e $({\mathcal C} , J , {\mathcal O})$ est un faisceau de groupes ab\'eliens ${\mathcal M}$ muni d'un morphisme de pr\'efaisceaux
$$
{\mathcal O} \times {\mathcal M} \longrightarrow {\mathcal M}
$$
tel que, pour tout objet $X$ de ${\mathcal C}$, l'application induite
$$
{\mathcal O} (X) \times {\mathcal M} (X) \longrightarrow {\mathcal M} (X)
$$
fasse du groupe ab\'elien ${\mathcal M} (X)$ un module sur l'anneau ${\mathcal O} (X)$.

\smallskip

Un morphisme de faisceaux de modules
$$
{\mathcal M}_1 \longrightarrow {\mathcal M}_2
$$
sur un site annel\'e $({\mathcal C} , J , {\mathcal O})$ est un morphisme de faisceaux de groupes ab\'eliens ${\mathcal M}_1 \to {\mathcal M}_2$ tel que le carr\'e induit dans la cat\'egorie des pr\'efaisceaux
$$
\xymatrix{
{\mathcal O} \times {\mathcal M}_1 \ar[d] \ar[r] &{\mathcal M}_1 \ar[d] \\
{\mathcal O} \times {\mathcal M}_2 \ar[r] &{\mathcal M}_2
}
$$
soit commutatif.

\end{listeimarge}
\end{defn}

\begin{remarksqed}
\begin{listeisansmarge}
\item Autrement dit, un faisceau de groupes ab\'eliens sur un site $({\mathcal C} , J)$ est un faisceau d'ensembles ${\mathcal M}$ muni d'une structure de groupe ab\'elien de ${\mathcal M} (X)$ pour tout objet $X$ de ${\mathcal C}$ et tel que, pour tout morphisme $X \to Y$ de ${\mathcal C}$, l'application
$$
{\mathcal O} (Y) \longrightarrow {\mathcal O} (X)
$$
soit un morphisme de groupes ab\'eliens.

\medskip

\item De m\^eme, un faisceau de modules sur un site annel\'e $({\mathcal C} , J , {\mathcal O})$ est un faisceau d'ensembles ${\mathcal M}$ muni pour tout objet $X$ de ${\mathcal C}$ d'une structure de module de ${\mathcal M}(X)$ sur l'anneau ${\mathcal O}(X)$ et tel que, pour tout morphisme $X \to Y$ de ${\mathcal C}$, l'application
$$
{\mathcal M}(Y) \longrightarrow {\mathcal M}(X)
$$
soit un morphisme de modules sur ${\mathcal O}(Y)$ via le morphisme d'anneaux ${\mathcal O} (Y) \to {\mathcal O} (X)$.

\medskip

\item Les faisceaux de groupes ab\'eliens sur un site $({\mathcal C},J)$ et leurs morphismes forment une cat\'egorie que l'on peut noter ${\mathcal A}b_{{\mathcal C},J}$. C'est une sous-cat\'egorie pleine de la cat\'egorie $[{\mathcal C}^{\rm op} , {\rm Ab}]$ des foncteurs contravariants ${\mathcal C}^{\rm op} \to {\rm Ab}$.

\medskip

\item Un morphisme de faisceaux de modules sur un site annel\'e $({\mathcal C},J,{\mathcal O})$
$$
{\mathcal M}_1 \longrightarrow {\mathcal M}_2
$$
est un morphisme de faisceaux de groupes ab\'eliens sur $({\mathcal C},J)$ tel que, pour tout objet $X$ de ${\mathcal C}$, le carr\'e induit dans la cat\'egorie des ensembles
$$
\xymatrix{
{\mathcal O}(X) \times {\mathcal M}_1(X) \ar[d] \ar[r] &{\mathcal M}_1(X) \ar[d] \\
{\mathcal O}(X) \times {\mathcal M}_2(X) \ar[r] &{\mathcal M}_2(X)
}
$$
soit commutatif.

\medskip

\item Les faisceaux de modules sur un site annel\'e $({\mathcal C},J,{\mathcal O})$ et leurs morphismes forment une cat\'egorie que l'on peut noter ${\mathcal M}od_{{\mathcal C},J,{\mathcal O}}$ ou, s'il n'y a pas d'ambigu{\"\i}t\'e sur le site sous-jacent, ${\mathcal M}od_{\mathcal O}$.

\smallskip

Elle est munie du foncteur d'oubli de la structure de multiplication par les scalaires
$$
{\mathcal M}od_{{\mathcal C},J,{\mathcal O}} \longrightarrow {\mathcal A}b_{{\mathcal C},J}
$$
lequel est un foncteur fid\`ele.

\smallskip

Les objets de ${\mathcal M}od_{{\mathcal C} , J , {\mathcal O}}$ sont parfois appel\'es les Modules sur ${\mathcal O}$ ou les ${\mathcal O}$-Modules.

\medskip

\item Toutes les notions cat\'egoriques usuelles s'appliquent dans les cat\'egories ${\mathcal A}b_{{\mathcal C},J}$ ou ${\mathcal M}od_{{\mathcal C},J,{\mathcal O}}$.

\smallskip

En particulier, un sous-Module d'un ${\mathcal C}$-Module ${\mathcal M}$ est un monomorphisme
$$
{\mathcal N} \xhookrightarrow{ \ \ \ } {\mathcal M}
$$
consid\'er\'e \`a  isomorphisme pr\`es, et un Module quotient d'un ${\mathcal O}$-Module ${\mathcal M}$ est un \'epimorphisme
$$
{\mathcal M} -\!\!\!-\!\!\!\twoheadrightarrow {\mathcal Q}
$$
consid\'er\'e \`a  isomorphisme pr\`es.

\smallskip

Les sous-Modules d'un ${\mathcal O}$-Module ${\mathcal M}$ sont les collections de sous-modules sur les anneaux ${\mathcal O}(X)$
$$
{\mathcal N}(X) \subset {\mathcal M} (X)
$$
index\'ees par les objets $X$ de ${\mathcal C}$ telles que, pour tout morphisme $X \to Y$ de ${\mathcal C}$, l'application
$$
{\mathcal M} (Y) \longrightarrow {\mathcal M} (X)
$$
envoie ${\mathcal N}(Y)$ dans ${\mathcal N} (X)$.

\smallskip

En revanche, si ${\mathcal M} \twoheadrightarrow {\mathcal Q}$ est un quotient d'un ${\mathcal O}$-Module ${\mathcal M}$, les applications induites
$$
{\mathcal M}(X) \longrightarrow {\mathcal Q}(X) \, , \quad X \in {\rm Ob} ({\mathcal C}) \, ,
$$
ne sont pas n\'ecessairement surjectives.

\medskip

\item Sur un site annel\'e $({\mathcal C},J,{\mathcal O})$, le faisceau de groupes ab\'eliens ${\mathcal O}$ muni de la loi de multiplication ${\mathcal O} \times {\mathcal O} \to {\mathcal O}$ d\'efinit un objet de ${\mathcal M}od_{{\mathcal C},J,{\mathcal O}}$.

\smallskip

Si ${\mathcal O}$ est un faisceau d'anneaux commutatifs, les sous-Modules du Module ${\mathcal O}$ sont appel\'es les faisceaux d'id\'eaux de $({\mathcal C},J,{\mathcal O})$ ou encore les Id\'eaux de ${\mathcal O}$.

\medskip

\item Sur un site annel\'e $({\mathcal C},J,{\mathcal O})$, on dispose pour toute paire de ${\mathcal O}$-Modules ${\mathcal M}$, ${\mathcal N}$ du faisceau de groupes ab\'eliens
$$
{\mathcal H}om_{\mathcal O} ({\mathcal M}, {\mathcal N})
$$
qui associe \`a  tout objet $X$ de ${\mathcal C}$ le groupe ab\'elien des morphismes de ${\mathcal O}_X$-Modules
$$
{\mathcal M}_X \longrightarrow {\mathcal N}_X
$$
o\`u ${\mathcal O}_X , {\mathcal M}_X$ et ${\mathcal N}_X$ d\'esignent les restrictions des faisceaux ${\mathcal O} , {\mathcal M}$ et ${\mathcal N}$ \`a  la cat\'egorie relative ${\mathcal C}/X$ munie de la topologie induite par $J$.

\smallskip

Si ${\mathcal O}$ est un faisceau d'anneaux commutatifs, le faisceau de groupes ab\'eliens ${\mathcal H}om_{\mathcal O} ({\mathcal M}, {\mathcal N})$ a une structure naturelle de ${\mathcal O}$-Module.

\smallskip

D'autre part, si ${\mathcal M} = {\mathcal N}$, le faisceau ${\mathcal H}om_{\mathcal O} ({\mathcal M}, {\mathcal M})$ est muni d'une loi de composition
$$
{\mathcal H}om_{\mathcal O} ({\mathcal M}, {\mathcal M}) \times {\mathcal H}om_{\mathcal O} ({\mathcal M}, {\mathcal M}) \longrightarrow {\mathcal H}om_{\mathcal O} ({\mathcal M}, {\mathcal M})
$$
et, si de plus ${\mathcal O}$ est un faisceau d'anneaux commutatifs, d'un morphisme naturel
$$
{\mathcal O} \longrightarrow {\mathcal H}om_{\mathcal O} ({\mathcal M}, {\mathcal M})
$$
qui respecte les lois d'addition et de multiplication. 
\end{listeisansmarge}
\end{remarksqed}

\bigskip

Sur un espace annel\'e commutatif, on dispose encore des Modules de d\'erivations:

\begin{defn}\label{defII48}

Soit $(X,{\mathcal O}_X) \to (S,{\mathcal O}_S)$ un morphisme d'espaces annel\'es commutatifs.

\smallskip

On appelle faisceau des d\'erivations de $X$ sur $S$ \`a  valeurs dans un ${\mathcal O}_X$-Module ${\mathcal M}$ et on note
$$
{\mathcal D}er_{X/S} ({\mathcal M})
$$
le ${\mathcal O}_X$-Module qui associe \`a  tout ouvert $U$ de $X$ le ${\mathcal O}_X(U)$-module des morphismes de faisceaux sur $U$
$$
{\rm d} : {\mathcal O}_U \longrightarrow {\mathcal M}_U
$$
qui

$
\left\{ \begin{matrix}
\bullet &\mbox{respectent l'addition des fonctions} \hfill \\
{ \ } \\
&{\rm d}(\varphi_1 + \varphi_2) = {\rm d}\varphi_1 + {\rm d} \varphi_2 \, , \\
{ \ } \\
\bullet &\mbox{envoient sur $0$ les images des morphismes} \hfill \\
{ \ } \\
&{\mathcal O}_S (V) \longrightarrow {\mathcal O}_X (U') \\
{ \ } \\
&\mbox{pour tout ouvert $U'$ de $U$ se projetant dans un ouvert $V$ de $S$,} \hfill \\
{ \ } \\
\bullet &\mbox{satisfont la r\`egle de Leibniz} \hfill \\
{ \ } \\
&{\rm d}(\varphi_1 \cdot \varphi_2) = \varphi_1 \cdot {\rm d} \varphi_2 + \varphi_2 \cdot {\rm d} \varphi_1 \, .
\end{matrix}
\right.
$
\end{defn}

\begin{remarksqed}
\begin{listeisansmarge}
\item Ainsi, ${\mathcal D}er_{X/S} ({\mathcal M})$ est un sous-Module du ${\mathcal O}_X$-Module
$$
{\mathcal H}om_{\mathbb Z} ({\mathcal O}_X , {\mathcal M}) 
$$
qui associe \`a  tout ouvert $U$ de $X$ le module des morphismes de faisceaux ${\mathcal O}_U \to {\mathcal M}_U$ qui respectent l'addition.

\medskip

\item En particulier, si ${\mathcal M} = {\mathcal O}_X$, ${\mathcal D}er_{X/S} ({\mathcal O}_X)$ est un sous-Module du ${\mathcal O}_X$-Module
$$
{\mathcal H}om_{\mathbb Z} ({\mathcal O}_X , {\mathcal O}_X)
$$
lequel est muni d'une loi de composition
$$
{\mathcal H}om_{\mathbb Z} ({\mathcal O}_X , {\mathcal O}_X) \times {\mathcal H}om_{\mathbb Z} ({\mathcal O}_X , {\mathcal O}_X) \longrightarrow {\mathcal H}om_{\mathbb Z} ({\mathcal O}_X , {\mathcal O}_X)
$$
et d'un morphisme naturel respectant l'addition et la multiplication
$$
{\mathcal O}_X \longrightarrow {\mathcal H}om_{\mathbb Z} ({\mathcal O}_X , {\mathcal O}_X) \, .
$$

On note ${\mathcal D}_X$ et on appelle faisceau des op\'erateurs diff\'erentiels lin\'eaires de $X$ sur $S$ le plus petit sous-Module de ${\mathcal H}om_{\mathbb Z} ({\mathcal O}_X , {\mathcal O}_X)$ qui contient ${\mathcal O}_X$ et ${\mathcal D}er_{X/S} ({\mathcal O}_X)$ et qui est stable par la loi de composition.

\smallskip

Ainsi, ${\mathcal D}_X$ est un faisceau d'anneaux (non commutatifs) sur l'espace topologique $X$. 
\end{listeisansmarge}
\end{remarksqed}

\pagebreak

On pose la d\'efinition naturelle:

\begin{defn}\label{defII49}

Soit $({\mathcal C} , J , {\mathcal O})$ un site annel\'e.

\begin{listeimarge}

\item Deux ${\mathcal O}$-Modules ${\mathcal M}$ et ${\mathcal N}$ sont dits ``localement isomorphes'' si tout objet $X$ de ${\mathcal C}$ admet un crible couvrant $S$ tel que les restrictions de ${\mathcal M}$ et ${\mathcal N}$ \`a  $S$ (vu comme une sous-cat\'egorie pleine de ${\mathcal C}/X$ et muni de la topologie induite par $J$ et de la restriction ${\mathcal O}_S$ du faisceau d'anneaux ${\mathcal O}$) soient isomorphes comme ${\mathcal O}_S$-Modules.

\medskip

\item Un ${\mathcal O}$-Module ${\mathcal M}$ est dit ``localement libre de rang $n$'' s'il est localement isomorphe au ${\mathcal O}$-Module ``libre de rang $n$''
$$
{\mathcal O}^n = \underbrace{{\mathcal O} \times \cdots \times {\mathcal O}}_{n \ {\rm fois}} \, .
$$
\end{listeimarge}
\end{defn}

Un exemple tr\`es important de Module localement libre est le faisceau tangent d'une vari\'et\'e diff\'erentielle ou analytique, ou encore d'un sch\'ema lisse sur un sch\'ema de base:

\begin{prop}\label{propII410}
\begin{listeimarge}
\item Soit $p : X \to S$ un morphisme de vari\'et\'es diff\'erentielles de classe $C^{\infty}$ [resp. analytiques] qui est une submersion de dimension relative $n$.

\smallskip

Alors le ${\mathcal O}_X$-Module
$$
{\mathcal D}er_{X/S} ({\mathcal O}_X) = \Omega_{X/S}^{\vee} \, ,
$$
qui est appel\'e le faisceau tangent de $X$ relativement \`a  $S$, est localement libre de rang $n$.

\medskip

\item Soit $p : X \to S$ un morphisme de sch\'emas qui est lisse de dimension relative $n$.

\smallskip

Alors, pour tout ${\mathcal O}_X$-Module ${\mathcal M}$, le ${\mathcal O}_X$-Module
$$
{\mathcal D}er_{X/S} ({\mathcal M})
$$
est localement isomorphe \`a  ${\mathcal M}^n = {\mathcal M} \times \cdots \times {\mathcal M}$.

\smallskip

En particulier, le faisceau tangent de $X$ relativement \`a  $S$ d\'efini comme
$$
\Omega_{X/S}^{\vee} = {\mathcal D}er_{X/S} ({\mathcal O}_X)
$$
est un ${\mathcal O}_X$-Module localement libre de rang $n$.
\end{listeimarge}
\end{prop}

\begin{remarks}
\begin{listeisansmarge}
\item Si $X$ est une vari\'et\'e diff\'erentielle ou analytique de dimension $n$ et $S$ est la vari\'et\'e r\'eduite \`a  un point, l'unique projection $p : X \to S$ est une submersion.

\smallskip

Dans ce cas, le ${\mathcal O}_X$-Module
$$
\Omega_{X/S}^{\vee} = \Omega_X^{\vee}
$$
est simplement appel\'e le faisceau tangent de $X$.

\smallskip

Il est localement libre de rang $n = {\rm dim} \, (X)$.

\medskip

\item L'\'enonc\'e de la proposition ne s'applique plus si $(X,{\mathcal O}_X)$ est une vari\'et\'e diff\'erentielle de classe $C^k$, si $k < \infty$ est un entier $\geq 1$.
\end{listeisansmarge}
\end{remarks}

\begin{demo}

L'\'enonc\'e est local sur $X$ et sur $S$.

\begin{listeisansmarge}

\item Dans le cas des vari\'et\'es diff\'erentielles [resp. analytiques], on peut supposer que $S$ est un ouvert $V$ de ${\mathbb R}^m$ [resp. ${\mathbb C}^m$], que $X$ est le produit de $V$ et d'un ouvert $W$ de ${\mathbb R}^n$ [resp. ${\mathbb C}^n$] et que la submersion $p : X \to S$ est la projection $W \times V \to V$.

\smallskip

Notons $y_1 , \cdots , y_m$ les fonctions de coordonn\'ees cart\'esiennes de l'ouvert $V$ de ${\mathbb R}^m$ [resp. ${\mathbb C}^m$] et $z_1 , \cdots , z_n$ celles de l'ouvert $W$ de ${\mathbb R}^n$ [resp. ${\mathbb C}^n$].

\smallskip

Il suffit de montrer que l'application ${\mathcal O}_X (X)$-lin\'eaire
$$
\begin{matrix}
{\mathcal D}er_{X/S} ({\mathcal O}_X)(X) &\longrightarrow &{\mathcal O}_X (X)^n \, , \hfill \\
\hfill {\rm d} &\longmapsto &({\rm d} z_1 , \cdots , {\rm d} z_n) 
\end{matrix}
$$
est un isomorphisme.

\smallskip

Ce morphisme est surjectif car toute famille de fonctions
$$
(\varphi_1 , \cdots , \varphi_n) \in {\mathcal O}_X(X)^n
$$
admet pour ant\'ec\'edent la d\'erivation
$$
{\rm d} : {\mathcal O}_X \longrightarrow {\mathcal O}_X
$$
d\'efinie par la formule
$$
\varphi \longmapsto {\rm d}\varphi = \sum_{1 \leq i \leq n} \frac{\partial \varphi}{\partial z_i} \cdot \varphi_i \, .
$$

Pour montrer qu'il est injectif, consid\'erons une d\'erivation ${\rm d} \in {\mathcal D}er_{X/S} ({\mathcal O}_X)(X)$ telle que ${\rm d} z_1 = 0 , \cdots , {\rm d} z_n = 0$ et montrons que ${\rm d}=0$.

\smallskip

Il suffit de prouver que pour toute fonction $\varphi \in {\mathcal O}_X (U)$ d\'efinie sur un ouvert convexe $U$ de $X = W \times V$ et tout point $x^0 = (z_1^0 , \cdots , z_n^0 , y_1^0 , \cdots , y_m^0)$ de $U$, on a
$$
({\rm d} \varphi)(x^0)=0 \, .
$$
Or la fonction $\varphi$ des points $x = (z_1 , \cdots ,z_n , y_1 , \cdots , y_m) \in U$ s'\'ecrit
$$
\varphi (x) = \varphi (x^0) + \sum_{1 \leq i \leq n} (z_i - z_i^0) \cdot \varphi_i (x) + \sum_{1 \leq j \leq m} (y_j - y_j^0) \cdot \psi_j(x)
$$
o\`u les $\varphi_i$ et $\psi_j$ sont les fonctions de classe $C^{\infty}$ [resp. analytiques] d\'efinies par les int\'egrales
$$
\varphi_i (x) = \int_0^1 {\rm d}t \cdot \frac{\partial \varphi}{\partial z_i} (t \cdot x + (1 - t) \cdot x^0) \, ,
$$
$$
\psi_j (x) = \int_0^1 {\rm d}t \cdot \frac{\partial \varphi}{\partial y_j} (t \cdot x + (1 - t) \cdot x^0) \, .
$$

Comme l'image par ${\rm d}$ de la constante $\varphi (x^0)$ et des fonctions $z_i - z_i^0$ et $y_j - y_j^0$ est $0$ et que les fonctions $z_i - z_i^0$ et $y_j - y_j^0$ s'annulent en $0$, on conclut comme voulu
$$
{\rm d}\varphi (x_0) = 0 \, .
$$

\medskip

\item Par d\'efinition des morphismes lisses de sch\'emas $X \to S$ dans la d\'efinition \ref{defII317}, on peut supposer que $X = {\rm Spec} (B)$ et $S = {\rm Spec} (A)$ sont deux sch\'emas affines et que le $B$-module des diff\'erentielles
$$
\Omega_{B/A}
$$
est libre de la forme $\underset{1 \leq i \leq n}{\bigoplus} B \cdot {\rm d}f_i$ pour des \'el\'ements $f_i \in B$.

\smallskip

D'autre part, on sait d'apr\`es la remarque (iv) qui suit cette d\'efinition que pour tout \'el\'ement $f \in B$, $\Omega_{B_f / A}$ s'identifie \`a $B_f \otimes_B \Omega_{B/A}$.

\smallskip

On en d\'eduit comme voulu que pour tout ${\mathcal O}_X$-Module ${\mathcal M}$, l'application $B$-lin\'eaire
$$
\begin{matrix}
{\mathcal D}er_{X/S} ({\mathcal M})(X) &\longrightarrow &{\mathcal M} (X)^n \, , \hfill \\
\hfill {\rm d} &\longmapsto &({\rm d}f_1 , \cdots , {\rm d}f_n)
\end{matrix}
$$
est un isomorphisme. 
\end{listeisansmarge}
\end{demo}

\subsection{Faisceaux coh\'erents}\label{subsec245}

\medskip

Consid\'erons un site annel\'e $({\mathcal C} , J , {\mathcal O})$.

\smallskip

Se donner $n$ sections $m_1 , \cdots , m_n$ d'un ${\mathcal O}$-Module ${\mathcal M}$ sur un objet $U$ de ${\mathcal C}$ \'equivaut \`a se donner un morphisme de ${\mathcal O}_U$-Modules
$$
{\mathcal O}_U^n \longrightarrow {\mathcal M}_U
$$
entre les restrictions ${\mathcal O}^n_U$ et ${\mathcal M}_U$ de ${\mathcal O}^n$ et ${\mathcal M}$ \`a la cat\'egorie relative ${\mathcal C}/U$ munie de la topologie $J$.

\smallskip

Si ce morphisme ${\mathcal O}_U^n \to {\mathcal M}_U$ est un \'epimorphisme, on dit que les sections $m_1 , \cdots , m_n \in {\mathcal M} (U)$ engendrent le ${\mathcal O}$-Module ${\mathcal M}$ au-dessus de $U$.

\smallskip

On peut poser:

\begin{defn}\label{defII411}

Soit $({\mathcal C} , J , {\mathcal O})$ un site annel\'e.

\begin{listeimarge}

\item Un ${\mathcal O}$-Module ${\mathcal M}$ est dit ``fini'' si tout objet $X$ de ${\mathcal C}$ admet un crible couvrant constitu\'e de morphismes $U \to X$ tels que la restriction ${\mathcal M}_U$ \`a chaque ${\mathcal C} / U$ soit engendr\'ee par un nombre fini de sections dans ${\mathcal M} (U)$.

\medskip

\item Un ${\mathcal O}$-Module ${\mathcal M}$ est dit ``coh\'erent'' s'il est fini et si, pour tout objet $X$ de ${\mathcal C}$ et tout \'epimorphisme
$$
{\mathcal O}_X^n \longrightarrow {\mathcal M}_X \, ,
$$
le noyau de celui-ci est un ${\mathcal O}_X$-Module fini. 
\end{listeimarge}
\end{defn}

Les Modules coh\'erents sur les sch\'emas affines sont d\'ecrits par le th\'eor\`eme suivant:

\begin{thm}\label{thmII412}

Soit $(X , {\mathcal O}) = {\rm Spec} (A)$ le sch\'ema affine associ\'e \`a un anneau commutatif $A$.

\smallskip

Alors:

\begin{listeimarge}

\item Pour tout $A$-module $M$, il existe un ${\mathcal O}$-Module $\widetilde M$, unique \`a unique isomorphisme pr\`es, tel que pour tout \'el\'ement $f$ de $A$, on ait
$$
\widetilde M ({\rm Spec} (A)_f) = M_f = A_f \otimes_A M \, .
$$

\item Le foncteur
$$
\begin{matrix}
{\rm Mod}_A &\longrightarrow &{\mathcal M}od_{\mathcal O} \, , \\
\hfill M &\longmapsto &\widetilde M \hfill
\end{matrix}
$$
est adjoint \`a gauche du foncteur des sections globales
$$
\begin{matrix}
{\mathcal M}od_{\mathcal O} &\longrightarrow &{\rm Mod}_A \, , \\
\hfill {\mathcal M} &\longmapsto &{\mathcal M} (X) \, , \hfill
\end{matrix}
$$
et son compos\'e avec celui-ci est canoniquement isomorphe au foncteur
$$
{\rm id} : {\rm Mod}_A \longrightarrow {\rm Mod}_A \, .
$$
Il est pleinement fid\`ele, et il respecte les colimites arbitraires ainsi que les limites finies.

\medskip

\item Un ${\mathcal O}$-Module ${\mathcal M}$ est coh\'erent si et seulement si il est isomorphe \`a un ${\mathcal O}$-Module de la forme
$$
\widetilde M
$$
pour un $A$-module $M$ de pr\'esentation finie.
\end{listeimarge}
\end{thm}

\begin{demosansqed}
\begin{listeisansmarge}
\item L'unicit\'e du faisceau $\widetilde M$ associ\'e \`a un $A$-module $M$ r\'esulte de ce que, par d\'efinition de la topologie de Zariski, tout ouvert de ${\rm Spec} (A)$ est r\'eunion d'ouverts de la forme ${\rm Spec} (A)_f$.

\smallskip

Son existence r\'esulte du lemme \ref{lemI46} (iii) en d\'ecidant que, pour tout ouvert $U$ de ${\rm Spec} (A)$, $\widetilde M (U)$ est le ${\mathcal O} (U)$-Module constitu\'e des familles d'\'el\'ements
$$
m_f \in M_f = A_f \otimes_A M_f \quad \mbox{index\'es par les $f \in A$ v\'erifiant ${\rm Spec} (A)_f \subset U$,}
$$
telles que, pour tout multiple $f' = fa$ d'un tel $f$, $m_f$ s'envoie sur $m_{f'}$ par le morphisme $M_f \to M_{f'}$ induit par le morphisme canonique $M_f \to M_{f'}$.

\medskip

\item Tout morphisme de $A$-modules $M \to N$ induit pour tout \'el\'ement $f \in A$ un morphisme de $A_f$-modules
$$
M_f \longrightarrow N_f \, .
$$
De plus, pour tout multiple $f' = fa$ d'un tel \'el\'ement $f$, le carr\'e
$$
\xymatrix{
M_f \ar[d] \ar[r] &N_f \ar[d] \\
M_{f'} \ar[r] &N_{f'}
}
$$
est commutatif.

\smallskip

Cela d\'efinit un foncteur
$$
\begin{matrix}
{\rm Mod}_A &\longrightarrow &{\mathcal M}od_{\mathcal O} \, , \\
\hfill M &\longmapsto &\widetilde M \, . \hfill
\end{matrix}
$$
Montrons qu'il est adjoint \`a gauche du foncteur ${\mathcal M} \to {\mathcal M} (X)$. Consid\'erons pour cela un $A$-module $M$ et un ${\mathcal O}$-Module ${\mathcal M}$. Tout morphisme de ${\mathcal O}$-Modules
$$
\widetilde M \longrightarrow {\mathcal M}
$$
induit un morphisme de $A$-modules
$$
M = \widetilde M(X) \longrightarrow {\mathcal M} (X) \, .
$$
R\'eciproquement, tout morphisme de $A$-modules
$$
M \longrightarrow {\mathcal M} (X)
$$
induit pour tout \'element $f \in A$ un morphisme de $A_f$-modules
$$
M_f = A_f \otimes_A M \longrightarrow A_f \otimes_A {\mathcal M} (X) \longrightarrow {\mathcal M} ({\rm Spec} (A)_f)
$$
par composition avec le morphisme de restriction
$$
{\mathcal M} (X) \longrightarrow {\mathcal M} ({\rm Spec} (A)_f) \, .
$$

Pour tout multiple $f' = fa$ d'un tel $f$, le carr\'e
$$
\xymatrix{
M_f \ar[d] \ar[r] &{\mathcal M} ({\rm Spec} (A)_f) \ar[d] \\
M_{f'} \ar[r] &{\mathcal M} ({\rm Spec} (A)_{f'})
}
$$
est commutatif, et donc les morphismes $M_f \to {\mathcal M} ({\rm Spec} (A)_f)$ d\'efinissent un morphisme de ${\mathcal O}$-Modules
$$
\widetilde M \longrightarrow {\mathcal M} \, .
$$
Le compos\'e du foncteur $M \mapsto \widetilde M$ et de son adjoint \`a droite ${\mathcal M} \mapsto {\mathcal M}(X)$ s'identifie au foncteur ${\rm id} : {\rm Mod}_A \to {\rm Mod}_A$ puisque, pour tout $A$-module $M$, $\widetilde M (X)$ s'identifie \`a $M$ par construction.

\smallskip

Il en r\'esulte que le foncteur $M \mapsto \widetilde M$ est pleinement fid\`ele.

\smallskip

Il respecte les colimites arbitraires puisque c'est un adjoint \`a gauche.

\smallskip

Enfin, il respecte les limites finies puisque, pour tout \'el\'ement $f \in A$, $A_f$ est plat sur $A$ d'apr\`es la remarque (viii) qui suit la d\'efinition \ref{defII315}. Cela signifie en effet que, pour tout tel \'el\'ement $f$, le foncteur
$$
\begin{matrix}
{\rm Mod}_A &\longrightarrow &{\rm Mod}_{A_f} \, , \hfill\\
\hfill M &\longmapsto &M_f = A_f \otimes_A M
\end{matrix}
$$
respecte les limites finies.

\medskip

\item Montrons d'abord que si $M$ est un $A$-module de pr\'esentation finie, alors le ${\mathcal O}$-Module $\widetilde M$ est coh\'erent.

\smallskip

Le $A$-module $M$ est a fortiori de type fini, c'est-\`a-dire admet un \'epimorphisme de $A$-modules
$$
A^n \longrightarrow M
$$
et celui-ci induit un \'epimorphisme de ${\mathcal O}$-Modules
$$
\widetilde{A^n} = {\mathcal O}^n \longrightarrow \widetilde M \, .
$$
Il reste \`a montrer que pour tout \'el\'ement $f \in A$ et tout \'epimorphisme de ${\mathcal O}$-Modules sur ${\rm Spec} (A)_f$
$$
p : {\mathcal O}^m \longrightarrow \widetilde M \, ,
$$
le noyau de $p$ est un ${\mathcal O}$-Module fini. Or ce noyau est de la forme
$$
{\rm Ker} (p) = \widetilde N
$$
si $\widetilde N$ est le ${\mathcal O}$-Module sur ${\rm Spec} (A)_f$ d\'efini par le $A_f$-module
$$
N = {\rm Ker} \left( A_f^m \xrightarrow{ \ p \ } M_f \right).
$$
On est donc ramen\'e au lemme suivant:
\end{listeisansmarge}
\end{demosansqed}

\begin{lem}\label{lemII413}

Soient $A$ un anneau commutatif et $M$ un $A$-module de pr\'esentation finie.

\smallskip

Alors le noyau de tout \'epimorphisme de $A$-modules
$$
p : A^m \longrightarrow M
$$
est de type fini.
\end{lem}


\begin{demolem}

Par hypoth\`ese, le $A$-module $M$ s'inscrit dans une suite exacte de $A$-modules
$$
0 \longrightarrow L \longrightarrow A^n \xrightarrow{ \ q \ } M \longrightarrow 0
$$
dans laquelle $L$ est engendr\'e par un nombre fini d'\'el\'ements.

\smallskip

Il s'agit de prouver que le noyau $K$ de $p : A^m \to M$ est \'egalement engendr\'e par un nombre fini d'\'el\'ements.

\smallskip

Comme $q$ et $p$ sont des \'epimorphismes, il existe deux morphismes de $A$-modules
$$
p' : A^m \longrightarrow A^n \quad \mbox{et} \quad q' : A^n \longrightarrow A^m
$$
tels que
$$
q \circ p' = p \qquad \mbox{et} \qquad p \circ q'=q \, .
$$
On a un carr\'e commutatif
$$
\xymatrix{
A^m \oplus A^n \ar[d]_{q_1} \ar[r]^-{p_1} &A_n \ar[d]^q \\
A^m \ar[r]^p &M
}
$$
avec
$$
\begin{matrix}
&&p_1 (a,b) = p'(a) + b \, , \\
&&q_1 (a,b) = a + q' (b) \hfill
\end{matrix}
$$
et
$$
r = q \circ p_1 = p \circ q_1 : (a,b) \longmapsto p(a) + q(b) \, .
$$

On a deux suites exactes
$$
\begin{matrix}
0 \longrightarrow A^m &\longrightarrow &A^m \oplus A^n \xrightarrow{ \ p_1 \ } A^n \longrightarrow 0 \\
\hfill a &\longmapsto &(a,-p'(a)) \hfill
\end{matrix}
$$
et
$$
\begin{matrix}
0 \longrightarrow A^n &\longrightarrow &A^m \oplus A^n \xrightarrow{ \ q_1 \ } A^m \longrightarrow 0 \\
\hfill b &\longmapsto &(-q'(b),b) \hfill
\end{matrix}
$$
si bien que le noyau $J$ de $r = q \circ p_1 = p \circ q_1$ s'inscrit dans deux suites exactes
$$
0 \longrightarrow A^m \longrightarrow J \xrightarrow{ \ p_1 \ } L \longrightarrow 0
$$
et
$$
0 \longrightarrow A^n \longrightarrow J \xrightarrow{ \ q_1 \ } K \longrightarrow 0 \, .
$$
Comme $L$ est engendr\'e par un nombre fini d'\'el\'ements, il en est de m\^eme de $J$ et donc aussi de son quotient $K$.

\smallskip

C'est ce que l'on voulait. 
\end{demolem}


\noindent {\bf Fin de la d\'emonstration du th\'eor\`eme \ref{thmII412} (iii):}

\smallskip

Il reste \`a prouver que si ${\mathcal M}$ est un ${\mathcal O}$-Module coh\'erent sur ${\rm Spec} (A)$, il est de la forme $\widetilde M$ pour un $A$-module $M$ de pr\'esentation finie.

\smallskip

Il existe un recouvrement de $X = {\rm Spec} (A)$ par les ouverts affines ${\rm Spec} (A)_{f_i}$, $1 \leq i \leq n$, telle que la restriction ${\mathcal M}_i$ de ${\mathcal M}$ \`a chaque ${\rm Spec} (A)_{f_i}$ admette un \'epimorphisme
$$
{\mathcal O}^{m_i} \longrightarrow {\mathcal M}_i
$$
dont le noyau ${\mathcal N}_i$ admet lui-m\^eme un \'epimorphisme
$$
{\mathcal O}^{n_i} \longrightarrow {\mathcal N}_i \, .
$$
D'apr\`es les propri\'et\'es g\'en\'erales des cat\'egories de faisceaux lin\'eaires que l'on verra au chapitre III, cela signifie que
$$
{\mathcal M}_i = {\rm coker} ({\mathcal O}^{n_i} \longrightarrow {\mathcal O}^{m_i})
$$
et donc
$$
{\mathcal M}_i = \widetilde M_i
$$
o\`u $M_i$ est le $A_{f_i}$-module de pr\'esentation finie
$$
M_i = {\rm coker} (A_{f_i}^{n_i} \longrightarrow A_{f_i}^{m_i}) \, .
$$
Pour tous indices $i$ et $j$, les $A_{f_i f_j}$-modules
$$
A_{f_i f_j} \otimes_{A_{f_i}} M_i \qquad \mbox{et} \qquad A_{f_i f_j} \otimes_{A_{f_j}} M_j
$$
s'identifient tous deux \`a ${\mathcal M} ({\rm Spec} (A)_{f_i f_j})$.

\smallskip

Comme ${\mathcal M}$ est un faisceau, le $A$-module ${\mathcal M}(X)$ s'identifie \`a l'\'egalisateur
$$
M = {\rm eg} \left( \prod_{1 \leq i \leq n} M_i \rightrightarrows \prod_{1 \leq i,j \leq n} A_{f_i f_j} \otimes_{A_{f_i}} M_i \right).
$$
Pour tout indice $k$, on a
\begin{eqnarray}
A_{f_k} \otimes_A M &= &{\rm eg} \left( \prod_{1 \leq i \leq n} A_{f_i f_k} \otimes_{A_{f_i}} M_i \rightrightarrows \prod_{1 \leq i,j \leq n} A_{f_i f_j f_k} \otimes_{A_{f_i}} M_i \right) \nonumber \\
&= &{\rm eg} \left( \prod_{1 \leq i \leq n} A_{f_i f_k} \otimes_{A_{f_k}} M_k \rightrightarrows \prod_{1 \leq i,j \leq n} A_{f_i f_j f_k} \otimes_{A_{f_k}} M_k \right) \nonumber \\
&= &M_k = {\mathcal M} ({\rm Spec} (A)_{f_k}) \, . \nonumber
\end{eqnarray}

Cela prouve que le morphisme canonique
$$
\widetilde M \longrightarrow {\mathcal M}
$$
est un isomorphisme.

\smallskip

Il reste \`a prouver que le $A$-module $M$ est de pr\'esentation finie.

\smallskip

Chaque $M_{f_i} = A_{f_i} \otimes_A M = M_i$ est un $A_{f_i}$-module de type fini. Donc on peut trouver dans $M$ une famille finie d'\'el\'ements qui engendre chaque $M_{f_i}$ comme $A_{f_i}$-module. Alors cette famille engendre $M$ c'est-\`a-dire d\'efinit un \'epimorphisme
$$
A^m \longrightarrow M \, .
$$
Soit $K$ le noyau de cet \'epimorphisme.

\smallskip

Pour tout $i$, $1 \leq i \leq n$, $A_{f_i} \otimes_A K$ est le noyau du morphisme induit
$$
A_{f_i}^m \longrightarrow M_{f_i} = M_i
$$
donc, d'apr\`es le lemme \ref{lemII413}, il est engendr\'e par un nombre fini d'\'el\'ements.

\smallskip

On peut alors trouver dans $K$ une famille finie d'\'el\'ements qui engendre chaque $A_{f_i} \otimes_A K$ sur $A_{f_i}$ donc engendrent $K$ comme $A$-module.

\smallskip

Cela termine la d\'emonstration du th\'eor\`eme. \hfill $\Box$

\bigskip

Le th\'eor\`eme \ref{thmII412} implique:

\begin{cor}\label{corII414}

Sur un sch\'ema $(X,{\mathcal O})$, un ${\mathcal O}$-Module ${\mathcal M}$ est coh\'erent si et seulement si sa restriction \`a tout ouvert affine ${\rm Spec} (A)$ de $X$ a la forme
$$
\widetilde M
$$
pour un $A$-module $M$ de pr\'esentation finie.
\end{cor}

\begin{remarksqed}
\begin{listeisansmarge}
\item La propri\'et\'e d'\^etre un ${\mathcal O}$-Module coh\'erent est locale.

\smallskip

Pour qu'un ${\mathcal O}$-Module ${\mathcal M}$ soit coh\'erent, il suffit donc qu'existe un recouvrement ouvert de $X$ par des sch\'emas affines ${\rm Spec} (A)$ tels que la restriction de ${\mathcal M}$ \`a chaque ${\rm Spec} (A)$ soit de la forme $\widetilde M$ pour un $A$-module $M$ de pr\'esentation finie.

\medskip

\item Un ${\mathcal O}$-Module ${\mathcal M}$ est dit ``quasi-coh\'erent'' si sa restriction \`a tout ouvert affine ${\rm Spec} (A)$ de $X$ a la forme
$$
\widetilde M
$$
pour un $A$-module $M$.

\smallskip

Il r\'esulte de la d\'emonstration du th\'eor\`eme \ref{thmII412} (iii) que cette propri\'et\'e est locale: pour qu'un ${\mathcal O}$-Module ${\mathcal M}$ sur un sch\'ema $X$ soit quasi-coh\'erent, il suffit qu'il le soit sur les \'el\'ements d'un recouvrement ouvert de $X$. 
\end{listeisansmarge}
\end{remarksqed}

\bigskip

Voici une famille importante de Modules quasi-coh\'erents (au sens de la derni\`ere remarque) ou coh\'erents:

\begin{prop}\label{propII415}

Consid\'erons un morphisme de sch\'emas
$$
(X , {\mathcal O}_X) \longrightarrow (S , {\mathcal O}_S) \, .
$$
Alors:

\begin{listeimarge}

\item Il existe sur $X$ un ${\mathcal O}_X$-Module quasi-coh\'erent
$$
\Omega_{X/S} \, ,
$$
unique \`a unique isomorphisme pr\`es et appel\'e le Module des diff\'erentielles relatives de $X$ sur $S$, tel que pour tout ouvert affine ${\rm Spec} (B)$ de $X$ s'envoyant dans un ouvert affine ${\rm Spec} (A)$ de $S$, on ait
$$
\Omega_{X/S} ({\rm Spec} (B)) = \Omega_{B/A} \, .
$$

\item Si le morphisme
$$
X \longrightarrow S
$$
est localement de pr\'esentation finie, le ${\mathcal O}_X$-Module $\Omega_{X/S}$ des diff\'erentielles relatives de $X$ sur $S$ est coh\'erent.
\end{listeimarge}
\end{prop}

\begin{demo}
\begin{listeisansmarge}
\item L'unicit\'e r\'esulte de ce que les ouverts affines de $X$ qui s'envoient dans des ouverts affines de $S$ forment une base de la topologie de $X$.

\smallskip

L'existence de $\Omega_{X/S}$ comme ${\mathcal O}_X$-Module quasi-coh\'erent r\'esulte des deux faits suivants:

\smallskip

D'une part, si un morphisme d'anneaux commutatifs
$$
A \longrightarrow B
$$
se factorise en
$$
A \longrightarrow A_f \longrightarrow B
$$
pour un \'element $f \in A$, alors le morphisme canonique de $B$-modules
$$
\Omega_{B/A_f} \longrightarrow \Omega_{B/A}
$$
est un isomorphisme. En effet, une diff\'erentielle $d : B \to M$ \`a valeurs dans un $B$-module $M$ s'annule sur les \'el\'ements de $A$ si et seulement si elle s'annule sur ceux de $A_f$.

\smallskip

D'autre part, pour tout morphisme d'anneaux commutatifs
$$
A \longrightarrow B
$$
et tout \'el\'ement $f \in B$, les $B_f$-modules
$$
\Omega_{B_f/A} \qquad \mbox{et} \qquad B_f \otimes_B \Omega_{B/A}
$$
s'identifient, comme on l'avait d\'ej\`a not\'e dans la remarque (iv) qui suit la d\'efinition \ref{defII317}.

\medskip

\item r\'esulte de ce que, si $B$ est le quotient d'une alg\`ebre de polynomes $A [X_1 , \cdots ,X_n]$ par un id\'eal engendr\'e par un nombre fini d'\'el\'ements $P_1 , \cdots , P_k$, alors $\Omega_{B/A}$ est un $B$-module de pr\'esentation finie d'apr\`es le corollaire~\ref{corII310} (ii). 
\end{listeisansmarge}
\end{demo}

\section{Faisceautisation et foncteurs canoniques}\label{sec25}

\subsection{Le foncteur de faisceautisation}\label{subsec251}

\medskip

Pour tout site $({\mathcal C},J)$, la cat\'egorie $\widehat{\mathcal C}_J$ des faisceaux sur $({\mathcal C},J)$ est par d\'efinition la sous-cat\'egorie pleine de $\widehat{\mathcal C} = [{\mathcal C}^{\rm op} , {\rm Ens}]$ constitu\'ee des pr\'efaisceaux d'ensembles sur ${\mathcal C}$ qui sont des faisceaux pour la topologie $J$.

\smallskip

Elle est donc munie du foncteur pleinement fid\`ele
$$
j_* : \widehat{\mathcal C}_J \longrightarrow \widehat{\mathcal C}
$$
qui consiste \`a voir tout faisceau comme un pr\'efaisceau.

\smallskip

Nous allons montrer que ce foncteur de plongement admet un adjoint \`a gauche
$$
j^* : \widehat{\mathcal C} \longrightarrow \widehat{\mathcal C}_J
$$
qui donc transforme tout pr\'efaisceau en un faisceau.

\smallskip

Pour cela, on a besoin de la d\'efinition suivante:

\begin{defn}\label{defII51}

Soient $({\mathcal C},J)$ un site et
$$
F : {\mathcal C}^{\rm op} \longrightarrow {\rm Ens}
$$
un pr\'efaisceau sur ${\mathcal C}$.

\begin{listeimarge}

\item On dit que deux sections de $F$ sur un objet $X$ de ${\mathcal C}$
$$
a,b \in F(X)
$$
co{\"\i}ncident localement pour la topologie $J$ s'il existe un crible $J$-couvrant $S$ de $X$ tel que, pour tout \'el\'ement $u : U \to X$ de $S$, on ait l'\'egalit\'e
$$
F(u)(a) = F(u)(b) \quad \mbox{dans} \quad F(U) \, .
$$

\item Pour tout objet $X$ de ${\mathcal C}$, on note $F_s (X)$ l'ensemble quotient de $F(X)$ par la relation de co{\"\i}ncidence locale.

\medskip

\item On dit que le pr\'efaisceau $F$ est s\'epar\'e pour la topologie $J$ si, pour tout objet $X$ de ${\mathcal C}$, toutes sections de $F$ sur $X$ qui co{\"\i}cident localement sont \'egales ou, ce qui revient au m\^eme, la surjection canonique
$$
F(X) \longrightarrow F_s (X)
$$
est une bijection.
\end{listeimarge}
\end{defn}

\begin{remarksqed}
\begin{listeisansmarge}
\item La relation de co{\"\i}cidence locale est une relation d'\'equivalence sur chaque ensemble $F(X)$ car si une section $b \in F(X)$ co{\"\i}ncide avec une section $a \in F(X)$ sur les \'el\'ements $u : U \to X$ d'un crible couvrant $S \in J(X)$ et une section $c \in F(X)$ co{\"\i}ncide avec $b$ sur les \'el\'ements d'un crible couvrant $S' \in J(X)$, alors $a$ et $c$ co{\"\i}ncident sur tous les \'el\'ements de n'importe quel crible couvrant de $X$ qui raffine \`a la fois $S$ et $S'$.

\medskip

\item Il r\'esulte de l'axiome de stabilit\'e v\'erifi\'e par la topologie de Grothendieck $J$ que, pour tout morphisme $f : X \to Y$ de ${\mathcal C}$, l'application
$$
F(f) : F(Y) \longrightarrow F(X)
$$
d\'efinit par passage aux quotients une application
$$
F_s (f) : F_s (Y) \longrightarrow F_s (X) \, .
$$
Ainsi, $F_s$ est d\'efini comme pr\'efaisceau sur ${\mathcal C}$ reli\'e \`a $F$ par un \'epimorphisme canonique
$$
F -\!\!\!-\!\!\!\twoheadrightarrow F_s \, .
$$
Associer \`a tout pr\'efaisceau $F$ le pr\'efaisceau quotient $F_s$ d\'efinit un foncteur
$$
\begin{matrix}
F &\longmapsto &F_s \, , \\
\widehat{\mathcal C} &\longrightarrow &\widehat{\mathcal C} \hfill
\end{matrix}
$$
reli\'e au foncteur identit\'e par une transformation naturelle
$$
{\rm id} \longrightarrow (F \mapsto F_s) \, .
$$
\end{listeisansmarge}
\end{remarksqed}

On a:

\begin{lem}\label{lemII52}

Soit $({\mathcal C},J)$ un site.

\begin{listeimarge}

\item Pour tout pr\'efaisceau $F$ sur ${\mathcal C}$, le pr\'efaisceau $F_s$ est s\'epar\'e pour la topologie $J$.

\medskip

\item Le foncteur
$$
F \longmapsto F_s
$$
est adjoint \`a gauche du foncteur de plongement
$$
\widehat{\mathcal C}_J^+ \xhookrightarrow { \ { \ } \ } \widehat{\mathcal C}
$$
de la sous-cat\'egorie pleine $\widehat{\mathcal C}_J^+$ de $\widehat{\mathcal C}$ constitu\'ee des pr\'efaisceaux qui sont s\'epar\'es pour la topologie $J$.

\medskip

\item Le foncteur
$$
F \longmapsto F_s
$$
respecte les colimites arbitraires.
\end{listeimarge}
\end{lem}

\begin{demo}
\begin{listeisansmarge}
\item Consid\'erons deux sections $\overline a , \overline b \in F_s (X)$ de $F_s$ sur un objet $X$ de ${\mathcal C}$ et un crible couvrant $S \in J(X)$ dont tout \'el\'ement $u : U \to X$ v\'erifie
$$
F_s (u)(\overline a) = F_s (u) (\overline b) \quad \mbox{dans} \quad F_s (U) \, .
$$
Les \'el\'ements $\overline a , \overline b \in F_s (X)$ se rel\`event en deux \'el\'ements $a,b \in F(X)$.

\smallskip

Pour tout \'el\'ement $u : U \to X$ de $S$, les \'el\'ements $F(u)(a)$ et $F(u)(b)$ co{\"\i}ncident localement, ce qui signifie qu'existe un crible couvrant $S_u \in J(U)$ dont tout \'el\'ement $v : V \to U$ v\'erifie
$$
F(v) (F(u)(a)) = F(v) (F(u)(b)) \quad \mbox{dans} \quad F(V) \, .
$$
Alors la famille des compos\'es
$$
u \circ v : V \xrightarrow{ \ v \ } U \xrightarrow{ \ u \ } X
$$
d'un \'el\'ement $(U \xrightarrow{ \ u \ } X) \in S$ et d'un \'el\'ement $(V \xrightarrow{ \ v \ } U) \in S_u$ est un crible $J$-couvrant de $X$.

\smallskip

Cela montre que les sections $a,b \in F(X)$ co{\"\i}ncident localement, autrement dit ont la m\^eme image $\overline a = \overline b$ dans $F_s (X)$.

\medskip

\item En effet, si $F$ est un pr\'efaisceau et $G$ est un pr\'efaisceau s\'epar\'e pour la topologie $J$, tout morphisme $F \to G$ induit un carr\'e commutatif
$$
\xymatrix{
F \ar[d] \ar[r] &G \ar[d] \\
F_s \ar[r] &G_s
}
$$
dont la seconde fl\`eche verticale $G \to G_s$ est l'identit\'e de $G$, ce qui signifie qu'il se factorise canoniquement en
$$
F \longrightarrow F_s \longrightarrow G \, .
$$

\item Le foncteur $F \mapsto F_s$ respecte les colimites arbitraires car c'est un adjoint \`a gauche. 
\end{listeisansmarge}
\end{demo}


Puis on a: 

\begin{prop}\label{propII53}

Soit $({\mathcal C},J)$ un site.

\begin{listeimarge}

\item La formule
$$
F^+ (X) = \varinjlim_{S \in J(X)} \ \varprojlim_{(U \xrightarrow{u} X) \in S} F(U)
$$
appliqu\'ee \`a tout pr\'efaisceau $F$ sur ${\mathcal C}$ et \`a tout objet $X$ de ${\mathcal C}$ d\'efinit un foncteur
$$
\begin{matrix}
F &\longrightarrow &F^+ \, , \\
\widehat{\mathcal C} &\longrightarrow &\widehat{\mathcal C} \hfill
\end{matrix}
$$
muni d'une transformation naturelle
$$
{\rm id} \longrightarrow (F \mapsto F^+) \, .
$$

\item Pour tout pr\'efaisceau $F$, $F^+$ est un pr\'efaisceau s\'epar\'e pour la topologie $J$.

\medskip

\item Pour tout pr\'efaisceau s\'epar\'e $F$, $F^+$ est un faisceau pour la topologie $J$.

\medskip

\item Pour tout faisceau $F$, le morphisme canonique
$$
F \longrightarrow F^+
$$
est un isomorphisme.

\medskip

\item Le foncteur
$$
F \longmapsto F^+
$$
respecte les limites finies.
\end{listeimarge}
\end{prop}

\medskip

\begin{demosansqed}
\begin{listeisansmarge}
\item[(i)] Tout morphisme $X \xrightarrow{ \ f \ } Y$ de ${\mathcal C}$ induit une application
$$
F^+ (Y) \longmapsto F^+ (X)
$$
car l'image r\'eciproque de tout crible $J$-couvrant de $Y$ est un crible $J$-couvrant de $X$.

\smallskip

Donc $F^+$ est bien d\'efini comme pr\'efaisceau sur ${\mathcal C}$.

\smallskip

Pour tout objet $X$ de ${\mathcal C}$, la formule de d\'efinition de $F^+ (X)$ est fonctorielle en $F$, et donc $F \mapsto F^+$ d\'efinit un foncteur
$$
\widehat{\mathcal C} \longrightarrow \widehat{\mathcal C} \, .
$$
Pour tout $S \in J(X)$, on a une application canonique
$$
F(X) \longrightarrow \varprojlim_{(U \xrightarrow{u} X) \in S} F(U) \, .
$$
D'o\`u une application naturelle $F(X) \to F^+ (X)$.

\medskip

\item[(iv)] Si $F$ est un faisceau et $X$ un objet de ${\mathcal C}$, les applications
$$
F(X) \longrightarrow \varprojlim_{(U \xrightarrow{u} X) \in S} F(U) \, , \quad S \in J(X) \, ,
$$
sont des bijections.

\smallskip

Par cons\'equent, chaque
$$
F(X) \longrightarrow F^+(X)
$$
est un bijection et
$$
F \longrightarrow F^+
$$
est un isomorphisme.

\medskip

\item[(ii)] Soient $F$ un pr\'efaisceau quelconque, $X$ un objet de ${\mathcal C}$, $S$ un crible $J$-couvrant de $X$ et
$$
a^+ , b^+ \in F^+ (X)$$
deux \'el\'ements qui ont m\^eme image dans $\underset{(U \xrightarrow{u} X)}{\varprojlim} F^+(U)$.

\smallskip

Il existe un crible $J$-couvrant $S'$ de $X$ tel que $a^+ , b^+$ soient les images de deux \'el\'ements
$$
a,b \in \varprojlim_{(U \xrightarrow{u} X) \in S'} F(U) \, .
$$
Puis il existe un crible $J$-couvrant $S''$ de $X$ contenu \`a la fois dns $S$ et dans $S'$ et tel que $a,b$ aient la m\^eme image dans
$$
\varprojlim_{(U \xrightarrow{u} X) \in S''} F(U) \, .
$$
Cela montre que $F$ est s\'epar\'e pour la topologie $J$.

\medskip

\item[(iii)] Soit $F$ un pr\'efaisceau s\'epar\'e pour la topologie $J$.

\smallskip

D'apr\`es (ii) qui est d\'ej\`a d\'emontr\'e, on sait que $F^+$ est s\'epar\'e.

\smallskip

Il reste \`a prouver que pour tout crible $J$-couvrant $S$ d'un objet $X$ de ${\mathcal C}$, l'application
$$
F^+ (X) \longrightarrow \varprojlim_{(U \xrightarrow{u} X) \in S} F^+ (U)
$$
est surjective.

\smallskip

Consid\'erons donc un \'el\'ement
$$
a^+ \in \varprojlim_{(U \xrightarrow{u} X) \in S} F^+ (U) \, .
$$

Pour tout \'el\'ement $U \xrightarrow{ \ u \ } X$ de $S$, il existe un crible couvrant $S_u$ de $U$ tel que l'image $a_u^+$ de $a^+$ dans $F^+ (U)$ provienne d'un \'el\'ement
$$
a_u \in \varprojlim_{(V \xrightarrow{v} U) \in S_u} F(V) \, .
$$
Soit $S'$ le crible de $X$ constitu\'e de tous les compos\'es d'une fl\`eche $U \xrightarrow{ \ u \ } X$ de $S$ et d'une fl\`eche $V \xrightarrow{ \ v \ } U$ de $S_u$. Pour tout \'el\'ement $U \xrightarrow{ \ u \ } X$ de $S$, le crible $u^* S'$ contient $S_u$ et donc le crible $S'$ de $X$ est $J$-couvrant.

\smallskip

Pour tout diagramme commutatif de ${\mathcal C}$
$$
\xymatrix{
V_1 \ar[dd]_w \ar[rr]^{v_1} &&U_1 \ar[rd]^{u_1} \\
&&&X \\
V_2 \ar[rr]^{v_2} &&U_2 \ar[ru]^{u_2}
}
$$
avec $u_1 , u_2 \in S$, $v_1 \in S_{u_1}$, $v_2 \in S_{u_2}$, les \'el\'ements
$$
a^+_{u_1} \in F^+ (U_1) \quad \mbox{et} \quad a^+_{u_2} \in F^+ (U_2)
$$
ont m\^eme image dans $F^+ (V_1)$.

\smallskip

Comme $F$ est s\'epar\'e, on en d\'eduit que l'image de $a_{u_2}$ dans $F(V_2)$ est envoy\'ee par $F(w)$ sur l'image de $a_{u_1}$ dans $F(V_1)$.

\smallskip

Donc les $a_u$, $(U \xrightarrow{ \ u \ } X) \in S$, d\'efinissent un \'el\'ement
$$
a \in \varprojlim_{(V \rightarrow X) \in S'} F(V)
$$
dont l'image dans $\underset{(U \rightarrow X) \in S}{\varprojlim} F^+ (U)$ est \'egale \`a $a^+$.

\smallskip

Cela montre que $F$ est un $J$-faisceau.

\medskip

\item[(v)] Pour tout objet $X$ de ${\mathcal C}$, le foncteur
$$
F \longmapsto \varinjlim_{S \in J(X)} \ \varprojlim_{(U \xrightarrow{u} X)\in S} F(U)
$$
respecte les limites finies. En effet, pour tout $S \in J(X)$, le foncteur de limite
$$
\varprojlim_{(U \xrightarrow{u} X)}
$$
respecte les limites arbitraires et, d'autre part, le foncteur de colimite sur l'ensemble ordonn\'e filtrant $S$
$$
\varinjlim_{S \in J(X)}
$$
respecte les limites finies d'ensembles, comme il r\'esulte du lemme suivant qui ach\`eve la preuve de la proposition:
\end{listeisansmarge}
\end{demosansqed}

\begin{lem}\label{lemII54}

Soit ${\mathcal D}$ une petite cat\'egorie filtrante au sens que

\medskip

$\left\{ \begin{matrix}
\bullet &\mbox{pour tous objets $d_1 , d_2$ de ${\mathcal D}$, il existe dans ${\mathcal D}$ au moins deux morphismes vers un m\^eme objet $d$} \hfill \\
&\begin{matrix}
d_1 { \ }_{\displaystyle\searrow} \\
&\!\!\!\!d \\
d_2 { \ }^{\displaystyle\nearrow}
\end{matrix} \\
\bullet & \mbox{pour toute paire de fl\`eches $\xymatrix{d \dar[r]^{^{\mbox{\footnotesize$\alpha$}}}_{{\mbox{\footnotesize$\beta$}}} &d'}$ de ${\mathcal D}$, il existe dans ${\mathcal D}$ au moins une fl\`eche $d' \xrightarrow{ \ \gamma \ } d''$} \hfill \\
&\mbox{telle que} \hfill \\
&\gamma \circ \alpha = \gamma \circ \beta \, . 
\end{matrix} \right.
$

\medskip

\noindent Alors le foncteur de colimite
$$
\varinjlim_{\mathcal D} : [{\mathcal D} , {\rm Ens}] \longrightarrow {\rm Ens}
$$
respecte les limites finies.
\end{lem}


\begin{remark}

Le lemme s'applique en particulier lorsque la cat\'egorie ${\mathcal D}$ est un ensemble ordonn\'e et filtrant au sens que, pour tous \'el\'ements $d_1 , d_2 \in {\mathcal D}$, il existe $d \in {\mathcal D}$ v\'erifiant
$$
d_1 \leq d \quad \mbox{et} \quad d_2 \leq d \, .
$$

\end{remark}

\begin{demo}

Pour tout objet $X_{\bullet}$ de $[{\mathcal D} , {\rm Ens}]$
$$
\begin{matrix}
{\mathcal D} &\longrightarrow &{\rm Ens} \, , \\
\hfill d &\longmapsto &X_d \, , \hfill
\end{matrix}
$$
tout \'el\'ement de l'ensemble
$$
\varinjlim_{\mathcal D} \, X_{\bullet}
$$
est repr\'esent\'e par au moins un \'el\'ement $x_d \in X_d$ pour un indice $d \in {\rm Ob} (D)$. Deux \'el\'ements $x_d \in X_d$, $x_{d'} \in X_{d'}$ ont la m\^eme image dans l'ensemble colimite si et seulement si il existe deux morphismes de ${\mathcal D}$
$$
\xymatrix{
d \ar[rd]^{\alpha} \\
&d'' \\
d' \ar[ru]_{\beta}
}
$$
qui transforment $x_d$ et $x_{d'}$ en le m\^eme \'el\'ement de $X_{d''}$.

\smallskip

On en d\'eduit facilement que le foncteur $\underset{\mathcal D}{\varinjlim}$ respecte les produits finis et les produits fibr\'es, donc respecte les limites finies arbitraires. 
\end{demo}


On peut maintenant montrer:

\begin{thm}\label{thmII55}

Soit $({\mathcal C} , J)$ un site.

\begin{listeimarge}

\item Le foncteur de plongement des faisceaux dans les pr\'efaisceaux
$$
j_* : \widehat{\mathcal C}_J \xhookrightarrow{ \ { \ } \ } \widehat{\mathcal C}
$$
admet un adjoint \`a gauche
$$
j^* : \widehat{\mathcal C} \longrightarrow \widehat{\mathcal C}_J
$$
appel\'e le foncteur de faisceautisation.

\medskip

\item Ce foncteur s'\'ecrit comme le compos\'e
$$
F \longmapsto F^{++} = (F^+)^+
$$
du foncteur $F \mapsto F^+$ avec lui-m\^eme, ou encore comme le compos\'e
$$
F \longmapsto F_s^+ = (F_s)^+
$$
du foncteur $F \mapsto F_s$ suivi du foncteur $F \mapsto F^+$.

\medskip

\item Pour tout pr\'efaisceau $F$, son faisceautis\'e $j^* F$ est donc donn\'e par la formule en tout objet $X$ de ${\mathcal C}$
$$
j^* F(X) = \varinjlim_{S \in J (X)} \, \varprojlim_{(U \to X) \in S} F_s (U) \, .
$$

\item Le foncteur $j^*$ respecte les colimites arbitraires et les limites finies.
\end{listeimarge}
\end{thm}

\begin{remarks}
\begin{listeisansmarge}
\item Pour tout pr\'efaisceau $F$ et tout objet $X$ de ${\mathcal C}$, l'application $F(X) \to F_s (X)$ est surjective et l'application $F_s (X) \to (F_s)^+ (X) = j_* j^* F(X)$ est injective. Donc $F \to F_s \to j_*j^* F$ est la factorisation du morphisme canonique $F \to j_* j^* F$ en un \'epimorphisme suivi d'un monomorphisme.

\medskip

\item La transformation naturelle d'adjonction
$$
{\rm id} \longrightarrow j_* \circ j^*
$$
est un isomorphisme de foncteurs $\widehat{\mathcal C}_J \to \widehat{\mathcal C}_J$. D'apr\`es le lemme \ref{lemI82}, cette propri\'et\'e \'equivaut formellement \`a la pleine fid\'elit\'e du foncteur $j_* : \widehat{\mathcal C}_J \hookrightarrow \widehat{\mathcal C}$.
\end{listeisansmarge}
\end{remarks}

\begin{demo}

D'apr\`es la proposition \ref{propII53}, le foncteur
$$
F \longmapsto F^{++} = (F^+)^+
$$
transforme les pr\'efaisceaux en $J$-faisceaux. Il est muni d'une transformation naturelle
$$
{\rm id} \longrightarrow (F \longmapsto F^{++})
$$
telle que, si $F$ est un $J$-faisceau, le morphisme canonique
$$
F \longrightarrow F^{++}
$$
est un isomorphisme.

\smallskip

Ainsi, tout morphisme de pr\'efaisceaux $F \to G$ induit un carr\'e commutatif
$$
\xymatrix{
F \ar[d] \ar[r] &G \ar[d] \\
F^{++} \ar[r] &G^{++}
}
$$
qui, si $G$ est un $J$-faisceau, se lit comme une factorisation canonique
$$
F \longrightarrow F^{++} \longrightarrow G
$$
de $F \to G$.

\smallskip

Cela prouve (i) et la premi\`ere partie de (ii).

\smallskip

De m\^eme, d'apr\`es le lemme \ref{lemII52} et la proposition \ref{propII53}, le foncteur
$$
F \longmapsto F_s^+ = (F_s)^+
$$
transforme les pr\'efaisceaux en faisceaux. Il est muni d'une transformation naturelle
$$
{\rm id} \longrightarrow (F \longmapsto F_s^+)
$$
telle que, si $F$ est un faisceau, le morphisme canonique
$$
F \longrightarrow F_s^+
$$
est un isomorphisme.

\smallskip

On en d\'eduit que tout morphisme $F \to G$ d'un pr\'efaisceau $F$ dans un faisceau $G$ se factorise canoniquement en
$$
F \longrightarrow F_s^+ \longrightarrow G
$$
et donc que $F \mapsto F_s^+$ est un adjoint \`a gauche de
$$
j_* : \widehat{\mathcal C}_J \xhookrightarrow{ \ { \ } \ } \widehat{\mathcal C} \, .
$$
Par unicit\'e des adjoints, les deux foncteurs
$$
\begin{matrix}
&F &\longmapsto &F^{++} \\
\mbox{et} &F &\longmapsto &F_s^+
\end{matrix}
$$
sont canoniquement isomorphes.

\smallskip

Cela prouve la seconde partie de (ii).

\smallskip

(iii) en r\'esulte.

\smallskip

D'apr\`es la proposition \ref{propII53} (v), le foncteur
$$
F \longmapsto F^+
$$
respecte les limites finies, donc aussi le foncteur
$$
j^* : F \longmapsto F^{++} \, .
$$
Cela prouve (iv).
\end{demo}


On a le lemme utile suivant:

\begin{lem}\label{lemII56}

Soient $({\mathcal C},J)$ un site et $j^* : \widehat{\mathcal C} \to \widehat{\mathcal C}_J$ le foncteur de faisceautisation associ\'e.

\smallskip

Soit $P \to F$ un morphisme de pr\'efaisceaux.

\smallskip

Alors:

\begin{listeimarge}

\item Le morphisme induit $j^* P \to j^* F$ est un \'epimorphisme de faisceaux si et seulement si, pour tout objet $X$ de ${\mathcal C}$ et tout \'el\'ement $x \in F(X)$, il existe une famille $J$-couvrante de morphismes de ${\mathcal C}$
$$
X_i \longrightarrow X \, , \qquad i \in I \, ,
$$
telle que la restriction de $x$ \`a chaque $F(X_i)$ est dans l'image de $P(X_i)$.

\medskip

\item Le morphisme induit $j^* P \to j^* F$ est un monomorphisme de faisceaux si et seulement si, pour tout objet $X$ de ${\mathcal C}$, deux \'el\'ements $x_1 , x_2 \in P(X)$ qui ont m\^eme image dans $F(X)$ co{\"\i}ncident localement pour la topologie~$J$.
\end{listeimarge}
\end{lem}

\begin{remarks}
\begin{listeisansmarge}
\item On verra au paragraphe \ref{subsec331} qu'un morphisme de faisceaux est un isomorphisme si (et seulement si) il est \`a la fois un \'epimorphisme et un monomorphisme. 

\smallskip

Donc un morphisme de pr\'efaisceaux $P \to F$ induit un isomorphisme de faisceaux $j^* P \to j^* F$ si et seulement si il satisfait les conditions de (i) et (ii).

\medskip

\item Si $P$ est un faisceau ou plus g\'en\'eralement un pr\'efaisceau s\'epar\'e, (ii) dit que $j^* P \to j^* F$ est un monomorphisme si et seulement si l'application $P(X) \to F(X)$ est injective pour tout objet $X$ de ${\mathcal C}$ c'est-\`a-dire si $P \to F$ est un monomorphisme de $\widehat{\mathcal C}$.

\medskip

\item Le crit\`ere de (i) vaut en particulier quand $P$ et $F$ sont des faisceaux: un morphisme $P \to F$ de $\widehat{\mathcal C}_J$ est un \'epimorphisme si et seulement si toute section de $F$ sur un objet $X$ de ${\mathcal C}$ se rel\`eve localement en des sections de $P$.
\end{listeisansmarge}
\end{remarks}

\begin{demo}
\begin{listeisansmarge}
\item Formons la somme amalgam\'ee $F \amalg_P F$ de $F$ et $F$ sous $P$ dans $\widehat{\mathcal C}$. C'est le pr\'efaisceau qui associe \`a tout objet $X$ de ${\mathcal C}$ la somme amalgam\'ee des ensembles $F(X)$ et $F(X)$ sous $P(X)$. Le foncteur $j^*$ respecte les colimites, donc $j^* P \to j^* F$ est un \'epimorphisme si et seulement si
$$
j^* (F \amalg_P F) \longrightarrow j^* F
$$
est un isomorphisme.

\smallskip

Si cette condition est v\'erifi\'ee,
$$
(F \amalg_P F)_s \longrightarrow F_s
$$
est un monomorphisme de pr\'efaisceaux, ce qui \'equivaut \`a la condition de (i).

\smallskip

R\'eciproquement, si cette condition est v\'erifi\'ee et $(F \amalg_P F)_s \to F_s$ est un monomorphisme, alors
$$
j^* (F \amalg_P F) \longrightarrow j^* F
$$
est aussi un monomorphisme, donc un isomorphisme puisqu'il admet une section $j^* F \to j^* (F \amalg_P F)$.

\medskip

\item Formons le produit fibr\'e $P \times_F P$ de $P$ et $P$ sur $F$ dans $\widehat{\mathcal C}$. C'est le pr\'efaisceau qui associe \`a tout objet $X$ de ${\mathcal C}$ le produit fibr\'e d'ensembles $P(X) \times_{F(X)} P(X)$.

\smallskip

Le foncteur $j^*$ respecte les limites finies, donc $j^* P \to j^* F$ est un monomorphisme si et seulement si le morphisme diagonal
$$
j^* (P) \longrightarrow j^* (P \times_F P)
$$
est un isomorphisme. Comme ses compos\'es avec les projections $j^* (P \times_F P) \rightrightarrows j^* (P)$ sont l'identit\'e de $j^* (P)$, cela revient \`a demander qu'il soit un \'epimorphisme.

\smallskip

On conclut en appliquant le crit\`ere de (i). 
\end{listeisansmarge}
\end{demo}

\subsection{Le foncteur canonique}\label{subsec252}

\medskip

Toute cat\'egorie essentiellement petite ${\mathcal C}$ s'envoie canoniquement dans la cat\'egorie $\widehat{\mathcal C}_J$ des faisceaux sur ${\mathcal C}$ pour tout choix d'une topologie $J$:

\begin{defn}\label{defII57}

Soit $({\mathcal C} , J)$ un site.

\smallskip

On appelle foncteur canonique et on note
$$
\ell : {\mathcal C} \longrightarrow \widehat{\mathcal C}_J
$$
le compos\'e du foncteur de Yoneda
$$
\begin{matrix}
y : {\mathcal C} &\longrightarrow &\widehat{\mathcal C} \, , \hfill \\
\hfill S &\longmapsto &{\rm Hom} (\bullet , S)
\end{matrix}
$$
et du foncteur de faisceautisation
$$
j^* : \widehat{\mathcal C} \longrightarrow \widehat{\mathcal C}_J \, .
$$
\end{defn}

\begin{remarksqed}
\begin{listeisansmarge}
\item Pour tout site $({\mathcal C},J)$, le foncteur canonique
$$
\ell : {\mathcal C} \xrightarrow{ \ y \ } \widehat{\mathcal C} \xrightarrow{ \ j^* \ } \widehat{\mathcal C}_J
$$
respecte les limites finies.

\smallskip

En effet, le foncteur de Yoneda $y$ respecte les limites arbitraires et, d'apr\`es le th\'eor\`eme \ref{thmII55} (iv), le foncteur de faisceautisation $j^*$ respecte les limites finies.

\medskip

\item En revanche, le foncteur $\ell : {\mathcal C} \to \widehat{\mathcal C}_J$ ne respecte pas en g\'en\'eral les colimites, bien que le foncteur $j^* : \widehat{\mathcal C} \to \widehat{\mathcal C}_J$ respecte les colimites arbitraires. Cela r\'esulte de ce que le foncteur de Yoneda $y$ ne respecte pas en g\'en\'eral les colimites.

\medskip

\item Plus g\'en\'eralement, tout foncteur
$$
\rho : {\mathcal C} \longrightarrow {\mathcal D}
$$
d'une cat\'egorie essentiellement petite ${\mathcal C}$ munie d'une topologie $J$ vers une cat\'egorie localement petite ${\mathcal D}$ induit un foncteur
$$
\ell_{\rho} : {\mathcal D} \longrightarrow \widehat{\mathcal C}_J \, .
$$
Il est d\'efini comme le compos\'e du foncteur de Yoneda
$$
y : {\mathcal D} \longrightarrow \widehat{\mathcal D} \, ,
$$
du foncteur de composition avec $\rho$
$$
\begin{matrix}
\rho^* : \widehat{\mathcal D} &\longrightarrow &\widehat{\mathcal C} \, , \hfill \\
\hfill F &\longmapsto &F \circ \rho
\end{matrix}
$$
et du foncteur de faisceautisation
$$
j^* : \widehat{\mathcal C} \longrightarrow \widehat{\mathcal C}_J \, .
$$
Le foncteur $\ell_{\rho} = j^* \circ \rho^* \circ y$ respecte les limites finies puisqu'il en est ainsi de $\rho^*$ comme de $y$ et $j^*$. En revanche il ne respecte pas les colimites en g\'en\'eral.

\medskip

\item En particulier, si ${\mathcal C}$ est une cat\'egorie localement petite, tout choix d'une sous-cat\'egorie essentiellement petite ${\mathcal C}'$ de ${\mathcal C}$ et d'une topologie $J$ de ${\mathcal C}'$ d\'efinit un foncteur canonique
$$
{\mathcal C} \longrightarrow \widehat{\mathcal C}'_J \, .
$$
\end{listeisansmarge}
\end{remarksqed}

\subsection{Topologies sous-canoniques et topologie canonique}\label{subsec253}

\medskip

Tout choix d'une topologie $J$ sur une cat\'egorie essentiellement petite ${\mathcal C}$ d\'efinit un foncteur
$$
\ell : {\mathcal C} \longrightarrow \widehat{\mathcal C}_J
$$
qui envoie ${\mathcal C}$ dans une cat\'egorie de faisceaux.

\smallskip

Particuli\`erement int\'eressants sont les cas o\`u ce foncteur est pleinement fid\`ele, c'est-\`a-dire permet de voir ${\mathcal C}$ comme une sous-cat\'egorie pleine de la cat\'egorie de faisceaux $\widehat{\mathcal C}_J$.

\smallskip

Voici une double caract\'erisation des topologies $J$ qui poss\`edent cette propri\'et\'e:

\begin{prop}\label{propII58}

Soit $({\mathcal C} , J)$ un site.

\smallskip

Alors les propri\'et\'es suivantes sont \'equivalentes:

\begin{enumerate}[label=(\arabic*)]

\item Le foncteur canonique
$$
\ell : {\mathcal C} \xrightarrow{ \ y \ } \widehat{\mathcal C} \xrightarrow{ \ j^* \ } \widehat{\mathcal C}_J
$$
est pleinement fid\`ele.

\medskip

\item Pour tout objet $X$ de ${\mathcal C}$, le pr\'efaisceau
$$
y(X) = {\rm Hom} (\bullet , X)
$$
est un faisceau pour la topologie $J$.

\medskip

\item Pour tout crible $J$-couvrant $S$ d'un objet $X$ de ${\mathcal C}$ et pour tout morphisme $f : X' \to X$ de ${\mathcal C}$, on a
$$
X' = \varinjlim_{(U' \to X') \in f^* S} U'
$$
o\`u la colimite est prise sur la sous-cat\'egorie pleine de ${\mathcal C} / X'$ dont les objets sont les \'el\'ements du crible $f^* S$.
\end{enumerate}
\end{prop}

\begin{demosansqed}

L'implication $(2) \Rightarrow (1)$ r\'esulte de ce que

\medskip

$
\left\{ \begin{matrix}
\bullet &\mbox{d'apr\`es le lemme de Yoneda, le foncteur} \hfill \\
&y : {\mathcal C} \longrightarrow \widehat{\mathcal C} \\
&\mbox{est pleinement fid\`ele,} \hfill \\
\bullet &\mbox{le foncteur de faisceautisation $j^*$ laisse invariant tout pr\'efaisceau qui est un faisceau } \hfill \\
&\mbox{pour la topologie $J$,} \hfill \\
\bullet &\mbox{par d\'efinition, $\widehat{\mathcal C}_J$ est une sous-cat\'egorie pleine de $\widehat{\mathcal C}$.} \hfill
\end{matrix} \right.
$

\medskip

Pour montrer l'implication $(3) \Rightarrow (2)$, consid\'erons deux objets $X$ et $Y$ de ${\mathcal C}$ ainsi qu'un crible $J$-couvrant $S$ de $X$. Sous l'hypoth\`ese (3), on a
$$
{\rm Hom} (X,Y) = \varprojlim_{(U \to X) \in S} {\rm Hom} (U,Y)
$$
c'est-\`a-dire
$$
y(Y)(X) = \varprojlim_{(U \to X) \in S} y(Y)(U) \, .
$$
Cela signifie que le pr\'efaisceau $y(Y)$ est un faisceau pour la topologie $J$.

\smallskip

Comme nous allons le voir, l'implication $(1) \Rightarrow (3)$ r\'esulte du lemme suivant:
\end{demosansqed}

\begin{lem}\label{lemII59}

Soient $({\mathcal C},J)$ un site et $\ell : {\mathcal C} \to \widehat{\mathcal C}_J$ le foncteur canonique associ\'e.

\begin{listeimarge}

\item Pour tout objet $X$ de ${\mathcal C}$ et tout $J$-faisceau $F$ sur ${\mathcal C}$, l'ensemble des sections
$$
F(X)
$$
s'identifie \`a l'ensemble de morphismes de $\widehat{\mathcal C}_J$
$$
{\rm Hom} (\ell (X),F) \, .
$$

\item Pour tout crible $J$-couvrant $S$ d'un objet $X$ de ${\mathcal C}$, on a dans la cat\'egorie $\widehat{\mathcal C}_J$ la formule
$$
\ell(X) = \varinjlim_{(U \to X) \in S} \ell (U) \, .
$$
\end{listeimarge}
\end{lem}

\begin{demolem}
\begin{listeisansmarge}
\item On sait d\'ej\`a d'apr\`es le lemme \ref{lemI72} (ii) que l'ensemble $F(X)$ s'identifie \`a l'ensemble de morphismes de $\widehat{\mathcal C}$
$$
{\rm Hom} (y(X) , j_* F) \, .
$$
Or, celui-ci s'identifie \`a l'ensemble de morphismes de $\widehat{\mathcal C}_J$
$$
{\rm Hom} (j^* \circ y(X),F)
$$
puisque le foncteur de faisceautisation
$$
j^* : \widehat{\mathcal C} \longrightarrow \widehat{\mathcal C}_J
$$
est adjoint \`a gauche du foncteur de plongement $j_* : \widehat{\mathcal C}_J \hookrightarrow \widehat{\mathcal C}$.

\medskip

\item En effet, pour tout objet $F$ de $\widehat{\mathcal C}_J$, la formule qui caract\'erise la propri\'et\'e de $F$ d'\^etre un faisceau
$$
F(X) = \varprojlim_{(U \to X) \in S} F(U)
$$
s'\'ecrit encore d'apr\`es (i)
$$
{\rm Hom} (\ell (X),F) = \varprojlim_{(U \to X) \in S} {\rm Hom} (\ell (U),F) \, .
$$
Le fait qu'elle soit v\'erifi\'ee par tout objet $F$ de $\widehat{\mathcal C}_J$ signifie exactement
$$
\ell (X) = \varinjlim_{(U \to X) \in S} \ell (U) \, .
$$
\end{listeisansmarge}
\end{demolem}

\noindent {\bf Fin de la d\'emonstration de la proposition \ref{propII58}.:}

\medskip

Il reste \`a montrer l'implication $(1) \Rightarrow (3)$.

\smallskip

Comme l'image r\'eciproque par un morphisme $f : X' \to X$ de ${\mathcal C}$ de tout crible $J$-couvrant de $X$ est un crible $J$-couvrant de $X'$, il suffit de montrer que, pour tout crible $J$-couvrant $S$ d'un objet $X$ de ${\mathcal C}$, on a dans la cat\'egorie ${\mathcal C}$
$$
X = \varinjlim_{(U \to X) \in S} U \, .
$$
D'apr\`es le lemme \ref{lemII59} (ii) ci-dessus, on sait d\'ej\`a que dans la cat\'egorie $\widehat{\mathcal C}_J$ est v\'erifi\'ee la formule
$$
\ell (X) = \varinjlim_{(U \to X) \in S} \ell (U)
$$
et donc aussi, pour tout objet $Y$ de ${\mathcal C}$,
$$
{\rm Hom} (\ell (X) , \ell (Y)) = \varprojlim_{(U \to X)\in S} {\rm Hom} (\ell (U) , \ell (Y)) \, .
$$
Comme le foncteur $\ell$ est pleinement fid\`ele par hypoth\`ese, cela s'\'ecrit encore
$$
{\rm Hom} (X,Y) = \varprojlim_{(U \to X) \in S} {\rm Hom} (U,Y) \, .
$$
Cette identit\'e signifie comme voulu que, dans la cat\'egorie ${\mathcal C}$,
$$
X = \varinjlim_{(U \to X)\in S} U \, .
$$
\hfill $\Box$

On d\'eduit de la proposition:

\begin{cor}\label{corII510}

Soit ${\mathcal C}$ une cat\'egorie essentiellement petite.

\begin{listeimarge}

\item Il existe une topologie $J_c$ de ${\mathcal C}$, appel\'ee la ``topologie canonique'' de ${\mathcal C}$, pour laquelle un crible $S$ d'un objet $X$ de ${\mathcal C}$ est couvrant si et seulement si on a pour tout morphisme $f : X' \to X$ de ${\mathcal C}$ la formule dans la cat\'egorie ${\mathcal C}$
$$
X' = \varinjlim_{(U' \to X') \in f^* S} U' \, .
$$

\item Pour que le foncteur canonique
$$
\ell : {\mathcal C} \longrightarrow \widehat{\mathcal C}_J
$$
associ\'e \`a une topologie $J$ de ${\mathcal C}$ soit pleinement fid\`ele, il faut et il suffit que $J$ soit ``sous-canonique'' au sens que
$$
J \subseteq J_c \, .
$$
\end{listeimarge}
\end{cor}

\begin{remark}

Ainsi, la topologie canonique d'une cat\'egorie essentiellement petite ${\mathcal C}$ est la topologie la plus fine pour laquelle la sous-cat\'egorie pleine de $\widehat{\mathcal C}$ constitu\'ee des faisceaux pour cette topologie contient l'image du foncteur de Yoneda $y = {\mathcal C} \to \widehat{\mathcal C}$.
\end{remark}

\begin{demo}
\begin{listeisansmarge}
\item Il faut v\'erifier que $J_c$ satisfait les trois axiomes de maximalit\'e, stabilit\'e et transitivit\'e qui d\'efinissent la notion de topologie de Grothendieck.

\smallskip

L'axiome de stabilit\'e est v\'erifi\'e par d\'efinition de $J_c$.

\smallskip

L'axiome de maximalit\'e r\'esulte de ce que le crible maximal $S$ d'un objet $X$ engendr\'e par le morphisme ${\rm id} : X \to X$ v\'erifie automatiquement la formule
$$
X = \varinjlim_{(U \to X) \in S} U \, .
$$

Enfin, l'axiome de transitivit\'e est satisfait puisque, si $S$ et $S'$ sont deux cribles d'un objet $X$ tels que
$$
X = \varinjlim_{(U \to X) \in S} U
$$
et
$$
U = \varinjlim_{(U' \to U) \in f^* S'} U' \, , \quad \forall \left( U \xrightarrow{ \ f \ } X \right) \in S \, ,
$$
on a n\'ecessairement
$$
X = \varinjlim_{(U' \to X) \in S'} U' \, .
$$
\end{listeisansmarge}
\end{demo}

\subsection{Exemples de topologies sous-canoniques}\label{subsec254}

\medskip

\noindent $\bullet$ {\bf La topologie discr\`ete d'une cat\'egorie essentiellement petite:}

\smallskip

Il r\'esulte du lemme de Yoneda que la topologie discr\`ete est toujours sous-canonique.

\medskip

\noindent $\bullet$ {\bf La topologie ordinaire des espaces topologiques:}

\smallskip

On a:

\begin{lem}\label{lemII511}

Soit ${\mathcal C}$ une sous-cat\'egorie pleine et essentiellement petite de la cat\'egorie ${\rm Top}$ des espaces topologiques, qui contient tout ouvert de tout objet de ${\mathcal C}$.

\smallskip

Alors la topologie ordinaire de ${\mathcal C}$ (pour laquelle une famille d'applications continues $U_i \to X$ est couvrante quand elle contient une sous-famille constitu\'ee d'immersions ouvertes dont la r\'eunion des images est $X$) est sous-canonique.
\end{lem}

\begin{remarks}
\begin{listeisansmarge}
\item En revanche, la topologie de la densit\'e de ${\mathcal C}$ (pour laquelle une famille d'applications continues $U_i \to X$ est couvrante quand elle contient une sous-famille constitu\'ee d'immersions ouvertes dont la r\'eunion des images est dense dans $X$) n'est pas sous-canonique en g\'en\'eral.

\smallskip

En effet, si $U$ est un ouvert dense d'un espace topologique $X$, une application continue $U \to Y$ vers un espace topologique $Y$ ne se prolonge pas en g\'en\'eral en une application continue $X \to Y$.

\medskip

\item On observe cependant que si $Y$ est un espace topologique s\'epar\'e, toute application continue $X \to Y$ est uniquement d\'etermin\'ee par sa restriction \`a n'importe quel ouvert dense de $X$.

\smallskip

Cela signifie que si $Y$ est un espace topologique s\'epar\'e, le pr\'efaisceau ${\rm Hom} (\bullet , Y)$ est s\'epar\'e pour la topologie de la densit\'e.
\end{listeisansmarge}
\end{remarks}


\begin{demo}

D'apr\`es la proposition \ref{propII58}, c'est le contenu du lemme \ref{lemII44}. 

\end{demo}

\medskip

Ce r\'esultat s'\'etend aux sous-cat\'egories de la cat\'egorie ${\rm Top}_{\rm an}$ des espaces annel\'es qui sont ``g\'eom\'etriques'' au sens de la d\'efinition \ref{defI58}:

\begin{prop}\label{propII512}

Soit ${\mathcal G}$ une sous-cat\'egorie ``g\'eom\'etrique'' de la cat\'egorie ${\rm Top}_{\rm an}$ des espaces annel\'es.

\smallskip

Soit ${\mathcal C}$ une sous-cat\'egorie pleine et essentiellement petite de ${\mathcal G}$ telle que, pour tout objet $(X,{\mathcal O}_X)$ de ${\mathcal G}$, ses ouverts $(U \hookrightarrow X$, ${\mathcal O}_U = i^* {\mathcal O}_X)$ sont encore des objets de ${\mathcal C}$.

\smallskip

Alors la topologie usuelle de ${\mathcal C}$ (pour laquelle un crible sur un objet $(X,{\mathcal O}_X)$ de ${\mathcal C}$ est couvrant s'il contient une famille d'immersions ouvertes $(U_i , {\mathcal O}_{U_i}) \hookrightarrow (X,{\mathcal O}_X)$ dont les images recouvrent $X$ au sens que $X = \underset{i}{\bigcup} \, U_i$) est sous-canonique.
\end{prop}


\begin{demo}

D'apr\`es la proposition \ref{propII58}, c'est le contenu de la proposition \ref{propII45}. 

\end{demo}


\noindent $\bullet$ {\bf La topologie canonique des treillis:}

\smallskip

On a aussi

\begin{lem}\label{lemII513}

Soit $O$ un treillis (au sens de la d\'efinition \ref{defII23}) consid\'er\'e comme une cat\'egorie.

\smallskip

Alors la topologie naturelle de $O$ (pour laquelle une famille de fl\`eches $u_i \leq u$, $i \in I$, est couvrante si $\underset{i \in I}{\bigvee} \, u_i = u$) est la topologie canonique de $O$.
\end{lem}

\medskip

\begin{remark}

En particulier, si $O(X)$ est le treillis des ouverts d'un espace topologique $X$, la topologie de $X$ est la topologie canonique de $O(X)$.

\end{remark}

\begin{demo}

En effet, pour tous \'el\'ements $v$ et $u_i$, $i \in I$, du treillis $O$, on a par d\'efinition du supremum $\underset{i \in I}{\bigvee} \, u_i$ une \'equivalence entre les relations
$$
u_i \leq v \, , \quad \forall \, i \in I \, ,
$$
et
$$
\bigvee_{i \in I} \, u_i \leq v \, .
$$
\end{demo}

\bigskip

\noindent $\bullet$ {\bf La topologie ordinaire d'un espace muni d'une mesure:}

\smallskip

Soit $X$ un ensemble muni d'une mesure $\mu$.

\smallskip

Soit $O$ l'ensemble des parties mesurables de $X$.

\smallskip

Soit $\overline O$ l'ensemble quotient de $O$ par la relation d'\'equivalence pour laquelle deux parties mesurables $Y$ et $Y'$ sont \'equivalentes lorsque
$$
\mu (Y - Y \cap Y') = 0 \quad \mbox{et} \quad \mu (Y' - Y \cap Y') = 0 \, .
$$

Soit $\leq$ la relation d'ordre sur $\overline O$ pour laquelle les classes $\overline Y$ et $\overline Z$ de deux parties mesurables $Y$ et $Z$ satisfont la relation
$$
\overline Y \leq \overline Z
$$
lorsque
$$
\mu (Y - Y \cap Z) = 0 \, .
$$

On observe que toute famille d\'enombrable $(\overline Y_{\!\!n})_{n \in {\mathbb N}}$ d'\'el\'ements de $\overline O$ admet un supremum $\underset{n \in {\mathbb N}}{\bigvee} \, \overline Y_{\!\!n}$ caract\'eris\'e par l'\'equivalence
$$
\bigvee_{n \in {\mathbb N}} \, \overline Y_{\!\!n} \leq \overline Z \Longleftrightarrow \overline Y_{\!\!n} \leq \overline Z \, , \quad \forall \, n \in {\mathbb N} \, .
$$

Consid\'erant l'ensemble ordonn\'e $\overline O$ comme une cat\'egorie, on a:

\begin{prop}\label{propII514}

Dans le contexte ci-dessus, et si $X$ est r\'eunion d\'enombrable de parties mesurables de mesure finie, la topologie canonique de $\overline O$ est celle pour laquelle une famille de fl\`eches
$$
\overline Y_{\!\!i} \leq \overline Z \, , \quad i \in I \, ,
$$
est couvrante si et seulement si elle compte une sous-famille d\'enombrable $(\overline Y_{\!\!i_n})_{n \in {\mathbb N}}$ telle que
$$
\bigvee_{n \in {\mathbb N}} \, \overline Y_{\!\!i_n} = \overline Z \, .
$$
\end{prop}


\begin{demo}

Si une famille de fl\`eches $(\overline Y_{\!\!i} \leq \overline Z)_{i \in I}$ poss\`ede une sous-famille d\'enombrable $(\overline Y_{\!\!i_n})_{n \in {\mathbb N}}$ v\'erifiant
$$
\mu \left(\overline Z - \bigcup_{n \in {\mathbb N}} \overline Y_{\!\!i_n}\right) \, ,
$$
alors pour tout \'el\'ement $\overline X$ de $\overline O$, les relations
$$
\overline Y_{\!\!i} \leq \overline X \, , \quad \forall \, i \in I \, ,
$$
et
$$
\overline Z \leq \overline X
$$
sont \'equivalentes.

\smallskip

Cela prouve que la topologie consid\'er\'ee est sous-canonique.

\smallskip

R\'eciproquement, consid\'erons une famille de fl\`eches
$$
(\overline Y_{\!\!i} \leq \overline Z)_{i \in I}
$$
qui est couvrante pour la topologie canonique de $\overline O$.

\smallskip

Comme $\overline Z$ est le supremum d'une famille d\'enombrable d'\'el\'ements de mesure finie, on peut supposer que
$$
\mu (\overline Z) < + \infty \, .
$$

Soit alors $M \in [0,+\infty[$ l'infimum des mesures
$$
\mu \left( \overline Z - \bigvee_{n \in {\mathbb N}} \overline Y_{\!\!i_n} \right)
$$
associ\'ees \`a tous les choix possibles de sous-familles d\'enombrables $(\overline Y_{\!\!i_n})_{n \in {\mathbb N}}$ de la famille $(\overline Y_{\!\!i})_{i \in I}$.

\smallskip

Comme toute r\'eunion d\'enombrable d'ensembles d\'enombrables est d\'enombrable, il existe une sous-famille d\'enombrable $(\overline Y_{\!\!i_n})_{n \in {\mathbb N}}$ de $(\overline Y_{\!\!i})_{i \in I}$ telle que
$$
\mu \left( \overline Z - \bigvee_{n \in {\mathbb N}} \overline Y_{\!\!i_n} \right) = M \, .
$$

Il en r\'esulte que pour tout \'el\'ement $i \in I$ on a n\'ecessairement
$$
\overline Y_{\!\!i} \leq \bigvee_{n \in {\mathbb N}} \overline Y_{\!\!i_n} \, .
$$

Alors l'\'equivalence des relations
$$
\overline Y_{\!\!i} \leq \overline X \, , \quad \forall \, i \in I \, ,
$$
et
$$
\overline Z \leq \overline X
$$
pour tout \'el\'ement $\overline X$ de $\overline O$ impose
$$
\mu \left( \overline Z - \bigvee_{n \in {\mathbb N}} \overline Y_{\!\!i_n} \right) = 0
$$
c'est-\`a-dire
$$
\overline Z = \bigvee_{n \in {\mathbb N}} \overline Y_{\!\!i_n} \, .
$$

Cela montre que la topologie canonique est contenue dans la topologie consid\'er\'ee, donc qu'elles se confondent. 

\end{demo}

\bigskip

\noindent $\bullet$ {\bf La topologie atomique d'une cat\'egorie codirig\'ee:}

\smallskip

D'apr\`es la d\'efinition \ref{defII24}, la topologie atomique d'une cat\'egorie codirig\'ee essentiellement petite est celle dont les cribles couvrants sont les cribles non vides.

\newpage

Afin de d\'eterminer quand elle est sous-canonique, on a besoin de la d\'efinition suivante:

\begin{defn}\label{defII515}

Dans une cat\'egorie ${\mathcal C}$, un \'epimorphisme [resp. un monomorphisme]
$$
X \xrightarrow{ \ f \ } Y
$$
est dit strict s'il existe une famille de paires de morphismes
$$
\qquad\quad \xymatrix{X_i \dar[r]^-{^{\mbox{\footnotesize$p_i$}}}_-{{\mbox{\footnotesize$q_i$}}} &X \, , \quad i \in I \, ,}
$$
$$
\mbox{[resp.} \quad \xymatrix{Y \dar[r]^-{^{\mbox{\footnotesize$p_i$}}}_-{{\mbox{\footnotesize$q_i$}}} &Y_i \, , \quad i \in I \, ,} \mbox{]}
$$
telle que, pour tout objet $Z$ de ${\mathcal C}$, la composition avec $f$ identifie
$$
{\rm Hom} (Y,Z) \qquad \mbox{[resp.} \quad {\rm Hom} (Z,X) \, \mbox{]}
$$
au sous-ensemble
$$
\qquad\quad \{ g \in {\rm Hom} (X,Z) \mid g \circ p_i = g \circ q_i \, , \quad \forall \, i \}
$$
$$
\mbox{[resp.} \quad \{ g \in {\rm Hom} (Z,Y) \mid p_i \circ g = q_i \circ g \, , \quad \forall i \} \, . \, \mbox{]}
$$

\end{defn}

Disposant de cette d\'efinition, on peut \'enoncer:

\begin{prop}\label{propII516}

La topologie atomique d'une cat\'egorie codirig\'ee essentiellement petite ${\mathcal C}$ est sous-canonique si et seulement si toute fl\`eche de ${\mathcal C}$ est un \'epimorphisme strict.
\end{prop}

\begin{demo}

D'apr\`es la proposition \ref{propII58}, la topologie atomique de ${\mathcal C}$ est sous-canonique si et seulement si, pour tout morphisme de ${\mathcal C}$
$$
X \xrightarrow{ \ f \ } Y \, , 
$$
on a
$$
Y = \varinjlim_{(X' \to Y) \in S_f} X'
$$
o\`u $S_f$ d\'esigne le crible de $Y$ engendr\'e par le morphisme $f$.

\smallskip

C'est \'equivalent \`a demander que, pour tout objet $Z$ de ${\mathcal C}$, l'ensemble
$$
{\rm Hom} (Y,Z)
$$
s'identifie au sous-ensemble de
$$
{\rm Hom} (X,Z)
$$
constitu\'e des morphismes $g : X \to Z$ tels que, pour toute paire de morphismes
$$
\xymatrix{X' \dar[r]^-{^{\mbox{\footnotesize$p$}}}_-{{\mbox{\footnotesize$q$}}} &X}
$$
v\'erifiant $f \circ p = f \circ q$, on ait $g \circ p = g \circ q$.

\smallskip

On passe de cette formulation \`a la d\'efinition des ``\'epimorphismes stricts'' en se souvenant que la cat\'egorie ${\mathcal C}$ est essentiellement petite. 

\end{demo}

\newpage

\noindent $\bullet$ {\bf Les sites g\'eom\'etriques:}

\smallskip

Consid\'erons une cat\'egorie essentiellement petite ${\mathcal C}$, une classe ${\mathcal M}$ de morphismes de ${\mathcal C}$ qui est ``g\'eom\'etrique'' au sens de la d\'efinition \ref{defII31} (i), et une propri\'et\'e (R) des familles de morphismes de ${\mathcal M}$ de m\^eme but qui est une ``notion g\'eom\'etrique de recouvrement'' au sens de la d\'efinition \ref{defII31} (ii).

\smallskip

Suivant la d\'efinition \ref{defII32}, on consid\`ere la topologie $J_R$ de ${\mathcal C}$ pour laquelle un crible sur un objet $X$ de ${\mathcal C}$ est couvrant s'il contient une sous-famille constitu\'ee d'\'el\'ements $U_i \to X$ de ${\mathcal C}$ et qui poss\`ede la propri\'et\'e (R).

\smallskip

On d\'eduit de la proposition \ref{propII58}:

\begin{cor}\label{corII517}

Dans les conditions ci-dessus, la topologie $J_R$ est sous-canonique si et seulement si pour tous objets $X,Y$ de ${\mathcal C}$, et toute famille de morphismes de ${\mathcal M}$
$$
U_i \longrightarrow X
$$
qui poss\`ede la propri\'et\'e {\rm (R)}, se donner un morphisme
$$
f : X \longrightarrow Y
$$
\'equivaut \`a se donner une famille de morphismes
$$
f_i : U_i \longrightarrow Y
$$
telle que, pour tous indices $i,j$, les deux morphismes d\'eduits de $f_i$ et $f_j$ par composition avec les deux projections de $U_i \times_X U_j$ sur $U_i$ et $U_j$
$$
U_i \times_X U_j \rightrightarrows Y
$$
sont \'egaux.
\end{cor}

\bigskip

\begin{remarks}
\begin{listeisansmarge}
\item Si ${\mathcal C}$ est une sous-cat\'egorie pleine et essentiellement petite de ${\rm Top}$ qui contient toutes les immersions ouvertes vers les objets de ${\mathcal C}$, que ${\mathcal M}$ est la classe g\'eom\'etrique des immersions ouvertes et que (R) est la propri\'et\'e d'\^etre globalement surjective, alors $J_R$ est la topologie ordinaire des espaces topologiques.

\smallskip

Elle est sous-canonique d'apr\`es le lemme \ref{lemII511}.

\medskip

\item Si ${\mathcal G}$ est une sous-cat\'egorie ``g\'eom\'etrique'' de la cat\'egorie ${\rm Top}_{\rm an}$ des espaces annel\'es, que ${\mathcal C}$ est une sous-cat\'egorie pleine et essentiellement petite de ${\mathcal G}$ qui contient toutes les immersions ouvertes vers les objets de ${\mathcal C}$, que ${\mathcal M}$ est la classe g\'eom\'etrique des immersions ouvertes et que (R) est la propri\'et\'e d'\^etre globalement surjective, alors $J_R$ est la topologie ordinaire des espaces annel\'es.

\smallskip

Elle est sous-canonique d'apr\`es la proposition \ref{propII512}.

\medskip

\item Si ${\mathcal G}$ est la cat\'egorie g\'eom\'etrique et essentiellement petite des vari\'et\'es diff\'erentielles de classe $C^k$ [resp. des vari\'et\'es analytiques] qui sont d\'enombrables \`a l'infini, que ${\mathcal M}$ est la classe g\'eom\'etrique des morphismes submersifs et que (R) est la propri\'et\'e d'\^etre globalement surjective, alors $J_R$ co{\"\i}ncide avec la topologie ordinaire comme on a vu dans la remarque (iii) suivant le corollaire \ref{corII35}.

\smallskip

Donc cette topologie est sous-canonique et, pour toutes vari\'et\'es de classe $C^k$ [resp. analytiques] $X,Y$ et toute famille globalement surjective de submersions
$$
U_i \longrightarrow X \, ,
$$
se donner un morphisme
$$
X \longrightarrow Y
$$
\'equivaut \`a se donner une famille de morphismes
$$
U_i \longrightarrow Y
$$
qui co{\"\i}ncident sur les $U_i \times_X U_j$.
\end{listeisansmarge}
\end{remarks}

\begin{demo}

C'est l'\'equivalence des propri\'et\'es (1) et (2) de la proposition \ref{propII58}. 

\end{demo}
\bigskip

\noindent $\bullet$ {\bf La topologie (fpqc) des sch\'emas:}

\smallskip

Les conditions \'equivalentes du corollaire pr\'ec\'edent sont satisfaites par la topologie (fpqc) des sch\'emas:

\begin{cor}\label{corII518}
\begin{listeimarge}
\item Soient $X,Y$ deux sch\'emas au-dessus d'un sch\'ema de base $S$ et
$$
(U_i \longrightarrow X)_{i \in I}
$$
une famille de morphismes plats de sch\'emas qui est (fpqc) c'est-\`a-dire globalement surjective et quasi-compacte.

\smallskip

Alors se donner un morphisme au-dessus de $S$
$$
X \longrightarrow Y
$$
\'equivaut \`a se donner une famille de morphismes au-dessus de $S$
$$
U_i \longrightarrow Y
$$
telle que, pour tous indices $i,j$, les deux morphismes induits
$$
U_i \times_X U_j \rightrightarrows Y
$$
soient \'egaux.

\medskip

\item Si $S$ est un sch\'ema de base, que ${\mathcal C}$ est une sous-cat\'egorie pleine de ${\rm Sch}/S$ qui est essentiellement petite, stable par limites finies et contient les immersions ouvertes vers des objets de ${\mathcal C}$, alors la topologie (fpqc) sur ${\mathcal C}$ est sous-canonique.

\medskip

\item A fortiori, la topologie (fppf) sur ${\mathcal C}$ et la topologie \'etale sur ${\mathcal C}$ sont sous-canoniques.
\end{listeimarge}
\end{cor}

\begin{demo}
\begin{listeisansmarge}
\item Il n'y a pas de restriction \`a supposer que les sch\'emas $S$, $X,Y$ et $U_i$, $i \in I$, sont affines et que l'ensemble d'indices $I$ est fini.

\smallskip

Ecrivons
$$
\begin{matrix}
\hfill S &= &{\rm Spec} (A) \, , \hfill \\
\hfill Y &= &{\rm Spec} (B) \, , \hfill \\
\hfill X &= &{\rm Spec} (C) \hfill \\
et \quad \hfill U_i &= &{\rm Spec} (D_i) \, , \quad i \in I \, .
\end{matrix}
$$
Ainsi, chaque morphisme
$$
C \longrightarrow D_i
$$
est plat et le morphisme
$$
C \longrightarrow D = \prod_{i \in I} \, D_i
$$
est fid\`element plat.

\smallskip

Il s'agit de prouver que se donner un morphisme au-dessus de $A$
$$
B \longrightarrow C
$$
\'equivaut \`a se donner un morphisme au-dessus de $A$
$$
B \longrightarrow D
$$
dont les compos\'es avec les deux morphismes
$$
D \rightrightarrows D \otimes_C D
$$
sont \'egaux.

\smallskip

Cela r\'esulte de ce que, d'apr\`es le lemme \ref{lemII320} (i), l'anneau $C$ s'identifie \`a l'\'egalisateur
$$
{\rm eg} \left( D \rightrightarrows D \otimes_C D \right) \, .
$$

\item est cons\'equence de (i) d'apr\`es le corollaire \ref{corII517}.

\medskip

\item r\'esulte de (ii) puisque la topologie (fpqc) est plus fine que la topologie (fppf), qui est elle-m\^eme plus fine que la topologie \'etale. 
\end{listeisansmarge}
\end{demo}


%% file: Chapitre3_num.tex







\vglue 15mm

\chapter{D\'efinition et propri\'et\'es cat\'egoriques des topos}\label{chap3}

\section{La notion de topos}\label{sec31}

\subsection{D\'efinition}\label{subsec311}

Grothendieck a propos\'e d'appeler ``topos'' les cat\'egories qu'il est possible de repr\'esenter comme des cat\'egories de faisceaux sur des sites:

\begin{defn}\label{defIII11}

On appelle topos les cat\'egories qui sont \'equivalentes \`a la cat\'egorie $\widehat{\mathcal C}_J$ des faisceaux sur au moins un site $({\mathcal C} , J)$.
\end{defn}

\begin{remarksqed}
\begin{listeisansmarge}
\item Ainsi, un topos est une cat\'egorie d'un genre particulier, sans autre structure que celle de cat\'egorie.

\smallskip

Cela implique que toute propri\'et\'e d'un topos ou toute construction g\'en\'erale dans les topos doit \^etre formul\'ee dans le seul langage des cat\'egories et \^etre invariante par \'equivalence de cat\'egories.

\medskip

\item Nous allons voir un peu plus loin que, en fait, tout topos peut \^etre pr\'esent\'e comme la cat\'egorie des faisceaux non pas seulement sur un site mais sur une infinit\'e de sites diff\'erents. 
\end{listeisansmarge}
\end{remarksqed}

\bigskip

Le but principal de ce chapitre est de montrer que les topos poss\`edent chacun toutes les propri\'et\'es cat\'egoriques ``constructives'' de la cat\'egorie des ensembles.

\smallskip

En effet, ces propri\'et\'es passent de la cat\'egorie ${\rm Ens}$ des ensembles aux cat\'egories de pr\'efaisceaux $\widehat{\mathcal C} = [{\mathcal C}^{\rm op} , {\rm Ens}]$ et de l\`a aux cat\'egories de faisceaux $\widehat{\mathcal C}_J$ par le foncteur de faisceautisation $j^* : \widehat{\mathcal C} \to \widehat{\mathcal C}_J$ ou le foncteur de plongement $j_* : \widehat{\mathcal C}_J \hookrightarrow \widehat{\mathcal C}$. Enfin, elles passent aux topos car elles sont invariantes par les \'equivalences de cat\'egories.

\smallskip

On verra aussi que, r\'eciproquement, toute cat\'egorie qui poss\`ede une partie de ces propri\'et\'es ``constructives'' est n\'ecessairement un topos et donc poss\`ede la liste compl\`ete de ces propri\'et\'es.

\smallskip

En revanche, les topos ne poss\`edent pas en g\'en\'eral les propri\'et\'es ``non constructives'' de la cat\'egorie des ensembles que sont ``l'axiome du choix'' (tout \'epimorphisme $p : X \to Y$ admet une section $s : Y \to X$ telle que $p \circ s = {\rm id}_Y$) et le ``principe du tiers exclu'' (la r\'eunion d'un sous-objet $Z$ d'un objet $X$ et de son compl\'ementaire $X-Z$, d\'efini comme le plus grand sous-objet de $X$ disjoint de $Z$, est toujours \'egale \`a $X$ tout entier).

\smallskip

Voici la premi\`ere propri\'et\'e que les topos h\'eritent de la cat\'egorie des ensembles:

\begin{lem}\label{lemIII12}

On a la propri\'et\'e suivante:

\smallskip

\begin{enumerate}[label=(\arabic*)]
\item[(0)] Tout topos est une cat\'egorie localement petite. Autrement dit, les morphismes $X \to Y$ entre deux objets $X,Y$ d'un topos ${\mathcal E}$ forment un ensemble ${\rm Hom} (X,Y)$.
\end{enumerate}
\end{lem}

\begin{remark}

En revanche, en dehors du topos trivial \`a un seul objet et un seul morphisme (constitu\'e des pr\'efaisceaux sur la cat\'egorie vide), un topos n'est jamais une cat\'egorie petite ou essentiellement petite.
\end{remark}

\bigskip

\begin{demo}

Il suffit de montrer que si ${\mathcal C}$ est une petite cat\'egorie munie d'une topologie $J$, alors la cat\'egorie $\widehat{\mathcal C}_J$ est localement petite.

\smallskip

Comme $\widehat{\mathcal C}_J$ est une sous-cat\'egorie pleine de $\widehat{\mathcal C}$, il suffit de traiter le cas de $\widehat{\mathcal C}$.

\smallskip

Or, pour tous objets $F_1 , F_2$ de $\widehat{\mathcal C}$, ${\rm Hom} (F_1 , F_2)$ est le sous-ensemble de l'ensemble produit
$$
\prod_{X \in {\rm Ob} (X)} {\rm Hom} (F_1 (X) , F_2 (X))
$$
constitu\'e des familles
$$
(u_X : F_1 (X) \longrightarrow F_2 (X))_{X \in {\rm Ob} (X)}
$$
telles que, pour tout morphisme $f : X \to Y$ de ${\mathcal C}$, le carr\'e d'applications entre ensembles
$$
\xymatrix{
F_1 (Y) \ar[d]_{F_1(f)} \ar[r]^{u_Y} &F_2 (Y) \ar[d]^{F_2 (f)} \\
F_1(X) \ar[r]^{u_X} &F_2(X)
}
$$
est commutatif, soit
$$
u_X \circ F_1 (f) = F_2 (f) \circ u_Y \, .
$$

\end{demo}

\subsection{Le ``lemme de comparaison'' de Grothendieck}\label{subsec312}

Nous allons introduire un proc\'ed\'e g\'en\'eral qui permet, \`a partir d'une repr\'esentation d'un topos ${\mathcal E} \cong \widehat{\mathcal C}_J$ comme cat\'egorie des faisceaux sur un site $({\mathcal C} , J)$, d'en d\'eduire une infinit\'e d'autres repr\'esentations par des sites diff\'erents.

\smallskip

Pour cela, on a besoin de la d\'efinition suivante:

\begin{defn}\label{III13}

Soit $({\mathcal C},J)$ un site.

\smallskip

Une sous-cat\'egorie pleine ${\mathcal C}'$ de ${\mathcal C}$ est dite ``dense'' si tout objet $X$ de ${\mathcal C}$ admet une famille $J$-couvrante de morphismes
$$
(U_i \longrightarrow X)_{i \in I}
$$
dont les sources $U_i$ sont des objets de ${\mathcal C}'$.

\smallskip

Dans ce cas, on appelle topologie induite par $J$ dans ${\mathcal C}'$ la topologie $J'$ de ${\mathcal C}'$ pour laquelle une famille de morphismes de ${\mathcal C}'$
$$
(U_i \longrightarrow U)_{i \in I}
$$ 
est $J'$-couvrante si elle est $J$-couvrante.

\end{defn}

\begin{remarkqed}

Tout crible $J$-couvrant d'un objet $U$ de ${\mathcal C}'$ contient un crible $J$-couvrant engendr\'e par des morphismes $U_i \to U$ de ${\mathcal C}'$.

\smallskip

Il en r\'esulte aussit\^ot que $J'$ est bien une topologie de Grothendieck.

\end{remarkqed}

\smallskip

On peut maintenant \'enoncer:

\begin{prop}\label{propIII14}

Soit $({\mathcal C},J)$ un site.

\smallskip

Soit ${\mathcal C}'$ une sous-cat\'egorie pleine de ${\mathcal C}$ qui est dense.

\smallskip

Soit $J'$ la topologie de ${\mathcal C}'$ qui est induite par la topologie $J$ de ${\mathcal C}$.

\smallskip

Alors:

\begin{listeimarge}

\item La restriction \`a ${\mathcal C}'$ de tout $J$-faisceau sur ${\mathcal C}$ est un $J'$-faisceau.

\medskip

\item Le foncteur de restriction ainsi d\'efini
$$
\widehat{\mathcal C}_J \longrightarrow \widehat{\mathcal C}'_{J'}
$$
est une \'equivalence de cat\'egories.

\medskip

\item Une \'equivalence en sens inverse
$$
\widehat{\mathcal C}'_{J'} \longrightarrow \widehat{\mathcal C}_J
$$
est d\'efinie en associant \`a tout faisceau $F$ sur $({\mathcal C}',J')$ le faisceau $\overline F$ sur $({\mathcal C},J)$ donn\'e par la formule
$$
\overline F(X) = \varprojlim_{(U \to X) \in {\mathcal C}' / X} F(U)
$$
si, pour tout objet $X$ de ${\mathcal C}$, ${\mathcal C}'/X$ d\'esigne la cat\'egorie ``relative'' des morphismes $U \to X$ d'un objet $U$ de ${\mathcal C}'$ vers $X$.
\end{listeimarge}
\end{prop}

\begin{demo}

On peut supposer que ${\mathcal C}$ est une petite cat\'egorie.

\begin{listeisansmarge}

\item[(i)] Consid\'erons un $J$-faisceau $F$ sur ${\mathcal C}$, un objet $U$ de ${\mathcal C}'$ et un crible $J'$-couvrant $S$ de $U$ dans ${\mathcal C}'$.

\smallskip

Par d\'efinition de $J'$, le crible $\overline S$ de $U$ dans ${\mathcal C}$ engendr\'e par $S$ est couvrant et on a
$$
F(U) = \varprojlim_{(U' \to U) \in \overline S} F(U') \, .
$$
Cela signifie qu'existe une famille de carr\'es commutatifs de ${\mathcal C}$
$$
\xymatrix{
V_i \ar[d] \ar[r] &U''_i \ar[d] \\
U'_i \ar[r] &U
}
$$
telle que tous les $U'_i \to U$ et $U''_i \to U$ soient \'el\'ements de $S$ et que
$$
F(U) = {\rm eg} \left[ \prod_{(U' \to U) \in S} F(U') \rightrightarrows \prod_i F(V_i) \right].
$$
Chaque $V_i$ poss\`ede une famille $J$-couvrante de morphismes
$$
(V_{i,j} \longrightarrow V_i)_{j \in I_i} 
$$
dont les sources $V_{i,j}$ sont des objets de ${\mathcal C}'$.

\smallskip

On en d\'eduit
$$
F(U) = {\rm eg} \left[ \prod_{(U' \to U) \in S} F(U') \rightrightarrows \prod_i \prod_{j \in I_i} F(V_{i,j}) \right]
$$
et a fortiori
$$
F(U) = \varprojlim_{(U' \to U) \in S} F(U') \, .
$$
Ainsi, la restriction de $F$ \`a ${\mathcal C}'$ est un $J'$-faisceau.

\medskip

\item[(iii)] D'apr\`es la proposition \ref{propI106}, le foncteur d'\'evaluation 
$$
[{\mathcal C}^{\rm op} , {\rm Ens}] = \widehat{\mathcal C} \longrightarrow \widehat{\mathcal C}' = [{\mathcal C}'^{\rm op} , {\rm Ens}] \, ,
$$
d\'efini par composition des pr\'efaisceaux ${\mathcal C}^{\rm op} \to {\rm Ens}$ avec le plongement ${\mathcal C}'^{\rm op} \to {\mathcal C}^{\rm op}$, admet pour adjoint \`a droite le foncteur d'extension de Kan
$$
\begin{matrix}
\widehat{\mathcal C}' &\longrightarrow &\widehat{\mathcal C} \, , \hfill \\
\hfill F &\longmapsto &\overline F = \displaystyle \left[ X \longmapsto \varprojlim_{(U \to X) \in {\mathcal C}' / X} F(U) \right]. 
\end{matrix}
$$
Le compos\'e de ces deux foncteurs s'identifie \`a
$$
{\rm id} : \widehat{\mathcal C}' \longrightarrow \widehat{\mathcal C}' \, ,
$$
ce qui signifie que le foncteur $F \mapsto \overline F$ est pleinement fid\`ele.

\smallskip

Montrons que si $F$ est un $J'$-faisceau sur ${\mathcal C}'$, alors $\overline F$ est un $J$-faisceau sur ${\mathcal C}$.

\smallskip

Consid\'erons pour cela un objet $X$ de ${\mathcal C}$ et un crible $J$-couvrant $S$ de $X$. Pour tout morphisme $p : U \to X$ de ${\mathcal C}$ dont la source $U$ est un objet de ${\mathcal C}'$, le crible $p^* S$ de $U$ est $J$-couvrant donc contient un crible $J'$-couvrant $S_p$ de $U$ dans ${\mathcal C}'$. 

\smallskip

Si $F$ est un $J'$-faisceau, on en d\'eduit que la limite
$$
\varprojlim_{(U \to X) \in {\mathcal C}' / X} F(U)
$$
s'identifie \`a la limite
$$
\varprojlim_{(U \to X) \in ({\mathcal C}' / X)_S} F(U)
$$
calcul\'ee sur la sous-cat\'egorie pleine $({\mathcal C}'/X)_S$ de ${\mathcal C}'/X$ constitu\'ee des morphismes $U \to X$ qui sont \'el\'ements de $S$.

\smallskip

Cela implique que $\overline F$ est un $J$-faisceau sur ${\mathcal C}$.

\smallskip

Ainsi, le foncteur $F \mapsto \overline F$ se restreint en un foncteur
$$
\begin{matrix}
\widehat{\mathcal C}'_{J'} &\longrightarrow &\widehat{\mathcal C}_J \, , \\
\hfill F&\longmapsto &\overline F \, .
\end{matrix}
$$
Comme $\widehat{\mathcal C}_J$ et $\widehat{\mathcal C}'_{J'}$ sont des sous-cat\'egories pleines de $\widehat{\mathcal C}$ et $\widehat{\mathcal C}'$, il est adjoint \`a droite du foncteur d'\'evaluation
$$
\widehat{\mathcal C}_J \longrightarrow \widehat{\mathcal C}'_{J'} \, .
$$
Il est pleinement fid\`ele, et le compos\'e
$$
\widehat{\mathcal C}'_{J'} \longrightarrow \widehat{\mathcal C}_J \longrightarrow \widehat{\mathcal C}'_{J'}
$$ 
s'identifie au foncteur ${\rm id} : \widehat{\mathcal C}'_{J'} \to \widehat{\mathcal C}'_{J'}$.

\smallskip

Il reste \`a montrer que le foncteur compos\'e
$$
\widehat{\mathcal C}_J \longrightarrow \widehat{\mathcal C}'_{J'} \longrightarrow \widehat{\mathcal C}_J
$$
est canoniquement isomorphe \`a ${\rm id} : \widehat{\mathcal C}_J \to \widehat{\mathcal C}_J$ ou, ce qui revient au m\^eme, que $\widehat{\mathcal C}_J \to \widehat{\mathcal C}'_{J'}$ est pleinement fid\`ele.

\smallskip

Or, tout objet $X$ de ${\mathcal C}$ admet une famille $J$-couvrante de morphismes $U_i \to X$ dont les sources $U_i$ sont des objets de ${\mathcal C}'$ et pour tout carr\'e commutatif de ${\mathcal C}$
$$
\xymatrix{
V_{i,j} \ar[d] \ar[r] &U_i \ar[d] \\
U_j \ar[r] &X
}
$$
il existe une famille $J$-couvrante de morphismes
$$
U_{i,j,k} \longrightarrow V_{i,j}
$$
dont les sources $U_{i,j,k}$ sont des objets de ${\mathcal C}'$.

\smallskip

Si $F$ est un $J$-faisceau, on en d\'eduit comme voulu que l'application canonique
$$
F(X) \longrightarrow \varprojlim_{(U \to X) \in {\mathcal C}'/X} F(U)
$$
est un isomorphisme.

\smallskip

Cela prouve (iii) et donc aussi (ii). 

\end{listeisansmarge}
\end{demo}

\newpage

Cette proposition comprend comme cas particulier le corollaire suivant:

\begin{cor}\label{corIII15}

Soit $X$ un espace topologique.

\smallskip

Soit ${\mathcal V}$ une base d'ouverts de $X$, consid\'er\'ee comme une sous-cat\'egorie pleine de l'ensemble ordonn\'e $O(X)$ des ouverts de $X$.

\smallskip

Soit $J$ la topologie de ${\mathcal V}$ induite par la topologie canonique de $O(X)$.

\smallskip

Alors le foncteur de restriction des faisceaux sur $X$ aux ouverts de la base ${\mathcal V}$
$$
{\mathcal E}_X \longrightarrow \widehat{\mathcal V}_J
$$
d\'efinit une \'equivalence de la cat\'egorie ${\mathcal E}_X$ des faisceaux sur l'espace topologique $X$ \`a la cat\'egorie $\widehat{\mathcal V}_J$ des faisceaux sur le site $({\mathcal V} , J)$.

\end{cor}

\begin{demo}

En effet, dire qu'une famille ${\mathcal V}$ d'ouverts de $X$ est une base signifie exactement que, consid\'er\'ee comme une sous-cat\'egorie pleine de $O(X)$, elle est dense.

\end{demo}

\section{Propri\'et\'es de compl\'etude et de cocompl\'etude des topos}\label{sec32}

\subsection{Compl\'etude des topos}\label{subsec321}

\medskip

On a:

\begin{prop}\label{propIII21}
\begin{listeimarge}
\item Tout topos ${\mathcal E}$ est une cat\'egorie compl\`ete au sens que les limites y sont toujours d\'efinies:

\begin{enumerate}
\item[(1)] Pour tout carquois $D$, le foncteur diagonal
$$
\Delta_D : {\mathcal E} \longrightarrow D\mbox{\rm -diag} ({\mathcal E})
$$
admet un adjoint \`a droite
$$
\varprojlim_D : D\mbox{\rm -diag} ({\mathcal E}) \longrightarrow {\mathcal E} \, .
$$
\end{enumerate}

\item De plus, si ${\mathcal E} = \widehat{\mathcal C}_J$ est le topos des faisceaux sur un site $({\mathcal C},J)$ et $F_{\bullet}$ est un $D$-diagramme de faisceaux, on a pour tout objet $X$ de ${\mathcal C}$ une identification canonique
$$
\left( \varprojlim_D F_{\bullet} \right) (X) = \varprojlim_D F_{\bullet} (X) \, .
$$
\end{listeimarge}
\end{prop}

\begin{remark}

En particulier, un topos ${\mathcal E}$ poss\`ede toujours un objet terminal que l'on note habituellement $1_{\mathcal E}$ ou $1$.

\smallskip

Si ${\mathcal E} = \widehat{\mathcal C}_J$, $1(X)$ est l'ensemble \`a un \'el\'ement pour tout objet $X$ de ${\mathcal C}$.

\end{remark}

\begin{demo}
\begin{listeisansmarge}
\item[(ii)] On sait d\'ej\`a d'apr\`es la proposition \ref{propI912} que la cat\'egorie $\widehat{\mathcal C} = [{\mathcal C}^{\rm op} , {\rm Ens}]$ est compl\`ete et plus pr\'ecis\'ement que tout $D$-diagramme de pr\'efaisceaux $P_{\bullet}$ admet pour limite le pr\'efaisceau
$$
\varprojlim_D P_{\bullet}
$$
d\'efini en tout objet $X$ de ${\mathcal C}$ par la formule
$$
\left( \varprojlim_D P_{\bullet} \right)(X) = \varprojlim_D P_{\bullet} (X) \, .
$$

De plus, si tous les $P_d$, $d \in {\rm Ob} (D)$, sont des faisceaux pour la topologie $J$, il en est de m\^eme de leur limite $\underset{D}{\varprojlim} \, P_{\bullet}$. En effet, un pr\'efaisceau $F$ sur ${\mathcal C}$ est un $J$-faisceau si et seulement si pour tout crible $J$-couvrant $S$ d'un objet $X$ de ${\mathcal C}$, le morphisme canonique
$$
F(X) \longrightarrow \varprojlim_{(U \to X) \in S} F(U)
$$
est un isomorphisme. Cette condition est pr\'eserv\'ee par les foncteurs $\underset{D}{\varprojlim}$ car ils respectent les limites.

\smallskip

La conclusion de (ii) r\'esulte alors de ce que $\widehat{\mathcal C}_J$ est par d\'efinition une sous-cat\'egorie pleine de $\widehat{\mathcal C}$.

\medskip

\item[(i)] r\'esulte de (ii) puisque les limites sont toujours pr\'eserv\'ees par les \'equivalences de cat\'egories. 
\end{listeisansmarge}
\end{demo}


La structure cat\'egorique de produit dans un topos ${\mathcal E}$ permet de d\'efinir les notions de mono{\"\i}de interne, de groupe interne ou d'anneau interne d'un topos:

\begin{defn}\label{defIII22}

Soit ${\mathcal E}$ un topos ou plus g\'en\'eralement une cat\'egorie localement petite qui poss\`ede des produits finis, en particulier un objet terminal $1_{\mathcal E}$.

\begin{listeimarge}

\item Un mono{\"\i}de interne de ${\mathcal E}$ consiste en un objet $M$ de ${\mathcal E}$ muni de

\medskip

$\left\{\begin{matrix}
\bullet &\mbox{un morphisme dit de multiplication} \hfill \\
{ \ } \\
&M \times M \xrightarrow{ \ m \ } M \\
{ \ } \\
&\mbox{qui est associatif au sens que commute le carr\'e:} \hfill \\
{ \ } \\
&\xymatrix{
M \times M \times M \ar[d]_{{\rm id}_M \times m} \ar[rr]^{m \times {\rm id}_M} &&M \times M \ar[d]^m \\
M \times M \ar[rr]^{m} &&M
} \\
{ \ } \\
\bullet &\mbox{un morphisme dit d'\'el\'ement unit\'e} \hfill \\
{ \ } \\
&1_{\mathcal E} \xrightarrow{ \ e \ } M \\
{ \ } \\
&\mbox{qui est neutre pour $m$ au sens que commute le diagramme:} \hfill \\
{ \ } \\
&\xymatrix{
M \times 1_{\mathcal E} \ar[d]_{{\rm id}_M \times e} \ar[r]^-{\sim} &M \ar[d]^-{{\rm id}_M} &1_{\mathcal E} \times M \ar[d]^{e \times {\rm id}_M} \ar[l]_-{\sim} \\
M \times M \ar[r]^m &M &M \times M \ar[l]_m
}
\end{matrix} \right.$

\bigskip

Un morphisme de mono{\"\i}des internes
$$
(M,m,e) \longrightarrow (M',m',e')
$$
est un morphisme de ${\mathcal E}$
$$
\rho : M \longrightarrow M'
$$
qui respecte les lois de multiplication et les \'el\'ements neutres au sens que commutent les deux diagrammes:
$$
\xymatrix{
M \times M \ar[d]_-m \ar[rr]^{\rho \times \rho} &&M' \times M' \ar[d]^-{m'} \\
M \ar[rr]^{\rho} &&M'
} \qquad\qquad \xymatrix{
&1_{\mathcal E} \ar[ld]_e \ar[rd]^{e'} \\
M \ar[rr] &&M'
}
$$

\item Un groupe interne de ${\mathcal E}$ est un mono{\"\i}de interne $(M,m,e)$ muni de surcro{\^\i}t d'un morphisme de passage \`a l'inverse
$$
i : M \longrightarrow M
$$
tel que commutent les carr\'es:
$$
\xymatrix{
M \ar[d] \ar[rr]^-{({\rm id} , i )} &&M \times M \ar[d]^m \\
1_{\mathcal E} \ar[rr]^{e} &&M
} \qquad\qquad \xymatrix{
M \ar[d] \ar[rr]^-{(i,{\rm id})} &&M \times M \ar[d]^m \\
1_{\mathcal E} \ar[rr]^e &&M
}
$$

Un morphisme de groupes internes est un morphisme de mono{\"\i}des internes.

\medskip

\item Un groupe [resp. mono{\"\i}de] interne $(M,m,e)$ est dit commutatif ou ab\'elien si, notant $M \times M \xrightarrow{ \ \sigma \ } M \times M$ le morphisme de permutation des deux composantes, commute le triangle:
$$
\xymatrix{
M \times M \ar[rr]^{\sigma} \ar[rd]_m &&M \times M \ar[ld]^m \\
&M
}
$$

\item Un anneau [resp. anneau commutatif] interne est constitu\'e d'un groupe ab\'elien interne $(A,+,0,-)$ et d'un mono{\"\i}de [resp. mono{\"\i}de commutatif] interne $(A,\cdot,1)$ de m\^eme objet sous-jacent $A$, tels que commutent les carr\'es:
$$
\xymatrix{
A \times A \times A \ar[d]_{{\rm id} \times (+)} \ar[rrr]^-{(p_1 \cdot p_2 , p_1 \cdot p_3)} &&&A \times A \ar[d]^+ \\
A \times A \ar[rrr]^{\cdot} &&&A
} \qquad\qquad \xymatrix{
A \times A \times A \ar[d]_{(+) \times {\rm id}} \ar[rrr]^-{(p_1 \cdot p_3 , p_2 \cdot p_3)} &&&A \times A \ar[d]^+ \\
A \times A \ar[rrr]^{\cdot} &&&A
}
$$

Un morphisme d'anneaux internes
$$
(A,+,0,-,\cdot , 1) \longrightarrow (A' , +,0,-,\cdot , 1)
$$
est un morphisme de groupes internes
$$
(A , + , 0 , -) \longrightarrow (A' , +,0,-)
$$
qui est aussi un morphisme de mono{\"\i}des internes
$$
(A , \cdot , 1) \longrightarrow (A',\cdot ,1) \, .
$$
\end{listeimarge}
\end{defn}

\begin{remarksqed}
\begin{listeisansmarge}
\item D'apr\`es le lemme de Yoneda, un triplet $(M,M \times M \xrightarrow{ \ m \ } M , 1_{\mathcal E} \xrightarrow{ \ e \ } M)$ est un mono{\"\i}de [resp. groupe, resp. mono{\"\i}de commutatif, resp. groupe commutatif] interne de ${\mathcal E}$ si et seulement si, pour tout objet $E$ de ${\mathcal E}$, l'ensemble ${\rm Hom} (E,M)$ muni de l'application induite ${\rm Hom} (E,M) \times {\rm Hom} (E,M) \to {\rm Hom} (E,M)$ et de l'\'el\'ement $E \longrightarrow 1_{\mathcal E} \xrightarrow{ \ e \ } M$ est un mono{\"\i}de [resp. un groupe, resp. un mono{\"\i}de commutatif, resp. un groupe commutatif].

\smallskip

De m\^eme, un sextuplet $(A,+,0,-,\cdot , 1)$ est un anneau [resp. anneau commutatif] interne de ${\mathcal E}$ si et seulement si, pour tout objet $E$ de ${\mathcal E}$, l'ensemble ${\rm Hom} (E,M)$ muni des lois et des \'el\'ements unit\'es induits est un anneau [resp. anneau commutatif].

\medskip

\item Si ${\mathcal E} = \widehat{\mathcal C}_J$ est le topos des faisceaux sur un site $({\mathcal C},J)$, il suffit de consid\'erer dans (i) les objets $E$ de ${\mathcal E}$ qui sont les images par le foncteur canonique $\ell : {\mathcal C} \to \widehat{\mathcal C}_J = {\mathcal E}$ d'objets $X$ de ${\mathcal C}$. Autrement dit, un mono{\"\i}de [resp. groupe, resp. anneau] interne de ${\mathcal E}$ est un faisceau de mono{\"\i}des [resp. de groupes, resp. d'anneaux] sur le site $({\mathcal C},J)$, et il est commutatif s'il en est ainsi de ses \'evaluations en les objets $X$ de ${\mathcal C}$. Cela r\'esulte de ce que la formation des produit dans $\widehat{\mathcal C}_J$ est respect\'ee par les foncteurs d'\'evaluation en les objets de ${\mathcal C}$. 
\end{listeisansmarge}
\end{remarksqed}

\newpage

On d\'eduit de la proposition \ref{propIII21}:

\begin{cor}\label{corIII23}
\begin{listeimarge}
\item La cat\'egorie des mono{\"\i}des [resp. groupes, resp. groupes commutatifs, resp. mono{\"\i}des commutatifs, resp. anneaux, resp. anneaux commutatifs] internes d'un topos ${\mathcal E}$ est compl\`ete, c'est-\`a-dire poss\`ede des limites arbitraires.

\medskip

\item De plus, le foncteur d'oubli des structures, qui va de cette cat\'egorie vers le topos ${\mathcal E}$, respecte les limites.

\medskip

\item Si ${\mathcal E} = \widehat{\mathcal C}_J$ est le topos des faisceaux sur un site $({\mathcal C} , J)$, le foncteur d'\'evaluation en un objet $X$ de ${\mathcal C}$ de cette cat\'egorie vers celle des mono{\"\i}des [resp. groupes, resp. groupes commutatifs, resp. mono{\"\i}des commutatifs, resp. anneaux, resp. anneaux commutatifs] respecte les limites.
\end{listeimarge}
\end{cor}

\begin{demo}

Compte tenu de la proposition \ref{propIII21}, cela r\'esulte de ce que la d\'efinition \ref{defIII22} ne fait intervenir que des produits d'objets reli\'es par des morphismes, puisque tous les foncteurs de limites respectent les produits. 

\end{demo}

\subsection{Cocompl\'etude des topos}\label{subsec322}

On a:

\begin{prop}\label{propIII24}
\begin{listeimarge}
\item Tout topos ${\mathcal E}$ est une cat\'egorie cocompl\`ete au sens que les colimites y sont toujours d\'efinies:

\begin{enumerate}
\item[(2)] Pour tout carquois $D$, le foncteur diagonal
$$
\Delta_D : {\mathcal E} \longrightarrow D\mbox{\rm -diag} ({\mathcal E})
$$
admet un adjoint \`a gauche
$$
\varinjlim_D : D\mbox{\rm -diag} ({\mathcal E}) \longrightarrow {\mathcal E} \, .
$$
\end{enumerate}

\item De plus, si ${\mathcal E} = \widehat{\mathcal C}_J$ est le topos des faisceaux sur un site $({\mathcal C},J)$ et $F_{\bullet}$ est un $D$-diagramme de faisceaux, le faisceau $\underset{D}{\varinjlim} \, F_{\bullet}$ est le transform\'e par le foncteur de faisceautisation $j^* : \widehat{\mathcal C} \to \widehat{\mathcal C}_J$ du pr\'efaisceau
$$
X \longmapsto \varinjlim_D F_{\bullet} (X) \, .
$$
\end{listeimarge}
\end{prop}

\begin{remark}

En particulier, un topos ${\mathcal E}$ poss\`ede toujours un objet initial que l'on note habituellement $\emptyset_{\mathcal E}$ ou $\emptyset$.

\smallskip

Si ${\mathcal E} = \widehat{\mathcal C}_J$, $\emptyset(X)$ est l'ensemble vide pour tout objet $X$ de ${\mathcal C}$.

\end{remark}

\begin{demo} 
\begin{listeisansmarge}
\item[(ii)] On sait d\'ej\`a d'apr\`es la proposition \ref{propI912} que la cat\'egorie $\widehat{\mathcal C} = [{\mathcal C}^{\rm op} , {\rm Ens}]$ est cocompl\`ete et plus pr\'ecis\'ement que tout $D$-diagramme de pr\'efaisceau $P_{\bullet}$ admet pour colimite le pr\'efaisceau
$$
\varinjlim_D P_{\bullet}
$$
d\'efini en tout objet $X$ de ${\mathcal C}$ par la formule
$$
\left( \varinjlim_D P_{\bullet} \right)(X) = \varinjlim_D P_{\bullet} (X) \, .
$$
Le foncteur de faisceautisation $j^* : \widehat{\mathcal C} \to \widehat{\mathcal C}_J$ est adjoint \`a gauche du foncteur $j_* : \widehat{\mathcal C}_J \hookrightarrow \widehat{\mathcal C}$ donc respecte les colimites arbitraires.

\smallskip

Il en r\'esulte que pour tout $D$-diagramme de faisceaux $F_{\bullet}$, le faisceau
$$
j^* \left( \varinjlim_D j_* F_{\bullet} \right)
$$
est colimite du $D$-diagramme $j^* j_* F_{\bullet}$ lequel s'identifie \`a $F_{\bullet}$.

\medskip

\item[(i)] r\'esulte de (ii) puisque les colimites sont toujours pr\'eserv\'ees par les \'equivalences de cat\'egories. 
\end{listeisansmarge}
\end{demo}


On d\'eduit de cette proposition:

\begin{cor}\label{corIII25}
\begin{listeimarge}
\item La cat\'egorie ${\mathcal M}_{\mathcal E}$ des mono{\"\i}des [resp. groupes, resp. groupes commutatifs, resp. mono{\"\i}des commutatifs, resp. anneaux, resp. anneaux commutatifs] internes d'un topos ${\mathcal E}$ est cocompl\`ete c'est-\`a-dire poss\`ede des colimites arbitraires.

\medskip

\item Si ${\mathcal E} = \widehat{\mathcal C}$ est le topos des pr\'efaisceaux sur une cat\'egorie essentiellement petite ${\mathcal C}$, la cat\'egorie ${\mathcal M}_{\mathcal E}$ s'identifie \`a la cat\'egorie des foncteurs ${\mathcal C}^{\rm op} \to {\mathcal M}_{\rm Ens}$ et, pour tout objet $X$ de ${\mathcal C}$, le foncteur d'\'evaluation
$$
\begin{matrix}
{\mathcal M}_{\mathcal E} = [{\mathcal C}^{\rm op} , {\mathcal M}_{\rm Ens}] &\longrightarrow &{\mathcal M}_{\rm Ens} \, , \\
\hfill M &\longmapsto &M(X) \hfill
\end{matrix}
$$
respecte les colimites.

\medskip

\item Si ${\mathcal E} = \widehat{\mathcal C}_J$ est le topos des faisceaux sur un site $({\mathcal C},J)$, les foncteurs $j^* : \widehat{\mathcal C} \to \widehat{\mathcal C}_J$ et $j_* : \widehat{\mathcal C}_J \hookrightarrow \widehat{\mathcal C}$ induisent des foncteurs adjoints
$$
j^* : {\mathcal M}_{\widehat{\mathcal C}} \to {\mathcal M}_{\mathcal E} \quad \mbox{et} \quad j_* : {\mathcal M}_{\mathcal E} \hookrightarrow {\mathcal M}_{\widehat{\mathcal C}}
$$
et tout $D$-diagramme $M_{\bullet}$ de ${\mathcal M}_{\mathcal E}$ admet pour colimite
$$
j^* \left( \varinjlim_D j_* M_{\bullet} \right) .
$$
\end{listeimarge}
\end{cor}

\begin{demo}
\begin{listeisansmarge}
\item Traitons d'abord le cas o\`u ${\mathcal E} = {\rm Ens}$ est le topos des ensembles.

\smallskip

Le morphisme d'oubli de la structure
$$
O : {\mathcal M}_{\rm Ens} \longrightarrow {\rm Ens}
$$
admet un adjoint \`a gauche
$$
I \longmapsto M_I
$$
qui associe \`a tout ensemble $I$ la ``structure libre'' engendr\'ee par les \'el\'ements de $I$.

\smallskip

Toute partie $I$ de l'ensemble sous-jacent \`a un objet $M$ de ${\mathcal M}_{\rm Ens}$ d\'efinit par adjonction un morphisme de ${\mathcal M}_{\rm Ens}$
$$
M_I \longrightarrow M \, .
$$
S'il est surjectif, on dit que $I$ engendre $M$ qui est alors le quotient de $M_I$ par une famille de relations.

\smallskip

R\'eciproquement, toute famille de relations dans une structure libre $M_I$ d\'efinit un objet quotient de ${\mathcal M}_{\rm Ens}$.

\smallskip

Si $M_{\bullet}$ est un $D$-diagramme de ${\mathcal M}_{\rm Ens}$, il admet pour colimite le quotient de la structure libre
$$
M_I
$$
associ\'ee \`a la somme disjointe des ensembles sous-jacents
$$
I = \coprod_{d \in {\rm Ob} (D)} O \, (M_d)
$$
par la r\'eunion des relations internes \`a chaque $M_d$ et des relations de la forme
$$
m' = u(m) \, , \qquad m' \in M_{d'} \, , \ m \in M_d \, ,
$$
associ\'ees \`a chaque morphisme $u : d \to d'$ de $D$.

\smallskip

Le cas o\`u ${\mathcal E} = {\rm Ens}$ \'etant d\'emontr\'e, il suffit de prouver (ii) et (iii) pour prouver aussi (i).

\medskip

\item La premi\`ere assertion est un cas particulier de la remarque (ii) qui suit la d\'efinition \ref{defIII22}.

\smallskip

Tout $D$-diagramme $M_{\bullet}$ de ${\mathcal M}_{\widehat{\mathcal C}}$ admet pour colimite le foncteur
$$
\begin{matrix}
{\mathcal C}^{\rm op} &\longrightarrow &{\mathcal M}_{\rm Ens} \, , \hfill \\
\hfill X &\longmapsto &\displaystyle \varinjlim_D M_{\bullet} (X)
\end{matrix}
$$
car $\underset{D}{\varinjlim} \, : D\mbox{-diag} ({\mathcal M}_{\rm Ens}) \to {\mathcal M}_{\rm Ens}$ est un foncteur.

\medskip

\item Les foncteurs adjoints
$$
j^* : \widehat{\mathcal C} \longrightarrow \widehat{\mathcal C}_J \qquad \mbox{et} \qquad j_* : \widehat{\mathcal C}_J \xhookrightarrow{ \ { \ } \ } \widehat{\mathcal C}
$$
induisent des foncteurs adjoints
$$
j^* : {\mathcal M}_{\widehat{\mathcal C}} \longrightarrow {\mathcal M}_{\mathcal E} \qquad \mbox{et} \qquad j_* :  {\mathcal M}_{\mathcal E} \longrightarrow {\mathcal M}_{\widehat{\mathcal C}}
$$
car ils respectent les produits finis et les structures des objets de ${\mathcal M}_{\widehat{\mathcal C}}$ ou ${\mathcal M}_{\mathcal E}$ sont d\'efinies en termes de produits finis et de morphismes entre eux.

\smallskip

Comme $j^* : {\mathcal M}_{\widehat{\mathcal C}} \to {\mathcal M}_{\mathcal E}$ est un adjoint \`a gauche, il respecte les colimites arbitraires.

\smallskip

Si donc $M_{\bullet}$ est un $D$-diagramme de ${\mathcal M}_{\mathcal E}$, l'objet de ${\mathcal M}_{\mathcal E}$
$$
j^* \left( \varinjlim_D j_* M_{\bullet} \right)
$$
est colimite du $D$-diagramme $j^* j_* M_{\bullet}$ lequel s'identifie \`a $M_{\bullet}$.

\smallskip

Cela termine la d\'emonstration. 
\end{listeisansmarge}
\end{demo}

\subsection{Cat\'egories de modules dans un topos}\label{subsec323}

\medskip

On compl\`ete la d\'efinition \ref{defIII22} par la suivante:

\begin{defn}\label{defIII26}

Soit $({\mathcal E} , A)$ un ``topos annel\'e'' constitu\'e d'un topos ${\mathcal E}$ et d'un anneau interne $A$ de ${\mathcal E}$ ou, plus g\'en\'eralement, soient ${\mathcal E}$ une cat\'egorie localement petite qui poss\`ede des produits finis et $A$ un anneau interne de ${\mathcal E}$.

\smallskip

Un $A$-module dans ${\mathcal E}$ consiste en un groupe ab\'elien interne $M$ de ${\mathcal E}$ muni d'un morphisme de ${\mathcal E}$
$$
A \times M \longrightarrow M
$$
tels que commutent les carr\'es:
$$
\xymatrix{
A \times A \times M \ar[d]_{{\rm id}_A \times \cdot} \ar[rr]^-{\cdot \times {\rm id}_M} &&A \times M \ar[d] \\
A \times M \ar[rr] &&M
} \qquad \qquad \xymatrix{
A \times M \times M \ar[d]_{{\rm id}_A \times (+)} \ar[rrr]^-{(p_1 \cdot p_2 , p_1 \cdot p_3)} &&&M \times M \ar[d]^+ \\
A \times M \ar[rrr] &&&M
}
$$

Un morphisme de $A$-modules est un morphisme de groupes ab\'eliens internes
$$
M_1 \xrightarrow{ \ u \ } M_2
$$
qui rend commutatif le carr\'e:
$$
\xymatrix{
A \times M_1 \ar[d] \ar[rr]^-{{\rm id}_A \times u} &&A \times M_2 \ar[d] \\
M_1 \ar[rr]^-u &&M_2
}
$$

On note $A$-{\rm mod} ou ${\mathcal M}od_A$ la cat\'egorie des $A$-modules internes de ${\mathcal E}$.
\end{defn}

\begin{remarkqed}

Si ${\mathcal E} = \widehat{\mathcal C}_J$ est le topos des faisceaux sur un site $({\mathcal C},J)$, $A$ est un faisceau d'anneaux.

\smallskip

Alors un $A$-module interne de ${\mathcal E}$ est un faisceau de modules sur le site annel\'e $({\mathcal C} , J , A)$ au sens de la d\'efinition \ref{defII47} (ii). 

\smallskip 

Dans ce cas, on note aussi ${\mathcal M}od_A = {\mathcal M}od_{{\mathcal C},J,A}$ et ${\mathcal M}od_A = {\mathcal M}od_{{\mathcal C},A}$ si $J$ est la topologie discr\`ete de ${\mathcal C}$ et ${\mathcal E} = \widehat{\mathcal C}$. 

\end{remarkqed}


On d\'emontre \`a la suite du corollaire \ref{corIII23}:

\begin{cor}\label{corIII27}

Soit $({\mathcal E} , A)$ un topos annel\'e.

\begin{listeimarge}

\item La cat\'egorie ${\mathcal M}od_A$ des $A$-modules de ${\mathcal E}$ est compl\`ete, c'est-\`a-dire poss\`ede des limites arbitraires, et le foncteur d'oubli des morphismes de structure $A \times M \to M$
$$
{\mathcal M}od_A \longrightarrow {\mathcal A}b_{\mathcal E}
$$
vers la cat\'egorie ${\mathcal A}b_{\mathcal E}$ des groupes ab\'eliens internes de ${\mathcal E}$ respecte toutes les limites.

\medskip

\item Si ${\mathcal E} = \widehat{\mathcal C}_J$ est le topos des faisceaux sur un site $({\mathcal C},J)$, le foncteur d'\'evaluation en un objet $X$ de ${\mathcal C}$
$$
\begin{matrix}
\hfill M &\longmapsto &M(X) \, , \hfill \\
{\mathcal M}od_A &\longrightarrow &{\rm Mod}_{A(X)}
\end{matrix}
$$
vers la cat\'egorie des $A(X)$-modules respecte les limites.
\end{listeimarge}
\end{cor}

\newpage

\begin{demo}
\begin{listeisansmarge}
\item Consid\'erons un carquois $D$ et un $D$-diagramme $M_{\bullet}$ de ${\mathcal M}od_A$. Soit $M$ la limite du $D$-diagramme induit dans ${\mathcal A}b_{\mathcal E}$. Comme le foncteur d'oubli
$$
{\mathcal A}b_{\mathcal E} \longrightarrow {\mathcal E}
$$
respecte les limites, les morphismes de structure
$$
A \times M_d \longrightarrow M_d \, , \quad d \in {\rm Ob} (D) \, ,
$$
induisent un morphisme
$$
A \times M \longrightarrow M \, .
$$
Ce morphisme rend commutatifs les deux carr\'es de la d\'efinition \ref{defIII26} puisqu'il en est ainsi des morphismes $A \times M_d \to M_d$, $d \in {\rm Ob} (D)$.

\smallskip

Ainsi, $M$ a une structure de $A$-module dans ${\mathcal E}$.

\smallskip

Enfin, si $N$ est un groupe ab\'elien interne de ${\mathcal E}$, se donner un morphisme
$$
N \longrightarrow M
$$
\'equivaut \`a se donner une famille compatible de morphismes
$$
N \longrightarrow M_d \, , \quad d \in {\rm Ob} (D) \, .
$$
Un morphisme $A \times N \to N$ rend commutatif le carr\'e
$$
\xymatrix{
A \times N \ar[d] \ar[r] &A \times M \ar[d] \\
N \ar[r] &M
}
$$
si et seulement si il rend commutatifs les carr\'es:
$$
\xymatrix{
A \times N \ar[d] \ar[r] &A \times M_d \ar[d] \\
N \ar[r] &M_d
}
$$
Cela implique que $M$ muni de son morphisme de structure $A \times M \to M$ est une limite du $D$-diagramme $M_{\bullet}$ dans ${\mathcal M}od_A$.

\medskip

\item Comme le foncteur d'oubli ${\mathcal M}od_A \to {\mathcal A}b_{\mathcal E}$ et le foncteur d'\'evaluation
$$
\begin{matrix}
{\mathcal A}b_{\mathcal E} &\longrightarrow &{\rm Ab} \, , \hfill \\
\hfill M &\longmapsto &M(X) \hfill
\end{matrix}
$$
respectent les limites, il en va de m\^eme de leur compos\'e qui est aussi le compos\'e
$$
{\mathcal M}od_A \longrightarrow {\rm Mod}_{A(X)} \longrightarrow {\rm Ab} \, .
$$
Si $M$ est la limite d'un $D$-diagramme $M_{\bullet}$ de ${\mathcal M}od_A$, le morphisme de structure
$$
A \times M \longrightarrow M
$$
est la limite dans ${\mathcal E}$ du syst\`eme compatible de morphismes
$$
A \times M_d \longrightarrow M_d \, , \quad d \in {\rm Ob} (D) \, .
$$
Comme les limites dans ${\mathcal E} = \widehat{\mathcal C}_J$ se calculent composante par composante, l'application de structure
$$
A(X) \times M(X) \longrightarrow M(X)
$$
est la limite dans ${\rm Ens}$ du syst\`eme compatible d'applications
$$
A_d(X) \times M_d(X) \longrightarrow M_d(X) \, .
$$
Cela signifie comme voulu que le $A(X)$-module $M(X)$ est la limite du $D$-diagramme des $A(X)$-modules $M_d (X)$. 

\end{listeisansmarge}
\end{demo}


On d\'emontre d'autre part \`a la suite du corollaire \ref{corIII25}:

\begin{cor}\label{corIII28}

 Soit $({\mathcal E} , A)$ un topos annel\'e.

\begin{listeimarge}

\item La cat\'egorie ${\mathcal M}od_A$ des $A$-modules de ${\mathcal E}$ est cocompl\`ete, c'est-\`a-dire poss\`ede des colimites arbitraires, et le foncteur d'oubli
$$
{\mathcal M}od_A \longrightarrow {\mathcal A}b_{\mathcal E}
$$
respecte toutes les colimites.

\medskip

\item Si ${\mathcal E} = \widehat{\mathcal C}$ est le topos des pr\'efaisceaux sur une cat\'egorie essentiellement petite ${\mathcal E}$, le foncteur d'\'evaluation
$$
\begin{matrix}
{\mathcal M}od_A &\longrightarrow &{\rm Mod}_{A(X)} \, , \\
\hfill M &\longmapsto &M(X) \hfill
\end{matrix}
$$
respecte les colimites.

\medskip

\item Si ${\mathcal E} = \widehat{\mathcal C}_J$ est le topos des faisceaux sur un site $({\mathcal C},J)$, les foncteurs $j^* : \widehat{\mathcal C} \to \widehat{\mathcal C}_J$ et $j_* : \widehat{\mathcal C}_J \to \widehat{\mathcal C}$ induisent des foncteurs adjoints
$$
j^* : {\mathcal M}od_{{\mathcal C},A} \longrightarrow {\mathcal M}od_{{\mathcal C},J,A} \quad \mbox{et} \quad j_* : {\mathcal M}od_{{\mathcal C},J,A} \longrightarrow {\mathcal M}od_{{\mathcal C},A}
$$
et tout $D$-diagramme $M_{\bullet}$ de ${\mathcal M}_{{\mathcal C},J,A}$ admet pour colimite
$$
j^* \left( \varinjlim_D j_* M_{\bullet} \right).
$$
\end{listeimarge}
\end{cor}

\begin{remark}

En particulier, toute famille $(M_i)_{i \in I}$ de $A$-modules dans un topos ${\mathcal E}$ admet une somme not\'ee $\underset{i \in I}{\bigoplus} \, M_i$.

\end{remark}

\begin{demo}
\begin{listeisansmarge}
\item[] Traitons d'abord le cas o\`u ${\mathcal E} = \widehat{\mathcal C}$ est un topos de pr\'efaisceaux.

\smallskip

La cat\'egorie ${\mathcal M}od_{{\mathcal C},A}$ des $A$-modules dans $\widehat{\mathcal C}$ peut \^etre comme la cat\'egorie des pr\'efaisceaux de groupes ab\'eliens $M$ sur ${\mathcal C}$ munis d'une famille de morphismes de groupes ab\'eliens 
$$
m_a : M(X) \longrightarrow M(X)
$$
index\'es par les objets $X$ de ${\mathcal C}$ et les \'el\'ements $a \in A(X)$, telle que

\medskip

$\left\{ \begin{matrix}
\bullet &\mbox{pour tout objet $X$ de ${\mathcal C}$ et tous \'el\'ements $a,b \in A(X)$,} \hfill \\
{ \ } \\
&m_a \circ m_b = m_{ab} \, , \\
{ \ } \\
\bullet &\mbox{pour tout morphisme $X \xrightarrow{ \ f \ } Y$ de ${\mathcal C}$ et tout \'el\'ement $b \in A(Y)$} \\ &\mbox{d'image $a = A(f)(b) \in A(X)$, le carr\'e} \hfill \\
{ \ } \\
&\xymatrix{
M(Y) \ar[d]_{M(f)} \ar[r]^{m_b} &M(Y) \ar[d]^{M(f)} \\
M(X) \ar[r]^{m_a} &M(X)
} \\
{ \ } \\
&\mbox{est commutatif.} \hfill
\end{matrix} \right.$

\bigskip

\noindent Pour tout $D$-diagramme $M_{\bullet}$ de ${\mathcal M}od_{{\mathcal C},A}$, leur colimite $M$ dans ${\mathcal A}b_{\widehat{\mathcal C}}$ se trouve munie naturellement d'une telle famille de morphismes v\'erifiant ces conditions.

\smallskip

Si $N$ est un $A$-module dans $\widehat{\mathcal C}$, se donner un morphisme de $A$-modules
$$
M \longrightarrow N
$$
\'equivaut \`a se donner une famille compatible de morphismes de $A$-modules
$$
M_d \longrightarrow N \, , \quad d \in {\rm Ob} (D) \, ,
$$
ce qui signifie que $M$ est une colimite du $D$-diagramme $M_{\bullet}$ dans ${\mathcal M}od_{{\mathcal C},A}$.

\smallskip

Cela prouve (ii) ainsi que (i) dans le cas o\`u ${\mathcal E}$ est un topos de pr\'efaisceaux $\widehat{\mathcal C}$.

\item[(iii)] Les foncteurs adjoints
$$
j^* : \widehat{\mathcal C} \longrightarrow \widehat{\mathcal C}_J \quad \mbox{et} \quad j_* : \widehat{\mathcal C}_J \xhookrightarrow{ \ { \ } \ } \widehat{\mathcal C}
$$
induisent des foncteurs adjoints
$$
j^* : {\mathcal M}od_{{\mathcal C},A} \longrightarrow {\mathcal M}od_{{\mathcal C},J,A} \quad \mbox{et} \quad j_* : {\mathcal M}od_{{\mathcal C},J,A} \longrightarrow {\mathcal M}od_{{\mathcal C},A}
$$
car ils respectent les produits finis et les structures des objets de ${\mathcal M}od_{{\mathcal C},A}$ ou ${\mathcal M}od_{{\mathcal C},A,J}$ sont d\'efinies en termes de produits finis et de morphismes entre eux.

\smallskip

Comme $j^* : {\mathcal M}od_{{\mathcal C},A} \to {\mathcal M}od_{{\mathcal C},J,A}$ est un adjoint \`a gauche il respecte les colimites arbitraires.

\smallskip

Si donc $M_{\bullet}$ est un $D$-diagramme de ${\mathcal M}od_{{\mathcal C},J,A}$, l'objet de ${\mathcal M}od_{{\mathcal C},J,A}$
$$
j^* \left( \varinjlim_D j_* M_{\bullet} \right)
$$
est une colimite du diagramme $j^* j_* \, M_{\bullet}$, lequel s'identifie \`a $M_{\bullet}$.

\medskip

\item[(i)] Les propri\'et\'es consid\'er\'ees sont invariantes par \'equivalence de cat\'egories, donc il suffit de traiter le cas o\`u ${\mathcal E} = \widehat{\mathcal C}_J$ est le topos des faisceaux sur un site $({\mathcal C},J)$.

\smallskip

Comme (iii) est d\'ej\`a d\'emontr\'e, il reste \`a montrer que le foncteur d'oubli
$$
{\mathcal M}od_{{\mathcal C},J,A} \longrightarrow {\mathcal A}b_{\widehat{\mathcal C}_J}
$$
respecte les colimites.

\smallskip

On sait d\'ej\`a qu'il en est ainsi du foncteur d'oubli
$$
{\mathcal M}od_{{\mathcal C},A} \longrightarrow {\mathcal A}b_{\widehat{\mathcal C}} \, .
$$

La conclusion r\'esulte de ce que le carr\'e
$$
\xymatrix{
{\mathcal M}od_{{\mathcal C},A} \ar[d]_{j^*} \ar[r] &{\mathcal A}b_{\widehat{\mathcal C}} \ar[d]^{j^*} \\
{\mathcal M}od_{{\mathcal C},J,A} \ar[r] &{\mathcal A}b_{\widehat{\mathcal C}_J}
}
$$

\noindent est commutatif, puisque les deux foncteurs $j^*$ verticaux respectent les colimites. 

\end{listeisansmarge}
\end{demo}

\bigskip

Les cat\'egories ${\mathcal M}od_{A}$ des $A$-modules d'un topos annel\'e $({\mathcal E} , A)$ partagent encore la propri\'et\'e fondamentale d'\^etre des cat\'egories additives au sens suivant:

\begin{defn}\label{defIII29}
\begin{listeimarge}
\item Une cat\'egorie localement petite ${\mathcal M}$ est dite ``semi-additive'' si
\begin{enumerate}
\item[$\bullet$] elle poss\`ede des produits finis et des sommes finies (not\'ees $\oplus$) arbitraires, en particulier un objet initial $0$ et un objet terminal $1$,
\item[$\bullet$] l'unique morphisme
$$
0 \longrightarrow 1
$$
est un isomorphisme,
\item[$\bullet$] pour tous objets $M$ et $N$ de ${\mathcal M}$, le morphisme
$$
M \oplus N \longrightarrow M \times N
$$
somme de
$$
M \xrightarrow{ \ \sim \ } M \times 1 \xrightarrow{ \ \sim \ } M \times 0 \longrightarrow M \times N
$$
et
$$
N \xrightarrow{ \ \sim \ } 1 \times N \xrightarrow{ \ \sim \ } 0 \times N \longrightarrow M \times N \, ,
$$
est un isomorphisme.
\end{enumerate}

\medskip

\item Une telle cat\'egorie semi-additive ${\mathcal M}$ est dite additive si, pour tout objet $M$ de ${\mathcal M}$, le groupo{\"\i}de commutatif interne d\'efini par la loi d'addition
$$
M \times M \xrightarrow{ \ \sim \ } M \oplus M \longrightarrow M
$$
et le morphisme d'\'el\'ement neutre
$$
1 \xrightarrow{ \ \sim \ } 0 \longrightarrow M
$$
est un groupe ab\'elien interne.
\end{listeimarge}
\end{defn}

\begin{remarksqed}
\begin{listeisansmarge}
\item Une cat\'egorie localement petite ${\mathcal M}$ est semi-additive [resp. additive] si elle poss\`ede des produits finis arbitraires et qu'il est possible de munir les ensembles
$$
{\rm Hom} (N,M)
$$
de structures de mono{\"\i}des commutatifs [resp. de groupes ab\'eliens] telles que, pour tout morphisme
$$
u : M \longrightarrow M' \qquad \mbox{ou} \qquad v : N' \longrightarrow N \, ,
$$
les applications de composition avec $u$ ou $v$
$$
{\rm Hom} (N,M) \longrightarrow {\rm Hom} (N,M')
$$
ou
$$
{\rm Hom} (N,M) \longrightarrow {\rm Hom} (N',M)
$$
respectent ces structures.

\smallskip

Dans ce cas, ces structures sont uniquement d\'etermin\'ees par la structure cat\'egorique de ${\mathcal M}$.

\medskip

\item Il r\'esulte de (i) que si ${\mathcal C}$ est une cat\'egorie additive [resp. semi-additive], toute sous-cat\'egorie pleine de ${\mathcal C}$ qui admet des produits finis est additive [resp. semi-additive].

\medskip

\item Une cat\'egorie ${\mathcal M}$ est additive [resp. semi-additive] si et seulement si sa cat\'egorie oppos\'ee ${\mathcal M}^{\rm op}$ l'est. 

\end{listeisansmarge}
\end{remarksqed}


Cette d\'efinition \'etant pos\'ee, on peut donc \'enoncer:

\begin{lem}\label{lemIII210}

Pour tout topos annel\'e $({\mathcal E},A)$ ou plus g\'en\'eralement pour toute cat\'egorie localement petite ${\mathcal E}$ qui admet des produits finis et pour tout anneau interne $A$ de ${\mathcal E}$, la cat\'egorie ${\mathcal M}od_A$ des $A$-modules de ${\mathcal E}$ est additive.
\end{lem}

\begin{remark}

En particulier, pour tout topos ${\mathcal E}$ ou plus g\'en\'eralement pour toute cat\'egorie localement petite ${\mathcal E}$ qui admet des produits finis, la cat\'egorie ${\mathcal A}b_{\mathcal E}$ des groupes ab\'eliens internes de ${\mathcal E}$ est additive.

\end{remark}

\begin{demo}

Cela r\'esulte de la remarque (i) qui suit la d\'efinition \ref{defIII29}.

\smallskip

En effet, pour tous objets $M,N$ de ${\mathcal M}od_A$, le morphisme de structure $M \times M \to M$ munit l'ensemble ${\rm Hom} (N,M)$ d'une structure de groupe ab\'elien. Ces structures sont respect\'ees par les morphismes $N' \to N$ ou $M \to M'$ de ${\mathcal M}od_A$.

\end{demo}

\section{Correspondance entre quotients et relations d'\'equivalence}\label{sec33}

\subsection{La propri\'et\'e de balancement}\label{subsec331}

\medskip

Les topos h\'eritent encore de ${\rm Ens}$ l'importante propri\'et\'e de balancement:

\begin{prop}\label{propIII31}

Tout topos ${\mathcal E}$ est une cat\'egorie ``balanc\'ee'' au sens qu'elle poss\`ede la propri\'et\'e suivante:

\begin{enumerate}
	\item[(3)] Tout morphisme de ${\mathcal E}$ qui est \`a la fois un monomorphisme et un \'epimorphisme est un isomorphisme.
\end{enumerate}
\end{prop}

\begin{demo}

C'est vrai si ${\mathcal E} = {\rm Ens}$ est le topos des ensembles car, dans cette cat\'egorie, un monomorphisme est une injection, un \'epimorphisme est une surjection, et un isomorphisme est une bijection.

\smallskip

Ce cas entra{\^\i}ne celui des cat\'egories ${\mathcal E} = \widehat{\mathcal C}$ de pr\'efaisceaux sur une cat\'egorie essentiellement petite ${\mathcal C}$ car un morphisme de $\widehat{\mathcal C}$
$$
u : F \longrightarrow G
$$
est un monomorphisme, un \'epimorphisme ou un isomorphisme si et seulement si chacune de ses composantes
$$
u_X : F(X) \longrightarrow G(X) \, , \quad X \in {\rm Ob} ({\mathcal C}) \, ,
$$
est un monomorphisme, un \'epimorphisme ou un isomorphisme de ${\rm Ens}$.

\smallskip

Consid\'erons alors le cas o\`u ${\mathcal E} = \widehat{\mathcal C}_J$ est le topos des faisceaux sur un site $({\mathcal C},J)$, reli\'e \`a $\widehat{\mathcal C}$ par les deux foncteurs adjoints $j^* : \widehat{\mathcal C} \to \widehat{\mathcal C}_J$ et $j_* : \widehat{\mathcal C}_J \hookrightarrow \widehat{\mathcal C}$.

\smallskip

Tout monomorphisme de $\widehat{\mathcal C}_J$
$$
u : F \longrightarrow G
$$
induit un monomorphisme de $\widehat{\mathcal C}$
$$
j_* F \longrightarrow j_* G \, .
$$
Soit $H$ la colimite du diagramme
$$
\xymatrix{
j_* F \ar[d] \ar[r] &j_* G \\
j_* G
}
$$
dans $\widehat{\mathcal C}$.

\smallskip

On observe que, $j_* F \to j_* G$ \'etant un monomorphisme de $\widehat{\mathcal C}$, le carr\'e commutatif
$$
\xymatrix{
j_* F \ar[d] \ar[r] &j_* G \ar[d]^-{p_1} \\
j_* G \ar[r]_{p_2} &H
}
$$
est cart\'esien.

\smallskip

Comme $j^*$ pr\'eserve les limites finies, le transform\'e de ce carr\'e par $j^*$
$$
\xymatrix{
F \ar[d]_u \ar[r]^u &G \ar[d]^-{j^* p_2} \\
G \ar[r]_{j^* p_1} &j^* H
}
$$
est un carr\'e cart\'esien de $\widehat{\mathcal C}_J$.

\smallskip

Si $u$ est un \'epimorphisme, l'\'egalit\'e
$$
j^* p_1 \circ u = j^* p_2 \circ u
$$
impose
$$
j^* p_1 = j^* p_2
$$
et donc le morphisme
$$
u : F \longrightarrow G
$$
est un isomorphisme.

\smallskip

C'est ce que l'on voulait. 
\end{demo}

\subsection{Effectivit\'e des relations d'\'equivalence dans les topos}\label{subsec332}

\medskip

Les notions de relation d'\'equivalence et de quotient par une relation d'\'equivalence ont un sens dans n'importe quelle cat\'egorie:

\begin{defn}\label{defIII32}

Soit ${\mathcal C}$ une cat\'egorie localement petite.

\begin{listeimarge}

\item Si $X$ est un objet de ${\mathcal C}$ dont le carr\'e $X \times X$ est repr\'esentable dans ${\mathcal C}$, une relation d'\'equivalence sur $X$ est un sous-objet
$$
(i_1 , i_2) : R \xhookrightarrow{ \ { \ } \ } X \times X
$$
tel que, pour tout objet $Y$ de ${\mathcal C}$, le sous-ensemble
$$
{\rm Hom} (Y,R) \subset {\rm Hom} (Y,X) \times {\rm Hom} (Y,X)
$$
satisfait les trois axiomes de

\medskip

$\left\{ \begin{matrix}
\bullet &\mbox{r\'eflexivit\'e: il contient la diagonale ${\rm Hom} (Y,X)$,} \hfill \\
{ \ } \\
\bullet &\mbox{sym\'etrie: il est respect\'e par l'\'echange des deux facteurs de ${\rm Hom} (Y,X) \times {\rm Hom} (Y,X)$,} \hfill \\
{ \ } \\
\bullet &\mbox{transitivit\'e: si trois \'el\'ements $x_1 , x_2 , x_3$ de ${\rm Hom} (Y,X)$ v\'erifient les conditions} \hfill \\
&\mbox{$(x_1 , x_2) , (x_2 , x_3) \in {\rm Hom} (Y,R)$, alors on a aussi $(x_1 , x_3) \in {\rm Hom} (Y,R)$.} \hfill
\end{matrix} \right.
$

\medskip

\item Un quotient d'une telle relation d'\'equivalence est un \'epimorphisme
$$
p : X \longrightarrow Q
$$
tel que le carr\'e
$$
\xymatrix{
R \ar[d]_-{i_2} \ar[r]^{i_1} &X \ar[d]^-p \\
X \ar[r]^p &Q
}
$$

\noindent soit \`a la fois commutatif, cocart\'esien et cart\'esien.
\end{listeimarge}
\end{defn}

\begin{remarksqed}
\begin{listeisansmarge}
\item Pour tout morphisme $X \to Q$ de ${\mathcal C}$ tel que les limites finies $X \times X$ et $X \times_Q X$ soient repr\'esentables,
$$
R = X \times_Q X \xhookrightarrow{ \ { \ } \ } X \times X
$$
est une relation d'\'equivalence de $X$.

\medskip

\item Les conditions de (ii) signifient que
$$
p \circ i_1 = p \circ i_2 \, ,
$$
que $Q$ est le co\'egalisateur
$$
Q = {\rm coeg} \left( \raisebox{.7ex}{\xymatrix{ R  \dar[r]^-{^{^{\mbox{\scriptsize$i_i$}}}}_-{i_2} &X}} \right),
$$
donc est uniquement d\'etermin\'e par $R$ \`a unique isomorphisme pr\`es, et que $R$ est uniquement d\'etermin\'e par $p : X \to Q$ comme le produit fibr\'e
$$
Q = X \times _Q X \, .
$$

\item Si une relation d'\'equivalence sur un objet $X$ d'une cat\'egorie localement petite admet un quotient $X \xrightarrow{ \ p \ } Q$, on dit qu'elle est effective. 

\end{listeisansmarge}
\end{remarksqed}

\bigskip

La propri\'et\'e d'effectivit\'e des relations d'\'equivalence s'\'etend elle aussi de la cat\'egorie ${\rm Ens}$ \`a tous les topos:

\begin{prop}\label{propIII33}
On a la propri\'et\'e suivante:

\begin{enumerate}

\item[(4)] Dans un topos ${\mathcal E}$, toute relation d'\'equivalence sur un objet $F$
$$
R \xhookrightarrow{ \ { \ } \ } F \times F
$$
est effective, c'est-\`a-dire d\'efinit un quotient $F \to Q$ qui s'inscrit dans un carr\'e cart\'esien et cocart\'esien:
$$
\xymatrix{
R \ar[d] \ar[r] &F \ar[d] \\
F \ar[r] &Q
}
$$

\end{enumerate}
\end{prop}

\begin{demo}

Dans la cat\'egorie ${\rm Ens}$ des ensembles, une relation d'\'equivalence $R \hookrightarrow F \times F$ d\'efinit un quotient $F \xrightarrow{ \ p \ } Q$. L'ensemble $Q$ est celui des classes d'\'equivalence de $F$ pour la relation $R$ et l'application $p$ associe \`a tout \'el\'ement de $F$ la classe \`a laquelle il appartient.

\smallskip

Cette construction est fonctorielle et donc toute relation d'\'equivalence
$$
R \xhookrightarrow{ \ { \ } \ } F \times F
$$
dans une cat\'egorie $\widehat{\mathcal C}$ de pr\'efaisceaux sur une cat\'egorie essentiellement petite ${\mathcal C}$ d\'efinit un quotient $Q$. C'est le pr\'efaisceau qui associe \`a tout objet $X$ de ${\mathcal C}$ le quotient $Q(X)$ de $F(X)$ par la relation d'\'equivalence $R(X) \subset F(X) \times F(X)$.

\smallskip

Si ${\mathcal E} = \widehat{\mathcal C}_J$ est le topos des faisceaux sur un site $({\mathcal C},J)$, toute relation d'\'equivalence sur un faisceau $F$
$$
R \xhookrightarrow{ \ { \ } \ } F \times F
$$
induit via le foncteur de plongement $j_* : \widehat{\mathcal C}_J \hookrightarrow \widehat{\mathcal C}$ une relation d'\'equivalence sur le pr\'efaisceau $j_* F$
$$
j_* R \xhookrightarrow{ \ { \ } \ } j_* F \times j_* F \, .
$$

Celle-ci est effective c'est-\`a-dire d\'efinit un carr\'e de $\widehat{\mathcal C}$
$$
\xymatrix{
j_* R \ar[d] \ar[r] &j_* F \ar[d]^-{p'} \\
j_* F \ar[r]^{p'} &Q'
}
$$
\`a la fois cart\'esien et cocart\'esien.

\smallskip

Comme $j^*$ pr\'eserve les colimites et les limites finies, il transforme ce carr\'e en un carr\'e cart\'esien et cocart\'esien
$$
\xymatrix{
R \ar[d] \ar[r] &F \ar[d]^-p \\
F \ar[r]^p &Q
}
$$
dans lequel $Q = j^* Q'$ et $p = j^* p'$.

\smallskip

Le cas d'un topos arbitraire s'en d\'eduit car la propri\'et\'e consid\'er\'ee est invariante par \'equivalence de cat\'egories. 

\end{demo}

\subsection{La correspondance entre objets quotients et relations d'\'equivalence}\label{subsec333}

\medskip

Les topos v\'erifient aussi la r\'eciproque de la propri\'et\'e (4):

\begin{prop}\label{propIII34}

On a la propri\'et\'e suivante:

\medskip
\begin{enumerate}
\item[(5)] Dans un topos ${\mathcal E}$, tout \'epimorphisme
$$
p : F \longrightarrow Q
$$
est le quotient de la relation d'\'equivalence
$$
R = F \times_Q F \xhookrightarrow{ \ { \ } \ } F \times F
$$
qu'il d\'efinit.
\end{enumerate}
\end{prop}

\begin{remarks}
\begin{listeisansmarge}
\item Autrement dit, pour tout \'epimorphisme $p : F \to Q$ dans un topos, le carr\'e cart\'esien
$$
F \times_Q F = \xymatrix{R \ar[d] \ar[r] &F \ar[d] \\
F \ar[r] &Q
}
$$
est aussi cocart\'esien.

\medskip

\item Les propri\'et\'es (4) et (5) signifient que, pour tout objet $F$ d'un topos ${\mathcal E}$, il y a correspondance bijective entre les objets quotients de $F$ et les relations d'\'equivalence sur $F$.

\medskip

\item Il en r\'esulte que tout morphisme $u : F \to G$ de ${\mathcal E}$ se factorise comme le compos\'e d'un \'epimorphisme $F \twoheadrightarrow {\rm Im} (u)$ suivi d'un monomorphisme ${\rm Im} (u) \hookrightarrow G$, o\`u ${\rm Im} (u)$ est le quotient de $F$ par la relation d'\'equivalence $F \times_G F$.

\smallskip

En effet, si le morphisme canonique ${\rm Im} (u) \to G$ n'\'etait pas un monomorphisme, il se factoriserait \`a travers un quotient non trivial de ${\rm Im} (u)$, donc un quotient de $F$ plus petit que ${\rm Im} (u)$. C'est impossible.
\end{listeisansmarge}
\end{remarks}

\begin{demo}

C'est \'evident si ${\mathcal E} = {\rm Ens}$.

\smallskip

Le cas des topos de pr\'efaisceaux ${\mathcal E} = \widehat{\mathcal C}$ s'en d\'eduit puisque les foncteurs d'\'evaluation en les objets $X$ de ${\mathcal C}$ respectent les limites et les colimites.

\smallskip

Si ${\mathcal E} = \widehat{\mathcal C}_J$ est le topos des faisceaux sur un site $({\mathcal C},J)$, la relation d'\'equivalence dans $\widehat{\mathcal C}_J$
$$
R = F \times_Q F \xhookrightarrow{ \ { \ } \ } F \times F
$$
induit par le foncteur de plongement $j_* : \widehat{\mathcal C}_J \hookrightarrow \widehat{\mathcal C}$ une relation d'\'equivalence dans $\widehat{\mathcal C}$
$$
j_* R = j_* F \times_{j_* Q} j_* F \xhookrightarrow{ \ { \ } \ } j_* F \times j_* F \, .
$$

Celle-ci d\'efinit un quotient $Q'$ dans $\widehat{\mathcal C}$ et le morphisme $j_* p : j_* F \to j_* Q$ se factorise canoniquement en
$$
Q' \longrightarrow j_* Q \, .
$$

En tout objet $X$ de ${\mathcal C}$, $Q'(X)$ est l'ensemble quotient de $F(X)$ par la relation d'\'equivalence $F(X) \times_{Q(X)} F(X)$ et donc l'application 
$$
Q'(X) \longrightarrow j_* Q(X)
$$
est injective.

\smallskip

Donc $Q' \to j_* Q$ est un monomorphisme de $\widehat{\mathcal C}$ et son transform\'e par $j^*$
$$
j^* Q' \longrightarrow Q
$$
est un monomorphisme de $\widehat{\mathcal C}_J$.

\smallskip

L'\'epimorphisme $F \to Q$ de $\widehat{\mathcal C}_J$ se factorise en
$$
F \longrightarrow j^* Q' \longrightarrow Q
$$
et donc $j^* Q' \to Q$ est aussi un \'epimorphisme.

\smallskip

Comme le topos ${\mathcal E}$ est balanc\'e, le morphisme de $\widehat{\mathcal C}_J$
$$
j^* Q' \longrightarrow Q
$$
est un isomorphisme.

\smallskip

Le carr\'e cart\'esien de $\widehat{\mathcal C}$
$$
\xymatrix{
j_* R \ar[d] \ar[r] &j_* F \ar[d] \\
j_* F \ar[r] &Q'
}
$$
\'etant aussi cocart\'esien, son transform\'e par $j^*$
$$
\xymatrix{
R \ar[d] \ar[r] &F \ar[d] \\
F \ar[r] &Q
}
$$
est un carr\'e cocart\'esien de $\widehat{\mathcal C}_J$.

\smallskip

Le cas d'un topos arbitraire en r\'esulte puisque la propri\'et\'e consid\'er\'ee est invariante par \'equivalence de cat\'egories. 

\end{demo}

\subsection{Application aux cat\'egories d'objets lin\'eaires d'un topos}\label{subsec334}

\medskip

Nous allons montrer que la correspondance entre relations d'\'equivalences et quotients d'un topos entra{\^\i}ne que les cat\'egories ${\mathcal M}od_A$ des $A$-modules d'un topos annel\'e $({\mathcal E},A)$ sont des cat\'egories ab\'eliennes au sens suivant:

\begin{defn}\label{defIII35}

Une cat\'egorie additive ${\mathcal A}$ est dite ``ab\'elienne'' si
\begin{enumerate}
\item[$\bullet$] elle poss\`ede des limites finies et des colimites finies arbitraires, autrement dit tout morphisme de ${\mathcal M}$
$$
M \xrightarrow{ \ u \ } N
$$
poss\`ede un ``noyau''
$$
\ker (u) = \varprojlim \left( \raisebox{.7ex}{\xymatrix{ M  \dar[r]^-{^{^{\mbox{\scriptsize$u$}}}}_-{0} &N}} \right)
$$
et un ``conoyau''
$$
{\rm coker} (u) = \varinjlim 
\left( \raisebox{.7ex}{\xymatrix{ M  \dar[r]^-{^{^{\mbox{\scriptsize$u$}}}}_-{0} &N}} \right),
$$
\item[$\bullet$] d\'efinissant l'image d'un tel morphisme $u : M \to N$ comme
$$
{\rm Im} (u) = {\rm coker} (\ker (u) \xhookrightarrow{ \ { \ } \ } M) \, ,
$$
le morphisme canonique
$$
{\rm Im} (u) \longrightarrow \ker (N \longrightarrow {\rm coker} (u))
$$
est un isomorphisme.
\end{enumerate}
\end{defn}

\begin{remarksqed}
\begin{listeisansmarge}
\item Ainsi, pour tout morphisme $M \xrightarrow{ \ u \ } N$ d'une cat\'egorie ab\'elienne, ${\rm Im} (u) = M / \ker (u)$ est un sous-objet de $N$ et $u$ se factorise canoniquement comme le compos\'e d'un \'epimorphisme $M \twoheadrightarrow {\rm Im} (u)$ suivi d'un monomorphisme ${\rm Im} (u) \hookrightarrow N$.

\medskip

\item Toute sous-cat\'egorie pleine d'une cat\'egorie ab\'elienne qui est stable par limites finies et colimites finies est encore une cat\'egorie ab\'elienne.

\medskip

\item Dans une cat\'egorie ab\'elienne, une suite de morphismes
$$
\cdots \longrightarrow M_{i-1} \xrightarrow{ \ u_i \ } M_i \xrightarrow{ \ u_{i+1} \ } M_{i+1} \longrightarrow \cdots
$$
est appel\'ee un complexe si, pour tout indice $i$,
$$
u_{i+1} \circ u_i = 0
$$
ou, ce qui revient au m\^eme,
$$
{\rm Im} (u_i) \subset \ker (u_{i+1}) \, .
$$
Elle est appel\'ee un complexe exact si, de plus, pour tout indice $i$
$$
{\rm Im} (u_i) = \ker (u_{i+1}) \, .
$$

\item Une suite longue est un complexe index\'e par ${\mathbb Z}$.

\smallskip

Une suite courte est un complexe de la forme
$$
0 \longrightarrow M \xrightarrow{ \ u \ } N \xrightarrow{ \ v \ } Q \longrightarrow 0 \, .
$$
Si un tel complexe est exact, on parle de suite exacte longue ou de suite exacte courte.

\medskip

\item Pour toute cat\'egorie ab\'elienne ${\mathcal A}$ et tout carquois $D$ [resp. toute cat\'egorie essentiellement petite ${\mathcal D}$], la cat\'egorie $D\mbox{-diag} ({\mathcal A})$ des $D$-diagrammes de ${\mathcal A}$ [resp. $[{\mathcal D} , {\mathcal A}]$ des foncteurs ${\mathcal D} \to {\mathcal A}$] est encore une cat\'egorie ab\'elienne.

\smallskip

Si $D$ est le carquois de cha{\^\i}ne
$$
\cdots \longrightarrow \underset{-1}{\bullet} \longrightarrow \underset{0}{\bullet} \longrightarrow \underset{1}{\bullet}
\longrightarrow \underset{2}{\bullet} \longrightarrow \cdots \, ,
$$
les suites longues forment une sous-cat\'egorie pleine $K({\mathcal A})$ de $D\mbox{-diag} ({\mathcal A})$ qui est stable par limites finies et colimites finies. Donc $K({\mathcal A})$ est aussi une cat\'egorie ab\'elienne.

\smallskip

Elle est munie de foncteurs dits ``de cohomologie''
$$
H^i : K({\mathcal A}) \longrightarrow {\mathcal A} \, ,
$$
$$
\begin{matrix}
\left( \cdots \longrightarrow M_{i-1} \xrightarrow{ \ u_i \ } M_i \xrightarrow{ \ u_{i+1} \ } M_{i+1} \longrightarrow \cdots \right) &\longmapsto &\ker (u_{i+1}) / {\rm Im} (u_i) \, . \\
\Vert &&\Vert \\
M_{\bullet} &&H^i (M_{\bullet})
\end{matrix}
$$

\item Deux morphismes de $K({\mathcal A})$
$$
\raisebox{.7ex}{\xymatrix{ M_{\bullet}  \dar[rrr]^-{^{^{\mbox{\scriptsize$u_{\bullet}$}}}}_-{v_{\bullet}} &&&N_{\bullet}}} 
$$
\vglue-3mm
$$
\Vert \qquad\qquad\qquad\qquad\qquad \Vert
$$
\vglue-3mm
$$
\left( \cdots \longrightarrow M_{i-1} \xrightarrow{ \ d \ } M_i \xrightarrow{ \ d \ } M_{i+1} \longrightarrow \cdots \right) \quad \left(\cdots \longrightarrow N_{i-1} \xrightarrow{ \ d \ } N_i \xrightarrow{ \ d \ } N_{i+1} \longrightarrow \cdots \right) 
$$
sont dit ``homotopes'' s'il existe une suite de morphismes de ${\mathcal A}$
$$
h_i : M_i \longrightarrow N_{i-1}
$$
tels que
$$
u_i - v_i = d \circ h_i + h_{i+1} \circ d \, , \qquad \forall \, i \, .
$$
La relation d'homotopie est une relation d'\'equivalence entre morphismes $M_{\bullet} \to N_{\bullet}$.

\smallskip

On observe que deux morphismes homotopes
$$
\raisebox{.7ex}{\xymatrix{ M_{\bullet} \dar[r] &N_{\bullet}}} 
$$
induisent les m\^emes morphismes en cohomologie
$$
H^i (M_{\bullet}) \longrightarrow H^i (N_{\bullet}) \, , \qquad i \in {\mathbb Z} \, .
$$

\item Un morphisme de $K({\mathcal A})$
$$
M_{\bullet} \xrightarrow{ \ u_{\bullet} \ } N_{\bullet}
$$
est appel\'e un ``quasi-isomorphisme'' s'il induit des isomorphismes en cohomologie
$$
H^i (M_{\bullet}) \xrightarrow{ \ \sim \ } H^i (N_{\bullet}) \, .
$$

Il suffit pour cela qu'il soit une ``\'equivalence d'homotopie'' au sens qu'existe un morphisme en sens inverse
$$
N_{\bullet} \xrightarrow{ \ v_{\bullet} \ } M_{\bullet}
$$
tel que $v_{\bullet} \circ u_{\bullet}$ soit homotope \`a ${\rm id}_{M_{\bullet}}$ et $u_{\bullet} \circ v_{\bullet}$ soit homotope \`a ${\rm id}_{N_{\bullet}}$. 

\end{listeisansmarge}
\end{remarksqed}


On a:

\begin{prop}\label{propIII36}

Pour tout topos annel\'e $({\mathcal E} , A)$, la cat\'egorie ${\mathcal M}od_A$ des $A$-modules de ${\mathcal E}$ est ab\'elienne.
\end{prop}

\begin{demo}

On sait d\'ej\`a que la cat\'egorie additive ${\mathcal M}od_A$ a des limites et colimites arbitraires, ce qui permet d'associer \`a tout morphisme $u : M \to N$
\begin{eqnarray}
\ker (u) &= &\varprojlim \left(\raisebox{.7ex}{\xymatrix{ M \dar[r]^-{^{^{\mbox{\scriptsize$u$}}}}_-{0} &N}}\right) ,\nonumber \\
{\rm coker} (u) &= &\varinjlim \left(\raisebox{.7ex}{\xymatrix{ M \dar[r]^-{^{^{\mbox{\scriptsize$u$}}}}_-{0} &N}}\right) , \nonumber \\
{\rm Im} (u) &= &{\rm coker} (\ker (u) \longrightarrow M) \, , \nonumber \\
{\rm coIm} (u) &= &\ker (N \longrightarrow {\rm coker} (u)) \, . \nonumber
\end{eqnarray}
Il s'agit de prouver que le morphisme canonique
$$
{\rm Im} (u) \longrightarrow {\rm coIm} (u)
$$
est un isomorphisme.

\smallskip

On observe qu'une relation d'\'equivalence sur un objet $M$ de ${\mathcal M}od_A$
$$
R \xhookrightarrow{ \ { \ } \ } M \times M
$$
est un sous-objet de $M \times M$ dans ${\mathcal M}od_A$ qui contient la diagonale.

\smallskip

Se donner une telle relation d'\'equivalence \'equivaut \`a se donner le sous-objet
$$
K = R \times_{M \times M} (0 \times M) \, .
$$
En effet, pour tous morphismes
$$
h_1 , h_2 : P \longrightarrow M \, ,
$$
on a
$$
(h_1 , h_2) \in {\rm Hom} (P,R)
$$
si et seulement si
$$
h_2 - h_1 \in {\rm Hom} (P,K) \, .
$$
D'autre part, un morphisme $M\to Q$ de ${\mathcal M}od_A$ est un \'epimorphisme si et seulement si c'est un \'epimorphisme de ${\mathcal E}$. ll y a donc correspondance bijective entre les sous-objets de $M$ et ses quotients.

\smallskip

Revenant au morphisme $u : M \to N$, on voit que le morphisme
$$
{\rm Im} (u) \longrightarrow N
$$
est un monomorphisme. S'il en \'etait autrement, $M \xrightarrow{ \ u \ } N$ se factoriserait \`a travers un quotient de $M$ plus petit que ${\rm Im} (u)$ et $M \xrightarrow{ \ u \ } N$ s'annulerait sur un sous-objet de $M$ plus grand que ${\rm Ker} (u)$.

\smallskip

Comme ${\rm Im} (u) \to N$ est un monomorphisme, il en est a fortiori de m\^eme de sa factorisation
$$
{\rm Im} (u) \longrightarrow {\rm coIm} (u) \, .
$$

De m\^eme, le morphisme
$$
M \longrightarrow {\rm coIm} (u)
$$
est un \'epimorphisme. S'il en \'etait autrement, $M \xrightarrow{ \ u \ } N$ se factoriserait \`a travers un sous-objet de $N$ plus petit que ${\rm coIm} (u)$ et il existerait un quotient $N'$ de $N$ plus grand que ${\rm coker} (u)$ tel que le compos\'e $M \xrightarrow{ \ u \ } N \longrightarrow N'$ s'annulerait.

\smallskip

Comme $M \to {\rm coIm} (u)$ est un \'epimorphisme, il en est a fortiori de m\^eme de sa factorisation
$$
{\rm Im} (u) \longrightarrow {\rm coIm} (u) \, .
$$

Ainsi, le morphisme
$$
{\rm Im} (u) \longrightarrow {\rm coIm} (u)
$$
est \`a la fois un monomorphisme et un \'epimorphisme dans ${\mathcal M}od_A$.

\smallskip

Comme le topos ${\mathcal E}$ est une cat\'egorie balanc\'ee, il en est de m\^eme de la cat\'egorie ${\mathcal M}od_A$ et le morphisme
$$
{\rm Im} (u) \longrightarrow {\rm coIm} (u)
$$
est un isomorphisme comme annonc\'e. 

\end{demo}

\section{Compatibilit\'es entre limites et colimites dans un topos}\label{sec34}

\subsection{Changement de base et colimites}\label{subsec341}

\medskip

Les topos h\'eritent de la cat\'egorie des ensembles que les produits fibr\'es respectent les colimites:

\begin{prop}\label{propIII41}

On a la propri\'et\'e suivante:

\medskip
\begin{enumerate}
\item[(6)] Dans un topos ${\mathcal E}$, les foncteurs de changement de base associ\'es aux morphismes $S' \to S$ de ${\mathcal E}$
$$
\begin{matrix}
\hfill {\mathcal E}/S &\longrightarrow &{\mathcal E}/S' \, , \hfill \\
(F \to S) &\longmapsto &(F \times_S S' \to S')
\end{matrix}
$$
respectent les colimites.
\end{enumerate}
\end{prop}

\bigskip

\begin{remark}

En particulier, pour tout objet $S$ d'un topos ${\mathcal E}$, le foncteur $S \times \bullet$ de produit avec $S$ respecte les sommes~$\amalg$.

\smallskip

On remarque que la propri\'et\'e duale (qui voudrait que les foncteurs $S \amalg \bullet$ respectent les produits) n'est absolument pas v\'erifi\'ee dans la cat\'egorie ${\rm Ens}$.

\smallskip

Cela signifie que la cat\'egorie duale ${\rm Ens}^{\rm op}$ de ${\rm Ens}$ n'est pas un topos.

\smallskip

En fait, le seul topos ${\mathcal E}$ dont la cat\'egorie duale ${\mathcal E}^{\rm op}$ est aussi un topos est le topos \`a un seul objet et un seul morphisme.
\end{remark}
\bigskip

\begin{demo}

Il s'agit de prouver que pour tout carquois $D$ et tout $D$-diagramme $F_{\bullet}$ d'objets de ${\mathcal E}/S$, le morphisme canonique
$$
\varinjlim_D \, (F_{\bullet} \times_S S') \longrightarrow (\varinjlim_D F_{\bullet}) \times_S S'
$$
est un isomorphisme.

Si ${\mathcal E} = {\rm Ens}$ est le topos des ensembles, cela r\'esulte de la description concr\`ete des limites et colimites dans ${\rm Ens}$ donn\'ee par le lemme \ref{lemI911}.

\smallskip

Le cas o\`u ${\mathcal E} = \widehat{\mathcal C}$ est le topos des pr\'efaisceaux sur une cat\'egorie essentiellement petite ${\mathcal C}$ s'en d\'eduit puisque les foncteurs d'\'evaluation en les objets $X$ de ${\mathcal C}$
$$
F \longmapsto F(X)
$$
respectent \`a la fois les limites et les colimites.

\smallskip

Le cas o\`u ${\mathcal E} = \widehat{\mathcal C}_J$  est le topos des faisceaux sur un site $({\mathcal C},J)$ se d\'eduit du cas de $\widehat{\mathcal C}$ gr\^ace aux deux foncteurs adjoints $j^* : \widehat{\mathcal C} \to \widehat{\mathcal C}_J$ et $j_* : \widehat{\mathcal C}_J \hookrightarrow \widehat{\mathcal C}$.

\smallskip

En effet, l'isomorphisme
$$
\varinjlim_D \left( j_* F_{\bullet} \times_{j_* S} j_* S' \right) \xrightarrow{ \ \sim \ } \left( \varinjlim_D j_* F_{\bullet} \right) \times_{j_* S} j_* S'
$$
est transform\'e par $j^*$ en un isomorphisme
$$
\varinjlim_D \left( F_{\bullet} \times_S S' \right) \xrightarrow{ \ \sim \ } \left( \varinjlim_D F_{\bullet} \right) \times_S S'
$$
puisque $j^*$ respecte les colimites arbitraires et les limites finies et que le compos\'e $j^* \circ j_*$ s'identifie au foncteur ${\rm id}$ de $\widehat{\mathcal C}_J$.

\smallskip

Enfin, le cas g\'en\'eral r\'esulte de ce que la propri\'et\'e consid\'er\'ee est invariante par \'equivalence de cat\'egories. 

\end{demo}

\subsection{Exactitude des colimites filtrantes}\label{subsec342}

\medskip

Voici encore une propri\'et\'e qui s'\'etend de la cat\'egorie des ensembles \`a tous les topos:

\begin{prop}\label{propIII42}

On a la propri\'et\'e suivante:

\medskip
\begin{enumerate}
\item[(7)] Pour tout topos ${\mathcal E}$ et toute petite cat\'egorie filtrante ${\mathcal D}$, le foncteur de colimite
$$
\varinjlim_{\mathcal D} \, [{\mathcal D} , {\mathcal E}] \longrightarrow {\mathcal E}
$$
respecte les limites finies.
\end{enumerate}
\end{prop}

\begin{demo}

Le cas o\`u ${\mathcal E} = {\rm Ens}$ est la cat\'egorie des ensembles est le contenu du lemme \ref{lemII54}.

\smallskip

Le cas des topos de pr\'efaisceaux ${\mathcal E} = \widehat{\mathcal C}$ s'en d\'eduit puisque les limites et colimites dans $\widehat{\mathcal C}$ se calculent composante par composante.

\smallskip

Le cas du topos ${\mathcal E} = \widehat{\mathcal C}_J$ des faisceaux sur un site $({\mathcal C},J)$ se d\'eduit de celui de $\widehat{\mathcal C}$ via les foncteurs adjoints $j^* : \widehat{\mathcal C} \to \widehat{\mathcal C}_J$ et $j_* : \widehat{\mathcal C}_J \hookrightarrow \widehat{\mathcal C}$. On utilise encore une fois la propri\'et\'e du foncteur $j^*$ de respecter les colimites arbitraires et les limites finies, et l'identification de $j^* \circ j_*$ avec le foncteur ${\rm id}$ de $\widehat{\mathcal C}_J$.

\smallskip

Enfin, le cas g\'en\'eral r\'esulte de ce que la propri\'et\'e consid\'er\'ee est invariante par \'equivalences de cat\'egories.

\end{demo}

\subsection{Caract\`ere disjoint des sommes}\label{subsec343}

On obtient toujours de la m\^eme fa\c con \`a partir du cas ensembliste:

\begin{prop}\label{propIII43}

On a la propri\'et\'e suivante:

\medskip
\begin{enumerate}
\item[(8)] Dans un topos ${\mathcal E}$, le carr\'e commutatif associ\'e \`a deux objets $F_1 , F_2$ et \`a l'objet initial $\emptyset$
$$
\xymatrix{
\emptyset \ar[d] \ar[r] &F_1 \ar[d] \\
F_2 \ar[r] &F_1 \amalg F_2
}
$$

est cart\'esien.
\end{enumerate}
\end{prop}

\begin{remark}

Cela implique que, plus g\'en\'eralement, pour toute famille $(F_i)_{i \in I}$ d'objets de ${\mathcal E}$ et tous indices $i_1 , i_2  \in I$ distincts, on a
$$
F_{i_1} \times_{\underset{i \in I}{\amalg} F_i} F_{i_2} = \emptyset \, .
$$
\end{remark}

\begin{demo}

Le cas de ${\rm Ens}$ est \'evident et implique celui des topos de pr\'efaisceaux $\widehat{\mathcal C}$.

\smallskip

Le passage de $\widehat{\mathcal C}$ \`a $\widehat{\mathcal C}_J$ utilise le fait que le foncteur $j^* : \widehat{\mathcal C} \to \widehat{\mathcal C}_J$ respecte les colimites et les produits fibr\'es ainsi que l'identification de $j^* \circ j_*$ avec ${\rm id} : \widehat{\mathcal C}_J \to \widehat{\mathcal C}_J$.

\smallskip

Le cas g\'en\'eral r\'esulte de ce que la propri\'et\'e consid\'er\'ee est invariante par \'equivalence de cat\'egories. 

\end{demo}

\subsection{Le cas des cat\'egories de modules dans un topos}\label{subsec344}

\medskip

On observe d'abord que la propri\'et\'e (8) est v\'erifi\'ee par n'importe qu'elle cat\'egorie semi-additive:

\begin{lem}\label{lemIII44}

Soit ${\mathcal A}$ une cat\'egorie semi-additive.

\smallskip

Alors, pour tous objets $M_1$ et $M_2$ de ${\mathcal A}$, le carr\'e commutatif
$$
\xymatrix{
0 \ar[d] \ar[r] &M_1 \ar[d] \\
M_2 \ar[r] &M_1 \oplus M_2 = M_1 \times M_2
}
$$

est cart\'esien.
\end{lem}

\begin{demo}

Pour tout objet $M$ de ${\mathcal A}$, les compos\'es de deux morphismes $u_1 : M \to M_1$ et $u_2 : M \to M_2$ avec $M_1 \to M_1 \oplus M_2$ et $M_2 \to M_1 \oplus M_2$ sont les paires
$$
(u_1,0) , (0,u_2) \in {\rm Hom} (M,M_1 \times M_2) = {\rm Hom} (M,M_1) \times {\rm Hom} (M,M_2) \, .
$$
Elles sont \'egales si et seulement si $u_1 = 0$ et $u_2 = 0$. 

\end{demo}

\bigskip

D'autre part, les cat\'egories de modules d'un topos annel\'e poss\`edent toujours la propri\'et\'e (7):

\begin{prop}\label{propIII45}

Soit $({\mathcal E},A)$ un topos annel\'e.

\smallskip

Soit ${\mathcal M}od_A$ la cat\'egorie des $A$-modules de ${\mathcal E}$.

\smallskip

Alors, pour toute petite cat\'egorie filtrante ${\mathcal D}$, le foncteur de colimite
$$
\varinjlim_{\mathcal D} \, [{\mathcal D} , {\mathcal M}od_A] \longrightarrow {\mathcal M}od_A
$$
respecte les limites finies.
\end{prop}

\begin{demosansqed}

On sait d\'ej\`a que la cat\'egorie ${\mathcal M}od_A$ poss\`ede des limites et des colimites arbitraires.

\smallskip

Compte-tenu de la proposition \ref{propIII42}, la propri\'et\'e consid\'er\'ee de ${\mathcal M}od_A$ r\'esulte du lemme suivant:
\end{demosansqed}

\begin{lem}\label{lemIII46}

Soit $({\mathcal E},A)$ un topos annel\'e.

\smallskip

Soit ${\mathcal M}od_A$ la cat\'egorie des $A$-modules de ${\mathcal E}$.

\smallskip

Alors le foncteur d'oubli de la structure de module
$$
{\mathcal M}od_A \longrightarrow {\mathcal E}
$$
poss\`ede les propri\'et\'es suivantes:

\begin{listeimarge}

\item Il respecte les limites.

\medskip

\item Il respecte les colimites filtrantes.

\medskip

\item Un morphisme de ${\mathcal M}od_A$
$$
M_1 \longrightarrow M_2
$$
est un isomorphisme [resp. un monomorphisme, resp. un \'epimorphisme] si et seulement si il en est ainsi de son image dans ${\mathcal E}$ par le foncteur d'oubli.
\end{listeimarge}
\end{lem}

\begin{remark}

Ce lemme vaut aussi en rempla\c cant la cat\'egorie ${\mathcal M}od_A$ par celle des mono{\"\i}des [resp. groupes, resp. mono{\"\i}des commutatifs, resp. groupes ab\'eliens, resp. anneaux, resp. anneaux commutatifs] internes d'un topos ${\mathcal E}$.

\smallskip

Il en va de m\^eme de la proposition \ref{propIII45}.

\end{remark}

\begin{demolem}
\begin{listeisansmarge}
\item[(i)] r\'esulte de ce que les foncteurs d'oubli
$$
\begin{matrix}
\hfill {\mathcal M}od_A &\longrightarrow &{\mathcal A}b_{\mathcal E} \\
\mbox{et} \qquad\qquad {\mathcal A}b_{\mathcal E} &\longrightarrow &{\mathcal E} \hfill
\end{matrix}
$$
respectent les limites d'apr\`es le corollaire \ref{corIII28} (i) et le corollaire \ref{corIII23} (i).

\medskip

\item[(iii)] Dans le cas du topos ${\mathcal E} = {\rm Ens}$ des ensembles, cela r\'esulte de ce qu'un morphisme de modules sur un anneau $A$
$$
M_1 \longrightarrow M_2
$$
est un isomorphisme [resp. un monomorphisme, resp. un \'epimorphisme] si et seulement si il est bijectif [resp. injectif, resp. surjectif].

\smallskip

Le cas des topos ${\mathcal E} = \widehat{\mathcal C}$ de pr\'efaisceaux s'en d\'eduit puisqu'un morphisme de ${\mathcal M}od_A$ ou de ${\mathcal E}$
$$
M_1 \longrightarrow M_2
$$
est un isomorphisme [resp. un monomorphisme, resp. un \'epimorphisme] si et seulement si ses composantes
$$
M_1 (X) \longrightarrow M_2 (X) \, , \qquad X \in {\rm Ob} ({\mathcal C}) \, ,
$$
le sont.

\smallskip

Si ${\mathcal E} = \widehat{\mathcal C}_J$ est le topos des faisceaux sur un site $({\mathcal C},J)$, un morphisme de ${\mathcal M}od_A = {\mathcal M}od_{{\mathcal C},J,A}$ ou de ${\mathcal E}$
$$
M_1 \longrightarrow M_2
$$
est un isomorphisme [resp. un monomorphisme] si et seulement si son image
$$
j_* M_1 \longrightarrow j_* M_2
$$
par le foncteur $j_* : {\mathcal M}od_{{\mathcal C},J,A} \hookrightarrow {\mathcal M}od_{{\mathcal C},A}$ ou $j_* : \widehat{\mathcal C}_J \hookrightarrow \widehat{\mathcal C}$ est un isomorphisme [resp. un monomorphisme].

\smallskip

Enfin, un morphisme de ${\mathcal M}od_{{\mathcal C},J,A}$
$$
u : M_1 \longrightarrow M_2
$$
est un \'epimorphisme si et seulement si l'image
$$
j^* \, M
$$
par le foncteur de faisceautisation $j^* : {\mathcal M}od_{{\mathcal C},A} \to {\mathcal M}od_{{\mathcal C},J,A}$ du conoyau dans ${\mathcal M}od_{{\mathcal C},A}$
\begin{eqnarray}
M &= &{\rm coker} (j_* u : j_* M_1 \longrightarrow j_* M_2)  \nonumber \\
&= &\left[{\rm Ob} ({\mathcal C}) \ni X \longmapsto {\rm coker} (M_1 (X) \longrightarrow M_2(X))\right] \nonumber
\end{eqnarray}
est $0$.

\smallskip

C'est \'equivalent \`a demander que pour tout objet $X$ de ${\mathcal C}$ et tout \'el\'ement $m_2 \in M_2 (X)$, $X$ admette un crible $J$-couvrant $S$ tel que pour tout \'el\'ement $(U \xrightarrow{ \ p \ } X) \in S$, la restriction
$$
M_2 (p) (m_2) \in M_2 (U)
$$
soit dans l'image de $M_1 (U)$.

\smallskip

Cela revient encore \`a demander que $M_1 \xrightarrow{ \ u \ } M_2$ induise un \'epimorphisme de ${\mathcal E} = \widehat{\mathcal C}_J$.

\smallskip

Le cas d'un topos arbitraire ${\mathcal E}$ s'en d\'eduit puisque la propri\'et\'e consid\'er\'ee est invariante par \'equivalence de cat\'egories.

\medskip

\item[(ii)] Le cas g\'en\'eral se r\'eduit \`a celui des topos de faisceaux ${\mathcal E} = \widehat{\mathcal C}_J$ puisque la propri\'et\'e consid\'er\'ee est invariante par \'equivalence de cat\'egories.

\smallskip

Le cas du topos ${\mathcal E} = \widehat{\mathcal C}_J$ des faisceaux sur un site $({\mathcal C},J)$ se r\'eduit \`a celui de $\widehat{\mathcal C}$ via les foncteurs adjoints
$$
j^* : {\mathcal M}od_{{\mathcal C},A} \longrightarrow {\mathcal M}od_{{\mathcal C},J,A} \quad \mbox{et} \quad j_* : {\mathcal M}od_{{\mathcal C},J,A} \longrightarrow {\mathcal M}od_{{\mathcal C},A}
$$
puisque $j^*$ respecte les colimites et que le compos\'e $j^* \circ j_*$ s'identifie au foncteur ${\rm id}$ de ${\mathcal M}od_{{\mathcal C},J,A}$.

\smallskip

Le cas des topos $\widehat{\mathcal C}$ de pr\'efaisceaux se r\'eduit \`a celui de ${\mathcal E} = {\rm Ens}$ puisque les colimites dans ${\mathcal M}od_{{\mathcal C},A}$ se calculent composante par composante.

\smallskip

Enfin, si ${\mathcal E} = {\rm Ens}$ est le topos des ensembles, consid\'erons le foncteur d'oubli des structures de $A$-module
$$
O : {\rm Mod}_A \longrightarrow {\rm Ens} \, .
$$
Si ${\mathcal D}$ est une petite cat\'egorie filtrante et $M_{\bullet} : {\mathcal D} \to {\rm Mod}_A$ un foncteur, consid\'erons la colimite 
$$
M = \varinjlim_{\mathcal D} \, O \, (M_{\bullet}) \qquad \mbox{dans} \quad {\rm Ens} \, .
$$

Comme ${\mathcal D}$ est filtrante, il existe pour toute famille finie
$$
\overline m_1 , \cdots , \overline m_k
$$
d'\'el\'ements de $M$ un objet $d'$ de ${\mathcal D}$ et des rel\`evements
$$
m'_1 , \cdots , m'_k
$$
de $\overline m_1 , \cdots , \overline m_k$ dans $M_{d'}$.

\smallskip

De plus, si $d''$ est un objet de ${\mathcal D}$ tel que $\overline m_1 , \cdots , \overline m_k$ admettent aussi des rel\`evements
$$
m''_1 , \cdots , m''_k \qquad \mbox{dans} \quad M_{d''} \, ,
$$
il existe deux morphismes de ${\mathcal D}$
$$
\xymatrix{
d' \ar[rd] \\
&d \\
d'' \ar[ru]
}
$$
qui transforment les deux familles $m'_1 , \cdots , m'_k$ et $m''_1 , \cdots , m''_k$ en la m\^eme famille
$$
m_1 , \cdots , m_k \qquad \mbox{dans} \quad M_d \, .
$$
On en d\'eduit que $M$ h\'erite des $M_d$ une structure naturelle de $A$-module qui en fait une colimite de $M_{\bullet}$ dans ${\rm Mod}_A$.

\smallskip

Cela termine la d\'emonstration du lemme, donc aussi de la proposition \ref{propIII45}. 

\end{listeisansmarge}
\end{demolem}

\bigskip

En revanche, une cat\'egorie de modules d'un topos annel\'e ou plus g\'en\'eralement une cat\'egorie semi-additive ne v\'erifie jamais la propri\'et\'e (6), \`a moins d'\^etre triviale:

\begin{lem}\label{lemIII47}

Soit ${\mathcal A}$ une cat\'egorie semi-additive dans laquelle les foncteurs de produits fibr\'es
$$
N \longmapsto N \times_M M'
$$
sont bien d\'efinis et respectent les sommes $\oplus$.

\smallskip

Alors ${\mathcal A}$ est triviale au sens que tous ses objets sont $0$.
\end{lem}

\begin{remark}

Il r\'esulte de ce lemme et de la proposition \ref{propIII41} qu'une cat\'egorie semi-additive non triviale n'est jamais un topos.

\end{remark}

\begin{demo}

Consid\'erons un objet $M$ de ${\mathcal A}$, le morphisme diagonal
$$
M \xrightarrow{ \ \Delta \ } M \times M = M \oplus M
$$
et les deux morphismes canoniques $i_1 : M \to M \oplus M$ et $i_2 : M \to M \oplus M$.

\smallskip

Les produits fibr\'es de $i_1$ ou $i_2$ avec $\Delta$ sont $0$.

\smallskip

En revanche, le produit fibr\'e de $i_1 \oplus i_2 = {\rm id} : M \oplus M \to M \oplus M$ avec $\Delta$ est $M$.

\smallskip

Si donc le foncteur $M \times_{\Delta , M \oplus M} \bullet$ respecte les sommes, on a $M=0$. 

\end{demo}

\section{Exponentielles dans les topos}\label{sec35}

\subsection{La notion g\'en\'erale d'exponentielle dans une cat\'egorie}\label{subsec351}

\medskip

On pose la d\'efinition naturelle suivante:

\begin{defn}\label{defIII51}

Soit ${\mathcal C}$ une cat\'egorie localement petite.

\smallskip

Soit $X$ un objet de ${\mathcal C}$ qui est ``carrable'' au sens qu'il admet un produit $X \times Y$ avec tout objet $Y$ de ${\mathcal C}$.

\smallskip

Alors l'objet $X$ est dit ``exponentiable'' si le foncteur
$$
\begin{matrix}
{\mathcal C} &\longrightarrow &{\mathcal C} \, , \hfill \\
Y &\longmapsto &X \times Y
\end{matrix}
$$
admet un adjoint \`a droite
$$
\begin{matrix}
{\mathcal C} &\longrightarrow &{\mathcal C} \, , \hfill \\
Z &\longmapsto &Z^X = {\mathcal H}om (X,Z) \, .
\end{matrix}
$$
\end{defn}

\begin{remarksqed}
\begin{listeisansmarge}
\item Par d\'efinition de l'adjonction, si $X$ est un objet exponentiable de ${\mathcal C}$, se donner un morphisme de ${\mathcal C}$
$$
X \times Y \longrightarrow Z
$$
\'equivaut \`a se donner un morphisme
$$
Y \longrightarrow Z^X \, .
$$

\item Si ${\mathcal C}$ est la cat\'egorie ${\rm Ens}$ des ensembles, tout objet $X$ de ${\mathcal C}$ est exponentiable car le foncteur
$$
Y \longmapsto X \times Y
$$
admet pour adjoint \`a droite
$$
Z \longmapsto \prod_{x \in X} Z = {\rm Hom} (X,Z) \, .
$$
C'est pourquoi dans le cas g\'en\'eral on emploie les notations
$$
Z^X = {\mathcal H}om (X,Z)
$$
et on dit qu'un tel objet $X$ est ``exponentiable''.

\smallskip

La notation $Z^X$ se lit ``$Z$ exponentielle $X$'' et la notation ${\mathcal H}om (X,Z)$ se lit ``le ${\mathcal H}om$ interne de $X$ dans $Z$'', l'adjectif ``interne'' renvoyant au fait que ${\mathcal H}om (X,Z)$ est un objet de la cat\'egorie ${\mathcal C}$.

\medskip

\item Tout morphisme
$$
X' \xrightarrow{ \ u \ } X
$$
entre deux objets exponentiables $X$ et $X'$ de ${\mathcal C}$ d\'efinit par adjonction une transformation naturelle du foncteur
$$
Z \longmapsto Z^X
$$
vers le foncteur
$$
Z \longmapsto Z^{X'} \, .
$$
Ainsi, on a pour tout morphisme de ${\mathcal C}$ 
$$
v : Z_1 \longrightarrow Z_2
$$
un carr\'e commutatif induit par $u$ et $v$:
$$
\xymatrix{
Z_1^X \ar[d] \ar[r] &Z_2^X \ar[d] \\
Z_1^{X'} \ar[r] &Z_2^{X'}
}
$$

\end{listeisansmarge}
\end{remarksqed}

\bigskip

Nous allons montrer un peu plus loin que tous les objets d'un topos sont exponentiables et donner une description concr\`ete de leurs exponentielles.

\subsection{Repr\'esentation des faisceaux et caract\'erisation des cribles couvrants}\label{subsec352}

\medskip

Afin de montrer que les objets des topos sont toujours exponentiables, on a besoin de deux lemmes g\'en\'eraux qui sont utiles dans de nombreux contextes.

\smallskip

Voici le premier qui exprime tout faisceau sur un site en termes d'objets repr\'esentables:

\begin{lem}\label{lemIII52}

Soit $({\mathcal C},J)$ un site, et soit $\ell$ le foncteur canonique
$$
\ell : {\mathcal C} \xrightarrow{ \ y \ } \widehat{\mathcal C} \xrightarrow{ \ j^* \ } \widehat{\mathcal C}_J
$$
compos\'e du foncteur de Yoneda $y : {\mathcal C} \to \widehat{\mathcal C}$ et du foncteur de faisceautisation $j^* : \widehat{\mathcal C} \to \widehat{\mathcal C}_J$.

\smallskip

Alors tout faisceau $F$ sur $({\mathcal C} , J)$ s'\'ecrit dans $\widehat{\mathcal C}_J$ comme la colimite
$$
F = \varinjlim_{(X,x) \in \int\!\!F} \ell (X)
$$
calcul\'ee sur la cat\'egorie relative
$$
{\mathcal C} / F = \textstyle{\int\!\!F}
$$
dont les objets sont les paires
$$
(X,x)
$$
constitu\'ees d'un objet $X$ de ${\mathcal C}$ et d'un \'el\'ement
$$
x \in {\rm Hom} (\ell (X) , F) = F(X)
$$
et dont les morphismes
$$
(X_1 , x_1),\longrightarrow (X_2 , x_2)
$$
sont les morphismes de ${\mathcal C}$
$$
u : X_1 \longrightarrow X_2
$$
tels que commute le triangle
$$
\xymatrix{
\ell (X_1) \ar[rd]_{x_1} \ar[rr]^{\ell (u)} &&\ell (X_2 \ar[ld]^{x_2}) \\
&F
}
$$
c'est-\`a-dire $F(u)(x_2) = x_1$.
\end{lem}

\begin{remark}

La cat\'egorie
$$
\textstyle{\int\!\!F} = {\mathcal C}/F
$$
est appel\'ee la ``cat\'egorie des \'el\'ements de $F$''.

\smallskip

Elle est essentiellement petite, comme ${\mathcal C}$, ce qui donne un sens au foncteur de colimite
$$
\varinjlim_{(X,x) \in \int\!\!F} \, .
$$
\end{remark}

\begin{demo}

Le foncteur de faisceautisation $j^*$ respecte les colimites arbitraires donc il suffit de prouver que si $F$ est un pr\'efaisceau sur ${\mathcal C}$, on a dans $\widehat{\mathcal C}$ une identification
$$
F = \varinjlim_{(X,x) \in \int\!\!F} y(X) \, .
$$

En effet, pour tout pr\'efaisceau $G$ sur ${\mathcal C}$, se donner un morphisme de $\widehat{\mathcal C}$
$$
F \longrightarrow G
$$
revient \`a se donner une famille d'applications
$$
F(X) \longrightarrow G(X) = {\rm Hom} (y(X),G)
$$
telle que, pour tout morphisme de ${\mathcal C}$
$$
u : X_1 \longrightarrow X_2 \, ,
$$
commute le carr\'e:
$$
\xymatrix{
F(X_2) \ar[d]_{F(u)} \ar[r] &G(X_2) = {\rm Hom} (y(X_2),G) \ar[d] \\
F(X_1) \ar[r] &G(X_1) = {\rm Hom} (y(X_1),G)
}
$$

Cela revient \`a demander que pour tous \'el\'ements $x_1 \in F(X_1)$, $x_2 \in F(X_2)$ v\'erifiant $x_1 = F(x_2)$, le triangle associ\'e
$$
\xymatrix{
y(X_2) \ar[rd] \ar[rr]^{y(u)} &&y(X_1) \ar[ld] \\
&G
}
$$
commute. 

\end{demo}

\bigskip

Voici le second lemme qui donne plusieurs caract\'erisations des cribles couvrants:

\begin{lem}\label{lemIII53}

Soient $({\mathcal C},J)$ un site et $\ell : {\mathcal C} \to \widehat{\mathcal C}_J$ le foncteur canonique.

\smallskip

Alors, pour toute famille de morphismes $(U_i \to X)_{i \in I}$ de ${\mathcal C}$ engendrant un crible $S$, les propri\'et\'es suivantes sont \'equivalentes:

\begin{enumerate}

\item[(A)] Cette famille $(U_i \to X)_{i \in I}$ est couvrante, autrement dit le crible $S$ est \'el\'ement de $J(X)$.

\medskip

\item[(B)] On a dans $\widehat{\mathcal C}_J$ la formule
$$
\ell (X) = \varinjlim_{(U \to X)\in S} \ell (U) \, .
$$

\medskip

\item[(C)] La famille de morphismes de $\widehat{\mathcal C}_J$
$$
(\ell (U_i) \longrightarrow \ell (X))_{i \in I}
$$
est globalement \'epimorphique au sens que, pour tout objet $F$ de $\widehat{\mathcal C}_J$, l'application
$$
{\rm Hom} (\ell (X),F) \longrightarrow \prod_{i \in I} {\rm Hom} (\ell (U_i) , F)
$$
est injective.

\medskip

\item[(D)] L'objet $\ell (X)$ de $\widehat{\mathcal C}_J$ est colimite du diagramme constitu\'e par les objets $\ell (U_i)$, $i \in I$, et les objets $\ell (U')$ munis des morphismes $\ell (U') \to \ell (U_i)$ et $\ell (U') \to \ell (U_j)$ associ\'es aux carr\'es commutatifs de ${\mathcal C}$:
$$
\xymatrix{
U' \ar[d] \ar[r] &U_i \ar[d] \\
U_j \ar[r] &X
}
$$

\end{enumerate}
\end{lem}

\begin{remark}

L'\'equivalence de (A) et (C) permet de retrouver la topologie $J$ de ${\mathcal C}$ \`a partir de la structure cat\'egorique de $\widehat{\mathcal C}_J$, via le foncteur canonique $\ell : {\mathcal C} \to \widehat{\mathcal C}_J$.

\end{remark}

\begin{demo}

Pour tout faisceau $F$ sur $({\mathcal C},J)$ et tout objet $X$ de ${\mathcal C}$, on a
$$
F(X) = {\rm Hom} (y(X), j_* F) = {\rm Hom} (\ell (X),F) \, .
$$

Si $S$ est un crible $J$-couvrant d'un objet $X$ de ${\mathcal C}$, la formule v\'erifi\'ee par d\'efinition par tout faisceau $F$
$$
F(X) = \varprojlim_{(U \to X) \in S} F(U)
$$
se r\'e\'ecrit donc
$$
{\rm Hom} (\ell (X),F) = \varprojlim_{(U \to X) \in S} {\rm Hom} (\ell (U),F) \, .
$$

Elle signifie que l'on a dans la cat\'egorie $\widehat{\mathcal C}_J$
$$
\ell (X) = \varinjlim_{(U\to X) \in S} \ell (U) \, .
$$

Cela montre que (A) implique (B).

\smallskip

Il est \'evident que (B) implique (C).

\smallskip

Si (C) est v\'erifi\'ee, il r\'esulte de la proposition \ref{propIII34} que $\ell (X)$ est la colimite dans $\widehat{\mathcal C}_J$ du diagramme constitu\'e par les $\ell (U_i)$ et les $\ell (U_i) \times_{\ell (X)} \ell (U_j)$ reli\'es par les projections
$$
\ell (U_i) \times_{\ell (X)} \ell (U_j) \longrightarrow \ell (U_i) \quad \mbox{et} \quad \ell (U_i) \times_{\ell (X)} \ell (U_j) \longrightarrow \ell (U_j) \, .
$$

Or, pour tous indices $i,j$, on a
$$
\ell (U_i) \times_{\ell (X)} \ell (U_j) = j^* (y(U_i) \times _{y(X)} y (U_j)) \, .
$$

Donc $\ell (U_i) \times_{\ell (X)} \ell (U_j)$ admet une famille globalement \'epimorphique de fl\`eches
$$
\ell (U') \longrightarrow \ell (U_i) \times_{\ell (X)} \ell (U_j)
$$
induites par des objets $U'$ de ${\mathcal C}$ munis de fl\`eches de $\widehat{\mathcal C}$
$$
y(U') \longrightarrow y(U_i) \times_{y(X)} y (U_j)
$$
c'est-\`a-dire par des carr\'es commutatifs de ${\mathcal C}$:
$$
\xymatrix{
U' \ar[d] \ar[r] &U_i \ar[d] \\
U_j \ar[r] &X
}
$$
Cela montre que (C) implique (D).

\smallskip

Enfin, (D) implique (A) puisque, pour tout faisceau $F$ sur $({\mathcal C},J)$, l'identit\'e
$$
{\rm Hom} (\ell (X),F) = {\rm eg} \Biggl( \prod_i {\rm Hom} (\ell (U_i) , F) \rightrightarrows \prod_{\begin{pmatrix} U' &\to &U_i \\ \downarrow &&\downarrow \\ U_j &\to &X \end{pmatrix}} {\rm Hom} (\ell (U'),F) \Biggl)
$$
se r\'e\'ecrit
$$
F(X) = {\rm eg} \Biggl( \prod_i F(U_i) \rightrightarrows \prod_{\begin{pmatrix} U' &\to &U_i \\ \downarrow &&\downarrow \\ U_j &\to &X \end{pmatrix}} F(U') \Biggl).
$$

Cela ach\`eve la d\'emonstration du lemme. 

\end{demo}

\subsection{Les exponentielles dans un topos}\label{subsec353}

\medskip

Nous pouvons maintenant d\'emontrer:

\begin{prop}\label{propIII54}
\begin{listeimarge}
\item Dans un topos, tout objet est exponentiable:

\medskip
\begin{enumerate}
\item[(9)] Pour tout objet $F$ d'un topos ${\mathcal E}$, le foncteur
$$
\begin{matrix}
{\mathcal E} &\longrightarrow &{\mathcal E} \, , \hfill \\
G &\longmapsto &F \times G
\end{matrix}
$$
admet un adjoint \`a droite
$$
\begin{matrix}
{\mathcal E} &\longrightarrow &{\mathcal E} \, , \hfill \\
H &\longmapsto &{\mathcal H}om (F,H) \, .
\end{matrix}
$$
\end{enumerate}
\item De plus, si ${\mathcal E} = \widehat{\mathcal C}_J$ est le topos des faisceaux sur un site $({\mathcal C},J)$, le faisceau
$$
{\mathcal H}om (F,H)
$$
associe \`a tout objet $X$ de ${\mathcal C}$ l'ensemble
$$
{\rm Hom} (F_X , H_X)
$$
des morphismes entre les restrictions $F_X$ et $H_X$ des faisceaux $F$ et $H$ \`a la cat\'egorie relative ${\mathcal C}/X$.
\end{listeimarge}
\end{prop}

\begin{remark}

\noindent (ii) signifie en particulier que pour tous faisceaux $F,H$ sur $({\mathcal C} , J)$, le pr\'efaisceau
$$
\begin{matrix}
{\mathcal C}^{\rm op} &\longrightarrow &{\rm Ens} \, , \hfill \\
X &\longmapsto &{\rm Hom} (F_X,H_X)
\end{matrix}
$$
est un faisceau pour la topologie $J$.

\end{remark}

\begin{demo}
\begin{listeisansmarge}
\item[(ii)] Pour tout objet $X$ de ${\mathcal C}$, se donner un morphisme de faisceaux
$$
F \times \ell (X) \longrightarrow H
$$
\'equivaut \`a se donner un morphisme de pr\'efaisceaux
$$
F \times y(X) \longrightarrow H \, .
$$
Un tel morphisme consiste \`a associer \`a tout morphisme de ${\mathcal C}$
$$
Y \xrightarrow{ \ v \ } X
$$
une application
$$
F(Y) \xrightarrow{ \ f_v \ } H(Y)
$$
de telle sorte que tout triangle commutatif
$$
\xymatrix{
Y_1 \ar[rd]_{v_1} \ar[rr]^u &&Y_2 \ar[ld]^{v_2} \\
&X
}
$$
induise un carr\'e commutatif:
$$
\xymatrix{
F(Y_2) \ar[d]_{F(u)} \ar[r]^{f_{v_2}} &H (Y_2) \ar[d]^{H(u)} \\
F(Y_1) \ar[r]^{f_{v_1}} &H(Y_1)
}
$$

Autrement dit, les morphismes de faisceaux
$$
F \times \ell (X) \longrightarrow H
$$
s'identifient aux morphismes de pr\'efaisceaux sur ${\mathcal C}/X$
$$
F_X \longrightarrow H_X \, .
$$

Le pr\'efaisceau
$$
\begin{matrix}
{\mathcal H}om (F,H) : {\mathcal C}^{\rm op} &\longrightarrow &{\rm Ens} \hfill \\
\hfill X &\longmapsto &{\rm Hom} (F_X , H_X) = {\rm Hom} (F \times \ell (X),H)
\end{matrix}
$$
est un $J$-faisceau. En effet, si $S$ est un crible couvrant d'un objet $X$ de ${\mathcal C}$, on a d'apr\`es le lemme \ref{lemIII53}
$$
\ell (X) = \varinjlim_{(U \to X) \in S} \ell (U)
$$
puis d'apr\`es la proposition \ref{propIII41}
$$
F \times \ell (X) = \varinjlim_{(U \to X) \in S} F \times \ell (U)
$$
et donc
$$
{\mathcal H}om (F,H)(X) = \varprojlim_{(U \to X) \in S} {\mathcal H}om (F,H)(U) \, .
$$

Enfin, tout faisceau $G$ s'\'ecrit d'apr\`es le lemme \ref{lemIII52}
$$
G = \varinjlim_{(X,x) \in \int\!\!G} \ell (X) \, .
$$
On en d\'eduit d'apr\`es la proposition \ref{propIII41}
$$
F \times G = \varinjlim_{(X,x) \in \int\!\!G} F \times \ell (X)
$$
et donc
\begin{eqnarray}
{\rm Hom} (F \times G , H) &= &\varprojlim_{(X,x) \in \int\!\!G} {\rm Hom} (F \times \ell (X),H) \nonumber \\
&= &\varprojlim_{(X,x) \in \int\!\!G} {\mathcal H}om (F,H)(X) \nonumber \\
&= &\varprojlim_{(X,x) \in \int\!\!G} {\rm Hom} (\ell (X),{\mathcal H}om (F,H)) \nonumber \\
&= &{\rm Hom} (G,{\mathcal H}om (F,H)) \, . \nonumber
\end{eqnarray}
Cela prouve que le foncteur
$$
H \longmapsto {\mathcal H}om (F,H)
$$
est adjoint \`a droite de
$$
G \longmapsto F \times G \, .
$$

\item[(i)] r\'esulte de (ii) puisque la propri\'et\'e consid\'er\'ee est invariante par \'equivalence de cat\'egories. 

\end{listeisansmarge}
\end{demo}

\subsection{Les faisceaux de morphismes lin\'eaires entre faisceaux lin\'eaires}\label{subsec354}

\medskip

Posons la d\'efinition suivante:

\begin{defn}\label{defIII55}

Soit ${\mathcal E}$ une cat\'egorie localement petite qui poss\`ede des produits finis arbitraires.

\smallskip

Soit $A$ un anneau interne de ${\mathcal E}$.

\smallskip

Soit ${\mathcal M}od_A$ la cat\'egorie des $A$-modules internes de ${\mathcal E}$.

\smallskip

Alors:

\begin{listeimarge}

\item Pour tous objets $M,N$ de ${\mathcal M}od_A$ et tout objet $E$ de ${\mathcal E}$, on dit qu'un morphisme de ${\mathcal E}$
$$
M \times E \longrightarrow N
$$
est $A$-lin\'eaire en le premier facteur si, pour tout objet $X$ de ${\mathcal E}$, l'application
$$
{\rm Hom} (X,M) \times {\rm Hom} (X,E) \longrightarrow {\rm Hom} (X,N)
$$
est ${\rm Hom} (X,A)$-lin\'eaire en la premi\`ere variable.

\medskip

\item On dit qu'un objet $M$ de ${\mathcal M}od_A$ est exponentiable s'il existe un foncteur
$$
\begin{matrix}
{\mathcal M}od_A &\longrightarrow &{\mathcal E} \, , \hfill \\
\hfill N &\longmapsto &{\mathcal H}om_A(M,N)
\end{matrix}
$$
tel que, pour tout objet $N$ de ${\mathcal M}od_A$, l'objet
$$
{\mathcal H}om_A (M,N) \qquad \mbox{de} \quad {\mathcal E}
$$
repr\'esente le foncteur
$$
\begin{matrix}
{\mathcal E} &\longrightarrow &{\rm Ens} \, , \hfill \\
E &\longmapsto &\{\mbox{morphismes $M \times E \to N$ qui sont $A$-lin\'eaires en le premier facteur}\} \, .
\end{matrix}
$$
\end{listeimarge}
\end{defn}

\begin{remarksqed}
\begin{listeisansmarge}
\item Si $M$ est un objet exponentiable de ${\mathcal M}od_A$, $N$ un objet de ${\mathcal M}od_A$ et $E$ un objet de ${\mathcal E}$, l'ensemble
$$
{\rm Hom} (E,{\mathcal H}om_A (M,N))
$$
des morphismes
$$
M \times E \longrightarrow N
$$
qui sont $A$-lin\'eaires en le premier facteur, a une structure naturelle de groupe ab\'elien [resp. de module sur l'anneau ${\rm Hom} (E,A)$ si $A$ est un anneau commutatif interne de ${\mathcal E}$].

\smallskip

Par cons\'equent, le foncteur
$$
N \longmapsto {\mathcal H}om_A (M,N)
$$
peut \^etre vu comme un foncteur
$$
{\mathcal M}od_A \longrightarrow {\mathcal A}b_{\mathcal E}
$$
vers la cat\'egorie ${\mathcal A}b_{\mathcal E}$ des groupes ab\'eliens internes de ${\mathcal E}$ [resp. peut \^etre vu comme un foncteur
$$
{\mathcal M}od_A \longrightarrow {\mathcal M}od_A
$$
si $A$ est un anneau commutatif interne de ${\mathcal E}$].

\medskip

\item Si $M$ et $M'$ sont deux objets exponentiables de ${\mathcal M}od_A$, tout morphisme de ${\mathcal M}od_A$
$$
M' \longrightarrow M
$$
induit une transformation naturelle du foncteur
$$
N \longmapsto {\mathcal H}om_A (M,N)
$$
vers le foncteur
$$
N \longmapsto {\mathcal H}om_A (M',N) \, .
$$
Cela signifie que tout morphisme de ${\mathcal M}od_A$
$$
N_1 \longrightarrow N_2
$$
induit un carr\'e commutatif
$$
\xymatrix{
{\mathcal H}om_A (M,N_1) \ar[d] \ar[r] &{\mathcal H}om_A (M,N_2) \ar[d] \\
{\mathcal H}om_A (M',N_1) \ar[r] &{\mathcal H}om_A (M',N_2)
}
$$
de la cat\'egorie ${\mathcal A}b_{\mathcal E}$ [resp. ${\mathcal M}od_A$ si l'anneau interne $A$ est commutatif]. 

\end{listeisansmarge}
\end{remarksqed}


Dans un topos, les modules internes sont toujours exponentiables:

\begin{prop}\label{propIII56}

Soit $({\mathcal E},A)$ un topos annel\'e.

\smallskip

Alors:

\begin{listeimarge}

\item Tout objet $M$ de ${\mathcal M}od_A$ est exponentiable, c'est-\`a-dire d\'efinit un foncteur
$$
N \longmapsto {\mathcal H}om_A (M,N)
$$
tel que chaque ${\mathcal H}om_A (M,N)$ repr\'esente le foncteur qui associe \`a tout objet $E$ de ${\mathcal E}$ l'ensemble des morphismes $A$-lin\'eaires en le premier facteur
$$
M \times E \longrightarrow N \, .
$$

\item De plus, si ${\mathcal E} = \widehat{\mathcal C}_J$ est le topos des faisceaux sur un site $({\mathcal C} , J)$, le faisceau
$$
{\mathcal H}om_A (M,N)
$$
associe \`a tout objet $X$ de ${\mathcal C}$ le groupe ab\'elien [resp. le $A(X)$-module, si $A$ est un faisceau d'anneaux commutatifs]
$$
{\rm Hom}_{A_X} (M_X , N_X)
$$
des morphismes de $A_X$-modules entre les restrictions $M_X$ et $N_X$ des faisceaux $M$ et $N$ \`a la cat\'egorie relative ${\mathcal C}/X$.
\end{listeimarge}
\end{prop}

\bigskip

\begin{demo}

Supposons que ${\mathcal E} = \widehat{\mathcal C}_J$ est le topos des faisceaux sur un site $({\mathcal C} , J)$.

\smallskip

Le pr\'efaisceau
$$
{\mathcal H}om_A (M,N) : X \longmapsto {\rm Hom}_{A_X} (M_X , N_X)
$$
est un sous-objet dans $\widehat{\mathcal C}$ du faisceau
$$
{\mathcal H}om (M,N) : X \longmapsto {\rm Hom} (M_X , N_X) = {\rm Hom} (M \times \ell (X),N) \, .
$$

De plus, un morphisme de faisceaux sur ${\mathcal C}/X$
$$
M_X \longrightarrow N_X
$$
vu comme un morphisme de faisceaux sur ${\mathcal C}$
$$
M \times \ell (X) \longrightarrow N
$$
est $A_X$-lin\'eaire si et seulement si les carr\'es
$$
\xymatrix{
M \times M \times \ell (X) \ar[d] \ar[r] &N \times N \ar[d] \\
M \times \ell (X) \ar[r] &N
}
$$
et 
$$
\xymatrix{
A \times M \times \ell (X) \ar[d] \ar[r] &A \times N \ar[d] \\
M \times \ell (X) \ar[r] &N
}
$$
sont commutatifs.

\smallskip

Or, si $S$ est un crible couvrant d'un objet $X$ de ${\mathcal C}$, on a d'apr\`es le lemme \ref{lemIII53}
$$
\ell (X) = \varinjlim_{(U \to X) \in S} \ell (U) \, .
$$

Il en r\'esulte que le sous-pr\'efaisceau
$$
{\mathcal H}om_A (M,N) \qquad \mbox{de} \quad {\mathcal H}om (M,N)
$$
est en fait un faisceau.

\smallskip

Si $E$ est un objet de ${\mathcal E} = \widehat{\mathcal C}_J$, consid\'erons un morphisme de ${\mathcal E}$
$$
M \times E \longrightarrow N \, .
$$

Il lui correspond un morphisme $E \to {\mathcal H}om (M,N)$.

\smallskip

Le morphisme $M \times E \to N$ est $A$-lin\'eaire en le premier facteur si et seulement si, pour tout objet $F$ de ${\mathcal E} = \widehat{\mathcal C}_J$ et tout morphisme $F \to E$, l'application induite
$$
{\rm Hom} (F,M) \longrightarrow {\rm Hom} (F,N)
$$
est ${\rm Hom} (F,A)$-lin\'eaire.

\smallskip

Or un tel faisceau $F$ s'\'ecrit d'apr\`es le lemme \ref{lemIII52} comme une colimite
$$
F = \varinjlim_{(X,x) \in \int\!\!F} \ell (X) \, .
$$

Il suffit donc que pour tout objet $X$ de ${\mathcal C}$ et tout \'el\'ement $e \in E(X)$, l'application induite
$$
M(X) \longrightarrow N(X)
$$
soit $A(X)$-lin\'eaire.

\smallskip

C'est \'equivalent \`a demander que pour tout $X$ et tout \'el\'ement $e \in E(X) = {\rm Hom} (\ell (X),E)$, le morphisme induit
$$
M_X \longrightarrow N_X
$$
soit $A_X$-lin\'eaire, autrement dit que le morphisme correspondant
$$
\ell (X) \longrightarrow {\mathcal H}om (M,N)
$$
se factorise en
$$
\ell (X) \longrightarrow {\mathcal H}om_A (M,N) \, .
$$

Comme d'apr\`es le lemme \ref{lemIII52},
$$
E = \varinjlim_{(X,e) \in \int\!\!E} \ell (X) \, ,
$$
cela revient \`a demander que le morphisme
$$
E \longrightarrow {\mathcal H}om (M,N)
$$
se factorise en
$$
E \longrightarrow {\mathcal H}om_A (M,N) \, .
$$

Cela d\'emontre (ii) et (i) dans le cas o\`u ${\mathcal E} = \widehat{\mathcal C}_J$. Le cas g\'en\'eral de (i) s'en d\'eduit puisque la propri\'et\'e consid\'er\'ee est invariante par \'equivalence de cat\'egories. 

\end{demo}

\subsection{Les produits tensoriels de faisceaux lin\'eaires}\label{subsec355}

\medskip

Par d\'efinition, les foncteurs
$$
H \longmapsto {\mathcal H}om (F,H)
$$
admettent pour adjoints \`a gauche les foncteurs
$$
G \longmapsto F \times G \, .
$$

Dans le cas de faisceaux lin\'eaires, c'est-\`a-dire de modules internes d'un topos annel\'e, les foncteurs
$$
N \longmapsto {\mathcal H}om_A (M,N)
$$
admettent encore des adjoints \`a gauche qui sont les foncteurs de produits tensoriels:

\begin{prop}\label{propIII57}

Soit ${\mathcal E}$ un topos.

\smallskip

Soit $A$ un anneau [resp. un anneau commutatif] interne de ${\mathcal E}$.

\smallskip

Alors:

\begin{listeimarge}

\item Pour tout $A$-module interne $M$, le foncteur
$$
\begin{matrix}
\hfill N &\longmapsto & {\mathcal H}om_A (M,N) \, , \hfill \\
{\mathcal M}od_A &\longrightarrow &{\mathcal A}b_{\mathcal E} \quad \mbox{[resp. ${\mathcal M}od_A$]}
\end{matrix}
$$
poss\`ede un adjoint \`a gauche not\'e
$$
\begin{matrix}
\otimes : &\hfill L &\longmapsto &M \otimes L \, , \\
&\hfill {\mathcal A}b_{\mathcal E} &\longrightarrow &{\mathcal M}od_A \hfill
\end{matrix}
$$
[resp.
$$
\begin{matrix}
\otimes_A : &\hfill L &\longmapsto &M \otimes_A L \, , \\
&\hfill {\mathcal M}od_A &\longrightarrow &{\mathcal M}od_A \ \mbox{]}. \hfill
\end{matrix}
$$

\item De plus, si ${\mathcal E} = \widehat{\mathcal C}_J$ est le topos des faisceaux sur un site $({\mathcal C},J)$, le foncteur
$$
L \longmapsto M \otimes L \qquad \mbox{[resp.} \quad L \longmapsto M \otimes_A L \ \mbox{]}
$$
peut \^etre construit comme le compos\'e du foncteur de faisceautisation
$$
j^* : {\mathcal A}b_{\mathcal C} \to {\mathcal A}b_{{\mathcal C},J} \qquad \mbox{[resp.} \quad j^* : {\mathcal M}od_{{\mathcal C},J} \to {\mathcal M}od_{{\mathcal C},J,A} \ \mbox{]}
$$
avec le foncteur qui associe \`a tout faisceau de groupes ab\'eliens [resp. de $A$-modules] $L$ sur $({\mathcal C},J)$ le pr\'efaiseau
$$
\begin{matrix}
&X &\longmapsto &M(X) \otimes L(X) \hfill \\
\mbox{[resp.} &X &\longmapsto &M(X) \otimes_{A(X)} L(X) \ \mbox{]}.
\end{matrix}
$$
\end{listeimarge}
\end{prop}

\begin{remarks}
\begin{listeisansmarge}
\item Tout morphisme de $A$-modules internes
$$
M \longrightarrow M'
$$
induit une transformation naturelle du foncteur
$$
L \longmapsto M \otimes L \qquad \mbox{[resp.} \quad M \otimes_A L \ \mbox{]}
$$
vers le foncteur
$$
L \longmapsto M' \otimes L \qquad \mbox{[resp.} \quad M' \otimes_A L \ \mbox{]}.
$$

\item Si $A$ est un anneau commutatif interne, les deux foncteurs
$$
\begin{matrix}
{\mathcal M}od_A \times {\mathcal M}od_A &\longrightarrow &{\mathcal M}od_A \, , \hfill \\
\hfill (M,L) &\longmapsto &M \otimes_A L \, , \\
\hfill (M,L) &\longmapsto &L \otimes_A M \hfill
\end{matrix}
$$
sont canoniquement isomorphes.

\smallskip

En effet, pour tout objet $N$ de ${\mathcal M}od_A$, se donner un morphisme de ${\mathcal M}od_A$
$$
M \otimes_A L \longrightarrow N
$$
\'equivaut \`a se donner un morphisme de ${\mathcal E}$
$$
M \times L \longrightarrow N
$$
qui est $A$-lin\'eaire en chacun des deux facteurs.

\end{listeisansmarge}
\end{remarks}

\begin{demo}
\begin{listeisansmarge}
\item se r\'eduit \`a (ii) puisque la propri\'et\'e consid\'er\'ee est invariante par \'equivalence de cat\'egories.

\medskip

\item Le cas o\`u ${\mathcal E} = {\rm Ens}$ est le topos des ensembles est la propri\'et\'e de d\'efinition des produits tensoriels de groupes ab\'eliens ou de modules sur un anneau commutatif.

\smallskip

Le cas o\`u ${\mathcal E} = \widehat{\mathcal C}$ est le topos des pr\'efaisceaux sur une cat\'egorie essentiellement petite s'en d\'eduit par fonctorialit\'e des produits tensoriels
$$
M \otimes L \qquad \mbox{[resp.} \quad M \otimes_A L \ \mbox{]}
$$
par rapport aux deux facteurs $M,L$ [resp. et par rapport \`a l'anneau commutatif $A$].

\smallskip

Si ${\mathcal E} = \widehat{\mathcal C}_J$ est le topos des faisceaux sur un site $({\mathcal C},J)$, muni des deux foncteurs adjoints $j^* : \widehat{\mathcal C} \to \widehat{\mathcal C}_J$ et $j_* : \widehat{\mathcal C}_J \to \widehat{\mathcal C}$, on commence par remarquer que pour tous faisceaux de $A$-modules $M$ et $N$, on a
$$
j_* {\mathcal H}om_A (M,N) = {\mathcal H}om_{j_* A} (j_* M , j_* N) \, .
$$
Si donc $L$ est un faisceau de groupes ab\'eliens [resp. de $A$-modules], se donner un morphisme de faisceaux de groupes ab\'eliens [resp. de $A$-modules]
$$
L \longrightarrow {\mathcal H}om_A (M,N)
$$
\'equivaut \`a se donner un morphisme de pr\'efaisceaux de $j_* A$-modules
$$
j_* M \otimes j_* L \longrightarrow j_* N \qquad \mbox{[resp.} \quad j_* M \otimes_{j_* A} j_* L \longrightarrow j_* N \ \mbox{]}
$$
et donc \`a se donner un morphisme de faisceaux de $A$-modules
$$
j^*  (j_* M \otimes j_* L) \longrightarrow N
$$
[resp.
$$
j^* (j_* M \otimes_{j_* A} j_* L) \longrightarrow N \ \mbox{]}.
$$
C'est ce que l'on voulait. 
\end{listeisansmarge}
\end{demo}

\section{Classification des sous-objets dans les topos}\label{sec36}

\subsection{Ensembles de sous-objets et de quotients}\label{subsec361}

\medskip

D\'emontrons encore \`a partir du cas du topos des ensembles:

\begin{prop}\label{propIII61}

Tout topos ${\mathcal E}$ poss\`ede les deux propri\'et\'es suivantes:

\medskip
\begin{enumerate}
\item[(10)] Les sous-objets de tout objet $F$ de ${\mathcal E}$ forment un ensemble.

\medskip

\item[(11)] Les objets quotients de tout objet $F$ de ${\mathcal E}$ forment un ensemble.
\end{enumerate}
\end{prop}

\begin{demo}

La propri\'et\'e (10) est connue lorsque ${\mathcal E} = {\rm Ens}$ est le topos des ensembles: en effet, les parties d'un ensemble $E$ forment un ensemble ${\mathcal P} (E)$.

\smallskip

Si ${\mathcal E} = \widehat{\mathcal C}$ est le topos des pr\'efaisceaux sur une petite cat\'egorie ${\mathcal C}$, les sous-objets d'un objet $F$ de ${\mathcal E}$ sont les familles de parties
$$
P_X \subset F(X) \, , \qquad X \in {\rm Ob} (X) \, ,
$$
telles que, pour tout morphisme $f : X \to Y$ de ${\mathcal C}$, l'application induite
$$
F(f) : F(Y) \longrightarrow F(X)
$$
satisfait la condition
$$
F(f)(P_Y) \subset P_X \, .
$$

Par cons\'equent, les sous-objets d'un tel objet $F$ de $\widehat{\mathcal C}$ forment un ensemble.

\smallskip

Si ${\mathcal E} = \widehat{\mathcal C}_J$ est le topos des faisceaux sur un petit site $({\mathcal C},J)$, les sous-objets d'un objet $F$ de ${\mathcal E}$ sont les sous-objets de l'objet $j_* F$ de $\widehat{\mathcal C}$ qui sont des faisceaux c'est-\`a-dire satisfont la condition
$$
F(X) \xrightarrow{ \ \sim \ } \varprojlim_{(U \to X) \in S} F(U)
$$
pour tout crible $J$-couvrant $S$ de tout objet $X$ de ${\mathcal C}$.

\smallskip

Ces conditions d\'efinissent un sous-ensemble de l'ensemble des sous-objets de $j_* F$ dans $\widehat{\mathcal C}$.

\smallskip

La propri\'et\'e (10) dans le cas g\'en\'eral d'un topos ${\mathcal E}$ arbitraire s'en d\'eduit puisqu'elle est invariante par \'equivalence de cat\'egories. 

\smallskip

D'apr\`es les propositions \ref{propIII33} et \ref{propIII34}, se donner un quotient d'un objet $F$ d'un topos ${\mathcal E}$ \'equivaut \`a se donner un sous-objet
$$
R \xhookrightarrow{ \ { \ } \ } F \times F
$$ 
qui est une relation d'\'equivalence.

\smallskip

Cette condition d\'efinit un sous-ensemble dans l'ensemble des sous-objets de $F \times F$.

\smallskip

Cela d\'emontre la propri\'et\'e (11). 

\end{demo}

\subsection{Op\'erations sur les sous-objets d'un objet}\label{subsec362}

\medskip

Dans un topos, l'ensemble ordonn\'e des sous-objets d'un objet est un treillis au sens de la d\'efinition II.2.3:

\begin{prop}\label{propIII62}

Pour tout objet $F$ d'un topos ${\mathcal E}$, l'ensemble ordonn\'e $\Omega (F)$ des sous-objets de $F$ est un treillis.

\smallskip

Autrement dit:

\begin{listeimarge}

\item Toute famille $(F_i)_{i \in I}$ de sous-objets de $F$ a une r\'eunion
$$
\bigvee_{i \in I} F_i
$$
caract\'eris\'ee par la propri\'et\'e que pour les sous-objets $F'$ de $F$, il y a \'equivalence entre les conditions
$$
F_i \subset F' \, , \qquad \forall \, i \in I \, ,
$$
et
$$
\bigvee_{i \in I} F_i \subset F' \, .
$$

\item Toute famille finie $(F_i)_{1 \leq i \leq n}$ de sous-objets de $F$ a une intersection
$$
F_1 \wedge \cdots \wedge F_n
$$
caract\'eris\'ee par la propri\'et\'e que pour les sous-objets $F'$ de $F$ il y a \'equivalence entre les conditions
$$
F' \subset F_i \, , \qquad 1 \leq i \leq n \, ,
$$
et
$$
F' \subset F_1 \wedge \cdots \wedge F_n \, .
$$

\item Pour toute famille $(F_i)_{i \in I}$ de sous-objets de $F$ et tout sous-objet $F'$ de $F$, on a
$$
F' \wedge \left( \bigvee_{i \in I} F_i \right) = \bigvee_{i \in I} (F' \cap F_i) \, .
$$
\end{listeimarge}
\end{prop}

\begin{remarks}
\begin{listeisansmarge}
\item En particulier, tout objet $F$ de ${\mathcal E}$ admet un sous-objet minimal. Il est repr\'esent\'e par l'unique morphisme
$$
\emptyset \longrightarrow F
$$
d\'efini sur l'objet initial $\emptyset$ de ${\mathcal E}$, car c'est un monomorphisme.

\medskip

\item Tout sous-objet $F'$ de $F$ admet un compl\'ementaire
$$
\neg \, F'
$$
caract\'eris\'e par la propri\'et\'e que pour les sous-objets $F''$ de $F$ il y a \'equivalence entre les conditions
$$
F'' \wedge F' = \emptyset
$$
et
$$
F'' \subset \ \neg \, F' \, .
$$

En effet, $\neg \, F'$ peut \^etre construit comme la r\'eunion des sous-objets $F''$ de $F$ tels que
$$
F'' \wedge F' = \emptyset \, .
$$
Il r\'esulte de (iii) que l'intersection de cette r\'eunion et de $F'$ est le sous-objet minimal $\emptyset$.

\medskip

\item Plus g\'en\'eralement, pour tout sous-objet $F'$ de $F$, le foncteur
$$
\begin{matrix}
\Omega (F) &\longrightarrow &\Omega (F) \hfill \\
\hfill G &\longmapsto &F' \wedge G
\end{matrix}
$$
admet un adjoint \`a droite not\'e
$$
H \longmapsto (F' \Rightarrow H) \, .
$$
Il est caract\'eris\'e par la propri\'et\'e que pour tous sous-objets $G$ et $H$ de $F$
$$
G \subset (F' \Rightarrow H) \quad \mbox{si et seulement si} \quad F' \wedge G \subseteq H \, .
$$
Le sous-objet $(F' \Rightarrow H)$ est construit comme la r\'eunion des sous-objets $G$ tels que $F' \wedge G \subseteq H$.
\end{listeisansmarge}
\end{remarks}

\begin{demosansqed}
\begin{listeisansmarge}
\item Consid\'erons l'ensemble de tous les sous-objets de $F$
$$
F' \xhookrightarrow{ \ { \ } \ } F
$$
tels que chaque $F_i \hookrightarrow F$, $i \in I$, se factorise en
$$
F_i \xhookrightarrow{ \ { \ } \ } F' \xhookrightarrow{ \ { \ } \ } F \, .
$$
Alors la limite $\underset{i \in I}{\bigvee} \, F_i$ du diagramme constitu\'e de $F$, de ces objets $F'$ et des fl\`eches $F' \hookrightarrow F$ r\'epond \`a la question pos\'ee.

\smallskip

En effet, le morphisme
$$
\bigvee_{i \in I} F_i \longrightarrow F
$$
est un monomorphisme, chaque $F_i \hookrightarrow F$ se factorise en 
$$
F_i \xhookrightarrow{ \ { \ } \ } \bigvee_{i \in I} F_i \, ,
$$
et tout $F' \hookrightarrow F$ \`a travers lequel se factorisent les $F_i \hookrightarrow F$ s'inscrit dans un triangle commutatif:
$$
\xymatrix{
\underset{i \in I}{\bigvee} F_i \ar@{^{(}->}[rd] \ar@{^{(}->}[rr] &&F' \ar@{^{(}->}[ld] \\
&F
}
$$

\item La limite du diagramme constitu\'e de $F$, des objets $F_i$ et des fl\`eches $F_i \hookrightarrow F$, $1 \leq i \leq n$, est un sous-objet
$$
F_1 \wedge \cdots \wedge F_n \xhookrightarrow{ \ { \ } \ } F
$$
qui r\'epond \`a la question pos\'ee.

\medskip

\item est le cas particulier des monomorphismes dans la partie (iii) du lemme suivant:
\end{listeisansmarge}
\end{demosansqed}
\begin{lem}\label{lemIII63}

Soit ${\mathcal E}$ un topos.

\smallskip

Consid\'erons pour tout objet $F$ de ${\mathcal E}$ l'ensemble $\Omega (F)$ de ses sous-objets.

\smallskip

Alors, pour tout morphisme $G \xrightarrow{ \ f \ } F$ de ${\mathcal E}$, le foncteur
$$
(F' \xhookrightarrow{ \ { \ } \ } F) \longmapsto (F' \times_F G \xhookrightarrow{ \ { \ } \ } G)
$$
d\'efinit une application
$$
f^{-1} = \Omega (f) : \Omega (F) \longrightarrow \Omega (G)
$$
qui poss\`ede les propri\'et\'es suivantes:

\medskip

\begin{listeimarge}
\item[(i)] Elle pr\'eserve la relation d'ordre.
\item[(ii)] Elle commute avec les intersections.
\item[(iii)] Elle commute avec les r\'eunions.
\end{listeimarge}
\end{lem}

\begin{remark}

L'application $f^{-1} = \Omega (f) : \Omega (F) \to \Omega (G)$ respecte aussi les op\'erations $\neg$ et $\Rightarrow$ au sens que pour tous sous-objets $F'$ et $F''$ de $F$ on a
$$
f^{-1} (\neg \, F') = \neg \, f^{-1} (F') \qquad \mbox{et} \qquad f^{-1} (F' \Rightarrow F'') = (f^{-1} (F') \Rightarrow f^{-1} (F'')) \, .
$$
Cette propri\'et\'e est laiss\'ee en exercice au lecteur.

\end{remark}

\begin{demolem}

Le foncteur de changement de base $\bullet \times_F G$ d\'efinit une application
$$
f^{-1} = \Omega (f) : \Omega (F) \longrightarrow \Omega (G)
$$
car il transforme les monomorphismes en monomorphismes.

\begin{listeisansmarge}

\item Comme c'est un foncteur, l'application induite pr\'eserve la relation d'ordre.

\medskip

\item Ce foncteur $\bullet \times_F G$ pr\'eserve les limites, donc l'application induite pr\'eserve les intersections.

\medskip

\item Si ${\mathcal E} = {\rm Ens}$ est le topos des ensembles, $G \xrightarrow{ \ f \ } F$ est une application et on a pour toute famille de parties $F_i \subset F$, $i \in I$, de $F$ la formule
$$
f^{-1} \left( \bigcup_{i \in I} F_i \right) = \bigcup_{i \in I} f^{-1} (F_i) \, .
$$

Si ${\mathcal E} = \widehat{\mathcal C}$ est le topos des pr\'efaisceaux sur une cat\'egorie essentiellement petite ${\mathcal C}$, $G \xrightarrow{ \ f \ } F$ est une famille compatible d'applications
$$
f_X : G(X) \longrightarrow F(X) \, .
$$
Chaque $F_i$, $i \in I$, est une famille compatible de parties
$$
F_i (X) \subset F(X) \, , \qquad i \in I \, ,
$$
et $f^{-1} \left( \underset{i \in I}{\bigvee} \, F_i \right)$ est la famille compatible des parties
$$
f_X^{-1} \left( \bigcup_{i \in I} F_i(X) \right)
$$
tandis que $\underset{i \in I}{\bigvee} \, f^{-1} F_i$ est la famille compatible des parties 
$$
\bigcup_{i \in I} f_X^{-1} (F_i (X)) \, .
$$

Le r\'esultat dans ce cas o\`u ${\mathcal E} = \widehat{\mathcal C}$ se d\'eduit donc du cas du topos des ensembles.

\smallskip

Si ${\mathcal E} = \widehat{\mathcal C}_J$ est le topos des faisceaux sur un site $({\mathcal C},J)$, le foncteur $j_* : \widehat{\mathcal C}_J \hookrightarrow \widehat{\mathcal C}$ transforme les sous-objets en sous-objets et pr\'eserve les produits fibr\'es tandis que son adjoint \`a gauche $j^* : \widehat{\mathcal C} \to \widehat{\mathcal C}_J$ transforme les sous-objets en sous-objets et pr\'eserve les produits fibr\'es et les r\'eunions.

\smallskip

Si les $F_i$, $i \in I$, sont une famille de sous-objets de $F$ dans $\widehat{\mathcal C}_J$, l'identit\'e dans $\widehat{\mathcal C}$
$$
j_* G \times_{j_* F} \left( \bigvee_{i \in I} j_* F_i \right) = \bigvee_{i \in I} (j_* G \times_{j_* F} j_* F_i )
$$
induit donc l'identit\'e dans $\widehat{\mathcal C}_J$
$$
G \times_F \left( \bigvee_{i \in I} F_i \right) = \bigvee_{i \in I} (G \times_F F_i) \, .
$$

Le cas g\'en\'eral d'un topos arbitraire ${\mathcal E}$ en r\'esulte puisque la propri\'et\'e consid\'er\'ee est invariante par \'equivalence de cat\'egories.

\smallskip

Cela termine la d\'emonstration du lemme \ref{lemIII63}, donc aussi celle de la proposition \ref{propIII62}. 

\end{listeisansmarge}
\end{demolem}

\subsection{Le classificateur des sous-objets}\label{subsec363}

\medskip

D'apr\`es le lemme \ref{lemIII63}, associer \`a tout objet $F$ d'un topos ${\mathcal E}$ l'ensemble $\Omega (F)$ de ses sous-objets d\'efinit un pr\'efaisceau
$$
\begin{matrix}
\Omega : &\hfill {\mathcal E}^{\rm op} &\longrightarrow &{\rm Ens} \, , \hfill \\
&\hfill F &\longmapsto &\Omega (F) \, , \hfill \\
&\left( G \xrightarrow{ \ f \ } F \right) &\longmapsto &\left( \Omega (F) \xrightarrow{ \ \Omega (f) \ } \Omega (G) \right) .
\end{matrix}
$$

Ce pr\'efaisceau est toujours repr\'esentable:

\begin{prop}\label{propIII64}

Soit ${\mathcal E}$ un topos.

\begin{listeimarge}

\item Il poss\`ede la propri\'et\'e suivante:

\medskip
\begin{enumerate}
\item[(12)] Le pr\'efaisceau
$$
\begin{matrix}
\Omega : &\hfill {\mathcal E}^{\rm op} &\longrightarrow &{\rm Ens} \, , \hfill \\
&\hfill F &\longmapsto &\Omega (F) = \{\mbox{sous-objets de $F$}\}
\end{matrix}
$$
est repr\'esentable.\end{enumerate}

\medskip

\item De plus, si ${\mathcal E} = \widehat{\mathcal C}_J$ est le topos des faisceaux sur un site $({\mathcal C},J)$, le pr\'efaisceau $\Omega$ sur ${\mathcal E}$ est repr\'esentable par l'objet de ${\mathcal E} = \widehat{\mathcal E}_J$
$$
\begin{matrix}
\Omega_J : &\hfill {\mathcal C}^{\rm op} &\longrightarrow &{\rm Ens} \, , \hfill \\
&\hfill X &\longmapsto &\Omega_J (X) = \{\mbox{cribles $J$-ferm\'es de $X$}\}
\end{matrix}
$$
introduit au paragraphe \ref{subsec242}.
\end{listeimarge}
\end{prop}

\begin{demosansqed}
\begin{listeisansmarge}
\item Comme la propri\'et\'e (12) est invariante par \'equivalence de cat\'egories, il suffit de d\'emontrer que si ${\mathcal E} = \widehat{\mathcal C}_J$, le pr\'efaisceau $\Omega : {\mathcal E}^{\rm op} \to {\rm Ens}$ est repr\'esentable par l'objet $\Omega_J$ de $\widehat{\mathcal C}_J = {\mathcal E}$.

\smallskip

Ainsi, (i) est r\'eduit \`a (ii).

\medskip

\item Commen\c cons par d\'emontrer que, pour tout objet $X$ de ${\mathcal C}$ et son image $\ell(X)$ par le foncteur canonique
$$
\ell : {\mathcal C} \xrightarrow{ \ y \ } \widehat{\mathcal C} \xrightarrow{ \ j^* \ } \widehat{\mathcal C}_J \, ,
$$
l'ensemble $\Omega (\ell (X))$ des sous-objets de $\ell(X)$ dans $\widehat{\mathcal C}_J$ s'identifie \`a l'ensemble $\Omega_J (X)$ des cribles $J$-ferm\'es de~$X$.

\smallskip

C'est le contenu du lemme suivant:
\end{listeisansmarge}
\end{demosansqed}

\begin{lem}\label{lemIII65}

Soient $({\mathcal C},J)$ un site et $\ell : {\mathcal C} \xrightarrow{ \ y \ } \widehat{\mathcal C} \xrightarrow{ \ j^* \ } \widehat{\mathcal C}_J$ le foncteur canonique associ\'e.

\smallskip

Alors les deux applications en sens inverse
$$
( F \xhookrightarrow{ \ { \ } \ } \ell (X)) \longmapsto (F \times_{\ell (X)} y(X) \xhookrightarrow{ \ { \ } \ } y(X))
$$
et
$$
(S \xhookrightarrow{ \ { \ } \ } y(X)) \longmapsto (j^* S \xhookrightarrow{ \ { \ } \ } j^* y(X)=\ell(X))
$$
d\'efinissent deux bijections r\'eciproques l'une de l'autre entre l'ensemble des sous-objets de $\ell (X)$ dans $\widehat{\mathcal C}_J$ et l'ensemble des cribles $J$-ferm\'es de $X$.
\end{lem}

\begin{demolem}

On rappelle qu'un cribre $S$ d'un objet $X$ de ${\mathcal C}$ n'est pas autre chose qu'un sous-objet $S \hookrightarrow y(X)$ dans $\widehat{\mathcal C}$.

\smallskip

D'autre part, on a pour tout sous-objet $F$ de $\ell (X)$ une identification
$$
j^* (F \times_{\ell(X)} y(X)) = F \times_{\ell(X)} j^* \circ y(X) = F 
$$
puisque le foncteur $j^*$ respecte les limites finies et qu'il fixe les faisceaux.

\smallskip 

Il en r\'esulte que les deux applications consid\'er\'ees d\'efinissent deux bijections r\'eciproques l'une de l'autre entre l'ensemble des sous-objets de $\ell (X)$ dans $\widehat{\mathcal C}_J$ et l'ensemble des cribles $S$ de $X$ tels que
$$
j^* S \times_{\ell (X)} y(X) = S
$$
si $S$ est vu comme un sous-objet $S \hookrightarrow y(X)$ dans $\widehat{\mathcal C}$.

\smallskip

Or, le foncteur $j^* : \widehat{\mathcal C} \to \widehat{\mathcal C}_J$ est construit comme le compos\'e
$$
P \longmapsto (P^+)^+
$$
avec lui-m\^eme du foncteur
$$
P \longmapsto P^+ = \left[ X' \longmapsto \varinjlim_{S' \in J(X')} \ \varprojlim_{(U' \to X') \in S'} P(U') \right].
$$

Comme le foncteur $P \mapsto P^+$ respecte les limites finies et fixe les faisceaux, l'identit\'e
$$
j^* S \times_{\ell (X)} y(X) = S
$$
est \'equivalente \`a l'identit\'e
$$
S^+ \times_{y(X)^+} y(X)=S
$$
laquelle signifie exactement que $S$ est un crible $J$-ferm\'e de $X$.

\smallskip

Cela ach\`eve la d\'emonstration du lemme. 

\end{demolem}

\bigskip

\noindent {\bf Suite de la d\'emonstration de la proposition \ref{propIII64} (ii):}

\smallskip

Consid\'erons un objet $F$ de $\widehat{\mathcal C}_J$.

\smallskip

Nous devons montrer que les ensembles $\Omega (F)$ et ${\rm Hom} (F,\Omega_J)$ s'identifient.

\smallskip

On sait d'apr\`es le lemme \ref{lemIII52} que $F$ s'\'ecrit comme la colimite dans $\widehat{\mathcal C}_J$
$$
F = \varinjlim_{(X,x) \in \int\!\!F} \ell (X) \, .
$$

On en d\'eduit d'apr\`es le lemme \ref{lemIII65}
\begin{eqnarray}
{\rm Hom} (F,\Omega_J) &= &\varprojlim_{(X,x) \in \int\!\!F} {\rm Hom} (\ell(X),\Omega_J) \nonumber \\
&= &\varprojlim_{(X,x) \in \int\!\!F} \Omega (\ell (X)). \nonumber
\end{eqnarray}

La conclusion cherch\'ee
$$
{\rm Hom} (F,\Omega_J) = \Omega (F)
$$
r\'esulte alors du lemme suivant:

\begin{lem}\label{lemIII66}

Dans un topos ${\mathcal E}$, le pr\'efaisceau des sous-objets
$$
\Omega : {\mathcal E}^{\rm op} \longrightarrow {\rm Ens}
$$
transforme les colimites en limites, au sens que pour tout carquois $D$ et tout $D$-diagramme $F_{\bullet}$ de ${\mathcal E}$, l'application canonique
$$
\Omega \left( \varinjlim_D F_{\bullet} \right) \longrightarrow \varprojlim_D \Omega (F_{\bullet})
$$
est une bijection.
\end{lem}

\begin{demolem}

Si ${\mathcal E} = {\rm Ens}$ est le topos des ensembles, la propri\'et\'e r\'esulte de ce que le pr\'efaisceau $\Omega$ est alors repr\'esentable par l'ensemble $\{0,1\}$.

\smallskip

Si ${\mathcal E} = \widehat{\mathcal C}$ est le topos des pr\'efaisceaux sur une cat\'egorie essentiellement petite ${\mathcal C}$, la propri\'et\'e r\'esulte de ce que les colimites dans $\widehat{\mathcal C}$ se calculent composante par composante et que les sous-objets d'un pr\'efaisceau $F$ sont les familles compatibles de sous-ensembles des $F(X)$, $X \in {\rm Ob} ({\mathcal C})$.

\smallskip

Si ${\mathcal E} = \widehat{\mathcal C}_J$ est la cat\'egorie des faisceaux sur un site $({\mathcal C},J)$, consid\'erons donc un $D$-diagramme $F_{\bullet}$ de $\widehat{\mathcal C}_J$ et sa colimite $F$. Si $F'$ est un sous-objet de $F$, il r\'esulte de la proposition \ref{propIII41} que
$$
F' = \varinjlim_D F' \times_F F_{\bullet} \, .
$$

R\'eciproquement, si les $F'_d$, $d \in {\rm Ob}(D)$, sont une famille de sous-objets des $F_d$ tels que, pour tout morphisme $d_1 \to d_2$ de $D$, on ait
$$
F'_{d_1} = F'_{d_2} \times_{F_{d_2}} F_{d_1} \, ,
$$
alors la colimite $F'$ des $j_* F'_d$ dans $\widehat{\mathcal C}$ est un sous-objet de $j_* F$ et v\'erifie la propri\'et\'e
$$
j_* F'_d = F' \times_{j_* F} j_* F_d \, , \qquad \forall \, d \in {\rm Ob}(D) \, .
$$

Donc la colimite $j^* F'$ des $F'_d$ dans $\widehat{\mathcal C}_J$ est un sous-objet de $F$ et v\'erifie la propri\'et\'e
$$
F'_d = j^* F' \times_F F_d \, , \qquad \forall \, d \in {\rm Ob} (D) \, .
$$

Cela prouve comme voulu que l'application
$$
\Omega (F) \longrightarrow \varprojlim_D \Omega (F_{\bullet})
$$
est une bijection.

\smallskip

Ainsi, le lemme est d\'emontr\'e et avec lui la proposition \ref{propIII64}. 

\end{demolem}

\bigskip

On d\'eduit de la proposition \ref{propIII64}, du lemme \ref{lemIII63} et de la proposition \ref{propIII62}:

\begin{cor}\label{corIII67}

Soit ${\mathcal E}$ un topos.

\smallskip

Soit $\Omega$ l'objet de ${\mathcal E}$ qui classifie les sous-objets.

\smallskip

Alors la r\'eunion et l'intersection des sous-objets d\'efinissent deux op\'erations associatives et commutatives
$$
\vee : \Omega \times \Omega \longrightarrow \Omega
$$
et
$$
\wedge : \Omega \times \Omega \longrightarrow \Omega \, .
$$

Elles rendent commutatif le carr\'e:
$$
\xymatrix{
\Omega \times \Omega \times \Omega \ar[d]_{{\rm id} \times \vee} \ar[rrr]^{(p_1 \wedge p_2 , p_1 \wedge p_3)} &&&\Omega \times \Omega \ar[d]^{\vee} \\
\Omega \times \Omega \ar[rrr]^{\wedge} &&&\Omega
}
$$

\end{cor}

\begin{remark}

Il r\'esulte de la remarque qui suit le lemme \ref{lemIII63} que l'objet $\Omega$ de ${\mathcal E}$ est \'egalement muni de deux op\'erations
$$
\neg : \Omega \longrightarrow \Omega
$$
et
$$
\Rightarrow \, : \Omega \times \Omega \longrightarrow \Omega \, .
$$

\end{remark}

\begin{demo}

Comme le pr\'efaisceau des ensembles de sous-objets
$$
\Omega : {\mathcal E}^{\rm op} \longrightarrow {\rm Ens}
$$
est repr\'esentable par un objet de ${\mathcal E}$ que l'on note encore $\Omega$, l'existence des morphismes $\vee$ et $\wedge$ de $\Omega \times \Omega$ dans $\Omega$ r\'esulte d'apr\`es le lemme de Yoneda des parties (ii) et (iii) du lemme \ref{lemIII63}.

\smallskip

La commutativit\'e du carr\'e r\'esulte alors de la proposition \ref{propIII62} (iii). 
\end{demo}

\subsection{Le cas des cat\'egories d'objets lin\'eaires des topos}\label{subsec364}

\medskip

Dans les cat\'egories ab\'eliennes, on a de mani\`ere g\'en\'erale:

\begin{lem}\label{lemIII68}

Soit ${\mathcal A}$ une cat\'egorie ab\'elienne.

\smallskip

Alors, pour tout objet $M$, les deux applications
$$
(M' \xhookrightarrow{ \ { \ } \ } M) \longmapsto {\rm coker} (M' \longrightarrow M) \, ,
$$
$$
(M \, -\!\!\!\twoheadrightarrow M'') \longmapsto {\rm ker} (M \longrightarrow M'')
$$
d\'efinissent deux bijections r\'eciproques l'une de l'autre entre sous-objets de $M$ et quotients de $M$.
\end{lem}

\begin{remark}

En particulier, les sous-objets d'un objet $M$ forment un ensemble si et seulement si ses quotients forment un ensemble.
\end{remark}


\begin{demo}

En effet, il r\'esulte de la d\'efinition de la notion de cat\'egorie ab\'elienne que:
\begin{enumerate}
\item[$\bullet$] si $M' \hookrightarrow M$ est un monomorphisme et $M'' = {\rm coker} (M' \to M)$, alors 

$M' \to \ker (M \to M'')$ est un isomorphisme,
\item[$\bullet$] si $M \twoheadrightarrow M''$ est un \'epimorphisme et $M' = \ker (M \to M'')$, alors 

${\rm coker} (M' \to M) \to M''$ est un isomorphisme.
\end{enumerate}
\end{demo}

\medskip

On d\'eduit de ce lemme et de la proposition \ref{propIII61}:

\begin{cor}\label{corIII69}

Soit $({\mathcal E},A)$ un topos annel\'e.

\smallskip

Alors la cat\'egorie ab\'elienne ${\mathcal M}od_A$ des $A$-modules internes de ${\mathcal E}$ poss\`ede les deux propri\'et\'es:

\medskip
\begin{enumerate}
\item[(10)] Les sous-objets de tout objet $M$ de ${\mathcal M}od_A$ forment un ensemble.

\medskip

\item[(11)] Les objets quotients de tout objet $M$ de ${\mathcal M}od_A$ forment un ensemble.
\end{enumerate}
\end{cor}

\begin{demo}

D'apr\`es le lemme \ref{lemIII68}, il suffit de prouver que ${\mathcal M}od_A$ poss\`ede la propri\'et\'e (10).

\smallskip

Or tout sous-objet $M'$ d'un objet $M$ de ${\mathcal M}od_A$ peut \^etre vu comme un sous-objet de $M$ dans ${\mathcal M}od_A$ puisque le foncteur d'oubli de la structure de $A$-module
$$
{\mathcal M}od_A \longrightarrow {\mathcal E}
$$
respecte les limites, en particulier les monomorphismes.

\smallskip

R\'eciproquement, un sous-objet $M'$ de $M$ dans ${\mathcal E}$ d\'efinit un (unique) sous-objet de $M$ dans ${\mathcal M}od_A$ si et seulement si il est respect\'e par le morphisme d'addition
$$
M \times M \longrightarrow M
$$
et par celui de multiplication par les scalaires
$$
A \times M \longrightarrow M \, .
$$

Ainsi, les sous-objets de $M$ dans ${\mathcal M}od_A$ forment un sous-ensemble de l'ensemble des sous-objets de $M$ dans ${\mathcal E}$.

\end{demo}

\section{Familles s\'eparantes dans les topos et crit\`eres de repr\'esentabilit\'e}\label{sec37}

\subsection{Familles s\'eparantes et familles cos\'eparantes dans les topos}\label{subsec371}

\medskip

Commen\c cons par la d\'efinition suivante:

\begin{defn}\label{defIII71}

Soit ${\mathcal C}$ une cat\'egorie localement petite.

\begin{listeimarge}

\item Une famille $(X_i)_{i \in I}$ d'objets de ${\mathcal C}$ est dite ``s\'eparante'' si, pour tous objets $X,Y$ de ${\mathcal C}$, l'application

\begin{eqnarray}
{\rm Hom} (X,Y) &\longrightarrow &\displaystyle \prod_{i \in I} \ \prod_{f_i \in {\rm Hom} (X_i,X)} {\rm Hom} (X_i,Y) \, , \nonumber \\
 f &\longmapsto &(f \circ f_i)_{i \in I , f_i \in {\rm Hom} (X_i,X)} \nonumber
\end{eqnarray}
est injective.

\medskip

\item Une famille $(X_i)_{i \in I}$ d'objets de ${\mathcal C}$ est dite ``cos\'eparante'' si elle est s\'eparante dans la cat\'egorie oppos\'ee ${\mathcal C}^{\rm op}$ c'est-\`a-dire si, pour tous objets $X,Y$ de ${\mathcal C}$, l'application
\begin{eqnarray}
{\rm Hom} (Y,X) &\longrightarrow &\displaystyle \prod_{i \in I} \ \prod_{f_i \in {\rm Hom} (X,X_i)} {\rm Hom} (Y,X_i) \, , \nonumber \\
g &\longmapsto &(f_i \circ g)_{i \in I , f_i \in {\rm Hom} (X,X_i)} \nonumber
\end{eqnarray}
est injective.
\end{listeimarge}
\end{defn}

On a:

\begin{prop}\label{propIII72}

Soit ${\mathcal E}$ un topos.

\begin{listeimarge}

\item Il poss\`ede les deux propri\'et\'es duales:

\medskip
\begin{enumerate}
\item[(13)] La cat\'egorie ${\mathcal E}$ admet des familles s\'eparantes.

\medskip

\item[(14)] La cat\'egorie ${\mathcal E}$ admet des familles cos\'eparantes.
\end{enumerate}
\medskip

\item Si ${\mathcal E} = \widehat{\mathcal C}_J$ est le topos des faisceaux sur un petit site $({\mathcal C},J)$, une famille s\'eparante de ${\mathcal E}$ est constitu\'ee des images
$$
\ell (X)
$$ par le foncteur canonique $\ell : {\mathcal C} \xrightarrow{ \ y \ } \widehat{\mathcal C} \xrightarrow{ \, j^* \ } \widehat{\mathcal C}_J$ des objets $X$ de la petite cat\'egorie ${\mathcal C}$.

\smallskip

D'autre part, une famille cos\'eparante de ${\mathcal E}$ est constitu\'ee des faisceaux
$$
{\mathcal H}om (\ell (X) , \Omega_J) = \Omega_J^{\ell (X)} \, , \qquad X \in {\rm Ob} ({\mathcal C}) \, ,
$$
o\`u $\Omega_J$ d\'esigne le faisceau classificateur des sous-objets.
\end{listeimarge}
\end{prop}


\begin{demo}
\begin{listeisansmarge}
\item Comme les propri\'et\'es (13) et (14) sont invariantes par \'equivalence de cat\'egories, il suffit de d\'emontrer~(ii).

\medskip

\item La premi\`ere assertion r\'esulte de ce qu'un morphisme de faisceau sur ${\mathcal C}$
$$
F \longrightarrow G
$$
est enti\`erement d\'etermin\'e par les applications associ\'ees
$$
{\rm Hom} (\ell (X),F) = F(X) \longrightarrow G(X) = {\rm Hom} (\ell (X),G)
$$
index\'ees par l'ensemble des objets $X$ de ${\mathcal C}$.

\smallskip

Puis montrons que la famille des
$$
\Omega_J^{\ell (X)} = {\mathcal H}om (\ell (X) , \Omega_J)
$$
est cos\'eparante.

\smallskip

Consid\'erons pour cela deux morphismes de faisceaux
$$
\raisebox{.7ex}{\xymatrix{ F  \dar[r]^-{^{^{\mbox{\scriptsize$u$}}}}_-{v} &G}} \, .
$$
S'ils sont distincts, il existe un objet $X$ de ${\mathcal C}$ et un \'el\'ement $x \in F(X)$ tels que
$$
u(x) \ne v(x) \qquad \mbox{dans} \quad G(X) \, .
$$
L'\'el\'ement $y = u(x) \in G(X)$ peut \^etre vu comme un morphisme
$$
\ell (X) \longrightarrow G \, .
$$
Son graphe est un sous-objet de $\ell (X) \times G$ qui correspond \`a un morphisme
$$
\ell (X) \times G \longrightarrow \Omega_J
$$
ou, ce qui revient au m\^eme, \`a un morphisme
$$
G \longrightarrow \Omega_J^{\ell (X)} \, .
$$
Les compos\'es de ce morphisme avec les deux morphismes $u$ et $v$ sont deux morphismes
$$
F \rightrightarrows \Omega_J^{\ell (X)}
$$
qui correspondent \`a deux morphismes
$$
F \times \ell (X) \rightrightarrows \Omega_J
$$
et donc \`a deux sous-objets de $F \times \ell (X)$.

\smallskip

Le premier est le graphe de $x \in F(X)$ vu comme un morphisme
$$
\ell (X) \longrightarrow F
$$
tandis que le second diff\`ere de ce graphe sans quoi on aurait
$$
y = v(x) \, .
$$
Cela montre que la famille des $\Omega_J^{\ell (X)}$ est cos\'eparante comme annonc\'e. 

\end{listeisansmarge}
\end{demo}

\bigskip

Nous allons voir au paragraphe suivant que les propri\'et\'es (1) et (2) des topos combin\'ees avec les propri\'et\'es (10) et (11) et avec les propri\'et\'es (13) et (14) permettent d'\'etablir un remarquable crit\`ere de repr\'esentabilit\'e dans les topos.

\subsection{Un crit\`ere de repr\'esentabilit\'e}\label{subsec372}

\medskip

Nous allons d\'emontrer le th\'eor\`eme tr\`es g\'en\'eral suivant:

\begin{thm}\label{thmIII73}
\begin{listeimarge}
\item Soit ${\mathcal E}$ un topos ou plus g\'en\'eralement une cat\'egorie localement petite qui poss\`ede les propri\'et\'es suivantes:

\medskip

$
\left\{\begin{matrix}
{\rm (2)} &\mbox{Pour tout carquois $D$ existe dans ${\mathcal E}$ un foncteur de colimite} \hfill \\
{ \ } \\
&\displaystyle\varinjlim_D : D\mbox{\rm -diag} \, ({\mathcal E}) \longrightarrow {\mathcal E} \, . \\
{ \ } \\
{\rm (11)} &\mbox{Les objets quotients de tout objet $E$ de ${\mathcal E}$ forment un ensemble.} \hfill \\
{ \ } \\
{\rm (13)} &\mbox{La cat\'egorie ${\mathcal E}$ admet une famille s\'eparante d'objets $(X_i)_{i \in I}$.} \hfill
\end{matrix} \right.
$

\bigskip

Alors un foncteur contravariant 
$$
F : {\mathcal E}^{\rm op} \longrightarrow {\rm Ens}
$$
est repr\'esentable si et seulement si il transforme les colimites en limites.

\smallskip

Dans ce cas, l'objet qui le repr\'esente peut \^etre construit \`a partir des $X_i$ comme la colimite du diagramme constitu\'e de tous les quotients de la somme $\underset{i \in I}{\coprod} \ \underset{x \in F(x_i)}{\coprod} (X_i , x)$ dans la cat\'egorie ${\mathcal E} / F = \int\!F$ des \'el\'ements de~$F$.

\medskip

\item Soit ${\mathcal E}$ un topos ou plus g\'en\'eralement une cat\'egorie localement petite qui poss\`ede les propri\'et\'es suivantes:

\medskip

$
\left\{\begin{matrix}
{\rm (1)} &\mbox{Pour tout carquois $D$ existe dans ${\mathcal E}$ un foncteur de limite} \hfill \\
{ \ } \\
&\displaystyle \varprojlim_D : D\mbox{\rm -diag} \, ({\mathcal E}) \longrightarrow {\mathcal E} \, . \\
{ \ } \\
{\rm (10)} &\mbox{Les sous-objets de tout objet $E$ de ${\mathcal E}$ forment un ensemble.} \hfill \\
{ \ } \\
{\rm (14)} &\mbox{La cat\'egorie ${\mathcal E}$ admet une famille cos\'eparante d'objets $(X_i)_{i \in I}$.} \hfill
\end{matrix} \right.
$

\bigskip

Alors un foncteur covariant
$$
F : {\mathcal E} \longrightarrow {\rm Ens}
$$
est repr\'esentable si et seulement si il pr\'eserve les limites.

\smallskip

Dans ce cas, l'objet qui le repr\'esente peut \^etre construit \`a partir des $X_i$ par une formule qui utilise seulement des foncteurs de limites.
\end{listeimarge}
\end{thm}

\begin{demosansqed}

Quitte \`a remplacer ${\mathcal E}$ par ${\mathcal E}^{\rm op}$, il suffit de traiter le cas (i) d'un foncteur contravariant
$$
F : {\mathcal E}^{\rm op} \longrightarrow {\rm Ens} \, .
$$

Consid\'erons la cat\'egorie ${\mathcal E}/F = \int\!F$ des \'el\'ements de $F$. Ses objets sont les paires $(X,x)$ constitu\'ees d'un objet $X$ de ${\mathcal E}$ et d'un \'el\'ement $x \in F(X)$, et ses morphismes
$$
(X,x) \longrightarrow (Y,y)
$$
sont les morphismes $f : X \to Y$ de ${\mathcal E}$ tels que $F(f)(y) = x$.

\smallskip

On observe que $F$ est repr\'esentable si et seulement si la cat\'egorie $\int\!F$ admet un objet terminal.

\smallskip

Comme la cat\'egorie ${\mathcal E}$ a des colimites arbitraires et le foncteur $F$ transforme les colimites en limites, la cat\'egorie $\int\!F$ a des colimites arbitraires et le foncteur d'oubli
$$
\begin{matrix}
\hfill \textstyle\int\!F &\longrightarrow &{\mathcal E} \, , \\
(X,x) &\longmapsto &X \hfill
\end{matrix}
$$
respecte les colimites.

\smallskip

En particulier, pour tout \'epimorphisme de $\int\!F$
$$
(X,x) \longrightarrow (Y,y) \, ,
$$
le morphisme induit $X \to Y$ est un \'epimorphisme de ${\mathcal E}$ et l'application $F(Y) \to F(X)$ est injective. Il en r\'esulte que si la cat\'egorie ${\mathcal E}$ poss\`ede la propri\'et\'e (11), il en est de m\^eme de la cat\'egorie $\int\!F$.

\smallskip

Enfin, si une famille $(X_i)_{i \in I}$ d'objets de ${\mathcal E}$ est s\'eparante, alors la famille des paires $(X_i ,x)$, $i \in I$, $x \in F(X_i)$, est s\'eparante dans la cat\'egorie $\int\!F$.

\smallskip

Ainsi, le th\'eor\`eme est ramen\'e au lemme suivant:
\end{demosansqed}

\begin{lem}\label{lemIII74}

Soit ${\mathcal E}$ une cat\'egorie localement petite qui poss\`ede les propri\'et\'es {\rm (2), (11)} et {\rm (13)}.

\smallskip

Alors ${\mathcal E}$ poss\`ede un objet terminal.

\smallskip

De plus, si $(X_i)_{i \in I}$ est une famille s\'eparante d'objets de ${\mathcal E}$, l'objet terminal peut \^etre construit \`a partir des $X_i$ comme la colimite du diagramme constitu\'e de tous les quotients de $\underset{i \in I}{\coprod} X_i$.
\end{lem}

\begin{demo}

Partant d'une famille s\'eparante $(X_i)_{i \in I}$ d'objets de $E$, formons leur somme
$$
S = \coprod_{i \in I} X_i
$$
munie des morphismes canoniques $s_i : X_i \to S$.

\smallskip

Par hypoth\`ese, les quotients $Q$ de $S$ forment un ensemble.

\smallskip

Notons $1$ la colimite du diagramme constitu\'e de l'objet $S$ et de toutes ses fl\`eches de quotients
$$
S \longrightarrow Q \, .
$$

Elle est munie d'un morphisme canonique
$$
s : S \longrightarrow 1
$$
qui est un \'epimorphisme, et tout \'epimorphisme $S \to Q$ s'inscrit dans un triangle commutatif:
$$
\xymatrix{
S \ar[rr] \ar[rd] &&Q \ar[ld] \\
&1
}
$$

Autrement dit, $1$ est le plus petit quotient de $S$.

\smallskip

On pr\'etend que $1$ est un objet terminal de ${\mathcal E}$, c'est-\`a-dire que tout objet $X$ de ${\mathcal C}$ admet une unique fl\`eche
$$
X \longrightarrow 1 \, .
$$

Pour l'unicit\'e, consid\'erons deux fl\`eches
$$
\raisebox{.7ex}{\xymatrix{ X  \dar[r]^-{^{^{\mbox{\scriptsize$f$}}}}_-{g} &1}}.
$$

Leur co\'egalisateur est un quotient de $1$ et donc de $S$. Comme $1$ est le plus petit quotient de $S$, ce co\'egalisateur est n\'ecessairement \'egal \`a $1$, ce qui signifie comme voulu
$$
f=g \, .
$$

Pour l'existence, consid\'erons la somme
$$
T = \coprod_{i \in I} \ \coprod_{f \in {\rm Hom} (X_i , X)} X_i
$$
munie des morphismes canoniques
$$
s_{i,f} : X_i \longrightarrow T \, , \qquad i \in I \, , \quad f \in {\rm Hom} (X_i,X)
$$
et du morphisme
$$
t:T \longrightarrow X
$$
tel que, pour tout $i \in I$ et tout $f \in {\rm Hom} (X_i,X)$,
$$
t \circ s_{i,f} = f \, .
$$

Comme la famille $(X_i)_{i \in I}$ est s\'eparante, le morphisme
$$
t:T \longrightarrow X
$$
est un \'epimorphisme.

\smallskip

Consid\'erons alors le morphisme
$$
p : T \longrightarrow 1
$$
tel que
$$
p \circ s_{i,f} = s \circ s_i \, , \qquad \forall \, i \in I \, , \quad \forall \, f \in {\rm Hom} (X_i,X) \, .
$$

Puis formons le carr\'e cocart\'esien:
$$
\xymatrix{
T \ar[d]_-t \ar[r]^-p &1 \ar[d] \\
X \ar[r] &Q
}
$$

Comme $t : T \to X$ est un \'epimorphisme, il en est de m\^eme de la fl\`eche
$$
1 \longrightarrow Q \, .
$$

Or, comme $1$ est le plus petit quotient de $S$, un tel \'epimorphisme
$$
1 \longrightarrow Q
$$
est n\'ecessairement un isomorphisme, et donc on a construit un morphisme
$$
X \longrightarrow 1 \, .
$$

Cela termine la d\'emonstration du lemme, donc aussi celle du th\'eor\`eme \ref{thmIII73}. 

\end{demo}

\newpage

On d\'eduit de ce th\'eor\`eme:

\begin{cor}\label{corIII75}

Soit
$$
F : {\mathcal E} \longrightarrow {\mathcal D}
$$
un foncteur entre deux cat\'egories localement petites.

\begin{listeimarge}

\item Supposons que ${\mathcal E}$ est un topos, ou plus g\'en\'eralement qu'elle poss\`ede les propri\'et\'es {\rm (2), (11)} et {\rm (13)}.

\smallskip

Alors $F$ admet un adjoint \`a droite si (et seulement si) il pr\'eserve les colimites.

\medskip

\item Supposons que ${\mathcal E}$ est un topos, ou plus g\'en\'eralement qu'elle poss\`ede les propri\'et\'es {\rm (1), (10)} et {\rm (14)}.

\smallskip

Alors $F$ admet un adjoint \`a gauche si (et seulement si) il pr\'eserve les limites.
\end{listeimarge}
\end{cor}

\begin{demo}
\begin{listeisansmarge}
\item On sait d\'ej\`a que la condition est n\'ecessaire.

\smallskip

Pour la suffisance, on sait d'apr\`es la remarque (ii) qui suit la d\'efinition \ref{defI81} que
$$
F : {\mathcal E} \longrightarrow {\mathcal D}
$$
a un adjoint \`a droite $G$ si (et seulement si), pour tout objet $Y$ de ${\mathcal D}$, le foncteur
$$
\begin{matrix}
{\mathcal E}^{\rm op} &\longrightarrow &{\rm Ens} \, , \hfill \\
\hfill X &\longmapsto &{\rm Hom} (F(X),Y)
\end{matrix}
$$
est repr\'esentable par un objet $G(Y)$ de ${\mathcal E}$ et qu'il est possible de choisir les $G(Y)$ uniform\'ement en $Y$.

\smallskip

Or, si $F$ pr\'eserve les colimites, chaque foncteur
$$
X \longmapsto {\rm Hom} (F(X),Y)
$$
transforme les colimites en limites, et le th\'eor\`eme \ref{thmIII73} (i) s'applique. Ces foncteurs sont repr\'esentables par des objets $G(Y)$ construits par des formules, donc de mani\`ere uniforme en $Y$.

\medskip

\item r\'esulte de (i) en voyant $F$ comme un foncteur ${\mathcal E}^{\rm op} \to {\mathcal D}^{\rm op}$. 

\end{listeisansmarge}
\end{demo}

\subsection{Application aux cat\'egories d'objets lin\'eaires des topos}\label{subsec373}

\medskip

Montrons encore comme cons\'equence du th\'eor\`eme \ref{thmIII73}:

\begin{cor}\label{corIII76}

Soit $({\mathcal E},A)$ un topos annel\'e.

\smallskip

Soit ${\mathcal M}od_A$ la cat\'egorie des $A$-modules internes de ${\mathcal E}$.

\smallskip

Alors:

\begin{listeimarge}

\item Un foncteur contravariant
$$
F : {\mathcal M}od_A^{\rm op} \longrightarrow {\rm Ens}
$$
est repr\'esentable si (et seulement si) il transforme les colimites en limites.

\medskip

\item Un foncteur vers une cat\'egorie localement petite ${\mathcal A}$
$$
F : {\mathcal M}od_A \longrightarrow {\mathcal A}
$$
admet un adjoint \`a droite si (et seulement si) il respecte les colimites.
\end{listeimarge}
\end{cor}

\begin{demosansqed}

On sait d\'ej\`a que la cat\'egorie ${\mathcal M}od_A$ poss\`ede la propri\'et\'e (2) -- l'existence de colimites arbitraires -- et la propri\'et\'e (11): les quotients de tout objet forment un ensemble.

\smallskip

Pour conclure en appliquant le th\'eor\`eme \ref{thmIII73} (i), il suffit de conna{\^\i}tre le lemme suivant:

\end{demosansqed}

\begin{lem}\label{lemIII77}

Soit $({\mathcal E},A)$ un topos annel\'e.

\smallskip

Alors la cat\'egorie ${\mathcal M}od_A$ des $A$-modules internes de ${\mathcal E}$ poss\`ede la propri\'et\'e:

\medskip
\begin{enumerate}
\item[(13)] Elle admet une famille s\'eparante d'objets.
\end{enumerate}
\end{lem}

\begin{demo}

Il suffit de traiter le cas o\`u ${\mathcal E} = \widehat{\mathcal C}_J$ est le topos des faisceaux sur un petit site $({\mathcal C},J)$.

\smallskip

On observe que, pour tout objet $X$ de ${\mathcal C}$, le foncteur sur la cat\'egorie ${\mathcal M}od_{{\mathcal C},A}$ des pr\'efaisceaux de $A$-modules
$$
\begin{matrix}
{\mathcal M}od_{{\mathcal C},A} &\longrightarrow &{\rm Ens} \, , \hfill \\
\hfill P &\longmapsto &P(X) \hfill
\end{matrix}
$$
est repr\'esentable par le pr\'efaisceau de $A$-modules
$$
\begin{matrix}
A_X : &\hfill U &\longmapsto &\displaystyle \bigoplus_{f \in {\rm Hom} (U,X)} A(U) \, , \hfill \\
{ \ } \\
&\left( V \xrightarrow{ \ h \ } U \right) &\longmapsto &\left[ \begin{matrix}
\displaystyle \bigoplus_{f \in {\rm Hom} (U,X)} A(U) &\longrightarrow &\displaystyle \bigoplus_{g \in {\rm Hom} (V,X)} A(V) \hfill \\
{ \ } \\
\hfill (a_f)_{f \in {\rm Hom} (U,X)} &\longmapsto &\Biggl( \displaystyle \sum_{f \in {\rm Hom} (U,X) \atop f \circ h = g} A(h) (a_f) \Biggl)_{\!\!g \in {\rm Hom} (V,X)}
\end{matrix}\right].
\end{matrix}
$$

Le foncteur de faisceautisation
$$
j^* : {\mathcal M}od_{{\mathcal C},A} \longrightarrow {\mathcal M}od_{{\mathcal C},J,A} = {\mathcal M}od_A
$$
est adjoint \`a gauche du foncteur de plongement
$$
j_* : {\mathcal M}od_A = {\mathcal M}od_{{\mathcal C},J,A} \longrightarrow {\mathcal M}od_{{\mathcal C},A}
$$
et donc le foncteur
$$
\begin{matrix}
{\mathcal M}od_A &\longrightarrow &{\rm Mod}_{A(X)} \, , \\
\hfill F &\longmapsto &F(X) \hfill
\end{matrix}
$$
est repr\'esentable par le faisceau de $A$-modules
$$
j^* A_X \, .
$$

Alors la famille des objets
$$
(j^* A_X)_{X \in {\rm Ob} ({\mathcal C})}
$$
est s\'eparante dans la cat\'egorie ${\mathcal M}od_A$. 

\smallskip

Cela d\'emontre le lemme et donc aussi le corollaire \ref{corIII76}. 

\end{demo}


On sait d'autre part que les cat\'egories ${\mathcal M}od_A$ de modules internes sur un anneau interne $A$ d'un topos poss\`edent la propri\'et\'e (1) -- l'existence de limites arbitraires -- et la propri\'et\'e (10): les sous-objets de tout objet forment un ensemble.

\smallskip

Afin de montrer qu'elles satisfont aussi la propri\'et\'e (14) -- l'existence de familles cos\'eparantes d'objets -- on a besoin de montrer qu'elles ont ``assez d'injectifs'' au sens de la d\'efinition suivante:

\begin{defn}\label{defIII78}

Soit ${\mathcal C}$ une cat\'egorie localement petite.

\begin{listeimarge}

\item Un objet $I$ de ${\mathcal C}$ est dit ``injectif'' si, pour tout monomorphisme de ${\mathcal C}$,
$$
Z \xhookrightarrow{ \ { \ } \ } X \, ,
$$
l'application induite
$$
{\rm Hom} (X,I) \longrightarrow {\rm Hom} (Z,I)
$$
est surjective.

\medskip

\item On dit que la cat\'egorie ${\mathcal C}$ a ``assez d'injectifs'' si tout objet $X$ de ${\mathcal C}$ admet au moins un monomorphisme
$$
X \xhookrightarrow{ \ { \ } \ } I
$$
dans un objet injectif $I$ de ${\mathcal C}$.
\end{listeimarge}
\end{defn}

\begin{remarksqed}
\begin{listeisansmarge}
\item On dit qu'un objet $P$ de ${\mathcal C}$ est ``projectif'' s'il est injectif dans ${\mathcal C}^{\rm op}$ c'est-\`a-dire si tout \'epimorphisme de ${\mathcal C}$
$$
X \, -\!\!\!\twoheadrightarrow Y
$$
induit une application surjective
$$
{\rm Hom} (P,X) \longrightarrow {\rm Hom} (P,Y) \, .
$$

\item On dit que ${\mathcal C}$ a ``assez de projectifs'' si ${\mathcal C}^{\rm op}$ a assez d'injectifs, c'est-\`a-dire si tout objet $X$ de ${\mathcal C}$ s'\'ecrit comme un quotient
$$
P \, -\!\!\!\twoheadrightarrow X
$$
d'un objet projectif $P$ de ${\mathcal C}$.

\medskip

\item Un objet $I$ d'une cat\'egorie additive ${\mathcal A}$ est injectif si et seulement si le foncteur
$$
{\rm Hom} (\bullet,I) : {\mathcal A}^{\rm op} \longrightarrow {\rm Ab}
$$
transforme les limites finies en colimites finies.

\smallskip

Si ${\mathcal A}$ est une cat\'egorie ab\'elienne, c'est \'equivalent \`a demander que ${\rm Hom} (\bullet , I)$ transforme toute suite exacte de ${\mathcal A}$
$$
\cdots \longrightarrow M_{i-1} \longrightarrow M_i \longrightarrow M_{i+1} \longrightarrow \cdots
$$
en une suite exacte de ${\rm Ab}$:
$$
\cdots \longrightarrow {\rm Hom} (M_{i+1} , I) \longrightarrow {\rm Hom} (M_i,I) \longrightarrow {\rm Hom} (M_{i-1} , I) \longrightarrow \cdots
$$

\item De m\^eme, un objet $P$ d'une cat\'egorie additive ${\mathcal A}$ est projectif si et seulement si le foncteur
$$
{\rm Hom} (P,\bullet) : {\mathcal A} \longrightarrow {\rm Ab}
$$
respecte les colimites finies.

\smallskip

Si ${\mathcal A}$ est une cat\'egorie ab\'elienne, c'est \'equivalent \`a demander que ${\rm Hom} (P,\bullet)$ transforme toute suite exacte de ${\mathcal A}$
$$
\cdots \longrightarrow M_{i-1} \longrightarrow M_i \longrightarrow M_{i+1} \longrightarrow \cdots
$$
en une suite exacte de ${\rm Ab}$:
$$
\cdots \longrightarrow {\rm Hom} (P,M_{i-1}) \longrightarrow {\rm Hom} (P,M_i) \longrightarrow {\rm Hom} (P,M_{i+1}) \longrightarrow \cdots
$$

\end{listeisansmarge}
\end{remarksqed}

On a le r\'esultat important:

\begin{thm}\label{thmIII79}

Soit ${\mathcal A}$ la cat\'egorie des modules internes d'un topos annel\'e $({\mathcal E}, A)$ ou, plus g\'en\'eralement, une cat\'egorie ab\'elienne qui poss\`ede les propri\'et\'es suivantes:

\bigskip

$\left\{\begin{matrix}
{\rm (2)} &\mbox{Elle admet des colimites arbitraires.} \hfill \\
{ \ } \\
{\rm (6')} &\mbox{Pour tout morphisme $N \to M$ de ${\mathcal A}$, le foncteur $N \times_M \bullet$ respecte les colimites filtrantes.} \hfill \\
{ \ } \\
{\rm (10)} &\mbox{Les sous-objets de tout objet de ${\mathcal A}$ forment un ensemble.} \hfill \\
{ \ } \\
{\rm (13)} &\mbox{Elle admet une famille s\'eparante d'objets $(U_i)_{i \in I}$.} \hfill
\end{matrix} \right.$

\bigskip

Alors la cat\'egorie ${\mathcal A}$ a assez d'injectifs.
\end{thm}

\begin{remarks}
\begin{listeisansmarge}
\item Comme ${\mathcal A}$ est une cat\'egorie ab\'elienne, la propri\'et\'e (2) revient \`a demander que ${\mathcal A}$ a des sommes arbitraires.

\medskip

\item La propri\'et\'e ($6'$) est un cas particulier de la propri\'et\'e suivante que poss\`edent les cat\'egories de modules de topos annel\'es:

\medskip

$\left\{\begin{matrix}
(7) &\mbox{Pour toute petite cat\'egorie filtrante ${\mathcal D}$, le foncteur $\underset{\mathcal D}{\varprojlim}$ respecte les limites finies.} \hfill
\end{matrix} \right.$

\medskip

\item En revanche, les cat\'egories ${\mathcal M}od_A$ de modules internes de topos annel\'es $({\mathcal E},A)$ n'ont en g\'en\'eral pas assez de projectifs.

\medskip

\item Cependant, ${\mathcal M}od_A$ a assez de projectifs si le topos ${\mathcal E}$ est \'equivalent au topos $\widehat{\mathcal C}$ des pr\'efaisceaux sur une petite cat\'egorie ${\mathcal C}$.

\smallskip

En effet, pour tout objet $X$ de ${\mathcal C}$, l'objet $A_X$ de ${\mathcal M}od_{{\mathcal C},A}$ qui repr\'esente le foncteur d'\'evaluation
$$
\begin{matrix}
{\mathcal M}od_{{\mathcal C},A} &\longrightarrow &{\rm Ab} \, , \hfill \\
\hfill P &\longmapsto &P(X)
\end{matrix}
$$
est projectif puisque ce foncteur respecte les colimites.

\smallskip

Alors toute somme d'objets de la forme $A_X$ est encore un objet projectif de ${\mathcal M}od_{{\mathcal C},A}$ et tout objet $P$ de ${\mathcal M}od_{{\mathcal C},A}$ s'\'ecrit comme un quotient de l'objet projectif
$$
\bigoplus_{X \in {\rm Ob} ({\mathcal C})} \, \bigoplus_{x \in P(X)} A_X \, .
$$
\end{listeisansmarge}
\end{remarks}

\begin{demosansqed}

Commen\c cons par d\'emontrer le lemme suivant:

\end{demosansqed}

\begin{lem}\label{lemIII710}

Soit ${\mathcal A}$ une cat\'egorie ab\'elienne qui poss\`ede les propri\'et\'es {\rm (2), ($6'$), (10)} et {\rm (13)} du th\'eor\`eme.

\smallskip

Soit $U = \underset{i \in I}{\bigoplus} \ U_i$ la somme des objets d'une famille s\'eparante $(U_i)_{i \in I}$.

\smallskip

Alors:

\begin{listeimarge}

\item Pour qu'un objet $I$ de ${\mathcal A}$ soit injectif, il suffit (et il faut) que pour tout sous-objet $U'$ de $U$, tout morphisme
$$
U' \longrightarrow I
$$
se prolonge en un morphisme
$$
U \longrightarrow I \, .
$$

\item Il existe un foncteur (non additif)
$$
\begin{matrix}
{\mathcal A} &\longrightarrow &{\mathcal A} \, , \\
M &\longmapsto &M_1
\end{matrix}
$$
tel que tout morphisme
$$
U' \longrightarrow M
$$
d'un sous-objet $U'$ de $U$ dans un objet $M$ de ${\mathcal A}$ se prolonge en un morphisme
$$
U \longrightarrow M_1 \, .
$$
\end{listeimarge}
\end{lem}

\begin{demo}
\begin{listeisansmarge}
\item Cette condition est n\'ecessaire par d\'efinition de la notion d'objet injectif.

\smallskip

R\'eciproquement, supposons qu'elle soit v\'erifi\'ee.

\smallskip

Consid\'erons un objet $M$ de ${\mathcal A}$, un sous-objet $M'$ de $M$ et un morphisme $M' \xrightarrow{ \ f \ } I$.

\smallskip

Soit $S$ l'ensemble des paires $(M_1 , f_1)$ constitu\'ees d'un sous-objet $M_1$ de $M$ qui contient $M'$ et d'un morphisme $M_1 \xrightarrow{ \ f_1 \ } I$ qui prolonge $f$.

\smallskip

On d\'efinit une relation d'ordre sur $S$ en d\'ecidant que
$$
(M_1 , f_1) \leq (M_2 , f_2)
$$
si $M_2$ contient $M_1$ et $f_2 : M_2 \to I$ prolonge $f_1 : M_1 \to I$.

\smallskip

Pour tout sous-ensemble totalement ordonn\'e $S'$ de $S$, la colimite
$$
M_2 = \varinjlim_{(M_1 , f_1) \in S'} M_1
$$
est un sous-objet de $M$ muni d'un morphisme
$$
f_2 : M_2 \longrightarrow I
$$
qui prolonge tous les $f_1 : M_1 \to I$, $(M_1 , f_1) \in S'$.

\smallskip

D'apr\`es le lemme de Zorn, $S$ poss\`ede un \'el\'ement maximal $(M_0 , f_0)$.

\smallskip

Montrons que $M_0 = M$.

\smallskip

Pour tout morphisme
$$
U \longrightarrow M \, ,
$$
$U_0 = U \times_M M_0$ est un sous-objet de $U$, et le compos\'e
$$
U_0 \longrightarrow M_0 \xrightarrow{ \ f_0 \ } I
$$
se prolonge par hypoth\`ese en un morphisme
$$
f : U \longrightarrow I \, .
$$
La paire $(f,f_0)$ d\'efinit un morphisme du sous-objet
$$
(U \oplus M_0) / (U \times_M M_0) \quad \mbox{de} \quad M
$$
vers $I$ qui prolonge $f_0$.

\smallskip

Par maximalit\'e de la paire $(M_0 , f_0)$ dans $S$, cela impose que le morphisme
$$
U \longrightarrow M
$$
se factorise n\'ecessairement par $M_0$.

\smallskip

Comme $U$ est la somme des objets de la famille s\'eparante $(U_i)_{i \in I}$, cela impose comme voulu
$$
M_0 = M \, .
$$

\item Par hypoth\`ese, les sous-objets $U'$ de l'objet $U = \underset{i \in I}{\bigoplus} \ U_i$ de ${\mathcal A}$ forment un ensemble $\Omega (U)$.

\smallskip

Pour tout objet $M$ de ${\mathcal A}$, notons $M_1$ le conoyau du morphisme
$$
\bigoplus_{U' \in \Omega (U)} \bigoplus_{f \in {\rm Hom} (U',M)} U' \longrightarrow M \oplus \left( \bigoplus_{U' \in \Omega (U)} \bigoplus_{f \in {\rm Hom} (U',M)} U \right)
$$
dont la premi\`ere composante
$$
\bigoplus_{U' \in \Omega (U)} \bigoplus_{f \in {\rm Hom} (U',M)} U' \longrightarrow M
$$
est la somme des morphismes
$$
f : U' \longrightarrow M \, , \qquad U' \in \Omega (U) \, , \quad f \in {\rm Hom} (U',M) \, ,
$$
et dont la seconde composante
$$
\bigoplus_{U' \in \Omega (U)} \bigoplus_{f \in {\rm Hom} (U',M)} U' \longrightarrow \bigoplus_{U' \in \Omega (U)} \bigoplus_{f \in {\rm Hom} (U',M)} U
$$
est la somme des monomorphismes
$$
U' \xhookrightarrow{ \ { \ } \ } U
$$
index\'es par les $U' \in \Omega (U)$ et les $f \in {\rm Hom} (U',M)$.

\smallskip

Par construction, $M_1$ est muni d'un monomorphisme
$$
M \xhookrightarrow{ \ { \ } \ } M_1
$$
et, pour tout $U' \in \Omega(U)$ et tout $f \in {\rm Hom} (U',M)$, d'un morphisme
$$
U \longrightarrow M_1
$$
qui rend commutatif le carr\'e:
$$
\xymatrix{
U' \, \ar[d]_-f \ar@{^{(}->}[r] &U \ar[d] \\
M \, \ar@{^{(}->}[r] &M_1
}
$$

Cela termine la preuve du lemme. 
\end{listeisansmarge}
\end{demo}

\bigskip

\noindent {\bf Fin de la d\'emonstration du th\'eor\`eme \ref{thmIII79}:}

\smallskip

Soit $M$ un objet de ${\mathcal A}$.

\smallskip

Construisons \`a partir de $M = M_0$ un syst\`eme inductif d'objets de ${\mathcal A}$
$$
M_i \quad \mbox{index\'es par les ordinaux $i$}
$$
et reli\'es par des monomorphismes $M_i \hookrightarrow M_j$ pour $i \leq j$.

\smallskip

La construction est par r\'ecurrence transfinie:

\medskip

$\left\{\begin{matrix}
\bullet &\mbox{si $j = i+1$, $M_j$ est d\'eduit de $M_i$ par le proc\'ed\'e du lemme \ref{lemIII710} (ii),} \hfill \\
{ \ } \\
\bullet &\mbox{si $j$ est la limite des $i < j$, on pose} \hfill \\
{ \ } \\
&M_j = \displaystyle \varinjlim_{i < j} M_i \, .
\end{matrix}\right.$

\medskip

Soit $k$ un ordinal dont la puissance est strictement plus grande que celle de l'ensemble $\Omega (U)$ des sous-objets de $U$ et qui est la limite des $i < k$.

\smallskip

Montrons que l'objet $M_k$ de ${\mathcal A}$ est injectif.

\smallskip

Pour cela, consid\'erons un sous-objet $U'$ de $U$ et un morphisme
$$
U' \longrightarrow M_k \, .
$$

La formule
$$
M_k = \varinjlim_{i < k} M_i
$$
implique
$$
U' = \varinjlim_{i < k} \ U' \times_{M_k} M_i \, .
$$

Comme la puissance de $k$ est strictement plus grande que celle de l'ensemble $\Omega (U)$ des sous-objets de $U'$, il existe un ordinal $i < k$ tel que
$$
U' \times_{M_i} M_k = U'
$$
c'est-\`a-dire que le morphisme
$$
U' \longrightarrow M_k
$$
se factorise en
$$
U' \longrightarrow M_i \, .
$$

Il se prolonge donc comme voulu en un morphisme
$$
U' \longrightarrow M_{i+1} \xhookrightarrow{ \ { \ } \ } M_k \, .
$$
\hfill $\Box$

On d\'eduit du th\'eor\`eme \ref{thmIII79}:

\begin{cor}\label{corIII711}

Soit ${\mathcal A}$ la cat\'egorie des modules internes ${\mathcal M}od_A$ d'un topos annel\'e $({\mathcal E},A)$ ou, plus g\'en\'eralement, une cat\'egorie ab\'elienne qui poss\`ede les propri\'et\'es {\rm (2), ($6'$), (10)} et {\rm (13)}.

\smallskip

Alors la cat\'egorie ${\mathcal A}$ admet une famille cos\'eparante d'objets.
\end{cor}

\begin{demo}

Soit $(U_i)_{i \in I}$ une famille s\'eparante d'objets de ${\mathcal A}$.

\smallskip

Alors l'objet
$$
U = \bigoplus_{i \in I} U_i
$$
est s\'eparant au sens que pour toute paire de morphismes de ${\mathcal A}$ suppos\'es distincts
$$
 \raisebox{.7ex}{\xymatrix{ M  \dar[r]^-{^{^{\mbox{\scriptsize$f$}}}}_-{g} &N}},
$$
il existe un morphisme
$$
u : U \longrightarrow M
$$
tel que $f \circ u \ne g \circ u$.

\smallskip

Soit $K$ l'ensemble des objets quotients $Q_k$ de $U \oplus U$.

\smallskip

D'apr\`es le th\'eor\`eme \ref{thmIII79}, pour tout $k \in K$, l'objet $Q_k$ admet un monomorphisme
$$
Q_k \xhookrightarrow{ \ { \ } \ } I_k
$$
dans un objet injectif $I_k$ de ${\mathcal A}$.

\smallskip

On pr\'etend que la famille des objets $I_k$, $k \in K$, est cos\'eparante.

\smallskip

Consid\'erons pour cela une paire de morphismes de ${\mathcal A}$
$$
 \raisebox{.7ex}{\xymatrix{ M  \dar[r]^-{^{^{\mbox{\scriptsize$f$}}}}_-{g} &N}} .
$$
S'ils sont distincts, il existe un morphisme
$$
u : U \longrightarrow M
$$
tel que $f \circ u \ne g \circ u$.

\smallskip

La paire $(f \circ u , g \circ u)$ d\'efinit un morphisme
$$
f \circ u + g \circ u : U \oplus U \longrightarrow N
$$
qui se factorise en un \'epimorphisme
$$
U \oplus U -\!\!\!\twoheadrightarrow Q_k
$$
suivi d'un monomorphisme
$$
Q_k \xhookrightarrow{ \ { \ } \ } N \, .
$$

Le morphisme
$$
Q_k \xhookrightarrow{ \ { \ } \ } I_k
$$
se prolonge en un morphisme
$$
n : N \longrightarrow I_k
$$
tel que
$$
n \circ f \ne n \circ g
$$
puisque 
$$
n \circ f \circ u \ne n \circ g \circ u \, .
$$

C'est ce que l'on voulait. 

\end{demo}

\medskip

Ce corollaire permet de d\'eduire du th\'eor\`eme \ref{thmIII73} (ii):

\begin{cor}\label{corIII712}

Soit $({\mathcal E},A)$ un topos annel\'e.

\smallskip

Soit ${\mathcal M}od_A$ la cat\'egorie des $A$-modules internes de ${\mathcal E}$.

Alors:

\begin{listeimarge}

\item Un foncteur covariant
$$
F : {\mathcal M}od_A \longrightarrow {\rm Ens}
$$
est repr\'esentable si (et seulement si) il respecte les limites.

\medskip

\item Un foncteur vers une cat\'egorie localement petite ${\mathcal A}$
$$
F : {\mathcal M}od_A \longrightarrow {\mathcal A}
$$
admet un adjoint \`a gauche si (et seulement si) il respecte les limites.
\end{listeimarge}
\end{cor}

\begin{demo}

On sait d\'ej\`a que la cat\'egorie ${\mathcal M}od_A$ poss\`ede la propri\'et\'e (1) -- l'existence de limites arbitraires -- et la propri\'et\'e (10): les sous-objets de tout objet forment un ensemble.

\smallskip

D'apr\`es le corollaire \ref{corIII711}, elle poss\`ede aussi la propri\'et\'e (14): elle admet une famille d'objets cos\'eparante.

\smallskip

Les hypoth\`eses du th\'eor\`eme \ref{thmIII79} (ii) sont v\'erifi\'ees, donc aussi sa conclusion.
\end{demo}

\section{Caract\'erisation des topos par leurs propri\'et\'es}\label{sec38}

\subsection{Enonc\'e du th\'eor\`eme de Giraud}\label{subsec381}

\medskip

Les topos ont \'et\'e d\'efinis de mani\`ere constructive: ce sont les cat\'egories qu'il est possible de construire, \`a \'equivalence pr\`es, comme cat\'egories des faisceaux sur un site.

\smallskip

Le th\'eor\`eme suivant montre que les topos ont aussi une d\'efinition axiomatique: ce sont les cat\'egories qui poss\`edent une certaine liste de propri\'et\'es.

\begin{thm}\label{thmIII81}

Soit ${\mathcal E}$ une cat\'egorie qui poss\`ede les propri\'et\'es suivantes:

\medskip

$
\left\{\begin{matrix}
(0) &\mbox{Elle est localement petite.} \hfill \\
{ \ } \\
(1') &\mbox{Elle admet des limites finies arbitraires.} \hfill \\
{ \ } \\
(2) &\mbox{Elle admet des colimites arbitraires.} \hfill \\
{ \ } \\
(4) &\mbox{Toute relation d'\'equivalence $R \hookrightarrow E \times E$ sur un objet $E$ de ${\mathcal E}$ est effective,} \hfill \\
&\mbox{c'est-\`a-dire d\'efinit un quotient $E \twoheadrightarrow Q$ tel que $R = E \times_Q E$.} \hfill \\
{ \ } \\
(5) &\mbox{R\'eciproquement, tout \'epimorphisme $E \twoheadrightarrow Q$ de ${\mathcal E}$ fait de $Q$ le quotient de $E$} \hfill \\ 
&\mbox{par la relation d'\'equivalence $R = E \times_Q E$.} \hfill \\
{ \ } \\
(6) &\mbox{Pour tout morphisme $E' \to E$ de ${\mathcal E}$, le foncteur de changement de base} \hfill \\
&\mbox{$E' \times_E \bullet$ respecte les colimites.} \hfill \\
{ \ } \\
(13) &\mbox{Elle admet une famille s\'eparante d'objets.} \hfill
\end{matrix} \right.
$

\bigskip

Alors:

\begin{listeimarge}

\item ${\mathcal E}$ est un topos.

\medskip

\item Plus pr\'ecis\'ement, si ${\mathcal C}$ est une petite sous-cat\'egorie pleine de ${\mathcal E}$ dont les objets forment une famille s\'eparante, et que $J$ est la topologie de ${\mathcal C}$ pour laquelle un crible $S$ d'un objet $X$ est couvrant si et seulement si
$$
X' = \varinjlim_{(U' \to X') \in f^* S} U' \quad \mbox{dans ${\mathcal E}$, pour toute fl\`eche $X' \xrightarrow{ \ f \ } X$ de ${\mathcal C}$,}
$$
alors la cat\'egorie ${\mathcal E}$ est \'equivalente au topos $\widehat{\mathcal C}_J$ des faisceaux sur le site $({\mathcal C},J)$.

\medskip

\item Plus pr\'ecis\'ement encore, l'\'equivalence de {\rm (ii)} est d\'efinie par les deux foncteurs
$$
\begin{matrix}
H : &{\mathcal E} &\longrightarrow &\widehat{\mathcal C}_J \, , \hfill \\
&E &\longmapsto &{\rm Hom} (\bullet , E) = \left\{ \begin{matrix}
{\mathcal C}^{\rm op} &\longrightarrow &{\rm Ens} \, , \hfill \\
\hfill X &\longmapsto &{\rm Hom}_{\mathcal E} (X,E) \, ,
\end{matrix} \right.
\end{matrix}
$$
et

\hglue34mm$
\begin{matrix}
G : &\widehat{\mathcal C}_J &\longrightarrow &{\mathcal E}, \hfill \\
&\hfill F &\longmapsto &\displaystyle \varinjlim_{(X,x) \in \int\!\!F} X \, .
\end{matrix}
$

\medskip

\noindent En particulier, le foncteur $H$ pr\'eserve les colimites.
\end{listeimarge}
\end{thm}

\begin{remarks}
\begin{listeisansmarge}
\item Dans l'\'enonc\'e originel du th\'eor\`eme de Giraud au d\'ebut du chapitre IV de SGA~4, les hypoth\`eses donn\'ees sont l\'eg\`erement diff\'erentes: en effet, la propri\'et\'e (5) est remplac\'ee par la propri\'et\'e (8) que dans ${\mathcal E}$ les sommes sont disjointes.

\smallskip

Nous avons choisi cet \'enonc\'e plut\^ot que celui de Giraud car il a une d\'emonstration naturelle plus concr\`ete.

\smallskip

De toute fa\c con, le plus important est que les topos aient, en plus de leur d\'efinition constructive, une d\'efinition axiomatique qui les caract\'erise par une liste de propri\'et\'es cat\'egoriques.

\medskip

\item Ce th\'eor\`eme montre en particulier que tout topos non trivial a tellement de pr\'esentations diff\'erentes qu'elles ne forment m\^eme pas un ensemble.

\end{listeisansmarge}
\end{remarks}

\begin{demo}

Elle va occuper les paragraphes suivants. 
\end{demo}

\subsection{D\'efinition de la topologie et du foncteur vers les faisceaux}\label{subsec382}

\medskip

Sous les hypoth\`eses du th\'eor\`eme \ref{thmIII81}, on consid\`ere une petite sous-cat\'egorie pleine ${\mathcal C}$ de ${\mathcal E}$ dont les objets forment une famille s\'eparante de ${\mathcal E}$.

\smallskip

On v\'erifie d'abord:

\begin{lem}\label{lemIII82}

Pour tout objet $X$ de ${\mathcal C}$, notons $J(X)$ l'ensemble des cribles $S$ de $X$ tels que, pour toute fl\`eche $X' \xrightarrow{ \ f \ } X$ de ${\mathcal C}$ de but $X$, on ait
$$
X' = \varprojlim_{(U' \to X') \in f^* S} U' \qquad \mbox{dans} \quad {\mathcal E} \, .
$$

Alors $J$ est une topologie de Grothendieck de ${\mathcal C}$.
\end{lem}

\begin{demo}

Il faut montrer que $J$ v\'erifie les trois axiomes de maximalit\'e (M), de stabilit\'e (S) et de transitivit\'e (T) de la d\'efinition \ref{defII21}.

\smallskip

C'est \'evident pour la maximalit\'e.

\smallskip

L'axiome de stabilit\'e est satisfait par d\'efinition de $J$ puisque, pour toutes fl\`eches $X'' \xrightarrow{ \ g \ } X' \xrightarrow{ \ f \ } X$ de ${\mathcal C}$ et tout crible $S$ de $X$, on a
$$
g^* (f^* S) = (f \circ g)^* S \, .
$$

Pour la transitivit\'e, consid\'erons deux cribles $S$ et $S'$ d'un objet $X$ de ${\mathcal C}$ tels que $S \in J(X)$ et
$$
h^* S' \in J(U) \, , \qquad \forall \, (h : U \longrightarrow X) \in S \, .
$$

Consid\'erant un morphisme $f : X' \to X$ et un objet $Y$ de ${\mathcal C}$, il faut montrer que toute famille compatible de morphismes de ${\mathcal E}$
$$
U' \longrightarrow Y \quad \mbox{index\'es par les} \ (U' \longrightarrow X') \in f^* S'
$$
provient d'un unique morphisme
$$
X' \longrightarrow Y \, .
$$

Or, pour tout morphisme $U' \xrightarrow{ \ g \ } X'$ de $f^* S$, on a dans ${\mathcal E}$
$$
\varinjlim_{(U'' \to U') \in (f \circ g)^* S'} U'' = U'
$$
et donc les morphismes
$$
U'' \longrightarrow Y \, , \qquad (U'' \longrightarrow U') \in (f \circ g)^* S' \, ,
$$
proviennent d'un unique morphisme
$$
U' \longrightarrow Y \, .
$$

Comme on a aussi
$$
\varinjlim_{(U' \xrightarrow{g} X') \in f^* S} U' = X' \qquad \mbox{dans} \quad {\mathcal E} \, ,
$$
la famille des morphismes
$$
U' \longrightarrow Y \, , \qquad (U' \xrightarrow{ \ g \ } X') \in f^* S
$$
provient d'un unique morphisme
$$
X' \longrightarrow Y \, .
$$
Pour tout $(U' \xrightarrow{ \ g \ } X') \in f^* S'$, le triangle
$$
\xymatrix{
U' \ar[rd] \ar[rr]^-g &&X' \ar[ld] \\
&Y
}
$$
est commutatif car pour tout $(U'' \xrightarrow{ \ h \ } U') \in (f \circ g)^* S$, le triangle induit
$$
\xymatrix{
U'' \ar[rd] \ar[rr]^-{g\circ h} &&X' \ar[ld] \\
&Y
}
$$
est commutatif.

\smallskip

Ainsi, l'axiome de transitivit\'e est bien satisfait. 

\end{demo}

\bigskip

La topologie $J$ sur ${\mathcal C}$ \'etant bien d\'efinie, on a aussit\^ot:

\begin{cor}\label{corIII83}

Pour tout objet $E$ de ${\mathcal E}$, le pr\'efaisceau sur ${\mathcal C}$
$$
\begin{matrix}
{\rm Hom} (\bullet , E) : {\mathcal C}^{\rm op} &\longrightarrow &{\rm Ens} \, , \hfill \\
\hfill X &\longmapsto &{\rm Hom}_{\mathcal E} (X,E)
\end{matrix}
$$
est un faisceau pour la topologie $J$.

\smallskip

Cela d\'efinit un foncteur
$$
\begin{matrix}
H : {\mathcal E} &\longrightarrow &\widehat{\mathcal C}_J \, , \hfill \\
\hfill E &\longmapsto &H(E) = {\rm Hom} (\bullet,E) \, .
\end{matrix}
$$
\end{cor}

\bigskip

\begin{demo}

Si $S$ est un crible $J$-couvrant d'un objet $X$ de ${\mathcal C}$, on a
$$
X = \varinjlim_{(U \to X) \in S} U \qquad \mbox{dans} \quad {\mathcal E} \, ,
$$
et donc pour tout objet $E$ de ${\mathcal E}$
$$
{\rm Hom} (X,E) = \varprojlim_{(U \to X) \in S} {\rm Hom} (U,E) \, .
$$
Cela signifie que les pr\'efaisceaux ${\rm Hom} (\bullet , E)$ sont des faisceaux pour la topologie $J$. 

\end{demo}

\subsection{Construction d'un adjoint \`a gauche}\label{subsec383}

\medskip

On note comme d'habitude $j^* : \widehat{\mathcal C} \to \widehat{\mathcal C}_J$ le foncteur de faisceautisation adjoint \`a gauche du foncteur de plongement $j_* : \widehat{\mathcal C}_J \hookrightarrow \widehat{\mathcal C}$.

\smallskip

Le foncteur $H : {\mathcal E} \to \widehat{\mathcal C}_J$ a un adjoint \`a gauche explicite:

\begin{lem}\label{lemIII84}

Consid\'erons le foncteur
$$
\begin{matrix}
G : \widehat{\mathcal C} &\longrightarrow &{\mathcal E} \, , \hfill \\
\hfill P &\longmapsto &\displaystyle \varinjlim_{(X,x) \in \int\!\!P} X
\end{matrix}
$$
et sa restriction \`a la sous-cat\'egorie pleine $\widehat{\mathcal C}_J \hookrightarrow \widehat{\mathcal C}$ \'egalement not\'ee $G$.

\smallskip

Alors:

\begin{listeimarge}

\item Le triangle
$$
\xymatrix{
\widehat{\mathcal C} \ar[rd]_-{j^*} \ar[rr]^-G &&{\mathcal E} \\
&\widehat{\mathcal C}_J \ar[ru]_-G
}
$$
est commutatif \`a unique isomorphisme pr\`es.

\medskip

\item Le foncteur
$$
G : \widehat{\mathcal C} \longrightarrow {\mathcal E}
$$
est adjoint \`a gauche du foncteur compos\'e
$$
H : {\mathcal E} \longrightarrow \widehat{\mathcal C}_J \xhookrightarrow{ \ { \ } \ } \widehat{\mathcal C} \, .
$$

\item Le foncteur
$$
G : \widehat{\mathcal C}_J \longrightarrow {\mathcal E}
$$
est adjoint \`a gauche du foncteur
$$
H : {\mathcal E} \longrightarrow \widehat{\mathcal C}_J \, .
$$
\end{listeimarge}
\end{lem}

\begin{remark}

En particulier, $G$ respecte les colimites arbitraires et $H$ respecte les limites arbitraires.

\end{remark}

\begin{demo}
\begin{listeisansmarge}
\item[(ii)] Pour tout objet $P$ de $\widehat{\mathcal C}$ et son faisceautis\'e $F = j^* P$, et pour tout objet $E$ de ${\mathcal E}$, on a
\begin{eqnarray}
{\rm Hom} (G(F),E) &= &\varprojlim_{(X,x) \in \int\!\!F} {\rm Hom} (X,E) \nonumber \\
&= &\varprojlim_{(X,x) \in \int\!\!F} H(E)(X) \nonumber \\
&= &\varprojlim_{(X,x) \in \int\!\!F} {\rm Hom}_{\widehat{\mathcal C}} \, (y(X),H(E)) \nonumber \\
&= &{\rm Hom}_{\widehat{\mathcal C}} \, (P,H(E)) \nonumber
\end{eqnarray}
puisque, par d\'efinition de $G$,
$$
G(F) = \varinjlim_{(X,x) \in \int\!\!F} X \qquad \mbox{dans} \quad {\mathcal E}
$$
et que l'on a d'autre part
$$
P = \varinjlim_{(X,x) \in \int\!\!P} y(X) \qquad \mbox{dans} \quad \widehat{\mathcal C} \, .
$$

\item[(iii)] r\'esulte de (ii) puisque le foncteur
$$
\begin{matrix}
{\mathcal E} &\longrightarrow &\widehat{\mathcal C} \, , \hfill \\
E &\longmapsto &{\rm Hom} (\bullet , E)
\end{matrix}
$$
se factorise en
$$
H : {\mathcal E} \longrightarrow \widehat{\mathcal C}_J
$$
et que le plongement
$$
j_* : \widehat{\mathcal C}_J \xhookrightarrow{ \ { \ } \ } \widehat{\mathcal C}
$$
est pleinement fid\`ele.

\medskip

\item[(i)] r\'esulte de (ii) et (iii) puisque les deux foncteurs
$$
\widehat{\mathcal C} \xrightarrow{ \ G \ } {\mathcal E}
$$
et
$$
\widehat{\mathcal C} \xrightarrow{ \ j^* \ } \widehat{\mathcal C}_J \xrightarrow{ \ G \ } {\mathcal E}
$$
sont adjoints \`a gauche du m\^eme foncteur
$$
{\mathcal E} \xrightarrow{ \ H \ } \widehat{\mathcal C}_J \xrightarrow{ \ j_* \ } \widehat{\mathcal C} \, .
$$
\end{listeisansmarge}
\end{demo}

\subsection{Pleine fid\'elit\'e du foncteur vers les faisceaux}\label{subsec384}

Montrons:

\begin{lem}\label{lemIII85}
\begin{listeimarge}
\item Le foncteur
$$
\begin{matrix}
H:{\mathcal E} &\longrightarrow &\widehat{\mathcal C}_J \, , \hfill \\
\hfill E &\longmapsto &{\rm Hom} (\bullet , E)
\end{matrix}
$$
est pleinement fid\`ele.

\medskip

\item Le morphisme d'adjonction
$$
G \circ H \longrightarrow {\rm id}_{\mathcal E}
$$
est un isomorphisme de foncteurs.
\end{listeimarge}
\end{lem}

\begin{demo}
\begin{listeisansmarge}
\item[] Comme $G$ est adjoint \`a gauche de $H$, (i) et (ii) sont \'equivalentes.

\medskip

\item[(ii)] signifie que si $E$ est un objet de ${\mathcal E}$ et $F = {\rm Hom} (\bullet , E)$ est le faisceau associ\'e sur $({\mathcal C},J)$, le morphisme de ${\mathcal E}$
$$
\varinjlim_{(X,x) \in \int\!\!F} X \longrightarrow E
$$
est un isomorphisme.

\smallskip

Notons $\left( X_i , X_i \xrightarrow{ \ x_i \ } E \right)$ les paires constitu\'ees d'un objet $X_i$ de ${\mathcal C}$ et d'un morphisme $X_i \xrightarrow{ \ x_i \ } E$ de ${\mathcal E}$.

\smallskip

Comme la famille des objets de ${\mathcal C}$ est s\'eparante, le morphisme
$$
\coprod_i X_i \longrightarrow E
$$
est un \'epimorphisme.

\smallskip

D'apr\`es la propri\'et\'e (5) de ${\mathcal E}$, $E$ s'\'ecrit comme la colimite
$$
E = \varinjlim \left( R \rightrightarrows \coprod_i X_i \right)
$$
si $R$ est la relation d'\'equivalence
$$
R = \left( \coprod_i X_i \right) \times_E \left( \coprod_j X_j \right)
$$
qui, d'apr\`es la propri\'et\'e (6), s'\'ecrit encore
$$
R = \coprod_{i,j} X_i \times_E X_j \, .
$$

Pour tous indices $i,j$, notons $\left( Y_k , Y_k \xrightarrow{ \ x_k \ } X_i \times_E X_j \right)$ les paires constitu\'ees d'un objet $Y_k$ de ${\mathcal C}$ et d'un morphisme
$$
Y_k \xrightarrow{ \ x_k = (x'_k , x''_k) \ } X_i \times_E X_j \, .
$$
Pour tous indices $i,j$, le morphisme
$$
\coprod_k Y_k \longrightarrow X_i \times_E X_j
$$
est un \'epimorphisme, et donc $E$ est la colimite du diagramme constitu\'e par les objets de ${\mathcal C}$
$$
X_i \, , \ X_j \, , \ Y_k
$$
reli\'es par les morphismes:
$$
\xymatrix{
Y_k \ar[d]_-{x'_k} \ar[r]^-{x''_k} &X_j \\
X_i
}
$$

Cela prouve comme voulu que
$$
E = \varinjlim_{(X,x) \in \int\!\!E} X \, .
$$
\end{listeisansmarge}
\end{demo}

\subsection{Correspondance des familles globalement \'epimorphiques}\label{subsec385}

\medskip

On note:

\begin{lem}\label{lemIII86}

Une famille de morphismes de ${\mathcal C}$
$$
X_i \longrightarrow X \, , \qquad i \in I \, ,
$$
est couvrante pour la topologie $J$ si et seulement si elle est globalement \'epimorphique dans ${\mathcal E}$.
\end{lem}

\begin{demo}

Si $X$ est la colimite dans ${\mathcal E}$ du crible engendr\'e par les $X_i \to X$, a fortiori la famille des $X_i \to X$ est globalement \'epimorphique.

\smallskip

R\'eciproquement, si
$$
\coprod_i X_i \longrightarrow X
$$
est un \'epimorphisme, $X$ est la colimite du diagramme
$$
R \rightrightarrows \coprod_i X_i
$$
pour
\begin{eqnarray}
R &= &\left( \coprod_i X_i \right) \times_X \left( \coprod_j X_j \right) \nonumber \\
&= &\coprod_{i,j} X_i \times_X X_j \, . \nonumber
\end{eqnarray}

Pour tous indices $i,j$, la famille des paires $\left( Y_k , Y_k \xrightarrow{ \ y_k \ } X_i \times_X X_j \right)$ constitu\'ees d'un objet $Y_k$ de ${\mathcal C}$ et d'un morphisme
$$
Y_k \xrightarrow{ \ y_k = (y'_k , y''_k) \ } X_i \times_X X_j
$$
est globalement \'epimorphique au-dessus de $X_i \times_X X_j$.

\smallskip

Donc $X$ est la colimite du diagramme constitu\'e des objets $X_i$ et $Y_k$ reli\'es par les morphismes
$$
\xymatrix{
Y_k \ar[d]_{y'_k} \ar[r]^{y''_k} &X_j \\
X_i
}
$$
ce qui signifie
$$
X = \varinjlim_{(U \to X) \in S} U
$$
si $S$ est le crible de $X$ engendr\'e par les morphismes $X_i \to X$.

\smallskip

De plus, pour tout morphisme $X' \xrightarrow{ \ f \ } X$, le foncteur $X' \times_X \, \bullet$ respecte les familles globalement \'epimorphiques de ${\mathcal E}$.

\smallskip

Il en r\'esulte que la famille des morphismes de ${\mathcal C}$
$$
X'_j \longrightarrow X'
$$
qui s'inscrivent dans au moins un carr\'e commutatif
$$
\xymatrix{
X'_j \ar[d] \ar[r] &X_i \ar[d] \\
X' \ar[r]^-f &X
}
$$
est globalement \'epimorphique.

\smallskip

Or cette famille est celle des \'el\'ements du crible $f^* S$.

\smallskip

Comme elle est globalement \'epimorphique, on a d'apr\`es ce qui pr\'ec\`ede
$$
X' = \varinjlim_{(U' \to X') \in f^* S} U' \, .
$$

Autrement dit, le crible $S$ est $J$-couvrant.

\end{demo}

\medskip

On d\'eduit de ce lemme:

\begin{cor}\label{corIII87}

Une famille de morphismes de ${\mathcal E}$
$$
E_i \to E \, , \qquad i \in I \, ,
$$
est globalement \'epimorphique si et seulement si son image par le foncteur
$$
H : {\mathcal E} \longrightarrow \widehat{\mathcal C}_J
$$
est globalement \'epimorphique.
\end{cor}

\begin{demo}

La suffisance de la condition r\'esulte de ce que le foncteur $H$ est pleinement fid\`ele.

\smallskip

Pour la n\'ecessit\'e, supposons que la famille des $E_i \to E$ est globalement \'epimorphique.

\smallskip

La famille des morphismes
$$
\ell(X) = H(X) \xrightarrow{ \ H(x) \ } H(E)
$$
index\'es par les objets $X$ de ${\mathcal C}$ et les $x \in {\rm Hom} (X,E) = H(E)(X)$ est globalement \'epimorphique et, pour tout objet $X$, la famille des
$$
X \times_E E_i \longrightarrow X
$$
reste globalement \'epimorphique dans ${\mathcal E}$.

\smallskip

De plus, chaque $X \times_E E_i$ admet une famille globalement \'epimorphique de morphismes
$$
X_k \longrightarrow X \times_E E_i
$$
dont les sources $X_k$ sont des objets de ${\mathcal C}$.

\smallskip

Donc on peut supposer que $E = X$ est un objet de ${\mathcal C}$, de m\^eme que les $E_i = X_i$.

\smallskip

Dire que la famille des morphismes
$$
\ell (X_i) \longrightarrow \ell (X)
$$
est globalement \'epimorphique dans $\widehat{\mathcal C}_J$ \'equivaut d'apr\`es le lemme \ref{lemIII53} \`a dire que le crible de $X$ engendr\'e par les $X_i \to X$ est $J$-couvrant.

\smallskip

On est ramen\'e au lemme \ref{lemIII86}, ce qui conclut la d\'emonstration. 

\end{demo}

\subsection{Balancement de ${\mathcal E}$ et factorisation de ses morphismes}\label{subsec386}

\medskip

Nous pouvons maintenant d\'emontrer:

\begin{lem}\label{lemIII88}

La cat\'egorie ${\mathcal E}$ est balanc\'ee.
\end{lem}

\begin{demo}

Soit $u : E' \to E$ une fl\`eche de ${\mathcal E}$ qui est \`a la fois un monomorphisme et un \'epimorphisme. Son image
$$
H(u) : H(E') \longrightarrow H(E)
$$
est encore un monomorphisme puisque $H$, qui admet un adjoint \`a gauche, pr\'eserve les limites.

\smallskip

D'autre part, $H(u)$ est aussi un \'epimorphisme puisque, d'apr\`es le corollaire \ref{corIII87}, $H$ pr\'eserve les \'epimorphismes.

\smallskip

Comme le topos $\widehat{\mathcal C}_J$ est une cat\'egorie balanc\'ee, $H(u)$ est un isomorphisme.

\smallskip

Cela implique comme voulu que $u$ est un isomorphisme puisque le foncteur $H$ est pleinement fid\`ele. 

\end{demo}

On d\'eduit de ce lemme:

\begin{cor}\label{corIII89}

Tout morphisme de ${\mathcal E}$
$$
u : E' \longrightarrow E
$$
se factorise canoniquement en un \'epimorphisme
$$
u_1 : E' \ -\!\!\!\twoheadrightarrow {\rm Im} (u)
$$
suivi d'un monomorphisme
$$
u_2 : {\rm Im} (u) \xhookrightarrow{ \ { \ } \ } E \, .
$$
\end{cor}

\begin{demo}

Si une telle factorisation existe, ${\rm Im} (u)$ doit n\'ecessairement s'identifier au quotient de $E'$ par la relation d'\'equivalence $E' \times_{{\rm Im}(u)} E' = E' \times_E E'$.

\smallskip

Soit donc ${\rm Im} (u)$ le quotient de $E'$ par la relation d'\'equivalence
$$
R = E' \times_E E' \, .
$$

Alors $u$ se factorise canoniquement en
$$
E' \xrightarrow{ \ u_1 \ } {\rm Im} (u) \xrightarrow{ \ u_2 \ } E
$$
et $u_1$ est par construction un \'epimorphisme.

\smallskip

Il faut prouver que $u_2$ est un monomorphisme, c'est-\`a-dire que le morphisme diagonal
$$
\Delta : {\rm Im} (u) \longrightarrow {\rm Im} (u) \times_{E'} {\rm Im} (u)
$$
est un isomorphisme.

\smallskip

On sait d\'ej\`a que $\Delta$ est un monomorphisme.

\smallskip

De plus, on a un carr\'e cart\'esien
$$
\xymatrix{
R = E' \times_{{\rm Im}(u)} E' \ar[d] \ar[r]^-{\sim} &E' \times_E E' = R \ar[d]^-{u_1 \times u_1} \\
{\rm Im} (u) \ar[r]^-{\Delta} &{\rm Im} (u) \times_E {\rm Im} (u)
}
$$
dont la fl\`eche verticale $u_1 \times u_1$ est un \'epimorphisme.

\smallskip

Donc $\Delta$ est un \'epimorphisme et, comme c'est aussi un monomorphisme, $\Delta$ est un isomorphisme. 

\end{demo}

\subsection{Respect des colimites et fin de la d\'emonstration}\label{subsec387}

\medskip

On peut maintenant d\'emontrer:

\begin{lem}\label{lemIII810}

Le foncteur
$$
H : {\mathcal E} \longrightarrow \widehat{\mathcal C}_J
$$
respecte les colimites arbitraires.
\end{lem}

\begin{demo}

Consid\'erons un carquois $D$, un $D$-diagramme $E_{\bullet}$ de ${\mathcal E}$ et sa colimite
$$
E = \varinjlim_D E_{\bullet} \, .
$$

Il faut montrer que le morphisme canonique
$$
\varinjlim_D H(E_{\bullet}) \longrightarrow H(E)
$$
est un isomorphisme.
 
\smallskip
 
D'apr\`es le corollaire \ref{corIII87}, on sait d\'ej\`a que c'est un \'epimorphisme.
 
\smallskip

Soit
$$
E' = \coprod_{d \in {\rm Ob} (D)} E_d
$$
et
$$
R = E' \times_E E' \, .
$$

Pour tout diagramme fini
$$
D_n = (d_0 \leftrightarrow d_1 \leftrightarrow \cdots \leftrightarrow d_n)
$$
constitu\'e d'objet $d_0 , d_1 , \cdots , d_n$ de $D$ et de fl\`eches $d_{i-1} \to d_i$ ou $d_i \to d_{i-1}$ de $D$ not\'ees $d_{i-1} \leftrightarrow d_i$, on a un morphisme
$$
\varprojlim_{D_n} E_{\bullet} \longrightarrow E_{d_0} \times E_{d_n}
$$
d'o\`u une fl\`eche
$$
\coprod_{D_n} \left( \varprojlim_{D_n} E_{\bullet} \right) \longrightarrow E' \times E' \, .
$$

Son image est une relation d'\'equivalence de $E'$, n\'ecessairement \'egale \`a $R$ puisqu'elle d\'efinit le quotient $E$ de $E'$.

\smallskip

Donc la famille des fl\`eches
$$
\varprojlim_{D_n} E_{\bullet} \longrightarrow R
$$
est globalement \'epimorphique dans ${\mathcal E}$.

\smallskip

Comme $H$ pr\'eserve les limites ainsi que les familles globalement \'epimorphiques, $H(R)$ est la relation d'\'equivalence de $H(E')$ qui d\'efinit son quotient $H(E)$, et la famille des fl\`eches
$$
\varprojlim_{D_n} H(E_{\bullet}) \longrightarrow H(R)
$$
est globalement \'epimorphique.

\smallskip

Cela entra{\^\i}ne que le morphisme
$$
\varinjlim_{D} H(E_{\bullet}) \longrightarrow H(E)
$$
est un isomorphisme. 

\end{demo}

On d\'eduit de ce lemme le corollaire suivant qui termine la d\'emonstration du th\'eor\`eme \ref{thmIII81}:

\begin{cor}\label{corIII811}

Le morphisme d'adjonction
$$
{\rm id}_{\widehat{\mathcal C}_J} \longrightarrow H \circ G
$$
est un isomorphisme.
\end{cor}

\begin{demo}

Les deux foncteurs $H$ et $G$ respectent les colimites et tout objet $F$ de $\widehat{\mathcal C}_J$ s'\'ecrit comme la colimite
$$
F = \varinjlim_{(X,x) \in \int\!\!F} \ell (X)
$$
des images $\ell (X)$ des objets $X$ de ${\mathcal C}$ par le foncteur canonique
$$
\ell : {\mathcal C} \xrightarrow{ \ y \ } \widehat{\mathcal C} \xrightarrow{ \ j^* \ } \widehat{\mathcal C}_J \, .
$$

Il suffit donc de prouver que pour tout objet $X$ de ${\mathcal C}$, le morphisme canonique
$$
\ell (X) \longrightarrow H \circ G (\ell (X))
$$
est un isomorphisme.

\smallskip

Or le foncteur
$$
\widehat{\mathcal C} \longrightarrow \widehat{\mathcal C}_J \xrightarrow{ \ G \ } {\mathcal E}
$$
prolonge le foncteur de plongement
$$
{\mathcal C} \xhookrightarrow{ \ { \ } \ } {\mathcal E}
$$
ce qui signifie que pour tout objet $X$ de ${\mathcal C}$
$$
G(\ell (X)) = X
$$
et donc
$$
H \circ G (\ell (X)) = {\rm Hom} (\bullet , X) = \ell (X)
$$
comme voulu. 

\end{demo}

Le th\'eor\`eme \ref{thmIII81} est d\'emontr\'e puisque, d'apr\`es le lemme \ref{lemIII85} (ii) et le corollaire \ref{corIII811} ci-dessus, les deux morphismes d'adjonction
$$
G \circ H \longrightarrow {\rm id}_{\mathcal E}
$$
et
$$
{\rm id}_{\widehat{\mathcal C}_J} \longrightarrow H \circ G
$$
sont des isomorphismes.

\section{Application aux probl\`emes de repr\'esentabilit\'e dans les cat\'egories}\label{sec39}

\subsection{Conditions faisceautiques de repr\'esentabilit\'e}\label{subsec391}

\medskip

D'apr\`es la proposition \ref{propII58}, une topologie $J$ sur une cat\'egorie essentiellement petite ${\mathcal C}$ est sous-canonique, au sens que le foncteur canonique
$$
\ell : {\mathcal C} \xrightarrow{ \ y \ } \widehat{\mathcal C} \xrightarrow{ \ j^* \ } \widehat{\mathcal C}_J
$$
est pleinement fid\`ele, si et seulement si les pr\'efaisceaux repr\'esentables
$$
y(X) = {\rm Hom} (\bullet , X) : {\mathcal C}^{\rm op} \longrightarrow {\rm Ens} \, , \qquad X \in {\rm Ob} ({\mathcal C}) \, ,
$$
sont des faisceaux pour la topologie $J$.

\smallskip

D'apr\`es le corollaire \ref{corII510}, cette condition est satisfaite si et seulement si $J$ est contenue dans une certaine topologie $J_c$ de ${\mathcal C}$, appel\'ee la topologie canonique.

\smallskip

Ces r\'esultats peuvent \^etre vus comme exprimant des conditions n\'ecessaires pour qu'un pr\'efaisceau
$$
F : {\mathcal C}^{\rm op} \longrightarrow {\rm Ens}
$$
soit repr\'esentable:

\begin{cor}\label{corIII91}

Soit ${\mathcal C}$ une cat\'egorie essentiellement petite.

\smallskip

Pour qu'un foncteur contravariant
$$
F : {\mathcal C}^{\rm op} \longrightarrow {\rm Ens}
$$
soit repr\'esentable, il est n\'ecessaire qu'il soit un faisceau pour la topologie canonique de ${\mathcal C}$, et a fortiori pour toute topologie sous-canonique $J$ de ${\mathcal C}$.
\end{cor}

\begin{remarksqed}
\begin{listeisansmarge}
\item La topologie canonique d'une cat\'egorie essentiellement petite ${\mathcal C}$ est souvent difficile \`a expliciter en termes concrets.

\smallskip

Dans ce cas, le crit\`ere ci-dessus s'appliquera plut\^ot en prenant pour $J$ une topologie sous-canonique facile \`a expliciter concr\`etement et aussi proche que possible de la topologie canonique de ${\mathcal C}$.

\smallskip

Par exemple, si ${\mathcal C} = {\rm Sch} /_{\!{\rm pf}} \, S$ est la cat\'egorie essentiellement petite des sch\'emas de pr\'esentation finie sur un sch\'ema de base $S$, la topologie canonique de ${\mathcal C}$ est difficile \`a expliciter mais la topologie fppf est sous-canonique et proche de la topologie canonique de ${\mathcal C}$. C'est celle-l\`a que l'on choisira g\'en\'eralement pour appliquer le crit\`ere du corollaire.

\medskip

\item Le corollaire dit que si un pr\'efaisceau
$$
F : {\mathcal C}^{\rm op} \longrightarrow {\rm Ens}
$$
n'est pas un faisceau pour une certaine topologie sous-canonique $J$ de ${\mathcal C}$, il ne peut \^etre repr\'esentable.

\smallskip

Dans ce cas, on peut cependant remplacer $F$ par son image
$$
j^* F : {\mathcal C}^{\rm op} \longrightarrow {\rm Ens}
$$
par le foncteur de $J$-faisceautisation
$$
j^* : \widehat{\mathcal C} \longrightarrow \widehat{\mathcal C}_J
$$
et se demander si le faisceau $j^* F$ est repr\'esentable.

\smallskip

Cela demande bien s\^ur que $j^* F$ soit un faisceau non seulement pour la topologie $J$ mais pour la topologie canonique de ${\mathcal C}$ et a fortiori pour toute topologie sous-canonique.

\medskip

\item Par exemple, si $P$ est un sch\'ema projectif et lisse sur un sch\'ema de base $S$, le pr\'efaisceau
$$
\begin{matrix}
({\rm Sch} /_{\!{\rm pf}} \, S)^{\rm op} &\longrightarrow &{\rm Ens} \, , \hfill \\
\hfill (X \to S) &\longmapsto &{\mathcal P}ic_{P/S} (X)
\end{matrix}
$$
qui associe \`a tout sch\'ema de pr\'esentation finie $X$ sur $S$ l'ensemble ${\mathcal P}ic_{P/S} (X)$ des classes d'isomorphie de ${\mathcal O}_{P \times_S X}$-Modules localement libres de rang 1 (pour la topologie de Zariski) sur le sch\'ema $P \times_S X$ n'est en g\'en\'eral pas un faisceau pour la topologie \'etale.

\smallskip

A fortiori, il n'est pas repr\'esentable en g\'en\'eral.

\smallskip

On montre en revanche que son transform\'e
$$
j^* {\mathcal P}ic_{P/S}
$$
par le foncteur $j^*$ de faisceautisation pour la topologie \'etale est repr\'esentable par un sch\'ema ${\rm Pic}_{P/S}$, appel\'e le sch\'ema de Picard de $X$ sur $S$, qui est une somme disjointe de sch\'emas projectifs lisses sur $S$. A fortiori, c'est un faisceau non seulement pour la topologie \'etale mais aussi pour la topologie fppf et pour la topologie canonique de ${\mathcal C} = {\rm Sch} /_{\!{\rm pf}} \, S$.

\medskip

\item Pour toute topologie sous-canonique $J$, donc contenue dans la topologie canonique $J_c$ d'une cat\'egorie essentiellement petite ${\mathcal C}$, le foncteur de Yoneda
$$
y : {\mathcal C} \xhookrightarrow{ \ { \ } \ } \widehat{\mathcal C}
$$
sa factorise en une suite de plongements pleinement fid\`eles
$$
{\mathcal C} \xhookrightarrow{ \ { \ } \ } \widehat{\mathcal C}_{J_c} \xhookrightarrow{ \ { \ } \ } \widehat{\mathcal C}_J \xhookrightarrow{ \ { \ } \ } \widehat{\mathcal C} \, .
$$

Mais les topos $\widehat{\mathcal C}_J$ et $\widehat{\mathcal C}_{J_c}$ sont \'enorm\'ement plus gros que la cat\'egorie essentiellement petite ${\mathcal C}$, puisque contrairement \`a elle ils ne sont pas des cat\'egories essentiellement petites d\`es lors qu'ils sont non triviaux.

\smallskip

Cela signifie que la condition n\'ecessaire du corollaire est tr\`es loin d'\^etre suffisante.

\end{listeisansmarge}
\end{remarksqed}

A paragraphe suivant, nous allons expliciter en termes de la structure cat\'egorique de ${\mathcal C}$ ce que signifie la condition n\'ecessaire du corollaire \ref{corIII91}.

\subsection{Description des cribles couvrants de la topologie canonique et application}\label{subsec392}

\medskip

Le corollaire \ref{corII510} (i) permet de d\'ecrire en termes de la structure cat\'egorique de ${\mathcal C}$ ce que sont les cribles couvrants de la topologie canonique d'une cat\'egorie essentiellement petite ${\mathcal C}$.

\smallskip

Pour ce faire, on a besoin de la g\'en\'eralisation suivante de la d\'efinition \ref{defII515}:

\begin{defn}\label{defIII92}

Soit ${\mathcal C}$ une cat\'egorie localement petite.

\begin{listeimarge}

\item Une famille globalement \'epimorphique de morphismes de ${\mathcal C}$
$$
U_i \xrightarrow{ \ p_i \ } X \, , \qquad i \in I \, ,
$$
est dite ``stricte'' s'il existe une famille d'objets $U_{i,j,k}$ de ${\mathcal C}$ s'inscrivant dans des carr\'es commutatifs
$$
\xymatrix{
U_{i,j,k} \ar[d]_-{q_{i,k}} \ar[r]^-{q_{j,k}} &U_j \ar[d]^-{p_j} \\
U_i \ar[r]^-{p_i} &X
}
$$
tels que, pour tout objet $Z$ de ${\mathcal C}$, la composition avec les $p_i$ identifie l'ensemble
$$
{\rm Hom} (X,Z)
$$
au sous-ensemble de
$$
\prod_i {\rm Hom} (U_i , Z)
$$
constitu\'e des familles $(f_i)_{i \in I}$ qui satisfont les \'equations
$$
f_i \circ q_{i,k} = f_j \circ q_{j,k} \, , \qquad \forall \, (i,j,k) \, .
$$

\item Une telle famille
$$
\left( U_i \xrightarrow{ \ p_i \ } X \right)_{i \in I}
$$
est dite ``\'epimorphique stricte universelle'' si, pour tout morphisme de ${\mathcal C}$ de but $X$
$$
X' \xrightarrow { \ x \ } X \, ,
$$
il existe une famille de morphismes
$$
U'_j \xrightarrow{ \ p'_j \ } X'
$$
qui soit ``\'epimorphique stricte'' et dont chaque \'el\'ement $p'_j$ s'inscrit dans au moins un carr\'e commutatif de la forme:
$$
\xymatrix{
U'_j \ar[d]_-{p'_j} \ar[r] &U_i \ar[d]^-{p_i} \\
X' \ar[r]^-{x} &X
}
$$

\end{listeimarge}
\end{defn}

\begin{remarksqed}
\begin{listeisansmarge}
\item La condition de (i) signifie que $X$ est la colimite du diagramme constitu\'e des objets $U_i$ et $U_{i,j,k}$ reli\'es par les morphismes:
$$
\xymatrix{
U_{i,j,k} \ar[d]_-{q_{i,k}} \ar[r]^-{q_{j,k}} &U_j \\
U_i
}
$$

\item Une famille de morphismes carrables de ${\mathcal C}$
$$
U_i \longrightarrow X \, , \qquad i \in I \, ,
$$
est ``\'epimorphique stricte'' si $X$ est la colimite du diagramme constitu\'e des objets $U_i$ et $U_i \times_X U_j$ reli\'es par les projections des $U_i \times_X U_j$ sur leurs facteurs $U_i$ et $U_j$.

\smallskip

Elle est ``\'epimorphique stricte universelle'' si, pour tout morphisme $X' \to X$ de ${\mathcal C}$, $X'$ est la colimite du diagramme constitu\'e des objets $U_i \times_X X'$ et $U_i \times_X U_j \times_X X'$ reli\'es par les projections des $U_i \times_X U_j \times_X X'$ sur leurs facteurs $U_i \times_X X'$ et $U_j \times_X X'$. 

\end{listeisansmarge}
\end{remarksqed}

\medskip

On d\'eduit du corollaire \ref{corII510} (i):

\begin{lem}\label{lemIII93}

Dans une cat\'egorie essentiellement petite ${\mathcal C}$, une famille de morphismes
$$
U_i \longrightarrow X \, , \qquad i \in I \, ,
$$
est couvrante pour la topologie canonique de ${\mathcal C}$ si et seulement si elle est ``\'epimorphique stricte universelle'' au sens de la d\'efinition pr\'ec\'edente.
\end{lem}

\begin{demo}

Il r\'esulte de la d\'efinition qu'une famille de morphismes
$$
U_i \longrightarrow X \, , \qquad i \in I \, ,
$$
est \'epimorphique stricte si et seulement si le crible $S$ de $X$ qu'elle engendre v\'erifie la condition
$$
X = \varinjlim_{(U \to X) \in S} U \, .
$$

Elle est \'epimorphique stricte universelle si et seulement si on a pour tout morphisme $x : X' \to X$ de ${\mathcal C}$
$$
X' = \varinjlim_{(U' \to X') \in x^* S} U' \, .
$$

La conclusion r\'esulte du lemme \ref{lemII511} (i). 

\end{demo}

Compte tenu de ce lemme, le corollaire \ref{corIII91} se r\'e\'ecrit:

\begin{cor}\label{corIII94}

Soit ${\mathcal C}$ une cat\'egorie essentiellement petite.

\smallskip

Pour qu'un foncteur contravariant
$$
F : {\mathcal C}^{\rm op} \longrightarrow {\rm Ens}
$$
soit repr\'esentable, il est n\'ecessaire qu'il satisfasse la condition
$$
F(X) = \varprojlim_{(U \to X) \in S} F(U)
$$
pour tout crible $S$ d'un objet $X$ engendr\'e par une famille \'epimorphique stricte universelle
$$
U_i \longrightarrow X \, .
$$
\end{cor}

\begin{remarks}
\begin{listeisansmarge}
\item A fortiori la condition doit \^etre v\'erifi\'ee par tout crible $J$-couvrant, pour n'importe quelle topologie sous-canonique $J$ de ${\mathcal C}$.

\medskip

\item Si les $U_i \to X$ sont carrables, la condition s'\'ecrit
$$
F(X) = {\rm eg} \left( \prod_{i \in I} F(U_i) \rightrightarrows \prod_{i,j} F(U_i \times_X U_j) \right) .
$$
\end{listeisansmarge}
\end{remarks}

\begin{demo}

C'est une simple traduction, r\'ealis\'ee par le lemme \ref{lemIII93}. 

\end{demo}

\subsection{Calcul de limites et d'exponentielles par les faisceaux}\label{subsec393}

\medskip

On observe que chercher des limites de diagrammes ou des exponentielles dans une cat\'egorie essentiellement petite ${\mathcal C}$ \'equivaut \`a les calculer dans une cat\'egorie de faisceaux dans laquelle ${\mathcal C}$ est plong\'ee puis \`a se demander si les faisceaux obtenus sont repr\'esentables:

\begin{prop}\label{propIII95}

Soit ${\mathcal C}$ une cat\'egorie essentiellement petite munie d'une topologie sous-canonique $J$ et donc du foncteur pleinement fid\`ele associ\'e
$$
\ell : {\mathcal C} \xrightarrow{ \ y \ } \widehat{\mathcal C} \xrightarrow{ \ j^* \ } \widehat{\mathcal C}_J \, .
$$
Alors:

\begin{listeimarge}

\item Pour tout carquois $D$ et tout $D$-diagramme $X_{\bullet}$ de ${\mathcal C}$, la limite
$$
\varprojlim_D X_{\bullet}
$$
est bien d\'efinie dans ${\mathcal C}$ si et seulement si l'objet de $\widehat{\mathcal C}_J$
$$
\varprojlim_D \ell (X_{\bullet})
$$
est repr\'esentable comme pr\'efaisceau par un objet de ${\mathcal C}$ qui n'est alors autre que
$$
\varprojlim_D X_{\bullet} \, .
$$

\item Pour tout objet $X$ de ${\mathcal C}$ qui est carrable et pour tout objet $Z$ de ${\mathcal C}$, le foncteur contravariant
$$
\begin{matrix}
{\mathcal C}^{\rm op} &\longrightarrow &{\rm Ens} \, , \hfill \\
\hfill Y &\longmapsto &{\rm Hom} (X \times Y,Z)
\end{matrix}
$$
est repr\'esentable par un objet $Z^X$ de ${\mathcal C}$ si et seulement si l'objet de $\widehat{\mathcal C}_J$
$$
{\mathcal H}om (\ell (X),\ell(Z))
$$
est repr\'esentable comme pr\'efaisceau par un objet de ${\mathcal C}$ qui n'est alors autre que $Z^X$.
\end{listeimarge}
\end{prop}

\begin{demo}
\begin{listeisansmarge}
\item Par d\'efinition, la limite
$$
\varprojlim_D X_{\bullet}
$$
est bien d\'efinie dans ${\mathcal C}$ si et seulement si le pr\'efaisceau
$$
\varprojlim_D \, y(X_{\bullet})
$$
est repr\'esentable par un objet de ${\mathcal C}$ not\'e alors $\underset{D}{\varprojlim} \, X_{\bullet}$.

\smallskip

La conclusion r\'esulte alors de ce que, $J$ \'etant sous-canonique, le foncteur
$$
y : {\mathcal C} \xhookrightarrow{ \ { \ } \ } \widehat{\mathcal C}
$$
sa factorise en
$$
{\mathcal C} \xhookrightarrow{ \ \ell \ } \widehat{\mathcal C}_J \xhookrightarrow{ \ j_* \ } \widehat{\mathcal C}
$$
et que le foncteur pleinement fid\`ele $j_* : \widehat{\mathcal C}_J \hookrightarrow \widehat{\mathcal C}$ respecte les limites.

\medskip

\item S'il existe un objet $Z^X$ de ${\mathcal C}$ tel que
$$
\ell (Z^X) \cong {\mathcal H}om (\ell (X), \ell(Z)) \, ,
$$
alors on a pour tout objet $Y$ de ${\mathcal C}$ la suite d'identifications
\begin{eqnarray}
{\rm Hom}_{\mathcal C} (Y , Z^X) &= &{\rm Hom}_{\widehat{\mathcal C}_J} (\ell (Y) , \ell (Z^X)) \nonumber \\
&\cong &{\rm Hom}_{\widehat{\mathcal C}_J} (\ell (Y) , {\mathcal H}om (\ell (X), \ell (Z)) \nonumber \\
&= &{\rm Hom} (\ell (X) \times \ell (Y) , \ell (Z)) \nonumber \\
&= &{\rm Hom} (\ell (X \times Y) , \ell (Z)) \nonumber \\
&= &{\rm Hom}_{\mathcal C} (X \times Y , Z) \nonumber
\end{eqnarray}
qui montre que $Z^X$ repr\'esente le foncteur
$$
Y \longmapsto {\rm Hom} (X \times Y,Z) \, .
$$

R\'eciproquement, si $Z^X$ est un objet de ${\mathcal C}$ qui repr\'esente ce foncteur, on a pour tout objet $Y$ de ${\mathcal C}$ un isomorphisme canonique
$$
{\rm Hom} (\ell (X) \times \ell (Y),\ell (Z)) = {\rm Hom} (\ell (Y) , \ell (Z^X)) \, .
$$
Or tout objet $F$ de $\widehat{\mathcal C}_J$ s'\'ecrit canoniquement comme la colimite
$$
F = \varinjlim_{(X,x) \in \int\!\!F} \ell (X)
$$
calcul\'ee dans $\widehat{\mathcal C}_J$ sur la cat\'egorie $\int\!\!F$ des \'el\'ements $(X,x)$ de $F$.

\smallskip

Comme le foncteur $\ell (X) \times \bullet$ de produit avec $\ell (X)$ dans $\widehat{\mathcal C}_J$ respecte les colimites, on en d\'eduit pour un tel objet $F$ de $\widehat{\mathcal C}_J$ un isomorphisme canonique
$$
{\rm Hom} (\ell (X) \times F , \ell (Z)) = {\rm Hom} (F,\ell (Z^X)) \, .
$$
Autrement dit, l'objet $\ell (Z^X)$ de $\widehat{\mathcal C}_J$ repr\'esente le foncteur contravariant
$$
\begin{matrix}
(\widehat{\mathcal C}_J)^{\rm op} &\longrightarrow &{\rm Ens} \, , \hfill \\
\hfill F &\longmapsto &{\rm Hom} (\ell (X) \times F , \ell (Z))
\end{matrix}
$$
donc est canoniquement isomorphe \`a l'objet ${\mathcal H}om (\ell (X) , \ell (Z))$ de $\widehat{\mathcal C}_J$.

\smallskip

Cela termine la d\'emonstration. 

\end{listeisansmarge}
\end{demo}

\subsection{Faisceaux et d\'efinition d'une notion de colimite g\'eom\'etrique}\label{subsec394}

Le foncteur de Yoneda
$$
y : {\mathcal C} \xhookrightarrow{ \ { \ } \ } \widehat{\mathcal C}
$$
sur une cat\'egorie essentiellement petite ${\mathcal C}$ ne respecte pas les colimites, non plus que les foncteurs
$$
{\mathcal C} \xhookrightarrow{ \ \ell \ } \widehat{\mathcal C}_J \xhookrightarrow{ \ j_* \ } \widehat{\mathcal C}
$$
associ\'ees \`a n'importe quelle topologie sous-canonique $J$.

\smallskip

C'est pourquoi calculer des colimites dans $\widehat{\mathcal C}_J$ et se demander si elles sont repr\'esentables d\'efinit des notions de colimites dans ${\mathcal C}$ qui sont plus fortes que la notion interne \`a ${\mathcal C}$ et qui d\'ependent du choix de $J$:

\begin{thm}\label{thmIII96}

Soit ${\mathcal C}$ une cat\'egorie essentiellement petite munie d'une topologie sous-canonique $J$ et donc du foncteur pleinement fid\`ele
$$
\ell : {\mathcal C} \xhookrightarrow{ \ { \ } \ } \widehat{\mathcal C}_J
$$
\`a travers lequel se factorise $y : {\mathcal C} \hookrightarrow \widehat{\mathcal C}$.

\smallskip

Soient $D$ un carquois et $X_{\bullet}$ un $D$-diagramme de ${\mathcal C}$.

\smallskip

Alors:

\begin{listeimarge}

\item Si la colimite calcul\'ee dans $\widehat{\mathcal C}_J$
$$
\varinjlim_D \ell (X_{\bullet})
$$
est repr\'esentable comme pr\'efaisceau par un objet $X$ de ${\mathcal C}$, celui-ci satisfait les conditions suivantes:

\medskip

$\left\lmoustache \begin{matrix}
{\rm (A)} &\mbox{L'objet $X$ est une colimite du $D$-diagramme $X_{\bullet}$ de ${\mathcal C}$.} \hfill \\
{ \ } \\
{\rm (B)} &\mbox{La famille des morphismes de ${\mathcal C}$} \hfill \\
{ \ } \\
&X_d \longrightarrow X \, , \quad d \in {\rm Ob} (D) \, , \\
{ \ } \\
&\mbox{est $J$-couvrante.} \hfill \\
{ \ } \\
{\rm (C)} &\mbox{Pour tout carr\'e commutatif de ${\mathcal C}$} \hfill \\
{ \ } \\
&\xymatrix{
U \ar[d] \ar[r] &X_{d'} \ar[d] \\
X_d \ar[r] &X
} \\
{ \ } \\
&\mbox{il existe une famille $J$-couvrante de morphismes} \hfill \\
{ \ } \\
&U' \longrightarrow U 
\end{matrix} \right.
$

$
\left\rmoustache \begin{matrix}
&\mbox{tels que les compos\'es $U' \to X_d$ et $U' \to X_{d'}$ s'inscrivent dans un diagramme commutatif de ${\mathcal C}$} \hfill \\
{ \ } \\
&\def\troispoints{\ar@{}[rd]|{\displaystyle\cdots}}
\xymatrix{
&U' \ar[ld] \ar[d] \troispoints \ar[rrd] \\
X_d = X_{d_0} \ar@{<->}[r]&X_{d_1} \ar@{<->}[r] &\cdots \ar@{<->}[r]&X_{d_n} = X_{d'}
} \\
{ \ } \\
&\mbox{o\`u $d_0 = d$, $d_n = d'$ et chaque $X_{d_{k-1}} \leftrightarrow X_{d_k}$ est un morphisme} \hfill \\
{ \ } \\
&X_{d_{k-1}} \longrightarrow X_{d_k} \qquad \mbox{ou} \qquad X_{d_k} \longrightarrow X_{d_{k-1}} \\
{ \ } \\
&\mbox{associ\'e \`a une fl\`eche de $D$} \hfill \\
&d_{k-1} \longrightarrow d_k \qquad \mbox{ou} \qquad d_k \longrightarrow d_{k-1} \, .
\end{matrix} \right.
$

\medskip

\item R\'eciproquement, si le $D$-diagramme $X_{\bullet}$ de ${\mathcal C}$ admet une colimite $X$ dans ${\mathcal C}$ qui satisfait les conditions {\rm (B)} et {\rm (C)} de {\rm (i)}, l'objet
$$
\ell (X) \qquad \mbox{dans} \qquad \widehat{\mathcal C}_J
$$
est une colimite du $D$-diagramme $\ell (X_{\bullet})$.

\smallskip

Autrement dit, le pr\'efaisceau
$$
{\rm Hom} (\bullet , X)
$$
est le transform\'e par le foncteur de faisceautisation
$$
j^* : \widehat{\mathcal C} \longrightarrow \widehat{\mathcal C}_J
$$
du pr\'efaisceau
$$
\begin{matrix}
\displaystyle\varinjlim_D y(X_{\bullet}) : &{\mathcal C}^{\rm op} &\longrightarrow &{\rm Ens} \, , \hfill \\
&\hfill Y &\longmapsto &\displaystyle\varinjlim_D {\rm Hom} (Y,X_{\bullet}) \, .
\end{matrix}
$$
\end{listeimarge}
\end{thm}

\begin{remarks}
\begin{listeisansmarge}
\item Si un $D$-diagramme $X_{\bullet}$ de ${\mathcal C}$ admet une colimite $X$ qui satisfait les conditions (B) et (C) relativement \`a une certaine topologie sous-canonique $J$ de ${\mathcal C}$, alors elle les v\'erifie a fortiori relativement \`a toute topologie sous-canonique plus fine que $J$.

\smallskip

En particulier, elle les satisfait relativement \`a la topologie canonique de ${\mathcal C}$.

\medskip

\item Si $X_{\bullet}$ est un $D$-diagramme de ${\mathcal C}$ dont toutes les fl\`eches $X_{d'} \to X_d$ sont carrables et qui admet une colimite $X$ telle que toutes les projections
$$
X_d \longrightarrow X \, , \qquad d \in {\rm Ob} ({\mathcal C}) \, , 
$$
sont carrables, la condition (C) s'\'ecrit encore:

\medskip

$
\left\{\begin{matrix}
&\mbox{Pour tous objets $d,d'$ de $D$, consid\'erant la famille des sous-carquois finis $D_n$ de $D$ de la forme} \hfill \\
{ \ } \\
&d=d_0 \longleftrightarrow d_1 \longleftrightarrow \cdots \longleftrightarrow d_n = d' \\
{ \ } \\
&\mbox{o\`u $d=d_0 , d_1 , \cdots , d_{n-1}$, $d_n = d'$ sont des objets de $D$ et chaque $d_{i-1} \leftrightarrow d_i$ est une fl\`eche} \hfill \\
&\mbox{$d_{i-1} \to d_i$ ou $d_i \to d_{i-1}$ de $D$, la famille des morphismes} \hfill \\
{ \ } \\
&\displaystyle\varprojlim_{D_n} X_{\bullet} \longrightarrow X_d \times_X X_{d'} \\
&\mbox{est $J$-couvrante.} \hfill
\end{matrix} \right.
$
\end{listeisansmarge}
\end{remarks}

\medskip

\begin{demo}

Si
$$
\ell (X) = \varinjlim_D \ell (X_{\bullet})
$$
dans $\widehat{\mathcal C}_J$, on a pour tout objet $Y$ de ${\mathcal C}$
$$
{\rm Hom}_{\widehat{\mathcal C}_J} (\ell (X) , \ell (Y)) = \varprojlim_D {\rm Hom}_{\widehat{\mathcal C}_J} (\ell (X_{\bullet}) , \ell (Y)) \, .
$$

Comme le foncteur $\ell : {\mathcal C} \to \widehat{\mathcal C}_J$ est pleinement fid\`ele par hypoth\`ese, cela se r\'e\'ecrit
$$
{\rm Hom}_{\mathcal C} (X,Y) = \varprojlim_D {\rm Hom}_{\mathcal C} (X_{\bullet} , Y)
$$
qui signifie
$$
X = \varinjlim_D X_{\bullet} \qquad \mbox{dans} \quad {\mathcal C} \, .
$$

R\'eciproquement, si $X$ est une colimite dans ${\mathcal C}$ du $D$-diagramme $X_{\bullet}$, le morphisme canonique de $\widehat{\mathcal C}_J$
$$
\varinjlim_D \ell (X_{\bullet}) \longrightarrow \ell (X)
$$
est un isomorphisme si et seulement si sont v\'erifi\'ees les deux conditions suivantes, comme il r\'esulte des propositions \ref{propIII33} et \ref{propIII34}:

\medskip

$
\left\{ \begin{matrix}
({\rm B}') &\mbox{La famille des morphismes de $\widehat{\mathcal C}_J$} \hfill \\
{ \ } \\
&\ell (X_d) \longrightarrow \ell (X) \, , \qquad d \in {\rm Ob} (D) \, , \\
&\mbox{est globalement \'epimorphique.} \hfill \\
{ \ } \\
({\rm C}') &\mbox{Pour tous objets $d$ et $d'$ de $D$, consid\'erant la famille des sous-carquois finis $D_n$ de $D$ de la forme} \hfill \\
{ \ } \\
&d = d_0 \longleftrightarrow d_1 \longleftrightarrow \cdots \longleftrightarrow d_n = d' \\
{ \ } \\
&\mbox{o\`u $d=d_0 , d_1 , \cdots , d_{n-1}$, $d_n = d'$ sont des objets de $D$ et chaque $d_{i-1} \leftrightarrow d_i$ est une fl\`eche} \hfill \\
&\mbox{$d_{i-1} \to d_i$ ou $d_i \to d_{i-1}$ de $D$, la famille des morphismes} \hfill \\
{ \ } \\
&\displaystyle \varprojlim_{D_n} \ell (X_{\bullet}) \longrightarrow \ell (X_d) \times_{\ell (X)} \ell (X_{d'}) \\
&\mbox{est globalement \'epimorphique.} \hfill
\end{matrix} \right.
$

\bigskip

Or, d'apr\`es le lemme \ref{lemIII53}, la condition (B$'$) ci-dessus \'equivaut \`a la condition (B) de l'\'enonc\'e.

\smallskip

Quant \`a la condition (C$'$), elle signifie que pour tout objet $U$ de ${\mathcal C}$ muni d'un morphisme de $\widehat{\mathcal C}$
$$
y(U) \longrightarrow y(X_d) \times_{y(X)} y(X_{d'}) \, ,
$$
il existe une famille $J$-couvrante de morphismes de ${\mathcal C}$
$$
U' \longrightarrow U
$$
dont les images par $y : {\mathcal C} \hookrightarrow \widehat{\mathcal C}$ s'inscrivent dans des carr\'es commutatifs de $\widehat{\mathcal C}$ de la forme:
$$
\xymatrix{
y(U') \ar[d] \ar[r] &y(U) \ar[d] \\
\displaystyle \varprojlim_{D_n} \, y(X_{\bullet}) \ar[r] &y(X_d) \times_{y(X)} y(X_{d'})
}
$$
C'est la condition (C) de l'\'enonc\'e. 

\end{demo}

\bigskip

Ce th\'eor\`eme conduit \`a proposer la d\'efinition suivante:

\begin{defn}\label{defIII97}

Soit ${\mathcal C}$ une cat\'egorie essentiellement petite.

\smallskip

Etant donn\'e un carquois $D$, une colimite $X$ dans ${\mathcal C}$ d'un $D$-diagramme $X_{\bullet}$ de ${\mathcal C}$ sera appel\'ee ``g\'eom\'etrique'' s'il existe une topologie sous-canonique $J$ de ${\mathcal C}$ telle que le foncteur $\ell$
$$
\ell : {\mathcal C} \xhookrightarrow{ \ { \ } \ } \widehat{\mathcal C}_J
$$
induise un isomorphisme de $\widehat{\mathcal C}_J$
$$
\varinjlim_D \ell (X_{\bullet}) \xrightarrow{ \ \sim \ } \ell (X) \, .
$$
\end{defn}

\begin{remarksqed}
\begin{listeisansmarge}
\item D'apr\`es le th\'eor\`eme \ref{thmIII96}, la colimite $X$ du $D$-diagramme $X_{\bullet}$ de ${\mathcal C}$ est ``g\'eom\'etrique'' si et seulement si les conditions (B) et (C) de ce th\'eor\`eme sont v\'erifi\'ees par une certaine topologie sous-canonique $J$.

\smallskip

Elles le sont alors par toute topologie sous-canonique plus fine que $J$, en particulier par la topologie canonique de ${\mathcal C}$.

\medskip

\item Si les morphismes $X_{d'} \to X_d$ de $X_{\bullet}$ sont carrables ainsi que les projections $X_d \to X$, la condition (C) \'equivaut \`a demander que pour tous objets $d,d'$ de $D$, la famille des morphismes de ${\mathcal C}$
$$
\varprojlim_{D_n} X_{\bullet} \longrightarrow X_d \times_X X_{d'}
$$
index\'es par les sous-diagrammes finis de $D$ de la forme
$$
d = d_0 \longleftrightarrow d_1 \longleftrightarrow \cdots \longleftrightarrow d_{n-1} \longleftrightarrow d_n = d'
$$
est $J$-couvrante.

\medskip

\item Si $X$ est une colimite g\'eom\'etrique d'un $D$-diagramme $X_{\bullet}$ de ${\mathcal C}$, le foncteur
$$
y(X) = {\rm Hom} (\bullet , X)
$$
peut \^etre calcul\'e \`a partir des foncteurs
$$
y(X_d) = {\rm Hom} (\bullet , X_d) \, , \qquad d \in {\rm Ob} (D) \, ,
$$
comme le transform\'e par le foncteur de faisceautisation
$$
j^* : \widehat{\mathcal C} \longrightarrow \widehat{\mathcal C}_J
$$
du pr\'efaisceau
$$
\begin{matrix}
\displaystyle\varinjlim_D \, y(X_{\bullet}) : {\mathcal C}^{\rm op} &\longrightarrow &{\rm Ens} \, , \hfill \\
\hfill Y &\longmapsto &\displaystyle\varinjlim_D {\rm Hom} (Y,X_{\bullet})
\end{matrix}
$$
pour n'importe quelle topologie sous-canonique $J$ de ${\mathcal C}$ qui satisfait les conditions (B) et (C) du th\'eor\`eme~\ref{thmIII96}.

\smallskip

Cela rend la notion de colimite ``g\'eom\'etrique'' d'autant plus int\'eressante que, dans le cas de simples colimites dans ${\mathcal C}$,
$$
X = \varinjlim_D X_{\bullet} \, ,
$$
il n'existe pas de formule g\'en\'erale qui exprime le foncteur ${\rm Hom} (\bullet , X)$ \`a partir du $D$-diagramme des foncteurs ${\rm Hom} (\bullet , X_d)$. 

\end{listeisansmarge}
\end{remarksqed}

\bigskip

La notion de colimite g\'eom\'etrique s'applique en particulier \`a celle de quotient d'une action interne d'un groupe interne au sens suivant:

\begin{defn}\label{defIII98}

Soit ${\mathcal C}$ une cat\'egorie localement petite qui poss\`ede des produits finis.

\begin{listeimarge}

\item Une action (\`a gauche) interne d'un groupe interne $G$ de ${\mathcal C}$ consiste en un objet $X$ de ${\mathcal C}$ et un morphisme
$$
G \times X \xrightarrow{ \ a \ } X
$$
qui rend commutatif les deux carr\'es:
$$
\xymatrix{
G \times G \times X \ar[d]_-{m \times {\rm id}_X} \ar[rr]^-{{\rm id}_G \times a} &&G \times X \ar[d]^-a \\
G \times X \ar[rr]^-a &&X
} \qquad\qquad \xymatrix{
1 \times X \ar[d]^-{\wr} \ar[rr]^-{e \times {\rm id}_X} &&G \times X \ar[d]^-a \\
X \ar[rr]^-= &&X
}
$$

\item Un quotient d'une telle action est une colimite dans ${\mathcal C}$ du diagramme
$$
\raisebox{.7ex}{\xymatrix{ G \times X  \dar[r]^-{^{^{\mbox{\scriptsize$a$}}}}_-{p_2} &X}}
$$
constitu\'e du morphisme d'action
$$
a : G \times X \longrightarrow X
$$
et de la seconde projection
$$
p_2 : G \times X \longrightarrow X \, .
$$

\end{listeimarge}
\end{defn}

\begin{remarkqed}

Une action interne d'un groupe interne $G$ de ${\mathcal C}$ sur un objet $X$ est donc un morphisme
$$
a : G \times X \longrightarrow X
$$
tel que, pour tout objet $Y$ de ${\mathcal C}$, l'application induite
$$
{\rm Hom} (Y,G) \times {\rm Hom} (Y,X) \longrightarrow {\rm Hom} (Y,X)
$$
d\'efinisse une action du groupe ${\rm Hom} (Y,G)$ sur l'ensemble ${\rm Hom} (Y,X)$. 
\end{remarkqed}

\medskip

On d\'eduit du th\'eor\`eme \ref{thmIII96} et de la d\'efinition \ref{defIII97}:

\begin{cor}\label{corIII99}

Soit ${\mathcal C}$ une cat\'egorie essentiellement petite qui poss\`ede des produits finis.

\smallskip

Un quotient
$$
G \backslash X = \varinjlim \, (G \times X  \rightrightarrows X)
$$
d'un objet $X$ de ${\mathcal C}$ par l'action d'un groupe interne $G$ est g\'eom\'etrique si et seulement si il existe une topologie sous-canonique $J$ de ${\mathcal C}$ pour laquelle:

\medskip

$
\left\{ \begin{matrix}
({\rm B}) &\mbox{Le morphisme} \hfill \\
&\qquad\qquad\qquad\qquad X \longrightarrow G \backslash X \\
&\mbox{est $J$-couvrant.} \hfill \\
{ \ } \\
({\rm C}) &\mbox{Pour tout carr\'e commutatif de ${\mathcal C}$} \hfill \\
{ \ } \\
&\qquad\qquad\qquad\qquad \xymatrix{
U \ar[d] \ar[r] &X \ar[d] \\
X \ar[r] &G \backslash X
} \\
{ \ } \\
&\mbox{il existe une famille $J$-couvrante de morphismes} \hfill \\
{ \ } \\
&\qquad\qquad\qquad\qquad U' \longrightarrow U \\
&\mbox{tels que les compos\'es} \hfill \\
{ \ } \\
&\qquad\qquad\qquad\qquad U' \longrightarrow U \longrightarrow X \times X \\
{ \ } \\
&\mbox{se factorisent \`a travers le morphisme} \hfill \\
{ \ } \\
&\qquad\qquad\qquad\qquad G \times X \xrightarrow{ \ (a,p_2) \ } X \times X \, .
\end{matrix} \right.
$

\bigskip

Dans ce cas, le foncteur ${\rm Hom} (\bullet , G \backslash X)$ est le transform\'e du pr\'efaisceau
$$
\begin{matrix}
{\mathcal C}^{\rm op} &\longrightarrow &{\rm Ens} \, , \hfill \\
\hfill Y &\longmapsto &{\rm Hom} (Y,G) \backslash {\rm Hom} (Y,X)
\end{matrix}
$$
par le foncteur de faisceautisation 
$$
j^* : \widehat{\mathcal C} \to \widehat{\mathcal C}_J \, .
$$
\end{cor}

\begin{remark}

Si le morphisme $X \to G \backslash X$ est carrable, la condition (C) \'equivaut \`a:

\medskip

$
\left\{ \begin{matrix}
({\rm C}') &\mbox{Le morphisme} \hfill \\
&\qquad\qquad\qquad\qquad\qquad G \times X \xrightarrow{ \ (a,p_2) \ } X \times_{G \backslash X} X \\
&\mbox{est $J$-couvrant.} \hfill
\end{matrix} \right.
$
\end{remark}

\bigskip

\begin{demo}

Cela r\'esulte du th\'eor\`eme \ref{thmIII96} puisque, pour tout objet $Y$ de ${\mathcal C}$, l'image de l'application
$$
{\rm Hom} (Y,G) \times {\rm Hom} (Y,X) \longrightarrow {\rm Hom} (Y,X) \times {\rm Hom} (Y,X)
$$
est une relation d'\'equivalence de l'ensemble ${\rm Hom} (Y,X)$. 

\end{demo}

\section{Enrichissements g\'eom\'etriques de cat\'egories}\label{sec310}

\subsection{Une pr\'esentation alternative des cat\'egories g\'eom\'etriques}\label{subsec3101}

\medskip

La notion de sous-cat\'egorie g\'eom\'etrique de la cat\'egorie ${\rm Top}_{\rm an}$ des espaces annel\'es a \'et\'e introduite dans la d\'efinition \ref{defI58} du paragraphe \ref{subsec153}.

\smallskip

On d\'eduit du ``lemme de comparaison'' de Grothendieck:

\begin{prop}\label{propIII101}

Soit ${\mathcal G}$ une sous-cat\'egorie g\'eom\'etrique de la cat\'egorie ${\rm Top}_{\rm an}$ des espaces annel\'es.

\smallskip

Soit ${\mathcal C}$ une sous-cat\'egorie pleine de ${\mathcal G}$ qui est essentiellement petite et telle que tout objet de ${\mathcal G}$ admette un recouvrement par des ouverts qui sont isomorphes \`a des objets de ${\mathcal C}$.

\smallskip

Soit $J$ la topologie de ${\mathcal C}$ pour laquelle une famille de morphismes est couvrante si elle admet une sous-famille compos\'ee d'immersions ouvertes
$$
U_i \xhookrightarrow{ \ { \ } \ } U
$$
dont la r\'eunion des images est \'egale \`a $U$.

\smallskip

Alors:

\begin{listeimarge}

\item Le foncteur
$$
\begin{matrix}
{\mathcal G} &\longrightarrow &\widehat{\mathcal C} \, , \hfill \\
X &\longmapsto &{\rm Hom} (\bullet , X)
\end{matrix}
$$
identifie ${\mathcal G}$ \`a une sous-cat\'egorie pleine de $\widehat{\mathcal C}_J$.

\medskip

\item Un objet $F$ de $\widehat{\mathcal C}_J$ est isomorphe \`a l'image par le foncteur pleinement fid\`ele
$$
{\mathcal G} \xhookrightarrow{ \ { \ } \ } \widehat{\mathcal C}_J
$$
d'un objet de ${\mathcal G}$ si et seulement si il existe une famille d'objets $U_i$ de ${\mathcal C}$ et de monomorphismes de faisceaux
$$
{\rm Hom} (\bullet , U_i) \xhookrightarrow{ \ { \ } \ } F
$$
tels que
\begin{enumerate}
\item[$\bullet$] chaque 
$$
{\rm Hom} (\bullet , U_i) \xhookrightarrow{ \ { \ } \ } F
$$
est une immersion ouverte au sens que, pour tout objet $V$ de ${\mathcal C}$ muni d'un \'el\'ement de $F(V)$, le sous-objet
$$
{\rm Hom} (\bullet , U_i) \times_F {\rm Hom} (\bullet , V) \xhookrightarrow{ \ { \ } \ } {\rm Hom} (\bullet , V)
$$
est une r\'eunion de sous-objets
$$
{\rm Hom} (\bullet , V') \xhookrightarrow{ \ { \ } \ } {\rm Hom} (\bullet , V)
$$
associ\'es \`a des immersions ouvertes
$$
V' \xhookrightarrow{ \ { \ } \ }  V
$$
de ${\mathcal C}$,
\item[$\bullet$] le faisceau $F$ est la r\'eunion de ses sous-objets
$$
{\rm Hom} (\bullet , U_i) \, .
$$
\end{enumerate}
\end{listeimarge}
\end{prop}

\begin{remarks}
\begin{listeisansmarge}
\item Si ${\mathcal G}$ est la cat\'egorie des vari\'et\'es diff\'erentielles de classe $C^k$, $k \geq 1$ [resp. des vari\'et\'es analytiques], on peut prendre pour ${\mathcal C}$ la cat\'egorie des ouverts des ${\mathbb R}^n$ [resp. des $C^n$] et des applications de classe $C^k$ [resp. des applications holomorphes] entre ces ouverts, ou bien n'importe quelle sous-cat\'egorie pleine qui est dense.

\smallskip

Par exemple, on peut prendre pour ${\mathcal C}$ la cat\'egorie constitu\'ee de la boule unit\'e de chaque ${\mathbb R}^n$ [resp. de chaque ${\mathbb C}^n$] et des applications de classe $C^k$ [resp. des applications holomorphes] entre ces boules.

\smallskip

Les vari\'et\'es diff\'erentielles de classe $C^k$ [resp. les vari\'et\'es holomorphes] peuvent donc \^etre d\'efinies comme les faisceaux sur une telle cat\'egorie ${\mathcal C}$ munie de la topologie des recouvrements ouverts, qui admettent un recouvrement par des sous-faisceaux ouverts et repr\'esentables.

\medskip

\item Si ${\mathcal G}$ est la cat\'egorie des sch\'emas localement de pr\'esentation finie sur un sch\'ema affine de base $S$, on peut prendre pour ${\mathcal C}$ la cat\'egorie des sch\'emas affines de pr\'esentation finie sur $S$.

\smallskip

Ainsi, les sch\'emas localement de pr\'esentation finie sur $S$ apparaissent comme des faisceaux sur la cat\'egorie essentiellement petite des sch\'emas affines de pr\'esentation finie sur $S$.

\smallskip

La proposition recoupe dans ce cas le lemme \ref{lemI924} dont la preuve avait \'et\'e laiss\'ee au lecteur.
\end{listeisansmarge}
\end{remarks}
\medskip

\begin{demo}
\begin{listeisansmarge}
\item Soit ${\mathcal C}'$ n'importe quelle sous-cat\'egorie pleine de ${\mathcal G}$ qui contient ${\mathcal C}$, est essentiellement petite et contient tout ouvert de n'importe quel objet de ${\mathcal C}'$.

\smallskip

Soit $J'$ la topologie des recouvrements ouverts sur ${\mathcal C}'$.

\smallskip

D'apr\`es la proposition \ref{propII512}, la topologie $J'$ de ${\mathcal C}'$ est sous-canonique. Autrement dit, le foncteur
$$
\ell : {\mathcal C}' \longrightarrow \widehat{\mathcal C}'_{J'}
$$
est pleinement fid\`ele.

\smallskip

Or, d'apr\`es la proposition \ref{propIII14}, le foncteur de restriction
$$
\widehat{\mathcal C}'_{J'} \longrightarrow \widehat{\mathcal C}_J
$$
est une \'equivalence de cat\'egories.

\smallskip

Donc le foncteur compos\'e
$$
\begin{matrix}
{\mathcal C}' &\longrightarrow &\widehat{\mathcal C}_J \, , \hfill \\
X &\longmapsto &{\rm Hom} (\bullet , X)
\end{matrix}
$$
est pleinement fid\`ele.

\smallskip

Comme ${\mathcal C}'$ peut \^etre choisie assez grande pour contenir n'importe quelle paire d'objets de ${\mathcal G}$ que l'on consid\`ere, on conclut comme voulu que le foncteur
$$
\begin{matrix}
{\mathcal G} &\longrightarrow &\widehat{\mathcal C}_J \hfill \\
X &\longmapsto &{\rm Hom} (\bullet , X)
\end{matrix}
$$
est pleinement fid\`ele.

\medskip

\item Si $X$ est un objet de la sous-cat\'egorie g\'eom\'etrique ${\mathcal G}$ de ${\rm Top}_{\rm an}$, il peut \^etre recouvert par des immersions ouvertes
$$
U_i \xhookrightarrow{ \ { \ } \ } X
$$
d'objets $U_i$ de ${\mathcal C}$.

\smallskip

Le faisceau ${\rm Hom} (\bullet , X)$ est r\'eunion des sous-faisceaux ${\rm Hom} (\bullet , U_i)$.

\smallskip

De plus, pour tout morphisme $V \to X$ dont la source $V$ est un objet de ${\mathcal C}$, les faisceaux ${\rm Hom} (\bullet , U_i) \times_{{\rm Hom} (\bullet , X)} {\rm Hom} (\bullet , V)$ sont repr\'esentables par des ouverts de $V$, lesquels s'\'ecrivent comme des r\'eunions d'ouverts isomorphes \`a des objets de ${\mathcal C}$.

\smallskip

R\'eciproquement, consid\'erons un objet $F$ de $\widehat{\mathcal C}_J$ qui est r\'eunion de sous-objets repr\'esentables ${\rm Hom} (\bullet , U_i)$ qui sont des ouverts au sens de l'\'enonc\'e.

\smallskip

Pour tous indices $i,j$, le produit fibr\'e
$$
{\rm Hom} (\bullet , U_i) \times_F {\rm Hom} (\bullet , U_j)
$$
est un sous-objet de $F$ qui est r\'eunion de sous-objets repr\'esentables par des objets de ${\mathcal C}$ munis d'immersions ouvertes dans $U_i$ et dans $U_j$. Donc il est repr\'esentable par un ouvert $U_{i,j}$ de $U_i$ et de $U_j$.

\smallskip

Le carr\'e de $\widehat{\mathcal C}_J$
$$
\xymatrix{
\displaystyle\coprod_{i,j} {\rm Hom} (\bullet , U_{i,j}) \ar[d] \ar[r] &\displaystyle\coprod_j {\rm Hom} (\bullet , U_j) \ar[d] \\
\displaystyle\coprod_i {\rm Hom} (\bullet , U_i) \ar[r] &F
}
$$
est \`a la fois cart\'esien et cocart\'esien.

\smallskip

Soit $X$ la colimite dans Top du diagramme
$$
\coprod_{i,j} U_{i,j} \rightrightarrows \coprod_i U_i \, .
$$
Elle s'inscrit dans un carr\'e de Top
$$
\xymatrix{
\displaystyle\coprod_{i,j} U_{i,j} \ar[d] \ar[r] &\displaystyle\coprod_j U_i \ar[d] \\
\displaystyle\coprod_i U_i \ar[r] &X
}
$$
qui est \`a la fois cart\'esien et cocart\'esien, et les applications
$$
U_i \longrightarrow X
$$
sont des immersions ouvertes.

\smallskip

Les faisceaux de structure ${\mathcal O}_{U_i}$ des objets $U_i$ de ${\rm Top}_{\rm an}$ co{\"\i}ncident sur les intersections $U_{i,j} = U_i \times_X U_j$, donc d\'efinissent un unique faisceau d'anneaux ${\mathcal O}_X$ sur $X$ dont la restriction \`a chaque $U_i$ est ${\mathcal O}_{U_i}$.

\smallskip

Alors l'objet $(X,{\mathcal O}_X)$ de ${\rm Top}_{\rm an}$ est un objet de ${\mathcal G}$ puisqu'il admet un recouvrement ouvert par des objets de ${\mathcal C}$. Il repr\'esente le faisceau $F$ sur le site $({\mathcal C},J)$.

\smallskip

C'est ce que l'on voulait. 
\end{listeisansmarge}
\end{demo}

\subsection{Construction de cat\'egories g\'eom\'etriques}\label{subsec3102}

\medskip

R\'eciproquement, on a le r\'esultat suivant qui fournit un proc\'ed\'e g\'en\'eral de construction de sous-cat\'egories g\'eom\'etriques de ${\rm Top}_{\rm an}$:

\begin{thm}\label{thmIII102}

Soit ${\mathcal C}$ une cat\'egorie essentiellement petite munie d'un foncteur fid\`ele
$$
G : {\mathcal C} \longrightarrow {\rm Top}_{\rm an}
$$
tel que:

\medskip

$
\left\{ \begin{matrix}
\bullet &\mbox{un morphisme $U \xrightarrow{ \ f \ } V$ de ${\mathcal C}$ est un isomorphisme si $G(f)$ est un isomorphisme de ${\rm Top}_{\rm an}$,} \hfill \\
{ \ } \\
\bullet &\mbox{pour tout objet $U$ de ${\mathcal C}$, tout ouvert de $G(U)$ est r\'eunion des images d'immersions ouvertes} \hfill \\
{ \ } \\
&G(U_i) \longrightarrow G(U) \\
{ \ } \\
&\mbox{qui sont les transform\'ees par $G$ de morphismes de ${\mathcal C}$} \hfill \\
{ \ } \\
&U_i \longrightarrow U \, , \\
{ \ } \\
\bullet &\mbox{un morphisme $U \to V$ de ${\mathcal C}$ se factorise \`a travers un morphisme $V' \to V$ qui induit une} \hfill \\
&\mbox{immersion ouverte $G(V') \to G(V)$ si et seulement si $G(U) \to G(V)$ se factorise \`a travers $G(V')$,} \hfill \\
{ \ } \\
\bullet &\mbox{pour tous objets $U,V$ de ${\mathcal C}$ et tous morphismes} \hfill \\
{ \ } \\
&U_i \longrightarrow U \, , \qquad i \in I \, , \\
{ \ } \\
&\qquad \qquad \qquad U_{i,j,k} \longrightarrow U_i \qquad \mbox{et} \qquad U_{i,j,k} \longrightarrow U_j \\
{ \ } \\
&\mbox{dont les transform\'ees par $G$ sont des immersions ouvertes v\'erifiant} \hfill \\
{ \ } \\
&\displaystyle\bigcup_{i \in I} G(U_i) = G(U) \, , \\
{ \ } \\
&\qquad \qquad \qquad \quad \displaystyle\bigcup_k G(U_{i,j,k}) = G(U_i) \times_{G(U)} G(U_j) \, , \\
&\mbox{se donner un morphisme de ${\mathcal C}$} \hfill \\
&U \longrightarrow V \\
{ \ } \\
&\mbox{\'equivaut \`a se donner une famille de morphismes} \hfill \\
{ \ } \\
&U_i \longrightarrow V \\
{ \ } \\
&\mbox{dont les compos\'es avec les $U_{i,j,k} \to U_i$ et $U_{i,j,k} \to U_j$ co{\"\i}ncident.} \hfill
\end{matrix} \right.
$

\medskip

Alors:

\begin{listeimarge}

\item Il existe une topologie $J$ de ${\mathcal C}$ pour laquelle une famille de morphismes vers un objet $U$ est couvrante si elle contient une sous-famille de morphismes
$$
U_i \longrightarrow U \, , \qquad i \in I \, ,
$$
dont les transform\'es par $G$ sont des immersions ouvertes dont la r\'eunion est $G(U)$.

\medskip

\item Cette topologie est sous-canonique.

\medskip

\item Soit ${\mathcal G}$ la sous-cat\'egorie pleine de ${\mathcal C}$ constitu\'ee des faisceaux $F$ pour lesquels existe une famille d'objets $U_i$ de ${\mathcal C}$ et de monomorphismes de faisceaux
$$
{\rm Hom} (\bullet , U_i) \xhookrightarrow{ \ { \ } \ } F
$$
tels que

\medskip

$
\left\{ \begin{matrix}
\bullet &\mbox{chaque ${\rm Hom} (\bullet , U_i) \hookrightarrow F$ est une immersion ouverte au sens que, pour tout objet $V$ de ${\mathcal C}$} \hfill \\
&\mbox{muni d'un \'el\'ement de $F(V)$, le sous-objet} \hfill \\
{ \ } \\
&{\rm Hom} (\bullet , U_i) \times_F {\rm Hom} (\bullet , V) \xhookrightarrow{ \ { \ } \ } {\rm Hom} (\bullet , V) \\
{ \ } \\
&\mbox{est r\'eunion de sous-objets repr\'esentables} \hfill \\
{ \ } \\
&{\rm Hom} (\bullet , V') \xhookrightarrow{ \ { \ } \ } {\rm Hom} (\bullet , V) \\
&\mbox{tels que les morphismes} \hfill \\
&G(V') \longrightarrow G(V) \\
&\mbox{soient des immersions ouvertes,} \hfill \\
{ \ } \\
\bullet &\mbox{le faisceau $F$ est la r\'eunion de ses sous-objets} \hfill \\
{ \ } \\
&{\rm Hom} (\bullet , U_i) \, .
\end{matrix} \right.
$

\bigskip

Alors le foncteur fid\`ele
$$
G : {\mathcal C} \longrightarrow {\rm Top}_{\rm an}
$$
se prolonge en une \'equivalence de cat\'egories
$$
{\mathcal G} \xrightarrow{ \ \sim \ } {\mathcal G}'
$$
vers une sous-cat\'egorie g\'eom\'etrique ${\mathcal G}'$ de ${\rm Top}_{\rm an}$ dont tout objet est r\'eunion d'ouverts isomorphes aux images par $G$ d'objets de ${\mathcal C}$.
\end{listeimarge}
\end{thm}

\bigskip

\begin{demo}
\begin{listeisansmarge}
\item L'axiome de maximalit\'e est v\'erifi\'e car, pour tout objet $U$ de ${\mathcal C}$, $G$ transforme ${\rm id}_U : U \to U$ en ${\rm id}_{G(U)} : G(U) \to G(U)$.

\smallskip

Pour l'axiome de stabilit\'e, consid\'erons un morphisme de ${\mathcal C}$
$$
f : U \longrightarrow V
$$
et une famille de morphismes
$$
V_i \longrightarrow V
$$
dont les images par $G$ sont des immersions ouvertes qui recouvrent $G(V)$. Chaque $G(f)^{-1} (G(V_i))$ est un ouvert de $G(U)$ donc est recouvert par une famille d'immersions ouvertes $G(U_{i,j}) \to G(U)$ qui sont les transform\'ees par $G$ de morphismes $U_{i,j} \to U$. Chaque compos\'e $U_{i,j} \to U \to V$ se factorise en $U_{i,j} \to V_i \to V$ puisque son transform\'e $G(U_{i,j}) \to G(V)$ se factorise \`a travers l'ouvert $G(V_i)$. Ainsi, les $U_{i,j} \to U$ sont \'el\'ements de l'image r\'eciproque par $f$ du crible de $V$ engendr\'e par les $V_i \to V$ et ils induisent des immersions ouvertes $G(U_{i,j}) \hookrightarrow G(U)$ qui recouvrent $G(U)$.

\smallskip

Enfin, pour l'axiome de transitivit\'e, consid\'erons des morphismes de ${\mathcal C}$
$$
U_i \longrightarrow U
$$
qui induisent des immersions ouvertes $G(U_i) \hookrightarrow G(U)$ recouvrant $G(U)$ et, pour tout indice $i$, des morphismes de ${\mathcal C}$
$$
U_{i,j} \longrightarrow U_i
$$
qui induisent des immersions ouvertes $G(U_{i,j}) \hookrightarrow G(U_i)$ recouvrant $G(U_i)$. Alors les morphismes compos\'es
$$
U_{i,j} \longrightarrow U_i \longrightarrow U
$$
induisent des immersions ouvertes $G(U_{i,j}) \hookrightarrow G(U)$ recouvrant $G(U)$.

\medskip

\item En effet, la derni\`ere hypoth\`ese de l'\'enonc\'e signifie que si $U,V$ sont deux objets de ${\mathcal C}$ et $S$ est le crible sur $U$ engendr\'e par une famille de morphismes
$$
U_i \longrightarrow U \, , \qquad i \in I \, ,
$$
alors on a
$$
{\rm Hom} (U,V) = \varprojlim_{(U' \to U) \in S} {\rm Hom} (U',V) \, .
$$

\item Comme le foncteur $G : {\mathcal C} \to {\rm Top}_{\rm an}$ est fid\`ele, il existe une sous-cat\'egorie ${\mathcal C}'$ de ${\rm Top}_{\rm an}$ dont les objets sont les images par $G$ d'objets de ${\mathcal C}$ et dont les morphismes sont les images par $G$ de morphismes de ${\mathcal C}$. Le foncteur $G$ se factorise en une \'equivalence de cat\'egories
$$
{\mathcal C} \xrightarrow{ \ \sim \ } {\mathcal C}' \, .
$$
La topologie $J$ sur ${\mathcal C}$ correspond \`a la topologie des recouvrements ouverts $J'$ sur ${\mathcal C}'$ et on a encore une \'equivalence induite
$$
\widehat{\mathcal C}'_{J'} \xrightarrow{ \ \sim \ } \widehat{\mathcal C}_J \, .
$$

Soit alors ${\mathcal G}'$ la sous-cat\'egorie de ${\rm Top}_{\rm an}$ dont

\medskip

$
\left\{\begin{matrix}
\bullet &\mbox{les objets sont les espaces annel\'es $(X,{\mathcal O}_X)$ qui admettent un recouvrement ouvert par des espaces} \hfill \\
&\mbox{annel\'es $(U_i , {\mathcal O}_{U_i})$ qui sont isomorphes \`a des objets de ${\mathcal C}'$,} \hfill \\
{ \ } \\
\bullet &\mbox{les morphismes entre objets de ${\mathcal G}'$} \hfill \\
&(X,{\mathcal O}_X) \longrightarrow (Y,{\mathcal O}_Y) \\
{ \ } \\
&\mbox{sont les morphismes d'espaces annel\'es tels qu'existe une famille de carr\'es commutatifs} \hfill \\
{ \ } \\
&\xymatrix{
(U_{i,j} , {\mathcal O}_{U_{i,j}}) \ar@{_{(}->}[d] \ar[r] &(V_j , {\mathcal O}_{V_j}) \ar@{_{(}->}[d] \\
(X,{\mathcal O}_X) \ar[r] &(Y,{\mathcal O}_Y)
} \\
{ \ } \\
&\mbox{dont les fl\`eches verticales sont des immersions ouvertes qui recouvrent $X$ et $Y$, et les morphismes} \hfill \\
{ \ } \\
&(U_{i,j} , {\mathcal O}_{U_{i,j}}) \longrightarrow (V_j , {\mathcal O}_{V_j}) \\
&\mbox{sont des morphismes de ${\mathcal C}'$.} \hfill
\end{matrix} \right.
$

\bigskip

Il r\'esulte des hypoth\`eses faites que ${\mathcal G}'$ est une sous-cat\'egorie g\'eom\'etrique de ${\rm Top}_{\rm an}$ et que ${\mathcal C}'$ est une sous-cat\'egorie pleine de ${\mathcal G}'$.

\smallskip

Par construction, tout objet de ${\mathcal G}'$ admet un recouvrement ouvert par des objets isomorphes \`a des objets de ${\mathcal C}'$.

\smallskip

Donc la proposition \ref{propIII101} s'applique et le foncteur d'\'evaluation
$$
X \longmapsto {\rm Hom} (\bullet , X)
$$
d\'efinit une \'equivalence de cat\'egories
$$
{\mathcal G}' \xrightarrow{ \ \sim \ } \widehat{\mathcal C}'_{J'} \, .
$$
Son compos\'e avec l'\'equivalence
$$
\widehat{\mathcal C}'_{J'} \xrightarrow{ \ \sim \ } \widehat{\mathcal C}_J
$$
est encore une \'equivalence de cat\'egories
$$
{\mathcal G}' \xrightarrow{ \ \sim \ } \widehat{\mathcal C}_J \, .
$$
\end{listeisansmarge}
\end{demo}

\subsection{Enrichissement de cat\'egories par des quotients d'un type g\'eom\'etrique}\label{subsec3103}

\medskip

Rappelons que nous avons introduit dans la d\'efinition \ref{defII31} du paragraphe \ref{subsec231} les notions de ``classe g\'eom\'etrique de morphismes'' et de ``propri\'et\'e g\'eom\'etrique de recouvrement''. Elles permettent d'introduire un proc\'ed\'e g\'en\'eral d'enrichissement g\'eom\'etrique de cat\'egories par ajouts de quotients dont le type est sp\'ecifi\'e par la ``classe de morphismes'' et la ``propri\'et\'e de recouvrement'' choisies:

\begin{defn}\label{defIII103}

Soit ${\mathcal C}$ une cat\'egorie localement petite.

\smallskip

Soit ${\mathcal M}$ une ``classe g\'eom\'etrique'' de morphismes de ${\mathcal C}$.

\smallskip

Soit {\rm (R)} une ``notion g\'eom\'etrique de recouvrement'' par des \'el\'ements de la classe ${\mathcal M}$.

\smallskip

Alors:

\begin{listeimarge}

\item On dit qu'un morphisme de $\widehat{\mathcal C}$
$$
P' \longrightarrow P
$$
appartient \`a la classe ${\mathcal M}$ si, pour tout objet $X$ de ${\mathcal C}$ et tout \'el\'ement $x \in P(X) = {\rm Hom} (y(X),P)$, le produit fibr\'e
$$
P' \times_P y(X)
$$
est repr\'esentable par un objet $X'$ de ${\mathcal C}$ tel que le morphisme
$$
X' \longrightarrow X
$$
soit \'el\'ement de la classe ${\mathcal M}$ de ${\mathcal C}$.

\medskip

\item On dit qu'une famille de morphismes de la classe ${\mathcal M}$ de $\widehat{\mathcal C}$
$$
P_i \longrightarrow P \, , \qquad i \in I \, ,
$$
poss\`ede la propri\'et\'e {\rm (R)} si, pour tout objet $X$ de ${\mathcal C}$ et tout \'el\'ement $x \in P(X)$, les morphismes
$$
P_i \times_P y(X) \longrightarrow y(X)
$$
sont repr\'esent\'es par une famille de morphismes de ${\mathcal M}$
$$
X_i \longrightarrow X \, , \qquad i \in I \, ,
$$
qui poss\`ede la propri\'et\'e {\rm (R)}.

\medskip

\item On dit qu'un objet $P$ de $\widehat{\mathcal C}$ est $({\mathcal M},R)$-constructible s'il existe dans $\widehat{\mathcal C}$ une famille de morphismes de la classe ${\mathcal M}$
$$
P_i \longrightarrow P
$$
dont les sources $P_i$ sont repr\'esentables par des objets $X_i$ de ${\mathcal C}$, et qui poss\`ede la propri\'et\'e (R).

\medskip

\item On appelle cat\'egorie enrichie de ${\mathcal C}$ par $({\mathcal M},R)$ la sous-cat\'egorie pleine de $\widehat{\mathcal C}$ constitu\'ee des objets de $\widehat{\mathcal C}$ qui sont $({\mathcal M},R)$-constructibles.
\end{listeimarge}
\end{defn}

\begin{remarkqed}

Si la cat\'egorie ${\mathcal C}$ est essentiellement petite, la cat\'egorie $\widehat{\mathcal C}$ est localement petite.

\smallskip

La classe ${\mathcal M}$ dans $\widehat{\mathcal C}$ est alors une classe g\'eom\'etrique de morphismes de $\widehat{\mathcal C}$, et la propri\'et\'e (R) des familles de morphismes $(P_i \to P)_{i \in I}$ \'el\'ements de ${\mathcal M}$ est une notion g\'eom\'etrique de recouvrement.

\smallskip

Enfin, la cat\'egorie enrichie de ${\mathcal C}$ par $({\mathcal M},R)$ est essentiellement petite comme ${\mathcal C}$, et elle est close sous l'op\'eration d'enrichissement, au sens qu'elle est \'egale \`a sa propre cat\'egorie enrichie par $({\mathcal M} , R)$. 
\end{remarkqed}

\medskip

On d\'eduit du lemme de comparaison de Grothendieck:

\begin{prop}\label{propIII104}

Soit ${\mathcal C}$ une cat\'egorie essentiellement petite.

\smallskip

Soit ${\mathcal M}$ une ``classe g\'eom\'etrique'' de morphismes de ${\mathcal C}$.

\smallskip

Soit {\rm (R)} une ``notion g\'eom\'etrique de recouvrement'' par des \'el\'ements de la classe ${\mathcal M}$.

\smallskip

Enfin, soit ${\mathcal C}'$ la cat\'egorie enrichie de ${\mathcal C}$ par $({\mathcal M},R)$ et soit $J'_R$ la topologie de ${\mathcal C}'$ d\'efinie par la notion g\'eom\'etrique de recouvrement {\rm (R)} par des \'el\'ements de ${\mathcal M}$ dans ${\mathcal C}'$.

\smallskip

Alors:

\begin{listeimarge}

\item La sous-cat\'egorie pleine ${\mathcal C}$ de ${\mathcal C}'$ est dense pour la topologie $J_R$.

\medskip

\item La topologie de ${\mathcal C}$ induite par la topologie $J'_R$ de ${\mathcal C}'$ n'est autre que la topologie $J_R$ de ${\mathcal C}$ d\'efinie par la notion g\'eom\'etrique de recouvrement {\rm (R)} par des \'el\'ements de ${\mathcal M}$ dans ${\mathcal C}$.

\medskip

\item Le foncteur de restriction des faisceaux
$$
\widehat{\mathcal C}'_{J'_R} \longrightarrow \widehat{\mathcal C}_{J_R}
$$
est une \'equivalence de cat\'egories.
\end{listeimarge}
\end{prop}

\begin{remarks}
\begin{listeisansmarge}
\item En particulier, si le site $({\mathcal C} , J_R)$ est muni d'un faisceau d'anneaux $A$, celui-ci s'\'etend de mani\`ere unique en un faisceau d'anneaux $A$ sur le site $({\mathcal C}',J')$.

\smallskip

Dans ce cas, le petit ou le gros topos associ\'e \`a tout objet $P$ de ${\mathcal C}'$ devient un topos annel\'e, tout comme le petit ou le gros topos associ\'e \`a un objet $X$ de ${\mathcal C}$.

\medskip

\item Si par exemple ${\mathcal C}$ est la cat\'egorie des sch\'emas de pr\'esentation finie sur un sch\'ema de base $S$, ${\mathcal M}$ est la classe des morphismes \'etales et (R) est la propri\'et\'e d'une famille de morphismes \'etales d'\^etre surjective, la cat\'egorie enrichie ${\mathcal C}'$ contient des objets qui ne sont plus n\'ecessairement des sch\'emas. On les appelle des espaces alg\'ebriques.
\end{listeisansmarge}
\end{remarks}

\medskip

\begin{demo}
\begin{listeisansmarge}
\item est \'evident sur la d\'efinition des objets $({\mathcal M},R)$-constructibles.

\medskip

\item est aussi \'evident sur les d\'efinitions.

\medskip

\item r\'esulte alors de la proposition \ref{propIII14}. 
\end{listeisansmarge}
\end{demo}


%% file: Chapitre4_num.tex







\vglue 15mm

\chapter{G\'eom\'etrie des topos}\label{chap4}

\section{Morphismes de topos}\label{sec41}

\subsection{Des applications continues aux morphismes de topos}\label{subsec411}

\medskip

On a vu au chapitre \ref{chap1} qu'\`a tout espace topologique $X$ est associ\'ee la cat\'egorie ${\mathcal E}_X$ des faisceaux sur $X$. Avec le vocabulaire du chapitre \ref{chap2}, c'est la cat\'egorie des faisceaux sur le site constitu\'e de la cat\'egorie $O(X)$ des ouverts de $X$ et de sa topologie canonique.  

\smallskip

On a vu \'egalement que toute application continue entre espaces topologiques
$$
f : X \to Y
$$
d\'efinit un ``foncteur d'image directe'' par $f$
$$
f_* : {\mathcal E}_X \longrightarrow {\mathcal E}_Y
$$
qui transforme les faisceaux sur $X$ en faisceaux sur $Y$.

\smallskip

Montrons maintenant:

\begin{prop}\label{propIV11}

Soit $f : X \to Y$ une application continue entre deux espaces topologiques $X$ et $Y$.

\smallskip

Alors:

\begin{listeimarge}

\item Le foncteur d'image directe par $f$
$$
f_* : {\mathcal E}_X \longrightarrow {\mathcal E}_Y
$$
admet pour adjoint \`a gauche
$$
f^* : {\mathcal E}_Y \longrightarrow {\mathcal E}_X
$$
le compos\'e du foncteur
$$
\begin{matrix}
{\mathcal E}_Y &\longrightarrow &\widehat{O(X)} \, , \hfill \\
\hfill G &\longmapsto &\left[\begin{matrix}
O(X)^{\rm op} &\longrightarrow &{\rm Ens} \hfill \\
\hfill U &\longmapsto &\!\!\!\!\!\!\displaystyle \varinjlim_{V \in O(Y) \atop f^{-1} (V) \supset \, U} \!\!\!\!\!G(V)
\end{matrix} \right]
\end{matrix}
$$
et du foncteur de faisceautisation des pr\'efaisceaux sur $O(X)$ 
$$
j^* : \widehat{O(X)} \longrightarrow {\mathcal E}_X \, .
$$

\item Ce foncteur adjoint \`a gauche
$$
f^* : {\mathcal E}_Y \longrightarrow {\mathcal E}_X
$$
respecte non seulement les colimites arbitraires mais aussi les limites finies.
\end{listeimarge}
\end{prop}

\begin{demo}
\begin{listeisansmarge}
\item L'application continue $f : X \to Y$ d\'efinit un foncteur
$$
\begin{matrix}
f^{-1} : O(Y) &\longrightarrow &O(X) \, , \hfill \\
\hfill V &\longmapsto &f^{-1} V \, . \hfill
\end{matrix}
$$
La composition avec $\rho = f^{-1}$ d\'efinit un foncteur
$$
\rho^* : \widehat{O(X)} \longrightarrow \widehat{O(Y)}
$$
et le foncteur $f_* : {\mathcal E}_X \to {\mathcal E}_Y$ est la restriction de $\rho^*$ \`a la sous-cat\'egorie pleine ${\mathcal E}_X$ de $\widehat{O(X)}$. Autrement dit, on a un carr\'e commutatif:
$$
\xymatrix{
{\mathcal E}_X \ar@{_{(}->}[d] \ar[r]^{f_*} &{\mathcal E}_Y \ar@{_{(}->}[d] \\
\widehat{O(X)} \ar[r]^{\rho^*} &\widehat{O(Y)}
}
$$
D'apr\`es le corollaire \ref{corI108}, le foncteur
$$
\rho^* : \widehat{O(X)} \longrightarrow \widehat{O(Y)}
$$
admet pour adjoint \`a gauche le foncteur
$$
\begin{matrix}
\rho_! : \widehat{O(Y)} &\longrightarrow &\widehat{O(X)} \, , \hfill \\
\hfill P &\longmapsto &\rho_! P = \left[\begin{matrix}
O(X)^{\rm op} &\longrightarrow &{\rm Ens} \hfill \\
\hfill U &\longmapsto &\!\!\!\!\!\!\displaystyle \varinjlim_{V \in O(Y) \atop f^{-1} (V) \supset \, U} \!\!\!\!\!P(V)
\end{matrix} \right].
\end{matrix}
$$
D'autre part, le foncteur de plongement
$$
j_* : {\mathcal E}_X \xhookrightarrow{ \ { \ } \ } \widehat{O(X)}
$$
admet pour adjoint \`a gauche le foncteur de faisceautisation
$$
j^* : \widehat{O(X)} \longrightarrow {\mathcal E}_X \, .
$$
Donc le foncteur compos\'e
$$
{\mathcal E}_X \xrightarrow{ \ f_* \ } {\mathcal E}_Y \xhookrightarrow{ \ { \ } \ } \widehat{O(Y)}
$$
admet pour adjoint \`a gauche le foncteur compos\'e
$$
\widehat{O(Y)} \xrightarrow{ \ \rho_! \ } \widehat{O(X)} \xrightarrow{ \ j^* \ } {\mathcal E}_X \, .
$$
A fortiori, le foncteur
$$
{\mathcal E}_X \xrightarrow{ \ f_* \ } {\mathcal E}_Y
$$
admet pour adjoint \`a gauche la restriction du foncteur
$$
j^* \circ \rho_!
$$
\`a la sous-cat\'egorie pleine ${\mathcal E}_Y$ de $\widehat{O(Y)}$.

\medskip

\item Le foncteur
$$
\rho_! : \widehat{O(Y)} \longrightarrow \widehat{O(X)}
$$
respecte les limites finies car, pour tout objet $U$ de $O(X)$, la cat\'egorie des objets $V$ de $O(Y)$ tels que $f^{-1} (V) \supset U$ est filtrante.

\smallskip

En effet, cela entra{\^\i}ne d'apr\`es le lemme II.5.4 que pour tout tel ouvert $U$ de $X$, le foncteur
$$
\varinjlim_{V \in O(Y) \atop f^{-1} (V) \supset U}
$$
respecte les limites finies.

\smallskip

D'autre part, on sait que le foncteur de faisceautisation
$$
j^* : \widehat{O(X)} \longrightarrow {\mathcal E}_X
$$
respecte les limites finies.

\smallskip

Il en est donc de m\^eme du foncteur compos\'e
$$
j^* \circ \rho_! : \widehat{O(X)} \longrightarrow {\mathcal E}_X
$$
et de sa restriction \`a ${\mathcal E}_Y$
$$
f^* : {\mathcal E}_Y \longrightarrow {\mathcal E}_X \, .
$$
\end{listeisansmarge}
\end{demo}

\medskip

Cette proposition conduit \`a poser la d\'efinition suivante:

\begin{defn}\label{defIV12}

Etant donn\'es deux topos ${\mathcal E},{\mathcal E}'$, on appelle morphisme de topos
$$
f : {\mathcal E}' \longrightarrow {\mathcal E}
$$
une paire de foncteurs adjoints
$$
\left( {\mathcal E} \xrightarrow{ \ f^* \ } {\mathcal E}' , {\mathcal E}' \xrightarrow{ \ f_* \ } {\mathcal E} \right)
$$
dont la composante de gauche $f^*$ respecte non seulement les colimites arbitraires mais aussi les limites finies.
\end{defn}

\begin{remarksqed}
\begin{listeisansmarge}
\item Dans un morphisme de topos $f = (f^* , f_*)$, la composante de gauche $f^*$ est appel\'ee le foncteur d'image r\'eciproque et la composante de droite $f_*$ est appel\'ee le foncteur d'image directe.

\medskip

\item Dans un morphisme de topos $f = (f^* , f_*)$, chaque composante d\'etermine l'autre \`a unique isomorphisme pr\`es.

\medskip

\item Un foncteur entre deux topos
$$
{\mathcal E} \longrightarrow {\mathcal E}'
$$
est la composante d'image r\'eciproque $f^*$ d'un morphisme de topos
$$
f : {\mathcal E}' \longrightarrow {\mathcal E}
$$
si et seulement si il respecte \`a la fois les colimites arbitraires et les limites finies.

\smallskip

En effet, d'apr\`es le corollaire \ref{corIII75}, un tel foncteur ${\mathcal E} \to {\mathcal E}'$ admet un adjoint \`a droite si et seulement si il respecte les colimites. 

\medskip

\item Pour tout site $({\mathcal C} , J)$, l'adjoint \`a gauche $j^*$ du foncteur de plongement $j_* : \widehat{\mathcal C}_J \hookrightarrow \widehat{\mathcal C}$ respecte les limites finies.

\smallskip

Donc la paire $(j^* , j_*)$ est un morphisme de topos $\widehat{\mathcal C}_J \to \widehat{\mathcal C}$.

\medskip

\item Pour tout foncteur $\rho : {\mathcal C} \to {\mathcal D}$ entre deux cat\'egories essentiellement petites, le foncteur $\rho^* : \widehat{\mathcal D} \to \widehat{\mathcal C}$ de composition avec $\rho$ admet d'apr\`es le corollaire \ref{corI108} un adjoint \`a droite $\rho_*$ et un adjoint \`a gauche $\rho_!$. Il respecte les limites arbitraires et la paire $(\rho^* , \rho_*)$ est un morphisme de topos $\widehat{\mathcal C} \to \widehat{\mathcal D}$. 
\end{listeisansmarge}
\end{remarksqed}

\smallskip

Les morphismes de topos se composent naturellement:

\begin{defn}\label{defIV13}
\begin{listeimarge}
\item Le compos\'e de deux morphismes de topos
$$
f = (f^* , f_*) : {\mathcal E}_1 \longrightarrow {\mathcal E}_2
$$
et
$$
g = (g^* , g_*) : {\mathcal E}_2 \longrightarrow {\mathcal E}_3
$$
est la paire de foncteurs compos\'es
$$
g \circ f = (f^* \!\circ g^* , g_* \circ f_*) \, .
$$

\item Le morphisme identique d'un topos ${\mathcal E}$ est la paire de foncteurs
$$
{\rm id}_{\mathcal E} = \left( {\mathcal E} \xrightarrow{ \ {\rm id} \ } {\mathcal E} , {\mathcal E} \xrightarrow{ \ {\rm id} \ } {\mathcal E} \right) .
$$
\end{listeimarge}
\end{defn}

\begin{remarkqed}

La loi de composition des morphismes de topos ainsi d\'efinie est associative et on a pour tout morphisme de topos
$$
f : {\mathcal E}' \longrightarrow {\mathcal E}
$$
les identit\'es
$$
{\rm id}_{\mathcal E} \circ f = f = f \circ {\rm id}_{{\mathcal E}'} \, .
$$

Ainsi, les topos et leurs morphismes forment une cat\'egorie. 
\end{remarkqed}

\bigskip

La cat\'egorie des topos n'est pas localement petite. Autrement dit, si ${\mathcal E}$ et ${\mathcal E}'$ sont deux topos, les morphismes
$$
{\mathcal E}' \longrightarrow {\mathcal E}
$$
ne forment pas en g\'en\'eral un ensemble.

\smallskip

Cependant, nous allons voir au paragraphe suivant que ces morphismes
$$
{\mathcal E}' \longrightarrow {\mathcal E}
$$
forment naturellement une cat\'egorie et que cette cat\'egorie est localement petite.

\subsection{La cat\'egorie des morphismes entre deux topos}\label{subsec412}

\medskip

Il est naturel de poser:

\begin{defn}\label{defIV14}

Soient
$$
f,g : {\mathcal E}' \rightrightarrows {\mathcal E}
$$
deux morphismes
$$
f = (f^* , f_*) \qquad \mbox{et} \qquad g = (g^* , g_*)
$$
d'un topos ${\mathcal E}'$ dans un topos ${\mathcal E}$.

\smallskip

On appelle transformation de $f$ dans $g$ la donn\'ee d'une transformation naturelle de foncteurs
$$
f^* \longrightarrow g^*
$$
ou, ce qui revient au m\^eme, d'une transformation naturelle de foncteurs dans l'autre sens
$$
g_* \longrightarrow f_* \, .
$$
\end{defn}

\begin{remarkqed}

Se donner une transformation naturelle
$$
f^* \longrightarrow g^*
$$
\'equivaut \`a se donner une transformation naturelle
$$
g_* \longrightarrow f_*
$$
puisque les paires $(f^* , f_*)$ et $(g^* , g_*)$ sont adjointes.
 
\end{remarkqed}

\medskip

Cette d\'efinition est compl\'et\'ee par la suivante:

\begin{defn}\label{defIV15}

Soient ${\mathcal E}$ et ${\mathcal E}'$ deux topos.

\begin{listeimarge}

\item L'identit\'e d'un morphisme
$$
f = (f^* , f_*) : {\mathcal E}' \longrightarrow {\mathcal E}
$$
est la transformation naturelle identique
$$
{\rm id} : f^* \longrightarrow f^*
$$
ou, ce qui revient au m\^eme, la transformation naturelle adjointe
$$
{\rm id} : f_* \longrightarrow f_* \, .
$$

\item La compos\'ee de deux transformations entre morphismes de topos de ${\mathcal E}'$ dans ${\mathcal E}$
$$
(f^* , f_*) \longrightarrow (g^* , g_*) \longrightarrow (h^* , h_*)
$$
est la compos\'ee des transformations naturelles
$$
f^* \longrightarrow g^* \longrightarrow h^*
$$
ou, ce qui revient au m\^eme, la compos\'ee des transformations naturelles adjointes
$$
h_* \longrightarrow g_* \longrightarrow f_* \, .
$$
\end{listeimarge}
\end{defn}

Avec ces d\'efinitions, on a imm\'ediatement:

\begin{prop}\label{propIV16}

Soient ${\mathcal E}$ et ${\mathcal E}'$ deux topos.

\smallskip

Alors:

\begin{listeimarge}

\item Les morphismes de topos de ${\mathcal E}'$ dans ${\mathcal E}$ et leurs transformations forment une cat\'egorie que l'on notera
$$
[{\mathcal E}' , {\mathcal E}]_T \, .
$$

\item Cette cat\'egorie est localement petite.

\medskip

\item Le foncteur
$$
\begin{matrix}
\hfill [{\mathcal E}' , {\mathcal E}]_T &\longrightarrow &[{\mathcal E} , {\mathcal E}'] \, , \\
f = (f^* , f_*) &\longmapsto &f^* \hfill
\end{matrix}
$$
est pleinement fid\`ele.

\smallskip

Il d\'efinit une \'equivalence de la cat\'egorie des morphismes de topos ${\mathcal E}' \to {\mathcal E}$ vers la sous-cat\'egorie pleine de $[{\mathcal E} , {\mathcal E}']$ constitu\'ee des foncteurs
$$
\rho : {\mathcal E} \longrightarrow {\mathcal E}'
$$
qui respectent les colimites arbitraires et les limites finies.

\medskip

\item De m\^eme, le foncteur
$$
\begin{matrix}
\hfill [{\mathcal E}' , {\mathcal E}]_T^{\rm op} &\longrightarrow &[{\mathcal E}', {\mathcal E}] \\
f = (f^* , f_*) &\longmapsto &f_* \hfill
\end{matrix}
$$
est pleinement fid\`ele.

\smallskip

Il d\'efinit une \'equivalence de l'oppos\'ee de la cat\'egorie des morphismes de topos ${\mathcal E}' \to {\mathcal E}$ vers la sous-cat\'egorie pleine de $[{\mathcal E}',{\mathcal E}]$ constitu\'ee des foncteurs
$$
{\mathcal E}' \longrightarrow {\mathcal E}
$$
qui admettent un adjoint \`a gauche respectant les limites finies.
\end{listeimarge}
\end{prop}

\begin{demo}
\begin{listeisansmarge}
\item[(i)] En effet, la loi de composition des transformations telle qu'elle a \'et\'e pos\'ee dans la d\'efinition \ref{defIV15} est associative.

\smallskip

D'autre part, la compos\'ee de toute transformation
$$
(f^* , f_*) \longrightarrow (g^* , g_*)
$$
avec l'identit\'e de $(f^* , f_*)$ ou de $(g^* , g_*)$ est elle-m\^eme.

\smallskip

Cela suffit \`a faire de $[{\mathcal E}' , {\mathcal E} ]_T$ une cat\'egorie.

\medskip

\item[(iii)] La premi\`ere assertion est \'evidente sur la d\'efinition des transformations de morphismes de topos.

\smallskip

La seconde r\'esulte des remarques (ii) et (iii) qui suivent la d\'efinition \ref{defIV12}.

\medskip

\item[(iv)] La premi\`ere assertion r\'esulte de ce que, pour deux morphismes de topos $(f^* , f_*)$ et $(g^* , g_*)$ de ${\mathcal E}'$ vers ${\mathcal E}$, se donner une transformation naturelle $f^* \to g^*$ \'equivaut \`a se donner  une transformation naturelle $g_* \to f_*$.

\smallskip

La seconde assertion r\'esulte de la remarque (ii) qui suit la d\'efinition \ref{defIV12}.

\medskip

\item[(ii)] Il existe dans le topos ${\mathcal E}$ une petite sous-cat\'egorie pleine ${\mathcal C}$ dont les objets forment une famille s\'eparante de ${\mathcal E}$.

\smallskip

D'apr\`es le th\'eor\`eme \ref{thmIII81} de Giraud, le topos ${\mathcal E}$ est \'equivalent \`a la cat\'egorie $\widehat{\mathcal C}_J$ des faisceaux sur la cat\'egorie ${\mathcal C}$ munie de la topologie $J$ induite par ${\mathcal E}$ et tout objet $E$ de ${\mathcal E}$ vu comme un faisceau s'\'ecrit comme la colimite d'objets de ${\mathcal C}$
$$
E = \varinjlim_{(X,x) \in {\mathcal C}/E} X \, .
$$
Il en r\'esulte que tout foncteur respectant les colimites
$$
f^* : {\mathcal E} \longrightarrow {\mathcal E}'
$$
est d\'etermin\'e \`a unique isomorphisme pr\`es par sa restriction \`a ${\mathcal C}$ et qu'une transformation naturelle entre deux tels foncteurs $f^*$ et $g^*$ respectant les colimites
$$
f^* \longrightarrow g^*
$$
est d\'etermin\'ee par sa restriction \`a ${\mathcal C}$.

\smallskip

R\'eciproquement, toute transformation naturelle entre les restrictions de $f^*$ et $g^*$ \`a ${\mathcal C}$ se prolonge en une transformation naturelle $f^* \to g^*$.

\smallskip

La conclusion r\'esulte de ce que la cat\'egorie $[{\mathcal C} , {\mathcal E}']$ est localement petite puisque ${\mathcal C}$ est une petite cat\'egorie et que ${\mathcal E}'$ est localement petite. 
\end{listeisansmarge}
\end{demo}

\medskip

On d\'eduit de cette proposition:

\begin{cor}\label{corIV17}
\begin{listeimarge}
\item Pour tout topos ${\mathcal E}'$ et tout morphisme de topos
$$
g = (g^* , g_*) : {\mathcal E}_1 \longrightarrow {\mathcal E}_2 \, ,
$$
la composition avec $g$ d\'efinit un foncteur
$$
[{\mathcal E}' , {\mathcal E}_1]_T \longrightarrow[{\mathcal E}' , {\mathcal E}_2]_T \, .
$$
De plus, si
$$
h = (h^* , h_*) : {\mathcal E}_2 \longrightarrow {\mathcal E}_3
$$
est un autre morphisme de topos, le foncteur de composition avec $h \circ g$
$$
[{\mathcal E}' , {\mathcal E}_1]_T \longrightarrow [{\mathcal E}' , {\mathcal E}_3]_T
$$
est le compos\'e du foncteur de composition avec $g$ et de celui de composition avec $h$
$$
[{\mathcal E}' , {\mathcal E}_1]_T \longrightarrow [{\mathcal E}' , {\mathcal E}_2] \longrightarrow [{\mathcal E}' , {\mathcal E}_3]_T \, .
$$

\item Pour tout topos ${\mathcal E}$ et tout morphisme de topos
$$
g = (g^* , g_*) : {\mathcal E}'_2 \longrightarrow {\mathcal E}'_1 \, ,
$$
la composition avec $g$ d\'efinit un foncteur
$$
[{\mathcal E}'_1 , {\mathcal E}]_T \longrightarrow [{\mathcal E}'_2 , {\mathcal E}]_T \, .
$$
De plus, si
$$
h = (h^* , h_*) : {\mathcal E}'_3 \longrightarrow {\mathcal E}'_2
$$
est un autre morphisme de topos, le foncteur de composition avec $g \circ h$
$$
[{\mathcal E}'_1 , {\mathcal E}]_T \longrightarrow [{\mathcal E}'_3 , {\mathcal E}]_T
$$
est le compos\'e du foncteur de composition avec $h$ et de celui de composition avec $g$
$$
[{\mathcal E}'_1 , {\mathcal E}]_T \longrightarrow [{\mathcal E}'_2,{\mathcal E}]_T \longrightarrow [{\mathcal E}'_3,{\mathcal E}]_T \, .
$$

\item Pour tous morphismes de topos
$$
\begin{matrix}
g : {\mathcal E}_1 &\longrightarrow &{\mathcal E}_2 \, , \\
h : {\mathcal E}'_2 &\longrightarrow &{\mathcal E}'_1 \hfill
\end{matrix}
$$
le carr\'e de foncteurs de composition \`a gauche avec $g$ et \`a droite avec $h$
$$
\xymatrix{
[{\mathcal E}'_1 , {\mathcal E}_1]_T \ar[d] \ar[r] &[{\mathcal E}'_1 , {\mathcal E}_2]_T \ar[d] \\
[{\mathcal E}'_2 , {\mathcal E}_1]_T \ar[r] &[{\mathcal E}'_2 , {\mathcal E}_2]_T
}
$$
est commutatif.

\end{listeimarge}
\end{cor}

\begin{remark}

On r\'esume ces propri\'et\'es en disant que les topos forment une 2-cat\'egorie.
\end{remark}

\begin{demo}
\begin{listeisansmarge}
\item La loi de composition
$$
\left( {\mathcal E}' \xrightarrow{ \ f \ } {\mathcal E}_1 \right) \longmapsto \left( {\mathcal E}' \xrightarrow{ \ g \circ f \ } {\mathcal E}_2 \right)
$$
se compl\`ete en un foncteur
$$
[{\mathcal E}' , {\mathcal E}_1]_T \longrightarrow [{\mathcal E}' , {\mathcal E}_2]_T
$$
qui associe \`a toute transformation de morphismes de topos
$$
f = (f^* , f_*) \longrightarrow i=(i^*,i_*)
$$
vue comme une transformation naturelle
$$
\alpha : \left( {\mathcal E}' \xrightarrow{ \ i_* \ } {\mathcal E}_1 \right) \longrightarrow \left( {\mathcal E}' \xrightarrow{ \ f_* \ } {\mathcal E}_1 \right)
$$
la transformation naturelle
$$
g_* (\alpha) : \left( {\mathcal E}' \xrightarrow{ \ g_* \circ i_* \ } {\mathcal E}_2 \right) \longrightarrow \left( {\mathcal E}' \xrightarrow{ \ g_* \circ f_* \ } {\mathcal E}_2 \right)
$$
qui est la famille des morphismes de ${\mathcal E}'$
$$
g_* \circ i_* (E') \xrightarrow{ \ g_* (\alpha_{E'}) \ } g_* \circ f_* (E')
$$
images par le foncteur $g_*$ des morphismes
$$
i_* (E') \xrightarrow{ \ \alpha_{E'} \ } f_* (E')
$$
index\'es par les objets $E'$ de ${\mathcal E}'$ qui constituent $\alpha$.

\smallskip

Cette d\'efinition \'etant pos\'ee, la compatibilit\'e avec la composition des morphismes de topos
$$
{\mathcal E}_1 \xrightarrow{ \ g \ } {\mathcal E}_2 \xrightarrow{ \ h \ } {\mathcal E}_3
$$
est \'evidente.

\medskip

\item La loi de composition
$$
\left( {\mathcal E}'_1 \xrightarrow{ \ f \ } {\mathcal E} \right) \longmapsto \left( {\mathcal E}'_2 \xrightarrow{ \ f \circ g \ } {\mathcal E} \right)
$$
se compl\`ete en un foncteur
$$
[{\mathcal E}'_1 , {\mathcal E}]_T \longrightarrow [{\mathcal E}'_2 , {\mathcal E}]_T
$$
qui associe \`a toute transformation de morphismes de topos
$$
f = (f^* , f_*) \longrightarrow i = (i^* , i_*)
$$
vue comme une transformation naturelle
$$
\beta : f^* \longrightarrow i^*
$$
la transformation naturelle
$$
g^* (\beta) : \left( {\mathcal E} \xrightarrow{ \ g^* \circ f^* \ } {\mathcal E}'_2 \right) \longrightarrow \left( {\mathcal E} \xrightarrow{ \ g^* \circ i^* \ } {\mathcal E}'_2 \right)
$$
qui est la famille des morphismes de ${\mathcal E}'_2$
$$
g^* \circ f^* (E) \xrightarrow{ \ g^* (\beta_E) \ } g^* \circ i^* (E)
$$
images par le foncteur $g^*$ des morphismes
$$
f^* (E) \xrightarrow{ \ \beta_E \ } i^* (E)
$$
index\'es par les objets $E$ de ${\mathcal E}$ qui constituent $\beta$.

\smallskip

Cette d\'efinition \'etant pos\'ee, la compatibilit\'e avec la composition des morphismes de topos
$$
{\mathcal E}'_3 \xrightarrow{ \ h \ } {\mathcal E}'_2 \xrightarrow{ \ g \ } {\mathcal E}'_1
$$
est \'evidente.

\medskip

\item Ce carr\'e est commutatif au niveau des objets $f = (f^* , f_*)$ de $[{\mathcal E}'_1 , {\mathcal E}_1]_T$ car on a
$$
(g_* \circ f_*) \circ h_* = g_* \circ (f_* \circ h_*)
$$
et
$$
h^* \circ (f^* \!\circ g^*) = (h^* \circ f^*) \circ g^* \, .
$$

On remarque d'autre part que pour toute transformation de morphismes de topos ${\mathcal E}'_1 \to {\mathcal E}_1$
$$
f = (f^* , f_*) \longrightarrow (i^* , i_*)
$$
consistant en une transformation naturelle
$$
\beta : f^* \longrightarrow i^*
$$
ou, ce qui revient au m\^eme, la transformation adjointe
$$
\alpha : i_* \longrightarrow f_* \, ,
$$
alors la transformation naturelle
$$
h^* (\beta) : h^* \circ f^* \longrightarrow h^* \circ i^*
$$
admet pour adjointe la transformation naturelle
$$
i_* \circ h_* \longrightarrow f_* \circ h_*
$$
qui consiste en la famille des morphismes
$$
i_* \circ h_* (E) \xrightarrow{ \ \alpha_{h_* (E)} \ } f_* \circ h_* (E)
$$
index\'es par les objets $E$ de ${\mathcal E}'_2$.

\smallskip

La commutativit\'e du carr\'e au niveau des morphismes s'en d\'eduit. 
\end{listeisansmarge}
\end{demo}

\subsection{La cat\'egorie des points d'un topos}\label{subsec413}

\medskip

D'apr\`es la proposition \ref{propIV11}, toute application continue entre deux espaces topologiques
$$
f : X \longrightarrow Y
$$
induit un morphisme entre les topos associ\'es
$$
(f^* , f_*) : {\mathcal E}_X \longrightarrow {\mathcal E}_Y \, .
$$

En particulier, tout \'el\'ement $x$ d'un espace topologique $X$, vu comme une application continue
$$
x : \{ \bullet \} \longrightarrow X
$$
de l'ensemble \`a un \'el\'ement $\{\bullet\}$ vers $X$, induit un morphisme de topos
$$
{\rm Ens} \longrightarrow {\mathcal E}_X \, .
$$

Cette observation conduit \`a poser la d\'efinition g\'en\'erale suivante:

\begin{defn}\label{defIV18}

Soit ${\mathcal E}$ un topos.

\begin{listeimarge}

\item On appelle points de ${\mathcal E}$ les morphismes de topos
$$
{\rm Ens} \longrightarrow {\mathcal E}
$$
qui vont du topos {\rm Ens} des ensembles vers ${\mathcal E}$.

\medskip

\item On appelle cat\'egorie des points de ${\mathcal E}$ et on note ${\rm pt} ({\mathcal E})$ la cat\'egorie
$$
{\rm pt} ({\mathcal E}) = [{\rm Ens} , {\mathcal E}]_T
$$
des morphismes de topos ${\rm Ens} \to {\mathcal E}$.
\end{listeimarge}
\end{defn}

\begin{remarksqed}
\begin{listeisansmarge}
\item Pour tout point d'un topos ${\mathcal E}$
$$
x = (x^* , x_*) : {\rm Ens} \longrightarrow {\mathcal E} \, ,
$$
sa composante d'image r\'eciproque
$$
x^* : {\mathcal E} \longrightarrow {\rm Ens}
$$ 
est appel\'ee le foncteur fibre du point $x$.

\smallskip

En particulier, si $x$ est un \'el\'ement d'un espace topologique $X$, le foncteur
$$
x^* : {\mathcal E}_X \longrightarrow {\rm Ens}
$$
est appel\'e le foncteur fibre en le point $x$ de $X$. Il transforme tout faisceau sur $X$ en un ensemble appel\'e la fibre de ce faisceau en le point $x$.

\medskip

\item Un foncteur d\'efini sur un topos ${\mathcal E}$
$$
{\mathcal E} \longrightarrow {\rm Ens}
$$
est le foncteur fibre d'un point de ${\mathcal E}$ si et seulement si il respecte les colimites arbitraires et les limites finies.
\end{listeisansmarge}
\end{remarksqed}

\newpage

On d\'eduit du corollaire \ref{corIV17} (i):

\begin{cor}\label{defIV19}
\begin{listeimarge}
\item Tout morphisme de topos
$$
f = (f^* , f_*) : {\mathcal E}_1 \longrightarrow {\mathcal E}_2
$$
d\'efinit par composition avec $f$ un foncteur
$$
{\rm pt} ({\mathcal E}_1) \longrightarrow {\rm pt} ({\mathcal E}_2) \, .
$$

\item De plus, pour tous morphismes de topos
$$
{\mathcal E}_1 \xrightarrow{ \ f \ } {\mathcal E}_2 \xrightarrow{ \ g \ } {\mathcal E}_3 \, ,
$$
le foncteur associ\'e \`a $g \circ f$
$$
{\rm pt} ({\mathcal E}_1) \longrightarrow {\rm pt} ({\mathcal E}_3)
$$
est le compos\'e des foncteurs associ\'es \`a $f$ et $g$
$$
{\rm pt} ({\mathcal E}_1) \longrightarrow {\rm pt} ({\mathcal E}_2) \longrightarrow {\rm pt} ({\mathcal E}_3) \, .
$$
\end{listeimarge}
\end{cor}

\begin{demo}

C'est le cas particulier du corollaire \ref{corIV17} (i) o\`u l'on prend pour ${\mathcal E}'$ le topos ${\rm Ens}$ des ensembles. 

\end{demo}

\bigskip

La d\'efinition de la notion de point d'un topos est justifi\'ee mieux encore par le th\'eor\`eme suivant et son corollaire:

\begin{thm}\label{thmIV110}

Soit $X$ un espace topologique.

\begin{listeimarge}

\item Soit $Y$ un sous-espace ferm\'e de $X$ qui est ``irr\'eductible'' au sens que l'intersection de deux ouverts non vides de $Y$ est toujours non vide.

\smallskip

Alors le foncteur d\'efini sur la cat\'egorie ${\mathcal E}_X$ des faisceaux sur $X$
$$
\begin{matrix}
p_Y^* : {\mathcal E}_X &\longrightarrow &{\rm Ens} \, , \hfill \\
\hfill F &\longmapsto &\displaystyle \varinjlim_{U \in O(X) \atop U \cap Y \ne \emptyset} F(U)
\end{matrix}
$$
est la composante d'image r\'eciproque d'un point
$$
p_Y : {\rm Ens} \longrightarrow {\mathcal E}_X
$$
du topos ${\mathcal E}_X$.

\medskip

\item R\'eciproquement, tout point du topos ${\mathcal E}_X$ est isomorphe dans la cat\'egorie ${\rm pt} ({\mathcal E}_X)$ \`a un point de la forme
$$
p_Y : {\rm Ens} \longrightarrow {\mathcal E}_X
$$
associ\'e \`a un ferm\'e irr\'eductible $Y$ de $X$.

\medskip

\item Deux points de ${\mathcal E}_X$ associ\'es \`a deux ferm\'es irr\'eductibles $Y_1$ et $Y_2$ sont reli\'es par un morphisme
$$
p_{Y_1} \longrightarrow p_{Y_2}
$$
si et seulement si
$$
Y_1 \subset Y_2 \, .
$$

Dans ce cas, ils sont reli\'es par l'unique morphisme
$$
p_{Y_1} \longrightarrow p_{Y_2}
$$
qui correspond \`a la transformation naturelle
$$
p_{Y_1}^* \longrightarrow p_{Y_2}^*
$$
d\'efinie par l'inclusion
$$
\{ U \in O(X) \mid U \cap Y_1 \ne 0 \} \subset \{ U \in O(X) \mid U \cap Y_2 \ne 0 \} \, .
$$
\end{listeimarge}
\end{thm}

\begin{remark}

Autrement dit, pour tout espace topologique $X$, la cat\'egorie ${\rm pt} ({\mathcal E}_X)$ des points du topos associ\'e ${\mathcal E}_X$ est \'equivalence \`a l'ensemble ordonn\'e des ferm\'es irr\'eductibles de $X$.

\smallskip

En particulier, elle est essentiellement petite.

\end{remark}

\begin{demo}
\begin{listeisansmarge}
\item Associons \`a tout ensemble $E$ le pr\'efaisceau sur $X$ d\'efini par la formule
$$
\begin{matrix}
P_E : O(X)^{\rm op} &\longrightarrow &{\rm Ens} \, , \hfill \\
\hfill U &\longmapsto &\left\{ \begin{matrix}
E &\mbox{si} &U \cap Y \ne \emptyset \hfill \\
\{ \bullet \} &\mbox{si} &U \cap Y = \emptyset \, .
\end{matrix} \right.
\end{matrix}
$$
L'hypoth\`ese d'irr\'eductibilit\'e de $Y$ entra{\^\i}ne que $P_E$ est un faisceau.

\smallskip

En effet, pour tout recouvrement ouvert
$$
U = \bigcup_{i \in I} U_i
$$
d'un ouvert $U$ de $X$, on a
$$
U \cap Y = \emptyset \quad \mbox{si et seulement si} \quad U_i \cap Y = \emptyset \, , \quad \forall \, i \, .
$$
Dans le cas contraire, on a pour tous indices $i,j$
$$
U_i \cap U_j \cap Y \ne \emptyset \quad \mbox{si et seulement si} \quad U_i \cap Y \ne 0 \quad \mbox{et} \quad U_j \cap Y \ne \emptyset \, .
$$
On en d\'eduit que l'application canonique
$$
P_E (U) \longrightarrow {\rm eg} \left( \prod_i P_E (U_i) \rightrightarrows \prod_{i,j} P_E (U_i \cap U_j) \right)
$$
est bijective.

\smallskip

Ainsi, on a d\'efini un foncteur
$$
\begin{matrix}
{\rm Ens} &\longrightarrow &{\mathcal E}_X \, , \\
\hfill E &\longmapsto &P_E \, .
\end{matrix}
$$

On observe qu'il est adjoint \`a droite du foncteur
$$
\begin{matrix}
p_Y^* : {\mathcal E}_X &\longrightarrow &{\rm Ens} \, , \hfill \\
\hfill F &\longmapsto &\displaystyle \varinjlim_{U \in O(X) \atop U \cap Y = \emptyset} F(U) \, .
\end{matrix}
$$

En effet, pour tout ensemble $E$ et tout faisceau $F$ sur $X$, se donner un morphisme de faisceaux
$$
F \longrightarrow P_E
$$
\'equivaut \`a se donner une famille compatible d'applications
$$
F(U) \longrightarrow E
$$
index\'ees par les ouverts $U$ de $X$ tels que $U \cap Y \ne \emptyset$.

\smallskip

Enfin, le foncteur $p_Y^* : {\mathcal E}_X \to {\rm Ens}$ respecte les limites finies car, le ferm\'e $Y$ \'etant irr\'eductible, l'ensemble ordonn\'e
$$
\{ U \in O(X) \mid U \cap Y \ne 0 \}
$$
est filtrant.

\medskip

\item Commen\c cons par observer que dans la cat\'egorie $O(X)$ des ouverts de $X$, l'objet $X$ est terminal et tout morphisme $U \to X$ est un monomorphisme.

\smallskip

Notant
$$
\ell : O(X) \xrightarrow{ \ y \ } \widehat{O(X)} \xrightarrow{ \ j^* \ } {\mathcal E}_X
$$
le foncteur canonique associ\'e, on en d\'eduit que $\ell (X)$ est l'objet terminal de ${\mathcal E}_X$ et que les morphismes induits
$$
\ell (U) \longrightarrow \ell (X)
$$
sont des monomorphismes.

\smallskip

Si
$$
(x^* , x_*) : {\rm Ens} \longrightarrow {\mathcal E}_X
$$
est un point de ${\mathcal E}_X$,
$$
x^* (\ell (X))
$$
est n\'ecessairement l'objet terminal $\{\bullet\}$ de ${\rm Ens}$ et, pour tout ouvert $U$,
$$
x^* (\ell (U))
$$
est un sous-ensemble de $\{\bullet\}$, n\'ecessairement \'egal \`a $\emptyset$ ou \`a tout $\{\bullet\}$.

\smallskip

Si $x^* (\ell (U)) = \emptyset$, on a n\'ecessairement
$$
x^* (\ell (U')) = \emptyset
$$
pour tout ouvert $U'$ contenu dans $U$.

\smallskip

Comme $x^*$ respecte les colimites, il existe un plus grand ouvert $U_{\emptyset}$ tel que
$$
x^* (\ell (U_{\emptyset})) = \emptyset \, .
$$

Le sous-espace ferm\'e $Y \hookrightarrow X$ compl\'ementaire de $U_{\emptyset}$ est n\'ecessairement irr\'eductible puisque, pour tous ouverts $U_1$ et $U_2$ qui rencontrent $Y$, le carr\'e cart\'esien
$$
\xymatrix{
x^* \circ \ell (U_1 \cap U_2) \ar[d] \ar[r] &x^* \circ \ell (U_1) = \{\bullet\} \ar[d] \\
\{\bullet\} = x^* \circ \ell (U_2) \ar[r] &x^* \circ \ell (U_1 \cup U_2) = \{\bullet\}
}
$$
impose que
$$
x^* \circ \ell (U_1 \cap U_2) = \{\bullet\} \, ,
$$
autrement dit que $U_1 \cap U_2$ rencontre $Y$.

\smallskip

On observe que les deux foncteurs
$$
x^* : {\mathcal E}_X \longrightarrow {\rm Ens}
$$
et
$$
p_Y^* : {\mathcal E}_X \longrightarrow {\rm Ens}
$$
ont le m\^eme compos\'e
$$
x^* \circ \ell = p_Y^* \circ \ell
$$
avec le foncteur canonique $\ell : O(X) \to {\mathcal E}_X$.

\smallskip

Comme tous deux respectent les colimites, ils sont canoniquement isomorphes.

\medskip

\item Pour tous ferm\'es irr\'eductibles $Y_1 , Y_2$ de $X$, le foncteur
$$
p_{Y_1}^* : {\mathcal E}_X \longrightarrow {\rm Ens} \qquad \mbox{[resp.} \quad p_{Y_2}^* : {\mathcal E}_X \longrightarrow {\rm Ens} \, \mbox{]}
$$
associe \`a tout objet de la forme $\ell (U)$, $U \in O(X)$, l'ensemble

\medskip

\hglue 1cm $\{\bullet\}$ si $U$ rencontre $Y_1$ [resp. $Y_2$],

\medskip

\hglue 1cm $\emptyset$ dans le cas contraire.

\medskip

Donc il ne peut exister de transformation naturelle
$$
p_{Y_1}^* \longrightarrow p_{Y_2}^*
$$
que si $Y_1 \subset Y_2$ et, dans ce cas, elle est uniquement d\'etermin\'ee et consiste \`a associer \`a tout faisceau $F$ sur $X$ l'application canonique
$$
\varinjlim_{U \cap Y_1 \ne \emptyset} F(U) \longrightarrow \varinjlim_{U \cap Y_2 \ne \emptyset} F(U) \, .
$$
Cela termine la d\'emonstration du th\'eor\`eme. 
\end{listeisansmarge}
\end{demo}

\bigskip

Revenant aux points des espaces topologiques, on d\'eduit de ce th\'eor\`eme:

\begin{cor}\label{corIV111}

Soit $X$ un espace topologique.

\begin{listeimarge}

\item Pour tout \'el\'ement $x \in X$, son adh\'erence $\overline{\{ x \}} = Y_x$ est un ferm\'e irr\'eductible de $X$ et le foncteur fibre
$$
x^* : {\mathcal E}_X \longrightarrow {\rm Ens}
$$
s'identifie au foncteur
$$
\begin{matrix}
p_{Y_x}^* : {\mathcal E}_X &\longrightarrow &{\rm Ens} \, , \hfill \\
\hfill F &\longmapsto &\displaystyle \varinjlim_{U \in O(X) \atop x \in U} F(U) \, .
\end{matrix}
$$

\item L'application
$$
x \longmapsto p_{Y_x}^*
$$
de $X$ vers l'ensemble des classes d'isomorphie de points de $X$ est une bijection si et seulement si $X$ est un espace topologique ``sobre'' au sens que tout ferm\'e irr\'eductible de $X$ est l'adh\'erence d'un unique \'el\'ement.

\medskip

\item Si $X$ est un espace topologique sobre, la cat\'egorie
$$
{\rm pt} ({\mathcal E}_X)
$$
est \'equivalente \`a l'ensemble
$$
X
$$
ordonn\'e par la relation
$$
x_1 \leq x_2 \quad \mbox{si et seulement si} \quad \overline{\{x_1\}} \subset \overline{\{x_2\}} \, .
$$
\end{listeimarge}
\end{cor}

\begin{remarks}
\begin{listeisansmarge}
\item Si $X$ est un espace topologique s\'epar\'e, ses ferm\'es irr\'eductibles sont ses sous-ensembles \`a un \'el\'ement $\{x\}$.

\smallskip

Donc un espace s\'epar\'e est toujours sobre.

\medskip

\item Un espace $X$ \'ecrit comme une r\'eunion d'ouverts
$$
X = \bigcup_i U_i
$$
est sobre si et seulement si chaque $U_i$ est sobre.

\smallskip

Autrement dit, la propri\'et\'e de sobri\'et\'e est locale.

\medskip

\item L'espace sous-jacent $X$ d'un sch\'ema $(X , {\mathcal O}_X)$ est toujours sobre.

\smallskip

D'apr\`es (ii), il suffit en effet de le v\'erifier dans le cas du sch\'ema affine ${\rm Spec} (A)$ associ\'e \`a un anneau commutatif $A$.

\smallskip

Si $Y$ est un ferm\'e irr\'eductible, soit $I$ l'intersection de tous les id\'eaux premiers de $A$ qui sont \'el\'ements de~$Y$.

\smallskip

Par d\'efinition de l'irr\'eductibilit\'e de $Y$, on a pour tous \'el\'ements $a,b$ de $A$ l'implication
$$
a \notin I \wedge b \notin I \Rightarrow ab \notin I \, .
$$

Donc $I$ est un id\'eal premier et $Y$ est l'adh\'erence du point $I$. Il y a unicit\'e car pour tous id\'eaux premiers $p,q$ de $A$ on a
$$
p \in \overline{\{q\}} \qquad \mbox{dans} \quad {\rm Spec} (A)
$$
si et seulement si
$$
q \subseteq p \, .
$$
\end{listeisansmarge}
\end{remarks}

\begin{demo}
\begin{listeisansmarge}
\item L'adh\'erence $\overline{\{x\}} = Y_x$ d'un point $x$ de $X$ est irr\'eductible puisque, pour tout ouvert $U$ de $X$, on a
$$
U \cap \overline{\{x\}} \ne \emptyset
$$
si et seulement si $x \in U$.

\smallskip

L'identit\'e des deux foncteurs $x^*$ et $p_{Y_x}^*$ vient de ce que, d'apr\`es la proposition \ref{propIV11}, le foncteur fibre $x^*$ est d\'efini par la formule
$$
F \longmapsto \varinjlim_{U \in O(X) \atop x \in U} F(U) \, .
$$

\item est cons\'equence imm\'ediate de (i) et du th\'eor\`eme \ref{thmIV110}.

\medskip

\item r\'esulte de (ii) et du th\'eor\`eme \ref{thmIV110} (iii). 
\end{listeisansmarge}
\end{demo}

\subsection{Espaces topologiques, applications continues, topos et morphismes de topos}\label{subsec414}

\medskip

La d\'efinition g\'en\'erale des morphismes de topos est encore davantage justifi\'ee par le th\'eor\`eme suivant qui compl\`ete le corollaire \ref{corIV111}:

\begin{thm}\label{thmIV112}

Soit $X$ un espace topologique.

\begin{listeimarge}

\item L'ensemble $O(X)$ des ouverts de $X$, avec sa relation d'ordre, ses r\'eunions arbitraires et ses intersections finies s'identifie \`a l'ensemble des sous-objets
$$
S \xhookrightarrow{ \ { \ } \ } 1_X
$$
de l'objet terminal $1_X$ du topos ${\mathcal E}_X$ avec sa relation d'ordre, ses r\'eunions arbitraires et ses intersections finies.

\smallskip

Dans cette correspondance, un \'el\'ement $x \in X$ est dans l'ouvert $U_S$ associ\'e \`a un sous-objet $S \hookrightarrow 1_X$ si et seulement si l'image $x^* S$ par le foncteur fibre
$$
x^* : {\mathcal E}_X \longrightarrow {\rm Ens}
$$
est le sous-objet $\{\bullet\}$ de l'objet terminal $\{\bullet\} = x^* 1_X$ de ${\rm Ens}$.

\medskip

\item Si $Y$ est un espace topologique sobre, se donner une application continue
$$
f : X \to Y
$$
\'equivaut \`a se donner, \`a unique isomorphisme pr\`es, un morphisme de topos
$$
(f^* , f_*) : {\mathcal E}_X \longrightarrow {\mathcal E}_Y \, .
$$
\end{listeimarge}
\end{thm}

\begin{demo}
\begin{listeisansmarge}
\item Dans la cat\'egorie $O(X)$, $X$ est un objet terminal et les morphismes $U \to X$ sont des monomorphismes.

\smallskip

Le foncteur canonique
$$
\ell : O(X) \xrightarrow{ \ y \ } \widehat{O(X)} \xrightarrow{ \ j^* \ } {\mathcal E}_X
$$
respecte les limites finies, donc il transforme $X$ en un objet terminal
$$
\ell (X) = 1_X \qquad \mbox{de} \qquad {\mathcal E}_X
$$
et les ouverts $U$ en des sous-objets
$$
\ell (U) \xhookrightarrow{ \ { \ } \ } 1_X \, .
$$

Il respecte la relation d'ordre et les intersections finies.

\smallskip

Si un ouvert $U$ est r\'eunion d'ouverts $U_i$, $i \in I$, la famille des
$$
U_i \longrightarrow U
$$
est couvrante et donc la famille des morphismes ${\mathcal E}_X$
$$
\ell (U_i) \longrightarrow \ell (U)
$$
est globalement \'epimorphique, ce qui signifie que le sous-objet $\ell (U)$ de $1_X$ est r\'eunion des sous-objets $\ell (U_i)$.

\smallskip

Tout sous-objet $S$ de $1_X$ vu comme un faisceau s'\'ecrit
$$
S = \varinjlim_{(U,u) \in \int\!\!S} \ell (U) \, .
$$
Si $U_S$ est l'ouvert de $X$ r\'eunion des ouverts $U$ tels que
$$
\ell (U) \xhookrightarrow{ \ { \ } \ } 1_X
$$
se factorise \`a travers $S \hookrightarrow 1_X$, on a donc
$$
S = \ell (U_S) \, .
$$

Pour tout \'el\'ement $x$ de $X$, $x^* 1_X$ est n\'ecessairement l'objet terminal $\{\bullet\}$ de ${\rm Ens}$ et $x$ est contenu dans un ouvert $U$ si et seulement si
$$
x^* \ell (U) = \{\bullet\} \, .
$$

Comme le foncteur $x^*$ respecte les colimites, un tel \'el\'ement $x$ est contenu dans l'ouvert $U_S$ associ\'e \`a un sous-objet $S \hookrightarrow 1_X$ si et seulement si
$$
x^* S = \{ \bullet \} \, .
$$

\item Toute application continue
$$
f : X \to Y
$$
d\'efinit un morphisme de topos
$$
(f^* , f_*) : {\mathcal E}_X \longrightarrow {\mathcal E}_Y
$$
qui respecte la formation des points.

\smallskip

R\'eciproquement, un morphisme de topos
$$
(f^* , f_*) : {\mathcal E}_X \longrightarrow {\mathcal E}_Y
$$
d\'efinit un foncteur
$$
{\rm pt} ({\mathcal E}_X) \longrightarrow {\rm pt} ({\mathcal E}_Y)
$$
et donc une application
$$
f : X \longrightarrow Y
$$
puisque tout \'el\'ement de $X$ d\'efinit un point de ${\mathcal E}_X$ et que, $Y$ \'etant sobre, les points de ${\mathcal E}_Y$ s'identifient aux \'el\'ements de $Y$.

\smallskip

De plus, le foncteur
$$
f^* : {\mathcal E}_Y \longrightarrow {\mathcal E}_X
$$
respecte les colimites arbitraires et les limites finies. Il transforme donc l'objet terminal $1_Y$ de $Y$ en l'objet terminal $1_X$ de $X$, et les sous-objets de $1_Y$ en sous-objets de $1_X$, en respectant la relation d'ordre, les r\'eunions arbitraires et les intersections finies.

\smallskip

D'apr\`es (i), cela signifie que $f^*$ d\'efinit une application
$$
O(Y) \longrightarrow O(X) \, .
$$
Celle-ci se confond avec l'application
$$
f^{-1} : O(Y) \longrightarrow O(X)
$$
induite par
$$
f : X \longrightarrow Y
$$
car pour tout \'el\'ement $x$ de $X$ et tout sous-objet $S \hookrightarrow 1_Y$, on a les \'equivalences
$$
\begin{matrix}
f(x) \in U_S &\Longleftrightarrow &f(x)^* S = \{\bullet\} \hfill \\
&\Longleftrightarrow &x^* \circ f^* S = \{\bullet\} \hfill \\
&\Longleftrightarrow &x \in U_{f^* S} \, . \hfill
\end{matrix}
$$
Cela signifie que $f$ est une application continue
$$
X \longrightarrow Y
$$
et que $(f^* , f_*)$ est canoniquement isomorphe au morphisme de topos d\'efini par $f$.

\smallskip

En effet, pour tout ouvert $U_S$ de $Y$ correspondant \`a un sous-objet $S \hookrightarrow 1_Y$, le foncteur $f^*$ transforme $\ell (U_S) = S$ en le sous-objet $\ell (f^{-1} (U_S)) = f^* S$.

\smallskip

Cela termine la d\'emonstration du th\'eor\`eme.
\end{listeisansmarge}
\end{demo}

\section{Morphismes de topos et structures alg\'ebriques}\label{sec42}

\subsection{Images directes et images r\'eciproques de structures alg\'ebriques}\label{subsec421}

\medskip

On observe:

\begin{prop}\label{propIV21}

Soit ${\mathbb T}$ un type de structure alg\'ebrique tel que celui de groupe, de mono{\"\i}de, de groupe commutatif, de mono{\"\i}de commutatif, d'anneau ou d'anneau commutatif.

\smallskip

Pour tout topos ${\mathcal E}$, consid\'erons la cat\'egorie ${\mathcal M}_{\mathbb T} ({\mathcal E})$ des objets de ${\mathcal E}$ munis d'une structure de type ${\mathbb T}$ c'est-\`a-dire des groupes internes, mono{\"\i}des internes, anneaux internes [resp. commutatifs] de ${\mathcal E}$.

\smallskip

Alors, pour tout morphisme de topos
$$
f = (f^* , f_*) : {\mathcal E}' \longrightarrow {\mathcal E} \, ,
$$
on a:

\begin{listeimarge}

\item Le foncteur d'image directe $f_* : {\mathcal E}' \to {\mathcal E}$ induit un foncteur
$$
f_* : {\mathcal M}_{\mathbb T} ({\mathcal E}') \longrightarrow {\mathcal M}_{\mathbb T} ({\mathcal E})
$$
qui rend commutatif le carr\'e
$$
\xymatrix{
{\mathcal M}_{\mathbb T} ({\mathcal E}') \ar[d] \ar[r]^-{f_*} &{\mathcal M}_{\mathbb T} ({\mathcal E}) \ar[d] \\
{\mathcal E}' \ar[r]^-{f_*} &{\mathcal E}
}
$$
dont les fl\`eches verticales sont les foncteurs d'oubli de la structure alg\'ebrique de type ${\mathbb T}$.

\medskip

\item Le foncteur d'image r\'eciproque $f^* : {\mathcal E} \to {\mathcal E}'$ induit un foncteur
$$
f^* : {\mathcal M}_{\mathbb T} ({\mathcal E}) \longrightarrow {\mathcal M}_{\mathbb T} ({\mathcal E}')
$$
qui rend commutatif le carr\'e
$$
\xymatrix{
{\mathcal M}_{\mathbb T} ({\mathcal E}) \ar[d] \ar[r]^-{f^*} &{\mathcal M}_{\mathbb T} ({\mathcal E}') \ar[d] \\
{\mathcal E} \ar[r]^-{f^*} &{\mathcal E}'
}
$$
dont les fl\`eches verticales sont les foncteurs d'oubli.

\medskip

\item Le foncteur
$$
f^* : {\mathcal M}_{\mathbb T} ({\mathcal E}) \longrightarrow {\mathcal M}_{\mathbb T} ({\mathcal E}')
$$
est adjoint \`a gauche du foncteur
$$
f_* : {\mathcal M}_{\mathbb T} ({\mathcal E}') \longrightarrow {\mathcal M}_{\mathbb T} ({\mathcal E}) \, .
$$
\end{listeimarge}
\end{prop}

\begin{remarks}
\begin{listeisansmarge}
\item Pour tous morphismes de topos
$$
{\mathcal E}_1 \xrightarrow{ \ f = (f^* , f_*) \ } {\mathcal E}_2 \xrightarrow{ \ g = (g^* , g_*) \ } {\mathcal E}_3 \, ,
$$
les foncteurs induits
$$
(g \circ f)_* : {\mathcal M}_{\mathbb T} ({\mathcal E}_1) \longrightarrow {\mathcal M}_{\mathbb T} ({\mathcal E}_3)
$$
et
$$
(g \circ f)^* : {\mathcal M}_{\mathbb T} ({\mathcal E}_3) \longrightarrow {\mathcal M}_{\mathbb T} ({\mathcal E}_1)
$$
sont les foncteurs compos\'es
$$
{\mathcal M}_{\mathbb T} ({\mathcal E}_1) \xrightarrow{ \ f_* \ } {\mathcal M}_{\mathbb T} ({\mathcal E}_2) \xrightarrow{ \ g_* \ } {\mathcal M}_{\mathbb T} ({\mathcal E}_3)
$$
et
$$
{\mathcal M}_{\mathbb T} ({\mathcal E}_3) \xrightarrow{ \ g^* \ } {\mathcal M}_{\mathbb T} ({\mathcal E}_2) \xrightarrow{ \ f^* \ } {\mathcal M}_{\mathbb T} ({\mathcal E}_1) \, .
$$

\item Si ${\mathbb T}$ est un type de structure alg\'ebrique et ${\mathbb S}$ un autre type qui se d\'eduit de ${\mathbb T}$ en oubliant une partie de la structure (comme la structure de groupe par rapport \`a celle de groupe ab\'elien, ou la structure de groupe ab\'elien par rapport \`a celle d'anneau), les foncteurs d'oubli
$$
{\mathcal M}_{\mathbb T} ({\mathcal E}) \longrightarrow {\mathcal M}_{\mathbb S} ({\mathcal E})
$$
s'inscrivent dans des carr\'es commutatifs
$$
\xymatrix{
{\mathcal M}_{\mathbb T} ({\mathcal E}') \ar[d] \ar[r]^-{f_*} &{\mathcal M}_{\mathbb T} ({\mathcal E}) \ar[d] \\
{\mathcal M}_{\mathbb S} ({\mathcal E}') \ar[r]^-{f_*} &{\mathcal M}_{\mathbb S} ({\mathcal E})
}
$$
et 
$$
\xymatrix{
{\mathcal M}_{\mathbb T} ({\mathcal E}) \ar[d] \ar[r]^-{f^*} &{\mathcal M}_{\mathbb T} ({\mathcal E}') \ar[d] \\
{\mathcal M}_{\mathbb S} ({\mathcal E}) \ar[r]^-{f^*} &{\mathcal M}_{\mathbb S} ({\mathcal E}')
}
$$
pour tout morphisme de topos
$$
f = (f^* , f_*) : {\mathcal E}' \longrightarrow {\mathcal E} \, .
$$

\end{listeisansmarge}
\end{remarks}

\begin{demo}

Une structure de type ${\mathbb T}$ sur un objet $E$ d'un topos ${\mathcal E}$ consiste en une famille de morphismes reliant $E$, ses produits finis et l'objet terminal $1_{\mathcal E}$ de ${\mathcal E}$, soumis \`a des axiomes de commutativit\'e de diagrammes compos\'es de ces morphismes et de produits de ces morphismes.

\smallskip

Un morphisme entre deux objets $E_1$ et $E_2$ de la cat\'egorie ${\mathcal M}_{\mathbb T} ({\mathcal E})$ est un morphisme de ${\mathcal E}$
$$
e : E_1 \longrightarrow E_2
$$
qui rend commutatifs les diagrammes compos\'es des morphismes de d\'efinition des structures de type ${\mathbb T}$ sur $E_1$ et $E_2$ et du morphisme $e$ ou de puissances de $e$.

\begin{listeisansmarge}

\item[(i) et (ii)] r\'esultent donc de ce que les deux foncteurs
$$
f_* : {\mathcal E}' \longrightarrow {\mathcal E}
$$
et
$$
f^* : {\mathcal E} \longrightarrow {\mathcal E}'
$$
respectent les limites finies, en particulier les objets terminaux et les produits finis.

\medskip

\item[(iii)] r\'esulte alors de ce que le foncteur $f^* : {\mathcal E} \to {\mathcal E}'$ est adjoint \`a gauche du foncteur $f_* : {\mathcal E}' \to {\mathcal E}$.

\smallskip

En effet, \'etant donn\'es des objets $E$ de ${\mathcal M}_{\mathbb T} ({\mathcal E})$ et $E'$ de ${\mathcal M}_{\mathbb T} ({\mathcal E}')$ un morphisme de ${\mathcal E}$
$$
E \longrightarrow f_* E'
$$
respecte les structures de type ${\mathbb T}$ de $E$ et de $f_* (E')$ si et seulement si le morphisme de ${\mathcal E}'$ qui lui correspond
$$
f^* E \longrightarrow E'
$$
respecte les structures de type ${\mathbb T}$ de $f^* E$ et de $E'$. 

\end{listeisansmarge}
\end{demo}

\subsection{Morphismes de topos annel\'es}\label{subsec422}

\medskip

La proposition \ref{propIV21} permet de g\'en\'eraliser la notion de morphisme d'espaces annel\'es en celle de morphisme  de topos annel\'es:

\begin{defn}\label{defIV22}

On appelle morphisme de topos annel\'es
$$
({\mathcal E}' , A') \longrightarrow ({\mathcal E},A)
$$
la donn\'ee d'un morphisme de topos
$$
f = (f^* , f_*) : {\mathcal E}' \longrightarrow {\mathcal E}
$$
et d'un morphisme d'anneaux internes de ${\mathcal E}$
$$
A \longrightarrow f_* A' \, .
$$
\end{defn}

\begin{remarksqed}
\begin{listeisansmarge}
\item D'apr\`es la proposition \ref{propIV21} (iii), se donner un morphisme d'anneaux internes de ${\mathcal E}$
$$
A \longrightarrow f_* A'
$$
\'equivaut \`a se donner un morphisme d'anneaux internes de ${\mathcal E}'$
$$
f^* A \longrightarrow A' \, .
$$

\item Le compos\'e de deux morphismes de topos annel\'es
$$
({\mathcal E}_1 , A_1) \longrightarrow ({\mathcal E}_2 , A_2) \longrightarrow ({\mathcal E}_3 , A_3)
$$
est le compos\'e des morphismes de topos sous-jacents
$$
{\mathcal E}_1 \xrightarrow{ \ f = (f^* , f_*) \ } {\mathcal E}_2 \xrightarrow{ \ g = (g^* , g_*) \ } {\mathcal E}_3
$$
compl\'et\'e par le compos\'e de morphismes d'anneaux internes de ${\mathcal E}_3$
$$
A_3 \longrightarrow g_* A_2 \longrightarrow g_* f_* A_1 = (g \circ f)_* A_1 \, .
$$

Cette loi de composition est associative et d\'efinit la cat\'egorie des topos annel\'es.

\medskip

\item Pour tous topos annel\'es $({\mathcal E}',A')$ et $({\mathcal E},A)$ reli\'es par deux morphismes
$$
(f^* , f_* , A \longrightarrow f_* A')
$$
et
$$
(g^* , g_* , A \longrightarrow g_* A') \, ,
$$
une transformation du premier morphisme vers le second est une transformation naturelle
$$
f^* \longrightarrow g^*
$$
telle que le morphisme de ${\mathcal E}'$ induit
$$
f^* A \longrightarrow g^* A
$$
respecte les structures d'anneaux et que le triangle
$$
\xymatrix{
f^* A \ar[rd] \ar[rr] &&g^*A \ar[ld] \\
&A'
}
$$
soit commutatif dans la cat\'egorie des anneaux internes de ${\mathcal E}'$.

\smallskip

La composition des transformations naturelles d\'efinit la cat\'egorie des morphismes de topos annel\'es de $({\mathcal E}',A')$ dans $({\mathcal E},A)$.

\medskip

\item Les topos annel\'es, leurs morphismes et les transformations entre ces morphismes forment une 2-cat\'egorie. 

\end{listeisansmarge}
\end{remarksqed}

\subsection{Images directes et images r\'eciproques de modules internes}\label{subsec423}

\medskip

Les morphismes de topos annel\'es induisent des morphismes additifs entre les cat\'egories de modules internes associ\'ees:

\begin{prop}\label{propIV23}

Soient deux topos annel\'es $({\mathcal E},A)$ et $({\mathcal E}',A')$ reli\'es par un morphisme constitu\'e de
$$
f = (f^* , f_*) : {\mathcal E}' \longrightarrow {\mathcal E}
$$
et
$$
A \longrightarrow f_* A' \, .
$$

Alors:

\begin{listeimarge}

\item Le foncteur d'image directe $f_* : {\mathcal E}' \to {\mathcal E}$ induit un foncteur additif
$$
f_* : {\mathcal M}od_{A'} \longrightarrow {\mathcal M}od_A
$$
qui transforme la loi de multiplication par $A'$ des objets $E'$ de ${\mathcal M}od_{A'}$
$$
A' \times E' \longrightarrow E'
$$
en la loi de multiplication par $A$ de $f_* E'$
$$
A \times f_* E' \longrightarrow f_* A' \times f_* E' = f_* (A' \times E') \longrightarrow f_* E' \, .
$$
Ce foncteur rend commutatif le carr\'e
$$
\xymatrix{
{\mathcal M}od_{A'} \ar[d] \ar[r]^-{f_*} &{\mathcal M}od_A \ar[d] \\
{\mathcal A}b_{{\mathcal E}'} \ar[r]^-{f_*} &{\mathcal A}b_{\mathcal E}
}
$$
dont les fl\`eches verticales sont les foncteurs d'oubli vers les cat\'egories des groupes ab\'eliens internes de ${\mathcal E}'$ et~${\mathcal E}$.

\medskip

\item Si le morphisme
$$
f^* A \longrightarrow A'
$$
est un isomorphisme, le foncteur $f^* : {\mathcal E} \to {\mathcal E}'$ d\'efinit un adjoint \`a gauche de $f_* : {\mathcal M}od_{A'} \to {\mathcal M}od_A$
$$
f^* : {\mathcal M}od_A \longrightarrow {\mathcal M}od_{A'}
$$
qui rend commutatif le carr\'e
$$
\xymatrix{
{\mathcal M}od_A \ar[d] \ar[r]^-{f^*} &{\mathcal M}od_{A'} \ar[d] \\
{\mathcal A}b_{\mathcal E} \ar[r]^-{f^*} &{\mathcal A}b_{{\mathcal E}'}
}
$$
dont les fl\`eches verticales sont les foncteurs d'oubli.

\smallskip

Il transforme la loi de multiplication par $A$ des objets $E$ de ${\mathcal M}od_A$
$$
A \times E \longrightarrow E
$$
en la loi de multiplication par $f^* A = A'$ de $f^* E$
$$
f^* A \times f^* E = f^* (A \times E) \longrightarrow f^* E \, .
$$

\item Si ${\mathcal E}' = {\mathcal E}$ et $f = {\rm id}_{\mathcal E}$, le foncteur d'oubli
$$
{\mathcal M}od_{A'} \longrightarrow {\mathcal M}od_A
$$
admet pour adjoint \`a gauche le foncteur
$$
\begin{matrix}
{\mathcal M}od_A &\longrightarrow &{\mathcal M}od_{A'} \, , \\
\hfill E &\longmapsto &A' \otimes_A E
\end{matrix}
$$
o\`u, pour tout objet $E$ de ${\mathcal M}od_A$, $A' \otimes_A E$ d\'esigne le conoyau dans ${\mathcal A}b_{\mathcal E}$ du morphisme
$$
\begin{matrix}
A' \times A \times E &\longrightarrow &A' \otimes E \, , \hfill \\
\hfill (a' , a , e) &\longmapsto &(a'a \otimes e - a' \otimes ae) \, ,
\end{matrix}
$$
muni de la loi de multiplication
$$
A' \times A' \otimes_A E \longrightarrow A' \otimes_A E
$$
induite par celle de $A'$
$$
A' \times A' \longrightarrow A' \, .
$$

\item Dans le cas g\'en\'eral, le foncteur
$$
f_* : {\mathcal M}od_{A'} \longrightarrow {\mathcal M}od_A
$$
admet un adjoint \`a gauche
$$
f^* : {\mathcal M}od_A \longrightarrow {\mathcal M}od_{A'}
$$
qui est le compos\'e du foncteur
$$
f^* : {\mathcal M}od_A \longrightarrow {\mathcal M}od_{f^* A}
$$
et du foncteur
$$
\begin{matrix}
{\mathcal M}od_{f^* A} &\longrightarrow &{\mathcal M}od_{A'} \, , \hfill \\
\hfill E &\longmapsto &A' \otimes_{f^* A} E \, .
\end{matrix}
$$
\end{listeimarge}
\end{prop}

\begin{remarks}
\begin{listeisansmarge}
\item Dans la situation de (iii) et si $A$ est un anneau commutatif interne de ${\mathcal E}$, le foncteur
$$
E \longmapsto A' \otimes_A E
$$
se confond avec celui introduit dans la proposition \ref{propIII57} du paragraphe \ref{subsec355}

\medskip

\item Dans le cas g\'en\'eral, l'existence d'un adjoint \`a gauche $f^*$ du foncteur
$$
f_* : {\mathcal M}od_{A'} \longrightarrow {\mathcal M}od_A
$$
est connue a priori d'apr\`es le corollaire \ref{corIII712} (ii).

\smallskip

En effet, comme le foncteur $f_* : {\mathcal E}' \to {\mathcal E}$ respecte les limites arbitraires, il en est de m\^eme du foncteur
$$
f_* : {\mathcal M}od_{A'} \longrightarrow {\mathcal M}od_A \, .
$$

\item Pour tous morphismes de topos annel\'es
$$
({\mathcal E}_1 , A_1) \xrightarrow{ \ f \ } ({\mathcal E}_2 , A_2) \xrightarrow{ \ g \ } ({\mathcal E}_3 , A_3) \, ,
$$
le foncteur d\'efini par leur compos\'e
$$
(g \circ f)_* : {\mathcal M}od_{A_1} \longrightarrow {\mathcal M}od_{A_3}
$$
est le compos\'e des deux foncteurs
$$
f_* : {\mathcal M}od_{A_1} \longrightarrow {\mathcal M}od_{A_2}
$$
et
$$
g_* : {\mathcal M}od_{A_2} \longrightarrow {\mathcal M}od_{A_3} \, .
$$

Par cons\'equent, l'adjoint \`a gauche de $(g \circ f)_*$
$$
(g \circ f)^* : {\mathcal M}od_{A_3} \longrightarrow {\mathcal M}od_{A_1}
$$
est canoniquement isomorphe au compos\'e des adjoints \`a gauche
$$
g^* : {\mathcal M}od_{A_3} \longrightarrow {\mathcal M}od_{A_2} \, ,
$$
$$
f^* : {\mathcal M}od_{A_2} \longrightarrow {\mathcal M}od_{A_1}
$$
de $g_*$ et $f_*$.
\end{listeisansmarge}
\end{remarks}

\medskip

\begin{demo}
\begin{listeisansmarge}
\item r\'esulte de ce que le foncteur $f_* : {\mathcal E}' \to {\mathcal E}$ respecte les produits finis.

\medskip

\item r\'esulte de ce que le foncteur $f^* : {\mathcal E}' \to {\mathcal E}$ adjoint \`a gauche de $f_* : {\mathcal E}' \to {\mathcal E}$ respecte les produits finis.

\medskip

\item Il suffit de traiter le cas o\`u ${\mathcal E} = \widehat{\mathcal C}_J$ est le topos des faisceaux sur un site $({\mathcal C},J)$. Comme le foncteur de faisceautisation $j^* : \widehat{\mathcal C} \to \widehat{\mathcal C}_J$ adjoint \`a gauche de $j_* : \widehat{\mathcal C}_J \hookrightarrow \widehat{\mathcal C}$ respecte les colimites et les limites finies, on peut m\^eme supposer que ${\mathcal E} = \widehat{\mathcal C}$ est le topos des pr\'efaisceaux sur une cat\'egorie essentiellement petite ${\mathcal C}$.

\smallskip

Les calculs de colimites et de limites se faisant composante par composante, on est ramen\'e au cas o\`u ${\mathcal E} = {\rm Ens}$ est le topos des ensembles.

\smallskip

Dans ce cas, pour tout objet $E$ de ${\rm Mod}_A$, la multiplication \`a gauche par les \'el\'ements de $A'$ d\'efinit des endomorphismes de $A' \otimes_A E$ qui font de celui-ci un $A'$-module.

\smallskip

Alors le foncteur
$$
\begin{matrix}
{\rm Mod}_A &\longrightarrow &{\rm Mod}_{A'} \, , \hfill \\
\hfill E &\longmapsto &A' \otimes _A E
\end{matrix}
$$
est adjoint \`a gauche du foncteur d'oubli ${\rm Mod}_{A'} \to {\rm Mod}_{A}$.

\medskip

\item r\'esulte de (ii) et (iii) puisque le foncteur $f_* : {\mathcal M}od_{A'} \to {\mathcal M}od_{A}$ est le compos\'e du foncteur d'oubli ${\mathcal M}od_{A'} \to {\mathcal M}od_{f^* A}$ et du foncteur $f_* : {\mathcal M}od_{f^* A} \to {\mathcal M}od_{A}$. 

\end{listeisansmarge}
\end{demo}

\section{L'\'equivalence de Diaconescu}\label{sec43}

\subsection{Prolongement d'un foncteur d'une cat\'egorie \`a ses pr\'efaisceaux}\label{subsec431}

\medskip

On d\'eduit du lemme \ref{lemIII52}:

\begin{prop}\label{propIV31}

Soit ${\mathcal C}$ une cat\'egorie essentiellement petite.

\smallskip

Soit ${\mathcal E}$ un topos ou, plus g\'en\'eralement, une cat\'egorie cocompl\`ete. 

\smallskip

Alors, pour tout foncteur
$$
\rho : {\mathcal C} \longrightarrow {\mathcal E} \, ,
$$
il existe un foncteur
$$
\widehat\rho : \widehat{\mathcal C} \longrightarrow {\mathcal E} \, ,
$$
unique \`a unique isomorphisme pr\`es, tel que

\medskip

$
\left\{\begin{matrix}
\bullet &\mbox{ce foncteur prolonge $\rho$ au sens que} \hfill \\
{ \ } \\
&\widehat\rho \circ y = \rho \\
{ \ } \\
&\mbox{si $y : {\mathcal C} \hookrightarrow \widehat{\mathcal C}$ d\'esigne le foncteur de Yoneda,} \hfill \\
{ \ } \\
\bullet &\mbox{il respecte les colimites.} \hfill
\end{matrix} \right.
$

\bigskip

Ce foncteur associe \`a tout pr\'efaisceau $P$ sur ${\mathcal C}$ la colimite
$$
\widehat\rho (P) = \varinjlim_{(X,x) \in \int\!\!F} \rho (X) \, .
$$
\end{prop}

\begin{demo}

On sait d'apr\`es le lemme \ref{lemIII52} que tout objet $P$ de $\widehat{\mathcal C}$ s'\'ecrit comme la colimite
$$
F = \varinjlim_{(X,x) \in \int\!\!P} y(X) \, .
$$
Un foncteur $\widehat{\rho} : \widehat{\mathcal C} \to {\mathcal E}$ qui prolonge $\rho$ et respecte les colimites est donc n\'ecessairement donn\'e par la formule
$$
\widehat\rho (P) = \varinjlim_{(X,x) \in \int\!\!P} \rho (X) \, .
$$

R\'eciproquement, le foncteur d\'efini par cette formule prolonge $\rho$.

\smallskip

Nous devons prouver qu'il respecte les colimites.

\smallskip

Consid\'erons donc un carquois $D$, un $D$-diagramme $P_{\bullet}$ de $\widehat{\mathcal C}$ et sa colimite
$$
P = \varinjlim_D P_{\bullet} \, .
$$

Le morphisme canonique
$$
\varinjlim_D \widehat\rho (P_{\bullet}) \longrightarrow \widehat\rho (P)
$$
s'\'ecrit par d\'efinition de $\widehat{\rho}$
$$
\varinjlim_D \varinjlim_{(X,x) \in \int\!\!P_{\bullet}} \rho (X) \longrightarrow \varinjlim_{(X,x) \in \int\!\!P} \rho (X) \, .
$$

C'est un isomorphisme car on a pour tout objet $X$ de ${\mathcal C}$
$$
P(X) = \varinjlim_D P_{\bullet} (X) \, .
$$

Ainsi, le foncteur
$$
\rho : \widehat{\mathcal C} \longrightarrow {\mathcal E}
$$
respecte comme annonc\'e les colimites. 

\end{demo}

\subsection{Morphismes de topos, foncteurs plats et pr\'eservation des recouvrements}\label{subsec432}

\medskip

On d\'emontre \`a partir de la proposition \ref{propIV31}:

\begin{thm}\label{thmIV32}

Soient $({\mathcal C},J)$ un site et
$$
\ell : {\mathcal C} \xrightarrow{ \ y \ } \widehat{\mathcal C} \xrightarrow{ \ j^* \ } \widehat{\mathcal C}_J
$$
le foncteur canonique associ\'e.

\smallskip

Soit
$$
\rho : {\mathcal C} \longrightarrow {\mathcal E}
$$
un foncteur \`a valeurs dans un topos ${\mathcal E}$.

\smallskip

Alors il existe un morphisme de topos
$$
f = (f^* , f_*) : {\mathcal E} \longrightarrow \widehat{\mathcal C}_J \, ,
$$
n\'ecessairement unique \`a unique isomorphisme pr\`es, tel que
$$
\rho = f^* \!\circ \ell
$$
si et seulement si

\medskip

$\left\{\begin{matrix}
\bullet &\mbox{le prolongement canonique de $\rho$} \hfill \\
&\widehat{\rho} : \widehat{\mathcal C} \longrightarrow {\mathcal E} \\
&\mbox{respecte les limites finies,} \hfill \\
{ \ } \\
\bullet &\mbox{le foncteur} \hfill \\
&\rho : {\mathcal C} \longrightarrow {\mathcal E} \\
{ \ } \\
&\mbox{transforme les familles $J$-couvrantes de morphismes de ${\mathcal C}$ en des familles (globalement)} \hfill \\
&\mbox{\'epimorphiques de ${\mathcal E}$.} \hfill
\end{matrix} \right.
$
\end{thm}

\bigskip

\begin{demo}

Si un tel morphisme de topos $f = (f^* , f_*) : {\mathcal E} \to \widehat{\mathcal C}_J$ existe, il est d\'etermin\'e \`a unique isomorphisme pr\`es par sa composante d'image r\'eciproque $f^*$.

\smallskip

Comme $f^* \!\circ j^* \circ y = \rho$ et $f^*$ pr\'eserve les colimites, $f^* \!\circ j^* = \widehat\rho$ est le prolongement canonique de $\rho$ en un foncteur $\widehat\rho : \widehat{\mathcal C} \to {\mathcal E}$ qui respecte les colimites. Le foncteur de faisceautisation $j^* : \widehat{\mathcal C} \to \widehat{\mathcal C}_J$ respecte les limites finies et son compos\'e $\ell$ avec $y$ transforme les familles $J$-couvrantes en familles globalement \'epimorphiques. D'autre part, le foncteur $f^* : \widehat{\mathcal C}_J \to {\mathcal E}$ respecte les limites finies et les familles globalement \'epimorphiques.

\smallskip

Donc $\widehat\rho$ respecte les limites finies et $\rho$ transforme les familles $J$-couvrantes en familles globalement \'epimorphiques.

\smallskip

R\'eciproquement, supposons que ces deux conditions soient satisfaites.

\smallskip

Montrons d'abord qu'elles impliquent que le foncteur $\widehat\rho : \widehat{\mathcal C} \to {\mathcal E}$ se factorise \`a isomorphisme pr\`es en
$$
\widehat{\mathcal C} \xrightarrow{ \ j^* \ } \widehat{\mathcal C}_J \longrightarrow {\mathcal E} \, .
$$

Pour cela, il suffit de prouver que pour tout pr\'efaisceau $P$ et son faisceautis\'e $F = j_* \circ j^* P$, $\widehat\rho$ transforme le morphisme canonique
$$
P \longrightarrow F
$$
en un isomorphisme de ${\mathcal E}$.

\smallskip

Le transform\'e par $j^*$ de $P \to F$ est un \'epimorphisme, et il en est de m\^eme du morphisme
$$
P \times_F y(X) \longrightarrow y(X)
$$
pour tout objet $X$ de ${\mathcal C}$ et tout \'el\'ement $y(X) \to F$ de ${\rm Hom} (y(X),F) = F(X)$. Donc il existe une famille $J$-couvrante de morphismes $U_i \to X$ de ${\mathcal C}$ s'inscrivant dans des carr\'es commutatifs de $\widehat{\mathcal C}$:
$$
\xymatrix{
y(U_i) \ar[d] \ar[r] &y(X) \ar[d] \\
P \ar[r] &F
}
$$

La famille des $\rho (U_i) \to \rho (X)$ est n\'ecessairement (globalement) \'epimorphique dans ${\mathcal E}$ et
$$
\widehat\rho (P) \longrightarrow \widehat\rho (F)
$$
est un \'epimorphisme de ${\mathcal E}$.

\smallskip

D'autre part, le transform\'e par $j^*$ du monomorphisme diagonal
$$
P \longrightarrow P \times_F P = P'
$$
est un isomorphisme. Pour tout objet $X$ de ${\mathcal C}$ et tout \'el\'ement $y(X) \to P'$ de ${\rm Hom} (y(X),P') = P'(X)$, il en est de m\^eme du monomorphisme
$$
P \times_{P'} y(X) \longrightarrow y(X) \, .
$$

Donc il existe une famille $J$-couvrante de morphismes $U_i \to X$ de ${\mathcal C}$ s'inscrivant dans des carr\'es commutatifs de $\widehat{\mathcal C}$:
$$
\xymatrix{
y(U_i) \ar[d] \ar[r] &y(X) \ar[d] \\
P \ar[r] &P \times_F P
}
$$

La famille des $\rho (U_i) \to \rho (X)$ est n\'ecessairement (globalement) \'epimorphique dans ${\mathcal E}$ et le monomorphisme diagonal de ${\mathcal E}$
$$
\widehat\rho (P) \longrightarrow \widehat{\rho} (P \times_F P) = \widehat\rho (P) \times_{\widehat\rho (F)} \widehat\rho (F)
$$
est un \'epimorphisme. Comme le topos ${\mathcal E}$ poss\`ede la propri\'et\'e de balancement, c'est un isomorphisme et
$$
\widehat\rho (P) \longrightarrow \widehat\rho (F)
$$
est un monomorphisme. Or on sait d\'ej\`a que c'est aussi un \'epimorphisme. On conclut comme voulu que
$$
\widehat\rho (P) \longrightarrow \widehat\rho (F)
$$
est un isomorphisme.

\smallskip

Ainsi, le foncteur $\widehat\rho : \widehat{\mathcal C} \to {\mathcal E}$ se factorise en un triangle
$$
\xymatrix{
\widehat{\mathcal C} \ar[rd]_-{\widehat\rho} \ar[rr]^{j^*} &&\widehat{\mathcal C}_J \ar[ld]^-{f^*} \\
&{\mathcal E}
}
$$
commutatif \`a isomorphisme pr\`es.

\smallskip

Comme le foncteur $\widehat\rho$ respecte les limites finies et les colimites arbitraires, il en est de m\^eme de $f^*$.

\smallskip

Enfin, $f^*$ respectant les colimites, il admet un adjoint \`a droite $f_*$ qui d\'efinit avec lui un morphisme de topos
$$
f = (f^* , f_*) : {\mathcal E} \longrightarrow \widehat{\mathcal C}_J \, .
$$
\end{demo}

\bigskip

Ce th\'eor\`eme conduit \`a poser:

\begin{defn}\label{defIV33}

Soit
$$
\rho : {\mathcal C} \longrightarrow {\mathcal E}
$$
un foncteur d'une cat\'egorie essentiellement petite ${\mathcal C}$ vers un topos ${\mathcal E}$.

\smallskip

Alors:

\begin{listeimarge}

\item Ce foncteur est dit ``plat'' si son prolongement canonique respectant les colimites
$$
\widehat\rho : \widehat{\mathcal C} \longrightarrow {\mathcal E}
$$
respecte aussi les limites finies.

\medskip

\item Il est dit ``pr\'eservant les recouvrements'' s'il transforme toute famille $J$-couvrante de morphismes de ${\mathcal C}$ en une famille (globalement) \'epimorphique de ${\mathcal E}$. 
\end{listeimarge}
\end{defn}

Avec cette d\'efinition, on d\'eduit du th\'eor\`eme \ref{thmIV32} le r\'esultat suivant appel\'e ``\'equivalence de Diaconescu'':

\begin{cor}\label{corIV34}

Soient $({\mathcal C},J)$ un site et
$$
\ell : {\mathcal C} \xrightarrow{ \ y \ } \widehat{\mathcal C} \xrightarrow{ \ j^* \ } \widehat{\mathcal C}_J
$$
le foncteur canonique associ\'e.

\smallskip

Alors, pour tout topos ${\mathcal E}$, le foncteur
$$
\begin{matrix}
[{\mathcal E} , \widehat{\mathcal C}_J]_T &\longrightarrow &[{\mathcal C} , {\mathcal E}] \, , \\
\hfill (f^* , f_*) &\longmapsto &f^* \!\circ \ell \hfill
\end{matrix}
$$
d\'efinit une \'equivalence de la cat\'egorie
$$
[{\mathcal E} , \widehat{\mathcal C}_J]_T
$$
des morphismes de topos de ${\mathcal E}$ dans $\widehat{\mathcal C}_J$ sur la sous-cat\'egorie pleine de
$$
[{\mathcal C} , {\mathcal E}]
$$
constitu\'ee des foncteurs
$$
\rho : {\mathcal C} \longrightarrow {\mathcal E}
$$
qui sont ``plats'' et ``pr\'eservent les recouvrements''.
\end{cor}

\begin{demo}

Compte tenu du th\'eor\`eme \ref{thmIV32}, il suffit de remarquer que pour toute paire de morphismes de topos $(f^* , f_*)$ et $(g^* , g_*)$ de ${\mathcal E}$ dans $\widehat{\mathcal C}_J$, toute transformation naturelle de foncteurs de ${\mathcal C}$ dans ${\mathcal E}$
$$
f^* \!\circ \ell \longrightarrow g^* \circ \ell
$$
se rel\`eve en une unique transformation naturelle
$$
f^* \longrightarrow g^* \, .
$$

\end{demo}

\subsection{Description des points d'un topos de faisceaux}\label{subsec433}

Le corollaire \ref{corIV34} comprend le cas particulier suivant:

\begin{cor}\label{corIV35}

Soit ${\mathcal E} \cong \widehat{\mathcal C}_J$ pr\'esent\'e comme le topos des faisceaux sur un site $({\mathcal C},J)$.

\smallskip

Associer \`a tout point $(x^* , x_*)$ de ${\mathcal E}$ le compos\'e du foncteur fibre $x^* : {\mathcal E} \to {\rm Ens}$ avec l'\'equivalence $\widehat{\mathcal C}_J \cong {\mathcal E}$ et le foncteur canonique $\ell : {\mathcal C} \to \widehat{\mathcal C}_J$ d\'efinit une \'equivalence de la cat\'egorie des points de ${\mathcal E}$
$$
{\rm pt} ({\mathcal E})
$$
sur la sous-cat\'egorie pleine de
$$
[{\mathcal C} , {\rm Ens}] = \widehat{{\mathcal C}^{\rm op}}
$$
constitu\'ee des foncteurs covariants
$$
{\mathcal C} \longrightarrow {\rm Ens}
$$
qui sont ``plats'' et ``pr\'eservent les recouvrements''.
\end{cor}

\begin{demo}

Cela r\'esulte du corollaire \ref{corIV34} puisque, par d\'efinition,
$$
{\rm pt} ({\mathcal E}) = [{\rm Ens} , {\mathcal E}]_T \, .
$$

\end{demo}

\newpage

La propri\'et\'e de platitude d'un foncteur ${\mathcal C} \to {\rm Ens}$ s'exprime de plusieurs mani\`eres diff\'erentes:

\begin{lem}\label{lemIV36}

Soit
$$
\rho : {\mathcal C} \longrightarrow {\rm Ens}
$$
un foncteur d\'efini sur une cat\'egorie essentiellement petite ${\mathcal C}$.

\smallskip

Alors les trois propri\'et\'es suivantes sont \'equivalentes:

\begin{enumerate}

\item[(A)] Le foncteur $\rho$ est plat.

\medskip

\item[(B)] Dans la cat\'egorie $[{\mathcal C} , {\rm Ens}]$, l'objet $\rho$ s'\'ecrit comme une colimite filtrante d'un diagramme de foncteurs repr\'esentables, c'est-\`a-dire de la forme ${\rm Hom} (X,\bullet)$.

\medskip

\item[(C)] La cat\'egorie $\int\!\!\rho = {\mathcal C}^{\rm op} / \rho$ des \'el\'ements de $\rho$ est filtrante.
\end{enumerate}
\end{lem}

\begin{demo}

L'implication (C) $\Rightarrow$ (B) r\'esulte de ce que, d'apr\`es le lemme \ref{lemIII52}, tout objet $\rho$ de $[{\mathcal C} , {\rm Ens}] = \widehat{{\mathcal C}^{\rm op}}$ s'\'ecrit comme la colimite
$$
\rho = \varinjlim_{(X,x) \in \int\!\!\rho} y(X)
$$
o\`u $y$ d\'esigne le foncteur
$$
\begin{matrix}
{\mathcal C}^{\rm op} &\longrightarrow &\widehat{{\mathcal C}^{\rm op}} = [{\mathcal C} , {\rm Ens}] \, , \\
\hfill X &\longmapsto &{\rm Hom} (X,\bullet) \, . \hfill
\end{matrix}
$$

Pour l'implication (B) $\Rightarrow$ (A), supposons qu'existe une petite cat\'egorie filtrante ${\mathcal D}$ et un foncteur
$$
\begin{matrix}
{\mathcal D} &\longrightarrow &{\mathcal C}^{\rm op} \, , \\
\hfill d &\longmapsto &X_d \hfill
\end{matrix}
$$
tel que
$$
\rho = \varinjlim_{d \in {\mathcal D}} {\rm Hom} (X_d , \bullet)
$$
dans $\widehat{{\mathcal C}^{\rm op}} = [{\mathcal C} , {\rm Ens}]$.

\smallskip

Pour tout objet $d$ de ${\mathcal D}$, le foncteur d'\'evaluation en $X_d$
$$
\begin{matrix}
\widehat{\mathcal C} &\longrightarrow &{\rm Ens} \, , \\
\hfill P &\longmapsto &P(X_d) \hfill
\end{matrix}
$$
pr\'eserve les colimites et son compos\'e avec
$$
y : {\mathcal C} \longrightarrow \widehat{\mathcal C}
$$
est le foncteur ${\rm Hom} (X_d , \bullet)$. Il n'est donc autre que le prolongement canonique du foncteur ${\rm Hom} (X_d , \bullet)$.

\smallskip

Comme les colimites pr\'eservent les colimites, le foncteur
$$
\begin{matrix}
\widehat{\mathcal C} &\longrightarrow &{\rm Ens} \, , \hfill \\
\hfill P &\longmapsto &\displaystyle \varinjlim_{d \in {\mathcal D}} P(X_d) \hfill
\end{matrix}
$$
respecte les colimites. Il est le prolongement canonique du foncteur $\rho$.

\smallskip

L'implication (B) $\Rightarrow$ (A) r\'esulte alors de ce que, dans la cat\'egorie ${\rm Ens}$, les colimites filtrantes respectent les limites finies.

\smallskip

Montrons enfin l'implication (A) $\Rightarrow$ (C). Supposons que $\rho$ est plat, c'est-\`a-dire que son prolongement canonique $\widehat\rho : \widehat{\mathcal C} \to {\rm Ens}$ respecte les limites finies.

\smallskip

La cat\'egorie $\int\!\!\rho$ est non vide car $\widehat\rho (1)$ est l'objet terminal $\{\bullet\}$ de ${\rm Ens}$.

\smallskip

Pour tous objets $X_1 , X_2$ de ${\mathcal C}$ munis d'\'el\'ements $e_1 \in \rho (X_1)$, $e_2 \in \rho (X_2)$, l'identit\'e
$$
\rho (X_1) \times \rho (X_2) = \widehat\rho (y(X_1) \times y (X_2)) = \varinjlim_{(X,x) \in \int\!\!y (X_1) \times y(X_2)} \rho (X)
$$
montre qu'il existe un objet $X$ de ${\mathcal C}$ muni de deux morphismes
$$
x_1 : X \longrightarrow X_1 \, , \qquad x_2 : X \longrightarrow X_2
$$
et d'un \'el\'ement $e \in \rho (X)$ tels que
$$
\rho (x_1)(e) = e_1 \, , \qquad \rho (x_2)(e) = e_2 \, .
$$

Puis, pour tous objets $X_1,X_2$ de ${\mathcal C}$ reli\'es par deux morphismes $X_1 \raisebox{.7ex}{\xymatrix{\dar[r]^-{^{^{\mbox{\scriptsize$f$}}}}_-{g} &X_2}}$ et pour tout \'el\'ement $e \in \rho (X_1)$ tel que $\rho (f)(e) = \rho (g)(e)$ dans $\rho (X_2)$, l'identit\'e
\begin{eqnarray}
{\rm eg} \left( \rho(X_1) \raisebox{.7ex}{\xymatrix{\dar[r]^-{^{^{\mbox{\scriptsize$\rho(f)$}}}}_-{\rho(g)} &\rho(X_2)}} \right) &= &\widehat\rho \left( {\rm eg} \left( y(X_1) \raisebox{.7ex}{\xymatrix{\dar[r]^-{^{^{\mbox{\scriptsize$y(f)$}}}}_-{y(g)} &y(X_2)}} \right)\right) \nonumber \\
&= &\varinjlim_{(X,x) \in \int\!\!{\rm eg} \bigl( y(X_1) \raisebox{.7ex}{\xymatrix{\dar[r]^-{^{^{\mbox{\scriptsize$y(f)$}}}}_-{y(g)} &y(X_2)}} \bigl)} \rho (X) \nonumber
\end{eqnarray}
montre qu'il existe un objet $X$ de ${\mathcal C}$ muni d'un morphisme
$$
x : X \longrightarrow X_1
$$
et d'un \'el\'ement $e \in \rho (X)$ tel que
$$
f \circ x = g \circ x
$$
et
$$
e_1 = \rho (x)(e) \, .
$$

Cela ach\`eve de montrer comme voulu que la cat\'egorie $\int\!\!\rho = {\mathcal C}^{\rm op} / \rho$ est filtrante. 

\end{demo}

\bigskip

On donne un nom aux pr\'efaisceaux qui s'\'ecrivent comme des colimites filtrantes de pr\'efaisceaux repr\'esen\-tables:

\begin{defn}\label{defIV37}

Soit ${\mathcal C}$ une cat\'egorie essentiellement petite.

\smallskip

On appelle ``cat\'egorie des ind-objets'' de ${\mathcal C}$ et on note
$$
{\rm Ind} ({\mathcal C})
$$
la sous-cat\'egorie pleine de
$$
\widehat{\mathcal C} = [{\mathcal C}^{\rm op} , {\rm Ens}]
$$
constitu\'ee des pr\'efaisceaux qui s'\'ecrivent comme des colimites filtrantes de pr\'efaisceaux repr\'esentables~${\rm Hom} (\bullet , X)$.
\end{defn}

\newpage

\begin{remarksqed}
\begin{listeisansmarge}
\item Le foncteur de Yoneda $y : {\mathcal C} \hookrightarrow \widehat{\mathcal C}$ se factorise en
$$
{\mathcal C} \xhookrightarrow{ \ { \ } \ } {\rm Ind} ({\mathcal C}) \xhookrightarrow{ \ { \ } \ } \widehat{\mathcal C} \, .
$$

\item Dualement, on appelle ``cat\'egorie des pro-objets'' de ${\mathcal C}$ et on note
$$
{\rm Pro} ({\mathcal C})
$$
la sous-cat\'egorie pleine
$$
\bigl({\rm Ind} ({\mathcal C}^{\rm op})\bigl)^{\rm op}
$$
de $[{\mathcal C} , {\rm Ens}]^{\rm op}$.

\smallskip

Le foncteur canonique
$$
\begin{matrix}
{\mathcal C} &\xhookrightarrow{ \ { \ } \ } &[{\mathcal C} , {\rm Ens}]^{\rm op} \hfill \\
X &\longmapsto &{\rm Hom} (X,\bullet)
\end{matrix}
$$
se factorise en
$$
{\mathcal C} \xhookrightarrow{ \ { \ } \ } {\rm Pro} ({\mathcal C}) \xhookrightarrow{ \ { \ } \ } [{\mathcal C} , {\rm Ens}]^{\rm op} \, .
$$
\end{listeisansmarge}
\end{remarksqed}

Avec cette d\'efinition, on d\'eduit du corollaire \ref{corIV35} et du lemme \ref{lemIV36}:

\begin{thm}\label{thmIV38}

Soit ${\mathcal C}$ une cat\'egorie essentiellement petite, munie du foncteur de Yoneda $y : {\mathcal C} \hookrightarrow \widehat{\mathcal C}$.

\smallskip

Alors le foncteur
$$
(x^* , x_*) \longmapsto x^* \circ y
$$
d\'efinit une \'equivalence de la cat\'egorie
$$
{\rm pt} (\widehat{\mathcal C})
$$
des points de $\widehat{\mathcal C}$ sur la cat\'egorie
$$
{\rm Ind} ({\mathcal C}^{\rm op})
$$
des ind-objets de ${\mathcal C}^{\rm op}$.
\end{thm}

\bigskip

\begin{demo}

C'est une simple traduction du corollaire \ref{corIV35}, compte tenu de l'\'equivalence des conditions (A) et (B) du lemme \ref{lemIV36}.
\end{demo}

\subsection{Caract\'erisation des foncteurs plats}\label{subsec434}

D'apr\`es l'\'equivalence de Diaconescu, conna{\^\i}tre les morphismes d'un topos ${\mathcal E}$ dans un topos de faisceaux $\widehat{\mathcal C}_J$ sur un site $({\mathcal C},J)$ revient \`a conna{\^\i}tre les foncteurs ${\mathcal C} \to {\mathcal E}$ qui sont plats et pr\'eservent les recouvrements.

\smallskip

La propri\'et\'e de platitude peut \^etre explicit\'ee concr\`etement (le résultat suivant se trouve déjà dans \cite{Sheaves} comme Théorème VII.9.1):

\begin{prop}\label{propIV39}

Pour un foncteur
$$
\rho : {\mathcal C} \longrightarrow {\mathcal E}
$$
d'une cat\'egorie essentiellement petite ${\mathcal C}$ dans un topos ${\mathcal E}$, les propri\'et\'es suivantes sont \'equivalentes:

\begin{enumerate}

\item[{\rm (A)}] Le foncteur $\rho$ est plat.

\medskip

\item[{\rm (B)}] Les trois conditions suivantes sont satisfaites:

\medskip

\adjustbox{scale=0.865}{$\left\{ \begin{matrix}
{\rm (B1)} &\mbox{Pour tout objet $E$ de ${\mathcal E}$, il existe une famille \'epimorphique de morphismes de ${\mathcal E}$} \hfill \\
{ \ } \\
&E_i \longrightarrow E \, , \qquad i \in I \, , \\
{ \ } \\
&\mbox{et pour tout indice $i$ un objet $X_i$ de ${\mathcal C}$ qui admet un morphisme} \hfill \\
{ \ } \\
&E_i \longrightarrow \rho (X_i) \, . \\
{ \ } \\
{\rm (B2)} &\mbox{Pour tous objets $X,Y$ de ${\mathcal C}$ et tous morphismes de ${\mathcal E}$} \hfill \\
{ \ } \\
&\rho (X) \longleftarrow E \longrightarrow \rho (Y) \, , \\
{ \ } \\
&\mbox{il existe une famille \'epimorphique de morphismes de ${\mathcal E}$} \hfill \\
{ \ } \\
&E_i \longrightarrow E \, , \qquad i \in I \, , \\
&\mbox{et pour chaque indice $i$ deux morphismes de ${\mathcal C}$} \hfill \\
{ \ } \\
&X \longleftarrow W_i \longrightarrow Y \\
&\mbox{et un morphisme de ${\mathcal E}$} \hfill \\
&E_i \longrightarrow \rho (W_i) \\
{ \ } \\
&\mbox{s'inscrivant dans un diagramme commutatif:} \hfill \\
{ \ } \\
&\xymatrix{
&E \ar[ld] \ar[rd] \\
\rho(X) &E_i \ar[u] \ar[d] &\rho (Y) \\
&\rho(W_i) \ar[lu] \ar[ru]
} \\

{\rm (B3)} &\mbox{Pour tous morphismes $X \raisebox{.7ex}{\xymatrix{\dar[r]^-{^{^{\mbox{\scriptsize$f$}}}}_-{g} &Y}}$ de ${\mathcal C}$ et tout morphisme de ${\mathcal E}$} \hfill \\
&E \xrightarrow{ \ e \ } \rho (X) \\
{ \ } \\
&\mbox{tel que $\rho (f) \circ e = \rho (g) \circ e$, il existe une famille \'epimorphique de morphismes de ${\mathcal E}$} \hfill \\
{ \ } \\
&E_i \longrightarrow E \, , \qquad i \in I \, , \\
{ \ } \\
&\mbox{et pour chaque indice $i$ un morphisme de ${\mathcal C}$} \hfill \\
{ \ } \\
&w_i : W_i \longrightarrow X \\
{ \ } \\
&\mbox{tel que $f \circ w_i = g \circ w_i$ et un morphisme de ${\mathcal E}$} \hfill \\
{ \ } \\
&E_i \longrightarrow \rho (W_i) \\
{ \ } \\
&\mbox{s'inscrivant dans un carr\'e commutatif:} \hfill \\
{ \ } \\
&\xymatrix{
E_i \ar[d] \ar[r] &\rho (W_i) \ar[d] \\
E \ar[r] &\rho (X)
} 
\end{matrix} \right.
$}

\medskip

\item[{\rm (C)}] Pour n'importe quelle  repr\'esentation de ${\mathcal E}$ comme topos des faisceaux sur un site $({\mathcal D},K)$
$$
\widehat{\mathcal D}_K \cong {\mathcal E} \, ,
$$
avec donc un foncteur canonique associ\'e
$$
\ell : {\mathcal D} \xrightarrow{ \ y \ } \widehat{\mathcal D} \xrightarrow{ \ j^* \ } \widehat{\mathcal D}_K \cong {\mathcal E} \, ,
$$
sont v\'erifi\'ees les propri\'et\'es (C1), (C2), (C3) d\'eduites de (B1), (B2), (B3) en prenant pour objets $E$ de ${\mathcal E}$ ceux de la forme
$$
\ell (D) \, , \qquad D \in {\rm Ob} ({\mathcal D}) \, ,
$$
et pour familles \'epimorphiques
$$
(E_i \longrightarrow E)_{i \in I}
$$
les images par $\ell$ de familles $K$-couvrantes de ${\mathcal D}$
$$
(D_i \longrightarrow D)_{i \in I} \, .
$$
\end{enumerate}
\end{prop}

\begin{remarks}
\begin{listeisansmarge}
\item Dans le cas o\`u ${\mathcal E} = {\rm Ens}$ est le topos des ensembles, on peut prendre pour ${\mathcal D}$ la cat\'egorie \`a un seul objet et un seul morphisme.

\smallskip

Alors la propri\'et\'e (C) exprime exactement que la cat\'egorie $\int\!\!\rho = {\mathcal C}^{\rm op} / \rho$ des \'el\'ements de $\rho$ est filtrante. C'est la condition (C) du lemme \ref{lemIV36}.

\medskip

\item Plus g\'en\'eralement, si ${\mathcal E} = \widehat{\mathcal D}$ est la cat\'egorie des pr\'efaisceaux sur une cat\'egorie essentiellement petite ${\mathcal D}$, (C) signifie qu'un foncteur
$$
\rho : {\mathcal C} \longrightarrow \widehat{\mathcal D}
$$
est plat si et seulement si, pour tout objet $D$ de ${\mathcal D}$, la cat\'egorie $\int\!\!\rho_D = {\mathcal C}^{\rm op} / \rho_D$ des \'el\'ements du foncteur
$$
\begin{matrix}
\rho_D : {\mathcal C} &\longrightarrow &{\rm Ens} \, , \hfill \\
\hfill X &\longmapsto &\rho (X)(D)
\end{matrix}
$$
est filtrante.

\medskip

\item Dans le cas g\'en\'eral, la platitude de $\rho : {\mathcal C} \to {\mathcal E}$ signifie par d\'efinition que son unique prolongement respectant les colimites
$$
\begin{matrix}
\widehat\rho : \widehat{\mathcal C} &\longrightarrow &{\mathcal E} \, , \hfill \\
\hfill P &\longmapsto &\displaystyle \varinjlim_{(X,x) \in \int\!\!P} \rho (X)
\end{matrix}
$$
respecte les limites finies, c'est-\`a-dire respecte les objets terminaux, les produits de deux objets et les \'egalisateurs de deux morphismes.

\smallskip

Autrement dit, cela signifie que $\rho$ doit satisfaire les trois conditions suivantes:

\bigskip

$\left\{\begin{matrix}
{\rm (A1)} &\mbox{Le morphisme canonique de ${\mathcal E}$} \hfill \\
{ \ } \\
&\displaystyle\varinjlim_{X \in {\mathcal C}} \rho (X) \longrightarrow 1_{\mathcal E} \\
&\mbox{est un isomorphisme.} \hfill \\
{ \ } \\
{\rm (A2)} &\mbox{Pour tous pr\'efaisceaux $P,Q$ sur ${\mathcal C}$, le morphisme canonique de ${\mathcal E}$} \hfill \\
{ \ } \\
&\displaystyle \varinjlim_{\mbox{\footnotesize$\begin{matrix} W \in {\mathcal C} \\ p \in P(W) \\ q \in Q(W) \end{matrix}$}} \rho (W) \longrightarrow \left( \varinjlim_{\mbox{\footnotesize$\begin{matrix} X \in {\mathcal C} \\ p \in P(X) \end{matrix}$}} \rho (X) \right) \times \left( \varinjlim_{\mbox{\footnotesize$\begin{matrix} Y \in {\mathcal C} \\ q \in Q(Y) \end{matrix}$}} \rho (Y) \right)  \\
{ \ } \\
&\mbox{est un isomorphisme.} \hfill \\
{\rm (A3)} &\mbox{Pour tous morphismes de pr\'efaisceaux $P \raisebox{.7ex}{\xymatrix{\dar[r]^-{^{^{\mbox{\scriptsize$f$}}}}_-{g} &Q}}$ sur ${\mathcal C}$, le morphisme canonique de ${\mathcal E}$} \hfill \\
{ \ } \\
&\displaystyle \varinjlim_{\mbox{\footnotesize$\begin{matrix} W \in {\mathcal C} \\ p \in P(W) , f_W(p)=g_W(p) \end{matrix}$}} \rho (W) \longrightarrow {\rm eg} \left[ \left( \varinjlim_{\mbox{\footnotesize$\begin{matrix} X \in {\mathcal C} \\ p \in P(X) \end{matrix}$}} \rho (X) \right) \rightrightarrows \left( \varinjlim_{\mbox{\footnotesize$\begin{matrix} Y \in {\mathcal C} \\ q \in Q(Y) \end{matrix}$}} \rho (Y) \right) \right] \\
{ \ } \\
&\mbox{est un isomorphisme.} \hfill
\end{matrix} \right.
$
\end{listeisansmarge}
\end{remarks}
\bigskip

\begin{demo}

Montrons d'abord que les propri\'et\'es (B) et (C) de l'\'enonc\'e sont \'equivalentes.

\smallskip

L'implication (B) $\Rightarrow$ (C) r\'esulte de ce que, si ${\mathcal E} \cong \widehat{\mathcal D}_K$ et $\ell : {\mathcal D} \xrightarrow{ \ y \ } \widehat{\mathcal D} \xrightarrow{ \ j^* \ } \widehat{\mathcal D}_K \cong {\mathcal E}$ est le foncteur compos\'e, alors pour tout objet $D$ de ${\mathcal D}$ et toute famille \'epimorphique de morphismes de ${\mathcal E}$
$$
E_i \longrightarrow \ell (D) \, , \qquad i \in I \, ,
$$
il existe pour chaque indice $i$ une famille d'objets $D_{i,j}$ de ${\mathcal D}$ et une famille \'epimorphique de morphismes de~${\mathcal E}$
$$
\ell (D_{i,j}) \longrightarrow E_i
$$
telles que les morphismes compos\'es
$$
\ell (D_{i,j}) \longrightarrow \ell (D)
$$
soient les images par $\ell$ de morphismes de ${\mathcal D}$
$$
D_{i,j} \longrightarrow D \, .
$$

En effet, la famille des morphismes
$$
\ell (D_{i,j}) \longrightarrow \ell (D)
$$
est \'epimorphique, ce qui signifie que la famille des morphismes de ${\mathcal D}$
$$
D_{i,j} \longrightarrow D
$$
est $K$-couvrante.

\smallskip

L'implication en sens inverse (C) $\Rightarrow$ (B) r\'esulte de ce que tout objet $E$ de ${\mathcal E}$ admet une famille \'epimorphique de morphismes de la forme
$$
\ell (D) \longrightarrow E
$$
pour des objets $D$ de ${\mathcal D}$.

\smallskip

Comme (C) est \'equivalente \`a (B), la propri\'et\'e (C) ne d\'epend pas du choix d'une pr\'esentation de ${\mathcal E}$
$$
\widehat{\mathcal D}_K \cong {\mathcal E}
$$
comme cat\'egorie des faisceaux sur un site $({\mathcal D},K)$.

\smallskip

On peut donc supposer que ${\mathcal D}$ est une sous-cat\'egorie pleine de ${\mathcal E}$ essentiellement petite, qui contient les images par $\rho$ des objets de ${\mathcal C}$ ainsi qu'une famille s\'eparante d'objets de ${\mathcal E}$. On sait en effet d'apr\`es le th\'eor\`eme \ref{thmIII81} que si $K$ d\'esigne la topologie de ${\mathcal D}$ pour laquelle une famille est couvrante quand elle est \'epimorphique dans ${\mathcal E}$, on a une \'equivalence
$$
\begin{matrix}
\widehat{\mathcal D}_K &\xrightarrow{ \ \sim \ } &{\mathcal E} \, , \hfill \\
\hfill F &\longmapsto &\displaystyle\varinjlim_{(D,d) \in \int\!\!F} D \, .
\end{matrix}
$$

Avec ce choix de ${\mathcal D}$, le foncteur $\rho : {\mathcal C} \to {\mathcal E}$ se factorise en
$$
{\mathcal C} \xrightarrow{ \ \rho \ } {\mathcal D} \xhookrightarrow{ \ { \ } \ } {\mathcal E}
$$
et le foncteur compos\'e
$$
{\mathcal C} \xrightarrow{ \ \rho \ } {\mathcal D} \xhookrightarrow{ \ y \ } \widehat{\mathcal D}
$$
se prolonge en un foncteur respectant les colimites
$$
\rho_! : \widehat{\mathcal C} \longrightarrow \widehat{\mathcal D} \, .
$$

Comme le foncteur de faisceautisation
$$
j^* : \widehat{\mathcal D} \longrightarrow \widehat{\mathcal D}_K
$$
respecte les colimites, le diagramme
$$
\xymatrix{
\widehat{\mathcal C} \ar[d]_{\rho_!} \ar[rrd]^{\widehat\rho} \\
\widehat{\mathcal D} \ar[r]^{j^*} &\widehat{\mathcal D}_K \ar[r]^{\sim} &{\mathcal E}
}
$$
est commutatif \`a isomorphisme canonique pr\`es.

\smallskip

Autrement dit, pour tout pr\'efaisceau $P$ sur ${\mathcal C}$, son image dans ${\mathcal E} \cong \widehat{\mathcal D}_K$
$$
\widehat\rho (P)
$$
s'identifie \`a la faisceautisation
$$
j^* (\rho_! (P))
$$
du pr\'efaisceau
$$
\begin{matrix}
\rho_! (P) : {\mathcal D}^{\rm op} &\longrightarrow &{\rm Ens} \, , \hfill \\
\hfill D &\longmapsto &\displaystyle\varinjlim_{(X,x) \in \int\!\!P} {\rm Hom} (D,\rho (X)) \, .
\end{matrix}
$$

Prenant pour $P$ l'objet terminal de $\widehat{\mathcal C}$ qui est le pr\'efaisceau constant \'egal \`a $\{\bullet\}$, ou bien le produit $y(X) \times y(Y)$ de deux pr\'efaisceaux repr\'esentables, ou encore l'\'egalisateur de deux morphismes $y(X) \rightrightarrows y(Y)$ reliant deux pr\'efaisceaux repr\'esentables, on voit que les conditions (C1), (C2) et (C3) se r\'e\'ecrivent sous la forme:

\medskip

$\left\{ \begin{matrix}
({\rm C1'}) &\mbox{Le morphisme canonique de ${\mathcal E} \cong \widehat{\mathcal D}_K$} \hfill \\
{ \ } \\
&\widehat\rho (1) \longrightarrow 1_{\mathcal E} \\
&\mbox{est un \'epimorphisme.} \hfill \\
{ \ } \\
({\rm C2'})  &\mbox{Pour tous objets $X,Y$ de ${\mathcal C}$, le morphisme canonique de ${\mathcal E}$} \hfill \\
{ \ } \\
&\widehat\rho (y(X) \times y(Y)) \longrightarrow \rho (X) \times \rho (Y) \\
&\mbox{est un \'epimorphisme.} \hfill \\
{ \ } \\
({\rm C3'})  &\mbox{Pour toute paire de morphismes $X \raisebox{.7ex}{\xymatrix{\dar[r]^-{^{^{\mbox{\scriptsize$f$}}}}_-{g} &Y}}$ de ${\mathcal C}$, le morphisme canonique de ${\mathcal E}$} \hfill \\
{ \ } \\
&\widehat\rho \left({\rm eg} (y(X) \rightrightarrows y(Y)\right) \longrightarrow {\rm eg} \left(\rho (X) \rightrightarrows \rho (Y)\right) \\
{ \ } \\
&\mbox{est un \'epimorphisme.} \hfill 
\end{matrix} \right.
$

\bigskip

Comme (A) s'explicite en les trois conditions (A1), (A2) et (A3) de la remarque (iii), on voit d\'ej\`a que (A) implique (C).

\smallskip

R\'eciproquement, supposons que (B) est satisfaite.

\smallskip

Montrons d'abord que la propri\'et\'e (A1) est satisfaite, c'est-\`a-dire que le morphisme de ${\mathcal E}$
$$
\widehat\rho (1) = \varinjlim_{X \in {\mathcal C}} \rho (X) \longrightarrow 1_{\mathcal E}
$$
est un isomorphisme. C'est un \'epimorphisme d'apr\`es (B1) donc il suffit de prouver que c'est un monomorphisme.

\smallskip

Consid\'erons donc un objet $E$ de ${\mathcal E}$ muni de deux morphismes
$$
E \rightrightarrows \varinjlim_{X \in {\mathcal C}} \rho (X) = \widehat\rho (1) \, .
$$

Alors $E$ est la colimite des produits fibr\'es
$$
E \times_{\widehat\rho(1) \times \widehat\rho(1)} \rho (X_1) \times \rho (X_2) \, .
$$

Il suffit donc de montrer que si $E$ est un objet de ${\mathcal E}$ muni de deux morphismes
$$
\rho (X_1) \longleftarrow E \longrightarrow \rho (X_2) \quad \mbox{avec} \quad X_1 , X_2 \in {\rm Ob} ({\mathcal C}) \, ,
$$
alors leurs compos\'es avec les morphismes canoniques
$$
\rho (X_1) \longrightarrow \varinjlim_{X \in {\mathcal C}} \rho (X) \longleftarrow \rho (X_2)
$$
sont \'egaux.

\smallskip

Or, d'apr\`es la propri\'et\'e (B2), il existe une famille \'epimorphique
$$
E_i \longrightarrow E \, , \qquad i \in I \, ,
$$
et pour chaque indice $i$ deux morphismes de ${\mathcal C}$
$$
X_1 \longleftarrow W_i \longrightarrow X_2
$$
et un morphisme de ${\mathcal E}$
$$
E_i \longrightarrow \rho (W_i)
$$
s'inscrivant dans un diagramme commutatif:
$$
\xymatrix{
&E \ar[ld] \ar[rd] \\
\rho(X_1) &E_i \ar[u] \ar[d] &\rho (X_2) \\
&\rho(W_i) \ar[lu] \ar[ru]
} 
$$

Donc les morphismes compos\'es
$$
E \rightrightarrows \varinjlim_{X \in {\mathcal C}} \rho (X)
$$
sont \'egaux et la propri\'et\'e (A1) est satisfaite.

\smallskip

Puis montrons que la propri\'et\'e (A2) est satisfaite c'est-\`a-dire que pour tous pr\'efaisceaux $P,Q$ sur ${\mathcal C}$, le morphisme
$$
\widehat\rho (P \times Q) \longrightarrow \widehat\rho(P) \times \widehat\rho(Q)
$$
est un isomorphisme.

\smallskip

Comme on a
$$
P = \varinjlim_{(X,x) \in \int\!\!P} y(X)
$$
et
$$
Q = \varinjlim_{(X,y) \in \int\!\!Q} y(Y) \, ,
$$
et que d'autre part le foncteur $\widehat\rho$ respecte les colimites, de m\^eme que les foncteurs de produit dans $\widehat{\mathcal C}$ et ${\mathcal E}$, il suffit de prouver que, pour tous objets $X$ et $Y$ de ${\mathcal C}$, le morphisme de ${\mathcal E}$
$$
\widehat\rho (y(X) \times y(Y)) \longrightarrow \rho (X) \times \rho (Y)
$$
est un isomorphisme. D'apr\`es la propri\'et\'e (B2), c'est un \'epimorphisme donc il suffit de montrer que c'est un monomorphisme.

\smallskip

Consid\'erons donc un objet $E$ de ${\mathcal E}$ et deux morphismes
$$
E \rightrightarrows \widehat\rho (y(X) \times y(Y))
$$
induisant le m\^eme compos\'e $E \to \rho (X) \times \rho (Y)$.

\smallskip

Comme
$$
\widehat\rho (y(X) \times y(Y)) = \varinjlim_{\mbox{\footnotesize $\begin{pmatrix} W \to X \\ \ \downarrow \hfill \\ \ Y \hfill \end{pmatrix}$}} \rho (W) \, ,
$$
on peut supposer qu'existent des objets $W_1$ et $W_2$ de ${\mathcal C}$ munis de morphismes
$$
X \longleftarrow W_1 \longrightarrow Y \, ,
$$
$$
X \longleftarrow W_2 \longrightarrow Y \, ,
$$
tels que les deux morphismes $E \rightrightarrows \widehat\rho (y(X) \times y(Y))$ se factorisent en
$$
E \longrightarrow \rho (W_1) \, ,
$$
$$
E \longrightarrow \rho (W_2) \, .
$$

Quitte \`a remplacer $E$ par les sources d'une famille \'epimorphique de morphismes de but $E$, il existe d'apr\`es la propri\'et\'e (B2) deux morphismes de ${\mathcal C}$
$$
W_1 \longleftarrow W \longrightarrow W_2 
$$
et un morphisme de ${\mathcal E}$
$$
E \longrightarrow \rho (W)
$$
s'inscrivant dans un diagramme commutatif:
$$
\xymatrix{
&E \ar[ld] \ar[d] \ar[rd] \\
\rho (W_1) &\rho(W) \ar[l] \ar[r] &\rho (W_2)
}
$$

Les compos\'es du morphisme $E \to \rho (W)$ avec les morphismes
$$
\rho (W) \rightrightarrows \rho (X) \qquad \mbox{et} \qquad \rho (W) \rightrightarrows \rho (Y)
$$
d\'eduits de $W_1 \leftarrow W \rightarrow W_2$ et de $X \leftarrow W_1 \rightarrow Y$ et $X \leftarrow W_2 \rightarrow Y$ sont \'egaux.

\smallskip

Quitte \`a remplacer \`a nouveau $E$ par les sources d'une famille \'epimorphique de morphismes de but $E$, il existe d'apr\`es la propri\'et\'e (B3) deux morphismes de ${\mathcal C}$
$$
V_1 \longrightarrow W \qquad \mbox{et} \qquad V_2 \longrightarrow W
$$
ayant respectivement le m\^eme compos\'e avec
$$
W \rightrightarrows X \qquad \mbox{et} \qquad W \rightrightarrows Y
$$
et tels que $E \to \rho (W)$ se factorise en
$$
E \longrightarrow \rho (V_1) \longrightarrow \rho (W) \qquad \mbox{et} \qquad E \longrightarrow \rho (V_2) \longrightarrow \rho (W) \, .
$$

Puis, quitte \`a remplacer encore $E$ par les sources d'une famille \'epimorphique de morphismes de but $E$, il existe d'apr\`es la propri\'et\'e (B2) deux morphismes de ${\mathcal C}$
$$
V_1 \longleftarrow V \longrightarrow V_2
$$
et un morphisme de ${\mathcal E}$
$$
E \longrightarrow \rho (V)
$$
s'inscrivant dans un diagramme commutatif:
$$
\xymatrix{
&E \ar[ld] \ar[d] \ar[rd] \\
\rho (V_1) &\rho(V) \ar[l] \ar[r] &\rho (V_2)
}
$$

Les compos\'es de $V \to V_1 \to W$ avec les deux morphismes $W \rightrightarrows X$ sont \'egaux, de m\^eme que les compos\'es de $V \to V_2 \to W$ avec les deux morphismes $W \rightrightarrows Y$.

\smallskip

Donc les morphismes compos\'es
$$
E \rightrightarrows \widehat \rho (y(X) \times y(Y)) \longrightarrow \rho (X) \times \rho (Y) 
$$
sont \'egaux et la propri\'et\'e (A2) est satisfaite.

\smallskip

Il reste \`a montrer que la propri\'et\'e (A3) est satisfaite, c'est-\`a-dire que pour tous morphismes de pr\'efaisceaux sur ${\mathcal C}$
$$
P \rightrightarrows Q \, ,
$$
le morphisme de ${\mathcal E}$
$$
\widehat\rho \left({\rm eg} (P \rightrightarrows Q) \right) \longrightarrow {\rm eg} \left( \widehat\rho (P) \rightrightarrows \widehat\rho (Q) \right)
$$ 
est un isomorphisme.

\smallskip

Comme on peut \'ecrire
$$
Q = \varinjlim_{(Y,y) \in \int\!\!Q} y(Y)
$$
et que le foncteur $\widehat\rho$ respecte les colimites, de m\^eme que les foncteurs de produits fibr\'es dans $\widehat{\mathcal C}$ et ${\mathcal E}$, on peut supposer que $Q$ est repr\'esentable par un objet $Y$ de ${\mathcal C}$.

\smallskip

Pour les m\^emes raisons, on peut supposer que $P$ est repr\'esentable par un objet $X$ de ${\mathcal C}$ et donc que $P \rightrightarrows Q$ provient d'une paire de morphismes $X \raisebox{.7ex}{\xymatrix{\dar[r]^-{^{^{\mbox{\scriptsize$f$}}}}_-{g} &Y}}$ de ${\mathcal C}$.

\smallskip

La propri\'et\'e (B3) signifie que le morphisme
$$
\widehat\rho \left({\rm eg} (y(X) \rightrightarrows y(Y)) \right) \longrightarrow {\rm eg} \left( \rho (X) \rightrightarrows \rho (Y) \right)
$$ 
est un \'epimorphisme.

\smallskip

Comme les deux morphismes de $\widehat{\mathcal C}$ et ${\mathcal E}$
$$
{\rm eg} (y(X) \rightrightarrows y(Y)) \longrightarrow y(X)
$$
et
$$
{\rm eg} (\rho(X) \rightrightarrows \rho(Y)) \longrightarrow \rho(X)
$$
sont des monomorphismes, il suffit de d\'emontrer que $\widehat\rho$ transforme tout monomorphisme de $\widehat{\mathcal C}$ de la forme
$$
S \xhookrightarrow{ \ { \ } \ } y(X)
$$
en un monomorphisme de ${\mathcal E}$
$$
\widehat\rho (S) \xhookrightarrow{ \ { \ } \ } \rho (X) \, .
$$

Consid\'erons donc un objet de $E$ muni de deux morphismes
$$
E \rightrightarrows \widehat\rho (S)
$$
dont les compos\'es avec $\widehat\rho (S) \to \rho (X)$ sont \'egaux.

\smallskip

Comme on a
$$
\widehat\rho (S) = \varinjlim_{(W,w) \in \int\!\!S} \rho (W) \, ,
$$
on peut supposer que ces deux morphismes se factorisent en
$$
E \longrightarrow \rho (W_1) \longrightarrow \widehat\rho (S)
$$
et
$$
E \longrightarrow \rho (W_2) \longrightarrow \widehat\rho (S)
$$
pour deux objets $W_1$ et $W_2$ de ${\mathcal C}$ munis de morphismes
$$
y(W_1) \longrightarrow S \qquad \mbox{et} \qquad y(W_2) \longrightarrow S \, .
$$

Quitte \`a remplacer $E$ par les sources d'une famille \'epimorphique de morphismes de but $E$, il existe d'apr\`es la propri\'et\'e (B2) deux morphismes de ${\mathcal C}$
$$
W_1 \longleftarrow W \longrightarrow W_2
$$
et un morphisme de ${\mathcal E}$
$$
E \longrightarrow \rho (W)
$$
s'inscrivant dans un diagramme commutatif:
$$
\xymatrix{
&E \ar[ld] \ar[d] \ar[rd] \\
\rho (W_1) &\rho(W) \ar[l] \ar[r] &\rho (W_2)
}
$$

Les compos\'es du morphisme
$$
E \longrightarrow \rho (W)
$$
avec les deux morphismes $\rho (W) \rightrightarrows \rho (X)$ images par $\rho$ des deux compos\'es
$$
W \longrightarrow W_1 \longrightarrow X
$$
et
$$
W \longrightarrow W_2 \longrightarrow X
$$
sont \'egaux.

\smallskip

Quitte \`a remplacer \`a nouveau $E$ par les sources d'une famille \'epimorphique de morphismes de but $E$, il existe d'apr\`es la propri\'et\'e (B3) un morphisme de ${\mathcal C}$
$$
V \longrightarrow W
$$
ayant le m\^eme compos\'e avec les deux morphismes $W \rightrightarrows X$, tel que le morphisme
$$
E \longrightarrow \rho (W)
$$
se factorise en
$$
E \longrightarrow \rho (V) \longrightarrow \rho (W) \, .
$$

Or les deux compos\'es de $y(V) \to y(W)$ avec
$$
y(W) \longrightarrow y(W_1) \longrightarrow S
$$
et
$$
y(W) \longrightarrow y(W_2) \longrightarrow S
$$
sont \'egaux, puisqu'il en est ainsi de leurs compos\'es avec le monomorphisme
$$
S \xhookrightarrow{ \ { \ } \ } y(X) \, .
$$

Donc les deux morphismes $E \rightrightarrows \widehat\rho (S)$ sont \'egaux et la propri\'et\'e (A3) est satisfaite, ce qui termine la d\'emonstration. \end{demo}

On d\'eduit du th\'eor\`eme \ref{thmIV38} le cas particulier suivant:

\begin{cor}\label{corIV310}

Soit ${\mathcal C}$ une cat\'egorie essentiellement petite qui poss\`ede des limites finies arbitraires.

\smallskip

Alors un foncteur
$$
\rho : {\mathcal C} \longrightarrow {\mathcal E}
$$
de ${\mathcal C}$ dans un topos ${\mathcal E}$ est plat si et seulement si il respecte les produits finis.

\smallskip

Pour que $\rho$ soit plat, il suffit m\^eme qu'il satisfasse les conditions suivantes:

\medskip

$\left\{ \begin{matrix}
{\rm (B1')} &\mbox{Notant $1_{\mathcal C}$ et $1_{\mathcal E}$ les objets terminaux de ${\mathcal C}$ et ${\mathcal E}$, le morphisme canonique} \hfill \\
{ \ } \\
&\rho (1_{\mathcal C}) \longrightarrow 1_{\mathcal E} \\
&\mbox{est un \'epimorphisme.} \hfill \\
{ \ } \\
{\rm (B2')} &\mbox{Pour tous objets $X$ et $Y$ de ${\mathcal C}$, le morphisme canonique} \hfill \\
{ \ } \\
&\rho (X \times Y) \longrightarrow \rho(X) \times \rho (Y) \\
&\mbox{est un \'epimorphisme.} \hfill \\
{\rm (B3')} &\mbox{Pour tous morphismes $X \raisebox{.7ex}{\xymatrix{\dar[r]^-{^{^{\mbox{\scriptsize$f$}}}}_-{g} &Y}}$ de ${\mathcal C}$, le morphisme canonique} \hfill \\
{ \ } \\
&\rho \left({\rm eg} (X \rightrightarrows Y)\right) \longrightarrow {\rm eg} \left( \rho(X) \rightrightarrows \rho (Y) \right) \\
{ \ } \\
&\mbox{est un \'epimorphisme.} \hfill
\end{matrix} \right.
$
\end{cor}

\bigskip

\begin{demo}

Si $\rho : {\mathcal C} \to {\mathcal E}$ est plat, c'est-\`a-dire si son prolongement canonique pr\'eservant les colimites
$$
\widehat\rho : \widehat{\mathcal C} \longrightarrow {\mathcal E}
$$
respecte les limites finies, il en est a fortiori de m\^eme de
$$
\rho = \widehat\rho \circ y
$$
puisque le foncteur de Yoneda
$$
y : {\mathcal C} \xhookrightarrow{ \ { \ } \ } \widehat{\mathcal C}
$$
respecte toujours les limites.

\smallskip

Si $\rho$ respecte les limites finies, les morphismes de ${\mathcal E}$ envisag\'ees dans (B1$'$), (B2$'$) et (B3$'$) sont toujours des isomorphismes. A fortiori, ce sont des \'epimorphismes.

\smallskip

Enfin, supposons que $\rho$ satisfasse les conditions (B1$'$), (B2$'$) et (B3$'$).

\smallskip

Alors, pour tout objet $E$ de ${\mathcal E}$, la projection
$$
E' = E \times \rho (1_{\mathcal E}) \longrightarrow E
$$
est un \'epimorphisme et $E'$ est muni d'un morphisme
$$
E' = E \times \rho (1_{\mathcal E}) \longrightarrow \rho (1_{\mathcal E})
$$
si bien que la condition (B1) du th\'eor\`eme \ref{thmIV38} satisfaite.

\smallskip

De m\^eme, les conditions (B2) et (B3) de ce th\'eor\`eme sont satisfaites.

\smallskip

En effet, pour tous objets $X$ et $Y$ de ${\mathcal C}$ et tout objet $E$ de ${\mathcal E}$ muni d'un morphisme
$$
E \longrightarrow \rho (X) \times \rho (Y) \, ,
$$
le morphisme
$$
E' = E \times_{\rho (X) \times \rho (Y)} \rho (X \times Y) \longrightarrow \rho (X \times Y)
$$
est un \'epimorphisme, et il s'inscrit dans le carr\'e commutatif:
$$
\xymatrix{
E' \ar[d] \ar[r] &\rho (X \times Y) \ar[d] \\
E \ar[r] &\rho (X) \times \rho (Y)
}
$$

D'autre part, pour tous morphismes $X \raisebox{.7ex}{\xymatrix{\dar[r]^-{^{^{\mbox{\scriptsize$f$}}}}_-{g} &Y}}$ de ${\mathcal C}$ et tout objet $E$ de ${\mathcal E}$ muni d'un morphisme
$$
E \longrightarrow {\rm eg} \left( \rho (X) \rightrightarrows \rho (Y)\right) \, ,
$$
le morphisme
$$
E' = E \times_{{\rm eg} (\rho (X) \rightrightarrows \rho (Y))} \rho \left( {\rm eg} (X \rightrightarrows Y) \right) \longrightarrow E
$$
est un \'epimorphisme, et il s'inscrit dans le carr\'e commutatif:
$$
\xymatrix{
E' \ar[d] \ar[r] &\rho ({\rm eg} (X \rightrightarrows Y)) \ar[d] \\
E \ar[r] &{\rm eg} (\rho (X) \rightrightarrows \rho (Y))
}
$$

Le foncteur $\rho$ v\'erifiant les conditions (B1), (B2), (B3) du th\'eor\`eme \ref{thmIV38}, on conlut d'apr\`es ce th\'eor\`eme qu'il est un foncteur plat. 

\end{demo}

\subsection{Description des morphismes d'un topos dans un topos de faisceaux}\label{subsec435}

\medskip

Combinant le corollaire \ref{corIV34} et le th\'eor\`eme \ref{thmIV38}, on obtient:

\begin{thm}\label{thmIV311}

Soient $({\mathcal C},J)$ un site et
$$
\ell : {\mathcal C} \xhookrightarrow{ \ y \ } \widehat{\mathcal C} \xrightarrow{ \ j^* \ } \widehat{\mathcal C}_J
$$
le foncteur canonique associ\'e.

\smallskip

Alors, pour tout topos ${\mathcal E}$, le foncteur
$$
\begin{matrix}
[{\mathcal E} , \widehat{\mathcal C}_J]_T &\longrightarrow &[{\mathcal C},{\mathcal E}] \, , \\
\hfill (f^* , f_*) &\longmapsto &f^* \!\circ \ell \hfill
\end{matrix}
$$
d\'efinit une \'equivalence de la cat\'egorie
$$
[{\mathcal E} , \widehat{\mathcal C}_J]_T
$$
des morphismes de topos de ${\mathcal E}$ dans $\widehat{\mathcal C}_J$ sur la sous-cat\'egorie pleine de
$$
[{\mathcal C} , {\mathcal E}]
$$
constitu\'ee des foncteurs
$$
\rho : {\mathcal C} \longrightarrow {\mathcal E}
$$
qui 

\medskip

\begin{enumerate}
\item[$\bullet$] ``pr\'eservent les recouvrements'' c'est-\`a-dire transforment toute famille $J$-couvrante de morphismes de~${\mathcal C}$
$$
(X_i \longrightarrow X)_{i \in I}
$$
en une famille \'epimorphique de morphismes de ${\mathcal E}$,
\item[$\bullet$] sont ``plats'' c'est-\`a-dire satisfont les conditions suivantes:
\end{enumerate}

\adjustbox{scale=0.85}{
$\left\{\begin{matrix}
{\rm (B1)} &\mbox{Tout objet $E$ de ${\mathcal E}$ admet une famille \'epimorphique} \hfill \\
{ \ } \\
&E_i \longrightarrow E \, , \qquad i \in I \, , \\
{ \ } \\
&\mbox{dont chaque source admet un morphisme} \hfill \\
{ \ } \\
&E_i \longrightarrow \rho (X_i) \\
{ \ } \\
&\mbox{vers l'image par $\rho$ d'un objet $X_i$ de ${\mathcal C}$.} \hfill \\
{\rm (B2)} &\mbox{Tout objet $E$ de ${\mathcal E}$ muni d'un morphisme} \hfill \\
{ \ } \\
&E \longrightarrow \rho (X) \times \rho (Y) \\
{ \ } \\
&\mbox{vers le produit des images par $\rho$ de deux objets $X,Y$ de ${\mathcal C}$ admet une famille \'epimorphique} \hfill \\
{ \ } \\
&E_i \longrightarrow E \, , \qquad i \in I \, , \\
{ \ } \\
&\mbox{dont chaque source s'inscrit dans un carr\'e commutatif} \hfill \\
{ \ } \\
&\xymatrix{
E_i \ar[d] \ar[r] &\rho (W_i) \ar[d] \\
E \ar[r] &\rho (X) \times \rho (Y)
} \\
{ \ } \\
&\mbox{o\`u la seconde fl\`eche verticale est l'image par $\rho$ d'une paire de morphismes de ${\mathcal C}$} \hfill \\
{ \ } \\
&X \longleftarrow W_i \longrightarrow Y \, . \\

{\rm (B3)} &\mbox{Tout objet de $E$ muni d'un morphisme} \hfill \\
{ \ } \\
&E \longrightarrow {\rm eg} \left( \rho (X) \rightrightarrows \rho (Y) \right) \\
{ \ } \\
&\mbox{vers l'\'egalisateur des images de $\rho$ de deux morphismes $X \raisebox{.7ex}{\xymatrix{\dar[r]^-{^{^{\mbox{\scriptsize$f$}}}}_-{g} &Y}}$ de ${\mathcal C}$ admet une famille} \hfill \\
&\mbox{\'epimorphique} \hfill \\
&E_i \longrightarrow E \, , \qquad i \in I \, , \\
{ \ } \\
&\mbox{dont chaque source s'inscrit dans un carr\'e commutatif} \hfill \\
{ \ } \\
&\xymatrix{
E_i \ar[d] \ar[r] &\rho (W_i) \ar[d] \\
E \ar[r] &{\rm eg} (\rho (X) \rightrightarrows \rho(Y))
} \\
{ \ } \\
&\mbox{o\`u la seconde fl\`eche verticale est l'image par $\rho$ d'un morphisme de ${\mathcal C}$} \hfill \\
{ \ } \\
&w_i : W_i \longrightarrow X \\
&\mbox{tel que $f \circ w_i = g \circ w_i$.} \hfill
\end{matrix} \right.
$}
\end{thm}

\bigskip

\begin{demo}

C'est imm\'ediat. \end{demo}

\bigskip

Par combinaison avec le corollaire \ref{corIV310}, on obtient aussi le cas particulier de ce th\'eor\`eme:

\begin{cor}\label{corIV3.12}

Soient ${\mathcal C}$ une cat\'egorie essentiellement petite qui a des limites finies arbitraires, $J$ une topologie sur ${\mathcal C}$, et
$$
\ell : {\mathcal C} \xhookrightarrow{ \ y \ } \widehat{\mathcal C} \xrightarrow{ \ j^* \ } \widehat{\mathcal C}_J
$$
le foncteur canonique associ\'e.

\smallskip

Alors, pour tout topos ${\mathcal E}$, le foncteur
$$
\begin{matrix}
[{\mathcal E} , \widehat{\mathcal C}_J] &\longrightarrow &[{\mathcal C}, {\mathcal E}] \, , \\
(f^* , f_*) &\longmapsto &f^* \!\circ \ell
\end{matrix}
$$
d\'efinit une \'equivalence de la cat\'egorie
$$
[{\mathcal E} , \widehat{\mathcal C}_J]_T
$$
des morphismes de topos de ${\mathcal E}$ dans $\widehat{\mathcal C}_J$ sur la sous-cat\'egorie pleine de
$$
[{\mathcal C},{\mathcal E}]
$$
constitu\'ee des foncteurs
$$
\rho : {\mathcal C} \longrightarrow {\mathcal E}
$$
qui

\medskip

\begin{enumerate}
\item[$\bullet$] ``pr\'eservent les recouvrements'' c'est-\`a-dire transforment toute famille $J$-couvrante de morphismes de~${\mathcal C}$
$$
(X_i \longrightarrow X)_{i \in I}
$$
en une famille \'epimorphique de morphismes de ${\mathcal E}$,
\item[$\bullet$] sont ``plats'' c'est-\`a-dire respectent les limites finies, ce pour quoi il suffit que soient satisfaites les trois conditions suivantes:
\end{enumerate}

\medskip

$
\left\{ \begin{matrix}
{\rm (B1')} &\mbox{Le morphisme canonique} \hfill \\
&\rho (1_{\mathcal C}) \longrightarrow 1_{\mathcal E} \\
{ \ } \\
&\mbox{est un \'epimorphisme de ${\mathcal E}$.} \hfill \\
{ \ } \\
{\rm (B2')} &\mbox{Pour tous objets $X,Y$ de ${\mathcal C}$, le morphisme canonique} \hfill \\
{ \ } \\
&\rho (X \times Y) \longrightarrow \rho (X) \times \rho (Y) \\
{ \ } \\
&\mbox{est un \'epimorphisme de ${\mathcal E}$.} \hfill \\
{ \ } \\
{\rm (B3')} &\mbox{Pour tous morphismes $X \rightrightarrows Y$ de ${\mathcal C}$, le morphisme canonique} \hfill \\
{ \ } \\
&\rho \left( {\rm eg} (X \rightrightarrows Y)\right) \longrightarrow {\rm eg} \left(\rho (X) \rightrightarrows \rho (Y)\right) \\
{ \ } \\
&\mbox{est un \'epimorphisme de ${\mathcal E}$.} \hfill
\end{matrix} \right.
$
\end{cor}

\bigskip

\begin{demo}

C'est une cons\'equence imm\'ediate du corollaire \ref{corIV34} et du corollaire \ref{corIV310}. \end{demo}

\section{Morphismes de sites}\label{sec44}

\subsection{La notion de morphisme de sites}\label{subsec441}

\medskip

Il est naturel de poser la d\'efinition suivante:

\begin{defn}\label{defIV41}

Soient $({\mathcal C} , J)$ et $({\mathcal D} , K)$ deux sites, avec les morphismes canoniques associ\'es
$$
\ell : {\mathcal C} \xrightarrow{ \ y \ } \widehat{\mathcal C} \xrightarrow{ \ j^* \ } \widehat{\mathcal C}_J
$$
et
$$
\ell : {\mathcal D} \xrightarrow{ \ y \ } \widehat{\mathcal D} \xrightarrow{ \ j^* \ } \widehat{\mathcal D}_K \, .
$$

Un foncteur
$$
\rho : {\mathcal C} \longrightarrow {\mathcal D}
$$
sera appel\'e un morphisme de sites s'il existe un morphisme de topos, n\'ecessairement unique \`a isomorphisme canonique pr\`es,
$$
(\rho_! , \rho^!) : \widehat{\mathcal D}_K \longrightarrow \widehat{\mathcal C}_J
$$
dont la composante d'image r\'eciproque
$$
\rho_! : \widehat{\mathcal C}_J \longrightarrow \widehat{\mathcal D}_K
$$
s'inscrit dans un carr\'e commutatif:
$$
\xymatrix{
{\mathcal C} \ar[d]_-{\ell} \ar[r]^-{\rho} &{\mathcal D} \ar[d]^-{\ell} \\
\widehat{\mathcal C}_J \ar[r]^-{\rho_!} &\widehat{\mathcal D}_K
}
$$
\end{defn}
\bigskip

\begin{remarksqed}
\begin{listeisansmarge}
\item S'il existe un morphisme de topos
$$
(\rho_! , \rho^!) : \widehat{\mathcal D}_K \longrightarrow \widehat{\mathcal C}_J
$$
qui rende commutatif le carr\'e de l'\'enonc\'e, le compos\'e
$$
\widehat{\mathcal C} \xrightarrow{ \ j^* \ } \widehat{\mathcal C}_J \xrightarrow{ \ \rho_! \ } \widehat{\mathcal D}_K
$$
respecte les colimites et prolonge le foncteur compos\'e
$$
{\mathcal C} \xrightarrow{ \ \rho \ } {\mathcal D} \xrightarrow{ \ \ell \ } \widehat{\mathcal D}_K \, .
$$
Il s'identifie donc au foncteur
$$
\widehat\rho : {\mathcal C} \longrightarrow \widehat{\mathcal D}_K
$$
de la proposition \ref{propIV31}, et le foncteur
$$
\rho_! : \widehat{\mathcal C}_J \longrightarrow \widehat{\mathcal D}_K
$$
est canoniquement isomorphe au foncteur
$$
\widehat\rho \circ j_* : \widehat{\mathcal C}_J \xhookrightarrow{ \ { \ } \ } \widehat{\mathcal C} \longrightarrow \widehat{\mathcal D}_K
$$
d\'eduit de $\widehat\rho$ par composition avec le foncteur de plongement pleinement fid\`ele
$$
j_* : \widehat{\mathcal C}_J \xhookrightarrow{ \ { \ } \ } \widehat{\mathcal C} \, .
$$
Enfin, l'adjoint \`a droite $\rho^!$ de $\rho_!$ est uniquement d\'etermin\'e \`a isomorphisme canonique pr\`es.

\smallskip

C'est pourquoi le morphisme de topos $(\rho_! , \rho^!)$ est uniquement d\'etermin\'e par $\rho$ \`a isomorphisme canonique pr\`es, s'il existe.

\medskip

\item Un tel morphisme de topos est not\'e $(\rho_! , \rho^!)$ pour signifier que sa composante d'image r\'eciproque est covariante par rapport \`a $\rho$ et sa composante d'image directe est contravariante. On utilise le signe ``!'' plut\^ot que le signe ``$*$'' pour signifier que c'est l'adjoint \`a gauche qui est covariant et l'adjoint \`a droite qui est contravariant.

\medskip

\item Si $({\mathcal C}_1 , J_1)$, $({\mathcal C}_2 , J_2)$ et $({\mathcal C}_3 , J_3)$ sont trois sites reli\'es par deux morphismes de sites au sens de la d\'efinition
$$
{\mathcal C}_1 \xrightarrow{ \ \rho \ } {\mathcal C}_2 \xrightarrow{ \ \eta \ } {\mathcal C}_3 \, ,
$$
alors le foncteur compos\'e
$$
\eta \circ \rho : {\mathcal C}_1 \longrightarrow {\mathcal C}_3
$$
est encore un morphisme de sites et admet pour morphisme de topos associ\'e le compos\'e des morphismes de topos associ\'es \`a $\rho$ et $\eta$
$$
((\eta \circ \rho)_! , (\eta \circ \rho)^!) = (\eta_! \circ \rho_! , \rho^! \circ \eta^!) \, .
$$

\item Pour tous sites $({\mathcal C},J)$ et $({\mathcal D},K)$, et toute paire de morphismes de sites ${\mathcal C} \!\!\raisebox{.7ex}{\xymatrix{\dar[r]^-{^{^{\mbox{\scriptsize$\rho$}}}}_-{\rho'} &{\mathcal D}}}$, toute transformation naturelle $\alpha : \rho \to \rho'$ induit un morphisme de $(\rho_! , \rho^!)$ vers $(\rho'_! , \rho'^!)$ dans la cat\'egorie $[\widehat{\mathcal D}_K , \widehat{\mathcal C}_J]_T$.

\medskip

\item Pour toute cat\'egorie essentiellement petite ${\mathcal C}$ munie de deux topologies $J \supseteq K$, le foncteur ${\rm id} : {\mathcal C} \to {\mathcal C}$ est un morphisme du site $({\mathcal C},J)$ vers le site $({\mathcal C},K)$ qui induit le morphisme de topos $\widehat{\mathcal C}_J \to \widehat{\mathcal C}_K$.

\medskip

\item Pour tout morphisme de topos
$$
(f^* , f_*) : {\mathcal E} \longrightarrow {\mathcal B}
$$
et toute repr\'esentation du topos de base ${\mathcal B}$ comme un topos de faisceaux
$$
 \widehat{\mathcal C}_J \xrightarrow{ \ \sim \ } {\mathcal B}
$$
sur un site $({\mathcal C},J)$, tout choix d'une sous-cat\'egorie pleine essentiellement petite ${\mathcal D}$ de ${\mathcal E}$ qui contient une famille s\'eparante d'objets et les images des objets de ${\mathcal C}$ par le foncteur compos\'e
$$
{\mathcal C} \xrightarrow{ \ \ell \ } \widehat{\mathcal C}_J \xrightarrow{ \ \sim \ } {\mathcal B} \xrightarrow{ \ f^* \ } {\mathcal E}
$$
permet de repr\'esenter ${\mathcal E}$ comme le topos des faisceaux
$$
\widehat{\mathcal D}_K \xrightarrow{ \ \sim \ } {\mathcal E} 
$$
sur la cat\'egorie ${\mathcal D}$ munie de la topologie $K$ induite par la notion de famille \'epimorphique de ${\mathcal E}$ puis de voir le morphisme de topos $(f^* , f_*) : {\mathcal E} \to {\mathcal B}$ comme associ\'e au morphisme de sites
$$
{\mathcal C} \longrightarrow {\mathcal D}
$$
d\'eduit du foncteur ${\mathcal C} \to {\mathcal E}$.

\smallskip

En effet, l'\'equivalence $\widehat{\mathcal D}_K \xrightarrow{ \ \sim \ } {\mathcal E}$ r\'esulte du th\'eor\`eme \ref{thmIII81} puisque la sous-cat\'egorie pleine ${\mathcal D}$ de ${\mathcal E}$ contient une famille s\'eparante d'objets.

\smallskip

Comme ${\mathcal D}$ contient aussi les images des objets de ${\mathcal C}$ par le foncteur compos\'e ${\mathcal C} \xrightarrow{ \ \ell \ } \widehat{\mathcal C}_J \xrightarrow{ \ \sim \ } {\mathcal B} \xrightarrow{ \ f^* \ } {\mathcal E}$, celui-ci se factorise en un foncteur
$$
\rho : {\mathcal C} \longrightarrow {\mathcal D}
$$ 
qui est un morphisme de sites puisque le carr\'e
$$
\xymatrix{
{\mathcal C} \ar[d]_-{\ell} \ar[r]^-{\rho} &{\mathcal D} \ar[d]^-{\ell} \\
\widehat{\mathcal C}_J \ar[d]^-{\wr} &\widehat{\mathcal D}_K \ar[d]^-{\wr} \\
{\mathcal B} \ar[r]^-{f^*} &{\mathcal E} 
}
$$
est commutatif. 

\medskip

\item Tout morphisme entre deux topos de faisceaux sur des sites $({\mathcal D},K)$ et $({\mathcal C},J)$
$$
(f^* , f_*) : \widehat{\mathcal D}_K \longrightarrow \widehat{\mathcal C}_J
$$
se repr\'esente comme associ\'e \`a des morphismes du site $({\mathcal C},J)$ dans des sites $({\mathcal D}',K')$
$$
\rho : {\mathcal C} \longrightarrow {\mathcal D}'
$$
o\`u ${\mathcal D}'$ contient ${\mathcal D}$ comme sous-cat\'egorie pleine et dense pour la topologie $K'$ et o\`u $K'$ induit la topologie $K$ de ${\mathcal D}$, si bien que d'apr\`es le lemme de comparaison de Grothendieck la restriction des faisceaux de ${\mathcal D}'$ \`a ${\mathcal D}$ d\'efinit une \'equivalence
$$
\widehat{\mathcal D}'_{K'} \xrightarrow{ \ \sim \ } \widehat{\mathcal D}_K \, .
$$

En effet, il suffit de relever le foncteur compos\'e
$$
{\mathcal C} \xrightarrow{ \ \ell \ } \widehat{\mathcal C}_J \xrightarrow{ \ f^* \ } \widehat{\mathcal D}_K
$$
en n'importe quel foncteur $\rho : {\mathcal C} \to \widehat{\mathcal D}$ qui rende le carr\'e
$$
\xymatrix{
{\mathcal C} \ar[d]_-{\ell} \ar[r]^-{\rho} &\widehat{\mathcal D} \ar[d]^-{j^*} \\
\widehat{\mathcal C}_J \ar[r]^-{f^*} &\widehat{\mathcal D}_K
}
$$
commutatif \`a isomorphisme pr\`es, puis de choisir pour ${\mathcal D}'$ n'importe quelle sous-cat\'egorie pleine et essentiellement petite de $\widehat{\mathcal D}$ qui contienne l'image de $\rho$ et celle de $y : {\mathcal D} \to \widehat{\mathcal D}$.

\smallskip

On prend pour $K'$ la topologie de ${\mathcal D}'$ pour laquelle une famille de morphismes de m\^eme but est couvrante quand son image par $j^* : \widehat{\mathcal D} \to \widehat{\mathcal D}_K$ est globalement \'epimorphique. 
\end{listeisansmarge}
\end{remarksqed} 


Les composantes du morphisme de topos $(\rho_! , \rho^!)$ d\'efini par un morphisme de sites $\rho$ sont identifi\'ees par le lemme suivant:

\begin{lem}\label{lemIV42}

Soient $({\mathcal C} , J)$ et $({\mathcal D} , K)$ deux sites munis des foncteurs canoniques
$$
\ell : {\mathcal C} \xhookrightarrow{ \ y \ } \widehat{\mathcal C} \xrightarrow{ \ j^* \ } \widehat{\mathcal C}_J
$$
et
$$
\ell : {\mathcal D} \xhookrightarrow{ \ y \ } \widehat{\mathcal D} \xrightarrow{ \ j^* \ } \widehat{\mathcal D}_K \, .
$$

Alors:

\begin{listeimarge}

\item Un foncteur
$$
\rho : {\mathcal C} \longrightarrow {\mathcal D}
$$
s'inscrit dans un carr\'e commutatif
$$
\xymatrix{
{\mathcal C} \ar[d]_-{\ell} \ar[r]^-{\rho} &{\mathcal D} \ar[d]^-{\ell} \\
\widehat{\mathcal C}_J \ar[r]^-{\rho_!} &\widehat{\mathcal D}_K
}
$$
tel que $\rho_!$ admette un adjoint \`a droite $\rho^!$, si et seulement si le foncteur de composition avec $\rho$
$$
\begin{matrix}
\rho^* : \widehat{\mathcal D} = [{\mathcal D}^{\rm op} , {\rm Ens}] &\longrightarrow &[{\mathcal C}^{\rm op} , {\rm Ens}] = \widehat{\mathcal C} \\
\hfill G &\longmapsto &G \circ \rho \hfill
\end{matrix}
$$
transforme les faisceaux sur ${\mathcal D}$ en faisceaux sur ${\mathcal C}$.

\medskip

\item Dans ce cas, le foncteur
$$
\rho^! : \widehat{\mathcal D}_K \longrightarrow \widehat{\mathcal C}_J
$$
est la restriction \`a $\widehat{\mathcal D}_K$ du foncteur $\rho^*$ de composition avec $\rho$, et son adjoint \`a gauche le foncteur
$$
\rho_! : \widehat{\mathcal C}_J \longrightarrow \widehat{\mathcal D}_K
$$
associe \`a tout faisceau $F$ sur ${\mathcal C}$ le faisceau $\rho_! F$ sur ${\mathcal D}$ obtenu comme le faisceautis\'e du pr\'efaisceau
$$
\begin{matrix}
{\mathcal D}^{\rm op} &\longrightarrow &{\rm Ens} \, , \hfill \\
\hfill Y &\longmapsto &\displaystyle\varinjlim_{(Y \backslash {\mathcal C}) \ni X} F(X) 
\end{matrix}
$$
o\`u la colimite est calcul\'ee sur la cat\'egorie $Y \backslash {\mathcal C}$ des objets $X$ de ${\mathcal C}$ munis d'un morphisme $Y \to \rho (X)$.
\end{listeimarge}
\end{lem}

\begin{remark}

Le lemme s'applique en particulier aux morphismes de topos
$$
(\rho_! , \rho^!) : \widehat{\mathcal D}_K \longrightarrow \widehat{\mathcal C}_J
$$
induits par un morphisme de sites $\rho : {\mathcal C} \to {\mathcal D}$, mais il est plus g\'en\'eral car on ne demande pas que le foncteur $\rho_! : \widehat{\mathcal C}_J \to \widehat{\mathcal D}_K$ respecte les limites finies.
\end{remark}
\bigskip

\begin{demo}

Soit $\rho_! : \widehat{\mathcal C} \to \widehat{\mathcal D}$ l'unique prolongement pr\'eservant les colimites du compos\'e ${\mathcal C} \xrightarrow{ \ \rho \ } {\mathcal D} \xhookrightarrow{ \ y \ } \widehat{\mathcal D}$.

\smallskip

On pr\'etend qu'il est adjoint \`a gauche du foncteur $\rho^*$ de composition avec $\rho$
$$
\begin{matrix}
\rho^* : \widehat{\mathcal D} &\longrightarrow &\widehat{\mathcal C} \, , \hfill \\
\hfill G &\longmapsto &G \circ \rho \, ,
\end{matrix}
$$
autrement dit qu'il est le foncteur d'extension de Kan \`a gauche le long de $\rho$.

\smallskip

Cela r\'esulte de ce que pour tout pr\'efaisceau $F$ sur ${\mathcal C}$ les formules
$$
F = \varinjlim_{(X,x) \in \int\!\!F} y(X)
$$
et
$$
\rho_! F = \varinjlim_{(X,x) \in \int\!\!F} y \circ \rho (X)
$$
impliquent pour tout pr\'efaisceau $G$ sur ${\mathcal D}$
$$
{\rm Hom} (F,G \circ \rho) = \varprojlim_{(X,x) \in \int\!\!F} {\rm Hom} (y(X) , G \circ \rho) = \varprojlim_{(X,x) \in \int\!\!F} G \circ \rho (X)
$$
et
$$
{\rm Hom} (\rho_! F,G) = \varprojlim_{(X,x) \in \int\!\!F} {\rm Hom} (y \circ \rho (X) , G) = \varprojlim_{(X,x) \in \int\!\!F} G (\rho (X)) \, .
$$

Ce fait \'etant \'etabli, nous pouvons d\'emontrer (i) et (ii):

\medskip

(i) Un foncteur $\rho_! : \widehat{\mathcal C}_J \to \widehat{\mathcal D}_K$ qui prolonge $\rho : {\mathcal C} \to {\mathcal D}$ admet un adjoint \`a droite $\rho^!$ si et seulement si il pr\'eserve les colimites c'est-\`a-dire s'inscrit dans un carr\'e commutatif \`a isomorphisme canonique pr\`es:
$$
\xymatrix{
\widehat{\mathcal C} \ar[d]_-{j^*} \ar[r]^-{\rho_!} &\widehat{\mathcal D} \ar[d]^-{j^*} \\
\widehat{\mathcal C}_J \ar[r]^-{\rho_!} &\widehat{\mathcal D}_K
}
$$

En passant aux adjoints \`a droite, cela \'equivaut \`a demander que le foncteur compos\'e
$$
\widehat{\mathcal D}_K \xhookrightarrow{ \ j_* \ } \widehat{\mathcal D} \xrightarrow{ \ \rho^* \ } \widehat{\mathcal C}
$$
se factorise \`a travers $j_* : \widehat{\mathcal C}_J \hookrightarrow \widehat{\mathcal C}$ c'est-\`a-dire s'inscrit dans un carr\'e commutatif:
$$
\xymatrix{
\widehat{\mathcal D} \ar[r]^-{\rho^*} &\widehat{\mathcal C} \\
\!\!\widehat{\mathcal D}_K\ar@{^{(}->}[u]^-{j_*} \ar[r]^-{\rho^!} &\ \widehat{\mathcal C}_J \ar@{^{(}->}[u]_{j_*}
}
$$

C'est l'\'equivalence de (i).

\smallskip

La premi\`ere partie de (ii) reformule la propri\'et\'e de commutativit\'e de ce carr\'e.

\smallskip

La seconde partie de (ii) r\'esulte de la commutativit\'e du carr\'e
$$
\xymatrix{
\widehat{\mathcal C} \ar[d]_-{j^*} \ar[r]^-{\rho_!} &\widehat{\mathcal D} \ar[d]^-{j^*} \\
\widehat{\mathcal C}_J \ar[r]^-{\rho_!} &\widehat{\mathcal D}_K
}
$$
en rappelant la formule de la proposition \ref{propI106} pour le foncteur d'extension de Kan \`a gauche le long de $\rho$
$$
\rho_! : \widehat{\mathcal C} \longrightarrow \widehat{\mathcal D}
$$
d\'efini comme l'adjoint \`a gauche de $\rho^* : \widehat{\mathcal D} \to \widehat{\mathcal C}$. 

\end{demo}

\subsection{Caract\'erisation des morphismes de sites}\label{subsec442}

\medskip

On d\'eduit du th\'eor\`eme \ref{thmIV311} et de la proposition \ref{propIV39} le crit\`ere suivant pour qu'un foncteur soit un morphisme de sites:

\begin{thm}[Caramello (Définition 3.2 et Remarque 3.3 de \cite{Denseness})]\label{thmIV43}

Soient $({\mathcal C} , J)$ et $({\mathcal D} , K)$ deux sites.

\smallskip

Un foncteur
$$
\rho : {\mathcal C} \longrightarrow {\mathcal D}
$$
est un morphisme de sites, c'est-\`a-dire d\'efinit un morphisme de topos
$$
(\rho_! , \rho^!) : \widehat{\mathcal D}_K \longrightarrow \widehat{\mathcal C}_J
$$
dont la composante d'image r\'eciproque
$$
\rho_! : \widehat{\mathcal C}_J \longrightarrow \widehat{\mathcal D}_K
$$
prolonge $\rho : {\mathcal C} \to {\mathcal D}$, et dont la composante d'image directe
$$
\rho^! : \widehat{\mathcal D}_K \longrightarrow \widehat{\mathcal C}_J
$$
est d\'efinie par la composition avec $\rho$
$$
G \longmapsto G \circ \rho \, ,
$$ 
si et seulement si:

\medskip

\begin{enumerate}
\item[$\bullet$] le foncteur $\rho$ transforme toute famille $J$-couvrante de morphismes de ${\mathcal C}$
$$
(X_i \longrightarrow X)_{i \in I}
$$
en une famille $K$-couvrante de morphismes de ${\mathcal D}$
$$
(\rho (X_i) \longrightarrow \rho (X))_{i \in I} \, ,
$$
\item[$\bullet$] le compos\'e ${\mathcal C} \xrightarrow{ \ \rho \ } {\mathcal D} \longrightarrow \widehat{\mathcal D}_K$ est ``plat'', ce qui signifie que $\rho$ satisfait les trois conditions:
\end{enumerate}

\medskip

$
\left\lmoustache \begin{matrix}
{\rm (C1)} &\mbox{Pour tout objet $D$ de ${\mathcal D}$ existe une famille $K$-couvrante} \hfill \\
{ \ } \\
&D_i \longrightarrow D \, , \quad i \in I \, , \\
{ \ } \\
&\mbox{tel que chaque $D_i$ admette un morphisme de ${\mathcal D}$} \hfill \\
{ \ } \\
&D_i \longrightarrow \rho (X_i) \\
{ \ } \\
&\mbox{vers l'image par $\rho$ d'un objet $X_i$ de ${\mathcal C}$.} \hfill \\
{ \ } \\
{\rm (C2)} &\mbox{Pour tous objets $X,Y$ de ${\mathcal C}$ et tous morphismes de ${\mathcal D}$} \hfill \\
{ \ } \\
&\rho (X) \longleftarrow D \longrightarrow \rho (Y) \\
{ \ } \\
&\mbox{existent une famille $K$-couvrante de morphismes de ${\mathcal D}$} \hfill \\
{ \ } \\
&D_i \longrightarrow D \, , \quad i \in I \, , \\
&\mbox{et pour chaque indice $i$ deux morphismes de ${\mathcal C}$} \hfill \\
{ \ } \\
&X \longleftarrow W_i \longrightarrow Y \\
&\mbox{et un morphisme de ${\mathcal D}$} \hfill \\
&D_i \longrightarrow \rho (W_i) \\
{ \ } \\
&\mbox{s'inscrivant dans un carr\'e commutatif de ${\mathcal D}$:} \hfill \\
{ \ } \\
&\xymatrix{
&D \ar[ld] \ar[rd] \\
\rho(X) &D_i \ar[u] \ar[d] &\rho (Y) \\
&\rho (W_i) \ar[lu] \ar[ru]
} \\
{ \ } \\
{\rm (C3)} &\mbox{Pour tous morphismes $X \raisebox{.7ex}{\xymatrix{\dar[r]^-{^{^{\mbox{\scriptsize$f$}}}}_-{g} &Y}}$ de ${\mathcal C}$ et tout morphisme de ${\mathcal D}$} \hfill \\
&D \xrightarrow{ \ d \ } \rho (X) \\
{ \ } \\
&\mbox{tel que $\rho (f) \circ d = \rho (g) \circ d$, il existe une famille $K$-couvrante de morphismes de ${\mathcal D}$} \hfill \\
{ \ } \\
&D_i \longrightarrow D \, , \quad i \in I \, , \\
{ \ } \\
&\mbox{et, pour chaque indice $i$, un morphisme de ${\mathcal C}$} \hfill \\
{ \ } \\
&w_i : W_i \longrightarrow X \\
{ \ } \\
&\mbox{tel que $f \circ w_i = g \circ w_i$ et un morphisme de ${\mathcal D}$} \hfill \\
{ \ } \\
&D_i \longrightarrow \rho (W_i) \\
\end{matrix} \right.
$

$
\left\rmoustache \begin{matrix}
&\mbox{s'inscrivant dans un carr\'e commutatif:} \hfill \\
{ \ } \\
&\xymatrix{
D_i \ar[d] \ar[r] &\rho (W_i) \ar[d] \\
D \ar[r] &\rho (X)
}
\end{matrix} \right.
$
\end{thm}

\begin{demosansqed}

On sait d'apr\`es le lemme \ref{lemIII53} qu'une famille de morphismes de ${\mathcal D}$
$$
(D_i \longrightarrow D)_{i \in I}
$$
est $K$-couvrante si et seulement si son image par le foncteur canonique $\ell : {\mathcal D} \to \widehat{\mathcal D}_K$ est une famille \'epimorphique.

\smallskip

Le th\'eor\`eme est alors une simple traduction du th\'eor\`eme \ref{thmIV311}, compte tenu de l'\'equivalence des conditions (B) et (C) de la proposition \ref{propIV39} et du lemme suivant:
\end{demosansqed}

\begin{lem}\label{lemIV44}

Soit $({\mathcal D}, K)$ un site.

\smallskip

Soit $\ell : {\mathcal D} \xhookrightarrow{ \ y \ } \widehat{\mathcal D} \xrightarrow{ \ j^* \ } \widehat{\mathcal D}_K$ le foncteur canonique associ\'e.

\smallskip

Alors on a pour tous objets $D,D'$ de ${\mathcal D}$:

\begin{listeimarge}

\item Pour tout morphisme de $\widehat{\mathcal D}_K$
$$
\ell (D) \longrightarrow \ell (D')
$$
existe une famille $K$-couvrante de morphismes de ${\mathcal D}$
$$
D_i \longrightarrow D \, , \qquad i \in I \, ,
$$
telle que chaque morphisme compos\'e
$$
\ell (D_i) \longrightarrow \ell (D) \longrightarrow \ell (D')
$$
soit l'image par $\ell$ d'un morphisme de ${\mathcal D}$.

\medskip

\item Si deux morphismes de ${\mathcal D}$
$$
D \rightrightarrows D'
$$
induisent via $\ell$ le m\^eme morphisme de $\widehat{\mathcal D}_K$
$$
\ell (D) \longrightarrow \ell (D') \, ,
$$
il existe une famille $K$-couvrante de morphismes de ${\mathcal D}$
$$
D_i \longrightarrow D \, , \qquad i \in I \, ,
$$
telle que pour chaque indice $i$ les deux compos\'es
$$
D_i \to D \rightrightarrows D'
$$
soient \'egaux.
\end{listeimarge}
\end{lem}

\begin{demo}

Consid\'erons le pr\'efaisceau
$$
P = y(D') \, ,
$$
sa faisceautisation
$$
j^* P = \ell (D')
$$
et la m\^eme vue comme un pr\'efaisceau
$$
F = j_* \circ j^* P = j_* \circ \ell (D') \, .
$$

Pour tout objet $D$ de ${\mathcal D}$, on a
$$
{\rm Hom} (D,D') = P(D)
$$
et
$$
{\rm Hom} (\ell (D),\ell (D')) = {\rm Hom} (y(D) , j_* \circ \ell (D')) = F(D) \, .
$$

Le morphisme canonique de $\widehat{\mathcal D}$
$$
P \longrightarrow j_* \circ j^* P = F
$$
induit un isomorphisme de $\widehat{\mathcal D}_K$
$$
j^* P \xrightarrow{ \ \sim \ } j^* F \, .
$$

Alors les parties (i) et (ii) du lemme traduisent les crit\`eres (i) et (ii) du lemme \ref{lemII56}.

\smallskip

Cela prouve le lemme et donc aussi le th\'eor\`eme \ref{thmIV43}. \end{demo}

\medskip

On d\'eduit du th\'eor\`eme \ref{thmIV43} et du corollaire \ref{corIV310} le cas particulier suivant:

\begin{cor}\label{corIV45}

Soient $({\mathcal C} , J)$ et $({\mathcal D},K)$ deux sites.

\smallskip

Supposons que la cat\'egorie ${\mathcal C}$ poss\`ede des limites finies arbitraires.

\smallskip

Alors un foncteur
$$
\rho : {\mathcal C} \longrightarrow {\mathcal D}
$$
est un morphisme de sites si et seulement si:

\medskip

\begin{enumerate}
\item[$\bullet$] il transforme toute famille $J$-couvrante de morphismes de ${\mathcal C}$ en une famille $K$-couvrante de morphismes de ${\mathcal D}$,
\item[$\bullet$] le compos\'e ${\mathcal C} \xrightarrow{ \ \rho \ } {\mathcal D} \longrightarrow \widehat{\mathcal D}_K$ est ``plat'' c'est-\`a-dire respecte les limites finies, ce pour quoi il suffit que soient satisfaites les trois conditions suivantes:
\end{enumerate}

\medskip

$
\left\{\begin{matrix}
{\rm (C1')} &\mbox{Tout objet $D$ de ${\mathcal D}$ admet une famille $K$-couvrante} \hfill \\
{ \ } \\
&D_i \longrightarrow D \, , \quad i \in I \, , \\
{ \ } \\
&\mbox{telle que chaque $D_i$ poss\`ede un morphisme} \hfill \\
{ \ } \\
&D_i \longrightarrow \rho (1_{\mathcal C}) \, . \\
{ \ } \\
{\rm (C2')} &\mbox{Pour tous objets $X,Y$ de ${\mathcal C}$ et tous morphismes de ${\mathcal D}$} \hfill \\
{ \ } \\
&\rho (X) \longleftarrow D \longrightarrow \rho (Y) \, , \\
{ \ } \\
&\mbox{il existe une famille $K$-couvrante de morphismes de ${\mathcal D}$} \hfill \\
{ \ } \\
&D_i \longrightarrow D \, , \quad i \in I \, , \\
{ \ } \\
&\mbox{et pour chaque indice $i$ un morphisme} \hfill \\
{ \ } \\
&D_i \longrightarrow \rho (X \times Y) \\
{ \ } \\
&\mbox{s'inscrivant dans un diagramme commutatif de ${\mathcal D}$:} \hfill \\
{ \ } \\
&\xymatrix{
&D \ar[ld] \ar[rd] \\
\rho(X) &D_i \ar[u] \ar[d] &\rho(Y) \\
&\rho(X \times Y) \ar[lu] \ar[ru]
} \\
{\rm (C3')} &\mbox{Pour tous morphismes $X \raisebox{.7ex}{\xymatrix{\dar[r]^-{^{^{\mbox{\scriptsize$f$}}}}_-{g} &Y}}$ de ${\mathcal C}$ et tout morphisme de ${\mathcal D}$} \hfill \\
&D \xrightarrow{ \ d \ } \rho (X) \\
{ \ } \\
&\mbox{tel que $\rho (f) \circ d = \rho (g) \circ d$, il existe une famille $K$-couvrante de morphismes de ${\mathcal D}$} \hfill \\
{ \ } \\
&D_i \longrightarrow D \, , \quad i \in I \, , \\
{ \ } \\
&\mbox{s'inscrivant, pour chaque indice $i$, dans un carr\'e commutatif de ${\mathcal D}$:} \hfill \\
{ \ } \\
&\xymatrix{
D_i \ar[d] \ar[r] &\rho ({\rm eg} (X \rightrightarrows Y)) \ar[d] \\
D \ar[r] &\rho (X)
}
\end{matrix} \right.
$
\end{cor}

\begin{demo}

On sait d'apr\`es le corollaire \ref{corIV310} que, la cat\'egorie ${\mathcal C}$ poss\'edant des limites finies arbitraires, le foncteur
$$
{\mathcal C} \xrightarrow{ \ \rho \ } {\mathcal D} \longrightarrow \widehat{\mathcal D}_K
$$
est plat si et seulement si il pr\'eserve les limites finies.

\smallskip

Le reste du corollaire se ram\`ene au th\'eor\`eme \ref{thmIV43} puisque la cat\'egorie ${\mathcal C}$ poss\`ede un objet terminal $1_{\mathcal C}$, que toute paire d'objets $X$ et $Y$ y poss\`ede un produit $X \times Y$ et que toute paire de morphismes $X \rightrightarrows Y$ y poss\`ede un \'egalisateur ${\rm eg} (X \rightrightarrows Y)$.

\end{demo}

\subsection{Exemples de morphismes de sites}\label{subsec443}

\medskip

Les treillis au sens de la d\'efinition \ref{defII23} forment une cat\'egorie dont les morphismes sont d\'efinis de la mani\`ere suivante:

\begin{defn}\label{defIV46}

Un morphisme d'un treillis $O$ dans un treillis $O'$ est une application
$$
\rho : O \longrightarrow O'
$$
telle que:

\begin{enumerate}
\item[$\bullet$] pour toute famille $(v_i)_{i \in I}$ d'\'el\'ements de $O$, on a
$$
\rho \left(\bigvee_{i \in I} u_i \right) = \bigvee_{i \in I} \rho (u_i) \, ,
$$
\item[$\bullet$] pour toute famille finie $u_1 , \ldots , u_n$ d'\'el\'ements de $O$, on a
$$
\rho (u_1 \wedge \cdots \wedge u_n) = \rho (u_1) \wedge \cdots \wedge \rho (u_n) \, .
$$
\end{enumerate}
\end{defn}

\begin{remarkqed}

En particulier, un morphisme de treillis $O \to O'$ respecte la relation d'ordre au sens que pour tous \'el\'ements $u,v$ de $O$, on a
$$
u \leq v \Rightarrow \rho (u) \leq \rho (v) \, .
$$
\end{remarkqed}

\medskip

Tout morphisme de treillis munis de leur topologie canonique est un morphisme de sites:

\begin{prop}\label{propIV47}

Soient $O$ et $O'$ deux treillis vus comme des cat\'egories munies de leur topologie canonique.

\smallskip

Alors tout morphisme de treillis
$$
\rho : O \longrightarrow O'
$$
est un morphismes de sites.

\smallskip

Il d\'efinit un morphisme entre les topos associ\'es
$$
(\rho_! , \rho^!) : {\mathcal E}_{O'} \longrightarrow O_{\mathcal E}
$$ 
dont la composante d'image directe est le foncteur de composition avec $\rho$
$$
\begin{matrix}
\rho^! : {\mathcal E}_{O'} &\longrightarrow &{\mathcal E}_O \, , \hfill \\
\hfill G &\longmapsto &G \circ \rho
\end{matrix}
$$
et dont le foncteur d'image r\'eciproque
$$
\rho_! : {\mathcal E}_O \longrightarrow {\mathcal E}_{O'}
$$
associe \`a tout faisceau $F$ sur $O$ la faisceautisation du pr\'efaisceau sur $O'$
$$
O' \ni u' \longmapsto \varinjlim_{u \in O \atop u' \leq \rho (u)} F(u) \, .
$$
\end{prop}

\begin{remark}

En particulier, toute application continue entre espaces topologiques
$$
f : X \longrightarrow Y
$$
d\'efinit un morphisme de treillis
$$
\rho = f^{-1} : O(Y) \longrightarrow O(X) \, .
$$

Le morphisme de topos associ\'e
$$
(\rho_! , \rho^!) : {\mathcal E}_X \longrightarrow {\mathcal E}_Y
$$
n'est autre que celui d\'ej\`a associ\'e \`a $f$
$$
(f^* , f_*) : {\mathcal E}_X \longrightarrow {\mathcal E}_Y \, .
$$
\end{remark}

\begin{demo}

Un treillis vu comme une cat\'egorie admet des limites finies arbitraires: tout diagramme $u_{\bullet}$ d'un treillis index\'e par un carquois fini $D$ admet pour limite $\underset{d \in {\rm Ob} (D)}{\wedge} \, u(d)$.

\smallskip

D'autre part, une famille de morphismes $(u_i \leq u)_{i \in I}$ est couvrante si et seulement si $\underset{i \in I}{\bigvee} \, u_i = u$.

\smallskip

Ainsi, tout morphisme de treillis
$$
\rho : O' \longrightarrow O
$$
respecte les limites finies et transforme les familles couvrantes en familles couvrantes.

\smallskip

D'apr\`es le corollaire \ref{corIV45}, c'est un morphisme de sites et il d\'efinit un morphisme de topos
$$
(\rho_! , \rho^!) : {\mathcal E}_O \longrightarrow {\mathcal E}_{O'} \, .
$$

Les formules pour ses composantes $\rho_!$ et $\rho^!$ sont donn\'ees par le lemme \ref{lemIV42}. 

\end{demo}

\bigskip

Les foncteurs entre cat\'egories essentiellement petites munies chacune d'une ``classe g\'eom\'etrique de morphismes'' et d'une ``notion g\'eom\'etrique de recouvrement'' au sens de la d\'efinition \ref{defII31}, et qui respectent ces structures, fournissent une autre classe importante de morphismes de sites:

\bigskip

\begin{prop}\label{propIV48}
Soient ${\mathcal C}$ et ${\mathcal C}'$ deux cat\'egories essentiellement petites munies de ``classes g\'eom\'etriques de morphismes'' ${\mathcal M}$ et ${\mathcal M}'$ et de ``notions g\'eom\'etriques de recouvrement'' (R) et (R').

\smallskip

Supposons que ${\mathcal C}$ poss\`ede des limites finies arbitraires.

\smallskip

Soit un foncteur
$$
\rho : {\mathcal C} \longrightarrow {\mathcal C}'
$$
qui

\medskip

$
\left\{\begin{matrix}
\bullet &\mbox{respecte les limites finies,} \hfill \\
{ \ } \\
\bullet &\mbox{transforme les morphismes \'el\'ements de ${\mathcal M}$ en des morphismes \'el\'ements de ${\mathcal M}'$,} \hfill \\
{ \ } \\
\bullet &\mbox{transforme les familles d'\'el\'ements de ${\mathcal M}$} \hfill \\
&(U_i \longrightarrow X)_{i \in I} \\
{ \ } \\
&\mbox{qui poss\`edent la propri\'et\'e (R) en des familles d'\'el\'ements de ${\mathcal M}'$ qui poss\`edent la propri\'et\'e (R$'$).} \hfill
\end{matrix} \right.
$

\bigskip

Alors:

\begin{listeimarge}

\item Pour tout objet $X$ de ${\mathcal C}$, le foncteur $\rho$ d\'efinit un morphisme de sites

\medskip

$
\left\{\begin{matrix}
\bullet &\mbox{du ``gros site'' $({\mathcal C} / X , J_R)$ vers le ``gros site'' $({\mathcal C}' / \rho (X) , J_{R'})$,} \hfill \\
{ \ } \\
\bullet &\mbox{du ``petit site'' $(({\mathcal C}/X)_{\mathcal M} , J_R)$ vers le ``petit site'' $(({\mathcal C}'/\rho (X))_{{\mathcal M}'} , J_{R'})$} \hfill
\end{matrix} \right.
$

\medskip

\noindent et donc un morphisme en sens inverse entre les topos associ\'es
$$
(\rho_! , \rho^!) : [\widehat{{\mathcal C}' / \rho(X)}]_{J_{R'}} \longrightarrow [\widehat{{\mathcal C}/X}]_{J_R}
$$
ou
$$
(\rho_! , \rho^!) : [\widehat{({\mathcal C}' / \rho(X))}_{{\mathcal M}'}]_{J_{R'}} \longrightarrow [\widehat{({\mathcal C}/X)_{\mathcal M}}]_{J_R} \, .
$$

\item Sa composante d'image directe $\rho^!$ est d\'efinie par la composition avec $\rho$.

\medskip

\item Sa composante d'image r\'eciproque associe \`a tout faisceau $F$ sur ${\mathcal C}/X$ [resp. $({\mathcal C}/X)_{\mathcal M}$] le faisceautis\'e du pr\'efaisceau sur ${\mathcal C}'/\rho (X)$ [resp. $({\mathcal C}'/\rho(X))_{{\mathcal M}'}$] d\'efini par la formule
$$
\left( U' \xrightarrow{ \ u' \ } \rho (X)\right) \longmapsto \varinjlim_{U \to X \atop \begin{matrix} \mbox{\tiny$U'$} \to \mbox{\tiny$\rho(U)$} \\ \searrow \ \ \swarrow \\ \mbox{\tiny${ \ }^{\rho (X)}$} \end{matrix}} F(U \to X)
$$
o\`u la colimite est calcul\'ee sur la cat\'egorie des objets
$$
U \xrightarrow{ \ u \ } X \quad \mbox{de} \quad {\mathcal C}/X \qquad \mbox{[ resp.} \quad ({\mathcal C}/X)_{\mathcal M} \ \mbox{]}
$$
munis d'une factorisation
$$
\xymatrix{
U' \ar[rd]_-{u'} \ar[rr] &&\rho(U) \ar[ld]^-{\rho (u)} \\
&\rho (X)
}
$$
dans ${\mathcal C}' / \rho (X)$ [resp. $({\mathcal C}'/\rho(X))_{{\mathcal M}'}$].

\end{listeimarge}
\end{prop}

\bigskip

\begin{remarks}
\begin{listeisansmarge}
\item En particulier, pour tout objet $X$ de ${\mathcal C}$, le foncteur de plongement
$$
({\mathcal C} / X)_{\mathcal M} \xhookrightarrow{ \ { \ } \ } {\mathcal C}/X
$$
v\'erifie les conditions de la proposition.

\smallskip

C'est un morphisme de sites et il d\'efinit un morphisme de topos
$$
[\widehat{{\mathcal C}/X}]_{J_R} \longrightarrow [\widehat{({\mathcal C}/X)_{\mathcal M}}]_{J_R}
$$
dont la composante d'image directe est le foncteur de restriction \`a $({\mathcal C}/X)_{\mathcal M}$ des faisceaux sur ${\mathcal C}/X$.

\medskip

\item La proposition s'applique en particulier aussi \`a tout foncteur de changement de base
$$
\begin{matrix}
\rho_f &: &\hfill {\mathcal C}/X &\longrightarrow &{\mathcal C}/X' \, , \hfill \\
&&(U \to X) &\longmapsto &(U \times_X X' \to X')
\end{matrix}
$$
associ\'e \`a n'importe quel morphisme
$$
f : X' \longrightarrow X
$$
dans une cat\'egorie essentiellement petite ${\mathcal C}$ qui poss\`ede des produits fibr\'es et qui est munie d'une classe g\'eom\'etrique ${\mathcal M}$ de morphismes et d'une notion g\'eom\'etrique de recouvrement (R).

\smallskip

Le foncteur d'image directe associ\'e $f_*$ associe \`a tout faisceau $G$ sur ${\mathcal C}/X'$ [resp. $({\mathcal C}/X')_{\mathcal M}$] le faisceau $f_* G$ sur ${\mathcal C}/X$ [resp. $({\mathcal C}/X)_{\mathcal M}$] d\'efini comme le compos\'e
$$
(U \to X) \longmapsto G(U \times_X X' \to X') \, .
$$

Son adjoint \`a gauche le foncteur d'image r\'eciproque $f^*$ associe \`a tout faisceau $F$ sur ${\mathcal C}/X$ [resp. $({\mathcal C}/X)_{\mathcal M}$] le faisceau $f^* F$ sur ${\mathcal C}/X'$ [resp. $({\mathcal C}/X')_{\mathcal M}$] qui est le faisceautis\'e du pr\'efaisceau
$$
(U' \to X') \longmapsto \varinjlim_{\begin{pmatrix} U' \to U \\ \downarrow \qquad \downarrow \\ X' \xrightarrow{f} X \end{pmatrix}} F(U \to X) \, .
$$

\item La remarque pr\'ec\'edente s'applique par exemple si ${\mathcal C}$ est la cat\'egorie des sch\'emas de pr\'esentation finie sur un sch\'ema de base $S$, ${\mathcal M}$ est la classe des morphismes \'etales [resp. lisses, resp. plats] et (R) est la propri\'et\'e des familles de tels morphismes d'\^etre globalement surjectives.

\medskip

\item Plus g\'en\'eralement, la proposition s'applique aux foncteurs de changement de base
$$
(U \to X) \longmapsto (U \times_X X' \to X')
$$
associ\'es \`a n'importe quel morphisme de sch\'emas
$$
f : X' \longrightarrow X \, .
$$
On prend alors pour ${\mathcal C}$ et ${\mathcal C}'$ les cat\'egories des sch\'emas de pr\'esentation finie sur $X$ et $X'$, pour ${\mathcal M}$ et ${\mathcal M}'$ les classes des morphismes \'etales [resp. lisses, resp. plats] et pour (R) et (R$'$) les propri\'et\'es des familles de morphismes de ${\mathcal M}$ ou ${\mathcal M}'$ d'\^etre globalement surjectives.

\smallskip

Alors le morphisme
$$
f : X' \longrightarrow X
$$
d\'efinit un morphisme de topos
$$
(f^* , f_*)
$$
du gros ou petit topos \'etale [resp. lisse, resp. fppf] de $X'$ vers celui de $X$.

\smallskip

Sa composante d'image directe $f_*$ associe \`a tout faisceau $G$ sur le site de $X'$ le faisceau $f_* G$ sur le site de $X$ d\'efini par la formule
$$
(U \to X) \longmapsto G(U \times_X X' \to X') \, ,
$$
et sa composante d'image r\'eciproque $f^*$ associe \`a tout faisceau $F$ sur le site de $X$ le faisceau $f^* F$ sur le site de $X'$ obtenu comme le faisceautis\'e du pr\'efaisceau
$$
(U' \to X') \longmapsto \varinjlim_{\begin{pmatrix} U' \to U \\ \downarrow \qquad \downarrow \\ X' \xrightarrow{f} X \end{pmatrix}} F(U \to X) \, .
$$
\end{listeisansmarge}
\end{remarks}

\begin{demo}
\begin{listeisansmarge}
\item[(i)] est une cons\'equence imm\'ediate du corollaire \ref{corIV45}.

\smallskip

En effet, le foncteur induit par $\rho$
$$
{\mathcal C}/X \longrightarrow {\mathcal C}' / \rho (X)
$$
ou
$$
({\mathcal C}/X)_{\mathcal M} \longrightarrow ({\mathcal C}' / \rho (X))_{{\mathcal M}'}
$$
transforme les familles $J_R$-couvrantes de morphismes en familles $J_{R'}$-couvrantes.

\smallskip

D'autre part, il respecte les limites finies, donc aussi son compos\'e avec le foncteur canonique
$$
\ell : {\mathcal C}' / \rho (X) \xrightarrow{ \ y \ } \widehat{{\mathcal C}' / \rho (X)} \xrightarrow{ \ j^* \ } [\widehat{{\mathcal C}' / \rho (X)}]_{J_{R'}}
$$
ou
$$
\ell : ({\mathcal C}' / \rho (X))_{{\mathcal M}'} \xrightarrow{ \ y \ } (\widehat{{\mathcal C}' / \rho (X)})_{{\mathcal M}'} \xrightarrow{ \ j^* \ } [(\widehat{{\mathcal C}' / \rho (X))}_{{\mathcal M}'}]_{J_{R'}} \, .
$$

\item[(ii) et (iii)] r\'esultent alors du lemme \ref{lemIV42}. 
\end{listeisansmarge}
\end{demo}

\subsection{Foncteurs continus et topologies induites}\label{subsec444}

\medskip

Le lemme \ref{lemIV42} conduit \`a introduire une notion plus faible que celle de morphisme de sites:

\begin{defn}[Grothendieck (section III.1 of \cite{SGA4tome1})]\label{defIV49}

Soient $({\mathcal C},J)$ et $({\mathcal D},K)$ deux sites.

\smallskip

Un foncteur 
$$
\rho : {\mathcal C} \longrightarrow {\mathcal D}
$$
est appel\'e ``$(J,K)$-continu'' ou plus simplement ``continu'' s'il satisfait les deux conditions \'equivalentes suivantes:

\begin{enumerate}

\item[(A)] Le foncteur de composition avec $\rho$
$$
\begin{matrix}
\rho^* : \widehat{\mathcal D} &\longrightarrow &\widehat{\mathcal C} \, , \hfill \\
\hfill P &\longmapsto &P \circ \rho
\end{matrix}
$$
transforme les $K$-faisceaux sur ${\mathcal D}$ en $J$-faisceaux sur ${\mathcal C}$.

\medskip

\item[(B)] Son adjoint \`a gauche
$$
\begin{matrix}
\rho_! : \widehat{\mathcal C} &\longrightarrow &\widehat{\mathcal D} \, , \hfill \\
\hfill Q &\longmapsto &\displaystyle \varinjlim_{(X,x) \in \int\!\!Q} y \circ \rho (X) 
\end{matrix}
$$
s'inscrit dans un carr\'e commutatif \`a isomorphisme pr\`es:
$$
\xymatrix{
\widehat{\mathcal C} \ar[d]_-{j^*} \ar[r]^-{\rho_!} &\widehat{\mathcal D} \ar[d]^-{j^*} \\
\widehat{\mathcal C}_J \ar[r]^-{\rho_!} &\widehat{\mathcal D}_K
}
$$

\end{enumerate}
\end{defn}

\begin{remarkqed}

Le fait que $\rho_!$ est adjoint \`a gauche de $\rho^*$, et l'\'equivalence des conditions (A) et (B) qui en d\'ecoule, ont \'et\'e vus dans la d\'emonstration du lemme \ref{lemIV42}. 

\end{remarkqed}

\bigskip

On a:

\begin{prop}[Caramello (Proposition 4.13 de \cite{Denseness})]\label{propIV410}

Soient $({\mathcal C},J)$ et $({\mathcal D},K)$ deux sites.

\smallskip

Soit $\ell : {\mathcal D} \to \widehat{\mathcal D}_K$ le morphisme canonique associ\'e \`a $({\mathcal D},K)$.

\smallskip

Alors un foncteur
$$
\rho : {\mathcal C} \longrightarrow {\mathcal D}
$$
est continu si et seulement si tout crible $J$-couvrant $S$ d'un objet $X$ de ${\mathcal C}$ satisfait les deux conditions \'equivalentes suivantes:

\begin{enumerate}

\item[(A)] Si $S$ est vu comme une sous-cat\'egorie pleine de ${\mathcal C}/X$, on a dans $\widehat{\mathcal D}_K$
$$
\varinjlim_{(U \to X) \in S} \ell \circ \rho (U) = \ell \circ \rho (X) \, .
$$

\item[(B)] Il existe dans $S$ une famille d'objets
$$
(U_i \xrightarrow{ \ u_i \ } X)_{i \in I}
$$
dont les images par $\rho$
$$
\bigl( \rho (U_i) \xrightarrow{ \ \rho (u_i) \ } \rho (X)\bigl)_{i \in I}
$$
forment une famille $K$-couvrante de $\rho (X)$.

\smallskip

De plus, pour tous indices $i,j$ et tout carr\'e commutatif de ${\mathcal D}$
$$
\xymatrix{
V \ar[d]_-{v_i} \ar[r]^-{v_j} &\rho (U_i) \ar[d]^-{\rho (u_i)} \\
\rho (U_j) \ar[r]_-{\rho (u_j)} &\rho (X)
}
$$
il existe une famille de diagrammes de $S$ de la forme
$$
\def\troispoints{\ar@{}[ld]|{\displaystyle\cdots}}
\xymatrix{
U_i = U_{d_0} \ar[rd]_{u_i} \ar@{<->}[r] &U_{d_1} \ar[d] \ar@{<->}[r] &\ldots \ar@{<->}[r] \troispoints &U_{d_n} = U_j \ar[lld]^{u_j} \\
&X
}
$$
compl\'et\'es par des diagrammes commutatifs de ${\mathcal D}$

$$\def\troispoints{\ar@{}[rd]|{\displaystyle\cdots}}
\xymatrix{
&V' \ar[ld]_{v'_i} \ar[d] \troispoints \ar[rrd]^{v'_j} \\
\rho(U_i) = \rho(U_{d_0}) \ar@{<->}[r]&\rho(U_{d_1}) \ar@{<->}[r] &\cdots \ar@{<->}[r]&\rho (U_{d_n}) = \rho (U_j)
} 
$$
tels que les $v'_i$ et $v'_j$ se factorisent en
$$
v'_i : V' \xrightarrow{ \ v' \ } V \xrightarrow{ \ v_i \ } \rho (U_i) \, , 
$$
$$
v'_j : V' \xrightarrow{ \ v' \ } V \xrightarrow{ \ v_j \ } \rho (U_j)
$$
et que les
$$
V' \xrightarrow{ \ v' \ } V
$$
forment une famille $K$-couvrante de $V$.

\end{enumerate}
\end{prop}

\bigskip

\begin{remark}

Il r\'esulte du th\'eor\`eme \ref{thmIV43} que si le foncteur compos\'e
$$
{\mathcal C} \xrightarrow{ \ \rho \ } {\mathcal D} \xrightarrow{ \ \ell \ } \widehat{\mathcal D}_K
$$
est plat, $\rho$ est continu si et seulement si il transforme toute famille $J$-couvrante de morphismes de ${\mathcal C}$ en une famille $K$-couvrante de morphismes de ${\mathcal D}$.

\end{remark}

\bigskip

\begin{demo}

La condition (A) de la d\'efinition \ref{defIV49} est satisfaite si et seulement si on a pour tout faisceau $F$ sur $({\mathcal D},J)$ et tout crible $J$-couvrant $S$ d'un objet $X$ de ${\mathcal C}$
$$
F(\rho (X)) = \varprojlim_{(U\to X) \in S} F(\rho (U))
$$
soit
$$
{\rm Hom} (\ell \circ \rho (X),F) = \varprojlim_{(U \to X) \in S} {\rm Hom} (\ell \circ \rho (U),F) \, .
$$

Cela \'equivaut \`a demander que tout crible $J$-couvrant $S$ d'un objet $X$ de ${\mathcal C}$ satisfasse la condition (A) de l'\'enonc\'e
$$
\ell \circ \rho (X) = \varinjlim_{(U \to X) \in S} \ell \circ \rho (U) \, .
$$

Un tel crible $S$ satisfait cette condition (A) s'il contient une famille d'\'el\'ements
$$
U_i \longrightarrow X \, , \qquad i \in I \, ,
$$
telle que la famille des images
$$
\ell \circ \rho (U_i) \longrightarrow \ell \circ \rho (X)
$$
soit \'epimorphique, et que la relation d'\'equivalence
$$
\left(\coprod_i \ell \circ \rho (U_i) \right) \times_{\ell \circ \rho (X)} \left( \coprod_j \ell \circ \rho (U_j)\right) = \coprod_{i,j} \ell \circ \rho (U_i) \times_{\ell \circ \rho (X)} \ell \circ \rho (U_j)
$$
soit engendr\'ee par des morphismes de la sous-cat\'egorie pleine $S$ de ${\mathcal C} / X$.

\smallskip

Cela \'equivaut \`a la condition (B) de l'\'enonc\'e.

\smallskip

D'o\`u la conclusion. 

\end{demo}

\medskip

Cette proposition conduit \`a poser:

\begin{defn}\label{defIV411}

Soit
$$
\rho : {\mathcal C} \longrightarrow {\mathcal D}
$$
un foncteur entre deux cat\'egories essentiellement petites.

\smallskip

Pour toute topologie $K$ sur ${\mathcal D}$, on appelle topologie induite par $K$ via $\rho$, et on note
$$
\rho^{-1} K \, ,
$$
la topologie de ${\mathcal C}$ pour laquelle un crible $S$ d'un objet $X$ est couvrant si et seulement si, pour tout morphisme de ${\mathcal C}$
$$
f : X' \longrightarrow X \, ,
$$
le crible $f^{-1} S$ de $X'$ satisfait les conditions \'equivalentes {\rm (A)} et {\rm (B)} de la proposition \ref{propIV410}.
\end{defn}

\bigskip

\begin{remarksqed}
\begin{listeisansmarge}
\item La famille $\rho^{-1} K$ de cribles de ${\mathcal C}$ est bien une topologie de Grothendieck.

\smallskip

En effet, il est \'evident sur sa d\'efinition qu'elle satisfait les axiomes de maximalit\'e et de stabilit\'e.

\smallskip

Pour la transitivit\'e, consid\'erons deux cribles $S$ et $S'$ d'un objet $X$ de ${\mathcal C}$ tels que $S$ soit \'el\'ement de $\rho^{-1} K$ et que, pour tout morphisme $u : U \to X$ de $S$, $u^{-1} S'$ soit \'el\'ement de $\rho^{-1} K$.

\smallskip

Notant $\ell : {\mathcal D} \to \widehat{\mathcal D}_K$ le morphisme canonique, on a pour tout morphisme $X' \xrightarrow{ \ f \ } X$ de ${\mathcal C}$
$$
\ell \circ \rho (X') = \varinjlim_{(U' \to X') \in f^{-1} S} \ell \circ \rho (U')
$$
et, pour tout \'el\'ement $u' : U' \to X'$ de $f^{-1} S$,
$$
\ell \circ \rho (U') = \varinjlim_{(U'' \to U') \in (f \circ u')^{-1} S'} \ell \circ \rho (U'') \, .
$$
On en d\'eduit que, pour tout faisceau $G$ sur $({\mathcal D},K)$, tout morphisme
$$
\varinjlim_{(U'' \to X') \in f^{-1} S'} (\ell \circ \rho)(U') \longrightarrow G
$$
induit une famille compatible de morphismes
$$
(\ell \circ \rho)(U') \longrightarrow G
$$
index\'es par les \'el\'ements $(U' \to X')$ de $f^{-1} S$ et donc provient d'un unique morphisme
$$
\ell \circ \rho (X') \longrightarrow G \, .
$$

Cela montre comme voulu que
$$
\ell \circ \rho (X') = \varinjlim_{(U'' \to X') \in f^{-1} S'} \ell \circ \rho (U'')
$$
qui signifie que le crible $S'$ est \'el\'ement de $\rho^{-1} K$.

\medskip

\item Il r\'esulte de la remarque qui suit la proposition \ref{propIV410} que, si le foncteur compos\'e
$$
{\mathcal C} \xrightarrow{ \ \rho \ } {\mathcal D} \xrightarrow{ \ \ell \ } \widehat{\mathcal D}_K
$$
est plat, la topologie induite $\rho^{-1} K$ est constitu\'ee des cribles $S$ d'objets $X$ tels que, pour tout morphisme
$$
f : X' \longrightarrow X \, ,
$$
le crible de $\rho (X')$ engendr\'e par les images par $\rho$ des \'el\'ements $U' \to X'$ de $f^{-1} S$ est $K$-couvrant. 

\end{listeisansmarge}
\end{remarksqed}

\bigskip

Avec cette d\'efinition, la proposition \ref{propIV410} se reformule ainsi:

\begin{cor}\label{corIV412}

Soient $({\mathcal C},J)$ et $({\mathcal D},K)$ deux sites.

\smallskip

Alors un foncteur
$$
\rho : {\mathcal C} \longrightarrow {\mathcal D}
$$
est $(J,K)$-continu si et seulement si $J$ est contenue dans la topologie $\rho^{-1} K$ induite par $K$ via $\rho$.
\end{cor}

\begin{remarkqed}

Tout compos\'e de foncteurs continus est un foncteur continu.

\smallskip

Par cons\'equent, pour tous foncteurs entre des cat\'egories essentiellement petites
$$
{\mathcal C}' \xrightarrow{ \ \eta \ } {\mathcal C} \xrightarrow{ \ \rho \ } {\mathcal D}
$$
et toute topologie $K$ sur ${\mathcal D}$ on a
$$
(\rho \circ \eta)^{-1} K \subseteq \eta^{-1} (\rho^{-1} K) \, .
$$

En revanche, il n'y a pas n\'ecessairement \'egalit\'e. \end{remarkqed}

\section{Comorphismes de sites}\label{sec45}

\subsection{La notion de comorphisme de sites}\label{subsec451}

\medskip

Les morphismes d'un site $({\mathcal C} , J)$ dans un site $({\mathcal D},K)$ ont \'et\'e d\'efinis comme des foncteurs $\rho : {\mathcal C} \to {\mathcal D}$ qui induisent un morphisme de topos $(\rho_! , \rho^!) : \widehat{\mathcal D}_K \to \widehat{\mathcal C}_J$ dont la composante d'image r\'eciproque $\rho_! : \widehat{\mathcal C}_J \to \widehat{\mathcal D}_K$ se d\'eduit de $\rho$ par cocompl\'etion.

\smallskip

Dualement, on d\'efinit les ``comorphismes'' d'un site $({\mathcal C},J)$ dans un site $({\mathcal D},K)$ comme les foncteurs $\rho : {\mathcal C} \to {\mathcal D}$ qui induisent un morphisme de topos $(\rho^* , \rho_*) : \widehat{\mathcal C}_J \to \widehat{\mathcal D}_K$ dont la composante d'image r\'eciproque $\rho^* : \widehat{\mathcal D}_K \to \widehat{\mathcal C}_J$ se d\'eduit du foncteur $\rho^* : \widehat{\mathcal D} \to \widehat{\mathcal C}$ de composition avec $\rho$:

\begin{defn}\label{defIV51}

Soient $({\mathcal C},J)$ et $({\mathcal D},K)$ deux sites, avec les foncteurs de faisceautisation associ\'es
$$
j^* : \widehat{\mathcal C} \longrightarrow \widehat{\mathcal C}_J
$$
et
$$
j^* : \widehat{\mathcal D} \longrightarrow \widehat{\mathcal D}_K \, .
$$

Un foncteur
$$
\rho : {\mathcal C} \longrightarrow {\mathcal D}
$$
sera appel\'e un comorphisme de sites s'il existe un morphisme de topos, n\'ecessairement unique \`a unique isomorphisme pr\`es,
$$
(\rho^* , \rho_*) : \widehat{\mathcal C}_J \longrightarrow \widehat{\mathcal D}_K
$$
dont la composante d'image r\'eciproque
$$
\rho^* : \widehat{\mathcal D}_K \longrightarrow \widehat{\mathcal C}_J
$$
s'inscrit dans un carr\'e commutatif \`a isomorphisme pr\`es
$$
\xymatrix{
\widehat{\mathcal D} \ar[d]_-{j^*} \ar[r]^-{\rho^*} &\widehat{\mathcal C} \ar[d]^-{j^*} \\
\widehat{\mathcal D}_K \ar[r]^-{\rho^*} &\widehat{\mathcal C}_J
}
$$
dont la premi\`ere fl\`eche horizontale est le foncteur de composition avec $\rho$
$$
\begin{matrix}
\rho^* : \widehat{\mathcal D} &\longrightarrow &\widehat{\mathcal C} \, \hfill \\
\hfill G &\longmapsto &G \circ \rho \, .
\end{matrix}
$$
\end{defn}

\begin{remarksqed}
\begin{listeisansmarge}
\item Comme la transformation naturelle d'adjonction
$$
j^* \circ j_* \longrightarrow {\rm id}
$$ 
est un isomorphisme de foncteurs $\widehat{\mathcal D}_K \to \widehat{\mathcal D}_K$, la composante d'image r\'eciproque $\rho^*$ du morphisme de topos $(\rho^* , \rho_*) : \widehat{\mathcal C}_J \to \widehat{\mathcal D}_K$ est canoniquement isomorphe au foncteur compos\'e
$$
\widehat{\mathcal D}_K \xhookrightarrow{ \ j_* \ } \widehat{\mathcal D} \xrightarrow{ \ \rho^* \ } \widehat{\mathcal C} \xrightarrow{ \ j^* \ } \widehat{\mathcal C}_J \, .
$$
C'est pourquoi le morphisme de topos $(\rho^* , \rho_*) : \widehat{\mathcal C}_J \to \widehat{\mathcal D}_K$ est uniquement d\'etermin\'e par $\rho$ \`a unique isomorphisme pr\`es, s'il existe.

\medskip

\item Un tel morphisme de topos est not\'e $(\rho^* , \rho_*)$ pour signifier que sa composante d'image r\'eciproque est contravariante par rapport \`a $\rho$ et sa composante d'image directe est covariante. On utilise le signe ``$*$'' plut\^ot que le signe ``$!$'' pour signifier que c'est l'adjoint \`a gauche qui est contravariant et l'adjoint \`a droite covariant.

\medskip

\item Si $({\mathcal C}_1 , J_1)$, $({\mathcal C}_2 , J_2)$ et $({\mathcal C}_3,J_3)$ sont trois sites reli\'es par deux comorphismes de sites au sens de la d\'efinition
$$
{\mathcal C}_1 \xrightarrow{ \ \rho \ } {\mathcal C}_2 \xrightarrow{ \ \eta \ } {\mathcal C}_3 \, ,
$$
alors le foncteur compos\'e
$$
\eta \circ \rho : {\mathcal C}_1 \longrightarrow {\mathcal C}_3
$$
est encore un comorphisme de sites et admet pour morphisme de topos associ\'e le compos\'e des morphismes de topos associ\'es \`a $\rho$ et $\eta$
$$
((\eta \circ \rho)^* , (\eta \circ \rho)_*) = (\rho^* \circ \eta^* , \eta_* \circ \rho_*) \, .
$$

\item Toute transformation naturelle
$$
\rho \xrightarrow{ \ \alpha \ } \rho'
$$
entre deux foncteurs
$$
{\mathcal C} \raisebox{.7ex}{\xymatrix{\dar[r]^-{^{^{\mbox{\scriptsize$\rho$}}}}_-{\rho'} &{\mathcal D}}}
$$
qui sont des comorphismes d'un site $({\mathcal C},J)$ dans un site $({\mathcal D} ,K)$ induit un morphisme du morphisme de topos
$$
(\rho'^* , \rho'_*) : \widehat{\mathcal C}_J \longrightarrow \widehat{\mathcal D}_K
$$
dans le morphisme de topos
$$
(\rho^* , \rho_*) : \widehat{\mathcal C}_J \longrightarrow \widehat{\mathcal D}_K \, .
$$

Celui-ci est d\'efini par la transformation naturelle du foncteur
$$
\rho'^* : \widehat{\mathcal D} \longrightarrow \widehat{\mathcal C}
$$
dans le foncteur
$$
\rho^* : \widehat{\mathcal D} \longrightarrow \widehat{\mathcal C}
$$
qui associe \`a tout pr\'efaisceau $G$ sur ${\mathcal D}$ le morphisme de pr\'efaisceau sur ${\mathcal C}$
$$
G \circ \rho' \longrightarrow G \circ \rho
$$
qui consiste en la famille des applications
$$
G \circ \rho' (X) \longrightarrow G \circ \rho (X) \, , \qquad X \in {\rm Ob} ({\mathcal C}) \, ,
$$
associ\'ees aux morphismes de ${\mathcal D}$
$$
\alpha_X : \rho (X) \longrightarrow \rho'(X)
$$
qui composent la transformation naturelle
$$
\alpha : \rho \longrightarrow \rho' \, .
$$

\end{listeisansmarge}
\end{remarksqed}

On rappelle que d'apr\`es la proposition \ref{propI106}, le foncteur de composition
$$
\begin{matrix}
\rho^* : \widehat{\mathcal D} &\longrightarrow &\widehat{\mathcal C} \, , \hfill \\
\hfill G &\longmapsto &G \circ \rho
\end{matrix}
$$
d\'efini par un foncteur entre cat\'egories essentiellement petites
$$
\rho : {\mathcal C} \longrightarrow {\mathcal D}
$$
poss\`ede \`a la fois un adjoint \`a gauche
$$
\rho_! : \widehat{\mathcal C} \longrightarrow \widehat{\mathcal D}
$$
et un adjoint \`a droite
$$
\rho_* : \widehat{\mathcal C} \longrightarrow \widehat{\mathcal D} \, .
$$

De plus, on a pour tout pr\'efaisceau $F$ sur ${\mathcal C}$ et tout objet $D$ de ${\mathcal D}$ les formules
$$
\rho_! F(D) = \varinjlim_{D \backslash {\mathcal C} \ni X} F(X)
$$
et
$$
\rho_* F(D) = \varprojlim_{X \in {\mathcal C} / D} F(X)
$$
o\`u $D \backslash {\mathcal C}$ [resp. ${\mathcal C} / D$] d\'esigne la cat\'egorie dont les objets sont les paires $(X,d)$ constitu\'ees d'un objet $X$ de ${\mathcal C}$ et d'un morphisme de ${\mathcal D}$
$$
D \xrightarrow{ \ d \ } \rho (X) \qquad \mbox{[resp.} \quad \rho (X) \xrightarrow{ \ d \ } D \ \mbox{]}
$$
et dont les morphismes
$$
(X,d) \longrightarrow (Y,e)
$$
sont les morphismes $f : X \to Y$ de ${\mathcal C}$ qui rendent commutatif le triangle
$$
\xymatrix{
&D \ar[ld]_d \ar[rd]^e \\
\rho(X) \ar[rr]^{\rho (f)} &&\rho (Y)
} \qquad
\xymatrix{
\rho(X) \ar[rd]_{\mbox{$\Bigl[$resp. \qquad\qquad }} \ar[rr]^{\rho (f)} &&\rho(Y) \ar[ld]^{\mbox{\qquad\qquad$\Bigl].$}} \\
&D
}
$$

La formule qui d\'efinit l'adjoint \`a droite $\rho_* : \widehat{\mathcal C} \to \widehat{\mathcal D}$ permet d'identifier la composante d'image directe du morphisme de topos $(\rho^* , \rho_*) : \widehat{\mathcal C}_J \to \widehat{\mathcal D}_K$ associ\'e \`a un comorphisme de sites $\rho : {\mathcal C} \to {\mathcal D}$. En effet, on a:

\begin{lem}\label{lemIV52}

Soient $({\mathcal C},J)$ et $({\mathcal D},K)$ deux sites munis des paires de foncteurs adjoints
$$
\left(\widehat{\mathcal C} \xrightarrow{ \ j^* \ } \widehat{\mathcal C}_J \, , \quad \widehat{\mathcal C}_J \xhookrightarrow{ \ j_* \ } \widehat{\mathcal D}_K \right)
$$
et
$$
\left(\widehat{\mathcal D} \xrightarrow{ \ j^* \ } \widehat{\mathcal D}_K \, , \quad \widehat{\mathcal D}_K \xhookrightarrow{ \ j_* \ } \widehat{\mathcal D} \right).
$$

Soit $\rho : {\mathcal C} \to {\mathcal D}$ un foncteur.

\smallskip

Alors:

\begin{listeimarge}

\item Le foncteur de composition avec $\rho$
$$
\begin{matrix}
\rho^* : \widehat{\mathcal D} &\longrightarrow &\widehat{\mathcal C} \, , \hfill \\
\hfill G &\longmapsto &G \circ \rho
\end{matrix}
$$
s'inscrit dans un carr\'e commutatif \`a isomorphisme pr\`es
$$
\xymatrix{
\widehat{\mathcal D} \ar[d]_-{j^*} \ar[r]^-{\rho^*} &\widehat{\mathcal C} \ar[d]^-{j^*} \\
\widehat{\mathcal D}_K \ar[r]^-{\rho^*} &\widehat{\mathcal C}_J
}
$$
si et seulement si son adjoint \`a droite
$$
\rho_* : \widehat{\mathcal C} \longrightarrow \widehat{\mathcal D}
$$
s'inscrit dans un carr\'e commutatif
$$
\xymatrix{
\widehat{\mathcal C} \ar[r]^-{\rho_*} &\widehat{\mathcal D}  \\
\widehat{\mathcal C}_J \ar@{^{(}->}[u]^{j_*} \ar[r]^-{\rho_*} &\widehat{\mathcal D}_K \ar@{^{(}->}[u]_{j_*}
}
$$
c'est-\`a-dire si $\rho_* : \widehat{\mathcal C} \to \widehat{\mathcal D}$ transforme les $J$-faisceaux sur ${\mathcal C}$ en $K$-faisceaux sur ${\mathcal D}$.

\medskip

\item Dans ce cas, le foncteur 
$$
\rho_* :  \widehat{\mathcal C}_J \longrightarrow  \widehat{\mathcal D}_K
$$
est adjoint \`a droite de $\rho^* :  \widehat{\mathcal D}_K \to  \widehat{\mathcal C}_J$. Il associe \`a tout faisceau $F$ sur $({\mathcal C},J)$ le faisceau $\rho_* F$ sur $({\mathcal D},K)$ d\'efini par la formule
$$
\rho_* F(D) = \varprojlim_{X \in {\mathcal C}/D} F(X) \, .
$$

\item Le foncteur
$$
\rho^* : \widehat{\mathcal D}_K \longrightarrow \widehat{\mathcal C}_J
$$
associe \`a tout faisceau $G$ sur $({\mathcal D},K)$ le faisceautis\'e du pr\'efaisceau $G \circ \rho$ sur ${\mathcal C}$.
\end{listeimarge}
\end{lem}

\bigskip

\begin{demo}
\begin{listeisansmarge}
\item Si le foncteur $\rho^* : \widehat{\mathcal D} \to \widehat{\mathcal C}$ de composition avec $\rho$ [resp. son adjoint \`a droite $\rho_* : \widehat{\mathcal C} \to \widehat{\mathcal D}$] s'inscrit dans un carr\'e commutatif \`a isomorphisme pr\`es
$$
\xymatrix{
\widehat{\mathcal D} \ar[d]_-{j^*} \ar[r]^-{\rho^*} &\widehat{\mathcal C} \ar[d]^-{j^*} \\
\widehat{\mathcal D}_K \ar[r]^-{\rho^*} &\widehat{\mathcal C}_J
} \qquad \quad
\xymatrix{
\widehat{\mathcal C} \ar[r]^-{\rho_*} &\widehat{\mathcal D}  \\
\widehat{\mathcal C}_J \ar@{^{(}->}[u]^{\mbox{[resp.} \qquad j_*} \ar[r]^-{\rho_*} &\widehat{\mathcal D}_K \ar@{^{(}->}[u]_{j_* \ \mbox{]},}
}
$$
le foncteur $\rho^* : \widehat{\mathcal D}_K \to \widehat{\mathcal C}_J$ [resp. $\rho_* : \widehat{\mathcal C}_J \to \widehat{\mathcal D}_K$] respecte les colimites arbitraires [resp. les limites arbitraires] puisqu'il en est ainsi des foncteurs $\rho^* : \widehat{\mathcal D} \to \widehat{\mathcal C}$, $j^* : \widehat{\mathcal C} \to \widehat{\mathcal C}_J$ et $j^* : \widehat{\mathcal D} \to \widehat{\mathcal D}_K$ [resp. $\rho_* : \widehat{\mathcal C} \to \widehat{\mathcal D}$, $j_* : \widehat{\mathcal C}_J \hookrightarrow \widehat{\mathcal C}$ et $j_* : \widehat{\mathcal D}_K \hookrightarrow \widehat{\mathcal D}$] et que les deux foncteurs $j_* : \widehat{\mathcal C}_J \hookrightarrow \widehat{\mathcal C}$ et $j_* : \widehat{\mathcal D}_K \hookrightarrow \widehat{\mathcal D}$ sont pleinement fid\`eles.

\smallskip

Donc le foncteur $\rho^* : \widehat{\mathcal D}_K \to \widehat{\mathcal C}_J$ [resp. $\rho_* : \widehat{\mathcal C}_J \to \widehat{\mathcal D}_K$] admet un adjoint \`a droite
$$
\rho_* : \widehat{\mathcal C}_J \longrightarrow \widehat{\mathcal D}_K
$$
[resp. un adjoint \`a gauche
$$
\rho^* : \widehat{\mathcal D}_K \longrightarrow \widehat{\mathcal C}_J \ \mbox{]}.
$$

Alors les deux carr\'es se d\'eduisent l'un de l'autre par adjonction.

\smallskip

Dire que le foncteur $\rho_* : \widehat{\mathcal C} \to \widehat{\mathcal D}$ s'inscrit dans un carr\'e commutatif \`a isomorphisme pr\`es
$$
\xymatrix{
\widehat{\mathcal C} \ar[r]^-{\rho_*} &\widehat{\mathcal D}  \\
\widehat{\mathcal C}_J \ar@{^{(}->}[u]^{j_*} \ar[r]^-{\rho_*} &\widehat{\mathcal D}_K \ar@{^{(}->}[u]_{j_*}
}
$$
signifie qu'il transforme tout $J$-faisceau sur ${\mathcal C}$ en un $K$-faisceau sur ${\mathcal D}$.

\smallskip

Quitte \`a remplacer $\rho_* : \widehat{\mathcal C}_J \to \widehat{\mathcal D}_K$ par la restriction de $\rho_* : \widehat{\mathcal C} \to \widehat{\mathcal D}$ \`a $\widehat{\mathcal C}_J \hookrightarrow \widehat{\mathcal C}$, qui lui est canoniquement isomorphe, on peut supposer que le carr\'e
$$
\xymatrix{
\widehat{\mathcal C} \ar[r]^-{\rho_*} &\widehat{\mathcal D}  \\
\widehat{\mathcal C}_J \ar@{^{(}->}[u]^{j_*} \ar[r]^-{\rho_*} &\widehat{\mathcal D}_K \ar@{^{(}->}[u]
}
$$
est commutatif.

\medskip

\item Si les deux conditions \'equivalentes de (i) sont satisfaites, les deux carr\'es de (i) se d\'eduisent l'un de l'autre par adjonction.

\smallskip

Le foncteur $\rho_* : \widehat{\mathcal C}_J \to \widehat{\mathcal D}_K$ est la restriction du foncteur
$$
\rho_* : \widehat{\mathcal C} \longrightarrow \widehat{\mathcal D}
$$
\`a la sous-cat\'egorie pleine $\widehat{\mathcal C}_J \hookrightarrow \widehat{\mathcal C}$. Donc il est d\'efini par la m\^eme formule
$$
F \longmapsto \rho_* F = \left[ D \longmapsto \varprojlim_{X \in {\mathcal C}/D} F(X) \right]
$$
donn\'ee par la proposition \ref{propI106} pour le foncteur
$$
\rho_* : \widehat{\mathcal C} \longrightarrow \widehat{\mathcal D} \, .
$$

\item r\'esulte de ce que le foncteur
$$
\rho^* : \widehat{\mathcal D}_K \longrightarrow \widehat{\mathcal C}_J
$$
est canoniquement isomorphe au compos\'e
$$
\widehat{\mathcal D}_K \xhookrightarrow{ \ j_* \ } \widehat{\mathcal D} \xrightarrow{ \ \rho^* \ } \widehat{\mathcal C} \xrightarrow{ \ j^* \ } \widehat{\mathcal C}_J
$$ 
comme not\'e dans la remarque (ii) qui suit la d\'efinition \ref{defIV51}. 
\end{listeisansmarge}
\end{demo}

\subsection{Caract\'erisation des comorphismes de sites}\label{subsec452}

\medskip

Voici les conditions n\'ecessaires et suffisantes pour qu'un foncteur $\rho : {\mathcal C} \to {\mathcal D}$ soit un comorphisme d'un site $({\mathcal C},J)$ dans un site $({\mathcal D},K)$:

\begin{thm}\label{thmIV53}

Soient $({\mathcal C},J)$ et $({\mathcal D},K)$ deux sites.

\smallskip

Pour tout foncteur
$$
\rho : {\mathcal C} \longrightarrow {\mathcal D} \, ,
$$
les conditions suivantes sont \'equivalentes:

\begin{listeimarge}

\item[(1)] Le foncteur $\rho$ est un comorphisme de sites.

\medskip

\item[(2)] L'adjoint \`a droite
$$
\begin{matrix}
\rho_* : \widehat{\mathcal C} &\longrightarrow &\widehat{\mathcal D} \, , \hfill \\
\hfill F &\longmapsto &\rho_* F = \left[ D \longmapsto \displaystyle \varprojlim_{X \in {\mathcal C} / D} F(X) \right]
\end{matrix}
$$
du foncteur de composition avec $\rho$
$$
\begin{matrix}
\rho^* : \widehat{\mathcal D} &\longrightarrow &\widehat{\mathcal C} \, , \hfill \\
\hfill G &\longmapsto &G \circ \rho
\end{matrix}
$$
transforme tout $J$-faisceau sur ${\mathcal C}$ en un $K$-faisceau sur ${\mathcal D}$.

\medskip

\item[(3)] Pour que le foncteur compos\'e
$$
\widehat{\mathcal D} \xrightarrow{ \ \rho^* \ } \widehat{\mathcal C} \xrightarrow{ \ j^* \ } \widehat{\mathcal C}_J
$$
transforme un morphisme de pr\'efaisceaux sur ${\mathcal D}$
$$
G \longrightarrow G'
$$
en un isomorphisme de $\widehat{\mathcal C}_J$, il suffit que le morphisme induit entre les faisceautis\'es soit un isomorphisme dans $\widehat{\mathcal D}_K$.

\medskip

\item[(4)] Pour tout crible $K$-couvrant $S$ d'un objet $D$ de ${\mathcal D}$ vu comme une sous-cat\'egorie pleine de ${\mathcal D}/D$, le morphisme canonique de $\widehat{\mathcal C}_J$
$$
\varinjlim_{(V \to D) \in S} j^* \circ \rho^* \circ y(V) \longrightarrow j^* \circ \rho^* \circ y(D)
$$
est un isomorphisme.

\medskip

\item[(5)] Pour tout crible $K$-couvrant $S$ d'un objet $D$ de ${\mathcal D}$ et tout objet $X$ de ${\mathcal C}$ muni d'un morphisme
$$
d : \rho (X) \longrightarrow D \, ,
$$
les morphismes $u : U \to X$ de ${\mathcal C}$ tels que le compos\'e
$$
\rho (U) \xrightarrow{ \ \rho (u) \ } \rho (X) \xrightarrow{ \ d \ } D
$$
soit \'el\'ement de $S$, forment un crible $J$-couvrant de $X$.
\end{listeimarge}
\end{thm}

\bigskip

\begin{demo}

La condition (3) \'equivaut \`a l'existence d'un foncteur
$$
\rho^* : \widehat{\mathcal D}_K \longrightarrow \widehat{\mathcal C}_J
$$
qui rende le carr\'e
$$
\xymatrix{
\widehat{\mathcal D} \ar[d]_{j^*} \ar[r]^{\rho^*} &\widehat{\mathcal C} \ar[d]^{j^*} \\
\widehat{\mathcal D}_K \ar[r]^{\rho^*} &\widehat{\mathcal C}_J
}
$$
commutatif \`a isomorphisme pr\`es.

\smallskip

Donc les conditions (2) et (3) sont \'equivalentes d'apr\`es le lemme \ref{lemIV52}.

\smallskip

Pour v\'erifier qu'elles \'equivalent aussi \`a la condition (1), il suffit de montrer que tout foncteur
$$
\rho^* : \widehat{\mathcal D}_K \longrightarrow \widehat{\mathcal C}_J
$$
qui rend le carr\'e ci-dessus commutatif \`a isomorphisme pr\`es, respecte n\'ecessairement les limites finies.

\smallskip

Cela r\'esulte de ce qu'il en est ainsi des foncteurs
$$
\rho^* : \widehat{\mathcal D} \longrightarrow \widehat{\mathcal C} \, , \quad j^* : \widehat{\mathcal C} \longrightarrow \widehat{\mathcal C}_J \quad \mbox{et} \quad j^* : \widehat{\mathcal D} \longrightarrow \widehat{\mathcal D}_K
$$
et de ce que le foncteur $j_* : \widehat{\mathcal D}_K \hookrightarrow \widehat{\mathcal D}$ est pleinement fid\`ele.

\smallskip

La condition (2) signifie que pour tout faisceau $F$ sur $({\mathcal C},J)$ et tout crible $K$-couvrant $S$ d'un objet $D$ de ${\mathcal D}$, on a
$$
\rho_* (F)(D) = \varprojlim_{(V \to D) \in S} \rho_* (F)(V)
$$
soit
$$
{\rm Hom} (y(D),\rho_* F) = \varprojlim_{(V \to D) \in S} {\rm Hom} (y(V),\rho_* F)
$$
ou encore par adjonction
$$
{\rm Hom} (\rho^* \circ y(D),F) = \varprojlim_{(V \to D) \in S} {\rm Hom} (\rho^* \circ y(V),F)
$$
et, puisque $F$ est un $J$-faisceau sur ${\mathcal C}$,
$$
{\rm Hom} (j^* \circ \rho^* \circ y (D),F) = \varprojlim_{(V \to D) \in S} {\rm Hom} (j^* \circ \rho^* \circ y (V),F) \, .
$$

Cela signifie que l'on a dans la cat\'egorie $\widehat{\mathcal C}_J$
$$
j^* \circ \rho^* \circ y(D) = \varinjlim_{(V \to D)\in S} j^* \circ \rho^* \circ y (V) \, .
$$

Ainsi, les conditions (2) et (4) sont \'equivalentes.

\smallskip

Comme le foncteur $j^* : \widehat{\mathcal C} \to \widehat{\mathcal C}_J$ respecte les colimites, (4) \'equivaut \`a demander que le morphisme canonique de $\widehat{\mathcal C}$
$$
\varinjlim_{(V \to D) \in S} \rho^* \circ y(V) \longrightarrow \rho^* \circ y(D)
$$
devienne un isomorphisme apr\`es faisceautisation.

\smallskip

S'il en est ainsi, il existe pour tout objet $X$ de ${\mathcal C}$ et tout morphisme
$$
d : \rho (X) \longrightarrow D
$$
une famille $J$-couvrante de $X$
$$
U_i \xrightarrow{ \ u_i \ } X \, , \qquad i \in I \, ,
$$
et des \'el\'ements de $S$
$$
V_i \longrightarrow D
$$
qui s'inscrivent dans des carr\'es commutatifs de ${\mathcal D}$:
$$
\xymatrix{
\rho (U_i) \ar[d]_{\rho (u_i)} \ar[r] &V_i \ar[d] \\
\rho(X) \ar[r] &D
}
$$

Cela montre que (4) implique (5).

\smallskip

R\'eciproquement, si la condition (5) est v\'erifi\'ee, le morphisme de $\widehat{\mathcal C}_J$
$$
j^* \left( \varinjlim_{(V \to D) \in S} \rho^* \circ y(V) \right) \longrightarrow j^* \circ \rho^* \circ y(D)
$$
est un \'epimorphisme.

\smallskip

Pour conclure que (5) implique (4), il suffit de remarquer que c'est n\'ecessairement un monomorphisme.

\smallskip

En effet, le morphisme de $\widehat{\mathcal D}$
$$
\varinjlim_{(V \to D) \in S} y(V) = S \longrightarrow y(D)
$$
est un monomorphisme, et le foncteur $\widehat{\mathcal D} \xrightarrow{ \ \rho^* \ } \widehat{\mathcal C} \xrightarrow{ \ j^* \ } \widehat{\mathcal C}_J$ respecte les colimites et les monomorphismes.

\smallskip

Cela ach\`eve la preuve du th\'eor\`eme. 

\end{demo}

\subsection{Topologies coinduites et comorphismes de sites}\label{subsec453}

\medskip

Le crit\`ere du th\'eor\`eme \ref{thmIV53} conduit \`a poser la d\'efinition suivante:

\begin{defn}\label{defIV54}

Soit
$$
\rho : {\mathcal C} \longrightarrow {\mathcal D}
$$
un foncteur entre deux cat\'egories essentiellement petites.

\smallskip

Pour toute topologie $J$ sur ${\mathcal C}$, on appelle topologie coinduite par $J$ via $\rho$, et on note
$$
\rho_* J \, ,
$$
la topologie sur ${\mathcal D}$ pour laquelle un crible $S$ d'un objet $D$ est couvrant si et seulement si, pour tout objet $X$ de ${\mathcal C}$ muni d'un morphisme
$$
d : \rho (X) \longrightarrow D \, ,
$$
le crible de $X$ constitu\'e des morphismes
$$
u : U \longrightarrow X
$$
tels que le compos\'e
$$
\rho (U) \xrightarrow{\ \rho (u) \ } \rho (X) \xrightarrow{ \ d \ } D
$$
soit \'el\'ement de $S$, est un crible $J$-couvrant de $X$.
\end{defn}

\bigskip

\begin{remarkqed}

La collection ainsi d\'efinie $\rho_* J$ de cribles de ${\mathcal D}$ est bien une topologie de Grothendieck.

\smallskip

En effet, il est \'evident sur sa d\'efinition qu'elle satisfait les axiomes de maximalit\'e et de stabilit\'e.

\smallskip

Pour la transitivit\'e, consid\'erons deux cribles $S$ et $S'$ d'un objet $D$ de ${\mathcal D}$ tels que $S$ soit \'el\'ement de $j_* J$ et que, pour tout morphisme $v : V \to D$ de $S$, $v^{-1} S'$ soit \'el\'ement de $\rho_* J$.

\smallskip

Consid\'erons un objet $X$ de ${\mathcal C}$ muni d'un morphisme
$$
d : \rho (X) \longrightarrow D \, .
$$
Par hypoth\`ese, les morphismes $U \xrightarrow{ \ u \ } X$ de ${\mathcal C}$ tels que $\rho (u)$ s'inscrive dans un carr\'e commutatif
$$
\xymatrix{
\rho(U) \ar[d]_{\rho (u)} \ar[r] &V \ar[d]^v \\
\rho(X) \ar[r]^d &D
}
$$
avec $v \in S$, forment un crible $J$-couvrant de $X$.

\smallskip

Puis, pour de tels $U \xrightarrow{ \ u \ } X$, $V \xrightarrow{ \ v \ } D$ et $\rho (U) \to V$, les morphismes $U' \xrightarrow{ \ u' \ } U$ de ${\mathcal C}$ tels que $\rho (u')$ s'inscrive dans un carr\'e commutatif
$$
\xymatrix{
\rho(U') \ar[d]_{\rho (u')} \ar[r] &V' \ar[d]^{v'} \\
\rho(U) \ar[r]^d &V
}
$$
avec $v' \in v^{-1} S'$, forment un crible $J$-couvrant de $X$.

\smallskip

Donc la transitivit\'e de $\rho_* S$ se ram\`ene \`a celle de $S$.
 
\end{remarkqed}

\bigskip

Cette d\'efinition \'etant pos\'ee, l'\'equivalence des conditions (1) et (5) du th\'eor\`eme \ref{thmIV53} se reformule de la mani\`ere suivante:

\begin{cor}\label{corIV55}

Soient $({\mathcal C},J)$ et $({\mathcal D},K)$ deux sites.

\smallskip

Alors le foncteur
$$
\rho : {\mathcal C} \longrightarrow {\mathcal D}
$$
est un comorphisme de sites si et seulement si $K$ est contenue dans la topologie $\rho_* J$ coinduite par $J$ via $\rho$.
\end{cor}

\bigskip

\begin{remarkqed}

Tout compos\'e de comorphismes de sites est un comorphisme de sites.

\smallskip

Par cons\'equent, pour tous foncteurs entre des cat\'egories essentiellement petites
$$
{\mathcal C} \xrightarrow{ \ \rho \ } {\mathcal D} \xrightarrow{ \ \eta \ } {\mathcal D}'
$$
et toute topologie $J$ sur ${\mathcal C}$, on a
$$
\eta_* (\rho_* J) \subseteq (\eta \circ \rho)_* J \, .
$$
En revanche, il n'y a pas n\'ecessairement \'egalit\'e. \end{remarkqed}

\subsection{Exemples de comorphismes de sites}\label{subsec454}

\medskip

\noindent $\bullet$ Si ${\mathcal C}$ et ${\mathcal D}$ sont deux cat\'egories essentiellement petites, ${\mathcal C}$ est munie d'une topologie $J$ et ${\mathcal D}$ est munie de la topologie discr\`ete, tout foncteur
$$
\rho : {\mathcal C} \longrightarrow {\mathcal D}
$$
est un comorphisme de sites.

\smallskip

Il d\'efinit un morphisme de topos
$$
(\rho^* , \rho_*) : \widehat{\mathcal C}_J \longrightarrow \widehat{\mathcal D}
$$
dont la composante d'image r\'eciproque est la compos\'ee du foncteur de composition avec $\rho$
$$
\begin{matrix}
\rho^* : \widehat{\mathcal D} &\longrightarrow &\widehat{\mathcal C} \, , \hfill \\
\hfill G &\longmapsto &G \circ \rho
\end{matrix}
$$
et du foncteur de faisceautisation
$$
j^* : \widehat{\mathcal C} \longrightarrow \widehat{\mathcal C}_J \, .
$$

Sa composante d'image directe est donc la compos\'ee du foncteur de plongement
$$
j_* : \widehat{\mathcal C}_J \longrightarrow \widehat{\mathcal C}
$$
et du foncteur d'extension de Kan \`a droite
$$
\begin{matrix}
\rho_* : \widehat{\mathcal C} &\longrightarrow &\widehat{\mathcal D} \, , \hfill \\
\hfill F &\longmapsto &\rho_* F = \displaystyle \left[ D \longmapsto \varprojlim_{X \in {\mathcal C} / D} F(X) \right].
\end{matrix}
$$

\noindent $\bullet$ En particulier, pour tout site $({\mathcal C} , J)$, le foncteur
$$
{\rm id} : {\mathcal C} \longrightarrow {\mathcal C}
$$
est un comorphisme du site $({\mathcal C},J)$ vers le site discret de ${\mathcal C}$.

\smallskip

Il d\'efinit le morphisme de topos
$$
(j^* , j_*) : \widehat{\mathcal C}_J \longrightarrow \widehat{\mathcal C} \, .
$$

\noindent $\bullet$ Plus g\'en\'eralement, si $J$ et $K$ sont deux topologies sur une cat\'egorie essentiellement petite ${\mathcal C}$ avec $K \subseteq J$, le foncteur
$$
{\rm id} : {\mathcal C} \longrightarrow {\mathcal C}
$$
est un comorphisme du site $({\mathcal C},J)$ vers le site $({\mathcal C},K)$.

\smallskip

Il d\'efinit le morphisme de topos
$$
\widehat{\mathcal C}_J \longrightarrow \widehat{\mathcal C}_K
$$
dont la composante d'image directe est le foncteur de plongement des $J$-faisceaux dans les $K$-faisceaux
$$
\widehat{\mathcal C}_J \longrightarrow \widehat{\mathcal C}_K
$$
et dont la composante d'image r\'eciproque
$$
\widehat{\mathcal C}_K \longrightarrow \widehat{\mathcal C}_J
$$
est la restriction \`a la sous-cat\'egorie pleine $\widehat{\mathcal C}_K$ de $\widehat{\mathcal C}$ du foncteur de faisceautisation
$$
j^* : \widehat{\mathcal C} \longrightarrow \widehat{\mathcal C}_J \, .
$$

\medskip

\noindent $\bullet$ Pour tout site $({\mathcal C},J)$, l'unique foncteur
$$
p : {\mathcal C} \longrightarrow \un
$$
vers la cat\'egorie $\un$ \`a un objet et un morphisme est un comorphisme de site.

\smallskip

Il d\'efinit un morphisme de topos
$$
(p^* , p_*) : \widehat{\mathcal C}_J \longrightarrow {\rm Ens} \, .
$$
C'est en fait le seul morphisme de n'importe quel topos vers le topos des ensembles:

\begin{prop}\label{propIV56}

Pour tout topos ${\mathcal E}$, la cat\'egorie
$$
[{\mathcal E} , {\rm Ens}]_T
$$
des morphismes de topos
$$
{\mathcal E} \longrightarrow {\rm Ens}
$$
est \'equivalente \`a la cat\'egorie \`a un seul objet et un seul morphisme.

\smallskip

L'unique morphisme de topos
$$
(p^* , p_*) : {\mathcal E} \longrightarrow {\rm Ens}
$$
a pour composante d'image r\'eciproque
$$
\begin{matrix}
p^* : {\rm Ens} &\longrightarrow &{\mathcal E} \, , \hfill \\
\hfill I &\longmapsto &\displaystyle \coprod_I 1_{\mathcal E}
\end{matrix}
$$
et pour composante d'image directe
$$
\begin{matrix}
p_* : {\mathcal E} &\longrightarrow &{\rm Ens} \, , \hfill \\
\hfill E &\longmapsto &{\rm Hom} (1_{\mathcal E} , E) \, ,
\end{matrix}
$$
en notant $1_{\mathcal E}$ l'objet terminal de ${\mathcal E}$.
\end{prop}

\begin{remark}

Si ${\mathcal E} = \widehat{\mathcal C}_J$ est le topos des faisceaux sur un site $({\mathcal C},J)$ et ${\mathcal C}$ poss\`ede un objet terminal $X$, la composante d'image directe de $(p^* , p_*) : \widehat{\mathcal C}_J \to {\rm Ens}$ s'\'ecrit
$$
\begin{matrix}
p_* : \widehat{\mathcal C}_J &\longrightarrow &{\rm Ens} \, , \hfill \\
\hfill F &\longmapsto &F(X) \, .
\end{matrix}
$$
C'est le ``foncteur des sections globales''.

\smallskip

D'autre part, les objets de $\widehat{\mathcal C}_J$ qui sont isomorphes aux objets images du foncteur d'image r\'eciproque
$$
p^* : {\rm Ens} \longrightarrow \widehat{\mathcal C}_J
$$
sont appel\'es les ``faisceaux constants'' sur le site $({\mathcal C},J)$.

\smallskip

Ce sont les sommes de copies de l'objet terminal de $\widehat{\mathcal C}_J$.
\end{remark}
\bigskip

\begin{demo}

Pour tout morphisme de topos
$$
(p^* , p_*) : {\mathcal E} \longrightarrow {\rm Ens} \, ,
$$
le foncteur d'image r\'eciproque
$$
p^* : {\rm Ens} \longrightarrow {\mathcal E}
$$
respecte les limites finies, donc envoie l'objet terminal $\{\bullet\}$ de Ens sur un objet terminal $1_{\mathcal E}$ de ${\mathcal E}$.

\smallskip

Il respecte aussi les colimites donc envoie tout ensemble
$$
I = \coprod_I \ \{\bullet\}
$$
sur $\underset{I}{\coprod} \ 1_{\mathcal E}$.

\smallskip

Enfin, pour tout objet $E$ de ${\mathcal E}$, l'identit\'e d'adjonction
$$
{\rm Hom} (p^* \{\bullet\} , E) = {\rm Hom} (\{\bullet\} , p_* E)
$$
se r\'e\'ecrit
$$
{\rm Hom} (1_{\mathcal E} , E) = p_* E
$$
puisque tout ensemble $I$ s'identifie dans ${\rm Ens}$ \`a ${\rm Hom} (\{\bullet\},I)$. 

\end{demo}

\bigskip

\noindent $\bullet$ Pour tout site $({\mathcal C},J)$ et toute sous-cat\'egorie pleine et dense ${\mathcal C}'$ de ${\mathcal C}$ munie de la topologie induite $J'$, le foncteur de plongement
$$
{\mathcal C}' \xhookrightarrow{ \ { \ } \ } {\mathcal C}
$$
est \`a la fois un morphisme et un comorphisme de sites.

\smallskip

Il d\'efinit les deux \'equivalences $\widehat{\mathcal C}'_{J'} \xrightarrow{ \ \sim \ } \widehat{\mathcal C}_J$ et $\widehat{\mathcal C}_J \xrightarrow{ \ \sim \ } \widehat{\mathcal C}'_{J'}$, r\'eciproques l'une de l'autre.

\medskip

\noindent $\bullet$ On a vu dans la remarque (i) qui suit la proposition \ref{propIV48} que les foncteurs de plongement d'un ``petit site'' dans un ``gros site'' sont des comorphismes de sites (pourvu que la cat\'egorie sous-jacente poss\`ede des limites finies).

\smallskip

En fait, ce sont aussi des morphismes de sites:

\begin{prop}\label{propIV57}

Soit $C$ une cat\'egorie essentiellement petite munie d'une ``classe g\'eom\'etrique de morphismes'' ${\mathcal M}$ et d'une ``notion g\'eom\'etrique de recouvrement'' (R).

\smallskip

Alors, pour tout objet $X$ de ${\mathcal C}$, le foncteur de plongement pleinement fid\`ele
$$
({\mathcal C} / X)_{\mathcal M} \xhookrightarrow{ \ i \ } {\mathcal C}/X
$$
est un morphisme de sites d\'efinis par la topologie $J_R$ sur ${\mathcal C}/X$ ou $({\mathcal C}/X)_{\mathcal M}$.

\smallskip

Il d\'efinit un morphisme de topos
$$
(i^* , i_*) : [ \widehat{({\mathcal C}/X)_{\mathcal M}}]_{J_R} \longrightarrow [\widehat{{\mathcal C}/X}]_{J_R}
$$
dont la composante d'image r\'eciproque
$$
i^* : [\widehat{{\mathcal C}/X}]_{J_R} \longrightarrow [\widehat{({\mathcal C}/X)_{\mathcal M}}]_{J_R}
$$
est le foncteur de restriction \`a $({\mathcal C}/X)_{\mathcal M}$ des faisceaux sur ${\mathcal C}/X$ et dont la composante d'image directe
$$
i_* : [\widehat{({\mathcal C}/X)_{\mathcal M}}]_{J_R} \longrightarrow [\widehat{{\mathcal C}/X}]_{J_R}
$$
associe \`a tout faisceau $F$ sur $({\mathcal C}/X)_{\mathcal M}$ le faisceau sur ${\mathcal C}/X$ d\'efini par la formule
$$
(Z \to X) \longmapsto \varprojlim_{\mbox{\footnotesize$\begin{pmatrix} U \to Z \\ \searrow \ \ \swarrow \\ X\end{pmatrix}$}} F(U \to X) \, .
$$
\end{prop}

\begin{remark}

En particulier, la transformation naturelle canonique
$$
i^* \circ i_* \longrightarrow {\rm id}
$$
est un isomorphisme.

\smallskip

Or, si ${\mathcal C}$ a des limites finies arbitraires, $i$ induit en tant que morphisme de sites un morphisme de topos en sens inverse
$$
[\widehat{{\mathcal C}/X}]_{J_R} \longrightarrow [\widehat{({\mathcal C}/X)_{\mathcal M}}]_{J_R}
$$
dont la composante d'image directe n'est autre que le foncteur de restriction
$$
i^* : [\widehat{{\mathcal C}/X}]_{J_R} \longrightarrow [\widehat{({\mathcal C}/X)_{\mathcal M}}]_{J_R} \, .
$$
Cela signifie que le compos\'e des deux morphismes de topos
$$
[\widehat{({\mathcal C}/X)_{\mathcal M}}]_{J_R} \longrightarrow [\widehat{{\mathcal C}/X}]_{J_R} \longrightarrow [\widehat{({\mathcal C}/X)_{\mathcal M}}]_{J_R}
$$
est canonique isomorphe au foncteur identique.

\smallskip

On dit que le ``petit topos'' $[\widehat{({\mathcal C}/X)_{\mathcal M}}]_{J_R}$ est un r\'etracte du ``gros topos'' $[\widehat{{\mathcal C}/X}]_{J_R}$.
\end{remark}

\bigskip

\begin{demo}

En effet, il est \'evident sur la d\'efinition que le foncteur de restriction des pr\'efaisceaux sur ${\mathcal C}/X$ \`a la sous-cat\'egorie pleine $({\mathcal C}/X)_{\mathcal M}$
$$
i^* : \widehat{{\mathcal C}/X} \longrightarrow \widehat{({\mathcal C}/X)_{\mathcal M}}
$$
transforme les faisceaux en faisceaux et respecte les foncteurs de faisceautisation. 

\end{demo}

\section{Dualit\'e des morphismes et des comorphismes de sites}\label{sec46}

\subsection{Morphismes et comorphismes adjoints}\label{subsec461}

\medskip

On a le r\'esultat g\'en\'eral suivant:

\begin{lem}\label{lemIV61}

Soient $({\mathcal C},J)$ et $({\mathcal D},K)$ deux sites dont les cat\'egories sous-jacentes sont reli\'ees par une paire de foncteurs adjoints
$$
\left({\mathcal C} \xrightarrow{ \ \rho \ } {\mathcal D} \, , \ {\mathcal D} \xrightarrow{ \ \eta \ } {\mathcal C} \right).
$$

Alors:

\begin{listeimarge}

\item Le foncteur compos\'e ${\mathcal D} \xrightarrow{ \ \eta \ } {\mathcal C} \xrightarrow{ \ \ell \ } \widehat{\mathcal C}_J$ est n\'ecessairement plat.

\medskip

\item Le foncteur $\eta$ est un morphisme de sites si et seulement si son adjoint \`a gauche $\rho$ est un comorphisme de sites.

\medskip

\item Dans ces conditions, les morphismes de topos associ\'es
$$
(\eta_! , \eta^!) : \widehat{\mathcal C}_J \longrightarrow \widehat{\mathcal D}_K
$$
et
$$
(\rho^* , \rho_*) : \widehat{\mathcal C}_J \longrightarrow \widehat{\mathcal D}_K
$$
s'identifient.
\end{listeimarge}
\end{lem}

\begin{demo}

On observe que le foncteur
$$
\eta : {\mathcal D} \longrightarrow {\mathcal C}
$$
et le foncteur de composition avec $\rho$
$$
\begin{matrix}
\rho^* : \widehat{\mathcal D} &\longrightarrow &\widehat{\mathcal C} \, , \hfill \\
\hfill G &\longmapsto &G \circ \rho
\end{matrix}
$$
sont compatibles au sens que le carr\'e
$$
\xymatrix{
{\mathcal D}_{ \ } \ar@{_{(}->}[d]_y \ar[r]^{\eta} &{\mathcal C}_{ \ } \ar@{_{(}->}[d]^y \\
\widehat{\mathcal D} \ar[r]^{\rho^*} &\widehat{\mathcal C}
}
$$
est commutatif \`a isomorphisme canonique pr\`es.

\smallskip

Cela r\'esulte en effet de ce que, pour tout objet $X$ de ${\mathcal C}$ et tout objet $D$ de ${\mathcal D}$, les ensembles de morphismes
$$
{\rm Hom} (\rho (X),D) \qquad \mbox{et} \qquad {\rm Hom} (X,\eta (D))
$$
s'identifient par d\'efinition de l'adjonction.

\begin{listeisansmarge}

\item Il r\'esulte de cette observation que l'unique foncteur
$$
\widehat{\mathcal D} \longrightarrow \widehat{\mathcal C}_J
$$
qui prolonge le compos\'e
$$
{\mathcal D} \xrightarrow{ \ \eta \ } {\mathcal C} \xrightarrow{ \ \ell \ } \widehat{\mathcal C}_J
$$
et respecte les colimites est le compos\'e
$$
\widehat{\mathcal D} \xrightarrow{ \ \rho^* \ } \widehat{\mathcal C} \xrightarrow{ \ j^* \ } \widehat{\mathcal C}_J \, .
$$
Il respecte les limites finies puisque le foncteur $\rho^*$ respecte les limites arbitraires.

\smallskip

Cela signifie que ${\mathcal D} \xrightarrow{ \ \eta \ } {\mathcal C} \xrightarrow{ \ \ell \ } \widehat{\mathcal C}_J$ est plat.

\medskip

\item Ainsi, le foncteur $\eta$ est un morphisme de sites si et seulement si son unique prolongement respectant les colimites
$$
\eta_! = \rho^* : \widehat{\mathcal D} \longrightarrow \widehat{\mathcal C}
$$
s'inscrit dans un carr\'e commutatif \`a isomorphisme pr\`es
$$
\xymatrix{
\widehat{\mathcal D} \ar[d]_-{j^*} \ar[rr]^{\eta_! = \rho^*} &&\widehat{\mathcal C} \ar[d]^-{j^*} \\
\widehat{\mathcal D}_K \ar[rr] &&\widehat{\mathcal C}_J
}
$$
c'est-\`a-dire si et seulement si $\rho$ est un comorphisme de sites.

\end{listeisansmarge}
\smallskip

Dans ce cas, les deux foncteurs induits
$$
\widehat{\mathcal D}_K \!\! \raisebox{.7ex}{\xymatrix{\dar[r]^-{^{^{\mbox{\scriptsize$\eta_!$}}}}_-{\rho^*} & \ \widehat{\mathcal C}_J}}
$$
s'identifient, ce qui prouve (iii).

\end{demo}

\pagebreak

\begin{remarkqed}

Pour toute cat\'egorie essentiellement petite ${\mathcal C}$, ${\mathcal C} \xrightarrow{ \ {\rm id} \ } {\mathcal C}$ est son propre adjoint.

\smallskip

C'est ainsi que, pour toutes topologies $J \supseteq K$ sur ${\mathcal C}$, le morphisme de topos
$$
\widehat{\mathcal C}_J \longrightarrow \widehat{\mathcal C}_K
$$
peut \^etre vu comme associ\'e au comorphisme de sites $({\mathcal C},J) \to ({\mathcal C},K)$ aussi bien qu'au morphisme de sites $({\mathcal C},K) \to ({\mathcal C},J)$. 

\end{remarkqed}

\bigskip

Nous allons montrer un peu plus loin que tout morphisme de sites se factorise canoniquement comme le compos\'e de  deux morphismes de sites dont le premier a un adjoint \`a gauche -- auquel s'applique le lemme ci-dessus -- et le second induit une \'equivalence de topos.

\smallskip

Cela permettra de passer de tout morphisme de sites repr\'esentant un morphisme de topos \`a une repr\'esentation alternative par un comorphisme de sites.

\subsection{Les cat\'egories c\^ones d'un foncteur}\label{subsec462}

\medskip

Nous allons utiliser la factorisation canonique de tout foncteur \`a travers la cat\'egorie c\^one \`a gauche de ce foncteur d\'efinie de la mani\`ere suivante:

\begin{defn}\label{defIV62}

Soit
$$
\rho : {\mathcal C} \longrightarrow {\mathcal D}
$$
un foncteur entre deux cat\'egories.

\smallskip

On appelle c\^one \`a droite [resp. \`a gauche] de $\rho$ et on note
$$
\qquad \qquad {\mathcal C} /_{\!\!\rho} \, {\mathcal D} \qquad \mbox{ou plus bri\`evement} \qquad {\mathcal C}/{\mathcal D}
$$
$$
\mbox{[resp. \qquad ${\mathcal D}_{\rho}\!\backslash \, {\mathcal C}$ \qquad ou plus bri\`evement  \qquad ${\mathcal D} \backslash {\mathcal C}$ \ ]}
$$
la cat\'egorie dont

\medskip

$\left\lmoustache\begin{matrix}
\bullet &\mbox{les objets sont les triplets $(X,D,\rho (X) \to D)$ [resp. $(X,D,D \to \rho (X))$] consitu\'es d'un objet $X$ de ${\mathcal C}$,} \hfill \\
&\mbox{d'un objet $D$ de ${\mathcal D}$ et d'un morphisme de ${\mathcal D}$} \hfill \\
{ \ } \\
&\rho (X) \longrightarrow D \qquad \mbox{[resp.} \quad D \longrightarrow \rho (X) \ \mbox{]}, \\
{ \ } \\
\end{matrix}\right.$

$\left\rmoustache\begin{matrix}
\bullet &\mbox{les morphismes} \hfill \\
&(X',D',\rho (X') \to D') \longrightarrow (X,D,\rho (X) \to D) \\
&\mbox{[resp.} \hfill \\
&(X',D',D' \to \rho (X')) \longrightarrow (X,D,D \to \rho (X)) \ \mbox{]} \\
{ \ } \\
&\mbox{sont les paires $(X' \xrightarrow{ \ x \ } X , D' \xrightarrow{ \ d \ } D)$ constitu\'ees d'un morphisme de ${\mathcal C}$} \hfill \\
{ \ } \\
&X' \xrightarrow{ \ x \ } X \\
&\mbox{et d'un morphisme de ${\mathcal D}$} \hfill \\
&D' \xrightarrow{ \ d \ } D \\
&\mbox{qui rendent commutatif le carr\'e de ${\mathcal D}$} \hfill \\
{ \ } \\
&\xymatrix{
\rho (X') \ar[d]_-{\rho (x)} \ar[r] &D' \ar[d]^-{d} \\
\rho(X) \ar[r] &D
} \qquad \quad
\xymatrix{
D' \ar[d]_{\mbox{[resp.} \qquad d} \ar[r] &\rho(X') \ar[d]^-{\rho(x) \ \mbox{]}.} \\
D \ar[r] &\rho(X)
} 
\end{matrix} \right.
$

\bigskip

On note $\pi_{\mathcal D}$ et $\pi_{\mathcal C}$ les foncteurs de ``projection'' sur ${\mathcal D}$ et ${\mathcal C}$ qui associent

\bigskip

$\left\{\begin{matrix}
\bullet &\mbox{\`a tout objet $(X,D,\rho(X) \to D)$ [resp. $(X,D,D \to \rho(X))$] de ${\mathcal C}/{\mathcal D}$ [resp. ${\mathcal D} \backslash {\mathcal C}$] les objets $D$ et $X$} \hfill \\
&\mbox{de ${\mathcal D}$ et ${\mathcal C}$,} \hfill \\
{ \ } \\
\bullet &\mbox{\`a tout morphisme $(X' \xrightarrow{ \ x \ } X , D' \xrightarrow{ \ d \ } D)$ de ${\mathcal C}/{\mathcal D}$ [resp. ${\mathcal D} \backslash {\mathcal C}$], les morphismes} \hfill \\
{ \ } \\
&D' \xrightarrow{ \ d \ } D \quad \mbox{de} \quad {\mathcal D} \\
&\mbox{et} \hfill \\
&X' \xrightarrow{ \ x \ } X \quad \mbox{de} \quad {\mathcal C} \, .
\end{matrix} \right.
$

\bigskip

Enfin, on note $i_{\mathcal C}$ le foncteur de ${\mathcal C}$ dans ${\mathcal C} / {\mathcal D}$ [resp. ${\mathcal D} \backslash {\mathcal C}$] qui associe

\bigskip

$\left\{\begin{matrix}
\bullet &\mbox{\`a tout objet $X$ de ${\mathcal C}$ le triplet} \hfill \\
{ \ } \\
&\left(X,\rho (X),\rho (X) \xrightarrow{ \ {\rm id} \ } \rho (X)\right), \\
{ \ } \\
\bullet &\mbox{\`a tout morphisme $X' \xrightarrow{ \ x \ } X$ de ${\mathcal C}$ la paire de morphismes de ${\mathcal C}$ et ${\mathcal D}$} \hfill \\
{ \ } \\
&\left(X' \xrightarrow{ \ x \ } X , \rho (X') \xrightarrow{ \ \rho (x) \ } \rho (X) \right).
\end{matrix} \right.
$
\end{defn}

\bigskip

\begin{remarksqed}
\begin{listeisansmarge}
\item On a
$$
\pi_{\mathcal D} \circ i_{\mathcal C} = \rho
$$
et
$$
\pi_{\mathcal C} \circ i_{\mathcal C} = {\rm id}_{\mathcal C} \, .
$$
De plus, le foncteur
$$
i_{\mathcal C} : {\mathcal C} \longrightarrow {\mathcal C}/{\mathcal D} \qquad \mbox{[resp.} \quad i_{\mathcal C} : {\mathcal C} \longrightarrow {\mathcal D} \backslash {\mathcal C} \ \mbox{]}
$$ 
est adjoint \`a gauche [resp. \`a droite] du foncteur
$$
\pi_{\mathcal C} : {\mathcal C}/{\mathcal D} \longrightarrow {\mathcal C} \qquad \mbox{[resp.} \quad \pi_{\mathcal C} : {\mathcal D} \backslash {\mathcal C} \longrightarrow {\mathcal C} \ \mbox{]}.
$$ 

\item Si les cat\'egories ${\mathcal C}$ et ${\mathcal D}$ sont localement petites [resp. petites, resp. essentiellement petites], il en est de m\^eme des cat\'egories ${\mathcal C}/{\mathcal D}$ et ${\mathcal D} \backslash {\mathcal C}$.

\medskip

\item Dans la litt\'erature math\'ematique, les cat\'egories c\^ones de $\rho$
$$
{\mathcal C}/{\mathcal D} = {\mathcal C}/_{\!\!\rho} \, {\mathcal D} \qquad \mbox{et} \qquad {\mathcal D}\backslash{\mathcal C} = {\mathcal D}_{\rho}\!\backslash \, {\mathcal C} 
$$
sont souvent appel\'ees les cat\'egories ``comma'' et not\'ees
$$
(\rho \downarrow {\mathcal D}) = {\mathcal C}/{\mathcal D} \qquad \mbox{et} \qquad ({\mathcal D} \downarrow \rho) = 
{\mathcal D}\backslash{\mathcal C} \, .
$$
\end{listeisansmarge}
\end{remarksqed}

\subsection{Des morphismes de sites aux comorphismes}\label{subsec463}

\medskip

On a vu dans les remarques (vi) et (vii) qui suivent la d\'efinition \ref{defIV41} que tout morphisme de topos se pr\'esente de multiples fa\c cons comme induit par des morphismes de sites.

\smallskip

La factorisation d'un morphisme de sites \`a travers son c\^one \`a gauche permet de repr\'esenter n'importe quel morphisme de topos induit par un morphisme de sites, et donc n'importe quel morphisme de topos, comme induit par un comorphisme de sites:

\begin{thm}[Caramello (Théorème 3.16. \cite{Denseness})]\label{thmIV63}

Soient $({\mathcal C},J)$ et $({\mathcal D},K)$ deux sites reli\'es par un foncteur
$$
\rho : {\mathcal C} \longrightarrow {\mathcal D}
$$
qui est un morphisme de sites et induit donc un morphisme de topos
$$
(\rho_! , \rho^!) : \widehat{\mathcal D}_K \longrightarrow \widehat{\mathcal C}_J \, .
$$

Consid\'erons la factorisation de $\rho$ \`a travers son c\^one \`a gauche
$$
\rho : {\mathcal C} \xrightarrow{ \ i_{\mathcal C} \ } {\mathcal D} \backslash {\mathcal C} \xrightarrow{ \ \pi_{\mathcal D} \ } {\mathcal D}
$$
et munissons ${\mathcal D} \backslash {\mathcal C}$ de la topologie $\widetilde K = \pi_{\mathcal D}^{-1} K$ induite par la topologie $K$ de ${\mathcal D}$ via $\pi_{\mathcal D}$.

\smallskip

Alors:

\begin{listeimarge}

\item Une famille de morphismes de ${\mathcal D} \backslash {\mathcal C}$
$$
(X_i , D_i , D_i \to \rho (X_i)) \longrightarrow (X,D,D \to \rho (X)) \, , \quad i \in I \, ,
$$
est couvrante pour la topologie $\widetilde K$ si et seulement si son image par le foncteur $\pi_{\mathcal D}$
$$
D_i \longrightarrow D \, , \qquad i \in I \, ,
$$
est couvrante pour la topologie $K$.

\medskip

\item Le foncteur
$$
\pi_{\mathcal D} : {\mathcal D} \backslash {\mathcal C} \longrightarrow {\mathcal D}
$$
est \`a la fois un morphisme et un comorphisme du site $({\mathcal D} \backslash {\mathcal C} , \widetilde K)$ dans le site $({\mathcal D},K)$.

\smallskip

Il induit deux morphismes de topos
$$
\widehat{\mathcal D}_K \longrightarrow \widehat{({\mathcal D} \backslash {\mathcal C})}_{\widetilde K}
$$
et
$$
\widehat{({\mathcal D} \backslash {\mathcal C})}_{\widetilde K} \longrightarrow \widehat{\mathcal D}_K
$$
qui sont des \'equivalences r\'eciproques l'une de l'autre.

\medskip

\item Le foncteur
$$
i_{\mathcal C} : {\mathcal C} \longrightarrow {\mathcal D} \backslash {\mathcal C}
$$
est un morphisme de sites
$$
({\mathcal C},J) \longrightarrow ({\mathcal D} \backslash {\mathcal C} , \widetilde K) \, ,
$$
et son adjoint \`a gauche
$$
\pi_{\mathcal C} : {\mathcal D} \backslash {\mathcal C} \longrightarrow {\mathcal C}
$$
est un comorphisme de sites
$$
({\mathcal D} \backslash {\mathcal C} , \widetilde K) \longrightarrow ({\mathcal C},J) \, .
$$
Ils induisent un m\^eme morphisme de topos
$$
\widehat{({\mathcal D} \backslash {\mathcal C})}_{\widetilde K} \longrightarrow \widehat{\mathcal C}_J \, .
$$

\item Le morphisme de topos
$$
\widehat{\mathcal D}_K \longrightarrow \widehat{\mathcal C}_J
$$
induit par le morphisme de sites $\rho : ({\mathcal C},J) \to ({\mathcal D},K)$ se factorise comme le compos\'e de l'\'equivalence
$$
\widehat{\mathcal D}_K \xrightarrow{ \ \sim \ } \widehat{({\mathcal D} \backslash {\mathcal C})}_{\widetilde K}
$$
et du morphisme de topos
$$
\widehat{({\mathcal D} \backslash {\mathcal C})}_{\widetilde K} \longrightarrow \widehat{\mathcal C}_J
$$
induit par le comorphisme de sites
$$
\pi_{\mathcal C} : ({\mathcal D} \backslash {\mathcal C} , \widetilde K) \longrightarrow ({\mathcal C},J) \, .
$$
\end{listeimarge}
\end{thm}

\bigskip

\begin{demosansqed}
\begin{listeisansmarge}
\item Notons $\widetilde K$ la collection des cribles d'objets $(X,D,D \to \rho (X))$ de ${\mathcal D} \backslash {\mathcal C}$ qui contiennent une famille de morphismes
$$
(X_i , D_i , D_i \to (\rho (X_i)) \longrightarrow (X,D,D \to \rho(X)) \, , \quad i \in I \, ,
$$
tels que la famille induite
$$
D_i \longrightarrow D \, , \qquad i \in I \, ,
$$
soit $K$-couvrante.

\smallskip

Pour montrer que $\widetilde K$ est la topologie induite $\pi_{\mathcal D}^{-1} K$ par $K$ via $\pi_{\mathcal D}$, il suffit d'apr\`es la d\'efinition \ref{defIV411} de v\'erifier que, pour tout crible $S$ d'un objet $(X,D,D \to \rho (X))$ de ${\mathcal D} \backslash {\mathcal C}$ qui est \'el\'ement de $\widetilde K$, on a:

\medskip

$
\left\{ \begin{matrix}
(1) &\mbox{Pour tout morphisme de ${\mathcal D} \backslash {\mathcal C}$} \hfill \\
{ \ } \\
&(X',D',D' \to \rho (X')) \xrightarrow{ \ (x,d) \ } (X,D,D \to \rho (X)) \, , \\
{ \ } \\
&\mbox{le crible $(x,d)^{-1} S$ de $(X',D',D' \to \rho (X'))$ est encore \'el\'ement de $\widetilde K$.} \hfill \\
{ \ } \\
(2) &\mbox{Toute paire d'objets de $S$} \hfill \\
&\xymatrix{
&(X'' , D'' , D'' \to \rho (X'')) \ar[d] \\
(X' , D' , D' \to \rho (X')) \ar[r] &(X,D,D \to \rho (X))
} \\
{ \ } \\
&\mbox{peut \^etre compl\'et\'ee en une famille de carr\'es commutatifs de ${\mathcal D} \backslash {\mathcal C}$} \hfill \\
{ \ } \\
&\xymatrix{
(X_i , D_i , D_i \to \rho (X_i)) \ar[d] \ar[r] &(X'' , D'' , D'' \to \rho (X'')) \ar[d] \\
(X' , D' , D' \to \rho (X')) \ar[r] &(X,D,D \to \rho (X))
} \\
{ \ } \\
&\mbox{telle que, notant} \hfill \\
&\ell : {\mathcal D} \xrightarrow{ \ y \ } \widehat{\mathcal D} \xrightarrow{ \ j^* \ } \widehat{\mathcal D}_K \\
{ \ } \\
&\mbox{le foncteur canonique, la famille correspondante des} \hfill \\
{ \ } \\
&\ell (D_i) \longrightarrow \ell (D') \times_{\ell (D)} \ell (D'') \\
{ \ } \\
&\mbox{soit (globalement) \'epimorphique.} \hfill
\end{matrix} \right.
$

\medskip

V\'erifions d'abord (1).

\smallskip

Comme le crible de $D$ engendr\'e par les images par $\rho$ des \'el\'ements de $S$ est $K$-couvrant, son image r\'eciproque par $d : D' \to D$ est un crible $K$-couvrant de $D'$.

\smallskip

Cela signifie qu'existe une famille $K$-couvrante de morphismes de ${\mathcal D}'$
$$
D'_i \longrightarrow D' \, , \qquad i \in I \, ,
$$
qui s'inscrivent dans des carr\'es commutatifs de $D'$
$$
\xymatrix{
D'_i \ar[d] \ar[r] &D' \ar[d] \\
D_i \ar[r]^{d_i} &D
}
$$
dont les composantes
$$
D_i \xrightarrow{ \ d_i \ } D
$$
sont les images par $\pi_{\mathcal D}$ de morphismes de ${\mathcal D} \backslash {\mathcal C}$
$$
(X_i,D_i,D_i \to \rho (X_i)) \xrightarrow{ \ (x_i,d_i) \ } (X,D,D \to \rho (X))
$$
qui sont \'el\'ements du crible $S$.

\smallskip

Pour tous indices $i$, les compos\'es
$$
D'_i \longrightarrow D_i \longrightarrow \rho (X_i)
$$
et
$$
D'_i \longrightarrow D' \longrightarrow \rho (X')
$$
s'inscrivent dans un carr\'e commutatif de ${\mathcal D}$:
$$
\xymatrix{
D'_i \ar[d] \ar[r] &\rho(X') \ar[d] \\
\rho(X_i) \ar[r] &\rho(X)
}
$$
Or, le foncteur
$$
{\mathcal C} \xrightarrow{ \ \rho \ } {\mathcal D} \xrightarrow{ \ \ell \ } \widehat{\mathcal D}_K
$$
est plat.

\smallskip

On en d\'eduit que, quitte \`a remplacer chaque $D'_i$ par une famille $K$-couvrante au-dessus de $D'_i$, on peut supposer qu'existent des objets $X'_i$ de ${\mathcal C}$, des morphismes de ${\mathcal C}$
$$
X_i \longleftarrow X'_i \longrightarrow X'
$$
et des morphismes de ${\mathcal D}$
$$
D'_i \longrightarrow \rho (X'_i)
$$
tels que $D'_i \to \rho (X_i)$ et $D'_i \to \rho (X')$ se factorisent en
$$
D'_i \longrightarrow \rho (X'_i) \longrightarrow \rho (X_i)
$$
et
$$
D'_i \longrightarrow \rho (X'_i) \longrightarrow \rho (X') \, .
$$

Puis, toujours d'apr\`es la platitude de ${\mathcal C} \to \widehat{\mathcal D}_K$, on peut supposer, quitte \`a remplacer \`a nouveau chaque $D'_i$ par une famille $K$-couvrante au-dessus de $D'_i$ et chaque objet $X'_i$ par d'autres objets au-dessus de $X'_i$, que les carr\'es de ${\mathcal C}$
$$
\xymatrix{
X'_i \ar[d] \ar[r] &X' \ar[d] \\
X_i \ar[r] &X
}
$$
sont commutatifs.

\smallskip

On a ainsi construit une famille de carr\'es commutatifs de ${\mathcal D} \backslash {\mathcal C}$
$$
\xymatrix{
(X'_i , D'_i , D'_i \to \rho (X'_i)) \ar[d] \ar[r] &(X',D',D' \to \rho (X')) \ar[d] \\
(X_i , D_i , D_i \to \rho (X_i)) \ar[r] &(X,D,D \to \rho (X))
}
$$
dont les composantes
$$
(X_i , D_i , D_i \to \rho (X_i))\longrightarrow(X,D,D \to \rho (X))
$$
sont \'el\'ements du crible $S$, et qui induisent une famille $K$-couvrante de morphismes de ${\mathcal D}$
$$
D'_i \longrightarrow D' \, .
$$

Cela ach\`eve la v\'erification de (1).

\smallskip

Puis v\'erifions (2).

\smallskip

Le diagramme de ${\mathcal D}$
$$
\xymatrix{
&D'' \ar[d] \\
D' \ar[r] &D
}
$$
peut \^etre compl\'et\'e en une famille de carr\'es commutatifs de ${\mathcal D}$
$$
\xymatrix{
D_i \ar[d] \ar[r] &D'' \ar[d] \\
D' \ar[r] &D
}
$$
telle que la famille des morphismes induits de $\widehat{\mathcal D}_K$
$$
\ell (D_i) \longrightarrow \ell (D') \times_{\ell (D)} \ell (D'')
$$
soit (globalement) \'epimorphique.

\smallskip

Les compos\'es $D_i \to D' \to \rho (X')$ et $D_i \to D'' \to \rho (X'')$ s'inscrivent dans des carr\'es commutatifs:
$$
\xymatrix{
D_i \ar[d] \ar[r] &\rho (X'') \ar[d] \\
\rho (X') \ar[r] &\rho (X)
}
$$
Utilisant \`a nouveau la platitude du foncteur
$$
{\mathcal C} \xrightarrow{ \ \rho \ } {\mathcal D} \xrightarrow{ \ \ell \ } \widehat{\mathcal D}_K \, ,
$$
on peut supposer, quitte \`a remplacer chaque $D_i$ par une famille $K$-couvrante au-dessus de $D_i$, que le triangle de ${\mathcal C}$
$$
\xymatrix{
&X'' \ar[d] \\
X' \ar[r] &X
}
$$
se compl\`ete en des carr\'es commutatifs de ${\mathcal C}$
$$
\xymatrix{
X_i \ar[d] \ar[r] &X'' \ar[d] \\
X' \ar[r] &X
}
$$
tels que, pour tout indice $i$, les deux morphismes
$$
D_i \longrightarrow \rho (X')
$$
et
$$
D_i \longrightarrow \rho (X'')
$$
se factorisent comme les compos\'es d'un morphisme
$$
D_i \longrightarrow \rho (X_i)
$$
avec les images par $\rho$ des deux morphismes
$$
X' \longleftarrow X_i \longrightarrow X'' \, .
$$

Autrement dit, on a une famille de carr\'es commutatifs de ${\mathcal D} \backslash {\mathcal C}$
$$
\xymatrix{
(X_i , D_i , D_i \to \rho (X_i)) \ar[d] \ar[r] &(X'' , D'' , D'' \to \rho (X'')) \ar[d] \\
(X' , D' , D' \to \rho (X')) \ar[r] &(X,D,D \to \rho (X))
}
$$
qui r\'epond \`a la question pos\'ee.

\smallskip

Cela ach\`eve la v\'erification de (2) et donc de la partie (i) du th\'eor\`eme.

\medskip

\item Comme $\widetilde K = \pi_{\mathcal D}^{-1} K$ est la topologie induite par $K$ dans ${\mathcal D} \backslash {\mathcal C}$ via $\pi_{\mathcal D}$, on sait d\'ej\`a que le foncteur
$$
\pi_{\mathcal D} : ({\mathcal D} \backslash {\mathcal C} , \widetilde K) \longrightarrow ({\mathcal D},K) 
$$
est continu. Pour montrer que c'est un morphisme de sites, il reste seulement \`a v\'erifier la platitude du foncteur compos\'e
$$
{\mathcal D} \backslash {\mathcal C} \xrightarrow{ \ \pi_{\mathcal D} \ } {\mathcal D} \xrightarrow{ \ \ell \ } \widehat{\mathcal D}_K \, .
$$
Elle r\'esulte de la platitude du foncteur
$$
{\mathcal C} \xrightarrow{ \ \rho \ } {\mathcal D} \xrightarrow{ \ \ell \ } \widehat{\mathcal D}_K \, .
$$
En effet, celle-ci s'exprime d'apr\`es la proposition \ref{propIV39} par les trois conditions suivantes:

\medskip

$
\left\lmoustache\begin{matrix}
({\rm C1}) &\mbox{Tout objet $D$ de ${\mathcal D}$ admet une famille $K$-couvrante} \hfill \\
{ \ } \\
&D_i \longrightarrow D \, , \qquad i \in I \, , \\
{ \ } \\
&\mbox{telle que chaque $D_i$ soit l'image par $\pi_{\mathcal D}$ d'un objet} \hfill \\
{ \ } \\
&(X_i , D_i , D_i \to \rho (X_i)) \quad \mbox{de} \quad {\mathcal D} \backslash {\mathcal C} \, . \\
{ \ } \\
({\rm C2}) &\mbox{Pour tout objet $D$ de ${\mathcal D}$ muni de deux morphismes} \hfill \\
{ \ } \\
&D' \longleftarrow D \longrightarrow D'' \\
{ \ } \\
&\mbox{vers les images $D'$ et $D''$ par $\pi_{\mathcal D}$ de deux objets} \hfill \\
{ \ } \\
&(X' , D' , D' \to \rho (X')) \\
&\mbox{et} \hfill \\
&(X'' , D'' , D'' \to \rho (X'')) \\
{ \ } \\
&\mbox{de ${\mathcal D} \backslash {\mathcal C}$, il existe une famille $K$-couvrante de $D$} \hfill \\
{ \ } \\
&D_i \longrightarrow D \, , \qquad i \in I \, , \\
{ \ } \\
&\mbox{telle que chaque $D_i$ soit l'image par $\pi_{\mathcal D}$ d'un objet} \hfill \\
{ \ } \\
&(X_i , D_i , D_i \to \rho (X_i)) \quad \mbox{de} \quad {\mathcal D} \backslash {\mathcal C} \, ,\\

&\mbox{et que les morphismes compos\'es} \hfill \\
{ \ } \\
&D_i \longrightarrow D \longrightarrow D' \\
&\mbox{et} \hfill \\
&D_i \longrightarrow D \longrightarrow D'' \\
{ \ } \\
&\mbox{se rel\`event en des morphismes de ${\mathcal D} \backslash {\mathcal C}$} \hfill \\
{ \ } \\
&\qquad (X_i , D_i , D_i \to \rho (X_i)) \longrightarrow (X' , D' , D' \to \rho (X')) \\
&\mbox{et} \hfill  \\
&\qquad (X_i , D_i , D_i \to \rho (X_i)) \longrightarrow (X'' , D'' , D'' \to \rho (X'')) \, . \\
{ \ } \\
\end{matrix} \right.
$

$
\left\rmoustache \begin{matrix}
({\rm C3}) &\mbox{Pour toute paire de morphismes de ${\mathcal D} \backslash {\mathcal C}$} \hfill \\
{ \ } \\
&(X',D',D' \to \rho (X')) \rightrightarrows (X,D,D \to \rho (X)) \\
{ \ } \\
&\mbox{et tout morphisme de ${\mathcal D}$} \hfill \\
{ \ } \\
&D'' \longrightarrow D' \\
{ \ } \\
&\mbox{tel que les deux compos\'es} \hfill \\
{ \ } \\
&D'' \longrightarrow D' \rightrightarrows D \\
{ \ } \\ 
&\mbox{soient \'egaux, il existe une famille $K$-couvrante} \hfill \\
{ \ } \\
&\qquad \quad D_i \longrightarrow D'' \, , \qquad i \in I \, , \\
{ \ } \\
&\mbox{tel que chaque $D_i$ soit l'image par $\pi_{\mathcal D}$ d'un objet de ${\mathcal D} \backslash {\mathcal C}$} \hfill \\
{ \ } \\
&\qquad \quad (X_i , D_i , D_i \to \rho (X_i)) \\
{ \ } \\
&\mbox{muni d'un morphisme} \hfill \\
{ \ } \\
&\qquad \quad (X_i , D_i , D_i \to \rho (X_i)) \longrightarrow (X',D',D' \to \rho (X')) \\
{ \ } \\
&\mbox{dont les compos\'es avec les deux morphismes} \hfill \\
{ \ } \\
&\qquad \qquad (X',D',D' \to \rho (X')) \rightrightarrows (X,D,D \to \rho (X)) \\
&\mbox{sont \'egaux.} \hfill
\end{matrix} \right.
$

\bigskip

Ainsi, le foncteur 
$$
\pi_{\mathcal D} : {\mathcal D} \backslash {\mathcal C} \longrightarrow {\mathcal D}
$$
induit un morphisme de topos
$$
(\pi_{{\mathcal D}!} , \pi_{\mathcal D}^!) : \widehat{\mathcal D}_K \longrightarrow \widehat{({\mathcal D} \backslash {\mathcal C})}_{\widetilde K} \, .
$$

On pr\'etend d'autre part que la topologie $K$ de ${\mathcal D}$ est contenue dans la topologie coinduite par la topologie $\widetilde K$ de ${\mathcal D} \backslash {\mathcal C}$ via le foncteur $\pi_{\mathcal D}$.

\smallskip

Consid\'erons en effet un crible $K$-couvrant $S$ d'un objet de ${\mathcal D}$ et un objet $(X',D',D' \to \rho (X'))$ de ${\mathcal D} \backslash {\mathcal C}$ muni d'un morphisme
$$
\pi_{\mathcal D} (X',D',D' \to \rho (X')) = D' \xrightarrow{ \ d \ } D \, .
$$
Le crible $d^{-1} S$ de $D'$ est $K$-couvrant donc contient une famille $K$-couvrante de morphismes de ${\mathcal D}$
$$
D'_i \longrightarrow D' \, , \qquad i \in I \, .
$$
Les morphismes compos\'es
$$
D'_i \longrightarrow D' \longrightarrow \rho (X')
$$
d\'efinissent des objets de ${\mathcal D} \backslash {\mathcal C}$
$$
(X',D'_i,D'_i \longrightarrow \rho (X'))
$$
munis de morphismes
$$
(X',D'_i,D'_i \to \rho (X')) \longrightarrow (X',D',D' \to \rho (X'))
$$
qui rel\`event les morphismes de ${\mathcal D}'$
$$
D'_i \longrightarrow D' \, .
$$
Cela montre que le crible $S$ de $D$ est couvrant pour la topologie de ${\mathcal D}$ coinduite par $\widetilde K$ via $\pi_{\mathcal D} : {\mathcal D} \backslash {\mathcal C} \to {\mathcal D}$.

\smallskip

Comme la topologie $K$ est contenue dans la topologie coinduite par $\widetilde K$, le foncteur
$$
\pi_{\mathcal D} : {\mathcal D} \backslash {\mathcal C} \longrightarrow {\mathcal D}
$$
est aussi un comorphisme de sites.

\smallskip

Il induit un morphisme de topos
$$
(\pi_{\mathcal D}^* , \pi_{{\mathcal D}*}) : \widehat{({\mathcal D} \backslash {\mathcal C})}_{\widetilde K} \longrightarrow \widehat{\mathcal D}_K \, .
$$

Pour terminer la preuve de la partie (ii) du th\'eor\`eme, nous allons appliquer le r\'esultat g\'en\'eral suivant:
\end{listeisansmarge}
\end{demosansqed}

\begin{lem}\label{lemIV64}

Soient $({\mathcal A},J)$ et $({\mathcal B},K)$ deux sites.

\smallskip

Soit
$$
\pi : {\mathcal A} \longrightarrow {\mathcal B}
$$
un foncteur qui est \`a la fois un morphisme de sites et un comorphisme de sites.

\smallskip

Pour que les deux morphismes de topos induits
$$
(\pi_! , \pi^!) : \widehat{\mathcal B}_K \longrightarrow \widehat{\mathcal A}_J \, ,
$$
$$
(\pi^* , \pi_*) : \widehat{\mathcal A}_J \longrightarrow \widehat{\mathcal B}_K
$$
soient des \'equivalences, il faut [resp. il suffit] que soient satisfaites les conditions suivantes [resp. les conditions {\rm (2), (3)} et {\rm (4)}]:

\begin{listeimarge}

\item[(1)] La topologie $K$ de ${\mathcal B}$ est coinduite par la topologie $J$ de ${\mathcal A}$.

\medskip

\item[(2)] La topologie $J$ de ${\mathcal A}$ est induite par la topologie $K$ de ${\mathcal B}$.

\medskip

\item[(3)] Le foncteur $\pi$ est ``$J$-plein'' au sens que pour tous objets $X,X'$ de ${\mathcal A}$ et pour tout morphisme de ${\mathcal B}$
$$
\pi (X') \longrightarrow \pi(X) \, ,
$$
existe une famille $J$-couvrante de morphismes de ${\mathcal A}$
$$
X'_i \longrightarrow X' \, , \qquad i \in I \, ,
$$
telle que les morphismes compos\'es
$$
\pi (X'_i) \longrightarrow \pi (X') \longrightarrow \pi(X)
$$
soient les images par $\pi$ de morphismes de ${\mathcal A}$.

\medskip

\item[(4)] Le foncteur $\pi$ est ``$K$-dense'' au sens que tout objet $D$ de ${\mathcal B}$ admet une famille $K$-couvrante de morphismes
$$
\pi (X_i) \longrightarrow D \, , \qquad i \in I \, ,
$$
dont les sources sont les images par $\pi$ d'objets $X_i$ de ${\mathcal A}$.
\end{listeimarge}
\end{lem}

\bigskip

\begin{demolem}

Prouvons d'abord que la condition (1) est cons\'equence de la condition (4) sous l'hypoth\`ese que $\pi : {\mathcal A} \to {\mathcal B}$ est \`a la fois un morphisme et un comorphisme de sites.

\smallskip

En effet, comme $\pi$ est un comorphisme de sites, on sait d\'ej\`a que $K$ est contenue dans la topologie $\pi_* J$ coinduite par $J$.

\smallskip

R\'eciproquement, consid\'erons un crible $S$ d'un objet $D$ de ${\mathcal B}$ qui est couvrant pour la topologie coinduite $\pi_* J$.

\smallskip

D'apr\`es la condition (4), il existe une famille $K$-couvrante de $D$ de la forme
$$
d_i : \pi (X_i) \longrightarrow D \, , \qquad i \in I \, ,
$$
pour des objets $X_i$ de ${\mathcal A}$.

\smallskip

Puis, par d\'efinition de la topologie coinduite, il existe pour tout indice $i$ une famille $J$-couvrante de $X_i$ dans ${\mathcal A}$
$$
X_{i,j} \longrightarrow X_i
$$
dont les images $\pi (X_{i,j}) \to \pi (X_i)$ sont \'el\'ements de $d_i^{-1} S$.

\smallskip

Comme $\pi$ est un morphisme de sites, il pr\'eserve les familles couvrantes si bien que les familles de morphismes de ${\mathcal B}$
$$
\pi (X_{i,j}) \longrightarrow \pi (X_i)
$$
sont $K$-couvrantes.

\smallskip

Cela prouve comme voulu que le crible $S$ est \'el\'ement de $K$ et donc que $\pi_* J$ est contenue dans $K$.

\smallskip

On a bien $\pi_* J = K$.

\smallskip

Puis montrons que $(\pi_! , \pi^!)$ et $(\pi^* , \pi_*)$ sont deux \'equivalences r\'eciproques si et seulement si $\pi$ satisfait les conditions (1) \`a (4).

\smallskip

Comme $\pi$ est un morphisme de sites, il est continu ce qui signifie que $J$ est contenue dans la topologie $J' = \pi^{-1} K$ induite par $K$ via $\pi$.

\smallskip

D'autre part, le fait que $\pi$ soit un comorphisme de sites signifie que $K$ est contenue dans la topologie $K' = \pi_* J$ coinduite par $J$ via $\pi$.

\smallskip

Soit $\widehat\pi : \widehat{\mathcal A} \to \widehat{\mathcal B}$ l'unique prolongement de ${\mathcal A} \xrightarrow{ \ \pi \ } {\mathcal B} \xrightarrow{ \ y \ } \widehat{\mathcal B}$ qui respecte les colimites.

\smallskip

Le foncteur compos\'e
$$
\widehat{\mathcal A} \xrightarrow{ \ \widehat\pi \ } \widehat{\mathcal B} \xrightarrow{ \ j^* \ } \widehat{\mathcal B}_K
$$
se factorise \`a travers le quotient $\widehat{\mathcal A}_{J'}$ de $\widehat{\mathcal A}_J$. Si le foncteur
$$
\pi_! : \widehat{\mathcal A}_J \longrightarrow \widehat{\mathcal B}_K
$$
est une \'equivalence, cela impose $J = J' = \pi^{-1} K$.

\smallskip

De m\^eme, le foncteur compos\'e
$$
\widehat{\mathcal B} \xrightarrow{ \ \pi^* \ } \widehat{\mathcal A} \xrightarrow{ \ j^* \ } \widehat{\mathcal A}_J
$$
se factorise \`a travers le quotient $\widehat{\mathcal B} _{K'}$ de $\widehat{\mathcal B}_K$. Si le foncteur
$$
\pi^* : \widehat{\mathcal B}_K \longrightarrow \widehat{\mathcal A}_J
$$
est une \'equivalence, cela impose $K = K' = \pi_* J$.

\smallskip

R\'eciproquement, supposons que $J = \pi^{-1} K$ et $K = \pi_* J$.

\smallskip

Montrons qu'alors les deux foncteurs
$$
\pi_! : \widehat{\mathcal A}_J \longrightarrow \widehat{\mathcal B}_K \qquad \mbox{et} \qquad \pi_* : \widehat{\mathcal B}_K \longrightarrow \widehat{\mathcal A}_J
$$
sont fid\`eles.

\smallskip

Comme ils pr\'eservent les limites finies, en particulier les \'egalisateurs, il suffit de prouver que tout monomorphisme
$$
F' \xhookrightarrow{ \ { \ } \ } F \qquad \mbox{de} \qquad \widehat{\mathcal A}_J \qquad \mbox{[resp.} \quad \widehat{\mathcal B}_K \ \mbox{]}
$$
que $\pi_!$ [resp. $\pi_*$] transforme en isomorphisme de $\widehat{\mathcal B}_K$ [resp. $\widehat{\mathcal A}_J$] est un isomorphisme.

\smallskip

Comme tout faisceau $F$ est une colimite de faisceautis\'es de pr\'efaisceaux repr\'esentables, il suffit de prouver que pour tout objet $X$ de ${\mathcal A}$ [resp. ${\mathcal B}$] et tout crible
$$
S \xhookrightarrow{ \ { \ } \ } y(X) \, ,
$$
son image par
$$
\widehat{\mathcal A} \xrightarrow{ \ \widehat\pi \ } \widehat{\mathcal B} \xrightarrow{ \ j^* \ } \widehat{\mathcal B}_K \qquad 
\mbox{[resp.} \quad \widehat{\mathcal B} \xrightarrow{ \ \pi^* \ } \widehat{\mathcal A} \xrightarrow{ \ j^* \ } \widehat{\mathcal A}_J \ \mbox{]}
$$
est un isomorphisme (si et) seulement si son image par
$$
\widehat{\mathcal A} \xrightarrow{ \ j^* \ } \widehat{\mathcal A}_J \qquad \mbox{[resp.} \quad \widehat{\mathcal B} \xrightarrow{ \ j^* \ } \widehat{\mathcal B}_K \ \mbox{]}
$$
est un isomorphisme.

\smallskip

Il en est ainsi car $J = \pi^{-1} K$ [resp. $K = \pi_* J$].

\smallskip

La fid\'elit\'e du foncteur
$$
\pi_! : \widehat{\mathcal A}_J \longrightarrow \widehat{\mathcal B}_K \qquad \mbox{[resp.} \quad \pi^* : \widehat{\mathcal B}_K \longrightarrow \widehat{\mathcal A}_J \ \mbox{]}
$$
signifie que pour tous faisceaux $F,F'$ sur $({\mathcal A},J)$ [resp. $G,G'$ sur $({\mathcal B},K)$], l'application
$$
{\rm Hom} (F',F) \longrightarrow {\rm Hom} (\pi_! F' , \pi_! F) = {\rm Hom} (F' , \pi^! \circ \pi_! F)
$$
[resp.
$$
{\rm Hom} (G',G) \longrightarrow {\rm Hom} (\pi^* G' , \pi^* G) = {\rm Hom} (G' , \pi_* \circ \pi^* G) \ \mbox{]}
$$
est injective. Autrement dit, elle signifie que pour tout faisceau $F$ sur $({\mathcal A},J)$ [resp. $G$ sur $({\mathcal B},K)$], le morphisme canonique
$$
F \longrightarrow \pi^! \circ \pi_! F \qquad \mbox{[resp.} \quad G \longrightarrow \pi_* \circ \pi^* G \ \mbox{]}
$$
est un monomorphisme.

\smallskip

Pour que les foncteurs $\pi_!$, $\pi^! \cong \pi^*$ et $\pi_*$ soient des \'equivalences, il faut et il suffit que tous ces monomorphismes canoniques
$$
F \longrightarrow \pi^! \circ \pi_! F \qquad \mbox{et} \qquad G \longrightarrow \pi_* \circ \pi^* G
$$
soient aussi des \'epimorphismes, donc des isomorphismes.

\smallskip

Les foncteurs $\pi_!$ et $\pi^! \cong \pi^*$ respectent les colimites et tout faisceau $F$ sur $({\mathcal A},J)$ s'\'ecrit comme une colimite de faisceaux de la forme $\ell (X) = j^* {\rm Hom} (\bullet , X)$, $X \in {\rm Ob} ({\mathcal A})$.

\smallskip

Donc les morphismes
$$
F \longrightarrow \pi^! \circ \pi_! F
$$
sont tous des \'epimorphismes si et seulement si, pour tout objet $X$ de ${\mathcal A}$, le morphisme de $\widehat{\mathcal A}$
$$
{\rm Hom} (\bullet , X) \longrightarrow {\rm Hom} (\pi(\bullet) , \pi(X))
$$
est transform\'e par $j^* : \widehat{\mathcal A} \to \widehat{\mathcal A}_J$ en un \'epimorphisme de $\widehat{\mathcal A}$.

\smallskip

D'apr\`es le crit\`ere du lemme \ref{lemII56} (i), cela signifie exactement que pour tous objets $X,X'$ de ${\mathcal A}$ et tout morphisme de ${\mathcal B}$
$$
\pi(X') \longrightarrow \pi(X)
$$
existe une famille $J$-couvrante de $X'$
$$
X'_i \longrightarrow X' \, , \qquad i \in I \, ,
$$
telle que tous les compos\'es
$$
\pi (X'_i) \longrightarrow \pi(X') \longrightarrow \pi(X)
$$
soient les images par $\pi$ de morphismes $X'_i \to X$ de ${\mathcal A}$.

\smallskip

Autrement dit, cela signifie que le foncteur $\pi : {\mathcal A} \to {\mathcal B}$ est $J$-plein.

\smallskip

Enfin, les morphismes
$$
G \longrightarrow \pi_* \circ \pi^* G = \Biggl[ D \longmapsto \varprojlim_{(X , \pi(X) \to D) \atop \in {\mathcal A} / D} G(\pi (X)) \Biggl]
$$
sont tous des \'epimorphismes si et seulement si tout objet $D$ de ${\mathcal B}$ admet une famille $K$-couvrante
$$
\pi (X_i) \longrightarrow D
$$
dont les sources sont les images par $\pi$ d'objets $X_i$ de ${\mathcal A}$. Cela revient \`a demander que le foncteur $\pi : {\mathcal A} \to {\mathcal B}$ soit $K$-dense.

\smallskip

Cela ach\`eve la d\'emonstration du lemme. 

\end{demolem}

\bigskip

\noindent {\bf Suite de la d\'emonstration du th\'eor\`eme \ref{thmIV63}:}

\begin{listeisansmarge}

\item[] Pour achever de montrer (ii), il reste d'apr\`es ce lemme \`a v\'erifier que le foncteur
$$
\pi_{\mathcal D} : {\mathcal D} \backslash {\mathcal C} \longrightarrow {\mathcal D}
$$
est $\widetilde K$-plein et $K$-dense.

\smallskip

Il est $K$-dense d'apr\`es la propri\'et\'e (C1) \'enonc\'ee plus haut.

\smallskip

Pour v\'erifier qu'il est $\widetilde K$-plein, consid\'erons deux objets $(X,D,D \to \rho (X))$, $(X',D',D' \to \rho (X'))$ de ${\mathcal D} \backslash {\mathcal C}$ et un morphisme entre leurs images par $\pi_{\mathcal D}$ dans ${\mathcal D}$
$$
D' \longrightarrow D \, .
$$

D'apr\`es la propri\'et\'e (C2), il existe une famille $K$-couvrante de morphismes de ${\mathcal D}$
$$
D'_i \longrightarrow D' \, , \qquad i \in I \, ,
$$
telle que chaque $D'_i$ soit l'image par $\pi_{\mathcal D}$ d'un objet $(X'_i , D'_i , D'_i \to \rho (X'_i))$ de ${\mathcal D} \backslash {\mathcal C}$ et que les morphismes de ${\mathcal D}$
$$
D'_i \longrightarrow D' \qquad \mbox{et} \qquad D'_i \longrightarrow D
$$
se rel\`event en des morphismes de ${\mathcal D} \backslash {\mathcal C}$
$$
(X'_i , D'_i , D'_i \to \rho(X'_i)) \longrightarrow (X',D',D' \to \rho (X'))
$$
et
$$
(X'_i , D'_i , D'_i \to \rho(X'_i)) \longrightarrow (X,D,D \to \rho (X)) \, .
$$

Enfin, la famille des morphismes
$$
(X'_i , D'_i , D'_i \to \rho(X'_i)) \longrightarrow (X',D',D' \to \rho (X'))
$$
est $\widetilde K$-couvrante puisque la famille image par $\pi_{\mathcal D}$
$$
D'_i \longrightarrow D'
$$
est $K$-couvrante.

\smallskip

Cela prouve que $\pi_{\mathcal D}$ est $\widetilde K$-plein et ach\`eve de d\'emontrer (ii).

\medskip

\item[(iii)] Le foncteur compos\'e de
$$
i_{\mathcal C} : {\mathcal C} \longrightarrow {\mathcal D} \backslash {\mathcal C} \qquad \mbox{et de} \qquad {\mathcal D} \backslash {\mathcal C} \xrightarrow{ \ \ell \ } \widehat{({\mathcal D} \backslash {\mathcal C})}_{\widetilde K}
$$
est plat d'apr\`es le lemme \ref{lemIV61} (i) puisque le foncteur $i_{\mathcal C}$ a un adjoint \`a gauche.

\smallskip

De plus, si
$$
X_i \longrightarrow X \, , \qquad i \in I \, ,
$$
est une famille $J$-couvrante de morphismes de ${\mathcal C}$, la famille image par $\rho$
$$
\rho (X_i) \longrightarrow \rho (X) \, , \qquad i \in I \, ,
$$
est une famille $K$-couvrante de morphismes de ${\mathcal D}$, et donc la famille des
$$
i_{\mathcal C} (X_i) = (X_i , \rho (X_i) , \rho (X_i) \xrightarrow{ \ {\rm id} \ } \rho (X_i)) \longrightarrow (X,\rho(X),\rho (X) \xrightarrow{ \ {\rm id} \ } \rho(X)) = i_{\mathcal C} (X)
$$
est $\widetilde K$-couvrante.

\smallskip

Il r\'esulte alors du th\'eor\`eme \ref{thmIV43} que le foncteur
$$
i_{\mathcal C} : {\mathcal C} \longrightarrow {\mathcal D} \backslash {\mathcal C}
$$
est un morphisme du site $({\mathcal C},J)$ vers le site $({\mathcal D}\backslash{\mathcal C},\widetilde K)$.

\smallskip

D'apr\`es le lemme \ref{lemIV61}, son adjoint \`a gauche
$$
\pi_{\mathcal C} : {\mathcal D} \backslash {\mathcal C} \longrightarrow {\mathcal C}
$$
est un comorphisme de sites, et tous deux induisent le m\^eme morphisme de topos
$$
\widehat{({\mathcal D} \backslash {\mathcal C})}_{\widetilde K} \longrightarrow \widehat{\mathcal C}_J \, .
$$

\item[(iv)] r\'esulte de (ii) et (iii) puisque le morphisme de sites 
$$
\rho : {\mathcal C} \longrightarrow {\mathcal D}
$$
est le compos\'e des morphismes de sites 
$$
i_{\mathcal C} : {\mathcal C} \longrightarrow {\mathcal D} \backslash {\mathcal C}
$$
et
$$
\pi_{\mathcal D} : {\mathcal D} \backslash {\mathcal C} \longrightarrow {\mathcal D} \, .
$$

Cela ach\`eve la d\'emonstration du th\'eor\`eme.\end{listeisansmarge} \hfill $\Box$

\pagebreak 

\section{Sous-topos}\label{sec47}

\subsection{La notion de plongement de topos et celle de sous-topos}\label{subsec471}

\medskip

On pose:

\begin{defn}\label{defIV71}
\begin{listeimarge}
\item Un morphisme de topos
$$
f = (f^*,f_*) : {\mathcal E}' \longrightarrow {\mathcal E}
$$
est appel\'e un ``plongement'' si sa composante d'image directe
$$
f_* : {\mathcal E}' \longrightarrow {\mathcal E}
$$
est un foncteur pleinement fid\`ele.

\medskip

\item Deux plongements de topos dans un m\^eme topos
$$
{\mathcal E}'_1 \longrightarrow {\mathcal E}
$$
et
$$
{\mathcal E}'_2 \longrightarrow {\mathcal E}
$$
sont dits \'equivalents s'il existe une \'equivalence de topos
$$
{\mathcal E}'_1 \xrightarrow{ \ \sim \ } {\mathcal E}'_2
$$
qui rende le triangle
$$
\xymatrix{
{\mathcal E}'_1 \ar[dd]^{\wr} \ar[rd] \\
&{\mathcal E} \\
{\mathcal E}'_2 \ar[ru]
}
$$
commutatif \`a isomorphisme pr\`es.

\medskip

\item On appelle sous-topos d'un topos ${\mathcal E}$ les classes d'\'equivalence de plongements de topos dans ${\mathcal E}$
$$
{\mathcal E}' \longrightarrow {\mathcal E} \, .
$$
\end{listeimarge}
\end{defn}

\begin{remarksqed}
\begin{listeisansmarge}
\item Pour signifier qu'un morphisme de topos
$$
{\mathcal E}' \longrightarrow {\mathcal E} 
$$
est un plongement, on pourra le noter
$$
{\mathcal E}' \xhookrightarrow{ \ { \ } \ } {\mathcal E} \, .
$$

\item Le compos\'e de deux plongements de topos
$$
{\mathcal E}'' \xhookrightarrow{ \ { \ } \ } {\mathcal E}' \xhookrightarrow{ \ { \ } \ } {\mathcal E}
$$
est un plongement de topos.

\medskip

\item R\'eciproquement, si
$$
{\mathcal E}'_1 \xhookrightarrow{ \ { \ } \ } {\mathcal E} 
$$
et
$$
{\mathcal E}'_2 \xhookrightarrow{ \ { \ } \ } {\mathcal E} 
$$
sont deux plongements de topos, tout morphisme de topos
$$
{\mathcal E}'_1 \longrightarrow {\mathcal E}'_2
$$
qui rend le triangle
$$
\xymatrix{
{\mathcal E}'_1 \ar[dd] \ar@{^{(}->}[rd] \\
&{\mathcal E} \\
{\mathcal E}'_2 \ar@{^{(}->}[ru]
}
$$
commutatif \`a isomorphisme pr\`es, est un plongement de topos.

\smallskip

De plus, deux morphismes de topos
$$
{\mathcal E}'_1 \rightrightarrows {\mathcal E}'_2
$$
qui sont chacun compatibles en ce sens avec les plongements de ${\mathcal E}'_1$ et ${\mathcal E}'_2$ dans ${\mathcal E}$ sont canoniquement isomorphes.

\medskip

\item En particulier, pour tout plongement de topos
$$
{\mathcal E}' \xhookrightarrow{ \ { \ } \ } {\mathcal E} \, ,
$$
tout morphisme de topos
$$
{\mathcal E}' \longrightarrow {\mathcal E}'
$$
qui rend le triangle
$$
\xymatrix{
{\mathcal E}' \ar[dd] \ar@{^{(}->}[rd] \\
&{\mathcal E} \\
{\mathcal E}' \ar@{^{(}->}[ru]
}
$$
commutatif \`a isomorphisme pr\`es, est canoniquement isomorphe au morphisme
$$
{\rm id} : {\mathcal E}' \longrightarrow {\mathcal E}' \, .
$$

\item Toute \'equivalence de topos
$$
{\mathcal E}_1 \xrightarrow{ \ \sim \ } {\mathcal E}_2
$$
induit une correspondance bijective entre les sous-topos de ${\mathcal E}_1$ et les sous-topos de ${\mathcal E}_2$.

\medskip

\item Par exemple, pour tout site $({\mathcal C},J)$, le morphisme de topos
$$
(j^* , j_*) : \widehat{\mathcal C}_J \longrightarrow \widehat{\mathcal C}
$$
est un plongement.

\smallskip

Si $K$ est n'importe quelle topologie sur ${\mathcal C}$ contenue dans $J$, il se factorise en un plongement
$$
\widehat{\mathcal C}_J \xhookrightarrow{ \ { \ } \ } \widehat{\mathcal C}_K
$$
qui rend le triangle
$$
\xymatrix{
\widehat{\mathcal C}_J \ar@{_{(}->}[dd] \ar@{^{(}->}[rd] \\
&\widehat{\mathcal C} \\
\widehat{\mathcal C}_K \ar@{^{(}->}[ru]
}
$$
commutatif \`a isomorphisme canonique pr\`es.

\smallskip

Nous allons voir au paragraphe suivant qu'en fait tous les plongements de topos se repr\'esentent sous cette forme. 
\end{listeisansmarge}
\end{remarksqed}

\subsection{Sous-topos et topologies de Grothendieck}\label{subsec472}

\medskip

Le th\'eor\`eme suivant montre que la notion intrins\`eque de sous-topos d'un topos donn\'e correspond, du c\^ot\'e des pr\'esentations de ce topos par des sites, aux topologies de Grothendieck des cat\'egories sous-jacentes:

\begin{thm}\label{thmIV72}

Soit ${\mathcal E} = \widehat{\mathcal C}_J$ le topos des faisceaux sur un site $({\mathcal C},J)$. 

\smallskip

Alors:

\begin{listeimarge}

\item L'application
$$
K \longmapsto \widehat{\mathcal C}_K
$$
d\'efinit une bijection de l'ensemble des topologies $K$ de ${\mathcal C}$ contenant $J$ sur la classe des sous-topos de ${\mathcal E}$.

\medskip

\item La bijection r\'eciproque consiste \`a associer \`a  tout sous-topos de ${\mathcal E}$
$$
(j^* , j_*) : {\mathcal E}' \xhookrightarrow{ \ { \ } \ } {\mathcal E}
$$
la topologie $K$ de ${\mathcal C}$ pour laquelle un crible $S$ d'un objet $X$ de ${\mathcal C}$ est couvrant si et seulement si, pour tout objet $F$ de ${\mathcal E}'$, l'application
$$
(j_* F)(X) \longrightarrow \varprojlim_{(U \to X) \in S} (j_* F)(U)
$$
est bijective.

\medskip

\item Cette bijection renverse la relation d'ordre au sens que, pour toutes topologies $K_1$ et $K_2$ de ${\mathcal C}$ contenant $J$, on a
$$
K_1 \subseteq K_2
$$
si et seulement si le plongement
$$
\widehat{\mathcal C}_{K_2} \xhookrightarrow{ \ { \ } \ } \widehat{\mathcal C}_J
$$
se factorise en un triangle
$$
\xymatrix{
\widehat{\mathcal C}_{K_2} \ar@{_{(}->}[dd] \ar@{^{(}->}[rd] \\
&\widehat{\mathcal C}_J \\
\widehat{\mathcal C}_{K_1} \ar@{^{(}->}[ru]
}
$$
commutatif \`a isomorphisme canonique pr\`es.

\end{listeimarge}
\end{thm}

\bigskip

\begin{remarks}
\begin{listeisansmarge}
\item En particulier, les sous-topos d'un topos donn\'e forment un ensemble ordonn\'e.

\smallskip

C'est un ``invariant'' de ce topos au sens que toute \'equivalence de topos
$$
{\mathcal E}' \xrightarrow{ \ \sim \ } {\mathcal E}
$$
induit une bijection respectant l'ordre entre les ensembles de sous-topos associ\'es.

\medskip

\item Il r\'esulte de ce th\'eor\`eme que pour n'importe quelle paire de pr\'esentations d'un topos ${\mathcal E}$ comme topos de faisceaux sur des sites $({\mathcal C},J)$ et $({\mathcal D},K)$
$$
\widehat{\mathcal C}_J \cong {\mathcal E} \cong \widehat{\mathcal D}_K \, ,
$$
l'ensemble ordonn\'e des topologies de ${\mathcal C}$ contenant $J$ s'identifie \`a l'ensemble ordonn\'e des topologies de ${\mathcal D}$ contenant $K$.

\medskip

\item Par exemple, si on note ${\mathcal E}_X$ le topos des faisceaux sur un espace topologique $X$, tout sous-espace de $X$
$$
j : Y \xhookrightarrow{ \ { \ } \ } X
$$
muni de la topologie induite par celle de $X$, d\'efinit un sous-topos
$$
(j^* , j_*) : {\mathcal E}_Y \xhookrightarrow{ \ { \ } \ } {\mathcal E}_X \, .
$$
La topologie de $O(X)$ qui correspond \`a ce sous-topos est celle pour laquelle une famille d'immersions entre ouverts de $X$
$$
(U_i \xhookrightarrow{ \ { \ } \ } U)_{i \in I}
$$
est couvrante si et seulement si
$$
\bigcup_{i \in I} \ (U_i \cap Y) = U \cap Y \, .
$$
En revanche, le topos ${\mathcal E}_X$ poss\`ede en g\'en\'eral des sous-topos
$$
{\mathcal E}' \xhookrightarrow{ \ { \ } \ } {\mathcal E}_X
$$
qui ne sont associ\'es \`a aucun sous-espace de $X$.

\smallskip

Il en est ainsi en g\'en\'eral du sous-topos de ${\mathcal E}_X$ d\'efini par la ``topologie de la densit\'e'' de $O(X)$, pour laquelle une famille d'immersions
$$
(U_i \xhookrightarrow{ \ { \ } \ } U)_{i \in I}
$$
est couvrante si et seulement si la r\'eunion $\underset{i \in I}{\bigcup} \, U_i$ est dense dans $U$.
\end{listeisansmarge}
\end{remarks}

\bigskip

\begin{demosansqed}

D'apr\`es les remarques (ii) et (iii) qui suivent la d\'efinition \ref{defIV71}, il suffit de traiter le cas o\`u ${\mathcal E} = \widehat{\mathcal C}$ est le topos des faisceaux sur une petite cat\'egorie ${\mathcal C}$.

\smallskip

Toute topologie $J$ de ${\mathcal C}$ d\'efinit un sous-topos
$$
(j^* , j_*) : \widehat{\mathcal C}_J \xhookrightarrow{ \ { \ } \ } \widehat{\mathcal C} = {\mathcal E} \, .
$$

L'existence d'une application en sens inverse r\'esulte du lemme suivant:
\end{demosansqed}

\begin{lem}\label{lemIV73}

Soit ${\mathcal C}$ une cat\'egorie essentiellement petite munie du foncteur de Yoneda $y : {\mathcal C} \hookrightarrow \widehat{\mathcal C}$.

\smallskip

Soit
$$
(f^* , f_*) : {\mathcal E}' \longrightarrow \widehat{\mathcal C}
$$
un morphisme d'un topos ${\mathcal E}'$ dans $\widehat{\mathcal C}$.

\smallskip

Associons \`a tout objet $X$ de ${\mathcal C}$ l'ensemble $J(X)$ des cribles $S$ de $X$ qui v\'erifient les trois conditions \'equivalentes suivantes:

\medskip

$
\left\{\begin{matrix}
(1) &\mbox{Pour tout objet $F$ de ${\mathcal E}'$, l'application} \hfill \\
{ \ } \\
&f_* F (X) \longrightarrow \displaystyle\varprojlim_{(U \to X) \in S} f_* F(U) \\
&\mbox{est une bijection.} \hfill \\
{ \ } \\
(2) &\mbox{Le morphisme de ${\mathcal E}$} \hfill \\
{ \ } \\
&\displaystyle\varinjlim_{(U \to X) \in S} f^* \!\circ y(U) \longrightarrow f^* \!\circ y(X) \\
{ \ } \\
&\mbox{est un isomorphisme.} \hfill \\
{ \ } \\
(3) &\mbox{Si $S$ est vu comme un sous-objet de $y(X)$ dans $\widehat{\mathcal C}$, le foncteur $f^*$ transforme le monomorphisme} \hfill \\
{ \ } \\
&S \xhookrightarrow{ \ { \ } \ } y(X) \\
&\mbox{en un isomorphisme de ${\mathcal E}'$.} \hfill
\end{matrix} \right.
$

\medskip

Alors $J$ est une topologie de ${\mathcal C}$.
\end{lem}

\pagebreak

\begin{demolem}

L'\'equivalence de (2) et (3) r\'esulte de ce que $S$, vu comme un pr\'efaisceau sur ${\mathcal C}$, se repr\'esente dans $\widehat{\mathcal C}$ comme la colimite
$$
S = \varinjlim_{(U \to X) \in S} y(U)
$$
puisque le foncteur d'image r\'eciproque
$$
f^* : \widehat{\mathcal C} \longrightarrow {\mathcal E}'
$$
respecte les colimites.

\smallskip

L'\'equivalence de (1) et (2) r\'esulte de ce que l'identit\'e de ${\mathcal E}'$
$$
\varinjlim_{(U \to X) \in S} f^* \!\circ y(U) = f^* \!\circ y(X)
$$
\'equivaut, d'apr\`es le lemme de Yoneda, \`a ce qu'on ait pour tout objet $F$ de ${\mathcal E}'$
$$
{\rm Hom} (f^* \!\circ y(X),F) = \varprojlim_{(U \to X) \in S} {\rm Hom} (f^* \!\circ y(U),F)
$$
c'est-\`a-dire, par adjonction de la paire $(f^*,f_*)$,
$$
{\rm Hom} (y(X) , f_* F) = \varprojlim_{(U \to X) \in S} {\rm Hom} (y(U),f_* F)
$$
qui s'\'ecrit encore d'apr\`es le lemme de Yoneda
$$
f_* F (X) = \varprojlim_{(U \to X) \in S} f_* F(U) \, .
$$

Montrons maintenant que $J$ est une topologie.

\smallskip

Il est \'evident que $J$ satisfait l'axiome de maximalit\'e.

\smallskip

Pour la stabilit\'e, consid\'erons un morphisme $X' \xrightarrow{ \ x \ } X$ de ${\mathcal C}$ et un crible $S \in J(X)$ vu comme un sous-objet $S \hookrightarrow y(X)$.

\smallskip

Alors le crible $x^{-1} S$ est le sous-objet $S \times_{y(X)} y(X') \hookrightarrow y(X')$.

\smallskip

Si $f^*$ transforme $S \hookrightarrow y(X)$ en un isomorphisme de ${\mathcal E}'$, il fait de m\^eme avec $S \times_{y(X)} y(X') \hookrightarrow y(X')$ puisqu'il respecte les limites finies.

\smallskip

Pour la transitivit\'e, consid\'erons deux cribles $S,S'$ d'un objet $X$ de ${\mathcal C}$ tels que $S \in J(X)$ et $u^{-1} S' \in J(U)$ pour tout \'el\'ement $U \xrightarrow{ \ u \ } X$ de $S$.

\smallskip

Les identit\'es de ${\mathcal E}'$
$$
\varinjlim_{(U \xrightarrow{u} X) \in S} f^* \!\circ y (U) = f^* \!\circ y(X)
$$
et, pour tout \'el\'ement $U \xrightarrow{ \ u \ } X$ de $S$,
$$
f^* (S' \times_{y(X)} y(U)) \xrightarrow{ \ \sim \ } f^* \!\circ y (U)
$$
impliquent les identifications
\begin{eqnarray}
f^* S' &= &\varinjlim_{(U \xrightarrow{u} X) \in S} f^* S' \times_{f^* \!\circ y(X)} f^* \!\circ y(U) \nonumber \\
&= &\varinjlim_{(U \xrightarrow{u} X) \in S} f^* (S' \times_{y(X)} y(U)) \nonumber \\
&= &\varinjlim_{(U \xrightarrow{u} X) \in S} f^* \!\circ y(U) \nonumber \\
&= &f^* \!\circ y(X) \nonumber
\end{eqnarray}
puisque les produits fibr\'es respectent les colimites dans le topos ${\mathcal E}'$ et que le foncteur $f^*$ respecte les limites finies.

\smallskip

Cela ach\`eve de montrer que $J$ est une topologie de ${\mathcal C}$. 

\end{demolem}

\bigskip

\noindent {\bf Suite de la d\'emonstration du th\'eor\`eme \ref{thmIV72}:}

\smallskip

Le lemme ci-dessus associe \`a tout sous-topos
$$
(j^* , j_*) : {\mathcal E}' \xhookrightarrow{ \ { \ } \ } \widehat{\mathcal C}
$$
une topologie $J_{{\mathcal E}'}$ de ${\mathcal C}$.

\smallskip

Il est clair que les deux applications
$$
\begin{matrix}
\hfill J &\longmapsto &(\widehat{\mathcal C}_J \hookrightarrow \widehat{\mathcal C}) \\
\mbox{et} \qquad ({\mathcal E}' \hookrightarrow {\mathcal E}) &\longmapsto &J_{{\mathcal E}'} \hfill
\end{matrix}
$$
renversent les relations d'ordre.

\smallskip

Il reste \`a prouver qu'elles sont r\'eciproques l'une de l'autre.

\smallskip

D'apr\`es l'\'equivalence des conditions (A) et (B) du lemme \ref{lemIII53}, l'application compos\'ee
$$
J \longmapsto (\widehat{\mathcal C}_J \hookrightarrow \widehat{\mathcal C}) \longmapsto J_{\widehat{\mathcal C}_J}
$$
est l'identit\'e.

\smallskip

En sens inverse, consid\'erons un sous-topos
$$
(j^* , j_*) : {\mathcal E}' \xhookrightarrow{ \ { \ } \ } \widehat{\mathcal C}
$$
et la topologie $J = J_{{\mathcal E}'}$ qui lui est associ\'ee.

\smallskip

L'image par $j_*$ de tout objet de ${\mathcal E}'$ est un $J$-faisceau, ce qui signifie que le foncteur
$$
j_* : {\mathcal E}' \longrightarrow \widehat{\mathcal C}
$$
se factorise en
$$
{\mathcal E}' \xrightarrow{ \ f_* \ } \widehat{\mathcal C}_J \xrightarrow{ \ j_* \ } \widehat{\mathcal C}
$$
o\`u $f_*$ respecte les limites arbitraires, donc admet un adjoint \`a gauche $f^*$.

\smallskip

En passant aux adjoints, le foncteur
$$
j^* : \widehat{\mathcal C} \longrightarrow {\mathcal E}'
$$
s'inscrit dans un triangle
$$
\xymatrix{
\widehat{\mathcal C} \ar[d]_-{j^*} \ar[r]^-{j^*} &{\mathcal E}' \\
\widehat{\mathcal C}_J \ar[ru]_-{f^*}
}
$$
commutatif \`a isomorphisme canonique pr\`es.

\smallskip

Le foncteur
$$
f_* : {\mathcal E}' \longrightarrow \widehat{\mathcal C}_J
$$
est pleinement fid\`ele, et donc le morphisme de foncteurs ${\mathcal E}' \to {\mathcal E}'$
$$
f^* \!\circ f_* \longrightarrow {\rm id}_{{\mathcal E}'}
$$
est un isomorphisme.

\smallskip

Il reste \`a montrer que le morphisme de foncteurs $\widehat{\mathcal C}_J \to \widehat{\mathcal C}_J$
$$
{\rm id}_{\widehat{\mathcal C}_J} \longrightarrow f_* \circ f^*
$$
est un isomorphisme, autrement dit que, pour tout faisceau $F$ sur $({\mathcal C},J)$ et tout objet $X$ de ${\mathcal C}$, l'application
$$
F(X) \longrightarrow f_* \circ f^* F(X)
$$
est bijective.

\smallskip

Consid\'erons d'abord deux \'el\'ements $x_1 , x_2 \in F(X)$ qui ont m\^eme image dans $f_* \circ f^* F(X)$.

\smallskip

On peut les voir comme deux morphismes
$$
y(X) \raisebox{.7ex}{\xymatrix{\dar[r]^-{^{^{\mbox{\scriptsize$x_1$}}}}_-{x_2} & \, F}}
$$
qui ont les m\^emes transform\'es par le foncteur
$$
j^* : \widehat{\mathcal C} \longrightarrow {\mathcal E}' \, .
$$
Leur \'egalisateur est un crible
$$
S \xhookrightarrow{ \ { \ } \ } y(X) \, .
$$
Comme $j^*$ respecte les \'egalisateurs, il transforme le monomorphisme
$$
S \xhookrightarrow{ \ { \ } \ } y(X)
$$
en un isomorphisme de ${\mathcal E}'$, ce qui signifie que le crible $S$ est $J$-couvrant.

\smallskip

Comme $F$ est un $J$-faisceau, cela impose
$$
x_1 = x_2 \, .
$$
Ainsi, l'application
$$
F(X) \longrightarrow f_* \circ f^* F(X)
$$
est injective.

\smallskip

Enfin, consid\'erons un \'el\'ement $x$ de $f_* \circ f^* (X)$ vu comme un morphisme
$$
y(X) \longrightarrow f_* \circ f^* F \, .
$$
Le produit fibr\'e
$$
S = y(X) \times_{f_* \circ f^* F} F
$$ 
est un sous-objet de $y(X)$, c'est-\`a-dire un crible de $X$.

\smallskip

Son transform\'e par le foncteur $f^*$ est
$$
f^* \!\circ y(X) \times_{f^* \!\circ f_* \circ f^* F} f^* F
$$
qui s'identifie \`a
$$
f^* \!\circ y(X)
$$
puisque le compos\'e $f^* \!\circ f_*$ est isomorphe \`a l'identit\'e de ${\mathcal E}'$.

\smallskip

Cela signifie que $S$ est un crible $J$-couvrant de $X$.

\smallskip

Pour tout \'el\'ement $U \xrightarrow{ \ u \ } X$ de $S$, l'image de $x$ par l'application de restriction
$$
f_* \circ f^* F(X) \longrightarrow f_* \circ f^* F(U)
$$
se rel\`eve par construction en un \'el\'ement de $F(U)$.

\smallskip

Cela prouve que
$$
F \longrightarrow f_* \circ f^* F
$$
est un \'epimorphisme et donc un isomorphisme.

\smallskip

Ainsi, $f^* : \widehat{\mathcal C}_J \to {\mathcal E}'$ et $f_* : {\mathcal E}' \to \widehat{\mathcal C}_J$ sont bien des \'equivalences de cat\'egories. \hfill $\Box$

\bigskip

On observe que, \'etant donn\'ee une cat\'egorie essentiellement petite ${\mathcal C}$, le lemme \ref{lemIV73} permet d'associer \`a tout morphisme de topos
$$
{\mathcal E}' \longrightarrow \widehat{\mathcal C}
$$
une topologie sur ${\mathcal C}$.

\smallskip

Nous allons en d\'eduire une factorisation canonique de tout morphisme de topos.

\smallskip

Pour cela, nous avons besoin de la d\'efinition suivante:

\begin{defn}\label{defIV74}

Un morphisme de topos
$$
f = (f^* , f_*) : {\mathcal E}' \longrightarrow {\mathcal E}
$$
est appel\'e ``surjectif'' s'il poss\`ede les trois propri\'et\'es \'equivalentes suivantes:

\begin{listeimarge}

\item[(A)] Le foncteur $f^*$ est fid\`ele.

\medskip

\item[(B)] Le foncteur $f^*$ est ``conservatif'' au sens qu'un morphisme de ${\mathcal E}$
$$
u : G \longrightarrow F
$$
est un isomorphisme si et seulement si son image par $f^*$
$$
f^* (u) : f^* (G) \longrightarrow f^* (F)
$$
est un isomorphisme de ${\mathcal E}'$.

\medskip

\item[(C)] Pour tout objet $F$ de ${\mathcal E}$, le morphisme canonique
$$
F \longrightarrow f_* \circ f^* F
$$
est un monomorphisme.
\end{listeimarge}
\end{defn}

\begin{remarksqed}
\begin{listeisansmarge}
\item L'\'equivalence de (A) et (C) r\'esulte de ce que, pour tous objets $F$ et $G$ de ${\mathcal E}$,
$$
{\rm Hom} (f^* G , f^* F)
$$
s'identifie par adjonction \`a
$$
{\rm Hom} (G , f_* \circ f^* F) \, .
$$

Donc l'application
$$
{\rm Hom} (G,F) \longrightarrow {\rm Hom} (f^* G , f^* F)
$$
est injective pour tous objets $F$ et $G$ si et seulement si le morphisme canonique associ\'e \`a tout objet $F$
$$
F \longrightarrow f_* \circ f^* F
$$
est un monomorphisme.

\medskip

\item Le foncteur $f^*$ respecte les limites finies.

\smallskip

En particulier, pour toute paire de morphismes de ${\mathcal E}$
$$
G \raisebox{.7ex}{\xymatrix{\dar[r]^-{^{^{\mbox{\scriptsize$u$}}}}_-{v} & \, F}} ,
$$
on a
$$
f^* {\rm eg} \left(G \!\raisebox{.7ex}{\xymatrix{\dar[r]^-{^{^{\mbox{\scriptsize$u$}}}}_-{v} &F}}\right) = {\rm eg} \left( f^*G \!\raisebox{.7ex}{\xymatrix{\dar[r]^-{^{^{\mbox{\scriptsize$f^*(u)$}}}}_-{f^*(v)} &f^* F}} \right).
$$

Ainsi, la condition (A) signifie que le morphisme canonique
$$
{\rm eg} \left( G \!\raisebox{.7ex}{\xymatrix{\dar[r]^-{^{^{\mbox{\scriptsize$u$}}}}_-{v} &F}} \right) \longrightarrow G
$$
est un isomorphisme si et seulement si son image par $f^*$ est un isomorphisme.

\smallskip

Il en r\'esulte que (B) implique (A).

\medskip

\item R\'eciproquement, supposons que (A) est v\'erifi\'ee.

\smallskip

Montrons d'abord que tout monomorphisme de ${\mathcal E}$
$$
u : G \longrightarrow F
$$
que $f^*$ transforme en isomorphisme, est un isomorphisme de ${\mathcal E}$.

\smallskip

En effet, comme ${\mathcal E}$ est un topos, le sous-objet $G$ de $F$ s'\'ecrit comme l'\'egalisateur
$$
G = {\rm eg} \left( F \rightrightarrows F \amalg_G F \right).
$$

Comme
$$
f^* (F \amalg_G F ) = f^* F \amalg_{f^* G} f^*F \, ,
$$
on voit que $f^*$ transforme le morphisme
$$
G={\rm eg} \left( F \rightrightarrows F \amalg_G F \right) \longrightarrow F
$$
en isomorphisme, et donc c'est un isomorphisme d'apr\`es (A).

\smallskip

Si maintenant
$$
u : G \longrightarrow F
$$
est un morphisme quelconque de ${\mathcal E}$ que $f^*$ transforme en isomorphisme, on voit d'abord que le monomorphisme diagonal
$$
G \longrightarrow G \times_F G
$$
est transform\'e par $f^*$ en isomorphisme, donc est un isomorphisme, ce qui signifie que
$$
G \longrightarrow F
$$
est un monomorphisme.

\smallskip

Comme $f^*$ le transforme en isomorphisme, c'est un isomorphisme.

\smallskip

Ainsi, (A) implique (B). 

\end{listeisansmarge}
\end{remarksqed}
 
\bigskip
 
Cette d\'efinition \'etant pos\'ee, nous pouvons \'enoncer:
 
\begin{prop}\label{propIV75}

Soit
$$
f = (f^* , f_*) : {\mathcal E}' \longrightarrow {\mathcal E}
$$
un morphisme de topos.

\smallskip

Alors:

\begin{listeimarge}

\item Le morphisme s'inscrit dans un triangle de morphismes de topos
$$
\xymatrix{
{\mathcal E}' \ar[rd]_-{p} \ar[rr]^f &&{\mathcal E} \\
&{\rm Im} (f) \ar@{^{(}->}[ru]_-i
}
$$
commutatif \`a isomorphisme pr\`es, o\`u
$$
i : {\rm Im} (f) \xhookrightarrow{ \ { \ } \ } {\mathcal E}
$$
est un plongement de topos, et
$$
p : {\mathcal E}' \longrightarrow {\rm Im} (f)
$$
est un morphisme surjectif de topos.

\medskip

\item Cette factorisation est unique \`a unique isomorphisme pr\`es, si bien que le sous-topos repr\'esent\'e par
$$
i : {\rm Im} (f) \xhookrightarrow{ \ { \ } \ } {\mathcal E}
$$
peut \^etre appel\'e le sous-topos ``image'' de $f$.

\medskip

\item Si
$$
{\mathcal E} = \widehat{\mathcal C}_K
$$
est le topos des faisceaux sur un site $({\mathcal C},K)$, la topologie $J \supseteq K$ qui d\'efinit le sous-topos
$$
{\rm Im} (f) \xhookrightarrow{ \ { \ } \ } \widehat{\mathcal C}_K
$$
est constitu\'ee des cribles $S$ des objets $X$ de ${\mathcal C}$ qui satisfont les conditions \'equivalentes {\rm (1), (2)} et {\rm (3)} du lemme \ref{lemIV73} relativement au foncteur compos\'e
$$
\widehat{\mathcal C} \xrightarrow{ \ j^* \ } \widehat{\mathcal C}_J \xrightarrow{ \ f^* \ } {\mathcal E}' \, .
$$
\end{listeimarge}
\end{prop}

\begin{demo}

On peut supposer que
$$
{\mathcal E} = \widehat{\mathcal C}_K
$$
est le topos des faisceaux sur un site $({\mathcal C},K)$.

\smallskip

S'il existe une factorisation
$$
{\mathcal E} \xrightarrow{ \ p \ } {\rm Im} ({\mathcal E}) \xhookrightarrow{ \ i \ } \widehat{\mathcal C}_K \xhookrightarrow{ \ j \ } \widehat{\mathcal C} \, ,
$$
telle que $i$ est un plongement et $p$ un morphisme surjectif, on observe que, le foncteur $p^*$ \'etant conservatif, il est \'equivalent de demander, pour un crible d'un objet $X$ de ${\mathcal C}$
$$
S \xhookrightarrow{ \ { \ } \ } y(X) \, ,
$$
que
$$
i^* \circ j^* S = i^* \circ j^* \circ y(X)
$$
ou que
$$
f^* \!\circ j^* S = f^* \!\circ j^* \circ y(X) \, .
$$
 
Donc la topologie $J$ de ${\mathcal C}$ correspondant au sous-objet ${\rm Im} ({\mathcal E})$ de $\widehat{\mathcal C}_K$ et $\widehat{\mathcal C}$ est n\'ecessairement celle associ\'ee au foncteur compos\'e
$$
\widehat{\mathcal C} \xrightarrow{ \ j^* \ } \widehat{\mathcal C}_J \xrightarrow{ \ f^* \ } {\mathcal E}
$$
par le lemme \ref{lemIV73}.

\smallskip

D'o\`u l'unicit\'e de cette factorisation si elle existe.

\smallskip

R\'eciproquement, consid\'erons la topologie $J$ de ${\mathcal C}$ d\'efinie par le foncteur $f^* \!\circ j^* : \widehat{\mathcal C} \to {\mathcal E}$.

\smallskip

Elle contient la topologie $K$, donc d\'efinit un sous-topos
$$
\widehat{\mathcal C}_J \xhookrightarrow{ \ i \ } \widehat{\mathcal C}_K \xhookrightarrow{ \ j \ } \widehat{\mathcal C} \, .
$$

Pour tout crible $J$-couvrant $S$ d'un objet $X$ de ${\mathcal C}$, et pour tout objet $F$ de ${\mathcal E}'$, l'application
$$
j_* \circ f_* F (X) \longrightarrow \varprojlim_{(U \to X) \in S} j_* \circ f_* F(U)
$$
est bijective.

\smallskip

Cela signifie que le foncteur
$$
f_* : {\mathcal E}' \longrightarrow \widehat{\mathcal C}_K
$$
se factorise \`a travers la sous-cat\'egorie pleine $\widehat{\mathcal C}_J$ de $\widehat{\mathcal C}_K$ des $J$-faisceaux, en un foncteur
$$
p_* : {\mathcal E}' \longrightarrow \widehat{\mathcal C}_J \, .
$$

Ce foncteur respecte les limites, donc admet un adjoint \`a gauche
$$
p^* : \widehat{\mathcal C}_J \longrightarrow {\mathcal E}'
$$
et le triangle
$$
\xymatrix{
\widehat{\mathcal C}_K \ar[rd]_-{i^*} \ar[rr]^-{f^*} &&{\mathcal E}' \\
&\widehat{\mathcal C}_J \ar[ru]_-{p^*}
}
$$
est n\'ecessairement commutatif \`a isomorphisme pr\`es.

\smallskip

Le foncteur $p^*$ respecte les limites finies puisqu'il en est ainsi de $f^*$ et $i^*$, et que $i_*$ est pleinement fid\`ele.

\smallskip

Ainsi, $p=(p^*,p_*)$ est un morphisme de topos
$$
{\mathcal E}' \longrightarrow \widehat{\mathcal C}_J \, .
$$

Enfin, montrons que $p$ est surjectif, autrement dit que pour tout objet $F$ de $\widehat{\mathcal C}_J$, le morphisme canonique
$$
F \longrightarrow p_* \circ p^* F
$$
est un monomorphisme.

\smallskip

Consid\'erons donc un objet $X$ de ${\mathcal C}$ et deux sections
$$
x_1 , x_2 \in F(X)
$$
qui ont m\^eme image dans $(p_* \circ p^* F)(X)$.

\smallskip

Ces sections peuvent \^etre vues comme deux morphismes de $\widehat{\mathcal C}$
$$
y(X) \rightrightarrows F
$$
qui ont m\^eme compos\'e avec
$$
F \longrightarrow p_* \circ p^* F \, .
$$

Leur \'egalisateur est un crible
$$
S \xhookrightarrow{ \ { \ } \ } y(X)
$$
dont le transform\'e par $p^* \circ i^* \circ j^*$ ou $f^* \!\circ j^*$ est un isomorphisme puisque le transform\'e par $p^*$ de
$$
F \longrightarrow p_* \circ p^* F
$$
est un isomorphisme.

\smallskip

Cela signifie que le crible $S$ de $X$ est $J$-couvrant.

\smallskip

Comme $F$ est un $J$-faisceau, cela impose
$$
x_1 = x_2
$$
et $p : {\mathcal E}' \to \widehat{\mathcal C}_J$ est bien un morphisme surjectif. 

\end{demo}

\subsection{Sous-topos et op\'erateurs de cl\^oture}\label{subsec473}

\medskip

Si
$$
i = (i^* , i_*) : {\mathcal E}' \xhookrightarrow{ \ { \ } \ } {\mathcal E}
$$
est un plongement de topos, on peut associer \`a tout sous-objet $S \hookrightarrow F$ d'un objet $F$ de ${\mathcal E}$ le sous-objet de $F$
$$
\overline S = F \times_{i_* \circ i^* (F)} i_* \circ i^* (F) \, .
$$

Nous allons voir que la famille d'applications ainsi d\'efinie
$$
(S \hookrightarrow F) \longmapsto (\overline S \hookrightarrow F)
$$
est un ``op\'erateur de cl\^oture'' au sens de la d\'efinition suivante, et que se donner un op\'erateur de cl\^oture sur un topos ${\mathcal E}$ \'equivaut \`a se donner un sous-topos de ${\mathcal E}$.

\begin{defn}\label{defIV76}

Soit ${\mathcal E}$ un topos.

\smallskip

Soit $\Omega$ le classificateur des sous-objets de ${\mathcal E}$ muni du morphisme canonique $1_{\mathcal E} \hookrightarrow \Omega$.

\smallskip

Alors:

\begin{listeimarge}

\item On appelle op\'erateurs de cl\^oture de ${\mathcal E}$ les morphismes de ${\mathcal E}$
$$
p : \Omega \longrightarrow \Omega
$$
tels que

\medskip

$
\left\{\begin{matrix}
\bullet &\mbox{$p$ est idempotent au sens que $p \circ p = p$,} \hfill \\
{ \ } \\
\bullet &\mbox{$p$ respecte le morphisme $1_{\mathcal E} \hookrightarrow \Omega$,} \hfill \\
{ \ } \\
\bullet &\mbox{$p$ respecte la loi d'intersection $\Omega \times \Omega \xrightarrow{ \ \wedge \ } \Omega$.} \hfill
\end{matrix} \right.
$

\medskip

\item On dit que deux op\'erateurs de cl\^oture $p$ et $q$ satisfont la relation
$$
p \leq q
$$
si
$$
p \circ q = p = q \circ p \, .
$$
\end{listeimarge}
\end{defn}

\begin{remarksqed}
\begin{listeisansmarge}
\item Par d\'efinition du classificateur des sous-objets $\Omega$, l'ensemble ${\rm Hom} (E,\Omega)$ s'identifie pour tout objet $E$ de ${\mathcal E}$ \`a celui des sous-objets $E' \hookrightarrow E$.

\smallskip

Plus pr\'ecis\'ement, pour tout morphisme $E \xrightarrow{ \ \omega \ } \Omega$, le sous-objet de $E$ qui lui correspond est le produit fibr\'e
$$
E \times_{\Omega} 1_{\mathcal E}
$$
des morphismes $E \xrightarrow{ \ \omega \ } \Omega$ et $1_{\mathcal E} \hookrightarrow \Omega$. En particulier, le sous-objet total $E$ de $E$ correspond au morphisme compos\'e $E \to 1_{\mathcal E} \hookrightarrow \Omega$.

\smallskip

Se donner un morphisme
$$
\Omega \longrightarrow \Omega
$$
\'equivaut \`a associer \`a tout sous-objet $S \hookrightarrow E$ de tout objet $E$ un sous-objet $\overline S \hookrightarrow E$ de sorte que, pour tout morphisme de ${\mathcal E}$
$$
e : E' \longrightarrow E \, ,
$$
on ait pour tout objet $S$ de $E$
$$
\overline{e^{-1} (S)} = e^{-1} (\overline S) \, .
$$

La condition $p \circ p = p$ signifie que l'on a toujours
$$
\overline{\overline S} = \overline S \, .
$$

D'autre part, la compatibilit\'e de $p$ avec $1_{\mathcal E} \hookrightarrow \Omega$ signifie que pour tout sous-objet
$$
s : S \xhookrightarrow{ \ { \ } \ } E
$$
on a
$$
s^{-1} (\overline S) = S
$$
c'est-\`a-dire
$$
S \leq \overline S \, .
$$

\smallskip

Enfin, la compatibilit\'e de $p$ avec $\Omega \times \Omega \xrightarrow{ \ \wedge \ } \Omega$ signifie que pour tous sous-objets $S_1 , \cdots , S_n$ d'un objet $E$ de ${\mathcal E}$, on a
$$
\overline{S_1 \wedge \cdots \wedge S_n} = \overline S_1 \wedge \cdots \wedge \overline S_n \, .
$$

\item La relation $\leq$ sur les op\'erateurs de cl\^oture est une relation d'ordre. Elle est transitive car si $p \leq q$ et $q \leq r$, on a
$$
p \circ r = (p \circ q) \circ r = p \circ (q \circ r) = p \circ q = p
$$
et
$$
r \circ p = r \circ (q \circ p) = (r \circ q) \circ p = q \circ p = p \, ,
$$
ce qui signifie $p \leq r$.

\smallskip

Elle est r\'eflexive car tout op\'erateur de cl\^oture $p$ v\'erifie $p \circ p = p$.

\smallskip

Enfin, elle est antisym\'etrique car $p \leq q$ et $q \leq p$ entra{\^\i}nent
$$
p = p \circ q = q \, .
$$
\end{listeisansmarge}
\end{remarksqed}

Cette d\'efinition \'etant pos\'ee, nous pouvons \'enoncer:

\begin{thm}\label{thmIV77}

Soit ${\mathcal E}$ un topos muni de son classificateur des sous-objets $\Omega$.

\smallskip

Alors:

\begin{listeimarge}

\item Pour tout sous-topos
$$
i = (i^* , i_*) : {\mathcal E}' \xhookrightarrow{ \ { \ } \ } {\mathcal E} \, ,
$$
la famille d'applications
$$
(S \hookrightarrow E) \longmapsto (E \times_{i_* \circ i^* E} i_* \circ i^* S = \overline S \hookrightarrow E)
$$
d\'efinit un op\'erateur de cl\^oture de ${\mathcal E}$
$$
p_{{\mathcal E}'} : \Omega \longrightarrow \Omega \, .
$$

\item L'application
$$
({\mathcal E}' \hookrightarrow {\mathcal E}) \longmapsto p_{{\mathcal E}'}
$$
est une bijection de l'ensemble des sous-topos de ${\mathcal E}$ sur celui de ses op\'erateurs de cl\^oture.

\smallskip

Elle respecte la relation d'ordre, ainsi que sa bijection r\'eciproque.

\medskip

\item Si ${\mathcal E} = \widehat{\mathcal C}_K$ est le topos des faisceaux sur un site $({\mathcal C},K)$, la bijection r\'eciproque associe \`a tout op\'erateur de cl\^oture
$$
\begin{matrix}
\hfill \Omega &\longrightarrow &\Omega \, , \hfill \\
(S \hookrightarrow E) &\longmapsto &(\overline S \hookrightarrow E)
\end{matrix}
$$
le sous-topos de ${\mathcal E} = \widehat{\mathcal C}_K$ d\'efini par la topologie $J \supseteq K$ de ${\mathcal C}$ pour laquelle un crible d'un objet $X$ de ${\mathcal C}$
$$
S \xhookrightarrow{ \ { \ } \ } y(X)
$$
est couvrant si
$$
\overline{j^* S} = j^* \circ y(X)
$$
en notant $j^*$ le foncteur de faisceautisation
$$
j^* : \widehat{\mathcal C} \longrightarrow \widehat{\mathcal C}_K = {\mathcal E} \, .
$$
\end{listeimarge}
\end{thm}

\begin{remark}

Un op\'erateur de cl\^oture sur un topos ${\mathcal E}$ est aussi appel\'e une ``topologie de Lawvere'' sur ${\mathcal E}$.

\end{remark}

\bigskip

\begin{demo}
\begin{listeisansmarge}
\item[(i)] Les foncteurs $i^*$ et $i_*$ respectent les monomorphismes.

\smallskip

Si donc $S \hookrightarrow E$ est un monomorphisme de ${\mathcal E}$, il en est de m\^eme du morphisme
$$
\overline S \longrightarrow E
$$
d\'eduit de
$$
i_* \circ i^* S \xhookrightarrow{ \ { \ } \ } i_* \circ i^* E
$$
par le changement de base
$$
E \longrightarrow i_* \circ i^* E \, .
$$
Les applications
$$
(S \hookrightarrow E) \longmapsto (E \times_{i_* \circ i^* E} i_* \circ i^* S = \overline S \hookrightarrow E)
$$
respectent les intersections car les foncteurs $i^*$ et $i_*$ respectent les limites finies.

\smallskip

Tout sous-objet $(S \hookrightarrow E)$ est contenu dans le sous-objet $(\overline S \hookrightarrow E)$ car le carr\'e
$$
\xymatrix{
S \ \ar[d] \ar@{^{(}->}[r] &E \ar[d] \\
i_* \circ i^* S \ \ar@{^{(}->}[r] &i_* \circ i^* E
}
$$
est commutatif.

\smallskip

Enfin, on a pour tout morphisme de ${\mathcal E}$
$$
e : E' \longrightarrow E
$$
et tout sous-objet $S \hookrightarrow E$ les identifications
\begin{eqnarray}
\overline{e^{-1} S} &= &E' \times_{i_* \circ i^* E'} i_* \circ i^* (E' \times_E S) \nonumber \\
&= &E' \times_{i_* \circ i^* E'} (i_* \circ i^* E' \times_{i_* \circ i^* E} i_* \circ i^* S) \nonumber \\
&= &E' \times_{i_* \circ i^* E} i_* \circ i^* S \nonumber \\
&= &E' \times_E (E \times_{i_* \circ i^* E} i_* \circ i^* S) \nonumber \\
&= &e^{-1} \overline S \, . \nonumber
\end{eqnarray}

D'apr\`es le lemme de Yoneda, la famille d'applications
$$
(S \hookrightarrow E) \longmapsto (E \times_{i_* \circ i^* E} i_* \circ i^* S = \overline S \hookrightarrow E)
$$
d\'efinit un morphisme
$$
p_{{\mathcal E}'} : \Omega \longrightarrow \Omega
$$
qui est un op\'erateur de cl\^oture.

\medskip

\item[(iii)] Consid\'erons donc un topos de faisceaux sur un site $({\mathcal C},K)$
$$
{\mathcal E} = \widehat{\mathcal C}_K
$$
muni du foncteur de faisceautisation $j^* : \widehat{\mathcal C} \to \widehat{\mathcal C}_K$ adjoint \`a gauche du plongement $j_* : \widehat{\mathcal C}_K \hookrightarrow \widehat{\mathcal C}$ et de l'objet $\Omega_K$ classificateur des sous-objets.

\smallskip

On sait d'apr\`es le th\'eor\`eme \ref{thmIV72} que les sous-topos de ${\mathcal E} = \widehat{\mathcal C}_K$ correspondent aux topologies de ${\mathcal C}$ qui contiennent $K$.

\smallskip

D'apr\`es (i), toute topologie $J \supseteq K$ d\'efinit donc un op\'erateur de cl\^oture $p_J : \Omega_K \to \Omega_K$ de ${\mathcal E} = \widehat{\mathcal C}_K$.

\smallskip

R\'eciproquement, si $p = \Omega_K \to \Omega_K$ est un op\'erateur de cl\^oture consistant en une famille d'applications
$$
(S \hookrightarrow F) \longmapsto (\overline S \hookrightarrow F) \, ,
$$
notons $J_p (X)$ l'ensemble des cribles des objets $X$ de ${\mathcal C}$
$$
S \xhookrightarrow{ \ { \ } \ } y(X)
$$
tels que
$$
\overline{j^* S} = j^* \circ y(X)
$$
o\`u $j^* S$ est vu comme un sous-objet de $j^* \circ y(X)$ dans $\widehat{\mathcal C}_K$.

\smallskip

La famille des $J_p (X)$, $X \in {\rm Ob}({\mathcal C})$, est une topologie de ${\mathcal C}$.

\smallskip

En effet, elle satisfait \'evidemment l'axiome de maximalit\'e.

\smallskip

Elle satisfait aussi l'axiome de stabilit\'e car, pour tout morphisme de ${\mathcal C}$
$$
X' \xrightarrow{ \ x \ } X
$$
et tout crible $S \hookrightarrow y(X)$ de $X$, on a
\begin{eqnarray}
\overline{j^* (S \times_{y(X)} y(X'))} &= &\overline{j^* S \times_{j^* \circ y(X)} j^* \circ y(X')} \nonumber \\
&= &\overline{j^* S} \times_{j^* \circ y(X)} j^* \circ y(X') \, . \nonumber
\end{eqnarray}

Enfin, elle satisfait l'axiome de transitivit\'e car si $S$ et $S'$ sont deux cribles d'un objet $X$ de ${\mathcal C}$ v\'erifiant
$$
S \in J_p (X)
$$
et
$$
u^{-1} S' \in J_p (U) \, , \qquad \forall \ (U \xrightarrow{ \ u \ } X) \in S \, ,
$$
on a pour tout \'el\'ement $(U \xrightarrow{ \ u \ } X)$ de $S$
$$
j^* \circ y(u)^{-1} (\overline{j^* S'}) = \overline{j^* (u^{-1} S')} = j^* \circ y(U)
$$
si bien que chaque morphisme
$$
j^* \circ y(U) \longrightarrow j^* \circ y(X)
$$
se factorise \`a travers le sous-objet
$$
\overline{j^* S'} \xhookrightarrow{ \ { \ } \ } j^* \circ y(X)
$$
donc aussi le morphisme
$$
\varinjlim_{(U \xrightarrow{u} X) \in S} j^* \circ y(U) = j^* S \longrightarrow j^* \circ y(X) \, ,
$$
ce qui impose
$$
\overline{j^* S'} = j^* \circ y(X)
$$
puisque
$$
\overline{j^* S} = j^* \circ y(X) \, .
$$

La topologie $J_p$ contient $K$ car, si $S \hookrightarrow y(X)$ est $K$-couvrant, on a $j^* S = j^* \circ y(X)$.

\smallskip

Il reste \`a prouver que les deux applications
$$
J \longmapsto p_J
$$
et
$$
p \longmapsto J_p
$$
sont r\'eciproques l'une de l'autre.

\smallskip

Si $J \supseteq K$ est une topologie de ${\mathcal C}$ qui correspond \`a un sous-topos
$$
(i^* , i_*) : {\mathcal E}' = \widehat{\mathcal C}_J \xhookrightarrow{ \ { \ } \ } \widehat{\mathcal C}_K = {\mathcal E} \, ,
$$
un crible d'un objet $X$ de ${\mathcal C}$
$$
S \xhookrightarrow{ \ { \ } \ } y(X)
$$
est $J$-couvrant si et seulement si
$$
i^* \circ j^* S = i^* \circ j^* \circ y(X) \, .
$$
C'est \'equivalent \`a demander que
$$
j^* \circ y(X) \times_{i_* \circ i^* \circ j^* \circ y(X)} i_* \circ i^* \circ j^* S = j^* \circ y(X)
$$
puisque le foncteur $i^*$ respecte les produits fibr\'es et que $i^* \circ i_*$ s'identifie au foncteur ${\rm id}$.

\smallskip

Ainsi, $J$ est la topologie d\'efinie par l'op\'erateur de cl\^oture $p_J$ associ\'e au sous-topos $\widehat{\mathcal C}_J \hookrightarrow \widehat{\mathcal C}_K$.

\smallskip

R\'eciproquement, consid\'erons un op\'erateur de cl\^oture $p$ sur $\widehat{\mathcal C}_K$, la topologie $J_p \supseteq K$ qu'il d\'efinit et le sous-topos associ\'e
$$
i = (i^* , i_*) : \widehat{\mathcal C}_{J_p} \xhookrightarrow{ \ { \ } \ } \widehat{\mathcal C}_K \, .
$$
Il faut prouver que pour tout sous-objet d'un objet $F$ de $\widehat{\mathcal C}_K$
$$
S \xhookrightarrow{ \ { \ } \ }  F \, ,
$$
on a
$$
\overline S = F \times_{i_* \circ i^* F} i_* \circ i^* S \, .
$$

Pour cela, consid\'erons une section de $F$ sur un objet $X$ de ${\mathcal C}$
$$
x \in F(X) = {\rm Hom} (y(X) , j_* F) = {\rm Hom} (j^* \circ y(X),F) \, .
$$
C'est une section de $\overline S$ si et seulement si on a
$$
\overline{S \times_F j^* \circ y(X)} = j^* \circ y(X)
$$
c'est-\`a-dire si le crible de $X$
$$
j_* S \times_{j_* F} y(X) \xhookrightarrow{ \ { \ } \ } y(X)
$$
est couvrant pour la topologie $J_p$.

\smallskip

D'autre part, c'est une section de
$$
F \times_{i_* \circ i^* F} i_* \circ i^* S
$$
si et seulement si il existe un crible $J_p$-couvrant $S_x$ de $X$ dont tout \'el\'ement $(U \xrightarrow{ u } X)$ v\'erifie
$$
F(u)(x) \in S(U) \, .
$$

Ainsi, $p$ est l'op\'erateur de cl\^oture associ\'e au sous-topos $\widehat{\mathcal C}_{J_p} \hookrightarrow \widehat{\mathcal C}_K$ d\'efini par la topologie $J_p \supseteq K$ de ${\mathcal C}$.

\smallskip

Cela ach\`eve de prouver (iii).

\medskip

\item[(ii)] On sait d\'ej\`a d'apr\`es (iii) que cette application est une bijection.

\smallskip

Il reste \`a prouver qu'elle respecte la relation d'ordre, ainsi que sa r\'eciproque.

\smallskip

Pour une suite de sous-topos
$$
{\mathcal E}'_2 \xhookrightarrow{ \ i = (i^*,i_*) \ }  {\mathcal E}'_1 \xhookrightarrow{ \ j=(j^* , j_*) \ } {\mathcal E} \, ,
$$
le foncteur $j^* \circ j_*$ est canoniquement isomorphe \`a ${\rm id}_{{\mathcal E}'}$, d'o\`u des isomorphismes canoniques
$$
(j_* \circ i_* \circ i^* \circ j^*) \circ (j_* \circ j^*) \xrightarrow{ \ \sim \ } j_* \circ i_* \circ i^* \circ j^* \, ,
$$
$$
(j_* \circ j^*) \circ (j_* \circ i_* \circ i^* \circ j^*) \xrightarrow{ \ \sim \ } j_* \circ i_* \circ i^* \circ j^* \, .
$$
Il en r\'esulte que si $p_1$ et $p_2$ sont les op\'erateurs de cl\^oture de ${\mathcal E}$ associ\'es \`a ses sous-topos ${\mathcal E}'_1$ et ${\mathcal E}'_2$, on a
$$
p_2 \circ p_1 = p_2 = p_1 \circ p_2 \, .
$$

R\'eciproquement, on peut supposer que ${\mathcal E} = \widehat{\mathcal C}_K$ est le topos des faisceaux sur un site $({\mathcal C},K)$ et consid\'erer les topologies $J_1 , J_2 \supseteq K$ sur ${\mathcal C}$ associ\'ees \`a deux op\'erateurs de cl\^oture $p_1 , p_2$ sur $\widehat{\mathcal C}_K$.

\smallskip

Pour tout crible $S \hookrightarrow y(X)$ sur un objet $X$ de ${\mathcal C}$ et le sous-objet image $j^* S$ par le foncteur de faisceautisation $j^* : \widehat{\mathcal C} \to \widehat{\mathcal C}_K$, on a
$$
p_2 (j^* S) = j^* \circ y(X)
$$
si
$$
p_1 (j^* S) = j^* \circ y(X) \quad \mbox{et} \quad p_2 \circ p_1 = p_2 \, .
$$
Donc la relation $p_2 \leq p_1$ entra{\^\i}ne $J_1 \subseteq J_2$. 
\end{listeisansmarge}
\end{demo}

\subsection{Op\'erations sur les sous-topos}\label{subsec474}

\medskip

Dans l'ensemble ordonn\'e des sous-topos d'un topos ${\mathcal E}$, toute famille d'\'el\'ements a un supremum:

\begin{prop}\label{propIV78}

Soit ${\mathcal E}$ un topos.

\begin{listeimarge}

\item Pour toute famille $({\mathcal E}_i)_{i \in I}$ de sous-topos de ${\mathcal E}$, il existe un unique sous-topos de ${\mathcal E}$
$$
\bigvee_{i \in I} {\mathcal E}_i
$$
tel que, pour tout sous-topos ${\mathcal E}'$ de ${\mathcal E}$ on ait l'\'equivalence
$$
{\mathcal E}' \geq \bigvee_{i \in I} {\mathcal E}_i \Longleftrightarrow {\mathcal E}' \geq {\mathcal E}_i \, , \quad \forall \, i \in I \, .
$$

\item Si ${\mathcal E} = \widehat{\mathcal C}_K$ est le topos des faisceaux sur un site $({\mathcal C},K)$ et les ${\mathcal E}_i$, $i \in I$, sont des sous-topos de ${\mathcal E}$ d\'efinis par des topologies $J_i \supseteq K$, $i \in I$, le supremum
$$
\bigvee_{i \in I} {\mathcal E}_i
$$
est d\'efini par la topologie
$$
J = \bigcap_{i \in I} J_i \, .
$$
\end{listeimarge}
\end{prop}

\begin{remark}

Si ${\mathcal E} = {\mathcal E}_X$ est le topos des faisceaux sur un espace topologique $X$ et les ${\mathcal E}_i = {\mathcal E}_{Y_i}$, $i \in I$, sont des sous-topos de ${\mathcal E}$ associ\'es \`a des sous-espaces
$$
Y_i \xhookrightarrow{ \ { \ } \ } X \, ,
$$
alors le supremum
$$
\bigvee_{i \in I} {\mathcal E}_{Y_i}
$$
est le sous-topos de ${\mathcal E} = {\mathcal E}_X$ associ\'e au sous-espace
$$
Y = \bigcup_{i \in I} Y_i \, .
$$

En effet, chaque ${\mathcal E}_{Y_i}$ est d\'efini par la topologie $J_i$ de $O(X)$ pour laquelle une famille d'immersions entre ouverts
$$
(U_j \xhookrightarrow{ \ { \ } \ } U)_j
$$
est couvrante si et seulement si
$$
\bigcup_j \ (U_j \cap Y_i) = U \cap Y_i \, .
$$

Comme $Y = \underset{i \in I}{\bigcup} \ Y_i$, l'intersection des topologies $J_i$ est la topologie $J$ de $O(X)$ pour laquelle une famille
$$
(U_j \xhookrightarrow{ \ { \ } \ } U)_j
$$
est couvrante si et seulement si
$$
\bigcup_j \ (U_j \cap Y) = U \cap Y \, .
$$

\end{remark}

\begin{demo}

Il suffit de consid\'erer la situation de (ii) o\`u ${\mathcal E} = \widehat{\mathcal C}_K$ est le topos des faisceaux sur un site $({\mathcal C},K)$ et les sous-topos ${\mathcal E}_i$ de ${\mathcal E}$ sont d\'efinis par des topologies $J_i \supseteq K$ de ${\mathcal C}$.

\smallskip

Alors un sous-topos ${\mathcal E}'$ de ${\mathcal E}$ d\'efini par une topologie $J' \supseteq K$ contient tous les ${\mathcal E}_i$ comme sous-topos si et seulement si
$$
J' \subseteq J_i \, , \qquad \forall \, i \in I \, ,
$$
c'est-\`a-dire
$$
J' \subseteq \bigcap_{i \in I} J_i \, .
$$
On conclut en remarquant que toute intersection de topologies sur une cat\'egorie ${\mathcal C}$ est encore une topologie. Si toutes contiennent une topologie $K$, il en est encore ainsi de leur intersection. 

\end{demo}

\bigskip

De m\^eme dans l'ensemble ordonn\'e des sous-topos d'un topos ${\mathcal E}$, toute famille d'\'el\'ements a aussi un infimum:

\begin{prop}\label{propIV79}

Soit ${\mathcal E}$ un topos.

\begin{listeimarge}

\item Pour toute famille $({\mathcal E}_i)_{i \in I}$ de sous-topos de ${\mathcal E}$, il existe un unique sous-topos de ${\mathcal E}$
$$
\bigwedge_{i \in I} {\mathcal E}_i
$$
tel que, pour tout sous-topos ${\mathcal E}'$ de ${\mathcal E}$, on ait l'\'equivalence
$$
{\mathcal E}' \leq \bigwedge_{i \in I} {\mathcal E}_i \Longleftrightarrow {\mathcal E}' \leq {\mathcal E}_i \, , \quad \forall \, i \in I \, .
$$

\item Si ${\mathcal E} = \widehat{\mathcal C}_K$ est le topos des faisceaux sur un site $({\mathcal C},K)$ et les ${\mathcal E}_i = \widehat{\mathcal C}_{J_i}$, $i \in I$, sont des sous-topos de ${\mathcal E}$ d\'efinis par des topologies $J_i \supseteq K$, $i \in I$, l'infimum
$$
\bigwedge_{i \in I} {\mathcal E}_i = \widehat{\mathcal C}_J
$$
est d\'efini par la topologie $J$ de ${\mathcal C}$ ``engendr\'ee'' par les $J_i$, c'est-\`a-dire par la plus petite toologie de ${\mathcal C}$ qui contient toutes les topologies $J_i$, $i \in I$.

\smallskip

De plus, un faisceau $F$ sur $({\mathcal C},K)$ est un faisceau pour la topologie $J$ si et seulement si c'est un faisceau pour chacune des topologies $J_i$, $i \in I$.
\end{listeimarge}
\end{prop}

\begin{remark}

Si ${\mathcal E} = {\mathcal E}_X$ est le topos des faisceaux sur un espace topologique $X$ et les ${\mathcal E}_i = {\mathcal E}_{Y_i}$, $1 \leq i \leq n$, sont une famille finie de sous-topos associ\'es \`a des sous-espaces localement ferm\'es $Y_i$, alors l'infimum
$$
\bigwedge_{1 \leq i \leq n} {\mathcal E}_{Y_i} = {\mathcal E}_{Y_1} \wedge \cdots \wedge {\mathcal E}_{Y_n}
$$
est le sous-topos de ${\mathcal E} = {\mathcal E}_X$ associ\'e au sous-espace localement ferm\'e
$$
Y = \bigcap_{1 \leq i \leq n} Y_i = Y_1 \cap \cdots \cap Y_n \, .
$$

En effet, si $J_i$ d\'esigne la topologie sur $O(X)$ qui d\'efinit chaque sous-topos ${\mathcal E}_{Y_i}$, $1 \leq i \leq n$, alors un faisceau $F$ sur $X$ est un faisceau pour la topologie $J_i$ si, pour toute immersion d'ouverts de $X$
$$
V \xhookrightarrow{ \ { \ } \ } U \, ,
$$
l'application de restriction associ\'ee
$$
F(U) \longrightarrow F(V)
$$
est bijective d\`es que $V \cap Y_i = U \cap Y_i$.

\smallskip

Comme $Y_i$ est par hypoth\`ese l'intersection d'un ouvert $U_i$ de $X$ et d'un ferm\'e $Z_i$, cela revient \`a demander que pour tout ouvert $U$ de $X$ on ait
$$
F (U - Z_i \cap U) = \emptyset
$$
et l'application de restriction
$$
F(U) \longrightarrow F(U \cap U_i)
$$
soit bijective.

\smallskip

Cette condition est v\'erifi\'ee pour tous les indices $i$, $1 \leq i \leq n$, si et seulement si, pour tout ouvert $U$ de $X$, on a
$$
F \left( U - \left( \bigcap_{1 \leq i \leq n} Z_i \right) \cap U \right) = \emptyset
$$
et l'application de restriction
$$
F(U) \longrightarrow F \left( U \cap \bigcap_{1 \leq i \leq n} U_i \right)
$$
est bijective.

\smallskip

Cela revient \`a demander que $F$ soit un faisceau pour la topologie $J$ associ\'ee au sous-espace localement ferm\'e
$$
Y = \left( \bigcap_{1 \leq i \leq n} Z_i \right) \cap \left( \bigcap_{1 \leq i \leq n} U_i \right) = Y_1 \cap \cdots \cap Y_n \, .
$$
\end{remark}

\begin{demosansqed}

Il suffit de consid\'erer la situation de (ii) o\`u ${\mathcal E} = \widehat{\mathcal C}_K$ est le topos des faisceaux sur un site $({\mathcal C},K)$ et les sous-topos ${\mathcal E}_i$ de ${\mathcal E}$ sont d\'efinis par des topologies $J_i \supseteq K$ de ${\mathcal C}$.

\smallskip

Alors un sous-topos ${\mathcal E}'$ de ${\mathcal E}$ d\'efini par une topologie $J' \supseteq K$ est contenu comme sous-topos dans tous les ${\mathcal E}_i$ si et seulement si
$$
J' \supseteq J_i \, , \qquad \forall \, i \in I \, .
$$
Ainsi, l'intersection $J$ de toutes les topologies $J'$ de ${\mathcal C}$ qui satisfont cette condition d\'efinit un sous-topos de ${\mathcal E} = \widehat{\mathcal C}_K$ qui est un infimum des ${\mathcal E}_i = \widehat{\mathcal C}_{J_i}$.

\smallskip

Un objet $F$ de $\widehat{\mathcal C}_K$ qui est un faisceau pour la topologie $J$ est a fortiori un faisceau pour chaque topologie $J_i$, $i \in I$.

\smallskip

L'implication r\'eciproque r\'esulte du lemme suivant:

\end{demosansqed}

\begin{lem}\label{lemIV710}

Soit 
$$
p : {\mathcal C}^{\rm op} \longrightarrow {\rm Ens}
$$
un pr\'efaisceau sur une cat\'egorie essentiellement petite ${\mathcal C}$.

\smallskip

Pour tout objet $X$ de ${\mathcal C}$, notons $J(X)$ l'ensemble des cribles $S$ de $X$ tels que, pour tout morphisme $X' \xrightarrow{ \ x \ } X$ de ${\mathcal C}$, l'application canonique
$$
P(X') \longrightarrow \varprojlim_{(U' \to X') \in x^{-1} S} P(U')
$$
est bijective.

\smallskip

Alors $J_P$ est une topologie de ${\mathcal C}$.
\end{lem}

\begin{remark}

Par cons\'equent, pour toute famille $(P_i)_{i \in I}$ de pr\'efaisceaux sur une cat\'egorie essentiellement petite ${\mathcal C}$, il existe une plus grande topologie $J$ pour laquelle tous les $P_i$ sont des faisceaux.

\smallskip

C'est
$$
J = \bigcap_{i \in I} J_{P_i} \, .
$$

\end{remark}

\begin{demolem}

Il est \'evident sur la d\'efinition que $J_P$ satisfait les axiomes de maximalit\'e et de stabilit\'e.

\smallskip

Pour la transitivit\'e, consid\'erons deux cribles $S$ et $S'$ d'un objet $X$ de ${\mathcal C}$ tels que $S = J_p (X)$ et $u^{-1} S' \in J_p (U)$ pour tout \'el\'ement $(U \xrightarrow{ \ u \ } X)$ de $S$.

\smallskip

Pour tout morphisme $X' \xrightarrow{ \ x \ } X$ de ${\mathcal C}$, on a
$$
P(X') = \varprojlim_{(U' \to X') \in x^{-1} S} P(U')
$$
avec, pour tout \'el\'ement $(U' \xrightarrow{ \ u' \ } X')$ de $x^{-1} S$,
$$
P(U') = \varprojlim_{(U'' \to U') \in (x \circ u')^{-1} S'} P(U'') \, .
$$

On en d\'eduit
$$
P(X') = \varprojlim_{(U'' \to X') \in x^{-1} S'} P(U'') \, ,
$$
ce qui signifie que $S' \in J_P(X)$. 

\end{demolem}

\bigskip

\noindent {\bf Fin de la d\'emonstration de la proposition \ref{propIV79}:}

\smallskip

Pour tout objet $F$ de $\widehat{\mathcal C}_K$ qui est un faisceau pour chaque topologie $J_i$, $i \in I$, consid\'erons la topologie $J_F$ qui lui est associ\'ee par le lemme ci-dessus.

\smallskip

Cette topologie $J_F$ contient toutes les topologies $J_i$, $i \in I$, donc aussi la topologie $J$ qu'elles engendrent.

\smallskip

Cela signifie comme voulu que $F$ est un faisceau pour la topologie $J$. 

\hfill $\Box$

\bigskip

L'ensemble ordonn\'e des sous-topos d'un topos est encore muni de l'op\'eration suivante induite par celle de r\'eunion:

\begin{thm}\label{thmIV711}

Soit ${\mathcal E}$ un topos.

\begin{listeimarge}

\item Pour tout sous-topos ${\mathcal E}_1$ de ${\mathcal E}$, le foncteur sur l'ensemble ordonn\'e des sous-topos ${\mathcal E}'$ de ${\mathcal E}$
$$
{\mathcal E}' \longmapsto {\mathcal E}' \vee {\mathcal E}_1
$$
admet un adjoint \`a gauche not\'e
$$
{\mathcal E}_2 \longmapsto {\mathcal E}_2 \backslash {\mathcal E}_1 \, .
$$

Il est caract\'eris\'e par la propri\'et\'e que, pour tous sous-topos ${\mathcal E}_2$ et ${\mathcal E}'$ de ${\mathcal E}$, on a l'\'equivalence
$$
{\mathcal E}_2 \leq {\mathcal E}' \vee {\mathcal E}_1 \Longleftrightarrow {\mathcal E}_2 \backslash {\mathcal E}_1 \leq {\mathcal E}' \, .
$$

\item Si ${\mathcal E} = \widehat{\mathcal C}_K$ est le topos des faisceaux sur un site $({\mathcal C},K)$ et les sous-topos ${\mathcal E}_1 , {\mathcal E}_2$ de ${\mathcal E}$ sont d\'efinis par des topologies $J_1 , J_2 \supseteq K$ de ${\mathcal C}$, alors le sous-topos ${\mathcal E}_2 \backslash {\mathcal E}_1$ est d\'efini par la topologie
$$
(J_1 \Rightarrow J_2) = J_0
$$
pour laquelle un crible $S$ d'un objet $X$ de ${\mathcal C}$ est couvrant si et seulement si il satisfait la condition suivante:

\medskip

\noindent {\rm (C)} Pour tout morphisme $X' \xrightarrow{ \ x \ } X$ de ${\mathcal C}$ et tout crible $S'$ de $X'$ tel que

\medskip
$
\left\{\begin{matrix}
\bullet &\mbox{$S'$ est $J_1$-couvrant,} \hfill \\
\bullet &\mbox{$S'$ est $J_2$-ferm\'e,} \hfill \\
\bullet &\mbox{$S'$ contient $x^{-1} (S)$,} \hfill
\end{matrix} \right.
$

\medskip

\noindent alors $S'$ est le crible maximal de $X'$.
\end{listeimarge}
\end{thm}

\begin{remarks}
\begin{listeisansmarge}
\item En particulier, \`a tout sous-topos ${\mathcal E}_1$ de ${\mathcal E}$ est associ\'e le sous-topos
$$
{\mathcal E} \backslash {\mathcal E}_1
$$
que l'on peut appeler le compl\'ementaire de ${\mathcal E}$.

\medskip

\item Pour tous sous-topos ${\mathcal E}_1 , {\mathcal E}_2$ de ${\mathcal E}$, la relation
$$
{\mathcal E}_2 \leq ({\mathcal E}_2 \backslash {\mathcal E}_1) \vee {\mathcal E}_1
$$
impose
$$
{\mathcal E}_2 \backslash ({\mathcal E}_2 \backslash {\mathcal E}_1) \leq {\mathcal E}_1 \, .
$$

En particulier, on a pour tout sous-topos ${\mathcal E}_1$ de ${\mathcal E}$
$$
{\mathcal E} \backslash ({\mathcal E} \backslash {\mathcal E}_1) \leq {\mathcal E}_1
$$
mais il n'y a pas n\'ecessairement \'egalit\'e.

\medskip

\item Si ${\mathcal E} = {\mathcal E}_X$ est le topos des faisceaux sur un espace topologique $X$ et ${\mathcal E}_1 = {\mathcal E}_{Y_1}$, ${\mathcal E}_2 = {\mathcal E}_{Y_2}$ sont deux sous-topos associ\'es \`a deux sous-espaces $Y_1$ et $Y_2$ de $X$, la relation entre sous-espaces
$$
Y_2 \subseteq Y_2 \backslash Y_1 \cup Y_1
$$
entra{\^\i}ne
$$
{\mathcal E}_{Y_2} \leq {\mathcal E}_{Y_2 \backslash Y_1} \vee {\mathcal E}_{Y_1}
$$
et donc
$$
{\mathcal E}_{Y_2} \backslash {\mathcal E}_{Y_1} \leq {\mathcal E}_{Y_2 \backslash Y_1} \, .
$$
En g\'en\'eral, il n'y a pas \'egalit\'e.

\smallskip

On verra cependant en remarque au corollaire suivant qu'il y a \'egalit\'e si $Y_1$ et $Y_2 \backslash Y_1$ sont tous deux localement ferm\'es.

\end{listeisansmarge}
\end{remarks}

\bigskip

\begin{demo}

Il suffit de consid\'erer la situation de (ii) o\`u ${\mathcal E} = \widehat{\mathcal C}_K$ est le topos des faisceaux sur un site $({\mathcal C},K)$ et les sous-topos ${\mathcal E}_1 , {\mathcal E}_2$ de ${\mathcal E}$ sont d\'efinis par des topologies $J_1 , J_2 \supseteq K$.

\smallskip

Un sous-topos ${\mathcal E}'$ de ${\mathcal E} = \widehat{\mathcal C}_K$ correspondant \`a une topologie $J \supseteq K$ satisfait l'in\'egalit\'e
$$
{\mathcal E}_2 \leq {\mathcal E}' \vee {\mathcal E}_1
$$
si et seulement si la topologie $J$ v\'erifie la propri\'et\'e
$$
J \cap J_1 \subseteq J_2 \, .
$$

Pour tout objet $X$ de ${\mathcal C}$, notons $J_0 (X)$ l'ensemble des cribles $S$ de $X$ qui satisfont la condition (C) de l'\'enonc\'e.

\smallskip

Nous devons d\'emontrer que $J_0$ est une topologie et que l'in\'egalit\'e
$$
J \cap J_1 \subseteq J_2
$$
\'equivaut \`a
$$
J \subseteq J_0 \, .
$$

Il est \'evident sur la d\'efinition que $J_0$ satisfait les axiomes de maximalit\'e et de stabilit\'e.

\smallskip

Pour la transitivit\'e, consid\'erons deux cribles $S,S'$ d'un objet $X$ de ${\mathcal C}$ tels que $S \in S_0 (X)$ et $u^{-1} S' \in J_0(U)$ pour tout \'el\'ement $(U \xrightarrow{ \ u \ } X)$ de $S$.

\smallskip

Puis consid\'erons un morphisme $x : X' \to X$ de ${\mathcal C}$ et un crible $S''$ de $X'$ qui est $J_1$-couvrant, $J_2$-ferm\'e et contient $x^{-1} S'$.

\smallskip

Nous devons montrer que $S''$ est le crible maximal de $X'$.

\smallskip

Comme $J_0$ satisfait l'axiome de stabilit\'e, on peut supposer que $X'=X$ et $x={\rm id}$.

\smallskip

Pour tout \'el\'ement $(U \xrightarrow{ \ u \ } X)$ de $S$, le crible $u^{-1} S''$ de $U$ est $J_1$-couvrant et $J_2$-ferm\'e, et il contient le crible $u^{-1} S'$ qui est \'el\'ement de $J_0 (U)$. Donc $u^{-1} S''$ est le crible maximal de $U$.

\smallskip

Cela signifie que $S''$ contient le crible $S$, qui est \'el\'ement de $J_0(X)$. Comme $S''$ est $J_1$-couvrant et $J_2$-ferm\'e, cela impose comme voulu que $S''$ est le crible maximal de $X$.

\smallskip

Ainsi, $J_0$ est bien une topologie de ${\mathcal C}$.

\smallskip

Montrons maintenant que
$$
J_0 \cap J_1 \subseteq J_2 \, .
$$
Soit $S$ un crible d'un objet $X$ de ${\mathcal C}$ qui est $J_0$-couvrant et $J_1$-couvrant.

\smallskip

Soit $\overline S$ la $J_2$-fermeture de $S$: c'est le crible des morphismes $U \xrightarrow{ \ u \ } X$ dont la source $U$ admet un crible $J_2$-couvrant constitu\'e de morphismes $U' \xrightarrow{ \ u' \ } U$ tels que $u \circ u' \in S$. Alors le crible $\overline S$ de $X$ est $J_2$-ferm\'e et contient $S$, donc est aussi $J_1$-couvrant.

\smallskip

Par d\'efinition de $J_0$, c'est le crible maximal de $X$, ce qui signifie que $S$ est $J_2$-couvrant.

\smallskip

Enfin, consid\'erons une topologie $J$ de ${\mathcal C}$ telle que
$$
J \cap J_1 \subseteq J_2 \, .
$$

Consid\'erons un crible $J$-couvrant $S$ d'un objet $X$ de ${\mathcal C}$, un morphisme $x : X' \to X$ et un crible $S'$ de $X'$ qui est $J_1$-couvrant, $J_2$-ferm\'e et contient $x^{-1} S$.

\smallskip

Alors $S'$ est $(J \cap J_1)$-couvrant et a fortiori $J_2$-couvrant. Comme il est $J_2$-ferm\'e, cela impose comme voulu que $S'$ est le cible maximal de $X'$.

\smallskip

Cela termine la d\'emonstration du th\'eor\`eme. 

\end{demo}

\bigskip

On d\'eduit aussit\^ot de ce th\'eor\`eme:

\begin{cor}\label{corIV712}

Soit ${\mathcal E}$ un topos.

\begin{listeimarge}

\item Pour tous sous-topos ${\mathcal E}'$ et $({\mathcal E}_i)_{i \in I}$ de ${\mathcal E}$, on a
$$
{\mathcal E}' \vee \left( \bigwedge_{i \in I} {\mathcal E}_i \right) = \bigwedge_{i \in I} ({\mathcal E}' \vee {\mathcal E}_i)
$$
et
$$
\left( \bigvee_{i \in I} {\mathcal E}_i \right) \backslash \, {\mathcal E}' = \bigvee_{i \in I} ({\mathcal E}_i \backslash {\mathcal E}') \, .
$$

\item Pour tous sous-topos ${\mathcal E}'$ et ${\mathcal E}_1 , \ldots , {\mathcal E}_n$, on a
$$
{\mathcal E}' \wedge ({\mathcal E}_1 \vee \cdots \vee {\mathcal E}_n) = ({\mathcal E}' \wedge {\mathcal E}_1) \vee \cdots \vee ({\mathcal E}' \wedge {\mathcal E}_n) \, .
$$

\item Si ${\mathcal E}_1$ et ${\mathcal E}_2$ sont deux sous-topos de ${\mathcal E}$ tels que
$$
{\mathcal E}_1 \vee {\mathcal E}_2 = {\mathcal E}
$$
et
$$
{\mathcal E}_1 \wedge {\mathcal E}_2 = {\mathcal E}_{\emptyset}
$$
(o\`u ${\mathcal E}_{\emptyset}$ d\'esigne le sous-topos minimal de ${\mathcal E}$, dont tout objet est terminal), alors
$$
{\mathcal E}_2 = {\mathcal E} \backslash {\mathcal E}_1 \qquad \mbox{et} \qquad {\mathcal E}_1 = {\mathcal E} \backslash {\mathcal E}_2 \, .
$$
\end{listeimarge}
\end{cor}

\begin{remarks}
\begin{listeisansmarge}
\item Si ${\mathcal E} = {\mathcal E}_X$ est le topos des faisceaux sur un espace topologique $X$ et que $Y_1 \subseteq Y_2$ sont deux sous-espaces de $X$ tels que $Y_1$ et $Y_2 \backslash Y_1$ soient localement ferm\'es, on a
$$
{\mathcal E}_{Y_2} \vee {\mathcal E}_{Y_2 \backslash Y_1} = {\mathcal E}_{Y_2} \qquad \mbox{et} \qquad {\mathcal E}_{Y_1} \wedge {\mathcal E}_{Y_2 \backslash Y_1} = {\mathcal E}_{\emptyset} \, ,
$$
d'o\`u d'apr\`es (iii)
$$
{\mathcal E}_{Y_2 \backslash Y_1} = {\mathcal E}_{Y_2} \backslash {\mathcal E}_{Y_1} \qquad \mbox{et} \qquad {\mathcal E}_{Y_1} = {\mathcal E}_{Y_2} \backslash {\mathcal E}_{Y_2 \backslash Y_1} \, .
$$

\item En particulier, si $U$ est un ouvert d'un espace topologique $X$ et $Z$ son ferm\'e compl\'ementaire, on a
$$
{\mathcal E}_X \backslash {\mathcal E}_U = {\mathcal E}_Z \qquad \mbox{et} \qquad {\mathcal E}_X \backslash {\mathcal E}_Z = {\mathcal E}_U \, .
$$

\end{listeisansmarge}
\end{remarks}

\begin{demo}
\begin{listeisansmarge}
\item Le foncteur
$$
{\mathcal E}_1 \longmapsto {\mathcal E}' \vee {\mathcal E}_1
$$
respecte les limites arbitraires puisqu'il a un adjoint \`a gauche.

\smallskip

De m\^eme, le foncteur
$$
{\mathcal E}_2 \longmapsto {\mathcal E}_2 \backslash {\mathcal E}'
$$
respecte les colimites arbitraires puisqu'il a un adjoint \`a droite.

\medskip

\item Il suffit de traiter le cas o\`u $n=2$.

\smallskip

D'apr\`es (i), on a
\begin{eqnarray}
({\mathcal E}' \wedge {\mathcal E}_1) \vee ({\mathcal E}' \wedge {\mathcal E}_2) &= &{\mathcal E}' \wedge ({\mathcal E}' \vee {\mathcal E}_1) \wedge ({\mathcal E}' \vee {\mathcal E}_2) \wedge ({\mathcal E}_1 \vee {\mathcal E}_2) \nonumber \\
&= &{\mathcal E}' \wedge ({\mathcal E}_1 \vee {\mathcal E}_2) \, . \nonumber
\end{eqnarray}
C'est ce que l'on voulait. 

\medskip

\item Par sym\'etrie, il suffit de montrer que ${\mathcal E}_1 = {\mathcal E} \backslash {\mathcal E}_2$.

\smallskip

L'\'egalit\'e ${\mathcal E}_1 \vee {\mathcal E}_2 = {\mathcal E}$ implique ${\mathcal E}_1 \geq {\mathcal E} \backslash {\mathcal E}_2$.

\smallskip

D'autre part, appliquant le foncteur ${\mathcal E}_1 \wedge \bullet$ \`a l'\'egalit\'e ${\mathcal E}_2 \vee {\mathcal E} \backslash {\mathcal E}_2 = {\mathcal E}$, on obtient
$$
{\mathcal E}_1 \wedge {\mathcal E} \backslash {\mathcal E}_2 = {\mathcal E}_1
$$
qui signifie
$$
{\mathcal E} \backslash {\mathcal E}_2 \leq {\mathcal E}_1 \, .
$$

\end{listeisansmarge}
\end{demo}

\subsection{Fonctorialit\'e des sous-topos}\label{subsec475}

\medskip

Tout morphisme de topos $f : {\mathcal E}' \to {\mathcal E}$ d\'efinit une application $f^{-1}$ qui transforme les sous-topos de ${\mathcal E}$ en sous-topos de ${\mathcal E}'$:

\begin{prop}\label{propIV713}

Soit $f = (f^* , f_*) : {\mathcal E}' \to {\mathcal E}$ un morphisme de topos, et soit
$$
j = (j^* , j_*) : {\mathcal E}_1 \xhookrightarrow{ \ { \ } \ } {\mathcal E}
$$
un sous-topos de ${\mathcal E}$.

\smallskip

Alors:

\begin{listeimarge}

\item Il existe un unique sous-topos de ${\mathcal E}'$
$$
j' = (j'^* , j'_*) : f^{-1} {\mathcal E}_1 = {\mathcal E}'_1 \xhookrightarrow{ \ { \ } \ } {\mathcal E}'
$$
tel qu'un morphisme de topos
$$
g = (g^* , g_*) : {\mathcal E}'' \longrightarrow {\mathcal E}'
$$
se factorise \`a travers ${\mathcal E}'_1 \hookrightarrow {\mathcal E}'$ si et seulement si le morphisme compos\'e
$$
f \circ g : {\mathcal E}'' \longrightarrow {\mathcal E}
$$
se factorise \`a travers ${\mathcal E}_1 \hookrightarrow {\mathcal E}$.

\medskip

\item Supposons que ${\mathcal E} = \widehat{\mathcal C}_J$ et ${\mathcal E}' = \widehat{\mathcal D}_K$ sont les topos des faisceaux sur deux sites $({\mathcal C},J)$ et $({\mathcal D},K)$, et que le morphisme de topos $f : {\mathcal E}' \to {\mathcal E}$ est induit par un morphisme de sites
$$
\rho : {\mathcal C} \longrightarrow {\mathcal D} \, .
$$

Soit $J_1 \supseteq J$ la topologie de ${\mathcal C}$ qui d\'efinit le sous-topos ${\mathcal E}_1$ de ${\mathcal E} = \widehat{\mathcal C}_J$.

\smallskip

Alors le sous-topos ${\mathcal E}'_1$ de ${\mathcal E}' = \widehat{\mathcal D}_K$ est d\'efini par la plus petite topologie $K_1$ de ${\mathcal D}$ qui contient $K$ et pour laquelle les images par $\rho$ des familles $J_1$-couvrantes de ${\mathcal C}$ sont des familles couvrantes de ${\mathcal D}$.
\end{listeimarge}
\end{prop}

\begin{remarks}
\begin{listeisansmarge}

\item Pour tout morphisme $f : {\mathcal E}' \to {\mathcal E}$, l'application
$$
f^{-1} : ({\mathcal E}_1 \hookrightarrow {\mathcal E}) \longmapsto (f^{-1} {\mathcal E}_1 \hookrightarrow {\mathcal E}')
$$
respecte la relation d'ordre entre sous-topos.

\smallskip

Elle respecte les intersections au sens que pour toute famille $({\mathcal E}_i \hookrightarrow {\mathcal E})_{i \in I}$ de sous-topos de ${\mathcal E}$, on a
$$
f^{-1} \left( \bigwedge_{i \in I} {\mathcal E}_i \right) = \bigwedge_{i \in I} f^{-1} {\mathcal E}_i \, .
$$

Enfin, on a pour toute telle famille de sous-topos
$$
\bigvee_{i \in I} f^{-1} {\mathcal E}_i \leq f^{-1} \left( \bigvee_{i \in I} {\mathcal E}_i \right).
$$

\item Pour tous morphismes de topos ${\mathcal E}'' \xrightarrow{ \ g \ } {\mathcal E}' \xrightarrow{ \ f \ } {\mathcal E}$, on a
$$
(f \circ g)^{-1} = g^{-1} \circ f^{-1} \, .
$$
\end{listeisansmarge}
\end{remarks}

\begin{demo}

D'apr\`es la remarque (vi) qui suit la d\'efinition \ref{defIV41}, on peut supposer que l'on est dans la situation de (ii):

\smallskip

On a ${\mathcal E} = \widehat{\mathcal C}_J$, ${\mathcal E}' = \widehat{\mathcal D}_K$ et $f = (f^* , f_*) : {\mathcal E}' \to {\mathcal E}$ est induit par un morphisme de sites
$$
\rho : {\mathcal C} \longrightarrow {\mathcal D} \, .
$$

Pour tout morphisme de topos
$$
g = (g^* , g_*) : {\mathcal E}'' \longrightarrow {\mathcal E}' = \widehat{\mathcal D}_K \, ,
$$
les foncteurs compos\'es
$$
{\mathcal D} \xrightarrow{ \ y \ } \widehat{\mathcal D} \xrightarrow{ \ j^* \ } \widehat{\mathcal D}_K \xrightarrow{ \ g^* \ } {\mathcal E}''
$$
et
$$
{\mathcal C} \xrightarrow{ \ \rho \ } {\mathcal D} \xrightarrow{ \ y \ } \widehat{\mathcal D} \xrightarrow{ \ j^* \ } \widehat{\mathcal D} \xrightarrow{ \ g^* \ } {\mathcal E}''
$$
sont plats et ils transforment toute famille $K$-couvrante [resp. $J$-couvrante] de morphismes de ${\mathcal D}$ [resp. de ${\mathcal C}$] en une famille globalement \'epimorphique de ${\mathcal E}''$.

\smallskip

Si
$$
{\mathcal E}_1 \xhookrightarrow{ \ { \ } \ } {\mathcal E}
$$
et
$$
{\mathcal E}'_1 \xhookrightarrow{ \ { \ } \ } {\mathcal E}'
$$
sont deux sous-topos de ${\mathcal E} = \widehat{\mathcal C}_J$ et ${\mathcal E}' = \widehat{\mathcal D}_K$ d\'efinis par des topologies $J_1$ et $K_1$, le morphisme
$$
g : {\mathcal E}'' \longrightarrow {\mathcal E}'
$$
se factorise \`a travers ${\mathcal E}'_1$ si et seulement si le foncteur
$$
{\mathcal D} \xrightarrow{ \ \ell \ } \widehat{\mathcal D}_K \xrightarrow{ \ g^* \ } {\mathcal E}''
$$
transforme toute famille $K_1$-couvrante de ${\mathcal D}$ en une famille globalement \'epimorphique de ${\mathcal E}''$, et le morphisme
$$
f \circ g : {\mathcal E}'' \longrightarrow {\mathcal E}' \longrightarrow {\mathcal E}
$$
se factorise \`a travers ${\mathcal E}_1$ si et seulement si le foncteur
$$
{\mathcal C} \xrightarrow{ \ \rho \ } {\mathcal D} \xrightarrow{ \ \ell \ } \widehat{\mathcal D}_K \xrightarrow{ \ g^* \ } {\mathcal E}''
$$
transforme toute famille $J_1$-couvrante de ${\mathcal C}$ en une famille globalement \'epimorphique de ${\mathcal E}''$.

\smallskip

Prenons pour $K_1$ la plus petite topologie de ${\mathcal D}$ qui contient $K$ et pour laquelle les images par $\rho$ des familles $J_1$-couvrantes de ${\mathcal C}$ sont couvrantes.

\smallskip

On voit que le morphisme compos\'e
$$
\widehat{\mathcal D}_{K_1} \xhookrightarrow{ \ { \ } \ } \widehat{\mathcal D}_K = {\mathcal E}' \xrightarrow{ \ f \ } {\mathcal E} = \widehat{\mathcal C}_J
$$
se factorise en
$$
\widehat{\mathcal D}_{K_1} \longrightarrow \widehat{\mathcal C}_{J_1} \xhookrightarrow{ \ { \ } \ } \widehat{\mathcal C}_J = {\mathcal E} \, .
$$

R\'eciproquement, consid\'erons un morphisme de topos
$$
g= (g^* , g_*) : {\mathcal E}'' \longrightarrow {\mathcal E}' = \widehat{\mathcal D}_K
$$
tel que le compos\'e
$$
f \circ g : {\mathcal E}'' \longrightarrow {\mathcal E}' \longrightarrow {\mathcal E}
$$
se factorise \`a travers ${\mathcal E}_1 \hookrightarrow {\mathcal E}$.

\smallskip

Comme le foncteur
$$
{\mathcal D} \xrightarrow{ \ \ell \ } \widehat{\mathcal D}_K \xrightarrow{ \ g^* \ } {\mathcal E}''
$$
est plat, il existe une topologie $K'$ de ${\mathcal D}$ pour laquelle un crible $S$ d'un objet $D$ est couvrant si et seulement si il contient une famille de morphismes
$$
(D_i \longrightarrow D)_{i \in I}
$$
dont les images par $g^* \circ \ell$ forment une famille \'epimorphique de morphismes de ${\mathcal E}''$.

\smallskip

Cette topologie $K'$ contient n\'ecessairement $K_1$, et donc le morphisme
$$
g : {\mathcal E}'' \longrightarrow {\mathcal E}' = \widehat{\mathcal D}_K
$$
se factorise \`a travers
$$
{\mathcal E}'_1 = \widehat{\mathcal D}_{K_1} \xhookrightarrow{ \ { \ } \ } \widehat{\mathcal D}_K = {\mathcal E}' \, .
$$

Cela prouve \`a la fois (ii) et (i). 

\end{demo}

\section{Localisation des topos}\label{sec48}

\subsection{Stabilit\'e des topos par localisation}\label{subsec481}

\medskip

On observe:

\begin{prop}\label{propIV81}
\begin{listeimarge}
\item Pour tout objet $E$ d'un topos ${\mathcal E}$, la cat\'egorie relative ${\mathcal E}/E$ est encore un topos.

\medskip

\item Plus pr\'ecis\'ement, pour tout site $({\mathcal C},J)$ et pour tout $J$-faisceau $F$ sur ${\mathcal C}$, la cat\'egorie relative
$$
\widehat{\mathcal C}_J / F
$$
s'identifie \`a la cat\'egorie des faisceaux sur la cat\'egorie des \'el\'ements
$$
{\textstyle \int} P = {\mathcal C}/P
$$
de n'importe quel pr\'efaisceau $P$ dont la faisceautisation est $F$, munie de la topologie pour laquelle une famille de morphismes
$$
(X_i , x_i) \longrightarrow (X,x) \, , \qquad i \in I \, ,
$$
est couvrante si et seulement si la famille induite de morphismes de ${\mathcal C}$
$$
X_i \longrightarrow X \, , \qquad i \in I \, ,
$$
est $J$-couvrante.

\medskip

\item En particulier, pour tout pr\'efaisceau $P$ sur une cat\'egorie essentiellement petite ${\mathcal C}$, la cat\'egorie relative
$$
\widehat{\mathcal C}/P
$$
s'identifie \`a la cat\'egorie $\widehat{{\mathcal C}/P}$ des pr\'efaisceaux sur la cat\'egorie des \'el\'ements de $P$
$$
{\textstyle \int} P = {\mathcal C}/P \, .
$$
\end{listeimarge}
\end{prop}

\begin{demo}
\begin{listeisansmarge}
\item[(iii)] Pour tout objet $X$ de ${\mathcal C}$, se donner un ensemble
$$
Q(X)
$$
muni d'une application
$$
p_X : Q(X) \longrightarrow P(X)
$$
\'equivaut \`a se donner une famille d'ensembles
$$
Q' (X,x)
$$
index\'es par les \'el\'ements $x \in P(X)$.

\smallskip

Il suffit en effet de prendre
$$
Q' (X,x) = p_X^{-1} (x) \, , \qquad \forall \, x \in P(X) \, ,
$$
et dans l'autre sens
$$
Q(X) = \coprod_{x \in P(X)} Q' (X,x)
$$
muni de l'application vers $P(X)$ dont la fibre au-dessus de chaque \'el\'ement $x$ est $Q'(X,x)$.

\smallskip

Puis, pour tout morphisme de ${\mathcal C}$
$$
X \xrightarrow{ \ u \ } Y \, ,
$$
se donner une application
$$
Q(Y) \longrightarrow Q(X)
$$
qui rende commutatif le carr\'e
$$
\xymatrix{
Q(Y) \ar[d] \ar[rr] &&Q(X) \ar[d] \\
P(Y) \ar[rr]^-{P(u)} &&P(X)
}
$$
\'equivaut \`a se donner une famille d'applications
$$
Q(Y,y) \longrightarrow Q(X,x)
$$
index\'ees par les paires d'\'el\'ements
$$
x \in P(X) \, , \qquad y \in P(Y)
$$
tels que $x = P(u)(y)$.

\smallskip

Cela prouve que $\widehat{\mathcal C} / P$ s'identifie \`a $\widehat{{\mathcal C}/P}$.

\medskip

\item[(ii)] On sait d\'ej\`a d'apr\`es (iii) que $\widehat{{\mathcal C}/F}$ s'identifie \`a $\widehat{\mathcal C}/F$.

\smallskip

Il est imm\'ediat sur la d\'efinition de la topologie de ${\mathcal C}/P = \int\!\!P$ qu'un pr\'efaisceau sur ${\mathcal C}/P$ est un faisceau si et seulement si le pr\'efaisceau sur ${\mathcal C}$ qui lui correspond est un faisceau pour la topologie $J$.

\smallskip

Ainsi, $\widehat{\mathcal C}_J/F$ s'identifie \`a la cat\'egorie des faisceaux sur ${\mathcal C}/F$.

\smallskip

Il reste \`a d\'emontrer que les cat\'egories des faisceaux sur ${\mathcal C}/F$ et ${\mathcal C}/P$ s'identifient d\`es lors que $F$ est la faisceautisation de $P$.

\smallskip

Le morphisme canonique $P \xrightarrow{ \ a \ } F$ induit le foncteur
$$
{\mathcal C}/P \longrightarrow {\mathcal C}/F
$$
qui associe \`a toute paire $(X,x)$ constitu\'ee d'un objet $X$ de ${\mathcal C}$ et d'un \'el\'ement $x \in P(X)$ la paire constitu\'ee de $X$ et de l'image de $x$ par l'application $a_X : P(X) \to F(X)$.

\smallskip

On sait d'apr\`es le lemme \ref{lemII56} que, pour tout objet $X$ de ${\mathcal C}$ et tout \'el\'ement $x \in F(X)$ [resp. tous \'el\'ements $x_1 , x_2 \in P(X)$ qui ont m\^eme image dans $F(X)$], il existe une famille $J$-couvrante de morphismes de ${\mathcal C}$
$$
X_i \longrightarrow X \, , \qquad i \in I \, ,
$$
telle que l'image de $x \in F(X)$ dans chaque $F(X_i)$ se rel\`eve en un \'el\'ement de $P(X_i)$ [resp. telle que les \'el\'ements $x_1 , x_2 \in P(X)$ aient m\^eme image dans chaque $P(X_i)$].

\smallskip

On en d\'eduit que la topologie $J_F$ de ${\mathcal C}/F$ est coinduite par la topologie $J_P$ de ${\mathcal C}/P$ (si bien que ${\mathcal C}/P \to {\mathcal C}/F$ est un comorphisme de sites) et que le foncteur ${\mathcal C}/P \to {\mathcal C}/F$ est $J_F$-dense et $J_P$-plein au sens du lemme \ref{lemIV64}.

\smallskip

Comme le foncteur ${\mathcal C}/F \to \widehat{({\mathcal C}/F)}_{J_F}$ est plat, on en d\'eduit \'egalement que le foncteur compos\'e
$$
{\mathcal C}/P \longrightarrow {\mathcal C}/F \longrightarrow \widehat{({\mathcal C}/F)}_{J_F}
$$
est plat puis que la topologie $J_P$ de ${\mathcal C}/P$ est induite par la topologie $J_F$ de ${\mathcal C}/F$ (si bien que ${\mathcal C}/P \to {\mathcal C}/F$ est un morphisme de sites).

\smallskip

D'apr\`es le lemme \ref{lemIV64}, on conclut comme voulu que le morphisme et comorphisme de sites
$$
{\mathcal C}/P  \longrightarrow {\mathcal C}/F
$$
d\'efinit deux \'equivalences r\'eciproques l'une de l'autre entre les topos de faisceaux sur les sites $({\mathcal C}/P , J_P)$ et $({\mathcal C}/F , J_F)$. 
\end{listeisansmarge}
\end{demo}

\bigskip

On pose:

\begin{defn}\label{defIV82}

On appelle localisations d'un topos ${\mathcal E}$ les topos qui sont \'equivalents \`a la cat\'egorie relative
$$
{\mathcal E}/E
$$
associ\'ee \`a un objet $E$ de ${\mathcal E}$. 
\end{defn}

\subsection{Les morphismes entre localisations associ\'es aux morphismes d'un topos}\label{subsec482}

\medskip

La proposition \ref{propIV81} se compl\`ete par la proposition suivante:

\begin{prop}\label{propIV83}

Soit ${\mathcal E}$ un topos.

\smallskip

Pour tout morphisme de ${\mathcal E}$
$$
E_2 \xrightarrow{ \ f \ } E_1 \, ,
$$
le foncteur
$$
\begin{matrix}
f^* : &\hfill {\mathcal E} / E_1 &\longrightarrow &{\mathcal E}/E_2 \, , \hfill \\
&(E \to E_1) &\longmapsto &(E \times_{E_1} E_2 \to E_2)
\end{matrix}
$$
respecte les limites et les colimites arbitraires.

\smallskip

Il admet un adjoint \`a gauche qui n'est autre que le foncteur de composition avec $f$
$$
\begin{matrix}
f_! : &\hfill {\mathcal E} / E_2 &\longrightarrow &{\mathcal E}/E_1 \, , \hfill \\
&(E \to E_2) &\longmapsto &(E \to E_2 \xrightarrow{f} E_1)
\end{matrix}
$$
et un adjoint \`a droite
$$
f_* : {\mathcal E} / E_2 \longrightarrow {\mathcal E}/E_1
$$
avec lequel il forme un morphisme de topos
$$
(f^* , f_*) : {\mathcal E} / E_2 \longrightarrow {\mathcal E}/E_1 \, .
$$
\end{prop}

\begin{remark}

Pour tous morphismes de ${\mathcal E}$
$$
E_3 \xrightarrow{ \ g \ } E_2 \xrightarrow{ \ f \ } E_1 \, ,
$$
le morphisme de topos localis\'es d\'efini par $f \circ g$
$$
((f \circ g)^* , (f \circ g)_*) : {\mathcal E}/E_3 \longrightarrow {\mathcal E}/E_1
$$
est canoniquement isomorphe au compos\'e des morphismes d\'efinis par $g$ et $f$
$$
{\mathcal E}/E_3 \xrightarrow{ \ (g^* , g_*) \ } {\mathcal E}/E_2 \xrightarrow{ \ (f^*,f_*) \ } {\mathcal E}/E_1 \, .
$$
\end{remark}

\begin{demo}

Pour tous morphismes $E' \to E_2$ et $E \to E_1$ de ${\mathcal E}$, les compl\'eter en un carr\'e commutatif
$$
\xymatrix{
E' \ar[d] \ar[r] &E \ar[d] \\
E_2 \ar[r]^-f &E_1
}
$$
\'equivaut \`a se donner un morphisme
$$
E' \longrightarrow E_2 \times_{E_1} E
$$
dont le compos\'e avec la projection $E_2 \times_{E_1} E \to E_2$ est le morphisme $E' \to E_2$.

\smallskip

Donc le foncteur $f^*$ admet pour adjoint \`a gauche $f_!$ et, a fortiori, il respecte les limites.

\smallskip

D'autre part, il respecte les colimites d'apr\`es la proposition \ref{propIII41}.

\smallskip

Cela implique d'apr\`es le corollaire \ref{corIII75} que le foncteur
$$
f^* : {\mathcal E}/E_1 \longrightarrow {\mathcal E}/E_2
$$
admet un adjoint \`a droite $f_*$ puisque la cat\'egorie ${\mathcal E}/E_1$ est un topos. 

\end{demo}

\bigskip

Ainsi, tout morphisme d'un topos ${\mathcal E}$
$$
f : E_2 \longrightarrow E_1
$$
induit un morphisme entre les topos localis\'es associ\'es
$$
(f^* , f_*) : {\mathcal E} / E_2 \longrightarrow {\mathcal E}/E_1
$$
qui est ``essentiel'' au sens de la d\'efinition suivante:

\begin{defn}\label{defIV84}

On dit qu'un morphisme de topos
$$
f = (f^* , f_*) : {\mathcal E}' \longrightarrow {\mathcal E}
$$
est ``essentiel'' si sa composante d'image r\'eciproque $f^*$ admet non seulement un adjoint \`a droite qui est la composante d'image directe $f_*$, mais aussi un adjoint \`a gauche 
$$
f_! : {\mathcal E}' \longrightarrow {\mathcal E} \, .
$$
\end{defn}

\begin{remarkqed}

Il r\'esulte du corollaire \ref{corIII75} qu'un morphisme de topos est ``essentiel'' si et seulement si sa composante d'image r\'eciproque respecte les limites arbitraires. 

\end{remarkqed}

\bigskip

Voici une description concr\`ete des foncteurs $f^* , f_*$ et $f_!$ lorsque ${\mathcal E}$ est pr\'esent\'e comme le topos des faisceaux sur un site et que $f : E_2 \to E_1$ est l'image par le foncteur canonique d'un morphisme carrable de la cat\'egorie sous-jacente:

\begin{lem}\label{lemIV85}

Soit ${\mathcal E} = \widehat{\mathcal C}_J$ le topos des faisceaux sur un site $({\mathcal C},J)$.

\smallskip

Soit
$$
f : E_2 \longrightarrow E_1
$$
un morphisme de ${\mathcal E}$ qui est l'image par le foncteur canonique
$$
\ell : {\mathcal C} \xhookrightarrow{ \ y \ } \widehat{\mathcal C} \xrightarrow{ \ j^* \ } \widehat{\mathcal C}_J = {\mathcal E}
$$
d'un morphisme de ${\mathcal C}$
$$
x : X_2 \longrightarrow X_1
$$
suppos\'e carrable.

\smallskip

Soient $J_1$ et $J_2$ les topologies des cat\'egories relatives ${\mathcal C} / X_1$ et ${\mathcal C}/X_2$ induites par $J$, si bien que les topos localis\'es ${\mathcal E}/E_1$ et ${\mathcal E}/E_2$ s'identifient \`a $\widehat{({\mathcal C}/X_1)}_{J_1}$ et $\widehat{({\mathcal C}/X_2)}_{J_2}$.

\smallskip

Alors:

\begin{listeimarge}

\item Le foncteur
$$
\begin{matrix}
\bullet \times_{X_1} x : &\hfill {\mathcal C}/X_1 &\longrightarrow &{\mathcal C}/X_2 \, , \hfill \\
&(X \to X_1) &\longmapsto &(X \times_{X_1} X_2 \to X_2)
\end{matrix}
$$
est un morphisme du site $({\mathcal C}/X_1 , J_1)$ vers le site $({\mathcal C}/X_2,J_2)$, et le morphisme de topos qu'il d\'efinit
$$
\widehat{({\mathcal C}/X_2)}_{J_2} \longrightarrow \widehat{({\mathcal C}/X_1)}_{J_1}
$$
s'identifie au morphisme
$$
(f^* , f_*) : {\mathcal E}/E_2 \longrightarrow {\mathcal E}/E_1 \, .
$$

\item En particulier, le foncteur
$$
f^* : {\mathcal E} / E_1 \longrightarrow {\mathcal E} / E_2
$$
s'inscrit dans un carr\'e commutatif
$$
\xymatrix{
{\mathcal C}/X_1 \ar[d]_-{\ell} \ar[rr]^-{\bullet \times_{X_1} x} &&{\mathcal C}/X_2 \ar[d]^-{\ell} \\
\widehat{({\mathcal C}/X_1)}_{J_1} \ar@{=}[d] &&\widehat{({\mathcal C}/X_2)}_{J_2} \ar@{=}[d] \\
{\mathcal E}/E_1 \ar[rr]^-{f^*} &&{\mathcal E}/E_2
}
$$
et, pour tout objet $F$ de ${\mathcal E}/E_1$ vu comme un faisceau sur $({\mathcal C}/X_1,J_1)$, on a
$$
f^* F = \varinjlim_{(X \to X_1 , \, x) \in \int\!F} \ell (X \times_{X_1} X_2 \to X_2) \, .
$$

\item Son adjoint \`a droite
$$
f_* : {\mathcal E} / E_2 \longrightarrow {\mathcal E}/E_1
$$
est la restriction aux faisceaux
$$
\widehat{({\mathcal C}/X_2)}_{J_2} \longrightarrow \widehat{({\mathcal C}/X_1)}_{J_1}
$$
du foncteur d'\'evaluation
$$
\begin{matrix}
\widehat{{\mathcal C}/X_2} &\longrightarrow &\widehat{{\mathcal C}/X_1} \, , \hfill \\
\hfill P &\longmapsto &P \circ (\bullet \times_{X_1} x) = [(X \to X_1) \mapsto P(X \times_{X_1} X_2 \to X_2)] \, .
\end{matrix}
$$

\item Son adjoint \`a gauche
$$
f_! : {\mathcal E}/E_2 \longrightarrow {\mathcal E}/E_1
$$
prolonge le foncteur de composition avec $x$
$$
\begin{matrix}
x \circ \bullet : &\hfill {\mathcal C}/X_2 &\longrightarrow &{\mathcal C}/X_1 \, , \hfill \\
&(X \to X_2) &\longmapsto &(X \to X_2 \xrightarrow{x} X_1)
\end{matrix}
$$
au sens qu'ils s'inscrivent dans un carr\'e commutatif
$$
\xymatrix{
{\mathcal C}/X_2 \ar[d]_-{\ell} \ar[rr]^-{x \, \circ \, \bullet} &&{\mathcal C}/X_1 \ar[d]^-{\ell} \\
\widehat{({\mathcal C}/X_2)}_{J_2} \ar@{=}[d] &&\widehat{({\mathcal C}/X_1)}_{J_1} \ar@{=}[d] \\
{\mathcal E}/E_2 \ar[rr]^-{f_!} &&{\mathcal E}/E_1
}
$$
et, pour tout objet $G$ de ${\mathcal E}/E_2$ vu comme un faisceau sur $({\mathcal C}/X_2,J_2)$, on a
$$
f_! \, G = \varinjlim_{(X \to X_2, \, x) \in \int\!G} \ell (X \longrightarrow X_2 \xrightarrow{ \ x \ } X_1) \, .
$$

\end{listeimarge}
\end{lem}

\begin{demo}

Le carr\'e
$$
\xymatrix{
{\mathcal C}/X_1 \ar[d]_-{\ell} \ar[rr]^-{\bullet \times_{X_1} x} &&{\mathcal C}/X_2 \ar[d]^-{\ell} \\
\widehat{({\mathcal C}/X_1)}_{J_1} \ar@{=}[d] &&\widehat{({\mathcal C}/X_2)}_{J_2} \ar@{=}[d] \\
\widehat{\mathcal C}_J / \ell (X_1) \ar@{=}[d] &&\widehat{\mathcal C}_J / \ell (X_2) \ar@{=}[d] \\
{\mathcal E}/E_1 \ar[rr] &&{\mathcal E}/E_2
}
$$
est commutatif car le foncteur canonique
$$
\ell : {\mathcal C} \longrightarrow \widehat{\mathcal C}_J
$$
respecte les limites finies.

\smallskip

Comme le foncteur
$$
\begin{matrix}
f^* : &\hfill {\mathcal E}/E_1 &\longrightarrow &{\mathcal E}/E_2 \, , \hfill \\
&(E \to E_1) &\longmapsto &(E \times_{E_1} E_2 \to E_2)
\end{matrix}
$$
respecte les colimites arbitraires, le compos\'e
$$
\widehat{{\mathcal C}/X_1} \longrightarrow \widehat{({\mathcal C}/X_1)}_{J_1} = {\mathcal E}/E_1 \xrightarrow{ \ f^* \ } {\mathcal E} / E_2 = \widehat{({\mathcal C}/X_2)}_{J_2}
$$
est l'unique prolongement \`a $\widehat{{\mathcal C}/X_1}$ respectant les colimites du foncteur compos\'e
$$
{\mathcal C}/X_1 \xrightarrow{ \ \bullet \times_{X_1} x \ } {\mathcal C}/X_2 \longrightarrow \widehat{({\mathcal C}/X_2)}_{J_2} = {\mathcal E}/E_2 
$$
qui donc est plat.

\smallskip

Comme il se factorise en
$$
\widehat{{\mathcal C}/X_1} \longrightarrow \widehat{({\mathcal C}/X_1)}_{J_1} = {\mathcal E}/E_1 \xrightarrow{ \ f^* \ } {\mathcal E}/E_2 = \widehat{({\mathcal C}/X_2)}_{J_2} \, ,
$$
le foncteur $\bullet \times_{X_1} x$ est un morphisme de sites, et il d\'efinit le morphisme de topos
$$
(f^* , f_*) : {\mathcal E}/E_2 \longrightarrow {\mathcal E} / E_1 \, .
$$

Ainsi, (i) est d\'emontr\'e.

\smallskip

(ii) et (iii) s'en d\'eduisent par le lemme \ref{lemIV42}.

\smallskip

Le carr\'e de (iv) est commutatif simplement parce que $\ell : {\mathcal C} \to \widehat{\mathcal C}_J$ est un foncteur, qui donc respecte la composition des morphismes.

\smallskip

Comme le foncteur $f_!$ est un adjoint \`a gauche, il respecte les colimites.

\smallskip

La formule finale de (iv) se d\'eduit alors de ce que tout faisceau $G$ sur le site $({\mathcal C}/X_2,J_2)$ s'\'ecrit comme la colimite
$$
G = \varprojlim_{(X \to X_2, \, x) \in \int\!G} \ell (X \longrightarrow X_2) \, .
$$

Cela termine la d\'emonstration du lemme. \end{demo}

\bigskip

Cette description s'\'etend au cas plus g\'en\'eral ou ${\mathcal E}$ est une cat\'egorie de faisceaux sur un site mais le morphisme $f : E_2 \to E_1$ de ${\mathcal E}$ ne provient pas n\'ecessairement d'un morphisme de la cat\'egorie sous-jacente au site:

\begin{cor}\label{corIV86}

Soit ${\mathcal E} = \widehat{\mathcal C}_J$ le topos des faisceaux sur un site $({\mathcal C},J)$.

\smallskip

Soit
$$
f : E_2 \longrightarrow E_1
$$
un morphisme de ${\mathcal E}$ qui est l'image par le foncteur de faisceautisation
$$
j^* : \widehat{\mathcal C} \longrightarrow \widehat{\mathcal C}_J
$$
d'un morphisme de pr\'efaisceaux sur ${\mathcal C}$
$$
p : P_2 \longrightarrow P_1 \, .
$$

Soient $J_1$ et $J_2$ les topologies des cat\'egories $\int \! P_1 = {\mathcal C}/P_1$ et $\int \! P_2 = {\mathcal C}/P_2$ induites par $J$, si bien que les topos localis\'es ${\mathcal E}/E_1$ et ${\mathcal E}/E_2$ s'identifient \`a $\widehat{({\mathcal C}/P_1)}_{J_1}$ et $\widehat{({\mathcal C}/P_2)}_{J_2}$.

\smallskip

Alors:

\begin{listeimarge}

\item Le foncteur 
$$
f^*: \xymatrix{
{\mathcal E}_1 / E_1 \ar@{=}[d] \ar[r] &{\mathcal E}_2 / E_2 \ar@{=}[d] \\
\widehat{({\mathcal C}/P_1)}_{J_1} &\widehat{({\mathcal C}/P_2)}_{J_2}
}
$$
associe \`a tout faisceau $F$ sur $({\mathcal C}/P_1 , J_1)$ la faisceautisation du pr\'efaisceau sur ${\mathcal C}/P_2$
$$
(X' , y(X') \to P_2) \longmapsto \varinjlim_{\left(\begin{matrix} X , &\!\!\!\!y (X') \longrightarrow y(X) \\ &\downarrow \qquad\qquad \downarrow \\ &P_2 \ \ \longrightarrow \ \ P_1 \end{matrix} \right)} F(X,y(X) \to P_1) 
$$
o\`u la colimite est calcul\'ee sur la cat\'egorie des objets $X$ de ${\mathcal C}$ munis d'un morphisme $y(X) \to P_1$ et d'un morphisme $X' \to X$ qui rendent commutatif le carr\'e de $\widehat{\mathcal C}$:
$$
\xymatrix{
y(X') \ar[d] \ar[r] &y(X) \ar[d] \\
P_2 \ar[r] &P_1
}
$$

\item Le foncteur
$$
f_*: \xymatrix{
{\mathcal E} / E_2 \ar@{=}[d] \ar[r] &{\mathcal E} / E_1 \ar@{=}[d] \\
\widehat{({\mathcal C}/P_2)}_{J_2} &\widehat{({\mathcal C}/P_1)}_{J_1}
}
$$
associe \`a tout faisceau $G$ sur $({\mathcal C}/P_2 , J_2)$ le faisceau sur ${\mathcal C}/P_1$
$$
(X, y(X) \to P_1) \longmapsto \varprojlim_{\left(\begin{matrix} X' , &\!\!\!\!y (X') \longrightarrow y(X) \\ &\downarrow \qquad\qquad \downarrow \\ &P_2 \ \ \longrightarrow \ \ P_1 \end{matrix} \right)} G(X',y(X') \to P_2) 
$$
o\`u la limite est calcul\'ee sur la cat\'egorie des objets $X'$ de ${\mathcal C}$ munis d'un morphisme $y(X') \to P_2$ et d'un morphisme $X' \to X$ qui rendent commutatif le carr\'e de $\widehat{\mathcal C}$:
$$
\xymatrix{
y(X') \ar[d] \ar[r] &y(X) \ar[d] \\
P_2 \ar[r] &P_1
}
$$

\item Le foncteur
$$
f_! : \xymatrix{
{\mathcal E} / E_2 \ar@{=}[d] \ar[r] &{\mathcal E} / E_1 \ar@{=}[d] \\
\widehat{({\mathcal C}/P_2)}_{J_2} &\widehat{({\mathcal C}/P_1)}_{J_1}
}
$$
associe \`a tout faisceau $G$ sur $({\mathcal C}/P_2,J_2)$ la faisceautisation du pr\'efaisceau sur ${\mathcal C}/P_2$
$$
(X, y(X) \to P_1) \longmapsto \varinjlim_{\left(\begin{matrix} X' , &\!\!\!\!y (X') \longleftarrow y(X) \\ &\downarrow \qquad\qquad \downarrow \\ &P_2 \ \ \longrightarrow \ \ P_1 \end{matrix} \right)} G(X',y(X') \to P_2) 
$$
o\`u la colimite est calcul\'ee sur la cat\'egorie des objets $X'$ de ${\mathcal C}$ munis d'un morphisme $y(X') \to P_2$ et d'un morphisme $X \to X'$ tels que $y(X) \to P_1$ soit \'egal au compos\'e
$$
y(X) \longrightarrow y(X') \longrightarrow P_2 \longrightarrow P_1 \, .
$$

\end{listeimarge}
\end{cor}

\begin{demo}
\begin{listeisansmarge}
\item[(i)] Comme les colimites respectent les colimites et que le foncteur de faisceautisation $j^* : \widehat{{\mathcal C}/P_2} \to \widehat{({\mathcal C}/P_2)}_{J_2}$ les respecte \'egalement, le foncteur d\'efini par la formule de (i) respecte les colimites.

\smallskip

Si $F$ est l'image par le foncteur canonique
$$
{\mathcal C}/P_1 \longrightarrow \widehat{({\mathcal C}/P_1)}_{J_1}
$$
d'un objet $(X,y(X) \to P_1)$ de ${\mathcal C}/P_1$, ce foncteur transforme $F$ en la faisceautisation de 
$$P_2 \times_{P_1} y(X) \longrightarrow P_2
$$
qui est 
$$
E_2 \times_{E_1} j^* \circ y(X) \longrightarrow E_2 \, .
$$
Il co{\"\i}ncide donc bien avec le foncteur
$$
\begin{matrix}
f^* : &\hfill {\mathcal E}/E_1 &\longrightarrow &{\mathcal E}/E_2 \, , \hfill \\
&(E \to E_1) &\longmapsto &(E \times_{E_1} E_2 \to E_2) \, .
\end{matrix}
$$

\item[(iii)] De m\^eme, le foncteur d\'efini par la formule de (iii) respecte les colimites.

\smallskip

Si $G$ est l'image par le foncteur canonique
$$
{\mathcal C} / P_2 \longrightarrow \widehat{({\mathcal C}/P_2)}_{J_2}
$$
d'un objet $(X' , y(X') \to P_2)$ de ${\mathcal C}/P_2$, ce foncteur transforme $G$ en la faisceautisation de
$$
y(X') \longrightarrow P_2 \longrightarrow P_1
$$
qui est le compos\'e
$$
j^* \circ y (X') \longrightarrow E_2 \longrightarrow E_1 \, .
$$
Il co{\"\i}ncide donc bien avec le foncteur
$$
\begin{matrix}
f_! : &\hfill {\mathcal E}/E_2 &\longrightarrow &{\mathcal E}/E_1 \, , \hfill \\
&(E \to E_2) &\longmapsto &(E \to E_2 \to E_1) \, .
\end{matrix}
$$

\item[(ii)] Pour tout objet $(X , y(X) \to P_1)$ de ${\mathcal C}/P_1$ et tout objet $G$ de ${\mathcal E}/E_2 = \widehat{({\mathcal C}/P_2)}_{J_2}$, l'ensemble
$$
f_* G (X , y(X) \longrightarrow P_1)
$$
s'identifie \`a celui des morphismes de $\widehat{({\mathcal C}/P_1)}_{J_1}$
$$
\ell (X,y(X) \to P_1) \longrightarrow f_* G
$$
et donc \`a celui des morphismes de $\widehat{({\mathcal C}/P_2)}_{J_2}$
$$
f^* \circ \ell (X,y(X) \to P_1) \longrightarrow G \, .
$$
Or l'objet $f^* \circ \ell (X,y(X) \to P_1)$ de $\widehat{({\mathcal C}/P_2)}_{J_2}$ est la faisceautisation de l'objet de $\widehat{{\mathcal C}/P_2} = \widehat{\mathcal C} / P_2$
$$
y(X) \times_{P_1} P_2 \longrightarrow P_2 \, .
$$
Donc l'ensemble
$$
f_* G (X,y(X) \to P_1)
$$
s'identifie \`a celui des morphismes de $\widehat{{\mathcal C}/P_2}$
$$
y(X) \times_{P_1} P_2 \longrightarrow G \, .
$$
On conclut par l'identit\'e de $\widehat{{\mathcal C}/P_2}$
$$
y(X) \times_{P_1} P_2 = \varinjlim_{\left(\begin{matrix} X' , &\!\!\!\!y (X') \longrightarrow y(X) \\ &\downarrow \qquad\qquad \downarrow \\ &P_2 \ \ \longrightarrow \ \ P_1 \end{matrix} \right)} y(X') \, .
$$
Cela ach\`eve la preuve du corollaire.

\end{listeisansmarge}
\end{demo}

\subsection{Sous-topos ouverts et sous-topos ferm\'es}\label{subsecIV83}

\medskip

Caract\'erisons les morphismes d'un topos dont le morphisme de topos localis\'es associ\'e est un plongement:

\begin{lem}\label{lemIV87}

Soit 
$$
f : E_2 \longrightarrow E_1
$$
un morphisme d'un topos ${\mathcal E}$.

\smallskip

Alors le morphisme de topos localis\'es associ\'e
$$
(f^* , f_*) : {\mathcal E} / E_2 \longrightarrow {\mathcal E} / E_1
$$
est un plongement si et seulement si $f$ est un monomorphisme de ${\mathcal E}$.
\end{lem}

\begin{demo}

Par d\'efinition, le morphisme de topos $(f^* , f_*)$ est un plongement si et seulement si sa composante d'image directe $f_*$ est pleinement fid\`ele.

\smallskip

Comme $f_*$ est adjoint \`a droite de $f^*$ qui a aussi un adjoint \`a gauche $f_!$, cela \'equivaut d'apr\`es la remarque qui suit le lemme \ref{lemI82} \`a demander que le foncteur
$$
\begin{matrix}
f_! : &\hfill {\mathcal E}/E_2 &\longrightarrow &{\mathcal E}/E_1 \, , \hfill \\
&(E \xrightarrow{g} E_2) &\longmapsto &(E \xrightarrow{f \circ g} E_1) 
\end{matrix}
$$
soit pleinement fid\`ele.

\smallskip

Or, ce foncteur $f_!$ est toujours fid\`ele.

\smallskip

Il est plein si, pour tous morphismes de ${\mathcal E}$
$$
E \xrightarrow{ \ g \ } E_2 \, , \quad E' \xrightarrow{ \ g' \ } E_2 \quad \mbox{et} \quad E' \xrightarrow{ \ e \ } E \, ,
$$
la condition
$$
f \circ g \circ e = f \circ g'
$$
\'equivaut \`a la condition
$$
g \circ e = g' \, .
$$

Cela revient \`a demander que $f$ soit un monomorphisme de ${\mathcal E}$. 

\end{demo}

\bigskip

On pose alors:

\begin{defn}\label{defIV88}

On appelle sous-topos ouverts d'un topos ${\mathcal E}$ les sous-topos de ${\mathcal E}$ qui sont aussi des localisations de ${\mathcal E}$.

\smallskip

Autrement dit, ce sont les topos de la forme
$$
{\mathcal E}_U = {\mathcal E}/U
$$
associ\'es aux objets ``sous-terminaux'' $U$ de ${\mathcal E}$, c'est-\`a-dire aux sous-objets de l'objet terminal $1_{\mathcal E}$.
\end{defn}

\bigskip

\begin{remarkqed}

Si ${\mathcal E} = {\mathcal E}_X$ est le topos des faisceaux sur un espace topologique $X$, on a vu dans le th\'eor\`eme \ref{thmIV112} (i) que l'ensemble des ouverts de $X$ s'identifie \`a celui des objets sous-terminaux de ${\mathcal E}_X$, avec sa relation d'ordre, ses r\'eunions arbitraires et ses intersections finies.

\smallskip

De plus, pour tout ouvert $U$ de $X$ qui correspond \`a un objet sous-terminal $S_U$ de ${\mathcal E}_X$, le sous-topos ouvert 
$$
{\mathcal E}_X / S_U
$$
s'identifie au topos ${\mathcal E}_U$ des faisceaux sur $U$.

\smallskip

Ainsi, la notion de ``sous-topos ouvert'' g\'en\'eralise celle d'ouvert d'un espace topologique. 

\end{remarkqed}

\bigskip

D\'ecrivons les topologies qui d\'efinissent les sous-topos ouverts d'un topos de faisceaux sur un site:

\begin{lem}\label{lemIV89}

Soit ${\mathcal E} = \widehat{\mathcal C}_J$ le topos des faisceaux sur un site $({\mathcal C},J)$ muni du foncteur canonique $\ell : {\mathcal C} \xrightarrow{y} \widehat{\mathcal C} \xrightarrow{j^*} \widehat{\mathcal C}_J$.

\smallskip

Soient $U$ un objet sous-terminal de ${\mathcal E}$, ${\mathcal E}/U$ le sous-topos ouvert de ${\mathcal E}$ associ\'e, et $J_U \supseteq J$ la topologie de ${\mathcal C}$ qui d\'efinit ce sous-topos.

\smallskip

Alors:

\begin{listeimarge}

\item Une famille de morphismes de ${\mathcal C}$
$$
X_i \longrightarrow X \, , \qquad i \in I \, ,
$$
est $J_U$-couvrante si et seulement si elle satisfait les conditions \'equivalentes:

\medskip

$
\left\{\begin{matrix}
{\rm (1)} &\mbox{La famille de morphismes de $\widehat{\mathcal C}_J = {\mathcal E}$} \hfill \\
{ \ } \\
&\ell (X_i) \times U \longrightarrow \ell (X) \times U \, , \qquad i \in I \, , \\
{ \ } \\
&\mbox{est globalement \'epimorphique.} \hfill \\
{ \ } \\
{\rm (2)} &\mbox{Pour tout morphisme de ${\mathcal C}$} \hfill \\
&X' \xrightarrow{ \ x \ } X \, , \\
{ \ } \\
&\mbox{tel que $U(X') \ne \emptyset$, il existe une famille $J$-couvrante} \hfill \\
{ \ } \\
&X'_j \longrightarrow X' \\
{ \ } \\
&\mbox{tel que chaque compos\'e $X'_j \to X' \xrightarrow{x} X$ se factorise \`a travers au moins un $X_i \to X$.} \hfill
\end{matrix}\right.
$

\medskip

\item Si $P$ est un objet de $\widehat{\mathcal C}$ dont la faisceautisation est $U$, un objet $F$ de $\widehat{\mathcal C}_J = {\mathcal E}$ est un $J_U$-faisceau si et seulement si, pour tout objet $X$ de ${\mathcal C}$, l'application canonique
$$
F(X) \longrightarrow \varprojlim_{(X',x) \in \int\!y(X) \times P} F(X')
$$
est bijective.
\end{listeimarge}
\end{lem}

\begin{remark}

Si la cat\'egorie ${\mathcal C}$ est stable par produits finis et contient un objet $P$ dont l'image par $\ell : {\mathcal C} \to \widehat{\mathcal C}_J = {\mathcal E}$ est $U$, la condition de (ii) est que, pour tout objet $X$ de ${\mathcal C}$, l'application
$$
F(X) \longrightarrow F(X \times P)
$$
est bijective.
\end{remark}

\bigskip

\begin{demo}
\begin{listeisansmarge}
\item Une famille de morphismes de ${\mathcal C}$
$$
X_i \longrightarrow X \, , \qquad i \in I \, ,
$$
est $J_U$-couvrante si et seulement si son image par le foncteur
$$
{\mathcal C} \xrightarrow{ \ \ell \ } \widehat{\mathcal C}_J = {\mathcal E} \xrightarrow{ \ u^* \ } {\mathcal E}/U
$$
compos\'e de $\ell$ et de la composante d'image r\'eciproque $u^*$ du morphisme
$$
(u^* , u_*) : {\mathcal E}/U \longrightarrow {\mathcal E}
$$
associ\'e \`a $u : U \hookrightarrow 1_{\mathcal E}$, est une famille \'epimorphique.

\smallskip

Cette condition \'equivaut \`a (1) puisque $u^*$ est le foncteur
$$
\begin{matrix}
{\mathcal E} &\longrightarrow &{\mathcal E}/U \, , \\
E &\longmapsto &E \times U \, .
\end{matrix}
$$

L'\'equivalence de (1) et (2) provient de ce que, pour tout objet $X'$ de ${\mathcal C}$, on a
$$
U(X') = {\rm Hom} (\ell (X'),U) \, .
$$

\item Un objet $F$ de $\widehat{\mathcal C}_J = {\mathcal E}$ est un $J_U$-faisceau si et seulement si le morphisme canonique
$$
F \longrightarrow u_* \circ u^* F
$$
est un isomorphisme.

\smallskip

Or, pour tout objet $X$ de ${\mathcal C}$, on a
\begin{eqnarray}
u_* \circ u^* F (X) &= &{\rm Hom} (\ell (X) , u_* \circ u^* F) \nonumber \\
&= &{\rm Hom} (u^* \circ \ell (X), u^* F) \nonumber \\
&= &{\rm Hom}_U (\ell (X) \times U , F \times U) \nonumber \\
&= &{\rm Hom} (\ell (X) \times U , F) \nonumber \\
&= &\varprojlim_{(X',x) \in \int\!\!y(X) \times P} F(X') \nonumber
\end{eqnarray}
puisque par hypoth\`ese $U = j^* P$. 
\end{listeisansmarge}
\end{demo}

\bigskip

Les sous-topos ouverts des topos poss\`edent les propri\'et\'es habituelles des ouverts des espaces topologiques:

\begin{prop}\label{propIV810}
\begin{listeimarge}
\item Les r\'eunions de sous-topos ouverts sont des topos ouverts.

\smallskip

Plus pr\'ecis\'ement, si les $U_i$ sont des objets sous-terminaux d'un topos ${\mathcal E}$ et $\underset{i \in I}{\bigvee} \, U_i = U$ est leur r\'eunion, on a
$$
\bigvee_{i \in I} {\mathcal E}/U_i = {\mathcal E}/U \, .
$$

\item Les intersections finies de sous-topos ouverts sont des sous-topos ouverts.

\smallskip

Plus pr\'ecis\'ement, si $U_1 , \cdots , U_n$ sont des objets sous-terminaux d'un topos ${\mathcal E}$ et $U_1 \wedge \cdots \wedge U_n = U$ est leur intersection, on a
$$
{\mathcal E} / U_1 \wedge \cdots \wedge {\mathcal E} / U_n = {\mathcal E} / U \, .
$$

\item Les images r\'eciproques de sous-topos ouverts sont des sous-topos ouverts.

\smallskip

Plus pr\'ecis\'ement, si
$$
f = (f^* , f_*) : {\mathcal E}' \longrightarrow {\mathcal E}
$$
est un morphisme de topos et $U$ est un objet sous-terminal de ${\mathcal E}$, alors $f^* U$ est un objet sous-terminal de ${\mathcal E}'$ et on a
$$
f^{-1} ({\mathcal E}/U) = {\mathcal E}'/f^* U \, .
$$

\item Pour tout morphisme de topos
$$
f = (f^* , f_*) : {\mathcal E}' \longrightarrow {\mathcal E}
$$
et tout objet sous-terminal $U$ de ${\mathcal E}$, $f$ se factorise \`a travers le sous-topos ouvert
$$
{\mathcal E} / U \xhookrightarrow{ \ { \ } \ } {\mathcal E}
$$
si et seulement si le morphisme canonique
$$
f^* U \longrightarrow 1_{{\mathcal E}'}
$$
est un isomorphisme.
\end{listeimarge}
\end{prop}

\bigskip

\begin{demo}
\begin{listeisansmarge}
\item Chaque $U_i$ est contenu dans $U$ donc on a une in\'egalit\'e
$$
\bigvee_{i \in I} {\mathcal E}/U_i \leq {\mathcal E}/U \, .
$$
Pour montrer qu'il y a \'egalit\'e, on peut supposer que ${\mathcal E} = \widehat{\mathcal C}_J$ est le topos des faisceaux sur un site $({\mathcal C},J)$ et consid\'erer les topologies $J_U$ et $J_{U_i}$, $i \in I$, de ${\mathcal C}$ qui d\'efinissent les sous-topos ouverts ${\mathcal E}/U$ et ${\mathcal E}/U_i$ de ${\mathcal E}$.

\smallskip

D'apr\`es le lemme \ref{lemIV89}, il s'agit de montrer que si une famille de morphismes de ${\mathcal C}$
$$
(X_j \longrightarrow X)_j
$$
est telle que chaque famille
$$
(\ell (X_j) \times U_i \longrightarrow \ell (X) \times U_i)_j
$$
est globalement \'epimorphique, il en va de m\^eme de la famille
$$
(\ell (X_j) \times U \longrightarrow \ell (X) \times U)_j \, .
$$
Cela r\'esulte de ce que la famille de monomorphismes
$$
(U_i \xhookrightarrow{ \ { \ } \ } U)_{i \in I}
$$
est globalement \'epimorphique, donc aussi la famille
$$
(\ell (X) \times U_i \xhookrightarrow{ \ { \ } \ } \ell (X) \times U)_{i \in I} \, .
$$

\item On peut supposer que ${\mathcal E} = \widehat{\mathcal C}_J$ est le topos des faisceaux sur un site $({\mathcal C} , J)$ tel que ${\mathcal C}$ soit une sous-cat\'egorie pleine de ${\mathcal E}$, soit stable par limites finies et contiennent les objets $1_{\mathcal E} , U_1 , \cdots , U_n$ et donc aussi $U = U_1 \times \cdots \times U_n$.

\smallskip

Notant $J_U$ et $J_{U_i}$, $1 \leq i \leq n$, les topologies de ${\mathcal C}$ qui d\'efinissent les sous-topos ${\mathcal E}/U$ et ${\mathcal E}/U_i$, on sait qu'un objet $F$ de ${\mathcal E} = \widehat{\mathcal C}_J$ est un faisceau pour la topologie $J_U$ [resp. $J_{U_i}$, $1 \leq i \leq n$] si et seulement si, pour tout objet $X$ de ${\mathcal C}$, l'application
$$
\begin{matrix}
\hfill F(X) &\longrightarrow &F (X \times U) \hfill \\
\mbox{[resp.} \qquad F(X) &\longrightarrow &F(X \times U_i) \ \mbox{]}
\end{matrix}
$$
est bijective.

\smallskip

Les \'egalit\'es $U = U_1 \times \cdots \times U_n$ et $U_i \times U = U$, $1 \leq i \leq n$, montrent que $F$ est un faisceau pour la topologie $J_U$ si et seulement si c'est un faisceau pour chacune des topologies $J_{U_i}$, $1 \leq i \leq n$.

\smallskip

D'apr\`es la proposition \ref{propIV79} (ii), cela signifie que
$$
{\mathcal E}/U = {\mathcal E}/U_1 \wedge \cdots \wedge {\mathcal E}/U_n \, .
$$

\item Le foncteur $f^*$ respecte les limites finies donc il transforme l'objet sous-terminal $U$ de ${\mathcal E}$ en un objet sous-terminal de ${\mathcal E}'$.

\smallskip

Le carr\'e
$$
\xymatrix{
{\mathcal E} \ar[d]_-{f^*} \ar[rr]^-{\bullet \times U} &&{\mathcal E}/U \ar[d]^-{f^*} \\
{\mathcal E}' \ar[rr]^-{\bullet \times f^* U} &&{\mathcal E}' / f^* U
}
$$
est commutatif \`a isomorphisme canonique pr\`es, et d\'efinit un carr\'e de morphismes de topos
$$
\xymatrix{
{\mathcal E}' / f^* U \ \ar[d] \ar@{^{(}->}[r] &{\mathcal E}' \ar[d] \\
{\mathcal E}/U \ \ar@{^{(}->}[r] &{\mathcal E}
}
$$
commutatif \`a isomorphisme canonique pr\`es.

\smallskip

On a donc une in\'egalit\'e
$$
{\mathcal E}' / f^* U \leq f^{-1} ({\mathcal E}/U) \, .
$$

Pour montrer qu'il y a \'egalit\'e, on peut supposer que ${\mathcal E} = \widehat{\mathcal C}_J$ et ${\mathcal E}' = \widehat{\mathcal D}_K$ sont les topos des faisceaux sur deux sites $({\mathcal C},J)$ et $({\mathcal D},K)$ tels que ${\mathcal C}$ et ${\mathcal D}$ soient deux sous-cat\'egories pleines de ${\mathcal E}$ et ${\mathcal E}'$, qu'elles soient stables par limites finies et que le foncteur
$$
f^* : {\mathcal E} \longrightarrow {\mathcal E}'
$$
se restreigne en un foncteur
$$
\rho : {\mathcal C} \longrightarrow {\mathcal D}
$$
qui donc est un morphisme de sites et d\'efinit $f : {\mathcal E}' \to {\mathcal E}$.

\smallskip

La topologie $J_U$ [resp. $K_{f^* U}$] de ${\mathcal C}$ [resp. ${\mathcal D}$] qui d\'efinit le sous-topos ${\mathcal E} / U$ [resp. ${\mathcal E}' / f^* U$] de ${\mathcal E} = \widehat{\mathcal C}_J$ [resp. ${\mathcal E}' = \widehat{\mathcal D}_K$] est la plus petite topologie qui contient $J$ [resp. $K$] et pour laquelle le morphisme
$$
\begin{matrix}
\hfill U &\xhookrightarrow{ \ { \ } \ } &1_{\mathcal C} \hfill \\
\mbox{[resp.} \qquad f^* U &\xhookrightarrow{ \ { \ } \ } &1_{\mathcal D} \ \mbox{]}
\end{matrix}
$$
est couvrant.

\smallskip

Comme $f^* U \hookrightarrow 1_{\mathcal D}$ est image de $U \hookrightarrow 1_{\mathcal C}$ par le foncteur $\rho = f^*$, on conclut d'apr\`es la proposition \ref{propIV713}~(ii) que 
$$
{\mathcal E}' / f^* U = f^{-1} ({\mathcal E}/U) \, .
$$

\item Un morphisme de topos
$$
f = (f^* , f_*) : {\mathcal E}' \longrightarrow {\mathcal E}
$$
se factorise \`a travers le sous-topos ouvert
$$
{\mathcal E}/U \xhookrightarrow{ \ { \ } \ } {\mathcal E}
$$
associ\'e \`a un objet sous-terminal $U$ de ${\mathcal E}$, si et seulement si on a
$$
f^{-1} ({\mathcal E}/U) = {\mathcal E}' \, .
$$
Or, d'apr\`es (iii), on a
$$
f^{-1} ({\mathcal E}/U) = {\mathcal E}' / f^* U \, .
$$
Cela prouve qu'il y a factorisation si et seulement si le morphisme de ${\mathcal E}'$
$$
f^* U \xhookrightarrow{ \ { \ } \ } 1_{{\mathcal E}'}
$$
est un isomorphisme. 
\end{listeisansmarge}
\end{demo}

\bigskip

Il est naturel de poser \`a la suite de la d\'efinition des sous-topos ouverts:

\begin{defn}\label{defIV811}

Un sous-topos d'un topos ${\mathcal E}$
$$
{\mathcal E}' \xhookrightarrow{ \ { \ } \ } {\mathcal E}
$$
est appel\'e un ``sous-topos ferm\'e'' s'il est le compl\'ementaire
$$
{\mathcal E}' = {\mathcal E}_U^c = {\mathcal E} \backslash {\mathcal E}_U
$$
d'un sous-topos ouvert
$$
{\mathcal E}_U = {\mathcal E} / U \xhookrightarrow{ \ { \ } \ } {\mathcal E}
$$
associ\'e \`a un objet sous-terminal $U$ de ${\mathcal E}$.
\end{defn}

\bigskip

\begin{remarkqed}

Si ${\mathcal E} = {\mathcal E}_X$ est le topos des faisceaux sur un espace topologique $X$ et $U$ est un ouvert de $X$ compl\'ementaire d'un ferm\'e $Z$, il r\'esulte de la remarque (iii) qui suit le corollaire \ref{corIV712} que
$$
{\mathcal E}_Z = {\mathcal E} \backslash {\mathcal E}_U \, .
$$

Cela justifie le nom de ``sous-topos ferm\'e'' donn\'e aux compl\'ementaires des sous-topos ouverts de n'importe quel topos. \end{remarkqed}

\bigskip

Nous allons d\'emontrer:

\begin{prop}\label{propIV812}

Soit $U$ un objet terminal d'un topos ${\mathcal E}$.

\smallskip

Soient ${\mathcal E}_U = {\mathcal E}/U$ le sous-topos ouvert associ\'e et ${\mathcal E}_U^c = {\mathcal E} \backslash {\mathcal E}_U$ le sous-topos ferm\'e compl\'ementaire.

\smallskip

Alors:

\begin{listeimarge}

\item L'intersection de ${\mathcal E}_U$ et de ${\mathcal E}_U^c$ est le sous-topos minimal ${\mathcal E}_{\emptyset}$ dont tout objet est terminal.

\medskip

\item On a par cons\'equent
$$
{\mathcal E} \backslash {\mathcal E}_U^c = {\mathcal E}_U \, .
$$
\end{listeimarge}
\end{prop}

\begin{remark}

Dans le cas o\`u ${\mathcal E} = {\mathcal E}_X$ est le topos des faisceaux sur un espace topologique $X$, et $U$ est un ouvert de $X$ de compl\'ementaire le ferm\'e $Z$, on savait d\'ej\`a d'apr\`es la remarque qui suit la proposition \ref{propIV79} que
$$
{\mathcal E}_U \wedge {\mathcal E}_Z = {\mathcal E}_{\emptyset}
$$
d'o\`u il r\'esulte
$$
{\mathcal E}_U = {\mathcal E}_X \backslash {\mathcal E}_Z \, .
$$
\end{remark}

\begin{demosansqed}

D'apr\`es le corollaire \ref{corIV712} (iii), (i) et (ii) sont cons\'equences du lemme suivant:
\end{demosansqed}

\begin{lem}\label{lemIV813}

Soit ${\mathcal E} = \widehat{\mathcal C}_J$ le topos des faisceaux sur un site $({\mathcal C},J)$ muni du foncteur canonique $\ell : {\mathcal C} \xhookrightarrow{y} \widehat{\mathcal C} \xrightarrow{j^*} \widehat{\mathcal C}_J$.

\smallskip

Soit $U$ un objet sous-terminal de ${\mathcal E}$.

\smallskip

Alors:

\begin{listeimarge}

\item Il existe une topologie $J_U^c$ de ${\mathcal C}$ pour laquelle une famille de morphismes
$$
X_i \longrightarrow X \, , \qquad i \in I \, ,
$$
est couvrante si et seulement si elle satisfait les conditions \'equivalentes suivantes:

\medskip

$
\left\{\begin{matrix}
{\rm (1)} &\mbox{La famille des morphismes de $\widehat{\mathcal C}_J = {\mathcal E}$} \hfill \\
{ \ } \\
&(\ell (X_i) \longrightarrow \ell (X))_{i \in I} \quad \mbox{compl\'et\'ee par} \quad \ell (X) \times U \to \ell (X) \\
{ \ } \\
&\mbox{est globalement \'epimorphique.} \hfill \\
{ \ } \\
{\rm (2)} &\mbox{Il existe une famille de morphismes de ${\mathcal C}$} \hfill \\
{ \ } \\
&Y_j \longrightarrow X \\
{ \ } \\
&\mbox{telle que $U(Y_j) \ne \emptyset$, $\forall \, j$, et que la r\'eunion des deux familles} \hfill \\
{ \ } \\
&X_i \longrightarrow X \, , \quad i \in I \, , \qquad \mbox{et} \qquad Y_j \longrightarrow X \\
{ \ } \\
&\mbox{soit $J$-couvrante.} \hfill
\end{matrix} \right.
$

\medskip

\item L'intersection de la topologie $J_U^c$ et de la topologie $J_U$ de ${\mathcal C}$ qui d\'efinit le sous-topos ouvert
$$
{\mathcal E}_U = {\mathcal E} / U \xhookrightarrow{ \ { \ } \ } {\mathcal E} = \widehat{\mathcal C}_J
$$
est \'egale \`a $J$.

\medskip

\item Un objet $F$ de $\widehat{\mathcal C}_J = {\mathcal E}$ est un $J_U^c$-faisceau si et seulement si le morphisme canonique
$$
F \times U \longrightarrow U
$$
est un isomorphisme, c'est-\`a-dire si $F$ a une unique section sur tout objet $X$ de ${\mathcal C}$ tel que $U(X) \ne \emptyset$.

\medskip

\item Un objet $F$ de $\widehat{\mathcal C}_J = {\mathcal E}$ qui est \`a la fois un $J_U$-faisceau et un $J_U^c$-faisceau est l'objet terminal de ${\mathcal E}$.
\end{listeimarge}
\end{lem}

\begin{demolem}
\begin{listeisansmarge}
\item L'\'equivalence des conditions (1) et (2) r\'esulte du lemme \ref{lemIII53}.

\smallskip

Pour tout objet $X$ de ${\mathcal C}$, notons $J_U^c(X)$ l'ensemble des cribles de $X$ qui contiennent une famille de morphismes $X_i \to X$, $i \in I$, v\'erifiant les conditions (1) et (2).

\smallskip

Il faut v\'erifier que $J_U^c$ satisfait les trois axiomes des topologies de Grothendieck.

\smallskip

C'est \'evident pour la maximalit\'e.

\smallskip

Pour la stabilit\'e, consid\'erons un morphisme $X' \xrightarrow{x} X$ de ${\mathcal C}$ et une famille de morphismes
$$
X_i \longrightarrow X \, , \qquad i \in I \, ,
$$
qui devient $J$-couvrante une fois compl\'et\'ee par une famille
$$
Y_j \longrightarrow X
$$
telle que 
$$
U(Y_j) \ne \emptyset \, , \qquad \forall \, j \, .
$$

Alors $X'$ admet une famille $J$-couvrante constitu\'ee de morphismes
$$
X'_{i'} \longrightarrow X'
$$
qui s'inscrivent dans des carr\'es commutatifs
$$
\xymatrix{
X'_{i'} \ar[d] \ar[r] &X_i \ar[d] \\
X' \ar[r]^-x &X
}
$$
et de morphismes
$$
Y'_{j'} \longrightarrow X'
$$
qui s'inscrivent dans des carr\'es commutatifs
$$
\xymatrix{
Y'_{j'} \ar[d] \ar[r] &Y_j \ar[d] \\
X' \ar[r]^-x &X
}
$$
et v\'erifient a fortiori
$$
U(Y'_{j'}) \ne \emptyset \, .
$$

Ainsi, l'image r\'eciproque par $X' \xrightarrow{x} X$ du crible de $X$ engendr\'e par la famille des $X_i \to X$ est \'el\'ement de $J_U^c (X')$.

\smallskip

Cela montre la stabilit\'e.

\smallskip

Pour la transitivit\'e, consid\'erons une famille $J$-couvrante d'un objet $X$ de ${\mathcal C}$ constitu\'ee de morphismes
$$
X_i \longrightarrow X \, , \qquad i \in I \, ,
$$
et de morphismes
$$
Y_j \longrightarrow X
$$
tels que
$$
U(Y_j) \ne \emptyset \, , \qquad \forall \, j \, .
$$

Et consid\'erons pour chaque indice $i \in I$ une famille $J$-couvrante de $X_i$ constitu\'ee de morphismes
$$
X_{i,k} \longrightarrow X_i
$$
et de morphismes
$$
Y_{i,\ell} \longrightarrow X_i
$$
tels que
$$
U(Y_{i,\ell}) \ne \emptyset \, .
$$

Alors la famille des morphismes compos\'es
$$
X_{i,k} \longrightarrow X_i \longrightarrow X
$$
compl\'et\'ee par les morphismes
$$
Y_{i,\ell} \longrightarrow X_i \longrightarrow X
$$
et les morphismes
$$
Y_j \longrightarrow X
$$
est $J$-couvrante.

\smallskip

Comme $U(Y_j) \ne \emptyset$, $\forall \, j$, et $U(Y_{i,\ell}) \ne \emptyset$, $\forall \, i,\ell$, cela implique que $J_U^c$ satisfait l'axiome de transitivit\'e.

\smallskip

Ainsi, $J_U^c$ est bien une topologie de ${\mathcal C}$.

\medskip

\item Soit $(X_i \to X)_{i \in I}$ une famille de morphismes de ${\mathcal C}$ qui est couvrante \`a la fois pour la topologie $J_U$ et pour la topologie $J_U^c$.

\smallskip

Ainsi, la famille des
$$
\ell (X_i) \times U \longrightarrow \ell (X) \times U \, , \qquad i \in I \, ,
$$
est globalement \'epimorphique, ainsi que la famille des 
$$
\ell (X_i) \longrightarrow \ell (X) \, , \qquad i \in I \, ,
$$
compl\'et\'ee par le morphisme
$$
\ell (X) \times U \longrightarrow \ell (X) \, .
$$
Il en r\'esulte que la famille des
$$
\ell (X_i) \longrightarrow \ell (X)
$$
compl\'et\'ee par les
$$
\ell (X_i) \times U \longrightarrow \ell (X)
$$
est globalement \'epimorphique, donc la seule famille des
$$
\ell (X_i) \longrightarrow \ell (X) \, , \qquad i \in I \, ,
$$
est globalement \'epimorphique, ce qui signifie que les
$$
X_i \longrightarrow X \, , \qquad i \in I \, ,
$$
forment une famille $J$-couvrante de l'objet $X$.

\medskip

\item Tout objet $Y$ de ${\mathcal C}$ tel que $U(Y) \ne \emptyset$ admet le crible vide comme crible $J_U^c$-couvrant.

\smallskip

Cela impose que pour tout $J_U^c$-faisceau $F$, l'ensemble $F(Y)$ est constitu\'e d'un unique \'el\'ement.

\smallskip

R\'eciproquement, soit $F$ un objet de $\widehat{\mathcal C}_J = {\mathcal E}$ tel que, pour tout objet $Y$ de ${\mathcal C}$ v\'erifiant $U(Y) \ne \emptyset$, l'ensemble $F(Y)$ soit constitu\'e d'un unique \'el\'ement.

\smallskip

Comme $F$ est un $J$-faisceau, on a pour toute famille $J$-couvrante
$$
X_i \longrightarrow X \, , \qquad i \in I \, ,
$$
compl\'et\'ee par une famille de morphismes
$$
X_i \longleftarrow X_{i,j,k} \longrightarrow X_j
$$
compatibles avec les $X_i \to X$ et tels que chaque famille induite
$$
\ell (X_{i,j,k}) \longrightarrow \ell (X_i) \times_{\ell (X)} \ell (X_j)
$$
soit globalement \'epimorphique, la formule
$$
F(X) = {\rm eg} \left( \prod_i F (X_i) \rightrightarrows \prod_{i,j,k} F(X_{i,j,k}) \right).
$$
Or, d'apr\`es l'hypoth\`ese faite sur $F$, cette formule reste vraie apr\`es \'elimination de tous les $X_i$ tels que
$$
U(X_i) \ne \emptyset
$$
et de tous les $X_{i,j,k}$ tels que
$$
U(X_{i,j,k}) \ne \emptyset \, .
$$

Cela signifie que $F$ est un faisceau pour la topologie $J_U^c$.

\medskip

\item Consid\'erons le plongement du sous-topos ouvert ${\mathcal E}_U$
$$
(i^* , i_*) : {\mathcal E}_U = {\mathcal E}/U \xhookrightarrow{ \ { \ } \ } {\mathcal E} \, .
$$
Sa composante d'image r\'eciproque est le foncteur
$$
\begin{matrix}
i^* : {\mathcal E} &\longrightarrow &{\mathcal E}/U \, , \hfill \\
\hfill E &\longmapsto &E \times U \, .
\end{matrix}
$$
Si un objet $F$ de $\widehat{\mathcal C}_J = {\mathcal E}$ est un $J_U^c$-faisceau, le morphisme canonique
$$
F \times U \longrightarrow U
$$
est un isomorphisme d'apr\`es (iii).

\smallskip

Donc tout $J_U^c$-faisceau a pour image par $i^*$ l'objet terminal de ${\mathcal E}/U$ et pour image par $i_* \circ i^*$ l'objet terminal de ${\mathcal E}$.

\smallskip

Cela ach\`eve de prouver le lemme. 
\end{listeisansmarge}
\end{demolem}

\bigskip

\noindent {\bf Fin de la d\'emonstration de la proposition \ref{propIV812}:}

\smallskip

Comme tout topos ${\mathcal E}$ se repr\'esente comme le topos des faisceaux sur un site, on est ramen\'e au corollaire suivant du lemme ci-dessus:

\begin{cor}\label{corIV814}

Soient ${\mathcal E} = \widehat{\mathcal C}_J$ le topos des faisceaux sur un site $({\mathcal C},J)$, $U$ un objet sous-terminal de ${\mathcal E}$, ${\mathcal E}_U = {\mathcal E}/U$ le sous-topos ouvert associ\'e et ${\mathcal E}_U^c = {\mathcal E} \backslash {\mathcal E}_U$ son sous-topos ferm\'e compl\'ementaire.

\smallskip

Alors le sous-topos ${\mathcal E}_U^c$ est d\'efini par la topologie $J_U^c$ de ${\mathcal C}$ introduite dans le lemme \ref{lemIV813}.
\end{cor}

\bigskip

\begin{demo}

On sait d'apr\`es le lemme \ref{lemIV813} que la topologie $J_U^c$ et la topologie $J_U$ qui d\'efinit le sous-topos ouvert ${\mathcal E}_U = {\mathcal E} / U$ ont pour intersection la topologie $J$ et qu'elles engendrent la topologie maximale de ${\mathcal C}$ (pour laquelle tout crible est couvrant).

\smallskip

La conclusion r\'esulte du corollaire \ref{corIV712} (iii).

\end{demo}

\bigskip

On d\'eduit de la proposition \ref{propIV810}, de la proposition \ref{propIV812} et du corollaire \ref{corIV814} que les sous-topos ferm\'es des topos poss\`edent les propri\'et\'es habituelles des ferm\'es des espaces topologiques:

\begin{cor}\label{corIV815}
\begin{listeimarge}
\item Les intersections de sous-topos ferm\'es sont des sous-topos ferm\'es.

\smallskip

Plus pr\'ecis\'ement, si les $U_i$ sont des objets sous-terminaux d'un topos ${\mathcal E}$ et $\underset{i \in I}{\bigvee} \, U_i = U$ est leur r\'eunion, on a
$$
\bigwedge_{i \in I} {\mathcal E}_{U_i}^c = {\mathcal E}_U^c \, .
$$

\item Les r\'eunions finies de sous-topos ferm\'es sont des sous-topos ferm\'es.

\smallskip

Plus pr\'ecis\'ement, si $U_1 , \cdots , U_n$ sont des objets sous-terminaux d'un topos ${\mathcal E}$ et $U_1 \wedge \cdots \wedge U_n = U$ est leur intersection, on a
$$
{\mathcal E}_{U_1}^c \vee \cdots \vee {\mathcal E}_{U_n}^c = {\mathcal E}_U^c \, .
$$

\item Les images r\'eciproques de sous-topos ferm\'es sont des sous-topos ferm\'es.

\smallskip

Plus pr\'ecis\'ement, si
$$
f = (f^* , f_*) : {\mathcal E}' \longrightarrow {\mathcal E}
$$
est un morphisme de topos et $U$ est un objet sous-terminal de ${\mathcal E}$, on a
$$
f^{-1} \, {\mathcal E}_U^c = {\mathcal E}_{f^* U}^c \, .
$$

\item Pour tout morphisme de topos
$$
f = (f^* , f_*) : {\mathcal E}' \longrightarrow {\mathcal E}
$$
et tout objet sous-terminal $U$ de ${\mathcal E}$, $f$ se factorise \`a travers le sous-topos ferm\'e
$$
{\mathcal E}_U^c \xhookrightarrow{ \ { \ } \ } {\mathcal E}
$$
si et seulement si $f^* U$ est un objet initial de ${\mathcal E}'$.
\end{listeimarge}
\end{cor}

\begin{remarks}
\begin{listeisansmarge}
\item Un sous-topos d'un topos ${\mathcal E}$ peut \^etre dit localement ferm\'e s'il s'\'ecrit comme l'intersection
$$
{\mathcal E}_U \wedge {\mathcal E}_{U'}^c
$$
d'un sous-topos ouvert ${\mathcal E}_U$ et d'un sous-topos ferm\'e ${\mathcal E}_{U'}^c$.

\medskip

\item L'intersection d'une famille finie de sous-topos localement ferm\'es
$$
{\mathcal E}_{U_i} \wedge {\mathcal E}_{U'_i}^c \, , \qquad 1 \leq i \leq n \, ,
$$
est le sous-topos localement ferm\'e
$$
{\mathcal E}_U \wedge {\mathcal E}_{U'}^c
$$
avec $U = U_1 \wedge \cdots \wedge U_n$ et $U' = U'_1 \vee \cdots \vee U'_n$.

\medskip

\item L'image r\'eciproque par un morphisme de topos
$$
f = (f^*,f_*) : {\mathcal E}' \longrightarrow {\mathcal E}
$$
d'un sous-topos localement ferm\'e
$$
{\mathcal E}_U \wedge {\mathcal E}_{U'}^c
$$
est le sous-topos localement ferm\'e
$$
{\mathcal E}'_{f^* U} \wedge {\mathcal E}'^c_{f^* U'} \, .
$$
De plus, le morphisme $f$ se factorise \`a travers le sous-topos
$$
{\mathcal E}_U \wedge {\mathcal E}_{U'}^c
$$
si et seulement si $f^* U$ est un objet terminal de ${\mathcal E}'$ et $f^* U'$ est un objet initial de ${\mathcal E}'$.
\end{listeisansmarge}
\end{remarks}

\bigskip

\begin{demo}
\begin{listeisansmarge}
\item[(i) et (ii)] Pour (i) et (ii), il suffit de traiter le cas o\`u $n=2$.

\smallskip

En formant l'intersection des deux formules
$$
{\mathcal E}_{U_1} \vee {\mathcal E}_{U_1}^c = {\mathcal E} \, , \qquad {\mathcal E}_{U_2} \vee {\mathcal E}_{U_2}^c = {\mathcal E}
$$
et d\'eveloppant, on obtient
$$
({\mathcal E}_{U_1} \wedge {\mathcal E}_{U_2}) \vee ({\mathcal E}_{U_1} \wedge {\mathcal E}_{U_2}^c) \vee ({\mathcal E}_{U_1}^c \wedge {\mathcal E}_{U_2}) \vee ({\mathcal E}_{U_1}^c \wedge {\mathcal E}_{U_2}^c) = {\mathcal E}
$$
qui implique
$$
({\mathcal E}_{U_1} \wedge {\mathcal E}_{U_2}) \vee ({\mathcal E}_{U_1}^c \vee {\mathcal E}_{U_2}^c) = {\mathcal E}
$$
ainsi que
$$
({\mathcal E}_{U_1} \vee {\mathcal E}_{U_2}) \vee ({\mathcal E}_{U_1}^c \wedge {\mathcal E}_{U_2}^c) = {\mathcal E} \, .
$$

D'autre part, en formant la r\'eunion des deux formules
$$
{\mathcal E}_{U_1} \wedge {\mathcal E}_{U_1}^c = {\mathcal E}_{\emptyset} \, , \qquad {\mathcal E}_{U_2} \wedge {\mathcal E}_{U_2}^c = {\mathcal E}_{\emptyset}
$$
et d\'eveloppant, on obtient
$$
({\mathcal E}_{U_1} \vee {\mathcal E}_{U_2}) \wedge ({\mathcal E}_{U_1} \vee {\mathcal E}_{U_2}^c) \wedge ({\mathcal E}_{U_1}^c \vee {\mathcal E}_{U_2}) \wedge ({\mathcal E}_{U_1}^c \vee {\mathcal E}_{U_2}^c) = {\mathcal E}_{\emptyset}
$$
qui implique
$$
({\mathcal E}_{U_1} \vee {\mathcal E}_{U_2}) \wedge ({\mathcal E}_{U_1}^c \wedge {\mathcal E}_{U_2}^c) = {\mathcal E}_{\emptyset}
$$
ainsi que
$$
({\mathcal E}_{U_1} \wedge {\mathcal E}_{U_2}) \wedge ({\mathcal E}_{U_1}^c \vee {\mathcal E}_{U_2}^c) = {\mathcal E}_{\emptyset} \, .
$$

On sait aussi d'apr\`es la proposition \ref{propIV810} que
$$
{\mathcal E}_{U_1} \vee {\mathcal E}_{U_2} = {\mathcal E}_{U_1 \vee \, U_2}
$$
et
$$
{\mathcal E}_{U_1} \wedge {\mathcal E}_{U_2} = {\mathcal E}_{U_1 \wedge \, U_2} \, .
$$

On conclut d'apr\`es le corollaire \ref{corIV712} (iii) que 
$$
{\mathcal E}_{U_1}^c \wedge {\mathcal E}_{U_2}^c = {\mathcal E} \backslash {\mathcal E}_{U_1 \vee \, U_2} = {\mathcal E}_{U_1 \vee \, U_2}^c
$$
et
$$
{\mathcal E}_{U_1}^c \vee {\mathcal E}_{U_2}^c = {\mathcal E} \backslash {\mathcal E}_{U_1 \wedge \, U_2} = {\mathcal E}_{U_1 \wedge \, U_2}^c \, .
$$

Cela prouve (i) et (ii).

\medskip

\item[(iii)] La formule
$$
{\mathcal E}_U \wedge {\mathcal E}_U^c = {\mathcal E}_{\emptyset}
$$
implique en passant aux images r\'eciproques
$$
{\mathcal E}'_{f^* U} \wedge f^{-1} {\mathcal E}_U^c = {\mathcal E}'_{\emptyset} \, .
$$
D'o\`u une in\'egalit\'e
$$
f^{-1} {\mathcal E}_U^c \leq {\mathcal E}'^c_{f^* U} \, .
$$

Pour montrer qu'il y a \'egalit\'e, on peut supposer que ${\mathcal E} = \widehat{\mathcal C}_J$ et ${\mathcal E}' = \widehat{\mathcal D}_K$ sont les topos des faisceaux sur deux sites $({\mathcal C} , J)$ et $({\mathcal D},K)$ tels que ${\mathcal C}$ et ${\mathcal D}$ soient deux sous-cat\'egories pleines de ${\mathcal E}$ et ${\mathcal E}'$, qu'elles soient stables par limites finies et contiennent les objets $U$ et $f^*U$, et que le foncteur
$$
f^* : {\mathcal E} \longrightarrow {\mathcal E}'
$$
se restreigne en un foncteur
$$
\rho : {\mathcal C} \longrightarrow {\mathcal D}
$$
qui donc est un morphisme de sites et d\'efinit $f : {\mathcal E}' \to {\mathcal E}$.

\smallskip

La topologie $J_U^c$ [resp. $K_{f^*U}^c$] de ${\mathcal C}$ [resp. ${\mathcal D}$] qui d\'efinit le sous-topos ${\mathcal E}_U^c$ [resp. ${\mathcal E}'^c_{f^*U}$] de ${\mathcal E} = \widehat{\mathcal C}_J$ [resp. ${\mathcal E}' = \widehat{\mathcal D}_K$] est la plus petite topologie qui contient $J$ [resp. $K$] et pour laquelle l'objet $U$ de ${\mathcal C}$ [resp. l'objet $f^* U$ de ${\mathcal D}$] admette pour crible couvrant le crible vide.

\smallskip

Comme l'objet $f^* U$ de ${\mathcal D}$ est image de l'objet $U$ de ${\mathcal C}$ par le foncteur $\rho = f^*$, on conclut d'apr\`es la proposition \ref{propIV713} (ii).

\medskip

\item[(iv)] Un morphisme de topos
$$
f = (f^* , f_*) : {\mathcal E}' \longrightarrow {\mathcal E}
$$
se factorise \`a travers le sous-topos ferm\'e
$$
{\mathcal E}_U^c \xhookrightarrow{ \ { \ } \ } {\mathcal E}
$$
compl\'ementaire du sous-topos ouvert ${\mathcal E}_U$ associ\'e \`a un objet sous-terminal $U$ de ${\mathcal E}$, si et seulement si on a
$$
f^{-1} ({\mathcal E}_U^c) = {\mathcal E}' \, .
$$
Or on a d'apr\`es (iii)
$$
f^{-1} ({\mathcal E}_U^c) = {\mathcal E}'^c_{f^* U} = {\mathcal E}' \backslash {\mathcal E}'_{f^*U} \, .
$$
Ainsi, il y a factorisation si et seulement si
$$
{\mathcal E}'_{f^*U} = {\mathcal E}'_{\emptyset}
$$
c'est-\`a-dire si et seulement si $f^* U$ est un objet initial de ${\mathcal E}$.

\smallskip

Cela ach\`eve la d\'emonstration du corollaire. 

\end{listeisansmarge}
\end{demo}

\subsection{Le faisceau des sous-topos}\label{subsec484}

\medskip

Les images r\'eciproques de sous-topos par les morphismes de localisation d'un topos se calculent de la mani\`ere suivante:

\begin{lem}\label{lemIV816}

Soient
$$
f : E_2 \longrightarrow E_1
$$
un morphisme d'un topos ${\mathcal E}$,
$$
(f^*,f_*) : {\mathcal E}/E_2 \longrightarrow {\mathcal E}/E_1
$$
le morphisme de topos localis\'es associ\'e,
$$
j_1 : {\mathcal E}_1 \xhookrightarrow{ \ { \ } \ } {\mathcal E}/E_1
$$
un sous-topos et
$$
j_2 : {\mathcal E}_2 \xhookrightarrow{ \ { \ } \ } {\mathcal E}/E_2
$$
son image r\'eciproque par le morphisme $(f^*,f_*)$.

\smallskip

Supposons que ${\mathcal E}$ est le topos des faisceaux sur un site $({\mathcal C},J)$, de sorte que ${\mathcal C}/E_1$ et ${\mathcal C}/E_2$ s'identifient aux topos des faisceaux sur les cat\'egories ${\mathcal C}/E_1$ et ${\mathcal C}/E_2$ pour les topologies induites par $J$.

\smallskip

Soit $J_1$ la topologie de ${\mathcal C}/E_1$ qui d\'efinit le sous-topos ${\mathcal E}_1$ de ${\mathcal E}/E_1$.

\smallskip

Alors la topologie $J_2$ de ${\mathcal C}/E_2$ qui d\'efinit le sous-topos ${\mathcal E}_2$ de ${\mathcal E}/E_1$ est la topologie induite par $J_1$.
\end{lem}

\bigskip

\begin{demo}

La topologie $J'_2$ de ${\mathcal C}/E_2$ induite par la topologie $J_1$ de ${\mathcal C}/E_1$ d\'efinit un sous-topos
$$
{\mathcal E}'_2 \xhookrightarrow{ \ { \ } \ } {\mathcal E}/E_2
$$
qui s'inscrit dans un carr\'e commutatif \`a isomorphisme pr\`es
$$
\xymatrix{
{\mathcal E}'_{2} \ar[d] \ar@{^{(}->}[r] &{\mathcal E}/E_2 \ar[d] \\
{\mathcal E}_1 \ar@{^{(}->}[r] &{\mathcal E}/E_1
}
$$
donc ${\mathcal E}'_2$ est un sous-topos de ${\mathcal E}_2$.

\smallskip

Pour la r\'eciproque, consid\'erons un morphisme de ${\mathcal E}/E_2$
$$
F' \longrightarrow F
$$
dont l'image par le foncteur de composition avec $f$
$$
f_! : {\mathcal E}/E_2 \longrightarrow {\mathcal E}/E_1
$$
est transform\'ee par
$$
j_1^* : {\mathcal E}/E_1 \longrightarrow {\mathcal E}_1
$$
en un \'epimorphisme.

\smallskip

Alors le morphisme de ${\mathcal E}/E_2$
$$
F' \times_{E_1} E_2 \longrightarrow F \times_{E_1} E_2
$$
est transform\'e par
$$
j_2^* : {\mathcal E}/E_2 \longrightarrow {\mathcal E}_2
$$
en un \'epimorphisme.

\smallskip

Or, le foncteur de changement de base par la diagonale
$$
\Delta : E_2 \longrightarrow E_2 \times_{E_1} E_2
$$
dans ${\mathcal E}/E_2$ transforme le morphisme
$$
F' \times_{E_1} E_2 \longrightarrow F \times_{E_1} E_2
$$
en le morphisme
$$
F' \longrightarrow F \, .
$$

On en d\'eduit que le transform\'e de ce morphisme de ${\mathcal E}/E_2$
$$
F' \longrightarrow F
$$
par le foncteur
$$
j_2^* : {\mathcal E} / E_2 \longrightarrow {\mathcal E}_2
$$
est un \'epimorphisme.

\smallskip

Cela montre comme voulu que le morphisme de plongement
$$
{\mathcal E}'_2 \xhookrightarrow{ \ { \ } \ } {\mathcal E}_2
$$
est une \'equivalence de topos. \end{demo}

\bigskip

Nous pouvons maintenant montrer que tout topos ${\mathcal E}$ poss\`ede un objet qui classifie les sous-objets des topos localis\'es en les objets de ${\mathcal E}$:

\begin{thm}\label{thmIV817}

Soit ${\mathcal E}$ un topos.

\smallskip

Associons \`a tout objet $E$ de ${\mathcal E}$ l'ensemble $\Omega_T(E)$ des sous-topos du topos localis\'e ${\mathcal E}/E$.

\smallskip

Et associons \`a tout morphisme de ${\mathcal E}$
$$
f : E_2 \longrightarrow E_1
$$
l'application en sens inverse
$$
f^{-1} = \Omega_T(f) : \Omega_T (E_1) \longrightarrow \Omega_T(E_2)
$$
qui transforme tout sous-topos
$$
{\mathcal E}'_1 \xhookrightarrow{ \ { \ } \ } {\mathcal E}/E_1
$$
en son image r\'eciproque
$$
f^{-1} {\mathcal E}'_1 = {\mathcal E}'_2 \xhookrightarrow{ \ { \ } \ } {\mathcal E}/E_2
$$
par le morphisme de topos
$$
(f^*,f_*) : {\mathcal E} / E_2 \longrightarrow {\mathcal E}/E_1 \, .
$$

Alors le foncteur contravariant
$$
\Omega_T : {\mathcal E}^{\rm op} \longrightarrow {\rm Ens}
$$
est repr\'esentable par un objet de ${\mathcal E}$ que l'on peut encore noter $\Omega_T$.
\end{thm}

\bigskip

\begin{remarks}
\begin{listeisansmarge}
\item L'objet $\Omega_T$ de ${\mathcal E}$ contient comme sous-objet
$$
\Omega \xhookrightarrow{ \ { \ } \ } \Omega_T
$$
le ``classificateur des sous-objets'' $\Omega$ de ${\mathcal E}$.

\smallskip

Le morphisme de plongement
$$
\Omega \xhookrightarrow{ \ { \ } \ } \Omega_T
$$
consiste \`a associer \`a tout sous-objet $U$ d'un objet $E$ de ${\mathcal E}$ le sous-topos ouvert ${\mathcal E}_U = {\mathcal E}/U$ de ${\mathcal E}/E$.

\smallskip

C'est bien un morphisme de foncteurs car, pour tout morphisme de ${\mathcal E}$
$$
f : E_2 \longrightarrow E_1
$$
et tout sous-objet $U_1$ de $E_1$, l'image r\'eciproque par
$$
(f^*,f_*) : {\mathcal E}/E_2 \longrightarrow {\mathcal E}/E_1
$$
du sous-topos ouvert ${\mathcal E}/U_1$ de ${\mathcal E}/E_1$ est le sous-topos ouvert de ${\mathcal E}/E_2$
$$
{\mathcal E} / U_1 \times_{E_1} E_2
$$
d\'efini par le sous-objet $U_1 \times_{E_1} E_2 = f^* U_1$ de $E_2$.

\medskip

\item Il existe un second morphisme de plongement
$$
\Omega \xhookrightarrow{ \ { \ } \ } \Omega_T
$$
qui consiste \`a associer \`a tout sous-objet $U$ d'un objet $E$ de ${\mathcal E}$ le sous-topos ferm\'e ${\mathcal E}_U^c = ({\mathcal E}/E) \backslash {\mathcal E}_U$ compl\'ementaire du sous-topos ouvert ${\mathcal E}_U = {\mathcal E}/U$ de ${\mathcal E}/E$ d\'efini par $U$.

\smallskip

C'est bien un morphisme de foncteurs d'apr\`es le corollaire \ref{corIV815} (iii).
\end{listeisansmarge}
\end{remarks}

\bigskip

\begin{demo}

D'apr\`es le th\'eor\`eme \ref{thmIII73} (i), il suffit de montrer que le foncteur $\Omega$ transforme les colimites en limites.

\smallskip

Consid\'erons donc un carquois $D$ et un $D$-diagramme $E_{\bullet}$ de ${\mathcal E}$ de colimite
$$
E = \varinjlim_D E_{\bullet} \, .
$$

Il s'agit de prouver que l'application
$$
\Omega_T (E) \longrightarrow \, \varprojlim_{D} \, \Omega_T (E_{\bullet})
$$
est bijective.

\smallskip

Quitte \`a remplacer le topos ${\mathcal E}$ par ${\mathcal E}/E$, on peut supposer que $E$ est l'objet terminal de ${\mathcal E}$.

\smallskip

Consid\'erons une famille compatible de sous-topos
$$
i_{\rm d} : {\mathcal E}_d \xhookrightarrow{ \ { \ } \ } {\mathcal E}/E_d
$$
index\'es par les objets $d$ de $D$. Il s'agit de prouver qu'elle est induite par un unique sous-topos
$$
{\mathcal E}' \xhookrightarrow{ \ { \ } \ } {\mathcal E} \, .
$$

Montrons d'abord l'unicit\'e d'un tel ${\mathcal E}' \hookrightarrow {\mathcal E}$ s'il existe.

\smallskip

Notant $i^* : {\mathcal E} \to {\mathcal E}'$ la composante d'image r\'eciproque du plongement de topos ${\mathcal E}' \xhookrightarrow{i} {\mathcal E}$, on note d'apr\`es le lemme \ref{lemIV816} que chaque ${\mathcal E}_d \hookrightarrow {\mathcal E}/E_d$ doit s'identifier \`a
$$
{\mathcal E}'/i^* E_d \xhookrightarrow{ \ { \ } \ } {\mathcal E}/E_d \, .
$$

Il en r\'esulte qu'un monomorphisme de ${\mathcal E}$
$$
S \xhookrightarrow{ \ { \ } \ } Y
$$
est transform\'e par $i^* : {\mathcal E} \to {\mathcal E}'$ en un isomorphisme de ${\mathcal E}'$ si et seulement si, pour tout objet $d$ de $D$, le monomorphisme induit de ${\mathcal E}/E_d$
$$
S \times E_d \xhookrightarrow{ \ { \ } \ } Y \times E_d
$$
est transform\'e par le foncteur d'image r\'eciproque
$$
{\mathcal E}/E_d \longrightarrow {\mathcal E}'_d
$$
en un isomorphisme.

\smallskip

Cela r\'esulte en effet de ce que le foncteur $i^*$ respecte les colimites et les limites finies, et que les foncteurs $S \times \bullet$ et $Y \times \bullet$ respectent les colimites.

\smallskip

Ainsi, la connaissance des foncteurs ${\mathcal E}/E_d \to {\mathcal E}'_d$ d\'etermine ceux des monomorphismes $S \hookrightarrow Y$ de ${\mathcal E}$ que $i^* : {\mathcal E} \to {\mathcal E}'$ transforme en isomorphismes.

\smallskip

Cela montre que le sous-topos ${\mathcal E}' \hookrightarrow {\mathcal E}$ est enti\`erement d\'etermin\'e par les ${\mathcal E}'_d \hookrightarrow {\mathcal E}/E_d$ s'il existe.

\smallskip

Pour l'existence, on peut supposer que ${\mathcal E} = \widehat{\mathcal C}_J$ est le topos des faisceaux sur un site $({\mathcal C},J)$ muni du foncteur canonique
$$
\ell : {\mathcal C} \xhookrightarrow{ \ y \ } \widehat{\mathcal C} \xrightarrow{ \ j^* \ } \widehat{\mathcal C}_J = {\mathcal E} \, .
$$

Pour tout objet $X$ de ${\mathcal C}$, notons $J'(X)$ l'ensemble des cribles de $X$
$$
S \xhookrightarrow{ \ { \ } \ } y(X)
$$
tels que, pour tout objet $d$ de $D$, le foncteur d'image r\'eciproque
$$
i_d^* : {\mathcal E} / E_d \longrightarrow {\mathcal E}'_d
$$
transforme le monomorphisme de ${\mathcal E}/E_d$
$$
j^* S \times E_d \xhookrightarrow{ \ { \ } \ } \ell (X) \times E_d
$$
en un isomorphisme de ${\mathcal E}'_d$.

\smallskip

Alors $J'$ est une topologie de ${\mathcal C}$ qui contient $J$.

\smallskip

En effet, $J'(X)$ contient $J(X)$ pour tout objet $X$ de ${\mathcal C}$ car un crible $S \hookrightarrow J(X)$ de $X$ est $J$-couvrant si et seulement si $j^*$ le transforme en un isomorphisme $j^* S \xrightarrow{ \ \sim \ } \ell (X)$ de ${\mathcal E} = \widehat{\mathcal C}_J$. A fortiori, $J'$ satisfait l'axiome de maximalit\'e.

\smallskip

Elle satisfait l'axiome de stabilit\'e car tous les foncteurs $i^*_d : {\mathcal E} / E_d \to {\mathcal E}'_d$ respectent les produits finis.

\smallskip

Enfin, $J'$ satisfait l'axiome de transitivit\'e car tous les foncteurs $i_d^* : {\mathcal E}/E_d \to {\mathcal E}'_d$ respectent les colimites arbitraires et les limites finies.

\smallskip

La topologie $J' \supseteq J$ de ${\mathcal C}$ d\'efinit un sous-topos
$$
i = (i^*,i_*) : {\mathcal E}' \xhookrightarrow{ \ { \ } \ } {\mathcal E}
$$
tel que, pour tout objet $d$ de $D$, le plongement
$$
{\mathcal E}_d \xhookrightarrow{ \ { \ } \ } {\mathcal E}/E_d
$$
se factorise en
$$
{\mathcal E}_d \xhookrightarrow{ \ { \ } \ } {\mathcal E}' / i^* E_d \xhookrightarrow{ \ { \ } \ } {\mathcal E}/E_d \, .
$$

Il reste \`a montrer que, pour tout objet $d$ de $D$, le plongement de topos
$$
{\mathcal E}'_d \xhookrightarrow{ \ { \ } \ } {\mathcal E}' / i^* E_d
$$
est une \'equivalence.

\smallskip

Consid\'erons pour cela un objet $(X,\ell(X) \to E_d)$ de ${\mathcal C}/E_d$ et un crible $S \hookrightarrow y(X)$ de $X$ tel que
$$
i_d^* : {\mathcal E} / E_d \longrightarrow {\mathcal E}_d
$$
transforme $j^* S \hookrightarrow \ell (X)$ en isomorphisme.

\smallskip

Il s'agit de prouver que, pour tout objet $d'$ de $D$, le foncteur
$$
i^*_{d'} : {\mathcal E} / E_{d'} \longrightarrow {\mathcal E}_{d'}
$$
transforme $j^* S \times E_{d'} \hookrightarrow \ell (X) \times E_{d'}$ en isomorphisme. 

\smallskip

On le sait d\'ej\`a si $d=d'$.

\smallskip

Puis supposons que $d'$ soit reli\'e \`a $d$ par une cha{\^\i}ne
$$
d = d_0 \longleftrightarrow d_1 \longleftrightarrow \cdots \longleftrightarrow d_n = d'
$$
constitu\'ee de morphismes $d_{i-1} \to d_i$ ou $d_i \to d_{i-1}$ de $D$ not\'es $d_{i-1} \leftrightarrow d_i$.

\smallskip

Il r\'esulte du lemme \ref{lemIV816} que, pour tout indice $i$, le morphisme de ${\mathcal E}/E_{d_{i-1}}$
$$
j^* S \times E_{d_{i-1}} \xhookrightarrow{ \ { \ } \ } \ell (X) \times E_{d_{i-1}}
$$
est transform\'e par $i^*_{d_{i-1}} : {\mathcal E} / E_{d_{i-1}} \to {\mathcal E}_{d_{i-1}}$ en isomorphisme si et seulement si le morphisme de ${\mathcal E}/E_{d_i}$
$$
j^* S \times E_{d_i} \xhookrightarrow{ \ { \ } \ } \ell (X) \times E_{d_i}
$$
est transform\'e par $i_{d_i}^* : {\mathcal E} / E_{d_i} \to {\mathcal E}_{d_i}$ en isomorphisme.

\smallskip

Donc le morphisme de ${\mathcal E}/E_{d'}$
$$
j^* S \times E_{d'} \xhookrightarrow{ \ { \ } \ } \ell (X) \times E_{d'}
$$
est transform\'e par $i_{d'}^* : {\mathcal E} / E_{d'} \to {\mathcal E}_{d'}$ en isomorphisme.

\smallskip

Reste enfin le cas o\`u $d$ et $d'$ ne sont pas reli\'es par une telle suite.

\smallskip

Alors les deux produits
$$
j^* S \times E_{d'} \qquad \mbox{et} \qquad \ell (X) \times E_{d'}
$$
s'identifient \`a l'objet initial de ${\mathcal E} / E_{d'}$ et le morphisme
$$
j^* S \times E_{d'} \longrightarrow \ell (X) \times E_{d'}
$$
est automatiquement un isomorphisme, donc aussi son image par le foncterur $i_{d'}^* : {\mathcal E} / E_{d'} \to {\mathcal E}_{d'}$.

\smallskip

Cela termine la d\'emonstration du th\'eor\`eme. 

\end{demo}

\subsection{Le faisceau des transformations naturelles entre deux foncteurs}\label{subsec485}

\medskip

Toute paire de foncteurs \`a valeurs dans un topos d\'efinit un faisceau des transformations naturelles de l'un dans l'autre:

\begin{prop}\label{propIV818}

Soient ${\mathcal E}$ un topos et ${\mathcal C}$ une cat\'egorie essentiellement petite.

\smallskip

Alors:

\begin{listeimarge}

\item Pour toute paire de foncteurs
$$
\rho_1 , \rho_2 : {\mathcal C} \longrightarrow {\mathcal E} \, ,
$$
le pr\'efaisceau qui associe \`a tout objet $E$ de ${\mathcal E}$ l'ensemble des transformations naturelles
$$
(\bullet \times E) \circ \rho_1 \longrightarrow (\bullet \times E) \circ \rho_2
$$
entre les compos\'es de $\rho_1$ et $\rho_2$ avec le foncteur
$$
\begin{matrix}
\bullet \times E : {\mathcal E} &\longrightarrow &{\mathcal E}/E \, , \\
\hfill F &\longmapsto &F \times E
\end{matrix}
$$
est repr\'esentable par un objet de ${\mathcal E}$ que l'on peut noter
$$
{\mathcal H}om_{\mathcal E} (\rho_1,\rho_2) \qquad \mbox{ou} \qquad {\mathcal H}om (\rho_1,\rho_2) \, .
$$

\item La composition des transformations naturelles d\'efinit pour tout triplet de foncteurs
$$
\rho_1 , \rho_2 , \rho_3 : {\mathcal C} \longrightarrow {\mathcal E}
$$
un morphisme de ${\mathcal E}$
$$
{\mathcal H}om (\rho_1,\rho_2) \times {\mathcal H}om (\rho_2,\rho_3) \longrightarrow {\mathcal H}om (\rho_1,\rho_3) \, ,
$$
$$
(\alpha , \beta) \longmapsto \beta \circ \alpha \, .
$$

\item Ces morphismes d\'efinissent une loi de multiplication associative au sens que pour tout quadruplet de foncteurs
$$
\rho_1 , \rho_2 , \rho_3 , \rho_4 : {\mathcal C} \longrightarrow {\mathcal E}
$$
le carr\'e de ${\mathcal E}$
$$
\xymatrix{
{\mathcal H}om (\rho_1,\rho_2) \times {\mathcal H}om (\rho_2,\rho_3) \times {\mathcal H}om (\rho_3,\rho_4) \ar[d] \ar[r] &{\mathcal H}om (\rho_1,\rho_3) \times {\mathcal H}om (\rho_3,\rho_4) \ar[d] \\
{\mathcal H}om (\rho_1,\rho_2) \times {\mathcal H}om (\rho_2,\rho_4) \ar[r] &{\mathcal H}om (\rho_1,\rho_4)
}
$$
est commutatif.

\medskip

\item Pour tout foncteur
$$
\rho : {\mathcal C} \longrightarrow {\mathcal E} \, ,
$$
associer \`a tout objet $E$ de ${\mathcal E}$ l'identit\'e du foncteur $(\bullet \times E) \circ \rho$ d\'efinit un morphisme de ${\mathcal E}$
$$
{\rm id}_{\rho} : 1_{\mathcal E} \longrightarrow {\mathcal H}om (\rho , \rho) = {\mathcal E}nd (\rho) \, .
$$
Les ${\rm id}_{\rho}$ sont neutres pour la composition au sens que pour tous foncteurs
$$
\rho_1 , \rho_2 : {\mathcal C} \longrightarrow {\mathcal E} \, ,
$$
les morphismes compos\'es
$$
1_{\mathcal E} \times {\mathcal H}om (\rho_1,\rho_2) \longrightarrow {\mathcal H}om (\rho_1,\rho_1) \times {\mathcal H}om (\rho_1,\rho_2) \longrightarrow {\mathcal H}om (\rho_1,\rho_2) \, ,
$$
$$
{\mathcal H}om (\rho_1,\rho_2) \times 1_{\mathcal E} \longrightarrow {\mathcal H}om (\rho_1,\rho_2) \times {\mathcal H}om (\rho_2,\rho_2) \longrightarrow {\mathcal H}om (\rho_1,\rho_2) 
$$
sont les isomorphismes d'identit\'e.

\medskip

\item En particulier, pour tout foncteur
$$
\rho : {\mathcal C} \longrightarrow {\mathcal E} \, ,
$$
${\mathcal H}om (\rho , \rho) = {\mathcal E}nd (\rho)$ est un mono{\"\i}de interne de ${\mathcal E}$, et son groupe interne associ\'e
$$
{\mathcal A}ut (\rho) = 1_{\mathcal E} \times_{{\mathcal E}nd(\rho)} [{\mathcal E}nd (\rho) \times {\mathcal E}nd (\rho)] \times_{{\mathcal E}nd (\rho)} 1_{\mathcal E}
$$
repr\'esente le pr\'efaisceau qui associe \`a tout objet $E$ de ${\mathcal E}$ le groupe des automorphismes du foncteur $(\bullet \times E) \circ \rho$.
\end{listeimarge}
\end{prop}

\begin{remark}

En particulier, si ${\mathcal C}$ est munie d'une topologie $J$, on peut se restreindre aux foncteurs
$$
\rho : {\mathcal C} \longrightarrow {\mathcal E} 
$$
qui sont plats et $J$-continus c'est-\`a-dire repr\'esentent des morphismes de topos
$$
{\mathcal E} \longrightarrow \widehat{\mathcal C}_J \, .
$$

Cela signifie que les paires de morphismes de topos
$$
f_1 , f_2 : {\mathcal E} \longrightarrow {\mathcal E}'
$$
d\'efinissent des objets de ${\mathcal E}$
$$
{\mathcal H}om_{\mathcal E} (f_1,f_2) = {\mathcal H}om (f_1,f_2)
$$
tels que, pour tout objet $E$ de ${\mathcal E}$,
$$
{\mathcal H}om (f_1,f_2)(E)
$$
s'identifie \`a l'ensemble des transformations entre les deux morphismes de topos compos\'es
$$
{\mathcal E}/E \longrightarrow {\mathcal E} \raisebox{.7ex}{\xymatrix{\dar[r]^-{^{^{\mbox{\scriptsize$f_1$}}}}_-{f_2} & \ {\mathcal E}'}} \, .
$$

Ces objets de ${\mathcal E}$ sont reli\'es par des lois de compositions associatives
$$
{\mathcal H}om (f_1,f_2) \times {\mathcal H}om (f_2,f_3) \longrightarrow {\mathcal H}om (f_1,f_3)
$$
relativement auxquelles les morphismes naturels
$$
1_{\mathcal E} \longrightarrow {\mathcal H}om (f,f) = {\mathcal E}nd (f)
$$
sont neutres.

\smallskip

En particulier, chaque
$$
{\mathcal E}nd (f) = {\mathcal H}om (f,f) 
$$
est un mono{\"\i}de interne de ${\mathcal E}$, et son groupe interne associ\'e
$$
{\mathcal A}ut (f) = 1_{\mathcal E} \times_{{\mathcal E}nd (f)} [{\mathcal E}nd (f) \times {\mathcal E}nd (f)] \times_{{\mathcal E}nd (f)} 1_{\mathcal E}
$$
repr\'esente le pr\'efaisceau qui associe \`a tout objet $E$ de ${\mathcal E}$ le groupe des automorphismes du morphisme de topos compos\'e
$$
{\mathcal E}/E \longrightarrow {\mathcal E} \xrightarrow{ \ f \ } {\mathcal E}' \, .
$$
\end{remark}

\bigskip 

\begin{demo}
\begin{listeisansmarge}
\item[(i)] D'apr\`es le th\'eor\`eme \ref{thmIII73} (i), il suffit de montrer que le foncteur
$$
\begin{matrix}
{\mathcal E}^{\rm op} &\longrightarrow &{\rm Ens} \, , \hfill \\
\hfill E &\longmapsto &{\rm Hom} ((\bullet \times E) \circ \rho_1 , (\bullet \times E) \circ \rho_2)
\end{matrix}
$$
transforme les colimites en limites.

\smallskip

Consid\'erons donc un carquois $D$ et un $D$-diagramme $E_{\bullet}$ dans ${\mathcal E}$ de colimite
$$
E = \varinjlim_D \ E_{\bullet} \, .
$$
Un \'el\'ement de
$$
{\rm Hom} ((\bullet \times E) \circ \rho_1 , (\bullet \times E) \circ \rho_2)
$$
est une famille de morphismes de ${\mathcal E}/E$
$$
\rho_1 (X) \times E \longrightarrow \rho_2(X) \times E
$$
index\'es par les objets $X$ de ${\mathcal C}$ et tels que, pour tout morphisme $X \to Y$ de ${\mathcal C}$, le carr\'e
$$
\xymatrix{
\rho_1(X) \times E \ar[d] \ar[r] &\rho_2 (X) \times E \ar[d] \\
\rho_1(Y) \times E \ar[r] &\rho_2(Y) \times E
}
$$
soit commutatif.

\smallskip

C'est \'equivalent \`a se donner une famille de morphismes de ${\mathcal E}$
$$
\rho_1 (X) \times E \longrightarrow \rho_2 (X)
$$
index\'es par les objets $X$ de ${\mathcal C}$ et tels que, pour tout morphisme $X \to Y$ de ${\mathcal C}$, le carr\'e
$$
\xymatrix{
\rho_1(X) \times E \ar[d] \ar[r] &\rho_2 (X)  \ar[d] \\
\rho_2(Y) \times E \ar[r] &\rho_2(Y) 
}
$$
soit commutatif.

\smallskip

La conclusion r\'esulte de ce que, dans le topos ${\mathcal E}$, les foncteurs de produits $\rho_1 (X) \times \bullet$ respectent les colimites.

\medskip

\item[(ii)-(v)] (ii), (iii), (iv) et (v) r\'esultent, compte tenu de (i), de ce que les foncteurs de ${\mathcal C}$ dans chaque ${\mathcal E}/E$ forment une cat\'egorie localement petite $[{\mathcal C} , {\mathcal E}/E]$ et de ce que tout morphisme
$$
f : E_2 \longrightarrow E_1
$$
d\'efinit un foncteur
$$
f^* : {\mathcal E} / E_1 \longrightarrow {\mathcal E}/E_2
$$
puis un foncteur de composition avec $f^*$
$$
[{\mathcal C} , {\mathcal E}/E_1] \longrightarrow [{\mathcal C} , {\mathcal E}/E_2] \, .
$$
\end{listeisansmarge}
\end{demo}


%% file: Chapitre5_num.tex







\chapter{Th\'eories, mod\`eles, cat\'egories syntactiques et topos classifiants}\label{chap5}

Dans ce chapitre nous présentons les bases de la logique catégorique du prémier ordre. Les résultats exposés dans les sections 1 à 7 sont essentiellement tirés de \cite{CategoricalLogic}, tandis que les sections 8 et 9 s'appuient principalement sur \cite{TST}. 

\section{El\'ements de langage et leurs interpr\'etations}\label{sec51}

\subsection{La notion de langage du premier ordre}\label{subsec511}

Commen\c cons par introduire la notion g\'en\'erale de langage d'une th\'eorie du premier ordre:

\begin{defn}\label{defV11}

Un ``langage du premier ordre'', ou ``signature'', consiste en:
\begin{enumerate}
\item[$\bullet$] une famille de ``sortes'' (c'est-\`a-dire de ``noms d'objets'') $A$,
\item[$\bullet$] une famille de ``symboles de fonctions'' (c'est-\`a-dire de ``noms de morphismes'') \'ecrits formellement
$$
f : A_1 \cdots A_n \longrightarrow B
$$
o\`u $A_1 , \cdots , A_n$ et $B$ sont des sortes,
\item[$\bullet$] une famille de ``symboles de relations'' (c'est-\`a-dire de ``noms de sous-objets'') \'ecrits formellement
$$
\xymatrix{
R \ \ \ar@{>->}[r] &A_1 \cdots A_n
}
$$
o\`u $A_1 , \cdots , A_n$ sont des sortes.
\end{enumerate}
\end{defn}

\bigskip

\begin{remarksqed}
\begin{listeisansmarge}
\item Si $n=0$, un symbole de fonction
$$
f : \ \longrightarrow B
$$
est appel\'e un ``symbole de constante de sorte $B$'', et un symbole de relation
$$
\xymatrix{
R \ \ \ar@{>->}[r] &{ \ }
}
$$
est qualifi\'e de ``propositionnel''.

\medskip

\item Si $n = 1$ [resp. $n=2$, resp. $n=3$], un symbole de fonction
$$
f : A_1 \cdots A_n \longrightarrow B
$$
ou un symbole de relation
$$
\xymatrix{
R \ \ \ar@{>->}[r] &A_1 \cdots A_n
}
$$
est dit ``unaire'' [resp. ``binaire'', resp. ``ternaire'']. 

\end{listeisansmarge}
\end{remarksqed}

\bigskip

Il est facile de d\'efinir le langage des th\'eories alg\'ebriques les plus communes:

\begin{defn}\label{defV12}
\begin{listeimarge}
\item Le langage de la th\'eorie des mono{\"\i}des [resp. des groupes] consiste en

\medskip

$\left\{\begin{matrix}
\bullet &\mbox{une sorte $G$ (le mono{\"\i}de [resp. le groupe]),} \hfill \\
{ \ } \\
\bullet &\mbox{un symbole de fonction binaire} \hfill \\
{ \ } \\
&\qquad \cdot : GG \longrightarrow G \qquad \mbox{(la multiplication)}, \hfill \\
{ \ } \\
&\mbox{un symbole de constante} \hfill \\
{ \ } \\
&\qquad 1 : \ \longrightarrow G \qquad \mbox{(la constante unit\'e)}, \hfill \\
{ \ } \\
&\mbox{[resp. et un symbole de fonction unaire} \hfill \\
{ \ } \\
&\qquad (\bullet)^{-1} : G \longrightarrow G \qquad \mbox{(le passage \`a l'inverse) ]}. \hfill
\end{matrix} \right.$

\medskip

\item Le langage de la th\'eorie des actions de mono{\"\i}des [resp. de groupes] est constitu\'e du langage de la th\'eorie des mono{\"\i}des $(G,\cdot , 1)$ [resp. des groupes $(G,\cdot , 1 , (\bullet)^{-1})$] compl\'et\'e par

\medskip

$\left\{\begin{matrix}
\bullet &\mbox{une sorte $X$ (l'objet de l'action),} \hfill \\
{ \ } \\
\bullet &\mbox{un symbole de fonction binaire} \hfill \\
{ \ } \\
&\qquad \cdot : GX \longrightarrow X \qquad \mbox{(l'action)}. \hfill
\end{matrix} \right.$

\medskip

\item Le langage de la th\'eorie des anneaux consiste en

\medskip

$\left\{\begin{matrix}
\bullet &\mbox{une sorte $A$ (l'anneau),} \hfill \\
{ \ } \\
\bullet &\mbox{deux symboles de fonctions binaires} \hfill \\
{ \ } \\
&\qquad \begin{matrix}
+ &: &AA \longrightarrow A &\mbox{(l'addition)}, \hfill \\
\cdot &: &AA \longrightarrow A &\mbox{(la multiplication)},
\end{matrix} \hfill \\
{ \ } \\
&\mbox{deux symboles de constantes} \hfill \\
{ \ } \\
&\qquad \begin{matrix}
0 : \ \longrightarrow A &\mbox{(l'unit\'e de l'addition)}, \hfill \\
1 : \ \longrightarrow A &\mbox{(l'unit\'e de la multiplication)},
\end{matrix} \hfill \\
{ \ } \\
&\mbox{et un symbole de fonction unaire} \hfill \\
{ \ } \\
&\qquad -(\bullet) : A \longrightarrow A \qquad \mbox{(le passage \`a l'oppos\'e)}. \hfill
\end{matrix} \right.$

\medskip

\item Le langage de la th\'eorie des modules sur un anneau est constitu\'e du langage de la th\'eorie des anneaux $(A , + , \cdot , 0,1,-(\bullet))$ compl\'et\'e par

\medskip

$\left\{\begin{matrix}
\bullet &\mbox{une sorte $M$ (le module),} \hfill \\
{ \ } \\
\bullet &\mbox{deux symboles de fonctions binaires} \hfill \\
{ \ } \\
&\qquad \begin{matrix}
+ &: &M\!M \longrightarrow M &\mbox{(l'addition interne au module)}, \hfill \\
\cdot &: &AM \longrightarrow M &\mbox{(la multiplication scalaire externe)},
\end{matrix} \hfill \\
{ \ } \\
&\mbox{et un symbole de constante} \hfill \\
{ \ } \\
&\qquad 0 :  \ \longrightarrow M \qquad \mbox{(l'unit\'e de l'addition dans le module)}. \hfill
\end{matrix} \right.$
\end{listeimarge}
\end{defn}

\begin{remarksqed}
\begin{listeisansmarge}
\item Le langage de la th\'eorie des actions d'un mono{\"\i}de ou groupe $G$ fix\'e (constitu\'e d'un ensemble muni d'une structure de mono{\"\i}de ou de groupe) consiste en

\medskip

$\left\{\begin{matrix}
\bullet &\mbox{une sorte $X$ (l'objet de l'action),} \hfill \\
{ \ } \\
\bullet &\mbox{un symbole de fonction unaire} \hfill \\
{ \ } \\
&\qquad g \cdot \bullet : X \longrightarrow X \hfill \\
{ \ } \\
&\mbox{pour chaque \'el\'ement $g$ de $G$.} \hfill
\end{matrix} \right.$

\medskip

\item Le langage de la th\'eorie des modules sur un anneau $A$ fix\'e (constitu\'e d'un ensemble muni d'une structure d'anneau) consiste en

\medskip

$\left\{\begin{matrix}
\bullet &\mbox{une sorte $M$ (le module),} \hfill \\
{ \ } \\
\bullet &\mbox{une op\'eration binaire} \hfill \\
{ \ } \\
&\qquad + : MM \longrightarrow M \qquad \mbox{(l'addition)}, \hfill \\
{ \ } \\
&\mbox{une famille d'op\'erations unaires} \hfill \\
{ \ } \\
&\qquad a \cdot \bullet : M \longrightarrow M \qquad \mbox{(la multiplication par le scalaire $a$)} \hfill \\
{ \ } \\
&\mbox{index\'ees par les \'el\'ements $a$ de $A$, et un symbole de constante} \hfill \\
{ \ } \\
&\qquad 0 : \ \longrightarrow M \qquad \mbox{(l'unit\'e de l'addition)}. \hfill
\end{matrix} \right.$

\medskip

\item Le langage de la th\'eorie des mono{\"\i}des [resp. groupes, resp. anneaux] commutatifs est le m\^eme que celui de la th\'eorie des mono{\"\i}des [resp. groupes, resp. anneaux]. Ces th\'eories ne diff\`erent que par leurs axiomes.

\medskip

\item De m\^eme, le langage de la th\'eorie des corps [resp. des corps commutatifs] est le m\^eme que celui de la th\'eorie des anneaux. 

\end{listeisansmarge}
\end{remarksqed}

\medskip

Explicitons d'autre part le langage de la th\'eorie des cat\'egories:

\begin{defn}\label{defV13}
\begin{listeimarge}
\item Le langage de la th\'eorie des carquois consiste en

\medskip

$\left\{\begin{matrix}
\bullet &\mbox{deux sortes} \hfill \\
{ \ } \\
&\qquad \begin{matrix}
{\rm Ob} &\mbox{(les objets)}, \hfill \\
{\rm Fl} &\mbox{(les fl\`eches)},
\end{matrix} \hfill \\
{ \ } \\
\bullet &\mbox{deux symboles de fonctions unaires} \hfill \\
{ \ } \\
&\qquad \begin{matrix}
s : {\rm Fl} \longrightarrow {\rm Ob} &\mbox{(la source)}, \hfill \\
b : {\rm Fl} \longrightarrow {\rm Ob} &\mbox{(le but)}. \hfill
\end{matrix} \hfill
\end{matrix} \right.$

\medskip

\item Le langage de la th\'eorie des cat\'egories consiste en

\medskip

$\left\{\begin{matrix}
\bullet &\mbox{trois sortes} \hfill \\
&\qquad \begin{matrix}
\hfill {\rm Ob} &\mbox{(les objets)}, \hfill \\
\hfill {\rm Hom} &\mbox{(les fl\`eches ou morphismes)}, \hfill \\
{\rm Comp} &\mbox{(les paires de morphismes composables)},
\end{matrix} \hfill \\
{ \ } \\
\bullet &\mbox{six symboles de fonctions unaires} \hfill \\
{ \ } \\
&\qquad \begin{matrix}
\hfill s : {\rm Hom} &\longrightarrow &{\rm Ob} \hfill &\mbox{(la source)}, \hfill \\
\hfill b : {\rm Hom} &\longrightarrow &{\rm Ob} \hfill &\mbox{(le but)}, \hfill \\
\hfill {\rm id} : {\rm Ob} &\longrightarrow &{\rm Hom} &\mbox{(le morphisme identit\'e d'un objet)}, \hfill \\
\hfill p_1 : {\rm Comp} &\longrightarrow &{\rm Hom} &\mbox{(premi\`ere composante d'une paire composable)}, \hfill \\
\hfill p_2 : {\rm Comp} &\longrightarrow &{\rm Hom} &\mbox{(seconde composante d'une paire composable)}, \hfill \\
\hfill \circ : {\rm Comp} &\longrightarrow &{\rm Hom} &\mbox{(compos\'e d'une paire de morphismes composables)}. \hfill 
\end{matrix} \hfill \\
\end{matrix} \right.$

\medskip

\item Le langage de la th\'eorie des foncteurs consiste en deux copies $({\rm Ob}_1 , {\rm Hom}_1 , {\rm Comp}_1 , \cdots)$ et $({\rm Ob}_2 , {\rm Hom}_2$, ${\rm Comp}_2 , \cdots)$ du langage de la th\'eorie des cat\'egories compl\'et\'ees par

\medskip

$\left\{\begin{matrix}
\bullet &\mbox{trois symboles de fonctions unaires} \hfill \\
{ \ } \\
&\begin{matrix}
F : {\rm Ob}_1 \longrightarrow {\rm Ob}_2 \, , \hfill \\
F : {\rm Hom}_1 \longrightarrow {\rm Hom}_2 \, ,  \hfill \\
F : {\rm Comp}_1 \longrightarrow {\rm Comp}_2 \, .
\end{matrix}
\end{matrix} \right.$
\end{listeimarge}
\end{defn}

\begin{remarksqed}
\begin{listeisansmarge}
\item Le langage de la th\'eorie des diagrammes est la r\'eunion du langage $({\rm Ob}_c , {\rm Fl} , \cdots)$ de la th\'eorie des carquois et du langage $({\rm Ob}, {\rm Hom}, {\rm Comp}, \cdots)$ de la th\'eorie des cat\'egories compl\'et\'ee par

\medskip

$\left\{\begin{matrix}
\bullet &\mbox{deux symboles de fonctions unaires} \hfill \\
{ \ } \\
&\begin{matrix}
o:{\rm Ob}_c &\longrightarrow &{\rm Ob} \, , \hfill \\
h:{\rm Fl} \hfill &\longrightarrow &{\rm Hom} \, .
\end{matrix}
\end{matrix} \right.$

\bigskip

\item Si $D$ est un carquois fix\'e (constitu\'e de deux ensembles ${\rm Ob} (D)$ et ${\rm Fl} (D)$ reli\'es par deux applications ${\rm Fl} (D) \! \raisebox{.7ex}{\xymatrix{\dar[r]^-{^{^{\mbox{\scriptsize$s$}}}}_-{b} &{\rm Ob} (D)}}$), le langage de la th\'eorie des $D$-diagrammes est constitu\'e du langage $({\rm Ob}, {\rm Hom}, {\rm Comp}, \cdots)$ de la th\'eorie des cat\'egories compl\'et\'e par

\medskip

$\left\{\begin{matrix}
\bullet &\mbox{une famille de symboles de constantes} \hfill \\
{ \ } \\
&X_d : \ \longrightarrow {\rm Ob} \\
{ \ } \\
&\mbox{index\'es par les \'el\'ements $d$ de ${\rm Ob} (D)$}, \hfill \\
{ \ } \\
\bullet &\mbox{une famille de symboles de constantes} \hfill \\
{ \ } \\
&x_{\alpha} : \ \longrightarrow {\rm Hom} \\
{ \ } \\
&\mbox{index\'es par les \'el\'ements $\alpha$ de ${\rm Fl} (D)$.} \hfill
\end{matrix} \right.$ 

\end{listeisansmarge}
\end{remarksqed}

\bigskip

La th\'eorie des relations d'\'equivalence, la th\'eorie des relations d'ordre et la th\'eorie des ensembles sont \'ecrites dans le m\^eme langage tr\`es simple:

\begin{defn}\label{defV14}

Le langage de la th\'eorie des relations d'\'equivalence [resp. des relations d'ordre, resp. des ensembles] consiste en
\begin{enumerate}
\item[$\bullet$] une sorte $E$,
\item[$\bullet$] un symbole de relation binaire
$$
\xymatrix{
R \ \ \ar@{>->}[r] &EE
}
$$
not\'ee $R = \ \sim$ [resp. $R = \ \leq$, resp. $R = \ \in$].
\end{enumerate}
\end{defn}

\bigskip

\begin{remarkqed}

Bien s\^ur, la th\'eorie des relations d'\'equivalence et celle des relations d'ordre ont deux axiomes en commun (r\'eflexivit\'e et transitivit\'e) et elles diff\`erent par leur troisi\`eme axiome (sym\'etrie pour les relations d'\'equivalence, antisym\'etrie pour les relations d'ordre).

\smallskip 

Quant aux axiomes de la th\'eorie des ensembles, ils sont compl\`etement diff\'erents. 
\end{remarkqed}

\bigskip

Le langage de la th\'eorie des entiers de Peano est \'egalement tr\`es simple:

\begin{defn}\label{defV15}

Le langage de la th\'eorie des entiers de Peano consiste en
\begin{enumerate}
\item[$\bullet$] une sorte

\medskip

\qquad \qquad $N \qquad \mbox{(les entiers),}$

\item[$\bullet$] un symbole de fonction unaire

\medskip

\qquad \qquad $\bullet + 1 : N \longrightarrow N \qquad \mbox{(le passage au successeur)}$

\medskip

et un symbole de constante

\medskip

\qquad \qquad $0 : \ \longrightarrow N \qquad \mbox{(l'\'el\'ement initial).}$
\end{enumerate}
\end{defn}

Enfin, explicitons le langage de la th\'eorie des plans affines ou des plans euclidiens:

\begin{defn}\label{defV16}
\begin{listeimarge}
\item Le langage de la th\'eorie des plans affines consiste en

\medskip

$\left\{\begin{matrix}
\bullet &\mbox{deux sortes} \hfill \\
&\qquad \begin{matrix}
{\mathcal P} &\mbox{(les points)}, \hfill \\
{\mathcal D} &\mbox{(les droites)},
\end{matrix} \hfill\\
{ \ } \\
\bullet &\mbox{deux symboles de relations binaires} \hfill \\
{ \ } \\
&\qquad \begin{matrix}
\xymatrix{
\in \ \ \ar@{>->}[r] &{\mathcal P}{\mathcal D}
} &\mbox{(l'appartenance d'un point \`a une droite)}, \hfill \\
\xymatrix{
/\!/ \ \ \ar@{>->}[r] &{\mathcal D}{\mathcal D}
} &\mbox{(le parall\'elisme de deux droites)}, \hfill \\
\end{matrix} \hfill \\
{ \ } \\
\bullet &\mbox{trois symboles de constante} \hfill \\
{ \ } \\
&\qquad \begin{matrix}
O : \ \longrightarrow {\mathcal P} &\mbox{(l'origine)}, \hfill \\
I : \ \longrightarrow {\mathcal P} &\mbox{(le rep\`ere de l'axe des abscisses)}, \hfill \\
J : \ \longrightarrow {\mathcal P} &\mbox{(le rep\`ere de l'axe des ordonn\'ees)}. \hfill
\end{matrix} \hfill
\end{matrix} \right.$

\medskip

\item Le langage de la th\'eorie des plans euclidiens est constitu\'e de celui de la th\'eorie des plans affines compl\'et\'e par

\medskip

$\left\{\begin{matrix}
\bullet &\mbox{un symbole de relation binaire} \hfill \\
{ \ } \\
&\qquad \begin{matrix}
\xymatrix{\perp \ \ \ar@{>->}[r] &{\mathcal D}{\mathcal D}
}&\mbox{(l'orthogonalit\'e de deux droites)}. \hfill 
\end{matrix} \hfill
\end{matrix} \right.$
\end{listeimarge}
\end{defn}

\begin{remarkqed}

On peut enrichir le langage de la th\'eorie des plans euclidiens en lui ajoutant encore

\medskip

$\left\{\begin{matrix}
\bullet &\mbox{une sorte} \hfill \\
{ \ } \\
&\qquad \begin{matrix}
{\mathcal C} &\mbox{(les cercles)}, \hfill 
\end{matrix} \hfill \\
{ \ } \\
\bullet &\mbox{trois symboles de relations binaires} \hfill \\
{ \ } \\
&\qquad \begin{matrix}
\xymatrix{\in \ \ \ar@{>->}[r] &{\mathcal P}{\mathcal C}
}&\mbox{(l'appartenance d'un point \`a un cercle)}, \hfill \\
\xymatrix{{\rm Tg} \ \ \ar@{>->}[r] &{\mathcal D}{\mathcal C}
}&\mbox{(la tangence d'une droite \`a un cercle)}, \hfill \\
\xymatrix{{\rm Tg} \ \ \ar@{>->}[r] &{\mathcal C}{\mathcal C}
}&\mbox{(la tangence de deux cercles)}, \hfill \\
\end{matrix} \hfill \\
{ \ } \\
\bullet &\mbox{un symbole de fonction unitaire} \hfill \\
{ \ } \\
&\qquad C : {\mathcal C} \longrightarrow {\mathcal P} \qquad \mbox{(associer \`a tout cercle son centre).} \hfill
\end{matrix} \right.$

\end{remarkqed}

\bigskip

En revanche, les th\'eories math\'ematiques qui font appel \`a des constructions ``exponentielles'' $\Omega^A$ (classifiant les sous-objets d'un objet $A$) o\`u $B^A$ (classifiant les morphismes d'un objet $A$ dans un objet $B$) ne peuvent pas \^etre exprim\'ees dans un langage du premier ordre.

\smallskip

On dit que ce sont des th\'eories d'ordre sup\'erieur.

\smallskip

Il en est ainsi par exemple de la th\'eorie des espaces topologiques et de celle des sites.

\subsection{La notion de $\Sigma$-structure dans une cat\'egorie}\label{subsec512}

\medskip

Un langage du premier ordre est susceptible de s'exprimer dans le contexte de n'importe quelle cat\'egorie qui poss\`ede des produits finis:

\begin{defn}\label{defV17}

Soit $\Sigma$ une signature (ou ``langage du premier ordre'').

\smallskip

Soit ${\mathcal C}$ une cat\'egorie localement petite qui poss\`ede des produits finis, en particulier un objet terminal $1_{\mathcal C}$.

\smallskip

On d\'efinit alors:

\begin{listeimarge}

\item Une $\Sigma$-structure dans ${\mathcal C}$ est une fonction $M$ qui associe

\medskip

$
\left\lmoustache\begin{matrix}
\bullet &\mbox{\`a toute sorte $A$ de $\Sigma$} \hfill \\
{ \ } \\
&\qquad \mbox{un objet $M\!A$ de ${\mathcal C}$,} \hfill  \\
{ \ } \\
\bullet &\mbox{\`a tout symbole de fonction $f : A_1 \cdots A_n \to B$ de $\Sigma$} \hfill \\
{ \ } \\
&\qquad \mbox{un morphisme $M\!A_1 \times \cdots \times M\!A_n \to M\!B$ de ${\mathcal C}$} \hfill \\
{ \ } \\
&\mbox{(c'est-\`a-dire $1_{\mathcal C} \to M\!B$ si $n=0$),} \hfill 
\end{matrix} \right.
$

$
\left\rmoustache\begin{matrix}
\bullet &\mbox{\`a tout symbole de relation $R \rightarrowtail A_1 \cdots A_n$} \hfill \\
{ \ } \\
&\qquad \mbox{un sous-objet  $M\!R \hookrightarrow M\!A_1 \times \cdots \times M\!A_n$} \hfill \\
{ \ } \\
&\mbox{(c'est-\`a-dire $M\!R \hookrightarrow 1_{\mathcal C}$ si $n=0$).} \hfill
\end{matrix} \right.
$

\medskip

\item Un morphisme
$$
u : M \longrightarrow N
$$
entre deux $\Sigma$-structures $M$ et $N$ de ${\mathcal C}$ est une famille de morphismes de ${\mathcal C}$
$$
u_A : M\!A \longrightarrow N\!A
$$
index\'ee par les sortes $A$ de ${\mathcal C}$ et telle que:

\medskip

$
\left\{\begin{matrix}
{\rm (1)} &\mbox{pour tout symbole de fonction $f : A_1 \cdots A_n \to B$ de $\Sigma$, le carr\'e de ${\mathcal C}$} \hfill \\
{ \ } \\
&\xymatrix{
M\!A_1 \times \cdots \times M\!A_n \ar[d]_-{u_{A_1} \times \cdots \times u_{A_n}} \ar[rr]^-{M\!f} &&M\!B \ar[d]^-{u_B} \\
N\!A_1 \times \cdots \times N\!A_n \ar[rr]^-{N\!f} &&N\!B
} \\
&\mbox{est commutatif,} \hfill \\
{ \ } \\
{\rm (2)} &\mbox{pour tout symbole de relation $R \rightarrowtail A_1 \cdots A_n$ de $\Sigma$, il existe dans ${\mathcal C}$ un carr\'e commutatif} \hfill \\
&\mbox{(n\'ecessairement unique):} \hfill \\
{ \ } \\
&\xymatrix{
M\!R \, \ar[d] \ar@{^{(}->}[r] &M\!A_1 \times \cdots \times M\!A_n \ar[d]^-{u_{A_1} \times \cdots \times u_{A_n}} \\
N\!R \, \ar@{^{(}->}[r] &N\!A_1 \times \cdots \times N\!A_n
}
\end{matrix} \right.
$
\end{listeimarge}
\end{defn}

\bigskip

\begin{remarksqed}
\begin{listeisansmarge}
\item La condition (1) de (ii) appliqu\'ee \`a un symbole de constante $f : \ \to B$ signifie que le triangle
$$
\xymatrix{
&M\!B \ar[dd]^{u_B} \\
1_{\mathcal C} \ar[ru]^{M\!f} \ar[rd]_{N\!f} \\
&N\!B
}
$$
doit \^etre commutatif.

\smallskip

De m\^eme, la condition (2) de (ii) appliqu\'ee \`a une relation propositionnelle $R \rightarrowtail$ signifie que doit exister un triangle commutatif
$$
\xymatrix{
M\!R \ \ \ar[dd] \ar@{^{(}->}[rd] \\
&1_{\mathcal C} \\
N\!R \ar@{^{(}->}[ru]
}
$$
autrement dit que le sous-objet $M\!R$ de $1_{\mathcal C}$ doit \^etre contenu dans le sous-objet $N\!R$.

\medskip

\item Deux morphismes de $\Sigma$-structures dans ${\mathcal C}$
$$
M \xrightarrow{ \ u \ } N \xrightarrow{ \ v \ } P
$$
admettent un morphisme de $\Sigma$-structure compos\'e
$$
v \circ u : M \longrightarrow P
$$
qui consiste en la famille des morphismes de ${\mathcal C}$ compos\'es
$$
v_A \circ u_A : M\!A \longrightarrow N\!A \longrightarrow P\!A
$$
index\'es par les sortes $A$ de $\Sigma$.

\smallskip

D'autre part, on peut associer \`a toute $\Sigma$-structure $M$ dans ${\mathcal C}$ le morphisme
$$
{\rm id} : M \longrightarrow M
$$
qui consiste en la famille des morphismes d'identit\'e de ${\mathcal C}$
$$
{\rm id} : M\!A \longrightarrow M\!A
$$
index\'es par les sortes $A$ de $\Sigma$.

\smallskip

La loi de composition des morphismes de $\Sigma$-structures est associative, et le morphisme ${\rm id} : M \to M$ associ\'e \`a toute $\Sigma$-structure $M$ de ${\mathcal C}$ est neutre pour la composition \`a gauche ou \`a droite des morphismes.

\smallskip

Ainsi, les $\Sigma$-structures de ${\mathcal C}$ forment une cat\'egorie que l'on peut noter
$$
\Sigma\mbox{-str} \, ({\mathcal C}) \, .
$$
Elle est localement petite puisque ${\mathcal C}$ est localement petite par hypoth\`ese.

\medskip

\item Si ${\mathcal C}$ et ${\mathcal D}$ sont deux cat\'egories localement petites qui admettent des produits finis, tout foncteur
$$
F : {\mathcal C} \longrightarrow {\mathcal D}
$$
qui respecte les produits finis [resp. et les monomorphismes si $\Sigma$ a des symboles de relations] induit un foncteur
$$
F : \Sigma\mbox{-str} \, ({\mathcal C}) \longrightarrow \Sigma\mbox{-str} \, ({\mathcal D}) \, .
$$
Il associe \`a toute $\Sigma$-structure $M$ de ${\mathcal C}$ constitu\'ee d'objets $M\!A$, de morphismes $M\!f : M\!A_1 \times \cdots \times M\!A_n \to M\!B$ et de sous-objets $M\!R \hookrightarrow M\!A_1 \times \cdots \times M\!A_n$ la $\Sigma$-structure $F(M)$ de ${\mathcal D}$ constitu\'ee

\medskip

$\left\{\begin{matrix}
\bullet &\mbox{des objets} \hfill \\
&F(M)A = F(M\!A) \, , \\
\bullet &\mbox{des morphismes} \hfill \\
{ \ } \\
&\begin{matrix}
F(M)f = F (M\!f) : F (M\!A_1 \times \cdots \times M\!A_n) &\longrightarrow &F(M\!B) \, , \\
\qquad\qquad\Vert &&\Vert \\
\hfill F(M) A_1 \times \cdots \times F(M)A_n &&F(M)B \hfill
\end{matrix} \\
{ \ } \\
\bullet &\mbox{et des sous-objets} \hfill \\
{ \ } \\
&\begin{matrix}
F(M\!R) &\xhookrightarrow{ \ { \ } \ } &F(M\!A_1 \times \cdots \times M\!A_n) \, . \\
\Vert &&\Vert \\
F(M)R &&F(M) A_1 \times \cdots \times F(M) A_n
\end{matrix}
\end{matrix}\right.$

\bigskip

\noindent Il associe d'autre part \`a tout morphisme de $\Sigma$-structures de ${\mathcal C}$
$$
M \xrightarrow{ \ u \ } N \qquad \mbox{constitu\'e des} \qquad u_A : M\!A \longrightarrow N\!A
$$
le morphisme de $\Sigma$-structures de ${\mathcal D}$
$$
F(M) \xrightarrow{ \ F(u) \ } F(N) \qquad \mbox{constitu\'e des} \qquad F(u_A) : F(M\!A) \longrightarrow F(N\!A) \, .
$$
\end{listeisansmarge}
\end{remarksqed}

\section{Formules, axiomes et leurs interpr\'etations}\label{sec52}

\subsection{La notion de formule}\label{subsec521}

\medskip

Pour n'importe quelle signature $\Sigma$, les sortes $A$ de $\Sigma$ servent \`a nommer des objets $M\!A$ dans des cat\'egories ${\mathcal C}$ et les suites finies $A_1 \cdots A_n$ de sortes, \'eventuellement r\'ep\'et\'ees, servent \`a nommer les produits $M\!A_1 \times \cdots \times M\!A_n$ des objets correspondants dans la cat\'egorie ${\mathcal C}$. De plus, les symboles de fonctions $f : A_1 \cdots A_n \to B$ servent \`a nommer des morphismes $M\!f : M\!A_1 \times \cdots \times M\!A_n \to M\!B$ de ${\mathcal C}$, et les symboles de relations $R \rightarrowtail A_1 \cdots A_n$ servent \`a nommer des sous-objets $M\!R \hookrightarrow M\!A_1 \times \cdots \times M\!A_n$. Une famille de telles donn\'ees constitue ce que l'on a appel\'e au paragraphe pr\'ec\'edent une $\Sigma$-structure $M$ dans ${\mathcal C}$.

\smallskip

Nous allons d\'efinir une notion g\'en\'erale de formule \'ecrite dans le langage d'une signature $\Sigma$. Toute formule a un ``contexte'' qui consiste en une famille finie $\vec x = \vec x^{\vec A}$ de variables $x_i^{A_i}$, $1 \leq i \leq n$, dont chacune est affect\'ee \`a une sorte $A_i$ de $\Sigma$. Plusieurs variables peuvent \^etre affect\'ees \`a la m\^eme sorte.

\smallskip

Une formule de contexte $\vec x = (x_i^{A_i})_{1 \leq i \leq n}$ est un proc\'ed\'e de construction qui permet d'associer \`a toute $\Sigma$-structure $M$ dans une cat\'egorie ${\mathcal C}$ qui poss\`ede suffisamment de propri\'et\'es (et en particulier \`a toute $\Sigma$-structure $M$ dans n'importe quel topos) un sous-objet de $M\!A_1 \times \cdots \times M\!A_n$.

\smallskip

Voici les \'el\'ements constitutifs possibles d'une formule \'ecrite dans le langage d'une signature $\Sigma$:

\begin{defn}\label{defV21}

Soit $\Sigma$ une signature.

\begin{listeimarge}

\item Les possibles \'el\'ements constitutifs d'une formule $\varphi$ de $\Sigma$ sont
\begin{enumerate}
\item[$\bullet$] une famille $\vec x = (x_i^{A_i})$ finie ou infinie de variables $x_i^{A_i}$ dont chacune est affect\'ee \`a une sorte $A_i$ de $\Sigma$ (sans exclure des r\'ep\'etitions),
\item[$\bullet$] des \'ecritures
$$
R (x_1^{A_1} \cdots x_n^{A_n})
$$
qui consistent en un symbole de relation $R \rightarrowtail A_1 \cdots A_n$ de $\Sigma$ et certaines des variables de $\vec x$ affect\'ees aux sortes $A_1 \cdots A_n$,
\item[$\bullet$] des \'ecritures
$$
(x_1^{A_1} \cdots x_n^{A_n}) = (y_1^{A_1} \cdots y_n^{A_n}) 
$$
qui consistent en le symbole d'\'egalit\'e $=$ appliqu\'e \`a deux sous-familles finies de variables de $\vec x$ affect\'ees \`a la m\^eme suite finie de sortes $A_1 \cdots A_n$,
\item[$\bullet$] des \'ecritures
$$
\top (x_1^{A_1} \cdots x_n^{A_n})
$$
ou
$$
\perp (x_1^{A_1} \cdots x_n^{A_n})
$$
qui consistent en le symbole du vrai $\top$ ou celui du faux $\perp$ appliqu\'e \`a une sous-famille finie de variables de $\vec x$,
\item[$\bullet$] des op\'erations de substitution qui consistent \`a remplacer une variable affect\'ee \`a une sorte $B$
$$
y^B
$$
par une \'ecriture de la forme
$$
f(x_1^{A_1} \cdots x_n^{A_n})
$$
pour un symbole de fonction $f : A_1 \cdots A_n \to B$ de $\Sigma$ et une sous-famille $(x_1^{A_1} \cdots x_n^{A_n})$ de variables de $\vec x$ affect\'ees aux sortes $A_1 \cdots A_n$,
\item[$\bullet$] des quantificateurs
$$
\mbox{existentiels} \quad \exists
$$
ou
$$
\mbox{universels} \quad \forall
$$
portant sur une partie finie des variables de $\vec x$,
\item[$\bullet$] des symboles
$$
\mbox{de conjonction finie $\wedge$ ou infinie $\bigwedge$}
$$
ou
$$
\mbox{de disjonction finie $\vee$ ou infinie $\bigvee$}
$$
qui permettent de former des sous-formules
$$
\varphi_1 \wedge \cdots \wedge \varphi_k \, , \qquad \bigwedge_{i \in I} \varphi_i \, ,
$$
ou
$$
\varphi_1 \vee \cdots \vee \varphi_k \, , \qquad \bigvee_{i \in I} \varphi_i 
$$
\`a partir de sous-formules plus petites $\varphi_i$,
\item[$\bullet$] des symboles d'implication $\Rightarrow$ et de n\'egation $\neg$ qui permettent de former des sous-formules
$$
\varphi_1 \Rightarrow \varphi_2 \qquad \mbox{ou} \qquad \neg \ \varphi
$$
\`a partir de sous-formules plus petites $\varphi_1 , \varphi_2$ ou $\varphi$.
\end{enumerate}

\medskip

\item On appelle ``contexte'' d'une formule $\varphi$ de $\Sigma$ une famille finie de variables $(x_1^{A_1} , \cdots , x_n^{A_n})$ affect\'ees chacune \`a une sorte $A_i$ de $\Sigma$ et qui contient toutes les variables constitutives de $\varphi$ sur lesquelles ne porte aucun quantificateur $\exists$ ou $\forall$.
\end{listeimarge}
\end{defn}

\bigskip

\begin{remarksqed}
\begin{listeisansmarge}
\item Les variables constitutives d'une formule de $\Sigma$ sur lesquelles porte un quantificateur sont dites li\'ees. Celles sur lesquelles ne porte aucun quantificateur sont dites libres. Ainsi, le contexte d'une formule doit comprendre toutes ses variables libres mais il peut \^etre plus grand.

\medskip

\item S'il n'y a pas ambigu{\"\i}t\'e sur les variables impliqu\'ees, on \'ecrit habituellement
$$
\top \qquad \mbox{ou} \qquad \perp
$$
plut\^ot que
$$
\top (x_1^{A_1} \cdots x_n^{A_n}) \qquad \mbox{ou} \qquad \perp (x_1^{A_1} \cdots x_n^{A_n}) \, .
$$

\item On pr\'ecisera plus loin comment les formules sont constitu\'ees par l'usage it\'er\'e de symboles $\exists$, $\forall$, $\wedge$, $\bigwedge$, $\vee$, $\bigvee$, $\Rightarrow$ ou $\neg$ \`a partir de formules dites ``atomiques'' au sens qu'elles ne peuvent \^etre d\'ecompos\'ees en formules plus petites. 

\end{listeisansmarge}
\end{remarksqed}

\bigskip

Illustrons tout de suite cette d\'efinition abstraite par des exemples de formules qui entrent dans la d\'efinition de th\'eories du premier ordre d'usage courant:

\smallskip

Ainsi, l'axiome d'associativit\'e de la th\'eorie des mono{\"\i}des de signature $(G, \cdot , 1)$ dit que la formule
$$
g_1 \cdot (g_2 \cdot g_3) = (g_1 \cdot g_2) \cdot g_3
$$
doit \^etre v\'erifi\'ee universellement en les variables $g_1 , g_2 , g_3$ de sorte $G$.

\smallskip

L'axiome des inverses de la th\'eorie des groupes de signature $(G , \cdot , 1 , (\bullet)^{-1})$ dit que les formules
$$
g \cdot g^{-1} = 1 \qquad \mbox{et} \qquad g^{-1} \cdot g = 1
$$
doivent \^etre v\'erifi\'ees universellement en la variable $g$ de sorte $G$.

\smallskip

L'axiome de distributivit\'e de la th\'eorie des anneaux de signature $(A,+,\cdot ,0,1,-(\cdot))$ dit que les formules
$$
a \cdot (a' + a'') = a \cdot a' + a \cdot a'' \quad \mbox{et} \quad (a' + a'') \cdot a = a' \cdot a + a'' \cdot a
$$
doivent \^etre v\'erifi\'ees universellement en les variables $a,a',a''$ de sorte $A$.

\smallskip

L'axiome d'associativit\'e de la multiplication scalaire externe de la th\'eorie des modules sur les anneaux de signature $(A,M,\cdots)$ dit que la formule
$$
a_1 \cdot (a_2 \cdot m) = (a_1 \cdot a_2) \cdot m
$$
doit \^etre v\'erifi\'ee universellement en les variables $a_1 , a_2$ de sorte $A$ et la variable $m$ de sorte $M$.

\smallskip

L'axiome d'existence des inverses en dehors de $0$ de la th\'eorie des corps de signature $(K,+,\cdot ,0,1,-(\bullet))$ dit que la formule
$$
k = 0 \vee \exists \, k' \, (k \cdot k' = 1 \wedge k' \cdot k = 1)
$$
doit \^etre v\'erifi\'ee universellement en la variable $k$ de sorte $K$.

\smallskip

L'axiome de condition de composabilit\'e des morphismes de la th\'eorie des cat\'egories de signature $({\rm Ob}, {\rm Hom}$, ${\rm Comp}, s, b, {\rm id}, p_1, p_2, \circ)$ dit que la formule en les variables $f,g$ de sorte Hom
$$
b(f) = s(g)
$$
doit \^etre \'equivalente \`a la formule \'ecrite avec la variable li\'ee $c$ de sorte Comp
$$
\exists \, c \ (p_1 (c) = f \wedge p_2 (c)=g) \, .
$$

L'axiome de transitivit\'e de la th\'eorie des relations d'\'equivalence de sorte $(E,\sim)$ dit que la formule en les variables $e_1 , e_2 , e_3$ de sorte $E$
$$
e_1 \sim e_2 \wedge e_2 \sim e_3
$$
doit impliquer la formule
$$
e_1 \sim e_3 \, .
$$

L'axiome d'existence d'un pr\'ed\'ecesseur de tout \'el\'ement autre que l'\'el\'ement initial dans la th\'eorie des entiers de Peano de signature $(N,(\bullet) + 1 , 0)$ dit que la formule
$$
n = 0 \vee \exists \, n' (n'+1 = n)
$$
doit \^etre v\'erifi\'ee universellement en la variable $n$ de sorte $N$.

\smallskip

Dans la th\'eorie des plans affines de signature $({\mathcal P}, {\mathcal D}, \in , /\!/ , \cdots)$, l'axiome d'existence d'une droite parall\`ele \`a une droite donn\'ee passant par un point donn\'e dit que la formule
$$
\exists \, D' (D' /\!/ D \wedge P \in D')
$$
doit \^etre v\'erifi\'ee universellement en la variable $D$ de sorte ${\mathcal D}$ et la variable $P$ de sorte ${\mathcal P}$.

\smallskip

Dans la th\'eorie des plans euclidiens de signature $({\mathcal P}, {\mathcal D}, \in , /\!/ , \perp, \cdots)$, l'axiome de rencontre des hauteurs d'un triangle dit que la formule qui implique les variables li\'ees $D_1 , D_2 , D_3 , D'_1 , D'_2 , D'_3$ de sorte ${\mathcal D}$ et la variable li\'ee $P$ de sorte ${\mathcal P}$
$$
\exists \, D_1 , D_2 , D_3 , D'_1 , D'_2 , D'_3 , P
$$
$$
(D_1 \perp D'_1 \wedge D_2 \perp D'_2 \wedge D_3 \perp D'_3 \wedge A \in D'_1 \wedge B \in D'_2 \wedge C \in D'_3
$$
$$
\wedge \, A \in D_2 \wedge A \in D_3 \wedge B \in D_1 \wedge B \in D_3 \wedge C \in D_1 \wedge C \in D_2 \wedge P \in D'_1 \wedge P \in D'_2 \wedge P \in D'_3)
$$
doit \^etre v\'erifi\'ee universellement en les variables libres $A,B,C$ de sorte ${\mathcal P}$.

\subsection{L'interpr\'etation d'une formule dans une $\Sigma$-structure}\label{subsec522}

On pose la d\'efinition suivante qui, pour le moment, introduit surtout une notation.

\smallskip

Le sens de cette notation sera progressivement pr\'ecis\'e et explicit\'e dans les paragraphes suivants.

\begin{defn}\label{defV22}

Soit $\Sigma$ une signature.

\smallskip

Soit ${\mathcal C}$ une cat\'egorie localement petite qui admet des produits finis, en particulier un objet terminal.

\smallskip

Soit $\varphi$ une formule de $\Sigma$ de contexte $\vec x = (x_1^{A_1} \cdots x_n^{A_n})$.

\smallskip

Si ${\mathcal C}$ poss\`ede suffisamment de propri\'et\'es cat\'egoriques, en un sens qui sera pr\'ecis\'e plus tard et d\'ependra de la liste des symboles qui apparaissent dans l'\'ecriture de $\varphi$, on notera, pour toute $\Sigma$-structure $M$ dans ${\mathcal C}$,
$$
M\varphi (\vec x) \xhookrightarrow{ \ { \ } \ } M\!A_1 \times \cdots \times M\!A_n
$$
le sous-objet de $M\!A_1 \times \cdots \times M\!A_n$ d\'efini par la formule $\varphi$ en les variables $x_1^{A_1} \cdots x_n^{A_n}$ de sortes $A_1 \cdots A_n$.
\end{defn}

\bigskip

\begin{remarksqed}
\begin{listeisansmarge}
\item Si $\varphi$ est une formule de $\Sigma$ de contexte $\vec x = \vec x^{\vec A} = (x_1^{A_1} \cdots x_n^{A_n})$, les logiciens notent habituellement son interpr\'etation dans une $\Sigma$-structure $M$ d'une cat\'egorie ${\mathcal C}$ par
$$
\llbracket \vec x \cdot \varphi \rrbracket_M \qquad \mbox{ou} \qquad \llbracket \vec x^{\vec A} \cdot \varphi \rrbracket_M
$$
plut\^ot que par
$$
M \varphi (\vec x) \, .
$$

\item Nous pouvons d\'ej\`a donner une id\'ee de la mani\`ere dont une formule $\varphi$ dans un contexte $\vec x$ est interpr\'et\'ee dans une $\Sigma$-structure $M$ d'une cat\'egorie ${\mathcal C}$ qui poss\`ede suffisamment de propri\'et\'es cat\'egoriques.

\smallskip

Pour cela, reprenons les \'el\'ements constitutifs des formules $\varphi$ et annon\c cons comment chacun sera interpr\'et\'e:
\begin{enumerate}
\item[$\bullet$] Si $R \rightarrowtail A_1 \cdots A_n$ est un symbole de relation de $\Sigma$, une \'ecriture
$$
R(x_1^{A_1} \cdots x_n^{A_n})
$$
est interpr\'et\'ee comme le sous-objet
$$
M\!R \xhookrightarrow{ \ { \ } \ } M\!A_1 \times \cdots \times M\!A_n \, .
$$
Dans un contexte plus large constitu\'e de $x_1^{A_1} \cdots x_n^{A_n}$ et d'autres variables $y_1^{B_1} \cdots y_k^{B_k}$, elle est interpr\'et\'ee comme le sous-objet
$$
M\!R \times M\!B_1 \times \cdots \times M\!B_k \xhookrightarrow{ \ { \ } \ } M\!A_1 \times \cdots \times M\!A_n \times M\!B_1 \times \cdots \times M\!B_k
$$
image r\'eciproque de
$$
M\!R \xhookrightarrow{ \ { \ } \ } M\!A_1 \times \cdots \times M\!A_n 
$$
par le morphisme de projection sur les $n$ premiers facteurs
$$
M\!A_1 \times \cdots \times M\!A_n \times M\!B_1 \times \cdots \times M\!B_k \longrightarrow M\!A_1 \times \cdots \times M\!A_n \, .
$$
\item[$\bullet$] Une \'ecriture de la forme
$$
(x_1^{A_1} \cdots x_n^{A_n}) = (y_1^{A_1} \cdots y_n^{A_n})
$$
est interpr\'et\'ee comme le sous-objet diagonal
$$
M\!A_1 \times \cdots \times M\!A_n \xhookrightarrow{ \ { \ } \ } (M\!A_1 \times \cdots \times M\!A_n) \times (M\!A_1 \times \cdots \times M\!A_n) \, .
$$
\item[$\bullet$] Les \'ecritures
$$
\top (x_1^{A_1} \cdots x_n^{A_n})
$$
sont interpr\'et\'ees comme les plus grands sous-objets
$$
M\!A_1 \times \cdots \times M\!A_n \xhookrightarrow{ \ = \ } M\!A_1 \times \cdots \times M\!A_n \, .
$$
Quant aux \'ecritures
$$
\perp (x_1^{A_1} \cdots x_n^{A_n}) \, ,
$$
elles sont interpr\'et\'ees comme les plus petits sous-objets
$$
\emptyset_{M\!A_1 \times \cdots \times M\!A_n} \xhookrightarrow{ \ { \ } \ } M\!A_1 \times \cdots \times M\!A_n
$$
s'ils sont bien d\'efinis dans la cat\'egorie ${\mathcal C}$.
\item[$\bullet$] Si une formule $\varphi$ de $\Sigma$ dans un contexte $(y^B y_1^{B_1} \cdots y_k^{B_k})$ est interpr\'et\'ee comme un sous-objet
$$
M \varphi \xhookrightarrow{ \ { \ } \ } M\!B \times M\!B_1 \times \cdots \times M\!B_k \, ,
$$
alors la formule $\varphi'$ qui se d\'eduit de $\varphi$ par substitution de la variable $y^B$ par une \'ecriture de la forme
$$
f(x_1^{A_1} \cdots x_n^{A_n})
$$
pour un symbole de fonction $f : A_1 \cdots A_n \to B$, s'interpr\`ete comme le sous-objet
$$
M\varphi' \xhookrightarrow{ \ { \ } \ } M\!A_1 \times \cdots \times M\!A_n \times M\!B_1 \times \cdots \times M\!B_k
$$
d\'eduit de
$$
M\varphi \xhookrightarrow{ \ { \ } \ } M\!B \times M\!B_1 \times \cdots \times M\!B_k 
$$
par le morphisme de changement de base
$$
M\!A_1 \times \cdots \times M\!A_n \times M\!B_1 \times \cdots \times M\!B_k \xrightarrow{ \ M\!f \times {\rm id} \times \cdots \times  {\rm id} \ } M\!B \times M\!B_1 \times \cdots \times M\!B_k
$$
pourvu que ce produit fibr\'e soit bien d\'efini dans la cat\'egorie ${\mathcal C}$.
\item[$\bullet$] Si une formule $\varphi$ dans un contexte $(x_1^{A_1} \cdots x_n^{A_n} \, y_1^{B_1} \cdots y_k^{B_k})$ s'interpr\`ete comme un sous-objet
$$
M\varphi \xhookrightarrow{ \ { \ } \ } M\!A_1 \times \cdots \times M\!A_n \times M\!B_1 \times \cdots \times M\!B_k \, ,
$$
alors la formule $\varphi'$ d\'eduite de $\varphi$ en appliquant le quantificateur
$$
\exists \qquad \mbox{[resp.} \quad \forall \ \mbox{]}
$$
en les variables $y_1^{B_1} \cdots y_k^{B_k}$, s'interpr\`ete comme le sous-objet
$$
M\varphi' \xhookrightarrow{ \ { \ } \ } M\!A_1 \times \cdots \times M\!A_n
$$
image de
$$
M\varphi \xhookrightarrow{ \ { \ } \ } M\!A_1 \times \cdots \times M\!A_n \times M\!B_1 \times \cdots \times M\!B_k 
$$
par le foncteur adjoint \`a gauche [resp. \`a droite]
$$
\exists_p = p_! \qquad \mbox{[resp.} \quad \forall_p = p_* \ \mbox{]}
$$
du foncteur $p^{-1}$ qui associe \`a tout sous-objet
$$
S \xhookrightarrow{ \ { \ } \ } M\!A_1 \times \cdots \times M\!A_n
$$
son image r\'eciproque
$$
S \times M\!B_1 \times \cdots \times M\!B_k \xhookrightarrow{ \ { \ } \ } M\!A_1 \times \cdots \times M\!A_n \times M\!B_1 \times \cdots \times M\!B_k
$$
par la projection
$$
p : M\!A_1 \times \cdots \times M\!A_n \times M\!B_1 \times \cdots \times M\!B_k \longrightarrow M\!A_1 \times \cdots \times M\!A_n \, ,
$$
pourvu que cet adjoint \`a gauche $\exists_p = p_!$ [resp. cet adjoint \`a droite $\forall_p = p_*$] soit bien d\'efini dans la cat\'egorie ${\mathcal C}$.
\item[$\bullet$]  Si des formules $\varphi_i$ de m\^eme contexte $(x_1^{A_1} \cdots x_n^{A_n})$ sont interpr\'et\'ees comme des sous-objets
$$
M\varphi_i \xhookrightarrow{ \ { \ } \ } M\!A_1 \times \cdots \times M\!A_n \, ,
$$
alors leur conjonction [resp. disjonction] finie ou infinie
$$
\qquad \qquad \varphi_1 \wedge \cdots \wedge \varphi_k \qquad \mbox{ou} \qquad \bigwedge_i \varphi_i
$$
$$
\mbox{[resp.} \qquad \ \varphi_1 \vee \cdots \vee \varphi_k \qquad \mbox{ou} \qquad \bigvee_i \varphi_i \ \mbox{]}
$$
est interpr\'et\'ee comme l'intersection [resp. la r\'eunion] finie ou infinie des sous-objets
$$
M\varphi_i \xhookrightarrow{ \ { \ } \ } M\!A_1 \times \cdots \times M\!A_n \, ,
$$
pourvu que les intersections [resp. r\'eunions] finies ou infinies de sous-objets soient bien d\'efinies dans la cat\'egorie ${\mathcal C}$.
\item[$\bullet$] Supposons que pour tout sous-objet d'un objet $E$ de ${\mathcal C}$
$$
S \xhookrightarrow{ \ { \ } \ } E
$$
existe un sous-objet 
$$
\neg \, S \xhookrightarrow{ \ { \ } \ } E
$$
caract\'eris\'e par la propri\'et\'e que, pour tout sous-objet
$$
S' \xhookrightarrow{ \ { \ } \ } E \, ,
$$
on a $S' \leq \neg \, S$ si et seulement si $S$ et $S'$ ont pour intersection le plus petit sous-objet $\emptyset_S$ de $S$.

\smallskip

Alors, si une formule $\varphi$ de contexte $(x_1^{A_1} \cdots x_n^{A_n})$ est interpr\'et\'ee comme un sous-objet
$$
M \varphi \xhookrightarrow{ \ { \ } \ } M\!A_1 \times \cdots \times M\!A_n \, ,
$$
sa n\'egation $\neg \, \varphi$ est interpr\'et\'ee comme le sous-objet
$$
\neg \, (M\varphi) \xhookrightarrow{ \ { \ } \ } M\!A_1 \times \cdots \times M\!A_n \, .
$$
\item[$\bullet$] Plus g\'en\'eralement, supposons que pour toute paire de sous-objets d'un objet $E$ de ${\mathcal C}$
$$
S_1 \xhookrightarrow{ \ { \ } \ } E \qquad \mbox{et} \qquad S_2 \xhookrightarrow{ \ { \ } \ } E
$$
existe un sous-objet de $E$
$$
(S_1 \Rightarrow S_2) \xhookrightarrow{ \ { \ } \ } E
$$
caract\'eris\'e par la propri\'et\'e que, pour tout sous-objet
$$
S' \xhookrightarrow{ \ { \ } \ } E \, ,
$$
on a $S' \leq (S_1 \Rightarrow S_2)$ si et seulement si l'intersection $S' \wedge S_1$ est bien d\'efinie et v\'erifie
$$
S' \wedge S_1 \leq S_2 \, .
$$

Alors, si deux formules $\varphi_1$ et $\varphi_2$ de m\^eme contexte $(x_1^{A_1} \cdots x_n^{A_n})$ sont interpr\'et\'ees comme deux sous-objets
$$
M \varphi_1 \xhookrightarrow{ \ { \ } \ } M\!A_1 \times \cdots \times M\!A_n \qquad \mbox{et} \qquad  M \varphi_2 \xhookrightarrow{ \ { \ } \ } M\!A_1 \times \cdots \times M\!A_n \, ,
$$
la formule d'implication $\varphi_1 \Rightarrow \varphi_2$ est interpr\'et\'ee comme le sous-objet
$$
(M\varphi_1 \Rightarrow M\varphi_2) \xhookrightarrow{ \ { \ } \ } M\!A_1 \times \cdots \times M\!A_n \, .
$$
\end{enumerate}

\medskip

\item Comme on verra, toutes les propri\'et\'es cat\'egoriques de ${\mathcal C}$ dont on a besoin pour interpr\'eter les formules $\varphi$ de signature $\Sigma$ dans les $\Sigma$-structures de ${\mathcal C}$ sont v\'erifi\'ees en particulier par n'importe quel topos. 

\end{listeisansmarge}
\end{remarksqed}

\subsection{La notion de s\'equent et celle de th\'eorie du premier ordre}\label{subsec523}

On pose la d\'efinition formelle suivante:

\begin{defn}\label{defV23}

Soit $\Sigma$ une signature.

\begin{listeimarge}

\item Un ``s\'equent'' de signature $\Sigma$ est une \'ecrite formelle
$$
\varphi \vdash_{\vec x} \psi
$$
o\`u $\varphi$ et $\psi$ sont deux formules de $\Sigma$ et $\vec x = (x_1^{A_1} \cdots x_n^{A_n})$ est un contexte \`a la fois pour $\varphi$ et pour $\psi$.

\medskip

\item Une ``th\'eorie du premier ordre ${\mathbb T}$'' de signature $\Sigma$ est une famille de s\'equents de $\Sigma$
$$
\varphi \vdash_{\vec x} \psi
$$
appel\'es les ``axiomes'' de la th\'eorie ${\mathbb T}$.
\end{listeimarge}
\end{defn}

\begin{remarksqed}
\begin{listeisansmarge}
\item Si le contexte d'un s\'equent est clair, on peut \'ecrire simplement
$$
\varphi \vdash \psi \, .
$$

\item La signification d'une \'ecrite formelle
$$
\varphi \vdash_{\vec x} \psi \qquad \mbox{ou} \qquad \varphi \vdash \psi
$$
est de demander que la formule $\varphi$ implique la formule $\psi$.

\smallskip

En particulier, pour demander qu'une formule $\varphi$ soit toujours v\'erifi\'ee [resp. jamais v\'erifi\'ee], on \'ecrit
$$
\top \vdash \varphi \qquad \mbox{[resp.} \quad \varphi \vdash \ \perp \ \mbox{]} \, .
$$

\item Pour demander que deux formules $\varphi$ et $\psi$ soient \'equivalentes, on peut \'ecrire
$$
\varphi \dashv \, \vdash \psi
$$
au lieu de
$$
\varphi \vdash \psi \qquad \mbox{et} \qquad \psi \vdash \varphi \, .
$$

\item La famille des axiomes d'une th\'eorie ${\mathbb T}$ peut \^etre finie ou infinie, de m\^eme que sa signature $\Sigma$ peut consister en des familles finies ou infinies de sortes, de symboles de fonctions ou de symboles de relations.

\medskip

\item Une th\'eorie ${\mathbb T}$ peut consister en sa seule signature $\Sigma$, sans aucun axiome.

\smallskip

Par exemple, la th\'eorie des carquois de signature $({\rm Ob}, {\rm Fl}, s : {\rm Fl} \to {\rm Ob} , b : {\rm Fl} \to {\rm Ob})$ est une th\'eorie sans axiomes. 
\end{listeisansmarge}
\end{remarksqed}

\bigskip

Illustrons cette d\'efinition g\'en\'erale en donnant la liste des axiomes constitutifs de th\'eories alg\'ebriques courantes:

\pagebreak 

\begin{defn}\label{defV24}
\begin{listeimarge}
\item La th\'eorie des mono{\"\i}des [resp. des groupes] de signature $(G,\cdot : GG \to G, 1 : \ \to G)$ [resp. compl\'et\'ee par $(\bullet)^{-1} : G \to G$] est d\'efinie par les axiomes
$$
\begin{matrix}
\top &\vdash &g_1 \cdot (g_2 \cdot g_3) = (g_1 \cdot g_2) \cdot g_3 \, , \\
\top &\vdash &g \cdot 1 = g \, , \hfill \\
\top &\vdash &1 \cdot g = g \, , \hfill
\end{matrix}
$$
[resp. et
$$
\mbox{\hglue-2cm}\begin{matrix}
\top &\vdash &g \cdot g^{-1} = 1 \, , \hfill \\
\top &\vdash &g^{-1} \cdot g = 1 \ \mbox{].} \hfill
\end{matrix}
$$

\item La th\'eorie des actions de mono{\"\i}des [resp. de groupes], de signature celle de la th\'eorie des mono{\"\i}des [resp. des groupes] $(G,\cdots)$ compl\'et\'ee par $(X,\cdot : GX \to X)$, est d\'efinie par les axiomes de la th\'eorie des mono{\"\i}des [resp. des groupes] auxquels on ajoute les axiomes
$$
\begin{matrix}
\top &\vdash &g_1 \cdot (g_2 \cdot x) = (g_1 \cdot g_2) \cdot x  \, , \\
\top &\vdash &1 \cdot x = x \hfill
\end{matrix}
$$
en les variables $g_1 , g_2$ de sorte $G$ et la variable $x$ de sorte $X$.

\medskip

\item La th\'eorie des anneaux de signature $(A,+:AA \to A$, $\cdot : AA \to A$, $0 : \ \to A$, $1 : \ \to A$, $-(\bullet) : A \to A)$ est d\'efinie par les trois familles d'axiomes

\medskip

$
\left\{\begin{matrix}
\top &\vdash & a_1 + (a_2 + a_3) = (a_1 + a_2) + a_3 \, , \\
\top &\vdash &a_1 + a_2 = a_2 + a_1 \, , \hfill \\
\top &\vdash &a + 0 = a \, , \hfill \\
\top &\vdash &a + (-a) = 0 \, , \hfill
\end{matrix} \right.
$

\medskip

$
\left\{\begin{matrix}
\top &\vdash & a_1 \cdot (a_2 \cdot a_3) = (a_1 \cdot a_2) \cdot a_3 \, , \\
\top &\vdash &a \cdot 1 = a \, , \hfill \\
\top &\vdash &1 \cdot a = a \, , \hfill
\end{matrix} \right.
$

\medskip

$
\left\{\begin{matrix}
\top &\vdash & a \cdot (a_1 + a_2) = a \cdot a_1 + a \cdot a_2 \, , \hfill \\
\top &\vdash &(a_1 + a_2) \cdot a = a_1 \cdot a + a_2 \cdot a \, . \hfill
\end{matrix} \right.
$

\medskip

\item La th\'eorie des modules sur un anneau, dont la signature est constitu\'ee de celle de la th\'eorie des anneaux $(A,\cdots)$ compl\'et\'ee par $(M , + : MM \to M$, $0 : \ \to M$, $-(\bullet) : M \to M$, $\cdot : AM \to M)$, est d\'efinie par les axiomes de la th\'eorie des anneaux compl\'et\'es par les deux familles d'axiomes

\medskip

$
\left\{\begin{matrix}
\top &\vdash & m_1 + (m_2 + m_3) = (m_1 + m_2) + m_3 \, , \\
\top &\vdash &m_1 + m_2 = m_2 + m_1 \, , \hfill \\
\top &\vdash &m + 0 = m \, , \hfill \\
\top &\vdash &m + (-m) = 0 \hfill
\end{matrix} \right.
$

\medskip

\noindent en des variables de sorte $M$, et

\medskip

$
\left\{\begin{matrix}
\top &\vdash &a_1 \cdot (a_2 \cdot m) = (a_1 \cdot a_2) \cdot m \, , \hfill \\
\top &\vdash &1 \cdot m = m \, , \hfill \\
\top &\vdash &(a_1 + a_2) \cdot m = a_1 \cdot m + a_2 \cdot m \, , \hfill \\
\top &\vdash &a \cdot (m_1 + m_2) = a \cdot m_1 + a \cdot m_2  \hfill
\end{matrix} \right.
$

\medskip

\noindent en des variables $a_1 , a_2 , a$ de sorte $A$ et des variables $m,m_1 , m_2$ de sorte $M$.
\end{listeimarge}
\end{defn}

\begin{remarksqed}
\begin{listeisansmarge}
\item La th\'eorie des mono{\"\i}des [resp. des groupes] commutatifs a la m\^eme signature que celle des mono{\"\i}des [resp. des groupes]. Elle est d\'efinie en ajoutant \`a ses axiomes l'axiome suppl\'ementaire
$$
\top \vdash \ g_1 \cdot g_2 = g_2 \cdot g_1 \, .
$$

\item De m\^eme, la th\'eorie des anneaux commutatifs a la m\^eme signature que celle des anneaux. Elle est d\'efinie en ajoutant \`a ses axiomes l'axiome suppl\'ementaire
$$
\top \vdash \ a_1 \cdot a_2 = a_2 \cdot a_1 \, .
$$

\item La th\'eorie des corps a la m\^eme signature que celle des anneaux. Elle est d\'efinie en ajoutant \`a ses axiomes l'axiome suppl\'ementaire
$$
\top \vdash_a  a = 0 \vee \exists \, b \ (a \cdot b = 1 \wedge b \cdot a = 1) \, .
$$

\item La th\'eorie des corps commutatifs est d\'efinie en ajoutant \`a la th\'eorie des anneaux commutatifs l'axiome
$$
\top \vdash_a  a = 0 \vee \exists \, b \ (a \cdot b = 1) \, .
$$

\item La th\'eorie des actions d'un groupe [resp. d'un mono{\"\i}de] fix\'e $G$, dont la signature est constitu\'ee d'une sorte $X$ et de symboles de fonctions unaires $g \cdot \bullet : X \to X$ index\'es par les \'el\'ements $g \in G$, est d\'efinie par les axiomes en la variable $x$ de sorte $X$
$$
\begin{matrix}
\top &\vdash \ 1 \cdot x = x \, , \hfill \\
\top &\vdash \ g_1 \cdot g_2 \cdot x = (g_1 \cdot g_2) \cdot x 
\end{matrix} 
$$
pour tous \'el\'ements $g_1 , g_2$ de $X$.

\medskip

\item La th\'eorie des modules sur un anneau fix\'e $A$, de signature $(M$, $+ : MM \to M$, $\cdot : MM \to M$, $0 : \ \to M$) compl\'et\'ee par des symboles de fonctions unaires $a \cdot \bullet : M \to M$ index\'es par les \'el\'ements $a$ de $A$, est d\'efinie par les axiomes en des variables de sorte $M$

\medskip

$
\left\{\begin{matrix}
\top &\vdash & m_1 + (m_2 + m_3) = (m_1 + m_2) + m_3 \, , \\
\top &\vdash &m_1 + m_2 = m_2 + m_1 \, , \hfill \\
\top &\vdash &m + 0 = m \, , \hfill \\
\top &\vdash &m + (-m) = 0 \hfill
\end{matrix} \right.
$

\medskip

$
\left\{\begin{matrix}
\top &\vdash &a_1 \cdot (a_2 \cdot m) = (a_1 \cdot a_2) \cdot m \, , \hfill \\
\top &\vdash &1 \cdot m = m \, , \hfill \\
\top &\vdash &(a_1 + a_2) \cdot m = a_1 \cdot m + a_2 \cdot m \, , \hfill \\
\top &\vdash &a \cdot (m_1 + m_2) = a \cdot m_1 + a \cdot m_2  \hfill
\end{matrix} \right.
$

\medskip

\noindent pour tous \'el\'ements $a , a_1 , a_2$ de $A$. 

\end{listeisansmarge}
\end{remarksqed}

\medskip

Puis donnons la liste des axiomes de la th\'eorie des cat\'egories et de celle des foncteurs:

\begin{defn}\label{defV25}
\begin{listeimarge}
\item La th\'eorie des cat\'egories, dont la signature est constitu\'ee des trois sortes {\rm Ob}, {\rm Hom}, {\rm Comp} et des six symboles de fonctions $s : {\rm Hom} \to {\rm Ob}$, $b : {\rm Hom} \to {\rm Ob}$, ${\rm id} : {\rm Ob} \to {\rm Hom}$, $p_1 : {\rm Comp} \to {\rm Hom}$, $p_2 : {\rm Comp} \to {\rm Hom}$ et $\circ : {\rm Comp} \to {\rm Hom}$, est d\'efinie par les axiomes suivants en les variables $X$ de sorte {\rm Ob}, $f,g$ de sorte {\rm Hom}, et $(f,g),(f',g'),(g,h)$ de sorte {\rm Comp}

\medskip
$
\left\lmoustache\begin{matrix}
&\top \vdash_X s({\rm id}_X) = X \, , \hfill \\
{ \ } \\
&\top \vdash_X b({\rm id}_X) = X \hfill \\
{ \ } \\
&p_1 (f,g) = p_1 (f',g') \wedge p_2 (f,g) = p_2 (f',g') \vdash (f,g) = (f',g') \hfill \\
{ \ } \\
&\mbox{(ce qui permet de noter $p_1 (f,g) = f$, $p_2 (f,g) = g$ et $\circ \, (f,g) = g \circ f$),} \hfill \\
{ \ } \\
&\exists \, (f,g) (p_1 (f,g) = f \wedge p_2 (f,g) = g) \dashv \, \vdash b(f) = s(g) \, , \hfill 
\end{matrix} \right.
$

$
\left\rmoustache\begin{matrix}
&\top \vdash_{(f,g)} s (g \circ f) = s(f) \, , \hfill \\
{ \ } \\
&\top \vdash_{(f,g)} b (g \circ f) = b(g) \, , \hfill \\
{ \ } \\
&p_1 (f,g) = {\rm id}_X \vdash_{(f,g),X} g \circ f = g \, , \hfill \\
{ \ } \\
&p_2 (f,g) = {\rm id}_X \vdash_{(f,g),X} g \circ f = f \, , \hfill \\
{ \ } \\
&p_2 (f,g) = p_1 (g,h) \wedge p_1 (f,g) = p_1 (f,f') \wedge p_2 (g,h) = p_2 (h',h) \wedge h' = g \circ f \wedge f' = h \circ g \qquad \qquad \quad { \ } \hfill \\
{ \ } \\
&\hfill \vdash_{{(f,g),(g,h) \atop (h',h),(f,f')}} h \circ h' = f' \circ f \, .
\end{matrix} \right.
$

\medskip

\item La th\'eorie des foncteurs, dont la signature est constitu\'ee de deux copies $({\rm Ob}, {\rm Hom}, {\rm Comp}, \cdots)$ et $({\rm Ob}'$, ${\rm Hom}', {\rm Comp}', \cdots)$ de la signature de la th\'eorie des cat\'egories compl\'et\'ees par trois symboles de fonctions $F : {\rm Ob} \to {\rm Ob}'$, $F : {\rm Hom} \to {\rm Hom}'$ et $F : {\rm Comp} \to {\rm Comp}'$, est d\'efinie par les axiomes de la th\'eorie des cat\'egories portant \`a la fois sur les variables de sortes ${\rm Ob}, {\rm Hom}, {\rm Comp}$ et sur celles de sortes ${\rm Ob}', {\rm Hom}', {\rm Comp}'$, compl\'et\'es par les axiomes de commutation

\medskip

$
\left\{\begin{matrix}
\top &\vdash_f \hfill & s' (F(f)) = F(s(f)) \, , \hfill \\
\top &\vdash_f \hfill & b' (F(f)) = F(b(f)) \, , \hfill \\
\top &\vdash_X \hfill &{\rm id}'_{F(X)} = F({\rm id}_X) \, , \hfill \\
\top &\vdash_{(f,g)} &p'_1 (F(f,g)) = F(p_1 (f,g)) \, , \hfill \\
\top &\vdash_{(f,g)} &p'_2 (F(f,g)) = F(p_2 (f,g)) \, , \hfill \\
\top &\vdash_{(f,g)} &\circ' (F(f,g)) = F(\circ (f,g)) \, . \hfill
\end{matrix} \right.
$
\end{listeimarge}
\end{defn}

\begin{remarksqed}
\begin{listeisansmarge}
\item Comme d\'ej\`a not\'e, la th\'eorie des carquois est la th\'eorie sans axiomes de signature $({\rm Ob} , {\rm Fl} , s : {\rm Fl} \to {\rm Ob} , b : {\rm Fl} \to {\rm Ob})$.

\medskip

\item La th\'eorie des diagrammes, dont la signature est la r\'eunion de celles de la th\'eorie des carquois $({\rm Ob}_c , {\rm Fl} , \cdots)$ et de celle de la th\'eorie des cat\'egories $({\rm Ob}, {\rm Hom}, {\rm Comp}, \cdots)$ compl\'et\'ee par deux symboles de fonctions unaires
$$
\begin{matrix}
o : {\rm Ob}_c &\longrightarrow &{\rm Ob} \, , \hfill \\
h : {\rm Fl} \hfill &\longrightarrow &{\rm Hom} \, ,
\end{matrix}
$$
est d\'efinie par les axiomes de la th\'eorie des cat\'egories compl\'et\'es par les axiomes en la variable $f$ de sorte ${\rm Fl}$
$$
\begin{matrix}
\top &\vdash_f &s (h(f)) = o(s(f)) \, , \hfill \\
\top &\vdash_f &b (h(f)) = o(b(f)) \, . \hfill 
\end{matrix}
$$

\item Si $D$ est un carquois fix\'e, constitu\'e de deux ensembles ${\rm Ob}(D)$ et ${\rm Fl}(D)$ reli\'es par deux applications ${\rm Fl} (D) \! \raisebox{.7ex}{\xymatrix{\dar[r]^-{^{^{\mbox{\scriptsize$s$}}}}_-{b} &{\rm Ob} (D)}}$, la th\'eorie des $D$-diagrammes, dont la signature est celle de la th\'eorie des cat\'egories $({\rm Ob}, {\rm Hom}, {\rm Comp}, \cdots)$ compl\'et\'ee par des symboles de constantes
$$
X_d : \ \longrightarrow {\rm Ob} \, , \qquad d \in {\rm Ob}(D) \, ,
$$
et
$$
x_{\alpha} : \ \longrightarrow {\rm Hom} \, , \qquad \alpha \in {\rm Fl}(D) \, ,
$$
est d\'efinie par les axiomes de la th\'eorie des cat\'egories compl\'et\'es par les axiomes
$$
s(x_{\alpha}) = X_{s(\alpha)} \, ,
$$
$$
b(x_{\alpha}) = X_{b(\alpha)}
$$
index\'es par les \'el\'ements $\alpha$ de ${\rm Fl}(D)$. 

\end{listeisansmarge}
\end{remarksqed}

\bigskip

Donnons encore la liste des axiomes de la th\'eorie des relations d'\'equivalence [resp. des relations d'ordre]:

\begin{defn}\label{defV26}

La th\'eorie des relations d'\'equivalence [resp. des relations d'ordre], dont la signature est constitu\'ee d'une sorte $E$ et d'un symbole de relation binaire $\sim \ \rightarrowtail EE$ [resp. $\leq \ \rightarrowtail EE$], est d\'efinie par les trois axiomes de
\begin{enumerate}
\item[$\bullet$] r\'eflexivit\'e
$$
\top \vdash_e e \sim e \qquad \mbox{[resp.} \quad \top \vdash_e e \leq e \ \mbox{]},
$$
\item[$\bullet$] transitivit\'e
$$
e_1 \sim e_2 \wedge e_2 \sim e_3 \vdash e_1 \sim e_3
$$
[resp. \hglue 45mm $e_1 \leq e_2 \wedge e_2 \leq e_3 \vdash e_1 \leq e_3$ ],
\item[$\bullet$] sym\'etrie [resp. anti-sym\'etrie]
$$
e_1 \sim e_2 \vdash e_2 \sim e_1
$$
[resp. \hglue 37mm $e_1 \leq e_2 \wedge e_2 \leq e_1 \vdash e_1 = e_2$ ].
\end{enumerate}
\end{defn}

\bigskip

\begin{remarkqed}

Comme d\'ej\`a not\'e, la th\'eorie des ensembles a la m\^eme signature que la th\'eorie des relations d'\'equivalence ou celle des relations d'ordre. Elle est constitu\'ee d'une sorte $E$ et d'un symbole de relation binaire
$$
\xymatrix{
\in \ \ \ar@{>->}[r] &EE
}
$$
qui se lit ``appartient \`a'' ou ``est \'el\'ement de''.

\smallskip

On laisse en exercice au lecteur le soin d'\'ecrire en langage formel ses axiomes qui sont:
\begin{enumerate}
\item[$\bullet$] Axiome d'extensionalit\'e:

\smallskip

Deux ensembles sont \'egaux si et seulement si ils ont m\^emes \'el\'ements.

\smallskip

De plus, il existe un ensemble sans \'el\'ement, not\'e $\emptyset$.

\item[$\bullet$] Axiome de la paire:

\smallskip

Deux ensembles quelconques forment un nouvel ensemble dont ils sont les seuls \'el\'ements.

\item[$\bullet$] Axiome de r\'eunion:

\smallskip

On peut former la r\'eunion d'un ensemble d'ensembles.

\item[$\bullet$] Axiome de l'ensemble des parties:

\smallskip

Appelant parties d'un ensemble $A$ les ensembles $B$ dont tout \'el\'ement est \'el\'ement de $A$, on peut former l'ensemble des parties de $A$.

\item[$\bullet$] Axiome de l'infini:

\smallskip

Il existe un ensemble qui contient $\emptyset$ comme \'el\'ement et qui est stable par l'op\'eration qui consiste \`a associer \`a tout ensemble $A$ l'ensemble $A \cup \{A\}$ r\'eunion de $A$ et de l'ensemble $\{A\}$ dont l'unique \'el\'ement est $A$.

\item[$\bullet$] Axiome de remplacement:

\smallskip

L'image d'un ensemble par une application d\'efinie par une formule du premier ordre et d\'ependant de param\`etres est un ensemble d\'ependant de ces param\`etres.

\item[$\bullet$] Axiome de fondation:

\smallskip

Tout ensemble $A$ non vide a un \'el\'ement $B$ qui n'a aucun \'el\'ement commun avec $A$.

\item[$\bullet$] Axiome du choix (que l'on peut inclure ou non dans la liste):

\smallskip

Toute application surjective admet une section.
\end{enumerate}
\end{remarkqed}

\bigskip

La th\'eorie des entiers de Peano consiste \`a formaliser le principe de r\'ecurrence:

\begin{defn}\label{defV27}

La th\'eorie des entiers de Peano, dont la signature $\Sigma_N$ est constitu\'ee d'une sorte $N$, d'un symbole de constante $0 : \ \to N$ et d'un symbole de fonction unaire $\bullet + 1 : N \to N$, est d\'efinie par les axiomes suivants:
\begin{enumerate}
\item[$\bullet$] $n+1 = n'+1 \vdash n = n'$

(qui signifie que l'op\'eration de passage au successeur est injective),

\item[$\bullet$] $\top \vdash n = 0 \vee \exists \, n' (n'+1=n)$

(qui signifie que tout \'el\'ement autre que $0$ a un pr\'ed\'ecesseur),

\item[$\bullet$] $n+1=0 \vdash \perp$

(qui signifie que $0$ n'a pas de pr\'ed\'ecesseur),

\item[$\bullet$] pour toute formule $\varphi$ de $\Sigma_N$ dans le contexte d'une variable $n$ de sorte $N$, on a
$$
\varphi (0) \wedge \ \forall \, n \, (\varphi (n) \Rightarrow \varphi (n+1)) \vdash (\forall \, n ) \ \varphi (n) \, .
$$
\end{enumerate}
\end{defn}

\begin{remark}

On observe que le dernier point de la liste, qui formalise le principe de r\'ecurrence, ne consiste pas en un ou plusieurs axiomes explicites mais en un ``sch\'ema d'axiomes'': il \'enonce en fait un axiome pour chaque formule $\varphi$ de $\Sigma_N$ en une variable libre $n$ de sorte $N$.

\smallskip

De m\^eme, ``l'axiome de remplacement'' de la th\'eorie des ensembles est en fait un sch\'ema d'axiomes. \hfill $\Box$

\bigskip

On laisse au lecteur l'exercice qui consiste \`a formaliser les axiomes de d\'efinition de la th\'eorie des plans affines de signature $({\mathcal P}, {\mathcal D}, \in \ \rightarrowtail {\mathcal P}{\mathcal D}, /\!/ \rightarrowtail {\mathcal D}{\mathcal D}, O: \ \to {\mathcal P}, I : \ \to {\mathcal P}, J: \ \to {\mathcal P})$.

\smallskip

Voici ces axiomes:

\medskip

$
\left\{\begin{matrix}
\bullet &\mbox{Deux droites sont parall\`eles si et seulement si elles sont \'egales ou n'ont aucun point commun.} \hfill \\
{ \ } \\
\bullet &\mbox{Par deux points distincts $A$ et $B$ passe une unique droite (souvent not\'ee $(AB)$).} \hfill \\
{ \ } \\
\bullet &\mbox{Par un point donn\'e passe une unique droite parall\`ele \`a une droite donn\'ee.} \hfill \\
{ \ } \\
\bullet &\mbox{Les trois points marqu\'es $O,I,J$ sont non align\'es.} \hfill \\
{ \ } \\
\bullet &\mbox{(D\'esargues) Si $D_1 , D_2 , D_3$ sont trois droites distinctes qui sont parall\`eles} \hfill \\
&\mbox{[resp. se rencontrent en un point $P$] et} \hfill \\
{ \ } \\
&P_1 , P'_1 \in D_1 \, , \quad P_2 , P'_2 \in D_2 \, , \quad P_3 , P'_3 \in D_3 \\
{ \ } \\
&\mbox{sont six points [resp. distincts de $P$] tels que} \hfill \\
{ \ } \\
&(P_1 \, P_2) /\!/ (P'_1 \, P'_2) \qquad \mbox{et} \qquad (P_2 \, P_3) /\!/ (P'_2 \, P'_3) \, , \\
&\mbox{alors on a} \hfill \\
&(P_1 \, P_3) /\!/ (P'_1 \, P'_3) \, . 
\end{matrix} \right.
$

\medskip

On peut choisir d'ajouter ou non \`a cette liste l'axiome suppl\'ementaire:

\medskip

$
\left\{\begin{matrix}
\bullet &\mbox{(Pappus) Si $D$ et $D'$ sont deux droites distinctes qui se rencontrent en un point $P$ et} \hfill \\
{ \ } \\
&P_1 , P_2 , P_3 \in D \, , \qquad P'_1 , P'_2, P'_3 \in D' \\
{ \ } \\
&\mbox{sont six points distincts de $P$ tels que} \hfill \\
{ \ } \\
&(P_1 \, P'_2) /\!/ (P_2 \, P'_3) \qquad \mbox{et} \qquad (P'_1 \, P_2) /\!/ (P'_2 \, P_3) \, , \\
&\mbox{alors on a} \hfill \\
&(P_1 \, P'_1) /\!/ (P_3 \, P'_3) \, .
\end{matrix} \right.
$

\medskip

On laisse \'egalement en exercice la formalisation des axiomes de la th\'eorie des plans euclidiens, dont la signature est constitu\'ee de celle des plans affines $({\mathcal P}, {\mathcal D} , \in \ \rightarrowtail {\mathcal P}{\mathcal D}, /\!/ \rightarrowtail {\mathcal D}{\mathcal D},O,I,J)$ compl\'et\'ee par un second symbole de relation dit d'orthogonalit\'e $\perp \ \rightarrowtail {\mathcal D}{\mathcal D}$.

\smallskip

Les axiomes sont ceux de la th\'eorie des plans affines, auxquels on ajoute les suivants:

\medskip

$
\left\{\begin{matrix}
\bullet &\mbox{La relation d'orthogonalit\'e $\perp$ est sym\'etrique.} \hfill \\
{ \ } \\
\bullet &\mbox{Si deux droites sont orthogonales entre elles, une troisi\`eme droite est parall\`ele \`a l'une} \hfill \\
&\mbox{si et seulement si elle est orthogonale \`a l'autre.} \hfill \\
{ \ } \\
\bullet &\mbox{Toute droite admet au moins une droite orthogonale.} \hfill \\
{ \ } \\
\bullet &\mbox{Les droites $(OI)$ et $(OJ)$ sont orthogonales.} \hfill \\ 
{ \ } \\
\bullet &\mbox{Pour tous points $A,B,C$ non align\'es, les hauteurs du triangle $ABC$, c'est-\`a-dire les droites} \hfill \\
&\mbox{orthogonales \`a $(BC)$, $(AC)$ et $(AB)$ passant par les points $A,B,C$, ont un point commun.} \hfill
\end{matrix} \right.
$
\end{remark}

\subsection{Interpr\'etation des s\'equents et mod\`eles des th\'eories}\label{subsec524}

\medskip

De l'interpr\'etation des formules d'une signature $\Sigma$ dans les $\Sigma$-structures d'une cat\'egorie ${\mathcal C}$, qui a d\'ej\`a \'et\'e introduite et sera enti\`erement pr\'ecis\'ee au paragraphe suivant, d\'erivent l'interpr\'etation des s\'equents puis la notion de mod\`eles d'une th\'eorie dans une cat\'egorie qui poss\`ede assez de propri\'et\'es cat\'egoriques pour que les axiomes de la th\'eorie y soient interpr\'etables:

\begin{defn}\label{defV28}

Soit $\Sigma$ une signature.

\smallskip

Soit ${\mathcal C}$ une cat\'egorie localement petite qui poss\`ede des produits finis, en particulier un objet terminal.

\begin{listeimarge}

\item Soit $\varphi \vdash_{\vec x}  \psi$ un s\'equent de $\Sigma$ de contexte $\vec x = (x_1^{A_1} \cdots x_n^{A_n})$.

\smallskip

Supposons que ${\mathcal C}$ a assez de propri\'et\'es cat\'egoriques pour que les symboles qui apparaissent dans les formules $\varphi$ et $\psi$ soient interpr\'etables comme des constructions cat\'egoriques dans ${\mathcal C}$.

\smallskip

Alors on dit qu'une $\Sigma$-structure $M$ de ${\mathcal C}$ v\'erifie la propri\'et\'e
$$
\varphi \vdash_{\vec x} \psi
$$
si les interpr\'etations $M\varphi (\vec x)$ et $M\psi(\vec x)$ des formules $\varphi$ et $\psi$ dans $M$ v\'erifient la relation d'inclusion
$$
M\varphi (\vec x) \leq M\psi(\vec x)
$$
comme sous-objets de $M\!A_1 \times \cdots \times M\!A_n$.

\medskip

\item Soit ${\mathbb T}$ une th\'eorie du premier ordre de signature $\Sigma$.

\smallskip

Supposons que ${\mathcal C}$ a assez de propri\'et\'es cat\'egoriques pour que tous les symboles qui apparaissent dans les axiomes de ${\mathbb T}$ soient interpr\'etables comme des constructions cat\'egoriques dans ${\mathcal C}$.

\smallskip

Alors on dit qu'une $\Sigma$-structure $M$ de ${\mathcal C}$ est un ${\mathbb T}$-mod\`ele si $M$ v\'erifie tous les axiomes
$$
\varphi \vdash_{\vec x} \psi
$$
de la th\'eorie ${\mathbb T}$.
\end{listeimarge}
\end{defn}

\begin{remarksqed}
\begin{listeisansmarge}
\item Les propri\'et\'es exactes d'une cat\'egorie ${\mathcal C}$ dont on a besoin pour pouvoir interpr\'eter les formules $\varphi,\psi$ des axiomes $\varphi \vdash_{\vec x} \psi$ d'une th\'eorie ${\mathbb T}$ de  signature $\Sigma$ seront enti\`erement pr\'ecis\'ees au paragraphe suivant.

\smallskip

Toutefois, la remarque (ii) qui suit la d\'efinition V.2.1 permet d'annoncer quelles seront ces conditions:
\begin{enumerate}
\item[$\bullet$] Si une formule $\varphi$ ne comporte que le symbole d'\'egalit\'e $=$ ou un symbole de relation $R \rightarrowtail A_1 \cdots A_n$ en des variables $x_1^{A_1} \cdots x_n^{A_n}$, elle s'interpr\`ete sans autre hypoth\`ese sur ${\mathcal C}$ que l'existence de produits finis arbitraires.
\item[$\bullet$] Il en va de m\^eme du symbole $\top$ en des variables $x_1^{A_1} , \cdots , x_n^{A_n}$ qui s'interpr\`ete dans une $\Sigma$-structure $M$ comme le sous-objet total de $M\!A_1 \times \cdots \times M\!A_n$.
\item[$\bullet$] Si la formule $\varphi$ comporte des substitutions de variables $y^B$ par des symboles de fonctions $f(x_1^{A_1} \cdots$ $x_n^{A_n})$ en des variables $x_1^{A_1} , \cdots , x_n^{A_n}$, il faut supposer que les monomorphismes de ${\mathcal C}$ sont carrables. Pour cela, il suffit bien s\^ur que ${\mathcal C}$ admette des limites finies.
\item[$\bullet$] Si la formule $\varphi$ comporte des symboles de conjonctions finies $\wedge$, il faut faire la m\^eme hypoth\`ese que les monomorphismes de ${\mathcal C}$ sont carrables.

Si la formule $\varphi$ comporte des symboles de conjonctions infinies $\bigwedge$, il suffit de supposer que la cat\'egorie ${\mathcal C}$ poss\`ede des limites arbitraires.
\item[$\bullet$] Si la formule $\varphi$ comporte des quantificateurs existentiels $\exists$ [resp. des quantificateurs universels $\forall$], il faut supposer pour l'interpr\'eter que les monomorphismes de ${\mathcal C}$ sont carrables, si bien que tout morphisme de ${\mathcal C}$
$$
p : E_2 \longrightarrow E_1
$$
induit un foncteur
$$
p^{-1} : (S_1 \hookrightarrow E_1) \longmapsto (S_2 = S_1 \times_{E_1} E_2 \hookrightarrow E_2)
$$
de la cat\'egorie $\Omega (E_1)$ des sous-objets de $E_1$ vers la cat\'egorie $\Omega (E_2)$ des sous-objets de $E_2$, et que les foncteurs ainsi d\'efinis
$$
p^{-1} : \Omega (E_1) \longrightarrow \Omega (E_2)
$$
admettent des adjoints \`a gauche
$$
\exists_p : \Omega (E_2) \longrightarrow \Omega (E_1)
$$
[resp. des adjoints \`a droite
$$
\forall_p : \Omega (E_2) \longrightarrow \Omega (E_1) \ \mbox{]}.
$$
Pour des raisons qui apparaitront plus tard, on sera m\^eme amen\'e \`a supposer que ${\mathcal C}$ a des limites finies et que les foncteurs
$$
\exists_p \qquad \mbox{[resp.} \quad \forall_p \ \mbox{]}
$$
commutent aux changements de base au sens que, pour tout carr\'e cart\'esien de ${\mathcal C}$
$$
\xymatrix{
E'_2 \ar[d]_{p'} \ar[r]^{f_2} &E_2 \ar[d]^p \\
E'_1 \ar[r]^{f_1} &E_1
}
$$
on a
$$
\exists_{p'} \circ f_2^{-1} = f_1^{-1} \circ \exists_p
$$
[resp.
$$
\forall_{p'} \circ f_2^{-1} = f_1^{-1} \circ \forall_p \ \mbox{]}.
$$
\item[$\bullet$] Si la formule $\varphi$ comporte le signe $\perp$ ou des symboles de disjonctions finies $\vee$ [resp. infinies $\bigvee$], il faut supposer que, pour tout objet $E$ de ${\mathcal E}$, l'ensemble ordonn\'e $\Omega (E)$ de ses sous-objets admet un \'el\'ement minimal $\emptyset_E$ et un supremum de n'importe quelle famille finie [resp. infinie] d'\'el\'ements.

Il faut supposer de plus que les monomorphismes sont carrables dans ${\mathcal C}$ et que, pour tout morphisme $p : E_2 \to E_1$ de ${\mathcal C}$, le foncteur induit
$$
p^{-1} : \Omega (E_1) \longrightarrow \Omega (E_2)
$$
envoie $\emptyset_{E_1}$ sur $\emptyset_{E_2}$ et respecte le supremum de toute famille finie [resp. infinie].
\item[$\bullet$] Enfin, si la formule $\varphi$ comporte le symbole d'implication $\Rightarrow$ [resp. de n\'egation $\neg$], il faut supposer que, pour tout sous-objet $S$ d'un objet $E$ de ${\mathcal C}$, le foncteur
$$
\begin{matrix}
\Omega (E) &\longrightarrow &\Omega(E) \, , \hfill \\
\hfill S' &\longmapsto &S' \wedge S \hfill
\end{matrix}
$$
admet un adjoint \`a droite $S'' \mapsto (S \Rightarrow S'')$ [resp. que la collection des sous-objets $S'$ de $E$ tels que $S' \wedge S = \emptyset_E$ admet un plus grand \'el\'ement $\neg \, S$].

Il faut supposer de plus que les foncteurs
$$
\begin{matrix}
\hfill (S,S'') &\longmapsto &(S \Rightarrow S'') \hfill \\
\mbox{[resp.} \qquad S &\longmapsto &\neg \, S \ \mbox{]} \hfill
\end{matrix}
$$
commutent aux changements de base.
\end{enumerate}

\medskip

\item On verra que toutes les propri\'et\'es cat\'egoriques dont on a besoin pour interpr\'eter les formules du premier ordre sont v\'erifi\'ees par n'importe quel topos ${\mathcal E}$.

\smallskip

On pourra donc parler des mod\`eles de n'importe quelle th\'eorie du premier ordre ${\mathbb T}$ de signature $\Sigma$ dans n'importe quel topos.

\medskip

\item Dire qu'une $\Sigma$-structure $M$ dans une cat\'egorie ${\mathcal C}$ v\'erifie un s\'equent de la forme
$$
\top \vdash_{(x_1^{A_1} \cdots x_n^{A_n})} \varphi
$$
signifie que la formule $\varphi$ d\'efinit le sous-objet total de $M\!A_1 \times \cdots \times M\!A_n$, soit
$$
M\varphi (\vec x) = M\!A_1 \times \cdots \times M\!A_n \, .
$$

\item De m\^eme, dire qu'une $\Sigma$-structure $M$ dans une cat\'egorie ${\mathcal C}$ v\'erifie un s\'equent de la forme
$$
\varphi \vdash_{x_1^{A_1} \cdots x_n^{A_n}} \, \perp
$$
signifie que la formule $\varphi$ d\'efinit le sous-objet minimal $M\!A_1 \times \cdots \times M\!A_n$, soit
$$
M\varphi (\vec x) = \emptyset_{M\!A_1 \times \cdots \times M\!A_n} \, .
$$

\end{listeisansmarge}
\end{remarksqed}

\bigskip

On peut alors d\'efinir la cat\'egorie des mod\`eles d'une th\'eorie ${\mathbb T}$ dans n'importe quelle cat\'egorie ${\mathcal C}$ qui poss\`ede assez de propri\'et\'es pour que les axiomes de ${\mathbb T}$ y soient interpr\'etables:

\begin{defn}\label{defV29}

Soit ${\mathbb T}$ une th\'eorie du premier ordre de signature $\Sigma$.

\smallskip

Soit ${\mathcal C}$ une cat\'egorie localement petite qui admet des produits finis et poss\`ede suffisamment de propri\'et\'es cat\'egoriques pour que les formules de $\Sigma$ qui composent les axiomes de ${\mathbb T}$ soient interpr\'etables dans ${\mathcal C}$.

\smallskip

Alors on appelle cat\'egorie des mod\`eles de ${\mathbb T}$ dans ${\mathcal C}$ (ou des ${\mathbb T}$-mod\`eles de ${\mathcal C}$) et on note
$$
{\mathbb T}\mbox{-{\rm mod}} \, ({\mathcal C})
$$
la sous-cat\'egorie pleine de
$$
{\Sigma}\mbox{-{\rm str}} \, ({\mathcal C})
$$
constitu\'ee des $\Sigma$-structures de ${\mathcal C}$ qui sont des ${\mathbb T}$-mod\`eles c'est-\`a-dire v\'erifient les axiomes de ${\mathbb T}$.
\end{defn}


\begin{remarksqed}
\begin{listeisansmarge}
\item Ainsi, \'etant donn\'es deux mod\`eles $M$ et $M'$ d'une th\'eorie ${\mathbb T}$ de signature $\Sigma$ dans ${\mathcal C}$, un morphisme de mod\`eles
$$
u : M \longrightarrow N
$$
est une famille de morphismes de ${\mathcal C}$
$$
u_A : M\!A \longrightarrow N\!A
$$
index\'es par les sortes $A$ de $\Sigma$, telle que pour tout symbole de fonction $f : A_1 \cdots A_n \to B$ de $\Sigma$ le carr\'e
$$
\xymatrix{
M\!A_1 \times \cdots \times M\!A_n \ar[d] \ar[rrr]^-{u_{A_1} \times \cdots \times u_{A_n}} &&&N\!A_1 \times \cdots \times N\!A_n \ar[d] \\
M\!B \ar[rrr]^{u_B} &&&N\!B
}
$$
soit commutatif, et que pour tout symbole de relation $R \rightarrowtail A_1 \cdots A_n$ de $\Sigma$ le morphisme
$$
M\!A_1 \times \cdots \times M\!A_n \xrightarrow{ \ u_{A_1} \times \cdots \times u_{A_n} \ } N\!A_1 \times \cdots \times N\!A_n
$$
se compl\`ete en un carr\'e commutatif
$$
\xymatrix{
M\!R_{ \ } \ar@{_{(}->}[d] \ar[rr] &&N\!R_{ \ } \ar@{_{(}->}[d] \\
M\!A_1 \times \cdots \times M\!A_n \ar[rr] &&N\!A_1 \times \cdots \times N\!A_n
}
$$
n\'ecessairement unique.

\smallskip

Ainsi, les propri\'et\'es qui caract\'erisent les morphismes de ${\mathbb T}$-mod\`eles ne d\'ependent pas des axiomes de ${\mathbb T}$ mais seulement de sa signature $\Sigma$.

\medskip

\item Les propri\'et\'es requises pour la cat\'egorie ${\mathcal C}$ sont v\'erifi\'ees en particulier par la cat\'egorie Ens des ensembles.

\smallskip

Par cons\'equent, toute th\'eorie du premier ordre ${\mathbb T}$ d\'efinit une cat\'egorie
$$
{\mathbb T}\mbox{-mod} \, ({\rm Ens})
$$
des mod\`eles ensemblistes de ${\mathbb T}$.

\medskip

\item Plus g\'en\'eralement, on v\'erifiera que n'importe quel topos ${\mathcal E}$ poss\`ede les propri\'et\'es requises pour la cat\'egorie ${\mathcal C}$.

\smallskip

On pourra donc associer \`a toute th\'eorie du premier ordre ${\mathbb T}$ et \`a tout topos ${\mathcal E}$ la cat\'egorie
$$
{\mathbb T}\mbox{-mod} \, ({\mathcal E})
$$
des mod\`eles de ${\mathbb T}$ dans ${\mathcal E}$.

\medskip

\item Si ${\mathcal C}$ et ${\mathcal C}'$ sont deux cat\'egories localement petites qui admettent des produits finis et poss\`edent toutes les propri\'et\'es requises pour que les formules qui composent les axiomes d'une th\'eorie ${\mathbb T}$ de signature $\Sigma$ soient interpr\'etables dans ${\mathcal C}$ et ${\mathcal C}'$, et si
$$
F : {\mathcal C} \longrightarrow {\mathcal C}'
$$
est un foncteur qui respecte les produits finis, les monomorphismes et les structures cat\'egoriques impliqu\'ees dans l'interpr\'etation des formules qui composent les axiomes de ${\mathbb T}$, alors le foncteur
$$
F : \Sigma\mbox{-str} \, ({\mathcal C}) \longrightarrow \Sigma\mbox{-str} \, ({\mathcal C}')
$$
se restreint en un foncteur
$$
F : {\mathbb T}\mbox{-mod} \, ({\mathcal C}) \longrightarrow {\mathbb T}\mbox{-mod} \, ({\mathcal C}') \, .
$$

Si les formules des axiomes de ${\mathbb T}$ comportent des substitutions de variables ou des symboles de conjonctions fnies $\wedge$, il suffit de supposer que $F$ respecte les limites finies [resp. les limites infinies si ces formules comportent des symboles de conjonctions infinies $\bigwedge$].

\smallskip

Si ces formules comportent des quantificateurs existentiels $\exists$ [resp. des quantificateurs universels $\forall$], il faut supposer que $F$ respecte les foncteurs $\exists_p$ [resp. $\forall_p$] adjoints \`a gauche [resp. \`a droite] des foncteurs $p^{-1}$ d'images r\'eciproques des sous-objets par des morphismes $p$ de ${\mathcal C}$ ou ${\mathcal C}'$.

\smallskip

Si ces formules comportent le signe $\perp$ ou des symboles de disjonctions finies $\vee$ [resp. infinies $\bigvee$], il faut supposer que le foncteur $F$ respecte les r\'eunions finies [resp. infinies] de sous-objets. 

\smallskip

Enfin, si ces formules comportent le symbole $\Rightarrow$ ou $\neg$, il faut supposer que le foncteur $F$ respecte les foncteurs $\Rightarrow$ ou $\neg$ internes aux cat\'egories de sous-objets de ${\mathcal C}$ et ${\mathcal C}'$. 
\end{listeisansmarge}
\end{remarksqed}

\section{Les \'el\'ements de la logique du premier ordre et \\ leurs interpr\'etations}\label{sec53}

\subsection{Substitutions et relations}\label{subsec531}

\medskip

On commence par d\'efinir les notions de ``formules atomiques'' -- qui sont celles des formules d'une signature qui ne peuvent \^etre d\'ecompos\'ees en des formules plus courtes -- et de ``termes'' -- qui sont les expressions d\'eduites d'une variable d'une certaine sorte par une suite de substitutions de certaines variables par des symboles de fonctions en d'autres variables:

\begin{defn}\label{defV31}

Soit $\Sigma$ une signature.

\smallskip

Une formule atomique [resp. un terme] de $\Sigma$ est une expression formelle
$$
\varphi = \varphi (\vec x) = \varphi (x_1^{A_1} \cdots x_n^{A_n})
$$
en des variables $\vec x = (x_1^{A_1} \cdots x_n^{A_n})$ de sortes $A_ 1 , \cdots , A_n$ qui peuvent \^etre construites en un nombre fini de pas
$$
\varphi_0 , \varphi_1 , \cdots , \varphi_k = \varphi
$$
tels que
\begin{enumerate}
\item[$\bullet$] $\varphi_0$ est ou bien une formule d'\'egalit\'e
$$
(x_1^{B_1} , \cdots , x_m^{B_m}) = (y_1^{B_1} , \cdots , y_m^{B_m})
$$
entre deux familles finies de variables de m\^eme suite de sortes $B_1 , \cdots , B_m$, ou bien une formule de relation en des variables $y_1^{B_1} , \cdots , y_m^{B_m}$
$$
R (y_1^{B_1} \cdots y_m^{B_m})
$$
pour un symbole de relation $R \rightarrowtail B_1 \cdots B_m$ de $\Sigma$ [resp. $\varphi_0$ est une variable $y^B$ d'une certaine sorte $B$ de $\Sigma$],
\item[$\bullet$] chaque $\varphi_{i+1}$ est d\'eduit de $\varphi_i$ en rempla\c cant une variable $y^B$ de $\varphi_i$ par une expression en des variables $y_1^{B_1} , \cdots , y_m^{B_m}$
$$
f(y_1^{B_1} \cdots y_m^{B_m})
$$
pour un choix de symbole de fonction de $\Sigma$
$$
f : B_1 \cdots B_m \longrightarrow B \, .
$$
\end{enumerate}
\end{defn}

\bigskip

\begin{remarksqed}
\begin{listeisansmarge}
\item Tout terme $\varphi = \varphi (\vec x)$ est \`a valeur dans une certaine sorte. C'est la sorte à laquelle est affect\'ee la variable $y^B$ dont $\varphi$ est d\'eduit par une suite de substitutions de symboles de fonctions \`a des variables
$$
y^B = \varphi_0 , \varphi_1 , \cdots , \varphi_k = \varphi = \varphi (\vec x) \, .
$$

\item Les formules atomiques sont d\'eduites des formules d'\'egalit\'e
$$
(x_1^{B_1} , \cdots , x_m^{B_m}) = (y_1^{B_1} , \cdots , y_m^{B_m})
$$
ou des formules de relation
$$
R(y_1^{B_1} \cdots y_m^{B_m})
$$
en rempla\c cant certaines de leurs variables par des termes. 
\end{listeisansmarge}
\end{remarksqed}

\bigskip

L'interpr\'etation des termes dans une cat\'egorie ${\mathcal C}$ demande seulement que celle-ci poss\`ede des produits finis:

\begin{defn}\label{defV32}

Soit $\Sigma$ une signature.

\smallskip

Soit ${\mathcal C}$ une cat\'egorie localement petite qui poss\`ede des produits finis.

\smallskip

Soit $M$ une $\Sigma$-structure de ${\mathcal C}$.

\smallskip

Soit enfin $\varphi$ un terme de $\Sigma$ construit en une suite de pas
$$
\varphi_0 , \varphi_1 , \cdots, \varphi_k = \varphi \, .
$$
Alors:
\begin{enumerate}
\item[$\bullet$] Si $\varphi_0$ a la forme d'une variable $y^B$ de sorte $B$, $\varphi_0$ est interpr\'et\'ee dans $M$ comme le morphisme identit\'e
$$
M\varphi_0 = {\rm id}_{M\!B} : M\!B \longrightarrow M\!B \, .
$$
\item[$\bullet$] Si $\varphi_i = \varphi_i (z_1^{C_1} \cdots z_{\ell}^{C_{\ell}})$ est interpr\'et\'ee dans $M$ comme un morphisme de ${\mathcal C}$
$$
M\varphi_i : M\!C_1 \times \cdots \times M\!C_{\ell} \longrightarrow M\!B
$$
et que $\varphi_{i+1}$ est d\'eduite de $\varphi_i$ par une substitution
$$
z_j^{C_j} = f (y_1^{B_1} \cdots y_m^{B_m})
$$
pour un certain symbole de fonction de $\Sigma$
$$
f : B_1 \cdots B_m \longrightarrow C_j \, ,
$$
$\varphi_{i+1}$ est interpr\'et\'ee dans $M$ comme le morphisme
$$
M\varphi_{i+1} : M\!C_1\times \cdots \times M\!C_{j-1} \times M\!B_1 \times \cdots \times M\!B_m \times M\!C_{j+1} \times \cdots \times M\!C_{\ell} \longrightarrow M\!B
$$
compos\'e du morphisme
$$
{\rm id}_{M\!C_1} \times \cdots \times {\rm id}_{M\!C_{j-1}} \times M\!f \times {\rm id}_{M\!C_{j+1}} \times \cdots \times {\rm id}_{M\!C_{\ell}}
$$
d\'eduit de $M\!f : M\!B_1 \times \cdots \times M\!B_m \to M\!C_j$, avec le morphisme
$$
M\varphi_i : M\!C_1 \times \cdots \times M\!C_{\ell} \longrightarrow M\!B \, .
$$
\end{enumerate} 
\end{defn}

Les formules atomiques s'interpr\`etent de m\^eme dans une cat\'egorie ${\mathcal C}$ qui poss\`ede des limites finies:

\begin{defn}\label{defV33}

Soit $\Sigma$ une signature.

\smallskip

Soit ${\mathcal C}$ une cat\'egorie localement petite qui poss\`ede des limites finies.

\smallskip

Soit $M$ une $\Sigma$-structure de ${\mathcal C}$.

\smallskip

Soit $\varphi$ une formule atomique de ${\mathcal C}$ construite en une suite de pas
$$
\varphi_0 , \varphi_1 , \cdots , \varphi_k = \varphi \, .
$$

Alors:
\begin{enumerate}
\item[$\bullet$] Si $\varphi_0$ est une formule d'\'egalit\'e [resp. de relation]
$$
(x_1^{B_1} , \cdots , x_m^{B_m}) = (y_1^{B_1} , \cdots , y_m^{B_m})
$$
[resp. 
$$
R (y_1^{B_1} \cdots y_m^{B_m}) \mbox{]} 
$$
pour deux familles finies de variables affect\'ees \`a une m\^eme suite de sortes $B_1 \cdots B_m$ [resp. pour un symbole de relation $R \rightarrowtail B_1 \cdots B_m$], $\varphi_0$ s'interpr\`ete dans $M$ comme le sous-objet diagonal
$$
M\!B_1 \times \cdots \times M\!B_m \xhookrightarrow{ \ { \ } \ } M\!B_1 \times \cdots \times M\!B_m \times M\!B_1 \times \cdots \times M\!B_m
$$
[resp. comme le sous-objet
$$
M\!R \xhookrightarrow{ \ { \ } \ } M\!B_1 \times \cdots \times M\!B_m \ \mbox{].}
$$
\item[$\bullet$] Si $\varphi_i = \varphi_i (z_1^{C_1} \, z_2^{C_2}  \cdots  z_{\ell}^{C_{\ell}})$ est interpr\'et\'ee dans $M$ comme un sous-objet
$$
M\varphi_i \xhookrightarrow{ \ { \ } \ } M\!C_1 \times M\!C_2 \times \cdots \times M\!C_{\ell} \, ,
$$
et $\varphi_{i+1}$ est d\'eduite de $\varphi_i$ par une substitution
$$
z_j^{C_j} = f (y_1^{B_1} \cdots y_m^{B_m})
$$
pour un symbole de fonction $f : B_1 \cdots B_m \to C_j$, $\varphi_{i+1}$ s'interpr\`ete dans $M$ comme le sous-objet
$$
M\varphi_{i+1} \xhookrightarrow{ \ { \ } \ } M\!C_1 \times \cdots \times M\!C_{j-1} \times M\!B_1 \times \cdots \times M\!B_m \times M\!C_{j+1} \times \cdots \times M\!C_{\ell}
$$
d\'eduit de
$$
M\varphi_i \xhookrightarrow{ \ { \ } \ } M\!C_1 \times \cdots \times M\!C_{\ell}
$$
par produit fibr\'e avec le morphisme
$$
{\rm id}_{M\!C_1} \times \cdots \times {\rm id}_{M\!C_{j-1}} \times M\!f \times {\rm id}_{M\!C_{j+1}} \times \cdots \times {\rm id}_{M\!C_{\ell}}
$$
d\'eduit du morphisme
$$
M\!f : M\!B_1 \times \cdots \times M\!B_m \longrightarrow M\!C_j \, .
$$
\end{enumerate}
\end{defn}

\bigskip

\begin{remarkqed}

De la m\^eme fa\c con, si $\varphi = \varphi (x_{\sigma (1)}^{A_{\sigma (1)}} , \cdots , x_{\sigma (m)}^{A_{\sigma (m)}})$ est une formule atomique dans un contexte $\vec x = (x_1^{A_1} \cdots x_n^{A_n})$ pour une application injective $\sigma : \{1,\cdots ,m\} \to \{1,\cdots ,n\}$, alors $\varphi (\vec x)$ s'interpr\`ete comme le sous-objet
$$
M\varphi (\vec x) \xhookrightarrow{ \ { \ } \ }  M\!A_1 \times \cdots \times M\!A_n
$$
d\'eduit du sous-objet
$$
M\varphi \xhookrightarrow{ \ { \ } \ } M\!A_{\sigma(1)} \times \cdots \times M\!A_{\sigma (m)}
$$
par changement de base avec le morphisme de projection d\'efini par $\sigma$
$$
M\!A_1 \times \cdots \times M\!A_n \longrightarrow M\!A_{\sigma(1)} \times \cdots \times M\!A_{\sigma (m)} \, .
$$
\end{remarkqed}

On a le lemme imm\'ediat suivant:

\begin{lem}\label{lemV34}

Soit $\Sigma$ une signature.

\smallskip

Soient ${\mathcal C},{\mathcal C}'$ deux cat\'egories localement petites qui poss\`edent des limites finies, reli\'ees par un foncteur
$$
F : {\mathcal C} \longrightarrow {\mathcal C}'
$$
qui respecte les limites finies.

\smallskip

Alors le foncteur induit
$$
F : \Sigma\mbox{\rm -str} ({\mathcal C}) \longrightarrow \Sigma\mbox{\rm -str} ({\mathcal C}')
$$
respecte les interpr\'etations des formules atomiques. 

\smallskip

En d'autres termes, pour toute formule atomique $\varphi$ de $\Sigma$ de contexte $\vec x = (x_1^{A_1} \cdots x_n^{A_n})$ et pour toute $\Sigma$-structure $M$ de ${\mathcal C}$ transform\'ee par $F$ en une $\Sigma$-structure $F(M)$ de ${\mathcal C}'$, l'interpr\'etation
$$
F(M) \, \varphi (\vec x) \xhookrightarrow{ \ { \ } \ } F(M\!A_1) \times \cdots \times F(M\!A_n)
$$
est la transform\'ee par $F$ de l'interpr\'etation
$$
M\varphi (\vec x) \xhookrightarrow{ \ { \ } \ } M\!A_1 \times \cdots \times M\!A_n \, .
$$
\end{lem}

\begin{demo}

Si $\varphi$ est une formule atomique arbitraire, cela r\'esulte de ce que par hypoth\`ese $F$ respecte les limites finies, donc en particulier les monomorphismes et les produits fibr\'es. 
\end{demo}

\subsection{Les conjonctions et les disjonctions}\label{subsec532}

\medskip

D\'efinir formellement une formule par conjonction ou disjonction de formules plus petites est possible d\`es lors que ces formules plus petites partagent le m\^eme contexte:

\begin{defn}\label{defV35}

Soit $\Sigma$ une signature.

\smallskip

Si $\varphi_1 , \cdots , \varphi_k$ [resp. $(\varphi_i)_{i \in I}$] est une famille finie [resp. infinie] de formules de $\Sigma$ de m\^eme contexte $\vec x = (x_1^{A_1} \cdots x_n^{A_n})$, leur conjonction est l'expression formelle
$$
\varphi_1 \wedge \cdots \wedge \varphi_k \qquad \mbox{[resp.} \quad \bigwedge_{i \in I} \varphi_i \ \mbox{]}
$$
et leur disjonction est l'expression formelle
$$
\varphi_1 \vee \cdots \vee \varphi_k \qquad \mbox{[resp.} \quad \bigvee_{i \in I} \varphi_i \ \mbox{]}.
$$
\end{defn}

\begin{remarkqed}

Ainsi, les conjonctions et disjonctions apparaissent ici comme de pures \'ecritures formelles.

\smallskip

Elles acqui\`erent un sens par leurs interpr\'etations dans des cat\'egories. 

\end{remarkqed}

\bigskip

Pour que les conjonctions et disjonctions de formules puissent \^etre interpr\'et\'ees dans une cat\'egorie ${\mathcal C}$, il faut que dans ${\mathcal C}$ les intersections et r\'eunions de sous-objets soient bien d\'efinies au sens suivant:

\begin{defn}\label{defV36}

Soit ${\mathcal C}$ une cat\'egorie qui poss\`ede des limites finies.

\smallskip

Alors:

\begin{listeimarge}

\item L'intersection d'une famille de sous-objets d'un objet $E$
$$
S_1 \xhookrightarrow{ \ { \ } \ } E , \cdots , S_k \xhookrightarrow{ \ { \ } \ }  E
$$
est d\'efinie comme la limite dans ${\mathcal C}$
$$
S_1 \times_E \cdots \times_E S_k = S_1 \wedge \cdots \wedge S_k \xhookrightarrow{ \ { \ } \ } E \, .
$$
Elle est caract\'eris\'ee par l'\'equivalence, pour tout sous-objet $S'$ de $E$,
$$
S' \leq S_1 \wedge \cdots \wedge S_k \Longleftrightarrow S' \leq S_i \, , \quad 1 \leq i \leq k \, .
$$

\item On dit que ${\mathcal C}$ admet des intersections arbitraires de sous-objets si, pour toute famille de sous-objets d'un objet $E$
$$
(S_i  \xhookrightarrow{ \ { \ } \ } E)_{i \in I} \, ,
$$
le diagramme des monomorphismes $S_i \hookrightarrow E$ poss\`ede une limite dans ${\mathcal C}$
$$
\bigwedge_{i \in I} S_i \xhookrightarrow{ \ { \ } \ } E
$$
alors caract\'eris\'ee par l'\'equivalence, pour tout sous-objet $S'$ de $S$,
$$
S' \leq \bigwedge_{i \in I} S_i \Longleftrightarrow S' \leq S_i \, , \quad \forall \, i \in I \, .
$$

\noindent (iii) On dit que ${\mathcal C}$ admet des r\'eunions finies [resp. infinies] de sous-objets si

\medskip

$
\left\lmoustache \begin{matrix}
\bullet &\mbox{pour toute famille finie [resp. infinie] de sous-objets d'un objet $E$} \hfill \\
{ \ } \\
&(S_i \xhookrightarrow{ \ { \ } \ }  E)_{1 \leq i \leq k} \qquad \mbox{[resp.} \quad (S_i \xhookrightarrow{ \ { \ } \ }  E)_{i \in I} \ \mbox{],} \\
&\mbox{il existe un sous-objet} \hfill \\
&S_1 \vee \cdots \vee S_k \xhookrightarrow{ \ { \ } \ }  E \quad \mbox{[resp.} \quad \displaystyle\bigvee_{i \in I} S_i \xhookrightarrow{ \ { \ } \ }  E \ \mbox{]} \\
{ \ } \\
&\mbox{caract\'eris\'e par l'\'equivalence, pour tout sous-objet $S' \hookrightarrow E$,} \hfill \\
{ \ } \\
&\qquad S' \geq S_1 \vee \cdots \vee S_k \Longleftrightarrow S' \geq S_i \, , \quad 1 \leq i \leq k \\
{ \ } \\
&\mbox{[resp.} \qquad \ S' \geq \displaystyle\bigvee_{i \in I} S_i \Longleftrightarrow S' \geq S_i \, , \quad \forall \, i \in I \ \mbox{],} \\
{ \ } \\
\end{matrix} \right.
$

$
\left\rmoustache \begin{matrix}
\bullet &\mbox{en particulier, tout objet $E$ de ${\mathcal C}$ a un plus petit sous-objet $\emptyset_E$,} \hfill \\
{ \ } \\
\bullet &\mbox{pour toute telle famille et pour tout morphisme $f : E' \to E$ de ${\mathcal C}$, le monomorphisme canonique} \hfill \\
{ \ } \\
& \ f^{-1} S_1 \vee \cdots \vee f^{-1} S_k \xhookrightarrow{ \ { \ } \ } f^{-1} (S_1 \vee \cdots \vee S_k) \\
{ \ } \\
&\mbox{[resp.} \quad \displaystyle\bigvee_{i \in I} f^{-1} S_i \xhookrightarrow{ \ { \ } \ }  f^{-1} \biggl( \displaystyle\bigvee_{i \in I} S_i \biggl) \ \mbox{]} \\
{ \ } \\
&\mbox{(o\`u l'on note $f^{-1} S' = S' \times_E E'$ pour tout sous-objet $S' \hookrightarrow E$) est un isomorphisme,} \hfill \\
{ \ } \\
\bullet &\mbox{en particulier, pour tout morphisme $f : E' \to E$, on a} \hfill \\
{ \ } \\
&\emptyset_{E'} = f^{-1} \, \emptyset_E \, .
\end{matrix} \right.
$
\end{listeimarge}
\end{defn}

\bigskip

\begin{remarksqed}
\begin{listeisansmarge}
\item Pour toute famille de sous-objets $(S_i \hookrightarrow E)_{i \in I}$ telle que $\underset{i \in I}{\bigwedge} \, S_i$ est bien d\'efinie, et pour tout morphisme $f : E' \to E$, on a
$$
f^{-1} \left( \bigwedge_{i \in I} S_i \right) = \bigwedge_{i \in I} f^{-1} S_i \, .
$$
Cela r\'esulte en effet de ce que les limites respectent toujours les limites.

\medskip

\item Soient ${\mathcal C},{\mathcal C}'$ deux cat\'egories localement petites qui poss\`edent des limites finies, reli\'ees par un foncteur
$$
F : {\mathcal C} \longrightarrow {\mathcal C}'
$$
qui respecte les limites finies.

\smallskip

Alors $F$ respecte les intersections finies au sens que pour toute famille finie de sous-objets $S_i \hookrightarrow E$, $1 \leq i \leq k$, d'un objet $E$ de ${\mathcal C}$, le morphisme
$$
F(S_1 \wedge \cdots \wedge S_k)  \xhookrightarrow{ \ { \ } \ } F(S_1) \wedge \cdots \wedge F(S_k)
$$
est un isomorphisme.

\smallskip

Dans les m\^emes conditions, on dit que $F$ respecte les r\'eunions finies si, pour toute telle famille $(S_i \hookrightarrow E)_{1 \leq i \leq k}$, le morphisme
$$
F(S_1) \vee \cdots \vee F(S_k)  \xhookrightarrow{ \ { \ } \ } F(S_1 \vee \cdots \vee S_k)
$$
est un isomorphisme.

\medskip

\item On dit qu'un tel foncteur respecte les intersections arbitraires [resp. les r\'eunions arbitraires] si, pour toute famille de sous-objets d'un objet $E$ de ${\mathcal C}$
$$
(S_i  \xhookrightarrow{ \ { \ } \ } E)_{i \in I} \, ,
$$
le morphisme canonique
$$
\ \ \ \ \ F \left( \bigwedge_{i \in I} S_i \right)  \xhookrightarrow{ \ { \ } \ }  \bigwedge_{i \in I} F(S_i)
$$
$$
\mbox{[resp.} \qquad \bigvee_{i \in I} F(S_i) \xhookrightarrow{ \ { \ } \ } F \left( \bigvee_{i \in I} S_i \right) \ \mbox{]}
$$
est un isomorphisme. 

\end{listeisansmarge}
\end{remarksqed}

\bigskip

On observe:

\begin{lem}\label{lemV37}
\begin{listeimarge}
\item Tout topos ${\mathcal E}$ admet des intersections arbitraires et des r\'eunions arbitraires de sous-objets au sens de la d\'efinition \ref{defV36}.

\medskip

\item Pour tout morphisme de topos
$$
f = (f^* , f_*) : {\mathcal E}' \longrightarrow {\mathcal E} \, ,
$$
sa composante d'image r\'eciproque
$$
f^* : {\mathcal E} \longrightarrow {\mathcal E}'
$$
respecte les intersections finies et les r\'eunions arbitraires de sous-objets, et sa composante d'image directe
$$
f_* : {\mathcal E}' \longrightarrow {\mathcal E}
$$
respecte les intersections arbitraires de sous-objets.
\end{listeimarge}
\end{lem}

\begin{demo}
\begin{listeisansmarge}
\item Un topos ${\mathcal E}$ poss\`ede des intersections arbitraires car toutes les limites sont bien d\'efinies dans ${\mathcal E}$.

\smallskip

On a vu dans la proposition \ref{propIII62} (i) que toute famille de sous-objets d'un objet $E$ de ${\mathcal E}$
$$
(S_i \xhookrightarrow{ \ { \ } \ } E)_{i \in I}
$$
a une r\'eunion
$$
\bigvee_{i \in I} S_i
$$
qui est l'intersection de tous les sous-objets de $E$ contenant tous les $S_i$.

\smallskip

De plus, $\underset{i \in I}{\bigvee} \, S_i$ est \'egalement la colimite du diagramme constitu\'e des $S_i$, des $S_i \times_E S_j$ et des projections $S_i \leftarrow S_i \times_E S_j \rightarrow S_j$.

\smallskip

Pour tout morphisme $u : E' \to E$ de ${\mathcal E}$, le foncteur de changement de base $E' \times_E \bullet$ respecte \`a la fois les limites et les colimites, ce qui implique
$$
u^{-1} \left( \bigvee_{i \in I} S_i \right) = \bigvee_{i \in I} u^{-1} S_i \, .
$$

\item r\'esulte de ce que $f^*$ respecte \`a la fois les colimites arbitraires et les limites finies et que $f_*$ respecte les limites arbitraires. 

\end{listeisansmarge}
\end{demo}

\bigskip

Voici comment les conjonctions et disjonctions de formules sont interpr\'et\'ees dans des cat\'egories qui admettent des intersections ou r\'eunions de sous-objets:

\begin{defn}\label{defV38}

Soit $\Sigma$ une signature.

\smallskip

Soit ${\mathcal C}$ une cat\'egorie localement petite qui poss\`ede des limites finies.

\smallskip

Soit $(\varphi_i)_{i \in I}$ une famille de formules de $\Sigma$ de m\^eme contexte $\vec x = (x_1^{A_1} \cdots x_n^{A_n})$ qui sont suppos\'ees interpr\'etables dans la cat\'egorie ${\mathcal C}$.

\smallskip

Enfin, soit $M$ une $\Sigma$-structure de ${\mathcal C}$, avec les interpr\'etations des $\varphi_i (\vec x)$ comme des sous-objets
$$
M\varphi_i (\vec x) \xhookrightarrow{ \ { \ } \ } M\!A_1 \times \cdots \times M\!A_n \, .
$$
Alors:

\begin{listeimarge}
	
\item Si $I = \{1,\cdots , k \}$, la formule $\varphi_1 \wedge \cdots \wedge \varphi_k$ est interpr\'et\'ee comme l'intersection finie
$$
M(\varphi_1 \wedge \cdots \wedge \varphi_k) (\vec x) = M\varphi_1 (\vec x) \wedge \cdots \wedge M \varphi_k (\vec x) \, .
$$
Si $I$ est infini et ${\mathcal C}$ admet des intersections arbitraires de sous-objets, la formule $\underset{i \in I}{\bigwedge} \, \varphi_i$ est interpr\'et\'ee comme l'intersection
$$
M \left( \bigwedge_{i \in I} \varphi_i \right)(\vec x) = \bigwedge_{i \in I} M \varphi_i (\vec x) \, .
$$

\item Si ${\mathcal C}$ admet des r\'eunions finies [resp. arbitraires] de sous-objets, et si $I = \{ 1,\cdots ,k\}$ [resp. $I$ est arbitraire], la formule $\varphi_1 \vee \cdots \vee \varphi_k$ [resp. $\underset{i \in I}{\bigvee} \, \varphi_i$] s'interpr\`ete comme la r\'eunion
$$
\qquad \qquad \qquad \ \ \ M (\varphi_1 \vee \cdots \vee \varphi_k)(\vec x) = M\varphi_1 (\vec x) \vee \cdots \vee M\varphi_k (\vec x)
$$
$$
\mbox{[resp.} \qquad M \left( \bigvee_{i \in I} \varphi_i\right)(\vec x) = \bigvee_{i \in I} M\varphi_i (\vec x) \ \mbox{].}
$$
\end{listeimarge}
\end{defn}

\begin{remarksqed}
\begin{listeisansmarge}
\item On a d\'ej\`a annonc\'e que le symbole $\top$ du vrai dans le contexte $\vec x = (x_1^{A_1} \cdots x_n^{A_n})$ appliqu\'e \`a une $\Sigma$-structure $M$ de ${\mathcal C}$ s'interpr\`ete comme le sous-objet total
$$
M\top (\vec x) = M\!A_1 \times \cdots \times M\!A_n \, .
$$

Si les r\'eunions finies de sous-objets des objets $E$ de ${\mathcal C}$ sont bien d\'efinies et que, en particulier, chaque $E$ a un plus petit sous-objet $\emptyset_E$, le symbole $\perp$ du faux dans le contexte $\vec x$ appliqu\'e \`a $M$ s'interpr\`ete comme
$$
M\!\!\perp (\vec x) = \emptyset_{M\!A_1 \times \cdots \times M\!A_n} \, .
$$

\item D'apr\`es le lemme \ref{lemV37} (i), les conjonctions et disjonctions arbitraires de formules interpr\'etables dans un topos y sont encore interpr\'etables. 
\end{listeisansmarge}
\end{remarksqed}

\bigskip

On d\'eduit de cette d\'efinition:

\begin{lem}\label{lemV39}

Soit $\Sigma$ une signature.

\smallskip

Soient ${\mathcal C}, {\mathcal C}'$ deux cat\'egories localement petites qui poss\`edent des limites finies, reli\'ees par un foncteur
$$
F : {\mathcal C} \longrightarrow {\mathcal C}'
$$
qui respecte les limites finies.

\smallskip

Soit $(\varphi_i)_{i \in I}$ une famille de formules de $\Sigma$ de contexte $\vec x = (x_1^{A_1} \cdots x_n^{A_n})$ telles que le foncteur
$$
F : \Sigma\mbox{\rm -str} ({\mathcal C}) \longrightarrow \Sigma\mbox{\rm -str}({\mathcal C}')
$$
respecte l'interpr\'etation de chaque formule $\varphi_i$, $i \in I$, dans n'importe quelle $\Sigma$-structure $M$ de ${\mathcal C}$.

\smallskip

Alors:

\begin{listeimarge}

\item Si $I = \{1,\cdots ,k\}$, le foncteur $F$ respecte les interpr\'etations de $\varphi_1 \wedge \cdots \wedge \varphi_k$ au sens que, pour toute $\Sigma$-structure $M$ de ${\mathcal C}$, il transforme le sous-objet
$$
M (\varphi_1 \wedge \cdots \wedge \varphi_k) (\vec x) \xhookrightarrow{ \ { \ } \ } M\!A_1 \times \cdots \times M\!A_n
$$
en le sous-objet
$$
F(M) (\varphi_1 \wedge \cdots \wedge \varphi_k) (\vec x) \xhookrightarrow{ \ { \ } \ } F(M\!A_1) \times \cdots \times F(M\!A_n) \, .
$$

\item Si $I = \{1,\cdots ,k\}$ et $F$ respecte les r\'eunions finies de sous-objets, il respecte les interpr\'etations de $\varphi_1 \vee \cdots \vee \varphi_k$ au sens que, pour toute $M$, il transforme le sous-objet
$$
M(\varphi_1 \vee \cdots \vee \varphi_k) (\vec x) \xhookrightarrow{ \ { \ } \ } M\!A_1 \times \cdots \times M\!A_n
$$
en le sous-objet
$$
F(M) (\varphi_1 \vee \cdots \vee \varphi_k) (\vec x) \xhookrightarrow{ \ { \ } \ } M\!A_1 \times \cdots \times M\!A_n \, .
$$

\item Si $F$ respecte les intersections [resp. les r\'eunions] arbitraires, il respecte les interpr\'etations de $\underset{i \in I}{\bigwedge} \, \varphi_i$ [resp. $\underset{i \in I}{\bigvee} \, \varphi_i$] au sens que, pour toute $M$, il transforme le sous-objet
$$
M \left( \bigwedge_{i \in I} \varphi_i \right)(\vec x)  \qquad \mbox{[resp.} \quad M \left( \bigvee_{i \in I} \varphi_i \right)(\vec x)\mbox{]} \quad \mbox{de} \quad M\!A_1 \times \cdots \times M\!A_n
$$
en le sous-objet
$$
F(M) \left( \bigwedge_{i \in I} \varphi_i \right)(\vec x) \qquad \mbox{[resp.} \quad F(M) \left( \bigvee_{i \in I} \varphi_i \right)(\vec x)\mbox{]} \quad \mbox{de} \quad F(M\!A_1) \times \cdots \times F(M\!A_n) \, .
$$
\end{listeimarge}
\end{lem}

\pagebreak 

\begin{remarks}
\begin{listeisansmarge}
\item Le foncteur $F$ transforme toujours les interpr\'etations
$$
M\top (\vec x) = M\!A_1 \times \cdots \times M\!A_n
$$
en les interpr\'etations
$$
F(M) \top (\vec x) = F(M\!A_1) \times \cdots \times F(M\!A_n) \, .
$$

Si $F$ respecte les r\'eunions finies, il envoie en particulier chaque $\emptyset_E$ sur $\emptyset_{F(E)}$, $E \in {\rm Ob} ({\mathcal C})$. Cela implique qu'il transforme les interpr\'etations
$$
M\!\!\perp (\vec x) = \emptyset_{M\!A_1 \times \cdots \times M\!A_n}
$$
en les interpr\'etations
$$
F(M)\!\perp (\vec x) = \emptyset_{F(M\!A_1) \times \cdots \times F(M\!A_n)} \, .
$$

\item Il r\'esulte du lemme \ref{lemV37} (ii) que pour tout morphisme de topos
$$
f = (f^* , f_*) : {\mathcal E}' \longrightarrow {\mathcal E} \, ,
$$
le foncteur $f^*$ respecte les interpr\'etations des disjonctions arbitraires et des conjonctions finies de formules dont il respecte les interpr\'etations, tandis que le foncteur $f_*$ respecte les interpr\'etations de conjonctions arbitraires de telles formules.
\end{listeisansmarge}
\end{remarks}

\bigskip

\begin{demo}

C'est \'evident puisque, par d\'efinition, les conjonctions sont interpr\'et\'ees comme des intersections et les disjonctions comme des r\'eunions. \end{demo}

\subsection{Les quantificateurs existentiels et universels}\label{subsec533}

\medskip

En appliquant le quantificateur $\exists$ ou $\forall$ \`a une partie des variables d'une formule $\varphi$, on d\'efinit une nouvelle formule:

\begin{defn}\label{defV310}

Soit $\Sigma$ une signature.

\smallskip

Soit $\varphi = \varphi (\vec x , \vec y)$ une formule de $\Sigma$ dont le contexte $(\vec x , \vec y)$ est \'ecrit comme la r\'eunion de deux contextes disjoints $\vec x = (x_1^{A_1} \cdots x_n^{A_n})$ et $\vec y = (y_1^{B_1} \cdots y_m^{B_m})$.

\smallskip

Alors l'expression formelle
$$
(\exists \, \vec y) \, \varphi \qquad \mbox{[resp.} \quad (\forall \, \vec y) \, \varphi \ \mbox{]}
$$
est consid\'er\'ee comme une formule de contexte $\vec x$ et appel\'ee la formule d\'eduite de $\varphi$ par le quantificateur existentiel $\exists$ [resp. universel $\forall$] en les variables de $\vec y$.
\end{defn}

\begin{remarkqed}

Ici encore, les expressions $(\exists \, \vec y) \, \varphi$ ou $(\forall \, \vec y) \, \varphi$ sont introduites comme des \'ecritures purement formelles.

\smallskip

Elles acqui\`erent un sens par leurs interpr\'etations dans des cat\'egories. \end{remarkqed}

\bigskip

Rappelons que si ${\mathcal C}$ est une cat\'egorie localement petite qui admet des limites finies, tout morphisme de ${\mathcal C}$
$$
e : E_2 \longrightarrow E_1
$$
induit un foncteur d'images r\'eciproques
$$
\begin{matrix}
e^{-1} &: &\hfill \Omega (E_1) &\longrightarrow &\Omega (E_2) \, , \hfill \\
&&(S_1 \hookrightarrow E_1) &\longmapsto &(S_1 \times_{E_1} E_2 \hookrightarrow E_2)
\end{matrix}
$$
de la cat\'egorie $\Omega (E_1)$ des sous-objets de $E_1$ vers celle $\Omega (E_2)$ des sous-objets de $E_2$.

\smallskip

Les cat\'egories dans lesquelles les symboles $\exists$ et $\forall$ sont interpr\'etables sont celles dans lesquelles ces foncteurs d'images r\'eciproques admettent des adjoints \`a gauche [resp. \`a droite]:

\begin{defn}\label{defV311}
\begin{listeimarge}
\item Soit ${\mathcal C}$ une cat\'egorie localement petite qui poss\`ede des limites finies.

\smallskip

On dit que ${\mathcal C}$  ``permet les quantifications existentielles'' [resp. qu'elle ``permet les quantifications universelles''] si
\begin{enumerate}
\item[$\bullet$] pour tout morphisme de ${\mathcal C}$
$$
e : E_2 \longrightarrow E_1 \, ,
$$
le foncteur associ\'e d'images r\'eciproques des sous-objets
$$
e^{-1} : \Omega (E_1) \longrightarrow \Omega (E_2)
$$
admet un adjoint \`a gauche [resp. \`a droite]
$$
\qquad \ \ \ \ \exists_e : \Omega (E_2) \longrightarrow \Omega (E_1)
$$
$$
\mbox{[resp.} \qquad \forall_e : \Omega (E_2) \longrightarrow \Omega (E_1) \ \mbox{]},
$$
\item[$\bullet$] les foncteurs $\exists_e$ [resp. $\forall_e$] sont compatibles avec les changements de base au sens que, pour tout carr\'e cart\'esien de ${\mathcal C}$
$$
\xymatrix{
E'_2 \ar[d]_{e'} \ar[r]^{g_2} &E_2 \ar[d]^e \\
E'_1 \ar[r]^{g_1} &E_1
}
$$
la transformation naturelle canonique
$$
\qquad \ \ \exists_{e'} \circ g_2^{-1} \longrightarrow g_1^{-1} \circ \exists_e
$$
$$
\mbox{[resp.} \qquad g_1^{-1} \circ \forall_e \longrightarrow \forall_{e'} \circ g_2^{-1} \ \mbox{]}
$$
est un isomorphisme.

\end{enumerate}

\item Soit
$$
F : {\mathcal C} \longrightarrow {\mathcal C}'
$$
un foncteur entre deux cat\'egories localement petites ${\mathcal C}$ et ${\mathcal C}'$ qui permettent les quantifications existentielles [resp. permettent les quantifications universelles].

\smallskip

On dit que le foncteur $F$ respecte les quantifications existentielles [resp. qu'il respecte les quantifications universelles] s'il respecte les limites finies et si, pour tout morphisme de ${\mathcal C}$
$$
e : E_2 \longrightarrow E_1 \, ,
$$
le carr\'e
$$
\xymatrix{
\Omega (E_2) \ar[d]_{\exists_e} \ar[r]^-{F} &\Omega (F(E_2)) \ar[d]^{\exists_{F(e)}}\\
\Omega (E_1) \ar[r]^-{F} &\Omega (F(E_1)) 
} \qquad \qquad 
\xymatrix{
\Omega (E_2) \ar[d]_{\mbox{[resp.} \quad \forall_e} \ar[r]^-{F} &\Omega (F(E_2)) \ar[d]^{\forall_{F(e) \ \mbox{]}}}\\
\Omega (E_1) \ar[r]^-{F} &\Omega (F(E_1)) 
}
$$
est commutatif.

\end{listeimarge}
\end{defn}

\begin{remarksqed}
\begin{listeisansmarge}
\item Pour tout carr\'e cart\'esien comme dans (i)
$$
\xymatrix{
E'_2 \ar[d]_{e'} \ar[r]^{g_2} &E_2 \ar[d]^e \\
E'_1 \ar[r]^{g_1} &E_1
}
$$
on a une \'egalit\'e
$$
e'^{-1} \circ g_1^{-1} = g_2^{-1} \circ e^{-1}
$$
donc le morphisme d'adjonction
$$
{\rm id} \longrightarrow e^{-1} \circ \exists_e
$$
[resp.
$$
e^{-1} \circ \forall_e \longrightarrow {\rm id} \ \mbox{]}
$$
induit un morphisme
$$
g_2^{-1} \longrightarrow g_2^{-1} \circ e^{-1} \circ \exists_e = e'^{-1} \circ g_1^{-1} \circ \exists_e
$$
[resp.
$$
e'^{-1} \circ g_1^{-1} \circ \forall_e = g_2^{-1} \circ e^{-1} \circ \forall_e \longrightarrow g_2^{-1} \ \mbox{]}
$$
qui correspond par adjonction \`a un morphisme
$$
\exists_{e'} \circ g_2^{-1} \longrightarrow g_1^{-1} \circ \exists_e
$$
[resp.
$$
g_1^{-1} \circ \forall_e \longrightarrow \forall_{e'} \circ g_2^{-1} \ \mbox{]}.
$$

C'est la transformation naturelle dont il est question dans la seconde partie de (i).

\medskip

\item De m\^eme, si ${\mathcal C} , {\mathcal C}'$ sont deux cat\'egories qui permettent les quantifications existentielles [resp. qui permettent les quantifications universelles] et
$$
F : {\mathcal C} \longrightarrow {\mathcal C}'
$$
est un foncteur qui respecte les limites finies, on a une \'egalit\'e
$$
F \circ e^{-1} = F(e)^{-1} \circ F \, .
$$
Alors le morphisme d'adjonction
$$
{\rm id} \longrightarrow e^{-1} \circ \exists_e \qquad \mbox{[resp.} \quad e^{-1} \circ \forall_e \longrightarrow {\rm id} \ \mbox{]}
$$
induit un morphisme
$$
F \longrightarrow F \circ e^{-1} \circ \exists_e \qquad \mbox{[resp.} \quad F \circ e^{-1} \circ \forall_e \longrightarrow F \ \mbox{]}
$$
qui s'\'ecrit encore
$$
F \longrightarrow F (e)^{-1} \circ F \circ \exists_e \qquad \!\!\!\!\! \mbox{[resp.} \quad F(e)^{-1} \circ F \circ \forall_e \longrightarrow F \ \mbox{]}
$$
et correspond par adjonction \`a un morphisme
$$
\exists_{F(e)} \circ F \longrightarrow F \circ \exists_e \qquad \mbox{[resp.} \quad F \circ \forall_e \longrightarrow \forall_{F(e)} \circ F \ \mbox{]}.
$$

Le foncteur $F$ respecte les quantifications existentielles [resp. respecte les quantifications universelles] si cette transformation naturelle est un isomorphisme.

\medskip

\item Si ${\mathcal C}$ est la cat\'egorie des ensembles et
$$
f : E_2 \longrightarrow E_1
$$
est donc une application entre deux ensembles $E_2,E_1$, alors pour tout sous-ensemble
$$
S_2 \xhookrightarrow{ \ { \ } \ } E_2 \, ,
$$
le sous-ensemble
$$
\exists_f \, S_2 \xhookrightarrow{ \ { \ } \ } E_1 \qquad \mbox{[resp.} \quad \forall_f \, S_2 \xhookrightarrow{ \ { \ } \ } E_1 \ \mbox{]}
$$
est constitu\'e des \'el\'ements
$$
e_1 \in E_1
$$
dont la fibre
$$
f^{-1} (e_1) = \{ e_2 \in E_2 \mid f(e_2) = e_1 \}
$$
contient au moins un \'el\'ement de $S_1$ [resp. est enti\`erement contenue dans $S_1$].

\smallskip

Ce cas particulier justifie la d\'efinition g\'en\'erale.

\medskip

\item Si ${\mathcal C}$ est une cat\'egorie qui permet les quantifications existentielles [resp. qui permet les quantifications universelles] et
$$
e : E' \longrightarrow E
$$
est un morphisme de ${\mathcal C}$, on a:

\medskip

$
\left\{ \begin{matrix}
\bullet &\mbox{Pour toute famille $(S_i)_{i \in I}$ de sous-objets de $E$ telle que} \hfill \\
{ \ } \\
&\displaystyle \bigwedge_{i \in I} S_i \qquad \mbox{[resp.} \quad \displaystyle \bigvee_{i\in I} S_i \ \mbox{]} \\
{ \ } \\
&\mbox{est bien d\'efini dans $\Omega (E)$, on a} \hfill \\
{ \ } \\
&e^{-1} \biggl( \displaystyle \bigwedge_{i \in I} S_i \biggl) = \displaystyle \bigwedge_{i \in I} e^{-1} S_i \qquad \mbox{[resp.} \quad e^{-1} \biggl(\displaystyle \bigvee_{i\in I} S_i \biggl) = \displaystyle \bigvee_{i\in I} e^{-1} S_i \ \mbox{]}. \\
{ \ } \\
\bullet &\mbox{Pour toute famille $(S'_i)_{i \in I}$ de sous-objets de $E'$ telle que} \hfill \\
{ \ } \\
&\displaystyle\bigvee_{i\in I} S'_i \qquad \mbox{[resp.} \quad \displaystyle\bigwedge_{i \in I} S'_i \ \mbox{]} \\
{ \ } \\
&\mbox{est bien d\'efini dans $\Omega (E')$, on a} \hfill \\
{ \ } \\
&\exists_e \biggl( \displaystyle \bigvee_{i \in I} S'_i \biggl) = \displaystyle \bigvee_{i \in I} \exists_e \, S'_i  \qquad \mbox{[resp.} \quad \forall_e \biggl( \displaystyle \bigwedge_{i \in I} S'_i \biggl) = \displaystyle \bigwedge_{i \in I} \forall_e \, S'_i  \ \mbox{]}.
\end{matrix} \right.
$

\medskip

Cela r\'esulte de ce que les adjoints \`a droite [resp. \`a gauche] respectent toujours les limites [resp. les colimites]. 
\end{listeisansmarge}
\end{remarksqed}

\pagebreak

On observe:

\begin{lem}\label{lemV312}
\begin{listeimarge}
\item Tout topos ${\mathcal E}$ permet \`a la fois les quantifications existentielles et les quantifications universelles.

\medskip

\item Pour tout morphisme de topos
$$
f = (f^* , f_*) : {\mathcal E}' \longrightarrow {\mathcal E} \, ,
$$
sa composante d'image r\'eciproque
$$
f^* : {\mathcal E} \longrightarrow {\mathcal E}'
$$
respecte les quantifications existentielles.
\end{listeimarge}
\end{lem}

\begin{demo}
\begin{listeisansmarge}
\item Pour tout morphisme de ${\mathcal E}$
$$
e : E_2 \longrightarrow E_1 \, ,
$$
le foncteur
$$
\begin{matrix}
e^{-1} &: &\hfill \Omega (E_1) &\longrightarrow &\Omega (E_2) \, , \hfill \\
&&(S_1 \hookrightarrow E_1) &\longmapsto &(S_1 \times_{E_1} E_2 \hookrightarrow E_2)
\end{matrix}
$$
respecte les intersections (qui sont des limites) et les r\'eunions (qui sont des colimites).

\smallskip

Donc il admet pour adjoint \`a gauche [resp. \`a droite] le foncteur
$$
\exists_e : \Omega (E_2) \longrightarrow \Omega(E_1)
$$
[resp.
$$
\ \ \forall_e : \Omega (E_2) \longrightarrow \Omega(E_1) \ \mbox{]}
$$
qui associe \`a tout sous-objet
$$
S_2 \xhookrightarrow{ \ { \ } \ } E_2
$$
l'intersection [resp. la r\'eunion] des sous-objets
$$
S_1 \xhookrightarrow{ \ { \ } \ } E_1
$$
tels que
$$
S_2 \leq e^{-1} \, S_1 \qquad \mbox{[resp.} \quad e^{-1} \, S_1 \leq S_2 \ \mbox{]}.
$$

Puis consid\'erons un carr\'e cart\'esien de ${\mathcal E}$:
$$
\xymatrix{
E'_2 \ar[d]_{e'} \ar[r]^{e_2} &E_2 \ar[d]^e \\
E'_1 \ar[r]^{e_1} &E_1
}
$$

\pagebreak
Pour tout sous-objet $S_2 \hookrightarrow E_2$, le sous-objet
$$
\exists_e \, S_2 \xhookrightarrow{ \ { \ } \ } E_1
$$
est l'image du morphisme compos\'e
$$
S_2 \xhookrightarrow{ \ { \ } \ } E_2 \xrightarrow{ \ e \ } E_1
$$
donc la colimite du diagramme
$$
S_2 \times_{E_1} S_2 \! \raisebox{.7ex}{\xymatrix{\dar[r] &S_2}} \, .
$$

Comme le foncteur de changement de base
$$
E'_1 \times_{E_1} \bullet
$$ 
respecte \`a la fois les limites et les colimites, on peut conclure que
$$
e_1^{-1} \circ \exists_e = \exists_{e'} \circ e_2^{-1} \, .
$$

Echangeant les r\^oles de $(e_1,e_2)$ et de $(e,e')$, on a aussi
$$
e^{-1} \circ \exists_{e_1} = \exists_{e_2} \circ e'^{-1}
$$
puis, en passant aux adjoints \`a droite,
$$
e_1^{-1} \circ \forall_e = \forall_{e'} \circ e_2^{-1} \, .
$$
Cela ach\`eve de montrer (i).

\medskip

\item Consid\'erons un morphisme $e : E_2 \to E_1$ de ${\mathcal E}$ et un sous-objet $S_2 \hookrightarrow E_2$ de $E_2$.

\smallskip

Le foncteur
$$
f^* : {\mathcal E} \longrightarrow {\mathcal E}'
$$
respecte les limites finies et les colimites.

\smallskip

Donc il transforme la colimite $\exists_e \, S_2$ du diagramme
$$
S_2 \times_{E_1} S_2 \! \raisebox{.7ex}{\xymatrix{\dar[r] &S_2}}
$$
en la colimite $\exists_{f^*(e)} \, f^* \, S_2$ du diagramme
$$
f^* S_2 \times_{f^* E_1} f^* S_2 \! \raisebox{.7ex}{\xymatrix{\dar[r] &f^* S_2}} \, .
$$
Cela signifie que le foncteur $f^*$ respecte les quantifications existentielles. 
\end{listeisansmarge}
\end{demo}

\bigskip

Voici comment les quantificateurs existentiels $\exists$ [resp. universels $\forall$] sont interpr\'et\'es dans des cat\'egories qui permettent les quantifications existentielles [resp. qui permettent les quantifications universelles]:

\begin{defn}\label{defV313}

Soit $\Sigma$ une signature.

\smallskip

Soit ${\mathcal C}$ une cat\'egorie localement petite qui permet les quantifications existentielles [resp. permet les quantifications universelles].

\smallskip

Soit $\varphi = \varphi (\vec x , \vec y)$ une formule de $\Sigma$ dont le contexte est la r\'eunion de deux familles disjointes $\vec x = (x_1^{A_1} \cdots x_n^{A_n})$ et $\vec y = ( y_1^{B_1} \cdots y_m^{B_m})$ et qui est suppos\'ee interpr\'etable dans la cat\'egorie ${\mathcal C}$.

\smallskip

Enfin, soit $M$ une $\Sigma$-structure de ${\mathcal C}$, avec l'interpr\'etation de $\varphi$ comme un sous-objet
$$
M\varphi (\vec x , \vec y) \xhookrightarrow{ \ { \ } \ } M\!A_1 \times \cdots \times M\!A_n \times M\!B_1 \times \cdots \times M\!B_m \, .
$$

Alors la formule
$$
(\exists \, \vec y) \, \varphi \qquad \mbox{[resp.} \quad (\forall \, \vec y) \, \varphi \ \mbox{]}
$$
est interpr\'et\'ee dans $M$ comme le sous-objet
$$
\ \ M (\exists \, \vec y) \, \varphi (\vec x) \xhookrightarrow{ \ { \ } \ } M\!A_1 \times \cdots \times M\!A_n
$$
[resp.
$$
M (\forall \, \vec y) \, \varphi (\vec x) \xhookrightarrow{ \ { \ } \ } M\!A_1 \times \cdots \times M\!A_n \ \mbox{]}
$$
image du sous-objet
$$
M\varphi (\vec x , \vec y) \xhookrightarrow{ \ { \ } \ } M\!A_1 \times \cdots \times M\!A_n \times M\!B_1 \times \cdots \times M\!B_m
$$
par le foncteur
$$
\exists_p \qquad \mbox{[resp.} \quad \forall_p \ \mbox{]}
$$
associ\'e au morphisme de projection
$$
p : M\!A_1 \times \cdots \times M\!A_n \times M\!B_1 \times \cdots \times M\!B_m \longrightarrow M\!A_1 \times \cdots \times M\!A_n \, .
$$
\end{defn}

\begin{remarkqed}

Il r\'esulte du lemme \ref{lemV312} (i) que si une formule $\varphi (\vec x , \vec y)$ est interpr\'etable dans un topos ${\mathcal E}$, alors les formules de contexte $\vec x$
$$
(\exists \, \vec y) \, \varphi \qquad \mbox{et} \qquad (\forall \, \vec y) \, \varphi
$$
sont encore interpr\'etables dans ${\mathcal E}$. \end{remarkqed}

\bigskip

On d\'eduit de cette d\'efinition:

\begin{lem}\label{lemV314}

Soit $\Sigma$ une signature.

\smallskip

Soient ${\mathcal C},{\mathcal C}'$ deux cat\'egories localement petites qui permettent les quantifications existentielles [resp. permettent les quantifications universelles], reli\'ees par un foncteur
$$
F : {\mathcal C} \longrightarrow {\mathcal C}'
$$
qui respecte les quantifications existentielles [resp. respecte les quantifications universelles].

\smallskip

Soit $\varphi$ une formule de contexte $(\vec x , \vec y)$ qui est suppos\'ee interpr\'etable dans ${\mathcal C}$ et ${\mathcal C}'$ et telle que le foncteur
$$
F : \Sigma\mbox{\rm -str} ({\mathcal C}) \longrightarrow \Sigma\mbox{\rm-str} ({\mathcal C}')
$$
respecte l'interpr\'etation de $\varphi$ dans chaque $\Sigma$-structure $M$ de ${\mathcal C}$.

\smallskip

Alors ce foncteur $F$ respecte \'egalement l'interpr\'etation de la formule
$$
(\exists \, \vec y) \, \varphi \qquad \mbox{[resp.} \quad (\forall \, \vec y) \, \varphi) \ \mbox{]}
$$ 
au sens que, pour toute $\Sigma$-structure $M$ de ${\mathcal C}$, il transforme le sous-objet
$$
\ \ M (\exists \, \vec y) \, \varphi (\vec x) \xhookrightarrow{ \ { \ } \ } M\!A_1 \times \cdots \times M\!A_n
$$
[resp.
$$
M (\forall \, \vec y) \, \varphi (\vec x) \xhookrightarrow{ \ { \ } \ } M\!A_1 \times \cdots \times M\!A_n \ \mbox{]}
$$
en le sous-objet
$$
\ \ F(M) (\exists \, \vec y) \, \varphi (\vec x) \xhookrightarrow{ \ { \ } \ } F(M\!A_1) \times \cdots \times F(M\!A_n)
$$
[resp.
$$
F(M) (\forall \, \vec y) \, \varphi (\vec x) \xhookrightarrow{ \ { \ } \ } F(M\!A_1) \times \cdots \times F(M\!A_n) \ \mbox{]}.
$$
\end{lem}

\begin{remark}

Il r\'esulte du lemme \ref{lemV312} (ii) que, pour tout morphisme de topos
$$
f = (f^* , f_*) : {\mathcal E}' \longrightarrow {\mathcal E} \, ,
$$
le foncteur $f^* : {\mathcal E} \longrightarrow {\mathcal E}'$ respecte l'interpr\'etation d'une formule dans un contexte $\vec x$
$$
(\exists \, \vec y) \, \varphi 
$$
d\`es lors qu'il respecte l'interpr\'etation de la formule $\varphi$ dans un contexte $(\vec x , \vec y)$.
\end{remark}

\bigskip

\begin{demo}

C'est \'evident sur les d\'efinitions. \end{demo}

\subsection{Les implications et les n\'egations}\label{subsec534}

\medskip

En appliquant le symbole $\Rightarrow$ [resp. $\neg$] \`a une paire de formules de m\^eme contexte [resp. \`a une formule], on d\'efinit une nouvelle formule:

\begin{defn}\label{defV315}

Soit $\Sigma$ une signature.

\begin{listeimarge}

\item Pour toutes formules $\varphi$ et $\psi$ de $\Sigma$ de m\^eme contexte $\vec x = (x_1^{A_1} \cdots x_n^{A_n})$, on appelle implication de $\psi$ par $\varphi$ l'expression formelle
$$
\varphi \Rightarrow \psi
$$
consid\'er\'ee comme une formule de contexte $\vec x$.

\medskip

\item Pour toute formule $\varphi$ de $\Sigma$ de contexte $\vec x$, on appelle n\'egation de $\varphi$ l'expression formelle
$$
\neg \, \varphi
$$
consid\'er\'ee comme une formule de m\^eme contexte $\vec x$ que $\varphi$.
\end{listeimarge}
\end{defn}

\begin{remarkqed}

Une fois encore, les expressions $\varphi \Rightarrow \psi$ ou $\neg \, \varphi$ sont des \'ecritures purement formelles.

\smallskip

Elles acqui\`erent un sens par la mani\`ere dont elles sont interpr\'et\'ees dans des cat\'egories. \end{remarkqed}

\bigskip

D\'efinissons les propri\'et\'es cat\'egoriques qui permettent d'interpr\'eter les symboles $\Rightarrow$ ou $\neg$:

\begin{defn}\label{defV316}

Soit ${\mathcal C}$ une cat\'egorie qui poss\`ede des limites finies.

\begin{listeimarge}

\item On dit que ${\mathcal C}$ permet les implications si

\medskip

$
\left\{ \begin{matrix}
\bullet &\mbox{pour tout sous-objet $S$ d'un objet $E$ de ${\mathcal C}$, le foncteur interne \`a la cat\'egorie $\Omega (E)$ des sous-objets} \hfill \\
&\mbox{de $E$} \hfill \\
&\begin{matrix}
\Omega(E) &\longrightarrow &\Omega (E) \, , \\
\hfill S'' &\longmapsto &S'' \wedge S
\end{matrix} \\
&\mbox{admet un adjoint \`a droite} \hfill \\
&\begin{matrix}
\Omega(E) &\longrightarrow &\Omega (E) \, , \\
\hfill S' &\longmapsto &(S \Rightarrow S') \, ,
\end{matrix} \\
{ \ } \\
\bullet &\mbox{les foncteurs $\Rightarrow$ sont compatibles aux changements de base au sens que, pour tout morphisme de ${\mathcal C}$} \hfill \\
{ \ } \\
&e : E_2 \longrightarrow E_1 \\
&\mbox{et pour tous sous-objets $S,S'$ de $E_1$, on a} \hfill \\
{ \ } \\
&e^{-1} (S \Rightarrow S') = (e^{-1} S \Rightarrow e^{-1} S') \, .
\end{matrix} \right.
$

\medskip

\item On dit que ${\mathcal C}$ permet les n\'egations si

\medskip

$
\left\{\begin{matrix}
\bullet &\mbox{tout objet $E$ de ${\mathcal C}$ poss\`ede un plus petit sous-objet $\emptyset_E$ et on a pour tout morphisme $e : E_2 \to E_1$} \hfill \\
&\mbox{de ${\mathcal C}$} \hfill \\
&e^{-1} \, \emptyset_{E_1} = \emptyset_{E_2} \, , \\
{ \ } \\
\bullet &\mbox{pour tout sous-objet $S$ d'un objet $E$ de ${\mathcal C}$ existe un sous-objet $\neg \, S \hookrightarrow E$ caract\'eris\'e par} \hfill \\
&\mbox{l'\'equivalence, pour tout sous-objet $S' \hookrightarrow E$,} \hfill \\
&S' \leq \neg \, S \Leftrightarrow S' \wedge S = \emptyset_E \, , \\
{ \ } \\
\bullet &\mbox{les foncteurs $\neg$ sont compatibles aux changements de base au sens que pour tout morphisme} \hfill \\
&\mbox{$e : E_2 \to E_1$ de ${\mathcal E}$ et tout sous-objet $S$ de $E_1$, on a} \hfill \\
{ \ } \\
&e^{-1} (\neg \, S) = \neg \, (e^{-1} \, S) \, . 
\end{matrix} \right.
$

\medskip

\item Si ${\mathcal C},{\mathcal C}'$ sont deux cat\'egories qui poss\`edent des limites finies et permettent les implications [resp. les n\'egations], on dit qu'un foncteur respectant les limites finies
$$
F : {\mathcal C} \longrightarrow {\mathcal C}'
$$
respecte les implications [resp. les n\'egations] si on a pour tout objet $E$ de ${\mathcal C}$
$$
F(S \Rightarrow S') = ( F(S) \Rightarrow F(S') ) \, , \qquad \forall \, S,S' \in \Omega (E)
$$
[resp.
$$
F(\emptyset_E) = \emptyset_{F(E)} 
$$
et
$$
F(\neg \, S) = \neg \, F(S) \, , \qquad \forall \, S \in \Omega(E) \ \mbox{]}.
$$
\end{listeimarge}
\end{defn}

\begin{remarksqed}
\begin{listeisansmarge}
\item Si ${\mathcal C}$ est une cat\'egorie localement petite qui poss\`ede des limites finies et dont tout objet $E$ a un plus petit sous-objet $\emptyset_E$ avec
$$
e^{-1} \, \emptyset_{E_1} = \emptyset_{E_2}
$$
pour tout morphisme $e : E_2 \to E_1$, et que de plus ${\mathcal C}$ admet des implications, alors ${\mathcal C}$ admet a fortiori des n\'egations et celles-ci sont d\'efinies par la formule
$$
\neg \, S = (S \Rightarrow \emptyset_E) \, , \qquad \forall \, E \in {\rm Ob} ({\mathcal C}) \, , \qquad \forall \, S \in \Omega (E) \, .
$$

\item Si ${\mathcal C}$ et ${\mathcal C}'$ sont deux cat\'egories comme dans (i), tout foncteur $F$ respectant les limites finies
$$
F : {\mathcal C} \longrightarrow {\mathcal C}'
$$ 
respecte les n\'egations d\`es lors qu'il respecte les implications et v\'erifie
$$
F(\emptyset_E) = \emptyset_{F(E)} \, , \qquad \forall \, E \in {\rm Ob} ({\mathcal C}) \, .
$$
\end{listeisansmarge}
\end{remarksqed}

\bigskip

On observe:

\begin{lem}\label{lemV317}

Tout topos ${\mathcal E}$ permet les implications et les n\'egations.
\end{lem}

\begin{demo}

L'existence des foncteurs sur les cat\'egories $\Omega (E)$ des sous-objets des objets $E$ de ${\mathcal E}$
$$
\begin{matrix}
\neg &: &\hfill \Omega (E) &\longrightarrow &\Omega(E) \, , \\
\Rightarrow &: &\Omega(E) \times \Omega(E) &\longrightarrow &\Omega (E) \hfill
\end{matrix}
$$
a \'et\'e observ\'ee dans les remarques (ii) et (iii) qui suivent la proposition \ref{propIII62}.

\smallskip

Elle r\'esulte de ce que chaque $\Omega (E)$ est un ensemble ordonn\'e dans lequel les r\'eunions arbitraires sont bien d\'efinies, et que les foncteurs
$$
\begin{matrix}
\Omega (E) &\longrightarrow &\Omega (E) \, , \hfill \\
\hfill S'' &\longmapsto &S'' \wedge S \hfill
\end{matrix}
$$
respectent les r\'eunions arbitraires.

\smallskip

Pour tout morphisme de ${\mathcal E}$
$$
e : E_2 \longrightarrow E_1 \, ,
$$
le foncteur induit
$$
e^{-1} : \Omega (E_1) \longrightarrow \Omega (E_2)
$$
respecte les r\'eunions d'apr\`es le lemme \ref{lemIII63}. En particulier, il envoie $\emptyset_{E_1}$ sur $\emptyset_{E_2}$.

\smallskip

De plus, pour tout sous-objet $S_1$ de $E_1$, le carr\'e
$$
\xymatrix{
\Omega(E_2) \ar[d]_{\exists_e} \ar[rr]^{\bullet \wedge e^{-1} S_1} &&\Omega (E_2) \ar[d]^{\exists_e} \\
\Omega (E_1) \ar[rr]^{\bullet \wedge S_1} &&\Omega (E_1)
}
$$
est commutatif puisque ${\mathcal E}$ est une cat\'egorie r\'eguli\`ere d'apr\`es le lemme \ref{lemV312} et donc que le foncteur
$$
\exists_e : \Omega (E_2) \longrightarrow \Omega (E_1)
$$
est respect\'e par le foncteur de changement de base
$$
\bullet \times_{E_1} S_1
$$
par le morphisme $S_1 \hookrightarrow E_1$.

\smallskip

En passant aux adjoints \`a droite, on obtient que le carr\'e
$$
\xymatrix{
\Omega(E_1) \ar[d]_{e^{-1}} \ar[rr]^{\bullet \Rightarrow S_1} &&\Omega (E_1) \ar[d]^{e^{-1}} \\
\Omega (E_2) \ar[rr]^{\bullet \Rightarrow e^{-1} S_1} &&\Omega (E_2)
}
$$
est commutatif.

\smallskip

Cela signifie que les foncteurs $\Rightarrow$ sont compatibles avec les changements de base par les morphismes $e : E_2 \to E_1$ de ${\mathcal E}$.

\smallskip

Comme $e^{-1} (\emptyset_{E_1}) = \emptyset_{E_2}$ pour tout tel morphisme $e$, on conclut que les foncteurs $\neg$ sont \'egalement compatibles avec les changements de base.

\smallskip

Cela termine la d\'emonstration du lemme. 

\end{demo}

\bigskip

Dans les cat\'egories qui permettent les implications [resp. les n\'egations], les symboles logiques $\Rightarrow$ [resp. $\neg$] sont interpr\'et\'es comme les foncteurs $\Rightarrow$ [resp. $\neg$]:

\begin{defn}\label{defV318}

Soit $\Sigma$ une signature.

\smallskip

Soit ${\mathcal C}$ une cat\'egorie localement petite avec limites finies qui permet les n\'egations [resp. les implications].

\smallskip

Soit $\varphi$ une formule [resp. $(\varphi , \psi)$ une paire de formules] de $\Sigma$ dans un contexte $\vec x = (x_1^{A_1} \cdots x_n^{A_n})$ qui est suppos\'ee interpr\'etable dans ${\mathcal C}$.

\smallskip

Enfin, soit $M$ une $\Sigma$-structure de ${\mathcal C}$, avec l'interpr\'etation de $\varphi$ [resp. et $\psi$] comme le sous-objet
$$
M\varphi (\vec x) \xhookrightarrow{ \ { \ } \ } M\!A_1 \times \cdots \times M\!A_n
$$
[resp. et le sous-objet
$$
M\psi (\vec x) \xhookrightarrow{ \ { \ } \ } M\!A_1 \times \cdots \times M\!A_n \ \mbox{]}.
$$

Alors la formule
$$
\neg \, \varphi \qquad \mbox{[resp.} \quad \varphi \Rightarrow \psi \ \mbox{]}
$$
est interpr\'et\'ee dans $M$ comme le sous-objet
$$
M(\neg \, \varphi)(\vec x) \xhookrightarrow{ \ { \ } \ } M\!A_1 \times \cdots \times M\!A_n
$$
[resp.
$$
M(\varphi \Rightarrow \psi) (\vec x) \xhookrightarrow{ \ { \ } \ } M\!A_1 \times \cdots \times M\!A_n \ \mbox{]}
$$
 image du sous-objet
 $$
 M\varphi (\vec x) \xhookrightarrow{ \ { \ } \ } M\!A_1 \times \cdots \times M\!A_n
$$
[resp. et du sous-objet
$$
M\psi (\vec x) \xhookrightarrow{ \ { \ } \ } M\!A_1 \times \cdots \times M\!A_n \ \mbox{]}
$$
par le foncteur
$$
\neg \qquad \mbox{[resp.} \quad \Rightarrow \ \mbox{]}
$$
agissant sur les sous-objets [resp. sur les paires de sous-objets] de l'objet $M\!A_1 \times \cdots \times M\!A_n$ de ${\mathcal E}$.
\end{defn}

\bigskip

\begin{remarkqed}

Il r\'esulte du lemme \ref{lemV317} que si une formule $\varphi$ [resp. une paire de formules $\varphi$ et $\psi$] de contexte $\vec x$ est interpr\'etable dans un topos ${\mathcal E}$, alors la formule de m\^eme contexte $\vec x$
$$
\neg \, \varphi \qquad \mbox{[resp.} \quad \varphi \Rightarrow \psi \ \mbox{]}
$$
est encore interpr\'etable dans ${\mathcal E}$. 

\end{remarkqed}

\bigskip

On d\'eduit de cette d\'efinition:

\begin{cor}\label{corV319}

Soit $\Sigma$ une signature.

\smallskip

Soient ${\mathcal C},{\mathcal C}'$ deux cat\'egories localement petites avec limites finies qui permettent les n\'egations [resp. les implications], reli\'ees par un foncteur
$$
F : {\mathcal C} \longrightarrow {\mathcal C}'
$$
qui respecte les limites finies et les n\'egations [resp. les implications].

\smallskip

Soit $\varphi$ une formule [resp. $(\varphi,\psi)$ une paire de formules] de $\Sigma$ dans un contexte $\vec x = (x_1^{A_1} \cdots x_n^{A_n})$, qui est suppos\'ee interpr\'etable dans ${\mathcal C}$ et ${\mathcal C}'$ et telle que le foncteur
$$
F : \Sigma\mbox{\rm -str} ({\mathcal C}) \longrightarrow \Sigma\mbox{\rm -Str} ({\mathcal C}')
$$
respecte l'interpr\'etation de $\varphi$ [resp. et de $\psi$] dans chaque $\Sigma$-structure $M$ de ${\mathcal C}$.

\smallskip

Alors le foncteur $F$ respecte \'egalement les interpr\'etations de la formule 
$$
\neg \, \varphi \qquad \mbox{[resp.} \quad \varphi \Rightarrow \psi \ \mbox{]}
$$
au sens que, pour toute $\Sigma$-structure $M$ de ${\mathcal C}$, il transforme le sous-objet
$$
M(\neg \, \varphi)(\vec x) \xhookrightarrow{ \ { \ } \ } M\!A_1 \times \cdots \times M\!A_n
$$
[resp.
$$
M(\varphi \Rightarrow \psi) (\vec x) \xhookrightarrow{ \ { \ } \ } M\!A_1 \times \cdots \times M\!A_n \ \mbox{]}
$$
en le sous-objet
$$
F(M) (\neg \, \varphi)(\vec x) \xhookrightarrow{ \ { \ } \ } F(M\!A_1) \times \cdots \times F(M\!A_n)
$$
[resp.
$$
F(M)(\varphi \Rightarrow \psi) (\vec x) \xhookrightarrow{ \ { \ } \ } F(M\!A_1) \times \cdots \times F(M\!A_n) \ \mbox{]}.
$$
\end{cor}

\begin{demo}

C'est \'evident sur les d\'efinitions. 
\end{demo}

\bigskip

\section{Types de th\'eories du premier ordre et cat\'egories de mod\`eles}\label{sec54}

\subsection{Types de formules et types de s\'equents}\label{subsec541}

\medskip

On distingue les types de formules suivants en fonction des symboles qui entrent dans leur construction \`a partir des formules atomiques:

\begin{defn}\label{defV41}

Soit $\Sigma$ une signature.

\begin{listeimarge}

\item Une formule de $\Sigma$ est dite ``de Horn'' si elle est construite \`a partir de formules atomiques et du symbole $\top$ en appliquant le symbole de conjonction finie $\wedge$.

\medskip

\item Une formule de $\Sigma$ est dite ``r\'eguli\`ere'' si elle est construite \`a partir de formules atomiques et du symbole $\top$ en appliquant des symboles de conjonctions finies $\wedge$ ou de quantifications existentielles $\exists$.

\medskip

\item Une formule de $\Sigma$ est dite ``coh\'erente'' [resp. ``g\'eom\'etrique''] si elle est construite \`a partir de formules atomiques et des symboles $\top$ et $\perp$ en appliquant des symboles de conjonctions finies $\wedge$, de quantifications existentielles $\exists$ et de disjonctions finies $\vee$ [resp. de disjonctions arbitraires $\bigvee$].

\medskip

\item Une formule de $\Sigma$ est dite ``du premier ordre finitaire'' [resp. ``du premier ordre infinitaire''] si elle est construite \`a partir de formules atomiques et des symboles $\top$ et $\perp$ en appliquant des symboles de conjonctions finies $\wedge$ et de disjonctions finies $\vee$ [resp. de conjonctions arbitraires $\bigwedge$ et de disjonctions arbitraires $\bigvee$] ainsi que de quantifications existentielles $\exists$, de quantifications universelles $\forall$, d'implications $\Rightarrow$ et de n\'egations $\neg$.
\end{listeimarge}
\end{defn}

\begin{remarkqed}

Les formules d'une signature $\Sigma$ qui sont ``de Horn'', ``r\'eguli\`eres'', ``coh\'erentes'' ou ``du premier ordre finitaires'', forment un ensemble.

\smallskip

En revanche, la collection des formules g\'eom\'etriques et celle des formules du premier ordre infinitaires ne sont pas des ensembles. 

\end{remarkqed}

\bigskip

Les types de s\'equents correspondants sont d\'efinis suivant les types de formules qui les composent:

\begin{defn}\label{defV42}

Soit $\Sigma$ une signature.

\smallskip

Un segment de $\Sigma$
$$
\varphi \vdash_{\vec x} \psi
$$
est dit ``de Horn'' [resp. ``r\'egulier'', resp. ``coh\'erent'', resp. ``g\'eom\'etrique'', resp. ``du premier ordre finitaire'', resp. ``du premier ordre infinitaire''] si ses deux composantes $\varphi$ et $\psi$ sont des formules de Horn [resp. r\'eguli\`eres, resp. coh\'erentes, resp. g\'eom\'etriques, resp. du premier ordre finitaires, resp. du premier ordre infinitaires]. \hfill $\Box$
\end{defn}

\subsection{Types de th\'eories}\label{subsec542}

\medskip

On distingue d'abord les th\'eories r\'eduites \`a leur seul langage:

\begin{defn}\label{defV43}

Une th\'eorie de signature $\Sigma$ est appel\'ee ``la th\'eorie vide'' si elle ne comporte aucun axiome.
\end{defn}

\begin{remarkqed}

Par exemple, on a d\'ej\`a not\'e que la th\'eorie des carquois est la th\'eorie vide dans la signature constitu\'ee des deux sortes ${\rm Ob}$ et ${\rm Fl}$ et des deux symboles de fonctions $s : {\rm Fl} \to {\rm Ob}$ et $b : {\rm Fl} \to {\rm Ob}$ de source et de but. 

\end{remarkqed}

\bigskip

Puis on distingue les th\'eories alg\'ebriques:

\begin{defn}\label{defV44}

Une th\'eorie est dite alg\'ebrique si
\begin{enumerate}
\item[$\bullet$] sa signature ne comporte que des sortes et des symboles de fonctions, autrement dit ne comporte aucun symbole de relation,
\item[$\bullet$] tous ses axiomes ont la forme
$$
\top \vdash_{\vec x} \varphi
$$
o\`u $\varphi$ est une formule atomique, donc n\'ecessairement une formule d'\'egalit\'e
$$
f(\vec x) = g(\vec x)
$$
entre deux termes de m\^eme contexte $\vec x = (x_1^{A_1} \cdots x_n^{A_n})$.
\end{enumerate}
\end{defn}

\bigskip

\begin{remarksqed}
\begin{listeisansmarge}
\item Sont par exemple alg\'ebriques les th\'eories
\begin{enumerate}
\item[$\bullet$] des groupes (ou des mono{\"\i}des),
\item[$\bullet$] des groupes (ou des mono{\"\i}des) commutatifs,
\item[$\bullet$] des actions de groupes (ou de mono{\"\i}des),
\item[$\bullet$] des actions d'un groupe (ou d'un mono{\"\i}de) fix\'e,
\item[$\bullet$] des anneaux,
\item[$\bullet$] des anneaux commutatifs,
\item[$\bullet$] des modules sur un anneau,
\item[$\bullet$] des modules sur un anneau fix\'e (et en particulier des espaces vectoriels sur un corps fix\'e).
\end{enumerate}

\medskip

\item En revanche, avec cette d\'efinition, la th\'eorie des corps n'est pas alg\'ebrique, non plus que celle des anneaux locaux.

\smallskip

En effet, la th\'eorie des corps est d\'efinie en adjoignant aux axiomes de la th\'eorie alg\'ebrique des anneaux l'axiome suppl\'ementaire
$$
\top \vdash_x x = 0 \vee \exists \, y \, (x \cdot y = 1) \, .
$$

De m\^eme, la th\'eorie des anneaux locaux est d\'efinie en adjoignant aux axiomes de la th\'eorie alg\'ebrique des anneaux commutatifs l'axiome suppl\'emenaire
$$
\exists \, z \, ((x+y) \cdot z = 1) \vdash_{x , y} \exists \, x' \, (x \cdot x' = 1) \vee \exists \, y' \, (y \cdot y' = 1)
$$
qui dit que si une somme $x+y$ est inversible, alors l'un au moins des deux \'el\'ements $x,y$ qui la composent est inversible. 
\end{listeisansmarge}
\end{remarksqed}

\bigskip

On distingue encore diff\'erents types de th\'eories suivant les types de s\'equents qui constituent leurs axiomes:

\begin{defn}\label{defV45}

Une th\'eorie ${\mathbb T}$ de signature $\Sigma$ est dite ``de Horn'' [resp. ``r\'eguli\`ere'', resp. ``coh\'erente'', resp. ``g\'eom\'etrique'', resp. ``du premier ordre finitaire'', resp. ``du premier ordre infinitaire''] si tous ses axiomes sont des s\'equents ``de Horn'' [resp. ``r\'eguliers'', resp. ``coh\'erents'', resp. ``g\'eom\'etriques'', resp. ``du premier ordre finitaires'', resp. ``du premier ordre infinitaires''].
\end{defn}

\pagebreak

\begin{remarksqed}
\begin{listeisansmarge}

\item Toute th\'eorie alg\'ebrique est de Horn.

\smallskip

Toute th\'eorie de Horn est r\'eguli\`ere.

\smallskip

Toute th\'eorie r\'eguli\`ere est coh\'erente.

\smallskip

Toute th\'eorie coh\'erente est g\'eom\'etrique.

\smallskip

Toute th\'eorie g\'eom\'etrique, ainsi que toute th\'eorie du premier ordre finitaire, est une th\'eorie du premier ordre infinitaire.

\medskip

\item Des exemples de th\'eories de Horn qui ne sont pas alg\'ebriques sont
\begin{enumerate}
\item[$\bullet$] la th\'eorie des relations d'\'equivalence,
\item[$\bullet$] la th\'eorie des relations d'ordre.
\end{enumerate}

\medskip

\item Des exemples de th\'eories r\'eguli\`eres qui ne sont pas de Horn sont
\begin{enumerate}
\item[$\bullet$] la th\'eorie des cat\'egories,
\item[$\bullet$] la th\'eorie des suites exactes de groupes ab\'eliens.
\end{enumerate}

\medskip

\item Des exemples de th\'eories coh\'erentes qui ne sont pas r\'eguli\`eres sont
\begin{enumerate}
\item[$\bullet$] la th\'eorie des corps,
\item[$\bullet$] la th\'eorie des anneaux locaux,
\item[$\bullet$] la th\'eorie des plans affines,
\item[$\bullet$] la th\'eorie des plans euclidiens.
\end{enumerate}

\medskip

\item Des exemples de th\'eories g\'eom\'etriques qui ne sont pas coh\'erentes sont
\begin{enumerate}
\item[$\bullet$] la th\'eorie des \'el\'ements nilpotents [resp. des racines de l'unit\'e] d'un anneau,
\item[$\bullet$] la th\'eorie des extensions alg\'ebriques d'anneaux commutatifs.
\end{enumerate}

Pour d\'efinir ces th\'eories, on a besoin d'enrichir le langage $(A,+,\cdot ,0,1,-(\bullet))$ de la th\'eorie des anneaux en lui adjoignant des symboles de fonctions
$$
(\bullet)^n : A \longrightarrow A
$$
index\'es par les entiers $n \geq 1$ et soumis aux axiomes
$$
\top \vdash_x (x)^1 = x
$$
et
$$
\top \vdash_x (x)^n = x \cdot (x)^{n-1}
$$
pour tout $n \geq 2$.

\smallskip

La th\'eorie des \'el\'ements nilpotents [resp. des racines de l'unit\'e] d'un anneau est d\'efinie en adjoignant \`a la th\'eorie des anneaux $(A,+,\cdots)$ compl\'et\'ee par les symboles de fonctions $(\bullet)^n$ et leurs axiomes un symbole de relation
$$
N \rightarrowtail A \qquad \mbox{[resp.} \quad R \rightarrowtail A \ \mbox{]}
$$
et l'axiome en la variable libre $x$ de sorte $A$
$$
\qquad \quad \ N(x) \dashv \, \vdash \bigvee_{n \geq 1} (x)^n = 0
$$
$$
\mbox{[resp.} \qquad R(x) \dashv \, \vdash \bigvee_{n \geq 1} (x)^n = 1 \ \mbox{]}.
$$

La th\'eorie des extensions alg\'ebriques est d\'efinie \`a partir de celle des morphismes d'anneaux commutatifs, qui consiste en deux copies $(A,+,\cdot , 1_A , 0_A , \cdots)$ et $(B,+,\cdot , 1_B , 0_B , \cdots)$ de la th\'eorie des anneaux commutatifs compl\'et\'ees par un symbole de fonction
$$
u : A \longrightarrow B
$$
soumis aux axiomes
$$
\begin{matrix}
\top &\vdash &u(a_1 + a_2) = u(a_1) + u(a_2) \, , \\
\top &\vdash &u(a_1 \cdot a_2) = u(a_1) \cdot u(a_2)  \hfill \\
\top &\vdash &u(0_A) = 0_B \, , \hfill \\
\top &\vdash &u(1_A) = 1_B \, . \hfill
\end{matrix}
$$

Alors la th\'eorie des extensions alg\'ebriques s'obtient en ajoutant les symboles de fonctions
$$
(\bullet)^n : B \longrightarrow B
$$
et leurs axiomes, plus l'axiome suppl\'ementaire en la variable $x$ de sorte $B$ et les variables li\'ees $a_0 , a_1, \cdots , a_n , \cdots$ de sorte $A$
$$
\top \vdash_x \bigvee_{n \geq 1} \exists \, (a_0 , a_1 , \cdots , a_{n-1})(x^n + u(a_{n-1}) \cdot x^{n-1} + \cdots + u(a_1) \cdot x + u(a_0) = 0_B) \, .
$$

\item Des exemples de th\'eories du premier ordre infinitaires sont
\begin{enumerate}
\item[$\bullet$] la th\'eorie des entiers de Peano,
\item[$\bullet$] la th\'eorie des ensembles.
\end{enumerate}
\end{listeisansmarge}
\end{remarksqed}

\bigskip

Parmi les th\'eories r\'eguli\`eres, on distingue en particulier les th\'eories cart\'esiennes:

\begin{defn}\label{defV46}

Une th\'eorie de signature $\Sigma$ est dite cart\'esienne si elle est r\'eguli\`ere et si ses axiomes sont
\begin{enumerate}
\item[$\bullet$] ou bien de Horn,
\item[$\bullet$] ou bien de la forme
$$
f_1 (\vec x) = g_1 (\vec x) \wedge \cdots \wedge f_k (\vec x) = g_k (\vec x) \dashv \, \vdash_{\vec x} (\exists \, y) (u_1 (y) = x_1^{A_1} \wedge \cdots \wedge y_n (y) = x_n^{A_n})
$$
et accompagn\'es alors d'un axiome
$$
u_1 (y) = u_1 (y') \wedge \cdots \wedge y_n (y) = y_n (y') \vdash_{y,y'} y = y'
$$
pour des variables $\vec x = (x_1^{A_1} \cdots x_n^{A_n})$ de sortes $A_1 , \cdots , A_n$ et des variables $y,y'$ de sorte $B$, des symboles de fonctions
$$
u_1 : B \longrightarrow A_1 , \cdots , u_n : B \longrightarrow A_n \, ,
$$
et des formules atomiques d'\'egalit\'e
$$
f_i (\vec x) = g_i (\vec x) \, , \qquad 1 \leq i \leq k \, .
$$
\end{enumerate}
\end{defn}

\begin{remarksqed}
\begin{listeisansmarge}
\item Ainsi, une th\'eorie r\'eguli\`ere est cart\'esienne si chacun de ses axiomes o\`u appara{\^\i}t un symbole de quantification existentielle $\exists$ consiste \`a demander qu'un certain sous-objet image d'un objet par une fonction soit caract\'eris\'e par une famille finie d'\'equations, et si chaque tel axiome est accompagn\'e d'un autre qui affirme l'unicit\'e des ant\'ec\'edents des images par cette fonction.

\medskip

\item Par exemple, la th\'eorie des cat\'egories est cart\'esienne.

\smallskip

En effet, le seul de ses axiomes qui comporte le symbole $\exists$ est l'axiome
$$
b(f) = s(g) \dashv \, \vdash (\exists \, c) (p_1 (c) = f \wedge p_2 (c) = g)
$$
en les variables libres $f,g$ de sorte Hom et la variable li\'ee $c$ de sorte Comp, et il est accompagn\'e de l'axiome d'unicit\'e 
$$
p_1 (c) = p_1 (c') \wedge p_2 (c) = p_2 (c') \vdash c=c' 
$$
en les variables libres $c,c'$ de sorte Comp. 
\end{listeisansmarge}
\end{remarksqed}

\subsection{Types de cat\'egories et mod\`eles}\label{subsec543}

On introduit diff\'erents types de cat\'egories dans lesquels les diff\'erents types de cat\'egories que nous avons introduits s'interpr\`etent et d\'efinissent des cat\'egories de mod\`eles:

\begin{defn}\label{defV47}
\begin{listeimarge}
\item Une cat\'egorie localement petite ${\mathcal C}$ est dite
\begin{enumerate}
\item[$\bullet$] alg\'ebrique si elle poss\`ede des produits finis, en particulier un objet terminal $1_{\mathcal C}$,
\item[$\bullet$] cart\'esienne si elle poss\`ede des limites finies,
\item[$\bullet$] r\'eguli\`ere si elle est cart\'esienne et permet les quantifications existentielles au sens de la d\'efinition \ref{defV311}~(i),
\item[$\bullet$] coh\'erente si elle est r\'eguli\`ere et a des r\'eunions finies de sous-objets au sens de la d\'efinition \ref{defV36}~(iii),
\item[$\bullet$] g\'eom\'etrique si elle est r\'eguli\`ere et a des r\'eunions arbitraires de sous-objets au sens de la d\'efinition \ref{defV36}~(iii),
\item[$\bullet$] de Heyting si elle est coh\'erente et permet les implications au sens de la d\'efinition \ref{defV316}~(i).
\end{enumerate}

\medskip

\item Un foncteur
$$
F : {\mathcal C} \longrightarrow {\mathcal D}
$$
entre deux cat\'egories ${\mathcal C}$, ${\mathcal D}$ suppos\'ees alg\'ebriques [resp. cart\'esiennes, resp. r\'eguli\`eres, resp. coh\'erentes, resp. g\'eom\'etriques, resp. de Heyting] est dit
\begin{enumerate}
\item[$\bullet$] alg\'ebrique s'il respecte les produits finis,
\item[$\bullet$] resp. cart\'esien s'il respecte les limites finies,
\item[$\bullet$] resp. r\'egulier s'il est cart\'esien et respecte les quantifications existentielles au sens de la d\'efinition \ref{defV311}~(ii),
\item[$\bullet$] resp. coh\'erent s'il est r\'egulier et respecte les r\'eunions finies de sous-objets au sens de la remarque {\rm (ii)} apr\`es la d\'efinition \ref{defV36},
\item[$\bullet$] resp. g\'eom\'etrique s'il est r\'egulier et respecte les r\'eunions arbitraires de sous-objets au sens de la remarque {\rm (iii)} apr\`es la d\'efinition \ref{defV36},
\item[$\bullet$] resp. de Heyting s'il est coh\'erent et respecte les implications au sens de la d\'efinition \ref{defV316}~(iii).
\end{enumerate}
\end{listeimarge}
\end{defn}

\begin{remarkqed}

D'apr\`es la remarque (i) qui suit la d\'efinition \ref{defV316}, toute cat\'egorie de Heyting permet les n\'egations. Plus pr\'ecis\'ement, le foncteur de n\'egation des sous-objets de tout objet $E$ est donn\'e par la formule
$$
\begin{matrix}
\neg &: &\hfill \Omega (E) &\longrightarrow &\Omega (E) \, , \hfill \\
&&\hfill (S \hookrightarrow E) &\longmapsto &(S \Rightarrow \emptyset_E) \, .
\end{matrix}
$$
Il est respect\'e par tout morphisme de changement de base
$$
e : E' \longrightarrow E
$$
puisqu'il en est ainsi du foncteur $\Rightarrow$ et que
$$
e^{-1} \, \emptyset_E = \emptyset_{E'} \, .
$$
Pour les m\^emes raisons, tout foncteur de Heyting
$$
F : {\mathcal C} \longrightarrow {\mathcal D}
$$
entre deux cat\'egories de Heyting respecte les n\'egations. 

\end{remarkqed}

\bigskip

On d\'eduit des lemmes \ref{lemV37}, \ref{lemV312} et \ref{lemV317}:

\begin{cor}\label{corV48}
\begin{listeimarge}
\item Tout topos est \`a la fois une cat\'egorie g\'eom\'etrique et une cat\'egorie de Heyting. De plus, il admet des intersections arbitraires de sous-objets au sens de la d\'efinition \ref{defV36}.

\medskip

\item Pour tout morphisme de topos
$$
f = (f^* , f_*) : {\mathcal E}' \longrightarrow {\mathcal E} \, ,
$$
sa composante d'image r\'eciproque
$$
f^* : {\mathcal E} \longrightarrow {\mathcal E}' 
$$
est un foncteur g\'eom\'etrique, et sa composante d'image directe
$$
f_* : {\mathcal E} \longrightarrow {\mathcal E}'
$$
est un foncteur cart\'esien qui, de plus, respecte les intersections arbitraires de sous-objets.
\end{listeimarge}
\end{cor}

\begin{remark}

En revanche, la composante d'image r\'eciproque $f^*$ d'un morphisme de topos $f=(f^*,f_*) : {\mathcal E}' \to {\mathcal E}$ ne respecte pas en g\'en\'eral les quantifications universelles, les implications, les n\'egations ou les intersections arbitraires de sous-objets.

\smallskip

Quant \`a sa composante d'image directe $f_*$, elle ne respecte pas en g\'en\'eral les quantifications existentielles ou universelles, les implications, les n\'egations ou les r\'eunions m\^eme finies de sous-objets.
\end{remark}

\bigskip

\begin{demo}
\begin{listeisansmarge}
\item On sait d\'ej\`a que tout topos est une cat\'egorie cart\'esienne.

\smallskip

D'apr\`es le lemme \ref{lemV37} (i), il admet des intersections arbitraires et des r\'eunions arbitraires de sous-objets.

\smallskip

D'apr\`es le lemme \ref{lemV312} (ii), il permet les quantifications existentielles et les quantifications universelles.

\smallskip

D'apr\`es le lemme \ref{lemV317} (i), il permet les implications et les n\'egations.

\medskip

\item On sait d\'ej\`a que pour tout morphisme de topos
$$
f = (f^* , f_*) : {\mathcal E}' \longrightarrow {\mathcal E} \, ,
$$
ses deux composantes $f^*$ et $f_*$ sont des foncteurs cart\'esiens.

\smallskip

De plus, d'apr\`es le lemme \ref{lemV37} (ii), sa composante d'image r\'eciproque $f^*$ respecte les r\'eunions arbitraires de sous-objets, et sa composante d'image directe $f_*$ respecte les intersections arbitraires de sous-objets.

\smallskip

Enfin, d'apr\`es le lemme \ref{lemV312} (ii), la composante $f^*$ respecte les quantifications existentielles. 
\end{listeisansmarge}
\end{demo}

\bigskip

Toute th\'eorie alg\'ebrique d\'efinit une cat\'egorie de mod\`eles dans toute cat\'egorie alg\'ebrique:

\begin{defn}\label{defV49}

Soit ${\mathbb T}$ une th\'eorie alg\'ebrique de signature $\Sigma$.

\smallskip

Soit ${\mathcal C}$ une cat\'egorie alg\'ebrique.

\smallskip

Alors:

\begin{listeimarge}

\item On appelle mod\`ele de ${\mathbb T}$ (ou ${\mathbb T}$-mod\`ele) dans ${\mathcal C}$ toute $\Sigma$-structure de ${\mathcal C}$
$$
M
$$
telle que, pour tout axiome de ${\mathbb T}$ n\'ecessairement de la forme
$$
\top \vdash f(\vec x) = g(\vec x)
$$
pour deux termes $f$ et $g$ de m\^eme contexte $\vec x = (x_1^{A_1} \cdots x_n^{A_n})$ \`a valeurs dans une sorte $B$, les deux morphismes
$$
M\!f , M\!g : M\!A_1 \times \cdots \times M\!A_n \longrightarrow M\!B
$$
sont \'egaux.

\medskip

\item On appelle cat\'egorie des mod\`eles de ${\mathbb T}$ (ou des ${\mathbb T}$-mod\`eles) dans ${\mathcal C}$ la sous-cat\'egorie pleine
$$
{\mathbb T}\mbox{\rm -mod} ({\mathcal C})
$$
de la cat\'egorie des $\Sigma$-structures de ${\mathcal C}$
$$
\Sigma\mbox{\rm -str} ({\mathcal C})
$$
constitu\'ee des mod\`eles de ${\mathbb T}$.
\end{listeimarge}
\end{defn}

\begin{remarksqed}
\begin{listeisansmarge}
\item Les interpr\'etations de termes $f(\vec x) , g(\vec x)$ comme des morphismes
$$
M\!f , M\!g : M\!A_1 \times \cdots \times M\!A_n \longrightarrow M\!B
$$
ont \'et\'e introduites dans la d\'efinition V.3.2.

\medskip

\item Si $\Sigma$ est une signature sans symbole de relation, tout foncteur alg\'ebrique entre deux cat\'egories alg\'ebriques
$$
F : {\mathcal C} \longrightarrow {\mathcal D}
$$
induit un foncteur
$$
F : \Sigma\mbox{-str} ({\mathcal C}) \longrightarrow \Sigma\mbox{-str} ({\mathcal D}) \, .
$$

Si ${\mathbb T}$ est une th\'eorie alg\'ebrique de signature $\Sigma$, ce foncteur $F$ transforme les ${\mathbb T}$-mod\`eles de ${\mathcal C}$ en ${\mathbb T}$-mod\`eles de ${\mathcal D}$. Il induit donc par restriction un foncteur
$$
F : {\mathbb T}\mbox{-mod} ({\mathcal C}) \longrightarrow {\mathbb T}\mbox{-mod} ({\mathcal D}) \, .
$$

\item Si ${\mathbb T}$ est la th\'eorie vide dans une signature $\Sigma$ (qui comporte \'eventuellement des symboles de relations), on pose pour toute cat\'egorie alg\'ebrique ${\mathcal C}$
$$
{\mathbb T}\mbox{-mod} ({\mathcal C}) = \Sigma\mbox{-str} ({\mathcal C}) \, .
$$

Si $\Sigma$ a des symboles de relations, un foncteur alg\'ebrique entre cat\'egories alg\'ebriques
$$
F : {\mathcal C} \longrightarrow {\mathcal D}
$$
induit un foncteur
$$
{\mathbb T}\mbox{-mod} ({\mathcal C}) = \Sigma\mbox{-str} ({\mathcal C}) \longrightarrow \Sigma\mbox{-str} ({\mathcal D}) = {\mathbb T}\mbox{-mod} ({\mathcal D})
$$
s'il transforme les monomorphismes en monomorphismes. 
\end{listeisansmarge}
\end{remarksqed}

\bigskip

On d\'efinit de m\^eme les cat\'egories de mod\`eles des diff\'erents types de th\'eories dans des cat\'egories qui poss\`edent les propri\'et\'es n\'ecessaires \`a l'interpr\'etation des axiomes de ces th\'eories:

\begin{defn}\label{defV4109}

Soit ${\mathbb T}$ une th\'eorie de signature $\Sigma$ qui est ``de Horn'' [resp. r\'eguli\`ere, resp. coh\'erente, resp. g\'eom\'etrique, resp. du premier ordre finitaire, resp. du premier ordre infinitaire].

\smallskip

Soit ${\mathcal C}$ une cat\'egorie qui est cart\'esienne [resp. r\'eguli\`ere, resp. coh\'erente, resp. g\'eom\'etrique, resp. ``de Heyting'', resp. un topos].

\smallskip

Alors:

\begin{listeimarge}

\item On appelle mod\`ele de ${\mathbb T}$ (ou ${\mathbb T}$-mod\`ele) dans ${\mathcal C}$ toute $\Sigma$-structure de ${\mathcal C}$
$$
M
$$
telle que, pour tout axiome de ${\mathbb T}$
$$
\varphi \vdash_{\vec x} \psi
$$
de contexte $\vec x = (x_1^{A_1} \cdots x_n^{A_n})$, les interpr\'etations
$$
\begin{matrix}
M\varphi &\xhookrightarrow{ \ { \ } \ } M\!A_1 \times \cdots \times M\!A_n \, , \\
M\psi &\xhookrightarrow{ \ { \ } \ } M\!A_1 \times \cdots \times M\!A_n \hfill
\end{matrix}
$$
satisfont la relation d'inclusion entre sous-objets de $M\!A_1 \times \cdots \times M\!A_n$
$$
M\varphi \leq M\psi \, .
$$

\item On appelle cat\'egorie des mod\`eles de ${\mathbb T}$ (ou des ${\mathbb T}$-mod\`eles) dans ${\mathcal C}$ la sous-cat\'egorie pleine
$$
{\mathbb T}\mbox{\rm -mod} ({\mathcal C}) \qquad \mbox{de} \qquad \Sigma\mbox{\rm -str} ({\mathcal C})
$$
constitu\'ee des mod\`eles de ${\mathcal C}$.
\end{listeimarge}
\end{defn}

\begin{remarksqed}
\begin{listeisansmarge}
\item Si ${\mathbb T}$ est une th\'eorie de Horn [resp. r\'eguli\`ere, resp. coh\'erente, resp. g\'eom\'etrique, resp. du premier ordre finitaire], tout foncteur cart\'esien [resp. r\'egulier, resp. coh\'erent, resp. g\'eom\'etrique, resp. de Heyting]
$$
F : {\mathcal C} \longrightarrow {\mathcal D}
$$
entre deux cat\'egories cart\'esiennes [resp. r\'eguli\`eres, resp. coh\'erentes, resp. g\'eom\'etriques, resp. de Heyting] induit un foncteur
$$
F : {\mathbb T}\mbox{-mod} ({\mathcal C}) \longrightarrow {\mathbb T}\mbox{-mod} ({\mathcal D})
$$
qui est la restriction du foncteur
$$
F : \Sigma\mbox{-str} ({\mathcal C}) \longrightarrow \Sigma\mbox{-str} ({\mathcal D}) \, .
$$

\item En particulier, pour tout morphisme de topos
$$
f = (f^* , f_*) : {\mathcal E}' \longrightarrow {\mathcal E} \, ,
$$
sa composante d'image r\'eciproque
$$
f^* : {\mathcal E} \longrightarrow {\mathcal E}'
$$
d\'efinit un foncteur
$$
f^* : {\mathbb T}\mbox{-mod} ({\mathcal E}) \longrightarrow {\mathbb T}\mbox{-mod} ({\mathcal E}')
$$
pour toute th\'eorie g\'eom\'etrique ${\mathbb T}$. 
\end{listeisansmarge}
\end{remarksqed}

\bigskip

Dans une cat\'egorie cart\'esienne ne sont pas seulement d\'efinis les mod\`eles des th\'eories alg\'ebriques ou de Horn mais aussi ceux des th\'eories cart\'esiennes:

\begin{defn}\label{defV411}

Soit ${\mathbb T}$ une th\'eorie cart\'esienne de signature $\Sigma$.

\smallskip

Soit ${\mathcal C}$ une cat\'egorie cart\'esienne.

\smallskip

Alors:

\begin{listeimarge}

\item On appelle mod\`ele de ${\mathbb T}$ (ou ${\mathbb T}$-mod\`ele) dans ${\mathcal C}$ toute $\Sigma$-structure de ${\mathcal C}$
$$
M
$$
telle que,

\medskip

$\left\{\begin{matrix}
\bullet &\mbox{pour tout axiome de ${\mathbb T}$ qui est de Horn} \hfill \\
{ \ } \\
&\varphi \vdash_{\vec x} \psi \\
{ \ } \\
&\mbox{dans un contexte $\vec x = (x_1^{A_1} \cdots x_n^{A_n})$, on a la relation d'inclusion} \hfill \\
{ \ } \\
&M\varphi \leq M\psi \\
{ \ } \\
&\mbox{entre sous-objets de $M\!A_1 \times \cdots \times M\!A_n$,} \hfill \\
\bullet &\mbox{pour tout axiome de ${\mathbb T}$ qui a la forme} \hfill \\
{ \ } \\
&f_1 (\vec x) = g_1 (\vec x) \wedge \cdots \wedge f_k (\vec x) = g_k (\vec x) \dashv \, \vdash_{\vec x}  (\exists \, y) (u_1 (y) = x_1^{A_1} \wedge \cdots \wedge u_n (y) = x_n^{A_n}) \\
{ \ } \\
&\mbox{et est alors accompagn\'e par un axiome} \hfill \\
{ \ } \\
&u_1 (y) = u_1 (y') \wedge \cdots \wedge u_n (y) = u_n (y') \vdash_{y,y'} y = y' \\
{ \ } \\
&\mbox{(o\`u $y$ et $y'$ sont affect\'ees \`a une sorte $B$) comme dans la d\'efinition {\rm V.4.6}, le morphisme} \hfill \\
{ \ } \\
&M\!u_1 \times \cdots \times M\!u_n : M\!B \longrightarrow M\!A_1 \times \cdots \times M\!A_n \\
{ \ } \\
&\mbox{identifie $M\!B$ au sous-objet de $M\!A_1 \times \cdots \times M\!A_n$ d\'efini par les \'equations} \hfill \\
{ \ } \\
&M\!f_1 = M\!g_1 , \cdots , M\!f_k = M\!g_k \, .
\end{matrix} \right.
$

\medskip

\item On appelle cat\'egorie des mod\`eles de ${\mathbb T}$ (ou des ${\mathbb T}$-mod\`eles) dans ${\mathcal C}$ la sous-cat\'egorie pleine
$$
{\mathbb T}\mbox{\rm -mod} ({\mathcal C}) \qquad \mbox{de} \qquad \Sigma\mbox{\rm -str} ({\mathcal C})
$$
constitu\'ee des mod\`eles de ${\mathcal C}$.
\end{listeimarge}
\end{defn}

\begin{remarksqed}
\begin{listeisansmarge}
\item Si ${\mathbb T}$ est une th\'eorie cart\'esienne, tout foncteur cart\'esien
$$
F : {\mathcal C} \longrightarrow {\mathcal D}
$$
entre deux cat\'egories cart\'esiennes induit un foncteur
$$
F : {\mathbb T}\mbox{-mod} ({\mathcal C}) \longrightarrow {\mathbb T}\mbox{-mod} ({\mathcal D})
$$
qui est la restriction du foncteur
$$
F : \Sigma\mbox{-str} ({\mathcal C}) \longrightarrow \Sigma\mbox{-str} ({\mathcal D}) \, .
$$

\item En particulier, pour tout morphisme de topos
$$
f = (f^* , f_*) : {\mathcal E}' \longrightarrow {\mathcal E} \, ,
$$
sa composante d'image r\'eciproque $f^* : {\mathcal E} \to {\mathcal E}'$ et sa composante d'image directe $f_* : {\mathcal E}' \to {\mathcal E}$ d\'efinissent deux foncteurs
$$
\ f^* : {\mathbb T}\mbox{-mod} ({\mathcal E}) \longrightarrow {\mathbb T}\mbox{-mod} ({\mathcal E}') \, ,
$$
$$
f_* : {\mathbb T}\mbox{-mod} ({\mathcal E}') \longrightarrow {\mathbb T}\mbox{-mod} ({\mathcal E})
$$
pour toute th\'eorie cart\'esienne ${\mathbb T}$. 
\end{listeisansmarge}
\end{remarksqed}

\bigskip

Les propri\'et\'es d'adjonction entre foncteurs sont pr\'eserv\'ees en passant \`a des cat\'egories de mod\`eles:

\begin{cor}\label{corV412}

Soit $\Sigma$ une signature.

\smallskip

Soient deux cat\'egories alg\'ebriques ${\mathcal C}$ et ${\mathcal D}$ reli\'ees par une paire de foncteurs adjoints
$$
\left( {\mathcal C} \xrightarrow{ \ F \ } {\mathcal D} \, , \ {\mathcal D} \xrightarrow{ \ G \ } {\mathcal C} \right)
$$
qui sont alg\'ebriques et respectent les monomorphismes si $\Sigma$ a des symboles de relations.

\smallskip

Alors:

\begin{listeimarge}

\item Les deux foncteurs induits
$$
F : \Sigma\mbox{\rm -str} ({\mathcal C}) \longrightarrow \Sigma\mbox{\rm -str} ({\mathcal D})
$$
et
$$
G : \Sigma\mbox{\rm -str} ({\mathcal D}) \longrightarrow \Sigma\mbox{\rm-str} ({\mathcal C})
$$
forment une paire de foncteurs adjoints.

\medskip

\item Si ${\mathbb T}$ est une th\'eorie alg\'ebrique, les deux foncteurs induits
$$
\ F : {\mathbb T}\mbox{\rm -mod} ({\mathcal C}) \longrightarrow {\mathbb T}\mbox{\rm -mod} ({\mathcal D}) \, ,
$$
$$
G : {\mathbb T}\mbox{\rm -mod} ({\mathcal D}) \longrightarrow {\mathbb T}\mbox{\rm -mod} ({\mathcal C})
$$
forment une paire de foncteurs adjoints.

\medskip

\item Si ${\mathbb T}$ est une th\'eorie de Horn ou plus g\'en\'eralement cart\'esienne [resp. r\'eguli\`ere, resp. coh\'erente, resp. g\'eom\'etrique, resp. du premier ordre finitaire] et $F,G$ sont des foncteurs cart\'esiens [resp. r\'eguliers, resp. coh\'erents, resp. g\'eom\'etriques, resp. de Heyting] entre deux cat\'egories cart\'esiennes [resp. r\'eguli\`eres, resp. coh\'erentes, resp. g\'eom\'etriques, resp. de Heyting], les deux foncteurs induits
$$
\ F : {\mathbb T}\mbox{\rm -mod} ({\mathcal C}) \longrightarrow {\mathbb T}\mbox{\rm -mod} ({\mathcal D}) \, ,
$$
$$
G : {\mathbb T}\mbox{\rm -mod} ({\mathcal D}) \longrightarrow {\mathbb T}\mbox{\rm -mod} ({\mathcal C})
$$
forment une paire de foncteurs adjoints.
\end{listeimarge}
\end{cor}

\begin{remarks}
\begin{listeisansmarge}
\item Comme $G$ est un adjoint \`a droite, il respecte automatiquement toutes les limites, en particulier les produits et les monomorphismes.

\medskip

\item En particulier, pour tout morphisme de topos
$$
f = (f^* , f_*) : {\mathcal E}' \longrightarrow {\mathcal E}
$$
et pour toute signature $\Sigma$, les deux foncteurs
$$
\ f^* : \Sigma\mbox{-str} ({\mathcal E}) \longrightarrow \Sigma\mbox{-str} ({\mathcal E}') \, , 
$$
$$
f_* : \Sigma\mbox{-str} ({\mathcal E}') \longrightarrow \Sigma\mbox{-str} ({\mathcal E})
$$
sont adjoints.

\smallskip

Ils induisent une paire de foncteurs adjoints
$$
\ f^* : {\mathbb T}\mbox{-mod} ({\mathcal E}) \longrightarrow {\mathbb T}\mbox{-mod} ({\mathcal E}') \, ,
$$
$$
f_* : {\mathbb T}\mbox{-mod} ({\mathcal E}') \longrightarrow {\mathbb T}\mbox{-mod} ({\mathcal E})
$$
pour toute th\'eorie ${\mathbb T}$ de signature $\Sigma$ qui est alg\'ebrique, de Horn ou plus g\'en\'eralement cart\'esienne.
\end{listeisansmarge}
\end{remarks}

\bigskip

\begin{demo}
\begin{listeisansmarge}
\item[(i)] Consid\'erons une $\Sigma$-structure $M$ de ${\mathcal C}$ et une $\Sigma$-structure $N$ de ${\mathcal D}$. 

\smallskip

Pour toute sorte $A$ de $\Sigma$, se donner un morphisme
$$
F(M\!A) \longrightarrow N\!A
$$
\'equivaut \`a se donner un morphisme
$$
M\!A \longrightarrow G(N\!A) \, .
$$

Consid\'erons donc deux familles qui se correspondent de morphismes
$$
F(M\!A) \xrightarrow{ \ u_A \ } N\!A \qquad \mbox{et} \qquad M\!A \xrightarrow{ \ v_A \ } G(N\!A)
$$
index\'es par les sortes $A$ de $\Sigma$.

\smallskip

Les foncteurs $F$ et $G$ respectent les produit finis, donc pour tout symbole de fonction de $\Sigma$
$$
f : A_1 \cdots A_n \longrightarrow B \, ,
$$
le carr\'e
$$
\xymatrix{
F(M\!A_1) \times \cdots \times F(M\!A_n) \ar[d]_-{u_{A_1} \times \cdots \times u_{A_n}} \ar[rr]^-{F(M\!f)} &&F(M\!B) \ar[d]^-{u_B} \\
N\!A_1 \times \cdots \times N\!A_n \ar[rr]^-{N\!f} &&N\!B
}
$$
est commutatif si et seulement si le carr\'e
$$
\xymatrix{
M\!A_1 \times \cdots \times M\!A_n \ar[d]_-{v_{A_1} \times \cdots \times v_{A_n}} \ar[rr]^-{M\!f} &&F(M\!B) \ar[d]^-{v_B} \\
G(N\!A_1) \times \cdots \times G(N\!A_n) \ar[rr]^-{G(N\!f)} &&G(N\!B)
}
$$
est commutatif.

\smallskip

De m\^eme, pour tout symbole de relation de $\Sigma$
$$
\xymatrix{
R \ \ \ar@{>->}[r] & \, A_1 \cdots A_n \, ,
}
$$
les monomorphismes
$$
\ M\!R \xhookrightarrow{ \ { \ } \ } M\!A_1 \times \cdots \times M\!A_n \, ,
$$
$$
N\!R \xhookrightarrow{ \ { \ } \ } N\!A_1 \times \cdots \times N\!A_n 
$$
s'inscrivent dans un carr\'e commutatif
$$
\xymatrix{
F(M\!R) \ \ar[d] \ar@{^{(}->}[r] &F(M\!A_1) \times \cdots \times F(M\!A_n) \ar[d]^-{u_{A_1} \times \cdots \times u_{A_n}} \\
N\!R \ \ar@{^{(}->}[r] &N\!A_1 \times \cdots \times N\!A_n
}
$$
si et seulement si ils s'inscrivent dans un carr\'e commutatif
$$
\xymatrix{
M\!R \ \ar[d] \ar@{^{(}->}[r] &M\!A_1 \times \cdots \times M\!A_n \ar[d]^-{v_{A_1} \times \cdots \times v_{A_n}} \\
G(N\!R) \ \ar@{^{(}->}[r] &G(N\!A_1) \times \cdots \times G(N\!A_n)
}
$$
dont le morphisme
$$
M\!R \longrightarrow G(N\!R)
$$
correspond par adjonction au morphisme
$$
F(M\!R) \longrightarrow N\!R \, .
$$

\item[(ii) et (iii)] sont cons\'equences de (i) puisque les cat\'egories de mod\`eles
$$
{\mathbb T}\mbox{-mod} ({\mathcal C}) \qquad \mbox{et} \qquad {\mathbb T}\mbox{-mod} ({\mathcal D}) 
$$
sont par d\'efinition des sous-cat\'egories pleines de
$$
\Sigma\mbox{-str} ({\mathcal C}) \qquad \mbox{et} \qquad \Sigma\mbox{-str} ({\mathcal D}) 
$$
et que, d'apr\`es les diff\'erentes hypoth\`eses, les deux foncteurs adjoints
$$
\ F : \Sigma\mbox{-str} ({\mathcal C}) \longrightarrow \Sigma\mbox{-str} ({\mathcal D}) \, ,
$$
$$
G : \Sigma\mbox{-str} ({\mathcal D}) \longrightarrow \Sigma\mbox{-str} ({\mathcal C})
$$
se restreignent \`a ces cat\'egories de mod\`eles. 
\end{listeisansmarge}
\end{demo}

\bigskip

Dans le cas du topos des faisceaux sur un site, les $\Sigma$-structures ou les mod\`eles de th\'eories cart\'esiennes de ce topos ne sont autres que les faisceaux de $\Sigma$-structures ou de mod\`eles ensemblistes:

\begin{prop}\label{propV413}

Soit ${\mathcal E} = \widehat{\mathcal C}_J$ le topos des faisceaux sur un site $({\mathcal C},J)$.

\smallskip

Alors

\begin{listeimarge}

\item Pour toute signature $\Sigma$, la cat\'egorie des $\Sigma$-structures de ${\mathcal E}$ s'identifie \`a la cat\'egorie des foncteurs contravariants
$$
\begin{matrix}
M_{\bullet} : {\mathcal C}^{\rm op} &\longrightarrow &\Sigma\mbox{\rm -str} ({\rm Ens}) \, , \\
\hfill X &\longmapsto &M_X  \hfill
\end{matrix}
$$
tels que

\medskip

$\left\{\begin{matrix}
\bullet &\mbox{pour toute sorte $A$ de $\Sigma$, le pr\'efaisceau induit} \hfill \\
{ \ } \\
&\begin{matrix}
M_{\bullet} A : {\mathcal C}^{\rm op} &\longrightarrow &{\rm Ens} \, , \hfill \\
\hfill X &\longmapsto &M_X A \hfill
\end{matrix} \\
{ \ } \\
&\mbox{est un faisceau pour la topologie $J$,} \hfill \\
{ \ } \\
\bullet &\mbox{pour tout symbole de relation de $\Sigma$} \hfill \\
{ \ } \\
&\xymatrix{R \ \ \ar@{>->}[r] & \, A_1 \cdots A_n \, ,} \\
{ \ } \\
&\mbox{le sous-pr\'efaisceau de $M_{\bullet} A_1 \times \cdots \times M_{\bullet} A_n$} \hfill \\
{ \ } \\
&\begin{matrix}
M_{\bullet} R : {\mathcal C}^{\rm op} &\longrightarrow &{\rm Ens} \, , \hfill \\
\hfill X &\longmapsto &M_X R \hfill 
\end{matrix} \\
{ \ } \\
&\mbox{est un faisceau pour la topologie $J$.} \hfill
\end{matrix} \right.
$

\bigskip

\item Pour toute th\'eorie ${\mathbb T}$ de signature $\Sigma$ qui est alg\'ebrique, de Horn ou plus g\'en\'eralement cart\'esienne, la cat\'egorie des ${\mathbb T}$-mod\`eles de ${\mathcal E} = \widehat{\mathcal C}_J$ s'identifie \`a la cat\'egorie des foncteurs contravariants
$$
\begin{matrix}
M_{\bullet} : {\mathcal C}^{\rm op} &\longrightarrow &{\mathbb T}\mbox{\rm -mod} ({\rm Ens}) \, , \\
\hfill X &\longmapsto &M_X \hfill 
\end{matrix} 
$$
qui satisfont les deux propri\'et\'es de {\rm (i)} pour toute sorte $A$ et tout symbole de relation $R \rightarrowtail A_1 \cdots A_n$ de $\Sigma$.
\end{listeimarge}
\end{prop}

\begin{remark}

Ce r\'esultat g\'en\'eralise ce qui avait \'et\'e d\'ej\`a observ\'e au sujet des types de structures alg\'ebriques les plus courants:

\smallskip

Les mono{\"\i}des [resp. groupes, resp. anneaux, resp. modules sur un anneau] internes d'un topos de faisceaux ${\mathcal E} = \widehat{\mathcal C}_J$ sur un site $({\mathcal C},J)$ peuvent \^etre vus comme des faisceaux de mono{\"\i}des [resp. groupes, resp. anneaux, resp. modules sur un anneau] sur le site $({\mathcal C},J)$.
\end{remark}

\bigskip

\begin{demo}
\begin{listeisansmarge}
\item Si ${\mathcal E} = \widehat{\mathcal C}$ est la cat\'egorie des pr\'efaisceaux sur ${\mathcal C}$, cela r\'esulte de ce que
\begin{enumerate}
\item[$\bullet$] le produit de pr\'efaisceaux $M_1 , \ldots , M_n$ sur ${\mathcal C}$ est le pr\'efaisceau des produits
$$
\begin{matrix}
M_1 \times \cdots \times M_n : {\mathcal C}^{\rm op} &\longrightarrow &{\rm Ens} \, , \hfill \\
\hfill X &\longmapsto &M_1(X) \times \cdots \times M_n(X) \, , 
\end{matrix} 
$$
\item[$\bullet$] un sous-objet d'un pr\'efaisceau $M$ sur ${\mathcal C}$
$$
R \xhookrightarrow{ \ { \ } \ } M
$$
est une famille compatible de sous-ensembles
$$
R(X) \subseteq M(X) \, , \qquad X \in {\rm Ob} ({\mathcal C}) \, .
$$
\end{enumerate}

\medskip

\item r\'esulte de (i) puisque
\begin{enumerate}
\item[$\bullet$] si ${\mathbb T}$ est une th\'eorie alg\'ebrique, de Horn ou plus g\'en\'eralement cart\'esienne, ${\mathbb T}\mbox{-mod} ({\mathcal E})$ est la sous-cat\'egorie pleine de $\Sigma\mbox{-str} ({\mathcal E})$ constitu\'ee des objets qui satisfont un certain nombre de conditions s'exprimant en termes de limites (finies),
\item[$\bullet$] la cat\'egorie $\widehat{\mathcal C}_J$ des faisceaux sur $({\mathcal C},J)$ est d\'efinie comme la sous-cat\'egorie pleine de $\widehat{\mathcal C}$ constitu\'ee des pr\'efaisceaux qui satisfont un certain nombre de conditions -- celles qui d\'efinissent la notion de faisceau -- s'exprimant en termes de limites,
\item[$\bullet$] les foncteurs de limites respectent toujours les limites. 
\end{enumerate}
\end{listeisansmarge}
\end{demo}

\bigskip

Terminons ce paragraphe en exhibant un processus canonique qui permet d'associer \`a toute th\'eorie du premier ordre finitaire [resp. du premier ordre infinitaire] une th\'eorie coh\'erente [resp. g\'eom\'etrique] qui a les m\^emes mod\`eles ensemblistes et, plus g\'en\'eralement, les m\^emes mod\`eles dans n'importe quel topos ``bool\'een'' au sens suivant:

\begin{defn}\label{defV414}

Un topos ${\mathcal E}$ est dit ``bool\'een'' si, pour tout sous-objet $S$ d'un objet $E$ de ${\mathcal E}$, on a
$$
S \vee (\neg \, S) = E \, .
$$
\end{defn}

\begin{remarksqed}
\begin{listeisansmarge}
\item Autrement dit, un topos est ``bool\'een'' si la ``loi du tiers exclu'' y est satisfaite.

\medskip

\item Le topos Ens des ensembles est bool\'een.

\smallskip

Plus g\'en\'eralement, pour tout groupe $G$, le topos $BG$ des actions de $G$ est bool\'een. 
\end{listeisansmarge}
\end{remarksqed}

\bigskip

Enon\c cons donc:

\begin{prop}\label{prodV415}

Il existe un proc\'ed\'e g\'en\'eral de transformation des formules composant les axiomes d'une th\'eorie du premier ordre finitaire [resp. du premier ordre infinitaire] qui permet d'associer \`a toute telle th\'eorie ${\mathbb T}$

\medskip

$\left\{\begin{matrix}
\mbox{une th\'eorie coh\'erente [resp. g\'eom\'etrique] ${\mathbb T}^m$ de signature $\Sigma^m$} \hfill \\
\mbox{(compos\'ee de $\Sigma$ et de symboles de relations compl\'ementaires)} \hfill
\end{matrix} \right.$

\medskip

\noindent telle que:

\medskip

$\left\{\begin{matrix}
\mbox{pour tout topos bool\'een ${\mathcal E}$,} \hfill \\
\mbox{la cat\'egorie ${\mathbb T}^m${\rm-mod}$({\mathcal E})$} \hfill \\
\mbox{est une sous-cat\'egorie de ${\mathbb T}${\rm-mod}$({\mathcal E})$} \hfill \\
\mbox{qui a les m\^emes objets.} \hfill
\end{matrix} \right.$
\end{prop}

\begin{remarks}
\begin{listeisansmarge}
\item En revanche, et y compris dans le cas ${\mathcal E} = {\rm Ens}$, la cat\'egorie
$$
{\mathbb T}^m\mbox{-mod} ({\mathcal E})
$$
a en g\'en\'eral moins de morphismes que la cat\'egorie
$$
{\mathbb T}\mbox{-mod} ({\mathcal E}) \, .
$$

\item La th\'eorie ${\mathbb T}^m$ est appel\'ee la ``morleyisation'' de ${\mathbb T}$.
\end{listeisansmarge}
\end{remarks}

\bigskip

\begin{demo}

Les axiomes de ${\mathbb T}$ sont des s\'equents de $\Sigma$
$$
\varphi_1 \vdash_{\vec x} \varphi_2
$$
dans la composition desquels peuvent appara{\^\i}tre des symboles $\forall , \Rightarrow , \neg$ [resp. et $\underset{i \in I}{\bigwedge}$ si ${\mathbb T}$ est infinitaire] qu'il s'agit d'\'eliminer pour transformer ${\mathbb T}$ en une th\'eorie ${\mathbb T}^m$ coh\'erente [resp. g\'eom\'etrique].

\smallskip

Dans la cat\'egorie des ensembles ou plus g\'en\'eralement dans tout topos bool\'een, une formule de la forme
$$
(\forall \, y_1^{B_1} \cdots y_k^{B_k}) \, \varphi (x_1^{A_1} \cdots x_n^{A_n} \, y_1^{B_1} \cdots y_k^{B_k})
$$

\noindent a la m\^eme interpr\'etation que la formule
$$
\neg \, (\exists \, y_1^{B_1} \cdots y_k^{B_k}) (\neg \, \varphi) (x_1^{A_1} \cdots x_n^{A_n} \, y_1^{B_1} \cdots y_k^{B_k})
$$

\noindent et une formule de la forme
$$
\varphi_1 \Longrightarrow \varphi_2
$$

\noindent a la m\^eme interpr\'etation que la formule
$$
(\neg \, \varphi_1) \vee \varphi_2 \, .
$$

De m\^eme, une formule de conjonction infinie
$$
\bigwedge_{i \in I} \varphi_i
$$

\noindent a la m\^eme interpr\'etation que la formule
$$
\neg \left( \bigvee_{i \in I} \neg \, \varphi_i \right) .
$$

Op\'erant ces substitutions, on remplace ${\mathbb T}$ par une th\'eorie ${\mathbb T}'$ de m\^eme signature $\Sigma$ dans laquelle le seul symbole restant \`a \'ecarter est $\neg$ et telle que, pour tout topos bool\'een ${\mathcal E}$,
$$
{\mathbb T}'\mbox{-mod} ({\mathcal E}) = {\mathbb T}\mbox{-mod} ({\mathcal E}) \, .
$$

Pour \'eliminer \`a leur tour les symboles $\neg$, rempla\c cons chaque n\'egation
$$
\neg \, \varphi (\vec x)
$$

\noindent d'une formule $\varphi$ dans un contexte $\vec x = (x_1^{A_1} \cdots x_n^{A_n})$ par un symbole de relation
$$
\xymatrix{
N_{\varphi} \ \, \ar@{>->}[r] &A_1 \cdots A_n
}
$$

\noindent ajout\'e \`a la signature $\Sigma$ de ${\mathbb T}$ et ${\mathbb T}'$ et soumis aux axiomes suppl\'ementaires
$$
{\mathcal A}_{\varphi} = \left\{\begin{matrix}
\hfill \top &\vdash_{\vec x} &\varphi \vee N_{\varphi} \, , \\
\varphi \wedge N_{\varphi} &\vdash_{\vec x} &\perp \, . \hfill
\end{matrix}\right.
$$

Soit $\Sigma^m$ la signature d\'eduite de $\Sigma$ en lui adjoignant les symboles de relation $N_{\varphi}$ associ\'es \`a toutes les n\'egations
$$
\neg \, \varphi (\vec x)
$$
qui apparaissent dans la construction des axiomes de ${\mathbb T}'$.

\smallskip

Soit alors ${\mathbb T}^m$ la th\'eorie de signature $\Sigma^m$ dont les axiomes sont la r\'eunion des axiomes ${\mathcal A}_{\varphi}$ ci-dessus et des axiomes d\'eduits de ceux de ${\mathbb T}'$ en rempla\c cant chaque n\'egation
$$
\neg \, \varphi (\vec x)
$$
par le symbole de relation $N_{\varphi} (\vec x)$.

\smallskip

Toute $\Sigma$-structure $M$ dans un topos bool\'een ${\mathcal E}$ se compl\`ete en une unique $\Sigma^m$-structure telle que, pour toute $\varphi$ comme ci-dessus, $\neg \, \varphi (\vec x)$ et $N_{\varphi} (\vec x)$ aient la m\^eme interpr\'etation dans $M$.

\smallskip

Il en r\'esulte que les mod\`eles de ${\mathbb T}'$ dans ${\mathcal E}$ s'identifient aux mod\`eles de ${\mathbb T}^m$ dans ${\mathcal E}$.

\smallskip

Enfin, \'etant donn\'es deux mod\`eles $M$ et $M'$ dans ${\mathcal E}$, les morphismes
$$
M \longrightarrow M' \qquad \mbox{de} \qquad {\mathbb T}^m\mbox{-mod} ({\mathcal E})
$$
sont ceux des morphismes
$$
M \longrightarrow M' \qquad \mbox{de} \qquad {\mathbb T}'\mbox{-mod} ({\mathcal E})
$$
qui envoient chaque sous-objet
$$
M\!N_{\varphi} \xhookrightarrow{ \ { \ } \ } M\!A_1 \times \cdots \times M\!A_n
$$
dans le sous-objet
$$
M'N_{\varphi} \xhookrightarrow{ \ { \ } \ } M' A_1 \times \cdots \times M'A_n \, .
$$

Ainsi,
$$
{\mathbb T}^m\mbox{-mod} ({\mathcal E})
$$
est une sous-cat\'egorie de 
$$
{\mathbb T}'\mbox{-mod} ({\mathcal E}) = {\mathbb T}\mbox{-mod} ({\mathcal E})
$$
qui a les m\^emes objets mais a en g\'en\'eral moins de morphismes.

\smallskip

Cela termine la preuve. 

\end{demo}

\section{D\'emontrabilit\'e}\label{sec55}

\subsection{Les r\`egles d'inf\'erence de la logique du premier ordre}\label{subsec551}

Un syst\`eme logique consiste en une liste de r\`egles d'inf\'erence qui permettent de d\'emontrer des s\'equents \`a partir des axiomes d'une th\'eorie donn\'ee.

\smallskip

Voici la liste des r\`egles d'inf\'erence qui d\'efinit la logique du premier ordre:

\begin{prop}\label{propV51}

Soit ${\mathbb T}$ une th\'eorie du premier ordre de signature $\Sigma$.

\smallskip

Un s\'equent $ \ \varphi \vdash_{\vec x} \ \psi$ de $\Sigma$ est dit ``d\'emontrable'' dans la th\'eorie ${\mathbb T}$ s'il peut \^etre d\'eduit des axiomes de ${\mathbb T}$ en appliquant les r\`egles d'inf\'erence suivantes:

\begin{listeimarge}

\item[(1)] La r\`egle des coupures:

\medskip

$\left\{\begin{matrix}
\bullet &\mbox{Deux s\'equents de la forme} \hfill \\
{ \ } \\
&\varphi_1 \vdash_{\vec x} \varphi_2 \qquad \mbox{et} \qquad \varphi_2 \vdash_{\vec x} \varphi_3 \\
{ \ } \\
&\mbox{impliquent le s\'equent} \hfill \\
{ \ } \\
&\varphi_1 \vdash_{\vec x} \varphi_3 \, .
\end{matrix} \right.
$

\bigskip

\item[(2)] La r\`egle d'identit\'e:

\medskip

$\left\{\begin{matrix}
\bullet &\mbox{Pour tout terme $f$, le s\'equent} \hfill \\
{ \ } \\
&\top \vdash_{\vec x} (f=f)\\
{ \ } \\
&\mbox{est un axiome implicite de n'importe quelle th\'eorie.} \hfill 
\end{matrix} \right.
$

\bigskip

\item[(3)] Les r\`egles des \'egalit\'es:

\medskip

$\left\{\begin{matrix}
\bullet &\mbox{Un s\'equent de la forme} \hfill \\
{ \ } \\
&\top \vdash_{\vec x} f_1 = f_2 \\
{ \ } \\
&\mbox{implique le s\'equent} \hfill \\
{ \ } \\
&\top \vdash_{\vec x} f_2 = f_1 \, . \\
{ \ } \\
\bullet &\mbox{Deux s\'equents de la forme} \hfill \\
{ \ } \\
&\top \vdash_{\vec x} f_1 = f_2 \qquad \mbox{et} \qquad \top \vdash_{\vec x} f_2 = f_3 \\
{ \ } \\
&\mbox{impliquent le s\'equent} \hfill \\
{ \ } \\
&\top \vdash_{\vec x} f_1 = f_3 \, .
\end{matrix} \right.
$

\bigskip

\item[(4)] Les r\`egles de substitution:

\medskip

$\left\{\begin{matrix}
\bullet &\mbox{Si $f_1,f_2$ sont deux termes de m\^eme contexte, et $f'_1,f'_2$ sont deux termes d\'eduits de $f_1,f_2$} \hfill \\
&\mbox{par substitution d'un terme $f$ \`a une variable [resp. d\'eduits d'un terme $f$ par substitution} \hfill \\
&\mbox{d'une variable par $f_1$ et $f_2$], alors le s\'equent} \hfill \\
{ \ } \\
&\top \vdash \ f_1 = f_2 \\
&\mbox{implique le s\'equent} \hfill \\
&\top \vdash \ f'_1 = f'_2 \, . \\
{ \ } \\
\bullet &\mbox{Si $f_1,f_2$ sont deux termes de m\^eme contexte, $R$ est une relation et $R_1,R_2$ sont les deux} \hfill \\
&\mbox{relations d\'eduites de $R$ par substitution d'une variable par $f_1$ et $f_2$, alors le s\'equent} \hfill \\
{ \ } \\
&\top \vdash \ f_1 = f_2 \\
&\mbox{implique le double s\'equent} \hfill \\
&R_1 \dashv \, \vdash R_2 \, . 
\end{matrix} \right.
$

\bigskip

\item[(5)] Les r\`egles des conjonctions finitaires [resp. infinitaires]:

\medskip

$\left\lmoustache\begin{matrix}
\bullet &\mbox{Pour toute formule $\varphi$ dans un contexte $\vec x$, le s\'equent} \hfill \\
{ \ } \\
&\varphi \vdash_{\vec x} \top \\
{ \ } \\
&\mbox{est un axiome implicite de n'importe quelle th\'eorie.} \hfill \\
{ \ } \\
\end{matrix} \right.
$

$\left\rmoustache\begin{matrix}
\bullet &\mbox{Pour toute famille finie $\varphi_1 , \cdots , \varphi_k$ [resp. infinie $(\varphi_i)_{i \in I}$] de formules de m\^eme contexte $\vec x$,} \hfill \\
&\mbox{et pour toute formule $\varphi$ dans le contexte $\vec x$, la famille de s\'equents} \hfill \\
{ \ } \\
&\varphi \vdash_{\vec x} \varphi_i \, , \quad 1 \leq i \leq k \qquad \mbox{[resp.} \quad i \in I \mbox{]} \\
&\mbox{implique le s\'equent} \hfill \\
&\varphi \vdash_{\vec x} \varphi_1 \wedge \cdots \wedge \varphi_k \qquad \mbox{[resp.} \quad \varphi \vdash_{\vec x} \underset{i \in I}{\bigwedge} \, \varphi_i \mbox{]}. \\
{ \ } \\

\bullet &\mbox{Pour toute telle famille finie $\varphi_1 , \cdots , \varphi_k$ [resp. infinie $(\varphi_i)_{i \in I}$], les s\'equents} \hfill \\
{ \ } \\
&\varphi_1 \wedge \cdots \wedge \varphi_k \vdash_{\vec x} \varphi_{i_0} \, , \quad 1 \leq i_0 \leq k \, , \\
{ \ } \\
&\mbox{[resp.} \qquad \underset{i \in I}{\bigwedge} \, \varphi_i \vdash_{\vec x} \varphi_{i_0} \, , \quad i_0 \in I \ \mbox{]} \\
{ \ } \\
&\mbox{sont des axiomes implicites de n'importe quelle th\'eorie.} \hfill
\end{matrix} \right.
$

\bigskip

\item[(6)] Les r\`egles des disjonctions finitaires [resp. infinitaires]:

\medskip

$\left\lmoustache\begin{matrix}
\bullet &\mbox{Pour toute formule $\varphi$ dans un contexte $\vec x$, le s\'equent} \hfill \\
{ \ } \\
&\perp \, \vdash_{\vec x} \varphi \\
{ \ } \\
&\mbox{est un axiome implicite de n'importe quelle th\'eorie.} \hfill \\
{ \ } \\
\bullet &\mbox{Pour toute famille finie $\varphi_1 , \cdots , \varphi_k$ [resp. infinie $(\varphi_i)_{i \in I}$] de formules de m\^eme contexte $\vec x$,} \hfill \\
&\mbox{et pour toute formule $\varphi$ dans le contexte $\vec x$, la famille de s\'equents} \hfill \\
{ \ } \\
&\varphi_i \vdash_{\vec x} \varphi \, , \quad 1 \leq i \leq k \qquad \mbox{[resp.} \quad i \in I \mbox{]} 
\end{matrix} \right.
$

$
\left\rmoustache\begin{matrix}
&\mbox{implique le s\'equent} \hfill \\
&\varphi_1 \vee \cdots \vee \varphi_k \vdash_{\vec x} \varphi \qquad \mbox{[resp.} \quad \underset{i \in I}{\bigvee} \, \varphi_i \vdash_{\vec x} \varphi \mbox{]}. \\
{ \ } \\
\bullet &\mbox{Pour toute telle famille finie $\varphi_1 , \cdots , \varphi_k$ [resp. infinie $(\varphi_i)_{i \in I}$], les s\'equents} \hfill \\
{ \ } \\
&\varphi_{i_0} \vdash_{\vec x} \varphi_1 \vee \cdots \vee \varphi_k \, , \quad 1 \leq i_0 \leq k \, , \\
{ \ } \\
&\mbox{[resp.} \qquad \varphi_{i_0}  \vdash_{\vec x} \underset{i \in I}{\bigvee} \, \varphi_i  \, , \quad i_0 \in I \ \mbox{]} \\
{ \ } \\
&\mbox{sont des axiomes implicites de n'importe quelle th\'eorie.} \hfill
\end{matrix} \right.
$

\bigskip

\item[(7)] Les r\`egles de distributivit\'e:

\medskip

$\left\{\begin{matrix}
\bullet &\mbox{Pour toutes formules $\varphi$ et $\varphi_1 , \cdots , \varphi_k$ [resp. $(\varphi_i)_{i \in I}$] de m\^eme contexte $\vec x$, le double s\'equent} \hfill \\
{ \ } \\
&\varphi \wedge (\varphi_1 \vee \cdots \vee \varphi_k) \vdash_{\vec x} (\varphi \wedge \varphi_1) \vee \cdots \vee (\varphi \wedge \varphi_k) \\
{ \ } \\
&\mbox{[resp.} \qquad \varphi \wedge \underset{i \in I}{\bigvee} \, \varphi_i \vdash_{\vec x} \underset{i \in I}{\bigvee} \, \varphi \wedge \varphi_i \ \mbox{]} \\
{ \ } \\
&\mbox{est un axiome implicite de n'importe quelle th\'eorie.} \hfill 
\end{matrix} \right.
$

\bigskip

\item[(8)] Les r\`egles de quantification existentielle:

\medskip

$\left\{\begin{matrix}
\bullet &\mbox{Pour tous contextes disjoints $\vec x$ et $\vec y$, toute formule $\varphi$ de contexte $(\vec x , \vec y)$ et toute formule $\psi$} \hfill \\
&\mbox{de contexte $\vec x$, le s\'equent} \hfill \\
&\varphi \vdash_{\vec x , \vec y} \psi \\
&\mbox{\'equivaut au s\'equent} \hfill \\
&(\exists \, \vec y) \, \varphi \vdash_{\vec x} \psi \, .
\end{matrix} \right.
$

\bigskip

\item[(9)] La r\`egle de Frobenius :

\medskip

$\left\{\begin{matrix}
\bullet &\mbox{Pour toute formule $\varphi$ de contexte $(\vec x , \vec y)$ et toute formule $\psi$ de contexte $\vec x$ comme dans {\rm (8)},} \hfill \\
&\mbox{le s\'equent} \hfill \\
&(\exists \, \vec y) \, \varphi \wedge \psi \vdash_{\vec x} (\exists \, \vec y) (\varphi \wedge \psi) \\
{ \ } \\
&\mbox{est un axiome implicite de toute th\'eorie.} \hfill 
\end{matrix} \right.
$

\bigskip

\item[(10)] Les r\`egles de quantification universelle:

\medskip

$\left\{\begin{matrix}
\bullet &\mbox{Pour toute formule $\varphi$ de contexte $(\vec x , \vec y)$ et toute formule $\psi$ de contexte $\vec x$ comme dans {\rm (8)} et {\rm (9)},} \hfill \\
&\mbox{le s\'equent} \hfill \\
&\psi \vdash_{\vec x , \vec y} \varphi \\
&\mbox{\'equivaut au s\'equent} \hfill \\
&\psi \vdash_{\vec x} (\forall \, \vec y) \, \varphi \, .
\end{matrix} \right.
$

\bigskip

\item[(11)] Les r\`egles d'implication:

\medskip

$\left\{\begin{matrix}
\bullet &\mbox{Pour toutes formules $\varphi , \psi , \chi$ de contexte $\vec x$, le s\'equent} \hfill \\
{ \ } \\
&\varphi \wedge \psi \vdash_{\vec x} \chi \\
{ \ } \\
&\mbox{\'equivaut au s\'equent} \hfill \\
&\varphi \vdash_{\vec x} (\psi \Rightarrow \chi) \, .
\end{matrix} \right.
$

\bigskip

\item[(12)] Les r\`egles de n\'egation:

\medskip

$\left\{\begin{matrix}
\bullet &\mbox{Pour toutes formules $\varphi$ et  $\psi$ de contexte $\vec x$, le s\'equent} \hfill \\
{ \ } \\
&\varphi \wedge \psi \vdash_{\vec x} \perp \\
{ \ } \\
&\mbox{\'equivaut au s\'equent} \hfill \\
&\varphi \vdash_{\vec x} \neg \, \psi \, .
\end{matrix} \right.
$
\end{listeimarge}
\end{prop}

\begin{remarksqed}
\begin{listeisansmarge}
\item En revanche, on ne fait pas figurer dans la liste des r\`egles d'inf\'erence de la logique du premier ordre le ``principe du tiers exclu'' qui s'\'ecrirait:

\medskip

$\left\{\begin{matrix}
\bullet &\mbox{Pour toute formule $\varphi$ dans un contexte $\vec x$, le s\'equent} \hfill \\
{ \ } \\
&\top \vdash_{\vec x} \varphi \vee \neg \, \varphi \\
{ \ } \\
&\mbox{est un axiome implicite de n'importe quelle th\'eorie.} \hfill 
\end{matrix} \right.
$

\bigskip

\noindent La logique qui consiste en les r\`egles d'inf\'erence (1) \`a (12), sans principe du tiers exclu, est appel\'ee logique ``intuitioniste'' du premier ordre.

\medskip

\item Les r\`egles des disjonctions finitaires [resp. infinitaires] (6) assurent que pour toutes formules $\varphi$ et $\varphi_1 , \cdots , \varphi_k$ [resp. $(\varphi_i)_{i \in I}$] de m\^eme contexte $\vec x$, les s\'equents en sens inverse de ceux de (7)
$$
(\varphi \wedge \varphi_1) \vee \cdots \vee (\varphi \wedge \varphi_k) \vdash_{\vec x} \varphi \wedge (\varphi_1 \vee \cdots \vee \varphi_k)
$$
$$
\mbox{[resp.} \qquad \bigvee_{i \in I} \varphi \wedge \varphi_i \vdash_{\vec x} \varphi \wedge \bigvee_{i \in I} \varphi_i \ \mbox{]}
$$
sont d\'emontrables dans n'importe quelle th\'eorie.

\medskip

\item Les r\`egles de (6) et (7) assurent que pour toutes formules $\varphi$ et $\varphi_1 , \cdots , \varphi_k$ dans un contexte $\vec x$, le double s\'equent
$$
(\varphi \vee \varphi_1) \wedge \cdots \wedge (\varphi \vee \varphi_k) \, \dashv \, \vdash_{\vec x} \varphi \vee (\varphi_1 \wedge \cdots \wedge \varphi_k)
$$
est d\'emontrable dans n'importe quelle th\'eorie.

\smallskip

En effet, il suffit de le v\'erifier dans le cas o\`u $k=2$.

\smallskip

Dans ce cas, le s\'equent
$$
(\varphi \vee \varphi_1) \wedge (\varphi \vee \varphi_2) \, \dashv \, \vdash_{\vec x} (\varphi \wedge \varphi) \vee (\varphi_1 \wedge \varphi) \vee (\varphi \wedge \varphi_2) \vee (\varphi_1 \wedge \varphi_2)
$$
est d\'emontrable, donc aussi le s\'equent
$$
(\varphi \vee \varphi_1) \wedge (\varphi \vee \varphi_2) \, \dashv \, \vdash_{\vec x} \varphi \vee (\varphi_1 \wedge \varphi_2)
$$
puisque les s\'equents
$$
\varphi \, \dashv \, \vdash_{\vec x} \varphi \wedge \varphi \, ,
$$
$$
\varphi_1 \wedge \varphi \vdash_{\vec x} \varphi \, ,
$$
et
$$
\varphi \wedge \varphi_2 \vdash_{\vec x} \varphi
$$
sont d\'emontrables. 
\end{listeisansmarge}
\end{remarksqed}

\bigskip

En retenant seulement certaines parties de la liste des r\`egles d'inf\'erence de la logique (intuitioniste) du premier ordre, on d\'efinit diff\'erents fragments de cette logique:

\begin{defn}\label{defV52}
\begin{listeimarge}
\item Les r\`egles d'inf\'erence de base de tous les fragments de logique du premier ordre sont

\medskip

$
\left\{\begin{matrix}
\bullet &\mbox{la r\`egle des coupures {\rm (1)},} \hfill \\
\bullet &\mbox{la r\`egle d'identit\'e {\rm (2)},} \hfill \\
\bullet &\mbox{la r\`egle des \'egalit\'es {\rm (3)},} \hfill \\
\bullet &\mbox{les r\`egles de substitution {\rm (4)}.} \hfill 
\end{matrix}\right.
$

\medskip

\item La logique de Horn consiste en ces r\`egles d'inf\'erence de base compl\'et\'ees par

\medskip

$
\left\{\begin{matrix}
\bullet &\mbox{les r\`egles des conjonctions finitaires {\rm (5)}.} \hfill
\end{matrix}\right.
$

\medskip

\item La logique r\'eguli\`ere consiste en les r\`egles d'inf\'erence de la logique de Horn compl\'et\'ees par

\medskip

$
\left\{\begin{matrix}
\bullet &\mbox{les r\`egles de quantification existentielle {\rm (8)},} \hfill \\
\bullet &\mbox{la r\`egle de Frobenius {\rm (9)}.} \hfill 
\end{matrix}\right.
$

\medskip

\item La logique coh\'erente [resp. g\'eom\'etrique] consiste en les r\`egles d'inf\'erence de la logique r\'eguli\`ere compl\'et\'ees par

\medskip

$
\left\{\begin{matrix}
\bullet &\mbox{les r\`egles des disjonctions finitaires [resp. infinitaires] {\rm (6)},} \hfill \\
\bullet &\mbox{les r\`egles de distributivit\'e finitaires [resp. infinitaires] {\rm (7)}.} \hfill 
\end{matrix}\right.
$

\medskip

\item La logique du premier ordre finitaire consiste en les r\`egles de la logique coh\'erente compl\'et\'ees par

\medskip

$
\left\{\begin{matrix}
\bullet &\mbox{les r\`egles de quantification universelle {\rm (10)},} \hfill \\
\bullet &\mbox{les r\`egles d'implication {\rm (11)},} \hfill \\
\bullet &\mbox{les r\`egles de n\'egation {\rm (12)}.} \hfill
\end{matrix}\right.
$
\end{listeimarge}
\end{defn}

\bigskip

\begin{remarkqed}

On verra plus loin que pour toute th\'eorie ${\mathbb T}$ de signature $\Sigma$ qui est de Horn [resp. r\'eguli\`ere, resp. coh\'erente, resp. g\'eom\'etrique, resp. du premier ordre finitaire], un s\'equent de $\Sigma$ qui est de Horn [resp. r\'egulier, resp. coh\'erent, resp. g\'eom\'etrique, resp. du premier ordre finitaire] est d\'emontrable dans la logique du premier ordre si et seulement si il l'est dans la logique de Horn [resp. r\'eguli\`ere, resp. coh\'erente, resp. g\'eom\'etrique, resp. intuitioniste]. 

\end{remarkqed}

\subsection{V\'erification des propri\'et\'es d\'emontrables d'une th\'eorie dans ses mod\`eles}\label{subsecVV2}

\medskip

Les mod\`eles d'une th\'eorie alg\'ebrique dans une cat\'egorie alg\'ebrique v\'erifient toutes les propri\'et\'es de nature alg\'ebrique qui sont d\'emontrables \`a partir des axiomes de la th\'eorie par les r\`egles d'\'egalit\'e et de substitution:

\begin{prop}\label{propV53}

Soit ${\mathbb T}$ une th\'eorie alg\'ebrique de signature $\Sigma$.

\smallskip

Soit $M$ un mod\`ele de ${\mathbb T}$ dans une cat\'egorie alg\'ebrique ${\mathcal C}$.

\smallskip

Soit un s\'equent de $\Sigma$ de la forme
$$
\top \vdash \, f(\vec x) = g(\vec x)
$$
associ\'e \`a deux termes $f,g$ de $\Sigma$ dans un m\^eme contexte $\vec x$.

\smallskip

Supposons que ce s\'equent est d\'emontrable dans la th\'eorie ${\mathbb T}$ par les r\`egles d'\'egalit\'e {\rm (3)} et les r\`egles de substitution de termes {\rm (4)} de la d\'efinition \ref{defV52}.

\smallskip

Alors le mod\`ele $M$ v\'erifie ce s\'equent au sens que l'on a une \'egalit\'e de morphismes
$$
M\!f = Mg \, .
$$
\end{prop}

\begin{demo}

En effet, si $f_1,f_2$ sont deux termes dans un m\^eme contexte $\vec x$, l'\'egalit\'e de morphismes
$$
M\!f_1 = M\!f_2
$$
entra{\^\i}ne l'\'egalit\'e de morphismes
$$
M\!f_2 = M\!f_1 \, ,
$$
et si $f_1 , f_2 , f_3$ sont trois termes dans un m\^eme contexte $\vec x$, les \'egalit\'es de morphismes
$$
M\!f_1 = M\!f_2 \qquad \mbox{et} \qquad M\!f_2 = M\!f_3
$$
entra{\^\i}nent l'\'egalit\'e de morphismes
$$
M\!f_1 = M\!f_3 \, .
$$

D'autre part, si $f_1 , f_2$ sont deux termes de m\^eme contexte, et $f'_1 , f'_2$ sont deux termes d\'eduits de $f_1 , f_2$ par substitution d'un terme $f$ \`a une variable [resp. d\'eduits d'un terme $f$ par substitution de $f_1$ et $f_2$ \`a une variable], alors l'\'egalit\'e de morphismes
$$
M\!f_1 = M\!f_2
$$
entra{\^\i}ne l'\'egalit\'e de morphismes
$$
M\!f'_1 = M\!f'_2
$$
puisque $M\!f'_1$ et $M\!f'_2$ se d\'eduisent de $M\!f_1$ et $M\!f_2$ par une formule de la forme
$$
\left\{\begin{matrix}
M\!f'_1 &= &M\!f_1 \circ (p \times M\!f \times q) \, , \\
M\!f'_2 &= &M\!f_2 \circ (p \times M\!f \times q) \hfill
\end{matrix}\right.
$$
[resp.
$$
\left\{\begin{matrix}
M\!f'_1 &= &M\!f \circ (p \times M\!f_1 \times q) \, , \hfill \\
M\!f'_2 &= &M\!f \circ (p \times M\!f_2 \times q) \ \mbox{]}\hfill
\end{matrix}\right.
$$
o\`u $p$ et $q$ sont des morphismes de projection d'un produit
$$
M\!A_1 \times \cdots \times M\!A_n
$$
sur deux parties de ses facteurs. 

\end{demo}

\bigskip

De la m\^eme fa\c con, on a pour les diff\'erents types de th\'eories et de logiques que nous avons introduits:

\begin{prop}\label{propV54}

Soit ${\mathbb T}$ une th\'eorie de Horn [resp. r\'eguli\`ere, resp. coh\'erente, resp. g\'eom\'etrique, resp. du premier ordre finitaire, resp. du premier ordre infinitaire] de signature $\Sigma$.

\smallskip

Soit $M$ un mod\`ele de ${\mathbb T}$ dans une cat\'egorie ${\mathcal C}$ suppos\'ee cart\'esienne [resp. r\'eguli\`ere, resp. coh\'erente, resp. g\'eom\'etrique, resp. de Heyting, resp. un topos].

\smallskip

Soit
$$
\varphi \vdash_{\vec x} \psi
$$
un s\'equent de $\Sigma$ dans un contexte $\vec x = (x_1^{A_1} \cdots x_n^{A_n})$, qui est de Horn [resp. r\'egulier, resp. coh\'erent, resp. g\'eom\'etrique, resp. du premier ordre finitaire, resp. du premier ordre infinitaire].

\smallskip

Supposons que ce s\'equent est d\'emontrable dans la th\'eorie ${\mathbb T}$ par les r\`egles d'inf\'erence de la logique de Horn [resp. r\'eguli\`ere, resp. coh\'erente, resp. g\'eom\'etrique, resp. du premier ordre finitaire, resp. du premier ordre infinitaire].

\smallskip

Alors le mod\`ele $M$ v\'erifie ce s\'equent au sens que l'on a la relation d'inclusion entre sous-objets de $M\!A_1 \times \cdots \times M\!A_n$
$$
M\varphi (\vec x) \leq M\psi (\vec x) \, .
$$
\end{prop}

\begin{demo}

Il s'agit de prouver que si un mod\`ele $M$ v\'erifie une famille de s\'equents
$$
\varphi_i \vdash \, \psi_i \, ,
$$
alors tout s\'equent $\varphi \vdash \, \psi$ qui se d\'eduit de cette famille par l'une des r\`egles d'inf\'erence de la logique consid\'er\'ee est \'egalement v\'erifi\'e par le mod\`ele $M$.

\smallskip

Passons en revue les diff\'erentes r\`egles d'inf\'erence de la logique du premier ordre:

\begin{listeisansmarge}

\item[(1)] Pour la r\`egle des coupures, consid\'erons trois formules $\varphi_1 , \varphi_2$ et $\varphi_3$ de $\Sigma$ dans un m\^eme contexte $\vec x = (x_1^{A_1} \cdots x_n^{A_n})$. Si les trois sous-objets
$$
M\varphi_1 , M\varphi_2 \quad \mbox{et} \quad M \varphi_3 \qquad \mbox{de} \qquad M\!A_1 \times \cdots \times M\!A_n
$$
v\'erifient les relations d'inclusion
$$
M\varphi_1 \leq M\varphi_2 \qquad \mbox{et} \qquad M\varphi_2 \leq M\varphi_3 \, ,
$$
alors ils v\'erifient a fortiori la relation d'inclusion
$$
M\varphi_1 \leq M\varphi_3 \, .
$$

\item[(2)] Pour la r\`egle d'identit\'e, il suffit de noter que si $f$ est n'importe quel terme dans un contexte $\vec x$, son interpr\'etation $M\!f$ comme morphisme de ${\mathcal C}$ satisfait automatiquement la relation
$$
M\!f = M\!f \, .
$$

\item[(3)] Pour les r\`egles d'\'egalit\'e, c'est une partie de la proposition \ref{propV53}.

\medskip

\item[(4)] Pour les r\`egles de substitution de variables de termes par des termes, c'est l'autre partie de la proposition \ref{propV53}.

\smallskip

Pour les r\`egles de substitution d'une variable $y_j^{B_j}$ du contexte $\vec y = (y_1^{B_1} \cdots y_k^{B_k})$ d'une formule atomique $R$ par deux termes $f_1 , f_2$ de contexte $\vec x = (x_1^{A_1} \cdots x_n^{A_n})$, consid\'erons les deux formules atomiques associ\'ees
$$
R_1 \left( y_1^{B_1} \cdots y_{j-1}^{B_{j-1}} \ x_1^{A_1} \cdots x_n^{A_n} \ y_{j+1}^{B_{j+1}} \cdots y_k^{B_k} \right)
$$
et
$$
\ R_2 \left( y_1^{B_1} \cdots y_{j-1}^{B_{j-1}} \ x_1^{A_1} \cdots x_n^{A_n} \ y_{j+1}^{B_{j+1}} \cdots y_k^{B_k} \right) \, .
$$
Elles s'interpr\`etent dans le mod\`ele $M$ comme les sous-objets de
$$
M\!B_1 \times \cdots \times M\!B_{j-1} \times M\!A_1 \times \cdots \times M\!A_n \times M\!B_{j+1} \times \cdots \times M\!B_k
$$
d\'eduits du sous-objet
$$
M\!R \xhookrightarrow{ \ { \ } \ } M\!B_1 \times \cdots \times M\!B_k
$$
par produit fibr\'e avec les deux morphismes
$$
{\rm id}_{M\!B_1} \times \cdots \times {\rm id}_{M\!B_{j-1}} \times M\!f_1 \times {\rm id}_{M\!B_{j+1}} \times \cdots \times {\rm id}_{M\!B_k}
$$
et
$$
\ {\rm id}_{M\!B_1} \times \cdots \times {\rm id}_{M\!B_{j-1}} \times M\!f_2 \times {\rm id}_{M\!B_{j+1}} \times \cdots \times {\rm id}_{M\!B_k} \, .
$$

Donc l'\'egalit\'e des morphismes
$$
M\!f_1 = M\!f_2
$$
entra{\^\i}ne l'\'egalit\'e des sous-objets
$$
R_1 = R_2 \, .
$$

\item[(5)] Pour la r\`egle des conjonctions finitaires [resp. infinitaires], consid\'erons n'importe quelle formule $\varphi$ dans un contexte $\vec x = (x_1^{A_1} \cdots x_n^{A_n})$.

\smallskip

On note d'abord que, par d\'efinition de l'interpr\'etation $M\varphi$ de $\varphi$, c'est un sous-objet
$$
M\varphi (\vec x) \xhookrightarrow{ \ { \ } \ } M\!A_1 \times \cdots \times M\!A_n \, ,
$$
ce qui signifie que $M$ v\'erifie le s\'equent
$$
\varphi \vdash_{\vec x} \top \, .
$$

Puis, si les $\varphi_1$, $1 \leq i \leq k$ [resp. $i \in I$], sont une famille finie [resp. infinie] de formules de contexte $\vec x$, on a
$$
M\varphi \leq M\varphi_i
$$
pour tout indice $i$ si et seulement si
$$
M\varphi \leq M\varphi_1 \wedge \cdots \wedge M\varphi_k = M (\varphi_1 \wedge \cdots \wedge \varphi_k)
$$
[resp.
$$
M\varphi \leq \bigwedge_{i \in I} M\varphi_i = M \Bigl(\bigwedge_{i \in I} \varphi_i\Bigl) \ ].
$$

Autrement dit, $M$ v\'erifie les s\'equents
$$
\varphi \vdash_{\vec x} \varphi_i \, , \quad 1 \leq i \leq k \qquad \mbox{[resp.} \quad i \in I \mbox{]}
$$
si et seulement si il v\'erifie le s\'equent
$$
\varphi \vdash_{\vec x} \varphi_1 \wedge \cdots \wedge \varphi_k \qquad \mbox{[resp.} \quad \varphi \vdash_{\vec x} \bigwedge_{i \in I} \varphi_i \ \mbox{]}.
$$
En particulier, on a pour tout indice $i_0$
$$
\varphi_1 \wedge \cdots \wedge \varphi_k \vdash_{\vec x} \varphi_{i_0} \qquad \mbox{[resp.} \quad \bigwedge_{i \in I} \varphi_i \vdash_{\vec x} \varphi_{i_0} \ \mbox{]}.
$$

\item[(6)] De m\^eme, pour les r\`egles de disjonctions finitaires [resp. infinitaires], on note d'abord que l'interpr\'etation $M\varphi (\vec x)$ d'une formule $\varphi$ de contexte $\vec x = (x_1^{A_1} \cdots x_n^{A_n})$ v\'erifie la relation d'inclusion
$$
\emptyset_{M\!A_1 \times \cdots \times M\!A_n} \leq M\varphi (\vec x)
$$
entre sous-objets de $M\!A_1 \times \cdots \times M\!A_n$. Cela signifie que $M$ v\'erifie automatiquement le s\'equent
$$
\perp \, \vdash_{\vec x} \varphi \, .
$$
Puis, si les $\varphi_i$, $1 \leq i \leq k$ [resp. $i \in I$], sont une famille finie [resp. infinie] de formules de contexte $\vec x$, on a
$$
M\varphi_i \leq M\varphi
$$
pour tout indice $i$ si et seulement si
$$
M (\varphi_1 \vee \cdots \vee \varphi_k) = M\varphi_1 \vee \cdots \vee M\varphi_k \leq M\varphi
$$
[resp.
$$
M \Bigl( \bigvee_{i \in I} \varphi_i \Bigl) = \bigvee_{i \in I} M\varphi_i \leq M\varphi \ \mbox{]}.
$$

Autrement dit, $M$ v\'erifie les s\'equents
$$
\varphi_i \vdash_{\vec x} \varphi \, , \quad 1 \leq i \leq k \qquad \mbox{[resp.} \quad i \in I \ \mbox{]}
$$
si et seulement si il v\'erifie le s\'equent
$$
\varphi_1 \vee \cdots \vee \varphi_k \vdash_{\vec x} \varphi \qquad \mbox{[resp.} \quad \bigvee_{i \in I} \varphi_i \vdash_{\vec x} \varphi \ \mbox{]}.
$$
En particulier, on a pour tout indice $i_0$
$$
\varphi_{i_0} \vdash_{\vec x} \varphi_1 \vee \cdots \vee \varphi_k \qquad \mbox{[resp.} \quad \varphi_{i_0} \vdash_{\vec x} \bigvee_{i \in I} \varphi_i \ \mbox{]}.
$$

\item[(7)] Pour les r\`egles de distributivit\'e, consid\'erons encore des formules $\varphi$ et $\varphi_i$, $1 \leq i \leq k$ [resp. $i \in I$], de m\^eme contexte $\vec x = (x_1^{A_1} \cdots x_n^{A_n})$.

\smallskip

Les formules
$$
\varphi \wedge (\varphi_1 \vee \cdots \vee \varphi_k) \qquad \mbox{et} \qquad (\varphi \wedge \varphi_1) \vee \cdots \vee (\varphi \wedge \varphi_k)
$$
[resp.
$$
\varphi \wedge \bigvee_{i \in I} \varphi_i \qquad \mbox{et} \qquad \bigvee_{i \in I} \varphi \wedge \varphi_i \ \mbox{]}
$$
s'interpr\`etent comme les sous-objets de $M\!A_1 \times \cdots \times M\!A_n$
$$
M\varphi \wedge (M\varphi_1 \vee \cdots \vee M\varphi_k) \qquad \mbox{et} \qquad (M\varphi \wedge M\varphi_1) \vee \cdots \vee (M\varphi \wedge M\varphi_k)
$$
[resp.
$$
M\varphi \wedge \bigvee_{i \in I} \varphi_i \qquad \mbox{et} \qquad \bigvee_{i \in I} (M\varphi \wedge M\varphi_i) \ \mbox{]} .
$$

Si ${\mathcal C}$ est une cat\'egorie coh\'erente [resp. g\'eom\'etrique], ces sous-objets co{\"\i}ncident puisque les r\'eunions finies [resp. infinies] de sous-objets sont respect\'ees par les foncteurs de changement de base dans ${\mathcal C}$.

\smallskip

Autrement dit, le mod\`ele $M$ v\'erifie le s\'equent
$$
\varphi \wedge (\varphi_1 \vee \cdots \vee \varphi_k) \vdash_{\vec x} (\varphi \wedge \varphi_1) \vee \cdots \vee (\varphi \wedge \varphi_k)
$$
[resp.
$$
\varphi \wedge \bigvee_{i \in I} \varphi_i \vdash_{\vec x} \bigvee_{i \in I} \varphi \wedge \varphi_i \ \mbox{]}.
$$

\item[(8)] Pour les r\`egles de quantification existentielle, consid\'erons deux contextes disjoints $\vec x = (x_1^{A_1} \cdots x_n^{A_n})$ et $\vec y = (y_1^{B_1} \cdots y_k^{B_k})$.

\smallskip

Si $\varphi$ est une formule de contexte $(\vec x , \vec y)$ et $\psi$ une formule de contexte $\vec x$, le s\'equent
$$
\varphi \vdash_{\vec x , \vec y} \psi
$$
s'interpr\`ete dans le mod\`ele $M$ comme la relation d'inclusion
$$
M\varphi \, (\vec x , \vec y) \leq p^{-1} M\psi \, (\vec x)
$$
en notant $p$ le morphisme de projection
$$
p : M\!A_1 \times \cdots \times M\!A_n \times M\!B_1 \times \cdots \times M\!B_k \longrightarrow M\!A_1 \times \cdots \times M\!A_n
$$
et $p^{-1}$ le foncteur d'image r\'eciproque des sous-objets par $p$.

\smallskip

Par adjonction, cette relation d'inclusion \'equivaut \`a la relation d'inclusion entre sous-objets de $M\!A_1 \times \cdots \times M\!A_n$
$$
\exists_p \, M\varphi \, (\vec x , \vec y) \leq M\psi \, (\vec x)
$$
qui est l'interpr\'etation dans $M$ du s\'equent
$$
(\exists \, \vec y) \, \varphi \vdash_{\vec x} \psi \, .
$$
Autrement dit, $M$ v\'erifie le s\'equent
$$
\varphi \vdash_{\vec x, \vec y} \psi
$$
si et seulement si il v\'erifie le s\'equent
$$
(\exists \, \vec y) \, \varphi \vdash_{\vec x} \psi \, .
$$

\item[(9)] Pour la r\`egle de Frobenius, consid\'erons encore une formule $\varphi$ de contexte $(\vec x , \vec y)$ et une formule $\psi$ de contexte $\vec x$.

\smallskip

Les formules de contexte $\vec x$
$$
(\exists \, \vec y) \, \varphi \wedge \psi \qquad \mbox{et} \qquad (\exists \, \vec y) \, (\varphi \wedge \psi)
$$
s'interpr\`etent comme les sous-objets de $M\!A_1 \times \cdots \times M\!A_n$
$$
\exists_p \, M\varphi \wedge M\psi \qquad \mbox{et} \qquad \exists_p \, (M\varphi \wedge p^{-1} M\psi) \, .
$$

Comme ${\mathcal C}$ est une cat\'egorie r\'eguli\`ere, les foncteurs $\exists_p$ commutent avec les changements de base et donc ces deux sous-objets co{\"\i}ncident.

\smallskip

Cela signifie que $M$ v\'erifie le double s\'equent
$$
(\exists \, \vec y) \, \varphi \wedge \psi \, \dashv \, \vdash_{\vec x} (\exists \, \vec y)(\varphi \wedge \psi) \, .
$$

\item[(10)] Pour les r\`egles de quantification universelle, consid\'erons toujours deux formules $\varphi$ de contexte $(\vec x , \vec y)$ et $\psi$ de contexte $\vec x$.

\smallskip

Le s\'equent
$$
\psi \vdash_{\vec x , \vec y} \varphi
$$
s'interpr\`ete dans le mod\`ele $M$ comme la relation d'inclusion
$$
p^{-1} \, M\psi \, (\vec x) \subseteq M\varphi \, (\vec x , \vec y)
$$
si $p$ d\'esigne toujours la projection
$$
p : M\!A_1 \times \cdots \times M\!A_n \times M\!B_1 \times \cdots \times M\!B_k \longrightarrow M\!A_1 \times \cdots \times M\!A_n \, .
$$

Par adjonction, cette relation d'inclusion \'equivaut \`a la relation d'inclusion entre sous-objets de $M\!A_1 \times \cdots \times M\!A_n$
$$
M\psi \, (\vec x) \leq \forall_p \, M\varphi \, (\vec x , \vec y)
$$
qui est l'interpr\'etation dans $M$ du s\'equent
$$
\psi \vdash_{\vec x} (\forall \, \vec y) \, \varphi \, .
$$
Autrement dit, $M$ v\'erifie le s\'equent
$$
\psi \vdash_{\vec x , \vec y} \varphi
$$
si et seulement si il v\'erifie le s\'equent
$$
\psi \vdash_{\vec x} (\forall \, \vec y) \, \varphi \, .
$$

\item[(11)] Pour les r\`egles d'implication, consid\'erons trois formules $\varphi , \psi$ et $\chi$ dans un m\^eme contexte $\vec x = (x_1^{A_1} \cdots x_n^{A_n})$.

\smallskip

Les s\'equents
$$
\varphi \wedge \psi \vdash_{\vec x} \chi
$$
et
$$
\varphi \vdash_{\vec x} (\psi \Rightarrow \chi)
$$
s'interpr\`etent dans le mod\`ele $M$ comme les relations d'inclusion entre sous-objets $M\!A_1 \times \cdots \times M\!A_n$
$$
M\varphi \, (\vec x) \wedge M\psi \, (\vec x) \leq M\chi \, (\vec x)
$$
et
$$
M\varphi \, (\vec x) \leq (M\psi \, (\vec x) \Rightarrow M\chi \, (\vec x)) \, .
$$
Or, ces deux relations sont \'equivalentes puisque, par d\'efinition, le foncteur
$$
M\psi (\vec x) \Rightarrow \bullet
$$
est adjoint \`a droite du foncteur
$$
\bullet \wedge M\psi \, (\vec x)
$$
sur les sous-objets de $M\!A_1 \times \cdots \times M\!A_n$.

\smallskip

Ainsi, le mod\`ele $M$ v\'erifie le s\'equent
$$
\varphi \wedge \psi \vdash_{\vec x} \chi
$$
si et seulement si il v\'erifie le s\'equent
$$
\varphi \vdash_{\vec x} (\psi \Rightarrow \chi) \, .
$$

\item[(12)] Pour les r\`egles de n\'egation, consid\'erons deux formules $\varphi$ et $\psi$ dans un m\^eme contexte $\vec x = (x_1^{A_1} \cdots x_n^{A_n})$.

\smallskip

Les s\'equents
$$
\varphi \wedge \psi \vdash_{\vec x} \, \perp \qquad \mbox{et} \qquad \varphi \vdash_{\vec x} \neg \, \psi
$$
s'interpr\`etent dans le mod\`ele $M$ comme les relations d'inclusion entre sous-objets de $M\!A_1 \times \cdots \times M\!A_n$
$$
M\varphi \, (\vec x) \wedge M\psi \, (\vec x) = \emptyset_{M\!A_1 \times \cdots \times M\!A_n}
$$
et
$$
M\varphi \, (\vec x) \leq \neg \, M\psi \, (\vec x) \, .
$$
Or, ces deux relations sont \'equivalentes puisque, par d\'efinition, le sous-objet de $M\!A_1 \times \cdots \times M\!A_n$
$$
\neg \, M\psi \, (\vec x)
$$
repr\'esente le foncteur contravariant
$$
{\rm Hom} \, (\bullet \wedge M\psi \, (\vec x) , \emptyset_{M\!A_1 \times \cdots \times M\!A_n}) \, .
$$
Ainsi, le mod\`ele $M$ v\'erifie le s\'equent
$$
\varphi \wedge \psi \vdash_{\vec x} \, \perp
$$
si et seulement si il v\'erifie le s\'equent
$$
\varphi \vdash_{\vec x} \neg \, \psi \, .
$$
Cela ach\`eve la preuve de la proposition. 

\end{listeisansmarge}
\end{demo}

\bigskip

Enfin, cette proposition est compl\'et\'ee par la suivante:

\begin{prop}\label{propV55}

Soit ${\mathbb T}$ une th\'eorie cart\'esienne de signature $\Sigma$.

\smallskip

Soit $M$ un mod\`ele de ${\mathbb T}$ dans une cat\'egorie cart\'esienne ${\mathcal C}$.

\smallskip

Consid\'erons un s\'equent de Horn de $\Sigma$ dans un contexte $\vec x = (x_1^{A_1} \cdots x_n^{A_n})$
$$
\varphi \vdash_{\vec x} \psi
$$
ou bien un syst\`eme de s\'equents
$$
(S) \left\{\begin{matrix}
f_1 (\vec x) = g_1 (\vec x) \wedge \cdots \wedge f_k(\vec x) = g_k (\vec x) \, \dashv \, \vdash_{\vec x} (\exists \, y)(u_1 (y) = x_1^{A_1} \wedge \cdots \wedge u_n(y) = x_n^{A_n}) \, , \\
{ \ } \\
u_1 (y) = u_1 (y') \wedge \cdots \wedge u_n (y) = u_n (y') \vdash_{y,y'} y=y' \hfill \\
{ \ } \\
\mbox{(o\`u $y$ et $y'$ sont affect\'ees \`a une sorte $B$),} \hfill
\end{matrix}\right.
$$
comme dans la d\'efinition \ref{defV46}.

\smallskip

Supposons que ce s\'equent ou ce syst\`eme de s\'equents est d\'emontrable dans la th\'eorie ${\mathbb T}$ par les r\`egles d'inf\'erence de la logique r\'eguli\`ere (ou m\^eme par celles de la logique du premier ordre infinitaire).

\smallskip

Alors le mod\`ele $M$ v\'erifie ce s\'equent de Horn [resp. ce syst\`eme de s\'equents] au sens que
$$
M\varphi \, (\vec x) \leq M\psi \, (\vec x) \quad \mbox{comme sous-objets de} \quad M\!A_1 \times \cdots \times M\!A_n
$$
[resp. au sens que le morphisme
$$
Mu_1 \times \cdots \times Mu_n : M\!B \longrightarrow M\!A_1 \times \cdots \times M\!A_n
$$
identifie $M\!B$ au sous-objet de $M\!A_1 \times \cdots \times M\!A_n$ d\'efini par les \'equations
$$
M\!f_1 = Mg_1 , \cdots , M\!f_k = Mg_k \ \mbox{]}.
$$
\end{prop}

\begin{demo}

Il existe une petite sous-cat\'egorie pleine ${\mathcal C}'$ de ${\mathcal C}$ qui est stable par toutes les limites finies et qui contient les objets
$$
M\!A
$$
associ\'es aux sortes $A$ de $\Sigma$ ainsi que les sous-objets
$$
M\!R \xhookrightarrow{ \ { \ } \ } M\!A_1 \times \cdots \times M\!A_n
$$
associ\'es aux symboles de relations $R \rightarrowtail A_1 \cdots A_n$ de $\Sigma$.

\smallskip

Le mod\`ele $M$ de la th\'eorie cart\'esienne ${\mathbb T}$ dans ${\mathcal C}$ peut \^etre vu comme un mod\`ele dans ${\mathcal C}'$. Il satisfait un s\'equent de Horn ou un syst\`eme de s\'equents $(S)$ dans ${\mathcal C}$ si et seulement si il les satisfait dans ${\mathcal C}'$.

\smallskip

On peut donc supposer que ${\mathcal C}$ est une petite cat\'egorie.

\smallskip

Consid\'erons alors le plongement de Yoneda
$$
{\mathcal C} \xhookrightarrow{ \ { \ } \ } \widehat{\mathcal C}
$$
dans le topos ${\mathcal E} = \widehat{\mathcal C}$ des pr\'efaisceaux sur ${\mathcal C}$.

\smallskip

Comme ce foncteur est pleinement fid\`ele et respecte les limites, le mod\`ele $M$ de ${\mathbb T}$ v\'erifie un s\'equent de Horn ou un syst\`eme de s\'equents $(S)$ si et seulement si il en est ainsi de son image dans ${\mathcal E} = \widehat{\mathcal C}$.

\smallskip

Nous sommes donc ramen\'es au cas o\`u $M$ est un mod\`ele de ${\mathbb T}$ dans un topos ${\mathcal E}$.

\smallskip

Pour que le s\'equent ou le syst\`eme de s\'equents consid\'er\'e soit v\'erifi\'e par $M$, il suffit donc qu'il soit d\'emontrable \`a partir des axiomes de ${\mathbb T}$ par les r\`egles d'inf\'erence de la logique r\'eguli\`ere ou m\^eme par celles de la logique du premier ordre infinitaire.

\smallskip

Cela conclut la d\'emonstration. 

\end{demo}

\subsection{Th\'eories \'equivalentes et th\'eories quotients}\label{subsec553}

\medskip

Il est naturel d'ordonner les th\'eories suivant leurs ensembles ou collections de propri\'et\'es d\'emontrables:

\begin{defn}\label{defV56}

Soit $\Sigma$ une signature.

\smallskip

Soient ${\mathbb T}$ et ${\mathbb T}'$ deux th\'eories de signature $\Sigma$ suppos\'ees alg\'ebriques [resp. de Horn, resp. cart\'esiennes ou plus g\'en\'eralement r\'eguli\`eres, resp. coh\'erentes, resp. g\'eom\'etriques, resp. du premier ordre finitaires, resp. du premier ordre infinitaires].

\smallskip

Alors:

\begin{listeimarge}

\item On dit que la th\'eorie ${\mathbb T}'$ est un quotient de ${\mathbb T}$ si tout axiome de ${\mathbb T}$ est d\'emontrable \`a partir des axiomes de ${\mathbb T}'$ par les r\`egles d'\'egalit\'e et de substitution [resp. par les r\`egles de la logique de Horn, resp. de la logique r\'eguli\`ere, resp. de la logique coh\'erente, resp. de la logique g\'eom\'etrique, resp. de la logique du premier ordre finitaire, resp. de la logique du premier ordre infinitaire].

\medskip

\item On dit que les th\'eories ${\mathbb T}$ et ${\mathbb T}'$ sont \'equivalentes si ${\mathbb T}'$ est un quotient de ${\mathbb T}$ et ${\mathbb T}$ est un quotient de ${\mathbb T}'$.
\end{listeimarge}
\end{defn}

\begin{remarksqed}
\begin{listeisansmarge}
\item Ainsi, ${\mathbb T}'$ est un quotient de ${\mathbb T}$ si tout s\'equent (ou syst\`eme de s\'equents) qui est d\'emontrable dans ${\mathbb T}$ par les r\`egles de la logique consid\'er\'ee est a fortiori d\'emontrable dans ${\mathbb T}'$ par ces r\`egles.

\medskip

\item En particulier, deux th\'eories ${\mathbb T}$ et ${\mathbb T}'$ sont \'equivalentes si tout s\'equent (ou syst\`eme de s\'equents) qui est d\'emontrable dans l'une est d\'emontrable dans l'autre.

\end{listeisansmarge}
\end{remarksqed}

\smallskip

On d\'eduit des propositions \ref{propV53}, \ref{propV54} et \ref{propV55}:

\begin{cor}\label{corV57}

Soit $\Sigma$ une signature.

\smallskip

Soient ${\mathbb T}$ et ${\mathbb T}'$ deux th\'eories de signature $\Sigma$ suppos\'ees alg\'ebriques [resp. de Horn ou plus g\'en\'eralement cart\'esiennes, resp. r\'eguli\`eres, resp. coh\'erentes, resp. g\'eom\'etriques, resp. du premier ordre finitaires, resp. du premier ordre infinitaires].

\smallskip

Soit ${\mathcal C}$ une cat\'egorie localement petite suppos\'ee alg\'ebrique [resp. cart\'esienne, resp. r\'eguli\`ere, resp. coh\'erente, resp. g\'eom\'etrique, resp. de Heyting, resp. un topos].

\smallskip

Alors:

\begin{listeimarge}

\item Si ${\mathbb T}'$ est un quotient de ${\mathbb T}$, la cat\'egorie
$$
{\mathbb T}'\mbox{\rm -mod} ({\mathcal C})
$$
est une sous-cat\'egorie pleine de la cat\'egorie
$$
{\mathbb T}\mbox{\rm -mod} ({\mathcal C}) \, .
$$

\item En particulier, si ${\mathbb T}$ et ${\mathbb T}'$ sont \'equivalentes, les deux cat\'egories de mod\`eles
$$
{\mathbb T}\mbox{\rm -mod} ({\mathcal C}) \qquad \mbox{et} \qquad {\mathbb T}'\mbox{\rm -mod} ({\mathcal C})
$$
se confondent.
\end{listeimarge}
\end{cor}

\begin{demo}
\begin{listeisansmarge}
\item Les cat\'egories ${\mathbb T}\mbox{-mod} ({\mathcal C})$ et ${\mathbb T}'\mbox{-mod} ({\mathcal C})$ sont par d\'efinition les sous-cat\'egories pleines de la cat\'egorie
$$
\Sigma\mbox{-str} ({\mathcal C})
$$
constitu\'ees respectivement des mod\`eles dans ${\mathcal C}$ de la th\'eorie ${\mathbb T}$ et de la th\'eorie ${\mathbb T}'$.

\smallskip

Si les axiomes de ${\mathbb T}$ sont d\'emontrables \`a partir des axiomes de ${\mathbb T}'$, il r\'esulte des propositions \ref{propV53}, \ref{propV54} et \ref{propV55} que tout mod\`ele $M$ de ${\mathbb T}'$ est a fortiori un mod\`ele de ${\mathbb T}$.

\smallskip

D'o\`u la conclusion cherch\'ee.

\medskip

\item est cons\'equence imm\'ediate de (i). 
\end{listeisansmarge}
\end{demo}

\section{Cat\'egories syntactiques}\label{sec56}

\subsection{Repr\'esentation des foncteurs des mod\`eles par les cat\'egories syntactiques}\label{subsec561}

\medskip

Consid\'erons une th\'eorie du premier ordre ${\mathbb T}$ dans une signature $\Sigma$ qui est suppos\'ee alg\'ebrique [resp. de Horn ou plus g\'en\'eralement cart\'esienne, resp. r\'eguli\`ere, resp. coh\'erente, resp. g\'eom\'etrique, resp. du premier ordre finitaire].

\smallskip

Cette th\'eorie ${\mathbb T}$ permet d'associer \`a toute cat\'egorie localement petite ${\mathcal C}$ qui est alg\'ebrique [resp. cart\'esienne, resp. r\'eguli\`ere, resp. coh\'erente, resp. g\'eom\'etrique, resp. de Heyting] la cat\'egorie
$$
{\mathbb T}\mbox{-mod} ({\mathcal C})
$$
des mod\`eles de ${\mathbb T}$ dans ${\mathcal C}$ et \`a tout foncteur
$$
F : {\mathcal C} \longrightarrow {\mathcal D}
$$
qui est alg\'ebrique [resp. cart\'esien, resp. r\'egulier, resp. coh\'erent, resp. g\'eom\'etrique, resp. de Heyting] le foncteur de transformation par $F$ des mod\`eles
$$
F : {\mathbb T}\mbox{-mod} ({\mathcal C}) \longrightarrow {\mathbb T}\mbox{-mod} ({\mathcal D})
$$
entre les cat\'egories de mod\`eles correspondantes.

\smallskip

Nous allons montrer que ce foncteur, appel\'e le foncteur des mod\`eles de ${\mathbb T}$, est repr\'esentable au sens suivant:

\begin{thm}\label{thmV61}

Consid\'erons comme ci-dessus une th\'eorie ${\mathbb T}$ suppos\'ee alg\'ebrique [resp. de Horn ou plus g\'en\'eralement cart\'esienne, resp. r\'eguli\`ere, resp. coh\'erente, resp. g\'eom\'etrique, resp. du premier ordre finitaire].

\smallskip

Alors il existe une cat\'egorie alg\'ebrique ${\mathcal C}_{\mathbb T} = {\mathcal C}_{\mathbb T}^{\rm alg}$ [resp. cart\'esienne ${\mathcal C}_{\mathbb T} = {\mathcal C}_{\mathbb T}^{\rm cart}$, resp. r\'eguli\`ere ${\mathcal C}_{\mathbb T} = {\mathcal C}_{\mathbb T}^{\rm reg}$, resp. coh\'erente ${\mathcal C}_{\mathbb T} = {\mathcal C}_{\mathbb T}^{\rm coh}$, resp. g\'eom\'etrique ${\mathcal C}_{\mathbb T} = {\mathcal C}_{\mathbb T}^{\rm geo}$, resp. de Heyting ${\mathcal C}_{\mathbb T} = {\mathcal C}_{\mathbb T}^{\rm He}$], munie d'un mod\`ele $M_{\mathbb T}$ de ${\mathbb T}$ dans ${\mathcal C}_{\mathbb T}$, telle que, pour toute cat\'egorie ${\mathcal C}$ alg\'ebrique [resp. cart\'esienne, resp. r\'eguli\`ere, resp. coh\'erente, resp. g\'eom\'etrique, resp. de Heyting], le foncteur
$$
(F : {\mathcal C}_{\mathbb T} \to {\mathcal C}) \longmapsto F(M_{\mathbb T})
$$
est une \'equivalence de la sous-cat\'egorie pleine de
$$
[{\mathcal C}_{\mathbb T} , {\mathcal C}]
$$
constitu\'ee des foncteurs
$$
F : {\mathcal C}_{\mathbb T} \longrightarrow {\mathcal C}
$$
qui sont alg\'ebriques [resp. cart\'esiens, resp. r\'eguliers, resp. coh\'erents, resp. g\'eom\'etriques, resp. de Heyting] sur la cat\'egorie
$$
{\mathbb T}\mbox{\rm -mod} ({\mathcal C})
$$
des mod\`eles de ${\mathbb T}$ dans ${\mathcal C}$.

\smallskip

De plus, la cat\'egorie ${\mathcal C}_{\mathbb T}$, qui est appel\'ee la cat\'egorie syntactique alg\'ebrique [resp. cart\'esienne, resp. r\'eguli\`ere, resp. coh\'erente, resp. g\'eom\'etrique, resp. de Heyting] de ${\mathbb T}$, est caract\'eris\'ee \`a \'equivalence pr\`es par cette propri\'et\'e.
\end{thm}

\begin{commdemo}

Prouvons d'abord l'unicit\'e \`a \'equivalence pr\`es des cat\'egories syntactiques ${\mathcal C}_{\mathbb T}$ v\'erifiant la propri\'et\'e de l'\'enonc\'e.

\smallskip

Consid\'erons donc deux cat\'egories alg\'ebriques [resp. cart\'esiennes, resp. r\'eguli\`eres, resp. coh\'erentes, resp. g\'eom\'etriques, resp. de Heyting] ${\mathcal C}_{\mathbb T}$ et ${\mathcal C}'_{\mathbb T}$, munies de deux mod\`eles $M_{\mathbb T}$ et $M'_{\mathbb T}$ de ${\mathbb T}$, qui repr\'esentent le foncteur des mod\`eles de ${\mathbb T}$
$$
{\mathcal C} \longmapsto {\mathbb T}\mbox{-mod} ({\mathcal C}) 
$$
au sens de l'\'enonc\'e.

\smallskip

Par hypoth\`ese, il existe deux foncteurs alg\'ebriques [resp. cart\'esiens, resp. r\'eguliers, resp. coh\'erents, resp. g\'eom\'etriques, resp. de Heyting] 
$$
F : {\mathcal C}_{\mathbb T} \longrightarrow {\mathcal C}'_{\mathbb T} \qquad \mbox{et} \qquad G : {\mathcal C}'_{\mathbb T} \longrightarrow {\mathcal C}_{\mathbb T}
$$
tels que les ${\mathbb T}$-mod\`eles $F(M_{\mathbb T})$ et $G(M'_{\mathbb T})$ soient respectivement isomorphes \`a $M'_{\mathbb T}$ et $M_{\mathbb T}$.

\smallskip

Il en r\'esulte que les ${\mathbb T}$-mod\`eles dans ${\mathcal C}_{\mathbb T}$ et dans ${\mathcal C}'_{\mathbb T}$
$$
G \circ F (M_{\mathbb T}) \qquad \mbox{et} \qquad F \circ G(M'_{\mathbb T})
$$
sont isomorphes respectivement \`a $M_{\mathbb T}$ et $M'_{\mathbb T}$.

\smallskip

Par hypoth\`ese, les isomorphismes de ${\mathbb T}$-mod\`eles
$$
G \circ F (M_{\mathbb T}) \xrightarrow{ \ \sim \ } M_{\mathbb T} \qquad \mbox{et} \qquad F \circ G (M'_{\mathbb T}) \xrightarrow{ \ \sim \ } M'_{\mathbb T} 
$$
proviennent d'isomorphismes de foncteurs
$$
G \circ F \xrightarrow{ \ \sim \ } {\rm id}_{{\mathcal C}_{\mathbb T}} \qquad \mbox{et} \qquad F \circ G \xrightarrow{ \ \sim \ }  {\rm id}_{{\mathcal C}'_{\mathbb T}} \, .
$$

Ainsi, $F$ et $G$ sont deux \'equivalences r\'eciproques l'une de l'autre entre les cat\'egories ${\mathcal C}_{\mathbb T}$ et ${\mathcal C}'_{\mathbb T}$.

\smallskip

Cela prouve l'unicit\'e \`a \'equivalence pr\`es de la cat\'egorie syntactique ${\mathcal C}_{\mathbb T}$ munie de son mod\`ele universel $M_{\mathbb T}$.

\smallskip

Pour l'existence, nous allons construire explicitement dans les trois prochains paragraphes une cat\'egorie ${\mathcal C}_{\mathbb T}$ susceptible de r\'epondre \`a la question pos\'ee.

\smallskip

Puis nous montrerons au paragraphe suivant e) que la cat\'egorie ${\mathcal C}_{\mathbb T}$ ainsi construite v\'erifie bien la propri\'et\'e de repr\'esentation du foncteur des mod\`eles de ${\mathbb T}$ \'enonc\'ee dans le th\'eor\`eme. 

\end{commdemo}

\subsection{Construction des cat\'egories syntactiques alg\'ebriques}\label{subsec562}

\medskip

Commen\c cons par construire les cat\'egories syntactiques alg\'ebriques des th\'eories alg\'ebriques et leurs mod\`eles universels:

\begin{defn}\label{defV62}

Soit $\Sigma$ une signature sans symbole de relation.

\smallskip

Soit ${\mathbb T}$ une th\'eorie alg\'ebrique de signature $\Sigma$.

\smallskip

On appelle ``cat\'egorie syntactique alg\'ebrique'' de ${\mathbb T}$ et on note
$$
{\mathcal C}_{\mathbb T} = {\mathcal C}_{\mathbb T}^{\rm alg}
$$
la petite cat\'egorie alg\'ebrique munie d'un mod\`ele canonique $M_{\mathbb T}$ ainsi d\'efinie:

\begin{listeimarge}

\item Les objets de ${\mathcal C}_{\mathbb T}^{\rm alg}$ sont les familles finies de sortes de $\Sigma$
$$
A_1 \cdots A_n \, .
$$

\item Les morphismes de ${\mathcal C}_{\mathbb T}^{\rm alg}$
$$
A_1 \cdots A_n \longrightarrow B_1 \cdots B_m
$$
sont les classes d'\'equivalence de familles de termes
$$
(f_1 (\vec x) , \cdots , f_m (\vec x))
$$
o\`u

\medskip

$
\left\{\begin{matrix}
\bullet &\mbox{les termes $f_1 (\vec x) , \cdots , f_m (\vec x)$ ont le m\^eme contexte} \hfill \\
{ \ } \\
&\vec x = (x_1^{A_1} \cdots x_n^{A_n}) \\
{ \ } \\
&\mbox{de sortes associ\'ees $A_1 \cdots A_n$ et sont \`a valeurs dans les sortes $B_1 , \cdots , B_m$,} \hfill \\
{ \ } \\
\bullet &\mbox{deux telles familles de termes} \hfill \\
{ \ } \\
&(f_1 (\vec x) , \cdots , f_m (\vec x)) \qquad \mbox{et} \qquad (g_1 (\vec y) , \cdots , g_m (\vec y)) \\
{ \ } \\
&\mbox{sont dites \'equivalentes si, apr\`es substitutions des variables $\vec x = (x_1^{A_1} \cdots x_n^{A_n})$ \`a $\vec y = (y_1^{A_1} \cdots y_n^{A_n})$,} \hfill \\
&\mbox{les s\'equents} \hfill \\
&\top \vdash_{\vec x} f_1 = g_1 , \cdots , \top \vdash_{\vec x} f_m = g_m \\
{ \ } \\
&\mbox{sont d\'emontrables dans ${\mathbb T}$.} \hfill 
\end{matrix} \right.
$

\medskip

\item La composition de deux morphismes
$$
A_1 \cdots A_n \xrightarrow{ \ f \ } B_1 \cdots B_m \xrightarrow{ \ g \ } C_1 \cdots C_{\ell}
$$
repr\'esent\'es par des familles de termes
$$
(f_1 (\vec x) , \cdots , f_m (\vec x)) \qquad \mbox{et} \qquad (g_1 (\vec y) , \cdots , g_{\ell} (\vec y))
$$
de contextes
$$
\vec x = (x_1^{A_1} \cdots x_n^{A_n}) \qquad \mbox{et} \qquad \vec y = (y_1^{B_1} \cdots y_m^{B_m})
$$
est d\'efinie par la substitution des termes $f_1 (\vec x) , \cdots , f_m (\vec x)$ aux variables $y_1^{B_1} , \cdots , y_m^{B_m}$.

\medskip

\item Le mod\`ele $M_{\mathbb T}$ de ${\mathbb T}$ dans la cat\'egorie alg\'ebrique ${\mathcal C}_{\mathbb T}^{\rm alg}$ consiste \`a associer

\medskip

$\left\{\begin{matrix}
\bullet &\mbox{\`a toute sorte $A$ de $\Sigma$ l'objet $A$ de ${\mathcal C}_{\mathbb T}^{\rm alg}$,} \hfill \\
{ \ } \\
\bullet &\mbox{\`a tout symbole de fonction de $\Sigma$} \hfill \\
&f : A_1 \cdots A_n \longrightarrow B \\
&\mbox{le morphisme} \hfill \\
&A_1 \cdots A_n \longrightarrow B \\
&\mbox{de ${\mathcal C}_{\mathbb T}^{\rm alg}$ repr\'esent\'e par le terme} \hfill \\
&f(\vec x) \\
{ \ } \\
&\mbox{pour n'importe quel choix de contexte $\vec x = (x_1^{A_1} \cdots x_n^{A_n})$ de sortes associ\'ees $A_1 \cdots A_n$.} \hfill
\end{matrix} \right.
$
\end{listeimarge}
\end{defn}

\begin{remarksqed}
\begin{listeisansmarge}

\item Il r\'esulte de la ``r\`egle d'identit\'e'' et des ``r\`egles d'\'egalit\'e'' de la logique du premier ordre que la relation entre familles de termes introduite dans (ii) est bien une relation d'\'equivalence.

\medskip

\item De m\^eme, il r\'esulte des ``r\`egles de substitution'' de la logique du premier ordre que cette relation d'\'equivalence entre familles de termes est respect\'ee par la loi de composition de (iii) d\'efinie par substitution de termes \`a des variables.

\smallskip

Ainsi, cette loi de composition est bien d\'efinie.

\medskip

\item La loi de composition est associative car il en est ainsi des op\'erations de substitution de termes \`a des variables.

\smallskip

De plus, tout objet $A_1 \cdots A_n$ admet pour morphisme identique la classe d'\'equivalence de la famille de termes
$$
(x_1^{A_1} , \cdots , x_n^{A_n})
$$
pour n'importe quel choix de variables $x_1^{A_1} , \cdots , x_n^{A_n}$ affect\'ees aux sortes $A_1 , \cdots , A_n$.

\smallskip

Ainsi, ${\mathcal C}_{\mathbb T}^{\rm alg}$ est bien une cat\'egorie.

\smallskip

Elle est petite par construction.

\medskip

\item La cat\'egorie ${\mathcal C}_{\mathbb T}^{\rm alg}$ est alg\'ebrique, c'est-\`a-dire admet des produits finis.

\smallskip

En effet, ${\mathcal C}_{\mathbb T}^{\rm alg}$ admet d'une part pour objet terminal la famille vide de sortes, puisque tout objet $A_1 \cdots A_n$ de ${\mathcal C}_{\mathbb T}^{\rm alg}$ admet pour unique morphisme vers cet objet la famille vide de termes.

\smallskip

D'autre part, deux objets $A_1 \cdots A_n$ et $B_1 \cdots B_m$ de ${\mathcal C}_{\mathbb T}^{\rm alg}$ admettent pour produit l'objet $A_1 \cdots A_n \, B_1 \cdots B_m$. 

\end{listeisansmarge}
\end{remarksqed}

\subsection{Construction des cat\'egories syntactiques r\'eguli\`eres, coh\'erentes, g\'eom\'etriques et de Heyting}\label{subsec563}

\medskip

Pour construire la cat\'egorie syntactique r\'eguli\`ere [resp. coh\'erente, resp. g\'eom\'etrique, resp. de Heyting] d'une th\'eorie r\'eguli\`ere [resp. coh\'erente, resp. g\'eom\'etrique, resp. du premier ordre finitaire], on a d'abord besoin de poser les deux d\'efinitions suivantes:

\begin{defn}\label{defV63}

Deux formules d'une signature $\Sigma$
$$
\varphi (\vec x) \qquad \mbox{et} \qquad \psi (\vec y)
$$
en des variables $\vec x = (x_1^{A_1} \cdots x_n^{A_n})$ et $\vec y = (y_1^{A_1} \cdots y_n^{A_n})$ affect\'ees aux m\^emes sortes $A_1 \cdots A_n$, sont dites ``\'equivalentes par substitution'' si $\psi (\vec y)$ se d\'eduit de $\psi (\vec x)$ par substitution des variables $\vec y$ aux variables $\vec x$.  \hfill $\Box$
\end{defn}

Dans cette d\'efinition n'intervient aucune notion de d\'emontrabilit\'e relative \`a aucune th\'eorie.

\smallskip

Il en va autrement de la d\'efinition suivante:

\begin{defn}\label{defV64}

Soit ${\mathbb T}$ une th\'eorie r\'eguli\`ere [resp. coh\'erente, resp. g\'eom\'etrique, resp. du premier ordre finitaire].

\begin{listeimarge}

\item Si $\vec x = (x_1^{A_1} \cdots x_n^{A_n})$ et $\vec y = (y_1^{B_1} \cdots y_m^{B_m})$ sont deux contextes disjoints, et $\varphi (\vec x)$, $\psi (\vec y)$ sont deux formules r\'eguli\`eres [resp. coh\'erentes, resp. g\'eom\'etriques, resp. du premier ordre finitaires] dans ces contextes, une formule r\'eguli\`ere  [resp. coh\'erente, resp. g\'eom\'etrique, resp. du premier ordre finitaire] 
$$
\theta (\vec x , \vec y) \quad \mbox{de contexte} \quad (\vec x , \vec y) = (x_1^{A_1} \cdots x_n^{A_n} \, y_1^{B_1} \cdots y_m^{B_m})
$$ 
est dite ``d\'emontrablement fonctionnelle'' de $\varphi (\vec x)$ dans $\psi (\vec y)$ si les s\'equents

\medskip

$
\left\{\begin{matrix}
&\theta \vdash_{\vec x , \vec y} \varphi \wedge \psi \, , \\
{ \ } \\
&\varphi \vdash_{\vec x} (\exists \, \vec y) \, \theta \, , \\
{ \ } \\
&\theta (\vec x , \vec y) \wedge \theta (\vec x , \vec y') \vdash_{\vec x , \vec y , \vec y'} \vec y = \vec y'
\end{matrix} \right.
$

\bigskip

\noindent sont d\'emontrables dans la th\'eorie ${\mathbb T}$ selon les r\`egles d'inf\'erence de la logique r\'eguli\`ere [resp. coh\'erente, resp. g\'eom\'etrique, resp. du premier ordre finitaire].

\medskip

\item Etant donn\'ees deux telles formules r\'eguli\`eres [resp. coh\'erentes, resp. g\'eom\'etriques, resp. du premier ordre finitaires]
$$
\varphi (\vec x) \qquad \mbox{et} \qquad \psi (\vec y)
$$
dans deux contextes disjoints $\vec x$ et $\vec y$, deux formules d\'emontrablement fonctionnelles de $\varphi (\vec x)$ dans $\psi (\vec y)$
$$
\theta (\vec x , \vec y) \qquad \mbox{et} \qquad \theta' (\vec x , \vec y)
$$
sont dites \'equivalentes si le double s\'equent
$$
\theta \dashv \, \vdash_{\vec x , \vec y} \theta'
$$
est d\'emontrable dans la th\'eorie ${\mathbb T}$ suivant les r\`egles d'inf\'erence de la logique r\'eguli\`ere [resp. coh\'erente, resp. g\'eom\'etrique, resp. du premier ordre finitaire].
\end{listeimarge}
\end{defn}

Ces deux d\'efinitions \'etant pos\'ees, on peut construire les cat\'egories syntactiques et les mod\`eles universels des th\'eories r\'eguli\`eres, coh\'erentes, g\'eom\'etriques ou du premier ordre finitaires:

\begin{defn}\label{defV65}

Soit $\Sigma$ une signature.

\smallskip

Soit ${\mathbb T}$ une th\'eorie r\'eguli\`ere [resp. coh\'erente, resp. g\'eom\'etrique, resp. du premier ordre finitaire] de signature $\Sigma$.

\smallskip

On appelle ``cat\'egorie syntactique r\'eguli\`ere'' [resp. coh\'erente, resp. g\'eom\'etrique, resp. de Heyting] de ${\mathbb T}$ et on note
$$
{\mathcal C}_{\mathbb T} = {\mathcal C}_{\mathbb T}^{\rm reg} \qquad \mbox{[resp.} \ {\mathcal C}_{\mathbb T}^{\rm coh} , \ \mbox{resp.} \ {\mathcal C}_{\mathbb T}^{\rm geo} , \ \mbox{resp.} \ {\mathcal C}_{\mathbb T}^{\rm He} \ \mbox{]}
$$
la cat\'egorie r\'eguli\`ere [resp. coh\'erente, resp. g\'eom\'etrique, resp. de Heyting] munie d'un mod\`ele canonique $M_{\mathbb T}$ ainsi d\'efinie:

\begin{listeimarge}

\item Les objets de ${\mathcal C}_{\mathbb T}$ sont les formules r\'eguli\`eres [resp. coh\'erentes, resp. g\'eom\'etriques, resp. du premier ordre finitaires] de la signature $\Sigma$
$$
\varphi (\vec x) \quad \mbox{dans un contexte} \quad \vec x = (x_1^{A_1} \cdots x_n^{A_n}) \, ,
$$
consid\'er\'ees \`a \'equivalence par substitution pr\`es.

\item Etant donn\'es deux objets de ${\mathcal C}_{\mathbb T}$ repr\'esent\'es par deux formules
$$
\varphi (\vec x) \qquad \mbox{et} \qquad \psi (\vec y)
$$
dans des contextes disjoints $\vec x$ et $\vec y$, les morphismes
$$
\varphi (\vec x) \longrightarrow \psi (\vec y)
$$
de ${\mathcal C}_{\mathbb T}$ sont les classes d'\'equivalence de formules r\'eguli\`eres [resp. coh\'erentes, resp. g\'eom\'etriques, resp. du premier ordre finitaires]
$$
\theta (\vec x , \vec y)
$$
qui sont d\'emontrablement fonctionnelles de $\varphi (\vec x)$ dans $\psi (\vec y)$.

\medskip

\item Etant donn\'es trois objets de ${\mathcal C}_{\mathbb T}$ repr\'esent\'es par trois formules
$$
\varphi (\vec x) \, , \quad \psi (\vec y) \quad \mbox{et} \quad \chi (\vec z)
$$
dans des contextes disjoints $\vec x$, $\vec y$ et $\vec z$, le compos\'e de deux morphismes
$$
\varphi (\vec x) \longrightarrow \psi (\vec y)
$$
et
$$
\psi (\vec y) \longrightarrow \chi (\vec z)
$$
repr\'esent\'es par deux formules d\'emontrablement fonctionnelles
$$
\theta (\vec x , \vec y) \qquad \mbox{et} \qquad \theta' (\vec y , \vec z)
$$
est le morphisme
$$
\varphi (\vec x) \longrightarrow \chi (\vec z)
$$
repr\'esent\'e par la formule d\'emontrablement fonctionnelle
$$
(\exists \, \vec y) \, (\theta (\vec x , \vec y) \wedge \theta' (\vec y , \vec z)) \, .
$$

\medskip

\item Le mod\`ele $M_{\mathbb T}$ de ${\mathbb T}$ dans la cat\'egorie r\'eguli\`ere [resp. coh\'erente, resp. g\'eom\'etrique, resp. de Heyting] ${\mathcal C}_{\mathbb T}$ consiste \`a associer

\medskip

$
\left\{ \begin{matrix}
\bullet &\mbox{\`a toute sorte $A$ l'objet repr\'esent\'e par la formule} \hfill \\
{ \ } \\
&\top (x^A) \\
{ \ } \\
&\mbox{pour n'importe quel choix d'une variable $x^A$ de sorte $A$,} \hfill \\
{ \ } \\
\bullet &\mbox{\`a tout symbole de fonction} \hfill \\
{ \ } \\
&f : A_1 \cdots A_n \longrightarrow B \\
&\mbox{le morphisme} \hfill \\
{ \ } \\
&\top (x_1^{A_1} \cdots x_n^{A_n}) \longrightarrow \top (y^B) \\
{ \ } \\
&\mbox{repr\'esent\'e par la formule d\'emontrablement fonctionnelle} \hfill \\
{ \ } \\
&f(x_1^{A_1} \cdots x_n^{A_n}) = y^B \, , \\
{ \ } \\
\bullet &\mbox{\`a tout symbole de relation} \hfill \\
{ \ } \\
&R \rightarrowtail A_1 \cdots A_n \\
&\mbox{le monomorphisme} \hfill \\
{ \ } \\
&R(x_1^{A_1} \cdots x_n^{A_n}) \xhookrightarrow{ \ { \ } \ } \top (y_1^{A_1} \cdots y_n^{A_n}) \\
{ \ } \\
&\mbox{repr\'esent\'e par la formule d\'emontrablement fonctionnelle} \hfill \\
{ \ } \\
&R(x_1^{A_1} \cdots x_n^{A_n}) \wedge x_1^{A_1} = y_1^{A_1} \wedge \cdots \wedge x_n^{A_n} = y_n^{A_n} \, .
\end{matrix} \right.
$
\end{listeimarge}
\end{defn}

\bigskip

\begin{remarksqed}
\begin{listeisansmarge}
\item Ainsi, la liste des objets de ${\mathcal C}_{\mathbb T}$ ne d\'epend que de la signature $\Sigma$ de ${\mathbb T}$ et du fragment de la logique du premier ordre (r\'egulier, coh\'erent, g\'eom\'etrique ou du premier ordre finitaire) que l'on consid\`ere. Elle ne d\'epend pas des axiomes de la th\'eorie ${\mathbb T}$.

\smallskip

En revanche, les morphismes entre les objets de ${\mathcal C}_{\mathbb T}$ d\'ependent de ces axiomes.

\medskip

\item Une th\'eorie alg\'ebrique ${\mathbb T}$ poss\`ede non seulement une cat\'egorie syntactique alg\'ebrique ${\mathcal C}_{\mathbb T}^{\rm alg}$ mais aussi des cat\'egories syntactiques r\'eguli\`ere ${\mathcal C}_{\mathbb T}^{\rm reg}$, coh\'erente ${\mathcal C}_{\mathbb T}^{\rm coh}$, g\'eom\'etrique ${\mathcal C}_{\mathbb T}^{\rm geo}$ et de Heyting ${\mathcal C}_{\mathbb T}^{\rm He}$.

\smallskip

De m\^eme, une th\'eorie r\'eguli\`ere ${\mathbb T}$ poss\`ede non seulement une cat\'egorie syntactique r\'eguli\`ere ${\mathcal C}_{\mathbb T}^{\rm reg}$ mais aussi des cat\'egories syntactiques coh\'erente ${\mathcal C}_{\mathbb T}^{\rm coh}$, g\'eom\'etrique ${\mathcal C}_{\mathbb T}^{\rm geo}$ et de Heyting ${\mathcal C}_{\mathbb T}^{\rm He}$.

\smallskip

Enfin, une th\'eorie coh\'erente poss\`ede des cat\'egories syntactiques coh\'erente ${\mathcal C}_{\mathbb T}^{\rm coh}$, g\'eom\'etrique ${\mathcal C}_{\mathbb T}^{\rm geo}$ et de Heyting ${\mathcal C}_{\mathbb T}^{\rm He}$. 

\end{listeisansmarge}
\end{remarksqed}


V\'erifions que la d\'efinition \ref{defV65} ci-dessus est valide:

\begin{prop}\label{propV66}

Soit ${\mathbb T}$ une th\'eorie r\'eguli\`ere [resp. coh\'erente, resp. g\'eom\'etrique, resp. du premier ordre finitaire] de signature $\Sigma$.

\smallskip

Alors:

\begin{listeimarge}

\item Les objets (ou les classes d'\'equivalence d'objets dans le cas g\'eom\'etrique) et les morphismes de ${\mathcal C}_{\mathbb T}$ introduits dans la d\'efinition \ref{defV65} forment des ensembles.

\medskip

\item Ces ensembles d'objets (ou de classes d'objets dans le cas g\'eom\'etrique) et de morphismes munis de la loi de composition introduite dans la d\'efinition \ref{defV65} (iii) forment une cat\'egorie petite (ou essentiellement petite dans le cas g\'eom\'etrique)
$$
{\mathcal C}_{\mathbb T} = {\mathcal C}_{\mathbb T}^{\rm reg} \qquad \mbox{[resp. ${\mathcal C}_{\mathbb T}^{\rm coh}$, resp. ${\mathcal C}_{\mathbb T}^{\rm geo}$, resp. ${\mathcal C}_{\mathbb T}^{\rm He}$].}
$$

\item Cette cat\'egorie petite (ou essentiellement petite) est r\'eguli\`ere [resp. coh\'erente, resp. g\'eom\'etrique, resp. de Heyting].

\smallskip

De plus, les sous-objets d'un objet $\varphi (\vec x)$ de cette cat\'egorie correspondent aux formules $\varphi_1 (\vec x)$ de m\^eme contexte $\vec x$ telles que le s\'equent $\varphi_1 \vdash_{\vec x} \varphi$ soit d\'emontrable dans ${\mathbb T}$.

\medskip

\item La $\Sigma$-structure
$$
M_{\mathbb T} \qquad \mbox{de} \qquad {\mathcal C}_{\mathbb T}
$$
introduite dans la d\'efinition \ref{defV65} (iv) est un mod\`ele de la th\'eorie ${\mathbb T}$.
\end{listeimarge}
\end{prop}

\begin{demosansqed}
\begin{listeisansmarge}
\item Par d\'efinition d'une signature, les sortes de $\Sigma$, ses symboles de fonctions et ses symboles de relations forment des ensembles.

\smallskip

Il en r\'esulte que les formules atomiques de $\Sigma$ consid\'er\'ees \`a substitution pr\`es des variables forment \'egalement un ensemble.

\smallskip

Or les formules r\'eguli\`eres [resp. coh\'erentes, resp. du premier ordre finitaires] de $\Sigma$ se d\'eduisent chacune d'un nombre fini de formules atomiques en appliquant un nombre fini de fois des symboles $\wedge$ et $\exists$ [resp. $\wedge$, $\exists$ et $\vee$, resp. $\wedge$, $\exists$, $\vee$, $\forall$, $\Rightarrow$ et $\neg$]. Donc elles forment un ensemble.

\smallskip

A fortiori, si ${\mathbb T}$ est une th\'eorie r\'eguli\`ere [resp. coh\'erente, resp. du premier ordre finitaire], les objets et les morphismes de ${\mathcal C}_{\mathbb T} = {\mathcal C}_{\mathbb T}^{\rm reg}$ [resp. ${\mathcal C}_{\mathbb T}^{\rm coh}$, resp. ${\mathcal C}_{\mathbb T}^{\rm He}$] forment des ensembles.

\smallskip

Le cas des formules g\'eom\'etriques est un peu plus d\'elicat du fait que ces formules peuvent comprendre des disjonctions infinitaires $\underset{i \in I}{\bigvee}$ sur des familles de sous-formules index\'ees par des ensembles arbitraires $I$.

\smallskip

Cependant, ce probl\`eme va \^etre r\'esolu gr\^ace au lemme suivant:
\end{listeisansmarge}
\end{demosansqed}

\begin{lem}\label{lemV67}

Soit $\Sigma$ une signature.

\smallskip

Alors:

\begin{listeimarge}

\item Toute formule g\'eom\'etrique [resp. coh\'erente] de $\Sigma$ est d\'emontrablement \'equivalente \`a une formule de la forme
$$
\varphi = \bigvee_{i \in I} \varphi_i \qquad \mbox{[resp.} \quad \varphi = \varphi_1 \vee \cdots \vee \varphi_k \ \mbox{]}
$$
o\`u les $\varphi_i$ sont des formules r\'eguli\`eres dans un m\^eme contexte $\vec x = (x_1^{A_1} \cdots x_n^{A_n})$ et $I$ est une partie d'un ensemble fix\'e.

\medskip

\item Toute formule r\'eguli\`ere de $\Sigma$ est d\'emontrablement \'equivalente \`a une formule de la forme
$$
\psi = (\exists \, \vec y) \, \psi'
$$
o\`u $\psi'$ est une formule de Horn dans un contexte qui comprend $\vec y$.
\end{listeimarge}
\end{lem}

\begin{demosansqed}

Toute formule g\'eom\'etrique [resp. coh\'erente] de $\Sigma$ est construite \`a partir d'une famille [resp. d'une famille finie] de formules atomiques en appliquant un nombre fini de fois des symboles de conjonctions finitaires $\wedge$, de quantifications existentielles $\exists$ et de disjonctions arbitraires $\underset{i \in I}{\bigvee}$ [resp. de disjonctions finitaires $\vee$].

\smallskip

De m\^eme, toute formule r\'eguli\`ere est construite \`a partir d'une famille finie de formules atomiques en appliquant un nombre fini de fois des symboles $\wedge$ et $\exists$.

\smallskip

Pour conclure, il suffit de montrer qu'il est possible de permuter l'ordre dans lequel les symboles $\underset{i \in I}{\bigvee}$ [resp. $\vee$], $\exists$ et $\wedge$ sont appliqu\'es.

\smallskip

Les symboles $\underset{i \in I}{\bigvee}$ [resp. $\vee$] et $\wedge$ peuvent toujours \^etre permut\'es d'apr\`es la ``r\`egle de distributivit\'e''.

\smallskip

Les deux autres paires de symboles peuvent aussi \^etre permut\'ees. C'est le contenu du lemme suivant:
\end{demosansqed}

\begin{lem}\label{lemV68}

Soit $\Sigma$ une signature.

\begin{listeimarge}

\item Si les $\varphi_i$, $1 \leq i \leq k$ [resp. $i \in I$], sont une famille de formules coh\'erentes [resp. g\'eom\'etriques] dans un contexte $(\vec x , \vec y)$, le double s\'equent
$$
(\exists \, \vec y) \left( \bigvee_{i \in I} \varphi_i (\vec x , \vec y) \right) \dashv \, \vdash_{\vec x} \bigvee_{i \in I} (\exists \, \vec y) \, \varphi_i (\vec x , \vec y)
$$
[resp.
$$
(\exists \, \vec y) \, (\varphi_1 (\vec x , \vec y) \vee \cdots \vee \varphi_k (\vec x , \vec y)) \dashv \, \vdash_{\vec x} (\exists \, \vec y) \, \varphi_1 (\vec x , \vec y) \vee \cdots \vee (\exists \, \vec y) \, \varphi_k (\vec x , \vec y) \ \mbox{]}
$$
est d\'emontrable suivant les r\`egles d'inf\'erence de la logique g\'eom\'etrique [resp. coh\'erente].

\medskip

\item Si $\varphi$ est une formule r\'eguli\`ere de contexte $\vec x$ et $\psi$ est une formule r\'eguli\`ere de contexte $(\vec x , \vec y)$, le double s\'equent
$$
\varphi \wedge (\exists \, \vec y) \, \psi \dashv \, \vdash_{\vec x} (\exists \, \vec y) (\varphi \wedge \psi)
$$
est d\'emontrable suivant les r\`egles d'inf\'erence de la logique r\'eguli\`ere.
\end{listeimarge}
\end{lem}

\begin{demo}
\begin{listeisansmarge}
\item Pour toute famille $\psi$ de contexte $\vec x$, le s\'equent
$$
(\exists \, \vec y) \left( \bigvee_{i \in I} \varphi_i (\vec x , \vec y) \right) \, \vdash_{\vec x} \psi (\vec x)
$$
\'equivaut d'apr\`es les r\`egles de quantification existentielle au s\'equent
$$
\bigvee_{i \in I} \varphi_i (\vec x , \vec y) \vdash_{\vec x , \vec y} \psi (\vec x) \, .
$$
Celui-ci \'equivaut \`a la famille de s\'equents
$$
\varphi_i (\vec x , \vec y) \vdash_{\vec x , \vec y} \psi (\vec x) \, , \quad i \in I \, ,
$$
d'apr\`es les r\`egles de disjonctions infinitaires [resp. finitaires si $I$ est fini].

\smallskip

Chaque s\'equent
$$
\varphi_i (\vec x , \vec y) \vdash_{\vec x , \vec y} \psi (\vec x)
$$
\'equivaut \`a son tour au s\'equent
$$
(\exists \, \vec y) \, \varphi_i (\vec x , \vec y) \vdash_{\vec x} \psi (\vec x) \, ,
$$
et leur famille \'equivaut au s\'equent
$$
\bigvee_{i \in I} (\exists \, \vec y) \, \varphi_i (\vec x , \vec y) \vdash_{\vec x} \psi (\vec x) \, .
$$
Cela prouve l'\'equivalence
$$
(\exists \, \vec y) \left(\bigvee_{i \in I} \varphi_i (\vec x , \vec y) \right) \dashv \, \vdash_{\vec x} \bigvee_{i \in I} (\exists \, \vec y) \, \varphi_i (\vec x , \vec y) \, .
$$

\item La ``r\`egle de Frobenius'' est l'implication
$$
\varphi \wedge (\exists \, \vec y) \, \psi \vdash_{\vec x} (\exists \, \vec y) (\varphi \wedge \psi) \, .
$$

Dans l'autre sens, l'implication tautologique
$$
(\exists \, \vec y) \, \psi (\vec x , \vec y) \vdash_{\vec x} (\exists \, \vec y') \, \psi (\vec x , \vec y')
$$
\'equivaut \`a l'implication
$$
\psi \vdash_{\vec x , \vec y} (\exists \, \vec y') \, \psi (\vec x , \vec y') \, .
$$
Celle-ci entra{\^\i}ne l'implication
$$
\varphi \wedge \psi \vdash_{\vec x , \vec y} \varphi \wedge (\exists \, \vec y') \, \psi (\vec x , \vec y')
$$
qui \'equivaut \`a l'implication
$$
(\exists \, \vec y) \, (\varphi \wedge \psi) \vdash_{\vec x} \varphi \wedge (\exists \, \vec y') \, \psi (\vec x , \vec y') \, .
$$

Ainsi, l'\'equivalence
$$
\varphi \wedge (\exists \, \vec y) \, \psi \dashv \, \vdash_{\vec x} (\exists \, \vec y) \, (\varphi \wedge \psi)
$$
est d\'emontrable.

\smallskip

Cela ach\`eve de prouver le lemme \ref{lemV68} et donc aussi le lemme \ref{lemV67}. 

\end{listeisansmarge}
\end{demo}

\bigskip

\noindent {\bf Suite de la d\'emonstration de la proposition \ref{propV66}:}

\begin{listeisansmarge}

\item Terminons d'abord de prouver (i) dans le cas o\`u ${\mathbb T}$ est une th\'eorie g\'eom\'etrique de signature $\Sigma$.

\smallskip

D'apr\`es le lemme \ref{lemV67} (i), toute formule g\'eom\'etrique de $\Sigma$ est d\'emontrablement \'equivalente \`a une formule de la forme
$$
\varphi = \bigvee_{i \in I} \varphi_i
$$
o\`u les $\varphi_i$ sont des formules r\'eguli\`eres dans un m\^eme contexte $\vec x = (x_1^{A_1} \cdots x_n^{A_n})$.

\smallskip

Or, les formules r\'eguli\`eres, consid\'er\'ees \`a \'equivalence pr\`es par substitution des variables, forment un ensemble.

\smallskip

On peut donc supposer que $I$ est une partie de cet ensemble.

\smallskip

Ainsi, tout objet de ${\mathcal C}_{\mathbb T}$, c'est-\`a-dire toute formule g\'eom\'etrique de $\Sigma$, est d\'emontrablement \'equivalent \`a un objet \'el\'ement d'un ensemble fix\'e.

\smallskip

Et, pour tous objets $\varphi$ et $\psi$, les morphismes
$$
\varphi \longrightarrow \psi \, ,
$$
qui sont par d\'efinition les formules d\'emontrablement fonctionnelles consid\'er\'ees \`a \'equivalence pr\`es, forment un ensemble.

\medskip

\item Montrons d'abord que la loi de composition des morphismes de ${\mathcal C}_{\mathbb T}$ est bien d\'efinie.

\smallskip

Il s'agit de prouver que si
$$
\theta : \varphi (\vec x) \longrightarrow \psi (\vec y)
$$
et
$$
\theta' : \psi (\vec y) \longrightarrow \chi (\vec z)
$$
sont deux formules d\'emontrablement fonctionnelles, alors la formule
$$
(\exists \, \vec y) (\theta \wedge \theta') : \varphi (\vec x) \longrightarrow \chi (\vec z)
$$
et encore d\'emontrablement fonctionnelle.

\smallskip  

Les s\'equents
$$
\theta \vdash_{\vec x , \vec y} \varphi \wedge \psi
$$
et
$$
\theta' \vdash_{\vec y , \vec z} \psi \wedge \chi
$$
impliquent le s\'equent
$$
\theta \wedge \theta' \vdash_{\vec x , \vec y, \vec z} \varphi \wedge \chi
$$
qui \'equivaut d'apr\`es les ``r\`egles de quantification existentielle'' au s\'equent
$$
(\exists \, \vec y) \, (\theta \wedge \theta') \vdash_{\vec x , \vec z} \varphi \wedge \chi \, .
$$

D'autre part, la formule
$$
(\theta \wedge \theta') (\vec x , \vec y , \vec z) \wedge (\theta \wedge \theta') (\vec x , \vec y' , \vec z')
$$
se r\'e\'ecrit
$$
\theta (\vec x , \vec y) \wedge \theta (\vec x , \vec y') \wedge \theta' (\vec y , \vec z) \wedge \theta' (\vec y' , \vec z') \, .
$$

Donc le s\'equent
$$
(\theta \wedge \theta') (\vec x , \vec y , \vec z) \wedge (\theta \wedge \theta') (\vec x , \vec y' , \vec z') \vdash \vec y = \vec y' \wedge \theta' (\vec y , \vec z) \wedge \theta' (\vec y' , \vec z')
$$
est d\'emontrable.

\smallskip

Il en est de m\^eme du s\'equent
$$
\vec y = \vec y' \wedge \theta' (\vec y , \vec z) \wedge \theta' (\vec y' , \vec z') \vdash \vec z = \vec z'
$$
donc aussi du s\'equent
$$
(\theta \wedge \theta') (\vec x , \vec y , \vec z) \wedge (\theta \wedge \theta') (\vec x , \vec y' , \vec z') \vdash \vec z = \vec z'
$$
qui \'equivaut au s\'equent
$$
(\exists \, \vec y) (\theta \wedge \theta') (\vec x , \vec y , \vec z) \wedge (\exists \, \vec y') (\vec x , \vec y , \vec z') \vdash \vec z = \vec z' \, .
$$
Enfin, les s\'equents
$$
\varphi (\vec x) \vdash (\exists \, \vec y) \, \theta (\vec x , \vec y) \, ,
$$
$$
\theta (\vec x , \vec y) \vdash \psi (\vec y) \, ,
$$
$$
\psi (\vec y) \vdash (\exists \, \vec z) \, \theta' (\vec y , \vec z)
$$
impliquent le s\'equent
$$
\varphi (\vec x) \vdash (\exists \, \vec z) (\exists \, \vec y) \, (\theta (\vec x , \vec y) \wedge \theta' (\vec y , \vec z)) \, .
$$

Ainsi, la loi de composition des morphismes de ${\mathcal C}_{\mathbb T}$ est bien d\'efinie.

\smallskip

Elle est associative car le symbole $\wedge$ est associatif (d'apr\`es les r\`egles de conjonctions finitaires) et commute avec le symbole $\exists$ (d'apr\`es le lemme \ref{lemV68} (ii)).

\smallskip

Enfin, tout objet de ${\mathcal C}_{\mathbb T}$ repr\'esent\'e par une formule
$$
\varphi (\vec x)
$$
dans un contexte $\vec x = (x_1^{A_1} \cdots x_n^{A_n})$ admet pour morphisme identit\'e la classe de la formule d\'emontrablement fonctionnelle
$$
\varphi (\vec x) \wedge x_1^{A_1} = y_1^{A_1} \wedge \cdots \wedge x_n^{A_n} = y_n^{A_n}
$$
de $\varphi (\vec x)$ dans $\varphi (\vec y)$, pour n'importe quel choix de contexte $\vec y = (y_1^{A_1} \cdots y_n^{A_n})$ disjoint de $\vec x$ affect\'e aux m\^emes sortes $A_1 \cdots A_n$ que $\vec x$.

\smallskip

Ainsi, ${\mathcal C}_{\mathbb T}$ est bien une cat\'egorie.

\smallskip

D'apr\`es (i), elle est petite si ${\mathbb T}$ est une th\'eorie r\'eguli\`ere, coh\'erente ou du premier ordre finitaire, et elle est essentiellement petite si ${\mathbb T}$ est une th\'eorie g\'eom\'etrique.

\medskip

\item Montrons que si ${\mathbb T}$ est une th\'eorie r\'eguli\`ere, la cat\'egorie syntactique ${\mathcal C}_{\mathbb T} = {\mathcal C}_{\mathbb T}^{\rm reg} , {\mathcal C}_{\mathbb T}^{\rm coh} , {\mathcal C}_{\mathbb T}^{\rm geo}$ ou ${\mathcal C}_{\mathbb T}^{\rm He}$ est r\'eguli\`ere.

\smallskip

Tout d'abord, c'est une cat\'egorie cart\'esienne.

\smallskip

En effet, elle admet d'une part pour objet terminal la formule $\top$ dans le contexte vide. D'autre part, toute paire de morphismes vers un objet
$$
\varphi (\vec x) \longrightarrow \psi (\vec y) \longleftarrow \varphi' (\vec x')
$$
repr\'esent\'es par deux formules d\'emontrablement fonctionnelles
$$
\theta (\vec x , \vec y) \qquad \mbox{et} \qquad \theta' (\vec x' , \vec y)
$$
admet pour produit fibr\'e la formule de contexte $(\vec x , \vec x')$
$$
(\exists \, \vec y) (\theta (\vec x , \vec y) \wedge \theta' (\vec x' , \vec y)) = \chi (\vec x , \vec x')
$$
munie des deux morphismes
$$
\chi (\vec x , \vec x') \longrightarrow \varphi (\vec x) \qquad \mbox{et} \qquad \chi (\vec x , \vec x') \longrightarrow \varphi' (\vec x')
$$
repr\'esent\'es, pour n'importe quelle substitution \`a $\vec x$ et $\vec x'$ de contextes $\vec u$ et $\vec u'$ disjoints de $\vec x$ et $\vec x'$, par les formules d\'emontrablement fonctionnelles
$$
\chi (\vec x , \vec x') \wedge x_1^{A_1} = u_1^{A_1} \wedge \cdots \wedge x_n^{A_n} = u_n^{A_n} \qquad \mbox{si} \qquad \vec x = (x_1^{A_1} \cdots x_n^{A_n})
$$
et
$$
\chi (\vec x , \vec x') \wedge x'^{B_1}_1 = u'^{B_1}_1 \wedge \cdots \wedge x'^{B_{n'}}_{n'} = u'^{B_{n'}}_{n'} \qquad \mbox{si} \qquad \vec x' = (x'^{B_1}_1 \cdots x'^{B_{n'}}_{n'}) \, .
$$

Pour poursuivre la d\'emonstration, nous avons besoin de l'important lemme suivant:
\end{listeisansmarge}

\begin{lem}\label{lemV69}

Soit ${\mathbb T}$ une th\'eorie r\'eguli\`ere [resp. coh\'erente, resp. g\'eom\'etrique, resp. du premier ordre finitaire] de signature $\Sigma$.

\smallskip

Soit $\varphi (\vec x)$ une formule r\'eguli\`ere [resp. coh\'erente, resp. g\'eom\'etrique, resp. du premier ordre finitaire] de $\Sigma$, vue comme un objet de la cat\'egorie syntactique
$$
{\mathcal C}_{\mathbb T} = {\mathcal C}_{\mathbb T}^{\rm reg} \quad \mbox{[resp. ${\mathcal C}_{\mathbb T}^{\rm coh}$, resp. ${\mathcal C}_{\mathbb T}^{\rm geo}$, resp. ${\mathcal C}_{\mathbb T}^{\rm He}$].}
$$
Alors:

\begin{listeimarge}

\item Les sous-objets de l'objet $\varphi (\vec x)$ de ${\mathcal C}_{\mathbb T}$ correspondent aux formules $\varphi_1 (\vec x)$ de m\^eme contexte $\vec x$ telles que le s\'equent
$$
\varphi_1 \vdash_{\vec x} \varphi
$$
soit d\'emontrable dans ${\mathbb T}$ suivant les r\`egles d'inf\'erence de la logique r\'eguli\`ere [resp. coh\'erente, resp. g\'eom\'etrique, resp. du premier ordre finitaire].

\medskip

\item Deux sous-objets $\varphi_1 (\vec x)$, $\varphi_2 (\vec x)$ de $\varphi (\vec x)$ v\'erifient la relation d'inclusion
$$
\varphi_1 (\vec x) \leq \varphi_2 (\vec x)
$$
si et seulement si le s\'equent
$$
\varphi_1 \vdash_{\vec x} \varphi_2
$$
est d\'emontrable dans ${\mathbb T}$.

\medskip

\item En particulier, deux formules $\varphi_1 (\vec x) , \varphi_2 (\vec x)$ d\'efinissent le m\^eme sous-objet de $\varphi (\vec x)$ si et seulement si le double s\'equent
$$
\varphi_1 \dashv \, \vdash_{\vec x} \varphi_2
$$
est d\'emontrable dans ${\mathbb T}$.
\end{listeimarge}
\end{lem}

\begin{demolem}
\begin{listeisansmarge}
\item Un morphisme
$$
\chi (\vec u) \longrightarrow \varphi (\vec x)
$$
repr\'esent\'e par une formule d\'emontrablement fonctionnelle
$$
\theta (\vec u , \vec x)
$$
est un monomorphisme si et seulement si le morphisme diagonal
$$
\chi \longrightarrow \chi \times_{\varphi} \chi
$$
est un isomorphisme. Si $\vec u = (u_1^{A_1} \cdots u_n^{A_n})$ et $\vec u' = (u'^{A_1} \cdots u'^{A_n})$ est un contexte disjoint de $\vec u$ dont les variables sont affect\'ees aux m\^emes sortes, cela \'equivaut \`a demander que le s\'equent
$$
\theta (\vec u , \vec x) \wedge \theta (\vec u' , \vec x) \vdash \ u_1^{A_1} = u'^{A_1}_1 \wedge \cdots \wedge u_n^{A_n} = u'^{A_n}_n
$$
soit d\'emontrable dans la th\'eorie ${\mathbb T}$.

\smallskip

Alors le morphisme repr\'esent\'e par $\theta (\vec u , \vec x)$ induit un isomorphisme de $\chi (\vec u)$ sur l'objet repr\'esent\'e par la formule de contexte $\vec x$
$$
(\exists \, \vec u) \, \theta (\vec u , \vec x)
$$
laquelle s'inscrit dans le s\'equent d\'emontrable
$$
(\exists \, \vec u) \, \theta (\vec u , \vec x) \vdash_{\vec x} \varphi (\vec x) \, .
$$
R\'eciproquement, toute formule $\varphi (\vec x)$ de contexte $\vec x$ telle que le s\'equent
$$
\varphi_1 \vdash_{\vec x} \varphi
$$
soit d\'emontrable d\'efinit un sous-objet de $\varphi (\vec x)$.

\medskip

\item r\'esulte de (i) puisque deux sous-objets $\varphi_1 (\vec x)$, $\varphi_2 (\vec x)$ de $\varphi (\vec x)$ satisfont la relation d'inclusion
$$
\varphi_1 (\vec x) \leq \varphi_2 (\vec x)
$$
si et seulement si $\varphi_1 (\vec x)$ est un sous-objet de $\varphi_2 (\vec x)$.

\medskip

\item r\'esulte de (ii). 

\end{listeisansmarge}
\end{demolem}

\bigskip

\noindent {\bf Suite de la d\'emonstration de la proposition \ref{propV66}:}

\begin{listeisansmarge}

\item[(iii)] Montrons maintenant que la cat\'egorie syntactique ${\mathcal C}_{\mathbb T} =  {\mathcal C}_{\mathbb T}^{\rm reg} , {\mathcal C}_{\mathbb T}^{\rm coh} , {\mathcal C}_{\mathbb T}^{\rm geo}$ ou ${\mathcal C}_{\mathbb T}^{\rm He}$ permet les quantifications existentielles et donc est r\'eguli\`ere.

\smallskip

Consid\'erons pour cela un morphisme de ${\mathcal C}_{\mathbb T}$
$$
\varphi (\vec x) \longrightarrow \psi (\vec y)
$$
repr\'esent\'e par une formule d\'emontrablement fonctionnelle $\theta (\vec x , \vec y)$.

\smallskip

Le foncteur d'image r\'eciproque par ce morphisme des sous-objets de $\psi (\vec y)$ associe \`a toute formule
$$
\psi_1 (\vec y) \quad \mbox{v\'erifiant le s\'equent} \quad \psi_1 \vdash_{\vec y} \psi
$$
la formule de contexte $\vec x$
$$
(\exists \, \vec y) (\theta (\vec x , \vec y) \wedge \psi_1 (\vec y))
$$
qui v\'erifie le s\'equent
$$
(\exists \, \vec y) (\theta \wedge \psi_1) \vdash_{\vec x} \varphi \, .
$$

Si $\varphi_1$ est une formule de contexte $\vec x$ qui v\'erifie le s\'equent
$$
\varphi_1 \vdash_{\vec x} \varphi \, ,
$$
le s\'equent
$$
\varphi_1 \vdash_{\vec x} (\exists \, \vec y) (\theta \wedge \psi_1)
$$
\'equivaut au s\'equent
$$
\theta \wedge \varphi_1 \vdash_{\vec x , \vec y} \theta \wedge \psi_1
$$
et donc au s\'equent
$$
(\exists \, \vec x)(\theta \wedge \varphi_1) \vdash_{\vec y} \psi_1 \, .
$$

En effet, il r\'esulte des s\'equents d\'emontrables

\medskip

$
\left\{\begin{matrix}
\varphi \vdash_{\vec x} (\exists \, \vec y) \, \theta \, , \\
{ \ } \\
\theta (\vec x , \vec y) \wedge \theta (\vec x , \vec y') \vdash \vec y = \vec y'
\end{matrix} \right.
$

\medskip

\noindent que les deux applications
$$
\varphi_1 (\vec x) \longmapsto \varphi_1 (\vec x) \wedge \theta (\vec x , \vec y)
$$
et
$$
\theta_1 (\vec x , \vec y) \longmapsto (\exists \, \vec y) \, \theta_1 (\vec x , \vec y)
$$
d\'efinissent deux isomorphismes r\'eciproques entre l'ensemble ordonn\'e des sous-objets $\varphi_1 (\vec x)$ de $\varphi (\vec x)$ et celui des sous-objets $\theta_1 (\vec x , \vec y)$ de $\theta (\vec x , \vec y)$.

\smallskip

Ainsi, le foncteur d'image r\'eciproque des sous-objets de $\psi (\vec y)$
$$
\psi_1 (\vec y) \longmapsto (\exists \, \vec y)(\theta \wedge \psi_1)
$$
admet pour adjoint \`a gauche le foncteur
$$
\varphi_1 (\vec x) \longmapsto (\exists \, \vec x) (\theta \wedge \varphi_1) \, .
$$

Il r\'esulte du lemme \ref{lemV68} (ii) que ce foncteur est respect\'e par tout morphisme de changement de base
$$
\theta' (\vec z , \vec y) : \chi (\vec z) \longrightarrow \psi (\vec y) \, .
$$
Cela ach\`eve de montrer que la cat\'egorie ${\mathcal C}_{\mathbb T}$ est r\'eguli\`ere.

\smallskip

Si ${\mathbb T}$ est une th\'eorie coh\'erente ou du premier ordre finitaire [resp. g\'eom\'etrique], la cat\'egorie r\'eguli\`ere ${\mathcal C}_{\mathbb T}^{\rm coh}$ ou ${\mathcal C}_{\mathbb T}^{\rm He}$ [resp. ${\mathcal C}_{\mathbb T}^{\rm geo}$] admet des r\'eunions finies [resp. arbitraires] de sous-objets. En effet, pour tout objet repr\'esent\'e par une formule $\varphi (\vec x)$ et toute famille finie [resp. arbitraire] de sous-objets repr\'esent\'es par des formules $\varphi_i (\vec x)$, $1 \leq i \leq k$ [resp. $i \in I$], v\'erifiant le s\'equent
$$
\varphi_i \vdash_{\vec x} \varphi \, ,
$$
le sous-objet repr\'esent\'e par la formule
$$
\varphi_1 \vee \cdots \vee \varphi_k \qquad \mbox{[resp.} \quad \bigvee_{i \in I} \varphi_i \ \mbox{]}
$$
est la r\'eunion de ces sous-objets d'apr\`es les r\`egles des disjonctions finitaires [resp. infinitaires].

\smallskip

La formation des r\'eunions est respect\'ee par les foncteurs de changement de base par les morphismes
$$
\theta (\vec u , \vec x) : \chi (\vec u) \longrightarrow \varphi (\vec x)
$$
qui associent \`a toute formule $\varphi' (\vec x)$ v\'erifiant le s\'equent
$$
\varphi' \vdash_{\vec x} \varphi
$$
la formule de contexte $\vec u$
$$
(\exists \, \vec x) (\theta (\vec u , \vec x) \wedge \varphi' (\vec x)) \, .
$$

Cela r\'esulte de ce que le symbole $\vee$ [resp. $\bigvee$] commute avec le symbole $\wedge$ d'apr\`es les r\`egles de distributivit\'e et avec le symbole $\exists$ d'apr\`es le lemme \ref{lemV68} (i).

\smallskip

Ainsi, la cat\'egorie ${\mathcal C}_{\mathbb T}^{\rm coh}$ ou ${\mathcal C}_{\mathbb T}^{\rm He}$ [resp. ${\mathcal C}_{\mathbb T}^{\rm geo}$] est coh\'erente [resp. g\'eom\'etrique].

\smallskip

Montrons enfin que si ${\mathbb T}$ est une th\'eorie du premier ordre finitaire la cat\'egorie ${\mathcal C}_{\mathbb T}^{\rm He}$ est de Heyting, c'est-\`a-dire admet les quantifications universelles $\forall_{\theta}$ associ\'ees \`a des morphismes
$$
\theta (\vec x , \vec y) : \varphi (\vec x) \longrightarrow \psi (\vec y)
$$
ainsi que des foncteurs $\Rightarrow$ d'implications entre sous-objets et a fortiori $\neg$ de n\'egation des sous-objets.

\smallskip

Le foncteur d'image r\'eciproque par $\theta (\vec x , \vec y) : \varphi (\vec x) \to \psi (\vec y)$ des sous-objets repr\'esent\'es par des formules $\psi' (\vec y)$ v\'erifiant le s\'equent
$$
\psi' \vdash_{\vec y} \psi
$$
s'\'ecrit
$$
\psi' \longmapsto (\exists \, \vec y) (\theta (\vec x , \vec y) \wedge \psi' (\vec y)) \, .
$$

D'apr\`es les r\`egles des quantifications universelles et celles des implications, il admet pour adjoint \`a droite le foncteur
$$
\varphi ' (\vec x) \longmapsto (\forall \, \vec x) (\theta (\vec x , \vec y) \Rightarrow \varphi' (\vec x)) \, .
$$

En effet, pour toute formule $\varphi' (\vec x)$ v\'erifiant le s\'equent
$$
\varphi' \vdash_{\vec x} \varphi \, ,
$$
le s\'equent
$$
(\exists \, \vec y) (\theta (\vec x , \vec y) \wedge \psi' (\vec y)) \vdash_{\vec x} \varphi' (\vec x)
$$
\'equivaut au s\'equent
$$
\theta (\vec x , \vec y) \wedge \psi' (\vec y) \vdash_{\vec x , \vec y} \varphi' (\vec x)
$$
donc aussi au s\'equent
$$
\psi' (\vec y) \vdash_{\vec x , \vec y} (\theta (\vec x , \vec y) \Rightarrow \varphi' (\vec x))
$$
puis au s\'equent
$$
\psi' (\vec y) \vdash_{\vec y} (\forall \, \vec x) (\theta (\vec x , \vec y) \Rightarrow \varphi' (\vec x)) \, .
$$

Les foncteurs de quantifications universelles et ceux d'images r\'eciproques des sous-objets associ\'es \`a un carr\'e cart\'esien de ${\mathcal C}_{\mathbb T}^{\rm He}$ commutent entre eux puisqu'ils sont adjoints \`a droite des foncteurs d'images r\'eciproques et de quantifications existentielles et que ceux-ci commutent entre eux.

\smallskip

Enfin, il r\'esulte des r\`egles d'implication [resp. de n\'egation] que l'application
$$
\begin{matrix}
&(\varphi_1 , \varphi_2) &\longmapsto &(\varphi_1 \Rightarrow \varphi_2) \wedge \varphi \\
\mbox{[resp.} &\hfill \varphi_1 &\longmapsto &(\neg \, \varphi_1) \wedge \varphi \hfill
\end{matrix}
$$
d\'efinit un foncteur d'implication [resp. de n\'egation] sur les sous-objets $\varphi_1 (\vec x) , \varphi_2 (\vec x)$ [resp. $\varphi_1 (\vec x)$] de tout objet $\varphi (\vec x)$.

\smallskip

Par d\'efinition, le foncteur
$$
\varphi_2 (\vec x) \longmapsto (\varphi_1 (\vec x) \Rightarrow \varphi_2 (\vec x))
$$
est adjoint \`a droite du foncteur
$$
\varphi' (\vec x) \longmapsto \varphi_1 (\vec x) \wedge \varphi' (\vec x) \, .
$$

Il est respect\'e par le foncteur d'images r\'eciproques associ\'e \`a un morphisme de changement de base car celui-ci est adjoint \`a droite du foncteur de quantification existentielle, lequel respecte le foncteur $\varphi_1 (\vec x) \wedge \bullet$ d'apr\`es le lemme \ref{lemV68} (ii).

\smallskip

Pour finir, le foncteur $\varphi_1 \mapsto \neg \, \varphi_1$ est \'egalement respect\'e par tout morphisme de changement de base puisqu'il se confond avec le foncteur $\varphi_1 \mapsto (\varphi_1 \Rightarrow \, \perp)$.

\medskip

\item[(iv)] Par d\'efinition, $M_{\mathbb T}$ est une $\Sigma$-structure dans la cat\'egorie syntactique
$$
{\mathcal C}_{\mathbb T} = {\mathcal C}_{\mathbb T}^{\rm reg} \qquad \mbox{[resp. ${\mathcal C}_{\mathbb T}^{\rm coh}$, resp. ${\mathcal C}_{\mathbb T}^{\rm geo}$, resp. ${\mathcal C}_{\mathbb T}^{\rm He}$].}
$$
Le fait que ce soit un mod\`ele de la th\'eorie ${\mathbb T}$ r\'esulte du lemme \ref{lemV69} (ii) combin\'e avec le lemme suivant:
\end{listeisansmarge}

\begin{lem}\label{lemV610}

Soit ${\mathbb T}$ une th\'eorie r\'eguli\`ere [resp. coh\'erente, resp. g\'eom\'etrique, resp. du premier ordre finitaire] de signature $\Sigma$.

\smallskip

Soit $\varphi$ une formule r\'eguli\`ere [resp. coh\'erente, resp. g\'eom\'etrique, resp. du premier ordre finitaire] de $\Sigma$, et soit $\vec x = (x_1^{A_1} \cdots x_n^{A_n})$ un contexte de $\varphi$.

\smallskip

Alors l'interpr\'etation de $\varphi (\vec x)$ dans la $\Sigma$-structure $M_{\mathbb T}$ de ${\mathcal C}_{\mathbb T} = {\mathcal C}_{\mathbb T}^{\rm reg}$ [resp. ${\mathcal C}_{\mathbb T}^{\rm coh}$, resp. ${\mathcal C}_{\mathbb T}^{\rm geo}$, resp. ${\mathcal C}_{\mathbb T}^{\rm He}$] est le sous-objet
$$
\varphi (\vec x) \xhookrightarrow{ \ { \ } \ } \top (x_1^{A_1} \cdots x_n^{A_n}) \, .
$$
\end{lem}

\begin{demo}

La $\Sigma$-structure $M_{\mathbb T}$ associe \`a toute sorte $A_i$ l'objet $\top (x_i^{A_i})$ de ${\mathcal C}_{\mathbb T}$. De plus, le produit dans ${\mathcal C}_{\mathbb T}$ des objet $\top (x_1^{A_1}), \cdots , \top (x_n^{A_n})$ est l'objet $\top (x_1^{A_1} \cdots x_n^{A_n})$.

\smallskip

Comme $M_{\mathbb T}$ associe \`a tout symbole de fonction de $\Sigma$
$$
f : A_1 \cdots A_n \longrightarrow B
$$
le morphisme
$$
\top (x_1^{A_1} \cdots x_n^{A_n}) \xrightarrow{ \ y^B = f (x_1^{A_1} \cdots x_n^{A_n}) \ } \top (y^B)
$$
et \`a tout symbole de relation de $\Sigma$
$$
\xymatrix{
R \ \ \ar@{>->}[r] &A_1 \cdots A_n
}
$$
le sous-objet
$$
R (x_1^{A_1} \cdots x_n^{A_n}) \xhookrightarrow{ \ { \ } \ } \top (x_1^{A_1} \cdots x_n^{A_n}) \, ,
$$
on voit d\'ej\`a que toute formule atomique $\varphi (\vec x)$ de $\Sigma$ dans le contexte $\vec x = (x_1^{A_1} \cdots x_n^{A_n})$ s'interpr\`ete comme le sous-objet
$$
\varphi (\vec x) \xhookrightarrow{ \ { \ } \ } \top (x_1^{A_1} \cdots x_n^{A_n}) \, .
$$

La conclusion du lemme r\'esulte alors de ce que:

\medskip

$\left\{ \begin{matrix}
\bullet &\mbox{le symbole $\wedge$ de conjonction finie s'interpr\`ete dans ${\mathcal C}_{\mathbb T}$ comme une intersection finie de sous-objets,} \hfill \\
{ \ } \\
\bullet &\mbox{le symbole $\vee$ [resp. $\bigvee$] de disjonction finie [resp. arbitraire] s'interpr\`ete dans ${\mathcal C}_{\mathbb T}$ comme une r\'eunion} \hfill \\
&\mbox{finie [resp. arbitraire] de sous-objets,} \hfill \\
{ \ } \\
\bullet &\mbox{le symbole de quantification existentielle [resp. universelle]} \hfill \\
{ \ } \\
&(\exists \, \vec y) \qquad \mbox{[resp.} \quad (\forall \, \vec y) \ \mbox{]} \\
{ \ } \\
&\mbox{sur une partie $\vec y = (y_1^{B_1} \cdots y_m^{B_m})$ d'un contexte $(\vec x , \vec y) = (x_1^{A_1} \cdots x_n^{A_n} \, y_1^{B_1} \cdots y_m^{B_m})$ s'interpr\`ete} \hfill \\
&\mbox{comme le foncteur} \hfill \\
&\exists_p \qquad \mbox{[resp.} \quad \forall_p \ \mbox{]} \\
&\mbox{associ\'e au morphisme de projection} \hfill \\
{ \ } \\
&p : \top (x_1^{A_1} \cdots x_n^{A_n} \, y_1^{B_1} \cdots y_m^{B_m}) \longrightarrow \top (x'^{A_1}_1 \cdots x'^{A_n}_n) \\
{ \ } \\
&\mbox{d\'efini par la formule d\'emontrablement fonctionnelle} \hfill \\
{ \ } \\
&x_1^{A_1} = x'^{A_1}_1 \wedge \cdots \wedge x_n^{A_n} = x'^{A_n}_n \, , \\
{ \ } \\
\bullet &\mbox{le symbole $\Rightarrow$ d'implication [resp. $\neg$ de n\'egation] s'interpr\`ete comme le foncteur d'implication} \hfill \\
&\mbox{[resp. de n\'egation] sur la cat\'egorie des sous-objets de l'objet $\top (x_1^{A_1} \cdots x_n^{A_n})$.} \hfill
\end{matrix} \right.
$

\bigskip

Cela termine la d\'emonstration du lemme et donc aussi de la proposition \ref{propV66}. 

\end{demo}

\pagebreak

On note l'importante cons\'equence suivante des lemmes \ref{lemV69} et \ref{lemV610}:

\begin{cor}\label{corV611}

Soit ${\mathbb T}$ une th\'eorie r\'eguli\`ere [resp. coh\'erente, resp. g\'eom\'etrique, resp. du premier ordre finitaire] de signature $\Sigma$.

\smallskip

Alors un s\'equent r\'egulier [resp. coh\'erent, resp. g\'eom\'etrique, resp. du premier ordre finitaire] de $\Sigma$
$$
\varphi \vdash_{\vec x} \psi
$$
est d\'emontrable dans la th\'eorie ${\mathbb T}$ suivant les r\`egles d'inf\'erence de la logique r\'eguli\`ere [resp. coh\'erente, resp. g\'eom\'etrique, resp. du premier ordre finitaire] si et seulement si il est satisfait par le mod\`ele
$$
M_{\mathbb T}
$$
de ${\mathbb T}$ dans la cat\'egorie syntactique
$$
{\mathcal C}_{\mathbb T} = {\mathcal C}_{\mathbb T}^{\rm reg} \qquad \mbox{[resp. ${\mathcal C}_{\mathbb T}^{\rm coh}$, resp. ${\mathcal C}_{\mathbb T}^{\rm geo}$, resp. ${\mathcal C}_{\mathbb T}^{\rm He}$].}
$$
\end{cor}

\begin{demo}

Notant $\vec x = (x_1^{A_1} \cdots x_n^{A_n})$, les formules $\varphi$ et $\psi$ de contexte $\vec x$ s'interpr\`etent dans le mod\`ele $M_{\mathbb T}$ comme les sous-objets
$$
\varphi (\vec x) \xhookrightarrow{ \ { \ } \ } \top (x_1^{A_1} \cdots x_n^{A_n})
$$
et
$$
\psi (\vec x) \xhookrightarrow{ \ { \ } \ } \top (x_1^{A_1} \cdots x_n^{A_n}) \, .
$$
D'apr\`es le lemme \ref{lemV69} (ii), ces sous-objets satisfont la relation d'inclusion
$$
\varphi (\vec x) \leq \psi (\vec x)
$$
si et seulement si le s\'equent
$$
\varphi \vdash_{\vec x} \psi
$$
est d\'emontrable dans ${\mathbb T}$. 

\end{demo}

\subsection{Construction des cat\'egories syntactiques cart\'esiennes}\label{subsec564}

\medskip

Afin de construire la cat\'egorie syntactique cart\'esienne d'une th\'eorie alg\'ebrique, de Horn ou plus g\'en\'eralement cart\'esienne, on a d'abord besoin de poser la d\'efinition suivante:

\begin{defn}\label{defV612}

Soit $\Sigma$ une signature.

\smallskip

Soit ${\mathbb T}$ une th\'eorie cart\'esienne de signature $\Sigma$.

\smallskip

Une formule r\'eguli\`ere de $\Sigma$ dans un contexte $\vec x$
$$
\varphi (\vec x)
$$
est dite ${\mathbb T}$-cart\'esienne si elle est ${\mathbb T}$-d\'emontrablement \'equivalente, suivant les r\`egles d'inf\'erence de la logique r\'eguli\`ere, \`a une formule de la forme
$$
(\exists \, \vec y) \, \varphi' (\vec x , \vec y) \, ,
$$
o\`u $\varphi' (\vec x , \vec y)$ est une formule de Horn telle que le s\'equent
$$
\varphi' (\vec x , \vec y) \wedge \varphi' (\vec x , \vec y') \vdash \ y_1^{B_1} = y'^{B_1}_1 \wedge \cdots \wedge y_k^{B_k} = y'^{B_k}_k \, ,
$$
avec $\vec y = (y_1^{B_1} \cdots y_k^{B_k})$ et $\vec y' = (y'^{B_1}_1 \cdots y'^{B_k}_k)$, soit ${\mathbb T}$-d\'emontrable.
\end{defn}

\begin{remarksqed}
\begin{listeisansmarge}
\item Alors que les notions de formules de Horn, r\'eguli\`ere, coh\'erente, g\'eom\'etrique ou du premier ordre finitaire ou infinitaire ne d\'ependent que de la signature $\Sigma$, la notion de formule cart\'esienne d\'epend de la th\'eorie cart\'esienne ${\mathbb T}$ de signature $\Sigma$ que l'on consid\`ere.

\medskip

\item En particulier, toute formule de Horn de $\Sigma$ est ${\mathbb T}$-cart\'esienne.

\medskip

\item Cette d\'efinition s'applique a fortiori si ${\mathbb T}$ est une th\'eorie alg\'ebrique ou plus g\'en\'eralement de Horn. 
\end{listeisansmarge}
\end{remarksqed}

\bigskip

Voici la d\'efinition de la cat\'egorie syntactique cart\'esienne d'une th\'eorie cart\'esienne et celle de son mod\`ele canonique:

\begin{defn}\label{defV613}

Soit $\Sigma$ une signature.

\smallskip

Soit ${\mathbb T}$ une th\'eorie cart\'esienne de signature $\Sigma$.

\smallskip

On appelle cat\'egorie syntactique cart\'esienne de ${\mathbb T}$ et on note
$$
{\mathcal C}_{\mathbb T} = {\mathcal C}_{\mathbb T}^{\rm cart}
$$
la cat\'egorie cart\'esienne munie d'un mod\`ele canonique $M_{\mathbb T}$ ainsi d\'efinie:

\begin{listeimarge}

\item La cat\'egorie ${\mathcal C}_{\mathbb T}^{\rm cart}$ est la sous-cat\'egorie de la cat\'egorie syntactique r\'eguli\`ere de ${\mathbb T}$
$$
{\mathcal C}_{\mathbb T}^{\rm reg}
$$
dont

\medskip

$
\left\{ \begin{matrix}
\bullet &\mbox{les objets sont les formules ${\mathbb T}$-cart\'esiennes} \hfill \\
{ \ } \\
&\varphi (\vec x) \, , \\
{ \ } \\
\bullet &\mbox{les morphismes sont les formules ${\mathbb T}$-d\'emontrablement fonctionnelles} \hfill \\
{ \ } \\
&\varphi (\vec x) \xrightarrow{ \ \theta (\vec x , \vec y) \ } \psi (\vec y) \\
{ \ } \\
&\mbox{qui sont ${\mathbb T}$-cart\'esiennes.} \hfill 
\end{matrix}\right.
$

\bigskip

\item Le mod\`ele canonique $M_{\mathbb T}$ de ${\mathbb T}$ dans ${\mathcal C}_{\mathbb T}^{\rm cart}$ est l'unique $\Sigma$-structure de ${\mathcal C}_{\mathbb T}^{\rm cart}$ dont l'image par le foncteur de plongement
$$
{\mathcal C}_{\mathbb T}^{\rm cart} \xhookrightarrow{ \ { \ } \ } {\mathcal C}_{\mathbb T}^{\rm reg}
$$
est le mod\`ele canonique de ${\mathbb T}$ dans ${\mathcal C}_{\mathbb T}^{\rm reg}$.
\end{listeimarge}
\end{defn}

\begin{remarksqed}
\begin{listeisansmarge}
\item La liste des objets de ${\mathcal C}_{\mathbb T}^{\rm cart}$ tout comme celle de ses morphismes d\'epend de la th\'eorie cart\'esienne ${\mathbb T}$ et pas seulement de sa signature $\Sigma$.

\medskip

\item Cette d\'efinition s'applique a fortiori si ${\mathbb T}$ est une th\'eorie alg\'ebrique ou plus g\'en\'eralement de Horn.

\medskip

\item Ainsi, une th\'eorie cart\'esienne ${\mathbb T}$ poss\`ede une cat\'egorie syntactique cart\'esienne ${\mathcal C}_{\mathbb T}^{\rm cart}$ en plus de ses cat\'egories syntactiques r\'eguli\`ere ${\mathcal C}_{\mathbb T}^{\rm reg}$, coh\'erente ${\mathcal C}_{\mathbb T}^{\rm coh}$, g\'eom\'etrique ${\mathcal C}_{\mathbb T}^{\rm geo}$ et de Heyting ${\mathcal C}_{\mathbb T}^{\rm He}$ [et m\^eme alg\'ebrique ${\mathcal C}_{\mathbb T}^{\rm alg}$ si ${\mathbb T}$ est une th\'eorie alg\'ebrique]. 
\end{listeisansmarge}
\end{remarksqed}

\bigskip

On v\'erifie que cette d\'efinition est valide:

\begin{prop}\label{propV614}

Soit ${\mathbb T}$ une th\'eorie cart\'esienne de signature $\Sigma$.

\smallskip

Alors:

\begin{listeimarge}

\item La loi de composition des morphismes de ${\mathcal C}_{\mathbb T}^{\rm reg}$ respecte le sous-ensemble des morphismes de ${\mathcal C}_{\mathbb T}^{\rm cart}$, si bien que ${\mathcal C}_{\mathbb T}^{\rm cart}$ est une sous-cat\'egorie de ${\mathcal C}_{\mathbb T}^{\rm reg}$.

\smallskip

Elle est petite.

\medskip

\item La cat\'egorie ${\mathcal C}_{\mathbb T}^{\rm cart}$ est cart\'esienne et le foncteur de plongement
$$
{\mathcal C}_{\mathbb T}^{\rm cart} \xhookrightarrow{ \ { \ } \ } {\mathcal C}_{\mathbb T}^{\rm reg}
$$
respecte les limites finies.

\medskip

\item Le foncteur ${\mathcal C}_{\mathbb T}^{\rm cart} \hookrightarrow {\mathcal C}_{\mathbb T}^{\rm reg}$ est aussi conservatif.

\medskip

\item Le mod\`ele canonique $M_{\mathbb T}$ de ${\mathbb T}$ dans ${\mathcal C}_{\mathbb T}^{\rm reg}$ provient d'une $\Sigma$-structure $M_{\mathbb T}$ de ${\mathcal C}_{\mathbb T}^{\rm cart}$ qui est un mod\`ele de ${\mathbb T}$.
\end{listeimarge}
\end{prop}

\begin{demo}
\begin{listeisansmarge}
\item Pour tout objet de ${\mathcal C}_{\mathbb T}^{\rm cart}$ c'est-\`a-dire toute formule ${\mathbb T}$-cart\'esienne dans un contexte $\vec x = (x_1^{A_1} \cdots x_n^{A_n})$
$$
\varphi (\vec x) \, ,
$$
son morphisme identit\'e
$$
{\rm id} : \varphi (\vec x) \longrightarrow \varphi (\vec x') \, , \qquad \mbox{avec} \quad \vec x' = (x'^{A_1}_1 \cdots x'^{A_n}_n) \, ,
$$
est la formule d\'emontrablement fonctionnelle de contexte $(\vec x , \vec x')$
$$
\varphi (\vec x) \wedge x_1^{A_1} = x'^{A_1}_1 \wedge \cdots \wedge x_n^{A_n} = x'^{A_n}_n \, .
$$
Cette formule est ${\mathbb T}$-cart\'esienne puisqu'il en est ainsi de $\varphi (\vec x)$.

\smallskip

Puis consid\'erons des morphismes de ${\mathcal C}_{\mathbb T}^{\rm cart}$
$$
\varphi(\vec x) \xrightarrow{ \ \theta_1 (\vec x , \vec y) \ } \psi (\vec y) \xrightarrow{ \ \theta_2 (\vec y , \vec z) \ } \chi (\vec z) \, .
$$
Les formules ${\mathbb T}$-d\'emontrablement fonctionnelles
$$
\theta_1 (\vec x , \vec y) \qquad \mbox{et} \qquad \theta_2 (\vec y , \vec z)
$$
sont ${\mathbb T}$-cart\'esiennes c'est-\`a-dire ${\mathbb T}$-d\'emontrablement \'equivalentes \`a des formules de la forme
$$
(\exists \, \vec w_1) \, \theta'_1 (\vec x , \vec y , \vec w_1)
$$
et
$$
(\exists \, \vec w_2) \, \theta'_2 (\vec y , \vec z , \vec w_2)
$$
o\`u $\theta'_1$ et $\theta'_2$ sont des formules de Horn telles que les s\'equents
$$
\theta'_1 (\vec x , \vec y , \vec w_1) \wedge \theta'_1 (\vec x , \vec y , \vec w'_1) \vdash \vec w_1 = \vec w'_1
$$
et
$$
\theta'_2 (\vec y , \vec z , \vec w_2) \wedge \theta'_2 (\vec y , \vec z , \vec w'_2) \vdash \vec w_2 = \vec w'_2
$$
soient ${\mathbb T}$-d\'emontrables.

\smallskip

Alors la formule compos\'ee
$$
(\exists \, \vec y) (\theta_1 (\vec x , \vec y) \wedge \theta_2 (\vec y , \vec z))
$$
est ${\mathbb T}$-d\'emontrablement \'equivalente \`a la formule
$$
(\exists \, \vec y) (\exists \, \vec w_1) (\exists \, \vec w_2) (\theta'_1 (\vec x , \vec y , \vec w_1) \wedge \theta'_2 (\vec y , \vec z , \vec w_2)) \, .
$$
Or, la formule
$$
\theta'_1 (\vec x , \vec y , \vec w_1) \wedge \theta'_2 (\vec y , \vec z , \vec w_2)
$$
est de Horn, et le s\'equent
$$
\theta'_1 (\vec x , \vec y , \vec w_1) \wedge \theta'_2 (\vec y , \vec z , \vec w_2) \wedge \theta'_1 (\vec x , \vec y' , \vec w'_1) \wedge \theta'_2 (\vec y' , \vec z , \vec w'_2) \vdash \ \vec y = \vec y' \wedge \vec w_1 = \vec w'_1 \wedge \vec w_2 = \vec w'_2
$$
est ${\mathbb T}$-d\'emontrable puisqu'il en est ainsi des s\'equents
$$
(\exists \, \vec w_1) \, \theta'_1 (\vec x , \vec y , \vec w_1) \wedge (\exists \, \vec w'_1) \, \theta'_1 (\vec x , \vec y' , \vec w'_1) \vdash \ \vec y = \vec y' \, ,
$$
$$
\theta'_1 (\vec x , \vec y , \vec w_1) \wedge \theta'_1 (\vec x , \vec y , \vec w'_1) \vdash \vec w_1 = \vec w'_1 \, ,
$$
$$
\theta'_2 (\vec y , \vec z , \vec w_2) \wedge \theta'_2 (\vec y , \vec z , \vec w'_2) \vdash \vec w_2 = \vec w'_2 \, .
$$

Ainsi, la formule compos\'ee
$$
(\exists \, \vec y) (\theta_1 (\vec x , \vec y) \wedge \theta_2 (\vec y , \vec z))
$$
est ${\mathbb T}$-cart\'esienne, et ${\mathcal C}_{\mathbb T}^{\rm cart}$ est bien une sous-cat\'egorie de ${\mathcal C}_{\mathbb T}^{\rm reg}$.

\smallskip

Elle est petite puisque ${\mathcal C}_{\mathbb T}^{\rm reg}$ est petite.

\medskip

\item L'objet terminal ${\mathcal C}_{\mathbb T}^{\rm reg}$ est la formule $\top$ dans le contexte vide. C'est une formule de Horn et donc un objet de ${\mathcal C}_{\mathbb T}^{\rm cart}$.

\smallskip

Pour toute formule ${\mathbb T}$-cart\'esienne dans un contexte $\vec x$
$$
\varphi (\vec x) \, ,
$$
vue comme un objet de ${\mathcal C}_{\mathbb T}^{\rm cart}$ et a fortiori de ${\mathcal C}_{\mathbb T}^{\rm reg}$, l'unique morphisme de ${\mathcal C}_{\mathbb T}^{\rm reg}$
$$
\varphi (\vec x) \longrightarrow \top
$$
est la formule d\'emontrablement fonctionnelle $\varphi (\vec x)$.

\smallskip

Cette formule \'etant ${\mathbb T}$-cart\'esienne, c'est un morphisme de ${\mathcal C}_{\mathbb T}^{\rm cart}$.

\smallskip

Ainsi, ${\mathcal C}_{\mathbb T}^{\rm cart}$ admet pour objet terminal l'objet terminal de ${\mathcal C}_{\mathbb T}^{\rm reg}$.

\smallskip

Puis consid\'erons une paire de morphismes de ${\mathcal C}_{\mathbb T}^{\rm cart}$
$$
\varphi (\vec x) \longrightarrow \psi (\vec y) \longleftarrow \chi (\vec z)
$$
repr\'esent\'es par deux formules d\'emontrablement fonctionnelles ${\mathbb T}$-cart\'esiennes
$$
\theta_1 (\vec x , \vec y)
$$
et
$$
\theta_2 (\vec z , \vec y) \, .
$$
Leur produit fibr\'e dans ${\mathcal C}_{\mathbb T}^{\rm reg}$ est la formule de contexte $(\vec x , \vec z)$
$$
(\exists \, \vec y) (\theta_1 (\vec x , \vec y) \wedge \theta_2 (\vec z , \vec y)) = \Omega (\vec x , \vec z) \, .
$$
On v\'erifie comme dans la preuve de (i) que cette formule est ${\mathbb T}$-cart\'esienne puisqu'il en est ainsi des formules $\theta_1 (\vec x , \vec y)$ et $\theta_2 (\vec z , \vec y)$. Donc c'est un objet de la sous-cat\'egorie ${\mathcal C}_{\mathbb T}^{\rm cart}$ de ${\mathcal C}_{\mathbb T}^{\rm reg}$.

\smallskip

Les projections canoniques dans ${\mathcal C}_{\mathbb T}^{\rm reg}$ de cet objet sur les objets $\varphi (\vec x)$ et $\psi (\vec z)$ sont d\'efinies par les formules d\'emontrablement fonctionnelles
$$
\Omega (\vec x , \vec z) \xrightarrow{ \ \Omega (\vec x , \vec z) \, \wedge \, \vec x \, = \, \vec x' \ } \varphi (\vec x')
$$
et
$$
\Omega (\vec x , \vec z) \xrightarrow{ \ \Omega (\vec x , \vec z) \, \wedge \, \vec z \, = \, \vec z' \ } \chi (\vec z') \, .
$$
Or, les formules
$$
\Omega (\vec x , \vec z) \wedge \vec x = \vec x' \qquad \mbox{et} \qquad \Omega (\vec x , \vec z) \wedge \vec z = \vec z'
$$
sont ${\mathbb T}$-cart\'esiennes puisqu'il en est ainsi de la formule $\Omega (\vec x , \vec z)$.

\smallskip

Enfin, pour tout carr\'e commutatif de ${\mathcal C}_{\mathbb T}^{\rm cart}$
$$
\xymatrix{
\Delta (\vec w) \ar[d]_{\theta_4 (\vec w , \vec x)} \ar[rr]^{\theta_3 (\vec w , \vec z)} &&\chi (\vec z) \ar[d]^{\theta_2 (\vec z , \vec y)} \\
\varphi (\vec x) \ar[rr]^{\theta_1 (\vec x , \vec y)} &&\psi (\vec y)
}
$$
sa factorisation canonique dans ${\mathcal C}_{\mathbb T}^{\rm reg}$
$$
\Delta (\vec w) \longrightarrow \Omega (\vec x , \vec z)
$$
est d\'efinie par la formule d\'emontrablement fonctionnelle
$$
\theta_4 (\vec w , \vec x) \wedge \theta_3 (\vec w , \vec z) \, .
$$
Celle-ci est ${\mathbb T}$-cart\'esienne puisqu'il en est ainsi par hypoth\`ese des formules $\theta_4 (\vec w , \vec x)$ et $\theta_3 (\vec w , \vec z)$.

\smallskip

Cela montre que la cat\'egorie ${\mathcal C}_{\mathbb T}^{\rm cart}$ a des produits fibr\'es arbitraires et qu'ils sont respect\'es par le foncteur ${\mathcal C}_{\mathbb T}^{\rm cart} \hookrightarrow {\mathcal C}_{\mathbb T}^{\rm reg}$.

\medskip

\item Il s'agit de montrer que si un morphisme de ${\mathcal C}_{\mathbb T}^{\rm cart}$
$$
\varphi (\vec x) \longrightarrow \psi (\vec y)
$$
d\'efini par une formule ${\mathbb T}$-cart\'esienne d\'emontrablement fonctionnelle de $\varphi (\vec x)$ dans $\psi (\vec y)$
$$
\theta (\vec x , \vec y)
$$
est un isomorphisme dans ${\mathcal C}_{\mathbb T}^{\rm reg}$, alors il est un isomorphisme dans ${\mathcal C}_{\mathbb T}^{\rm cart}$.

\smallskip

Ce morphisme est un monomorphisme dans ${\mathcal C}_{\mathbb T}^{\rm reg}$, ce qui signifie que le s\'equent
$$
\theta (\vec x , \vec y) \wedge \theta (\vec x' , \vec y) \vdash \vec x = \vec x'
$$
est ${\mathbb T}$-d\'emontrable, et il est un \'epimorphisme, ce qui signifie que le s\'equent
$$
\psi (\vec y) \vdash (\exists \, \vec x) \, \theta (\vec x , \vec y)
$$
est ${\mathbb T}$-d\'emontrable.

\smallskip

Ainsi, la formule ${\mathbb T}$-cart\'esienne
$$
\theta (\vec x , \vec y)
$$
est d\'emontrablement fonctionnelle de $\psi (\vec y)$ dans $\varphi (\vec x)$. Autrement dit, elle d\'efinit un morphisme de ${\mathcal C}_{\mathbb T}^{\rm cart}$
$$
\psi (\vec y) \longrightarrow \varphi (\vec x)
$$
qui est un inverse du morphisme
$$
\varphi (\vec x) \longrightarrow \psi (\vec y)
$$
d\'efini par la m\^eme formule ${\mathbb T}$-cart\'esienne $\theta (\vec x , \vec y)$.

\smallskip

Cela montre que le foncteur ${\mathcal C}_{\mathbb T}^{\rm cart} \hookrightarrow {\mathcal C}_{\mathbb T}^{\rm reg}$ est conservatif.

\medskip

\item La $\Sigma$-structure $M_{\mathbb T}$ de ${\mathcal C}_{\mathbb T}^{\rm reg}$ consiste \`a associer

\medskip

$
\left\{\begin{matrix}
\bullet &\mbox{\`a toute sorte $A$ de $\Sigma$ l'objet d\'efini par la formule} \hfill \\
{ \ } \\
&\top (x^A) \\
{ \ } \\
&\mbox{pour n'importe quelle variable $x^A$ affect\'ee \`a la sorte $A$,} \hfill \\
{ \ } \\
\bullet &\mbox{\`a tout symbole de fonction $f : A_1 \cdots A_n \to B$ le morphisme d\'efini par la formule} \hfill \\
{ \ } \\
&y^B = f(x_1^{A_1} \cdots x_n^{A_n}) \, , \\
{ \ } \\
\bullet &\mbox{\`a tout symbole de relation $R \rightarrowtail A_1 \cdots A_n$ le sous-objet d\'efini par la formule} \hfill \\
{ \ } \\
&R(x_1^{A_1} \cdots x_n^{A_n}) \hookrightarrow \top (x_1^{A_1} \cdots x_n^{A_n}) = \top (x_1^{A_1}) \times \cdots \times \top (x_n^{A_n}) \, . 
\end{matrix}\right.
$

\bigskip

Toutes ces formules sont ${\mathbb T}$-cart\'esiennes et donc $M_{\mathbb T}$ peut \^etre vue comme une $\Sigma$-structure de ${\mathcal C}_{\mathbb T}^{\rm cart}$.

\smallskip

Enfin, la $\Sigma$-structure $M_{\mathbb T}$ dans ${\mathcal C}_{\mathbb T}^{\rm cart}$ v\'erifie les axiomes de la th\'eorie cart\'esienne ${\mathbb T}$.

\smallskip

Cela r\'esulte de ce que le foncteur de plongement
$$
{\mathcal C}_{\mathbb T}^{\rm cart} \xhookrightarrow{ \ { \ } \ } {\mathcal C}_{\mathbb T}^{\rm reg}
$$
respecte les limites finies, est conservatif et pr\'eserve les interpr\'etations des axiomes de ${\mathbb T}$ comme

\medskip

$
\left\{\begin{matrix}
\bullet &\mbox{des relations d'inclusions entre sous-objets d\'efinis par des formules de Horn,} \hfill \\
{ \ } \\
\bullet &\mbox{des isomorphismes entre des objets associ\'es \`a des sortes et } \hfill \\
&\mbox{des sous-objets d\'efinis par des formules d'\'egalit\'e entre termes.} \hfill
\end{matrix}\right.
$

\bigskip

En effet, comme ce foncteur respecte les limites finies, il pr\'eserve les interpr\'etations des formules de Horn, en particulier celles des formules d'\'egalit\'e entre termes.

\smallskip

Cela ach\`eve de montrer la proposition \ref{propV614}. 
\end{listeisansmarge}
\end{demo}

\bigskip

On a encore le parall\`ele suivant du lemme \ref{lemV69}, du lemme \ref{lemV610} et du corollaire \ref{corV611}:

\begin{cor}\label{corV615}

Soit ${\mathbb T}$ une th\'eorie cart\'esienne de signature $\Sigma$.

\smallskip

Soit ${\mathcal C}_{\mathbb T} = {\mathcal C}_{\mathbb T}^{\rm cart}$ sa cat\'egorie syntactique cart\'esienne munie de son mod\`ele canonique $M_{\mathbb T}$.

\smallskip

Alors:

\begin{listeimarge}

\item Les sous-objets d'un objet $\varphi (\vec x)$ de ${\mathcal C}_{\mathbb T}^{\rm cart}$ correspondent aux formules ${\mathbb T}$-cart\'esiennes de m\^eme contexte $\vec x$
$$
\varphi_1 (\vec x)
$$
telles que le s\'equent
$$
\varphi_1 \vdash_{\vec x} \varphi
$$
soit d\'emontrable dans ${\mathbb T}$ suivant les r\`egles d'inf\'erence de la logique r\'eguli\`ere.

\medskip

\item Deux sous-objets $\varphi_1 (\vec x)$, $\varphi_2 (\vec x)$ d'un objet $\varphi (\vec x)$ de ${\mathcal C}_{\mathbb T}^{\rm cart}$ v\'erifient la relation d'inclusion
$$
\varphi_1 (\vec x) \leq \varphi_2 (\vec x)
$$
si et seulement si le s\'equent
$$
\varphi_1 \vdash_{\vec x} \varphi_2
$$
est d\'emontrable dans ${\mathbb T}$.

\medskip

\item Pour toute formule ${\mathbb T}$-cart\'esienne $\varphi$ dans un contexte $\vec x = (x_1^{A_1} \cdots x_n^{A_n})$, son interpr\'etation dans le mod\`ele canonique $M_{\mathbb T}$ de ${\mathbb T}$ est le sous-objet
$$
\varphi (\vec x) \quad \mbox{de} \quad \top(\vec x) = \top (x_1^{A_1}) \times \cdots \times \top (x_n^{A_n}) \, .
$$

\item Etant donn\'ees deux formules ${\mathbb T}$-cart\'esiennes $\varphi_1$ et $\varphi_2$ dans un m\^eme contexte $\vec x$, le s\'equent
$$
\varphi_1 \vdash_{\vec x} \varphi_2
$$
est d\'emontrable dans la th\'eorie ${\mathbb T}$ si et seulement si il est v\'erifi\'e par le mod\`ele canonique $M_{\mathbb T}$ de ${\mathbb T}$ dans ${\mathcal C}_{\mathbb T}^{\rm cart}$.
\end{listeimarge}
\end{cor}

\begin{demo}
\begin{listeisansmarge}
\item Si $\varphi_1$ et $\varphi$ sont deux formules ${\mathbb T}$-cart\'esiennes dans un m\^eme contexte $\vec x = (x_1^{A_1} \cdots x_n^{A_n})$ telles que le s\'equent
$$
\varphi_1 \vdash_{\vec x} \varphi
$$
soit d\'emontrable dans la th\'eorie ${\mathbb T}$, alors la formule d\'emontrablement fonctionnelle de $\varphi_1 (\vec x')$ dans $\varphi (\vec x)$, pour $\vec x' = (x'^{A_1}_1 \cdots x'^{A_n}_n)$,
$$
\varphi_1 (\vec x') \wedge x_1^{A_1} = x'^{A_1}_1 \wedge \cdots \wedge x^{A_n}_n = x'^{A_n}_n
$$
est ${\mathbb T}$-cart\'esienne. Elle d\'efinit un monomorphisme de ${\mathcal C}_{\mathbb T}^{\rm cart}$
$$
\varphi_1 (\vec x') \xhookrightarrow{ \ { \ } \ } \varphi (\vec x) \, .
$$
R\'eciproquement, pour tout monomorphisme de ${\mathcal C}_{\mathbb T}^{\rm cart}$
$$
\theta (\vec y , \vec x) : \psi (\vec y) \xhookrightarrow{ \ { \ } \ } \varphi (\vec x) \, ,
$$
la formule induite
$$
(\exists \, \vec y) \ \theta (\vec y , \vec x)
$$
est ${\mathbb T}$-cart\'esienne puisqu'il en est ainsi de $\theta (\vec y , \vec x)$.

\smallskip

En effet, $\theta (\vec y , \vec x)$ est ${\mathbb T}$-d\'emontrablement \'equivalente \`a une formule de la forme
$$
(\exists \, \vec w) \, \theta' (\vec y , \vec x , \vec w)
$$
o\`u $\theta' (\vec y , \vec x , \vec w)$ est une formule de Horn telle que le s\'equent
$$
\theta' (\vec y , \vec x , \vec w) \wedge \theta' (\vec y , \vec x , \vec w') \vdash \vec w = \vec w'
$$
soit ${\mathbb T}$-d\'emontrable.

\smallskip

Comme $\theta (\vec y , \vec x) : \psi (\vec y) \hookrightarrow \varphi (\vec x)$ est un monomorphisme, le s\'equent
$$
(\exists \, \vec w) \, \theta' (\vec y , \vec x , \vec w) \wedge (\exists \, \vec w') \, \theta' ( \vec y' , \vec x , \vec w') \vdash \vec y = \vec y'
$$
est \'egalement ${\mathbb T}$-d\'emontrable.

\smallskip

On en d\'eduit que le s\'equent
$$
\theta' ( \vec y , \vec x , \vec w) \wedge \theta' ( \vec y' , \vec x , \vec w') \vdash \vec y = \vec y' \wedge \vec w = \vec w'
$$
est ${\mathbb T}$-d\'emontrable, et donc la formule
$$
(\exists \, \vec y) \, (\exists \, \vec w) \, \theta' (\vec y , \vec x , \vec w)
$$
est ${\mathbb T}$-cart\'esienne ainsi que la formule \'equivalente
$$
(\exists \, \vec y) \, \theta (\vec y , \vec x) \ .
$$
C'est ce que l'on voulait.

\medskip

\item r\'esulte de (i).

\medskip

\item r\'esulte du lemme \ref{lemV610} puisque le foncteur de plongement
$$
{\mathcal C}_{\mathbb T}^{\rm cart} \xhookrightarrow{ \ { \ } \ } {\mathcal C}_{\mathbb T}^{\rm reg}
$$
transforme le mod\`ele canonique de $M_{\mathbb T}$ de ${\mathbb T}$ dans ${\mathcal C}_{\mathbb T}^{\rm cart}$ en celui de ${\mathbb T}$ dans ${\mathcal C}_{\mathbb T}^{\rm reg}$, et qu'il respecte les limites finies donc les interpr\'etations des formules de Horn et les monomorphismes.

\medskip

\item r\'esulte de (ii) et (iii). 

\end{listeisansmarge}
\end{demo}

\subsection{V\'erification de l'universalit\'e des cat\'egories syntactiques}\label{subsec565}

\medskip

Nous allons maintenant montrer que la cat\'egorie syntactique alg\'ebrique ${\mathcal C}_{\mathbb T}^{\rm alg}$ [resp. r\'eguli\`ere ${\mathcal C}_{\mathbb T}^{\rm reg}$, resp. coh\'erente ${\mathcal C}_{\mathbb T}^{\rm coh}$, resp. g\'eom\'etrique ${\mathcal C}_{\mathbb T}^{\rm geo}$, resp. de Heyting ${\mathcal C}_{\mathbb T}^{\rm He}$, resp. cart\'esienne ${\mathcal C}_{\mathbb T}^{\rm cart}$] d'une th\'eorie alg\'ebrique [resp. r\'eguli\`ere, resp. coh\'erente, resp. g\'eom\'etrique, resp. du premier ordre finitaire, resp. cart\'esienne] ${\mathbb T}$, munie du mod\`ele canonique $M_{\mathbb T}$ de ${\mathbb T}$, poss\`ede la propri\'et\'e de repr\'esentation du foncteur des mod\`eles de ${\mathbb T}$ \'enonc\'ee dans le th\'eor\`eme \ref{thmV61}.

\smallskip

Commen\c cons par le cas des cat\'egories syntactiques alg\'ebriques:

\begin{prop}\label{propV616}

Soit $\Sigma$ une signature sans symbole de relation.

\smallskip

Soit ${\mathbb T}$ une th\'eorie alg\'ebrique de signature $\Sigma$.

\smallskip

Soit ${\mathcal C}_{\mathbb T} = {\mathcal C}_{\mathbb T}^{\rm alg}$ la cat\'egorie syntactique alg\'ebrique de ${\mathbb T}$ munie du mod\`ele canonique $M_{\mathbb T}$ introduit dans la d\'efinition \ref{defV62}.

\smallskip

Alors on a pour toute cat\'egorie alg\'ebrique ${\mathcal C}$:

\begin{listeimarge}

\item Le foncteur 
$$
(F : {\mathcal C}_{\mathbb T} \to {\mathcal C}) \longmapsto F(M_{\mathbb T})
$$
est une \'equivalence de la sous-cat\'egorie pleine de
$$
[{\mathcal C}_{\mathbb T} , {\mathcal C}]
$$
constitu\'ee des foncteurs
$$
F : {\mathcal C}_{\mathbb T} \longrightarrow {\mathcal C}
$$
qui sont alg\'ebriques, sur la cat\'egorie
$$
{\mathbb T}\mbox{\rm -mod} \, ({\mathcal C})
$$
des mod\`eles de ${\mathbb T}$ dans ${\mathcal C}$.

\medskip

\item Il admet pour \'equivalence r\'eciproque le foncteur
$$
M \longmapsto (F_M : {\mathcal C}_{\mathbb T} \to {\mathcal C})
$$
qui associe \`a tout mod\`ele $M$ de ${\mathbb T}$ dans ${\mathcal C}$ le foncteur alg\'ebrique
$$
F_M : {\mathcal C}_{\mathbb T} \longrightarrow {\mathcal C}
$$
qui associe

\medskip

$
\left\{ \begin{matrix}
\bullet &\mbox{\`a toute famille de sortes de $\Sigma$} \hfill \\
{ \ } \\
&A_1 \cdots A_n \\
{ \ } \\
&\mbox{l'objet $M\!A_1 \times \cdots \times M\!A_n$ de ${\mathcal C}$,} \hfill \\
{ \ } \\
\bullet &\mbox{\`a tout morphisme de ${\mathcal C}_{\mathbb T}$} \hfill \\
{ \ } \\
&A_1 \cdots A_n \longrightarrow B_1 \cdots B_m \\
{ \ } \\
&\mbox{d\'efini par une famille de termes} \hfill \\
{ \ } \\
&(f_1 (\vec x) , \cdots , f_m (\vec x)) \\
{ \ } \\
&\mbox{dans un contexte $\vec x = (x_1^{A_1} \cdots x_n^{A_n})$, la famille} \hfill \\
{ \ } \\
&M\!f_1 (\vec x) \times \cdots \times M\!f_m (\vec x) \\
{ \ } \\
&\mbox{des interpr\'etations de ces termes dans le mod\`ele $M$.} \hfill
\end{matrix} \right.
$
\end{listeimarge}
\end{prop}

\begin{demo}

Commen\c cons par v\'erifier que (ii) d\'efinit bien un foncteur alg\'ebrique
$$
F_M : {\mathcal C}_{\mathbb T} \longrightarrow {\mathcal C} \, .
$$

Par construction de ${\mathcal C}_{\mathbb T}$, ses objets sont les familles finies de sortes $A_1 \cdots A_n$ et ses morphismes
$$
A_1 \cdots A_n \longrightarrow B_1 \cdots B_m
$$
sont les classes d'\'equivalence de familles de termes
$$
(f_1 (\vec x) , \cdots , f_m (\vec x))
$$
\`a valeurs dans les sortes $B_1 \cdots B_m$ dans un m\^eme contexte $\vec x = (x_1^{A_1} \cdots x_n^{A_n})$ constitu\'e de variables affect\'ees aux sortes $A_1 \cdots A_n$.

\smallskip

Deux telles familles de termes
$$
(f_1 (\vec x) , \cdots , f_m (\vec x)) \quad \mbox{et} \quad (g_1 (\vec y) , \cdots , g_m (\vec y))
$$
sont \'equivalentes, c'est-\`a-dire d\'efinissent le m\^eme morphisme de ${\mathcal C}_{\mathbb T}$ si, apr\`es substitution des variables $\vec x = (x_1^{A_1} \cdots x_n^{A_n})$ \`a $\vec y = (y_1^{A_1} \cdots y_n^{A_n})$, les s\'equents
$$
\top \vdash_{\vec x} f_1 = g_1 , \cdots , \top \vdash_{\vec x} f_n = g_n
$$
sont d\'emontrables dans la th\'eorie alg\'ebrique ${\mathbb T}$.

\smallskip

S'il en est ainsi, les interpr\'etations
$$
M\!f_1 (\vec x) \times \cdots \times M\!f_m (\vec x) \quad \mbox{et} \quad M\!g_1 (\vec y) \times \cdots \times M\!g_m (\vec x)
$$
sont \'egales en tant que morphismes de ${\mathcal C}$
$$
M\!A_1 \times \cdots \times M\!A_n \longrightarrow M\!B_1 \times \cdots \times M\!B_m \, .
$$

De plus, par d\'efinition de l'interpr\'etation des termes dans le mod\`ele $M$, elle respecte la loi de composition. Ainsi on a bien un foncteur
$$
F_M : {\mathcal C}_{\mathbb T} \longrightarrow {\mathcal C} \, .
$$

Ce foncteur est alg\'ebrique, c'est-\`a-dire respecte les produits puisque le produit de deux objets ${\mathcal C}_{\mathbb T}$
$$
A_1 \cdots A_n \quad \mbox{et} \quad A'_1 \cdots A'_{n'}
$$
est l'objet constitu\'e de la juxtaposition des deux familles
$$
A_1 \cdots A_n \, A'_n \cdots A'_{n'} \, .
$$

Enfin, tout morphisme de mod\`eles de ${\mathbb T}$ dans ${\mathcal C}$
$$
u : M \longrightarrow N
$$
consistant en une famille de morphismes de ${\mathcal C}$
$$
u_A : M\!A \longrightarrow N\!A
$$
qui rend commutatifs les carr\'es
$$
\xymatrix{
M\!A_1 \times \cdots \times M\!A_n \ar[d]_{u_{A_1} \times \cdots \times u_{A_n}} \ar[r]^-{M\!f} &M\!B \ar[d]^{u_B} \\
N\!A_1 \times \cdots \times N\!A_n \ar[r]^-{N\!f} &N\!B
}
$$
associ\'es aux symboles de fonctions $f : A_1 \cdots A_n \to B$ de $\Sigma$, d\'efinit une transformation naturelle
$$
F_M \longrightarrow F_N
$$
qui consiste en la famille des morphismes de ${\mathcal C}$
$$
u_{A_1} \times \cdots \times u_{A_n} : M\!A_1 \times \cdots \times M\!A_n \longrightarrow N\!A_1 \times \cdots \times N\!A_n
$$
index\'es par les familles de sortes $A_1 \cdots A_n$ de $\Sigma$.

\smallskip

Ainsi, on a un foncteur bien d\'efini
$$
M \longmapsto F_M
$$
de la cat\'egorie ${\mathbb T}$-mod \!$({\mathcal C})$ sur la cat\'egorie des foncteurs alg\'ebriques
$$
F : {\mathcal C}_{\mathbb T} \longrightarrow {\mathcal C} \, .
$$

Il reste \`a d\'emontrer que les deux foncteurs
$$
(F : {\mathcal C}_{\mathbb T} \to {\mathcal C}) \longmapsto F(M_{\mathbb T})
$$
et
$$
M \longmapsto F_M
$$
sont deux \'equivalences de cat\'egories r\'eciproques l'une de l'autre, autrement dit que leurs compos\'es dans un sens et dans l'autre sont isomorphes aux foncteurs d'identit\'e de ${\mathbb T}$-mod \!$({\mathcal C})$ et de la cat\'egorie des foncteurs alg\'ebriques ${\mathcal C}_{\mathbb T} \to {\mathcal C}$.

\smallskip

Le mod\`ele $M_{\mathbb T}$ de ${\mathbb T}$ dans ${\mathcal C}_{\mathbb T}$ associe \`a toute sorte $A$ de $\Sigma$ elle-m\^eme consid\'er\'ee comme un objet de ${\mathcal C}_{\mathbb T}$, et \`a tout symbole de fonction $f : A_1 \cdots A_n \to B$ le terme $f(\vec x)$ pour n'importe quel choix de contexte $\vec x = (x_1^{A_1} \cdots x_n^{A_n})$ affect\'e aux sortes $A_1 \cdots A_n$.

\smallskip

Il en r\'esulte que, pour tout mod\`ele $M$ de ${\mathbb T}$ dans ${\mathcal C}$, le mod\`ele
$$
F_M (M_{\mathbb T})
$$
est \'egal \`a $M$ et que, pour tout morphisme de mod\`eles $M \to N$, le morphisme image
$$
F_M (M_{\mathbb T}) \longrightarrow F_N (M_{\mathbb T})
$$
lui est \'egal. Autrement dit, le compos\'e du foncteur
$$
M \longmapsto F_M
$$
suivi du foncteur
$$
(F : {\mathcal C}_{\mathbb T} \to {\mathcal C}) \longmapsto F(M_{\mathbb T})
$$
est \'egal au foncteur d'identit\'e de ${\mathbb T}$-mod \!$({\mathcal C})$.

\smallskip

Dans l'autre sens, partons d'un foncteur alg\'ebrique
$$
F : {\mathcal C}_{\mathbb T} \longrightarrow {\mathcal C} \, .
$$

Le mod\`ele $F(M_{\mathbb T})$ associe \`a toute sorte $A$ de $\Sigma$ l'objet
$$
F(A)
$$
et \`a tout symbole de fonction $f : A_1 \cdots A_n \to B$ le morphisme
$$
F(A_1) \times \cdots \times F(A_n) = F(A_1 \cdots A_n) \xrightarrow{ \ F(f(x_1^{A_1} \cdots x_n^{A_n})) \ } F(B)
$$
pour n'importe quel choix de contexte $\vec x = (x_1^{A_1} \cdots x_n^{A_n})$ affect\'e aux sortes $A_1 \cdots A_n$.

\smallskip

Les identifications
$$
F(A_1) \times \cdots \times F(A_n) = F(A_1 \cdots A_n)
$$
d\'efinissent un isomorphisme canonique du foncteur $F$ sur le foncteur ${\mathcal C}_{\mathbb T} \to {\mathcal C}$ associ\'e au mod\`ele $F(M_{\mathbb T})$ de ${\mathbb T}$ dans ${\mathcal C}$.

\smallskip

De plus, ces isomorphismes
$$
F \xrightarrow{ \ \sim \ } F_{F(M_{\mathbb T})}
$$
respectent les transformations naturelles
$$
F \longrightarrow G
$$
entre foncteurs alg\'ebriques ${\mathcal C}_{\mathbb T} \to {\mathcal C}$.

\smallskip

Cela ach\`eve de montrer la proposition. 

\end{demo}

\bigskip

Traitons de la m\^eme fa\c con le cas des cat\'egories syntactiques cart\'esiennes, r\'eguli\`eres, coh\'erentes, g\'eom\'etriques ou de Heyting:

\begin{prop}\label{propV617}

Soit ${\mathbb T}$ une th\'eorie cart\'esienne [resp. r\'eguli\`ere, resp. coh\'erente, resp. g\'eom\'etrique, resp. du premier ordre finitaire] de signature $\Sigma$.

\smallskip

Soit
$$
{\mathcal C}_{\mathbb T} = {\mathcal C}_{\mathbb T}^{\rm cart} \quad \mbox{[resp.} \ {\mathcal C}_{\mathbb T}^{\rm reg}, \mbox{resp.} \ {\mathcal C}_{\mathbb T}^{\rm coh}, \mbox{resp.} \ {\mathcal C}_{\mathbb T}^{\rm geo}, \mbox{resp.} \ {\mathcal C}_{\mathbb T}^{\rm He} \, \mbox{]}
$$
la cat\'egorie syntactique cart\'esienne [resp. r\'eguli\`ere, resp. coh\'erente, resp. g\'eom\'etrique, resp. de Heyting] de ${\mathbb T}$ munie du mod\`ele canonique $M_{\mathbb T}$ introduite dans la d\'efinition \ref{defV613} [resp. \ref{defV65}].

\smallskip

Alors on a pour toute cat\'egorie cart\'esienne [resp. r\'eguli\`ere, resp. coh\'erente, resp. g\'eom\'etrique, resp. de Heyting] ${\mathcal C}$:

\begin{listeimarge}

\item Le foncteur
$$
(F : {\mathcal C}_{\mathbb T} \to {\mathcal C}) \longmapsto F(M_{\mathbb T})
$$
est une \'equivalence de la sous-cat\'egorie pleine de
$$
[{\mathcal C}_{\mathbb T} , {\mathcal C}]
$$
constitu\'ee des foncteurs
$$
F : {\mathcal C}_{\mathbb T} \longrightarrow {\mathcal C}
$$
qui sont cart\'esiens [resp. r\'eguliers, resp. coh\'erents, resp. g\'eom\'etriques, resp. de Heyting], sur la cat\'egorie
$$
{\mathbb T}\mbox{\rm -mod} ({\mathcal C})
$$
des mod\`eles de ${\mathbb T}$ dans ${\mathcal C}$.

\medskip

\item Il admet pour \'equivalence r\'eciproque le foncteur
$$
M \longmapsto (F_M : {\mathcal C}_{\mathbb T} \to {\mathcal C})
$$
qui associe \`a tout mod\`ele $M$ de ${\mathbb T}$ dans ${\mathcal C}$ le foncteur cart\'esien [resp. r\'egulier, resp. coh\'erent, resp. g\'eom\'etrique, resp. de Heyting]
$$
F_M : {\mathcal C}_{\mathbb T} \longrightarrow {\mathcal C}
$$
qui associe

\medskip

$
\left\{ \begin{matrix}
\bullet &\mbox{\`a toute formule ${\mathbb T}$-cart\'esienne [resp. r\'eguli\`ere, resp. coh\'erente, resp. g\'eom\'etrique,} \hfill \\
&\mbox{resp. du premier ordre finitaire] de $\Sigma$ dans un contexte $\vec x = (x_1^{A_1} \cdots x_n^{A_n})$} \hfill \\
{ \ } \\
&\varphi (\vec x) \\
&\mbox{son interpr\'etation dans le mod\`ele $M$} \hfill \\
&M \varphi (\vec x) \\
{ \ } \\
&\mbox{comme sous-objet de $M\!A_1 \times \cdots \times M\!A_n$ dans ${\mathcal C}$,} \hfill \\
{ \ } \\
\bullet &\mbox{\`a tout morphisme de ${\mathcal C}_{\mathbb T}$} \hfill \\
&\varphi (\vec x) \longrightarrow \psi (\vec y) \\
{ \ } \\
&\mbox{repr\'esent\'e comme une formule ${\mathbb T}$-cart\'esienne [resp. r\'eguli\`ere, resp. coh\'erente, resp. g\'eom\'etrique,} \hfill \\
&\mbox{resp. du premier ordre finitaire] et ${\mathbb T}$-d\'emontrablement fonctionnelle de $\varphi (\vec x)$ dans $\psi (\vec y)$} \hfill \\
{ \ } \\
&\theta (\vec x , \vec y) \\
&\mbox{le morphisme de ${\mathcal C}$} \hfill \\
&M\varphi (\vec x) \longmapsto M\psi (\vec y) \\
{ \ } \\
&\mbox{dont le graphe est l'interpr\'etation dans le mod\`ele $M$} \hfill \\
{ \ } \\
&M\theta (\vec x , \vec y) \\
&\mbox{comme sous-objet de $M\varphi (\vec x) \times M \psi (\vec y)$.} \hfill
\end{matrix} \right.
$
\end{listeimarge}
\end{prop}

\begin{remarks}
\begin{listeisansmarge}
\item Si ${\mathbb T}$ est une th\'eorie alg\'ebrique [resp. cart\'esienne, resp. r\'eguli\`ere, resp. coh\'erente], ses diff\'erentes cat\'egories syntactiques sont reli\'ees par des foncteurs canoniques de plongement
$$
{\mathcal C}_{\mathbb T}^{\rm alg} \xhookrightarrow{ \ { \ } \ } {\mathcal C}_{\mathbb T}^{\rm cart} \xhookrightarrow{ \ { \ } \ } {\mathcal C}_{\mathbb T}^{\rm reg} \xhookrightarrow{ \ { \ } \ } {\mathcal C}_{\mathbb T}^{\rm coh} \xhookrightarrow{ \ { \ } \ } {\mathcal C}_{\mathbb T}^{\rm geo} 
$$
$$
\mbox{[resp.} \quad {\mathcal C}_{\mathbb T}^{\rm cart} \xhookrightarrow{ \ { \ } \ } {\mathcal C}_{\mathbb T}^{\rm reg} \xhookrightarrow{ \ { \ } \ } {\mathcal C}_{\mathbb T}^{\rm coh} \xhookrightarrow{ \ { \ } \ } {\mathcal C}_{\mathbb T}^{\rm geo} \, ,
$$
$$
\mbox{resp.} \quad {\mathcal C}_{\mathbb T}^{\rm reg} \xhookrightarrow{ \ { \ } \ } {\mathcal C}_{\mathbb T}^{\rm coh} \xhookrightarrow{ \ { \ } \ } {\mathcal C}_{\mathbb T}^{\rm geo} \, ,
$$
$$
\mbox{resp.} \quad {\mathcal C}_{\mathbb T}^{\rm coh} \xhookrightarrow{ \ { \ } \ } {\mathcal C}_{\mathbb T}^{\rm geo} \ \mbox{]}
$$
et
$$
{\mathcal C}_{\mathbb T}^{\rm coh} \xhookrightarrow{ \ { \ } \ }  {\mathcal C}_{\mathbb T}^{\rm He} 
$$
qui respectent les mod\`eles canoniques $M_{\mathbb T}$ de ${\mathbb T}$.

\medskip

\item Si ${\mathbb T}$ est une th\'eorie de Horn, c'est a fortiori une th\'eorie cart\'esienne \`a laquelle s'appliquent la proposition et la remarque (i).
\end{listeisansmarge}
\end{remarks}

\bigskip

\begin{demo}

Nous devons d'abord v\'erifier que pour tout mod\`ele $M$ dans ${\mathcal C}$ de la th\'eorie cart\'esienne [resp. r\'eguli\`ere, resp. coh\'erente, resp. g\'eom\'etrique, resp. du premier ordre finitaire] ${\mathbb T}$, (ii) d\'efinit bien un foncteur
$$
F_M : {\mathcal C}_{\mathbb T} \longrightarrow {\mathcal C}
$$
et que ce foncteur est cart\'esien [resp. r\'egulier, resp. coh\'erent, resp. g\'eom\'etrique, resp. de Heyting].

\smallskip

Si $\theta (\vec x , \vec y)$ est une formule ${\mathbb T}$-cart\'esienne [resp. r\'eguli\`ere, resp. coh\'erente, resp. g\'eom\'etrique, resp. de Heyting] et ${\mathbb T}$-d\'emontrablement fonctionnelle d'une formule $\varphi (\vec x)$ de contexte $\vec x = (x_1^{A_1} \cdots x_n^{A_n})$ dans une formule $\psi (\vec y)$ de contexte $\vec y = (y_1^{B_1} \cdots y_m^{B_m})$ disjoint de $\vec x$, l'interpr\'etation
$$
M\theta (\vec x , \vec y)
$$
est d\'efinie a priori comme sous-objet de
$$
M\!A_1 \times \cdots \times M\!A_n \times M\!B_1 \times \cdots \times M\!B_m \, .
$$

Comme c'est une formule ${\mathbb T}$-d\'emontrablement fonctionnelle, son interpr\'etation est m\^eme un sous-objet
$$
M\theta (\vec x , \vec y) \xhookrightarrow{ \ { \ } \ } M\varphi (\vec x) \times M\psi (\vec y)
$$
et la projection
$$
M\theta (\vec x , \vec y) \longrightarrow M\varphi (\vec x)
$$
est un isomorphisme de ${\mathcal C}$.

\smallskip

Autrement dit, $M\theta (\vec x , \vec y)$ est le graphe d'un unique morphisme
$$
M\varphi (\vec x) \longrightarrow M\psi (\vec y) \, .
$$

De plus, \'etant donn\'ees deux telles formules
$$
\varphi (\vec x) \xrightarrow{ \ \theta_1 (\vec x , \vec y) \ } \psi (\vec y) \xrightarrow{ \ \theta_2 (\vec y , \vec z) \ } \chi (\vec z) \, ,
$$
le graphe du morphisme compos\'e des morphismes dont les graphes sont
$$
M\theta_1 (\vec x , \vec y) \xhookrightarrow{ \ { \ } \ } M\varphi (\vec x) \times M\psi (\vec y)
$$
et
$$
M\theta_2 (\vec y , \vec z) \xhookrightarrow{ \ { \ } \ } M\psi (\vec y) \times M\chi (\vec z)
$$
est l'image par le morphisme de projection
$$
M\varphi (\vec x) \times M\psi (\vec y) \times M\chi (\vec z) \longrightarrow M\varphi (\vec x) \times M\chi (\vec z)
$$
de l'intersection des images r\'eciproques des sous-objets
$$
M\theta_1 (\vec x , \vec y) \qquad \mbox{et} \qquad M\theta_2 (\vec y , \vec z)
$$
de
$$
M\varphi (\vec x) \times M\psi (\vec y) \qquad \mbox{et} \qquad M\psi (\vec y) \times M\chi (\vec z) \, .
$$

Ce graphe du morphisme compos\'e est donc l'interpr\'etation dans le mod\`ele $M$ de la formule
$$
(\exists \, \vec y) (\theta_1 (\vec x , \vec y) \wedge \theta_2 (\vec y , \vec z)) \, .
$$

Cela d\'emontre que
$$
F_M : {\mathcal C}_{\mathbb T} \longrightarrow {\mathcal C}
$$
respecte la composition des morphismes.

\smallskip

Pour tout objet de ${\mathcal C}_{\mathbb T}$ repr\'esent\'e par une formule $\varphi (\vec x)$, $F_M$ transforme le morphisme d'identit\'e de $\varphi (\vec x)$
$$
\varphi (\vec x) \xrightarrow{ \ \varphi (\vec x) \wedge \vec x \, = \, \vec x' \ } \varphi (\vec x')
$$
en le morphisme d'identit\'e de $M\varphi (\vec x) = M\varphi (\vec x')$, et donc
$$
F_M : {\mathcal C}_{\mathbb T} \longrightarrow {\mathcal C}
$$
est bien un foncteur.

\smallskip

L'objet terminal de ${\mathcal C}_{\mathbb T}$ est la formule $\top$ dans le contexte vide. Il est transform\'e par $F_M$ en l'objet terminal $1_{\mathcal C}$ de ${\mathcal C}$.

\smallskip

De plus, le produit fibr\'e d'une paire de morphismes de ${\mathcal C}_{\mathbb T}$
$$
\varphi (\vec x) \longrightarrow \psi (\vec y) \longleftarrow \chi (\vec z)
$$
repr\'esent\'es par deux formules d\'emontrablement fonctionnelles
$$
\theta_1 (\vec x , \vec y) \qquad \mbox{et} \qquad \theta_2 (\vec z , \vec y)
$$
est l'objet de ${\mathcal C}_{\mathbb T}$ repr\'esent\'e par la formule
$$
(\exists \, \vec y) (\theta_1 (\vec x , \vec y) \wedge \theta_2 (\vec z , \vec y)) \, .
$$

Or, l'interpr\'etation dans $M$ de cette formule est le sous-objet
$$
M\varphi (\vec x) \times_{M\psi (\vec y)} M\chi (\vec z) \xhookrightarrow{ \ { \ } \ } M\varphi (\vec x) \times M \chi (\vec z) \, .
$$

Cela montre que le foncteur
$$
F_M : {\mathcal C}_{\mathbb T} \longrightarrow {\mathcal C}
$$
respecte les produits fibr\'es, donc les limites finies puisqu'il respecte aussi les objets terminaux.

\smallskip

C'est un foncteur cart\'esien.

\smallskip

Puis montrons que si ${\mathbb T}$ est une th\'eorie r\'eguli\`ere [resp. coh\'erente, resp. g\'eom\'etrique, resp. du premier ordre finitaire], $F_M$ est un foncteur r\'egulier c'est-\`a-dire respecte les images des sous-objets par des morphismes.

\smallskip

Consid\'erons donc un morphisme de ${\mathcal C}_{\mathbb T}$
$$
\theta (\vec x , \vec y) : \varphi (\vec x) \longrightarrow \psi (\vec y)
$$
et un sous-objet de $\varphi (\vec x)$ n\'ecessairement de la forme
$$
\varphi_1 (\vec x) \xhookrightarrow{ \ { \ } \ } \varphi (\vec x)
$$
pour un s\'equent ${\mathbb T}$-d\'emontrable $\varphi_1 \vdash_{\vec x} \varphi$ de $\Sigma$.

\smallskip

Son image par le morphisme $\theta (\vec x , \vec y)$ est d\'efinie par le s\'equent ${\mathbb T}$-d\'emontrable
$$
(\exists \, \vec x) (\varphi_1 (\vec x) \wedge \theta (\vec x , \vec y)) \vdash \psi (\vec y) \, .
$$

Or, l'interpr\'etation dans le mod\`ele $M$ de ${\mathbb T}$ de la formule
$$
(\exists \, \vec x) (\varphi_1 (\vec x) \wedge \theta (\vec x , \vec y))
$$
est l'image du sous-objet
$$
M\varphi_1 (\vec x) \xhookrightarrow{ \ { \ } \ } M\varphi (\vec x)
$$
par le morphisme
$$
M\varphi (\vec x) \longrightarrow M\psi (\vec x)
$$
dont le graphe est
$$
M\theta (\vec x , \vec y) \xhookrightarrow{ \ { \ } \ } M\varphi (\vec x) \times M\psi (\vec y) \, .
$$

Ainsi, le foncteur $F_M$ est bien r\'egulier.

\smallskip

Notons encore que si ${\mathbb T}$ est une th\'eorie coh\'erente ou finitaire du premier ordre [resp. g\'eom\'etrique], alors le foncteur $F_M$ respecte les r\'eunions finies [resp. arbitraires] de sous-objets donc est coh\'erent [resp. g\'eom\'etrique]. Cela r\'esulte de ce que, pour toute famille finie [resp. arbitraire] de sous-objets d'un objet $\varphi (\vec x)$ de ${\mathcal C}_{\mathbb T}$ repr\'esent\'es par des formules
$$
\varphi_i (\vec x) \, , \qquad 1 \leq i \leq k \quad \mbox{[resp.} \quad i \in I \ \mbox{]}
$$
pour des s\'equents ${\mathbb T}$-d\'emontrables $\varphi_i \vdash_{\vec x} \varphi$, la formule
$$
\varphi_1 (\vec x) \vee \cdots \vee \varphi_k (\vec x) \qquad \mbox{[resp.} \quad \bigvee_{i \in I} \varphi_i (\vec x) \ \mbox{]}
$$
s'interpr\`ete dans le mod\`ele $M$ de ${\mathbb T}$ comme la r\'eunion des sous-objets
$$
M\varphi_i (\vec x) \xhookrightarrow{ \ { \ } \ } M\varphi (\vec x) \, , \qquad 1 \leq i \leq k \quad \mbox{[resp.} \quad i \in I \ \mbox{]}.
$$

Enfin, si ${\mathbb T}$ est une th\'eorie du premier ordre finitaire, le foncteur
$$
F_M : {\mathcal C}_{\mathbb T} \longrightarrow {\mathcal C}
$$
respecte les foncteurs $\forall , \Rightarrow$ et $\neg$.

\smallskip

En effet, pour tout morphisme de ${\mathcal C}_{\mathbb T}$
$$
\theta = \theta (\vec x , \vec y) : \varphi (\vec x) \longrightarrow \psi (\vec y) \, ,
$$
le foncteur associ\'e $\forall_{\theta}$ sur les sous-objets de $\varphi (\vec x)$ transforme toute formule
$$
\varphi_1 (\vec x) \quad \mbox{telle que} \quad \varphi_1 \vdash_{\vec x} \varphi \quad \mbox{soit ${\mathbb T}$-d\'emontrable, }
$$
en la formule de contexte $\vec y$
$$
(\forall \, \vec x) \, (\theta (\vec x , \vec y) \Rightarrow \varphi_1 (\vec x)) \, .
$$

Or, celle-ci s'interpr\`ete dans le mod\`ele $M$ comme le sous-objet
$$
\forall_{M\theta} \, M\varphi_1 (\vec x)  \xhookrightarrow{ \ { \ } \ } M\psi (\vec y)
$$
image du sous-objet
$$
M\varphi_1 (\vec x) \xhookrightarrow{ \ { \ } \ } M\varphi (\vec x)
$$
par le foncteur
$$
\forall_{M\theta}
$$
associ\'e au morphisme
$$
M\theta : M\varphi (\vec x) \longrightarrow M\psi (\vec y)
$$
dont le graphe est
$$
M\theta (\vec x , \vec y) \xhookrightarrow{ \ { \ } \ } M\varphi (\vec x) \times M\psi (\vec y) \, .
$$

De la m\^eme fa\c con, pour tous sous-objets dans ${\mathcal C}_{\mathbb T}$
$$
\varphi_1 (\vec x) \xhookrightarrow{ \ { \ } \ } \varphi (\vec x) \qquad \mbox{et} \qquad \varphi_2 (\vec x) \xhookrightarrow{ \ { \ } \ } \varphi (\vec x) \, ,
$$
les formules
$$
(\varphi_1 (\vec x) \Rightarrow \varphi_2 (\vec x)) \wedge \varphi (\vec x)
$$
ou
$$
\neg \, \varphi_1 (\vec x) \wedge \varphi (\vec x)
$$
s'interpr\`etent dans le mod\`ele $M$ comme les images par les foncteurs $\Rightarrow$ ou $\neg$ op\'erant sur les sous-objets de $M\varphi (\vec x)$ dans ${\mathcal C}$ des sous-objets
$$
(M\varphi_1 (\vec x) , M\varphi_2 (\vec x))
$$
ou
$$
M\varphi_1 (\vec x) \, .
$$

Cela ach\`eve de montrer que, si ${\mathbb T}$ est une th\'eorie du premier ordre finitaire, le foncteur
$$
F_M : {\mathcal C}_{\mathbb T} \longrightarrow {\mathcal C}
$$
associ\'e \`a tout mod\`ele $M$ de ${\mathbb T}$ dans ${\mathcal C}$ est un foncteur de Heyting.

\smallskip

Il est clair d'autre part que, si ${\mathbb T}$ est une th\'eorie cart\'esienne [resp. r\'eguli\`ere, resp. coh\'erente, resp. g\'eom\'etrique, resp. du premier ordre finitaire], tout morphisme de mod\`eles de ${\mathbb T}$ dans ${\mathcal C}$
$$
M \longrightarrow N
$$
induit une famille de morphismes de ${\mathcal C}$
$$
M\varphi (\vec x) \longrightarrow N\varphi (\vec x)
$$
qui relient les interpr\'etations dans $M$ et $N$ des formules $\varphi (\vec x)$ repr\'esentant les objets de ${\mathcal C}_{\mathbb T}$. Cette famille de morphismes constitue une transformation naturelle
$$
F_M \longrightarrow F_N \, .
$$

Cette construction respecte la composition des morphismes, ce qui signifie que
$$
M \longmapsto F_M
$$
est bien d\'efini comme foncteur de ${\mathbb T}$-mod \!$({\mathcal C})$ dans la cat\'egorie des foncteurs cart\'esiens [resp. r\'eguliers, resp. coh\'erents, resp. g\'eom\'etriques, resp. de Heyting] ${\mathcal C}_{\mathbb T} \to {\mathcal C}$.

\smallskip

Il reste \`a d\'emontrer que les deux foncteurs
$$
\begin{matrix}
&(F : {\mathcal C}_{\mathbb T} \to {\mathcal C}) &\longmapsto &F(M_{\mathbb T}) \hfill \\
\mbox{et} &\hfill M &\longmapsto &(F_M : {\mathcal C}_{\mathbb T} \to {\mathcal C})
\end{matrix}
$$
sont deux \'equivalences r\'eciproques l'une de l'autre.

\smallskip

Par d\'efinition, le mod\`ele canonique $M_{\mathbb T}$ de ${\mathbb T}$ dans ${\mathcal C}_{\mathbb T}$ associe

\medskip

$
\left\{ \begin{matrix}
\bullet &\mbox{\`a toute sorte $A$ de $\Sigma$ l'objet repr\'esent\'e par} \hfill \\
{ \ } \\
&\top (x^A) \\
{ \ } \\
&\mbox{pour n'importe quelle variable $x^A$ affect\'ee \`a la sorte $A$,} \hfill \\
{ \ } \\
\bullet &\mbox{\`a tout symbole de fonction $f : A_1 \cdots A_n \to B$ le morphisme repr\'esent\'e par} \hfill \\
{ \ } \\
&\top (x_1^{A_1}) \times \cdots \times \top (x_n^{A_n}) = \top (x_1^{A_1} \cdots x_n^{A_n}) \xrightarrow{ \ y^B = f (x_1^{A_1} \cdots x_n^{A_n}) \ } \top (y^B) \\
{ \ } \\
&\mbox{pour des variables $x_1^{A_1} \cdots x_n^{A_n}$ et $y^B$ affect\'ees aux sortes $A_1 \cdots A_n$ et $B$,} \hfill \\
{ \ } \\
\bullet &\mbox{\`a tout symbole de relation $R \rightarrowtail A_1 \cdots A_n$ le sous-objet} \hfill \\
{ \ } \\
&R(x_1^{A_1} \cdots x_n^{A_n}) \xhookrightarrow{ \ { \ } \ } \top (x_1^{A_1} \cdots x_n^{A_n})  = \top (x_1^{A_1} )\times \cdots \times \top (x_n^{A_n}) \, .
\end{matrix} \right.
$

\bigskip

Il en r\'esulte que le compos\'e
$$
M \longmapsto F_M (M_{\mathbb T})
$$
du foncteur $M \mapsto F_M$ suivi du foncteur $F \mapsto F(M_{\mathbb T})$ n'est autre que le foncteur d'identit\'e de la cat\'egorie ${\mathbb T}$-mod \!$({\mathcal C})$
$$
M \longmapsto M \, .
$$

Dans l'autre sens, partons d'un foncteur cart\'esien [resp. r\'egulier, resp. coh\'erent, resp. g\'eom\'etrique, resp. de Heyting]
$$
F : {\mathcal C}_{\mathbb T} \longrightarrow {\mathcal C} \, .
$$

Le mod\`ele $F(M_{\mathbb T})$ de ${\mathbb T}$ dans ${\mathcal C}$ associe

\medskip

$
\left\{ \begin{matrix}
\bullet &\mbox{\`a toute sorte $A$ de $\Sigma$ l'objet} \hfill \\
{ \ } \\
&F(\top (x^A)) \, , \\
{ \ } \\
\bullet &\mbox{\`a tout symbole de fonction $f : A_1 \cdots A_n \to B$ le morphisme} \hfill \\
{ \ } \\
&F(\top (x_1^{A_1})) \times \cdots \times F(\top (x_n^{A_n})) = F(\top (x_1^{A_1} \cdots x_n^{A_n})) \longrightarrow F(\top (y^B)) \\
{ \ } \\
&\mbox{dont le graphe est le sous-objet} \hfill \\
{ \ } \\
&F(y^B = f(x_1^{A_1} \cdots x_n^{A_n})) \, , \\
{ \ } \\
\bullet &\mbox{\`a tout symbole de relation $R \rightarrowtail A_1 \cdots A_n$ le sous-objet} \hfill \\
{ \ } \\
&F(R(x_1^{A_1} \cdots x_n^{A_n})) \xhookrightarrow{ \ { \ } \ } F(\top (x_1^{A_1} \cdots x_n^{A_n})) \, .
\end{matrix} \right.
$

\bigskip

On en d\'eduit que toute formule ${\mathbb T}$-cart\'esienne [resp. r\'eguli\`ere, resp. coh\'erente, resp. g\'eom\'etrique, resp. du premier ordre finitaire] de $\Sigma$ dans un contexte $\vec x = (x_1^{A_1} \cdots x_n^{A_n})$
$$
\varphi (\vec x)
$$
s'interpr\`ete dans le mod\`ele $F(M_{\mathbb T})$ comme le sous-objet
$$
F(\varphi (x_1^{A_1} \cdots x_n^{A_n})) \xhookrightarrow{ \ { \ } \ } F(\top (x_1^{A_1} \cdots x_n^{A_n})) =  F(\top (x_1^{A_1})) \times \cdots \times F(\top (x_n^{A_n}))
$$
image par $F$ du sous-objet dans ${\mathcal C}_{\mathbb T}$
$$
\varphi (x_1^{A_1} \cdots x_n^{A_n}) \xhookrightarrow{ \ { \ } \ } \top (x_1^{A_1} \cdots x_n^{A_n}) = \top (x_1^{A_1}) \times \cdots \times \top (x_n^{A_n}) \, .
$$

Il en r\'esulte que le foncteur
$$
F_M : {\mathcal C}_{\mathbb T} \longrightarrow {\mathcal C} 
$$
associ\'e au mod\`ele $M = F(M_{\mathbb T})$ de ${\mathbb T}$ dans ${\mathcal C}$ est canoniquement isomorphe au foncteur $F$.

\smallskip

De plus, les isomorphismes canoniques
$$
F_{F(M_{\mathbb T})} \xrightarrow{ \ \sim \ } F
$$
sont compatibles avec les transformations naturelles
$$
F \longrightarrow G \, .
$$

Ainsi, le compos\'e du foncteur
$$
F \longmapsto F(M_{\mathbb T})
$$
suivi du foncteur
$$
M \longmapsto F_M
$$
est isomorphe au foncteur identit\'e de la cat\'egorie des foncteurs cart\'esiens [resp. r\'eguliers, resp. coh\'erents, resp. g\'eom\'etriques, resp. de Heyting] ${\mathcal C}_{\mathbb T} \to {\mathcal C}$.

\smallskip

Cela ach\`eve de montrer la proposition. 

\end{demo}

\subsection{La notion d'\'equivalence syntactique entre cat\'egories}\label{subsec566}

Il est naturel de poser:

\begin{defn}\label{defV618}

Deux th\'eories alg\'ebriques [resp. cart\'esiennes, resp. r\'eguli\`eres, resp. coh\'erentes, resp. g\'eom\'etriques, resp. du premier ordre finitaires] ${\mathbb T}_1$ et ${\mathbb T}_2$ dans des signatures $\Sigma_1$ et $\Sigma_2$ sont dites syntactiquement \'equivalentes si leurs cat\'egories syntactiques

\medskip

$
\begin{matrix}
&{\mathcal C}_{{\mathbb T}_1}^{\rm alg} &\mbox{et} &{\mathcal C}_{{\mathbb T}_2}^{\rm alg} \\
{ \ } \\
\mbox{[resp.} &{\mathcal C}_{{\mathbb T}_1}^{\rm cart} &\mbox{et} &{\mathcal C}_{{\mathbb T}_2}^{\rm cart} \, , \\
{ \ } \\
\mbox{resp.} &{\mathcal C}_{{\mathbb T}_1}^{\rm reg} &\mbox{et} &{\mathcal C}_{{\mathbb T}_2}^{\rm reg} \, , \\
{ \ } \\
\mbox{resp.} &{\mathcal C}_{{\mathbb T}_1}^{\rm coh} &\mbox{et} &{\mathcal C}_{{\mathbb T}_2}^{\rm coh} \, , \\
{ \ } \\
\mbox{resp.} &{\mathcal C}_{{\mathbb T}_1}^{\rm geo} &\mbox{et} &{\mathcal C}_{{\mathbb T}_2}^{\rm geo} \, , \\
{ \ } \\
\mbox{resp.} &{\mathcal C}_{{\mathbb T}_1}^{\rm He} &\mbox{et} &{\mathcal C}_{{\mathbb T}_2}^{\rm He} \ \mbox{]} 
\end{matrix}
$

\bigskip

\noindent sont \'equivalentes.
\end{defn}

\begin{remarksqed}
\begin{listeisansmarge}
\item Si ${\mathbb T}_1$ et ${\mathbb T}_2$ ont m\^eme signature $\Sigma_1 = \Sigma_2$ et que leurs cat\'egories syntactiques sont reli\'ees par des \'equivalences qui respectent leurs mod\`eles canoniques \`a isomorphismes pr\`es en tant que $\Sigma_1$-structures, alors ${\mathbb T}_1$ et ${\mathbb T}_2$ sont \'equivalentes au sens de la d\'efinition \ref{defV56}. Cela impose l'\'egalit\'e des cat\'egories syntactiques de ${\mathbb T}_1$ et ${\mathbb T}_2$.

\medskip

\item La notion d'\'equivalence syntactique est beaucoup plus g\'en\'erale et plus souple que la notion d'\'equivalence de th\'eories de m\^eme signature. 
\end{listeisansmarge}
\end{remarksqed}

\bigskip

On d\'eduit des propositions \ref{propV616} et \ref{propV617}:

\begin{cor}\label{corV619}

Soient ${\mathbb T}_1$ et ${\mathbb T}_2$ deux th\'eories alg\'ebriques [resp. cart\'esiennes, resp. r\'eguli\`eres, resp. coh\'erentes, resp. g\'eom\'etriques, resp. du premier ordre finitaires] dans des signatures $\Sigma_1$ et $\Sigma_2$.

\smallskip

Alors ces th\'eories sont syntactiquement \'equivalentes si et seulement si on peut associer \`a toute cat\'egorie alg\'ebrique [resp. r\'eguli\`ere, resp. coh\'erente, resp. g\'eom\'etrique, resp. de Heyting] ${\mathcal C}$ deux foncteurs
$$
F_{\mathcal C} : {\mathbb T}_1\mbox{\rm -mod} \, ({\mathcal C}) \longrightarrow {\mathbb T}_2\mbox{\rm -mod} \, ({\mathcal C}) \, ,
$$
$$
G_{\mathcal C} : {\mathbb T}_2\mbox{\rm -mod} \, ({\mathcal C}) \longrightarrow {\mathbb T}_1\mbox{\rm -mod} \, ({\mathcal C})
$$
tels que:

\begin{listeimarge}

\item[$\bullet$] ces foncteurs sont naturels au sens que, pour tout foncteur alg\'ebrique [resp. cart\'esien, resp. r\'egulier, resp. coh\'erent, resp. g\'eom\'etrique, resp. de Heyting]
$$
{\mathcal C} \longrightarrow {\mathcal D} \, ,
$$
les carr\'es induits
$$
\xymatrix{
{\mathbb T}_1\mbox{\rm -mod} \, ({\mathcal C}) \ar[d] \ar[r]^{F_{\mathcal C}} &{\mathbb T}_2\mbox{\rm -mod} \, ({\mathcal C}) \ar[d] \\
{\mathbb T}_1\mbox{\rm -mod} \, ({\mathcal D}) \ar[r]^{F_{\mathcal D}} &{\mathbb T}_2\mbox{\rm -mod} \, ({\mathcal D})
} \qquad \xymatrix{
{\mathbb T}_2\mbox{\rm -mod} \, ({\mathcal C}) \ar[d] \ar[r]^{G_{\mathcal C}} &{\mathbb T}_1\mbox{\rm -mod} \, ({\mathcal C}) \ar[d] \\
{\mathbb T}_2\mbox{\rm -mod} \, ({\mathcal D}) \ar[r]^{G_{\mathcal D}} &{\mathbb T}_1\mbox{\rm -mod} \, ({\mathcal D})
} 
$$
sont commutatifs \`a isomorphisme pr\`es,

\medskip

\item[$\bullet$] ces foncteurs
$$
F_{\mathcal C} \qquad \mbox{et} \qquad G_{\mathcal C}
$$
sont des \'equivalences de cat\'egories r\'eciproques l'une de l'autre pour toute telle cat\'egorie ${\mathcal C}$.
\end{listeimarge}
\end{cor}

\bigskip

\begin{remarks}
\begin{listeisansmarge}
\item Autrement dit, deux telles th\'eories sont syntactiquement \'equivalentes si leurs foncteurs des mod\`eles sont \'equivalents.

\medskip

\item Il r\'esulte de la caract\'erisation donn\'ee que
\begin{enumerate}
\item[$\bullet$] si ${\mathbb T}_1$ et ${\mathbb T}_2$ sont des th\'eories alg\'ebriques syntactiquement \'equivalentes, elles le sont aussi en tant que th\'eories cart\'esiennes,
\item[$\bullet$] si ${\mathbb T}_1$ et ${\mathbb T}_2$ sont des th\'eories cart\'esiennes [resp. r\'eguli\`eres, resp. coh\'erentes] syntactiquement \'equivalentes, elles le sont aussi en tant que th\'eories r\'eguli\`eres [resp. coh\'erentes, resp. g\'eom\'etriques ou du premier ordre finitaires].
\end{enumerate}
\end{listeisansmarge}
\end{remarks}

\bigskip

\begin{demo}

C'est une cons\'equence formelle de ce que, d'apr\`es les propositions \ref{propV616} et \ref{propV617}, les foncteurs des mod\`eles
$$
{\mathcal C} \longmapsto {\mathbb T}_1\mbox{\rm -mod} \, ({\mathcal C}) \qquad \mbox{et} \qquad {\mathcal C} \longmapsto {\mathbb T}_2\mbox{\rm -mod} \, ({\mathcal C})
$$
sont repr\'esent\'es par les cat\'egories syntactiques de ${\mathbb T}_1$ et ${\mathbb T}_2$, lesquelles sont caract\'eris\'ees \`a \'equivalence pr\`es par cette propri\'et\'e. \end{demo}

\section{Topos classifiants}\label{sec57}

\subsection{Repr\'esentation des foncteurs des mod\`eles par les topos}\label{subsec571}

\medskip

Consid\'erons une th\'eorie ${\mathbb T}$ qui est suppos\'ee g\'eom\'etrique. Ce peut \^etre en particulier une th\'eorie alg\'ebrique, ou plus g\'en\'eralement de Horn, ou plus g\'en\'eralement encore cart\'esienne, ou r\'eguli\`ere, ou coh\'erente.

\smallskip

Tout topos ${\mathcal E}$ est une cat\'egorie g\'eom\'etrique, et donc la th\'eorie g\'eom\'etrique ${\mathbb T}$ d\'efinit la cat\'egorie localement petite
$$
{\mathbb T}\mbox{-mod} \, ({\mathcal E})
$$
des mod\`eles de ${\mathbb T}$ dans ${\mathcal E}$.

\smallskip 

De plus, pour tout morphisme de topos
$$
f = (f^* , f_*) : {\mathcal E}' \longrightarrow {\mathcal E} \, ,
$$
sa composante d'image r\'eciproque
$$
f^* : {\mathcal E} \longrightarrow {\mathcal E}'
$$
est un foncteur g\'eom\'etrique. Elle induit par cons\'equent un foncteur
$$
f^* : {\mathbb T}\mbox{-mod} \, ({\mathcal E}) \longrightarrow {\mathbb T}\mbox{-mod} \, ({\mathcal E}')
$$
entre les cat\'egories de mod\`eles de ${\mathbb T}$.

\smallskip

Enfin, pour toute paire de morphismes $f = (f^* , f_*)$ et $g=(g^*,g_*)$ d'un topos ${\mathcal E}'$ dans un topos ${\mathcal E}$, tout morphisme $f \to g$ c'est-\`a-dire toute transformation naturelle
$$
f^* \longrightarrow g^*
$$
induit une transformation naturelle du foncteur
$$
f^* : {\mathbb T}\mbox{-mod} \, ({\mathcal E}) \longrightarrow {\mathbb T}\mbox{-mod} \, ({\mathcal E}')
$$
dans le foncteur
$$
g^* : {\mathbb T}\mbox{-mod} \, ({\mathcal E}) \longrightarrow {\mathbb T}\mbox{-mod} \, ({\mathcal E}') \, .
$$

Nous allons montrer que le 2-foncteur ainsi d\'efini, appel\'e le foncteur des mod\`eles topossiques de ${\mathbb T}$, est repr\'esentable au sens suivant:

\begin{thm}[Makkai et Reyes (\cite{CategoricalLogic})]\label{thmV71}

Consid\'erons comme ci-dessus une th\'eorie ${\mathbb T}$ suppos\'ee g\'eom\'etrique.

\smallskip

Alors il existe un topos ${\mathcal E}_{\mathbb T}$, muni d'un mod\`ele $U_{\mathbb T}$ de ${\mathbb T}$ dans ${\mathcal E}_{\mathbb T}$, tel que, pour tout topos ${\mathcal E}$, le foncteur
$$
(f : {\mathcal E} \longrightarrow {\mathcal E}_{\mathbb T}) \longmapsto f^*  U_{\mathbb T}
$$
est une \'equivalence de la cat\'egorie des morphismes de topos
$$
f = (f^* , f_*) : {\mathcal E} \longrightarrow {\mathcal E}_{\mathbb T}
$$
sur la cat\'egorie 
$$
{\mathbb T}\mbox{\rm -mod} \, ({\mathcal E})
$$
des mod\`eles de ${\mathbb T}$ dans ${\mathcal E}$.

\smallskip

De plus, le topos ${\mathcal E}_{\mathbb T}$, qui est appel\'e le topos classifiant de la th\'eorie g\'eom\'etrique ${\mathbb T}$, est caract\'eris\'e \`a \'equivalence pr\`es par cette propri\'et\'e.
\end{thm}

\begin{commdemo}

Prouvons d'abord l'unicit\'e \`a \'equivalence pr\`es des topos ${\mathcal E}_{\mathbb T}$ v\'erifiant la propri\'et\'e de l'\'enonc\'e.

\smallskip

Elle se d\'emontre comme la partie d'unicit\'e du th\'eor\`eme \ref{thmV61}.

\smallskip

Si ${\mathcal E}_{\mathbb T}$ et ${\mathcal E}'_{\mathbb T}$ sont deux topos munis de mod\`eles $U_{\mathbb T}$ et $U'_{\mathbb T}$ de ${\mathbb T}$ qui satisfont chacun la propri\'et\'e de l'\'enonc\'e, il existe deux morphismes de topos
$$
f = (f^* , f_*) : {\mathcal E}'_{\mathbb T} \longrightarrow {\mathcal E}_{\mathbb T}
$$
et
$$
g = (g^* , g_*) : {\mathcal E}_{\mathbb T} \longrightarrow {\mathcal E}'_{\mathbb T}
$$
tels que les ${\mathbb T}$-mod\`eles $f^* U_{\mathbb T}$ et $g^* U'_{\mathbb T}$ soient respectivement isomorphes \`a $U'_{\mathbb T}$ et \`a $U_{\mathbb T}$.

\smallskip

Alors les isomorphismes de ${\mathbb T}$-mod\`eles
$$
g^* \circ f^* U_{\mathbb T} \xrightarrow{ \ \sim \ } U_{\mathbb T}
$$
et
$$
f^* \circ g^* U'_{\mathbb T} \xrightarrow{ \ \sim \ } U'_{\mathbb T}
$$
proviennent d'isomorphismes de foncteurs
$$
g^* \circ f^* \xrightarrow{ \ \sim \ } {\rm id}_{{\mathcal E}_{\mathbb T}}
$$
et
$$
f^* \circ g^* \xrightarrow{ \ \sim \ } {\rm id}_{{\mathcal E}'_{\mathbb T}} \, .
$$

Ainsi, $f$ et $g$ sont deux \'equivalences r\'eciproques l'une de l'autre entre les topos ${\mathcal E}_{\mathbb T}$ et ${\mathcal E}'_{\mathbb T}$.

\smallskip

Cela prouve l'unicit\'e \`a \'equivalence pr\`es.

\smallskip

Pour l'existence, nous construirons un peu plus loin le topos ${\mathcal E}_{\mathbb T}$ comme associ\'e \`a un site $({\mathcal C}_{\mathbb T} , J_{\mathbb T})$ constitu\'e d'une cat\'egorie syntactique ${\mathcal C}_{\mathbb T}$ de ${\mathbb T}$ munie d'une certaine topologie $J_{\mathbb T}$ appel\'ee la topologie syntactique.

\smallskip

Nous montrerons alors que le topos ainsi construit ${\mathcal E}_{\mathbb T} = \widehat{({\mathcal C}_{\mathbb T})}_{J_{\mathbb T}}$ v\'erifie bien la propri\'et\'e du th\'eor\`eme. 

\end{commdemo}

\subsection{Topologies syntactiques}\label{subsec572}

\medskip

On d\'efinit sur les cat\'egories syntactiques g\'eom\'etriques [resp. coh\'erentes, resp. r\'eguli\`eres, resp. cart\'esiennes] des th\'eories g\'eom\'etriques [resp. coh\'erentes, resp. r\'eguli\`eres, resp. cart\'esiennes] des topologies appel\'ees leurs topologies syntactiques:

\begin{defn}\label{defV72}

Soit ${\mathbb T}$ une th\'eorie g\'eom\'etrique du premier ordre dans une signature $\Sigma$.

\smallskip

Alors:

\begin{listeimarge}

\item On munit la cat\'egorie syntactique g\'eom\'etrique de ${\mathbb T}$
$$
{\mathcal C}_{\mathbb T} = {\mathcal C}^{\rm geo}_{\mathbb T}
$$
de la topologie de Grothendieck $J_{\mathbb T}$ pour laquelle un crible sur un objet $\varphi (\vec x)$ est couvrant s'il contient une famille de morphismes
$$
\theta_i (\vec x_i , \vec x) : \varphi_i (\vec x_i) \longrightarrow \varphi (\vec x)
$$
dont la r\'eunion des images est \'egale \`a tout $\varphi (\vec x)$, c'est-\`a-dire tel que le s\'equent
$$
\varphi \vdash_{\vec x} \bigvee_{i \in I} (\exists \, \vec x_i) \, \theta_i
$$
soit d\'emontrable dans la th\'eorie ${\mathbb T}$ suivant les r\`egles d'inf\'erence de la logique g\'eom\'etrique.

\medskip

\item Si ${\mathbb T}$ est une th\'eorie coh\'erente, on munit sa cat\'egorie syntactique coh\'erente
$$
{\mathcal C}_{\mathbb T} = {\mathcal C}^{\rm coh}_{\mathbb T}
$$
de la topologie $J_{\mathbb T}$ pour laquelle un crible sur un objet $\varphi (\vec x)$ est couvrant s'il contient une famille finie de morphismes
$$
\theta_i (\vec x_i , \vec x) : \varphi_i (\vec x_i) \longrightarrow \varphi (\vec x) \, , \quad 1 \leq i \leq k \, ,
$$
dont la r\'eunion des images est \'egale \`a tout $\varphi (\vec x)$, c'est-\`a-dire tels que le s\'equent
$$
\varphi \vdash_{\vec x} (\exists \, \vec x_1) \, \theta_1 \vee \cdots \vee (\exists \, \vec x_k) \, \theta_k
$$
soit d\'emontrable dans la th\'eorie ${\mathbb T}$ suivant les r\`egles d'inf\'erence de la logique coh\'erente.

\medskip

\item Si ${\mathbb T}$ est une th\'eorie r\'eguli\`ere, on munit sa cat\'egorie syntactique r\'eguli\`ere
$$
{\mathcal C}_{\mathbb T} = {\mathcal C}^{\rm reg}_{\mathbb T}
$$
de la topologie $J_{\mathbb T}$ pour laquelle un crible sur un objet $\varphi (\vec x)$ est couvrant s'il contient  un morphisme
$$
\theta (\vec x' , \vec x) : \varphi' (\vec x') \longrightarrow \varphi (\vec x) 
$$
dont l'image est \'egale \`a tout $\varphi (\vec x)$, c'est-\`a-dire tel que le s\'equent
$$
\varphi \vdash_{\vec x} (\exists \, \vec x') \, \theta
$$
soit d\'emontrable dans la th\'eorie ${\mathbb T}$ suivant les r\`egles d'inf\'erence de la logique r\'eguli\`ere.

\medskip

\item Enfin, si ${\mathbb T}$ est une th\'eorie cart\'esienne, on munit sa cat\'egorie syntactique cart\'esienne
$$
{\mathcal C}_{\mathbb T} = {\mathcal C}^{\rm cart}_{\mathbb T}
$$
de la topologie discr\`ete $J_{\mathbb T}$ pour laquelle les seuls cribles couvrants sont les cribles maximaux.
\end{listeimarge}
\end{defn}

\begin{remarksqed}
\begin{listeisansmarge}
\item Si une th\'eorie g\'eom\'etrique ${\mathbb T}$ est coh\'erente [resp. r\'eguli\`ere, resp. cart\'esienne], elle poss\`ede non seulement un site syntactique g\'eom\'etrique
$$
({\mathcal C}^{\rm geo}_{\mathbb T} , J_{\mathbb T})
$$
mais aussi un site syntactique coh\'erent $({\mathcal C}^{\rm coh}_{\mathbb T} , J_{\mathbb T})$ [resp. et un site syntactique r\'egulier $({\mathcal C}^{\rm reg}_{\mathbb T} , J_{\mathbb T})$, resp. et un site syntactique cart\'esien $({\mathcal C}^{\rm cart}_{\mathbb T} , J_{\mathbb T})$].

\medskip

\item La topologie syntactique $J_{\mathbb T}$ de la cat\'egorie syntactique ${\mathcal C}_{\mathbb T} = {\mathcal C}^{\rm geo}_{\mathbb T}$ [resp. ${\mathcal C}^{\rm coh}_{\mathbb T}$, resp. ${\mathcal C}^{\rm reg}_{\mathbb T}$, resp. ${\mathcal C}^{\rm cart}_{\mathbb T}$] est d\'efinie par sa structure cat\'egorique.

\smallskip

En particulier, si deux th\'eories ${\mathbb T}_1$ et ${\mathbb T}_2$ sont syntactiquement \'equivalentes au sens que leurs cat\'egories syntactiques sont reli\'ees par des \'equivalences de cat\'egories
$$
{\mathcal C}^{\rm geo}_{{\mathbb T}_1} \cong {\mathcal C}^{\rm geo}_{{\mathbb T}_2} \, ,
$$
ou
$$
{\mathcal C}^{\rm coh}_{{\mathbb T}_1} \cong {\mathcal C}^{\rm coh}_{{\mathbb T}_2} \, ,
$$
ou
$$
{\mathcal C}^{\rm reg}_{{\mathbb T}_1} \cong {\mathcal C}^{\rm reg}_{{\mathbb T}_2} \, ,
$$
ou
$$
{\mathcal C}^{\rm cart}_{{\mathbb T}_1} \cong {\mathcal C}^{\rm cart}_{{\mathbb T}_2} \, ,
$$
alors ces \'equivalences respectent leurs topologies syntactiques. 
\end{listeisansmarge}
\end{remarksqed}

\bigskip

V\'erifions que la d\'efinition \ref{defV72} ci-dessus est valide:

\begin{prop}\label{propV73}

Soit ${\mathcal C}$ une cat\'egorie r\'eguli\`ere [resp. coh\'erente, resp. g\'eom\'etrique] essentiellement petite.

\smallskip

Alors:

\begin{listeimarge}

\item Il existe une topologie de Grothendieck $J$ de ${\mathcal C}$ pour laquelle un crible sur un objet $X$ est couvrant s'il contient un morphisme
$$
X' \longrightarrow X
$$
[resp. une famille finie de morphismes, resp. une famille de morphismes
$$
X_i \longrightarrow X \ \mbox{]}
$$
dont l'image [resp. dont la r\'eunion des images] est \'egale \`a $X$.

\medskip

\item Cette topologie est sous-canonique.
\end{listeimarge}
\end{prop}

\begin{remark}

On peut montrer que toute cat\'egorie r\'eguli\`ere [resp. coh\'erente, resp. g\'eom\'etrique] essentiellement petite ${\mathcal C}$ est \'equivalente \`a la cat\'egorie syntactique r\'eguli\`ere ${\mathcal C}^{\rm reg}_{\mathbb T}$ [resp. coh\'erente ${\mathcal C}^{\rm coh}_{\mathbb T}$, resp. g\'eom\'etrique ${\mathcal C}^{\rm geo}_{\mathbb T}$] d'une th\'eorie r\'eguli\`ere [resp. coh\'erente, resp. g\'eom\'etrique] ${\mathbb T}$.
\end{remark}

\bigskip

\begin{demo}
\begin{listeisansmarge}
\item Il faut v\'erifier que la collection $J$ des cribles v\'erifiant la propri\'et\'e de l'\'enonc\'e satisfait les axiomes de maximalit\'e, de stabilit\'e et de transitivit\'e.

\smallskip

C'est \'evident pour la maximalit\'e: tout crible maximal est couvrant.

\smallskip

La stabilit\'e r\'esulte de ce que, dans une telle cat\'egorie ${\mathcal C}$, la formation des images r\'eciproques de sous-objets par des morphismes respecte la formation des images ainsi que, dans le cas g\'eom\'etrique [resp. coh\'erent], la formation des r\'eunions [resp. des r\'eunions finies] de sous-objets.

\smallskip

Pour la transitivit\'e dans le cas r\'egulier, consid\'erons deux morphismes de ${\mathcal C}$
$$
x : X \longrightarrow S \qquad \mbox{et} \qquad y : Y \longrightarrow S
$$
avec le carr\'e cart\'esien associ\'e:
$$
\xymatrix{
X \times_S Y \ar[d]_{y'} \ar[r]^-{x'} &Y \ar[d]^y \\
X \ar[r]^x &S
}
$$

Si
$$
\exists_x (X) = S \qquad \mbox{et} \qquad \exists_{y'} (X \times_S Y) = X \, ,
$$
alors on a
$$
\exists_{x \circ y'} (X \times_S Y) = S \, .
$$
Comme $x \circ y' = y \circ x'$, on en d\'eduit comme voulu
$$
\exists_y (Y) = S \, .
$$

De m\^eme, pour la transitivit\'e dans le cas g\'eom\'etrique [resp. coh\'erent], consid\'erons deux familles [resp. deux familles finies] de morphismes
$$
x_i : X_i \longrightarrow S \qquad \mbox{et} \qquad y_j : Y_j \longrightarrow S
$$
telles que
$$
S = \bigvee_i \exists_{x_i} (X_i)
$$
et, pour tout indice $i$,
$$
X_i = \bigvee_j \exists_{y'_j} (X_i \times_S Y_j)
$$
en notant $y'_j : X_i \times_S Y_j \longrightarrow X_i$ le morphisme d\'eduit de $y_j : Y_j \to S$ par le changement de base $x_i : X_i \to S$.

\smallskip

Alors on a
$$
S = \bigvee_{i,j} \exists_{x_i \circ y'_j} (X_i \times_S Y_j)
$$
et a fortiori
$$
S = \bigvee_j \exists_{y_j} (Y_j)
$$
comme on voulait.

\medskip

\item Il s'agit de prouver que pour tout objet $Y$ de ${\mathcal C}$, le pr\'efaisceau repr\'esentable ${\rm Hom} (\bullet , Y)$ est un faisceau pour la topologie $J$ de ${\mathcal C}$.

\smallskip

Dans le cas g\'eom\'etrique [resp. coh\'erent, resp. r\'egulier], cela signifie que pour toute famille [resp. toute famille finie, resp. toute famille \`a un \'el\'ement] de morphismes
$$
X_i \longrightarrow X
$$
dont la r\'eunion des images est $X$, se donner un morphisme
$$
X \longrightarrow Y
$$
\'equivaut \`a se donner une famille de morphismes
$$
X_i \longrightarrow Y
$$
qui co{\"\i}ncident sur les $X_i \times_X X_j$.

\smallskip

Rempla\c cant les morphismes par leurs graphes, il suffit de prouver que se donner un sous-objet
$$
X' \xhookrightarrow{ \ { \ } \ } X
$$
\'equivaut \`a se donner une famille de sous-objets
$$
X'_i \xhookrightarrow{ \ { \ } \ } X_i
$$
dont les images r\'eciproques dans les $X_i \times_X X_j$ co{\"\i}ncident.

\smallskip

Tout sous-objet
$$
X' \xhookrightarrow{ \ { \ } \ } X
$$
d\'efinit une famille de sous-objets
$$
X' \times_X X_i \xhookrightarrow{ \ { \ } \ } X_i
$$
qui co{\"\i}ncident dans les $X_i \times_X X_j$.

\smallskip

Dans l'autre sens, toute famille de sous-objets
$$
X'_i \xhookrightarrow{ \ { \ } \ } X_i
$$
qui co{\"\i}ncident dans les $X_i \times_S X_j$ d\'efinit un sous-objet
$$
X' \xhookrightarrow{ \ { \ } \ } X
$$
qui est la r\'eunion des images des $X'_i \hookrightarrow X_i$ par les morphismes $X_i \to X$.

\smallskip

Ces deux applications sont inverses l'une de l'autre.

\smallskip

Si $X' \hookrightarrow X$ est un sous-objet, le fait que $X$ soit r\'eunion des images des $X_i$ implique, par compatibilit\'e des images avec les changements de base, que $X'$ est la r\'eunion des images des $X' \times_X X_i$ par les morphismes $X_i \to X$.

\smallskip

Si $(X'_i \hookrightarrow X_i)$ est une famille de sous-objets qui co{\"\i}ncident dans les $X_i \times_X X_j$, et $X' \hookrightarrow X$ est la r\'eunion de leurs images dans $X$, alors pour tout indice $j$, le produit fibr\'e
$$
X' \times_X X_j \xhookrightarrow{ \ { \ } \ } X_j
$$
est la r\'eunion des images par les morphismes $X_i \times_S X_j \to X_j$ des
$$
X'_i \times_X X_j = X_i \times_X X'_j \xhookrightarrow{ \ { \ } \ } X_i \times _X X_j
$$
et co{\"\i}ncide donc avec $X'_j \hookrightarrow X_j$.

\smallskip

C'est ce que l'on voulait. 
\end{listeisansmarge}
\end{demo}

\bigskip

On note le cas particulier suivant de la proposition \ref{propV73}:

\begin{cor}\label{corV74}

Soit ${\mathbb T}$ une th\'eorie g\'eom\'etrique [resp. coh\'erente, resp. r\'eguli\`ere].

\smallskip

Alors la topologie syntactique $J_{\mathbb T}$ de sa cat\'egorie syntactique g\'eom\'etrique ${\mathcal C}_{\mathbb T} = {\mathcal C}^{\rm geo}_{\mathbb T}$ [resp. coh\'erente ${\mathcal C}_{\mathbb T} = {\mathcal C}^{\rm coh}_{\mathbb T}$, resp. r\'eguli\`ere ${\mathcal C}_{\mathbb T} = {\mathcal C}^{\rm reg}_{\mathbb T}$] est sous-canonique, et le foncteur canonique
$$
\ell : {\mathcal C}_{\mathbb T} \longrightarrow \widehat{({\mathcal C}_{\mathbb T})}_{J_{\mathbb T}}
$$
est pleinement fid\`ele.
\end{cor}

\bigskip

\begin{remark}

Si ${\mathbb T}$ est une th\'eorie cart\'esienne, la topologie syntactique $J_{\mathbb T}$ sur la cat\'egorie syntactique cart\'esienne ${\mathcal C}^{\rm cart}_{\mathbb T}$ est la topologie discr\`ete, qui est sous-canonique.
\end{remark}

\bigskip

\begin{demo}

C'est une simple application de la proposition \ref{propV73} au cas particulier des cat\'egories syntactiques. 

\end{demo}

\subsection{Construction des topos classifiants}\label{subsec573}

\medskip

Le lemme suivant montre que la notion de foncteur cart\'esien [resp. r\'egulier, resp. coh\'erent, resp. g\'eom\'etri\-que] \`a valeurs dans un topos co{\"\i}ncide avec celle de morphisme plat et $J$-continu pour la topologie $J$ introduite au paragraphe pr\'ec\'edent:

\begin{lem}\label{lemV75}

Soit ${\mathcal C}$ une cat\'egorie cart\'esienne essentiellement petite.

\smallskip

Soit ${\mathcal E}$ un topos.

\smallskip

Alors:

\begin{listeimarge}

\item Un foncteur
$$
F : {\mathcal C} \longrightarrow {\mathcal E}
$$
est plat si et seulement si il est cart\'esien.

\medskip

\item Si ${\mathcal C}$ est une cat\'egorie r\'eguli\`ere [resp. coh\'erente, resp. g\'eom\'etrique], un foncteur
$$
F : {\mathcal C} \longrightarrow {\mathcal E}
$$
est plat et $J$-continu, pour la topologie $J$ de la proposition \ref{propV73}, si et seulement si il est r\'egulier [resp. coh\'erent, resp. g\'eom\'etrique].
\end{listeimarge}
\end{lem}

\begin{demo}
\begin{listeisansmarge}
\item La cat\'egorie ${\mathcal C}$ admettant des limites finies arbitraires, un foncteur
$$
F : {\mathcal C} \longrightarrow {\mathcal E}
$$
est plat si et seulement si il respecte les limites finies.

\medskip

\item Compte tenu de (i), il s'agit de prouver que si ${\mathcal C}$ est une cat\'egorie r\'eguli\`ere [resp. coh\'erente, resp. g\'eom\'etrique], un foncteur
$$
F : {\mathcal C} \longrightarrow {\mathcal E}
$$
respecte les images [resp. et les r\'eunions finies, resp. et les r\'eunions arbitraires] si et seulement si il transforme les familles $J$-couvrantes de morphismes
$$
x_i : X_i \longrightarrow X
$$
en familles globalement \'epimorphiques.

\smallskip

Or, dans le cas r\'egulier, les familles $J$-couvrantes sont celles qui comprennent un morphisme
$$
x : X' \longrightarrow X
$$
dont l'image est $X$.

\smallskip

Et, dans le cas g\'eom\'etrique [resp. r\'egulier], les familles $J$-couvrantes sont celles constitu\'ees de morphismes [resp. qui comprennent une sous-famille finie de morphismes]
$$
x_i : X_i \longrightarrow X
$$
dont la r\'eunion des images est $X$.

\smallskip

D'o\`u l'\'equivalence annonc\'ee. 
\end{listeisansmarge}
\end{demo}

\pagebreak

Gr\^ace \`a ce lemme, on d\'eduit de la proposition \ref{propV617} et de l'\'equivalence de Diaconescu que le topos ${\mathcal E}_{\mathbb T}$ associ\'e au site constitu\'e de la cat\'egorie syntactique ${\mathcal C}_{\mathbb T}$ d'une th\'eorie g\'eom\'etrique [resp. coh\'erente, resp. r\'eguli\`ere, resp. cart\'esienne] ${\mathbb T}$ et de sa topologie syntactique $J_{\mathbb T}$ v\'erifie la propri\'et\'e de repr\'esentation du foncteur des mod\`eles topossiques de ${\mathbb T}$ du th\'eor\`eme \ref{thmV71}:

\begin{prop}\label{propV76}

Soit ${\mathbb T}$ une th\'eorie g\'eom\'etrique [resp. coh\'erente, resp. r\'eguli\`ere, resp. cart\'esienne].

\smallskip

Soit la cat\'egorie syntactique g\'eom\'etrique ${\mathcal C}_{\mathbb T} = {\mathcal C}_{\mathbb T}^{\rm geo}$ [resp. coh\'erente ${\mathcal C}_{\mathbb T} = {\mathcal C}_{\mathbb T}^{\rm coh}$, resp. r\'eguli\`ere ${\mathcal C}_{\mathbb T} = {\mathcal C}_{\mathbb T}^{\rm reg}$, resp. cart\'esienne ${\mathcal C}_{\mathbb T} = {\mathcal C}_{\mathbb T}^{\rm cart}$] de ${\mathbb T}$, munie de sa topologie g\'eom\'etrique [resp. coh\'erente, resp. r\'eguli\`ere, resp. discr\`ete] $J_{\mathbb T}$.

\smallskip

Soit ${\mathcal E}_{\mathbb T}$ le topos des faisceaux sur le site $({\mathcal C}_{\mathbb T}, J_{\mathbb T})$.

\smallskip

Soit $U_{\mathbb T}$ le mod\`ele de ${\mathbb T}$ dans ${\mathcal E}_{\mathbb T}$ d\'efini comme l'image du mod\`ele canonique $M_{\mathbb T}$ de ${\mathbb T}$ dans ${\mathcal C}_{\mathbb T}$ par le foncteur canonique
$$
\ell : {\mathcal C}_{\mathbb T} \longrightarrow {\mathcal E}_{\mathbb T} \, .
$$

Alors, pour tout topos ${\mathcal E}$, le foncteur
$$
(f : {\mathcal E} \longrightarrow {\mathcal E}_{\mathbb T}) \longmapsto f^* U_{\mathbb T}
$$
est une \'equivalence de la cat\'egorie des morphismes de topos
$$
f = (f^* , f_*) : {\mathcal E}  \longrightarrow {\mathcal E}_{\mathbb T}
$$
sur la cat\'egorie
$$
{\mathbb T}\mbox{\rm -mod} \, ({\mathcal E})
$$
des mod\`eles de ${\mathbb T}$ dans ${\mathcal E}$.
\end{prop}


\begin{remark}

Il r\'esulte de cette proposition et de la propri\'et\'e d\'ej\`a d\'emontr\'ee d'unicit\'e \`a \'equivalence pr\`es des topos classifiants que
\begin{enumerate}
\item[$\bullet$] si ${\mathbb T}$ est une th\'eorie cart\'esienne (en particulier si c'est une th\'eorie alg\'ebrique ou plus g\'en\'eralement de Horn), ses cat\'egories syntactiques cart\'esienne ${\mathcal C}_{\mathbb T}^{\rm cart}$, r\'eguli\`ere ${\mathcal C}_{\mathbb T}^{\rm reg}$, coh\'erente ${\mathcal C}_{\mathbb T}^{\rm coh}$ et g\'eom\'etrique ${\mathcal C}_{\mathbb T}^{\rm geo}$ munies des topologies discr\`ete, r\'eguli\`ere, coh\'erente et g\'eom\'etrique, d\'efinissent des topos \'equivalents,
\item[$\bullet$] si ${\mathbb T}$ est une th\'eorie r\'eguli\`ere [resp. coh\'erente], ses cat\'egories syntactiques r\'eguli\`ere, coh\'erente et g\'eom\'etrique [resp. coh\'erente et g\'eom\'etrique] munies des topologies correspondantes d\'efinissent des topos \'equivalents.
\end{enumerate}
\end{remark}


\begin{demo}

D'apr\`es l'\'equivalence de Diaconescu, le foncteur
$$
(f : {\mathcal E} \longrightarrow {\mathcal E}_{\mathbb T}) \longmapsto (f^* \circ \ell : {\mathcal C}_{\mathbb T} \longrightarrow {\mathcal E})
$$
est une \'equivalence de la cat\'egorie des morphismes de topos
$$
f = (f^* , f_*) : {\mathcal E} \longrightarrow {\mathcal E}_{\mathbb T}
$$ 
sur la cat\'egorie des foncteurs
$$
F : {\mathcal C}_{\mathbb T} \longrightarrow {\mathcal E}
$$
qui sont plats et $J_{\mathbb T}$-continus.

\smallskip

Or, d'apr\`es le lemme \ref{lemV75}, un foncteur
$$
F : {\mathcal C}_{\mathbb T} \longrightarrow {\mathcal E}
$$
est plat et continu pour la topologie g\'eom\'etrique [resp. coh\'erente, resp. r\'eguli\`ere, resp. discr\`ete] $J_{\mathbb T}$ si et seulement si c'est un foncteur g\'eom\'etrique [resp. coh\'erent, resp. r\'egulier, resp. cart\'esien].

\smallskip

On conclut d'apr\`es la proposition \ref{propV617}. \end{demo}

\subsection{Compl\'etude du mod\`ele universel}\label{subsec574}

\medskip

On d\'eduit de la proposition \ref{propV76}, du corollaire \ref{corV74}, du corollaire \ref{corV611} et du corollaire \ref{corV615}~(iv) que les s\'equents d\'emontrables d'une th\'eorie sont ceux que v\'erifie son mod\`ele universel dans son topos classifiant:

\begin{cor}\label{corV77}

Soit ${\mathbb T}$ une th\'eorie g\'eom\'etrique du premier ordre dans une signature $\Sigma$.

\smallskip

Soit ${\mathcal E}_{\mathbb T}$ son topos classifiant muni du mod\`ele universel $U_{\mathbb T}$ de ${\mathbb T}$.

\smallskip

Alors:

\begin{listeimarge}

\item Un s\'equent g\'eom\'etrique de $\Sigma$
$$
\varphi \vdash_{\vec x} \psi
$$
est ${\mathbb T}$-d\'emontrable suivant les r\`egles d'inf\'erence de la logique g\'eom\'etrique si et seulement si il est v\'erifi\'e par le mod\`ele universel $U_{\mathbb T}$ de ${\mathbb T}$ dans ${\mathcal E}_{\mathbb T}$.

\medskip

\item Si ${\mathbb T}$ est une th\'eorie coh\'erente [resp. r\'eguli\`ere, resp. cart\'esienne], un s\'equent coh\'erent [resp. r\'egulier, resp. ${\mathbb T}$-cart\'esien] de $\Sigma$
$$
\varphi \vdash_{\vec x} \psi
$$
est ${\mathbb T}$-d\'emontrable suivant les r\`egles d'inf\'erence de la logique coh\'erente [resp. r\'eguli\`ere] si et seulement si il est v\'erifi\'e par le mod\`ele $U_{\mathbb T}$ de ${\mathbb T}$ dans ${\mathcal E}_{\mathbb T}$.
\end{listeimarge}
\end{cor}

\begin{remark}

Il r\'esulte de ce corollaire que si ${\mathbb T}$ est une th\'eorie coh\'erente [resp. r\'eguli\`ere, resp. cart\'esienne] de signature $\Sigma$, alors un s\'equent coh\'erent [resp. r\'egulier, resp. ${\mathbb T}$-cart\'esien] de $\Sigma$
$$
\varphi \vdash_{\vec x} \psi
$$
est ${\mathbb T}$-d\'emontrable suivant les r\`egles d'inf\'erence de la logique coh\'erente [resp. r\'eguli\`ere] si et seulement si il est ${\mathbb T}$-d\'emontrable suivant les r\`egles d'inf\'erence de la logique g\'eom\'etrique.
\end{remark}

\bigskip

\begin{demo}

Si ${\mathbb T}$ est une th\'eorie g\'eom\'etrique [resp. coh\'erente, resp. r\'eguli\`ere, resp. cart\'esienne], son topos classifiant ${\mathcal E}_{\mathbb T}$ s'identifie au topos des faisceaux sur la cat\'egorie g\'eom\'etrique ${\mathcal C}_{\mathbb T} = {\mathcal C}_{\mathbb T}^{\rm geo}$ [resp. coh\'erente ${\mathcal C}_{\mathbb T} = {\mathcal C}_{\mathbb T}^{\rm coh}$, resp. r\'eguli\`ere ${\mathcal C}_{\mathbb T} = {\mathcal C}_{\mathbb T}^{\rm reg}$, resp. cart\'esienne ${\mathcal C}_{\mathbb T} = {\mathcal C}_{\mathbb T}^{\rm cart}$] munie de la topologie g\'eom\'etrique [resp. coh\'erente, resp. r\'eguli\`ere, resp. discr\`ete] $J_{\mathbb T}$.

\smallskip

De plus, son mod\`ele universel $U_{\mathbb T}$ s'identifie \`a l'image par le foncteur canonique
$$
\ell : {\mathcal C}_{\mathbb T} \longrightarrow {\mathcal E}_{\mathbb T} 
$$
du mod\`ele canonique $M_{\mathbb T}$ de ${\mathbb T}$ dans ${\mathcal C}_{\mathbb T}$.

\smallskip

On sait d'apr\`es le corollaire \ref{corV611} [resp. \ref{corV615} (iv)] qu'un s\'equent g\'eom\'etrique, coh\'erent ou r\'egulier [resp. ${\mathbb T}$-cart\'esien] de $\Sigma$
$$
\varphi \vdash_{\vec x} \psi
$$
est ${\mathbb T}$-d\'emontrable suivant les r\`egles d'inf\'erence de la logique g\'eom\'etrique, coh\'erente ou r\'eguli\`ere si et seulement si, dans la cat\'egorie ${\mathcal C}_{\mathbb T}$, les deux sous-objets
$$
M_{\mathbb T} \, \varphi (\vec x) \qquad \mbox{et} \qquad M_{\mathbb T} \, \psi (\vec x)
$$
de l'objet $M_{\mathbb T} \top (\vec x)$ satisfont la relation d'inclusion
$$
M_{\mathbb T} \, \varphi (\vec x) \leq M_{\mathbb T} \, \psi (\vec x) \, .
$$

Or, d'apr\`es le corollaire \ref{corV74}, le foncteur canonique
$$
\ell : {\mathcal C}_{\mathbb T} \longrightarrow {\mathcal E}_{\mathbb T}
$$
est pleinement fid\`ele. De plus, il respecte les limites finies. Il en r\'esulte que deux sous-objets d'un objet de ${\mathcal C}_{\mathbb T}$ satisfont une relation d'inclusion si et seulement si leurs images par le foncteur $\ell$ satisfont cette relation.

\smallskip

D'o\`u la conclusion du corollaire. \end{demo}

\subsection{La notion d'\'equivalence de Morita entre th\'eories}\label{subsec575}

\medskip

Il est naturel de poser:

\begin{defn}\label{defV78}

Deux th\'eories g\'eom\'etriques ${\mathbb T}_1$ et ${\mathbb T}_2$ dans des signatures $\Sigma_1$ et $\Sigma_2$ sont dites \'equivalentes au sens de Morita, ou Morita-\'equivalentes, si leurs topos classifiants
$$
{\mathcal E}_{{\mathbb T}_1} \qquad \mbox{et} \qquad {\mathcal E}_{{\mathbb T}_2}
$$
sont \'equivalents.
\end{defn}

\begin{remarksqed}
\begin{listeisansmarge}
\item Comme les topologies syntactiques $J_{{\mathbb T}_1}$ et $J_{{\mathbb T}_2}$ sur les cat\'egories syntactiques ${\mathcal C}_{{\mathbb T}_1}^{\rm geo}$ et ${\mathcal C}_{{\mathbb T}_2}^{\rm geo}$ sont les topologies g\'eom\'etriques induites par leur structure cat\'egorique, on voit que deux th\'eories qui sont syntactiquement \'equivalentes sont a fortiori \'equivalentes au sens de Morita.

\medskip

\item En revanche, deux th\'eories peuvent \^etre Morita-\'equivalentes sans \^etre syntactiquement \'equivalentes.

\medskip

\item Le nom d'\'equivalence de Morita vient du cas particulier des th\'eories ${\mathbb T}_1$ et ${\mathbb T}_2$ des modules sur deux anneaux fix\'es $R_1$ et $R_2$.

\smallskip

Classiquement, de tels anneaux $R_1$ et $R_2$ sont dits \'equivalents au sens de Morita si les cat\'egories
$$
{\mathbb T}_1\mbox{-mod} \, ({\rm Ens}) \qquad \mbox{et} \qquad {\mathbb T}_2\mbox{-mod} \, ({\rm Ens})
$$
des modules sur $R_1$ et sur $R_2$ sont \'equivalentes.

On verra plus loin que, ces th\'eories \'etant alg\'ebriques et a fortiori cart\'esiennes, c'est \'equivalent \`a demander que leurs topos classifiants
$$
{\mathcal E}_{{\mathbb T}_1} \qquad \mbox{et} \qquad {\mathcal E}_{{\mathbb T}_2} 
$$
soient \'equivalents. 
\end{listeisansmarge}
\end{remarksqed}

\bigskip

On d\'eduit de la proposition \ref{propV76}:

\begin{cor}\label{corV79}

Soient ${\mathbb T}_1$ et ${\mathbb T}_2$ deux th\'eories g\'eom\'etriques dans des signatures $\Sigma_1$ et $\Sigma_2$.

\smallskip

Alors ces th\'eories sont Morita-\'equivalentes si et seulement si on peut associer \`a tout topos ${\mathcal E}$ deux foncteurs
$$
F_{\mathcal E} = {\mathbb T}_1\mbox{\rm -mod} \, ({\mathcal E}) \longrightarrow {\mathbb T}_2\mbox{\rm -mod} \, ({\mathcal E}) \, ,
$$
$$
G_{\mathcal E} = {\mathbb T}_2\mbox{\rm -mod} \, ({\mathcal E}) \longrightarrow {\mathbb T}_1\mbox{\rm -mod} \, ({\mathcal E})
$$
tels que:
\begin{enumerate}
\item[$\bullet$] ces foncteurs sont naturels au sens que pour tout morphisme de topos
$$
f = (f^* , f_*) : {\mathcal E}' \longrightarrow {\mathcal E} \, ,
$$
les carr\'es induits
$$
\xymatrix{
{\mathbb T}_1\mbox{\rm -mod} \, ({\mathcal E}) \ar[d]_{f^*} \ar[r]^{F_{\mathcal E}} &{\mathbb T}_2\mbox{\rm -mod} \, ({\mathcal E}) \ar[d]^{f^*} \\
{\mathbb T}_1\mbox{\rm -mod} \, ({\mathcal E}') \ar[r]^{F_{{\mathcal E}'}} &{\mathbb T}_2\mbox{\rm -mod} \, ({\mathcal E}')
} \qquad \qquad \xymatrix{
{\mathbb T}_2\mbox{\rm -mod} \, ({\mathcal E}) \ar[d]_{f^*} \ar[r]^{G_{\mathcal E}} &{\mathbb T}_1\mbox{\rm -mod} \, ({\mathcal E}) \ar[d]^{f^*} \\
{\mathbb T}_2\mbox{\rm -mod} \, ({\mathcal E}') \ar[r]^{G_{{\mathcal E}'}} &{\mathbb T}_1\mbox{\rm -mod} \, ({\mathcal E}')
}
$$
sont commutatifs \`a isomorphisme pr\`es,
\item[$\bullet$] ces foncteurs
$$
F_{\mathcal E} \qquad \mbox{et} \qquad G_{\mathcal E}
$$
sont des \'equivalences de cat\'egories r\'eciproques l'une de l'autre pour tout topos ${\mathcal E}$.

\end{enumerate}
\end{cor}

\begin{remarksqed}
\begin{listeisansmarge}
\item Autrement dit, deux th\'eories g\'eom\'etriques sont Morita-\'equivalentes si leurs foncteurs des mod\`eles topossiques sont \'equivalents.

\smallskip

Pour cette raison, l'\'equivalence de Morita des th\'eories g\'eom\'etriques peut aussi \^etre appel\'ee \'equivalence s\'emantique. 

\smallskip

Ainsi, l'\'equivalence syntactique de deux th\'eories g\'eom\'etriques entra{\^\i}ne leur \'equivalence s\'emantique, mais la r\'eciproque n'est pas vraie en g\'en\'eral.

\medskip

\item En particulier, si deux th\'eories g\'eom\'etriques ${\mathbb T}_1$ et ${\mathbb T}_2$ sont s\'emantiquement \'equivalentes, leurs cat\'egories de mod\`eles ensemblistes
$$
{\mathbb T}_1\mbox{-mod} \, ({\rm Ens}) \qquad \mbox{et} \qquad {\mathbb T}_2\mbox{-mod} ({\rm Ens})
$$
sont \'equivalentes.

\smallskip

La r\'eciproque n'est pas vraie en g\'en\'eral.

\medskip

\item La notion d'équivalence de Morita entre théories constitue le fondement de la méthodologie des topos comme ``ponts'' de Caramello (\cite{Memoire}, \cite{TST}), qui permet d'utiliser les topos comme des structures intermédiaires puissantes pour relier des théories distinctes, en exploitant les invariants définis sur eux et leur calcul à travers leurs différentes présentations. 
\end{listeisansmarge}
\end{remarksqed}

\subsection{Pr\'esentation des topos par des th\'eories}\label{subsec576}

\medskip

R\'eciproquement, tout topos s'\'ecrit comme le topos classifiant d'une infinit\'e de th\'eories g\'eom\'etriques du premier ordre.

\smallskip

Cela r\'esulte du th\'eor\`eme suivant:

\begin{thm}\label{thmV710}

Soit ${\mathcal C}$ une petite cat\'egorie munie d'une topologie de Grothendieck $J$.

\smallskip

Soit $\Sigma_{\mathcal C}$ la signature sans symbole de relation dont

\medskip

$\left\{\begin{matrix}
\bullet &\mbox{les sortes sont les objets $X$ de ${\mathcal C}$,} \hfill \\
{ \ } \\
\bullet &\mbox{les symboles de fonctions sont les morphismes de ${\mathcal C}$} \hfill \\
{ \ } \\
&f : X \longrightarrow Y \, .
\end{matrix} \right.
$

\bigskip

Alors:

\begin{listeimarge}

\item Il existe une th\'eorie alg\'ebrique ${\mathbb T}$ de signature $\Sigma$ telle que, pour toute cat\'egorie alg\'ebrique ${\mathcal D}$, la cat\'egorie de mod\`eles
$$
{\mathbb T}\mbox{\rm -mod} \, ({\mathcal D})
$$
s'identifie \`a la cat\'egorie
$$
[{\mathcal C},{\mathcal D}]
$$
des foncteurs ${\mathcal C} \to {\mathcal D}$ et, pour tout foncteur alg\'ebrique
$$
F : {\mathcal D} \longrightarrow {\mathcal D}' \, ,
$$
le foncteur
$$
{\mathbb T}\mbox{\rm -mod} \, ({\mathcal D}) \longrightarrow {\mathbb T}\mbox{\rm -mod} \, ({\mathcal D}')
$$
s'identifie au foncteur de composition avec $F$
$$
[{\mathcal C},{\mathcal D}] \longrightarrow [{\mathcal C},{\mathcal D}'] \, .
$$

\item Il existe une th\'eorie g\'eom\'etrique ${\mathbb T}_p$ quotient de ${\mathbb T}$ telle que, pour tout topos ${\mathcal E}$, la cat\'egorie de mod\`eles
$$
{\mathbb T}_{\!p}\mbox{\rm -mod} \, ({\mathcal E})
$$
s'identifie \`a la sous-cat\'egorie pleine de 
$$
[{\mathcal C},{\mathcal E}]
$$
constitu\'ee des foncteurs plats ${\mathcal C} \to {\mathcal E}$.

\medskip

\item Il existe une th\'eorie g\'eom\'etrique ${\mathbb T}_J$ quotient de ${\mathbb T}$ telle que, pour tout topos ${\mathcal E}$, la cat\'egorie de mod\`eles
$$
{\mathbb T}_J\mbox{\rm -mod} \, ({\mathcal E})
$$
s'identifie \`a la sous-cat\'egorie pleine de 
$$
[{\mathcal C},{\mathcal E}]
$$
constitu\'ee des foncteurs $J$-continus ${\mathcal C} \to {\mathcal E}$.

\item Si ${\mathbb T}_{p,J}$ d\'esigne la th\'eorie g\'eom\'etrique de signature $\Sigma$ dont les axiomes sont la r\'eunion de ceux de ${\mathbb T}_p$ et de ${\mathbb T}_J$, la cat\'egorie de mod\`eles dans n'importe quel topos ${\mathcal E}$
$$
{\mathbb T}_{p,J}\mbox{\rm -mod} \, ({\mathcal E})
$$
s'identifie \`a la cat\'egorie des morphismes de topos
$$
{\mathcal E} \longrightarrow \widehat{\mathcal C}_J \, .
$$
\end{listeimarge}
\end{thm}

\pagebreak
\begin{demo}
\begin{listeisansmarge}
\item Une $\Sigma$-structure dans une cat\'egorie alg\'ebrique ${\mathcal D}$ consiste en

\medskip

$\left\{\begin{matrix}
\bullet &\mbox{une famille d'objets $M\!X$ de ${\mathcal D}$ index\'es par les objets $X$ de ${\mathcal C}$,} \hfill  \\
{ \ } \\
\bullet &\mbox{une famille de morphismes de ${\mathcal D}$} \hfill \\
{ \ } \\
&M\!f : M\!X \longrightarrow MY \\
{ \ } \\
&\mbox{index\'es par les morphismes $f : X \to Y$ de ${\mathcal C}$.} \hfill
\end{matrix} \right.
$

\bigskip

\noindent C'est un foncteur si et seulement si

\medskip

$\left\{\begin{matrix}
\bullet &\mbox{$M {\rm id}_X = {\rm id}_{M\!X}$ pour tout objet $X$ de ${\mathcal C}$,} \hfill  \\
{ \ } \\
\bullet &\mbox{$M(g \circ f) = M\!g \circ M\!f$ pour tous morphismes de ${\mathcal C}$} \hfill \\
{ \ } \\
&X \xrightarrow{ \ f \ } Y \xrightarrow{ \ g \ } Z \, . 
\end{matrix} \right.
$

\bigskip

\noindent Cela revient \`a demander que $M$ soit un mod\`ele de la th\'eorie alg\'ebrique d\'efinie par les axiomes

\medskip

$\left\{\begin{matrix}
\bullet &\mbox{$\top \vdash_x {\rm id}_X(x) = x$,} \hfill  \\
{ \ } \\
\bullet &\mbox{$\top \vdash_x (g \circ f)(x) = g(f(x))$} \hfill 
\end{matrix} \right.
$

\medskip

\noindent pour toute sorte $X$ \`a laquelle on affecte une variable $x$, et pour tous morphismes $X \xrightarrow{ \ f \ } Y \xrightarrow{ \ g \ } Z$ de ${\mathcal C}$.

\medskip

\item Un foncteur \`a valeurs dans un topos ${\mathcal E}$
$$
{\mathcal C} \longrightarrow {\mathcal E} \, ,
$$
c'est-\`a-dire un mod\`ele $M$ dans ${\mathcal E}$ de la th\'eorie ${\mathbb T}$, est plat si et seulement si il poss\`ede les propri\'et\'es (B1), (B2) et (B3) de la proposition \ref{propIV39}.

\smallskip

Or, en notant $x^X , y^Y$ et $w^W$ des variables affect\'ees \`a des sortes $X,Y,Z$ de $\Sigma$ c'est-\`a-dire \`a des objets de ${\mathcal C}$, $M$ poss\`ede ces propri\'et\'es si et seulement si il satisfait les axiomes g\'eom\'etriques suivants:

\medskip
\begin{listeimarge}
\item[(B1)] est l'axiome g\'eom\'etrique sans variable libre
$$
\top \vdash \bigvee_{X \in {\rm Ob} ({\mathcal C})} (\exists \, x^X) \, \top (x^X) \, .
$$

\item[(B2)] est la famille d'axiomes g\'eom\'etriques index\'es par les paires d'objets $X,Y$ de ${\mathcal C}$
$$
\top \vdash_{x^X , y^Y} \bigvee_{\begin{matrix} \mbox{\footnotesize $W \in {\rm Ob}({\mathcal C})$} \\ \mbox{\footnotesize $(f,g) \in {\rm Hom} (W,X) \times {\rm Hom} (W,Y)$}\end{matrix}} (\exists \, w^W) (x^X = f(w^W) \wedge y^Y = g(w^W)) \, .
$$

\item[(B3)] est la famille d'axiomes g\'eom\'etriques index\'es par les paires de morphismes $X \!\! \raisebox{.7ex}{\xymatrix{\dar[r]^-{^{^{\mbox{\scriptsize$f$}}}}_-{g} &Y}}$ de ${\mathcal C}$
$$
f(x^X) = g(x^X) \vdash_{x^X} \bigvee_{\begin{matrix} \mbox{\footnotesize $W \in {\rm Ob} ({\mathcal C})$} \\  \mbox{\footnotesize $h \in {\rm Hom} (W,X)$} \\  \mbox{\footnotesize tel que $f \circ h = g \circ h$} \end{matrix}} (\exists \, w^W) (x^X = h(w^W)) \, .
$$
\end{listeimarge}

\item Un foncteur \`a valeurs dans un topos ${\mathcal E}$
$$
M : {\mathcal C} \longrightarrow {\mathcal E}
$$
est $J$-continu s'il transforme toute famille $J$-couvrante de morphismes de ${\mathcal C}$
$$
f_i : X_i \longrightarrow X \, , \qquad i \in I \, ,
$$
en une famille globalement \'epimorphique.

\smallskip

Cela revient \`a demander que, pour toute telle famille $J$-couvrante, $M$ consid\'er\'e comme un mod\`ele dans ${\mathcal E}$ de la th\'eorie ${\mathbb T}$ satisfasse l'axiome g\'eom\'etrique
$$
\top \vdash_{x^X} \bigvee_{i \in I} (\exists \, x_i^{X_i}) (x^X = f_i (x_i^{X_i}))
$$
o\`u les $x_i^{X_i}$ sont des variables affect\'ees aux sortes $X_i$.

\medskip

\item r\'esulte de (ii) et (iii) d'apr\`es l'\'equivalence de Diaconescu
$$
(f : {\mathcal E} \longrightarrow \widehat{\mathcal C}_J) \longmapsto (f^* \circ \ell : {\mathcal C} \xrightarrow{ \ \ell \ } \widehat{\mathcal C}_J \xrightarrow{ \ f^* \ } {\mathcal E})
$$
de la cat\'egorie des morphismes de topos
$$
f = (f^* , f_*) : {\mathcal E} \longrightarrow \widehat{\mathcal C}_J
$$
sur la cat\'egorie des foncteurs plats et $J$-continus
$$
M : {\mathcal C} \longrightarrow {\mathcal E} \, .
$$
\end{listeisansmarge}
\end{demo}

\section{Logique g\'eom\'etrique et topologies de Grothendieck}\label{sec58}

\subsection{Th\'eories quotients et sous-topos}\label{subsec581}

\medskip

On rappelle qu'une th\'eorie g\'eom\'etrique ${\mathbb T}'$ de m\^eme signature $\Sigma$ qu'une th\'eorie g\'eom\'etrique ${\mathbb T}$ est appel\'ee un quotient de ${\mathbb T}$ si tout s\'equent g\'eom\'etrique de $\Sigma$
$$
\varphi \vdash_{\vec x} \psi
$$
qui est d\'emontrable dans ${\mathbb T}$ est a fortiori d\'emontrable dans ${\mathbb T}'$.

\smallskip

Deux th\'eories quotients ${\mathbb T}_1$ et ${\mathbb T}_2$ d'une th\'eorie g\'eom\'etrique ${\mathbb T}$ sont dites \'equivalentes si chacune est un quotient de l'autre.

\smallskip

Montrons qu'il y a correspondance bijective entre les th\'eories quotients d'une th\'eorie g\'eom\'etrique ${\mathbb T}$, consid\'er\'ees \`a \'equivalence pr\`es, et les sous-topos du topos classifiant ${\mathcal E}_{\mathbb T}$ de ${\mathbb T}$:

\begin{thm}[Caramello (Théorème 3.2.5 \cite{TST})]\label{thmV81} 

Soit ${\mathbb T}$ une th\'eorie g\'eom\'etrique de signature $\Sigma$.

\smallskip

Soit ${\mathcal E}_{\mathbb T}$ son topos classifiant.

\smallskip

Alors:

\begin{listeimarge}

\item Pour toute th\'eorie quotient ${\mathbb T}'$ de ${\mathbb T}$, son topos classifiant ${\mathcal E}_{{\mathbb T}'}$ est muni d'un plongement canonique
$$
{\mathcal E}_{{\mathbb T}'} \xhookrightarrow{ \ { \ } \ } {\mathcal E}_{{\mathbb T}} \, .
$$

\item L'application
$$
{\mathbb T}' \longmapsto ({\mathcal E}_{{\mathbb T}'} \xhookrightarrow{ \ { \ } \ } {\mathcal E}_{{\mathbb T}})
$$
est une bijection de la collection des classes d'\'equivalence de th\'eories quotients ${\mathbb T}'$ de ${\mathbb T}$ sur l'ensemble des sous-topos de ${\mathcal E}_{\mathbb T}$.

\medskip

\item Etant donn\'ees deux th\'eories quotients ${\mathbb T}_1$ et ${\mathbb T}_2$ de ${\mathbb T}$, les sous-topos associ\'es ${\mathcal E}_{{\mathbb T}_1}$ et ${\mathcal E}_{{\mathbb T}_2}$ de ${\mathcal E}_{\mathbb T}$ satisfont la relation d'inclusion
$$
{\mathcal E}_{{\mathbb T}_1} \leq {\mathcal E}_{{\mathbb T}_2}
$$
si et seulement si ${\mathbb T}_1$ est une th\'eorie quotient de ${\mathbb T}_2$.
\end{listeimarge}
\end{thm}

\begin{remark}

Il r\'esulte en particulier de ce th\'eor\`eme que les th\'eories quotients d'une th\'eorie g\'eom\'etrique ${\mathbb T}$, consid\'er\'ees \`a \'equivalence pr\`es, forment un ensemble.
\end{remark}

\bigskip

\begin{demosansqed}

Si ${\mathbb T}$ et ${\mathbb T}'$ sont deux th\'eories g\'eom\'etriques de m\^eme signature $\Sigma$, alors pour tout topos ${\mathcal E}$, les cat\'egories de mod\`eles
$$
{\mathbb T}\mbox{-mod} \, ({\mathcal E}) \qquad \mbox{et} \qquad {\mathbb T}'\mbox{-mod} \, ({\mathcal E})
$$
sont par d\'efinition deux sous-cat\'egories pleines de la cat\'egorie de $\Sigma$-structures
$$
\Sigma\mbox{-str} \, ({\mathcal E}) \, .
$$

Si ${\mathbb T}'$ est un quotient de ${\mathbb T}$, alors pour tout topos ${\mathcal E}$, ${\mathbb T}'\mbox{-mod} \, ({\mathcal E})$ est une sous-cat\'egorie pleine de ${\mathbb T}\mbox{-mod} \, ({\mathcal E})$.

\smallskip

La famille des foncteurs de plongement
$$
{\mathbb T}'\mbox{-mod} \, ({\mathcal E}) \xhookrightarrow{ \ { \ } \ } {\mathbb T}\mbox{-mod} \, ({\mathcal E})
$$
correspond \`a un morphisme canonique de topos
$$
{\mathcal E}_{{\mathbb T}'} \longrightarrow {\mathcal E}_{\mathbb T}
$$
dont la composante d'image r\'eciproque prolonge le morphisme canonique entre les cat\'egories syntactiques
$$
{\mathcal C}_{\mathbb T} \longrightarrow {\mathcal C}_{{\mathbb T}'} \, .
$$

Pour montrer que le morphisme de topos
$$
{\mathcal E}_{{\mathbb T}'} \longrightarrow {\mathcal E}_{\mathbb T}
$$
est un plongement, puis que l'application
$$
{\mathbb T}' \longmapsto ({\mathcal E}_{{\mathbb T}'} \xhookrightarrow{ \ { \ } \ } {\mathcal E}_{\mathbb T})
$$
est une bijection pr\'eservant les relations d'ordre entre classes d'\'equivalence de th\'eories quotients de ${\mathbb T}$ et sous-topos de ${\mathcal E}_{\mathbb T}$, on a besoin du lemme suivant:
\end{demosansqed}

\begin{lem}\label{lemV82}

Soit ${\mathbb T}$ une th\'eorie g\'eom\'etrique de signature $\Sigma$.

\smallskip

Soit ${\mathcal C}_{\mathbb T}$ sa cat\'egorie syntactique g\'eom\'etrique, munie de sa topologie syntactique g\'eom\'etrique $J_{\mathbb T}$.

\smallskip

Associons \`a toute th\'eorie quotient ${\mathbb T}'$ de ${\mathbb T}$ la topologie sur ${\mathcal C}_{\mathbb T}$
$$
J_{{\mathbb T}'} \supseteq J_{\mathbb T}
$$
engendr\'ee par $J_{\mathbb T}$ et par les recouvrements
$$
(\varphi \wedge \psi)(\vec x) \xhookrightarrow{ \ { \ } \ } \varphi (\vec x)
$$
associ\'es aux axiomes
$$
\varphi \vdash_{\vec x} \psi
$$
de ${\mathbb T}'$.

\smallskip

En sens inverse, associons \`a toute topologie $J \supseteq J_{\mathbb T}$ de ${\mathcal C}_{\mathbb T}$ contenant $J_{\mathbb T}$ la th\'eorie ${\mathbb T}_J$ quotient de ${\mathbb T}$ d\'efinie par les axiomes de ${\mathbb T}$ compl\'et\'es par les axiomes
$$
\varphi \vdash_{\vec x} \bigvee_{i \in I} (\exists \, \vec x_i) \, \theta_i (\vec x_i , \vec x)
$$
associ\'es aux familles $J$-couvrantes
$$
\theta_i (\vec x_i , \vec x) : \varphi_i (\vec x_i) \longrightarrow \varphi (\vec x) \, , \qquad i \in I \, ,
$$
de morphismes de ${\mathcal C}_{\mathbb T}$.

\smallskip

Alors ces deux applications sont r\'eciproques l'une de l'autre au sens que:

\begin{listeimarge}

\item Pour toute th\'eorie quotient ${\mathbb T}'$ de ${\mathbb T}$, la th\'eorie
$$
{\mathbb T}_{J_{{{\mathbb T}'}}}
$$
associ\'ee \`a la topologie $J_{{\mathbb T}'}$ d\'efinie par ${\mathbb T}'$ lui est \'equivalente.

\medskip

\item Pour toute topologie $J \supseteq J_{\mathbb T}$ de ${\mathcal C}_{\mathbb T}$ contenant $J_{\mathbb T}$, la topologie
$$
J_{{\mathbb T}_J}
$$
d\'efinie par la th\'eorie quotient ${\mathbb T}_J$ de ${\mathbb T}$ associ\'ee \`a $J$ lui est \'egale.
\end{listeimarge}
\end{lem}

\begin{demolem}
\begin{listeisansmarge}
\item Montrons d'abord que tout axiome de la th\'eorie g\'eom\'etrique ${\mathbb T}'$ est d\'emontrable dans la th\'eorie ${\mathbb T}_{J_{{\mathbb T}'}}$.

\smallskip

Pour tout tel axiome
$$
\varphi \vdash_{\vec x} \psi \, ,
$$
le monomorphisme de ${\mathcal C}_{\mathbb T}$
$$
(\varphi \wedge \psi)(\vec x) \xhookrightarrow{ \ { \ } \ } \varphi (\vec x)
$$
est couvrant pour la topologie $J_{{\mathbb T}'}$ associ\'ee \`a ${\mathbb T}'$.

\smallskip

Il en r\'esulte que le s\'equent g\'eom\'etrique
$$
\varphi \vdash_{\vec x} \varphi \wedge \psi
$$
est un axiome de la th\'eorie ${\mathbb T}_{J_{{\mathbb T}'}}$.

\smallskip

Or il \'equivaut au s\'equent
$$
\varphi \vdash_{\vec x} \psi \, .
$$

Puis montrons que, r\'eciproquement, tout axiome de ${\mathbb T}_{J_{{\mathbb T}'}}$ est d\'emontrable dans la th\'eorie ${\mathbb T}'$.

\smallskip

Comme on vient de le voir, il en est ainsi des s\'equents
$$
\varphi \vdash_{\vec x} \varphi \wedge \psi
$$
associ\'es aux morphismes couvrants
$$
(\varphi \wedge \psi)(\vec x) \xhookrightarrow{ \ { \ } \ } \varphi (\vec x)
$$
index\'es par les axiomes $\varphi \vdash_{\vec x} \psi$ de ${\mathbb T}'$.

\smallskip

Or, ces morphismes couvrants engendrent la topologie $J_{{\mathbb T}'}$.

\smallskip

On est donc r\'eduit \`a prouver que la collection des familles de morphismes de ${\mathcal C}_{\mathbb T}$
$$
\theta_i (\vec x_i , \vec x) : \varphi_i (\vec x_i) \longrightarrow \varphi (\vec x) \, , \quad i \in I \, ,
$$
telles que le s\'equent
$$
\varphi \vdash_{\vec x} \bigvee_{i \in I} (\exists \, \vec x_i) \, \theta_i (\vec x_i , \vec x)
$$
soit ${\mathbb T}'$-d\'emontrable, est stable par changement de base et par transitivit\'e.

\smallskip

Pour la stabilit\'e par changement de base, consid\'erons donc un morphisme de ${\mathcal C}_{\mathbb T}$
$$
\theta (\vec y , \vec x) : \psi (\vec y) \longrightarrow \varphi (\vec x) \, .
$$
Si le s\'equent
$$
\varphi \vdash_{\vec x} \bigvee_{i \in I} (\exists \, \vec x_i) \, \theta_i (\vec x_i , \vec x)
$$
est ${\mathbb T}'$-d\'emontrable, il en est a fortiori de m\^eme du s\'equent
$$
\psi \vdash_{\vec y} \bigvee_{i \in I} (\exists \, \vec x_i) (\exists \, \vec x) (\theta_i (\vec x_i , \vec x) \wedge \theta (\vec y , \vec x))
$$
puisque le s\'equent
$$
\psi \vdash_{\vec y} (\exists \, \vec x) \, \theta (\vec y , \vec x)
$$
est ${\mathbb T}$-d\'emontrable.

\smallskip

Pour la stabilit\'e par transitivit\'e, consid\'erons aussi une famille de morphismes de ${\mathcal C}_{\mathbb T}$
$$
\theta' (\vec y_j , \vec x) : \psi_j (\vec y_j) \longrightarrow \varphi (\vec x) \, , \quad j \in I' \, ,
$$
telle que, pour tout indice $i \in I$, la famille de morphismes qui s'en d\'eduit par changement de base via $\theta_i (\vec x_i , \vec x) : \varphi_i (\vec x_i) \to \varphi (\vec x)$ poss\`ede la propri\'et\'e que le s\'equent associ\'e
$$
\varphi_i \vdash_{\vec x_i} \bigvee_{j \in I'} (\exists \, \vec y_j) (\exists \, \vec x) (\theta'_j (\vec y_j , \vec x) \wedge \theta_i (\vec x_i , \vec x))
$$
soit ${\mathbb T}'$-d\'emontrable.

\smallskip

Pour tout tel $i \in I$, le sous-objet $\theta_i (\vec x_i , \vec x)$ de $\varphi_i (\vec x_i) \times \varphi (\vec x)$ dans ${\mathcal C}_{\mathbb T}$ se projette sur $\varphi_i (\vec x_i)$ par un isomorphisme, et donc le s\'equent
$$
\theta_i \vdash_{\vec x_i , \vec x} \bigvee_{j \in I'} (\exists \, \vec y_j) (\theta'_j (\vec y_j , \vec x) \wedge \theta_i (\vec x_i , \vec x))
$$
est \'egalement ${\mathbb T}'$-d\'emontrable.

\smallskip

Comme il en est de m\^eme du s\'equent
$$
\varphi \vdash_{\vec x} \bigvee_{i \in I} (\exists \, \vec x_i) \, \theta_i (\vec x_i , \vec x) \, ,
$$
on obtient que le s\'equent
$$
\varphi \vdash_{\vec x} \bigvee_{i \in I} \ \bigvee_{j \in I'} (\exists \, \vec x_i) (\exists \, \vec y_j) (\theta'_j (\vec y_j , \vec x) \wedge \theta_i (\vec x_i , \vec x))
$$
est ${\mathbb T}'$-d\'emontrable, donc a fortiori le s\'equent
$$
\varphi \vdash_{\vec x} \bigvee_{j \in I'} (\exists \, \vec y_j) \, \theta'_j (\vec y_j , \vec x) \, .
$$
Cela ach\`eve la preuve de (i).

\medskip

\item Montrons d'abord que la topologie $J$ est contenue dans la topologie $J_{{\mathbb T}_J}$ associ\'ee \`a la th\'eorie ${\mathbb T}_J$ d\'efinie par $J$.

\smallskip

Pour toute famille $J$-couvrante de morphismes de ${\mathcal C}_{\mathbb T}$
$$
\theta_i (\vec x_i , \vec x) : \varphi_i (\vec x_i) \longrightarrow \varphi (\vec x) \, , \quad i \in I \, ,
$$
le s\'equent
$$
\varphi \vdash_{\vec x} \bigvee_{i \in I} (\exists \, \vec x_i) \, \theta_i (\vec x_i , \vec x)
$$
est un axiome de la th\'eorie ${\mathbb T}_J$, et donc le monomorphisme de ${\mathcal C}_{\mathbb T}$
$$
\varphi (\vec x) \wedge \bigvee_{i \in I} (\exists \, \vec x_i) \, \theta_i (\vec x_i , \vec x) \xhookrightarrow{ \ { \ } \ } \varphi (\vec x)
$$
est couvrant pour la topologie $J_{{\mathbb T}_J}$.

\smallskip

Comme la famille de morphismes
$$
\theta_{i'} (\vec x_{i'} , \vec x) : \varphi_{i'} (\vec x_{i'}) \longrightarrow \bigvee_{i \in I} (\exists \, \vec x_i) \, \theta_i (\vec x_i , \vec x) \, , \quad i' \in I \, ,
$$
est couvrante pour la topologie $J_{\mathbb T} \subseteq J_{{\mathbb T}_J}$, on conclut comme voulu que la famille des morphismes de ${\mathcal C}_{\mathbb T}$
$$
\theta_i (\vec x_i , \vec x) : \varphi_i (\vec x_i) \longrightarrow \varphi (\vec x) \, , \quad i \in I \, ,
$$
est couvrante pour la topologie $J_{{\mathbb T}_J}$.

\smallskip

Enfin, v\'erifions que la topologie $J_{{\mathbb T}_J}$ est contenue dans $J$.

\smallskip

Or, $J_{{\mathbb T}_J}$ est par d\'efinition la topologie engendr\'ee sur $J_{\mathbb T} \subseteq J$ par les monomorphismes couvrants de ${\mathcal C}_{\mathbb T}$
$$
\varphi (\vec x) \wedge \bigvee_{i \in I} (\exists \, \vec x_i) \, \theta_i (\vec x_i , \vec x) \xhookrightarrow{ \ { \ } \ } \varphi (\vec x)
$$
associ\'es aux familles de morphismes de ${\mathcal C}_{\mathbb T}$
$$
\theta_i (\vec x_i , \vec x) : \varphi_i (\vec x_i) \longrightarrow \varphi (\vec x) \, , \quad i \in I \, ,
$$
qui sont $J$-couvrantes.

\smallskip

On conclut en remarquant qu'une telle famille de morphismes de ${\mathcal C}_{\mathbb T}$
$$
\theta_i (\vec x_i , \vec x) : \varphi_i (\vec x_i) \longrightarrow \varphi (\vec x) \, , \quad i \in I \, ,
$$
est $J$-couvrante si et seulement si le monomorphisme associ\'e
$$
\varphi (\vec x) \wedge \bigvee_{i \in I} (\exists \, \vec x_i) \, \theta_i (\vec x_i , \vec x) \xhookrightarrow{ \ { \ } \ } \varphi (\vec x)
$$
est $J$-couvrant.

\smallskip

Cela ach\`eve de prouver le lemme. 
\end{listeisansmarge}
\end{demolem}

\bigskip

\noindent {\bf Fin de la d\'emonstration du th\'eor\`eme \ref{thmV81}:}

\smallskip

Le topos classifiant ${\mathcal E}_{\mathbb T}$ de ${\mathbb T}$ est le topos des faisceaux sur le site $({\mathcal C}_{\mathbb T} , J_{\mathbb T})$. Il y a donc une bijection renversant les relations d'ordre entre l'ensemble ordonn\'e des sous-topos de ${\mathcal E}_{\mathbb T}$ et celui des topologies $J \supseteq J_{\mathbb T}$ sur ${\mathcal C}_{\mathbb T}$.

\smallskip

On en d\'eduit d'apr\`es le lemme \ref{lemV82} ci-dessus une bijection pr\'eservant les relations d'ordre
$$
{\mathbb T}' \longmapsto {\mathcal E}'_{{\mathbb T}'} = \widehat{({\mathcal C}_{\mathbb T})_{J_{{\mathbb T}'}}} 
$$
entre l'ensemble (ordonn\'e par la relation de quotient) des classes d'\'equivalence de th\'eories g\'eom\'etriques quotients ${\mathbb T}'$ de ${\mathbb T}$ et celui des sous-topos de ${\mathcal E}_{\mathbb T}$.

\smallskip

Il reste \`a d\'emontrer que toute telle th\'eorie quotient ${\mathbb T}'$ de ${\mathbb T}$ admet pour topos classifiant le topos
$$
{\mathcal E}'_{{\mathbb T}'} = \widehat{({\mathcal C}_{\mathbb T})_{J_{{\mathbb T}'}}} \, .
$$

Consid\'erons un mod\`ele $M$ de ${\mathbb T}$ dans un topos ${\mathcal E}$. Il existe un morphisme de topos, unique \`a unique isomorphisme pr\`es,
$$
f = (f^* , f_*) : {\mathcal E} \longrightarrow {\mathcal E}_{\mathbb T}
$$
tel que $M$ s'identifie au transform\'e du mod\`ele canonique $M_{\mathbb T}$ de ${\mathbb T}$ dans ${\mathcal C}_{\mathbb T}$ par le foncteur compos\'e
$$
{\mathcal C}_{\mathbb T} \xrightarrow{ \ \ell \ } {\mathcal E}_{\mathbb T} \xrightarrow{ \ f^* \ } {\mathcal E} \, .
$$

Un tel mod\`ele $M$ de ${\mathbb T}$ est un mod\`ele de ${\mathbb T}'$ si et seulement si, pour tout axiome de ${\mathbb T}'$
$$
\varphi \vdash_{\vec x} \psi \, ,
$$
on a
$$
M\varphi (\vec x) \wedge M\psi (\vec x) = M\varphi (\vec x) \, .
$$

Cela revient \`a demander que, pour tout tel axiome, le monomorphisme de ${\mathcal C}_{\mathbb T}$ correspondant
$$
\varphi (\vec x) \wedge \psi (\vec x) \xhookrightarrow{ \ { \ } \ } \varphi (\vec x)
$$
soit transform\'e par le foncteur $f^* \circ \ell$ en un \'epimorphisme de ${\mathcal E}$.

\smallskip

Comme la topologie $J_{{\mathbb T}'}$ est engendr\'ee sur $J_{\mathbb T}$ par les monomorphismes couvrants
$$
\varphi (\vec x) \wedge \psi (\vec x) \xhookrightarrow{ \ { \ } \ } \varphi (\vec x)
$$
associ\'es aux axiomes $\varphi \vdash_{\vec x} \psi$ de ${\mathbb T}'$, cela signifie exactement que le foncteur
$$
f^* \circ \ell : {\mathcal C}_{\mathbb T} \longrightarrow {\mathcal E}
$$
est $J_{{\mathbb T}'}$-continu c'est-\`a-dire d\'efinit un morphisme de topos
$$
{\mathcal E} \longrightarrow {\mathcal E}'_{{\mathbb T}'} \, .
$$

Cela ach\`eve la d\'emonstration du th\'eor\`eme. \hfill $\Box$

\subsection{D\'emontrabilit\'e et topologies engendr\'ees}\label{subsec582}

\medskip

Le th\'eor\`eme \ref{thmV81} permet d'exprimer les probl\`emes de d\'emontrabilit\'e dans les th\'eories g\'eom\'etriques comme des probl\`emes de topologies de Grothendieck:

\begin{cor}\label{corV83}

Soit ${\mathbb T}$ une th\'eorie g\'eom\'etrique de signature $\Sigma$.

\smallskip

Soit $({\mathcal C},J)$ un site de pr\'esentation du topos classifiant ${\mathcal E}_{\mathbb T}$ de ${\mathbb T}$, au sens qu'existe une \'equivalence
$$
{\mathcal E}_{\mathbb T} \cong \widehat{\mathcal C}_J \, .
$$

Soit d'autre part un s\'equent g\'eom\'etrique de $\Sigma$
$$
\varphi \vdash_{\vec x} \psi \, .
$$

Consid\'erant la th\'eorie quotient ${\mathbb T}'$ de ${\mathbb T}$ d\'efinie en adjoignant cet axiome \`a ceux de ${\mathbb T}$ et l'unique topologie $J' \supseteq J$ de ${\mathcal C}$ qui d\'efinit le topos classifiant ${\mathcal E}_{{\mathbb T}'}$ de ${\mathbb T}'$ comme sous-topos $\widehat{\mathcal C}_{J'}$ de $\widehat{\mathcal C}_J \cong {\mathcal E}_{\mathbb T}$, on a:

\begin{listeimarge}

\item Le s\'equent
$$
\varphi \vdash_{\vec x} \psi
$$
contredit les axiomes de ${\mathbb T}$ si et seulement si la topologie $J'$ est maximale, c'est-\`a-dire si tout objet de ${\mathcal C}$ admet pour crible $J'$-couvrant le crible vide.

\medskip

\item Le s\'equent
$$
\varphi \vdash_{\vec x} \psi
$$
est d\'emontrable dans la th\'eorie ${\mathbb T}$ si et seulement si
$$
J' = J \, .
$$
\end{listeimarge}
\end{cor}

\begin{remark}

Cette \'equivalence entre ${\mathbb T}$-d\'emontrabilit\'e et les topologies $J'$ de ${\mathcal C}$ qui raffinent $J$ pose la question de savoir comment engendrer concr\`etement $J'$ \`a partir de $J$ et du s\'equent $\varphi \vdash_{\vec x} \psi$.

\smallskip

Comme nous allons le voir tout de suite, cela se d\'eduit du lemme \ref{lemV82} dans le cas o\`u $({\mathcal C},J)$ est le site syntactique $({\mathcal C}_{\mathbb T} , J_{\mathbb T})$ de ${\mathbb T}$.

\smallskip

Le cas g\'en\'eral sera trait\'e plus loin.
\end{remark}

\bigskip

\begin{demo}

D'apr\`es le th\'eor\`eme \ref{thmV81} combin\'e au th\'eor\`eme \ref{thmIV72}, il y a correspondance bijective entre les th\'eories g\'eom\'etriques quotients de ${\mathbb T}$ et les topologies de ${\mathcal C}$ qui raffinent $J$.

\smallskip

Dans la situation du corollaire, le s\'equent $\varphi \vdash_{\vec x} \psi$ contredit les axiomes de ${\mathbb T}$ si et seulement si le topos classifiant ${\mathcal E}_{{\mathbb T}'}$ est trivial c'est-\`a-dire si la topologie $J'$ de ${\mathcal C}$ qui le d\'efinit est la topologie maximale.

\smallskip

D'autre part, ce s\'equent $\varphi \vdash_{\vec x} \psi$ est d\'emontrable dans ${\mathbb T}$ si et seulement si ${\mathbb T}'$ est \'equivalente \`a ${\mathbb T}$ c'est-\`a-dire
$$
{\mathcal E}_{{\mathbb T}'} = {\mathcal E}_{\mathbb T} \qquad \mbox{et} \qquad J' = J \, .
$$

C'est ce que l'on voulait. \end{demo}

\bigskip

Dans le cas de sites syntactiques, on peut pr\'eciser d\`es \`a pr\'esent:

\begin{prop}\label{propV84}

Soit ${\mathbb T}$ une th\'eorie cart\'esienne [resp. r\'eguli\`ere, resp. coh\'erente, resp. g\'eom\'etrique] dans une signature $\Sigma$.

\smallskip

Soit ${\mathcal C}_{\mathbb T} = {\mathcal C}_{\mathbb T}^{\rm cart}$ [resp. ${\mathcal C}_{\mathbb T}^{\rm reg}$, resp. ${\mathcal C}_{\mathbb T}^{\rm coh}$, resp. ${\mathcal C}_{\mathbb T}^{\rm geo}$] sa cat\'egorie syntactique cart\'esienne [resp. r\'eguli\`ere, resp. coh\'erente, resp. g\'eom\'etrique] munie de la topologie discr\`ete [resp. r\'eguli\`ere, resp. coh\'erente, resp. g\'eom\'etrique] $J_{\mathbb T}$.

\smallskip

Soit ${\mathbb T}'$ la th\'eorie quotient de ${\mathbb T}$ d\'efinie en adjoignant aux axiomes de ${\mathbb T}$ des axiomes suppl\'ementaires constitu\'es de s\'equents ${\mathbb T}$-cart\'esiens [resp. r\'eguliers, resp. coh\'erents, resp. g\'eom\'etriques]
$$
\varphi_i \vdash_{\vec x_i} \psi_i \, .
$$
Soit $J'$ la topologie de ${\mathcal C}_{\mathbb T}$ engendr\'ee par $J_{\mathbb T}$ et par la famille de cribles engendr\'es par les monomorphismes de ${\mathcal C}_{\mathbb T}$
$$
\varphi_i (\vec x_i) \wedge \psi_i (\vec x_i) \xhookrightarrow{ \ { \ } \ } \varphi_i (\vec x_i) \, .
$$

Alors:

\begin{listeimarge}

\item La th\'eorie ${\mathbb T}'$ est contradictoire si et seulement si $J'$ est la topologie maximale de ${\mathcal C}_{\mathbb T}$, c'est-\`a-dire admet pour cribles couvrants tous les cribles vides.

\medskip

\item Un s\'equent ${\mathbb T}$-cart\'esien [resp. r\'egulier, resp. coh\'erent, resp. g\'eom\'etrique] de $\Sigma$
$$
\varphi \vdash_{\vec x} \psi
$$
est d\'emontrable dans la th\'eorie ${\mathbb T}'$ si et seulement si le monomorphisme de ${\mathcal C}_{\mathbb T}$
$$
\varphi (\vec x) \wedge \psi (\vec x) \xhookrightarrow{ \ { \ } \ } \varphi (\vec x)
$$
est $J'$-couvrant.
\end{listeimarge}
\end{prop}

\begin{demo}

On sait a priori que toute th\'eorie quotient ${\mathbb T}'$ de ${\mathbb T}$ admet pour topos classifiant ${\mathcal E}_{{\mathbb T}'}$ un sous-topos de
$$
{\mathcal E}_{\mathbb T} = \widehat{({\mathcal C}_{\mathbb T})}_{J_{\mathbb T}}
$$
d\'efini par une unique topologie $J'$ de ${\mathcal C}_{\mathbb T}$ qui contient $J_{\mathbb T}$.

\smallskip

Pour d\'emontrer (i) et (ii), il suffit de v\'erifier que si une telle th\'eorie quotient ${\mathbb T}'$ est d\'efinie en adjoignant aux axiomes de ${\mathbb T}$ une liste d'axiomes
$$
\varphi_i \vdash_{\vec x_i} \psi_i
$$
qui sont des s\'equents ${\mathbb T}$-cart\'esiens [resp. r\'eguliers, resp. coh\'erents, resp. g\'eom\'etriques] de $\Sigma$, alors $J'$ est la plus petite topologie de ${\mathcal C}_{\mathbb T}$ contenant $J_{\mathbb T}$ pour laquelle les monomorphismes de ${\mathcal C}_{\mathbb T}$
$$
\varphi_i (\vec x_i) \wedge \psi_i (\vec x_i)  \xhookrightarrow{ \ { \ } \ } \varphi (\vec x_i)
$$
sont couvrants.

\smallskip

On utilise pour cela le m\^eme argument qu'\`a la fin de la d\'emonstration du th\'eor\`eme \ref{thmV81}:

\smallskip

Les mod\`eles $M$ de ${\mathbb T}$ dans un topos ${\mathcal E}$ sont les images du mod\`ele canonique $M_{\mathbb T}$ de ${\mathbb T}$ dans ${\mathcal C}_{\mathbb T}$ par les foncteurs plats et $J_{\mathbb T}$-continus
$$
F : {\mathcal C}_{\mathbb T} \longrightarrow {\mathcal E} \, .
$$

Un tel mod\`ele $M$ de ${\mathbb T}$ est un mod\`ele de ${\mathbb T}'$ si et seulement si, pour tout axiome de la liste qui d\'efinit ${\mathbb T}'$
$$
\varphi_i \vdash_{\vec x_i} \psi_i \, ,
$$
le foncteur $F : {\mathcal C}_{\mathbb T} \to {\mathcal E}$ transforme le monomorphisme de ${\mathcal C}_{\mathbb T}$
$$
u_i : \varphi_i (\vec x_i) \wedge \psi_i (\vec x_i)  \xhookrightarrow{ \ { \ } \ } \varphi_i (\vec x_i)
$$
en un \'epimorphisme de ${\mathcal E}$.

\smallskip

Cela revient \`a demander que ce foncteur soit $J'$-continu pour la plus petite topologie $J'$ de ${\mathcal C}_{\mathbb T}$ contenant $J_{\mathbb T}$ pour laquelle les monomorphismes $u_i$ de ${\mathcal C}_{\mathbb T}$ sont couvrants.

\smallskip

C'est ce que l'on voulait. \end{demo}

\section{Th\'eories de type pr\'efaisceau}\label{sec59}

\subsection{La notion de th\'eorie de type pr\'efaisceau}\label{subsec591}

\medskip

Toute th\'eorie g\'eom\'etrique du premier ordre ${\mathbb T}$ admet un topos classifiant ${\mathcal E}_{\mathbb T}$.

\smallskip

Pour n'importe quelle pr\'esentation de celui-ci comme topos des faisceaux sur un site $({\mathcal C},J)$
$$
{\mathcal E}_{\mathbb T} \cong \widehat{\mathcal C}_J \, ,
$$
le topos classifiant de ${\mathbb T}$ appara{\^\i}t comme un sous-topos du topos $\widehat{\mathcal C}$ des pr\'efaisceaux sur une cat\'egorie essentiellement petite ${\mathcal C}$.

\smallskip

D'apr\`es le th\'eor\`eme \ref{thmV710}, il existe des th\'eories g\'eom\'etriques, par exemple la th\'eorie ${\mathbb T}_p$ des foncteurs plats sur ${\mathcal C}$, dont le topos classifiant est \'equivalent \`a $\widehat{\mathcal C}$, et il existe des th\'eories quotients de celles-ci, par exemple la th\'eorie ${\mathbb T}_{p,J}$ des foncteurs plats et $J$-continus sur ${\mathcal C}$, dont le topos classifiant s'identifie au sous-topos $\widehat{\mathcal C}_J \cong {\mathcal E}_{\mathbb T}$ de $\widehat{\mathcal C}$.

\smallskip

Cela fait appara{\^\i}tre la th\'eorie ${\mathbb T}$ comme s\'emantiquement \'equivalente \`a une th\'eorie quotient d'une th\'eorie dont le topos classifiant s'identifie au topos $\widehat{\mathcal C}$ des pr\'efaisceaux sur ${\mathcal C}$.

\smallskip

Ces consid\'erations suffisent \`a justifier l'importance des ``th\'eories de type pr\'efaisceau'' au sens de la d\'efinition suivante:

\begin{defn}\label{defV91}

Une th\'eorie g\'eom\'etrique du premier ordre ${\mathbb T}$ est dite ``de type pr\'efaisceau'' si son topos classifiant
$$
{\mathcal E}_{\mathbb T}
$$
est \'equivalent au topos $\widehat{\mathcal C}$ des pr\'efaisceaux sur une cat\'egorie (essentiellement) petite ${\mathcal C}$.
\end{defn}

\begin{remarkqed}

Les consid\'erations qui pr\'ec\`edent l'\'enonc\'e se r\'esument donc \`a dire que toute pr\'esentation du topos classifiant d'une th\'eorie g\'eom\'etrique ${\mathbb T}$ comme topos des faisceaux sur un site fait appara{\^\i}tre ${\mathbb T}$ comme s\'emantiquement \'equivalente \`a une th\'eorie quotient d'une th\'eorie de type pr\'efaisceau. \end{remarkqed}

\bigskip

La proposition \ref{propV76} fournit une classe tr\`es importante de th\'eories de type pr\'efaisceau:

\begin{cor}\label{corV92}

Toute th\'eorie cart\'esienne est de type pr\'efaisceau.

\smallskip

En particulier, toute th\'eorie alg\'ebrique et plus g\'en\'eralement toute th\'eorie de Horn est de type pr\'efaisceau.
\end{cor}

\begin{remarks}
\begin{listeisansmarge}
\item Ainsi, la th\'eorie des groupes, celle des mono{\"\i}des, celles des actions de groupes ou de mono{\"\i}des, celle des anneaux, celle des modules sur un anneau, celle des modules sur un anneau fix\'e ou celle des espaces vectoriels sur un corps fix\'e, sont des th\'eories de type pr\'efaiscau.

\medskip

\item En revanche, et comme on verra plus loin, la th\'eorie des corps n'est pas une th\'eorie de type pr\'efaisceau.

\medskip

\item La th\'eorie des relations d'\'equivalence et celle des relations d'ordre sont de type pr\'efaisceau. En effet, ce sont des th\'eories de Horn.

\medskip

\item La th\'eorie des cat\'egories -- qui est une th\'eorie cart\'esienne -- est une th\'eorie de type pr\'efaisceau.
\end{listeisansmarge}
\end{remarks}

\medskip

\begin{demo}

En effet, d'apr\`es la proposition \ref{propV76}, toute th\'eorie cart\'esienne ${\mathbb T}$ admet pour topos classifiant le topos des pr\'efaisceaux sur sa cat\'egorie syntactique cart\'esienne ${\mathcal C}_{\mathbb T}^{\rm cart}$. \end{demo}

\bigskip

Le cas particulier suivant du corollaire pr\'ec\'edent m\'erite d'\^etre \'enonc\'e \`a part:

\begin{cor}\label{corV93}

Pour toute signature $\Sigma$, la th\'eorie vide -- c'est-\`a-dire sans axiome -- de signature $\Sigma$ est une th\'eorie de type pr\'efaisceau.
\end{cor}

\begin{remark}

Ainsi, toute th\'eorie g\'eom\'etrique ${\mathbb T}$ appara{\^\i}t comme un quotient d'au moins une th\'eorie de type pr\'efaisceau.
\end{remark}

\bigskip

\begin{demo}

En effet, toute th\'eorie vide est une th\'eorie de Horn. \end{demo}

\subsection{Mod\`eles ensemblistes des th\'eories de type pr\'efaisceau}\label{subsec592}

\medskip

La propri\'et\'e de d\'efinition des topos classifiants et le th\'eor\`eme \ref{thmIV38} entra{\^\i}nent:

\begin{prop}\label{propV94}

Soit ${\mathbb T}$ une th\'eorie de type pr\'efaisceau.

\smallskip

Alors, pour toute cat\'egorie (essentiellement) petite ${\mathcal C}$ telle que
$$
{\mathcal E}_{\mathbb T} \cong \widehat{\mathcal C} \, ,
$$
la cat\'egorie des mod\`eles ensemblistes de ${\mathbb T}$
$$
{\mathbb T}\mbox{\rm -mod} \, ({\rm Ens})
$$
s'identifie \`a la cat\'egorie
$$
{\rm Ind} \, ({\mathcal C}^{\rm op})
$$
des {\rm ind}-objets de la cat\'egorie oppos\'ee ${\mathcal C}^{\rm op}$.
\end{prop}

\begin{remark}

En particulier, pour toutes cat\'egories (essentiellement) petites ${\mathcal C}$ et ${\mathcal D}$, toute paire d'\'equivalences de topos
$$
\widehat{\mathcal C} \cong {\mathcal E}_{\mathbb T} \cong \widehat{\mathcal D}
$$
induit une \'equivalence de cat\'egories
$$
{\rm Ind} \, ({\mathcal C}^{\rm op}) \cong {\rm Ind} \, ({\mathcal D}^{\rm op}) \, .
$$
\end{remark}

\bigskip

\begin{demo}

En effet, toute \'equivalence de topos
$$
{\mathcal E}_{\mathbb T} \cong \widehat{\mathcal C}
$$
induit en particulier par d\'efinition des topos classifiants une \'equivalence de cat\'egories
$$
{\mathbb T}\mbox{\rm -mod} \, ({\rm Ens}) \cong {\rm pt} (\widehat{\mathcal C}) \, .
$$
Or, d'apr\`es le th\'eor\`eme \ref{thmIV38}, la cat\'egorie des points du topos $\widehat{\mathcal C}$
$$
{\rm pt} (\widehat{\mathcal C})
$$
est canoniquement \'equivalente \`a la cat\'egorie
$$
{\rm Ind} \, ({\mathcal C}^{\rm op}) \, .
$$
\end{demo}

\bigskip

Afin de pr\'eciser dans quelle mesure exactement une th\'eorie de type pr\'efaisceau ${\mathbb T}$ d\'etermine les cat\'egories essentiellement petites ${\mathcal C}$ telles que
$$
{\mathcal E}_{\mathbb T} \cong \widehat{\mathcal C} \, ,
$$
introduisons la d\'efinition suivante:

\begin{defn}\label{defV95}

Soit ${\mathcal C}$ une cat\'egorie localement petite qui poss\`ede des colimites filtrantes arbitraires.

\smallskip

Alors un objet $X$ de ${\mathcal C}$ est dit ``compact'' si le foncteur
$$
{\rm Hom} (X,\bullet) : {\mathcal C} \longrightarrow {\rm Ens}
$$
respecte les colimites filtrantes.
\end{defn}

\begin{remarkqed}

Ainsi, un objet $X$ d'une telle cat\'egorie ${\mathcal C}$ est compact si, pour tout diagramme $X_{\bullet}$ de ${\mathcal C}$ index\'e par une petite cat\'egorie filtrante ${\mathcal I}$, tout morphisme
$$
X \longrightarrow \varinjlim_{i \in {\mathcal I}} X_i
$$
se factorise \`a travers l'un des $X_i$. \end{remarkqed}

\bigskip

Cette d\'efinition s'applique aux cat\'egories de mod\`eles des th\'eories g\'eom\'etriques du fait du lemme suivant:

\begin{lem}\label{lemV96}

Soit ${\mathbb T}$ une th\'eorie g\'eom\'etrique de signature $\Sigma$.

\smallskip

Soit ${\mathcal E}$ un topos.

\smallskip

Alors:

\begin{listeimarge}

\item La cat\'egorie $\Sigma\mbox{\rm -str} \, ({\mathcal E})$ des $\Sigma$-structures dans ${\mathcal E}$ poss\`ede des limites et des colimites filtrantes arbitraires.

\medskip

\item Les foncteurs de colimites filtrantes respectent les interpr\'etations des formules g\'eom\'etriques de $\Sigma$.

\medskip

\item La cat\'egorie ${\mathbb T}\mbox{\rm -mod} \, ({\mathcal E})$ des mod\`eles de ${\mathbb T}$ dans ${\mathcal E}$ poss\`ede des colimites filtrantes arbitraires.
\end{listeimarge}
\end{lem}

\begin{demo}
\begin{listeisansmarge}
\item La notion de $\Sigma$-structure dans un topos ${\mathcal E}$ utilise seulement la structure cat\'egorique de ${\mathcal E}$ et la notion de sous-objet.

\smallskip

Comme, dans un topos, les colimites filtrantes de sous-objets sont des sous-objets (puisque les foncteurs de colimites filtrantes sont exacts) et que les limites de sous-objets sont des sous-objets dans n'importe quelle cat\'egorie, on voit que $\Sigma\mbox{-str} \, ({\mathcal E})$ a des limites et des colimites filtrantes arbitraires.

\medskip

\item L'interpr\'etation de toute formule g\'eom\'etrique dans une $\Sigma$-structure $M$ de ${\mathcal E}$ se formule enti\`erement en termes de limites finies et de colimites.

\smallskip

Or, dans un topos, les foncteurs de colimites filtrantes respectent non seulement les colimites mais aussi les limites finies.

\smallskip

Donc ils respectent aussi les interpr\'etations des formules g\'eom\'etriques.

\medskip

\item est cons\'equence imm\'ediate de (ii). \end{listeisansmarge}\end{demo}

\bigskip

Ainsi, les cat\'egories de mod\`eles des th\'eories g\'eom\'etriques ont toujours des colimites filtrantes arbitraires et on peut poser:

\begin{defn}\label{defV97}

Soit ${\mathbb T}$ une th\'eorie g\'eom\'etrique.

\smallskip

Un mod\`ele ensembliste $M$ de ${\mathbb T}$ est dit ``finiment pr\'esentable'' s'il est compact en tant qu'objet de la cat\'egorie
$$
{\mathbb T}\mbox{\rm -mod}  ({\rm Ens})
$$
des mod\`eles de ${\mathbb T}$ dans ${\rm Ens}$.

\smallskip

On note
$$
{\mathbb T}\mbox{\rm -mod}  ({\rm Ens})_{\rm fp}
$$

la sous-cat\'egorie pleine de ${\mathbb T}\mbox{\rm -mod}  ({\rm Ens})$ constitu\'ee des mod\`eles finiment pr\'esentables. 
\end{defn}

Cette d\'efinition \'etant pos\'ee, on a:

\begin{lem}\label{lemV98}

Soit ${\mathcal C}$ une cat\'egorie essentiellement petite.

\smallskip

Alors les objets compacts de la ${\rm ind}$-compl\'etion de ${\mathcal C}$
$$
{\rm Ind} \, ({\mathcal C})
$$
sont les pr\'efaisceaux sur ${\mathcal C}$
$$
R : {\mathcal C}^{\rm op} \longrightarrow {\rm Ens}
$$
qui sont des r\'etractes d'objets repr\'esentables, c'est-\`a-dire sont tels qu'existent un objet $X$ de ${\mathcal C}$ et deux morphismes
$$
i : R \xhookrightarrow{ \ { \ } \ } {\rm Hom} (\bullet , X) \qquad \mbox{et} \qquad r : {\rm Hom} (\bullet , X) \longrightarrow R
$$
v\'erifiant
$$
r \circ i = {\rm id}_R \, .
$$
\end{lem}

\begin{remarks}
\begin{listeisansmarge}

\item Pour toute telle cat\'egorie essentiellement petite ${\mathcal C}$, la sous-cat\'egorie pleine de 
$$
{\rm Ind} \, ({\mathcal C}) \qquad \mbox{ou} \qquad \widehat{\mathcal C}
$$
constitu\'ee des r\'etractes d'objets repr\'esentables est appel\'ee la ``compl\'etion Karoubienne'' de ${\mathcal C}$ et peut \^etre not\'ee
$$
{\rm Kar} \, ({\mathcal C}) \, .
$$

\item Pour tout objet $R$ de $\widehat{\mathcal C}$ reli\'e \`a un objet repr\'esentable ${\rm Hom} (\bullet ,X)$ par deux morphismes
$$
i : R \xhookrightarrow{ \ { \ } \ } {\rm Hom} (\bullet , X) \qquad \mbox{et} \qquad r : {\rm Hom} (\bullet , X) \longrightarrow R
$$
v\'erifiant
$$
r \circ i = {\rm id}_R \, ,
$$
le compos\'e
$$
p = i \circ r : {\rm Hom} (\bullet , X) \longrightarrow {\rm Hom} (\bullet , X)
$$
est idempotent au sens que
$$
p \circ p = p \, .
$$
De plus, pour tout objet $P$ de $\widehat{\mathcal C}$, l'ensemble de morphismes de $\widehat{\mathcal C}$
$$
{\rm Hom} (R,P)
$$
s'identifie au sous-ensemble de
$$
P(X)
$$
constitu\'e des morphismes
$$
x : {\rm Hom} (\bullet , X) \longrightarrow P
$$
tels que
$$
x \circ p = p \, .
$$

\item Il r\'esulte de la remarque (ii) que le foncteur pleinement fid\`ele
$$
{\mathcal C} \xhookrightarrow{ \ { \ } \ } {\rm Kar} ({\mathcal C})
$$
induit une \'equivalence de cat\'egories de pr\'efaisceaux
$$
\widehat{{\rm Kar} ({\mathcal C})} \xrightarrow{ \ \sim \ } \widehat{\mathcal C} \, .
$$

\item On dit que ${\mathcal C}$ est Karoubi-compl\`ete si tout r\'etracte d'un objet repr\'esentable de $\widehat{\mathcal C}$ est repr\'esentable, c'est-\`a-dire si
$$
{\mathcal C} \xhookrightarrow{ \ { \ } \ } {\rm Kar} ({\mathcal C})
$$
est une \'equivalence de cat\'egories, ou encore si tout endomorphisme de ${\mathcal C}$
$$
p : X \longrightarrow X
$$
idempotent au sens que
$$
p \circ p = p \, ,
$$
est scind\'e au sens qu'existent un objet $R$ de ${\mathcal C}$ et deux morphismes
$$
i : R \xhookrightarrow{ \ { \ } \ } X \qquad \mbox{et} \qquad r : X \longrightarrow R
$$
v\'erifiant
$$
r \circ i = {\rm id}_R \qquad \mbox{et} \qquad i \circ r = p \, .
$$

Avec cette d\'efinition, il r\'esulte de la remarque (iii) que la compl\'etion Karoubienne d'une cat\'egorie essentiellement petite ${\mathcal C}$
$$
{\rm Kar} ({\mathcal C})
$$
est toujours Karoubi-compl\`ete.
\end{listeisansmarge}
\end{remarks}

\begin{demo}

Un objet $R$ de ${\rm Ind} ({\mathcal C})$ est par d\'efinition un objet de $\widehat{\mathcal C}$ qui s'\'ecrit comme une colimite d'objets repr\'esentables ${\rm Hom} (\bullet , X_i)$
$$
R = \varinjlim_{\mathcal I} {\rm Hom} (\bullet , X_i)
$$
index\'ee par une petite cat\'egorie filtrante ${\mathcal I}$.

\smallskip

Si $R$ est un objet compact de ${\rm Ind} ({\mathcal C})$, le morphisme
$$
{\rm id} : R \longrightarrow R = \varinjlim {\rm Hom} (\bullet , X_i)
$$
se factorise pour un certain objet $i_0$ de ${\mathcal I}$ en
$$
R \xrightarrow{ \ i \ } {\rm Hom} (\bullet , X_{i_0}) \xrightarrow{ \ r \ } R \, .
$$

Cela signifie que l'objet $R$ de $\widehat{\mathcal C}$ est un r\'etracte d'un objet repr\'esentable.

\smallskip

R\'eciproquement, consid\'erons un objet $R$ de $\widehat{\mathcal C}$ qui est un r\'etracte d'un objet $X$ de ${\mathcal C}$, c'est-\`a-dire tel que ${\rm id}_R$ se factorise en
$$
R \xrightarrow{ \ i \ } {\rm Hom} (\bullet , X) \xrightarrow{ \ r \ } R \, .
$$
Posant
$$
p = i \circ r : {\rm Hom} (\bullet , X) \longrightarrow {\rm Hom} (\bullet , X) \, ,
$$
on a
$$
p \circ p = p
$$
et donc la cat\'egorie ${\mathcal P}$ constitu\'ee de l'unique objet ${\rm Hom} (\bullet , X)$ et de ses deux endomorphismes ${\rm id}$ et $p$ est filtrante.

\smallskip

Alors $R$ est la colimite du foncteur de plongement
$$
{\mathcal P} \xhookrightarrow{ \ { \ } \ } \widehat{\mathcal C}
$$
\`a valeurs dans les objets repr\'esentables, ce qui montre que $R$ est un objet de ${\rm Ind} ({\mathcal C})$.

\smallskip

Enfin, pour toute colimite filtrante d'objets de ${\rm Ind} ({\mathcal C})$
$$
P = \varinjlim_{\mathcal I} P_i \, ,
$$
la formule
$$
P(X) = \varinjlim_{\mathcal I} P_i (X)
$$
entra{\^\i}ne
$$
{\rm Hom} (R,P) = \varinjlim_{\mathcal I} {\rm Hom} (R,P_i)
$$
puisque, pour tout objet $Q$ de $\widehat{\mathcal C}$, on a
$$
{\rm Hom} (R,Q) = \{ x \in Q(X) = {\rm Hom} ({\rm Hom} (\bullet , X) , Q) \mid x \circ p = x \} \, .
$$

Cela prouve comme annonc\'e que $R$ est un objet compact de ${\rm Ind} ({\mathcal C})$. 
\end{demo}

\bigskip

On d\'eduit de ce lemme et de la proposition \ref{propV94}:

\begin{cor}\label{corV99}

Soit ${\mathbb T}$ une th\'eorie de type pr\'efaisceau.

\smallskip

Soit ${\mathcal C}$ une cat\'egorie essentiellement petite munie d'une \'equivalence de topos
$$
{\mathcal E}_{\mathbb T} \xrightarrow{ \ \sim \ } \widehat{\mathcal C} \, .
$$

Alors:

\begin{listeimarge}

\item L'\'equivalence de cat\'egories induite
$$
{\mathbb T}\mbox{\rm -mod} ({\rm Ens}) \xrightarrow{ \ \sim \ } {\rm Ind} ({\mathcal C}^{\rm op})
$$
se restreint en une \'equivalence
$$
{\mathbb T}\mbox{\rm -mod} ({\rm Ens})_{\rm fp} \xrightarrow{ \ \sim \ } {\rm Kar} ({\mathcal C}^{\rm op})
$$
de la cat\'egorie des mod\`eles ensemblistes finiment pr\'esentables de ${\mathbb T}$ sur la compl\'etion Karoubienne ${\rm Kar} ({\mathcal C}^{\rm op})$ de ${\mathcal C}^{\rm op}$.

\medskip

\item Notant
$$
{\mathcal M} = {\mathbb T}\mbox{\rm -mod} ({\rm Ens})_{\rm fp} \, ,
$$
l'\'equivalence de topos
$$
{\mathcal E}_{\mathbb T} \xrightarrow{ \ \sim \ } \widehat{\mathcal C}
$$
se rel\`eve en une \'equivalence
$$
{\mathcal E}_{\mathbb T} \xrightarrow{ \ \sim \ } \widehat{{\mathcal M}^{\rm op}} = [{\mathcal M} , {\rm Ens}] \, .
$$
\end{listeimarge}
\end{cor}

\begin{remarks}
\begin{listeisansmarge}

\item On d\'eduit de (i) que si
$$
{\mathcal E}_{\mathbb T} \cong \widehat{\mathcal C} \, ,
$$
alors une cat\'egorie essentiellement petite ${\mathcal D}$ admet une \'equivalence
$$
{\mathcal E}_{\mathbb T} \cong \widehat{\mathcal D}
$$
si et seulement si les compl\'etions Karoubiennes
$$
{\rm Kar} ({\mathcal C}^{\rm op}) = {\rm Kar} ({\mathcal C})^{\rm op} \qquad \mbox{et} \qquad {\rm Kar} ({\mathcal D}^{\rm op}) = {\rm Kar} ({\mathcal D})^{\rm op}
$$
sont \'equivalentes.

\medskip

\item Il r\'esulte de (ii) que le foncteur des mod\`eles d'une th\'eorie de type pr\'efaisceau ${\mathbb T}$ dans les topos ${\mathcal E}$
$$
{\mathcal E} \longmapsto {\mathbb T}\mbox{-mod} ({\mathcal E})
$$
est d\'etermin\'e par la cat\'egorie des mod\`eles ensemblistes finiment pr\'esentables de ${\mathbb T}$.
\end{listeisansmarge}
\end{remarks}

\bigskip

\begin{demo}
\begin{listeisansmarge}
\item On calcule l'invariant de topos
$$
{\mathcal E} \longmapsto \{\mbox{cat\'egorie des objets compacts de ${\rm pt} ({\mathcal E})$}\}
$$
\`a travers les deux pr\'esentations du m\^eme topos
$$
\xymatrix {
	& & {\mathcal E}_{\mathbb T} \cong \widehat{\mathcal C} \ar@/^8pt/@{--}[drr] & & \\
	\mbox{${\mathbb T}$} \ar@/^8pt/@{--}[urr] & & & & \mbox{${\mathcal C}$}}
$$
comme topos classifiant de ${\mathbb T}$ et topos des pr\'efaisceaux sur ${\mathcal C}$.

\smallskip

On obtient ainsi une \'equivalence
$$
{\mathbb T}\mbox{\rm -mod} ({\rm Ens})_{\rm fp} \xrightarrow{ \ \sim \ } {\rm Kar} ({\mathcal C}^{\rm op}) \, .
$$

\item r\'esulte de (i) et de la remarque (iii) qui suit le lemme \ref{lemV98}. \end{listeisansmarge}\end{demo}

\subsection{Une caract\'erisation syntactique des th\'eories de type pr\'efaisceau}\label{subsec593}

\medskip

Nous allons donner une premi\`ere caract\'erisation, de nature purement syntactique, des th\'eories de type pr\'efaisceau.

\smallskip

Dans ce but, introduisons la d\'efinition suivante:

\begin{defn}\label{defV910}
\begin{listeimarge}
\item Un objet $E$ d'un topos ${\mathcal E}$ est dit ``irr\'eductible'' si toute famille globalement \'epimorphique de morphismes
$$
E_i \xrightarrow{ \ e_i \ } E \, , \quad i \in I \, ,
$$
admet une section
$$
E \xrightarrow{ \ s \ } E_{i_0} \, , \quad \mbox{avec} \quad e_{i_0} \circ s = {\rm id}_E \, ,
$$
de l'un de ses \'el\'ements.

\medskip

\item Un objet $X$ d'une cat\'egorie essentiellement petite ${\mathcal C}$ munie d'une topologie $J$ est dit ``$J$-irr\'eductible'' si son seul crible $J$-couvrant est le crible maximal.

\medskip

\item Si ${\mathbb T}$ est une th\'eorie g\'eom\'etrique de signature $\Sigma$, une formule g\'eom\'etrique de $\Sigma$ dans un contexte $\vec x$
$$
\varphi (\vec x)
$$
est dite ``$ \ {\mathbb T}$-irr\'eductible'' si elle est $J_{\mathbb T}$-irr\'eductible en tant qu'objet de la cat\'egorie syntactique g\'eom\'etrique de ${\mathbb T}$
$$
{\mathcal C}_{\mathbb T} = {\mathcal C}_{\mathbb T}^{\rm geo}
$$
munie de sa topologie g\'eom\'etrique $J_{\mathbb T}$.
\end{listeimarge}
\end{defn}

\begin{remarksqed}
\begin{listeisansmarge}
\item Pour tout site $({\mathcal C},J)$, le foncteur canonique
$$
\ell : {\mathcal C} \xhookrightarrow{ \ { \ } \ } \widehat{\mathcal C} \longrightarrow \widehat{\mathcal C}_J
$$
transforme tout objet $J$-irr\'eductible de ${\mathcal C}$ en un objet irr\'eductible du topos $\widehat{\mathcal C}_J$.

\smallskip

R\'eciproquement, si $J$ est une topologie sous-canonique, tout objet de ${\mathcal C}$ que $\ell$ transforme en objet irr\'eductible du topos $\widehat{\mathcal C}_J$ est $J$-irr\'eductible comme objet de ${\mathcal C}$.

\medskip

\item Pour toute th\'eorie g\'eom\'etrique ${\mathbb T}$, la topologie $J_{\mathbb T}$ de ${\mathcal C}_{\mathbb T}$ est sous-canonique et donc une formule g\'eom\'etrique
$$
\varphi (\vec x)
$$
est ${\mathbb T}$-irr\'eductible si et seulement si le foncteur canonique
$$
\ell : {\mathcal C}_{\mathbb T} \longrightarrow {\mathcal E}_{\mathbb T} = \widehat{({\mathcal C}_{\mathbb T})_{J_{\mathbb T}}}
$$
la transforme en un objet irr\'eductible du topos classifiant ${\mathcal E}_{\mathbb T}$ de ${\mathbb T}$.

\medskip

\item Par d\'efinition de la structure cat\'egorique de ${\mathcal C}_{\mathbb T}$ et de sa topologie g\'eom\'etrique $J_{\mathbb T}$, une formule g\'eom\'etrique de la signature $\Sigma$ de ${\mathbb T}$
$$
\varphi (\vec x)
$$
est ${\mathbb T}$-irr\'eductible si et seulement si, pour toute famille de formules d\'emontrablement fonctionnelles
$$
\theta_i (\vec x_i , \vec x) : \varphi_i (\vec x_i) \longrightarrow \varphi (\vec x) \, , \quad i \in I \, ,
$$
telle que le s\'equent g\'eom\'etrique
$$
\varphi \vdash_{\vec x} \bigvee_{i \in I} (\exists \, \vec x_i) \, \theta_i (\vec x_i , \vec x)
$$
soit ${\mathbb T}$-d\'emontrable, il existe un indice $i_0 \in I$ et une formule d\'emontrablement fonctionnelle
$$
\theta' (\vec x , \vec x_{i_0}) : \varphi (\vec x) \longrightarrow \varphi_{i_0} (\vec x_{i_0})
$$
telle que la formule de contexte $(\vec x , \vec x')$
$$
(\exists \, \vec x_{i_0}) (\theta (\vec x_{i_0} , \vec x) \wedge \theta' (\vec x' , \vec x_{i_0}))
$$
soit ${\mathbb T}$-d\'emontrablement \'equivalente \`a la formule
$$
\vec x = \vec x' \, .
$$
\end{listeisansmarge}
\end{remarksqed}

Cette d\'efinition \'etant pos\'ee, on peut \'enoncer le crit\`ere suivant:

\begin{thm}[Caramello (Corollaire 6.1.10 \cite{TST})]\label{thmV911}

Soit ${\mathbb T}$ une th\'eorie g\'eom\'etrique de signature $\Sigma$.

\smallskip

Soit ${\mathcal C}_{\mathbb T}$ sa cat\'egorie syntactique g\'eom\'etrique munie de sa topologie g\'eom\'etrique $J_{\mathbb T}$.

\smallskip

Alors la th\'eorie ${\mathbb T}$ est de type pr\'efaisceau si et seulement si toute formule g\'eom\'etrique de $\Sigma$ dans un contexte $\vec x$
$$
\varphi (\vec x) \, ,
$$
vue comme un objet de ${\mathcal C}_{\mathbb T}$, admet une famille $J_{\mathbb T}$-couvrante
$$
\varphi_i (\vec x_i) \longrightarrow \varphi (\vec x)
$$
constitu\'ee de formules $J_{\mathbb T}$-irr\'eductibles $\varphi_i (\vec x_i)$.
\end{thm}

\begin{demosansqed}

Supposons d'abord que ${\mathcal C}_{\mathbb T}$ v\'erifie cette propri\'et\'e.

\smallskip

Elle signifie que la sous-cat\'egorie pleine ${\mathcal C}_{\mathbb T}^{\rm ir}$ de ${\mathcal C}_{\mathbb T}$ constitu\'ee des formules $J_{\mathbb T}$-irr\'eductibles est dense pour la topologie $J_{\mathbb T}$.

\smallskip

Or, par d\'efinition de la notion de formule irr\'eductible, $J_{\mathbb T}$ induit sur ${\mathcal C}_{\mathbb T}^{\rm ir}$ la topologie discr\`ete.

\smallskip

Il r\'esulte alors du lemme de comparaison de Grothendieck que le foncteur de restriction
$$
{\mathcal E}_{\mathbb T} = \widehat{({\mathcal C}_{\mathbb T})}_{J_{\mathbb T}} \longrightarrow \widehat{{\mathcal C}_{\mathbb T}^{\rm ir}}
$$
est une \'equivalence de topos.

\smallskip

Ainsi, ${\mathbb T}$ est de type pr\'efaisceau.

\smallskip

Pour l'implication en sens inverse, on a besoin du lemme suivant:
\end{demosansqed}

\begin{lem}\label{lemV912}

Soit $({\mathcal C} , J)$ un site.

\smallskip

Alors:

\begin{listeimarge}

\item Tout objet irr\'eductible du topos $\widehat{\mathcal C}_J$ est un r\'etracte de l'image par le foncteur canonique
$$
\ell : {\mathcal C} \xhookrightarrow{ \ { \ } \ } \widehat{\mathcal C} \longrightarrow \widehat{\mathcal C}_J
$$
d'un objet de ${\mathcal C}$.

\medskip

\item Si la topologie $J$ est sous-canonique et si tout endomorphisme idempotent de ${\mathcal C}$ est scind\'e, alors le foncteur canonique
$$
\ell : {\mathcal C} \longrightarrow \widehat{\mathcal C}_J
$$
induit une \'equivalence de la sous-cat\'egorie pleine de ${\mathcal C}$ constitu\'ee des objets $J$-irr\'eductibles sur la sous-cat\'egorie pleine du topos $\widehat{\mathcal C}_J$ constitu\'ee des objets irr\'eductibles.
\end{listeimarge}
\end{lem}

\begin{remark}

Si ${\mathcal C}$ est une cat\'egorie cart\'esienne, elle v\'erifie a fortiori la condition que tout endomorphisme idempotent est scind\'e.
\end{remark}

\bigskip

\begin{demo}
\begin{listeisansmarge}
\item r\'esulte de ce que tout objet $E$ du topos $\widehat{\mathcal C}_J$ admet une famille globalement \'epimorphique
$$
\ell (X_i) \longrightarrow E \, , \quad i \in I \, ,
$$
compos\'ee d'images par $\ell : {\mathcal C} \to \widehat{\mathcal C}_J$ d'objet $X_i$ de la cat\'egorie ${\mathcal C}$.

\medskip

\item r\'esulte de (i) et de la remarque (i) qui suit la d\'efinition \ref{defV910}. \hfill $\Box$
\end{listeisansmarge}
\end{demo}

\bigskip

\noindent {\bf Fin de la d\'emonstration du th\'eor\`eme \ref{thmV911}:}

\smallskip

Supposons que ${\mathbb T}$ est de type pr\'efaisceau c'est-\`a-dire admet une \'equivalence
$$
{\mathcal E}_{\mathbb T} \cong \widehat{\mathcal C}
$$
pour une certaine cat\'egorie essentiellement petite ${\mathcal C}$.

\smallskip

D'apr\`es la remarque (iii) qui suit le lemme \ref{lemV98}, on peut supposer que ${\mathcal C}$ est Karoubi-compl\`ete, c'est-\`a-dire que tout endomorphisme idempotent de ${\mathcal C}$ est scind\'e.

\smallskip

Calculons alors l'invariant de topos
$$
{\mathcal E} \longmapsto \{\mbox{cat\'egorie des objets irr\'eductibles de ${\mathcal E}$}\}
$$
\`a travers la double pr\'esentation du topos classifiant de ${\mathbb T}$
$$
\xymatrix {
	& & {\mathcal E}_{\mathbb T} \cong \widehat{\mathcal C} \ar@/^8pt/@{--}[drr] & & \\
	\mbox{$({\mathcal C}_{\mathbb T} , J_{\mathbb T})$} \ar@/^8pt/@{--}[urr] & & & & \mbox{$\ \ {\mathcal C} \ \ $}}
$$
comme topos des faisceaux sur le site sous-canonique $(J_{\mathbb T} , {\mathcal C}_{\mathbb T})$ et comme topos des pr\'efaisceaux sur ${\mathcal C}$.

\smallskip

On d\'eduit de ce pont une \'equivalence de cat\'egories
$$
{\mathcal C}_{\mathbb T}^{\rm ir} \cong {\mathcal C}
$$
entre la cat\'egorie ${\mathcal C}_{\mathbb T}^{\rm ir}$ des objets $J_{\mathbb T}$-irr\'eductibles de ${\mathcal C}_{\mathbb T}$ et la cat\'egorie ${\mathcal C}$.

\smallskip

A fortiori, la sous-cat\'egorie pleine ${\mathcal C}_{\mathbb T}^{\rm ir}$ de ${\mathcal C}_{\mathbb T}$ est $J_{\mathbb T}$-dense, comme on voulait. \hfill $\Box$

\bigskip

On d\'eduit de ce th\'eor\`eme et du lemme:

\begin{cor}\label{corV913}

Soit ${\mathbb T}$ une th\'eorie g\'eom\'etrique de type pr\'efaisceau.

\smallskip

Soient ${\mathcal C}_{\mathbb T} = {\mathcal C}_{\mathbb T}^{\rm geo}$ sa cat\'egorie syntactique g\'eom\'etrique et ${\mathcal C}_{\mathbb T}^{\rm ir}$ la sous-cat\'egorie pleine de ${\mathcal C}_{\mathbb T}$ constitu\'ee des formules g\'eom\'etriques ${\mathbb T}$-irr\'eductibles.

\smallskip

Soit d'autre part
$$
{\mathcal M} = {\mathbb T}\mbox{\rm -mod} \, ({\rm Ens})_{\rm fp}
$$
la cat\'egorie des mod\`eles ensemblistes finiment pr\'esentables de ${\mathbb T}$.

\smallskip

Alors:

\begin{listeimarge}

\item Les cat\'egories
$$
{\mathcal C}_{\mathbb T}^{\rm ir} \qquad \mbox{et} \qquad {\mathcal M}^{\rm op}
$$
s'identifient toutes deux \`a la cat\'egorie des objets irr\'eductibles du topos classifiant ${\mathcal E}_{\mathbb T}$ de ${\mathbb T}$.

\medskip

\item Elles sont reli\'ees par une \'equivalence
$$
{\mathcal C}_{\mathbb T}^{\rm ir} \xrightarrow{ \ \sim \ } {\mathcal M}^{\rm op}
$$
qui associe \`a toute formule g\'eom\'etrique ${\mathbb T}$-irr\'eductible dans un contexte $\vec x = (x_1^{A_1} \cdots x_n^{A_n})$
$$
\varphi (\vec x)
$$
un mod\`ele ensembliste
$$
M_{\varphi}
$$
caract\'eris\'e par la propri\'et\'e que, pour tout mod\`ele ensembliste $M$ de ${\mathbb T}$, se donner un morphisme
$$
M_{\varphi} \longrightarrow M
$$
\'equivaut \`a se donner une famille d'\'el\'ements
$$
(m_1 , \cdots , m_n) \in M\!A_1 \times \cdots \times M\!A_n
$$
qui v\'erifie la formule $\varphi$ au sens que
$$
(m_1 , \cdots , m_n) \in M \varphi (\vec x) \, .
$$
\end{listeimarge}
\end{cor}

\begin{remark}

On dit qu'un mod\`ele $M_{\varphi}$ qui v\'erifie la condition de (ii) relativement \`a une formule g\'eom\'etrique $\varphi$ est ``finiment pr\'esent\'e'' par $\varphi$.

\smallskip

Ainsi, tout mod\`ele ensembliste finiment pr\'esentable d'une th\'eorie de type pr\'efaisceau ${\mathbb T}$ est finiment pr\'esent\'e par une formule g\'eom\'etrique ${\mathbb T}$-irr\'eductible.
\end{remark}

\bigskip

\begin{demosansqed}
\begin{listeisansmarge}
\item Cela r\'esulte du lemme \ref{lemV912} et de la double pr\'esentation
$$
\widehat{({\mathcal C}_{\mathbb T})}_{J_{\mathbb T}} \cong {\mathcal E}_{\mathbb T} \cong \widehat{{\mathcal M}^{\rm op}}
$$
puisque la topologie $J_{\mathbb T}$ est sous-canonique, que la cat\'egorie ${\mathcal C}_{\mathbb T}$ est cart\'esienne et que la cat\'egorie ${\mathcal M}^{\rm op}$ est Karoubi-compl\`ete.

\medskip

\item Les deux identifications de (i) d\'efinissent une \'equivalence de cat\'egories
$$
{\mathcal C}_{\mathbb T}^{\rm ir} \xrightarrow{ \ \sim \ } {\mathcal M}^{\rm op}
$$
qui fait commuter \`a isomorphisme pr\`es le carr\'e:
$$
\xymatrix{
{\mathcal C}_{\mathbb T}^{\rm ir} \, \ar[d]^{\wr} \ar@{^{(}->}[r] &{\mathcal C}_{\mathbb T} \ar[r]^{\ell} &{\mathcal E}_{\mathbb T} \ar[d]^-{\wr} \\
{\mathcal M}^{\rm op} \, \ar@{^{(}->}[rr]^y &&\widehat{{\mathcal M}^{\rm op}}
}
$$

Or on a:
\end{listeisansmarge}
\end{demosansqed}

\begin{lem}\label{lemV914}

Soient ${\mathbb T}$ une th\'eorie de type pr\'efaisceau de signature $\Sigma$, ${\mathcal C}_{\mathbb T} = {\mathcal C}_{\mathbb T}^{\rm geo}$ sa cat\'egorie syntactique g\'eom\'etrique munie de son mod\`ele canonique $M_{\mathbb T}$ et
$$
\ell : {\mathcal C}_{\mathbb T} \longrightarrow {\mathcal E}_{\mathbb T}
$$
le foncteur canonique de ${\mathcal C}_{\mathbb T}$ dans le topos classifiant ${\mathcal E}_{\mathbb T}$ de ${\mathbb T}$ muni de son mod\`ele universel $U_{\mathbb T} = \ell (M_{\mathbb T})$.

\smallskip

Soit d'autre part
$$
{\mathcal M} = {\mathbb T}\mbox{\rm -mod} \, ({\rm Ens})_{\rm fp}
$$
la cat\'egorie des mod\`eles ensemblistes finiment pr\'esentables de ${\mathbb T}$.

\smallskip

Alors:

\begin{listeimarge}

\item L'image du mod\`ele universel $U_{\mathbb T}$ de ${\mathbb T}$ dans ${\mathcal E}_{\mathbb T}$ par l'\'equivalence
$$
{\mathcal E}_{\mathbb T} \xrightarrow{ \ \sim \ } \widehat{{\mathcal M}^{\rm op}}
$$
est le mod\`ele $U'_{\mathbb T}$ de ${\mathbb T}$ dans $\widehat{{\mathcal M}^{\rm op}}$ qui associe

\medskip

$
\left\{\begin{matrix}
\bullet &\mbox{\`a toute sorte $A$ de $\Sigma$ le pr\'efaisceau} \hfill \\
&\begin{matrix}
{\mathcal M} &\longrightarrow &{\rm Ens} \, , \\
M &\longmapsto &M\!A \, ,
\end{matrix} \\
{ \ } \\
\bullet &\mbox{\`a tout symbole de fonction $f : A_1 \cdots A_n \to B$ de $\Sigma$ le morphisme de pr\'efaisceaux constitu\'e des} \hfill \\
{ \ } \\
&M\!A_1 \times \cdots \times M\!A_n \xrightarrow{ \, M\!f \ } M\!B \, , \quad M \in {\rm Ob} ({\mathcal M}) \, , \\
{ \ } \\
\bullet &\mbox{\`a tout symbole de relation $R \rightarrowtail A_1 \cdots A_n$ de $\Sigma$ le sous-pr\'efaisceau constitu\'e des} \hfill \\
{ \ } \\
&(M\!R \xhookrightarrow{ \ { \ } \ } M\!A_1 \times \cdots \times M\!A_n) \, , \quad M \in {\rm Ob} ({\mathcal M}) \, .
\end{matrix} \right.
$

\bigskip

\item Pour toute famille g\'eom\'etrique $\varphi (\vec x)$ de $\Sigma$ de contexte $\vec x = (x_1^{A_1} \cdots x_n^{A_n})$, consid\'er\'ee comme un objet de ${\mathcal C}_{\mathbb T}$, son image par le foncteur compos\'e
$$
{\mathcal C}_{\mathbb T} \xrightarrow{ \ \ell \ } {\mathcal E}_{\mathbb T} \xrightarrow{ \ \sim \ } \widehat{{\mathcal M}^{\rm op}}
$$
s'identifie au sous-pr\'efaisceau
$$
(M\varphi (\vec x) \xhookrightarrow{ \ { \ } \ } M\!A_1 \times \cdots \times M\!A_n) \, , \quad M \in {\rm Ob} ({\mathcal M}) \, .
$$
\end{listeimarge}
\end{lem}

\bigskip

\begin{demolem}
\begin{listeisansmarge}

\item Tout objet $M$ de ${\mathcal M}^{\rm op}$ d\'efinit un point du topos $\widehat{{\mathcal M}^{\rm op}}$ c'est-\`a-dire un morphisme de topos
$$
{\rm Ens} \longrightarrow \widehat{{\mathcal M}^{\rm op}} \, .
$$
Via l'\'equivalence ${\mathcal E}_{\mathbb T} \xrightarrow{ \ \sim \ } \widehat{{\mathcal M}^{\rm op}}$, il s'identifie au point du topos classifiant
$$
{\rm Ens} \longrightarrow {\mathcal E}_{\mathbb T}
$$
qui correspond \`a $M$ vu comme un mod\`ele ensembliste de ${\mathbb T}$.

\smallskip

Donc l'\'evaluation du mod\`ele $U'_{\mathbb T}$ de ${\mathbb T}$ dans $\widehat{{\mathcal M}^{\rm op}}$ sur un tel objet $M$ de ${\mathcal M}$ s'identifie \`a $M$ vu comme un mod\`ele de ${\mathbb T}$ dans ${\rm Ens}$.

\smallskip

Cela prouve (i).

\medskip

\item r\'esulte de (i) puisque le sous-objet $\varphi (\vec x) \hookrightarrow \top (x_1^{A_1} \cdots x_n^{A_n})$ de ${\mathcal C}_{\mathbb T}$ s'identifie \`a l'interpr\'etation $M_{\mathbb T} \, \varphi (\vec x)$ de la formule g\'eom\'etrique $\varphi (\vec x)$ dans le mod\`ele canonique $M_{\mathbb T}$ de ${\mathbb T}$ dans ${\mathcal C}_{\mathbb T}$ et que le foncteur compos\'e
$$
{\mathcal C}_{\mathbb T} \xrightarrow{ \ \ell \ } {\mathcal E}_{\mathbb T} \xrightarrow{ \ \sim \ } \widehat{{\mathcal M}^{\rm op}}
$$
respecte les interpr\'etations des formules g\'eom\'etriques de $\Sigma$ dans les $\Sigma$-structures. 
\end{listeisansmarge}
\end{demolem}

\bigskip

\noindent {\bf Fin de la d\'emonstration du corollaire \ref{corV913} (ii):}

\smallskip

Il r\'esulte du lemme ci-dessus que, pour toute formule g\'eom\'etrique ${\mathbb T}$-irr\'eductible $\varphi (\vec x)$ de contexte $\vec x = (x_1^{A_1} \cdots x_n^{A_n})$, l'\'equivalence de cat\'egories
$$
{\mathcal C}_{\mathbb T}^{\rm ir} \xrightarrow{ \ \sim \ } {\mathcal M}^{\rm op}
$$
transforme $\varphi (\vec x)$ en un objet $M_{\varphi}$ qui repr\'esente le pr\'efaisceau
$$
\begin{matrix}
{\mathcal M} &\longrightarrow &{\rm Ens} \, , \hfill \\
M &\longmapsto &M\varphi (\vec x) \, .
\end{matrix}
$$

Autrement dit, $M_{\varphi}$ est caract\'eris\'e par la propri\'et\'e que, pour tout objet $M$ de ${\mathcal M}$
$$
{\rm Hom} (M_{\varphi} , M) = M\varphi (\vec x) \xhookrightarrow{ \ { \ } \ } M\!A_1 \times \cdots \times M\!A_n \, .
$$

Cette propri\'et\'e s'\'etend \`a tous les mod\`eles ensemblistes de ${\mathbb T}$ puisque tout tel mod\`ele s'\'ecrit comme une colimite filtrante d'objets $M$ de ${\mathcal M} = {\mathbb T}\mbox{-mod} \, ({\rm Ens})_{\rm fp}$ et que le foncteur
$$
{\rm Hom} (M_{\varphi} , \bullet) : {\mathbb T}\mbox{-mod} \, ({\rm Ens}) \longrightarrow {\rm Ens}
$$
respecte les colimites filtrantes, tout comme le foncteur
$$
\begin{matrix}
{\mathbb T}\mbox{-mod} \, ({\rm Ens}) &\longrightarrow &{\rm Ens} \, , \hfill \\
\hfill M &\longmapsto &M\varphi (\vec x) \, .
\end{matrix}
$$

Cela termine la d\'emonstration. \hfill $\Box$

\subsection{Une autre caract\'erisation liant syntaxe et s\'emantique ensembliste}\label{subsec594}

\medskip

Les th\'eories de type pr\'efaisceau sont \'egalement caract\'eris\'ees par des propri\'et\'es d'\'equivalence entre leur syntaxe et leur s\'emantique ensembliste:

\begin{thm}[Caramello (Théorème 7.1.6 \cite{TST})]\label{thmV915} 

Soit ${\mathbb T}$ une th\'eorie g\'eom\'etrique du premier ordre de signature $\Sigma$.

\smallskip

Alors ${\mathbb T}$ est de type pr\'efaisceau si et seulement si elle satisfait les trois conditions suivantes:

\begin{listeimarge}

\item[(1)] Pour tout s\'equent g\'eom\'etrique de $\Sigma$
$$
\varphi \vdash_{\vec x} \psi \, ,
$$
ce s\'equent est ${\mathbb T}$-d\'emontrable si et seulement si il est v\'erifi\'e par tous les mod\`eles ensemblistes finiment pr\'esentables de ${\mathbb T}$.

\medskip

\item[(2)] Tout mod\`ele ensembliste finiment pr\'esentable $M$ de ${\mathbb T}$ est pr\'esent\'e par une formule g\'eom\'etrique
$$
\varphi (\vec x) = \varphi (x_1^{A_1} \cdots x_n^{A_n}) \, ,
$$
au sens que pour tout mod\`ele ensembliste $N$ de $A$, l'ensemble de morphismes de mod\`eles
$$
{\rm Hom} (M,N)
$$
s'identifie au sous-ensemble
$$
N\varphi (\vec x) \xhookrightarrow{ \ { \ } \ } N\!A_1 \times \cdots \times N\!A_n \, .
$$

\item[(3)] Pour toute suite de sortes $A_1 , \cdots , A_n$ de $\Sigma$, et pour toute famille de sous-ensembles
$$
P_M \xhookrightarrow{ \ { \ } \ } M\!A_1 \times \cdots \times M\!A_n
$$
index\'es par les mod\`eles ensemblistes finiment pr\'esentables $M$ de ${\mathbb T}$ et telle que tout morphisme de mod\`eles
$$
M \longrightarrow N
$$
envoie $P_M$ dans $P_N$ par l'application induite
$$
M\!A_1 \times \cdots \times M\!A_n \longrightarrow N\!A_1 \times \cdots \times N\!A_n \, ,
$$
il existe une formule g\'eom\'etrique de $\Sigma$
$$
\varphi (\vec x) \quad \mbox{de contexte} \quad \vec x = (x_1^{A_1} \cdots x_n^{A_n})
$$
telle que, pour tout mod\`ele ensembliste finiment pr\'esentable $M$,
$$
P_M = M\varphi (\vec x) \xhookrightarrow{ \ { \ } \ } M\!A_1 \times \cdots \times M\!A_n \, .
$$\end{listeimarge}
\end{thm}

\begin{remarks}
\begin{listeisansmarge}
\item Pour (1), on dit que la collection des mod\`eles ensemblistes finiment pr\'esentables de ${\mathbb T}$ est compl\`ete.

\medskip

\item Pour (2), on dit que tout mod\`ele ensembliste finiment pr\'esentable de ${\mathbb T}$ est pr\'esent\'e par une formule g\'eom\'etrique, n\'ecessairement ${\mathbb T}$-irr\'eductible d'apr\`es le corollaire \ref{corV913} (ii).

\medskip

\item Pour (3), on dit que toute propri\'et\'e fonctorielle des mod\`eles ensemblistes finiment pr\'esentables de ${\mathbb T}$ est d\'efinie par une formule g\'eom\'etrique.

\medskip

\item R\'eciproquement, pour toute th\'eorie g\'eom\'etrique ${\mathbb T}$ de signature $\Sigma$, toute formule g\'eom\'etrique de $\Sigma$
$$
\varphi (\vec x) \quad \mbox{de contexte} \quad \vec x = (x_1^{A_1} \cdots x_n^{A_n})
$$
d\'efinit une propri\'et\'e fonctorielle des mod\`eles ensemblistes de ${\mathbb T}$
$$
M\longmapsto (M\varphi (\vec x) \hookrightarrow M\!A_1 \times \cdots \times M\!A_n) \, .
$$
\end{listeisansmarge}
\end{remarks}

\begin{demo}

Si ${\mathbb T}$ est une th\'eorie de type pr\'efaisceau, elle satisfait les conditions (1), (2) et (3) comme cons\'equence du pont
$$
\xymatrix {
	& & {\mathcal E}_{\mathbb T} \cong \widehat{{\mathcal M}^{\rm op}} \ar@/^8pt/@{--}[drr] & & \\
	\mbox{$({\mathcal C}_{\mathbb T} , J_{\mathbb T})$} \ar@/^8pt/@{--}[urr] & & & & \mbox{$\ \ {\mathcal M} \ \ $}}
$$
qui relie le site syntactique g\'eom\'etrique $({\mathcal C}_{\mathbb T} , J_{\mathbb T})$ de ${\mathbb T}$ \`a la cat\'egorie 
$$
{\mathcal M} = {\mathbb T}\mbox{-mod} \, ({\rm Ens})_{\rm fp}
$$
de ses mod\`eles ensemblistes finiment pr\'esentables.

\smallskip

Pour (1), consid\'erons en effet un s\'equent g\'eom\'etrique de $\Sigma$
$$
\varphi \vdash_{\vec x} \psi \quad \mbox{de contexte} \quad \vec x = (x_1^{A_1} \cdots x_n^{A_n}) \, .
$$

Notant $U_{\mathbb T}$ le mod\`ele universel de ${\mathbb T}$ dans le topos classifiant ${\mathcal E}_{\mathbb T}$ de ${\mathbb T}$, ce s\'equent est ${\mathbb T}$-d\'emontrable si et seulement si les deux sous-objets dans 
${\mathcal E}_{\mathbb T}$
$$
U_{\mathbb T} \, \varphi (\vec x) \xhookrightarrow{ \ { \ } \ } U_{\mathbb T} \, A_1 \times \cdots \times U_{\mathbb T} \, A_n
$$
et
$$
U_{\mathbb T} \, \psi (\vec x) \xhookrightarrow{ \ { \ } \ } U_{\mathbb T} \, A_1 \times \cdots \times U_{\mathbb T} \, A_n
$$
satisfont la relation d'inclusion
$$
U_{\mathbb T} \, \varphi (\vec x) \leq U_{\mathbb T} \, \psi (\vec x) \, .
$$

Or, dans le topos $\widehat{{\mathcal M}^{\rm op}}$, ces deux sous-objets correspondent aux deux sous-pr\'efaisceaux
$$
M \longmapsto M \varphi (\vec x)
$$
et
$$
M \longmapsto M \psi (\vec x)
$$
du pr\'efaisceau
$$
\begin{matrix}
{\mathcal M} &\longrightarrow &{\rm Ens} \, , \hfill \\
\hfill M &\longmapsto &M\!A_1 \times \cdots \times M\!A_n \, .
\end{matrix}
$$

On constate comme voulu que le s\'equent
$$
\varphi \vdash_{\vec x} \psi
$$
est ${\mathbb T}$-d\'emontrable si et seulement si
$$
M \varphi (\vec x) \leq M \psi (\vec x)
$$
pour tout mod\`ele ensembliste finiment pr\'esentable $M$.

\smallskip

Pour (2), c'est une forme faible du corollaire \ref{corV913} (ii).

\smallskip

Pour (3), consid\'erons les sous-objets de 
$$
U_{\mathbb T} \, A_1 \times \cdots \times U_{\mathbb T} \, A_n
$$
dans les deux repr\'esentations de ${\mathcal E}_{\mathbb T} \cong \widehat{{\mathcal M}^{\rm op}}$.

\smallskip

Dans la repr\'esentation de ${\mathcal E}_{\mathbb T}$ comme topos des pr\'efaisceaux sur ${\mathcal M}^{\rm op}$, $U_{\mathbb T} \, A_1 \times \cdots \times U_{\mathbb T} \, A_n$ est le pr\'efaisceau
$$
\begin{matrix}
{\mathcal M} &\longrightarrow &{\rm Ens} \, , \hfill \\
\hfill M &\longmapsto &M\!A_1 \times \cdots \times M\!A_n
\end{matrix}
$$
et ses sous-objets sont les familles fonctorielles de sous-ensembles
$$
P_M \xhookrightarrow{ \ { \ } \ } M\!A_1 \times \cdots \times M\!A_n \, , \quad M \in {\rm Ob} ({\mathcal M}) \, .
$$

Dans la repr\'esentation de ${\mathcal E}_{\mathbb T}$ comme topos des faisceaux sur $({\mathcal C}_{\mathbb T} , J_{\mathbb T})$, $U_{\mathbb T} \, A_1 \times \cdots \times U_{\mathbb T} \, A_n$ est l'image de
$$
M_{\mathbb T} \, A_1 \times \cdots \times M_{\mathbb T} \, A_n = \top (x_1^{A_1} \cdots x_n^{A_n})
$$
par le foncteur pleinement fid\`ele
$$
\ell : {\mathcal C}_{\mathbb T} \longrightarrow {\mathcal E}_{\mathbb T} \, .
$$

On sait que les sous-objets de $\top (x_1^{A_1} \cdots x_n^{A_n})$ dans ${\mathcal C}_{\mathbb T}$ sont les formules g\'eom\'etriques
$$
\varphi (x_1^{A_1} \cdots x_n^{A_n})
$$
consid\'er\'ees \`a \'equivalence d\'emontrable pr\`es. Elles sont transform\'ees par $\ell$ en des sous-objets de $U_{\mathbb T} \, A_1 \times \cdots \times U_{\mathbb T} \, A_n$ puisque le foncteur $\ell$ respecte les limites finies.

\smallskip

R\'eciproquement, tout sous-objet de $U_{\mathbb T} \, A_1 \times \cdots \times U_{\mathbb T} \, A_n$ a cette forme.

\smallskip

En effet, notant $\vec x = (x_1^{A_1} \cdots x_n^{A_n})$, il est le transform\'e par $\ell$ de la r\'eunion des images d'une famille de morphismes de ${\mathcal C}_{\mathbb T}$
$$
\theta_i (\vec x_i , \vec x) : \varphi_i (\vec x_i) \longrightarrow \varphi (\vec x) \, , \quad i \in I \, ,
$$
donc est le transform\'e par $\ell$ du sous-objet de $\varphi (\vec x)$ dans ${\mathcal C}_{\mathbb T}$ d\'efini par la formule g\'eom\'etrique
$$
\bigvee_{i \in I} (\exists \, \vec x_i) \, \theta_i (\vec x_i , \vec x) \, .
$$

Ainsi, les sous-objets de
$$
U_{\mathbb T} \, A_1 \times \cdots \times U_{\mathbb T} \, A_n
$$
s'\'ecrivent dans la repr\'esentation de ${\mathcal E}_{\mathbb T}$ \`a partir de $({\mathcal C}_{\mathbb T} , J_{\mathbb T})$ comme les transform\'es par le foncteur canonique
$$
\ell : {\mathcal C}_{\mathbb T} \longrightarrow {\mathcal E}_{\mathbb T}
$$
des formules g\'eom\'etriques
$$
\varphi (\vec x) = \varphi (x_1^{A_1} \cdots x_n^{A_n}) \, .
$$

En combinant les deux repr\'esentations, on obtient comme voulu que les propri\'et\'es fonctorielles
$$
M \longmapsto (P_M \xhookrightarrow{ \ { \ } \ } M\!A_1 \times \cdots \times M\!A_n) \, , \quad M \in {\rm Ob} ({\mathcal M}) \, ,
$$
correspondent aux formules g\'eom\'etriques
$$
\varphi (\vec x) = \varphi (x_1^{A_1} \cdots x_n^{A_n}) \, .
$$

Enfin, si $(M \mapsto P_M)$ et $\varphi (\vec x)$ se correspondent, on a d'apr\`es le lemme \ref{lemV914} (ii) la formule
$$
P_M = M\varphi (\vec x) \xhookrightarrow{ \ { \ } \ } M\!A_1 \times \cdots \times M\!A_n
$$
pour tout objet $M$ de ${\mathcal M} = {\mathbb T}\mbox{-mod} ({\rm Ens})_{\rm fp}$.

\smallskip

Cela ach\`eve de prouver (3).

\smallskip

R\'eciproquement, consid\'erons une th\'eorie g\'eom\'etrique ${\mathbb T}$ qui satisfait les conditions (1), (2) et (3).

\smallskip

Notons $({\mathcal C}_{\mathbb T} , J_{\mathbb T})$ le site syntactique g\'eom\'etrique de ${\mathbb T}$ et
$$
{\mathcal M} =  {\mathbb T}\mbox{-mod} ({\rm Ens})_{\rm fp}
$$
la cat\'egorie des mod\`eles ensemblistes finiment pr\'esentables de ${\mathbb T}$.

\smallskip

On dispose du foncteur des interpr\'etations
$$
{\mathcal C}_{\mathbb T} \longrightarrow \widehat{{\mathcal M}^{\rm op}} = [{\mathcal M} , {\rm Ens}]
$$
qui envoie

\medskip

$
\left\{\begin{matrix}
\bullet &\mbox{tout objet de ${\mathcal C}_{\mathbb T}$, c'est-\`a-dire toute formule g\'eom\'etrique} \hfill \\
{ \ } \\
&\varphi (\vec x) = \varphi (x_1^{A_1} \cdots x_n^{A_n}) \, , \\
{ \ } \\
&\mbox{sur le pr\'efaisceau des interpr\'etations} \hfill \\
{ \ } \\
&\begin{matrix}
I_{\varphi (x)} &: &{\mathcal M} &\longrightarrow &{\rm Ens} \, , \hfill \\
&&\hfill M &\longmapsto &M\varphi (\vec x) \, ,
\end{matrix} \\
{ \ } \\
\bullet &\mbox{tout morphisme de ${\mathcal C}_{\mathbb T}$, c'est-\`a-dire toute formule d\'emontrablement fonctionnelle} \hfill \\
{ \ } \\
&\theta (\vec x , \vec y) : \varphi (\vec x) \longrightarrow \psi (\vec y) \, , \\
{ \ } \\
&\mbox{sur le morphisme de pr\'efaisceaux constitu\'e des applications} \hfill \\
{ \ } \\
&M \varphi (\vec x) \longrightarrow M\psi (\vec y) \, , \quad M \in {\rm Ob} ({\mathcal M}) \, , \\
{ \ } \\
&\mbox{dont les graphes sont les interpr\'etations} \hfill \\
{ \ } \\
&M\theta (\vec x , \vec y) \xhookrightarrow{ \ { \ } \ } M\varphi (\vec x) \times M\psi (\vec y) \, .
\end{matrix} \right.
$

\bigskip

Il r\'esulte des conditions (1) et (3) que ce foncteur des interpr\'etations
$$
\begin{matrix}
\hfill {\mathcal C}_{\mathbb T} &\longrightarrow &\widehat{{\mathcal M}^{\rm op}} \, , \hfill \\
\hfill \varphi (\vec x) &\longmapsto &I_{\varphi (\vec x)} \hfill
\end{matrix}
$$
est pleinement fid\`ele.

\smallskip

D'apr\`es la condition (2), il existe pour tout objet $M$ de ${\mathcal M}$ une formule g\'eom\'etrique
$$
\varphi_M = \varphi_M (\vec x)
$$
dont l'image $I_{\varphi_M}$ dans $\widehat{{\mathcal M}^{\rm op}} = [{\mathcal M} , {\rm Ens}]$ est isomorphe \`a 
$$
y(M) = {\rm Hom} (M,\bullet) \, .
$$

Cette formule $\varphi_M$ est n\'ecessairement ${\mathbb T}$-irr\'eductible.

\smallskip

En effet, toute famille $J_{\mathbb T}$-couvrante de morphismes de ${\mathcal C}_{\mathbb T}$
$$
\theta_i (\vec x_i , \vec x) : \varphi_i (\vec x_i) \longrightarrow \varphi_M (\vec x) \, , \quad i \in I \, ,
$$
induit une famille globalement surjective d'applications
$$
\xymatrix{
M\varphi_i \ar@{=}[d] \ar[r] &M\varphi_M \ar@{=}[d] \\
I_{\varphi_i} (M) &I_{\varphi_M} (M) 
}
$$
donc il existe un indice $i_0 \in I$ et un \'el\'ement de
$$
I_{\varphi_{i_0}} (M) = {\rm Hom} (I_{\varphi_M} , I_{\varphi_{i_0}})
$$
dont le compos\'e avec le morphisme $I_{\varphi_{i_0}} \to I_{\varphi_M}$ est
$$
{\rm id} \in {\rm Hom} (M,M) = {\rm Hom} (I_{\varphi_M} , I_{\varphi_M}) = I_{\varphi_M} (M) \, .
$$

Comme le foncteur des interpr\'etations est pleinement fid\`ele, cet \'el\'ement de
$$
{\rm Hom} (I_{\varphi_M} , I_{\varphi_{i_0}})
$$
provient d'un morphisme de ${\mathcal C}_{\mathbb T}$
$$
\varphi_M (\vec x) \longrightarrow \varphi_{i_0} (\vec x_{i_0})
$$
dont le compos\'e avec
$$
\varphi_{i_0} (\vec x_{i_0}) \longrightarrow \varphi_M (\vec x)
$$
est l'identit\'e de $\varphi_M (\vec x)$.

\smallskip

Montrons que tout objet $\varphi (\vec x)$ de ${\mathcal C}_{\mathbb T}$ admet une famille $J_{\mathbb T}$-couvrante compos\'ee de formules $J_{\mathbb T}$-irr\'eductibles.

\smallskip

En effet, l'image $I_{\varphi}$ d'une telle formule $\varphi (\vec x)$ dans $\widehat{{\mathcal M}^{\rm op}}$ admet une famille globalement \'epimorphique de morphismes
$$
{\rm Hom} (M_i , \bullet) = y(M_i) \longrightarrow I_{\varphi} \, , \quad i \in I \, ,
$$
dont les sources sont les pr\'efaisceaux associ\'es \`a des objets $M_i$ de ${\mathcal M}$, donc sont des images
$$
I_{\varphi_i} \, , \quad i \in I \, ,
$$
de formules g\'eom\'etriques ${\mathbb T}$-irr\'eductibles $\varphi_i (\vec x_i)$, $i \in I$.

\smallskip

Alors les morphismes de $\widehat{{\mathcal M}^{\rm op}}$ 
$$
I_{\varphi_i} \longrightarrow I_{\varphi}
$$
sont les images par le foncteur pleinement fid\`ele
$$
{\mathcal C}_{\mathbb T} \longrightarrow \widehat{{\mathcal M}^{\rm op}}
$$
de morphismes de ${\mathcal C}_{\mathbb T}$, c'est-\`a-dire de formules ${\mathbb T}$-d\'emontrablement fonctionnelles
$$
\theta_i (\vec x_i , \vec x) : \varphi_i (\vec x_i) \longrightarrow \varphi (\vec x) \, , \quad i \in I \, .
$$

Par hypoth\`ese, la famille des applications
$$
M\varphi_i (\vec x_i) \longrightarrow M \varphi (\vec x) \, , \quad i \in I \, ,
$$
est globalement surjective pour tout objet $M$ de ${\mathcal M}$.

\smallskip

D'apr\`es la condition (1), cela implique que le s\'equent
$$
\varphi \vdash_{\vec x} \bigvee_{i \in I} (\exists \, \vec x_i) \, \theta_i (\vec x_i , \vec x)
$$
est ${\mathbb T}$-d\'emontrable.

\smallskip

Autrement dit, la famille des morphismes de ${\mathcal C}_{\mathbb T}$
$$
\theta_i (\vec x_i , \vec x) : \varphi_i (\vec x_i) \longrightarrow \varphi (\vec x)
$$
est $J_{\mathbb T}$-couvrante, et l'objet $\varphi (\vec x)$ admet un $J_{\mathbb T}$-recouvrement par la famille des formules irr\'eductibles $\varphi_i (\vec x_i)$.

\smallskip

D'apr\`es le th\'eor\`eme \ref{thmV911}, cela implique comme voulu que la th\'eorie ${\mathbb T}$ est de type pr\'efaisceau. 
\end{demo}
